%% file: main.tex
\documentclass[12pt,leqno]{article}

\usepackage{amsfonts,amsthm,amssymb,amsmath}
\usepackage[all]{xy}

\newdimen\sz

\makeatletter
\def\@oddfoot{ \hfill \thepage}
\def\@evenfoot{\thepage \hfill }
\def\@oddhead{{\sl\sbsct}\hfill}
\def\@evenhead{\hfill{\sl\sct}}
\makeatother

\DeclareMathSymbol\lar{\mathord}{AMSa}{"4C}
\DeclareMathSymbol\sbn{\mathrel}{AMSb}{36}
\def\hk{\mathrel{\relbar\mkern-12mu\raise3.8pt\hbox
{$\scriptscriptstyle\lar$}\mkern-7mu\rhook}}

\def\Section#1[#2]{%
  \section{#1}%
  \setcounter{equation}{0}%
  \setcounter{thm}{0}%
  \def\sct{\thesection. #2}}

\def\Subsubsection#1{%
  \subsection{#1}%
  \setcounter{equation}{0}%
  \setcounter{thm}{0}%
  \def\sbsct{\thesubsection. #1}%
  \def\sss{\thesection.}}

\newtheorem{thm}{Theorem}
\renewcommand{\thethm}{\thesubsection.\arabic{thm}}
\renewcommand{\theequation}{\thesubsection.\arabic{equation}}

\newtheorem{tspc}[thm]{\namespec}
\newtheorem{lem}[thm]{Lemma}
\newtheorem{ldf}[thm]{Lemma and Definition}
\newtheorem{cor}[thm]{Corollary}
\newtheorem{prp}[thm]{Proposition}
\newtheorem{tbl}[thm]{Table}
\theoremstyle{definition}
\newtheorem{dfn}[thm]{Definition}
\newtheorem{dfnn}[thm]{Definition and Notation}
\newtheorem{dfns}[thm]{Definitions}
\newtheorem{dspc}[thm]{\namespec}
\newtheorem{ntt}[thm]{Notation}
\newtheorem{nts}[thm]{Notations}
\newtheorem{cnv}[thm]{Convention}
\newtheorem{exm}[thm]{Example}
\newtheorem{exms}[thm]{Examples}
\theoremstyle{remark}
\newtheorem{rems}[thm]{Remarks}
\newtheorem{rem}[thm]{Remark}
\newtheorem{exe}[thm]{Exercise}
\newtheorem{rspc}[thm]{\namespec}
\newtheorem*{cau}{Caution}

\newtheorem{que}[thm]{Question}
\newtheorem*{pnt}{The Prime Number Theorem}
\theoremstyle{definition}
\newtheorem{prb}[thm]{Problem}
\newtheorem{prt}[thm]{}
\theoremstyle{remark}

\def\blm#1 {\begin{lem}\label{\sss l#1}}\def\elm{\end{lem}}
\def\bld#1 {\begin{ldf}\label{\sss l#1}}\def\eld{\end{ldf}}
\def\bth#1 {\begin{thm}\label{\sss t#1}}\def\eth{\end{thm}}
\def\bco#1 {\begin{cor}\label{\sss c#1}}\def\eco{\end{cor}}
\def\bpr#1 {\begin{prp}\label{\sss p#1}}\def\epr{\end{prp}}
\def\btb#1 {\begin{tbl}\label{\sss tb#1}}\def\etb{\end{tbl}}
\def\bqu#1 {\begin{que}\label{\sss q#1}}\def\equ{\end{que}}
\def\bdf#1 {\begin{dfn}\label{\sss d#1}}\def\edf{\end{dfn}}
\def\bdn#1 {\begin{dfnn}\label{\sss d#1}}\def\edn{\end{dfnn}}
\def\bds#1 {\begin{dfns}\label{\sss d#1}}\def\eds{\end{dfns}}
\def\bnt#1 {\begin{ntt}\label{\sss n#1}}\def\ent{\end{ntt}}
\def\bns#1 {\begin{nts}\label{\sss n#1}}\def\ens{\end{nts}}
\def\bcn#1 {\begin{cnv}\label{\sss cv#1}}\def\ecn{\end{cnv}}
\def\brs#1 {\begin{rems}\label{\sss r#1}}\def\ers{\end{rems}}
\def\brm#1 {\begin{rem}\label{\sss r#1}}\def\erm{\end{rem}}
\def\bex#1 {\begin{exe}\label{\sss ex#1}}\def\eex{\end{exe}}
\def\bxa#1 {\begin{exm}\label{\sss xa#1}}\def\exa{\end{exm}}
\def\bxs#1 {\begin{exms}\label{\sss xa#1}}\def\exs{\end{exms}}
\def\bpb#1 {\begin{prb}\label{\sss pb#1}}\def\exb{\end{prb}}
\def\bca{\begin{cases}}\def\eca{\end{cases}}

\def\bq#1 #2\e{\refstepcounter{equation}\lb{#1}\leftline{\er{#1}: ``#2''}\noindent}
\def\up#1{\uppercase{#1}}
\def\E{\expandafter\up}
\def\ag{abelian group}
\def\am{abelian monoid}
\def\av{absolute value}
\def\ad{addition}
\def\asm{antisymmetric}
\def\art{arithmetic}
\def\asc{associativ}
\def\as{assumption}
\def\Cm{C-monoid}
\def\cnc{cancellativ}
\def\ca{cardinality}
\def\ch{character}
\def\cf{coefficient}
\def\cmt{commutativ}
\def\cpl{complement}
\def\cme{composite }
\def\cm{composition}
\def\cn{condition}

\def\csq{Consequently}
\def\cd{contradict}
\def\cv{convolution}
\def\co{coordinate}
\def\crs{corresponding}
\def\ct{countabl}
\def\cre{creasing}
\def\dr{decimal representation}
\def\ds{decimal system}
\def\dcm{decomposable}
\def\DRT{Dedekind \rc on theorem}
\def\df{definition}
\def\dv{derivative}
\def\dm{dimension}
\def\dsb{distributiv}
\def\el{element}
\def\eq{equation}
\def\et{equalit}
\def\ep{equipotent}
\def\ev{equivalen}
\def\es{establish}
\def\ex{existence}
\def\epc{exponentiation}
\def\ext{extension}
\def\fc{factorization}
\def\fcg{finite cyclic group}
\def\fg{finitely generated}
\def\fw{following}
\def\fa{for all }
\def\fe{for every }
\def\fs{for some }
\def\fp{formal polynomial}
\def\f{function}
\def\fd{fundamental}
\def\gn{generalization}
\def\Gn{generator}
\def\GCD{greatest common divisor}
\def\hg{homogeneous}
\def\hm{homomorphi}
\def\idp{idempotent}
\def\im{immediate }
\def\Ip{in particular, }
\def\In{induction}
\def\iv{inductive}
\def\icg{infinite cyclic group}
\def\ig{integer}
\def\ido{integral domain}
\def\isc{intersection}
\def\il{interval}
\def\inv{invariant}
\def\is{iso\-mor\-phi}
\def\If{It follows }
\def\IT{iterate}
\def\jc{jective}
\def\ji{jectivity}
\def\jn{jection}
\def\lt{lattice}
\def\LHS{left-hand side}
\def\lc{linear combination}
\def\len{linear equation}
\def\lf{linear functional}
\def\lic{linear independence}
\def\lop{linear operator}

\def\lss{linear subspace}
\def\li{linearly independent}
\def\lo{lower bound}
\def\Mo{Moreover, }
\def\mf{morphism}
\def\ml{multipl}
\def\mlc{multiplication}
\def\mlv{multiplicative}
\def\nn{natural number}
\def\nog{natural ordering}
\def\ng{negative}
\def\nel{neutral element}
\def\nf{nonempty finite set}
\def\nfs{nonempty finite subset}
\def\ns{nonempty set}
\def\nss{nonempty subset}
\def\nm{number}
\def\nme{number of elements}
\def\oh{on the other hand, }
\def\ooo{one and only one }
\def\op{operation}
\def\Op{operator}
\def\ois{order-isomorphi}
\def\opr{order-preserving}
\def\oam{ordered abelian monoid}
\def\os{ordered set}
\def\og{ordering}
\def\PM{P-monoid}
\def\pt{partition}
\def\Pm{permutation}
\def\pl{polynomial}
\def\po{positive}
\def\pd{predecessor}
\def\Pf{prime field}
\def\Pn{prime number}
\def\pf{prime subfield}
\def\pn{principal}
\def\pj{projection}
\def\pp{propert}
\def\Pr{proposition}

\def\qt{quotient}
\def\ra{rational}
\def\rc{recursi}
\def\rd{reducible}
\def\rfl{reflexive}
\def\rl{relation}
\def\rp{represent}
\def\rt{restrict}
\def\rsp{respectively}
\def\RHS{right-hand side}
\def\ru#1{root#1 of unity}
\def\sf{satisf}
\def\sq{sequence}
\def\sg{semigroup}
\def\slt{semilattice}
\def\sr{semi\-ring}
\def\so{solution}
\def\sqr{square root}
\def\sc{structure}
\def\sbm{submonoid}
\def\sbt{subtraction}
\def\su{successor}
\def\st{such that }
\def\sft{sufficient}
\def\sy{symmetr}
\def\te{there exist}
\def\Tf{therefore }
\def\tos{totally ordered set}
\def\tr{transitiv}
\def\tl{translation}
\def\tp{transposition}
\def\uq{uniqueness}
\def\ue{unity element}
\def\uk{unknown}
\def\ub{upper bound}

\def\vs{vector space}
\def\vsf#1{vector space\ifcase#1{}\or s\fi\ over a field~}
\def\Wanp{We are now in a position to }

\def\wrt{with respect to }
\def\wlg{without loss of generality }
\def\zd{zero divisor}
\def\dw{\/{\rm:}}
\def\sd{\/{\rm;}}
\long\def\nic#1{}

\def\1{1\kern-3.3pt\hbox{\rm I}}
\def\<#1>{\left<#1\right>}
\let\a\alpha
\let\ax\approx
\let\b\beta

\mathchardef\Bc="120D

\def\Bca{\bigcap\limits}
\def\bcl{\bigcup\limits}
\def\Bg(#1){\Bigl(#1\Bigr)}
\def\bg(#1){\bigl(#1\bigr)}
\def\Bgg(#1){\biggl(#1\biggr)}
\def\Bij{\mathop{\rm Bij}}
\def\bu{{\boldsymbol u}}
\def\bv{{\boldsymbol v}}
\def\bw{{\boldsymbol w}}
\def\bx{{\boldsymbol\xi}}
\def\bet{{\boldsymbol\eta}}
\def\cA{\mathcal A}
\def\cB{\mathcal B}
\def\cL{\mathcal L}
\def\cM{\mathcal M}
\def\cP{\mathcal P}
\def\cR{\mathcal R}
\def\cX{\mathcal X}
\def\C{\mathbb C}
\def\chr{\mathop{\rm char}}

\let\d\delta
\let\D\Delta
\let\di\diamond

\def\dile{\stackrel\di\le}
\let\dsl\displaylines
\let\e\varepsilon
\let\Eq\equiv
\let\F\phi

\def\Fr{\mathop{\rm frac}}
\def\frp#1#2{\frac{#1}{#2}{}{\scriptstyle\prime}}
\let\g\gamma

\def\gcd{\mathop{\rm gcd}}
\let\ge\geqslant

\def\glb{\mathop{\rm g.l.b.}}
\def\hph#1,#2,{{\rm (#1)\hphantom{#2}}}
\def\id{\mathop{\rm id}\nolimits}
\def\iff{{}\Longleftrightarrow{}}
\def\imp{{}\Longrightarrow{}}

\def\Inj{\mathop{\rm Inj}}
\def\Int{\mathop{\rm int}}
\def\Inv{^{\rm inv}}
\def\INV{\mathop{\rm Inv}}
\let\iy\infty
\def\jh{\wh\jmath}
\def\ker{\mathop{\rm Ker}}
\let\la\lambda
\let\La\Lambda
\mathchardef\land="325E
\mathchardef\lor="325F

\def\LOR{\bigvee\limits}
\def\LB{{\rm LB}}
\def\lhp#1,#2,{\leavevmode\hphantom{(#1)}{\rm(#2)} }
\def\mid{\,|\,}
\def\mhc{\mathrel{\hat\circ}}

\def\hpl{\mathbin{\hat+}}
\def\hcd{\mathbin{\hat\cdot}}
\def\lcm{\mathop{\rm lcm}}
\let\le\leqslant
\let\lea\leftarrow
\def\lub{\mathop{\rm l.u.b.}}
\def\mo{^{\,-1}}
\def\mod#1 {\kern1.5pt{\rm(mod}\,#1)}
\def\mr#1{\mathrel{\rm #1}}
\let\mt\mapsto
\def\N{\mathbb N}
\def\Na{\mathbb N^\ast}
\def\Nf{(\N,+,\cdot,0,1)}
\def\nmid{\setbox0=\hbox{$/$}\kern1.5pt/\kern-\wd0 \hbox to\wd0{\hfil\kern-0.3pt$|$\hfil}\kern1.5pt}
\let\non\nonumber
\let\O\Omega
\let\o\omega

\let\ov\overline
\def\P{\mathrel{\rm P}}

\def\prodl{\prod\limits}
\def\pz#1#2{\ifcase#1(#2,+,\cdot,0,1)\or(#2,+,0)\or(#2,\cdot,1)\fi}
\def\Q{\mathbb Q}

\def\QU{\mathop{\raise1pt\hbox{\Large\hbox to0pt{\Large$\sqcap$\hss}$\sqcup$}}\limits}
\def\tQU{\mathop{\raise1pt\hbox{\Large\hbox to0pt{\Large$\widetilde\sqcap$\hss}$\sqcup$}}\limits}
\def\qu{\mathrel{\raise1pt\hbox{\hbox to0pt{$\scriptstyle\sqcap$\hss}$\scriptstyle\sqcup$}}}
\def\qum{\mathrel{\raise1pt\hbox{\hbox to0pt{$\scriptscriptstyle\sqcap$\hss}$\scriptscriptstyle\sqcup$}}}
\def\bqu{\mathrel{\raise1pt\hbox{\hbox to0pt{$\bar{\scriptstyle\sqcap}$\hss}$\scriptstyle\sqcup$}}}
\def\hqu{\mathrel{\raise1pt\hbox{\hbox to0pt{$\hat{\scriptstyle\sqcap}$\hss}$\scriptstyle\sqcup$}}}
\def\tqu{\mathrel{\raise1pt\hbox{\hbox to0pt{$\tilde{\scriptstyle\sqcap}$\hss}$\scriptstyle\sqcup$}}}
\def\dqu{\mathrel{\lower1pt\hbox{$\stackrel{\lower7pt\hbox{\smash{\hbox to0pt{$\scriptstyle\sqcap$\hss}$\scriptstyle\sqcup$}}}\cdot$}}}
\def\dquz{\mathrel{\lower1pt\hbox{$\stackrel{\lower7pt\hbox{\smash{\hbox to0pt{$
\scriptstyle\sqcap$\hss}$\scriptstyle\sqcup$}$_0$}}\cdot$}}}
\def\dqum#1{\mathrel{\lower1pt\hbox{$\stackrel{\lower7pt\hbox{\smash{\hbox to0pt{$
\scriptstyle\sqcap$\hss}$\scriptstyle\sqcup$}$_{#1}$}}\cdot$}}}
\def\dqun#1{\mathrel{\lower1pt\hbox{$\stackrel{\lower7pt\hbox{\smash{\hbox to0pt{$
\scriptstyle\sqcap$\hss}$\scriptstyle\sqcup$}$_0^{(#1)}$}}\cdot$}}}
\def\dqumz{\mathrel{\lower2pt\hbox{$\stackrel{\lower1pt\hbox{\smash{\hbox to0pt{$
\scriptscriptstyle\sqcap$\hss}$\scriptscriptstyle\sqcup$}$_0$}}\cdot$}}}
\def\dhqu{\mathrel{\lower1pt\hbox{$\stackrel{\lower7pt\hbox{\smash{\hbox
 to0pt{$\hat{\scriptstyle\sqcap}$\hss}$\scriptstyle\sqcup$}}}\cdot$}}}
\def\dtqu{\mathrel{\lower1pt\hbox{$\stackrel{\lower7pt\hbox{\smash{\hbox
 to0pt{$\tilde{\scriptstyle\sqcap}$\hss}$\scriptstyle\sqcup$}}}\cdot$}}}
\def\mqu{\stackrel{\smash{\scriptstyle\qu}}{{\smash-}\vrule height3pt depth0pt width0pt}}
\def\dpl{\mathrel{\lower2pt\hbox{$\stackrel{\lower1pt\hbox{$\scriptstyle+$}}\cdot$}}}
\def\dpln{\mathrel{\lower2pt\hbox{$\stackrel{+_n}\cdot$}}}
\def\dpls#1{\mathrel{\lower2pt\hbox{$\stackrel{+_{#1}}\cdot$}}}
\def\ddt{\mathrel{\lower4pt\hbox{$\stackrel\cdot\cdot$}}}
\def\ddtn{\mathrel{\lower2pt\hbox{$\stackrel{\cdot_n}\cdot$}}}
\def\ddtt#1{\mathrel{\lower2pt\hbox{$\stackrel{\cdot_{#1}}\cdot$}}}
\def\dci{\mathrel{\lower4pt\hbox{$\stackrel\circ\cdot$}}}
\def\dtci{\mathrel{\lower4pt\hbox{$\stackrel{\tilde\circ}\cdot$}}}
\def\dhp{\mathrel{\lower2pt\hbox{$\stackrel{\lower1pt\hbox{$\scriptstyle\hpl$}}\cdot$}}}
\def\lequ{\stackrel\qum\le}
\def\wtle{\mathrel{\wt\le}}
\def\dmi{\stackrel{\raise3pt\hbox{.}}{\smash-}}

\def\qn#1{\qu^{(#1)}}
\def\en#1{^{(#1)}}
\def\R{\mathbb R}
\def\rhu{{\upharpoonright}}
\def\sbk#1{{\textstyle{\substack{#1}}}}
\def\sbn#1#2{_{(#1)#2}}
\let\sbs\subset
\def\sgn{\mathop{\rm sgn}\nolimits}
\let\sm\setminus

\def\sms#1{\setminus\{#1\}}
\let\si\sigma
\let\Si\Sigma
\let\ssk\smallskip
\def\spn{\mathop{\rm span}}
\def\suml{\sum\limits}
\def\supp{\mathop{\rm supp}}
\let\td\widetilde
\let\t\times
\let\th\theta
\def\tR{\mathrel{\widetilde{\rm R}}}
\def\UB{{\rm UB}}

\let\ve\varepsilon
\let\vf\varphi
\let\vn\varnothing
\let\wh\widehat
\let\wt\widetilde
\def\Z{\mathbb Z}
\def\Zf{(\Z,+,\cdot,0,1)}
\def\zb#1#2#3{\{#1_{#2}\}_{#2\in #3}}
\def\zo#1 {[#1)}
\def\oz#1 {(#1]}

\def\lb#1{\label{\sss #1}}
\def\beq#1 #2\e{\begin{equation}\lb{#1}#2\end{equation}}
\def\bea#1 #2\e{\begin{align}\lb{#1}#2\end{align}}
\def\beag#1\e{\begin{align*}#1\end{align*}}
\def\bga#1 #2\e{\begin{gather}\lb{#1}#2\end{gather}}
\def\bgg#1\e{\begin{gather*}#1\end{gather*}}
\def\bml#1 #2\e{\begin{multline}\lb{#1}#2\end{multline}}
\def\bmlg#1\e{\begin{multline*}#1\end{multline*}}
\let\bal\aligned \let\eal\endaligned
\let\ga\gathered \let\ega\endgathered
\def\bmn{\left(\begin{matrix}}\def\emn{\end{matrix}\right)}
\def\bmk{\left[\begin{matrix}}\def\emk{\end{matrix}\right]}
\def\bit{\begin{itemize}}\def\eit{\end{itemize}}
\def\ben{\begin{enumerate}}\def\een{\end{enumerate}}
\def\labelenumi{{\rm\roman{enumi})}}
\def\sdim#1{\settowidth\labelwidth{#1} \addtolength\labelwidth\labelsep
\leftmargini=\labelwidth }

\def\vls{\vskip-\lastskip}
\def\hen#1 {\vls \bgroup\par \smallskip \leftskip\leftmargini \parindent0pt
\leavevmode \llap{\rm #1) }\ignorespaces}
\def\ehen{\par\egroup}

\def\er#1{\eqref{\sss #1}}
\def\era#1#2{{\rm(\ref{#1.#2})}}
\def\erp#1{{\rm(\Ps.\ref{\PS #1})}}

\def\rf#1{{\rm\ref{\sss #1}}}
\def\rfa#1#2{{\rm\ref{#1.#2}}}

\def\nad#1#2{\buildrel{#1}\over{#2}}
\def\nde#1 {\nad{\er{#1}}}
\def\ndp#1 {\nad{\erp{#1}}}
\def\nda#1#2 {\nad{\era{#1}{#2}}}
\let\q\quad
\def\qh#1{\quad\mbox{#1}\ }
\def\qhq#1{\quad\mbox{#1}\quad}

\def\Eproof{\hfill\qedsymbol}

\let\ti\textit
\let\noi\noindent
\let\lvm\leavevmode

\def\6{\begin{picture}(3,7)(0,0)\put(1.5,0){\line(0,1)7}\end{picture}}
\def\2{\begin{picture}(5.5,7)(0,0)\put(1.5,0){\line(0,1)7}\put(4,0){\line(0,1)7}\end{picture}}
\def\3{\begin{picture}(8,7)(0,0)\put(1.5,0){\line(0,1)7}\put(4,0){\line(0,1)7}\put(6.5,0){\line(0,1)7}\end{picture}}
\def\4{\begin{picture}(10.5,7)(0,0)\put(1.5,0){\line(0,1)7}\put(4,0){\line(0,1)7}\put(6.5,0){\line(0,1)7}\put(9,0){\line(0,1)7}\end{picture}}
\def\5{\begin{picture}(13,7)(0,0)\put(1.5,0){\line(0,1)7}\put(4,0){\line(0,1)7}\put(6.5,0){\line(0,1)7}\put(9,0){\line(0,1)7}\put(11.5,0){\line(0,1)7}\end{picture}}
\def\kr{\begin{picture}(10.5,7)(0,0)\line(3,2){10.5}\put(-9,0){\line(0,1)7}\put(-4,0){\line(0,1)7}\put(-6.5,0){\line(0,1)7}\put(-1.5,0){\line(0,1)7}\end{picture}}
\def\Xk{\begin{picture}(32,7)(0,0)\line(5,1){32}\hskip-31\unitlength XXXX\end{picture}}
\def\Ck{\begin{picture}(32,7)(0,0)\line(5,1){32}\hskip-30\unitlength CCCC\end{picture}}
\def\ite[#1]#2{{\parindent0pt \hangindent\leftmargini \hangafter1 \par\lvm \hbox to\leftmargini{#1 \hfil}#2}}

\makeindex
\makeglossary

\makeatletter

 \newenvironment{thesymb}
{\if@twocolumn \@restonecolfalse
\else \@restonecoltrue \fi
\twocolumn[\section*{\symbname}]%
\@mkboth{\MakeUppercase\symbname}{\MakeUppercase\symbname}%
\parindent\z@ \parskip\z@ \@plus .3\p@\relax
\columnseprule \z@ \columnsep 35\p@
\let\item\@idxitem}
\makeatother

  \def\AppendixSetup{%
  \def\thesection{\Alph{section}}%
  \def\Section##1[##2]{%
    \section{##1}%
    \setcounter{equation}{0}%
    \setcounter{thm}{0}%
    \def\sct{\thesection. ##2}}%
  \def\Subsubsection##1{%
    \subsection{##1}%
    \setcounter{equation}{0}%
    \setcounter{thm}{0}%
    \def\sbsct{\thesubsection. ##1}%
    \def\sss{\thesection.}}%
  \renewcommand{\thethm}{\thesubsection.\arabic{thm}}%
  \renewcommand{\theequation}{\thesubsection.\arabic{equation}}%
  \def\bqu##1 {\begin{que}\label{\sss q##1}}%
  \def\equ{\end{que}}%

  \def\ap{approximat}%
  \def\Ar{Archi\-medean}%
  \def\ass{assertion}%
  \def\Bs{Banach space}%
  \def\ba{bounded above}%
  \def\bb{bounded below}%
  \def\Cs{Cauchy sequence}%
  \def\chz{characterization}%
  \def\cp{complete}%
  \def\cnt{connect}%
  \def\ctn{continuous}%
  \def\cg{converge}%
  \def\cvg{convergent}%
  \def\de{decimal expansion}%
  \def\dn{decimal number}%
  \def\dfb{differentiab}%
  \def\dis{directed set}%
  \def\Hts{Hausdorff topological space}%
  \def\Hc{H\"older-continuous}%
  \def\he{homeomorphi}%
  \def\ipt{interior point}%
  \def\itd{introduced }%
  \def\Lc{Lipschitz-continuous}%
  \def\lg{lower limit}%
  \def\ms{metric space}%
  \def\nos{nonempty open subset}%
  \def\oag{ordered abelian group}%
  \def\od{ordered domain}%
  \def\of{ordered field}%
  \def\pw{pointwise}%
  \def\Ps{Polish space}%
  \def\rv{relative}%
  \def\sql{sequentially}%
  \def\tm{terminating}%
  \def\tps{topological space}%
  \def\tc{transcendental}%
  \def\ut{upper limit}%
  \def\vl{valuation}%
  \def\we{whenever }%
  \def\bc{{\boldsymbol c}}%
  \def\be{{\boldsymbol e}}%
  \def\bn{{\boldsymbol n}}%
  \def\bX{{\boldsymbol x}}%
  \def\bY{{\boldsymbol y}}%
  \def\bZ{{\boldsymbol z}}%
  \def\bz{{\boldsymbol 0}}%
  \def\bO{{\boldsymbol O}}%
  \def\bA{{\boldsymbol \a}}%
  \def\bB{{\boldsymbol \b}}%
  \def\bG{{\boldsymbol \g}}%
  \def\bD{{\boldsymbol \d}}%
  \def\bE{{\boldsymbol\ve}}%
  \def\cC{\mathcal C}%
  \def\cO{\mathcal O}%
  \def\cS{\mathcal S}%
  \let\da\downarrow
  \def\ee{{\rm e}}%
  \def\h{\hskip2pt minus0.7pt}%
  \def\Id{\mathop{\rm Id}\nolimits}%
  \def\infl{\inf\limits}%
  \def\kl{\mathaccent"7017 }%
  \def\limu{\mathop{\ov{\rm lim}}\limits}%
  \def\liml{\mathop{\underline{\rm lim}}\limits}%
  \def\lml{\mathop{\lim\limits_<}\limits}%
  \def\lmr{\mathop{\lim\limits_>}\limits}%
  \let\lm\Rightarrow
  \def\mid{\,\big|\,}%
  \def\nbmid{\setbox0=\hbox{$/$}\kern1.5pt/\kern-\wd0 \hbox to\wd0{\hfil\kern-0.3pt$\big|$\hfil}\kern1.5pt}%
  \def\hpl{\mathrel{\hat+}}%
  \def\hcd{\mathrel{\hat\cdot}}%
  \def\mo{^{-1}}%
  \def\htge{\mathrel{\hat\ge}}%
  \def\htg{\mathrel{\hat>}}%
  \def\whg{\mathrel{\wh>}}%
  \def\whge{\mathrel{\wh\ge}}%
  \def\whle{\mathrel{\wh\le}}%
  \def\htle{\mathrel{\hat\le}}%
  \def\htl{\mathrel{\hat<}}%
  \def\whl{\mathrel{\wh<}}%
  \def\wtg{\mathrel{\wt>}}%
  \def\wtl{\mathrel{\wt<}}%
  \def\wtge{\mathrel{\wt\ge}}%
  \def\spn{\mathop{\rm span}\limits}%
  \def\supl{\sup\limits}%
  \let\ua\uparrow
  \let\tb\textbf
  \def\labelenumi{{\rm(\roman{enumi})}}%
}

\begin{document}

\thispagestyle{empty}
\rightline{Ad maiorem Dei gloriam}
\vfill
\begin{center}

{\LARGE \bf
\uppercase{From Natural Numbers\\to Prime Fields and Finite Fields}

}
\end{center}

\vfill
\centerline{\Large \bf Philippe Cl\'ement}
\vskip10pt
\centerline{\large \bf Emeritus Professor TU Delft}
\vfill\vfill
\setbox0=\hbox{\large To Wilhelmina, our children}

\rightline{\vbox{\hsize\wd0 \large
\baselineskip20pt
\copy0
\centerline{and grandchildren}}}
\vfill\vfill\vfill\vfill
\eject
\thispagestyle{empty}
\ \vfill

\eject
\def\sct{}
\def\sbsct{}
\def\thepage{\roman{page}}
\setcounter{page}1
\hbadness1800

\tableofcontents
\newpage

\thispagestyle{empty}
\ \vfill

\eject

\addcontentsline{toc}{section}{Preface}
\section*{Preface}\setcounter{equation}0
\def\sct{Preface}\def\sbsct{Preface}
\def\sss{}

The aim of this book is to introduce the reader to the beauty of Algebra,
through a journey from the \nn s to \Pf s and finite fields, with some
detours. The author is not an algebraist but a retired professor of \E\f al
Analysis.

The basic fields used in Analysis are the fields of real \nm s and of complex
\nm s. Many books are devoted to the construction of these fields from the
\nn s. Perhaps the first one is~\cite{Grund}, see also the book \ti{Numbers}
\cite{Nrs}, Lay~\cite{18} and the references therein. An important step in this
process is the construction of the field of \ra\ \nm s. It turns out that this
field is the only (in an appropriate sense) infinite \ti{prime} field.
Motivated by the reading of Loonstra~\cite{Loonstra}, I~decided to write an
essentially self-contained book devoted to the construction not only of the
field of
\ra\ \nm s but of \ti{all\/} \Pf s and, more interestingly, to the proof of
\ex\ and ``\uq'' of \ti{all finite} fields. Only a knowledge of basic set
theory and some familiarity with mathematical reasoning (see for example
\cite[Chapters 1 and~2]{18}) are assumed.

This work is based on the classical book by van der Waerden \cite{Alg},
E.~Artin's lectures~\cite{Artin}, and the textbooks by Rotman~\cite{Groups}
and Isaacs~\cite{Is}. For further reading we recommend the books in Analysis by Amann and
Escher~\cite{21}. The standard reference for finite fields is 
Lidl and Niederreiter~\cite{Lidl}.

The author wishes to thank the Analysis research group at the Delft Institute
of Applied Mathematics TU Delft and the Department of Mathematics of the
University of Pavia for hospitality, the Department of Mathematics of the
University of Fribourg (Switzerland) for according him the use of the
Mathematics Library. He also expresses his gratitude to Herbert Amann
(Universit\"at Z\"urich), 
Jan van Neerven (TUD) for their valuable comments,
Mark Veraar (TUD), Giuseppe Savar\'e (Bocconi and Pavia University), Jan Maas (IST Austria), Guido
Sweers (Uni K\"oln), Marco Picasso (EPFL) and Carolus Reinecke (London) for their
interest. He is indebted to the late Hieltje Rijnks (TUD), Mark Veraar who brought reference
\cite{Nrs}, respectively \cite{18}, to his attention. He also warmly thanks Dr.~Jan Kowalski (IMPAN
Warsaw) for his careful typing of the manuscript.

I thank Wilhelmina and our family for their patience and support during the long process
of writing.

Comments and reports of mistakes, misprints are welcome. They can be sent to
{\tt philippeclem@gmail.com}\,.

\newpage

\clearpage
\pagestyle{headings}
\pagenumbering{roman}

\addcontentsline{toc}{section}{Introduction}
\section*{Introduction}
\begingroup
\setcounter{equation}{0}
\def\sct{Introduction}%
\def\sbsct{Introduction}%
\def\sss{I.}%
\def\rfa#1#2{{\rm\ref{#1.#2}}}%
\def\era#1#2{{\rm(\ref{#1.#2})}}%
\input detour-i
\endgroup
\ \vfill

\eject
%
%
%
%

%
\setcounter{page}{1}
\def\thepage{\arabic{page}}

\input DETOUR1

\newpage
\input DETOUR2

\newpage
\input DETOUR3

\newpage
\input DETOUR4

\input DETOUR44

\newpage
\input DETOUR5


\clearpage
\begingroup
\AppendixSetup
\setcounter{section}{0}
\setcounter{subsection}{0}
\setcounter{equation}{0}
\setcounter{thm}{0}
\def\sct{Appendix}\def\sbsct{Appendix}\def\sss{}
\addcontentsline{toc}{section}{Appendix: Archimedean Fields}
\section*{Appendix: Archimedean Fields}
\input app1

\input app2

\clearpage

\def\sct{References}\def\sbsct{References}
\input references


%
\newpage

\def\sct{Index}\def\sbsct{Index}
\def\sss{}
\addcontentsline{toc}{section}{Index}
\InputIfFileExists{\jobname.ind}{}{%
  \typeout{No \jobname.ind found. Run makeindex \jobname.idx to build the combined index.}%
}

\def\symbname{List of symbols}
\def\sct{List of symbols}\def\sbsct{List of symbols}
\def\sss{}
\addcontentsline{toc}{section}{List of symbols}
\begin{thesymb}
\advance\baselineskip by1pt
\input detour.sym
\end{thesymb}

\end{document}

%% file: detour-i.tex
\def\theequation{I.\arabic{equation}}

The first three chapters of this book deal with \nn s and the last two with
prime and finite fields. In the first section of Chapter one
(Section~\ref{sss.Sets}) we present a \df\ of a ``set of \nn s'' (s.o.n.n.)
(\E\df\ \rfa1{d1.3}), the \E\rc on theorem of Dedekind (Theorem \rfa1{t1.6})
and a theorem about \uq\ of s.o.n.n.\ (Theorem \rfa1{t1.4}). This section is
based on \S2.1 and \S2.2 of Chapter~1 of~\cite{Nrs}. Our \df\ of a s.o.n.n.\ is
a slight \gn\ of the \df\ given in~\cite{Nrs}. We define a s.o.n.n.\ as a
triple $(E,e,S)$, where $E$~is a set, $e$~is a distinguished \el\ of~$E$ and\glossary{s.o.n.n.}
$S$~is a map from~$E$ into itself called the ``\su\ \f'' \sf ying axioms (N1),
(N2), \ev t to the Peano axioms~\cite{19}. In~\cite{Nrs}, $E$~is denoted by~$\N$ and
$e$~by~$0$. Some authors use the notation $\{0,1,2,\dots\}$ for~$\N$ and call
this set the set of \ti{non\ng\ \ig s\/}. The case $e:=1$ corresponds to the
\df\ given in~\cite{Grund} (Axioms 1--5). The set~$E$ is denoted by $1,\dots$
in the preface for the student of~\cite{Grund}, and is called the totality of
the \ti{\po\ \ig s\/}. Note that this set is denoted by~$\N$ in
\cite[p.~99]{18}. We prefer the notation $\{0,\dots\}$ for~$E$ if $e=0$,
$\{0,1,\dots\}$ if $e=0$ and $1=S(0)$, and $\{1,\dots\}$ if $e=1$. We reserve
the notation~$\N$ for the set $\{0,1,\dots\}$ when its \el s are identified
with their \dr\ (see Section \ref{sss.Decimal}), and $\Na$~for the set $\N
\sms0$. In this book we assume the \ex\ of such sets. For a proof we refer
the reader e.g.\ to~\cite[pp. 271--272]{Kelley}.

We now turn to the {\bf \DRT} (Theorem \rfa1{t1.6}), the cornerstone of the theory.
Given a s.o.n.n.\ $(E,e,S)$ as in \E\df\ \rfa1{d1.3}, a nonempty set~$X$, a
map $g:X\to X$ and an \el~$a$ of~$X$, we want to find a map $\psi:E\to X$
(also called a \sq\ of \el s of~$X$) \sf ying the \rc ve \rl
\beq 1
\psi(S(n)) = g(\psi(n)) \qh{\fa} n\in E,
\e
and the ``initial value''
\beq 2
\psi(e)=a.
\e
\If from axiom (N2) of \E\df\ \rfa1{d1.3} that \te s \ti{at most one\/}
map~$\psi$ \sf ying (1), (2). Indeed, suppose that $\psi,\psi':E\to X$
\sf y (1) and~(2). Set $A:=\{n\in E: \psi(n)=\psi'(n)\}$. Then
by~(2) $\psi(e)=a=\psi'(e)$, hence $e\in A$. Suppose that $n\in A$. Then
$\psi(S(n))\nde1 = g(\psi(n)) \nad{n\in A}= g(\psi'(n))\nde1 = \psi'(S(n))$,
hence $S(n)\in A$. Since $n$~is arbitrary in~$A$, we have $S(A)\sbs A$. \E\Tf
we obtain $A=E$ from (N2), that is, $\psi=\psi'$. Such a proof is called a
{\bf proof by \ti{\In}}. It turns out that \te s \ti{at least one\/} map $\psi:E\to
X$ \sf ying \er1, \er2. Our ``\ex'' proof is a straightforward adaptation of
the one given in \cite[pp.~16--17]{Nrs}.

Proceeding as in \cite[p.~17]{Nrs} we show in Theorem \rfa1{t1.4} that if
$(E,e,S)$ and $(E',e',S')$ are s.o.n.n., then \te s a unique bi\jn\ $\vf:E\to
E'$ \sf ying
\beq3
\vf(e)=e', \q S'\circ\vf = \vf\circ S.
\e
In contrast to many authors (e.g.~\cite{Nrs}, \cite{Grund}, \cite{Alg}) who
first introduce an \ad\ and then an \og\ on~$E$, we \ti{first introduce an \og}
$\le$ on~$E$. It turns out that this \og\ is the only \og\ on~$E$ for which
$e$~is the least \el\ of~$E$ and $S$~is an in\cre\ map (see Theorem
\rfa1{t3.23}). We call this \og\ the \ti{natural\/} \og\ of the s.o.n.n.\
$(E,e,S)$. To this end we introduce what we call the ``\ti{set of \IT s}'' of
the \su\ \f~$S$. If we view $S$ as a self-map of~$E$, the set of \IT s of~$S$
could be denoted by $\{\id_E,S,S\circ S,\dots\}$ where $\circ$~denotes the
\cm\ of self-maps of~$E$. Given a s.o.n.n.\ $(E',e',S')$, possibly equal to
$(E,e,S)$, we want to define a map~$\Phi$ from~$E'$ into~$E^E$ (the set of
self-maps of~$E$) \sf ying the \rc ve \rl
\beq4
\Phi(S'(n')) = S\circ \Phi(n') \qh{\fa}n'\in E',
\e
and
\beq5
\Phi(e'):=\id_E.
\e
\E\ex\ and \uq\ of such a map $\Phi$ is guaranteed by the \rc on theorem with
\beq6
(E,e,S):=(E',e',S'), \q X:=E^E,
\e
and
\beq7
g(h):=S\circ h \qh{\fe} h\in E^E.
\e
We use the notation $\Phi_{E'}$ if we want to emphasize the dependence on~$E'$.
We set
$$
I_{E'}(S) := \bigcup_{n'\in E'}\Phi_{E'}(n').
$$
Observe that if we replace $S\circ h$ by $h\circ S$ in~\er7 and denote the \crs\
map by~$\check \Phi_{E'}$, then we have in view of \E\Pr\ \rfa1{p2.6}:
\beq8
\check\Phi_{E'}(n') = \Phi_{E'}(n') \qh{\fa}n'\in E'.
\e
If $(E'',e'',S'')$ is a s.o.n.n.\ and $\vf:E'\to E''$ is the bi\jn\ introduced
in Theorem \rfa1{t1.4} \sf ying
\beq8a
\vf(e')=e'' \qh{and }S''\circ\vf = \vf\circ S',
\e
then we have the \fw\ diagram:
$$
\xymatrix{(E',e',S')\ar[rr]^\vf \ar[rd]_{\Phi_{E'}}&&(E'',e'',S'')
\ar[ld]^{\Phi_{E''}}\cr
&E^E&\cr}
$$
and we obtain by \E\Pr\ \rfa1{p2.8}:
\bea9
\Phi_{E'}(n')&= \Phi_{E''}(\vf(n')) \qh{\fa} n'\in E',\\
I_{E''}S & =\Phi_{E'}(S). \lb{10}
\e
We set
\beq11
I(S):=I_{E'}(S)(=I_E(S)),
\e
and call $I(S)$ the {\bf set of \IT s of} $S$. This set is \ti{independent\/}\glossary{$I(S)$}
of the choice of the s.o.n.n.\ used to ``label'' the \IT s of~$S$.
\E\pp ies of the \IT s of~$S$
are collected in \E\Pr s \rfa1{p2.7} and \rfa1{p2.8}.

These \pp ies imply that the \fw\ \rl\ on~$E$ is an {\bf\og} on~$E$ (see
\E\df\ \rfa1{d3.2} and \E\Pr\ \rfa1{p3.5}):
\beq12
\hbox{\E\fa} m,n\in E:\, n\le m \hbox{ if \te s }h\in I(S) \hbox{ \st}
m=h(n).
\e
This \og\ will be called the \ti{natural\/} \og\ on $(E,e,S)$, and we shall
use the notation $(E,e,S,\le)$.

This \og\ is investigated in detail in Section \ref{sss.Ordering}. The \rl~$<$
defined by
\beq13
n< m \qh{if $n\le m$ and $n\ne m$ \fa} n,m\in E,
\e
plays an important role. For the convenience of the reader we collect some
\pp ies of the \og\ defined in~\er{12}.

For all $x,y,z\in E$,
\bea14
&\era1{3.11}\hphantom{Lemma }\ x\le y \hbox{ and }y< z \hbox{ implies }x<z,\\
&\era1{3.12}\hphantom{Lemma }\ x<y \hbox{ and }y\le z \hbox{ implies }x<z. \lb{15}\\
\intertext{\E\fa $x,y\in E$,}
&\era1{3.20}\hphantom{Lemma }\ \hbox{either } x=y, \hbox{ or }x<y, \hbox{ or }y<x, \lb{16}\\
&\era1{3.15}\hphantom{Lemma }\ x<y \hbox{ implies } \Phi(n')x < \Phi(n')y \hbox{ \fa} n'\in E', \lb{17}\\
&\era1{3.17}\hphantom{Lemma }\ x<y \hbox{ implies }S(x)\le y. \lb{18}\\
\intertext{\E\fa $x\in E$,}
&\era1{3.13}\hphantom{Lemma }\ x<S(x), \lb{19} \\
&\hbox{(Lemma \rfa1{l3.27})}\ \hbox{there is no $y\in E$ \st }x<y<S(x).\hskip120pt \lb{20}\\
\intertext{If $\le'$ denotes the \nog\ of $(E',e',S')$, then, \fa $n',m'\in E'$
\sf ying $n'<' m'$,}
&\era1{3.16}\hphantom{Lemma }\ \Phi(n')x < \Phi(m')x \q \hbox{\fa}x\in E. \lb{21}
\e

We next mention some \pp ies of nonempty subsets of $E$.

Let $A$ be a \ti{nonempty} subset of $E$.
\bea22
&\hbox{(Theorem \rfa1{t3.22})}\ \hbox{\E\te s \ooo \el\ $a\in A$ \sf ying
$a\le x$}\\
&\hphantom{(Theorem 1.2.22)}\ \ \hbox{\fa $x\in A$.}\non\\
&\era1{3.19}\hphantom{Theorem } \ \ e\le x\ \hbox{\fa $x\in E$ and $e<x$ \fa} x\in E\sms e. \lb{23}
\e

A subset $A$ of $E$ is called {\bf bounded} if $A$ is not empty and \te s
$m\in A$ \st $x\le m$ \fa $x\in A$. The \fw\ holds:
\bea24
&\text{(Theorem \rfa1{t3.39}) If $A$ is a bounded subset of~$E$, then \te s
one and only}\\
&\text{one \el~$\wh a$ of $A$, called the greatest \el\ of~$A$, \sf ying $x\le\wh a$}\non \\
&\text{\fa $x\in A$.}\non\\
&\text{(\E\Pr\ \rfa1{p3.48}) Let $A$ be a nonempty subset of~$E$, and $f$ be a
map from~$A$} \lb{25} \\
&\text{into~$E$. If $A$ is bounded, then so is $f(A)$.}\non\\
&\text{(Theorem \rfa1{t3.49}) A nonempty subset $A$ of $E$ is not bounded
(or {\bf unbounded})} \lb{26}\\
&\text{iff \te s a bi\jn\ $\vf:E\to A$ \sf ying $\vf(x)\le
\vf(y)$ iff $x\le y$ \fa}\non \\
&\text{$x,y\in E$.}\non
\e

In Section \ref{sss.Card} we introduce the notions of finiteness and
infiniteness of a set. Our goal in this introduction is to motivate \E\df\
\rfa1{d4.17}. Let $A,B$ be nonempty sets (possibly equal). Then $A$ and~$B$
are called \ti{\ep\/} (see \cite[p.~370]{Nrs}, \cite[p.~3]{Alg}), equipollent
(\cite[p.~28]{Kelley}) or equinumerous (\cite[p.~78]{18}), if \te s a bi\jn\
from $A$ onto~$B$. We use then the notation $A\approx B$. In view of \pp ies of
bi\jn s we infer $A\approx A$, $A\approx B$ implies $B\approx A$, and $A\approx B$,
$B\approx C$ implies $A\approx C$. We first investigate this notion for the set
of all nonempty subsets of a s.o.n.n.\ $(E,e,S)$ which we denote by $\dot\cP
(E)$. This notion induces an \ti{\ev ce \rl\/} (see \E\df\ \rfa1{d4.3}) on
$\dot\cP(E)$. By \er{26}, a nonempty unbounded subset of~$E$ is \ep\ to~$E$.
The set $E$ is unbounded. Indeed, suppose for \cd ion that \te s $m\in E$ \st
$x\le m$ \fa $x\in E$. By~\er{25} \te s $\wh a\in A$ \st $x\le\wh a$ \fa
$x\in A$. However, $\wh a< S(\wh a)$ by~\er{19}. \E\Tf $S(\wh a)\le\wh a$
and $\wh a<S(\wh a)$ implies $S(\wh a)<S(\wh a)$ by~\er{14}, hence $S(\wh a)\ne
S(\wh a)$ by~\er{13}. A~\cd ion. Thus $E$ is unbounded. We claim that a
\ti{bounded\/} subset $A$ of~$E$ is \ti{not\/} \ep\ to~$E$ (in notation
$A\not\approx E$). Indeed, suppose, for \cd ion, that $A\approx E$.
Then \te s a bi\jn\ $g:A\to E$, thus $g(A)=E$. However, from
\er{25} we infer that $g(A)$ is bounded since $A$~is bounded. A~\cd ion. Thus
$A\not\approx E$. \If from Lemma \rfa1{l4.10} that a nonempty \ti{bounded\/}
subset of~$E$ is \ti{\ep\/} to $[e,n]$ \fs $n\in E$ where
\beq27
[e,m]:=\{x\in E: e\le x\le m\} \qh{\fa} m\in E.
\e

We now show that if $(E',e',S')$ is a s.o.n.n.\ and $A$ is a bounded subset
of~$E$, then $A\approx[e',n']'$ \fs $n'\in E'$. Let $\vf:E\to E'$ be the bi\jn\
introduced in Theorem \rfa1{t1.4}. \If from Theorem \rfa1{t3.36} that $\vf (x)
\le'\vf(y)$ iff $x\le y$ \fa $x,y\in E$ where $\le'$ denotes the \nog\ of
$(E',e',S')$. \E\Tf $\vf$~is a bi\jn\ from $[e,n]$ onto $[e',\vf(n)]'$,
hence $[e,n]\approx [e',\vf(n)]'$. Since $A\approx[e,n]$, we obtain $A\approx
[e',n']'$ \fs $n'\in E'$, where $[x',y']':=\{z'\in E': x'\le' z'\le'y'\}$ \fa
$x',y'\in E'$. Conversely, if $A$~is a \nss\ of~$E$ \ep\ to $[e',n']'$ \fs
$n'\in E'$, then $A$~is bounded in $(E,e,S,\le)$. Indeed, note that $[e',n']'$
is bounded in $(E',e',S',\le')$, hence by what precedes $[e',n']'\approx[e,n]$
\fs $n\in E$. Thus $A$~is \ep\ to~$[e,n]$. Let $g:[e,n]\to A$ be a bi\jn. Then
$A=g([e,n])$ is bounded in $(E,e,S,\le)$ by~\er{25} since $[e,n]$ is bounded
in $(E,e,S,\le)$. By \E\df\ \rfa1{d4.17} a \ns~$A$ is
called {\bf finite} if $A$~is \ep\ to some \il\ $[e',n']'$, $n'\in E'$, of
some s.o.n.n. $(E',e',S',\le')$. The empty set is finite, and a set~$A$ which
is not finite is called {\bf infinite}.
A~set is called {\bf \ct y infinite} or denumerable (see
\cite[p.~79]{18}) if $A\approx E$. If a set $A$ is finite or \ct y infinite,
it is called \ct e, and
un\ct e otherwise. Cantor showed that the set $\cP(E)$, the power set of~$E$
is {\bf un\ct e}. In Section \ref{sss.Card} the \fw\ important \pp y of
nonempty finite sets is \es ed (Theorem \rfa1{t4.18}\,(iii)): a self-map of
a nonempty finite set is in\jc\ iff it is sur\jc. Dedekind called a nonempty
set~$A$ finite if every in\jc\ self-map of~$A$ is sur\jc\ (see \cite[p.~15]{Nrs}).
Such a set is called nowadays {\bf Dedekind-finite}. Thus a finite set is
Dedekind-finite. The converse can be shown, for example, by making use of the
axiom of choice. However, this topic exceeds the scope of this book. For
further use we mention (see Lemma \rfa1{l4.11} and Corollary \rfa1{c4.13})
the \fw\ result.

Let $(E,e,S)$ be a s.o.n.n., let $m,n\in E$ and let $f$ be a map from $[e,m]$
into $[e,n]$.
\beq28
\aligned
\hbox{(i)}\q &\hbox{If $f$ is in\jc, then $m\le n$.}\\
\hbox{(ii)}\q &\hbox{If $f$ is sur\jc, then $n\le m$.}\hskip220pt
\endaligned
\e

In Chapter two the notion of ``{\bf binary \op}'' plays an important role. A~basic
binary \op\ on a s.o.n.n.\ is the \ti{\ad\/}. We first motivate its \df. We
recall that a binary \op\ on a nonempty set~$X$ is a map from $X\times X$
into~$X$. A~\fd\ example is the \ti{\cm\/} of two self-maps $f$ and~$g$ of
a nonempty set~$Y$. This binary \op, usually denoted by~$f\circ g$, \sf ies
$(f\circ g)\circ h = f\circ(g\circ h)$ \fa $f,g,h\in Y^Y$ (\ti{\asc ity\/})
and $\id_Y$, the identity map in~$Y$, \sf ies $f\circ\id_Y = \id_Y\circ f=f$
\fa $f\in Y^Y$. A~nonempty set~$M$ equipped with a binary \op\ denoted by~$\qu$
is called a {\bf monoid}~\cite{Lattice} if the \op~$\qu$ is \asc e and if
\te s an \el\ $e\in M$, called the \ti{\nel\/} or \ti{identity\/}, \sf ying
$a\qu e= e\qu a=a$ \fa $a\in M$. We denote a monoid by $(M,\qu,e)$ (see \E\df\
\rfa1{d2.2}). In contrast to \cite{Nrs}, \cite{Grund},~\cite{Alg}, this
algebraic notion will be used extensively in this book. If $(E,e,S)$ is a
s.o.n.n., then $(E^E,\circ,\id_E)$ is a monoid. \Mo $I(S)$, the set
of \IT s of~$S$, is a subset of~$E^E$ \sf ying $\id_E\in I(S)$ by~\er5, and
$h_1\circ h_2\in I(S)$ whenever $h_1,h_2\in I(S)$, in view of \E\Pr\
\rfa1{p2.6} \era1{2.27}. \csq, the \op~$\circ$ is also a map from $I(S)\t
I(S)$ into $I(S)$, and $(I(S),\circ,\id_E)$ is a monoid. The set $I(S)$ is then
called a {\bf\sbm} of the monoid $(E^E,\circ,\id_E)$ (see \E\df\ \rfa1{d1.2}).
Observe that the map $\Phi:E\to I(S)$ defined by \er4, \er5 is sur\jc\ in view
of the \df\ of~$I(S)$. It turns out that $\Phi$~is also in\jc. Indeed, we have
\beq29
\Phi(n)e = n \qh{\fa}n\in E,
\e
by \E\Pr\ \rfa1{p2.7} \era1{2.35}. Thus, if $\Phi(m)=\Phi(n)$, $m,n\in E$, then
$m=\Phi(m)e = \Phi(n)e=n$. \E\Tf $\Phi$~is bi\jc\ as well as the map $\Phi\Inv
:I(S)\to E$. \If from Example \rfa2{xa1.5}\,(v) that the ``monoid \sc'' of
$I(S)$ can be ``transported'' to the set~$E$ by the map $\Phi\Inv$. Using
$\Phi=(\Phi\Inv)\Inv$, the ``transported'' \op~$\qu'$ is defined by:
\beq30
m\qu' n:= \Phi\Inv(\Phi(m)\circ \Phi(n)) \qh{\fa }m,n\in E.
\e
By \E\Pr\ \rfa1{p2.7} \era1{2.37}, \te s \ooo $p\in E$ \st $\Phi(m)\circ
\Phi(n) = \Phi(p)$. Thus $m\qu' n=p$. The \nel\ $e':=\Phi\Inv(\id_E)$, thus
$e'=e$. In Lemma \rfa2{l1.1} it is shown that the \fw\ holds \fa $m,n\in E$:
\bga31
m\qu' e=m,\\
m\qu'(S(n)) = S(m\qu' n).\lb{32}
\e
In \cite[p.~17]{Nrs} formulae \er{31}, \er{32} are the defining \rl s for the
{\bf\ad} on~$E$ whenever $e=0$. \E\Tf we use the notation $+_E$ instead of~$\qu'$,
and $+$ whenever $e=0$. The monoid $(I(S),\circ,\id_E)$ enjoys important \pp
ies, namely \fa $f,g,h\in I(S)$:
\bea33
\hbox{``\cmt ity'' } & f\circ g = g\circ f, \\
\hbox{``\cnc ity'' } & f\circ h= g\circ h \hbox{ implies } f=g, \lb{34} \\
\hbox{``positivity'' } & f\circ g = \id_E \hbox{ implies } f=g=\id_E. \lb{35}
\e
\E\pp y \er{33} (resp.\ \er{35}) is a con\sq\ of \E\Pr\ \rfa1{p2.6} \era1{2.26}
(resp.\ \E\Pr\ \rfa1{p2.7} \era1{2.41}). For the convenience of the reader we
\es~\er{34}. Suppose $\Phi(m)\circ \Phi(k)=\Phi(n)\circ\Phi(k)$ \fs $m,n,k\in
E$. Then $\Phi(k)m \nde29 = \Phi(k)(\Phi(m)e) = (\Phi(k)\circ\Phi(m))e \nde33
=(\Phi(m)\circ\Phi(k))e = (\Phi(n)\circ\Phi(k))e \nde33
= (\Phi(k)\circ\Phi(n))e = \Phi(k)(\Phi(n)e) \nde29 = \Phi(k)n$. Since $S$~is
in\jc\ by~(N1), we have $\Phi(k)$ in\jc\ by the \E\Pr\ \rfa1{p2.6} \era1{2.28}.
Hence $m=n$ and $\Phi(m)=\Phi(n)$, which proves \er{34}.

A monoid \sf ying \er{33} is called abelian. We use the \ti{shorthand
notation\/} \Cm\ for an abelian monoid \sf ying~\er{34}, and \PM\ for a \Cm\
\sf ying~\er{35}. \If from Example \rfa2{xa1.5}\,(v) that $(E,+_E,e)$ is also
a \PM.

Observe that the map $\Phi\Inv:I(S)\to E$ \sf ies $\Phi\Inv(\id_E)=e$. \Mo by
\er{30} $\Phi\Inv(h_1\circ h_2) = \Phi\Inv(h_1)+_E \Phi\Inv(h_2)$ \fa $h_1,h_2
\in I(S)$, since $h_1=\Phi(m)$, $h_2=\Phi(n)$ \fs $m,n\in E$. If $(M,\qu,e)$
and $(M',\qu',e)$ are monoids and $f:M\to M'$ \sf ies
\beq36
f(e)=e' \hbox{ and } f(a\qu b)=f(a)\qu' f(b)\qh{\fa}a,b\in M,
\e
then $f$ is called a {\bf(monoid-)\hm sm}. If, in \ad, $f$ is bi\jc, then
$f$~is called a {\bf(monoid-)\is sm\/} (see \E\df\ \rfa2{d1.7} and Lemma
\rfa2{l1.8}). Thus, $\id_M$ is an \is sm, if $f$~is an \is sm, then so is
$f\Inv$, and the \cm\ of two \is sms is an \is sm. \Mo the monoids
$M$~and~$M'$ are called \ti{\is c\/} if \te s an \is sm between $M$ and~$M'$.
\csq, $(I(S),\circ,\id_E)$ and $(E,+_E,e)$ are \is c \fe s.o.n.n.\ $(E,e,S)$.
We next look for \pp ies of a monoid which makes it \is c to $(I(S),\circ,
\id_E)$, \ev tly to $(E,+_E,e)$. To this end we define the notion of {\bf\IT s}
of an \el\ of a monoid. This is done in \E\Pr\ \rfa2{p1.13}. Given a monoid
$(M,\qu,e)$, an \el\ $a\in M$, and a s.o.n.n.\ $(\wt E,\wt e,\wt S)$, using
the \E\rc on theorem, one shows that \te s \ooo \sq\ $\vf_a:\wt E\to M$ \sf
ying $\vf_a(\wt e)=e$ and
\beq37
{\vf_a(\wt S(\wt n))=a} \qu {\vf_a(\wt n)} \qh{\fa}\wt n\in \wt E.
\e
Setting $I_{\td E}(a):= \bcl_{\td n\in\td E}\vf_a(\wt n)$ we find that
$I_{\td E}(a)$ is a \sbm\ of $(M,\qu,e)$, and that the monoid $(I_{\td E},\qu,
e)$ is \ti{abelian\/}. \E\Ip $I_{\td E}(e)$ and $(\{e\},\qu,e)$ is called the
trivial monoid. It turns out that if $(E',e',S')$ is a s.o.n.n., then
$I_{E'}(a)=I_E(a)$, which we denote by~$I(a)$. In \E\df\ \rfa2{d1.18} we call
a nontrivial monoid $(M,\qu,e)$ a {\bf\pn\ monoid} if \te s $a\in M\sms e$
\st $M=I(a)$. In this case $a$ is called a {\bf\Gn\/} of the monoid~$M$. This
terminology is not standard. Note that a \pn\ monoid is abelian. If $(E,e,S)$
is a s.o.n.n.\ and $I(S)$ denotes the set of \IT s of the self-map~$S$ in
$E^E$ defined in~\er{11}, then $I(S)$ coincides with the set of \IT s of the
\el~$S$ of the monoid $(E^E,\circ,\id_E)$, also denoted by~$I(S)$. \If that
$(I(S),\circ,\id_E)$ is an infinite \pn\ monoid with \Gn~$S$. We now claim
that if a monoid $(M,\qu,e)$ is an \ti{infinite \pn\/} monoid, then it is \is c
to $(I(S),\circ,\id_E)$, hence also to $(E,+_E,e)$. Indeed, by Theorem
\rfa2{t1.42}, $M$~is a \pn\ \PM\ and by \E\Pr\ \rfa2{p1.26}, $M$~is \is c to
$(E,+_E,e)$. Since the \cm\ of two \is sms is an \is sm, the claim is proved.

It also follows from Theorem \rfa2{t1.42} that an \ti{infinite \pn\/} monoid
has \ti{at most one} \Gn. \E\Tf\ $S$ is the only \Gn\ of the monoid $(I(S),
\circ,\id_E)$. It is shown in Lemma \rfa2{l1.22} that if $M,M'$ are monoids,
$\vf:M\to M'$ is a \hm sm, and $a\in M$, then
\beq38
\vf(\vf_a(\wt n)) = \vf_{\vf(a)}(\wt n) \qh{\fa}\wt n\in\wt E,
\e
where $(\wt E,\wt e,\wt S)$ is a s.o.n.n., and $\vf_a$ (resp.\ $\vf_{\vf(a)}$)
is the \sq\ of \IT s of~$a$ in~$M$ (resp.\ of $\vf(a)$ in~$M'$). \E\Ip if
$\vf$~is an \is sm and $M=I(a)$, then $M'=I(\vf(a))$. Indeed, let $x'\in M'$
and $x\in M$ be \st $\vf(x)=x'$. Then \te s $\wt n\in\wt E$ \sf ying $x=\vf_a
(\wt n)$. From \er{38} we infer that $x'=\vf(x) = \vf(\vf_a(\wt n)) =
\vf_{\vf(a)}(\wt n)$, hence $x'\in I(\vf(a))$. Since $x'$~is arbitrary in~$M'$,
we obtain $M'\sbs I(\vf(a))$, which together with $I(\vf(a))\sbs M'$ implies
$M'=I(\vf(a))$. \E\Ip $\vf(a)$~is a \Gn\ of~$M'$. Since $\Phi\Inv$ is an \is
sm from $(I(S),\circ,\id_E)$ into $(E,e,S)$ and $S$~is the \Gn\ of~$I(S)$, we
infer from what precedes that $S(e)= \Phi\Inv(S)$ is the \Gn\ of the \pn\
monoid $(E,e,S)$. Note that \er{37} with $a:=S(e)$ and ${\qu}:=+_E$ becomes
$\vf_{S(e)}(S(\wt n)) = S(e)+_E \vf_{S(e)}(\wt n)$ \fa $\wt n\in\wt E$. \E\Tf
since $E=I(S)$ and $\vf_{S(e)}(\wt e)=e$, we could write
\[
E = \bigl\{ e, S(e), S(e)+_E S(e), S(e)+_E S(e)+_E S(e), \dots\bigr\}.
\]

In Section \ref{sss.Ordering} we introduce the \il s $\zo e,m :=\{x\in E: e\le x
<m\}$ \fa $m\in E$ (see Notation \rfa1{n3.42}) where $\le$~is the natural order
of $(E,e,S)$, and define in Theorem \rfa2{t3.6} the {\bf\ca} of these \il s by
setting
\beq39
\#_E(\zo e,m ):=m \qh{\fa}m\in E.
\e
Note that $\zo e,e $ is empty. By \df\ the \ca\ (\wrt $E$) of the empty set
is~$e$. It is customary to denote the \ca\ of the empty set by~$0$. \Mo if
$m:=S(e)$, that is, $S(0)$, then $\zo 0,S(0) =\{0\}$ by~\er{20}. The set~$\{0
\}$ is a singleton, and it is customary to denote the \ca\ of a singleton
by~$1$. \E\Tf $1=\#_E(\{0\}) = \#_E(\zo0,S(0) )=S(0)$. \csq, at the beginning
of Section~\ref{sss.Mult} we require in the sequel that $e$, the distinguished
\el\ of a s.o.n.n.\ $(E,e,S)$, is~$0$, and that $1:=S(0)$. We denote such a
s.o.n.n.\ by $(E,0,S)$. In Lemma \rfa2{l0.27} the \ex\ of such a s.o.n.n.\ is
derived from the \ex\ of a s.o.n.n.\ $(\wt E,\wt e,\wt S)$. As mentioned above,
we denote the \ad\ by~$+$ and the \ca\ by $\#(\ )$. The set $(E,0,S)$ could be
denoted by $\{0,1,1+1,1+1+1,\dots\}$. In the \dr\ of $(E,0,S)$, $1+1$~is
denoted by~$2$, and in the binary \rp ation $1+1=10$ (see Section
\ref{sss.Decimal}). If $A$~is a \nf, then $\#(A):=m\in E^\t$, where $E^\t:=
E\sms0$, $A\approx\zo0,m $ (see Theorem \rfa2{t3.6}). Observe that by~(N1) \te
s \ooo $p\in E$ \st $m=S(p)$ $(=p+1)$, and by Lemma \rfa2{l3.22} $\zo0,p+1
\approx \zo0+1,(p+1)+1 $. Using Lemma \rfa2{l3.4} we obtain $\zo0+1,(p+1)+1 =
[0+1,p+1] = [1,m]$. Hence $A\approx [1,m]$. In \cite[p.~8]{Alg} the \nm~$m$ is
called the \ti{\nm\ of \el s\/} of~$A$, and in \cite[p.~78]{18} is called
the \ti{cardinal \nm\/} of~$A$. Another notation for $\#(A)$ is $|A|$
(which is less musical!).

\advance\abovedisplayskip by-1pt
\advance\belowdisplayskip by-1pt
\advance\abovedisplayshortskip by-1pt
\advance\belowdisplayshortskip by-1pt
In Section \ref{sss.Mult} \E\Pr\ \rfa2{p2.2} \era2{2.3} we collect \pp ies of
\IT s of \el s of an \am\ $(M,\qu,e)$ \wrt the s.o.n.n.\ $(E,0,S)$. We
introduce the (non-standard) notation
\beq40
n \dqu a
\e
for the $n$-fold \IT\ of $a$ in the monoid $(M,\qu,e)$ \wrt $(E,0,S)$. \Mo we
use the notation $S^n$ instead of $\Phi_E(n)$. Thus $S^m\circ S^n=S^{m+n}$
\fa $m,n\in E$ (see \era2{2.4}--\era2{2.7}).

Using notation \er{40}, the {\bf\mlc\ $\cdot$} defined in \cite[p.~18]{Nrs}
becomes:
\beq41
m\cdot n:= n\dpl m, \q n,m\in E,
\e
where $n\dpl m$ is the $n$-fold \IT\ of $m$ \wrt $E$ in the monoid $(E,+,0)$.
We prefer the \ev t \df
\beq42
m\cdot n:= m\dpl n, \q m,n\in E,
\e
given in \era2{2.8}.

We prove in \E\Pr s \rfa2{p2.7} and \rfa2{p2.8} that $(E,\cdot,1)$ is an \am,
that $E^\t$ is a \sbm\ of $(E,\cdot,1)$ and that $(E^\t,+,1)$ is a \PM. \Mo
the \fw\ holds:
\bea43
&0\cdot m = m\cdot 0 = 0 \qh{\fa} m\in E, \\
&k\cdot (l+m) = k\cdot l +k\cdot m, \lb{44} \\
&(k+l)\cdot m = k\cdot m+l\cdot m \qh{\fa} k,l,m\in E. \lb{45}
\e

Finally, we define the {\bf\epc} as follows:
\beq46
m^n := n \ddt m \qh{\fa} m,n\in E,
\e
where $n\ddt m$ is the $n$-fold \IT\ of~$m$ \wrt $(E,0,S)$ in the monoid
$(E,\cdot,1)$. \E\pp ies of the \epc\ are collected in \E\Pr s \rfa2{p2.16}
and \rfa2{p2.17}.

In Section \ref{sss.Repr} the \fw\ formulae are \es ed. Let $A,B$ be finite
sets. Then $A\cup B$, $A\cap B$, $A\t B$ are finite sets and
\bea47
& \#(A\cup B) + \#(A\cap B) = \#(A) + \#(B),\\
& \#(A\t B) = \#(A)\cdot \#(B). \lb{48}
\e
If, in \ad, $A,B$ are nonempty, then $B^A$ is not empty and
\beq49
\#(B^A) = (\#(B))^{\#(A)}.
\e
We invite the reader to take a look at the introductory discussion of Sections
\ref{sss.Comp} and \ref{sss.Decimal}.

In Chapter \ref{s.3} and in the sequel, we use the notation $(\N,0,S)$ for the
s.o.n.n.\ $(E,0,S)$ identified with its ``{\bf\dr}'' $\N=\{0,1,2,3,\dots,9,10,11,
\dots\}$ and $\Na:=\N\sms0$. Using the \df\ of the \mlc\ \er{42}, we reformulate
Theorem \rfa2{t1.38}, called ``{\bf division algorithm}'', as follows.

\textit{Given $a\in\Na$ and $b\in\N$, \te s \ooo pair $(q,r)\in\N\t\N$ \sf
ying}
\beq50
b=q\cdot a+r \qh{\ti{and} } r\in\zo0,a .
\e

A \nm\ $b\in\Na$ is said to be a {\bf\ml e} of~$a$ if $r=0$ in \er{50}. In
this case, the \nm~$a$ is called a {\bf divisor} of~$b$. We recall that both
$(\N,+,0)$ and $(\Na,\cdot,1)$ are\break \PM s. In \E\Pr\ \rfa3{p1.4} it is shown
that if $(X,\qu,e)$ is a \PM, then the \rl\ on~$X$ defined by
\beq51
x\lequ y \qh{if \te s $z\in X$ \st $y=x\qu z(=z\qu x)$ \fa}x,y\in X,
\e
is an \og\ called the \nog\ of $X$. The \fw\ holds:
\bea52
&e\stackrel\qum< x \qh{\fe} x\in X\sms e,\\
&a\stackrel\qum< b \qh{implies} a\qu c \stackrel\qum< b\qu c
\qh{\fa}a,b,c \in X. \lb{53}
\e
If $(X,\qu,e):=(\N,+,0)$, then $y=x+z \nde29 = \Phi(x+z)0 \nde33 = \Phi(z+x)0
\nde30 = (\Phi(z)\circ \Phi(x))0 = \Phi(z)(\Phi(x)0) \nde29 = \Phi(z)x$ \fa
$x,y,z\in X$. Hence the \og\ $\stackrel+\le$ is the same as the \og\ $\le$
defined in~\er{12}. \E\oh if $(X,\qu,e):=(\Na,\cdot,1)$, then $x\stackrel\cdot
\le y$ iff $x$~is a divisor of~$y$ \fa $x,y\in\Na$. Nowadays, the usual
notation for $\stackrel\cdot\le$ is~$|$.

We now prove
\beq54
x|y \hbox{ implies }x\le y \qh{\fa} x,y\in\Na.
\e
Indeed, let $x,y\in\Na$ and $q\in\N$ be \st $y=q\cdot x$. Then $q\ne0$,
otherwise $y=0\cdot x\nde43 = 0\notin\Na$. If $q:=1$, then $y=1\cdot x=x$ since
$1$~is the \nel\ of $(\Na,\cdot,1)$. If $q>1$, then \te s $p>0$ \st $q=1+p$.
Hence $y=(1+p)\cdot x\nde45 = 1\cdot x+p\cdot x = x+z$ with $z\in\Na$, since
both $p$~and~$x$ belong to~$\Na$. \E\Tf\ $x<y$, which completes the proof. If
we denote by~$f$ the bi\jn\ from $(\Na,|)$ onto $(\Na,\le)$ defined by
$f(x):=x$ \fa $x\in\Na$, then $f$~is in\cre\ but $f\Inv$ is \ti{not\/} in\cre.
(We recall that, in contrast, if $f$~is a bi\jc\ monoid-\hm sm from a monoid
$(M,\qu,e)$ onto a monoid $(M',\qu',e')$, then $f\Inv$ is also a \hm sm.)

Chapter \ref{s.3} is devoted to the study of the \PM\ $(\Na,\cdot,1)$. We first
observe that every \nss~$A$ of~$\Na$ is bounded below. Indeed, $1$~is the least
\el\ of~$\Na$ since $x=1\cdot x$ \fa $x\in\Na$. Recall that a subset of
$(\N,\le)$ is bounded (\ev tly, bounded above) iff it is nonempty and finite.
The same holds in $(\Na,|)$. Indeed, if $A\sbs \Na$ is bounded above \wrt $|$,
then so is~$A$ \wrt $\le$ by~\er{54}, hence $A$~is nonempty and finite.
Conversely, if $A$~is nonempty and finite, then the (composite) product of all
\el s of~$A$, denoted by $\prodl_{a\in A}a$ (see Notation \rfa3{n1.9}) is an
\ti{\ub\/} of~$A$ in $(\Na,|)$ since every \el\ of~$A$ is a divisor of
$\prodl_{a\in A}a$ by Lemma \rfa3{l9.10}. The case $\#(A)=2$ is obvious since
$x|xy$, $y|yx$ and $xy=yx$ \fa $x,y\in\Na$. If $A$~is a bounded subset of~$\Na$,
and if $\UB(A)$ denotes the set of all \ub s of~$A$ in $(\Na,|)$, then
$\UB(A)\ne\vn$, hence by \er{22} the set $\UB(A)$ has a least \el\ \wrt $\le$,
which we denote by~$\check a$. \If from \E\Pr\ \rfa3{p1.19} if $\#(A)=2$, and
from \era3{3.23} in Section~\ref{sss.latt} in the general case, that $\check a$
is also the least \el\ of $\UB(A)$ \wrt the \og~$|$.
An \el\ of $\UB(A)$ is called a \ti{common \ml e\/} of the
\el s of~$A$ (or simply of~$A$) and $\check a$ is called the \ti{least common
\ml e\/} of~$A$ ($\lcm(A)$) if we refer to the \og~$\le$, and $\check a$~is
called the {\bf least \ub} of~$A$ ($\lub(A)$) or \ti{supremum\/} of~$A$ ($\sup(A)$)
if we refer to the \og~$|$. We use the notation $\sup A$ for the {\bf supremum}
of~$A$ and $\lcm A$ for the least common \ml e of~$A$, \Ip $\lcm(x,y)$ if
$\#(A)=2$ and $A=\{x,y\}$.

If $x,y\in\Na$ and $x\ne y$, then $1|x$ and $1|y$, hence $1$~is a \ti{common
divisor\/} of the set $\{x,y\}$ (or a \ti{\lo\/} of $\{x,y\}$ \wrt $|$). If
$a\in\Na$ is a common divisor of $(x,y)$, then $a\le x$ by~\er{54} since $a|x$.
As a bounded subset of~$\N$, the set of common divisors of $\{x,y\}$ has a
greatest \el~$\wh a$ \wrt $\le$, called the \textbf{\GCD\/} of $\{x,y\}$ and
denoted by $\gcd(x,y)$. It is shown in \E\Pr\ \rfa3{p1.22} that $\gcd(x,y)$
is also the \GCD\ (or \lo) of $\{x,y\}$ \wrt $|$, called the {\bf infimum\/} of
the set $\{x,y\}$ and denoted by $\inf\{x,y\}$. \Mo we have
\beq55
\gcd(x,y)\cdot \lcm(x,y) = x\cdot y \qh{\fa}x,y\in\Na
\e
(see \E\Pr\ \rfa3{p1.22} \era3{1.25}).

We present {\bf the Euclid algorithm} for computing the $\gcd$ of two \el s
of~$\Na$ in Theorem \rfa3{t8.49}. \Mo we \ch ize in \E\Pr\ \rfa3{p8.50} the
set of \IT s of $\gcd(x,y)$ in the monoid $(\N,+,0)$. This \ch ization allows
us to prove {\bf Euclid's lemma} (Lemma \rfa3{l1.29}) stating that if $a,b,q
\in\Na$ and $q|(a\cdot b)$ then $q|b$ provided that $\gcd(a,q)=1$. Note that
if $a,b,c\in\Na\sms1$ \sf y $a|c$ and $b|c$, then $(a\cdot b)|c$ provided
$\gcd(a,b)=1$ (see Lemma \rfa3{l1.28}). Two \nm s $a,b\in\Na\sms1$ are called
{\bf coprime} or {\bf relatively prime} if $\gcd(a,b)=1$ (see \E\df\
\rfa3{d1.30}).

In Section \ref{sss.Fundament} we prove the {\bf\fd\ theorem of \art}. In the \PM\
$(\N,+,0)$, $1$,~the least \el\ of $\N\sms0$, is the \Gn\ of $(\N,+,0)$. In the
\PM\ $(\Na,\cdot,1)$ there is no least \el\ \wrt $|$ of the set $\Na\sms1$.
Suppose, for \cd ion, that \te s $u\in\Na\sms1$ \st $u|n$ \fa $n\in\Na\sms1$.
Then $u+1=(1\cdot u)+1$ and $0\le 1<u$. Thus by the \uq\ part of the ``division
algorithm'' theorem (see~\er{50}), $u$~is not a divisor of $u+1$. A~\cd ion.
However, the set $\Na\sms1$ possesses {\bf minimal \el s} (see \E\df\
\rfa3{d9.1}). A~minimal \el\ of $\Na\sms1$ is called a {\bf\Pn}. For example,
$D(2)$, the set of divisors of~$2$, is $\{1,2\}$, $D(3)=\{1,3\}$, $D(4)=\{1,2,
4\}$. Thus $2$~and~$3$ are \Pn s and $4$~is called a {\bf\cme \nm} (see
\E\df\ \rfa3{d9.8}). We denote by~$P$ (see Notation \rfa3{n9.7}) the set of
\Pn s. Let $n\in\Na\sms1$ and let $D(n)$ be the set of divisors of~$n$. Then
$1,n\in D(n)$. \If from \E\Pr\ \rfa3{p9.11} that $P$~is \ct y infinite. By
Lemma \rfa1{l4.20} $P$~is unbounded and by Theorem \rfa1{t3.49} \te s a bi\jn\
$\a$ from $(\N,0,S)$ onto $(P,\le)$ which is in\cre. By Lemma \rfa1{l3.34}\,(i),
$\a\Inv$~is also in\cre, and by~(ii) of the same lemma, $\a$~is strictly in\cre.
Suppose that \te s $\b:\N\to P$ bi\jc\ and in\cre, then $\b$~is strictly in\cre\
and $\a\Inv\circ\b:\N\to\N$ is bi\jc\ and in\cre. In view of Exercise
\rfa1{ex3.52} $\a\Inv \circ\b=\id_\N$. Hence $\a=\b$. Thus \te s \ooo strictly
in\cre\ bi\jn~$\a$ from~$\N$ onto~$P$. We have $\a(0)=2$, $\a(1)=3$, $\a(2)=5$,
\dots, and $\bcl_{n\in\N}\a(n)=P$.

Let $n\in\Na\sms1$, and let $D(n)\cap P$ be the set of primes that divide~$n$.
\If from Lemma \rfa3{l9.9} that $D(n)\cap P\ne\vn$. \Mo $D(n)\cap P$ is
bounded in $(\N,\le)$, hence finite, since $p\in D(n)\cap P$ implies $p\le n$.
By  Lemma \rfa3{l9.22} \te s \ooo $m_p\in\Na$ \st $p^{m_p}|n$ and $p^{m_p+1}
\nmid n$ ($p^{m_p+1}$ does not divide~$n$) \fe $p\in D(n)\cap P$. The \fd\
theorem of \art\ states that $n$~is the (composite) product of all $p^{m_p}$
with $p\in D(n)\cap P$ (Theorem \rfa3{t9.22}) and that this \fc\ is unique
(Theorem \rfa3{t9.20}). In the ``\ex'' part of the proof we make use of a
``second form of mathematical \In'' based on the fact that the \nog\ of~$\N$
is a {\bf well-\og} (see \pp y~\er{22}). In the ``\uq'' part, we apply Euclid's
lemma.

In the \os\ $(\Na,|)$ every subset consisting of two \el s possesses a supremum
and an infimum. Such an \os\ is called a {\bf\lt} (see \E\df\ \rfa3{d3.1}).
We use the notation $x\lor y:=\sup\{x\}\cup\{y\}$ \fa $x,y$ belonging to a \lt~$L$
(\Ip $\Na$). Similarly, infimum $\{x\}\cup\{y\}$ is denoted by $x\land y$.
Clearly, $x\lor x=x$ and $x\lor y=y\lor x$ \fa $x,y\in L$. \Mo we have
$(x\lor y)\lor z={\rm supremum}\,\{x\}\cup\{y\}\cup\{z\}$ (see \era3{3.7}) \fa
$x,y,z\in L$, which implies that the binary \op~$\lor$ is \asc e. In Lemma
\rfa3{l8.37} it is shown that if $(X,\qu,e)$ is a \PM\ with \nog\ $\lequ$ for
which $x\lor y$ exists \fa $x,y\in X$, then $(X,\lequ)$ is a \lt\ and
\beq56
(x\lor y)\qu(x\land y) = x\qu y\qh{\fa}x,y\in X
\e
(see \er{55}). We call such a \lt\ an \ti{L-monoid\/} (see \E\df\ \rfa3{d8.34}
and Remark \rfa3{r8.35}). In Theorem \rfa3{t3.19} it is shown that an L-monoid
is \ti{\dsb e\/} (see \era3{3.52} and \era3{3.53}) and that, in particular
in~$\Na$, the \fw\ holds:
\bea57
a\cdot \lcm(x,y) &= \lcm(a\cdot x,a\cdot y),\\
a\cdot \gcd(x,y) &= \gcd(a\cdot x,a\cdot y) \qh{\fa} a,x,y\in\Na. \lb{58}
\e
Using \er{58}, one verifies that, if $x,y\in\Na$ and $u$ (resp.\ $v$) is the
only \el\ of~$\Na$ \st $x=\gcd(x,y)\cdot u$ (usually denoted by $\frac x
{\gcd(x,y)}$) (resp.\ $v=\frac y{\gcd(x,y)}$), then $(\gcd(u,v)=1$. Note that
$\lcm(1,x)=x$ \fa $x\in\Na$. Hence $(\Na,\lor,1)$ is an \am, and $\LOR_{a\in A}
a$ is well-defined by \era2{7.17} for a finite set $A$ of~$\Na$. By Lemma\glossary{$\LOR_{a\in A}a$}
\rfa4{l2.25} we have the \dsb e formula:
\beq59
\Bigl(\LOR_{i\in[1,N]}a_i\Bigr) \land b = \LOR_{i\in[1,N]}(a_i\land b)
\e
\fe finite \sq\ $a:[1,N]\to\Na$, and all $N,b\in\Na$. If in \ad\ $\gcd(a_i,
a_j)=1$ \fa $i,j\in[1,N]$ with $i\ne j$, then we have
\beq60
\prod_{i\in[1,N]}a_i = \lcm\Bigl(\bcl_{i\in[1,N]}a_i\Bigr),
\e
by Lemma \rfa4{l2.26}.

If, in \ad, \te s $c\in\Na$ \st $a_i|c$ \fa $i\in[1,N]$, then
\beq61
\Bg(\prod_{i\in[1,N]} a_i)\Big|c,
\e
by Lemma \rfa4{l2.30}.

We recall that in a \PM\ $(X,\qu,e)$ with \nog~$\lequ$, if $x,y\in X$ \sf y
$x\lequ y$, then \te s \ooo $p\in X$ \st $y=x\qu p$. The \ex\ of~$p$ follows
from \er{51} and the \uq\ from the ``\cnc ity'' \pp y of~$p$. In $(\Na,+,0)$,
the \op\ ``finding~$p$'' is called a {\bf\sbt} and in $(\Na,\cdot,1)$ a
{\bf division}.

There are many conjectures concerning \Pn s. We conclude this part devoted to
\nn s by mentioning a very easily formulated one, namely the {\bf Goldbach
conjecture}\index{Goldbach conjecture} (\cite[p.~377]{Nrs}), which asserts that every even \ig\ greater than~$2$ is the sum
of two primes. For example, $4=2+2$, $6=3+3$, $8=3+5$, $10=5+5$, $12=5+7$,
$14=7+7$, and so~on?

In the second part of this book devoted to prime and finite fields, the \fw\
algebraic \sc\ plays an important role. A~set~$X$ equipped with two binary
\op s $\qu_0,\qu_1$, and containing two distinguished \el s $e_0,e_1$ is called
a \cmt e \sr\ with unity or simply a \ti{\sr\/} if the \fw\ holds:
\bea62
&(X,\qu_0,e_0) \hbox{ is an \am,}\\
&(X,\qu_1,e_1) \hbox{ is an \am,} \lb{63}\\
&e_0\qu_1x = x\qu_1e_0 = e_0, \lb{64}\\
&x\qu_1(y\qu_0z) = (x\qu_1y) \qu_0 (x\qu_1 z) \qh{\fa}x,y,z \in X.\lb{65}
\e
A more general \df\ of a \sr\ can be found in \cite{Wiki}. Note that the \fw\
\cn\ holds:
\beq66
(x\qu_0y)\qu_1z = (x\qu_1z) \qu_0 (y\qu_1z) \qh{\fa}x,y,z\in X,
\e
since $(x\qu_0y)\qu_1z \nde63 = z\qu_1(x\qu_0y)\nde65 = (z\qu_1x)\qu_0(z\qu_1
y)\nde63 = (x\qu_1z)\qu_0(y\qu_1z)$ for $x,y,z\in X$.

The \am\ $(X,\qu_0,e_0)$ is called the \ti{additive\/} monoid of the \sr~$X$
and \am\ $(X,\qu_1,e_1)$ is called the \ti{\mlv\/} monoid of the \sr~$X$. We use
the notation $(X,\qu_0,\qu_1,e_0,e_1)$ for a \sr. In some cases, when no
confusion arises, we use the symbols $+$~and~$0$ (resp.\ $\cdot$~and~$1$)
instead of $\qu_0,e_0$ (resp.\ $\qu_1,e_1$). Usually, the symbol~$\cdot$ is
omitted. We write $ab+c$ (resp.\ $a+bc$) instead of $(ab)+c$ (resp.\ $a+(bc)$).
``\E\mlc\ is performed first.'' If $X,X'$ are \sr s (possibly $X=X'$) and
${f:X\to X'}$, then $f$~is called a (\ti{\sr-})\ti{\hm sm} (or simply a
{\bf (ring-)\hm sm}) if $f$~is a monoid-\hm sm from $(X,\qu_0,e_0)$ into
$(X',\qu_0',e_0')$ and a monoid-\hm sm from $(X,\qu_1,e_1)$ into
$(X',\qu_1',e_1')$. A bi\jc\ ring-\hm sm is called a {\bf ring-\is sm}. The \cm\
of \hm sms is a \hm sm, the inverse of an \is sm is an \is sm and the identity
$\id_X$ is an \is sm. If $X=X'$, then \hm sm is usually replaced by endo\mf\
and iso- by auto-. Note that if $X$~is a \sr, then the map $\d_a:X\to X$,
$a\in X$, defined by
\beq67
\d_a(x):= a\cdot x\ (=x\cdot a) \qh{\fa} x\in X,
\e
is a monoid-endo\mf\ of $(X,\qu_0,e_0)$.

In this part of the book, we are mainly interested in \sr s \st the additive
monoid is a \Cm\ (the ``cancellation law'' \er{34} holds). We call such a \sr\
a\break {\bf c-\sr} (non-standard terminology, see Notation \rfa4{n2.1}). The first
important example of a c-\sr\ is the \sr\ $(\N,+,\cdot,0,1)$. Note that the
additive monoid of~$\N$ is \pn\ with $1$ as \Gn. We call a \sr~$X$ \ti{\pn\/}
if its additive monoid $(X,\qu_0,e_0)$ is \pn\ with $e_1$ as a \Gn\ (i.e.,
$X=I(e_1)$). Thus obviously $\N$~is a \pn\ c-\sr. We recall that an infinite
\pn\ monoid is a \PM, hence a \Cm. It turns out that an \ti{infinite \pn\/}
\sr~$X$ is ring-\is c to $\Nf$. Indeed, by Theorem \rfa4{t1.9n} the map
$\vf:\N\to X$ defined by
\beq68
\vf_n:= n\dquz e_1 \qh{\fa} n\in\N,
\e
is the \ti{only\/} ring-\is sm from $\N$ to~$X$.

In Section \ref{sss.pr.smr} we construct \fa $n\ge2$ a \pn\ c-\sr\ $\N_n$ \st
$\#(\N_n)=n$. It turns out (see Theorem \rfa4{t1.28}) that a finite \pn\ c-\sr\
$X$ is ring-\is c to $\N_n$ if $n=\#(X)$. More precisely, if $n\in\N$ with
$n\ge2$, we define
\beq69
\N_n:=\zo0,n ,
\e
and $\F_n:\N\to\N_n$ by setting
\beq70
\F_n(x):= r \in\N_n,
\e
where $r$ is the only \el\ of~$\N_n$ \sf ying $x=q\cdot n+r$ \fs $q\in\N$ in
view of the ``division algorithm'', Theorem \rfa2{t1.38} and \era2{2.8}.
Clearly, we have
\bga71
\F_n(x)=x \qh{\fa} x\in\zo0,n , \\
\F_n(x+mn)= \F_n(x) \qh{\fa} x,m\in\N. \lb{72}
\e
\Mo we define
\beq72a
x+_ny := \F_n(x+y) \qh{\fa}x,y\in \zo0,n ,
\e
and
\beq73
x\cdot_ny := \F_n(x\cdot y) \qh{\fa}x,y\in \zo0,n ,
\e
In \E\Pr\ \rfa4{p3.15} we show that $(\N_n,+_n,\cdot_n,0,1)$ is a \pn\ c-\sr\
with\break ${\#(\N_n)=n}$ and $\F_n:\Nf \to(\N_n,+_n,\cdot_n,0,1)$ is a sur\jc\
ring-\hm sm. In Theorem \rfa4{t1.28} we show that if $(X,\qu_0,\qu_1,e_0,e_1)$
is a \ti{finite\/} \pn\ c-\sr, then the map $\ov\vf:\zo0,n \to X$ with
$n:=\#(X)$ defined by
\beq74
\ov\vf(m) := m\dpln e_1, \q m\in\zo0,n ,
\e
is the \ti{only\/} ring-\is sm from $\N_n$ onto~$X$. \If that if $X$~and~$X'$
are \ep\ \pn\ c-\sr s, then $X$~and~$X'$ are ring-\is c. At the end of Section
\ref{sss.pr.smr} some divisibility criteria are given. We know that a \po\ \ig\
$m$ is even (that is, $2|m$) iff the last digit of its \dr\ is even. We show,
\Ip that the sum of its digits is divisible by~$3$ iff $3|m$. For example,
$3|123456$ but $3\nmid 1234567$.

Section \ref{sss.indecomp} deals with the problem of factorizing a \pn\ c-\sr\
into a product of \sr s. We first consider the notion of (direct) product of
monoids. In Example \rfa2{xa1.5}\,(iv) we prove that if $(M_i,\qu_i,e_i)$,
$i=1,2$, are monoids (resp.\ \Cm s, \PM s), then $(M_1\t M_2,\qu,(e_1,e_2))$
is a monoid (resp.\ \Cm, \PM) where
\beq75
(a_1,a_2)\qu(b_1,b_2):= (a_1\qu_1b_1,a_2\qu_2b_2) \qh{\fa}a_1,b_1\in M_1
\hbox{ and }a_2,b_2\in M_2.
\e
This monoid is called the (direct) products of the monoids $M_1$ and $M_2$.
Note that if $M_1$~and~$M_2$ are \pn\ monoids, then $M_1\t M_2$ is not
necessarily \pn. The maps $\pi_i:M_1\t M_2\to M_i$, $i=1,2$, defined by
\beq76
\pi_i((a_1,a_2)):=a_i \hbox{ \fa} a_i\in M_i,\q i=1,2,
\e
are called \ti{\pj s\/} on $M_i$, $i=1,2$, and are sur\jc\ \hm sms. If
$M_1\t M_2$ is a \pn\ monoid, then $M_1$~and~$M_2$ are \pn. This is a direct
con\sq\ of \er{38} and of the sur\ji\ of~$\pi_i$, $i=1,2$.

In \E\df\ \rfa4{d2.4} we call a \sr\ {\bf\dcm} if \te\ \sr s $X_1$~and~$X_2$
\st $X$ is ring-\is c to $X_1\t X_2$. Otherwise $X$~is called \ti{in\dcm\/}.
In Lemma \rfa4{l2.5} it is shown that if $X,X'$ are ring-\is c, then $X'$~is
a c-\sr\ (resp.\ \pn, \dcm\ \sr) if $X$~is a c-\sr\ (resp.\ \pn, \dcm\ \sr).

An infinite \pn\ \sr\ (\Ip $\Nf$) is in\dcm\ by \E\Pr\ \rfa4{p2.7}. A~finite
\pn\ c-\sr~$X$ is shown to be in\dcm\ iff $\#(X)=p^m$ \fs \Pn~$p$ and some
\po\ \ig~$m$, in Theorem \rfa4{t2.17}. Its proof makes use of Corollary
\rfa4{c2.15}, stating that if $k,l\in\Na$, then the product of the monoids
$(\N_k,+_k,0)$ and $(\N_l,+_l,0)$ is \pn\ with $(1,1)$ as a \Gn\ iff $k\ne l$
and $\gcd(k,l)=1$. The second main result of Section \ref{sss.indecomp} is
Theorem \rfa4{t2.24}. As a corollary one obtains the Chinese remainder theorem
(see Remark \rfa4{r2.27}). From Theorem \rfa4{t2.24}, Lemma \rfa4{l2.16},
Theorem \rfa4{t2.17} and the \fd\ theorem of \art, we infer Theorem \rfa4{t2.27}:

\ssk
\noindent {\it Every finite \pn\ c-\sr\ is either indecomposable or \is c to
a finite product of indecomposable finite c-\sr s.}

\ssk
The additive monoid of a \pn\ c-\sr~$X$ is \is c to $(\N,+,0)$ where $X$~is
infinite (see Theorem \rfa4{t1.9n}\,(iv)), and to $(\N_n,+_n,0)$ with $n:=
\#(X)$ when $X$~is finite (see Theorem \rfa4{t1.28}). It turns out that
$(\N_n,+_n,0)$, $n\ge2$, is a \ti{\fcg\/} and that $(\N,+,0)$ can be embedded
in an \ti{in\fcg\/}.

The first part of Section \ref{sss.fcg} is devoted to the study of \pp ies of
finite cyclic groups and the second part to the ``embedding'' theorem and the study
of in\fcg s.

Let $(M,\qu,e)$ be a monoid. An \el\ $x\in M$ is called \ti{invertible\/} if
\te s $y\in M$ \st
\beq77
x\qu y = e = y\qu x
\e
(see \E\df\ \rfa4{d3.2}). If such an \el~$y$ exists it is unique and is called
the inverse of~$x$. Clearly $x$~is the inverse of~$y$ if $y$~is the inverse
of~$x$. We denote the set of invertible \el s of~$M$ by~$M^\t$. It turns out
(see Lemma \rfa4{l3.6}) that $M^\t$ is a \sbm\ (possibly trivial) of~$M$.
A~monoid~$X$ \st $X=X^\t$ is called a {\bf group} (see \cite[p.~12]{Groups}),
and an \ag\ if, in \ad, $x\qu y=y\qu x$ \fa $x,y\in X$. Thus we call $M^\t$ the
{\bf subgroup} of~$M$. In Example \rfa4{xa3.7} we consider the monoid $(F^F,\circ,
\id_F)$ where $F$~is a \ns. It is shown that $\Bij(F)$, the set of bi\jc\
self-maps of~$F$ is subgroup of $(F^F,\circ,\id_F)$. \E\oh $M^\t=\{e\}$ iff
$M$ \sf ies \cn~\er{35}. One verifies that if $M_1,M_2$ are monoids, then
\beq78
(M_1\t M_2)^\t = M_1^\t \t M_2^\t
\e
(see Lemma \rfa4{l3.36} and Corollary \rfa4{c3.37}). It turns out that a
finite \Cm\ $(X,\qu,e)$ is a finite \ag\ and conversely. Indeed, let $\th_a:
X\to X$ be the map defined by $\th_a(x):=a\qu x$ ($=x\qu a$), $a,x\in X$. Then
$\th_a$~is in\jc\ in view of~\er{34}, hence also sur\jc\ since $X$~is finite
(Theorem \rfa1{t4.18}\,(iii)). Thus every $x\in X$ has an inverse. Conversely, if
$X$~is an \ag, then $X$ \sf ies \er{33}. \Mo if $x\qu z = y\qu z$, $x,y,z\in X$,
then $x=x\qu e = x\qu(z\qu z\Inv) = (x\qu z)\qu z\Inv= (y\qu z)\qu z\Inv = {y\qu
(z\qu z\Inv)} = y\qu e = y$. A~\sr\ whose additive group is an \ag\ is called
a (\cmt e) {\bf ring} (with unity) (see \cite{Alg}). An invertible \el\ of the
\mlv\ monoid of a ring is called a~{\bf unit} (see \cite[p.~487]{Lattice}).
We use the notation $(\N_n^\t,\cdot_n,1)$ (or $(\N_n,\cdot_n,1)^\t)$ for the
group of units of the ring $(\N_n,+_n,\cdot_n,0,1)$. In contrast to an infinite
\pn\ monoid, a~finite \pn\ \Cm~$X$ with $\#(X)>2$ has more than one \Gn s.
Indeed, if $X=I(a)$, $a\in X\sms e$, then $X=I(a\Inv)$ (see \E\Pr\ \rfa4{p3.26}).
If we denote by ${\rm Gen}(\N_n,+_n,0)$ ($n\ge2$) the set of \Gn s of
$(\N_n,+_n,0)$, then we have
\beq79
{\rm Gen}(\N_n,+_n,0) = (\N_n,\cdot_n,1)^\t,
\e
by \E\Pr\ \rfa4{p3.38}, and
\beq80
{\rm Gen}(\N_n,+_n,0) = \{k\in\zo1,n : \gcd(k,n)=1\},
\e
by Lemma \rfa4{l3.27}. There is no formula for a \Gn\ of $(\N_n,+_n,0)$, but
a formula for the {\bf Euler totient \f}:
\beq89
\F(n):= \#(\{k\in\zo1,n :\gcd(k,n)=1\})
\e
(see Theorem \rfa4{t3.41}).

A nontrivial group $G$ (not necessarily finite) is called {\bf cyclic} if
\te s $a\in G\sms e$ \st
\beq90
G=I(a) \cup I(a\Inv).
\e
\If from \E\Pr\ \rfa2{p1.13} \era2{1.21}, \era2{1.23} that a cyclic group is
abelian. We already proved that a finite \pn\ \Cm~$X$ is a finite \ag\ \st
$X=I(a)$ \fs $a\in X$, that $I(a\Inv)=X$ by \E\Pr\ \rfa4{p3.26}, hence that
$X$~is a \fcg. Conversely, if $X$~is a \fcg, then $X$~is abelian and by~\er{90}
\te s $a\in X$ \st $I(a)$ is a \sbm\ of~$X$. By Lemma \rfa4{l3.60} $I(a)$ is
a group, hence $a\Inv \in I(a)$, as well as every $x\in I(a\Inv)$. Thus
$I(a\Inv)\sbs I(a)$ and $G=I(a)$. \E\Tf $X$ is a finite \pn\ \Cm\ iff
$X$~is a \fcg. We conclude this discussion on \fcg s by mentioning a criterion
for a finite group~$G$ to be cyclic. By Theorem \rfa4{t3.61} if $G$~is a finite
group and $H$ a finite \sbm\ then $\#(H)$ divides $\#(G)$, and by \E\Pr s
\rfa4{p3.20} and \rfa4{p3.24}, if $G$~is a \fcg, then for every divisor~$d$
of~$\#(G)$, $d\ge2$, \te s \ooo cyclic subgroup of~$G$. It turns out that if
a finite group~$G$ possesses for every divisor~$d$ of~$\#(G)$ \ti{at most one}
\pn\ \sbm\ of order~$d$, then $G$~is cyclic (see Corollary \rfa4{c3.72}). This
criterion will be used in Section \ref{sss.pr.fld}.

In Theorem \rfa4{t3.84} it is shown that an infinite \Cm\ can be embedded in
an infinite \ag. In Corollary \rfa4{c3.85}\,(ii) this \ag\ is proved to be
cyclic provided the infinite \Cm\ is \pn. This result is used in Section
\ref{sss.ipr.fld}. The cyclic group in which $(\N,+,0)$ is embedded is the
group $(\Z,+,0)$ of \ig s. This group allows us to define {\bf ``signed'' \IT s}
of an \el\ of a \ti{group} (see \E\df\ \rfa4{d3.94} and \E\Pr\ \rfa4{p3.142}).
Finally, we mention several \ev t \df s of a group (see Theorem \rfa4{t3.52},
\E\Pr s \rfa4{p3.54} and \rfa4{p3.55}).

In Section \ref{sss.pr.fld} we first define a {\bf field} as a (\cmt e) ring (with unity)
$(F,+,\cdot,0,1)$ \st every nonzero \el\ is a unit (see \E\df\ \rfa4{d4.1}).
\If that $F\sms0$ is a \sbm\ of the monoid $(F,\cdot,1)$ and that $(F\sms0,
\cdot,1)$ is a group, called the \ti{\mlv\ group\/} of the field~$F$. Using
\er{79}, \er{80} we find that the ring $(\N_n,+_n,\cdot_n,0,1)$, $n\ge2$, is
a field iff $n$~is a \Pn. We next consider an \ev t \df\ of a field. Let
$(X,+,\cdot,0,1)$ be a \sr. Then $X$~is a field iff \fe $a\in X\sms0$, $b\in X$,
the \eq\ $a\cdot x+b$ is uniquely solvable. We also investigate the solvability
of a system of $n$~\ti{\len s in $n$ \uk s\/} (see \era4{4.15}). The main
result is stated in Theorem \rfa4{t4.21} and Corollary \rfa4{c4.22}.

We introduce the notion of {\bf subfield} of a field in \E\df\ \rfa4{d4.23}.
If $E$~is a subfield of~$F$, then $E$~is a \sbm\ of both the additive group
of~$F$ and the \mlv\ monoid of~$F$, and $(E,+,\cdot,0,1)$ is also a field.
A~subfield~$E$ of a field~$F$ is called \ti{proper\/} if $E\ne F$. A~field with
no proper subfield is called a {\bf\Pf}. In Example \rfa4{xa4.26} we show
that $(\N_p,+_p,\cdot_p,0,1)$ is a \Pf\ \fa primes~$p$. It turns out that every
field~$F$ possesses a unique subfield~$\check F$ which is a \Pf. The field
$\check F$ is called the {\bf\pf} of~$F$ (see \E\Pr\ \rfa4{p4.29}). The main
result about \pf s is stated in Theorem \rfa4{t4.35}. \If that if $(F,+,\cdot,
0,1)$ is a field and the set $\{n\in\Na: n\dpl1=0\}$ is not empty, hence
possesses a least \el\ $\ov n\ge2$ by~\er{22}, then the \pf\ $\check F$ of~$F$
is ring-\is c to $(\N_p,+_p,\cdot_p,0,1)$ with $p=\ov n$. \E\Ip $\ov n$~is
a \Pn. The \nm~$p$ is called the {\bf\ch istic} of the field~$F$ (hence also
of the field~$\check F$), and is denoted by $\chr(F)$. If the set $\{n\in\Na:
n\dpl1=0\}$ of a field~$F$ is empty, then the field~$F$ is said to have \ch
istic~$0$. \E\Ip a \ti{finite field\/} has \ch istic~$p$ \fs prime~$p$.
In Corollary \rfa4{c4.42}\,(ii) we show that the \mlv\ group $(\N_p^\t,\cdot_p,
1)$ is cyclic. In Theorem \rfa4{t4.58} it is shown that the group $(\N_n^\t,
\cdot_n,1)$, $n\ge3$, is cyclic iff $n=4$, $p^k$ or $2p^k$ \fs odd \Pn~$p$ and
some $k\in\Na$. The proof is based on \cite{Rham}.

The aim of the last section of Chapter \ref{s.4} is to construct a field of
``fractions'', the {\bf field of \ra\ \nm s}. We proceed in two ways. In the
first one, using Theorem \rfa4{t3.84} we embed the \PM\ $(\Na,\cdot,1)$ in a
group of ``\po\ fractions'' denoted by $(\Q_{>0},\cdot,1)$, we adjoin~$0$ and
define $0\cdot q=q\cdot 0=0$. Then $\Q_{\ge0}:=\Q_{>0}\cup\{0\}$ becomes an
\am\ $(\Q_{\ge0},\cdot,1)$. We next define an \ad\ on $\Q_{\ge0}$ by \era4{5.11},
\era4{5.12} and deduce that $(\Q_{\ge0},+,\cdot,0,1)$ is a \sr\ whose additive
monoid is a \PM. \If that the \nog\ of $(\Q_{\ge0},+,0)$ is a {\bf total \og}
(see \E\Pr\ \rfa4{p5.4}). Next, we apply Corollary \rfa4{c3.85}\,(i) in order
to embed the \PM\ $(\Q_{\ge0},+,0)$ in a group denoted by $(\wh X,+,0)$.
Finally, we define a \mlc~$\cdot$ on~$\wh X$ by \era4{5.18}, \era4{5.19}. In
Theorem \rfa4{t5.6} we show that $(\wh X,+,\cdot,0,1)$ is a field.

In the second way, we start from the in\fcg\ $(\Z,+,0)$ introduced in \E\Pr\
\rfa4{p3.90}. Then we define, as in~$\N$, a \mlc~$\cdot$ by setting $m\cdot n
:=m\dpl n$ the $m$-th signed \IT\ of~$n$ defined in \era4{3.141}. It turns out
that $(\Z,+,\cdot,0,1)$ is a {\bf domain} (see \E\df\ \rfa4{d5.28}), that is,
a ring \st $\Z\sms0$ is a \sbm\ of the monoid $(\Z,\cdot,1)$ and \st the monoid
$(\Z\sms0,\cdot,1)$ is a \Cm. Finally, we embed this domain in the ``{\bf field of
\qt s}'' of~$D$ as for example in \cite[pp.~40--41]{Alg} (see Theorem
\rfa4{t5.35}). We show that the field~$\wh X$ and the field of \qt s of the domain
$\Zf$ are ring-\is c (see \E\Pr\ \rfa4{p5.37}). These fields are called fields
of \ra\ \nm s, denoted by~$\Q$. We next show that these fields are \ti{infinite
\Pf s\/} and that an infinite \Pf\ is ring-\is c to the field~$\Q$. \Mo we show
that $\Q$~is \ct y infinite and that the field~$\Q$ is an {\bf archimedean
ordered field}. An ordered field has \ch istic zero.

Chapter \ref{s.5} is devoted to the problem of \ex\ and \uq\ of finite fields.
In Section \ref{sss.fdvs} we introduce the important algebraic notion of
{\bf\vsf0} (see \E\df\ \rfa5{d1.1}). As a first example we consider the set
$K^{[1,n]}$ of all maps from the \il\ $[1,n]$, $n\in\Na$, onto a field
$(K,+,\cdot,0,1)$. We define on $K^{[1,n]}$ an \ti{\ad\/} by setting
\bea91
(\bx+\bet)(i)&:= \bx(i)+\bet(i) \qh{\fa} i\in[1,n],\\
{\bf0}(i)&:= 0  \qh{\fa} i\in[1,n]. \lb{92}
\e
Then $(K^{[1,n]},+,{\bf0})$ is an \ag\ (see \era4{2.11}, \era4{2.12},
\era4{2.13} and Corollary \rfa4{c3.37}\,(i)).

We define a \ti{\mlc\ by scalars\/} $(\la,\bx)\mt \la\bx$ from $K\t K^{[1,n]}$
into $K^{[1,n]}$ by setting
\beq93
(\la\bx)(i) := \la\cdot \bx(i) \qh{\fa} i\in[1,n].
\e
It turns out that axioms \era5{1.1} are \sf ied, hence $(K,K^{[1,n]})$ is a
\vs\ over~$K$. The notation $K^n$ is also used for $K^{[1,n]}$.

A map $L$ from a \vs\ $(K,V)$ into a \vs\ $(K,V')$ is called {\bf linear} if
the \fw\ holds:
\beq94
L(a+b) = La+' Lb \qh{\fa}a,b\in V,
\e
and
\beq95
L(\la a)= \la L(a) \qh{\fa}\la\in K \hbox{ and } a\in V.
\e
A \ti{bi\jc} linear map is called a {\bf linear \is sm}. \E\cm s of linear maps
are linear, the inverse of a linear map is linear and the identity map is an
\is sm. The \vs s $(K,V)$ and $(K,V')$ are called \ti{\is c} if \te s a linear
\is sm between $V$~and~$V'$. A~\vs\ $(K,V)$ is called \ti{finite-\dm al\/} if
\te\ $n\in\Na$ and a linear \is sm from $(K,K^n)$ onto $(K,V)$. In Theorem
\rfa5{t1.12} it is shown that $K^n$ and $K^m$ are \is c iff $m=n$. \If that if
$(K,V)$ is finite-\dm al then \te s \ti{exactly one} $n\in\Na$ \st $K^n$ and
$(K,V)$ are \is c. Then $n$~is called the {\bf\dm} of the \vs\ $(K,V)$. Since
$\id_{K^n}$ is an \is sm, the \vs\ $(K,K^n)$ is $n$-\dm al.

A nontrivial \vs\ $(K,V)$ ($V\ne\{0\}$) is called \ti{\fg\/} if \te s a
(nonempty) set of \el s of~$V$, $a_i$, $i\in[1,n]$, $n\in\Na$, \st \fe $x\in V$,
$x=\suml_{i=1}^n \la_ia_i$ \fs $\la:[1,n]\to K$ (see \E\df\ \rfa5{1.15}\,(ii)).
In Theorem \rfa5{t1.29} it is shown that a \fg\ \vs\ is finite-\dm al. \E\Ip
a \ti{finite\/} \vs\ is finite-\dm al.

If $F$ is a field and $E$ is a subfield of~$F$, then $F$~can be viewed as a
\vs\ over the field~$E$ (see Example \rfa5{xa1.2}\,(ii)). \E\Ip if $F$~is a
finite field, then its \pf\ $\check F$ is ring-\is c to the field $(\N_p,+_p,
\cdot_p,0,1)$ \fs prime~$p$. Since $F$~is linear \is c to $\check F{}^n$ \fs
$n\in\Na$, we obtain
$$
\#(F) = \#(\check F{}^n) \nde49 = \#(\check F)^n = p^n.
$$
The first \et y follows from the fact that $F$ and $\check F{}^n$ are \ep.

In the remaining part of Chapter \ref{s.5}, our goal is to show that \fe $p$
prime and $n\in\Na$, \te s a finite field~$F$ \sf ying $\#(F)=p^n$ and that
if $F'$ is a finite field with $\#(F')=p^n$, then $F$~and~$F'$ are ring-\is c.

In Section \ref{sss.formp} we introduce the notion of {\bf \fp}. This notion
is different from the notion of \pl\ used in Analysis. By a \pl\ map of degree
$n\in\N$ we mean a \f\ (map) of the form $p_n(z):=a_0z^n+a_1z^{n-1}+\ldots+
a_n$, $a_0\ne0$, defined on (a nonempty subset of) a field.
For example, the \pl\ map $f(x):=x^p-x$ defined on the field
$(\N_p,+_p,\cdot_p,0,1)$, $p$~prime, is \ti{identically\/} equal to zero.
Indeed, if $x=0$, then $0^p:=p\ddtt p0 =0$ since $p\ne0$. Hence $f(0)=0$.
Recall that the \mlv\ group of the field~$\N_p$ is $(\N_p\sms0,\cdot_p,1)$. \If
from \E\Pr\ \rfa4{p3.65}\,(i) that $\#(\N_p\sms0,\cdot_p,1)\ddtt p x=1$ \fa
$x\in\N_p\sms0$. \E\Tf $x^{p-1}=1$, hence $x^p=x$ \fa $x\in\N_p\sms0$, and
$f(x)=0$ \fa $x\in\zo0,p $. Similarly, one shows (see Section \ref{sss.Galois})
that the \pl\ map $f(x)=x^{(p^n)}-x$ is identically equal to zero on a finite
field~$F$ with $\#(F)=p^n$. The notion of \fp\ is a particular case of the
notion of \ti{formal power series\/} (see for example~\cite{Analysis}). If
$K$~is a field, we denote by $K^\N$ the \vs\ of all \sq s of \el s introduced
in Example \rfa5{xa1.2}\,(i) with $I:=\N$. We denote by $K_0^\N$ (non-standard
terminology) the \lss\ (see \E\df\ \rfa5{d2.6}) of all \lc s (see \E\df\ \rfa5
{d1.15}\,(i)) of \el s of $\{\ve_k\in K^\N: k\in\N\}$ where $\ve_k(i)=1$ for
$k=i$ and $0$~for $k\ne i$. We denote by~$0$ the \el\ of $K^\N$ \sf ying
$0(i)=0$ \fa $i\in\N$ and set $\supp(a):=\{i\in\N: a(i)\ne0\}$ for $a\in K^\N$.
It turns out that if $a\in K^\N_0\sms0$, then $\supp(a)$ is not empty and
bounded. In view of~\er{22} $\supp(a)$ has a greatest \el\ called the
{\bf degree} of~$a$ ($\deg(a)$). In Theorem \rfa5{t2.24} we define a binary
\op~$\cdot$ on~$K_0^\N$ which makes the \vs\ $K_0^\N$ an \asc e, \cmt e
\ti{$K$-algebra with unity\/}~$\ve_0$ (see \E\df\ \rfa5{d2.18}). \E\Ip we have
$\ve_k\cdot \ve_l=\ve_{k+l}$, $k,l\in\N$. \If that $\ve_k=k\ddt \ve_1$, $k\in\N$.
Usually,  $\ve_k$~is denoted by~$X^k$ if $k\in\Na\sms1$, by~$X$ if $k=1$ and
$1$~if $k=0$. The set of \fp s on a field~$K$ is denoted by $K[X]$. \E\Ip the
\fp\ $X^p-X\in K[X]$ with $K:=(\N_p,+_p,\cdot_p,0,1)$ is \ti{not\/} equal
to~$0$. If we ``forget'' the \mlc\ by scalars, then $(K[X],+,\cdot,0,1)$ is
a \cmt e ring with unity. Since for $a,b\in K[X]\sms0$ we have $a\cdot b\ne0$,
$K[X]$~is a~domain. \E\Ip $K[X]$ can be embedded in a field of \qt s by Theorem
\rfa4{t5.35}. If $K=(\N_p,+_p,\cdot_p,0,1)$, then this field of \qt s is an
example of an \ti{infinite\/} field of \ch istic~$p$.

In Theorem \rfa5{t3.1} of Section \ref{sss.Kronecker} we prove an analogue of
Theorem \rfa2{t1.38} (``division algorithm'') for the $K$-algebra $K[X]$.
Given $a\in K[X]\sms0$, $b\in K[X]$, \te s \ooo pair $(q,r)\in K[X]\t K[X]$
\sf ying
\beq96
b=q\cdot a+r, \qh{and either $r=0$ or $\deg(r)<\deg(a)$.}
\e
We define a map $P_a: K[X]\to K[X]$ by setting
\beq97
P_ab:=r,
\e
and we denote by $K_a$ the range of $P_a$:
\beq98
K_a:=\{r\in K[X]: r=0 \hbox{ or }r\ne0,\ \deg(r)<\deg(a)\}.
\e
If $\deg(a)=0$, then $K_a=\{0\}$. \E\Tf we assume
\beq99
a\in K[X]\sms0 \qh{and }\deg(a)\ge1.
\e
\Mo the \fw\ holds: $P_a$ is \idp\ (i.e., $P_a\circ P_a=P_a$), linear and
$K_a$~is a \vs\ over~$K$ with \ad~$+$. On~$K_a$ we define the binary
\op~$\cdot_a$ by
\beq100
u\cdot_a v:= P_a(u\cdot v),\q u,v\in K_a.
\e
It turns out that $(K,K_a)$ is a $\deg(a)$-\dm al \vs\ over~$K$ and that\break
$(K_a,+,\cdot_a,0,\ve_0)$ is a \cmt e ring with unity~$\ve_0$. \Mo
$\zb \ve i {\zo0,\deg(a) }$ forms a basis of the \vs~$K_a$. The first main
result of Section \ref{sss.Kronecker} is Theorem \rfa5{t3.13}, which states
that the ring~$K_a$ is a field iff the \fp~$a$ is {\bf ir\rd} (see \E\df\
\rfa5{d3.11}). The proof of the ``if'' part is taken from \cite[Lemma p.~24]
{Artin}. Note that the \cn\ of irreducibility of~$a$ is the analogue of the
\cn\ ``$n$~is a \Pn'' for the ring $(\N_n,+_n,\cdot_n,0,1)$ to be a field.

It is shown in Lemma \rfa5{l3.14} that if $a$ is the ir\rd\ \fp\ $\ve_1-\a
\ve_0$, $\a\in K$, then
\beq101
P_ab = \Bg(\sum_{i=0}^{\deg(b)} b_i\a^i)\ve_0,
\e
where $\a^i$ is the $i$-th \IT\ of $\a$ in the monoid $(K,\cdot,1)$ and
$b_i\a^i$ is the product of~$b_i$ and~$\a^i$ in the field~$K$. As a corollary,
we find that if $b\in K[X]$ and $\a\in K$, then \te s $q\in K[X]$ \st $b=q\cdot
(\ve_1-\a\ve_0)$ iff
\beq102
\sum_{i=0}^{\deg(b)}b_i\a^i = 0
\e
(see Corollary \rfa5{c3.14}).

The map $\Phi_\a : K[X]\to K$ defined by
\beq103
\Phi_\a(b) := \bca
0 & \hbox{for }b=0,\\
\displaystyle \sum_{i=0}^{\deg(b)}b_i\a^i & \hbox{for }b\in K[X]\sms0,
\eca
\e
is called {\bf evaluation map} at $\a$ \cite{Is} (see \E\df\ \rfa5{d3.15}). If
$b\ne0$ and $\Phi_\a(b)=0$ then $\a\in K$ is called a \ti{root\/} of the
\fp~$b$. If we use the notation $1:=\ve_0$, $X:=\ve_1$ and $X^k:=\ve_k$, $k>1$,
then the \fp\ $X-\a1$ divides~$b$ iff $\a$~is a root of~$b$. \If that an ir\rd\
\fp\ $b\in K[X]$ has no root in~$K$. The converse is not necessarily true.

We now consider the case $a\in K[X]\sms0$, $\deg(a)\ge2$ and $a$ ir\rd. Then
the map $j:K\to K_a$ defined by $j(\a):=\a\ve_0$ \fe $\a\in K$ is shown to be
an in\jc\ ring-\hm sm. \Mo if $\jh:K^\N \to K_a^\N$ is defined by
\beq104
\jh(b)(l) := j(b(l)), \q l\in\N,\ b\in K_0^\N,
\e
then $\jh(K[X])\sbs K_a[X]$ and $\jh|_{K[X]} : K[X]\to K_a[X]$ is an in\jc\
ring-\hm sm (see \E\Pr\ \rfa5{p3.23}). Finally, by Theorem \rfa5{t3.29},
$\ve_1\in K_a[X]$ is a {\bf root} of~$\jh(a)$ in~$K_a$.

Since by Lemma \rfa5{l4.21} every $b\in K[X]\sms0$ with degree greater than one
is the product of an ir\rd\ $a\in K[X]\sms0$ and $c\in K[X]\sms0$, we infer by
Theorem \rfa5{t3.34} that given $b\in K[X]\sms0$ of degree greater than one,
\te s a field~$F$ and an in\jc\ ring-\hm sm $j:K\to F$ \st the \fp\ $\jh(b)
\in F$ has a root in~$F$.

In Theorem \rfa5{t3.40} it is shown that if $K$ is a field of \ch istic zero,
then \te s a ring-\is sm between the ring $(K[X],+,\cdot,0,\ve_0)$ and the ring
of \pl\ maps from~$K$ into itself denoted by $P(K)$ equipped with the pointwise
\ad\ and \mlc.

In Section \ref{sss.Galois} we consider an \ext\ of Theorem \rfa5{t3.34}. We
need some \df s. Let $K_1,K_2$ be fields and $j:K_1\to K_2$ be an in\jc\
ring-\hm sm. Then $(K_2,j)$ is called a {\bf field \ext} of~$K_1$ and $K_1$~is
said to be \ti{embedded\/} in~$K_2$ (see \E\df\ \rfa5{d3.21}). We can
reformulate the con\sq\ of Theorem \rfa5{t3.34} mentioned above by saying that
if $b\in K[X]\sms0$ and $\deg(b)\ge2$, then \te s a field \ext\ $(F,j)$ \st the
\fp\ $\jh(b)\in F[X]$, defined in \er{104}, has a root $\a$ in~$F$. From
Corollary \rfa5{c3.14} we infer that \te s $\wt q\in F[X]$ \st $\jh(b) = \wt q
\cdot (\jh(\ve_1)-\a\jh(\ve_0))$. \Wanp formulate the promised \ext\ of Theorem
\rfa5{t3.34}.

\ssk \noindent
{\it Let $b\in K[X]\sms0$ with $\deg(b)=n\ge1$ and $b(\deg(b))=1$. Then \te s a
field \ext\ $(F,j)$ and a map $\a:[1,n]\to F$ \st
\beq105
\jh(b) = \prod_{i=1}^n (\jh(\ve_1)-\a_i\jh(\ve_0)),
\e
where $\jh$ is defined in \er{104} and $\prodl_{i=1}^n$ is the \cme product in
$(F,\cdot,1)$ $($see Theorem \rfa5{t4.20}$)$. The \fp\ $\jh(b)\in F[X]$ is said to
{\bf split over}~$F$.}

(The proof of Theorem \rfa5{t4.20} is an adaptation of the proof of Theorem
17.8 given in~\cite{Is}.)
\ssk

In a field $(E,+,\cdot,0,1)$ of \ch istic~$p$ the \fw\ formulae hold:
\bea106
(a+b)^{(p^n)} &= a^{(p^n)} + b^{(p^n)}, \\
(a\cdot b)^{(p^n)} &= a^{(p^n)}\cdot b^{(p^n)}, \q a,b\in F,\ p \hbox{ prime},\
n\in\Na \lb{107}
\e
(see \era4{4.16} and \era2{2.3}\,I5).

It is shown in Lemma \rfa5{l4.7} that
\beq108
G:= \{a\in F: a^{(p^n)}=a\}
\e
is a subfield of $F$. It is shown in \E\Pr\ \rfa5{p4.6} that $\#(G)=p^n$. The
proof uses the notion of {\bf\dv} of a \fp\ (see \E\df\ \rfa5{d4.11}, Lemma
\rfa5{l4.13} and \E\Pr\ \rfa5{p4.19}). Finally, in Theorem \rfa5{t4.21} we
prove the \ex\ of a finite field of order~$p^n$. Setting $K:=(\N_p,+_p,
\cdot_p,0,1)$ and $b:=\ve_{p^n}-\ve_1\in K[X]$, we apply Theorem \rfa5{t4.20}
and obtain a field \ext\ $(F,j)$ of~$K$ \st $\jh(\ve_{p^n}-\ve_1)$ splits
over~$F$. \If that $G$~is a subfield of~$F$ \sf ying $\#(G)=p^n$. It turns out
that $G$ is a {\bf splitting field} for the \fp\ $\ve_{p^n}-\ve_1\in K[X]$
(see \E\df\ \rfa5{d4.24}, Lemmae \rfa5{l4.22} and \rfa5{l4.23}). In Theorem
\rfa5{t4.25} it is shown that every nonzero \fp\ $c\in K[X]$, where $K$~is an
arbitrary field, has a splitting field. It is shown in \cite[Theorem 10]{Artin},
\cite[Theorem 17.12]{Is} that if $F,F'$ are splitting fields of a \fp\
$b\in K[X]\sms0$, $K$~arbitrary field, then $F$~and~$F'$ are ring-\is c. We
give a proof of this theorem for the special case of finite $F$~and~$F'$ using
the fact that the \mlv\ group of a finite field is cyclic (see Corollary
\rfa4{c4.42}\,(ii)).

At the end of Section \ref{sss.Galois} we give an exercise showing that the
field $(K_a,+_a,\cdot_a,0,\ve_0)$ introduced in Theorem \rfa5{t3.13} is
ring-\is c to a field of $n\t n$ matrices with entries in~$K$.

%% file: DETOUR1.TEX
\Section{Ordered sets of natural numbers}[Ordered sets of natural numbers]\label{s.1}
\Subsubsection{Sets of \nn s}\label{sss.Sets}
The aim of this section is to introduce axioms for \nn s and some
of their con\sq s, in particular the \fd\ \DRT. Informally, by a \nn\ is meant
an \el\ of the set one, two, three, \dots. Some authors prefer to include the
number zero in the set of \nn s. We shall adopt this terminology, use the
symbol~$\N$ for this set and $\Na$~for the set $\N\sms{\rm zero}$.\glossary{$\N$}\glossary{$\Na$}
There are other notations
for~$\N$ and~$\N^\ast$, \Ip there are authors using the symbols~$\N$ instead
of~$\Na$ and $\N_0,\N_{\ge0}$ instead of~$\N$.
When \nn s are used for ordering (first, second, third, and so
on), it seems preferable to start with the number one.
When \nn s are used for counting (zero, one, two, three, and so on), it
is sometimes preferable  to start with the number zero. Arithmetic
operations on \nn s such as addition, multiplication have been introduced and
used, relying solely on this informal ``definition'' of \nn s. It is only at
the end of the 19th century that mathematicians and logicians started to
provide \df s of the \nn s using the framework of what is called ``set
theory'' \cite[Chapter 1, \S1]{Nrs}. In these notes we shall use freely notions and results from set
theory. As a reference on this subject we mention for instance
\cite[Chapter~14]{Nrs}, \cite[Section 9]{18} and the references therein.

Both $\N$ and $\N^\ast$ are sets, and $\Na:=\N\sm\{{\rm zero}\}$ is a
subset of~$\N$. The \nm\ zero is a distinguished \el\ of~$\N$ and is usually
denoted by~$0$. The notation $:=$ is not a misprint\glossary{$:=$}
but is called the \ti{Rutishauser symbol\/} meaning ``is defined by'' (see
\cite{Analysis}, p.~xi). We shall (almost) systematically use this notation in
the sequel.
An important feature shared by~$\N$ and $\Na$ is the fact that each \el\
possesses an immediate \su. For example, in~$\N$, $1$~is the \su\ of~$0$, two
denoted by~$2$ is the \su\ of~$1$, and so on. Similarly for $\Na$. This can
be formalized in the following way: there exists a map $S:\N\to\Na$, called
the \ti{\su\ \f}\index{successor function}. If $n\in\N$, then $S(n)$, the immediate \su\ of~$n$, is,
informally, the \nm\ in~$\N$ obtained by ``adding'' one. When the addition
in~$\N$ will be defined, $S(n)$ will indeed correspond to $n+1$. We also
introduce the notion of ``\ti{\pd\ \f}'' $P:\Na \to\N$
which, informally, corresponds to the ``subtraction'' of one unit. Therefore
it is natural to assume that the composition of~$S$ and~$P$ is the identity
in~$\N$ and that the composition of~$P$ and~$S$ is the identity in~$\Na$.
In formula we have
\bea1.1
P\circ S&=\id_{\N},\\
S\circ P&=\id_{\Na}. \lb{1.2}
\e
It follows from \er{1.1} and \er{1.2} that $S$ is a bijection from $\N$
onto~$\Na$, $P$~is a bijection from $\Na$ onto~$\N$, and that $P$ is the
inverse of~$S$, $P=S^{\rm inv}$ and $S$ is the inverse of~$P$, $S=P^{\rm
inv}$.

Another way of formulating \er{1.1} and \er{1.2}, without using explicitly the
\pd\ \f~$P$, is the following. Since $\Na$ is a subset of~$\N$, we can view
the map~$S$ as a map from $\N$ into itself. From \er{1.2} it follows that the
`range' of~$S$ is~$\Na$, in notation $R(S)=\Na$. Indeed, by \df\ $R(S)\subset
\Na$ and by \er{1.2} for every $y\in\Na$, we have $S(Py)=y$. Setting
$x:=Py\in\N$, we get $Sx=y$, hence $\Na\subset R(S)$, and $\Na=R(S)$.

Next we use \er{1.1} to show that $S:\N\to\N$ is \ti{injective}. Let
$x,y\in\N$ be \st $Sx=Sy$. Then, since $Sx,Sy\in R(S)=\Na$, we can apply $P$
to $Sx$ and~$Sy$, and obtain from \er{1.1}: $x=P(Sx)=P(Sy)=y$.

\bex1.1
Show that if $S:\N\to\N$ is injective and satisfies $R(S)=\N\sm\{0\}$, then
\te s a unique map $P:\N\sm\{0\}\to \N$ \st $P(Sx)=x$ for every $x\in\N$ and
$S(Py)=y$ for every $y\in\N\sm\{0\}$.
\eex

From what precedes, we infer that \as s \er{1.1} and \er{1.2} are
\ev t to the \fw\ ones. The \su\ \f~$S$, viewed as map from $\N$ into itself,
satisfies
\bea 1.3
{}&S \hbox{ is injective,}\\
&R(S)=\N\sm\{0\}. \lb{1.4}
\e

We now formulate the \df\ of the natural numbers given in \cite[pp.~14--15]{Nrs}.

\bdf1.2
The \ti{\nn s} form a set $\N$,\index{set of \nn s} containing a distinguished \el~$0$, called
\ti{zero}, together with a \ti{\su\ \f} $S:\N\to\N$, which \sf ies the \fw\ axioms:
\begin{align}
& S \hbox{ is injective}, \tag{\rm S1}\\
& 0 \notin R(S), \tag{\rm S2}\\
& \hbox{If a subset $M\subset \N$ contains $0$ and is mapped
into itself by $S$}  \tag{\rm S3}\\
& \hskip100pt \hbox{(i.e. $S(M)$ is included in $M$), then $M=\N$.}\notag
\end{align}

Axioms (S1)--(S3) are \ev t to the Peano axioms \cite[pp.~18--19]{Nrs}\index{Peano's axioms}.
\edf

Clearly (S1) is \er{1.3} and (S2) is implied by \er{1.4}. It turns out
that \er{1.4} is a con\sq\ of (S2) and (S3). Indeed, let $M$ be
the subset of~$\N$ consisting of the union of~$\{0\}$ and~$R(S)$, in formula
$M:=\{0\}\cup R(S)$. Clearly $0\in M$. Moreover, if $x\in M$, then $Sx\in
R(S)\subset M$. As a con\sq\ of (S3) we have $M=\N$. Hence $\N=\{0\}\cup
R(S)$. By~(S2) $\{0\}\cap R(S)=\vn$, hence $\N\sm\{0\}=R(S)$.

\smallskip
We now give an example of a set $A$ including $\N$ and of a map $\td S:A\to A$ \sf
ying (S1) and~(S2) but not (S3).

Set $A:=\N\cup\{a\}$, where $a\notin\N$ and define $\tilde S:A\to A$ by
setting $\tilde Sx=Sx$ for $x\in\N$ and $\tilde Sa=a$. We claim that $\tilde
S$ \sf ies (S1), (S2) but not (S3). Injectivity of~$\tilde S$: let $x,y\in A$
be \st $x\ne y$. We have to show that $\tilde Sx\ne \tilde Sy$. If both $x$
and~$y$ belong to~$\N$,  then $Sx=\tilde Sx$, $Sy=\tilde Sy$, hence $\tilde
Sx\ne \tilde Sy$ follows from the injectivity of~$S$. If $x\in\N$ and $y=a$,
then $\tilde Sx=Sx\in\N$, $\tilde Sy=\tilde Sa=a\notin\N$, so $\tilde Sx\ne
\tilde Sy$.

(S2): Since $a\ne0$ and $0\notin R(S)$, we have $0\notin R(S)\cup\{a\}=
R(\td S)$.

We now show that
(S3) is not \sf ied. Indeed, if $M:=\N$, then $0\in M$ and for every $x\in
M$, we have $\tilde Sx=Sx\in\N=M$, hence $\tilde S$ maps $M$ into itself.
However, $M\ne A$.

\bex1.2
Show that in the above example $R(\tilde S)=A\sm\{0\}$, hence \te s $\tilde
P: A\sm\{0\}\to A$ \st $\tilde P\circ\tilde S=\id_A$ and $\tilde S\circ
\tilde P=\id_{A\sm\{0\}}$ hold. So \er{1.3}--\er{1.4} does not imply (S3).
\eex

Until now we have stated axioms (S1)--(S3) which should be \sf ied by the set
$\N$, the \el\ $0\in\N$ and the \su\ \f~$S$. We shall now \ti{postulate the
existence of such a set\/}. It is interesting to notice that the existence of
such a set follows from the existence of a set $\wt\N$ containing a
distinguished \el\ which we denote by~$0$, and a \su\ \f\ which we denote
by~$\wt S$, \sf ying (S1) and~(S2) alone, i.e.\ $\wt S$~is injective and
$0\notin R(\wt S)$ (in particular $\wt S$ is not surjective). Indeed, let $I$ be
the class of all subsets $\wt M\sbs\wt \N$ \sf ying $0\in \wt M$ and $\wt S(\wt
M)\subset \wt M$. Such subsets will be called \ti{\iv\ subsets} of~$\wt \N$.\index{inductive subset}
Clearly $\wt\N$ is an \iv\ subset of~$\wt\N$. We denote by~$\N$ their
intersection,
$\N:=\bigcap _{\wt M\in I}\wt M$.
One verifies that $\N$ is an \iv\ subset of~$\wt\N$. Therefore $\wt S$
maps $\N$ into itself. Thus its restriction to~$\N$ is well-defined and we
denote it by~$S$. Clearly $S$ \sf ies (S1). Since $0\notin R(\wt S)$, we also
have $0\notin R(S)$ and (S2) is \sf ied. We now show that (S3) is also \sf
ied.

Indeed, if $M\sbs \N$ \sf ies $0\in M$ and $S(M)\sbs M$, then $0\in M$ and
$\wt S(M)\sbs M$, hence $M$~is an \iv\ set of~$\wt\N$, therefore
$\N:=\bigcap _{\wt M\in I}\wt M \sbs M$. Since $M\sbs \N$ by \as, we have
$M=\N$. Therefore (S3) is \sf ied.

It follows that the existence of a set of \nn s is a con\sq\ of the existence
of a set and a map from this set into itself which is in\jc\ but not sur\jc.
We shall define later the notion of finite set and show that any in\jc\ map
from a finite set into itself is also sur\jc, therefore in this sense $\N$ is
not a finite set. For a deeper discussion on the existence of a set of \nn s
in the framework of set theory we refer the reader to the literature, in
particular to~\cite[Chapter~14]{Nrs}.

We now turn to the problem of ``uniqueness'' of a set of \nn s. Assuming that
a set $\N'$ contains a distinguished \el~$0'$ and that a \su\ map $S':
\N'\to \N'$ \sf ies (S1), (S2) and~(S3), which relations can be established
between $\N$ and~$\N'$, $0$~and~$0'$, $S$~and~$S'$\,?

We first consider the \fw\ example for~$\N'$. Let $\N$ be as in \E\df\
\rf{d1.2}. Set $1:=S(0)$. Since $0\notin
R(S)$, $1\ne0$. Set $\N':=\N\sm\{0\}$ (i.e.\ $\N'$ is $\Na$ considered
above), $0':=1$. Since $R(S)=\N'$, the restriction of~$S$ to~$\N'$ is a
well-defined map from $\N'$ into itself which we denote by~$S'$. As a
restriction of the in\jc\ map~$S$, $S'$~is also in\jc, hence \sf ies (S1). We
claim that $1\notin R(S')$. Indeed, suppose for \cd ion that \te s $x\in\N'$
\st $1=S'(x)$. By \df\ of~$S'$, $S'(x)=S(x)$ and by \df\ of~$1$, $1=S(0)$. So
$S(0)=S(x)$ and by injectivity of~$S$ we would get $x=0$, \cd ing $0\notin
R(S)$. This proves (S2). Next we prove (S3) for $(\N',0',S')$. Let $M'$ be an
\iv\ subset of~$\N'$, i.e., $0'\in M'$ and $S'(M')\sbs M'$.
We have to show that $M'=\N'$. Since $M'\sbs\N'$ and $\N'\cap\{0\}=\vn$, we
have $M'\cap\{0\}=\vn$. \Mo using $\N=\N'\cup\{0\}$ we obtain $M'=\N'$ iff\index{iff}
(abbreviation for ``if and only if'', see e.g.\ \cite[p.~5]{18}).
$M'\cup\{0\}=\N'\cup\{0\}=\N$. Hence, by (S3), $M'=\N'$ follows from $0\in
M'\cup\{0\}$ and $S(M'\cup\{0\})\sbs M'\cup\{0\}$, where $M'$ is viewed as
a~subset of~$\N$. But clearly $0\in M'\cup\{0\}$, and $S(M'\cup\{0\})=S(M')
\cup\{S(0)\}=S'(M')\cup\{0'\}\sbs M'\cup M'=M'\sbs M'\cup\{0\}$.
 Therefore (S3) is \sf ied. We have thus shown that $(\Na,1,S')$ also
\sf ies  axioms (S1)--(S3) where $\N$ is replaced by~$\Na$, $0$~by~$1$,
and $S$ by~$S'$.

Motivated by this observation we introduce the \fw\ ``enlarged'' \df\ of a set
of \nn s.

\sdim{(N1)}
\bdf1.3
A set of natural numbers\index{set of \nn s} $(E,e,S)$ is a set $E$ containing a distinguished element $e$
together with a successor function (map) $S: E \to E$ satisfying the following\glossary{$(E,e,S)$}
conditions:

\ssk
\ite[(N1)] {The set $E\setminus\{e\}$ is not empty and the map $S$ is a bijection from $E$ onto $E\setminus\{e\}$.}
\ite[(N2)] {If $M$ is a subset of $E$ containing $e$ and \sf ying $S(M)\sbs M$,
then $M=E$.}
\ssk

A subset $M$ of $E$ \sf ying $e\in M$ and $S(M)\sbs M$ will be called an
\ti{\iv\ subset\/} of $(E,e,S)$.
\edf

We are now in a position to formulate a
``uniqueness'' theorem for sets of \nn s.

\begin{thm}[{\cite[p.~17]{Nrs}}]\lb{t1.4}
Let $(E,e,S)$ and $(E',e',S')$ be two sets of \nn s. Then \te s one and only
one bi\jn\ $\vf:E\to E'$ \sf ying $\vf(e)=e'$ and $\vf\circ S=S'\circ\vf$.
\eth

Note that in the case where $E:=\N$, $e:=0$, $S:=S$ and $E':=\Na$, $e':=1$ and
$S'$~is the \rt ion of~$S$ to~$\Na$, the \f\ $\vf:\N\to\Na$ \sf ies $\vf(x)=Sx$
for every $x\in\N$.

\bex1.5
Let $(E,e,S)$ be a set of \nn s. Let $E'$ be a set and $\vf:E\to E'$ be a
bi\jc\ map. Set $e':=\vf(e)$ and $S':=\vf\circ S\circ\vf^{\rm inv}$. Show
that $(E',e',S')$ is a set of \nn s. \E\Ip if $E':=E\sm\{e\}$,
$e':=S(e)$ and $S'$~is the \rt ion of~$S$ to~$E'$, show that $(E',e',S')$ is a set
of \nn s.
\eex

Theorem \rf{t1.4} is a direct con\sq\ of the main theorem of
this section, namely the Dedekind Recursion Theorem
(1888) \cite[p.~16]{Nrs}.

In order to motivate this theorem we consider a set of \nn s $(E,e,S)$ and
a set~$F$, containing an
\el~$a$, together with a map $f:F\to F$. We are interested in the \sq\ $a$, $f(a)$,
$f(f(a))$, $f(f(f(a)))$, and so on. We want to define a map $\F:E\to F$
\sf ying $\F(e):=a$, $\F(S(e)):=f(a)$, $\F(S(S(e))):=f(f(a))$, and so on. Thus
we have
\[
\bal
\F(S(e))&=f(a)=f(\F(e)),\\
\F(S(S(e)))&=f(f(a))=f(\F(S(e))),
\eal
\]
and more generally we want to have
\beq1.5
\F(S(n)):=f(\F(n)) \qh{for all} n\in E.
\e
So our goal is to define a map $\F:E\to F$ \sf ying the ``recursion formula''
\er{1.5} together with the \df\ of $\F$ at~$e$:
\beq1.6
\F(e):=a.
\e

Dedekind's recursion theorem shows that \te s one and only one such map.\index{theorem!Dedekind}

\begin{thm}[Dedekind, 1888]\lb{t1.6}
Let $(E,e,S)$ be a set of \nn s and let $F$ be a set containing an \el~$a$,
together with a map $f$ of~$F$ into itself. Then there is one and only one
map $\F:E \to F$ \sf ying
\beq1.7
\F(S(n)) = f(\F(n)) \qh{for every} n\in E,
\e
and
\beq1.8
\F(e)=a.
\e
\eth

\brs1.7 \

\hph i,ii, The equalities in \er{1.7} and \er{1.8} are equalities in the set~$F$.

\ssk
\hph ii,i, Formula \er{1.7} can be rewritten as $(\F\circ S)(n)=(f\circ\F)(n)$ for every
$n\in E$, \ev tly
\beq1.9
\F\circ S=f\circ\F,
\e
where the equality in \er{1.9} is an equality in the set $F^E$ which is the
set of all maps from $E$ into~$F$.

\ssk
\hph iii,, Formula \er{1.9} is \ev t to saying that the \fw\ diagram is
``commutative'':
\beq1.10
\xymatrix{
(E,e) \ar[r]^S \ar[d]_\F & (E,e) \ar[d] ^\F\\
(F,a) \ar[r]^f & (F,a)}
\e
\ers

\proof[Proof of Theorem \rf{t1.6}]
\E\wlg we may assume that $F\ne\{a\}$.

\ti{Uniqueness.}
We suppose that $\F$ and $\psi$ are two maps from $E$ into~$F$ \sf ying
\er{1.7}, \er{1.8}. We set $M:=\{n\in E: \F(n)=\psi(n)\}$ and we want to
prove that $M=E$. In view of (N2), it suffices to show that $M$ is an
inductive subset of~$E$. By \er{1.8} we have $\F(e)=a=\psi(e)$, hence $e\in
M$. Let $n\in M$, we have to show that $S(n)\in M$. By \er{1.7}, we have
$\F(S(n))=f(\F(n))$ and $\psi(S(n))=f(\psi(n))$. Since $n\in M$,
$\F(n)=\psi(n)$, hence $\F(S(n))=\psi(S(n))$, and $S(n)\in M$.

\ti{Existence.}
If such a map $\F:E\to F$ exists, its graph is a subset of $E\t F$ which we
denote by~$G$ and which \sf ies the \fw\ conditions: for every $n\in E$

\ssk
\hph i,ii, \te s $x\in F$ \st $(n,x)\in G$;

\ssk
\hph ii,i, for all $y,z\in F$, if $(n,y)\in G$ and $(n,z)\in G$, then $y=z$;

\ssk
\hph iii,, $(e,a)\in G$;

\ssk
\hph iv,, if $x\in F$ is \st $(n,x)\in G$, then $(S(n),f(x))\in G$.

\ssk \noi
Conditions (i), (ii) say that $G$ is the graph of a map from $E$ into~$F$ and
conditions (iii), (iv) say that this map \sf ies \er{1.8} and \er{1.7}.

The idea of the proof is to consider the sets $\cA$ of all subsets of
$E\t F$ which \sf y (iii) and (iv). It turns out that the intersection of all
subsets belonging to~$\cA$ is a subset of $E\t F$ which \sf ies all
conditions (i)--(iv).

Set $\cA:=\{A \sbs E\t F: A$ \sf ies (iii) and (iv), where $G$ is replaced
by~$A\}$. Clearly $E\t F\in \cA$, hence $\cA\ne\vn$. Set
$G:=\bigcap_{A\in\cA} A$. Since $(e,a)\in A$ for every $A\in\cA$, $(e,a)\in
G$ and (iii) holds. \Mo for every $A\in\cA$ and for every $n\in E$, if $x\in F$
and $(n,x)\in A$, then $(S(n),f(x))\in A$. Therefore \fe $n\in E$, \fe $A\in
\cA$, if $x\in F$ and $(n,x)\in A$, then $(S(n),f(x))\in A$.
Thus \fe $n\in E$, if $x\in F$ and $(n,x)\in G$, then $(S(n),f(x))\in\cA$
\fe $A\in \cA$, hence $(S(n),f(x))\in G$. This shows that $G$ \sf ies
condition (iv). So far we have used the fact that properties (iii) and~(iv)
are preserved under intersections.

We now use the fact that $G$ is the intersection of \ti{all\/} subsets of
$E\t F$ \sf ying (iii) and~(iv). The proof that $G$ \sf ies \cn s (i) and~(ii)
will be done by \ti{induction}\index{proof by induction} on~$n$, that is, using~(N2). We define
\[
M:=\{n\in E: \hbox{\te s one and only one }x\in F \hbox{ \st}(n,x)\in G\}.
\]
We want to show that $M=E$, \ev tly that $M$ is an \iv\ set.

\goodbreak
\ti{Beginning of the induction}: $e\in M$.

We clearly have $(e,a)\in G$. It remains to show that $a\in F$ is the only
\el\ \st $(e,a)\in G$. Suppose for \ti{\cd ion} that \te s $b\in F$, $b\ne a$,
\st $(e,b)\in G$. Then $G\sm\{(e,b)\}$ also \sf ies (iii) and~(iv). Indeed,
since $b\ne a$, $(e,a)\in G\sm\{(e,b)\}$, hence $G\sm\{(e,b)\}$ \sf ies~(iii).
We now show that $G\sm\{(e,b)\}$ \sf ies~(iv). Let $(n,x)\in G\sm\{(e,b)\}$.
Since $G\sm\{(e,b)\}\sbs G$, $(S(n),f(x))\in G$. Observe that $(S(n),f(x))\ne
(e,b)$ since by~(N1) there is no $n\in E$ \st $S(n)=e$. (Recall that two
pairs $(u,v)$ and $(\wt u,\wt v)$ are equal iff $u=\wt u$ and $v=\wt v$.) It
follows that $(S(n),f(x))\in G\sm\{(e,b)\}$ and $G\sm\{(e,b)\}$ \sf ies~(iv).
Therefore $G\sm\{(e,b)\}\in\cA$, hence $G=\bigcap_{A\in\cA}A \subset
G\sm\{(e,b)\}$. This implies $(e,b)\notin G$. \ti{A~\cd ion}. \endproof

\ti{\E\iv\ argument\/}: $n\in M$ implies $S(n)\in M$.

Let $n\in M$. We have to prove that $S(n)\in M$, i.e., \te s one and only one
$y\in F$ \st $(S(n),y)\in G$. Since $n\in M$, \te s one and only one $x\in F$
\st $(n,x)\in G$. Since $G$ \sf ies~(iv), $(S(n),f(x))\in G$. It remains to
show that if $(S(n),y)\in G$, then $y=f(x)$. Suppose for \ti{\cd ion} that \te s
$y\in F$, $y\ne f(x)$ \st $(S(n),y)\in G$. As in the proof of the beginning
of the induction, we shall show that $G\sm\{(S(n),y)\}$ \sf ies (iii) and~(iv).
By~(iii), $(e,a)\in G$. By~(N1) there is no $n\in E$ \st $S(n)=e$, hence
$(S(n),y)\ne (e,a)$. It follows that $(e,a)\in G\sm\{(S(n),y)\}$ and that
(iii) is \sf ied by $G\sm\{(S(n),y)\}$. Next we show that $G\sm\{(S(n),y)\}$
\sf ies~(iv), i.e., if $(n,x)\in G\sm\{(S(n),y)\}$, then $(S(n),f(x))\in
G\sm\{(S(n),y)\}$. Since $G\in\cA$, $(S(n),f(x))\in G$. By \as\ $y\ne f(x)$,
hence $(S(n),f(x))\in G\sm\{(S(n),y)\}$. Consequently $G\sm\{(S(n),y)\}\in
\cA$, hence, as above, $G\sbs G\sm\{(S(n),y)\}$, which is possible only if
$(S(n),y)\notin G$. \ti{A \cd ion}. As a con\sq\ of (N2), $M=E$, which proves the
theorem. \Eproof

\bex1.9
Let $U,X,Y,Z$ be \ns s and let $f:U\to X$, $g:X\to Y$, $h:Y\to Z$ be maps.
Let $g\circ f$ denote the \ti{\cm\ of $f$ and $g$} defined by $(g\circ f)(u)
:=g(f(u))$ \fa $u\in U$. Prove that $h\circ(g\circ f)=(h\circ g)\circ f$.
\eex

We are now in a position to prove Theorem \rf{t1.4}.

\proof[Proof of Theorem \rf{t1.4}]
We use Theorem~\rf{t1.6} with $F:=E'$, $a:=e'$ and $f:=S'$ and conclude
that \te s one and only one mapping $\vf:E\to E'$ \st $\vf(e)=e'$ and
$\vf\circ S=S'\circ\vf$ in view of \er{1.9}. Interchanging $E$ and~$E'$ we obtain a
map $\psi:E'\to E$ \st $\psi(e')=e$ and $\psi\circ S'=S\circ\psi$. Next
observe that the map $\psi\circ\vf :E\to E$ \sf ies $(\psi\circ\vf)(e)=e$ and
$(\psi\circ\vf)\circ S=S\circ(\psi\circ\vf)$. Indeed, $(\psi\circ\vf)(e)
:=\psi(\vf(e))=\psi(e')=e$, and $(\psi\circ\vf)\circ S=
\psi\circ(\vf\circ S)=\psi\circ(S'\circ\vf)=(\psi\circ S')\circ\vf=
(S\circ\psi)\circ\vf=S\circ(\psi\circ\vf)$. Note that the map $\id_E$, the\glossary{$\id_E$}
identity in~$E$ also \sf ies $\id_E(e)=e$ and $\id_E\circ S=S\circ\id_E$.
From the uniqueness part of Theorem~\rf{t1.6} with $F:=E$, $e:=e$, and
$f:=S$ we get $\psi\circ\vf=\id_E$. Interchanging $E$ and~$E'$ we
obtain $\vf\circ\psi=\id_{E'}$. It follows that $\vf:E\to E'$ is a bi\jn,
which concludes the proof of Theorem~\rf{t1.4}.
\endproof

\brm1.10
We used several times the \fw\ method of proof. Let $P(n)$ be a statement
depending on an \el~$n$ of a~set of \nn s $(E,e,S)$. In order to prove that
$P(n)$ is true \fa $n\in E$, we set $M:=\{n\in E: P(n)$ is true$\}$. Then we
show that $e\in M$ and $S(M)\subset M$.

This method is usually called a method of proof by \ti{complete induction}
(see~\cite{Nrs}, p.~15) or \ti{mathematical induction}\index{mathematical induction} (see~\cite{Alg}, p.~4).
The proof of $e\in M$ is called the \ti{basis} for induction and the proof of
$S(M)\sbs M$, i.e., if $n\in M$, then $S(n)\sbs M$, \fa $n\in M$, is called
the \ti{induction step}. The \as\ $n\in M$ in the induction step is called the
\ti{induction hypothesis}. We shall introduce in Lemma \rfa3{l9.12} another form of
mathematical induction.
\erm

We conclude this section by giving another application of Theorem \rf{t1.6}.

\blm4.28
Let $X$ be a nonempty set and let $X_0$ be a nonempty subset of~$X$ \st $X_0\ne X$. If
\te s $f:X\to X_0$ in\jc, then \te s a bi\jn\ $h:X\to X_0$.
\elm

\proof
Let $(E,e,S)$ be a set of \nn s. Let $A:E\to \cP(X)$, the set of all subsets
of~$X$ (power set of~$X$)\index{power set}, be defined \rc vely by means of Theorem \rf{t1.6}:
\[
\bal
A(e)&:= X\sm X_0,\\
A(S(n))&:= f(A(n)),
\eal
\]
where $f(A(n))$ denotes the image of $A(n)$ under the map~$f$.
Set $Y:=\bigcup\limits_{n\in E}A(n)$. Define $h:X\to X$ by setting
\[
\bal
h(x)&:=f(x) &&\qh{if} x\in Y,\\
h(x)&:=x &&\qh{if} x\in X\sm Y.
\eal
\]
We first show that $h:X\to X_0$. If $x\in Y$, then $h(x)=f(x)\in X_0$
by \as. Clearly $X\sm X_0=A(e)\sbs Y$, hence $X\sm Y\sbs X_0$. Then,
if $x\in X\sm Y$, $h(x)=x\in X_0$.

We next show that $f(Y)\sbs Y$. Let $y\in Y$. Then there is
$n\in E$ \st $y\in A(n)$. By \df\ $f(y)\in f(A(n))=A(S(n))\sbs Y$.

$h$ \ti{is in\jc}: Let $x,y\in X$, $x\ne y$. If $x,y\in Y$, then $h(x)=f(x)
\ne f(y)=h(y)$ since $f$ is in\jc; if $x,y\in X\sm Y$, then $h(x)=x\ne
y=h(y)$; if $x\in X\sm Y$ and $y\in Y$, then $h(x)=x\in X\sm Y$ and
$h(y)=f(y) \in Y$, hence $h(x)\ne h(y)$.

$h$ \ti{is sur\jc}: Let $z\in X_0$. Since $X\sm Y\sbs X_0$, we have either
$z\in X\sm Y$ or $z\in Y\cap X_0$. If $z\in X\sm Y$, then $h(z)=z$. If $z\in Y\cap
X_0$, then $z\in Y\sm A(e)$. \E\te s $n\in E\sm\{e\}$ \st $z\in A(n)$. By (N1)
\te s $p\in E$ \st $S(p)=n$. Then \te s $x\in A(p)$ \st
$f(x)=z$, by \df\ of $A(S(p))$. Since $x\in A(p)\sbs Y$, $h(x)=f(x)=z$.
\endproof

\bex4.29
Let $(E,e,S)$ be a set of \nn s. Set $X:=E$, $X_0:=E\sms e$, $f:=S$. By (N1)
$f:X\to X_0$ is in\jc. We can apply Lemma \rf{l4.28}. What is $h:X\to X_0$
defined in the proof of Lemma \rf{l4.28}.
\eex

As a corollary of Lemma \rf{l4.28} we have

\begin{thm}[Schr\"oder--Bernstein, \cite{Kelley}, p.~28]\lb{t4.30}\index{theorem!Schr\"oder--Bernstein}
Let $A$ and $B$ be two nonempty sets. If \te\ an in\jc\ map $\vf:A\to B$ and
an in\jc\ map $\psi:B\to A$, then \te s a bi\jn\ from $A$ onto $B$.
\eth

\proof
If $\psi(B)=A$, then $\psi:B\to A$ is bi\jc. Then $h:=\psi\Inv:A\to B$ is a
bi\jn. Suppose $\psi(B)\ne A$. We may
use Lemma \rf{l4.28} with $X:=A$, $X_0:=\psi(B)$ and $f:=\psi\circ\vf$.
Then $f:X\to X_0$ is in\jc. Let $h:X\to X_0$ be as in the conclusion of
Lemma~\rf{l4.28}. Define $\wt\psi:B\to X_0$ by setting $\wt\psi(x):=\psi(x)$,
$x\in B$. Then $\wt\psi$ is a bi\jn\ and $\chi:=\wt\psi{}\Inv\circ h$ is a bi\jn\
from $A$ onto~$B$.
\endproof

\nic{\bex1.14
Let $(E,e,S)$ be a set of \nn s, and let $(E',e',S')$ \sf y {\rm(N1)}. Let $f:E
\to E'$ \sf y\dw

\hph i,ii, $f$ is bi\jc.

\hph ii,i, $f(e)=e'$.

\noi Show that $(E',e',S')$ \sf ies {\rm(N2)} iff it \sf ies

\hph iii,, $S'\circ f=f\circ S'$.
\eex}

\newpage
\Subsubsection{Iterates of the \su\ \f}\label{sss.Iterates}
Let $F$ be a nonempty set and let $f$ be a map from $F$ into itself. We
denote by $F^F$ the set of all maps from~$F$ into itself and we sometimes call an \el\
of~$F^F$ \ti{a self-map}\index{self-map}, notation $f:F\to F$ or $f:F\hk$. If $f,g\in F^F$ then
$f=g$ if $f(x)=g(x)$ \fe $x\in F$, and
the composition $g\circ f$ of $f$ and~$g$ is defined by $(g\circ f)(x):=
g(f(x))$, \fe $x\in F$. The map $g\circ f$ belongs to $F^F$. We recall that the
composition of self-maps is \ti{associative},\index{associativity} i.e.\ if $f,g,h\in F^F$, then
\beq 2.1
h\circ (g\circ f)=(h\circ g)\circ f.
\e
Indeed, for every $x\in F$ we have
\[
(h\circ (g\circ f))(x) = h((g\circ f)(x)) = h(g(f(x))) = (h\circ g)(f(x))
= ((h\circ g)\circ f)(x).
\]
We denote by $\id_F$ the identity map in~$F$, i.e.\ $\id_F(x):=x$, \fe $x\in
F$. Clearly $\id_F \in F^F$ and we have
\beq 2.2
\id_F \circ f=f\circ \id_F = f \hbox{ \fa} f\in F^F.
\e

Let $(E,e,S)$ be a set of \nn s. In view of the \rc on theorem \rf{t1.6}, \fe
$a\in F$, \te s one and only one map $\F_a:E\to F$ \st
\beq 2.3
\F_a(S(n))=f(\F_a(n)),\  n\in E, \qh{and } \F_a(e)=a.
\e
We have $\F_a(S(e))=f(a)$, $\F_a(S(S(e)))=f(f(a))$, and so on. Since
$f(f(a))=(f\circ f)(a)$, $f(f(f(a))) = (f\circ (f\circ f))(a)$, and so on.
Therefore we can define a map $\Phi:E\to F^F$ by setting\glossary{$\Phi(n)$}
\beq 2.4
\Phi(n)(a):=\F_a(n) \qh{\fe $n\in E$ and every} a\in F.
\e
We have $\Phi(e)(a)=\F_a(e)=a$, $a\in F$, hence
\bea 2.5
{}&\Phi(e)=\id_F,\\
&\Phi(S(e))=f, \lb{2.6}\\
&\Phi(S(S(e)))(a) = \F_a(S(S(e))= (f\circ f)(a),\ a\in F,\non
\e
hence
\beq 2.7
\Phi((S\circ S)(e)) = f\circ f.
\e
Note that in view of \er{2.1} $f\circ (f\circ f)=(f\circ f)\circ f$.

It is usual to define
\beq 2.8
f\circ f\circ f:= f\circ (f\circ f)\ ({}=(f\circ f)\circ f),
\e
i.e.\ we can ``omit'' parentheses. Sometimes $f\circ f$ is called the
\ti{two-fold iterate} of the (self-)map $f$, and $f\circ f\circ f$ the
\ti{three-fold iterate} of~$f$.

\bdf 2.1
Let $(E,e,S)$ be a set of \nn s and let
$n\in E$. Then\break $\Phi(n): F\to F$ defined by \er{2.4} is called the
\ti{$n$-fold iterate}\index{n-fold iterate@$n$-fold iterate} of~$f$ \wrt $E$, and $\{\Phi(n)\in F^F: n\in E\}$ is called
the set of \ti{iterates}\index{set of iterates} \wrt $E$ of the self-map~$f$. It will be denoted by
$I_E(f)$, or simply by $I(f)$ when no confusion arises, thus\glossary{$I(f)$}
\beq 2.9
I_E(f):=\{\Phi(n)\in F^F: n\in E\}.
\e
\edf

We now claim that the map $\Phi:E\to F^F$ can be defined without using the
maps $\{\F_a\}_{a\in F}$. Indeed, we apply the \rc on theorem with $F$
replaced by~$F^F$, $a$~by $\id_F$, and the map $f$ by
the map $\wh f:F^F\hk$ defined by
\beq2.10
\wh f(g):=f\circ g \qh{\fe} g\in F^F.
\e

The induction step is based on the \fw\ observation. \E\fa $n\in E$, $a\in F$,
\[
\Phi(S(n))(a) \nde2.4 = \F_a(S(n)) \nde2.3 = f(\F_a(n))
\nde2.4 = f(\Phi(n)(a)) = (f\circ \Phi(n))(a).
\]
In $\nde 2.4 =$ we adopt the notation
used in \cite{Nrs}, Chapter~1, which consists of indicating above the \et y
sign the axiom(s) or reference(s) implying the \et y. This notation will be
used in the sequel.
Hence \fe $n\in E$, $\Phi(S(n)) = f\circ \Phi(n) = \wh f(\Phi(n))$. Therefore
the map~$\Phi$ \sf ies the \rc ve formula
\beq2.11
\Phi(S(n)) = \wh f(\Phi(n)) \qh{for every} n\in E,
\e
where the equality in \er{2.11} is an equality in $F^F$.

Recall that the basis for induction is
\beq2.12
\Phi(e)=\id_F.
\e
We can now invoke the \rc on theorem to prove the existence and uniqueness of
a map $\wh\Phi:E \to F^F$ \sf ying
\beq 2.13
\wh\Phi(S(n)) = \wh f( \Phi(n)) \qh{\fe} n\in E,
\e
and
\beq2.14
\wh\Phi(e) = \id_F.
\e

The uniqueness part of the conclusion of Theorem \rf{t1.6} implies that
$\Phi=\wh\Phi$, which proves the claim \fw\ \er{2.9}. We also claim that
there is another way to define $\Phi$ by induction in the set~$F^F$. Define
$\check f:F^F \hk$ by setting
\beq2.15
\check f(g):=g\circ f \qh{\fe} g\in F^F,
\e
and apply the \rc on theorem to~$F$ replaced by~$F^F$, $a$~by~$\id_F$ and
$f$~by~$\check f$. Set\break $\check \Phi:E\to F^F$ the corresponding map. We now
show that $\wh\Phi=\check \Phi$. From \er{2.11} we have $\Phi(S(n))=\wh f
(\Phi(n)) = f\circ\Phi(n)$, $n\in E$. If we can prove
\beq2.16
f\circ\Phi(n) = \Phi(n)\circ f \qh{\fe}n\in E,
\e
then $\check f(\Phi(n))=\wh f(\Phi(n))$, $n\in E$, and $\wh\Phi=\check \Phi$ follows from
the uniqueness part of the \rc on theorem.

\proof[Proof of $\er{2.16}$]
We proceed by induction on $n$. Set
\[
M:=\{n\in E: f\circ \Phi(n)=\Phi(n)\circ f \}.
\]

\noindent
\ti{Basis for induction}: $e\in M$. Indeed, $\Phi(e)=\id_F$ by \er{2.5} and $e\in M$
by~\er{2.2}.

\ti{Induction step}: $n\in M$ implies $S(n)\in M$. We suppose $f\circ \Phi(n)
=\Phi(n)\circ f$.  We have
\[
f\circ \Phi(S(n))\nad{(1)}= f\circ (f\circ \Phi(n)) = f\circ (\Phi(n)\circ f)
\nad{(2)}= (f\circ \Phi(n))\circ f\nad{(3)}= \Phi(S(n))\circ f.
\]
\E\et ies (1), (3) follow from \er{2.10}, \er{2.11}, and \et y~(2) follows
from \er{2.1}. Hence $S(n)\in M$.

Since $M$ is an \iv\ set in $E$, $M=E$. \endproof

Summarizing, we have proved
\beq2.17
\Phi = \wh \Phi = \check \Phi.
\e
We next show that \ti{all iterates of $f$ commute}, i.e.
\beq2.18
\Phi(m)\circ \Phi(n) = \Phi(n)\circ \Phi(m) \qh{\fa $m,n\in E$.}
\e
Note that \er{2.18} with $m=e$ follows from \er{2.5}, \er{2.2} and
that \er{2.18} with $m=S(e)$ is \er{2.16} in view of~\er{2.6}.

\proof[Proof of $\er{2.18}$]
Let $m\in E$ be fixed and set $M:=\{n\in E: \Phi(m)\circ \Phi(n)=\Phi(n)\circ
\Phi(m)\}$. By \er{2.5} and \er{2.2}, $e \in M$.

Suppose $n\in M$. Then
\bmlg
\Phi(m)\circ \Phi(S(n)) = \Phi(m)\circ (f\circ \Phi(n)) \nad{(1)}=
\Phi(m)\circ (\Phi(n)\circ f)\\ {} \nde2.1 =
(\Phi(m)\circ \Phi(n))\circ f= (\Phi(n)\circ \Phi(m))\circ f \nde2.1 =
\Phi(n)\circ (\Phi(m)\circ f) \nad{(2)}= \Phi(n)\circ (f\circ \Phi(m)) \\{}
\nde2.1 =
(\Phi(n)\circ f)\circ \Phi(m) = (f\circ \Phi(n))\circ \Phi(m) =
\Phi(S(n))\circ \Phi(m).
\e
\E\et ies (1) and (2) follow from \er{2.16}. Hence $S(n)\in M$, and $M$ is an
\iv\ set, which implies $M=E$. \endproof

Finally, we show that \ti{the composition of two iterates of~$f$ is also an
iterate of~$f$\/}, i.e.
\beq2.19
\Phi(m)\circ \Phi(n)\in I(f) \qh{\fa $m,n\in E$.}
\e

\proof[Proof of \er{2.19}]
Let $m\in E$ be fixed and let $M:=\{n\in E: \Phi(m)\circ \Phi(n)\in I(f)\}$.
We have $e\in M$ by \er{2.5} and \er{2.2}. We claim that $n\in M$ implies
$S(n)\in M$. Indeed, suppose $n\in M$, then we have
\[
\Phi(m)\circ \Phi(S(n)) = \Phi(m)\circ (f\circ \Phi(n)) \nad{(1)}=
\Phi(m)\circ (\Phi(n)\circ f) = (\Phi(m)\circ \Phi(n))\circ f.
\]
Since $n\in M$, \te s $p\in E$ \st $\Phi(m)\circ \Phi(n)=\Phi(p)$. Hence
\[
(\Phi(m)\circ \Phi(n))\circ f = \Phi(p)\circ f \nad{(2)} =
f\circ \Phi(p) = \Phi(S(p)).
\]
\E\et ies (1) and (2) follow from \er{2.16}. Hence $S(n)\in M$, and $M$ is an
\iv\ set, which completes the proof. \endproof

We now introduce some \ti{algebraic notions}, which will play an important
role in the sequel.

\bds2.2
Let $X$ be a nonempty set. A map ${\qu}:X\t X\to X$ is called a \ti{law of \cm}
in $X$ (or a \ti{binary \op} on~$X$). The law $\qu$ is called \ti{\asc e} if
for all $a,b,c\in X$ we have
\beq2.20
(a\qu b)\qu c = a\qu (b\qu c).
\e
As a con\sq\ we can omit parentheses in \er{2.20}.

A nonempty set $X$ together with an \asc e \cm\ law~$\qu$ is called a
\ti{semigroup}, and is denoted by $(X,\qu)$.
An \el~$e$ of a \sg\ $(X,\qu)$ \sf ying for every $a\in X$
\beq2.21
a\qu e = e\qu a=a
\e
is called a \ti{\nel\/}\index{neutral element} (or an identity). Note that if such a \nel\ exists,
it is unique. Indeed, let $e$ and~$e'$ be \nel s of $(X,\qu)$. Then $e'=e\qu
e'=e$. A~\sg\ with a \nel\ is called a \ti{monoid}\index{monoid}, and is denoted\glossary{$(X,\qu,e)$}
by $(X,\qu,e)$. The law of \cm\ $\qu$
in~$X$ is called \ti{\cmt e} if for all $a,b\in X$
\beq2.22
a\qu b = b\qu a.
\e
A \sg\ (resp.\ monoid) with a \cmt e law of \cm\ is called an \ti{abelian}\index{monoid!abelian}, or
\ti{\cmt e \sg} (resp.\ \ti{monoid\/}). The monoid $(\{e\},\qu,e)$ is called
the \ti{trivial monoid\/}\index{monoid!trivial}. A nonempty subset $Y$ of a monoid
$(X,\qu,e)$ \sf ying $e\in Y$ and $x\qu y\in Y$ whenever $x,y\in Y$ is
called a \ti{submonoid\/}\index{submonoid} of $(X,\qu,e)$. The \sbm\ $\{e\}$ is called
the \ti{trivial\/} \sbm\ of $(X,\qu,e)$.
\eds

Let $Y$ be a nonempty subset of~$X$, and let $\bqu$ denote the \rt ion of the
binary \op~$\qu$ to $Y\t Y$. Then $\bqu$ is not necessarily a binary \op\ on~$Y$.
It is a binary \op\ on~$Y$ iff $Y$ is a \sbm\ of~$X$. In this case $(Y,\bqu,e)$
is a monoid, since the binary \op\ is \asc e and $e$~is a \nel. It is customary
to still denote $\bqu$ by~$\qu$ (abuse of notation!). With this convention we
can say that if $Y$~is a \sbm\ of the monoid $(X,\qu,e)$, then $(Y,\qu,e)$ is
a~monoid. If $(X,\qu,e)$ is an \am\ and $Y$~is a~\sbm\ of~$X$, then the monoid
$(Y,\qu,e)$ is clearly also abelian. A~\ti{\sbm} $Y$ of a monoid~$(X,\qu,e)$ is
called \ti{abelian if the monoid $(Y,\qu,e)$ is abelian}.
We shall see examples of abelian \sbm s of non\am s. A~monoid is a~\sbm\ of
itself. A~\sbm~$Y$ of a~monoid~$X$ \st $Y\ne X$ is called \ti{proper}.

\goodbreak
\bxa2.3 \

\hph i,ii, Let $F$ be a nonempty set and let $X:=F^F$, the set of all self-maps
of~$F$. Let $\circ$ denote the \cm\ of two maps in~$F^F$.
Then $(F^F,\circ,\id_F)$ is a monoid by \er{2.1}, \er{2.2}.
If $F$ contains at least two \el s, then
$\qu$ is \ti{not \cmt e}. Indeed, let $a,b\in F$, $a\ne b$. Set $f_a(x):=a$,
$f_b(x):=b$ for all $x\in F$. Then $(f_a\circ f_b)(x)=a\ne b=(f_b\circ
f_a)(x)$ for every $x\in F$.

\hph ii,i, Let $f\in F^F$ and let $Y:=I(f)$ be defined by \er{2.9}. The \cm\
of two maps in~$Y$ is again an \el\ of~$Y$ by \er{2.19},
and $\id_F\in Y$ by \er{2.5}. Therefore $I(f)$
is a submonoid of $(F^F,\circ,\id_F)$ and \er{2.18} implies that the monoid
$(I(f),\circ,\id_F)$ is \ti{abelian}.

(iii) Since the \cm\ of two in\jc\ (resp.\ sur\jc, bi\jc) maps in~$F^F$ is
in\jc\ (sur\jc, bi\jc), and since the
identity map in~$F$ is bi\jc, the subset of~$F^F$ consisting of in\jc\ (resp.\
sur\jc, bi\jc) maps is a \sbm\ of $(F^F,\circ,\id_F)$.
\exa

We mention two properties of a self-map which
are inherited by its iterates.

\blm2.4
Let $F$ be a nonempty set, let $f\in F^F$ and let $I(f)$ be the set of
iterates of~$f$ defined by \er{2.9}.

If $f$ is \emph{in\jc} $($resp.\ \emph{sur\jc}$)$, then so is every iterate
of~$f$.
\elm

\proof
First observe that if $f,g\in F^F$ and $f,g$ are both in\jc\ (resp.\ sur\jc),
then $f\circ g$ is also is in\jc\ (resp.\ sur\jc). We suppose that $f$ is
in\jc\ (resp.\ sur\jc). Let $M:=\{n\in E: \Phi(n)$ is in\jc\ (resp.\
sur\jc)$\}$. We have $\Phi(e)=\id_F$ which is both in\jc\ and sur\jc.
Therefore $e\in M$. We claim: $n\in M$ implies $S(n)\in M$. Suppose that $\Phi(n)$
is in\jc\ (resp.\ sur\jc). Then $\Phi(S(n))=f\circ\Phi(n)$ which is also
in\jc\ (resp.\ sur\jc). Therefore $S(n)\in M$. (N2) implies $M=E$.
\endproof

\blm2.5
Let $F$ be a nonempty set, let $f\in F^F$ and let $\Phi:E\to F^F$ be defined
as in \er{2.4}. Then

\hph i,ii, \fe $n\in E$
\beq2.23
R(\Phi(S(n))) \sbs R(\Phi(n)),
\e
where $R(\Phi(n))$ denotes the range of the map $\Phi(n)$ in~$F$, i.e.\
$\Phi(n)(F)$.
\E\Ip if $\Phi(S(n))$ is sur\jc\ for some $n\in E$ then $f$ is sur\jc.

\hph ii,i, If $A$ is a \nss\ of~$F$ which is \inv\ under~$f$, i.e.\ $f(A)\sbs A$, then $\Phi(n)(A)
\sbs A$ \fe $n\in E$. \E\Ip if $f(x)=x$ for some $x\in E$, then $\Phi(n)x=x$
\fe $n\in E$.

\ssk
\hph iii,i, $R(\Phi(S(n)))\sbs R(f)$ \fe $n\in E$.
\elm

\proof\

(i) By induction on $n$. Set $M:=\bigl\{n\in E: R(\Phi(S(n))) \sbs
R(\Phi(n))\bigr\}$. We have $e\in M$ since $R(\Phi(S(e)))\nde2.6 = R(f)\sbs
F=R(\id_F)\nde2.5 = R(\Phi(e))$.
$n\in M$ \ti{implies} $S(n)\in M$: Suppose $R(\Phi(S(n))) \sbs
R(\Phi(n))$. We have $R(\Phi(S(S(n)))) \nad{\er{2.10},\er{2.11}}= R(f\circ \Phi(S(n)))
=(f\circ \Phi(S(n)))(F) = f(\Phi(S(n))(F)) = f(R(\Phi(S(n))))
\nad{(1)}\sbs R(f\circ \Phi(n))= R(\Phi(S(n)))$. Hence $S(n)\in M$. In
inclusion (1) we used $f(A)\sbs f(B)$ if $A,B\sbs F$ and $A\sbs B$. It
follows that $M=E$.
\Eproof

\smallskip
(ii) Set $M:=\{n\in E: \Phi(n)(A)\sbs A\}$. Clearly $e\in M$ since $\Phi(e)
=\id_F$. Suppose $n\in M$. $\Phi(S(n))(A)=(f\circ \Phi(n))(A)=f(\Phi(n)(A))$.
Since $\Phi(n)(A)\sbs A$, we have $f(\Phi(n)(A))\sbs f(A)$. Since $f(A)\sbs A$,
we obtain $\Phi(S(n))(A)\sbs A$, hence $S(n)\in M$. \E\Tf $M=E$.
\Eproof

\smallskip
(iii) Let $y\in R(\Phi(S(n)))$. Then \te s $x\in F$ \st $y=\Phi(S(n))x=
f(\Phi(n)x) \in R(f)$.
\endproof

For convenience we restate the results obtained above in the
\fw\ proposition.

\bpr2.6
Let $F$ be a nonempty set and let $f$ be a map of~$F$ into itself
$($self-map$)$. Let $\circ$ denote the \cm\ of self-maps on~$F$, which we
recall, is \asc e. Let $(E,e,S)$ be a set of \nn s $($see Definition
\rf{d1.3}$)$.

Then \te s one and only one map $\Phi\ (\hbox{resp.\ }\check \Phi)
:E\to F^F$ $($set of self-maps of~$F)$ \sf ying
\bga2.24
\Phi(e)=\id_F\qh{$($resp.\ } \check\Phi(e):=\id_F),\\
\Phi(S(n))\ (\hbox{resp.\ }\check\Phi(S(n))):=f\circ \Phi(n)
\ (\hbox{resp.\ }\check\Phi(n)\circ f) \qh{\fe }n\in E.\lb{2.25}
\e
Moreover, the \fw\ holds{\rm:}
\bea2.26a
{}&\Phi=\check\Phi. \\
{}&\Phi(m)\circ \Phi(n)=\Phi(n)\circ \Phi(m) \qh{\fa} m,n\in E.\lb{2.26}\\
&\hbox{For all $m,n\in E$ \te s at least one $p\in E$ \st}\lb{2.27}\\
&\qquad \Phi(m)\circ \Phi(n)=\Phi(p).\non\\
&\hbox{If $f$ is in\jc\ $($resp.\ sur\jc$)$, so is $\Phi(n)$ \fe} n\in E. \lb{2.28}\\
&R(\Phi(S(n))) \sbs R(\Phi(n)) \qh{\fe}n\in E. \lb{2.30}\\
\noalign{\eject}
&\hbox{If $A$ is a \nss\ of~$F$ which is \inv\ under $f$, i.e.\ $f(A)\sbs A$, then }\lb{2.31}\\
&\qquad \Phi(n)(A)\sbs A \qh{\fe} n\in E.\non\\
&R(\Phi(S(n))) \sbs R(f) \hbox{ \fe \ $n\in E$. \E\Ip if $\Phi(S(n))$ is sur\jc}\lb{2.32}\\
&\qquad \hbox{for some $n\in E$, then $f$ is sur\jc.}\non
\e
\epr

We are now in a position to define the iterates of the \su\ \f\ in a set of
\nn s. Let $(E,e,S)$ be a set of \nn s. Then $S$ can be viewed as an in\jc\
self-map of~$E$ with range $R(S)=E\sm\{e\}$. We \ti{still denote} by $\Phi:E \to
E^E$ the map defined by Proposition \rf{p2.6} where
\bea2.33
\Phi(e)&:= \id_E,\\
\Phi(S(n))&:=S\circ \Phi(n), \q n\in E.\lb{2.34}
\e

In the next proposition we establish specific \pp ies of the iterates of the
\su\ \f~$\Phi$. In what follows we use the simpler notation $\Phi(n)x$ instead
of $\Phi(n)(x)$, ${n,x\in E}$. However, we keep using the notation $S(x)$,
$x\in E$.

\bpr2.7
Let $\Phi$ be the map defined by \er{2.33} and \er{2.34} and
let $m,n,q\in E$. Then
\bea2.35
{}&\Phi(n)e=n,\\
&\Phi \hbox{ is in\jc,} \lb{2.36}\\
&\hbox{\te s one and only one $p\in E$ \st}\lb{2.37}
\Phi(m)\circ \Phi(n)=\Phi(p), \\
&\Phi(m)n=\Phi(n)m, \lb{2.38}\\
&\Phi(m)q=\Phi(n)q \hbox{ implies }m=n, \lb{2.39}\\
&\Phi(S(n))q\ne q, \lb{2.40}\\
&\Phi(m)\circ \Phi(n)=\id_E \hbox{ implies }m=n=e, \lb{2.41}\\
&R(\Phi(S(n))) \cap \{n\}=\vn, \lb{2.42}\\
&S(R(\Phi(m)))=R(\Phi(S(m))), \lb{2.43}\\
&R(\Phi(n))=R(\Phi(S(n)))\cup\{n\}. \lb{2.44}
\e
\epr

\proof\

\er{2.35} by induction. Set $M:=\{n\in E: \hbox{\er{2.35} holds}\}$. Then
$e\in M$ since $\Phi(e)e=\id_E e=e$. Let $n\in M$. Then $\Phi(S(n))e=(S\circ
\Phi(n))e$ by \er{2.34}, hence $\Phi(S(n))e=S(\Phi(n)e)=S(n)$, and $S(n)\in M$.

\er{2.36} Let $m\ne n$. Then by \er{2.35} $\Phi(m)e\ne \Phi(n)e$, hence
$\Phi(m)\ne \Phi(n)$.

\er{2.37} The existence of such a $p\in E$ follows from \er{2.27} and the
uniqueness from \er{2.36}.

\er{2.38} 
$\Phi(m)n\nde2.35 = \Phi(m)(\Phi(n)e) = (\Phi(m)\circ \Phi(n))(e)
\nde2.26 = (\Phi(n)\circ \Phi(m))(e)=\break \Phi(n)(\Phi(m)e) =\Phi(n)m$.

\er{2.39} Suppose $\Phi(m)q=\Phi(n)q$. Then by \er{2.38} we have $\Phi(q)m =
\Phi(q)n$. Since $S$ is in\jc, $\Phi(q)$ is in\jc\ by \er{2.28}, hence $m=n$.

\er{2.40} $S(n)\ne e$ since $e\notin R(S)$, hence
$\Phi(q)S(n)\ne \Phi(q)e$ by \er{2.39}. But $\Phi(q)S(n)=\Phi(S(n))q$ and
$\Phi(q) e=\Phi(e)q$ by \er{2.38}. Then \er{2.40} follows from \er{2.24}.

\er{2.41} In view of \er{2.26} it suffices to prove $m=e$. It follows from
$\Phi(m)\circ\Phi(n)=\id_E$ that $\Phi(m)$ is sur\jc. Suppose for \cd ion
that $m\ne e$. Then $m\in R(S)$ by~(N1), i.e.\ there is $p\in E$ \st $m=S(p)$.
Hence $\Phi(S(p))$ is sur\jc. By \er{2.32},
$S$~is sur\jc, \cd ing (N1). Hence $m=e$.

\er{2.42} Let $M:=\{n\in E: n\notin R(\Phi(S(n)))\}$. We show that $M$ is
\iv. Clearly $e\in M$ since $\Phi(S(e))=S$ and $e\notin R(S)$ by~(N1). We
suppose $n\in M$ and show that $S(n)\in M$. Suppose for \cd ion that
$S(n)\notin M$, i.e.\ $S(n)=\Phi(S(S(n)))p$ for some $p\in E$. By \er{2.34} and \er{2.26a},
$\Phi(S(S(n)))p = \Phi(S(n))S(p)$, hence $S(n)=\Phi(S(n))S(p)$. By \er{2.35},
$S(n)=\Phi(S(n))e$. Since $S$~is in\jc\ by~(N1), we have $\Phi(S(n))$ in\jc\
by~\er{2.28}. Hence $e=S(p)$, which is impossible by~(N1). Hence $S(n)\in M$.

\er{2.43} Let $x\in S(R(\Phi(m)))$. Then \te s $p\in E$ \st $x=S(\Phi(m)p)$.
But $S(\Phi(m)p)=\Phi(S(m))p$ by~\er{2.34}. Hence $x=\Phi(S(m)p)$, and $x\in
R(\Phi(S(m)))$. Let $x\in R(\Phi(S(m)))$, then $x=\Phi(S(m))p$ for some $p\in
E$. Hence by \er{2.34}, $x=S(\Phi(m)p)$ and $x\in S(R(\Phi(m)))$.

\er{2.44} Let $M:=\{n\in E: \hbox{\er{2.44} holds}\}$. Then $e\in M$ since
$R(\Phi(e)) = R(\id_E) = E \nad{(1)}= R(S)\cup \{e\} = R(\Phi(S(e))) \cup
\{e\}$ where \et y~(1) follows from~(N1). We suppose $n\in M$ and claim that
$S(n)\in M$. We have $R(\Phi(n))=R(\Phi(S(n)))\cup\{n\}$. Using
$S(A\cup B) = S(A) \cup S(B)$ \fa $A,B\sbs E$ we get
$S(R(\Phi(n)))=S(R(\Phi(S(n))))\cup\{S(n)\}$. From \er{2.43} we obtain
\[
R(\Phi(S(n)))=R(\Phi(S(S(n))))\cup\{S(n)\}.
\]
Hence $S(n)\in M$, and $M=E$ by (N2).
\endproof

We conclude this section by \es ing a \rl\ between \IT s of a selfmap ${f:F\hk}$
\wrt \ti{two} sets of \nn s.

\bpr2.8
Let $(E,e,S)$ and $(E',e',S')$ be sets of \nn s and let\break ${\vf:E\to E'}$ be
the bi\jc\ map introduced in Theorem \rf{t1.4}. Let $f:F\to F$ where $F$~is
a~\ns. Let $\Phi_E:E\to F^F$ denote the map defined in \er{2.24}, \er{2.25}
and let $\Phi_{E'}:E' \to F^F$ denote the \crs\ map where $E$~is replaced
by~$E'$.

Then
\bga2.45
\Phi_{E'}(\vf(m)) = \Phi_E(m) \qh{\fa} m\in E,\\
I_{E'}(f)=I_E(f). \lb{2.46}
\e
\epr

\proof \

\er{2.45}: Set $M:=\{n\in E: \Phi_{E'}(\vf(n))=\Phi_E(n)\}$. We have $e\in M$.
Indeed, $\Phi_{E'}(\vf(e))=\Phi_{E'}(e')\nde2.24 = \id_F\nde2.24 = \Phi_E(e)$
since $e'=\vf(e)$ by Theorem \rf{t1.4}. Suppose $n\in M$. Then
$\Phi_{E'}(\vf(S(n))) = \Phi_{E'}(S'(\vf(n)))$ by Theorem~\rf{t1.4}. But
$\Phi_{E'}(S'(\vf(n)))\nde2.25 = f\circ \Phi_{E'}(\vf(n)) \nad{n\in M}=
f\circ \Phi_E(n) \nde2.25 = \Phi_E(S(n))$. Hence $S(n)\in M$. \If from (N2)
that $M=E$.

\er{2.46}: Let $g\in I_{E'}(f)$. Then, by~\er{2.9}, \te s $n'\in E'$ \st
$g=\Phi_{E'}(n')$. Since $\vf$~is bi\jc, \te s $n\in E$ \st $n'=\vf(n)$.
Hence $\Phi_{E'}(n') =\break \Phi_{E'}(\vf(n)) \nde2.45 = \Phi_E(n) \in I_E(f)$
by~\er{2.9}. \If that $I_{E'}(f)\sbs I_E(f)$. Interchanging $E$
and~$E'$ we obtain $I_E(f) \sbs I_{E'}(f)$. Then \er{2.46} follows.
\endproof

\brm2.9
\If from \E\Pr\ \rf{p2.8} that the set of \IT s of a self-map~$f$ of a set~$F$
\ti{does not depend\/} of the choice of the set of \nn s used to define the
\IT s of~$f$. \csq, we may use the notation $I(f)$ instead of $I_E(f)$ or
$I_{E'}(f)$.
\erm

\bnt2.10
Let $F$ be a \ns, let $f\in F^F$ and let $(E,e,S)$, $(E',e',S')$ be sets of
\nn s. Then
\beq2.47
I(f):= I_E(f) = I_{E'}(f).
\e
\ent

\newpage

\Subsubsection{Ordering}\label{sss.Ordering}
In this section a ``natural'' ordering on a set of \nn s is introduced.
When one informally describes the set $\N^\ast$ one says
``the set 1, 2, 3, 4 and so on'', not ``the set 1, 3, 4, 2, 5 and so on''.
One implicitly associates a certain ordering with the set~$\N^\ast$. The \IT
s of the \su\ \f\ plays an essential role in the definition of this ordering.

Let $(E,e,S)$ be a set of \nn s. Let $x$ be an \el\ of~$E$. Then $S(x)$ will
be called a \ti{\su} of~$x$.  A~\su\
of a~\su\ of~$x$ will also be called a~\su\ of~$x$.
It is convenient to call~$x$ a \su\ of itself. So every iterate of~$S$
applied to~$x$, will be called a~\su\ of~$x$.

In view of \er{2.46} and \er{2.47} we may introduce the \fw\ relation on the~set~$E$.

\bdf3.1
Let $(E,e,S)$ be a set of \nn s, and let $x,y\in E$.
\beq3.1
\ga
\hbox{$y$ is called a \su\ of~$x$, \ev tly $x$ is called a \pd\ of~$y$,}\\
\hbox{if \te s $T\in I(S)$ \st $y=Tx$. In this case we write $x\P y $.}
\ega
\e
\edf

In the next proposition we show that the relation $\P$ defined in \er{3.1}
possesses \ti{three} \pp ies which characterize an ordering on a set.

\bdf3.2
Let $X$ be a nonempty set and let $\mr R$ be a relation on~$X$ (i.e.\ a
subset of $X\t X$).
The relation R is called

\noi --- \ti{transitive} if \fa $x,y,z\in X$ we have
\beq3.2
x\mr Ry \qhq{and} y\mr Rz \qhq{implies} x\mr R z,
\e
--- \ti{reflexive} if \fe $x\in X$ we have
\beq3.3
x\mr Rx,
\e
--- \ti{\asm} if \fa $x,y\in X$ we have
\beq3.4
x\mr Ry \qhq{and} y\mr Rx \qhq{implies} x=y.
\e

A relation R on a set $X$ which is \tr e, \rfl\ and \asm\ is called an \ti{order
relation}, or an \ti{order}, or an \ti{ordering}\index{ordering}, or a \ti{partial ordering}
on~$X$. A~set~$X$ with an order relation is called a (\ti{partially})
\ti{ordered set\/}. An order relation R on a set~$X$ is usually denoted by\glossary{$\leq$}
$\le$ or $\leq$.
\edf

\brm3.3 \

\hph i,i,
If $(X,\le)$ is an ordered set and $A$ is a nonempty subset of~$X$, then the
\rt ion of the relation $\le$ to $A$ is an order relation on~$A$.

\hph ii,, The \et y \rl\ on a \ns\ is a trivial order \rl. An order \rl~$\le$
on~$X$ is non-trivial if \te\ $x,y\in X$, $x\ne y$, \st $x\le y$.
\erm

\bex3.4 \

\hph i,ii, Let $(X,\le)$ be an ordered set. Set $x\mr{\wt R}y$ if $y\le x$, $x,y\in X$.
Show that $\mr{\wt R}$ is also an order relation. It is sometimes called the
\ti{reverse} order\index{reverse ordering} relation of $\le$ and is denoted by~$\wtle$.

\hph ii,i, Let $Y$ be a nonempty set and let $\cP(Y)$ be the set of all subsets
of~$Y$, the so-called power set of~$Y$. Define
\beq3.5
A\mr RB \qhq{if} A\subset B
\e
where $\sbs$ denotes the inclusion of $A$ in $B$. Show that R is an order
relation on $\cP(Y)$.

\hph iii,, Find a set and a relation R \st
R is \tr e and \rfl\ but not \asm,
R is \tr e and \asm\ but not \rfl,
R is \rfl\ and \asm\ but not \tr e.
\eex

\blm3.5
Let $(E,e,S)$ be a set of \nn s, and let $x,y\in E$. Then
\beq3.6n
x\P y \hbox{ iff \te s $m\in E$ \st} y=\Phi(m)x,
\e
where $\Phi:E\to E$ is defined in \E\Pr\ \rf{p2.6}.
\elm

\proof
In view of \E\df\ \rf{d3.1} and Notation \rf{n2.10}, $x\P y$ iff $y=Tx$
\fs $T\in I_E(S)$ or $y=T'x$ \fs $T'\in I_{E'}(S)$ where $(E',e',S')$ is a set
of \nn s. In the former case \te s $m\in E$ \st $T=\Phi_E(m)$ (see \E\Pr\
\rf{p2.8}), and observe that $\Phi_E=\Phi$ defined in \E\Pr\ \rf{p2.6}.
In the latter case, \te s $m'\in E'$ \st $T'=\Phi_{E'}(m')$ (see
\E\Pr\ \rf{p2.8}), and by \er{2.45} $\Phi_{E'}(m') = \Phi_E(\vf\Inv(m'))$.
Setting $m:=\vf\Inv(m')$ we obtain $T'=\Phi_E(m) = \Phi(m)$.
\endproof

\bpr3.5
Let $(E,e,S)$ be a set of \nn s. Then the relation {\rm P} defined in~\er{3.1} is
an order relation on~$E$, which will be denoted by $\le$.
\epr

\proof \

``\ti{Transitivity}''. Let $x,y,z\in E$ be \st $x\P y$ and $y\P z$. Then \te\
$m,n\in E$ \st $y=\Phi(m)x$ and $z=\Phi(n)y$. It follows that
$z=\Phi(n)(\Phi(m)x) =(\Phi(n)\circ\Phi(m))x=\Phi(p)x$ for some $p\in E$ by
\er{2.37}. Hence $x\P z$.

``\ti{Reflexivity}''. Let $x\in E$. Then $x=\id_E x=\Phi(e)x$ by \er{2.33}.
Hence $x\P x$.

``\ti{Antisymmetry}''. Let $x,y\in E$ be \st $x\P y$ and $y\P x$. Then \te\
$m,n \in E$ \st $y=\Phi(m)x$ and $x=\Phi(n)y$. Hence $x=\Phi(n)(\Phi(m))x =
\Phi(n)\circ \Phi(m)x$. By~\er{2.37} \te s $p\in E$ \st $\Phi(n)\circ \Phi(m)
=\Phi(p)$. Hence $x=\Phi(p)x$. By
\er{2.35} we get $\Phi(e)x=x=\Phi(p)x$. Using
\er{2.39} we obtain $p=e$. Therefore $\Phi(n)\circ \Phi(m)=\Phi(e)=\id_E$
which implies $n=m=e$ by \er{2.41}. Hence $y=\Phi(e)x=x$.
\endproof

\bns3.7
Let $(E,e,S)$ be a set of \nn s.
We shall denote by $(E,e,S,\le)$ the set $(E,e,S)$ together with the \og~$\le$
introduced in \E\Pr\ \rf{p3.5}. \Mo we set\glossary{$<$}\glossary{$\not\le$}\glossary{$\not<$}
\bea3.7
x<y &\qhq{if} x\le y \hbox{ and }x\ne y,\\
x>y &\qhq{if} y\le x \hbox{ and }x\ne y,\lb{3.8}\\
x\not\le y &\qhq{if}\hbox{the proposition $x\le y$ is not true,}\lb{3.9}\\
x\not< y &\qhq{if}\hbox{the proposition $x< y$ is not true.}\lb{3.10}
\e
similarly for $x\not\ge y$ and $x\not>y$.
\ens

Observe that if $x,y,z\in X$, then
\bea3.11
x\le y \hbox{ and }y<z \qhq{implies} x<z,\\
x< y \hbox{ and }y\le z \qhq{implies} x<z.\lb{3.12}
\e

We prove \er{3.11}. \E\tr ity implies $x\le z$. Suppose for \cd ion that $x=z$.
Then $z\le y$, and from $y<z$ we infer $z=y$ by anti\sy y. A~\cd ion since
$y<z$.

Let $(E,e,S,\le)$ be as in Notations \rf{n3.7}.
Then the \el~$e$, the \su\ \f\ as well as its iterates enjoy
special \pp ies which can be described in terms of the ordering.

\blm3.10
Let $(E,e,S,\le)$ be as in Notations \rf{n3.7}. Then
\beq3.13
e<x \qh{\fe} x\in E\sm\{e\} \qh{and } x<S(x) \qh{\fe} x\in E.
\e
Let $x,y\in E$ be \st $x<y$. Then
\bea3.14
&S(x)<S(y),\\
&\Phi(m)x<\Phi(m)y, \qh{\fe} m\in E,\lb{3.15}\\
&\Phi(x)z<\Phi(y)z, \qh{\fe} z\in E,\lb{3.16}\\
&S(x)\le y.\lb{3.17}
\e
\elm

\proof \

\er{3.13}: $x=\Phi(x)e$ by \er{2.35}, hence $e\le x$. Since $x\ne e$, $e<x$ by
\er{3.7}. From \er{2.34}, \er{2.33} and \er{2.2} we obtain $S=\Phi(S(e))$.
Hence $x\le S(x)$, $x\in E$. But \er{2.40} implies $x\ne S(x)$, $x\in E$.

\er{3.14}: Since $x\le y$, there is $p\in E$ \st $y=\Phi(p)x$. Then $S(y) =
S(\Phi(p)x) =(S\circ\Phi(p))x
\nad{(1)}= (\Phi(p)\circ S)x=\Phi(p)S(x)$. Hence $S(x)\le S(y)$.
\E\et y~(1) follows from \er{2.26} and $S=\Phi(S(e))$. Since $x\ne y$ and and
$S$~is in\jc, we have $S(x)\ne S(y)$, hence $S(x)<S(y)$, by \er{3.7}.

\er{3.15}: Let $M:=\{m\in E: \Phi(m)x<\Phi(m)y\}$. $e\in M$ since
$\Phi(e)=\id_E$ and $x<y$. $m\in M$ \ti{implies} $S(m)\in M$: Suppose
$\Phi(m)x<\Phi(m)y$. By~\er{3.14} $S(\Phi(m)x)<S(\Phi(m)y)$. But $S\circ
\Phi(m)= \Phi(S(m))$ by~\er{2.34}, hence $S(m)\in M$. So $M=E$ by~(N2).

\er{3.16}: $\Phi(x)z=\Phi(z)x$, $\Phi(y)z=\Phi(z)y$ by \er{2.38}. Since $x<y$,
$\Phi(z)x<\Phi(z)y$ by  \er{3.15}. Hence \er{3.16} follows.

\er{3.17}: $x<y$ implies $x\le y$, hence there is $p\in E$ \st $y=
\Phi(p)x$. Since $x\ne y$, we have $p\ne e$, hence $p\in R(S)$ by (N1). Thus
there is $q\in E$ \st $p=S(q)$. Then $y=\Phi(p)x=\Phi(S(q))x=\Phi(q)S(x)$, by
\er{2.25}, \er{2.26a}, \er{2.34}. Hence $S(x)\le y$.
\endproof

We now give two definitions motivated by \er{3.13} and \er{3.14}.

\bdf3.11
Let $(X,\le)$ be an ordered set. An \el\ $a\in X$ \sf ying
\beq3.18
a\le x \qh{(resp.\ $x\le a$) \ \fe $x\in X$}
\e
is called a \ti{least \el\/}\index{least element} (resp.\ \ti{greatest \el\/})\index{greatest element} of $X$.
\edf

In view of the antisymmetry of the ordering there is \ti{at most one} least
(resp.\ greatest) \el\ in~$X$. Indeed if $a,b\in E$ are \st $a\le x$, $b\le x$
\fe $x\in E$, then we have $a\le b$ and $b\le a$, hence $a=b$.

\blm3.12
Let $(E,e,S,\le)$ be as in Notations \rf{n3.7}. Then
$(E,\le)$ has no greatest \el.
\elm

\proof
Follows from \er{3.13}.
\endproof

\bdf3.13
Let $(X,\le)$ and $(Y,\le')$ be ordered sets. A~map $f:X\to Y$ is called

\ssk
\hph i,ii, \ti{order-preserving} or \ti{in\cre}\index{increasing map} if $f(a)\le' f(b)$ \fa $a,b\in X$
\st $a\le b$,

\hph ii,i, \ti{strictly in\cre} if $f(a)<' f(b)$ \fa $a,b\in X$ \st $a< b$,

\hph iii,i, \ti{de\cre} (resp.\ \ti{strictly de\cre}) if $f(b)\le' f(a)$ (resp.\
$f(b)<'f(a)$) \fa $a,b\in X$ \st $a<b$.
\edf

\bxa3.14
Let $(E,e,S,\le)$ be as in Notations \rf{n3.7}. Then, by \er{3.15}, \ti{every
\IT} of the \ti{\su\ \f} (\Ip the identity and the \su\ \f\ itself) is
\ti{strictly in\cre}.
\exa

\brm3.15
An in\cre\ \f\ is \ti{strictly in\cre} if it is \ti{in\jc}.
\erm

One may raise the question whether \te s another ordering $\le'$ on $(E,e,S)$
for which $S$ is in\cre. The answer is yes.

\bex3.16
Let $(E,e,S)$ be a set of \nn s. For $x,y\in E$ set
\beq3.19
x\le' y \qh{if \te s $p\in E$ \st} y=(\Phi(p)\circ\Phi(p))x.
\e

\hph i,ii, Show that $\le'$ is an ordering on $E$ \st $S$ is in\cre.

\hph ii,i, Show that $e\not\le' S(e)$, hence $\le'$ is not the same ordering as
$\le$.

\hph iii,, Show that $S(e)\not\le' e$.
\eex

\brm3.16
Some authors call an in\cre\ \f\ nonde\cre\ and a strictly in\cre\ \f\ an
in\cre\ \f.
\erm

In the next lemma we show that the ordering defined in \er{3.6n} possesses an
additional \pp y which is not shared by the ordering $\le'$ of Exercise
\rf{ex3.16}.

\blm3.17
Let $(E,e,S,\le)$ be as in Notations \rf{n3.7}. Then \fa $x,y\in E$ we have
\beq3.20
\hbox{either}\q x=y \qhq{or} x< y \qhq{or} y< x.
\e
\elm

\brm3.18
Not all ordered sets \sf y \er{3.20}. The ordering on $\cP(Y)$ defined
in~\er{3.5} does not \sf y \er{3.20} if the set $Y$ contains at least two
distinct \el s $a$ and~$b$. Indeed, neither $\{a\}\sbs\{b\}$ nor
$\{b\}\sbs\{a\}$. The ordering on $(E,e,S)$ defined in \er{3.19} does not \sf
y \er{3.20} since $e\not\le' S(e)$ and $S(e)\not\le' e$.
\erm

\proof[Proof of Lemma \rf{l3.17}]
Let $a\in E$ and $M:=\{b\in E: a\le b\hbox{ or }b\le a\}$. We have to
show that $M$ is \iv.

$e\in M$ follows from \er{3.13} since $e\le a$.

\ti{$b\in M$ implies $S(b)\in M$}: Let $b\in M$. Then either $a\le b$ or $b<a$. If
$a\le b$, then $S(a)\le S(b)$ by \er{3.14}. \Mo $a\le S(a)$ by
\er{3.13}. Therefore $a\le S(a)$ and $S(a)\le S(b)$ implies $a\le S(b)$, hence
$S(b)\in M$. If $b<a$ then $S(b)\le a$ by \er{3.17}. Hence $S(b)\in M$.
\endproof

\bdf3.19
Let $(X,\le)$ be an ordered set. The ordering $\le$ is called a \ti{total\/}\index{ordering!total}
(or \ti{linear}) \ti{ordering} on~$X$ if \er{3.20} holds \fa $x,y\in X$. An
ordered set with a total ordering is called a \ti{totally ordered set\/} or a
\ti{chain}.
\edf

\bex3.20
Let $(X,\le)$ be a \tos.

\hph i,i, Let $A$ be a nonempty subset of $X$. Show that $(A,\le)$ is also
totally ordered.

\hph ii,, Let $\ge$ be the reverse ordering on $(X,\le)$ (see Exercise
\rf{ex3.4}). Show that $(X,\ge)$ is also totally ordered.

In a \tos\ every subset consisting of two \el s has a least \el. We shall
show in the next proposition that in $(E,e,S,\le)$ \ti{every} subset of~$E$
has a least \el.
\eex

\bdf3.21
An ordered set is said to be \ti{well-ordered\/} if every nonempty subset has
a least \el. The ordering of a well-ordered set is called a
\ti{well-ordering}\index{well-ordering}.
\edf

\bth3.22
Let $(E,e,S)$ be a set of \nn s and let $\le$ be the ordering on~$E$ defined
by~\er{3.6n}. Then $(E,\le)$ is well-ordered, the \el\ $e$ is the least \el\
of $(E,\le)$ and $S:(E,\le)\to(E,\le)$ is strictly in\cre.
\eth

\proof
In view of Lemma \rf{l3.10} it is sufficient to prove that $(E,\le)$ is
well-ordered. Let $A$ be a nonempty subset of~$E$. Set
$L_A:=\{x\in E: x\le y \hbox{ \fe} y\in A\}$. Then $e\in L_A$ by \er{3.13}.
We claim that $L_A\ne E$. Indeed, since $A\ne\vn$, there is $a\in A$. In view
of \er{3.13}, $a<S(a)$. Since the order is total, we have $S(a)\not\le
a$, hence $S(a)\notin L_A$. Therefore $L_A\ne E$. Since $e\in L_A$ and $L_A$
is not \iv, \te s $m\in L_A$ \st $S(m)\notin L_A$ in view of (N2). We claim that
$m\in A$. Suppose for \cd ion that $m\notin A$. Since $m\in L_A$ and $m\notin A$
we have $m<y$ \fe $y\in A$. From \er{3.17} we obtain $S(m)\le y$ \fe $y\in
A$, hence $S(m)\in L_A$, a \cd ion. Hence $m\in A\cap L_A$, that is, $m$ is a
(the) least \el\ of~$A$.
\endproof

\bth3.23
Let $(E,e,S)$ be a set of \nn s. Then the ordering on~$E$ defined by~\er{3.6n}
is the \emph{only} ordering on~$E$ \st\ $S$ is in\cre\ and $e$~is the least
\el\ of~$E$ in this ordering.
\eth

\proof \

(i) Let $\le'$ be an ordering on $E$. Suppose that $x\le y$ implies $x\le' y$
\fa $x,y\in E$. Then the ordering $\le'$ is the same as~$\le$, i.e., $x\le y$
iff $x\le'y$ \fa $x,y\in E$. Indeed, let $x,y\in E$ be \st $x\le' y$. We have
to show that $x\le y$. Without
loss of generality we may assume $x<'y$. Since the ordering $\le$ is total
and $x\ne y$, we have either $x<y$ or~$y<x$. We show that $y<x$ is
impossible. Indeed, if $y<x$, then by \as\ $y\le'x$, hence also $y<'x$ since
$x\ne y$. But this is incompatible with $x<'y$, since $x<'y$ and $y<'x$
implies $x<'x$ by \er{3.11}.

(ii) We show that $x\le y$ implies $x\le' y$. It suffices to prove that
$x\le'\Phi(n)x$ \fa $x,n\in E$. Indeed, in this case, if $x\le y$ with
$x,y\in E$, then by \er{3.6n} \te s $m\in E$ \st $y=\Phi(m)x$. Hence $x\le'
\Phi(m)x=y$. Note that $x=\Phi(x)e$ and $\Phi(n)x=\Phi(x)n$ \fa $x,n\in E$
by~\er{2.35} and~\er{2.38}. Hence, it suffices to show that $\Phi(x)e\le'\Phi(x)n$ \fa
$x,n\in E$. By \as, we have on the one hand, $e\le'n$ \fa $n\in E$, and
\oh $S$~is in\cre\ \wrt $\le'$. Since $S$~is in\jc, $S$~is also strictly
in\cre\ \wrt $\le'$ in view of Remark \rf{r3.15}.
Proceeding as in the proof of \er{3.15}, one
shows that all \IT s of~$S$ are in\cre\ \wrt $\le'$. \csq, $\Phi(x)e \le'
\Phi(x)n$ \fa $x,n\in E$.
\endproof

The uniqueness result established in Theorem \rf{t3.23} motivates the \fw

\bdf3.24
Let $(E,e,S)$ be a set of \nn s. Then the ordering defined by \er{3.6n} is
called the \ti{natural\/} ordering\index{ordering!natural} of $(E,e,S)$.
\edf

At the beginning  of this section we motivated the \df\ of the \og~$\le$ by
introducing the notion of ``\su''. We have $x\le y$ in $(E,e,S)$ iff $y$ is a
``\su'' of~$x$. We now introduce the notion of ``immediate \su''.

\bdf3.25
Let $(X,\le)$ be a chain and let $x\in X$. An \el\ $y\in X$ is called an
\ti{immediate \su}\index{immediate successor} of~$x$ if the \fw\ holds:
\bea3.21
{}&x<y \qh{($y$ is a ``strict \su'' of~$x$)},\\
&\hbox{there is no $z\in X$ \st $x<z<y$ (i.e.\ $x<z$ and $z<y$)}.\lb{3.22}
\e
\edf

\blm3.26
Let $(X,\le)$ be a chain and let $x\in X$. Then $x$ has at most one immediate
\su.
\elm

\proof
Let $x\in X$ and $y,y'\in X$ be immediate \su s of~$x$. Then either $y=y'$ or
$y<y'$ or $y'<y$. The last two cases are impossible since we would have
either $x<y<y'$, hence $y'$ would not be an immediate \su\ of~$x$, or
$x<y'<y$ and $y$ would not be an immediate \su\ of~$x$. Hence $y=y'$.
\endproof

\blm3.27
Let $(E,e,S)$ be a set of \nn s and let $\le$ be its natural order $($defined
by~\er{3.6n}$)$. Then, \fe $x\in E$, $S(x)$ is the immediate \su\ of~$x$.
\elm

\proof
Let $x\in E$. By \er{3.13} $x< S(x)$. Suppose that there is $y\in E$
\st $x<y<S(x)$. Then by \er{3.17} we have $S(y)\le S(x)$. Since $S$~is
strictly in\cre\ by \er{3.14}, we obtain $S(x)<S(y)$. Using \er{3.12} we
arrive at $S(x)<S(x)$, a~\cd ion. Therefore there is no $y\in E$ \st $x<y<S(x)$.
\endproof

In what follows we shall give an \ev t \df\ of a ``set of \nn s with
natural ordering'' formulated in terms of notions based on the ordering.

\bpr3.28
Let $(E,\le,e)$ be a chain \st $e$ is the least \el. Suppose in addition
the existence of a map $S':E\to E$ \st
$S'(x)$ is the \im\su\ of~$x$ \fe $x\in E$. Then
\ssk

\hph i,i, $S'$ is strictly in\cre\ and in\jc,

\hph ii,, $e\notin R(S')$.

\ssk
\noi If, moreover, condition {\rm(N2)} is \sf ied with $S$ replaced by~$S'$, then $(E,e,S')$
is a set of \nn s and its natural ordering is identical to~$\le$, that is,
\fa $x,y\in E$ we have $x\le y$ iff \te s an \IT\ of~$S'$, depending on $x$ and~$y$,
which we denote by~$S''$, \st $y=S''(x)$.
\epr

\brm3.29
If (N2) does not hold, then the intersection of all \iv\ sets of
$(E,e,S')$ \sf ies (N2) (see Section \ref{sss.Sets}).
\erm

\proof \

(i) Let $a,b\in E$ be \st $a<b$. We have to show that $S'(a)<S'(b)$. Since
$(E,\le)$ is a chain either $b<S'(a)$ or $S'(a)\le b$. The first case implies
$a<b<S'(a)$, \cd ing \er{3.22}. Hence $S'(a)\le b$. We have $b<S'(b)$ by
\er{3.21}. Thus $S'(a)<S'(b)$ by \er{3.11}.

Suppose $a\ne b$. Then either $a<b$ or $b<a$. Hence by what precedes either
$S(a)<S(b)$ or $S(b)<S(a)$. In both cases $S(a)\ne S(b)$, hence $S$ is in\jc.

(ii) Suppose \te s $x\in E$ \st $S'(x)=e$. Then by \er{3.21} $x<S'(x)=e$,
hence $x<e$ \cd ing \er{3.18} with $a:=e$.

Next we assume (N2) and prove
(N1). Set $M:=\{e\}\cup R(S')$. Clearly $e\in M$ and $S'(M)\sbs M$ hence $E=M
=\{e\}\cup R(S')$. It follows from~(ii) that $R(S')=E\sm\{e\}$, hence (N1)
holds. So $(E,e,S')$ is a set of \nn s by \E\df~\rf{d1.3}. Let us denote by
$\le'$ its natural ordering defined by \er{3.6n} where $\Phi(m)$ is replaced
by the $m$-fold \IT\ of~$S'$. It follows from Lemma~\rf{l3.10}, where $S$~is
replaced by~$S'$ and $\le$ by~$\le'$, that $e\le'x$ \fe $x\in E$ \er{3.13}
and that $S'$ is strictly in\cre\ \wrt $\le'$ \er{3.14}. Since $\le$ is an
ordering on~$E$ \st $S'$ is in\cre\ by~(i) and \st $e$~is by \as\ the least
\el\ of~$E$ \wrt $\le$, it follows from Theorem \rf{t3.23} that $\le$ and
$\le'$ are identical.
\endproof

Proposition \rf{p3.28} motivates the \fw\ \df.

\bdf3.30
An \ti{ordered set of \nn s}\index{ordered set of natural numbers} $(E,\le,e,S)$ is a chain $(E,\le)$ \st $E\ne
\{e\}$ with a least
\el~$e$ together with a map $S:E\to E$ \sf ying the \fw\ axioms:
\sdim{(OS1)}

\ssk
\ite[(OS1)] {$S(x)$ is the \im\su\ of $x$ \fe $x\in E$,}
\ite[(OS2)] {If a subset $M\sbs E$ contains $e$ and is \inv\ under~$S$ (i.e.\ $S(M)
\sbs M$), then $M=E$.}
\edf

We next summarize some of the previous results of this section.

\bth3.31 \

\hangindent20pt \hangafter1
\noindent\leavevmode\hbox to20pt{\hfil\rm (a)\ }Let $(E,e,S)$ be a set of \nn s $($in the sense of \E\df\
\rf{d1.3}$)$ with natural \og~$\le$ defined in \er{3.6n}. Then $(E,\le,e,S)$
is an \os\ of \nn s.

\hangindent20pt \hangafter1
\noindent\leavevmode\hbox to20pt{\hfil\rm (b)\ }Let $(E,\le,e,S)$ be an ordered set of \nn s.
Then the \fw\ holds{\rm:}

\hangindent20pt \hangafter0 \noi
\hph i,ii, $S$ is strictly in\cre, in\jc\ and $R(S)=E\sm\{e\}$.

\noi\hph ii,i, $(E,e,S)$ is a set of \nn s $($in the sense of \E\df\ \rf{d1.3}$)$ and
$\le$~is the natural ordering of $(E,e,S)$ defined in \er{3.6n}.

\noi\hph iii,, If $\le'$ is an ordering of $E$ \st $e$ is the least \el\ and $S$ is
in\cre\ \wrt $\le'$, then $\le'$ and $\le$ are identical.

\noi\hph iv,, $(E,\le)$ is well-ordered.
\eth

\proof \

(a): $(E,\le)$ is a chain in view of Lemma \rf{l3.17}, $e$~is the least \el\ by
\er{3.13}. Then (OS1) follows from Lemma \rf{l3.27} and (OS2) from~(N2).

(b): (i), (ii) follow from \E\Pr\ \rf{p3.28}, (iii) follows from Theorem
\rf{t3.23}, (iv)~follows from Theorem \rf{t3.22}.
\endproof

\bex3.32
Let $(E,\le)$ be a well-ordered set without greatest \el.

\hph i,i, Show that every \el\ $x\in E$ has a (unique) \im\su~$x'$.

\hph ii,, If $S:E\to E$ is defined by $Sx:=x'$, \fe $x\in E$, and if $e$ denotes
the least \el\ of~$E$, then $(E,\le,e,S)$ is an ordered set of \nn s provided
$S$ \sf ies (OS2).
\eex

We next show that all ordered sets of \nn s are \ois c.

\bdf3.33
Two ordered sets $X$ and $Y$ are called \ti{\ois c}\index{order-isomorphic} if \te s a bi\jn\
$\vf:X\to Y$ \st $\vf$ and its inverse $\vf^{\rm inv}$ are \opr. In that case $\vf$ is
called an \ti{order-\is sm}.
\edf

\goodbreak
\bex3.34
Let $X,Y,Z$ be ordered sets. Show that

\hph i,ii, $X$ is \ois c to itself,

\hph ii,i, if $X$ and $Y$ are \ois c then $Y$ and $X$ are also \ois c,

\hph iii,, if $X,Y$ and $Y,Z$ are \ois c then $X,Z$ are \ois c.
\eex

\blm3.34
Let $X$ be a chain and let $Y$ be an ordered set. Let $\vf:X\to Y$ be
in\cre.

\hph i,i, If $\vf$ is bi\jc, then both $\vf$ and $\vf\Inv$ are strictly in\cre.

\hph ii,, If $\vf$ is strictly in\cre, then $\vf$ is in\jc.
\elm

\proof \

(i)
Let $c,d\in Y$ be \st $c\le d$. Set $a:=\vf^{\rm inv}(c)$ and $b:=\vf^{\rm
inv}(d)$. We have to show that $a\le b$. Suppose for \cd ion that $a\not\le
b$. Then $b<a$ since $X$ is a chain. Since $\vf$ is in\cre\
$\vf(b)\le\vf(a)$. The injectivity of~$\vf$ implies $\vf(b)\ne\vf(a)$ since
$b\ne a$. Then $d=\vf(b)<\vf(a)=c$, hence $d<c$. Finally, $c\le d$ and $d<c$
implies $c<c$ by~\er{3.11}, a~\cd ion.

(ii) Let $a,b\in X$ be \st $\vf(a)=\vf(b)$. We have to show that $a=b$.
Suppose for \cd ion that $a\ne b$. Then either $a<b$ or $b<a$. Since $\vf$ is
strictly in\cre, either $\vf(a)<\vf(b)$ or $\vf(b)<\vf(a)$, a~\cd ion.
\endproof

We now give an example of an \os\ $(X,\le)$, a \ti{chain} $(X',\le')$ and
a~bi\jn\ $\vf:X\to X'$ \st $\vf$ is in\cre\ and $\vf\Inv$ is \ti{not}.

Consider the sets $X=X'=\{a,b,c,d\}$. We define on~$X$ the \rl\ $<$ by
setting $a<b$, $a<c$, $a<d$, $b<d$, $c<d$. Then the \rl\ $\le$ defined by
$x\le y$ iff either $x=y$ or $x<y$, for $x,y\in X$ is a (partial) ordering
on~$X$. We define on~$X'$ the \rl\ $<'$ by setting $a<'b$, $a<'c$, $a<'d$,
$b<'c$, $b<'d$ and $c<'d$. Then the \rl\ on~$X'$ defined by $x'\le'y'$ iff
either $x'=y'$ or $x'<'y'$, for $x',y'\in X'$ is a total ordering on~$X'$.
Let $\vf:X\to X'$ be the identity. Note that $a$ is the least \el\ of
$(X,\le)$ and of $(X',\le')$ and $d$~is the greatest \el\ of $(X,\le)$ and
$(X,\le')$. Then we have $x\le y$ iff $x\le' y$ whenever $x=a$ or $y=d$.
Since the pairs $(b,c)$ and $(c,b)$ do not belong to the \rl~$\le$, the map
$\vf$ is in\cre. However, since $b\le' c$ and $\vf\Inv(b)=b\not\le
c=\vf\Inv(c)$, the map $\vf\Inv$ is not in\cre.

We have the \fw\ diagram:
\[
\ga X\\
\xymatrix{&d& \\ b\ar[ur] && c\ar[ul] \\ & a\ar[ul] \ar[ur] \ar[uu] }\ega
\qquad \qquad
\ga X' \\ \xymatrix{d\\c \ar[u]\\b\ar[u]\\a\ar[u] }\ega
\]
where arrows mean $\le$ (resp.\ $\le'$).

\bth3.36
Let $X,Y$ be two ordered sets of \nn s. Then the bi\jn\ $\vf:X \to Y$
introduced in Theorem \rf{t1.4} is an \ois sm.
\eth

\proof
Let $e$ (resp.\ $e'$) be the least \el\ of~$X$ (resp.\ $Y$) and let~$S$ (resp.\
$S'$) be the \su\ \f\ in~$X$ (resp.\ $Y$). Then $(X,e,S)$ and $(Y,e',S')$ are
sets of \nn s. Therefore by Theorem~\rf{t1.4} \te s a bi\jn\ $\vf:X\to Y$
\st $\vf(e)=e'$ and $\vf\circ S=S'\circ\vf$. In view of Lemma
\rf{l3.34} it is sufficient to show that $\vf$ is in\cre\ since $(X,\le)$
is totally ordered. By induction one shows from $\vf\circ S=S'\circ\vf$ that
$\vf\circ \Phi(m)=\wt\Phi{}'(m)\circ\vf$ where $\Phi(m)$ ($\wt\Phi{}'(m)$) is the
$m$-fold \IT\ of~$S$ (resp.\ $S'$) \wrt $X$ (see \E\df~\rf{d2.1}).
Therefore if $a,b\in X$ are \st $a\le b$,
then, by \er{3.6n}, \te s $m\in X$ \st $b=\Phi(m)a$. Hence $\vf(b)=\vf
(\Phi(m)a)= \wt\Phi{}'(m)\vf(a)$. Since $I_X(S')=I_{X'}(S')$ by~\er{2.46},
\te s $m'\in X'$ \st $\wt\Phi{}'(m)=\Phi'(m')$ where $\Phi'(m')$ denotes the
$m'$-fold \IT\ of~$S'$ \wrt $X'$. \If that $\vf(b)=\Phi'(m')\vf(a)$.
Hence $\vf(a)\le'\vf(b)$ where $\le'$ is the ordering
in~$Y$. Therefore $\vf$ is \opr.
\endproof

We now investigate some con\sq s of the fact that an ordered set of \nn s is
well-ordered.

\bdf3.37
Let $(X,\le)$ be an ordered set and let $A$ be a nonempty subset of~$X$.
An~\el\ $l$ (resp.~$u$) of $X$ is called a \ti{lower bound\/} (resp.\
\ti{upper bound\/}) of~$A$ if $l\le a$ (resp.\ $a\le u$) \fe $a\in A$.
A~nonempty subset $A$ of~$X$ is called \ti{bounded below} (resp.\
\ti{above}) if it has a lower bound (resp.\ upper bound), \ti{bounded\/}\index{bounded set}
if it is bounded above and below, and \ti{unbounded\/} if it is not bounded.
\edf

\blm3.38
Let $(X,\le)$ be an \os\ and let $A$ be a nonempty subset of~$X$. If $m\in A$
is an upper bound $($resp.\ lower bound\/$)$ of~$A$, then $m$ is the greatest
$($resp.\ least\/$)$ \el\ of~$A$ and moreover, $m$ is the least upper bound
$($resp.\ greatest lower bound\/$)$ of~$A$.
\elm

\proof
Let $\vn\ne A\sbs X$ and let $m\in A$ be \st $a\le m$ \fe $a\in A$ (upper
bound). Then clearly $m$ is a (hence the) greatest \el\ of~$A$. Let $u\in X$
be \st $a\le u$ \fe $a\in A$ ($u$~upper bound), then $m\le u$. Hence $m$ is
the least upper bound of~$A$. Replacing the order \rl\ $\le$ by $\ge$
(see Exercise \rf{ex3.4}) one obtains the remaining part of the statement.
\endproof

\bth3.39
Let $(E,\le,e,S)$ be an \os\ of \nn s and let $A$ be a nonempty subset of~$E$
which is bounded above. Then $A$ possesses a greatest \el.
\eth

\proof
Let $\vn\ne A\sbs E$ and let $u\in X$ be an upper bound of~$A$. Let $B$ be
the set of upper bounds of~$A$. Since $B\ne\vn$ and $(E,\le)$ is
well-ordered, $B$~has a least \el\ which we denote by~$m$. In view of the
previous lemma it is sufficient to show that $m\in A$. Suppose for \cd ion
$m\notin A$. Clearly $m\ne e$, otherwise $A$ would be empty. \Mo $a<m$ \fe
$a\in A$. Since $m\ne e$, \te s $p\in E$ \st $S(p)=m$, in view of~(N1).
By Lemma \rf{l3.27} there is no $\ov a\in E$ \st $p<\ov a<S(p)$. Note that
if $a\in A$, then $a<m=S(p)$, hence $p\not< a$. Thus since the order~$\le$
is total, every $a\in A$ \sf ies $a\le p$. \If that $p$~is an \ub\ of~$A$,
hence $p\in B$. In view of $p<S(p)$ by \er{3.13}, $m=S(p)$ is \ti{not\/}
the least \el\ of~$B$, a~\cd ion. \E\Tf $m\in A$, and $m$ is the greatest
\el\ of~$A$.
\endproof

We recall that if $f$ is a self-map of a set $F$~($\ne\vn$) and $\vn\ne M\sbs
F$ is \st\ $f(M)\sbs M$ (i.e.\ \fe $x\in M$, $f(x)\in M$), then $M$~is called
\ti{\inv\ under}~$f$. \E\Ip if $f(x)=x$ for some $x\in F$, then $x$ is called
a \ti{fixed point of\/}~$f$, and $\{x\}$ is \inv\ under~$f$. In case $F$ is
an \os\ of \nn s~$E$, $f$~is the \su\ \f\ $S$ and $M\sbs E$ is \inv\
under~$S$, then $M=E$ iff $e\in M$ by~(N2). We now characterize all subsets
$M\sbs E$ \st $e\notin M$ and $S(M)\sbs M$. Note that if $S(M)\sbs M$, then
$S|_M$, the \rt ion of~$S$ to $M$, $S|_M:M\to M$, is well-defined.

\blm3.40
Let $(E,\le,e,S)$ be an \os\ of \nn s and let $M\sbs E$ be \st $M\ne\vn$,
$S(M)\sbs M$ and $e\notin M$. Let $m$ denote the least \el\ of~$M$. Then
\beq3.23
M=\{x\in E: m\le x\}=R(\Phi(m)).
\e
\elm

\proof
Let $m$ denote the least \el\ of~$M$ which exists since $E$ is well-ordered.
\Mo $e\le m$ since $e$ is the least \el\ of~$E$ and $e\ne m$ since $m\in M$ and
$e\notin M$. Hence $e<m$.
Set $A:=\{x\in E: e\le x<m\}$, then $A\ne\vn$ since $e\in A$.
\Mo $A\cap M=\vn$. Indeed, since $m$ is the least \el\ of~$M$ we have $m\le
x$, $x\in M$, hence $x\notin A$, since the order is total. We claim that $A\cup M$ is
\iv. Indeed, $e\in A\sbs A\cup M$. \Mo if $x\in M$ then $S(x)\in M$ by \as.
If $x\in A$ then $S(x)\le m$ by \er{3.17}. If $S(x)<m$ then $S(x)\in A\sbs
A\cup M$ and if $S(x)=m$ then $S(x)\in M\sbs A\cup M$. Therefore by (N2)
$E=A\cup M$ and since $A\cap M=\vn$, we have $M=E\sm A=\{x\in E: m\le x\}$,
since the order is total.

Finally, if $x\in R(\Phi(m))$, then there is $p\in E$ \st
$x=\Phi(m)p=\Phi(p)m$ by \er{2.38}. Then $m\le x$ by \er{3.6n}. Conversely, if
$m\le x$, then there is $p\in E$ \st $x=\Phi(p)m=\Phi(m)p$, hence $x\in R(
\Phi(m))$.
\endproof

\bco3.41
\E\fe $m\in E\sm\{e\}$
\beq3.24
\bca
\{x\in E: e\le x<m\} \cup R(\Phi(m)) = E,\\
\{x\in E: e\le x<m\} \cap R(\Phi(m)) = \vn.
\eca
\e
\eco

We now introduce some notation for a class of subsets (possibly empty) of an
\os. These subsets are called \ti{\il s}.\index{interval}

\bnt3.42
Let $(E,\le)$ be an \os\ and let $a,b\in E$ be \st $a\le b$.
\bea3.25
[a,{\to})&:= \{x\in E:a\le x\},\q ({\leftarrow},b]:=\{x\in E:x\le b\},\\
(a,{\to})&:= \{x\in E:a< x\},\q ({\leftarrow},b):=\{x\in E:x< b\},\lb{3.26}\\
[a,b]&:= \{x\in E: a\le x\le b\}, \lb{3.27}\\
(a,b]&:= \{x\in E: a< x\le b\}, \lb{3.28}\\
[a,b)&:= \{x\in E: a\le x< b\}, \lb{3.29}\\
(a,b)&:= \{x\in E: a< x< b\}, \lb{3.30}
\e
\ent

\Mo
\beq3.31
[a,b]=[a,b)\cup \{b\}.
\e
Indeed, $[a,b)\cup\{b\}\sbs [a,b]$. \E\oh $[a,b]\sbs [a,b)\cup\{b\}$. Indeed,
if $x\in [a,b]$ then $a\le x$ and $x\le b$. If $x<b$ then $x\in [a,b)\sbs
[a,b)\cup \{b\}$ and if $x=b$ then $x\in\{b\}\sbs [a,b)\cup \{b\}$.

Similarly,
\beq3.32
\left\{
\aligned{}
[a,{\to})&=(a,{\to})\cup \{a\},\\
({\leftarrow},b]&=({\leftarrow},b)\cup\{b\}\\
[a,b]&=(a,b]\cup \{a\},\\
[a,b]&=(a,b)\cup \{a\}\cup\{b\}.
\endaligned \right.
\e

We now mention a \pp y of some \il s of a set of \nn s.

\blm3.43
Let $(E,\le,e,S)$ be an \os\ of \nn s and let $a,b\in E$ be \st $a<b$. Then
\te s a unique \el\ $c\in[a,b)$ \st
\beq3.33
[a,b)=[a,c].
\e
\Mo
\beq3.34
S(c)=b.
\e
\elm

\proof \

\er{3.33}: Since the subset $[a,b)$ is nonempty ($a\in[a,b)$) and is bounded
above, it possesses a greatest \el~$c$, by Theorem~\rf{t3.39}. Therefore, if
$x\in [a,b)$ then $a\le x\le c$ hence $[a,b)\sbs [a,c]$. Conversely, if
$x\in[a,c]$, then $a\le x$, $x\le c$. Since $c\in[a,b)$, we have $c<b$. Then
$x\le c<b$ implies $x<b$ by \er{3.11}, and $[a,c]\sbs [a,b)$. If \er{3.33} holds \fs
$c\in[a,b)$, then $c$ is the greatest \el\ of~$[a,b)$.

\er{3.34}: Since $c<b$ by \er{3.33}, we get $S(c)\le b$ by \er{3.17}. Since
$c<S(c)$ by \er{3.13}, we have $c<S(c)\le b$. If $S(c)<b$, then
$c$~would not be the greatest \el\ of $[a,b)$. Hence $S(c)=b$.
\endproof

\brm3.44
In Section \ref{sss.Sets} we showed that \te s a unique map
$P:E\sm\{e\} \to E$ \st $P\circ S=\id_E$, $S\circ P=\id_{E\sm\{e\}}$. \E\fe
$x\in E\sm\{e\}$, $P(x)$ can be called the \ti{\im \pd} of~$x$ \wrt the
natural ordering of~$E$. With this \df, we have
\beq3.35
[a,b)=[a,P(b)]
\e
whenever $a,b\in E$, $a<b$.
\erm

\bex3.45
Let $(E,\le,e,S)$, $a,b$ be as in Lemma \rf{l3.43} and $P$ as in Remark
\rf{r3.44}. Show that the \rt ion of~$S$ to $[a,P(b)]$ is a bi\jn\ from
$[a,P(b)]$ onto $[S(a),b]$.
\eex

We now give a characterization of an \il\ of the form $[a,b]$ by means of a
\rc ve \rl.

\blm3.46
Let $E$, $a$, $b$ be as in Lemma \rf{l3.43}. Let $M$ be a subset of $[a,b]$.
Then $M=[a,b]$ iff the \fw\ \cn s hold\/{\rm:}
\bea3.36
{}&a \in M,\\
&S(x)\in M \hbox{ \fe }x\in M \cap [a,b).\lb{3.37}
\e
\elm

\proof[Proof \rm(if)]
Set
\beq3.38
\wt M:=[e,a)\cup M\cup (b,{\to}).
\e
Since $a<b$, $E$ is the disjoint union of $[e,a)$, $[a,b]$ and $(b,{\to})$.
Since $M\sbs[a,b]$, we have $M=[a,b]$ iff $\wt M=E$ iff $\wt M$ is \iv. We now
show that $\wt M$ is \iv.

$e\in \wt M$: If $e=a$, then $e\in M$ since $a\in M$. Then $e=a\in M\sbs\wt M$,
hence $e\in\wt M$.

If $e<a$ then $e\in[e,a)\sbs\wt M$.

$S(\wt M)\sbs \wt M$: We have $S(\wt M)=S([e,a))\cup S(M)\cup S((b,{\to}))$
if $e<a$, and $S(\wt M)= M\cup (b,{\to})$ if $e=a$.

(i) \ti{$S([e,a))\sbs \wt M$ if $a>e$}:
Indeed, $S([e,a))\nde3.14 \sbs [S(e),S(a))\nde3.17 = [S(e),a]\nde 3.31 =
[S(e),a)\cup\{a\}\nde 3.13 \sbs [e,a)\cup\{a\}$. But $[e,a)\sbs \wt M$, and
$\{a\}\sbs\wt M$ by \er{3.36}, hence (i)~holds.

(ii) $S(M)\sbs \wt M$: We have $M\sbs(M\cap [a,b))\cup\{b\}$, hence $S(M)\sbs
\bigl(S(M\cap [a,b))\bigr)\cup\{S(b)\}$. But $S(M\cap[a,b))\sbs M$ by \er{3.37},
and $\{S(b)\}\sbs(b,{\to})$ by \er{3.13}. Since both $M$ and $(b,{\to})$
are included in~$\wt M$, we obtain $M\cup (b,{\to})\sbs\wt M$, which implies
(ii).

(iii) $S((b,{\to}))\sbs\wt M$: If $x\in(b,{\to})$, then $b<x$ by \er{3.26},
hence $S(b)<S(x)$ by \er{3.14}. From \er{3.13} we have $b<S(b)$, hence
by \er{3.11} $b<S(x)$, i.e.\ $S(x)\in(b,{\to})$. \If that $S((b,{\to}))\sbs
(b,{\to})$, hence (iii) holds since $(b,{\to})\sbs\wt M$.

Then $S(\wt M)\sbs\wt M$ follows from (i), (ii), (iii). Thus $\wt M$ is \iv\ and
$\wt M=E$, $M=[a,b]$.

(only if) Suppose $M=[a,b]$. Then \er{3.36} is obvious and if $x\in[a,b)$,
then $a<S(a)\le S(x)\le b$ by \er{3.13}, \er{3.14} and \er{3.17}.
Hence $S(x)\in (a,b]\sbs [a,b]=M$.
\endproof

\bex3.47
Let $E$, $a$ be as in Lemma \rf{l3.46}. Let $M$ be a subset of~$[a,{\to})$ \st
\er{3.36} holds. 

Show that if $S(x)\in M$ \fe $x\in M$, then the conclusion
becomes $M= [a,{\to})$.
\eex

We recall that every subset $A$ of an \os\ $(X,\le)$ is itself an \os\ under
the ordering of~$X$ \rt ed to~$A$. \E\Ip if $(X,\le)$ is totally ordered so
is $(A,\le)$.

In the remaining part of this section we shall \ti{characterize} the subsets of an
\os\ of \nn s $(E,\le)$ which are \ois c to $(E,\le)$. We first show that
bounded subsets of $(E,\le)$ are not \ois c to $(E,\le)$. Note that a
\nss\ $A$ of~$E$ is bounded in $(E,\le)$ iff $A$ is bounded above, since $e$~is
always a lower bound. Indeed, if $A$ is bounded in $(E,\le)$, $A$ has a~greatest
\el\ by Theorem~\rf{t3.39}.

In what follows $E$ will always denote an \ti{\os\ of \nn s} $(E,\le,e,S)$ and a
subset $A$ of~$E$ will always be endowed with the ordering of~$E$ \rt ed
to~$A$.

\bpr3.48
Let $A$ be a nonempty bounded subset of~$E$ and let $f$ be a map from $A$
into~$E$. Then $f(A)$ is bounded, \Ip $f(A)\ne E$.
\epr

\proof
We define a map $g:E\to E$ by setting $g(x):=f(x)$ if $x\in A$ and $g(x):=e$
if $x\notin A$. We claim that \fe $n\in E$ \te s $p(n)\in E$ \st
\beq3.39
g(x)\le p(n) \qh{\fe}x\in[e,n].
\e
Indeed, let $M:=\{n\in E: \hbox{\er{3.39} holds for some }p(n)\in E\}$. We
show that $M$ is \iv. $e\in M$: setting $p(e):=g(e)$ we have $g(e)\le p(e)$.
$n\in M$ \ti{implies} $S(n)\in M$: Suppose ${n\in M}$. Let $p(S(n))$ be the greatest \el\
of the set $\{p(n),g(S(n))\}$, that is, $p(S(n)):=p(n)$ if $g(S(n))<p(n)$ and
$p(S(n)):=g(S(n))$ if $p(n)\le g(S(n))$. Hence $p(n)\le p(S(n))$ and $g(S(n))
\le p(S(n))$. Since $n\in M$, we have $g(x)\le p(n)\le p(S(n))$, hence $g(x)
\le p(S(n))$ \fe $x\in[e,n]$. Since $[e,S(n)]=[e,n]\cup\{S(n)\}$, we obtain
$g(x) \le p(S(n))$ \fe $x\in[e,S(n)]$. \E\Tf $S(n)\in M$, and $M=E$ by (N2).
Since $A$ is bounded and nonempty, it follows from Theorem~\rf{t3.39} that $A$
possesses a greatest \el\ which we denote by~$m$. \E\Tf $g(x)\le p(m)$ \fe
$x\in[e,m]$. Since $A\sbs[e,m]$, $f(x)=g(x)\le p(m)$ \fe $x\in A$, hence $f(A)$
is bounded above. By Theorem~\rf{t3.39} $f(A)$ has a greatest \el\ which we
denote by~$\wt m$. Since $\wt m<S(\wt m)$ by \er{3.13}, $S(\wt m)\notin
f(A)$ and $f(A)\ne E$.
\endproof

\bth3.49
A nonempty subset $A$ of an \os\ of \nn s $E$ is \ois c to~$E$ iff $A$ is not
bounded.
\eth

\proof[Proof \rm(only if)]
It follows from \E\Pr\ \rf{p3.48} that there is no sur\jc\ map\break ${\vf:A\to E}$,
whenever $A$~is bounded, hence in this case $A$ and~$E$ are \ti{not\/} \ois c.

(if)
Let $A$ be a nonempty subset of $E$ which is not bounded. Since $(E,\le)$ is
well-ordered by Theorem \rf{t3.22}, $(A,\le)$ is also well-ordered. Set
$e':=$ the least \el\ of~$A$. \E\fe $x\in A$ the set $\{y\in A: x<y\}$ is
\ti{not empty} since $A$ is \ti{not bounded\/}. Define a~map $S':A\to A$ by
setting $S'(x):=$ the least \el\ of the set $\{y\in A: x< y\}$, \fe
$x\in A$. By Proposition \rf{p3.28}, $S':A \to A$ is strictly in\cre\
and in\jc. In view of Theorem~\rf{t1.6} \te s one (and
only one) map $\vf:E\to A$ \st $\vf(e):=e'$ and $\vf(S(n)):=S'(\vf(n))$ \fe
$n\in E$. We next prove that

\leavevmode\llap{(i) }$\vf$ \ti{is strictly in\cre\ and in\jc}.

We proceed by induction. Let $x\in E$ be given. Set $M:=\{n\in E: \vf(x) <
\vf(\Phi(S(n))x)\}$ where $\Phi(S(n))$ is the $S(n)$-fold \IT\ of~$S$
\wrt $E$. We claim that $M$ is \iv. $e\in M$:
$\vf(x)<S'(\vf(x))=\vf(S(x))=\vf(\Phi(S(e))x)$, since $\Phi(S(e))=S$.
The in\et y follows from the \df\ of~$S'$. The first \et y follows
from the \df\ of~$\vf$.

\ti{$n\in M$ implies $S(n)\in M$}: Since $n\in M$, we have $\vf(x)<\vf(\Phi(S(n))x)$ and
since\break $S'$~is strictly in\cre, $S'(\vf(x))<S'(\vf(\Phi(S(n))x))$. We have
$\vf(x) < S'(\vf(x))$ by\break \df\ of~$S'$.
Hence, we obtain $\vf(x) < S'(\vf(\Phi(S(n))x))$ by \er{3.11}.
We have\break $S'(\vf(\Phi(S(n))x))=
\vf(S(\Phi(S(n))x))$, by \df\ of~$\vf$, and $\vf(S(\Phi(S(n))x))=\break
\vf(\Phi(S(S(n)))x)$ by \er{2.34}. Hence we obtain $\vf(x)<
\vf(\Phi(S(S(n)))x)$, which implies $S(n)\in M$. Since $M$ is
\iv, we have $M=E$ by~(N2). Finally, let $y\in E$ be \st $x<y$. By~\er{3.17},
$S(x)\le y$, hence by the \df\ of the ordering~$\le$, there is $p\in E$ \st
$y=\Phi(p)S(x)$. Observe that $\Phi(p)S(x)=\Phi(S(p))x$ by \er{2.25},
\er{2.26a}. Hence, $y=\Phi(S(p))x$. Since $M=E$, we obtain $\vf(x)
< \vf(\Phi(S(p))x)=\vf(y)$. Therefore, $\vf$ is strictly in\cre. Hence $\vf$~is
in\jc\ by Lemma \rf{l3.34}\,(ii) since $(E,\le)$ is a chain.

We next prove that
\beq3.40
n\le \vf(n) \qh{\fe}n\in E.
\e
Set $M:=\{n\in E: n\le \vf(n)\}$.

$e\in M$: clear since $e\le x$ \fe $x\in E$.

\ti{$n\in M$ implies $S(n)\in M$}: we have $n\le \vf(n)$. By the \df\ of~$S'$, $\vf(n)<
S'(\vf(n))$ and by the \df\ of~$\vf$, $S'(\vf(n))=\vf(S(n))$. Hence $n\le
\vf(n) < \vf(S(n))$, and $n<\vf(S(n))$ by \er{3.11}. By~\er{3.17}, $S(n)\le \vf(S(n))$. It
follows that $M$ is \iv\ and $M=E$.

We next prove

\leavevmode\llap{(ii) }$\vf$ \ti{is sur\jc}.

Let $m\in A$. Then by the \df\ of~$e'$ and by \er{3.40}, we have $\vf(e)=e'
\le m\le\vf(m)$. We introduce the notation $[a,b]':=[a,b]\cap A = \{y\in A:
a\le y\le b\}$ for $a,b\in A$ \st $a<b$. With this notation we have $m\in
[\vf(e), \vf(m)]'$. We want to show that $[\vf(e),\vf(m)]'=\vf([e,m])$. Set
$M:= \{n\in E: [\vf(e),\vf(n)]'=\vf([e,n])\}$. $e\in M$, since
$[\vf(e),\vf(e)]'=\{\vf(e)\}=\vf([e,e])$.

\ti{$n\in M$ implies $S(n)\in M$}: Let $n\in E$. First, note that $[\vf(e),\vf(S(n))]'
=\{y\in A: \vf(e)\le y\le \vf(S(n))\}=
\{y\in A: \vf(e)\le y\le \vf(n)\} \cup \{y\in A: \vf(n)< y\le \vf(S(n))\}$.
By the \df\ of~$\vf$, $\vf(S(n))=S'(\vf(n))$. \Mo by the \df\ of~$S'$,
we find $\{y\in A: \vf(n)< y\le S'(\vf(n))\}=\{S'(\vf(n))\}$. It follows that
$[\vf(e),\vf(S(n))]' = [\vf(e),\vf(n)]'\cup \{\vf(S(n))\}$. Using $n\in M$ we
infer that $[\vf(e),\vf(S(n))]'=\vf([e,n])\cup \{\vf(S(n))\}=
\vf([e,S(n)])$. This shows that $S(n)\in M$, hence $M$~is \iv\ and $M=E$.
\E\Ip $[\vf(e),\vf(m)]'=\vf([e,m])\sbs R(\vf)$. \E\Tf $m\in[\vf(e),\vf(m)]'
=\vf([e,m])\sbs R(\vf)$. Since $m\in M$ is arbitrary we have $A\sbs R(\vf)$,
which completes the proof.

Finally, we prove that $\vf$ is an \ti{\ois sm between $(E,\le)$ and\/}
$(A,\le)$. Since $\vf$ is in\jc, it is bi\jc\ by~(ii).
Since $\vf$ is in\cre\ and $(E,\le)$ is totally ordered, we infer in view of
Lemma~\rf{l3.34} that
also $\vf^{\rm inv}$ is in\cre. This completes the proof.
\endproof

\bex3.50
Show that if $A$, $S'$, $e'$ are as in the proof of Theorem~\rf{t3.49},
the \fw\ holds: if $M$ is a subset of~$A$ \st $e'\in M$ and $S'(M)\sbs M$,
then $M=A$.
\eex

We conclude this section by giving some important \rl s between the \og~$\le$
of~$E$ and the \og~$\sbs$ of $\cP(E)$.

\blm3.49
\allowdisplaybreaks
Let $(E,\le,e,S)$ be an \os\ of \nn s and let $a,b,c\in E$.

If $a<b$ and $a< c$, we have
\bea3.41
{}[a,b] \sbs [a,c] &\qhq{iff} b\le c, \\
[a,b) \sbs [a,c) &\qhq{iff} b\le c, \lb{3.42} \\
(a,b] \sbs (a,c] &\qhq{iff} b\le c, \lb{3.43} \\
(a,b) \sbs (a,c) &\qhq{iff} b\le c, \lb{3.44} \\
[a,b] \sbs [a,c) &\qhq{iff} b< c, \lb{3.45} \\
[a,b] = [a,c] &\qhq{iff} b= c, \lb{3.46} \\
[a,b) = [a,c) &\qhq{iff} b= c, \lb{3.47} \\
(a,b] = (a,c] &\qhq{iff} b= c, \lb{3.48} \\
(a,b) = (a,c) &\qhq{iff} b= c. \lb{3.49}
\e
\elm

\bex3.51
Prove Lemma \rf{l3.49}.
\eex

\bex3.52
Show that if $(E,e,S)$ is a set of \nn s and $\le$ is the \nog, then $\id_E$
is the only bi\jc\ in\cre\ self-map of~$E$.
\eex

\newpage
\Subsubsection{Cardinality}\label{sss.Card}
In this section we first introduce an \ev ce \rl\ on the power set $\cP(E)$
of an \os\ of \nn s~$E$. This \rl\ allows us to define the notions of
finite and infinite subsets of~$E$. These notions are then extended to
arbitrary sets.

It turns out that the notion of infinite set has to be refined. To this end
the concepts of countably infinite and uncountable sets are introduced.

Let $(E,\le,e,S)$ be an \os\ of \nn s and let~$A$ be a subset of~$E$, $A\in
\cP(E)$. It follows from Theorem~\rf{t3.49} that if $A$ is not bounded,
\te s a bi\jn\ from~$A$ onto~$E$. \Mo if $A$ is bounded, \te s
no bi\jn\ from~$A$ onto~$E$ in view of \E\Pr~\rf{p3.48}.

\bdf4.1
Let $A$ and $B$ be nonempty subsets of~$E$. Then $A$ and~$B$ are said to be\index{equipotent}
\ti{equipotent\/} or \ti{equipollent\/} or to have \ti{the same cardinality}
if there is a bi\jn\ from~$A$ onto~$B$. If $A$ and~$B$ are equipotent, we
write $A\approx B$.\glossary{$\approx$}
\edf

\blm4.2
Let $A,B,C$ be nonempty subsets of~$E$. Then the \fw\ holds{\rm:}
\bea4.1
{}&A\ax B \qh{and} B\ax C \qh{implies} A\ax C,\\
&A\ax A, \lb{4.2}\\
&A\ax B \qh{implies} B\ax A \qh{$($symmetry$)$.}\lb{4.3}
\e
\elm

\proof[Proof \rm(sketch)]\

\er{4.1}: The \cm\ of two bi\jn s is a bi\jn.

\er{4.2}: $\id_A$ is a bi\jn.

\er{4.3}: The inverse of a bi\jn\ is a bi\jn.
\endproof

Let us denote by $\dot \cP(E)$ the set of all nonempty subsets of~$E$. Then the
\rl\ $\ax$ on~$\dot \cP(E)$ is \ti{\tr e} (see~\er{3.2}), \ti{\rfl}
(see~\er{3.3}) and \ti{\sy ic} (see~\er{4.3}).

\bdf4.3
Let R be a \rl\ on a nonempty set~$X$. The \rl~R is called an \ti{\ev ce \rl}\index{equivalence relation}
if it is \tr e, \rfl\ and \sy ic.
\edf

\bex4.4
Let R be an \ev ce \rl\ on a nonempty set~$X$. \E\fe $x\in X$ set
\beq4.4
[x]:=\{y\in X: x\mr R y\}.
\e
Show the \fw:
\ssk

\hph i,i, if $x,y\in X$, then $y\in[x]$ iff $[x]=[y]$;

\hph ii,, if $x,y\in X$, then $y\notin[x]$ iff $[x]\cap [y]=\vn$.
\eex

\bdf4.5
Let $X$ be a nonempty set. A~nonempty subset $\cA$ of $\cP(X)$ is called a
\ti{partition}\index{partition} of~$X$, if the \fw\ holds:
\ssk

\hph i,ii, $A,B\in \cA$, $A\ne B$ implies $A\cap B=\vn$ \fa $A,B\in\cA$,

\hph ii,i, $\bigcup\limits_{A\in \cA}A=X$,

\hph iii,, $\vn\notin\cA$.
\edf

\bdf4.6
Given an \ev ce \rl\ R on a set~$X$, a nonempty subset $A$ of~$X$ is called
an R-\ti{\ev ce class} (or simply an \ti{\ev ce class}) if
\te s an \el\ $x\in X$ \st $A=[x]$.
\edf

\bex4.7 \

\hph i,i, Let R be an \ev ce \rl\ on a set~$X$. Set $\cA_{\rm R}$ the set of all
\ev ce classes. Show that $\cA_{\rm R}$ is a partition of~$X$.

\hph ii,, Conversely, let $\cA$ be a partition of $X$, and define $x\mr Ry$ if
\te s $A\in\cA$ \st $x\in A$ and $y\in A$.
Show that R is an \ev ce \rl\ and that
$\cA_{\rm R}=\cA$.
\eex

\brm4.8
Clearly the \et y \rl\ on a set $X$ is an \ev ce \rl. The \et y \rl\ on a set
$X$ is ``finer'' than any \ev ce \rl~R on~$X$, that is,
if R is an \ev ce \rl\ on~$X$ and $x,y\in X$, then $x=y$
implies $x\mr Ry$.
\erm

\bex4.9
Let $X$ be a nonempty set and let R be a \rl\ which is \tr e and \rfl. Define
\fa $x,y\in X$:
\beq4.5
x\tR y \qh{if }\, x\mr Ry \hbox{ and }y\mr Rx.
\e
Show that $\tR$ is an \ev ce \rl\ and let $\cA_{\tR}$ be the corresponding
set of \ev ce classes. Define on $\cA_{\tR}$ the \rl
\[
A\le B \qh{if} A,B\in\cA_{\tR} \hbox{ and \te\ }x\in A,y\in B \hbox{ \st}
x\mr Ry.
\]
Show that $\le$ is an ordering on $\cA_{\tR}$.
\eex

We now return to the \rl\ $\ax$, which is an \ev ce \rl\ on $\dot \cP(E)$
where $(E,\le,e,S)$ is an \os\ of \nn s. One
can extend this \rl\ to $\cP(E)$ by setting \fe $A\in \cP(E)$:
\bea4.6
A\ax \vn \q &\hbox{if}\q A=\vn,\\
\vn\ax A \q &\hbox{if}\q A=\vn.\lb{4.7}
\e

It is easily verified that the extended \rl~$\ax$ on~$\cP(E)$ is an \ev ce \rl\
and that $[\vn]$ consists of the \ti{empty subset\/} of~$E$.
It follows from Theorem~\rf{t3.49} that $[E]$ consists of all \ti{unbounded
subsets} of~$E$. We now consider nonempty bounded subsets of~$E$.

\blm4.10
Let $A$ be a nonempty bounded subset of $E$. Then \te s some $n\in E$ \st
$A\ax [e,n]$.
\elm

\proof
By Theorem \rf{t3.39}, $A$ possesses a greatest \el\ which we denote by~$m$.
Set $\wt A:=A\cup(m,\to)$. Note that $A\cap(m,\to)=\vn$. Clearly $\wt A$ is
not bounded since $(m,\to)$ is not bounded. By Theorem~\rf{t3.49}, \te s a bi\jn\
$\vf:E\to \wt A$ which is in\cre\ as well as $\vf^{\rm inv}$. Set
$n:=\vf^{\rm inv}(m)$. We have $\vf([e,n])\subset A$. Indeed, if $x\in[e,n]$,
i.e.\ $e\le x\le n$, then $e\le \vf(x)\le\vf(n)=m$, since $e$~is the least
\el\ of~$E$ and $\vf$~is in\cre. Hence
\[
\vf(x)\in \wt A\cap[e,m] = A\cap[e,m]=A.
\]
\Mo $A\sbs \vf([e,n])$. Indeed, let $y\in A\sbs \wt A$, and set
$x:=\vf^{\rm inv}(y)$. Since $y=\vf(x)$, it is \sft\ to show that $e\le x\le
n$. The first in\et y follows from \er{3.13}. Since $m$~is the greatest \el\
of~$A$, we have $\vf(x)=y\le m=\vf(n)$. Hence $\vf(x)\le \vf(n)$. Since
$\vf\Inv$ is in\cre, we obtain $x=\vf\Inv(\vf(x))\le \vf\Inv(\vf(n))=n$, and
$x\in [e,n]$. Thus $y\in\vf([e,n])$, hence $A\sbs \vf([e,n])$. It follows that
the image of $[e,n]$ under $\vf$ is~$A$. Let us denote by $\wt\vf$ the \rt
ion of~$\vf$ to~$[e,n]$. It is a bi\jn\ from $[e,n]$ onto~$A$, since $\wt\vf$~is
a~\rt ion of the in\jc\ map~$\vf$.
\endproof

We now consider the problem of \ti{uniqueness} of an \el\ $n\in E$ \sf ying the
conclusion of Lemma~\rf{l4.10}. The answer to this question is a con\sq\ of
the next lemma.

\blm4.11
Let $m,n\in E$. If \te s an in\jc\ map $f:[e,m]\to[e,n]$, then $m\le n$.
\elm

\proof
We proceed by induction on $n$. Let $m\in E$ and set
$$
M:=\{n\in E: \hbox{if \te s an in\jc\ map
$f:[e,m]\to[e,n]$, then $m\le n$}\}.
$$

$e\in M$: Let $f:[e,m]\to\{e\}$ be in\jc. Then $f(e)=e=f(m)$
implies $m=e\le e$. Hence $e\in M$.

\ti{$n\in M$ implies $S(n)\in M$}: Let $n\in M$ and let
$f:[e,m]\to[e,S(n)]$ be in\jc. We show that $m\le S(n)$. If $m=e$, then
$m=e\le S(n)$ by \E\df\ \rf{d3.30}. If $e<m$, then \te s one and only one
$p\in E$ \st $S(p)=m$, by Theorem \rf{t3.31}\,(b)(i). Then $f$ maps $[e,S(p)]$
into $[e,S(n)]$. If $R(f)$, the range of~$f$, does not contain $S(n)$, then
$R(f)\sbs [e,n]$ by (OS1). Since $n\in M$, we infer $m\le n$, hence $m\le
n<S(n)$, consequently $m<S(n)$ by \er{3.11}, hence $m\le S(n)$.

If $S(n)\in R(f)$ and $S(n)=f(S(p))$, then, by the injectivity of~$f$,
$f$~maps $[e,p]$ into $[e,n]$ and $f|_{[e,p]}$ is also in\jc. Since $n\in M$,
we infer that $p\le n$. By Theorem \rf{t3.31}\,(b)(i) we get $m=S(p)\le S(n)$. It remains to
consider the case where $S(n)\ne f(S(p))$, i.e.\ $S(n)=f(x)$ for some
$x\in[e,p]$. The idea is to modify the map~$f$ in such a way that the
modified map $\wt f:[e,m]\to[e,S(n)]$ remains in\jc\ and \sf ies $\wt
f(m) =S(n)$. To this end we define a map $h:[e,S(n)]\to[e,S(n)]$ by
setting
\beq4.8
h(k):=\bca
S(n) & \hbox{if }k=f(m),\\
f(m) &\hbox{if }k=S(n),\\
f(k) &\hbox{otherwise,}
\eca
\e
i.e.\ $h$ exchanges $S(n)$ and $f(m)$, and leaves fixed the complement of
$\{S(n),f(m)\}$ in $[e,S(n)]$.

Note that
\beq4.9
h\circ h=\id_{[e,S(n)]}
\e
hence $h$ is bi\jc\ and $h^{\rm inv}=h$.
Set
\beq4.10
\wt f:=h\circ f.
\e
Then $\wt f:[e,m]\to[e,S(n)]$ is in\jc\ since both $f$ and~$h$ are in\jc.
\Mo $\wt f(m)=h(f(m))=S(n)$. Thus we find as above $m\le S(n)$.
This completes the proof of ``$n\in M$ implies
$S(n)\in M$''. \If that $M$ is \iv, hence $M=E$.
\endproof

\bex4.12
Prove Lemma \rf{l4.11} by defining $\wt f$ as in \er{4.10} without
distinguishing cases.
\eex

\bco4.13
Let $m,n\in E$ and let $f:[e,m]\to[e,n]$.
\ssk

\hph i,i, If $f$ is sur\jc, then $n\le m$.

\hph ii,, If $f$ is bi\jc, then $n=m$.
\eco

\proof \

(i) We construct a map $g:[e,n]\to[e,m]$ \st $f\circ g=\id_{[e,n]}$. Given
$y\in[e,n]$, let $A(y):=\{x\in[e,m]:f(x)=y\}$. Since $f$ is sur\jc,
$A(y)\ne\vn$. Set $g(y):={}$the least \el\ of~$A(y)$. Then $f(g(y))=y$. The
map $g$ is in\jc\ since $g(y)=g(\wt y)$, $y,\wt y\in[e,n]$ implies
$y=f(g(y))=f(g(\wt y))=\wt y$. The conclusion follows from Lemma~\rf{l4.11}
applied to~$g$.

(ii) This is a con\sq\ of Lemma \rf{l4.11} and part~(i) or of Lemma \rf{l4.11}
applied to~$f$ and~$f^{\rm inv}$.
\endproof

Corollary~\rf{c4.13}~(ii) implies that in Lemma \rf{l4.10} ``some $n\in E$'' can
be replaced by ``exactly one $n\in E$''. Indeed, if $A\ax[e,n]$ and
$A\ax[e,m]$, then $[e,n]\ax[e,m]$, hence $m=n$.

We recall that every unbounded subset of~$E$ is equipotent to~$E$. We
now prove that \fe $n\in E$, the only \ti{subset\/} of $[e,n]$ which is equipotent
to~$[e,n]$ is $[e,n]$ itself.

\blm4.14
Let $n\in E$ and let $A\sbs [e,n]$ be nonempty. If $A\ax[e,n]$, then
$A=[e,n]$.
\elm

\proof
If $n=e$, the only nonempty subset of $[e,n]$ is~$\{e\}$ and \Tf $A=[e,n]$.
So we may suppose $n>e$. By Theorem \rf{t3.31}\,(b)(i) \te s one and only one $p\in E$
\st $S(p)=n$. \Mo we have $p<S(p)$ by (OS1).
Since $A\ax[e,n]$ \te s $f:[e,n]\to A$ bi\jc, hence in\jc.
We suppose for \cd ion that $A\ne[e,n]$. Let $c\in[e,n]$ be \st $c\notin A$.
Set $B:=[e,n]\sm\{c\}$. We have $A\sbs B$ and define $j:A\to B$ by setting
$j(x):=x$, $x\in A$. Clearly $j$ is in\jc. If $c=n$, then $B=[e,p]$ and
the map $\vf:=j\circ f:[e,n]\to[e,p]$ is in\jc. By Lemma~\rf{l4.11}, we have
$n\le p$. Since $p<S(p)=n$, we get a \cd ion. If $c\ne n$, then we
define a map $h:B\to[e,p]$ by setting $h(n):=c\in[e,p]$, and $h(x):=x$
otherwise. We also define a map $g:[e,p]\to B$ by setting $g(c):=n$ and
$g(x):=x$ otherwise. Observe that $g\circ h=\id_B$, hence $h:B\to[e,p]$
is in\jc. \E\Tf the map $\vf:=h\circ
j\circ f:[e,n]\to[e,p]$ is in\jc, and we also obtain a \cd ion. \E\Tf\
$A=[e,n]$.
\endproof

So far we considered subsets of a fixed \os\ of \nn s $(E,\le,e,S)$.
We now want to introduce \pp ies of sets which are independent of a
particular choice of \os\ of \nn s. The next lemma will justify the \df\ of
\ti{finite} and \ti{infinite} sets.

\blm4.15
Let $A$ be a nonempty set, let $(E,\le,e,S)$, $(E',\le',e',S')$ be \os s of
\nn s. Suppose \te\ $n\in E$ and a bi\jc\ map $f:[e,n]\to A$.
Then \te\ $m'\in E'$ and a bi\jc\ map $g:[e',m']'\to A$.
Here $[e',m']':=\{x'\in E': e'\le' x'\le' m'\}$.
\elm

\proof
In view of \E\df s \rf{d3.30}, \rf{d3.33}, Theorems \rf{t3.31},
\rf{t3.36}, \te s a bi\jn\ $\vf:E\to E'$ which is in\cre\ as well as
$\vf^{\rm inv}$. \Mo $\vf(e)=e'$ and $\vf([e,m])=[e',\vf(m)]'$ \fe $m\in E$.
Let  $n\in E$ and $f$ be as above. Set $n':=\vf(n)$. We now construct a
bi\jn\ $g:[e',n']'$ onto~$A$. We denote by $\wt\vf$ the \rt ion of~$\vf$
to~$[e,n]$. Then $\wt\vf$ is a bi\jn\ from $[e,n]$ onto $[e',n']'$. Set $h:=
\wt\vf\circ f^{\rm inv}:A\to[e',n']'$. Clearly $h$ is bi\jc\ and $h^{\rm inv}
=f\circ\wt\vf^{\rm inv}$ is a bi\jn\ from $[e',n']'\to A$.
\endproof

\bdf4.16
Let $A$ and $B$ be nonempty sets.
The sets $A$ and~$B$ are said to be \ti{equipotent\/}\index{equipotent} or
\ti{equipollent\/} or to have \ti{the same cardinality} if there is a bi\jc\
map from~$A$ onto~$B$. We write $A\approx B$ if $A$ and~$B$ are equipotent
and $A\not\approx B$ if they are not equipotent.\glossary{$\not\approx$}
\edf

\bdf4.17
A nonempty set $X$ is called \ti{finite}\index{set!finite} if \te s an \il\ $[e,n]$ of \ti{some} \os\
of \nn s $(E,\le,e,S)$ \st $X$ and $[e,n]$ are \ep. The \ti{empty set\/} is
also called a \ti{finite} set. A~set which is not finite is called
\ti{infinite}\index{set!infinite}.
\edf

In the \fw\ theorem we establish some basic \pp ies of finite sets and of
maps between finite sets.

\bth4.18 \

\hph i,ii, A subset of a finite set is finite.

\hph ii,i, If $f$ is a map from a nonempty finite set $X$ into a \ns~$Y$,
then $f(X)$ is finite.

\hph iii,, Let $f$ be a map from a nonempty finite set into itself.
Then $f$ is bi\jc\ iff $f$~is in\jc\ iff  $f$~is sur\jc.

\hph iv,, Every set of \nn s is infinite.
\eth

\proof \

(i) Let $A$ be a finite set and let $B$ be a subset of~$A$.
If $B=\vn$, \Ip if $A=\vn$, then $B$~is finite. Thus we may assume both
$A$~and~$B$ nonempty.
By \E\df\ \rf{d4.17} \te\ an \os\ of \nn s $(E,\le,e,S)$, an \el\ $n\in E$ and a bi\jc\ map
$f:A\to [e,n]$. Let $g:B\to f(B)$ be defined by $g(x):=f(x)$ \fe $x\in B$.
Note that $g$ is a bi\jn. The range $f(B)$ is a nonempty subset of $[e,n]$,
hence it is bounded in $(E,\le)$. By Lemma \rf{l4.10}, \te\ $m\in E$ and a
bi\jn\ $h:f(B)\to[e,m]$. It follows that $h\circ g:B\to[e,m]$ is a bi\jn,
hence $B$ is finite, since the inverse of a bi\jn\ is a bi\jn.

(ii) By \E\df\ \rf{d4.17} \te\ an \os\ of \nn s $(E,\le,e,S)$, $n\in E$, and
$g:[e,n]\to X$ bi\jc.  Set $B:=f(X)$. Then $h:=f\circ
g:[e,n]\to B$ is sur\jc. \E\fe $y\in B$ set $A(y):=\{a\in[e,n]: h(a)=y\}$.
Define $\vf:B\to[e,n]$ by setting $\vf(y):={}$the least \el\ of $A(y)$, which
exists since $A(y)\ne\vn$. Then $h\circ\vf=\id_B$. It follows that $\vf$ is
in\jc, hence $\vf:B\to\vf(B)\sbs[e,n]$ is bi\jc. Since $\vf(B)$ is a subset
of $[e,n]$ which is finite, $\vf(B)$ is also finite by part~(i). Hence \te\ an
\os\ of \nn s $(E',\le',e',S')$, $m'\in E'$, and $\psi':\vf(B)\to[e',m']'$
bi\jc\ in view of \E\df\ \rf{d4.17}. Using Lemma \rf{l4.15} with $E:=E'$,
$E':=E$, we find $m\in E$ and a bi\jn\ $\psi:\vf(B)\to [e,m]$. Consequently,
the map $\psi\circ\vf:B \to
[e,m]$ is a bi\jn. Since $f(X)=B$, we infer that $f(X)$ is finite.

\sdim{(b)}
(iii) Let $X$ be a nonempty finite set and let $f$ be a map from $X$ into~$X$.
Let $n\in E$ and $\vf:[e,n]\to X$ be bi\jc, as in \E\df\ \rf{d4.17}.
Set $\wt f:=\vf^{\rm inv} \circ
f\circ \vf:[e,n]\to[e,n]$. Then $f=\vf\circ\wt f\circ\vf\Inv :X\to X$,
\[
\xymatrix{
[e,n] \ar[r]^{\tilde f} \ar[d]_\vf & [e,n] \ar[d] ^\vf\\
X \ar[r]^f & X.}
\]
One verifies that $f$ is in\jc\ (resp.\ sur\jc) iff $\wt f$ is in\jc\ (resp.\
sur\jc). It is sufficient to show that
\ssk
\ite[(a)] {$\wt f$ in\jc\ implies $\wt f$ sur\jc\ and}
\ite[(b)] {$\wt f$ sur\jc\ implies $\wt f$ in\jc.}
\ssk

\ti{Proof of\/ }(a).
Set $A:=\wt f([e,n])\sbs[e,n]$. Since $\wt f$ is in\jc,
we have $A\ax [e,n]$. By Lemma~\rf{l4.14}, $A=[e,n]$, hence $\wt f$ is sur\jc.

(b) Since $\wt f:[e,n]\to[e,n]$ is sur\jc, we can find as in the proof of
Corollary \rf{c4.13}\,(i) an in\jc\ map $\wt g:[e,n]\to[e,n]$ \st $\wt f\circ
\wt g=\id_{[e,n]}$. From part~(a) we infer that $\wt g$~is bi\jc. Hence
$\wt f\nde2.2 = \wt f\circ(\id_{[e,n]}) = \wt f\circ(\wt g\circ \wt g{}\Inv)
\nde2.1 = (\wt f\circ\wt g)\circ\wt g{}\Inv = \id_{[e,n]}\circ \wt g{}\Inv
\nde2.2 = \wt g{}\Inv$. Since the inverse of a bi\jc\ map is bi\jc, we infer
that $\wt f=\wt g{}\Inv$ is bi\jc\ hence in\jc.

(iv) The ``\su'' \f\ is in\jc\ but \ti{not\/} sur\jc. Hence (iv) follows
from (iii).
\endproof

\goodbreak
\bex4.19 \

\hph i,i, Let $A,B,C$ be nonempty sets. Show that

\ssk
\ite[(a)] {$A\ax B$ and $B\ax C$ implies $A\ax C$,}
\ite[(b)] {$A\ax A$,}
\ite[(c)] {$A\ax B$ implies $B\ax A$.}
\ssk

\noindent (Since there is no such thing as the set of all nonempty sets one cannot
speak of an \ev ce \rl\ on the set of all nonempty sets.)

\hph ii,, Let $A,B$ be nonempty sets \st $A\ax B$. Show that if $A$ is finite, then
so is~$B$, and that if $A$  is infinite, then so is~$B$.
\eex

In the next chapter we shall pursue our study of finite sets. We now turn
to the case of infinite sets.

\blm4.20
Let $(E,\le,e,S)$ be an \os\ of \nn s. Then a nonempty subset $A$ of~$E$ is infinite
iff it is unbounded, \ev tly, it is finite iff it is bounded.
\elm

\proof
If $A$ is bounded, then it is finite by Lemma~\rf{l4.10}. If $A$ is unbounded,
then it is \ois c to~$E$ by Theorem~\rf{t3.49}, hence $A\ax E$. Since $E$ is
infinite by Theorem~\rf{t4.18}\,(iv), $A$~is also infinite by Exercise
\rf{ex4.19}\,(ii).
\endproof

The infinite subsets of an \os\ of \nn s $(E,\le)$ are all \ep\ to~$E$ since
they are \ois c to~$E$. \Mo all sets of \nn s are \ep\ by Theorem~\rf{t1.4}.
This justifies the \fw\ \df s.

\bdf4.21
A set is called \ti{\ct y infinite}\index{set!countably infinite} if it has the same cardinality as a set
of \nn s. A~set which is finite or \ct y infinite is called \ti{\ct e}\index{set!countable}
and a set which is not \ct e is called \ti{un\ct e}\index{set!uncountable}.
\edf

\bex4.22
Let $A,B$ be two nonempty sets \st $A\ax B$. Show that if $A$ is \ct e
(resp.\ un\ct e), then so is~$B$.
\eex

\brm4.23
A \ct y infinite set is \ep\ to a set of \nn s hence infinite by
Theorem~\rf{t4.18}\,(iii). Not all infinite sets are \ct y infinite.
\erm

\blm4.24
Let $E$ be a set of \nn s. Then the power set $\cP(E)$ is infinite.
\elm

\proof
Suppose for \cd ion that $\cP(E)$ is finite. Since $E\in\cP(E)$,
$\cP(E)$ is not empty. Then, by \E\df\ \rf{d4.17} and Lemma \rf{l4.15},
\te\ $n\in E$, an \il\
$[e,n]$ in $(E,\le)$, where $\le$ is the natural ordering of~$E$, and a
bi\jn\ $f:\cP(E)\to[e,n]$. Let $m:=S(n)>n$.
Then the map $h:[e,m]\to \cP(E)$ defined by
$h(k):=[e,k]$, $k\in[e,m]$, is in\jc\ by \er{3.46} with $a:=e$.
Hence the map $f\circ h: [e,m] \to
[e,n]$ is in\jc. By Lemma~\rf{l4.11}, we have $m\le n$, a~\cd ion.
Hence $\cP(E)$ is infinite.
\endproof

The fact that $\cP(E)$ is not \ct y infinite, i.e.\ is not \ep\ to~$E$ is a
con\sq\ of the \fw\ lemma.

\begin{lem}[Cantor]\lb{l4.25}
Let $A$ be a nonempty set. Then there is no sur\jn\ from $A$ onto $\cP(A)$.
\elm

\proof
Let $f$ be a map from $A$ into $\cP(A)$. Let $U$ be the subset of~$A$ defined
by
\beq4.11
U:=\{x\in A: x\notin f(x)\}.
\e
We claim that $U$ does not belong to the range of~$f$. Suppose for \cd ion
that \te s an \el~${a\in A}$ \st $f(a)=U$. Then either $a\in U$ or $a\notin
U$. If $a\in U$, then $a\notin f(a)=U$, hence $a\notin U$, a \cd ion. If
$a\notin U$, then $a\in f(a)=U$, hence $a\in U$, a \cd ion. \E\Tf there is no $a\in
A$ \st $f(a)=U$.
\endproof

\begin{thm}[Cantor] \lb{t4.26}
Let $(E,e,S)$ be a set of \nn s. Then $\cP(E)$ is un\ct e.
\eth

\proof
Indeed, $\cP(E)$ is infinite by
Lemma~\rf{l4.24}. There is no sur\jn\ from $E$ onto $\cP(E)$ by
Lemma~\rf{l4.25}. Hence $\cP(E)$ is not \ct y infinite.
\endproof

\sdim{[C]}
\bpr4.27
Let $X,Y$ be \ns s and let $f$ be a map from $X$ into~$Y$.

\ssk
\ite[\rm (A)] {Suppose $f$ is {\em in\jc}. Then}

\hph i,ii, if $X$ is infinite, then so is $Y$,

\hph ii,i, if $X$ is un\ct e, then so is $Y$,

\hph iii,, if $Y$ is finite, then so is $X$,

\hph iv,, if $Y$ is \ct e, then so is $X$.

\ite[\rm (B)] {Suppose $f$ is {\em sur\jc}. Then}

\hph i,ii, if $X$ is finite, then so is $Y$,

\hph ii,i, if $X$ is \ct e, then so is $Y$,

\hph iii,, if $Y$ is infinite, then so is $X$,

\hph iv,, if $Y$ is un\ct e, then so is $X$.
\epr

\proof
Using contraposition, we see that assertions (A)(i) and (A)(iii), (A)(ii) and (A)(iv),
(B)(i) and (B)(iii), (B)(ii) and (B)(iv) are \ev t. We prove (A)(iii), (A)(iv), (B)(i) and
(B)(ii). Let $(E,\le,e,S)$ be an \os\ of \nn s.
\ssk

(A)(iii): We have $f(X)\sbs Y$ and $Y$ finite. Hence $f(X)$ is finite by Theorem
\rf{t4.18}\,(i). Since $f:X\to X$ is a bi\jn, we have $X\ax f(X)$, hence $X$
is finite by Exercise \rf{ex4.19}.

\ssk
(B)(i): Since $X$ is finite, $f(X)$ is finite by Theorem \rf{t4.18}\,(ii). But
$Y=f(X)$.

\ssk
(A)(iv): In view of (A)(iii) we may assume $Y$ \ct y infinite. Let $\vf:Y\to E$
bi\jc\ be as in \E\df\ \rf{d4.17}. Since $f(X)\sbs Y$, we have $\vf(f(X))\sbs E$. From Lemma \rf{l4.20}
and Theorem \rf{t3.49}, we infer that either \te s $n\in E$ \st $\vf(f(X))
\ax [e,n]$, or $\vf(f(X))\ax E$. Since $f$ is in\jc, we have $X\ax f(X)$.
Using Exercise \rf{ex4.19}\,(i)(a), we find that $X$~is either finite or
\ct y infinite, hence $X$~is \ct e.

\ssk
(B)(ii): It is \sft\ to construct a map $g:Y\to X$ \st $f\circ g=\id_Y$. Indeed,
in this case $g$~is in\jc, since $g(a)=g(b)$, $a,b\in Y$ implies $a=f\circ
g(a)=f\circ g(b)=b$, $a,b\in Y$, and (B)(ii) follows from (A)(iv), with $f:=g$,
$X:=Y$ and $Y:=X$. Observe that if $X$ is finite, then (B)(ii) follows from
(B)(i), hence we may suppose $X$~\ct y infinite. Let $\psi:E\to X$ be bi\jc\
where $E$~is an \os\ of \nn s, and
set $\wt f:=f\circ\psi:E\to Y$. Observe that $\wt f$~is sur\jc\ since both
$\psi$~and~$f$ are sur\jc. In view of what precedes it is \sft\ to find a~map
$\wt g:Y\to E$ \st $\wt f\circ\wt g=\id_Y$. Given $y\in Y$, set $V(y):=
\{m\in E: \wt f(m)=y\}$. Since $\wt f$ is sur\jc, then $V(y)\ne\vn$ \fa
$y\in Y$. Let $y\in Y$. Then $V(y)$ possesses a~least \el\ in view of
Theorem \rf{t3.31}\,(iv). We denote this \el\ by $\wt g(y)$. The map $\wt g:Y
\to E$ defined by $y\mt \wt g(y)$, denoted by~$\wt g$, is well-defined and
\sf ies $\wt f(\wt g(y))=y$ since $\wt g(y)\in V(y)$, \fa $y\in Y$. Hence
$\wt f\circ\wt g=\id_Y$ and $Y$~is \ct y infinite.
\endproof

\brm4.27
Let $X$ and $Y$ be \ns s and let $f$ be a map from~$X$ into~$Y$. A~map $g:Y\to
X$ is called a \ti{right inverse} for~$f$ if it \sf ies $f\circ g=\id_Y$ and
a~\ti{left inverse} for~$f$ if it \sf ies $g\circ f=\id_X$. If $f$~has a~right
inverse, then $f$~is sur\jc. Indeed, in this case we have $f(g(y))=y$ \fa
$y\in Y$. If $f$~has a~left inverse, then $f$~is in\jc. Indeed, $f(a)=f(b)$
implies $a=g(f(a))=g(f(b))=b$, \fa $a,b\in X$. If $f$~is in\jc, then the map
$\wt f:X\to f(X)$ defined by $\wt f(x):=f(x)$, $x\in X$, is bi\jc. Hence \te s
\ooo $\wt g:f(X)\to X$ \sf ying $\wt g\circ\wt f=\id_X$, $\wt f\circ \wt g=
\id_{f(X)}$. It is called the \ti{inverse} of~$\wt f$. Now define a~map
$g:Y\to X$ by setting $g(y):=\wt g(y)$ if $y\in f(X)$ and $g(y):=x$ \fa
$y\in Y\sm f(X)$ (possibly empty) where $x$~is an arbitrary \el\ of~$X$
(recall that $X$~is not empty). Then we have $g\circ f=\id_X$. \E\Tf if
$f$~is in\jc, it has a left inverse (not necessarily unique). \E\oh if $f$~is
sur\jc, then the proof of (B)(ii) shows that $f$ has a~right inverse whenever
$X$~is \ct y infinite. The same holds if $X$~is finite (see the proof of
Theorem \rf{t4.18}\,(ii), (iii)). Hence if $X$~is \ct e and $f$~is sur\jc, then
$f$~has a~right inverse. There is no analogue for $X$ un\ct e. It can be shown
that the \Pr\ stating that every sur\jc\ map has a~right inverse is \ev t to
the so-called axiom of choice (see~\cite{Nrs}, Ch.~14). We now show that the
\ex\ of a~right inverse for a~sur\jc\ map is a con\sq\ of the axiom of choice.
Let R denote the \rl\ on~$X$ defined by $a \mr R b$ if $f(a)=f(b)$, $a,b\in
X$. One verifies that R is an \ti{\ev ce \rl} on~$X$ (see \E\df~\rf{d4.5}).
Then the set of \ev ce classes $\cA_{\rm R}$ is a \pt\ of~$X$. The axiom of
choice allows us to choose exactly one \el\ in each \ev ce class and to make
a~set~$Z$ consisting of these \el s. Since $f$ is sur\jc, every \ev ce class
is of the form $[x]$ where $f(x)=y$, \fs $y\in Y$. Define $g:Y\to Z\sbs X$ by
setting $g(y)$ is the unique \el\ $z\in Z$ \st $f(z)=y$. Then $f\circ g(y)=
f(z)=y$, $y\in Y$. Observe that in the proof of (B)(ii), $V(y):=\{m\in E:
\wt f(m)=y\}$, $y\in Y$, is the \ev ce class~$[m]$, where $\wt f(m)=y$, and
that $z$~is the least \el\ of~$[m]$.
\erm

%% file: DETOUR2.TEX
\Section{Arithmetic operations}[Arithmetic operations]\label{s.2}
\Subsubsection{Addition}\label{sss.Add}

\def\ees{(\wt E,\wt e,\wt S)}
\def\dqui#1#2{\mathrel{\lower1pt\hbox{$\stackrel{\lower7pt\hbox{\smash{\hbox
to0pt{$\scriptstyle\sqcap$\hss}$\scriptstyle\sqcup_{#1}$}}}{\cdot_{#2}}$}}}

In what follows $\ees$ denotes an arbitrary but \ti{fixed\/} set of \nn s
(see \E\df~\rfa1{d1.3}), and $\le$ denotes its \nog\ (see \E\df~\rfa1{d3.24}).
Recall that the \tr ity of the \rl~$\le$ on~$\wt E$ is a con\sq\ of the fact
that the \cm\ of \IT s of~$\wt S$ (\wrt $\wt E$) is also an \IT\ of~$\wt S$,
by \era1{2.27}, \era1{2.37}. The \uq\ part of \E\Pr~\rfa1{p2.7} \era1{2.37}
allows us to define a~binary \op\ on~$\wt E$, which we denote for the moment
by~$\tqu$, by setting \fa $m,n\in\wt E$:
\bml1.1
m\tqu n:=p \qh{where $p$ is the unique \el\ of $\wt E$}\\
\hbox{\st $\Phi(m)\circ \Phi(n)=\Phi(p)$ in \era1{2.37}}.
\e

\blm1.1
Let $\tqu$ be the binary \op\ on~$\wt E$ defined by~\er{1.1}. Then \fa
$m,n\in\wt E$ we have\dw
\bga1.2
m\tqu \wt e =m,\\
m\tqu \wt S(n)=\wt S(m\tqu n).\lb{1.3}
\e
\elm

\proof
Let $m,n\in \wt E$.

\er{1.2}: $\Phi(m\tqu \wt e)\nde1.1 = \Phi(m)\circ\Phi(\wt e) \nda12.33 =\Phi(m)\circ
\id_{\td E}\nda12.2 =\Phi(m)$. Hence $\Phi(m\tqu\wt e)=\Phi(m)$, and \er{1.2}
follows from \era1{2.35}.

Note that we have
\beq1.4
\wt S=\Phi(\wt S(\wt e)).
\e
Indeed, $\Phi(\wt S(\wt e))\nda12.34 = \wt S\circ\Phi(\wt e)\nda12.33 =
\wt S\circ\id_{\td E} \nda12.2 = \wt S$.

We now prove

\er{1.3}: $\Phi(m\tqu\wt S(n)) \nde1.1 = \Phi(m)\circ \Phi(\wt S(n)) \nda12.34
= \Phi(m)\circ(\wt S\circ\Phi(n)) \nde1.4 = \Phi(m)\circ(\Phi(\wt S(\wt e))
\circ\Phi(n)) \nda12.1 = (\Phi(m)\circ\Phi(\wt S(\wt e)))\circ \Phi(n)
\nda12.26 = (\Phi(\wt S(\wt e))\circ\Phi(m)\circ \Phi(n) \nde1.4 =
(\wt S\circ\Phi(m))\circ\Phi(n)
\nda12.1 = \wt S\circ(\Phi(m)\circ\Phi(n)) \nde1.1 = \wt S\circ(\Phi(m\tqu n
)) \nda12.34 = \Phi(\wt S(m\tqu n))$. Hence $\Phi(m\tqu\wt S(n))=
\Phi(\wt S(m\tqu n))$, and \er{1.3} follows from \era1{2.35}.
\endproof

It turns out that there is \ti{at most\/} one binary \op\ on~$\wt E$ \sf ying
\er{1.2}, \er{1.3}. Motivated by \cite[Chapter~1, \S2,3]{Nrs} we shall call
this \op\ \ti{\ad} on~$\wt E$.\index{addition}

\bpr1.2
\E\te s exactly one binary \op\ on $\ees$, called \emph{\ad\ on }$\wt E$,
denoted by $+_{\td E}$ $($or simply by~$+$ when no confusion arises$)$,
\sf ying \er{1.2}, \er{1.3}. If $a,b\in \wt E$, then $a+b$ is called the
\emph{sum} of $a$ and~$b$.
\Mo the \fw\ \pp ies hold\dw

\E\fa $a,b,c\in \wt E$ we have
\bea1.5
&(a\qu b)\qu c = a\qu(b\qu c),\\
&a\qu b = b\qu a, \lb{1.6}\\
&a\qu e= e\qu a=a, \lb{1.7}\\
&a\qu c=b\qu c \qh{implies} a=b, \lb{1.8}\\
&a\qu b=e \qh{implies} a=b=e, \lb{1.9}
\e
where $\qu$ is replaced by $+_{\td E}$, and $e$ by~$\wt e$.
\epr

\proof \

\ti{\E\ex}: follows from Lemma \rf{l1.1}.

\ti{\E\uq}: Let $\tqu$ and $\qu'$ be binary \op s on~$\wt E$ \sf ying \er{1.2}
and \er{1.3}, and let $m\in\wt E$. Define $\F$ (resp.~$\F'):\wt E\to \wt E$
by setting $\F(n):=m\tqu n$ (resp.\ $\F'(n):=m\qu'n$), $n\in \wt E$. Set
$A:=\{n\in\wt E:\F(n)=\F'(n)\}$. We have $\wt e\in A$. Indeed, $\F(\wt e)
\nde1.2 = m\nde1.2 = \F'(\wt e)$. We next show that if $n\in A$, then $\wt S
(n)\in A$. Let $n\in A$. Then $\F(\wt S(n)) \nde1.3 = \wt S(\F(n)) \nad{n\in A}
= \wt S(\F'(n))\nde1.3 = \F'(\wt S(n))$. Hence $\wt S(n)\in A$, and $A=E$
follows from (1.N2).

\er{1.5}: $\Phi((a+b)+c)\nde1.1 = \Phi(a+b)\circ\Phi(c) \nde1.1 =
(\Phi(a)\circ\Phi(b))\circ\Phi(c)\nda12.1 = \Phi(a)\circ(\Phi(b)\circ\Phi(c))
\nde1.1 = \Phi(a)\circ\Phi(b+c) \nde1.1 = \Phi(a+(b+c))$. Then \er{1.5}
follows from \era1{2.35}.

\er{1.6}: $\Phi(a+b)\nde1.1 = \Phi(a)\circ\Phi(b) \nda12.26 = \Phi(b)\circ
\Phi(a)\nde1.1 = \Phi(b+a)$. Then \er{1.6} follows from \era1{2.35}.

\er{1.7}: $\Phi(a+\wt e)\nde1.1 = \Phi(a)\circ\Phi(\wt e)\nda12.33 = \Phi(a)
\circ \id_{\td E}\nda12.2 = \Phi(a)$. Then \er{1.7} follows from \era1{2.35}
and \er{1.6}.

\er{1.8}: From $a+c=b+c$ we infer $\Phi(a)\circ\Phi(c)=\Phi(b)\circ\Phi(c)$
by~\er{1.1}. Hence $\Phi(a)c\nda12.35 = \Phi(a)(\Phi(c)\wt e)=(\Phi(a)\circ
\Phi(c))\wt e=(\Phi(b)\circ\Phi(c))\wt e= \Phi(b)(\Phi(c)\wt e) \nda12.35 =
\Phi(b)c$. \E\Tf $\Phi(a)c = \Phi(b)c$, and $a=b$ follows from \era1{2.39}.

\er{1.9}: We have $\Phi(a+b)=\Phi(\wt e)$. Using \er{1.1} and \era1{2.33} we
obtain $\Phi(a)\circ\Phi(b)=\id_{\td e}$, hence \er{1.9} follows from
\era1{2.41}.
\endproof

\brm1.3 \

\hph i,i, \If from \E\Pr\ \rf{p1.2} that $(\wt E,+_{\td E},\wt e)$ is an abelian
monoid (see \E\df~\rfa1{d2.2}).

\hph ii,, In the sequel we shall refer to \er{1.5} (resp.\ \er{1.6}, \dots,
\er{1.9}) to indicate that a binary \op\ \sf ies \er{1.5} (resp.\ \er{1.6}, \dots,
\er{1.9}). For example the reference to \er{1.5} means that the binary \op\ in
consideration is \asc e.
\erm

We now introduce two (non-standard) \df s.

\bds1.4
An \ti{abelian} monoid $(M,\qu,e)$ \sf ying \er{1.8} is called a \ti{\Cm\/}\break
(C~stands for \ti{\cnc e}). A~\Cm\ \sf ying \er{1.9} is called\index{C-monoid}\index{cancellativity}
a~\ti{P-monoid\/}\index{P-monoid} (P~stands for \ti{\po}).\index{positivity}
\eds

\bxs1.5 \

\allowdisplaybreaks
\hph i,ii, $(\wt E,+_{\td E},\wt e)$ is a P-monoid.

\hph ii,i, The trivial monoid $(\{e\},\qu,e)$ is a P-monoid.

\hph iii,, Let $\{\a,\b\}$ be the set consisting of two \el s $\a$~and~$\b$.
Suppose that $\qu$ is a~binary \op\ on $\{\a,\b\}$ \st $\a$~is a~\nel. We
thus have $\a\qu\a=\a$, $\a\qu\b=\b$ and $\b\qu\a=\b$. \E\Tf $\qu$~is
completely defined by either $\b\qu\b=\b$ or $\b\qu\b=\a$. We denote by
$\qu^{(1)}$ (resp.~$\qu^{(2)}$) the \op\ \crs\ to the first case (resp.\ the
second case):
\[
\arraycolsep5pt
\begin{array}{c|cc}\qu^{(1)}&\a&\b\\
\hline
\a&\a&\b\\
\b&\b&\b
\end{array}
\qquad
\begin{array}{c|cc}\qu^{(2)}&\a&\b\\
\hline
\a&\a&\b\\
\b&\b&\a
\end{array}
\]
Clearly both $\qu^{(1)}$ and $\qu^{(2)}$ are \ti{\cmt e} (see \era1{2.22}).
We next show that both are \ti{\asc e} (see \era1{2.20}):
\begin{align*}
(\a \qu \a) \qu\a &= \a \qu\a = \a \qu (\a \qu\a ),\\
(\b \qu \a) \qu\a &= \b \qu\a = \b \qu (\a \qu\a ),\\
(\a \qu \b) \qu\a &= \b \qu\a = \a \qu (\b \qu\a ),\\
(\a \qu \a) \qu\b &= \a \qu\b = \a \qu (\a \qu\b ),\\
(\b \qu \b) \qu\a &= \b \qu\b = \b \qu (\b \qu\a ),\\
(\b \qu \a) \qu\b &= \b \qu\b = \b \qu (\a \qu\b ),\\
(\a \qu \b) \qu\b &= \b \qu\b = \a \qu (\b \qu\b ),\\
(\b \qu^{(1)} \b) \qu^{(1)}\b &= \b \qu^{(1)}\b = \b \qu^{(1)} (\b \qu^{(1)}\b ),\\
(\b \qu^{(2)} \b) \qu^{(2)}\b &= \a \qu^{(2)}\b =
\b\qu^{(2)}\a = \b \qu^{(2)} (\b \qu^{(2)}\b ).
\end{align*}

It follows that both $(\{\a,\b\},\qu^{(1)},\a)$ and
$(\{\a,\b\},\qu^{(2)},\a)$ are \ti{abelian} monoids, with \nel~$\a$.

Since $\a\qu^{(1)}\b= \b\qu^{(1)}\b$ and $\a\ne \b$,
$(\{\a,\b\},\qu^{(1)},\a)$ is \ti{not\/} a \Cm, hence also \ti{not\/}
a~P-monoid. However, we have $a\qu^{(1)}b=\a$, $a,b\in\{\a,\b\}$ iff
$a=b=\a$, hence condition \er{1.9} holds.

We now show that $(\{\a,\b\},\qu^{(2)},\a)$ is a~\Cm. Indeed, let $a,b\in
\{\a,\b\}$ be \st (i)~$a\qu^{(2)}\a=b\qu^{(2)}\a$. Then $a=a\qu^{(2)}\a
=b\qu^{(2)}\a=b$; (ii)~$a\qu^{(2)}\b=b\qu^{(2)}\b$. Then, if $a:=\a$, we have
$b\qu^{(2)}\b= a\qu^{(2)}\b=\b$. From $b\qu^{(2)}\b=\b$, we infer $b=\a=a$.
If $a:=\b$, we have $b\qu^{(2)}\b = a\qu^{(2)}\b=\a$. From
$b\qu^{(2)}\b=\a$, we infer $b=\b=a$. \E\Tf the monoid
$(\{\a,\b\},\qu^{(2)},\a)$ is a \Cm. Since $\b\qu^{(2)}\b=\a$ and $\b\ne\a$,
it is \ti{not\/} a~P-monoid.

\ssk
\hph iv,, Let $(X_i,\qu_i,e_i)$, $i=1,2$, be \Cm s (resp.\ P-monoids). Define
on $X_1\t X_2$ a~binary \op~$\qu$ by setting $(a_1,a_2)\qu(b_1,b_2):=
(a_1\qu_1b_1,a_2\qu_2\b_2)$ where $a_1,b_1\in X_1$ and $a_2,b_2\in X_2$. We
claim that
\beq1.10
(X_1\t X_2,\qu,(e_1,e_2)) \hbox{ is a \Cm\ (resp.\ P-monoid).}
\e

\ti{\asc ity}: Let $a_1,b_1,c_1\in X_1$ and $a_2,b_2,c_2\in X_2$. Then
$((a_1,a_2)\qu(b_1,b_2))\qu(c_1,c_2)=((a_1\qu_1b_1,a_2\qu_2b_2))\qu
(c_1,c_2)=((a_1\qu_1b_1)\qu_1c_1,(a_2\qu_2b_2)\qu_2c_2)\nad\ast=\break
(a_1\qu_1(b_1\qu_1c_1), a_2\qu_2(b_2\qu_2c_2)) = (a_1,a_2)\qu
((b_1\qu_1c_1),(b_2\qu_2c_2)) = (a_1,a_2) \qu ((b_1,b_2)\qu(c_1,c_2))$.
In $\nad\ast=$ we used the \asc ity of $\qu_1$ and~$\qu_2$.

\ti{\nel\/}: Let $a_1\in X_1$, $a_2\in X_2$. Then $(a_1,a_2)\qu(e_1,e_2)=
(a_1\qu_1e_1,a_2\qu_2e_2)=(a_1,a_2) = (e_1\qu_1a_1, e_2\qu_2a_2)=(e_1,e_2)
\qu(a_1,a_2)$.

\ti{\cmt ity}: Let $a_1,b_1\in X_1$ and $a_2,b_2\in X_2$. Then $(a_1,a_2)\qu
(b_1,b_2)=(a_1\qu_1b_1,\break a_2\qu_2b_2) \nad\ast= (b_1\qu_1a_1,b_2\qu_2a_2)=
(b_1,b_2)\qu(a_1,a_2)$. In $\nad\ast=$ we used the \cmt ity of $\qu_1$
and~$\qu_2$.

\ti{\cnc ity}: Let $a_1,b_1,c_1\in X_1$ and $a_2,b_2,c_2\in X_2$ be \st
$(a_1,a_2)\qu(c_1,c_2)=(b_1,b_2)\qu(c_1,c_2)$. Then we have
$(a_1\qu_1c_1,a_2\qu_2c_2)=(b_1\qu_1c_1,b_2\qu_2c_2)$, hence $a_1\qu_1c_1
=b_1\qu_1c_1$, $a_2\qu_2c_2=b_2\qu_2c_2$, by the \df\ of the cartesian
product $X_1\t X_2$.

From the \cnc ity of $\qu_1$ and $\qu_2$, we infer that $a_1=b_1$ and $a_2
=b_2$, hence $(a_1,a_2)=(b_1,b_2)$.

\ti{positivity}: Let $a_1,b_1\in X_1$ and $a_2,b_2\in X_2$ be \st $(a_1,a_2)
\qu(b_1,b_2)=(e_1,e_2)$. Then $(a_1\qu_1b_1,a_1\qu_2b_2)=(e_1,e_2)$,
hence $a_1\qu_1b_1=e_1$ and $a_2\qu_2b_2=e_2$. From the positivity of $\qu_1$
and~$\qu_2$, we obtain $a_1=b_1=e_1$, $a_2=b_2=e_2$, hence $(a_1,a_2)
=(e_1,e_2)=(b_1,b_2)$.

This completes the proof of \er{1.10}.

\ssk
\hph v,i, Let $(M,\qu,e)$ be a monoid, let $M'$ be a~\ns. Suppose \te s
$f:M\to M'$ bi\jc. Set
\beq1.11
e':=f(e); \q a'\qu'b':=f(f\Inv(a')\qu f\Inv(b')) \qh{\fa}a',b'\in M'.
\e
Then
\bga1.12
(M',\qu',e') \qh{is a monoid.} \\
\hbox{If in \ad\ $M$ is an (abelian, resp.\ C-, P-) monoid,}
\hbox{ then so is $(M',\qu',e')$.}\lb{1.13}
\e

\ti{\asc ity}: Let $a',b',c'\in M'$. Set $a:=f\Inv(a')$, $b:=f\Inv(b')$, and
$c:=f\Inv(c')$. Observe that
\beq1.14
a'\qu'b' = f\bigl(f\Inv(a')\qu f\Inv(b')\bigr)=f(a\qu b).
\e
\looseness=-1
Hence $(a'\qu'b')\qu'c' \nde1.14 = f\bigl(\bigl(f\Inv(f(a\qu b))\bigr)
\qu f\Inv(c')\bigr) = f((a\qu b)\qu c)$. \E\oh $a'\qu'(b'\qu'c') \nde1.14 =
f\bigl((f\Inv(a'))\qu f\Inv(f(b\qu c))\bigr) = f(a\qu(b\qu c))$.
Using the \asc ity of~$\qu$, we find $f((a\qu b)\qu c)=
f(a\qu(b\qu c))$. Hence $(a'\qu'b')\qu'c' = a'\qu'(b'\qu'c')$.

\ti{\nel\/}: Let $a'\in M'$ and let $a:=f\Inv(a')$. Then $a'\qu'e' \nde1.14 =
f(a\qu e)=f(a)=a' =f(a)=f(e\qu a)\nde1.14 =
e'\qu'a'$. Hence $a'\qu'e'=a'=e'\qu'a'$.

\ti{\cmt ity}: Let $a',b'\in M'$ and let $a:=f\Inv(a')$, $b:=f\Inv(b')$.
Then $a'\qu'b'=f(a\qu b)\nad\ast= f(b\qu a)=b'\qu'a'$. In $\nad\ast=$ we used
the \cmt ity of~$\qu$.

\ti{\cnc ity}: Let $a',b',c'\in M'$ and let $a:=f\Inv(a')$, $b:=f\Inv(b')$,
$c:=f\Inv(c')$. We suppose that $(M,\qu,e)$ is a \Cm\ and that $a'\qu'c'=
b'\qu'c'$. Then by \er{1.14}, we have $a'\qu'c'=f(a\qu c)$ and $b'\qu'c'=
f(b\qu c)$. From $f(a\qu c)=f(b\qu c)$ and the bi\ji\ of~$f$, we infer
$a\qu c=b\qu c$, hence $a=b$. Thus $a'=f(a)=f(b)=b'$.

\ti{positivity}: Let $a',b'\in M'$ and let $a:=f\Inv(a')$, $b:=f\Inv(b')$.
If $a'\qu'b'=e'$, then $f(a\qu b)=f(e)$, hence $a\qu b=e$ in view of the
bi\ji\ of~$f$, and $a=b=e$ in view of the positivity of~$\qu$. Hence $a'=f(a)
=f(b)=f(e)=b'=e'$.
\exs

\brm1.6 \

\hph i,ii, \If from the proof of \er{1.10} that if $(X_i,\qu_i,e_i)$, $i=1,2$,
are monoids (not necessarily abelian), then $(X_1\t X_2,\qu,(e_1,e_2))$ is
a~monoid. It is called the (\ti{direct\/}) \ti{product of the monoids} $X_1$
and~$X_2$. If $X_1$ and $X_2$ are abelian (resp.\ C-, resp.~P-) monoids, then
so is their direct product.

\hph ii,i, The trivial monoid $(\{e\},\qu,e)$ is a P-monoid.

\hph iii,, Let $M,M'$ and $f:M\to M'$ be as in Example \rf{xa1.5}\,(v). \If
from \er{1.11} and \er{1.14} that $f(e)=e'$ and $f(a\qu b)=f(a)\qu f(b)$ \fa
$a,b\in M$. \Mo $f\Inv:M'\to M$ \sf ies $f\Inv(e')=e$ and $f\Inv(a'\qu'b')=
f\Inv(a')\qu f\Inv(b')$ \fa $a',b'\in M'$. Such a~map $f:M\to M'$ will be
called a~\ti{monoid-\is sm}\index{monoid-isomorphism}.
\erm

\bdf1.7
Let $(X,\qu,e)$ and $(\wt X,\tqu,\wt e)$ be monoids not necessarily distinct.
A~map $\vf:X\to \wt X$ \sf ying
\bga1.15
\vf(e)=\wt e,\\
\vf(x\qu y)=\vf(x)\tqu\vf(y) \qh{\fa} x,y\in X, \lb{1.16}
\e
is called a \ti{monoid-homo\mf}\index{monoid-homomorphism}. If, in addition, \te s a homo\mf\break $\psi:\wt
X\to X$ \st $\psi\circ\vf=\id_X$ (resp.\ $\vf\circ\psi=\id_{\td X}$) then
$\vf$~is called a \ti{monoid-mono\mf} (resp.\ \ti{monoid-epi\mf}).

If $\vf$ is both a monoid-mono\mf\ and a monoid-epi\mf\ then it is called a
\ti{monoid-iso\mf}. If \te s a monoid-iso\mf\ between the monoids $(X,\qu,e)$
and $(\wt X,\tqu,\wt e)$ then $X$ and~$\wh X$ are called \ti{monoid-\is c}.

If the context is clear, then the word \ti{monoid\/} is usually deleted.

A \hm sm (resp.\ \is sm) from a monoid into itself is usually called
\ti{an endo\mf} (resp.\ \ti{auto\mf}).\index{endomorphism}\index{automorphism}
\edf

\blm1.8 \

\hph i,ii, A monoid-\hm sm is a monoid-\is sm iff it is bi\jc.

\hph ii,i, The \cm\ of two monoid-\hm sms $($resp.\ \is sms$)$ is a\break monoid-\hm sm
$($resp.\ \is sm$)$.

\hph iii,, The inverse of a monoid-\is sm is a monoid-\is sm.

\hph iv,, The identity in a monoid is a monoid-\is sm.

\hph v,i, A monoid monoid-\is c to an abelian monoid is abelian.
\elm

\proof \

(i) \ti{If\/}: Let $(X,\qu,e)$ and $(X',\qu',e')$ be monoids (not necessarily distinct).
Let $\vf:X\to X'$ be a bi\jc\ monoid-\hm sm. Since $\vf$ is bi\jc, \te s one
(and only one) $\psi:X'\to X$ \st $\psi\circ\vf=\id_X$ and
$\vf\circ\psi=\id_{X'}$. It is \sft\ to show that $\psi:X'\to X$ is a \hm sm.
Since $\vf(e)=e'$ we have $\psi(e')=\psi\circ\vf(e)=e$. Let $x',y'\in X'$. We
have to show that $\psi(x'\qu'y')=\psi(x')\qu \psi(y')$. Set $x:=\psi(x')$ and
$y:=\psi(y')$. Then $\vf(x)=\vf\circ \psi(x')=x'$ and $\vf(y)=\vf\circ
\psi(y')=y'$. \Mo $x'\qu' y'=\vf(x)\qu'\vf(y)\nde 1.16 =\vf(x\qu y)$. Hence
$\psi(x'\qu'y')=\psi\circ\vf(x\qu y)=x\qu y=\psi(x')\qu\psi(y')$.

\ti{Only if\/}: Clear since $\psi\circ\vf=\id_X$ and $\vf\circ\psi= \id_{\td
X}$.

(ii) Let $(X_i,\qu_i,e_i)$, $i=1,2,3$, be monoids. Let $\vf_1:X_1\to X_2$ and
$\vf_2:X_2\to X_3$ be monoid-\hm sms. Let $\vf:=\vf_2\circ\vf_1$, and $a,b\in
X_1$. Then

\leavevmode\llap{1) }$\vf(e_1)=\vf_2\circ\vf_1(e_1)=\vf_2(\vf_1(e_1))=\vf_2(e_2)=e_3$.

{\hangindent\parindent\hangafter1
\leavevmode\llap{2) }$\vf_1(a\qu_1 b)=\vf_1(a)\qu_2\vf_1(b)$ and $\vf(a\qu_1 b)=
\vf_2(\vf_1(a)\qu_2\vf_1(b))=\vf_2\circ\vf_1(a) \qu_3 \vf_2\circ\vf_1(b) =
\vf(a)\qu_3\vf(b)$. If $\vf_1,\vf_2$ are \is sms, then $\vf_2\circ\vf_1$ is
bi\jc\ and (ii) follows from~(i).

}(iii) Let $\vf$ and $\psi$ be as in \E\df\ \rf{d1.7}. Then $\psi$ is
a \hm sm and $\psi$ is bi\jc, hence $\psi$ is an \is sm by~(i).  But $\vf$~is
bi\jc\ and $\vf\Inv=\psi$.

(iv) Clear.

(v) Let $\vf:(X_1,\qu_1,e_1) \to (X_2,\qu_2,e_2)$ be a monoid-\is sm and let
$(X_1,\qu_1,e_1)$ be abelian. Then for all $a,b\in X_2$ we have $a\qu_2 b=
\vf\circ\vf\Inv(a\qu_2 b)=\vf(\vf\Inv(a) \qu_1 \vf\Inv(b)) = \vf(\vf\Inv
(b)\qu_1\vf\Inv(a)) = \vf\circ\vf\Inv(b) \qu_2 \vf\circ\vf\Inv(a)
=b\qu_2 a$.
\endproof

\bex1.9 \

\hph i,i, Show that a monoid which is monoid-\is c to a \Cm\ (resp.\ P-monoid)
is a~\Cm\ (resp.\ P-monoid).

\hph ii,, Let $M_1$ be a monoid and $M_2$ be a \Cm. Let $\vf:M_1\to M_2$ \sf y
\er{1.16}. Show that \er{1.15} holds. Give an example of an abelian
monoid~$M_2$ which is not a~\Cm\ and of a~map $\vf:M_1\to M_2$ \sf ying
\er{1.16} but not \er{1.15}.
\eex

\blm1.10
Let $(E',e',S')$ be a set of \nn s $($see \E\df\ \rfa1{d1.3}$)$, let $(\wt E,
+_{\td E},\wt e)$ be the P-monoid introduced in \E\Pr\ \rf{p1.2} $($see Example
\rf{xa1.5}\,{\rm (i))} and let $(E',+_{E'},e')$ be the \crs\ P-monoid introduced
in \E\Pr\ \rf{p1.2}. Let $\vf:\wt E\to E'$ denote the unique bi\jn\ \sf ying
$\vf(\wt e)=e'$ and $\vf\circ\wt S=S'\circ\vf$ introduced in Theorem \rfa1{t1.4}.
Then $\vf$~is a monoid-\is sm from $(\wt E,+_{\td E},\wt e)$ onto
$(E',+_{E'},e')$.
\elm

\proof
In view of Lemma \rf{l1.8}\,(i), it is \sft\ to show that $\vf$~\sf ies
\er{1.15}, \er{1.16}, since $\vf$~is bi\jc. Clearly \er{1.15} holds. We next
show that $\vf(m+_{\td E}n)=\vf(m)+_{E'}\vf(n)$ \fa $m,n\in\wt E$. Let
$m\in\wt E$ be fixed but arbitrary and set
\[
A:=\bigl\{n\in\wt E: \vf(m+_{\td E} n)=\vf(m)+_{E'}\vf(n)\bigr\}.
\]
In view of (1.N2) it is \sft\ to show that $A$ is \iv\ in $\ees$. We have $\wt e
\in A$. Indeed, $\vf(m+_{\td E}e) \nad\ast= \vf(m) \nad{\ast\ast}= \vf(m)
+_{E'}e' \nde1.15 = \vf(m)+_{E'}\vf(e)$. In $\nad\ast=$ (resp.~$\nad{\ast\ast}
=$) we used the fact that $\wt e$ (resp.~$e'$) is the \nel\ of $(\wt E,
+_{\td E},\wt e)$ (resp.\ $(E',+_{E'},e')$). Hence $\wt e\in A$.

We now show that if $n\in A$, then $S(n)\in A$. We suppose that $n\in A$. Then
$\vf(m+_{\td E}\wt S(n)) \nde1.3 = \vf(\wt S(m+_{\td E}n)) \nad\ast=
S'\vf(m+_{\td E}n) \nad{n\in A}= S'(\vf(m)+_{E'}\vf(n)) \nde1.3 = \vf(m)
+_{E'}S'\vf(n) \nad\ast= \vf(m)+_{E'}\vf(\wt S(n))$. In $\nad\ast=$ we used
$\vf\circ\wt S = S'\circ\vf$. \If that $\wt S(n)\in A$, hence $A$~is \iv\ and
by (1.N2) $A=\wt E$.
\endproof

\brm1.11
Lemma \rf{l1.10} tells us that all P-monoids of the form $(E,+_E,e)$ where
$(E,e,S)$ is a set of \nn s are {\it\is c}. Clearly the trivial P-monoid\break $(\{\wt e\},
+_{\td E},\wt e)$ is not \is c to $(\wt E, +_{\td E},\wt e)$ since $\{\wt e\}$
and~$\wt E$ are not \ep. It turns out that every \ti{nontrivial\/} P-monoid
is \ti{infinite}. Indeed, let $(X,\qu,e)$ be a \ti{nontrivial\/}\break P-monoid and
let $a\in X\sms e$. Define a~map $\th_a:X\to X$ by setting $\th_a(x):=a\qu x$,
$x\in X$. Then $\th_a$~is in\jc\ by \cnc ity \er{1.8} and not sur\jc\ by
``positivity'' \er{1.9}. There is no $x\in X$ \st $a\qu x=e$, since $a\ne e$.
\E\Tf $\th_a$~is \ti{not sur\jc}.

By Theorem \rfa1{t4.18}\,(iii), $X$ is infinite. However, not all infinite
P-monoids are \is c. For example, if $(X_1,\qu_1,e_1):=(\wt E,+_{\td E},\wt e)$
and $(X_2,\qu_2,e_2):=(X_1,\qu_1,e_1)$, then the product $(\wt E\t\wt E,\qu,
(e,e))$, introduced in Examples \rf{xa1.5}\,(iv) (see Remark~\rf{r1.6}\,(i)),
is an \ti{infinite} \hbox{P-monoid} which is \ti{not\/} \is c to $(\wt E,+_{\td E},
\wt e)$. Indeed, since $(X_i,\qu_i,e_i)$, $i=1,2$, are \hbox{P-monoids} by \E\Pr\
\rf{p1.2}, $(\wt E\t\wt E,\qu,(\wt e,\wt e))$ is also a P-monoid by \er{1.10}.
\Mo since $(\wt S(\wt e),\wt e)\ne(\wt e,\wt e)$ by \era1{3.13},
$(\wt E\t\wt E,\qu,(\wt e,\wt e))$ is nontrivial, hence infinite by what
precedes. Now, suppose for \cd ion that \te s an \is sm $f:\wt E\to \wt E\t
\wt E$. Let $m:=f\Inv((\wt S(\wt e),\wt e))$ and $n:=f\Inv((\wt e,\wt S
(\wt e)))$. Since $(\wt S(\wt e),\wt e)\ne (\wt e,\wt S(\wt e))$, and $f$~is
bi\jc, we have $m\ne n$. In view of \era1{3.20} we have either $m<n$ or $m>n$.
Suppose $m<n$, then by Lemma \rfa1{l3.5} \te s $p\in\wt E$ \st $n=\Phi(p)m$.
Using \era1{2.35} we have $\Phi(n)\wt e=\Phi(p)(\Phi(m)\wt e)=\bigl(\Phi(p)
\circ\Phi(m))\wt e$. Using \er{1.1} and \E\Pr\ \rf{p1.2} we obtain
$\Phi(n)\wt e=\Phi(p+m)\wt e$. Hence by \era1{2.35} again, we arrive at
$n=p+m$. \csq\ $(\wt e,\wt S(\wt e))=f(n)=f(p+m)\nde1.16 = f(p)\qu f(m)=
(a,b)\qu(\wt S(\wt e),\wt e)$ where $(a,b):=f(p)$ \fs $a,b\in\wt E$. \If that
$a+\wt S(\wt e)=\wt e$ and $\wt S(\wt e)=b+\wt e=b$. Therefore, $(a+b)+\wt e
\nde1.5 = a+(b+\wt e)=a+\wt S(\wt e)=\wt e=\wt e+\wt e$. By \er{1.8} we obtain
$a+b=\wt e$, hence by \er{1.9} $a=b=\wt e$. Thus $\wt e=b=\wt S(\wt e)$.
A~\cd ion since $\wt e\ne \wt S(\wt e)$ by \era1{3.13}. The case $m>n$ is similar. \csq\
$(\wt E,+_{\td E},\wt e)$ is \ti{not\/} \is c to $(\wt E\t\wt E,\qu,(\wt e,
\wt e))$.
\erm

Our next goal is to find an \ad al \cn\ on a nontrivial P-monoid which makes
it \is c to $(\wt E,+_{\td E},\wt e)$. To this end we introduce the notion of
\ti{\IT} of an \el\ of a~monoid.

Informally, if $(M,\qu,e)$ is a monoid (not necessarily abelian) and if
$a\in M$, then $a$, $a\qu a$, $a\qu (a\qu a)=(a\qu a)\qu a=:a\qu a\qu a$,
$a\qu(a\qu a\qu a)=(a\qu a\qu a)\qu a=:a\qu a\qu a\qu a$ will be called \IT s
of~$a$ in $(M,\qu,e)$. It will be convenient to call~$e$ an \IT\ of~$a$.

We shall see in Exercise \rf{ex2.19} an example of a~set~$M$ together with
a~non\asc e binary \op~$\qu$ \st \te s $a\in M$ \sf ying $(a\qu a)\qu a
\ne a\qu(a\qu a)$. Of course in that case $a\qu a\ne a$, hence $\qu$~is also
not \cmt e. If $\qu$~is an \asc e binary \op\ on a nonempty set~$M$,
then $(M,\qu)$ is a~\sg\ (see \E\df\ \rfa1{d2.2}). If $(M,\qu)$ is not a~monoid
and $A$~is a~set consisting of only one \el\ $a\notin M$, i.e.\ $A=\{a\}$ and
$M\cap A=\vn$, then we can define an \ext\ of $(M,\qu)$ which we denote by
$(\wt M,\tqu)$ where $\wt M:=M\cup\{a\}$, $x\tqu y:=x\qu y$ whenever $x,y\in
M$, $a\tqu x=x$, $x\tqu a=x$ \fa $x\in M$ and $a\tqu a:=a$.

\bex1.12 \

\hph i,i, Show that $(\wt M,\tqu,a)$ defined above is a monoid.

\hph ii,, Let $M$ be a \ns\ containing an~\el~$e$, \st $M\sms e\ne\vn$, and
let $\qu$~be a binary \op\ \sf ying $a\qu e=e\qu a=a$ \fa $a\in M$. Show that
$(a\qu b)\qu c=a\qu(b\qu c)$ \fa $a,b,c\in M$ whenever $a$ or~$b$ or~$c$
is equal to~$e$.
\eex

\bpr1.13
Let $\ees$ be a set of \nn s, let $(M,\qu,e)$ be a monoid not necessarily
abelian, and let $a$~be an \el\ of~$M$. Then \te s \ooo map $\vf_a$
$($resp.~${}_a\vf):\wt E\to M$ \sf ying
\bea1.17
{}&\vf_a(\wt e):=e \qh{$($resp.\ }{}_a\vf(\wt e):=e) \\
&\vf_a(\wt S(n)):=a\qu\vf_a(n) \q (\hbox{resp.\ }{}_a\vf(\wt S(n)):={}_a\vf
(n)\qu a) \qh{\fa}n\in\wt E. \lb{1.18}
\e
\Mo \fa $m,n\in\wt E$ we have
\bea1.19
{}&\vf_e(n)=e, \\
&\vf_a(m)\qu\vf_a(n) = (\Phi_a(m)\circ\Phi_a(n))(e), \lb{1.20} \\
&\vf_a(m)\qu\vf_a(n) = \vf_a(n)\qu\vf_a(m), \lb{1.21} \\
&\vf_a(n) = {}_a\vf(n), \lb{1.22}
\e
where $\Phi_a(m),\Phi_a(n)$ are defined in \er{1.30}, \er{1.31}.
If, in \ad, $b$ is an \el\ of~$M$ \sf ying $a\qu b=b\qu a$ $(a$~and~$b$
{\rm commute)}, then we have \fa $m,n\in \wt E$\dw
\bea1.23
{}&\vf_a(m) \qu \vf_b(n) = \vf_b(n)\qu\vf_a(m), \\
&\vf_{a\qum b}(n) = \vf_a(n) \qu \vf_b(n). \lb{1.24}
\e
\epr

\def\namespec{Definition and Notation}
\begin{dspc} \lb{d1.14}
Let $\wt E,a$ and $\vf_a$ be as in \E\Pr\ \rf{p1.13}, and let $n\in\wt E$.
Then $\vf_a(n)$ is called the $n$-th (or~$n$-fold) \ti{\IT\ of~$a$ \wrt
$\wt E$ in $(M,\qu,e)$}, and is denoted by (not standard notation)\glossary{$n\dqu_{\td E}a$}
\beq1.25
n\dqu_{\td E}a:=\vf_a(n)
\e
(or simply $n\dqu a$ when no confusion arises).\glossary{$n\dqu a$}

We shall use the notation
\beq1.26
I_{\td E}(a):= \bigl\{n\dqu_{\td E}a \in M: n\in\wt E\bigr\}.
\e
\end{dspc}

\bpr1.15
Let $\ees$ and $(E',e',S')$ be sets of \nn s and let $\vf:\wt E\to E'$ be the
bi\jc\ map introduced in Theorem \rfa1{t1.4}. Let $(M,\qu,e)$, $a\in M$, be
as in \E\Pr\ \rf{p1.13} and let $I_{\td E}(a)$, $I_{E'}(a)$ be as in
\E\df s \rf{d1.14}. Let $\vf_a':E'\to M$ be the analogue of~$\vf_a$ where
$\wt E$~is replaced by~$E'$ in \E\Pr\ \rf{p1.13}. Then we have
\bga1.27
\vf_a'(\vf(m)) = \vf_a(m) \qh{\fa} m\in\wt E, \\
I_{E'}(a) = I_{\td E}(a). \lb{1.28}
\e
\epr

\bdf1.16
Let $(M,\qu,e)$ be a monoid and let $a\in M$. Then in view of \er{1.28} we
may define the \ti{set of \IT s of~$a$ in} $(M,\qu,e)$ (notation: $I(a)$) by\glossary{$I(a)$}
\beq1.29
I(a):=I_{\td E}(a) \qh{where $\wt E$ is any set of \nn s.}
\e
\edf

\proof[Proof of \E\Pr\ \rf{p1.13}]
The \ex\ and \uq\ of the map $\vf_a$ (resp.~${}_a\vf):\wt E\to M$ \sf ying
\er{1.17}, \er{1.18} is a direct con\sq\ of Theorem \rfa1{t1.6} with $(E,e,S)
:=\ees$, $F:=M$, $a:=a$, and $f(x):=a\qu x$ (resp.\ $x\qu a$), $x\in M$. For
the proof of \er{1.19}--\er{1.22} we introduce the maps
$\Phi_a$ (resp.\ ${}_a\Phi):\wt E\to M^M$ \rc vely defined by
\bga1.30
\Phi_a(\wt e):=\id_M \qh{(resp.\ ${}_a\Phi(\wt e):=\id_M)$}, \\
\Phi_a(\wt S(n)):=\th_a\circ\Phi_a(n) \qh{(resp.\
${}_a\Phi(\wt S(n)):={}_a\Phi(n)\circ{}_a\th)$ \fe }n\in\wt E, \lb{1.31}
\e
where $\th_a$ (resp.\ ${}_a\th):M\to M$ are defined by
\beq1.32
\th_a(x):=a\qu x \qh{(resp.\ ${}_a\th(x):=x\qu a$),} x\in M.
\e
\E\ex\ and \uq\ of such maps follows from \E\Pr\ \rfa1{p2.6} with $(E,e,S):=
\ees$, $F:=M$, and $f:=\th_a$ (resp.\ ${}_a\th$).

We claim that the \fw\ holds \fa $n\in\wt E$ and $x\in M$:
\bea1.33
\Phi_a(n)(x)&=\vf_a(n)\qu x,\\
{}_a\Phi(n)(x)&=x\qu {}_a\vf(n). \lb{1.34}
\e
We prove \er{1.33} by induction on $n\in\wt E$.

Set $A:=\{n\in\wt E: \hbox{\er{1.33} holds}\}$. We have $\wt e\in A$, since
$\Phi_a(\wt e)(x)\nde1.30 = \id_M(x) = x\nde1.17 = \vf_a(\wt e)\qu x$,
$x\in M$. Suppose $n\in A$ and $x\in M$. Then $\Phi_a(\wt S(n))(x)\nde1.31 =
(\th_a\circ\Phi_a(n))(x) = \th_a(\Phi_a(n)(x)) \nad{n\in A}= \th_a(\vf_a(n)
\qu x) \nde1.32 = a\qu(\vf_a(n)\qu x) \nde1.5 = (a\qu \vf_a(n))\qu x
\nde1.18 = \break\vf_a(\wt S(n))\qu x$. Hence $S(n)\in A$. Then $A=\wt E$ follows
from (1.N2). The proof of \er{1.34} is similar (and is left as an exercise).

\er{1.19}: We first show $\Phi_e(n)=\id_M$ \fa $n\in\wt E$. If $a:=e$ in
\er{1.31}, then $\th_a=\id_M$. Then $\id_M$ \sf ies \er{1.30}, \er{1.31}
\fa $n\in\wt E$. In view of the \uq\ part of \E\Pr\ \rfa1{p2.6}, we have
$\Phi_e(n)=\id_M$, $n\in\wt E$. Then $\vf_e(n)=\vf_e(n)\qu e\nde1.33 =
\Phi_e(n)(e)=\id_M(e)=e$, $n\in\wt E$.

\er{1.20}: From \er{1.33} with $x:=e$ we obtain $\vf_a(n)=\Phi_a(n)e$,
$n\in\wt E$. Hence $\vf_a(m)\qu\vf_a(n) \nde1.33 = \Phi_a(m)(\vf_a(n))
=\Phi_a(m)(\Phi_a(n)e) = (\Phi_a(m)\circ\Phi_a(n))(e)$, $m,n\in\wt E$.

\er{1.21} follows from \er{1.20} and \era1{2.26}.

\er{1.22}: By induction on $n\in\wt E$. Set $A:=\{n\in\wt E: \hbox{\er{1.22}
holds}\}$. We have $\wt e\in A$ since ${}_a\vf(\wt e)\nde1.17 = e\nde1.17 =
\vf_a(\wt e)$. Suppose $n\in A$. Then ${}_a\vf(\wt S(n)) \nde1.18 = {}_a\vf(n)
\qu a \nad{n\in A}= {\vf_a(n)\qu a} \nad{\er{1.17},\er{1.18}}= \vf_a(n)\qu
\vf_a(\wt S(\wt e)) \nde1.21 = \vf_a(\wt S(\wt e))\qu \vf_a(n) =a\qu\vf_a(n)
\nde1.18 = \vf_a(\wt S(n))$.
Hence $\wt S(n)\in A$. \If from (1.N2) that $A=\wt E$.

\er{1.23}: We first prove the case $m\in\wt E$, $n:=\wt S(\wt e)$. Note
that $\vf_b(\wt S(\wt e))\nad{\er{1.17},\er{1.18}}=\break b\qu \vf_b(\wt e)\nde1.17
= b\qu e=b$. Set $A:=\{m\in\wt E:\vf_a(m)\qu b=b\qu\vf_a(m)\}$. We have $\wt
e\in A$ since $\vf_a(\wt e)\nde1.17 = e$. Suppose $m\in A$. Then $\vf_a(\wt
S(m))\qu b\nde1.18 = (a\qu\vf_a(m))\qu b\nad\ast=\break a\qu(\vf_a(m)\qu
b)\nad{m\in A}= a\qu(b\qu\vf_a(m)) \nad\ast= (a\qu b)\qu\vf_a(m) = (b\qu
a)\qu \vf_a(m)\nad\ast= b\qu(a\qu\vf_a(m)) \nde1.18 = b\qu\vf_a(\wt S(m))$.
In~$\nad\ast=$ we used the \asc ity of~$\qu$. Hence $\wt S(m)\in A$, and
$A=\wt E$ follows from (1.N2). We now prove \er{1.23}. Let $m\in\wt E$ and
set $B:=\{n\in\wt E:\vf_a(m)\qu\vf_b(n)=\vf_b(n)\qu\vf_a(m)\}$. We have $\wt
e\in B$ since $\vf_b(\wt e)=e$. Suppose $n\in B$. Then $\vf_a(m)\qu \vf_b(\wt
S(n)) = \vf_a(m)\qu(b\qu\vf_b(n)) = (\vf_a(m)\qu b)\qu\vf_b(n) \nad{A=\td E}=
(b\qu\vf_a(m))\qu\vf_b(n) = b\qu(\vf_a(m)\qu\vf_b(n)) \nad{m\in B}= b\qu
(\vf_b(n)\qu\vf_a(m)) = (b\qu\vf_b(n))\qu \vf_a(m) = \vf_b(\wt S(n))\qu
\vf_a(m)$. Hence $\wt S(n)\in B$, and $B=\wt E$ by (1.N2).

In the proof of \er{1.24} we shall use the \fw\ observations. Let $X$~be
a~\ns\ and let $\qu$~be an \asc e binary \op\ on~$X$, i.e.\ $(X,\qu)$ is
a~\sg\ (see \E\df s \rfa1{d2.2}). Let $\a,\b,\g\in X$, then $(\a\qu\b)\qu\g
=\a\qu(\b\qu\g)$. One usually ``omits parentheses'', i.e.
\beq1.35
\a\qu\b\qu\g: =\a\qu(\b\qu\g)=(\a\qu\b)\qu\g.
\e
Let $x,y,u,v\in X$, then we have
\beq1.36
(x\qu y)\qu(u\qu v)=x\qu(y\qu u)\qu v.
\e
Indeed, $(x\qu y)\qu(u\qu v)=x\qu (y\qu(u\qu v))=x\qu((y\qu u)\qu v)\nde1.35 =
x\qu(y\qu u)\qu v$.

\er{1.24}: Set $A:=\{n\in\wt E: \hbox{\er{1.24} holds}\}$. We have $\wt e\in
A$ since $\vf_{a\qum b}(\wt e) \nde1.17 = e =\break e\qu e = \vf_a(\wt e)\qu
\vf_b(\wt e)$. We suppose $\wt n\in A$. Then $\vf_a(\wt S(n))\qu
\vf_b(\wt S(n))\nde1.18 = (a\qu\vf_a(n))\qu(b\qu\vf_b(n)) \nde1.36 = a\qu
(\vf_a(n)\qu b)\qu\vf_b(n) = a\qu(\vf_a(n)\qu\vf_b(\wt S(\wt e)))\qu\vf_b(n)
\nde1.23 = a\qu(\vf_b(\wt S(\wt e))\qu\vf_a(n))\qu\vf_b(n) = a\qu
(b\qu\vf_a(n))\qu\vf_b(n) \nde1.36 = (a\qu b)\qu(\vf_a(n)\qu\vf_b(n))
\nad{n\in A}= (a\qu b)\qu\vf_{a\qum b}(n)\nde1.18 = \vf_{a\qum b}(\wt S(n))$.

Hence $\wt S(n)\in A$, and $A=\wt E$ by (1.N2).
\endproof

\proof[Proof of \E\Pr\ \rf{p1.15}]
We first apply \E\Pr\ \rfa1{p2.8} with $(E,e,S)$ replaced by $\ees$, $F:=M$
and $f:=\th_a$, where $\th_a$~is defined in \er{1.32}. Then we have
\bga1.37
\Phi_{E'}(\vf(m)) \nda1{2.45} = \Phi_{\td E}(m) \qh{\fa}m\in \wt E, \\
I_{E'}(\th_a)\nda1{2.46} = I_{\td E}(\th_a). \lb{1.38}
\e
We observe that $\Phi_{\td E}:\wt E\to M^M$ is equal to the map $\Phi_a$
defined by \er{1.30}, \er{1.31} in the proof of \E\Pr\ \rf{p1.13}. Similarly
$\Phi_{E'}:E' \to M^M$ is equal to the map $\Phi'_a$ defined by \er{1.30},
\er{1.31} where $\ees$ is replaced by $(E',e',S')$. We now use \er{1.33}
with $x:=e$ and find that $\vf_a(m)=\Phi_a(m)e$.
Similarly $\vf_a'(\vf(m)) = \Phi_a'(\vf(m))e$. Using \er{1.37} we obtain
\er{1.27}, and from \er{1.38}, \er{1.33} with $x:=e$, we arrive at \er{1.28}.
\endproof

\blm1.17
Let $\ees$ be a set of \nn s, and let $(\wt E,+_{\td E},\wt e)$ be the\break P-monoid
introduced in \E\Pr\ \rf{p1.2}. Then, we have
\beq1.39
\wt S(x)=x+_{\td E}\wt S(\wt e)=\wt S(\wt e)+_{\td E}x, \q x\in\wt E,
\e
and
\beq1.40
\wt E=I(\wt S(\wt e)).
\e
\elm

\proof \

\er{1.39}: $\wt S(x)\nde1.7 = \wt S(x+_{\td E}\wt e) \nde1.3 = x+_{\td E}
\wt S(\wt e)\nde1.6 = \wt S(\wt e)+_{\td E}x$, $x\in \wt E$.

\er{1.40}: Let $\vf_{\td S(\td e)}:\wt E\to \wt E$ be the map introduced in
\E\Pr\ \rf{p1.13} with $a:=\wt S(\wt e)$. We claim that $\vf
_{\td S(\td e)}=\id_{\td E}$. Observe that $\id_{\td E}$ \sf ies \er{1.17}
since $\id_{\td E}(\wt e)=\wt e$, and \er{1.18} since $\id_{\td E}(\wt S(n))
=\wt S(n)\nde1.39 = \wt S(\wt e)+_{\td E}n=\wt S(\wt e)+_{\td E}\id_{\td E}
(n)$, $n\in\wt E$.

Then the claim follows from the \uq\ part in \E\Pr\ \rf{p1.13}. \E\Tf
\er{1.40} holds since the range of $\id_{\td E}$ is~$\wt E$.
\endproof

\bdf1.18
A \ti{nontrivial\/} monoid $(M,\qu,e)$, i.e.\ $X\ne \{e\}$, is called
\ti{\pn\/}\index{monoid!principal} (not standard terminology), if \te s $u\in M\sms e$ \st
\beq1.41
M=I(u).
\e
An \el\ $u\in M\sms e$ \st \er{1.41} holds is called a \ti{\Gn}\index{generator} of the \pn\
monoid $(M,\qu,e)$.
\edf

For convenience we recall that if $\ees$ is a~set of \nn s (see \E\df\
\rfa1{d1.3}) and if $+_{\td E}$ is the binary \op\ on~$\wt E$ defined by
\er{1.1}, then $(\wt E,+_{\td E},\wt e)$ is a~\pn\ P-monoid with \Gn\ $\wt
S(\wt e)$ (see \E\Pr\ \rf{p1.2}, Lemma \rf{l1.17} and \E\df\ \rf{d1.18}).
Our next goal is to show that if $(M,\qu,e)$ is a~\pn\ P-monoid with \Gn~$a$,
then $(\wt E,+_{\td E},\wt e)$ and $(M,\qu,e)$ are monoid-\is c and that the
map $\vf_a:\wt E\to M$ defined in \E\Pr\ \rf{1.13} is the \ti{only} \is sm
from $\wt E$ onto~$M$, and $a$~is the only \Gn\ of $(M,\qu,e)$. To this end
we need some \Pr\ and lemmata.

In \E\Pr\ \rf{p1.2} we defined the \ad\ $+_{\td E}$ on $\ees$ by stating that
$\Phi(m)\circ \Phi(n)=\Phi(m+_{\td E}n)$ where $\Phi(m),\Phi(n)$ are \IT s of
the selfmap~$\wt S$ of~$\wt E$. It turns out that this \rl\ holds for the \IT
s of \ti{all\/} selfmaps of a \ns~$F$.

\bpr1.19
Let $F$, $f$ and $(E,e,S)$ be as in \E\Pr\ \rfa1{p2.6}. Let $\Phi(m)$
$($resp.\ $\Phi(n))$ be the $m$-fold $($resp.\ $n$-fold\/$)$ \IT\ of~$f$ \wrt
$E$. Then
\beq1.42
\Phi(m)\circ\Phi(n)= \Phi(m+_En), \q m,n\in E,
\e
where $+_E$ is defined in \er{1.1} with $\ees:=(E,e,S)$.
\epr

\proof
By induction on $n\in E$. Let $m\in E$ be fixed. Set $A:=\{n\in E: \er{1.42}
\hbox{ holds}\}$. We have $e\in A$ since $\Phi(e)=\id_F$ and $\Phi(m+_E e)=
\Phi(m)=\Phi(m)\circ\id_F = \Phi(m)\circ\Phi(e)$. We suppose $n\in A$. Then
$\Phi(m)\circ\Phi(S(n))\nad\ast= \Phi(m)\circ(f\circ\Phi(n)) = (\Phi(m)
\circ f)\circ\Phi(n) = (\Phi(m)\circ\Phi(S(e)))\circ\Phi(n) \nad{\ast\ast}=
(\Phi(S(e))\circ\Phi(m))\circ\Phi(n) = \Phi(S(e))\circ(\Phi(m)\circ\Phi(n))
\nad{n\in A}= f\circ(\Phi(m+_En))\nad\ast=\Phi(S(m+_En))\nde1.3 = \Phi(m
+_ES(n))$, where we used \era1{2.25} in~$\nad\ast=$, \era1{2.18} in $\nad{
\ast\ast}=$. Hence $S(n)\in A$. \If from (1.N2) that $A=E$.
\endproof

\blm1.20
Let $\vf_a:\wt E\to M$ be the map introduced in \E\Pr\ \rf{p1.13}. Then
\beq1.43
\vf_a(m)\qu \vf_a(n) = \vf_a(m+_{\td E}n), \q m,n\in \wt E.
\e
\elm

\proof
Using \er{1.33} with $x:=\wt e$, we obtain
\beq1.44
\vf_a(n)=\Phi_a(n)(e) \qh{\fa}n\in\wt E.
\e
Now let $m,n\in\wt E$. Then $\vf_a(m)\qu\vf_a(n) \nde1.20 = (\Phi_a(m)\qu
\Phi_a(n))(e) \nde1.42 = \Phi_a(m+_{\td E}n)(e) \nde1.44 =
\vf_a(m+_{\td E}n)$.
\endproof

\blm1.21
Let $(M,\qu,e)$ be a nontrivial monoid and let $M_0$ be a nontrivial \sbm\
of~$M$, possibly  equal to~$M$. Let $a\in M_0\sms e$, let $\vf_a:\wt E
\to M$ $($resp.\ $\psi_a:\wt E\to M)$ denote the maps of \IT s of~$a$ in the
monoid $(M_0,\qu,e)$ $($resp.\ $M,\qu,e))$ introduced in
\E\Pr\ \rf{p1.13}. Then

\hph i,ii, $\vf_a(n)=\psi_a(n)$ \fa $n\in\wt E$ where the \et y take place in
$(M,\qu,e)$. \E\Ip $I(a)$ the set of \IT s of~$a$ in $(M,\qu,e)$ \sf ies
$I(a)\sbs M_0$.

\hph ii,i, $I(a)$ is a \sbm\ of $(M_0,\qu,e)$, and $(I(a),\qu,e)$ is an
abelian monoid. \E\Ip a~\pn\ monoid is abelian.

\hph iii,, The map $\vf_a:\wt E\to(I(a),\qu,e)$ is a sur\jc\ \hm sm.
\elm

\proof \

(i) Since $M_0$ is a \sbm\ of the monoid $(M,\qu,e)$, the \rt ion of~$\qu$ to
$M_0\t M_0$ is a binary \op\ on~$M_0$ which we denote (only in this proof)
by~$\tqu$. Hence $(M_0,\tqu,e)$ is a monoid. Let $\vf_a$ denote the \sq\ of
\IT s of~$a$ in $(M_0,\tqu,e)$ and let $\psi_a$ denote the \sq\ of \IT s in
$(M,\qu,e)$. We claim that $\vf_a(n)=\psi_a(n)$ in~$M$ \fa $n\in\wt E$. Set
$A:= \{k\in\wt E: \vf_a(k)=\psi_a(k)\}$. We use \In\ to show that $A=\wt E$.

$\wt e\in A$: since $\vf_a(\wt e)=e=\psi_a(\wt e)$.

\ti{$k\in A$ implies $\wt S(k)\in A$}: Suppose $k\in A$. Then $\psi_a(\wt
S(k))\nde1.18 = a\qu \psi_a(k) \nad{k\in A} = a\qu \vf_a(k)$. Since $a$
and~$\vf_a(k)$ belong to~$M_0$, $a\qu \vf_a(k)= a\tqu \vf_a(k) \nde1.18 =
\vf_a(\wt S(k))$. Hence $\wt S(k)\in A$. \If that $A$~is \iv\ in~$\wt
E$, hence $A=\wt E$. Clearly  $I(a)=\bcl_{n\in \wt E}\{\psi_a(n)\} =
\bcl_{n\in\wt E}\{\vf_a(n)\}\sbs M_0$.

(ii) We have $e\in I(a)$ since $\vf_a(\wt e)=e$. Let $x,y\in I(a)$. \E\te\
$m,n\in\wt E$ \st $x=\vf_a(m)$, $y=\vf_a(n)$. Hence $x\qu y=\vf_a(m)\qu
\vf_a(n) \nde1.43 = \vf_a(m+_{\td E}n)\in I(a)$. \Mo $x\qu y=\vf_a(m)\qu
\vf_a(n) \nde1.21 = \vf_a(n)\qu \vf_a(m) = y\qu x$. \If that $(I(a),\qu,e)$
is an abelian monoid.

(iii) $\vf_a$ is sur\jc, since the range of $\vf_a$ is equal to $I(a)$ by the
\df\ of~$I(a)$. \Mo $\vf_a(\wt e)=e$ and $\vf_a(m+_{\td E}n)\nde1.43 =
\vf_a(m)\qu\vf_a(n)$, $m,n\in\wt E$. Hence $\vf_a$ is a \hm sm.
\endproof

In the  next lemma we consider images of ``\IT s'' under a \hm sm.

\blm1.22
Let $(E,e,S)$ be a set of \nn s. Let $(M_i,\qu_i,e_i)$, $i=1,2$, be monoids
$($not necessarily abelian$)$ and let $\vf:M_1\to M_2$ be a \hm sm. Let\break
$a\in M $, $n\in E$, then the \fw\ holds.
\beq1.45
\vf(n \dqui 1E a) = n\dqui 2E \vf(a).
\e
\elm

\proof
We use \In\ on $n\in E$. Set $A:=\{n\in E: \er{1.45}\hbox{ holds}\}$. Then
$e\in A$, since $e\dqui1E a=e_1$, $e\dqui 2E\vf(a)=e_2$ in view of \er{1.17},
and $\vf(e_1)=e_2$ by \er{1.15}. We now suppose that $n\in A$. We
have $\vf(S(n)\dqui1Ea) \nde1.18 = \vf(a\qu_1(n\dqui1E a)) \nde1.16 = {\vf(a)
\qu_2 \vf(n\dqui1E a)} \nad{n\in A}= \vf(a)\qu_2(n\dqui2E \vf(a)) \nde1.18 =
S(n)\dqui2E \vf(a)$. Hence $S(n)\in A$, and $A=E$ by (1.N2).
\endproof

In what follows it will be useful to introduce the \nog\ of the set of \nn s
$\ees$ (see \E\df\ \rfa1{d3.24}). In the next lemma we \es\ the \rl\ between
this \og, which we denote by~$\le$, and the monoid $(\wt E,+_{\td E},\wt e)$
introduced in \E\Pr\ \rf{p1.2}.

\blm1.23
Let $\ees$ and $(\wt E,+_{\td E},\wt e)$ be as in Lemma \rf{l1.17}, and let
$\le$~denote the \nog\ of $\ees$. Let $x,y\in\wt E$. Then
\bml1.46
x\le y\ (\hbox{resp.\ }x<y)\hbox{ iff \te s }p\in\wt E \ (\hbox{resp.\ }
p\in\wt E\sms{\wt e}) \\ \hbox{ \st }y=x+_{\td E} p.
\e
\elm

\proof
We first prove
\beq1.47
x+_{\td E}p = \Phi(p)x \qh{\fa} x,p\in\wt E.
\e
Let $x,p\in E$. Then $x+_{\td E}p\nda1{2.35} = \Phi(x+_{\td E}p)\wt e
\nde1.1 = (\Phi(x)\circ\Phi(p))\wt e \nda1{2.26} = (\Phi(p)
\circ\Phi(x))\wt e = \Phi(p)(\Phi(x)\wt e) \nda1{2.35} =
\Phi(p)x$.

Let $x,y\in\wt E$. Then $x\le y$ iff \te s $p\in\wt E$ \st $y=\Phi(p)x$, by
\era1{3.6n}. Then \er{1.46} with $\le$ follows from \er{1.47}. Now if $x<y$,
then $x\le y$ and $x\ne y$ by \era1{3.7}. Then \te s $p\in\wt E$ \st $y=x
+_{\td E}p$, and $p\ne \wt e$ for otherwise $y=x+_{\td E}\wt e\nde1.7 = x$.
Conversely, if $y=x+_{\td E}p$ with $p\in\wt E\sms{\wt e}$, then $y\ne x$,
for otherwise $x=x+_{\td E}p$, hence $x+_{\td E}
\wt e\nde1.7 = x=x+_{\td E}p$, which implies $p=\wt e$ by \er{1.8}, \er{1.6}, a~\cd ion.
\endproof

\brm1.24
If $x,y\in\wt E$ and $x\le y$, then \te s \ti{\ooo}$p\in\wt E$ \st $y+_{\td E} p$.
Indeed, from $x+_{\td E} p=x+_{\td E} q$, $p,q\in\wt E$ we infer $p=q$ by \er{1.8}, \er{1.6}.
Set $\D_{\td E}:=\{(x,y)\in \wt E\t\wt E: x\le y\}$. Then we can define a~map~$\psi$
from $\D_{\td E}$ into $\wt E$ by setting $\psi(x,y):=p\in\wt E$ \st $y=x+_{\td E} p$.
The map $\psi$ is called \ti{\sbt} and will be denoted by
$-_{\td E}$. We thus have
\beq1.48
p=y-_{\td E} x \qh{iff $x\le y$ and} y=x+_{\td E} p,\ x,y,p\in\wt E.
\e
\erm

\bth1.25
Let $(M,\qu,e)$ be a \pn\ P-monoid and let $a$~be a~\Gn\ of ${(M,{\qu},e)}$.
Then\dw

\hph i,ii, The map $\vf_a:\wt E\to M$ introduced in \E\Pr\ \rf{p1.13} is an
\ti{\is sm}.

\hph ii,i, The \el\ $a\in M\sms e$ is the only \Gn\ of $(M,\qu,e)$.

\hph iii,, The map $\vf_a$ is the only \is sm from $(\wt E,+_{\td E},\wt E)$
onto $(M,\qu,e)$.
\eth

\proof \

(i) In view of Lemma \rf{l1.21}\,(iii) it is \sft\ to show that $\vf_a$ is
\ti{in\jc}. Suppose for \cd ion that \te\ $m,n\in \wt E$, $m\ne n$, \st
$\vf_a(m)=\vf_a (n)$. Let $\le$ denote the \nog\ of $\ees$. Then, by \er{1.46}
\te\ $p,q\in\wt E\sms{\wt e}$ \st $m=n+_{\td E}p$ if $n<m$, and $n=m+_{\wt E}q$
if $m<n$. Recall that the \nog\ is total. We obtain in the first case
$\vf_a(n)=\vf_a(m)=\vf_a(n+_{\td E}p)\nde1.43 = \vf_a(n)\qu \vf_a(p)$. Since
$\vf_a(n)=\vf_a(n)\qu e$, we obtain $\vf_a(p)=e$ from \er{1.6}, \er{1.8}, since $M$~is
a~\Cm. Similarly in the second case we obtain $\vf_a(m)=\vf_a(n)=\vf_a(m
+_{\td E}q)=\vf_a(m)\qu \vf_a(q)$, hence $\vf_a(q)=e$. We find a \cd ion by
showing that if $r\ne\wt e$, then $\vf_a(r)\ne e$, $r\in \wt E$. Indeed, if $r\ne\wt e$, we
have $\wt e<r$ by \era1{3.13}, and $\wt S(\wt e)\le r$ by \era1{3.17}. \E\Tf
by \er{1.46} \te s $s\in \wt E$ \st $r=\wt S(\wt e)+_{\td E}s$. \If from \er{1.43}
that $\vf_a(r)=\vf_a(\wt S(\wt e)) \qu \vf_a(s)$. Note that $\vf_a(\wt S
(\wt e))=a$ by \er{1.17}, \er{1.18}. Thus $\vf_a(r)=a\qu \vf_a(s)$. If
$\vf_a(r)=e$, then $a=e$ and $\vf_a(s)=e$ by \er{1.9} since $(M,\qu,e)$ is
a~P-monoid. But this is impossible since $a\ne e$. Hence $\vf_a(r)\ne e$,
a~\cd ion.

(ii) Let $b\in M\sms e$ be \st $M=I(b)$. Since $b\in I(a)\sms e$ and $\vf_a
(\wt e)=e$, \te s $m\in\wt E\sms{\wt e}$ \st $b=\vf_a(m)$. Since $m\in\wt E
\sms{\wt e}$, \te s $p\in\wt E$ \st $m=\wt S(p)$ by (1.N1). \E\Tf $b=\vf_a(\wt S(p))
\nde1.18 = a\qu \vf_a(p)$. Similarly $a\in I(b)$ and $a\ne e$, hence \te s
$q\in\wt E$ \st $a=\vf_b(\wt S(q))=b\qu \vf_b(q)$. \If that $b = a\qu\vf_a(p) =
(b\qu \vf_b(q))\qu \vf_a(p) = b\qu(\vf_b(q)\qu\vf_a(p))$, hence $b\qu e
= b = b\qu(\vf_b(q)\qu\vf_a(p))$. Since $M$~is a \Cm, we infer that $\vf_b(q)
\qu \vf_a(p)=e$. Since $M$~is a P-monoid, we have $\vf_b(q)=\vf_a(p)=e$.
Hence $b=a\qu e=a$.

(iii) We first show that the identity is the \ti{only auto\mf} of $(\wt E,
+_{\td E},\wt e)$. We claim that $\id_{\td E}=\vf_{\td S(\td e)}$. Indeed, we
have $\id_{\td E}(\wt e)=\wt e$ and $\id_{\td E}(\wt S(n))= \wt S(n)\nde1.7 =\wt S
(n+_{\td E}\wt e) \nde1.3 = n+_{\td E}\wt S(\wt e) \nde1.6 = \wt S(\wt e)
+_{\td E}n = \wt S(\wt e)+_{\td E}\id_{\td E}(n)$ \fa $n\in \wt E$. Hence
$\id_{\td M}$ \sf ies \er{1.17}, \er{1.18} with $a:=\wt S(\wt e)$. \E\Tf we
obtain
\beq1.49
\id_{\td E}=\vf_{\td S(\td e)}.
\e
We now suppose that $\vf$ is an auto\mf\ (see \E\df\ \rf{d1.7}). We claim that
$\vf=\vf_b$ with $b:=\vf(\wt S(\wt e))$. Indeed, $\vf(\wt e)=\wt e$ and
$\vf(\wt S(n))\nde1.39  = \vf(\wt S(\wt e)+_{\td E}n) = \vf(\wt S(\wt e))
+_{\td E}\vf(n)$ \fa $n\in\wt E$. Hence $\vf$~\sf ies \er{1.17}, \er{1.18}
with $a:=\vf(\wt S(\wt e))$, \Tf $\vf=\vf_b$. Since the range of~$\vf$ is
equal to~$\wt E$, we have $\wt E=I(b)$. In view of \E\Pr\ \rf{p1.2} and
\er{1.40} $(\wt E,+_{\td E},\wt e)$ is a \pn\ P-monoid with \Gn\ $\wt S(\wt
e)$. Hence by~(ii), we obtain $b=\wt S(\wt e)$. \E\Tf $\vf=\vf_b=\vf_{\td S
(\td e)}\nde1.49 = \id_{\td E}$, and (iii) holds when $(M,\qu,e) := (\wt E,
+_{\td E},\wt e)$.

Finally, let $\vf$ be an \is sm from $(\wt E,+_{\td E},\wt e)$ onto $(M,\qu,e)$.
Then $\vf\Inv:M\to E$ is also an \is sm as well as $\psi:=\vf\Inv\qu\vf_a$,
in view of Lemma \rf{l1.8}\,(iii) and (ii). Since $\psi$~is an auto\mf\ of
$(\wt E,+_{\td E},\wt e)$, it is equal to~$\id_{\td E}$. \E\Tf $\vf =\vf
\circ\id_{\td E} = \vf\circ(\vf\Inv\circ\vf_a) = (\vf\circ\vf\Inv) \circ
\vf_a = \id_M\circ\vf_a=\vf_a$, and (iii) follows.
\endproof

\bpr1.26
Let $(M,\qu,e)$ and $(M',\qu',e')$ be \pn\ P-monoids $($not necessarily
distinct\/$)$. Let $a$ $($resp.~$a')$ denote the \Gn\ of~$M$ $($resp.~$M')$. Let
$\vf_a:\wt E\to M$ denote the map introduced in \E\Pr\ \rf{p1.13} and let
$\vf'_{a'}:\wt E \to M'$ denote the \crs\ map with $M$~replaced by~$M'$. Then
$\vf:=\vf'_{a'}\circ\vf_a\Inv$ is the \emph{only \is sm} from $M$ onto~$M'$,
and $\vf(a)=a'$. \E\Ip the identity is the \emph{only auto\mf} of a~\pn\
P-monoid.
\epr

\proof
We first prove the last assertion. Let $(M,\qu,e)$ be a~\pn\ P-monoid with
\Gn~$a$. In view of Theorem \rf{t1.25}\,(i) the map $\vf_a:\wt E\to M$ is an
\is sm. Let $\phi$~be an auto\mf\ of~$M$. Then the map $\vf_a\Inv\circ \phi
\circ\vf_a$ is an auto\mf\ of~$\wt E$ by Lemma \rf{l1.8}\,(ii). \If from the
proof of Theorem \rf{t1.25}\,(iii) that $\vf_a\Inv\circ \phi\circ\vf_a =
\id_{\td E}$. \E\Tf $\phi\circ\vf_a = (\vf_a\circ\vf_a\Inv)\circ\phi\circ
\vf_a = \vf_a\circ(\vf_a\Inv \circ\phi\circ\vf_a) = \vf_a\circ\id_{\td E}
=\vf_a$. Hence $\phi=\phi\circ(\vf_a\circ\vf_a\Inv) = (\phi\circ\vf_a)\circ
\vf_a\Inv = \vf_a\circ\vf_a\Inv = \id_{\td E}$. Thus $\phi=\id_{\td E}$.

Clearly $\vf:=\vf_a'\circ\vf_a\Inv$ is an \is sm from $M$ onto~$M'$, by Lemma
\rf{l1.8}\,(ii). Suppose that $\wh\vf$ is an \is sm from $M$ onto~$M'$. Then
$\wh\vf{}\Inv\circ\vf$ is an auto\mf\ of~$M$, hence $\wh\vf{}\Inv\circ\vf
=\id_M$. Then $\vf = (\wh\vf \circ \wh\vf{}\Inv)\circ \vf = \wh\vf\circ
(\wh\vf{}\Inv\circ\vf) = \wh\vf\circ\id_{M'}=\wh\vf$. Thus $\vf=\wh\vf$.
\Mo $\vf(a) = (\vf_{a'}'\circ\vf_a\Inv)(a) = \vf'_{a'}(\vf_a\Inv(a)) =
\vf'_{a'}(\wt S(\wt e))=a'$.
\endproof

\bpr1.27
Let $(M,\qu,e)$ be a monoid. Then $M$ is \is c to $(\wt E,+_{\td E},\wt e)$
introduced in \E\Pr\ \rf{p1.2} iff $M$ is a \pn\ P-monoid.
\epr

\bex1.28
Prove \E\Pr\ \rf{p1.27}.
\eex

\bnt1.29
In what follows we shall denote by $(M,\qu,e;u)$ a~\pn\ P-monoid $(M,\qu,e)$
with \Gn~$u$ which is unique by Theorem \rf{t1.25}\,(ii).
\ent

We recall that in this section we started with a set of \nn s $\ees$ and we
introduced in \E\Pr\ \rf{p1.2} a binary \op\ on~$\wt E$ denoted by $+_{\td E}$
making~$\wt E$ a~\pn\ P-monoid with \nel~$\wt e$ and \Gn\ $\wt S(\wt e)$,
denoted by $(\wt E,+_{\td E},\wt e;\wt S(\wt e))$. Our next goal is to show
that given $(M,\qu,e;u)$, a~\pn\ P-monoid with \Gn~$u$, \te s a~selfmap $S$
of~$M$ which makes $(M,e,S)$ a~set of \nn s \st $(M,+_M,e;u) = (M,\qu,e;u)$.

\bpr1.30
Let $(M,\qu,e;u)$ be a \pn\ P-monoid with \Gn~$u$. Define a~map $S:M\to M$ by
setting
\beq1.50
S(x):= u\qu x, \q x\in M.
\e
Then $(M,e,S)$ is a set of \nn s. Let $\vf_u:\wt E\to M$ denote the map
introduced in \E\Pr\ \rf{p1.13} with $a:=u$. Let $\vf:\wt E\to M$ be the
unique bi\jn\ introduced in Theorem \rfa1{t1.4}, \sf ying
\beq1.51
\vf(\wt e)=e \qh{and }\vf\circ\wt S = S\circ \vf.
\e
Then
\beq1.52
\vf=\vf_u.
\e

Let $+_M$ denote the unique binary \op\ on~$M$ \sf ying \er{1.2}, \er{1.3}
with $\wt E:=M$, $\wt e:=e$ and $\wt S:=S$. Then this \op~$+_M$ is identical
to the \op~$\qu$, and the \Gn~$u$ \sf ies
\beq1.53
u = \vf(\wt S(\wt e)).
\e
\epr

\proof
Let $\vf_u:\wt E\to M$ denote the map uniquely defined by \er{1.17} and
\er{1.18} with $a:=u\in M\sms e$. Let $n\in\wt E$. We have $(\vf_u\circ
\wt S)(n) = \vf_u(\wt S(n)) \nde1.18 = u\qu \vf_u(n) = S(\vf_u(n)) = (S\circ
\vf_u)(n)$. Since $n$~is arbitrary in~$\wt E$, we obtain
\beq1.54
\vf_u\circ \wt S = S\circ \vf_u.
\e
In view of Theorem \rf{t1.25}\,(i), $\vf_u$ is bi\jc, hence $S=\vf_u\circ
\wt S\circ\vf_u\Inv$. The in\ji\ of~$S$ is a direct con\sq\ of \er{1.6}, \er{1.8} since
$M$~is a~\Cm. We now show that the range of~$S$, $R(S)=M\sms e$. Since $M$~is
a~P-monoid, there is no $x\in M$ \st $u\qu x=e$. Indeed, if $u\qu x=e$
\fs $x\in M$, then $u=x=e$ by \er{1.9}, \cd ing $u\ne e$.
Let $y\in M\sms e$. Since $M=I(u)$,
\te s $n\in\wt E$ \st $\vf_u(n)=y$. We have $n\ne\wt e$, since $\vf_u(\wt e)
=e$ and $y\ne e$. \E\Tf $n\in R(\wt S)$ by (1.N1). Let $p\in\wt E$ be \st
$n=\wt S(p)$. Then $\vf_u(\wt S(p))=y$, i.e.\ $(\vf_u\circ\wt S)(p)=y$. From
\er{1.54} we obtain $(S\circ\vf_u)(p)=y$, hence $y=S(\vf_u(p))\in R(S)$. This
shows that $(M,e,S)$ \sf ies (1.N1).

We now prove (1.N2) for $(M,e,S)$. Let $A\sbs M$ be \iv, i.e.\ $e\in A$
and ${S(A)\sbs A}$. We have to show that $A=M$. It is \sft\ to show that
$B:=\vf_u\Inv(A)$ is \iv\ in $\ees$, since in this case $\vf_u\Inv(A)=
\wt E$, hence $A=
\vf_u(\wt E)$ by the bi\ji\ of~$\vf_u$ and $\vf_u(\wt E)=I(u)=M$. Clearly
$\wt e\in B$, since $e\in A$ and $\wt e=\vf_u\Inv(e)$. Now let $n\in B$, i.e.\
$\vf_u(n)\in A$. Then $\vf_u(\wt S(n))=(\vf_u\circ\wt S)(n)\nde1.54 =
(S\circ \vf_u)(n)=(S\circ
\vf_u(n))\in A$, since $\vf_u(n)\in A$ and $S(A)\sbs A$. \E\Tf $\wt S(n)\in B$,
hence $B$~is \iv\ in $\ees$. \If that $(M,e,S)$ is a~set of \nn s. In view of
Theorem \rfa1{t1.4} \te s \ooo $\vf:\wt E\to M$ \st \er{1.51} holds. \E\Tf
$\vf=\vf_u$, in view of \er{1.54} and \er{1.17}. This shows \er{1.51},
\er{1.52}. We have $\vf_u(\wt S(\wt e))=u$ by \er{1.17}, \er{1.18}, hence
\er{1.53} follows from \er{1.52}.

In view of \E\Pr\ \rf{p1.2},
it remains to show that the \op~$\qu$ in~$M$ \sf ies \er{1.2}, \er{1.3}
with $\wt E:=E$, $\wt e:=e$ and $\wt S:=S$. \E\fa $x,y\in M$ we have $x\qu e
=x$, hence \er{1.2} holds. Now $x\qu S(y) = x\qu(u\qu y)\nde1.5 = (x\qu u)
\qu y\nde1.6 = (u\qu x)\qu y \nde1.5 = u\qu(x\qu y)=\break S(x\qu y)$. Hence
\er{1.3} holds.
\endproof

\bex1.31
Let $(M,\qu,e)$ be a P-monoid. Define a \rl\ R on~$M$ by setting
\beq1.55
x\mathrel{\rm R}y \qh{if \te s $z\in M$ \st $y=x\qu z$}.
\e
Show that R is an \og\ on~$M$. Show that if $x,x',y,y'\in M$ \sf y
$x\mathrel{\rm R}y$ and $x'\mathrel{\rm R}y'$, then $(x\qu y)\mathrel{\rm R}
(x'\qu y')$.

Show that R is a well-\og\ if $(M,\qu,e)$ is a \pn\ P-monoid. Show that the
least \el\ of $M\sms e$ is the \Gn\ of $(M,\qu,e)$.
\eex

In \E\Pr\ \rf{p1.2} we showed that $(\wt E,+_{\td E},\wt e)$ is a~\PM\ and
in Lemma \rf{l1.17} that $(\wt E,+_{\td E},\wt e)$ is \pn\ with \Gn\
$\wt S(\wt e)$. In Remark \rf{r1.11} we observed that a~nontrivial \PM\ is
infinite. It turns out that an \ti{infinite} \pn\ monoid is a~\PM. An
important ingredient of the proof of this fact is what is sometimes called
the ``division algorithm'' theorem\index{division algorithm}, for reasons which will be clear later.
This theorem is in turn based on the \fw\ lemma.

\blm1.32
Let $\ees$ be a set of \nn s and let $\le$ denote its \nog. Let $c:\wt E \to
\wt E$ be a \emph{strictly in\cre} map $($see \E\df\ \rfa1{d3.13}$)$. We shall
denote by~$c_n$ the value of~$c$ at $n\in\wt E$, and by~$I_n$ the \il\
$[c_n,c_{\wt S(n)})$, $n\in\wt E$. \Mo we assume
\beq1.56
\wt c_{\td e}=\wt e.
\e
Then, \fe $x\in\wt E$, \te s $m\in\wt E$ \st $x\in I_m$. \Mo $I_n\cap I_m=\vn$
whenever $n,m\in\wt E$ with $m\ne n$.
\elm

\proof \

\ti{\E\ex}. Let $x\in \wt E$.
We have $\wt e\le x$ by \era1{3.13}. Set $B:=c(\wt E)\cap[\wt e,x]$. Since
$c_{\td e}\nde1.56 = \wt e\le x$, we also have $c_{\td e}\in[\wt e,x]$.
Clearly $c_{\td e}\in c(\wt E)$, hence $c_{\td e}\in B$. \E\Tf $B$~is
a~nonempty bounded subset of~$\wt E$ with $x$ as \ub. In view of Theorem
\rfa1{t3.39} it has a~greatest \el\ which we denote by~$y$. Since $y\in
c(\wt E)$, \te s $m\in \wt E$ \st $y=c_m$. Clearly $c_m\le x$. \E\oh since
$m<\wt S(m)$ by \era1{3.13} and $c:\wt E\to \wt E$ is strictly in\cre, we have
$c_m<c_{\td S(m)}$. Thus $c_{\td S(m)}\not\le x$, for otherwise, $c_{\td S(m)}$
would belong to~$B$, and $c_m$ would not be the greatest \el\ of~$B$. \E\Tf
we have $x<c_{\td S(m)}$. Hence $c_m\le x<c_{\td S(m)}$, i.e. $x\in I_m$.

\ti{Disjointness}. Suppose for \cd ion that \te\ $m,n\in \wt E$ with $m\ne n$,
and $y\in I_m\cap I_n$. Since $I_m\cap I_n=I_n\cap I_m$ and $\le$ is total,
we may assume that $n<m$. In view of \era1{3.17} we have $\wt S(n)\le m$,
hence by in\cre ness of~$c$, $c_{\td S(n)}\le c_m$. Since $y\in I_n$, we infer
$y<c_{\td S(n)}$, hence $y<c_{\td S(n)}\le c_m$. Using \era1{3.12} we obtain
$y<c_m$. Since $y\in I_m$, we have $c_m\le y$. Hence $c_m\le y<c_m$, and by
\era1{3.11}, $c_m<c_m$, a~\cd ion.
\endproof

\brm1.33 \

\hph i,i, Disjointness implies \uq\ of $m\in\wt E$ in the first part of the
conclusion of the lemma.

\hph ii,, The conclusion of the lemma can be rephrased by saying that the
subsets $I_n$, $n\in\wt E$, form a~\pt\ of~$\wt E$ (since $I_n\ne\vn$,
$n\in \wt E$).
\erm

\bex1.34
Let $\ees$ be a set of \nn s and let $\le$ denote its \nog. Let $f:\wt E \to
\wt E$. Show that $f$~is in\cre\ (resp.\ strictly in\cre) iff it \sf ies $f(n)
\le\hbox{(resp.~$<$)}f(\wt S(n))$ for all $n\in\wt E$.
\eex

\begin{dspc}\lb{d1.35}
A map from a set of \nn s $(E,e,S)$ into a set~$X$ is usually called
a~\ti{\sq}\index{sequence} of \el s of~$X$. The value of~$f$ at $n\in E$ is usually denoted
by~$f_n$ and it is customary to call $f_n$ a~\ti{term} of the \sq. The range
of~$f$ is the set $\{f_n\in X: n\in E\}$. Another notation for the \sq~$f$ is
$(f_n)_{n\in E}$. A~map~$f$ (resp.\ \sq\ $(f_n)_{n\in E}$) \st $f_n=f_m$ \fa
$m,n\in E$ is called a \ti{constant\/} map (resp.\ \sq).
\end{dspc}

\bnt1.36
The map $\vf_a:\wt E\to M$ introduced in \E\Pr\ \rf{p1.13} is a~\sq, the \sq\
of \IT s of $a\in M$. We shall (and we did) \ti{not\/} use the subscript
notation $(\vf_a)_n$ for $\vf_a(n)$. \Mo if $(M,\qu,e):=(\wt E,+_{\td E},
\wt e)$, we shall use the simpler notation $n\dpl a$ instead of
$\mathrel{\lower2pt\hbox{$\stackrel{+_{\tilde E}}{\cdot_{\td E}}$}}$, i.e.\
$n\dpl a:=\vf_a(n)$.
\ent

\blm1.37
Let $\ees$, $\le$ be as in Lemma \rf{l1.32}, and let $(\wt E,+_{\td E},\wt
e)$ denote the monoid introduced in \E\Pr\ \rf{p1.2}. Let
$\vf: (\wt E,+_{\td E},\wt e)\to (\wt E,+_{\td E},\wt e)$ be
a~\hm sm $($i.e.\ an endo\mf$)$. Set
\beq1.57
a:=\vf(\wt S(\wt e)).
\e
Then the \fw\ holds.

\hph i,ii, $\vf=\vf_a$ where $\vf_a$ is introduced in \E\Pr\ \rf{p1.13}.

\hph ii,i, If $a=\wt e$, then $\vf$ is constant.

\hph iii,, If $a=\wt S(\wt e)$, then $\vf=\id_{\td E}$.

\hph iv,, If $a\ge \wt S(\wt e)$, then $\vf$ is strictly in\cre.

\hph v,i, If $a>\wt S(\wt e)$, then $\vf$ is in\jc, but not sur\jc.
\elm

\proof \

(i) \er{1.17} holds since $\vf(\wt e)=\wt e$. Let $n\in\wt E$. Then $\wt S(n)
\nde1.39 = \wt S(\wt e)+_{\td E} n$, hence $\vf(\wt S(n))= \vf(\wt S(\wt e)
+_{\td E}n)=\vf(\wt S(\wt e))+_{\td E}\vf(n)$. Hence \er{1.18} holds with $a
=\vf(\wt S(\wt e))$. Hence $\vf=\vf_a$ by \E\Pr\ \rf{p1.13}. Therefore, using
Notation \rf{n1.36} we have
\beq1.58
\vf(n)=n\dpl a, \q n\in\wt E.
\e

(ii) If $a=\wt e$, then $n\dpl\wt e=\wt e$ \fa $n\in\wt E$ by \er{1.19}.

(iii) If $a=\wt S(\wt e)$, then $\vf=\id_{\td E}$ by \er{1.49}.

(iv) Let $m,n\in\wt E$ with $m<n$. In view of \er{1.46} \te s $p\in\wt E
\sms{\wt e}$ \st $n=m+_{\td E}p$. Hence by \er{1.43}, $\vf_a(n)=\vf_a(m
+_{\td E}p)=\vf_a(m)+_{\td
E}\vf_a(p)$. \If that $\vf_a(m)\le \vf_a(n)$ by \er{1.46}. It remains to show
that $\vf_a(p)\ne\wt e$. Since $p\in\wt E\sms{\wt e}$, there is $q\in\wt E$
\st $p=\wt S(q)$ by (1.N1).  Hence $\vf_a(\wt S(q))\nde1.18 = a+_{\td E}\vf_a(q)
\ge a$ by \er{1.46}. \If that $\vf_a(p)=\vf_a(\wt S(q))\ge a\ge \wt S(\wt e)$,
hence $\vf_a(p)\ge \wt S(\wt e)>\wt e$ by \era1{3.13}, hence $\vf_a(p)
>\wt e$ by \era1{3.12}. \E\Tf $\vf_a(p)\ne\wt e$, hence $\vf_a(m)<\vf_a(n)$.

(v) $\vf_a$ is in\jc\ by Lemma \rfa1{l3.34} since the natural order of
$(\wt E,\wt e,\wt S)$ is total and $\vf_a$ is strictly in\cre\ by~(iv).

\Mo $\vf_a$ is not sur\jc, otherwise $\wt E=I(a)$, i.e.\ $a$~would be a~\Gn\
of~$\wt E$. Since $\wt E$ is a~\pn\ \PM, it has only one \Gn\ $\wt S(\wt e)$
by Theorem \rf{t1.25}\,(ii).
But $a>\wt S(\wt e)$, hence $\vf_a$ is not sur\jc.
\endproof

\begin{thm}[``division algorithm'']\lb{t1.38}
Let $\ees$ be a set of \nn s, let $\le$ denote its \nog\ and let $(\wt
E,+,\wt e)$ denote the monoid introduced in \E\Pr\ \rf{p1.2}. Using
Notation \rf{n1.36}, the \fw\ holds.

Given $a\in\wt E\sms{\wt e}$, $b\in\wt E$, \te s \ooo pair $(n,r)\in\wt E\t
\wt E$ \sf ying
\bga1.59
b=n\dpl a+r,\\
r\in[\wt e,a). \lb{1.60}
\e
The \nm\ $b$ is called the \emph{dividend}, $a$~the \emph{divisor},\index{divisor} $n$~the
\emph{\qt}\index{quotient} and $r$~the \emph{remainder}. If $r=0$ in \er{1.59}, then $b$~is
also called a~\emph{\ml e}\index{multiple} of~$a$.
\eth

\proof
We use the simpler notation $n\cdot a$ instead of $n\dpl a$.

\ti{\E\ex}: Set
\beq1.61
A_n:= [n\cdot a,(n+1)\cdot a), \q n\in \wt E.
\e
In view of Lemma \rf{l1.37}\,(i), (iv) together with Notation \rf{n1.36}, the
\sq\ $(n\cdot a)_{n\in \wt E}$ is strictly in\cre, since $a\ge \wt S(\wt e)$ by
\era1{3.13}, \era1{3.17}. Observe that $\wt e\cdot a=\wt e$ by \df. Hence
by Lemma \rf{l1.32} with $I_n:=A_n$, $n\in
\wt E$, \te s $n\in\wt E$ \st $b\in A_n$, i.e.\ $n\cdot a\le b<
(n+1)\cdot a$. From the first in\et y and \er{1.46} we deduce the \ex\ of
$r\in \wt E$ \st $b=n\cdot a+r$, which is \er{1.59}. Similarly from the second
in\et y we infer the \ex\ of $s\in \wt E\sms{\wt e}$, \st $(n+1)\cdot a=b+s$. Since
$(n+1)\cdot a\nde1.18 = a+(n\cdot a)\nde1.6 = (n\cdot a)+a$, we obtain
$(n\cdot a)+a = (n\cdot a+r)+s \nde1.5 = n\cdot a+(r+s)$. Using \er{1.8} we
arrive at $r+s=a$ which implies $r<a$ by \er{1.46}, since $s\in \wt E\sms{\wt e}$.
Hence \er{1.60} holds.

\ti{\E\uq}: Using Lemma \rf{l1.32} with $I_n:=A_n$ as in the ``\E\ex'' part
of the proof, we obtain
\beq1.62
A_n\cap A_m=\vn \qh{whenever}m\ne  n.
\e
We next show that if $(n,r)\in \wt E\t \wt E$ \sf ies \er{1.59}, \er{1.60} then
$b\in A_n$. Clearly $a\cdot n\le b$ follows from \er{1.59} and \er{1.46}.
Since $r$~\sf ies \er{1.60}, \te s $s\in \wt E\sms{\wt e}$ \st $a=r+s$ in
view of \er{1.46}. \E\Tf $(n+1)\cdot a\nde1.18 = a+(n\cdot a) = (r+s)+
(n\cdot a)\nde1.6 = n\cdot a+(r+s) \nde1.5 = (n\cdot a+r)+s \nde1.59 = b+s$.
Hence $b<(n+1)\cdot a$ by \er{1.46}. Thus $b\in A_n$. Thus if $b=n'\cdot a
+r'$ with $n',r'\in \wt E\t \wt E$ together with $r'\in[e,a)$, then $b\in A_{n'}$. \If
from \er{1.62} that $n'=n$, hence $b=n\cdot a+r'$. Finally from $n\cdot a
+r=n\cdot a+r'$ we infer $r=r'$ using \er{1.6}, \er{1.8}.
\endproof

\bex1.39
Let $a,b,n,r$ be as in Theorem \rf{t1.38}. Show that $n$ is the \ti{greatest
\el\/} of the set $C:=\{m\cdot a\in\wt E: m\cdot a\le b \}$.
\eex

In the proof of Theorem \rf{t1.38} we used several times the fact that the
maps $x\mt x+y$, $y\mt x+y$, $x,y\in\wt E$, are strictly in\cre. Thus we have
\fa $x,y,z\in\wt E$:
\beq1.63
x\le z\ \hbox{(resp.\ }x<z) \qh{implies }x+y\le z+y
\ (\hbox{resp.\ }x+y<z+y)
\e
and
\beq1.64
y\le z\ \hbox{(resp.\ }y<z) \qh{implies }x+y\le  x+z
\ (\hbox{resp.\ }x+y<x+z).
\e
\Mo the \fw\ holds.
\bea1.65
x+y&\le \hbox{(resp.\,}<)\,z+y \hbox{ implies }x\le \hbox{(resp.\,}<)\,z,
\\
x+y&\le \hbox{(resp.\,}<)\,x+z \hbox{ implies }y\le \hbox{(resp.\,}<)\,z.
\lb{1.66}
\e

\bex1.40
Prove \er{1.63}--\er{1.66}.
\eex

In the next
lemma we show among other things that the range of a non-in\jc\ \hm sm from
$(\wt E,+_{\td E},\wt e)$ into a monoid is finite. The proof of this lemma
uses Theorem \rf{t1.38}.

\blm1.41
Let $\ees$ be a set of \nn s, let $\le$ be its \nog\ and let $(\wt E,+_{\td
E},\wt e)$ be the \PM\ introduced in \E\Pr\ \rf{p1.13}. Let $(M,\qu,e)$ be
a~monoid and let $\vf:\wt E\to M$ be a~\ti{non-in\jc} \hm sm. Then \te\
$m,n\in \wt E$ with $m<n$ \st the \fw\ holds.

\hph i,ii, The \rt ion of $\vf$ to $[\wt e,m)$ is in\jc.

\hph ii,i, The range of $\vf$ is equal to $\vf([\wt e,n))$.

\hph iii,, The range of $\vf$ is finite.
\elm

\proof
We use the simpler notation $n\cdot a$ instead of $n\dpl a$.

(i) Since $\vf$ is not in\jc, \te\ $i,j\in \wt E$ \st $\vf(i)=\vf(j)$ and
$i\ne j$. \E\wlg  we may assume $i<j$ since the \og~$\le$ is total. Set
\[
B:=\bigl\{i\in \wt E: \hbox{\te s $j\in\wt E$ (depending on $i$) \st $i<j$ and }
\vf(i)=\vf(j)\}.
\]
From what precedes $B$ is a nonempty subset of~$\wt E$. Since $\le$ is
a~well-\og, \te s a~least \el\ of~$B$ which we denote by~$m$. If $m=\wt e$,
then (i)~is trivially \sf ied since $[\wt e,m)$ is empty. If $m=\wt S(\wt e)$,
then $[\wt e,m)=\{\wt e\}$ and (i)~is also trivially \sf ied. Suppose $m>\wt
S(\wt e)$ and $r,s\in\wt E$ with $\wt e\le r<s<m$, then we have $\vf(r)\ne
\vf(s)$, for otherwise $r\in B$, $r<m$, \cd ing the fact that $m$~is the
least \el\ of~$B$. \If that $\vf|_{[\td e,m)}$ is in\jc.

(ii) Since $m$ is an \el\ of~$B$, \te s $n\in \wt E$ \st $m<n$ and $\vf(m)
=\vf(n)$. In view of \er{1.46} \te s $p\in\wt E\sms{\wt e}$ \st $n=p+m$,
hence
\beq1.66a
\vf(m+p)=\vf(m).
\e
Let $l\in[\wt e,p)$. Using \In\ and Notation \rf{n1.36}  we prove
\beq1.67
\vf(m+k\cdot p+l) = \vf(m+l), \q k\in\wt E.
\e
Set $A:=\{k\in\wt E: \hbox{\er{1.67} holds}\}$. We have $\wt e\in A$,
since $\vf(m+\wt e\cdot p+l)\nad{\rm\er{1.25}\er{1.17}}=
{\vf(m+\wt e+l)}\nde1.7 = \vf(m+l)$. Suppose $k\in A$. Then $\vf(m+\wt S(k)\cdot p
+l)\nde1.18 = \vf(m+p+k\cdot p+l)\break \nad{\er{1.5},\er{1.6}} = \vf(m+k\cdot p+p+l) \nad{\er{1.5},\er{1.6}} =
\vf(m+k\cdot p+l+p) \nde1.16 = \vf(m+k\cdot p+l)\qu\vf(p) \nad{k\in A}=
\vf(m+l)\qu\vf(p) = \vf(m+l+p)\nad{\er{1.5},\er{1.6}} = \vf(m+p+l)= \vf(m+p)\qu\vf(l) \nde1.66a = \vf(m)
\qu\vf(l)=\vf(m+l)$. Hence $\wt S(k)\in A$. \If from (1.N2) that $A=\wt E$,
so \er{1.67} holds.

Since $\vf([\wt e,n))\sbs\vf(\wt E)$, it suffices to
prove $\vf(\wt E)\sbs\vf([\wt e,n))$. Since the \og~$\le$ is total, we have
$\wt E=[\wt e,m)\cup [m,{\to})$ using \era1{3.25}. Hence $\vf(\wt E)=
\vf([\wt e,m))\cup \vf([m,{\to}))$. Recall that $m<n$, hence $[\wt e,m)\sbs
[\wt e,n)$ by \era1{3.42}. Thus $\vf([\wt e,m))\sbs\vf([\wt e,n))$. It
remains to show that $\vf([m,{\to}))\sbs\vf([\wt e,n))$. Let $s\in
[m,{\to})$. In view of \er{1.46}, \te s $b\in\wt E$ \st $s=m+b$, $b\in\wt E$.
Recall that $p\in \wt E\sms{\wt e}$. Using Theorem \rf{t1.38}
with $b:=b$, $a:=p$, we find $k\in \wt E$ and $l\in[\wt e,p)$ \st $b=k\cdot
p+l$. \E\Tf $\vf(s)=\vf(m+k\cdot p+l)\nde1.67 = \vf(m+l)$. Since $l<p$, we
have $m+l<m+p$ by \er{1.64}. Hence $m+l<n$, and $m+l\in[\wt e,n)$ by
\era1{3.13} and \era1{3.29}. \E\Tf $s\in[\wt e,n)$, hence
$\vf(s)\in\vf([\wt e,n))$.

\advance\baselineskip by-0.3pt
(iii) Since $\wt e\le m$ by \era1{3.13} and $m<n$, we have $\wt e<n$ by
\era1{3.11}. In view of Lemma \rfa1{l3.43} \te s $c\in\wt E$ \st $[\wt e,n)
=[\wt e,c]$. \E\Tf $[\wt e,n)$ is finite by \E\df\ \rfa1{d4.17}. Then
(iii) follows from (ii) and Theorem \rfa1{t4.18}\,(ii) since $[\wt e,n)$ is
finite.
\endproof

\Wanp prove one of the main results of this section.

\bth1.42
Let $(M,\qu,e)$ be an infinite \pn\ monoid with \Gn~$a$. Then $M$ is a \PM\
and $a$ is the only \Gn\ of~$M$.
\eth

\proof
Let $\ees$ be a set of \nn s and let $\vf_a:\wt E\to M$ denote the \sq\ of
\IT s of~$a$ introduced in \E\Pr\ \rf{p1.13}. Viewed as a~map from $(\wt E,
+_{\td E},\wt e)$ introduced in \E\Pr\ \rf{p1.2} into the monoid $(M,\qu,e)$,
$\vf_a$ is a~sur\jc\ \hm sm by Lemma \rf{l1.21}\,(iii). By part~(ii) of the same
lemma, $(M,\qu,e)$ is abelian. \Mo $\vf_a$ is in\jc\ by Lemma \rf{l1.41} since $M=\vf_a
(\wt E)$ is infinite. \E\Tf $\vf_a$ is an \is sm from $\wt E$ onto~$M$ by
Lemma \rf{l1.8}\,(i). We recall that $(\wt E,+_{\td E},\wt e)$ is a \PM\ by
\E\Pr\ \rf{p1.2}.

We now show that \cn s \er{1.8}, \er{1.9} are \sf ied by the \op~$\qu$.

\er{1.8}: Let $x,y,z\in M$ be \st $x\qu z=y\qu z$. Let $k,l,m\in\wt E$ be
\st $\vf_a(k)=x$, $\vf_a(l)=y$ and $\vf_a(m)=z$, and let $\psi:M\to \wt E$
denote the inverse of~$\vf_a$, which is also a~\hm sm by Lemma \rf{l1.8}\,(iii).
Then $k+_{\td E}m = \psi(x)+_{\td E}\psi(z) = \psi(x\qu z) = \psi(y\qu z)=
\psi(\vf_a(l)\qu\vf_a(m)) = \psi(\vf_a(l))+_{\td E}\psi(\vf_a(m)) =l+_{\td E}
m$. Using \er{1.8} for $+_{\td E}$, we obtain $k=l$, hence $x=\vf_a(k) =
\vf_a(l)=y$.

\er{1.9}: Let $x,y\in M$ be \st $x\qu y=e$. Then $\wt e=\psi(e)=\psi(x\qu y)
=\psi(x)+_{\td E}\psi(y)$. In view of \er{1.9} we have $\psi(x)=\psi(y)=\wt e$.
Hence $k=l=\wt e$, and $x=y=e$.

Finally, $a$ is the only \Gn\ of~$M$ by Theorem \rf{t1.25}\,(ii).
\endproof

\bex1.43
Show that the direct product of two infinite \pn\ monoids is not a~\pn\
monoid.
\eex

We conclude this section by considering \ti{submonoids} of an infinite \pn\ monoid.
We recall (see \E\df s \rfa1{d2.2}) that if $Y$~is a \sbm\ of a~monoid
$(X,\qu,e)$, then $(Y,\qu,e)$ is a~monoid.
One verifies that a~\sbm\ of a~C- (resp.~P-) \sbm\ is a~C- (resp.~P-) monoid.
However, a~\sbm\ of
a~\pn\ monoid does not need to be \pn. In the next exercise we consider such an
example.

We first observe that if $(X,\qu,e;u)$ denotes an infinite \pn\ monoid with
\Gn~$u$, and if $Y$~is a~\pn\ \sbm\ with \Gn\ $a\in Y\sms e$, then the set
of \IT s of~$a$ in $(Y,\qu,e)$ is the same as the set of \IT s of~$a$ in
$(X,\qu,e)$ (see \E\df s and Notation \rf{d1.14}) in view of Lemma \rf{l1.21}.
So we may denote this set by $I(a)$ and we have $Y=I(a)$. \E\Ip if
$(X,\qu,e;u):=(\wt E,+_{\td E},\wt e;\wt S(\wt e))$, then $Y=\{\vf_a(n):
n\in\wt E\}$ where $\vf_a:\wt E\to \wt E$ is the \sq\ of \IT s of~$a$
introduced in \E\Pr\ \rf{p1.13}. \If from Lemma \rf{l1.37}\,(iv) that the \Gn\
$a$ of~$Y$ is the least \el\ of the set $Y\sms e$, \wrt the \nog\ $\le$ of
$\ees$.

\Wanp give an example of a \sbm\ $Y$ of $X$ which is not \pn.

\bex1.44
Let $(X,\qu,e;u):=(\wt E,+_{\td E},\wt e;\wt S(\wt e))$ be as above and
let $\le$ denote the \nog\ of $\ees$. Set $\b:=\wt S(\wt S(\wt e))$ and
$\g:=\wt S(\b)$. Let $Y:=\{m{\dpl}\b+n{\dpl}\g: m,n\in\wt E\}$ (see Notation
\rf{n1.36}). Show the \fw:

\hph i,ii, $Y$ is a nontrivial \sbm\ of $X$.

\hph ii,i, $\b$ is the least \el\ of $Y\sms e$ \wrt the \og~$\le$.

\hph iii,, There is no $k\in\wt E$ \st $\g=k\dpl\b$.

Conclude in view of the discussion preceding the exercise that $(Y,+_{\td E},
\wt e)$ is not a \pn\ monoid.
\eex

\advance\baselineskip by0.3pt
\newpage

\Subsubsection{Multiplication and exponentiation}\label{sss.Mult}

The purpose of this section is to introduce two basic \art\ \op s: the \mlc\
and the \epc. We first continue our study of the \IT s of an \el\ of an
\ti{abelian} monoid. Applying these results to the abelian monoid $(E,+,0)$ where\glossary{$(E,0,S)$}
$(E,0,S)$ is a set of \nn s with distinguished \el~$0$ and where the \ad~$+$
is defined in \E\Pr\ \rf{p1.2}, we naturally arrive at the notion of \mlc.
Various \pp ies of this \op\ are summarized in Theorem \rf{t2.14}. The notion
of \epc\ follows from the notion of \IT s of an \el\ of the abelian monoid
$(E,\cdot,1)$ where $\cdot$~denotes the \mlc\ on~$E$, and $1:=S(0)$ (see
Lemma \rf{l0.27} below). 
Basic \pp ies of the \epc\ are \es ed in \E\Pr s \rf{p2.16} and \rf{p2.17}.

In Section \ref{sss.Add} we used the notation $(M,\qu,e)$ for a monoid with
\nel~$e$, and we used the notation $\ees$ for a set of \nn s with
distinguished \el~$\wt e$. In \E\Pr\ \rf{p1.2} we introduced a~binary \op\
$+_{\td E}$ on~$\wt E$, called \ad, which makes $\wt E$ an abelian monoid with
\nel~$\wt e$: $(\wt E,+_{\td E},\wt e)$. We showed in Lemma \rf{l1.17} that
this monoid is an infinite \pn\ monoid with \Gn\ $\wt S(\wt e)$. Note that
the notion of \pn\ monoid does not depend on a~particular choice of a~set of
\nn s for defining the \IT s of \el s of the monoid. In Theorem \rf{t1.42}
we showed that an infinite \pn\ monoid has only one \Gn\ and is a~\PM.

In this section we shall denote by $(E,+,0;1)$ a~\ti{fixed\/} but arbitrary
infinite \pn\ monoid with \nel~$0$ and \Gn~$1$ (see Notation \rf{n1.29} and
Theorem \rf{t1.42}). \If from \E\Pr\ \rf{p1.30} that if the map $S:E\to E$ is
defined by
\beq2.1
S(n):=n+1 \ (=1+n), \q n\in E,
\e
then $(E,0,S)$ is a set of \nn s with \su\ \f~$S$, \Ip $1=S(0)$. \Mo the
\ad~$+$ coincides with the binary \op\ introduced in \E\Pr\ \rf{p1.2} with
$\ees:=(E,0,S)$. \csq\ we shall say that a~subset $A$ of~$E$ is \ti{\iv} if
$0\in A$ and if \fe $n\in A$, we have $n+1\in A$. Then Axiom (1.N2) becomes
\beq2.2
E\hbox{ is the only \iv\ subset of }E.
\e

\blm0.27
Let $(\wt E,\wt e,\wt S)$ be a set of \nn s $($s.o.n.n.$)$. Then \te s a
s.o.n.n.\ $(E,0,S)$ \st $S(0)=1\ne0$.
\elm

\proof
We first construct a s.o.n.n.\ $(E',0,S')$. We consider three cases:
(i)~$\wt e=0$, (ii)~$\wt e\ne0$ and $0\notin\wt E$, (iii)~$\wt e\ne0$ and
$0\in\wt E$. Case~(i): We set $E':= \wt E$, $e':=0$, $S':=\wt S$. Then
$(E',0,S')$ is a s.o.n.n. Case~(ii): We set $E':=(\wt E\sms{\wt e})\cup\{0\}$
(i.e.\ we replace $\wt e$ by~$0$), $e':=0$, and $S'(0):=\wt S(\wt e)$,
$S'(x):=\wt S(x)$ \fa $x\in \wt E \sms{\wt e}$. Then $(E',0,S')$ is a s.o.n.n.
Case~(iii): We set $E':=\wt E$. Define $f:\wt E\to E'$ by setting $f(\wt e)
:=0$, $f(0):=\wt e$, $f(x):=x$ \fa $x\in \wt E\sms{\wt e,0}$, and $g:E'\to\wt E$
by setting $g(0):=\wt e$, $g(\wt e):=0$, $g(x):=x$ \fa $x\in\wt E\sms{0,\wt e}$.
Then $g\circ f=\id_{\td E}$, $f\circ g=\id_{E'}$, $f$~is a bi\jn\ and $g=f\Inv$.
Set $S':=f\circ S\circ f\Inv$ and note that $f(\wt e)=0$. Then $(E',0,S')$ is
a s.o.n.n.\ in view of Exercise \rfa1{ex1.5}. If $1\in E'$ and $S'(0)=1$
then set $E:=E'$, $S:=S'$ and $(E,0,S)$ is the required s.o.n.n. If $1\notin
E'$, then set $E:=(E'\sms{S'(0)})\cup\{1\}$ and $S(0):=1$, $S(1):=S'(S'(0))$
and $S(x):=S'(x)$ \fa $x\in E\sms{0,1}$. Then $(E,0,S)$ is the required s.o.n.n.
Finally, if $1\in E'\sms0$ and $S'(0)\ne 1$, we proceed as in Case~(iii) above.
Set $E:=E'$, define $f:E'\to E$ by setting $f(S'(0)):=1$, $f(1):=S'(0)$, and
$f(x):=x$ otherwise. Define $g:E\to E'$ by setting $g(1):=S'(0)$, $g(S'(0))
:=1$, and $g(x):=x$ otherwise. Set $S:=f\circ S'\circ f\Inv$. Observe that
$f(0)=0$ since $0\notin\{1,S'(0)\}$. Indeed, $1\ne0$ by \as\ and $S'(0)\ne0$
by \era1{3.13}. Hence $(E,0,S)$ is the required s.o.n.n.\ in view of Exercise
\rfa1{ex1.5}.
\endproof

\bnt2.1
Given $m,n\in E$, we shall use the notation $m\dpl n$ for the $m$-fold \IT\
of~$n$ \wrt $(E,0,S)$ in the monoid $(E,+,0)$ (see \E\df s and Notation
\rf{d1.14}).
\ent

In the next \Pr\ we collect some \pp ies of the ``\IT\ map'' introduced in
\E\Pr\ \rf{p1.13}.

\setbox9=\hbox{I0}
\def\lca#1 #2{\hbox to0.6\textwidth{\hbox to\wd9{\rm #1}\quad $#2$\hfil}}
\bpr2.2
Let $(M,\qu,e)$ be an \emph{abelian} monoid. Given $n\in E$ and $x\in M$, let
$n\dqu x$ $($or simply $n\cdot x)$ denote the $n$-fold \IT\ of~$x$ \wrt
$(E,0,S)$ in $(M,\qu,e)$. Then the \fw\ holds.

\E\fa $m,n\in E$ and all $x,y\in M$ we have
\begin{gather}
\begin{cases}
\lca I0 {0\cdot x=e,}\\
\lca I1 {1\cdot x=x,}
\end{cases}\non\\
\begin{cases}
\lca I2 {(m+n)\cdot x=m\cdot x \qu n\cdot x,}\\
\lca I3 {(m\dpl n)\cdot x=m\cdot(n\cdot x),}
\end{cases} \lb{2.3}\\
\begin{cases}
\lca I4 {m\cdot e=e,}\\
\lca I5 {m\cdot(x\qu y)=m\cdot x\qu m\cdot y.}\non
\end{cases}
\end{gather}
Here $m\cdot x\qu n\cdot x:=(m\cdot x)\qu(n\cdot x)$ and
$m\cdot x\qu m\cdot y:=(m\cdot x)\qu(m\cdot y)$.
\epr

\proof
It is \sft\ to prove I2 and I3. Indeed, I0 follows from \er{1.17}, I1~follows
from \er{1.17}, \er{1.18}, I4 follows from \er{1.19} and I5 from \er{1.24}.

I2: In view of \er{1.20} we have $m\cdot x\qu n\cdot x=(\Phi_x(m)\circ
\Phi_x(n))(e)$, where $\Phi_a(m)$ is defined in \E\Pr\ \rf{p1.13}. From
\er{1.42} we obtain $\Phi_x(m)\circ \Phi_x(n)=\Phi_x(m+n)$. \Mo we have
$\Phi_x (m +n)(e)\nde1.33 = ((m+n)\cdot x)\qu e \nde1.7 = (m+n)\cdot x$.

I3: Let $n\in E$ be fixed. We prove I3 by \In\ on $m\in E$. Let $A:= \{m\in
E: $\ I3~holds$\}$. We claim that $A$ is \iv. We have $0\in A$. Indeed, $(0\dpl
n)\cdot x=0\cdot x$ from~I0 for $(E,+,0)$. Then $0\cdot x\nad{\rm I0}= e
\nad{\rm I0}= 0\cdot (n\cdot x)$, hence $0\in A$. We suppose $m\in A$. We have
$((m+1)\dpl n)\cdot x= (n+(m\dpl n))\cdot x$ by \df\ of $(m+1)\dpl n$.
\Mo $(n+(m\dpl n))\cdot x \nad{\rm I2}= n\cdot x+(m\dpl n)\cdot x\nad{m\in A}=
n\cdot x+m\cdot (n\cdot x)\nde1.18 = (m+1)\cdot (n\cdot x)$. \E\Tf $((m+1)
\dpl n)\cdot x=(m+1)\cdot(n\cdot x)$, hence $m+1\in
A$. \If that $A=E$ by \er{2.2} since $A$~is \iv.
\endproof

\bxa2.2
Let $F$ be a \ns\ and let $f$ be a selfmap of~$F$. We denote by $\Phi(n)\in
F^F$ the $n$-fold \IT\ of $f$ \wrt $(E,0,S)$ (see \E\df\ \rfa1{d2.1}). Let
$m\in E$. We \ti{claim} that the $m$-fold \IT\ of
$\Phi(n)$ \wrt $(E,0,S)$ is
$\Phi(m\dpl n)$ where $m\dpl n$ is defined in Notation \rf{n2.1}. We first
observe that $\Phi(n)$ defined in \era1{2.4} is the $n$-fold \IT\ of~$f$ \wrt
$(E,0,S)$ (see \E\df s and Notation \rf{d1.14}) in the monoid $(F^F,\circ,
\id_F)$ introduced in Example \rfa1{xa2.3}\,(i).
In view of \er{1.25} we have $\Phi(n)=n\dci f$, hence by \er{1.26},
$\Phi(n)\in I(f)$, the \sbm\ of $(F^F,\circ,\id_F)$ introduced in Example
\rfa1{xa2.3}\,(ii). \If from Lemma \rf{1.21}\,(i) that $\Phi(n)$ is also the
$n$-fold \IT\ of~$f$ \wrt $(E,0,S)$ in $(I(f),\circ,\id_F)$, which we also
denote by $n\dci f$. Since $\Phi(n)$ belongs to $I(f)$, the $m$-fold \IT\ of
$\Phi(n)$ in the sense of \E\df\ \rfa1{d2.1}, which is also the $m$-fold \IT\
of $\Phi(n)$ in the monoid $(F^F,\circ,\id_F)$, is again by Lemma
\rf{l1.21}\,(i), the $m$-fold \IT\ of $\Phi(n)$ in the monoid $(I(f),\circ,
\id_F)$. \E\Tf we have $m\dci\Phi(n)=m\dci(n\dci f)$, where in the \LHS\
$\dci$~means the \IT\ in $(F^F,\circ,\id_F)$ and in the \RHS\ means the \IT\
in $(I(f),\circ,\id_F)$. Since $(I(f),\circ,\id_F)$ is an \ti{abelian} monoid
(see Example \rfa1{xa2.3}\,(ii)), we can apply \E\Pr\ \rf{p2.2} with $M:=I(f)$,
${\qu}:=\circ$ and $e:=\id_F$. Then by \er{2.3}~I3 we obtain $m\dci\Phi(n) =
(m\dpl n)\circ f$. Now, since the $(m\dpl n)$-fold \IT\ of~$f$ in
$(F^F,\circ,\id_F)$ is $\Phi(m\dpl n)$ by \era1{2.4}, we have $m\dci
\Phi(n)=\Phi(m\dpl n)$, which proves the claim.
\exa

We see from Example \rf{xa2.2} that the notation $\Phi(n)$, introduced in
\E\df\ \rfa1{d2.2}, for the $n$-fold \IT\ of a selfmap of a~set $F\ne\vn$,
is not very convenient when dealing with ``\IT s of \IT s'' of a~map~$f$.
\E\Tf \ti{from now on}, we shall use the \fw\ standard notation:
\beq2.4
f^n:=\Phi(n), \q n\in E,\ f\in F^F, \ F\ne\vn.
\e
Using this notation we obtain from Example \rf{xa2.2}:
\beq2.5
(f^n)^m = f^{m\dpl n}, \q m,n\in E,\ f\in F^F,\ F\ne\vn.
\e
Formula \er{1.42} becomes
\beq2.6
f^m \circ f^n = f^{m+n}, \q m,n\in E,\ f\in F^F,\ F\ne\vn.
\e
We also have
\beq2.7
f^0 = \id_F \qh{and }f^1=f,\ f\in F^F,\ F\ne\vn.
\e

We recall that if $F:=E$ and $f:=S$ the \su\ \f\ of $(E,0,S)$, then the map
$n\mt S^n$ from $E$ into $E^E$ is \ti{in\jc} by \era1{2.36}. \E\Tf given
$m,n\in E$, $m+n$ is the \ti{only} \el\ $p\in E$ \st $f^m\circ f^n=f^p$. This
fact was used for the motivation of the \df\ of the \ad\ in $(E,0,S)$ (see \E\Pr\
\rf{p1.2}). We now use the fact that given $m,n\in E$, $m\dpl n$ is the only
\el\ $p\in E$ \st $(S^n)^m = S^p$, in order to motivate the \df\ of the binary
\op\ on~$E$ called \ti{\mlc}.

\bdn2.3
Let $(E,+,0;1)$ be an infinite \pn\ monoid with \nel~$0$ and \Gn~$1$. Given
$m,n\in E$ we set
\beq2.8
m\cdot n:= m\dpl n \qh{(see Notation \rf{n2.1}).}
\e
The binary \op\ $\cdot$ on $E$ defined by \er{2.8} is called \ti{\mlc}.\index{multiplication} If
$m,n\in E$, then $m\cdot n$ is called the \ti{product\/} of $m$ and~$n$, and
it is also denoted by~$mn$ when no confusion arises. The \el s $m$~and~$n$ are
called \ti{factors}.\index{factors}
\edn

\bpr2.4
The \mlc\ $\cdot$ on $E$ is the only binary \op\ on~$E$ \sf ying\dw
\bea2.9
{}&m\cdot0 = 0 \qh{\fa $m\in E$,}\\
&m\cdot(n+1) = m\cdot n +m := (m\cdot n)+m \qh{\fa $m,n\in E$.} \lb{2.10}
\e
\Mo the \fw\ holds. \E\fa $a,b,c\in E$ we have\dw
\bea2.11
(ab)c&=a(bc), \\
ab&=ba, \lb{2.12} \\
0c&=0=c0 \hbox{ and }1c=c=c1, \lb{2.13} \\
a(b+c)&=ab+ac, \lb{2.14} \\
(a+b)c&=ac+bc, \lb{2.15}
\e
where $ab+ac:=(ab)+(ac)$.
\epr

\proof
In what follows we use \E\Pr\ \rf{p2.2} with $(M,\qu,e):=(E,+,0)$ and
$n\dqu x:=n\dpl x=nx$. We first show that the \mlc~$\cdot$ defined in \er{2.8}
\sf ies \er{2.9} and~\er{2.10}. Observe that \er{2.9} directly follows from
\er{2.3}\,I4. Now, let $m,n\in E$. Then $m\cdot(n+1) = m\cdot n+m\cdot1$ by
\er{2.3}\,I5. It remains to show that $k\cdot1=k$ \fa $k\in E$. We define
$\vf:E\to E$ by setting
$\vf(k):=k\cdot1$, $k\in E$. Note that $\vf(0)=0\dpl1=0$ by \er{2.3}\,I0,
$\vf(1)=1$ by \er{2.3}\,I1 and $\vf(k+l)=(k+l)\dpl1 = (k\dpl1)+(l\dpl1)=
\vf(k)+\vf(l)$ by \er{2.3}\,I2. Hence $\vf$~is an endo\mf\ of $(E,+,0)$ \sf ying $\vf(S(0))=1
=S(0)$. Thus by Lemma \rf{l1.37}\,(iii) with $(\wt E,\wt e,\wt S) := (E,0,S)$
we obtain $\vf=\id_E$. \csq, $m\cdot1=\vf(m)=m$. Thus $\cdot$~\sf ies \er{2.9},
\er{2.10}.

\er{2.13}: $c\cdot1=c$ by what precedes, and $c=1\cdot c$ by \er{2.3}\,I1;
$c\cdot 0=0$ by \er{2.9} and $0\cdot c=0$ by \er{2.3}\,I0.

\er{2.14}: $a\cdot(b+c) \nad{\er{2.3}\,\rm I5}= a\cdot b+a\cdot c$.

\er{2.15}: $(a+b)\cdot c \nad{\er{2.3}\,\rm I2} = a\cdot c+b\cdot c$.

\er{2.11}: $(a\cdot b)\cdot c \nad{\er{2.3}\,\rm I3}= a\cdot(b\cdot c)$.

\er{2.12}: Let $b\in E$. Set $A:=\{a\in E: a\cdot b=b\cdot a\}$. We show that
$A$~is \iv.

``$0\in A$'': Since $0\cdot b=0=b\cdot0$ by \er{2.13}.

``\ti{$n\in A$ implies $n+1\in A$}'': Let $n\in A$. Then $(n+1)\cdot a\nde2.15 =
n\cdot a+1\cdot a\nde2.13 = n\cdot a+a\cdot1 \nad{n\in A}= a\cdot n+a\cdot1
\nde2.14 = a\cdot(n+1)$. Hence $n+1\in A$. Thus $A$~is \iv, and $A=E$ by
(1.N2).

It remains to prove that if $\cdot'$ is a binary \op\ on~$E$ \sf ying \er{2.9},
\er{2.10} where $\cdot$~is replaced by~$\cdot'$, then $m\cdot n=m\cdot' n$ \fa
$m,n\in E$. Let $m\in E$ and let $A:=\{n\in E: m\cdot n=m\cdot'n\}$. We show
that $A$~is \iv.

``$0\in A$'': By \er{2.9} $m\cdot0=0=m\cdot'0$.

``\ti{$n\in A$ implies $n+1\in A$}'': Let $n\in A$. Then $m\cdot(n+1)\nde2.10 =
m\cdot n+m \nad{n\in A}= m\cdot'n+ m\nde2.10 = m\cdot'(n+1)$. Hence $n+1\in A$.
Thus $A$~is \iv\ and $A=E$.
\endproof

\brm2.5 \

\hph i,i, \If from \er{2.11}--\er{2.13} that $(E,\cdot,1)$ is an abelian
monoid.

\hph ii,, \If from \er{2.9} and \er{2.14} that \fe\ $m\in E$, the map $n\mt
m\cdot n$ from~$E$ into itself is an endo\mf\ of $(E,+,0)$.
\erm

\bex2.6
Let $\F:E\to E$ be an endo\mf\ of $(E,+,0)$. Show that \te s \ooo $m\in E$ \st
$\F(n)=mn$ \fa $n\in E$. What are the auto\mf s of $(E,+,0)$?
\eex

The abelian monoid $(E,\cdot,1)$ does \ti{not\/} \sf y \er{1.8} since $0\cdot
0=0=0\cdot 1$ with $0\ne1$. However, it \sf ies \er{1.9}. In what follows it
will be convenient to use the notation
\beq2.17
E^*:=E\sms0.
\e
We recall that the map $S$ defined in \er{2.1} is a \ti{bi\jn} from~$E$
onto~$E^*$ by (1.N1).

\bpr2.7
\beq2.16
\left\{
\vcenter{\hsize=0.85\textwidth
\hph i,ii, If $a,b\in E^*$, then $ab\in E^*$.\endgraf
\hph ii,i, If $a,b\in E$ and $ab=0$, then $a=0$ or $b=0$.\endgraf
\hph iii,, If $a\in E^*$, $b,c\in E$, then $ab=ac$ implies $b=c$.\endgraf
\hph iv,, If $a,b\in E$, then $ab =1$ implies $a=b=1$.}\right.
\e
\epr

\proof \

(i) Since $a,b\in E^*$ which is the range of $S$, \te\ $x,y\in E$ \st
$a=x+1$ and $b=y+1$. $ab=(x+1)(y+1)\nde2.15 = x(y+1)+1(y+1)\nad{\rm \er{2.14},
\er{2.13}}= (xy+x1)+(y+1) \nad{\er{2.13}}= (xy+x)+(y+1)\nde1.5 =
((xy+x)+y)+1 = S((xy+x)+y)\in E^*$.

(ii) Follows from (i).

(iii) Using the \nog\ $\le$ of $(E,0,S)$ we first suppose $b\le c$, i.e.\
$c=b+p$ \fs $p\in E$. Then $ab=ac$ implies $ab+0\nad{\rm\er{1.7},\er{2.14}}=
ab+ap$. Using \er{1.6}, \er{1.8} we infer $ap=0$. Since $a\in E^*$, we have $p=0$
by~(ii). Hence $b=c$. Interchanging $b$ and~$c$ we reduce the case $c\le b$ to
the case $b\le c$.

(iv) In view of \er{2.13} we have $a0=0=0b$, hence $a,b\in E^*$. As
in the proof of~(i), we have $a=x+1$, $b=y+1$, with $x,y\in E$, and $ab=
((xy+x)+y)+1$. Since $ab=1=0+1$, we obtain $(xy+x)+y=0$ by \er{1.8}. Using
\er{1.9} we have $xy+x=0$ and $y=0$. But $xy=0$ by \er{2.13}, hence
$x\nde1.7 = 0+x = xy+x = 0$. From $x=y=0$, we infer $a=b=1$.
\endproof

\bpr2.8 \

\hph i,ii, $(E,\cdot,1)$ is an \am.

\hph ii,i, $E^\ast$ is a \sbm\ of the monoid $(E,\cdot,1)$.

\hph iii,, $(E^*,\cdot,1)$ is a~\PM.
\epr

\proof \

(i) follows from \er{2.11}--\er{2.13}.

(ii) follows from \er{2.16}\,(i) and $1\in E^\ast$.

(iii) follows from \er{2.16}\,(iii) and (iv).
\endproof

We shall prove later that $(E^*,\cdot,1)$ is not \pn,
hence not \is c to $(E,+,0)$.

As already observed, using \er{2.13}, we have
\beq2.18
m\cdot 1=m \qh{\fa $m\in E$}.
\e
It turns out that \er{2.18} and \er{2.10} can be used to define the \mlc\
on~$E^*$.

\bex2.9
Using Theorem \rfa1{t1.6} show that \te s \ooo binary \op\ on~$E^*$, which we
denote by~$\qu$, \sf ying \er{2.18} and \er{2.10} where $\cdot$ is replaced
by~$\qu$. Show that $\qu$ coincides with the \rt ion to~$E^*$ of the \mlc\
on~$E$ defined in \er{2.8}.
\eex

\def\namespec{Historical remark}
\begin{rspc}
In E. Landau's \cite{Grund} \nn s are the \el s of what we denote by
$(E^*,1,S)$. \E\Ip \mlc\ is defined by \er{2.18} and \er{2.10}, and \asc ity,
\cmt ity and \cnc ity of this \op\ are deduced from these axioms and \pp ies
of the \ad\ on~$E^*$. At the end of his \ti{Preface for the Student\/}, he
made the \fw\ comment about his daughters: ``\ti{,
my daughters have been studying $($Chemistry$)$ at the University for several
semesters already and think that they learned the differential and integral calculus
in College{\rm;} and yet they don't know why $x\cdot y=y\cdot x$\,.}''
\end{rspc}

\bex2.11
Let $\ees:=(E^*,1,S)$ in Section \ref{sss.Add}. Show that the \ad\ $+_{\td
E}$ defined in \E\Pr\ \rf{p1.2} \sf ies $a+_{\td E}b= 1+((a-1)+(b-1))$ where
$a-1$ is defined in Remark \rf{1.24}. \E\Ip $1+_{\td E}1=1$. What is the
\crs\ \mlc\ defined by $m\cdot_{\td
E}n=m\mathrel{\lower2pt\hbox{$\stackrel{+_{\tilde E}}\cdot$}}n$, the $m$-fold
\IT\ of~$n$ \wrt $\ees$ in $(\wt E,+_{\td E},\wt e)$\,?
\eex

In the next \Pr\ we consider analogues of \er{1.64} and \er{1.65} for the
\mlc.

\bpr2.12
Let $x,y\in E$ and $z\in E^*$. Then
\bea2.19
x\le y \ (\hbox{resp.\ $x<y)$ implies }xz\le yz\ (\hbox{resp.\ $xz<yz)$},\\
xz\le yz\ (\hbox{resp.\ $xz<yz)$ implies }x\le y\ (\hbox{resp.\ $x<y)$}.
\lb{2.20}
\e
\epr

\bex2.13
Prove \E\Pr\ \rf{p2.12}.
\eex

In the next theorem we collect some of the results obtained in Section
\ref{sss.Iterates} and in this section, about the \ad\ and the \mlc.

\bth2.14
Let $(E,+,0;1)$ be an infinite \pn\ monoid with \Gn~$1$. Then the \fw\ holds.

\sdim{III\,(ii)}
\ssk
\ite[\rm \hphantom{II}I\,(i)\hphantom{i}]
{$(E,+,0)$ is a \PM.}
\ite[\rm \hphantom{III\,}(ii)]
{Let $\cdot$ denote the \mlc\ on~$E$ defined in \E\Pr\ \rf{p2.4}. Then
$(E,\cdot,1)$ is an abelian monoid.}
\ite[\rm \hphantom{I}II\,(i)\hphantom{i}]
{\E\fe $m\in E$, the map $\vf_m:E\to E$ defined by $\vf_m(n):=m\cdot n$ \fa
$n\in E$, is an endo\mf\ of the monoid $(E,+,0)$.}
\ite[\rm \hphantom{III\,}(ii)]
{\E\fa $m,n\in E$ we have $\vf_m(n)=\vf_n(m)$.}
\ite[\rm III\,(i)\hphantom{i}]
{$E^*\ (:=E\sms0)$ is a \sbm\ of the monoid $(E,\cdot,1)$.}
\ite[\rm \hphantom{III\,}(ii)]
{The monoid $(E^*,\cdot,1)$ is a \PM.}
\ssk

Let $S:E\to E$ be the map defined in \er{2.1} and let $(E,0,S)$ be the \crs\
set of \nn s. Let $\le$ denote the \og\ of $(E,0,S)$. Then\dw

\ssk
\ite[\rm IV\,(i)\hphantom i]
{\E\fe $m\in E$, the map $\th_m:E\to E$ defined by
$\th_m(n):=m+n$, \fa $n\in E$, is in\jc\ and strictly
in\cre. \Mo if $\th_m(n)<\th_m(p)$ with $m,n,p\in E$, then $n<p$.}
\ite[\rm \hphantom{IV\,}(ii)]
{Let $\vf_m$, $m\in E$, be as in {\rm II\,(i)}. Then $\vf_0(n)=0$ \fa $n\in E$,
$\vf_1=\id_E$. \Mo \fa $m\in E^*$, $\vf_m$ is in\jc\ and strictly in\cre\ and
if $\vf_m(n)<\vf_m(p)$, $n,p\in E$, then $n<p$.}
\eth

\proof \
{\setbox0=\hbox{I}
\sz=\wd0 \parindent0pt

\lvm \kern2\sz I\,\hph i,i, follows from Theorem \rf{t1.42}.

\lvm \kern3\sz \,\hph ii,, see Remark \rf{r2.5}\,(i).

\lvm \kern\sz II\,\hph i,i, see Remark \rf{r2.5}\,(ii).

\lvm \kern3\sz \,\hph ii,, follows from \er{2.12}.

\lvm III\,\hph i,i, follows from \er{2.16}\,(i).

\lvm \kern3\sz \,\hph ii,, follows from \er{2.16}\,(iii) and (iv).

\lvm IV\,\hph i,i, In\ji\ of $\th_m$ follows from \er{1.6}, \er{1.8} and strict in\cre
ness from \er{1.66}.

\lvm \kern3\sz\,\hph ii,, $\vf_0(n)=0$ and $\vf_1(n)=n$ for $n\in E$ follow from
\er{2.13}.

}Now  let $m\in E^*$. Then $\vf_m:E\to E$ is strictly in\cre\ and in\jc\ by
Lemma \rf{l1.37}\,(iv), (v). Finally, let $n,p\in E$ be \st
$\vf_m(n)<\vf_n(p)$. Then $n\ne p$ in view of the \df\ of~$<$. Suppose for
\cd ion that $p<n$. Then by what precedes we obtain $\vf_m(p)<\vf_m(n)$. From
$\vf_m(p)<\vf_m(n)$ and $\vf_m(n)<\vf_n(p)$, we obtain $\vf_m(p)<\vf_m(p)$ by
\era1{3.11}, a~\cd ion. Hence $n<p$ since the \og~$\le$ is total.
\endproof

We now introduce the \op\ ``\epc'' which is nothing but the notion of \IT\ of
an~\el\ of the monoid $(E,\cdot,1)$. Using notation \er{1.25}, we have

\bdf2.15
Let $m,n\in E$ and let $n\ddt m$ denote the $n$-fold \IT\ of~$m$ \wrt
$(E,0,S)$ in the abelian monoid $(E,\cdot,1)$. Set
\beq2.21
m^n := n\ddt m \qh{\fa $m,n\in E$.}
\e
The \nm\ $m^n$ is called the \ti{$n$-th power} of~$m$\index{n-th power of $m$@$n$-th power of $m$} (or $m$~to the
power~$n$), and $n$ is called the \ti{exponent\/}. Given $m\in E$, the map
$n\mt m^n$ from~$E$ into itself is called \ti{\epc}.\index{exponentiation} Note that in view
of \er{2.3}\,I0 we have
\beq2.22
0^0:=1.
\e
\edf

\bpr2.16
\E\fa $m,n,p$ and $q$, \el s of $E$ we have\dw
\bea2.23
{}&m^0=1,\\
&m^1=m, \lb{2.24}\\
&m^{p+q} = m^p\cdot m^q, \lb{2.25}\\
&m^{p\cdot q}=(m^p)^q, \lb{2.26}\\
&1^p = 1, \lb{2.27}\\
&(m\cdot n)^p = m^p\cdot n^p, \lb{2.28}\\
&0^p=1\hbox{ if }p=0 \qh{and } 0^p=0 \hbox{ if }p\in E^*. \lb{2.29}\\
&\hbox{If $m\in E^*$, then }m^p\in E^*. \lb{2.30}\\
&\hbox{If $m,p\in E$ and $m^p=0$, then $m=0$ and }p\in E^*. \lb{2.31}\\
&\hbox{If $m\in E\sms1$ and $p\in E$, then $m^p=1$ implies }p=0. \lb{2.32}
\e
\epr

\proof  \

\er{2.23} follows from \er{2.3}\,I0.

\er{2.24} follows from \er{2.3}\,I1.

\er{2.25} follows from \er{2.3}\,I2.

\er{2.26} follows from \er{2.3}\,I3, \er{1.6}.

\er{2.27} follows from \er{2.3}\,I4.

\er{2.28} follows from \er{2.3}\,I5.

\er{2.29}: $0^0=1$ by \er{2.22}. We prove $0^p=0$, $p\in E^*$, by \In\ on
$p\in E^*$. Recall that $(E^*,1,S)$ is a set of \nn s (see e.g.\ Exercise
\rfa1{ex1.5}). Set $A:=\{p\in E^*: 0^p=0\}$. Clearly $1\in A$ since $0^1\nde
2.24 = 0$. Suppose $p\in A$. Then $S(p)=p+1$ and $0^{p+1}\nde2.25 = 0^p\cdot
0^1= 0^p\cdot 0\nde2.9 = 0$. Hence $A$ is \iv\ in $(E^*,1,S)$ and  $A=E^*$.

\er{2.30}: Follows from Lemma \rf{l1.21}\,(i) with $(M,\qu,e):=(E,\cdot,1)$,
$M_0:=E^*$ $(=E\sms0)$ and $a:=m\in M_0\sms 1$. If $m:=1$, then \er{2.30}
follows from \er{2.27}.

\er{2.31}: If $m^p=0$, then $m=0$, since otherwise $m\in E^*$ and $m^p\in E^*$ by
\er{2.30}, which is impossible since $0\notin E^*$. \Mo $0^p=0$ only if $p\in
E^*$ by \er{2.29}.

\er{2.32}: Suppose for \cd ion that $p\in E^*$. Then \te s $q\in E$ \st
$p=q+1$. Hence $m^p=m^{q+1}\nde2.25 = m^q\cdot m^1\nde2.24 = m^q\cdot m$. But
if $m^q\cdot m=1$, then $m^q=m=1$ by \E\Pr\ \rf{p2.7}\,(iv). A~\cd ion since
$m\ne 1$.
\endproof

\bpr2.17 \

\hph i,i, Let $m\in E\sms{0,1}$ and let $p,q\in E$. Then we have
\bea2.33
p<q &\hbox{ implies }m^p<m^q, \\
m^p=m^q &\hbox{ implies }p=q, \lb{2.34} \\
m^p<m^q &\hbox{ implies }p<q. \lb{2.35}
\e

\hph ii,, Let $p\in E^*$ and let $m,n\in E$. Then we have
\bea2.36
m<n &\hbox{ implies }m^p<n^p, \\
m^p=n^p &\hbox{ implies }m=n, \lb{2.37} \\
m^p<n^p &\hbox{ implies }m<n. \lb{2.38}
\e
\epr

\proof \

\er{2.33}: Since $p<q$, \te s $r\in E^*$ \st $q=p+r$. Hence $m^q=m^{p+r}
\nde2.25 = m^p\cdot m^r$. Then $m^r\in E^*$ by \er{2.30} since $m\ne0$.
Suppose for \cd ion that $m^r=1$. Then $r=0$ by \er{2.32} since $m\ne1$.
Hence $q=p+0=p$, a~\cd ion. Since $m^r\in E\sms{0,1}$, we have $1<m^r$. \Mo
$m^p\in E^*$ by \er{2.30} since $m\ne0$. \E\Tf we have $m^p\nde2.13 = 1\cdot
m^p\nde2.19 < m^r\cdot m^p\nde2.12 = m^p\cdot m^r\nde2.25 = m^{p+r}= m^q$. Hence
$m^p<m^q$.

\er{2.34}: Suppose for \cd ion $p<q$. Then by \er{2.33} $m^p<m^q$, a~\cd ion.
If $q<p$, then $m^q<m^p$, a~\cd ion. Hence $p=q$, since the \og\ is total.

\er{2.35}: $p\ne q$ since otherwise $m^p=m^q$; $q\not<p$ since otherwise $m^q
<m^p$ by \er{2.33} and $m^p< m^q$ by \as. Hence $m^q<m^q$, a~\cd ion. Hence
$p<q$.

\er{2.36}: Let $m,n\in E$ and $p\in E^*$ be \st $m<n$. We show $m^p<n^p$. If
$m=0$, then $n\in E^*$, hence $m^p\nde2.29 = 0<n^p$ since $n^p\in E^*$ by
\er{2.30}. If $m\in E^*$, then $n\in E^*$. We prove \er{2.36} by \In\ on
$p\in E^*$. Set $A:=\{p\in E^*: m^p<n^p \hbox{ \fa}m,n\in E^*\}$. Then $1\in
A$ by \as\ and \er{2.24}. Suppose $p\in A$. Let $m,n\in E^*$ be \st $m<n$. Then
$m^{p+1}\nad{\er{2.25},\er{2.24}} = m^p\cdot m \nad{\er{2.19},p\in A} <n^p\cdot
m\nde2.12 = m\cdot n^p\nde2.19 < n\cdot n^p = n^{p+1}$. Hence $p+1\in A$, by
\era1{3.11}. Thus $A=E^*$, and \er{2.36} holds.

\er{2.37}: The proof of \er{2.37} is similar to the proof of \er{2.34}, \Tf
it is omitted.

\er{2.38}: The proof of \er{2.38} is similar to the proof of \er{2.35}.
\endproof

\brm2.18
If we view the map $(m,n)\mt m^n$ from $E\t E$ into~$E$ as a binary \op\
on~$E$, then this binary \op\ is neither \asc e nor \cmt e.

\ti{Non-\asc ity}: \E\fe $n\in E$ we have $(n^1)^n\nde2.26 = n^{1\cdot n}
\nde2.13 = n^n$ and $n^{(1^n)}\nde2.27 = {n^1\nde2.24 = n}$. If $n:=S(1)\nde2.1
= 1+1$, then $(1+1)^{1+1}\nde2.25 = (1+1)^1 \cdot (1+1)^1 \nde2.24 =
(1+1)\cdot(1+1)\nde2.19 > (1+1)\cdot1$ since $1+1>1$ by \er{1.46}. \Mo
$(1+1)\cdot1\nde2.13 = 1+1 = S(1)=n$. Hence $n^n>n$ by \era1{3.12}. \E\Tf
$(n^1)^n = n^n > n=n^{(1^n)}$, thus $(n^1)^n \ne n^{(1^n)}$.

\ti{Non-\cmt ity}: $0^1\nde2.24 = 0\ne 1\nde2.23 = 1^0$, hence $0^1\ne1^0$.
\erm

We conclude this section by showing as promised that the \PM\ $(E^*,\cdot,1)$
is \ti{not\/} \pn. Suppose for \cd ion that \te s $a\in E^*\sms1$ \st $E^*=\break
\{a^n\in E^*: n\in E\}$. Since $a^0=1$ by \er{2.23} we have $\{a^n\in E^*:
n\in E\} = \{1\}{\cup} \{a^n\in E^*:\break n\in E^*\}$.

\bex2.19
Set $a:=S(S(S(e)))$. Show that $(a^a)^a \ne a^{(a^a)}$. See Remark \rf{r1.11}.
\eex

\bex2.21
Show that the \PM\ $(E^*,\cdot,1)$ is not \pn.
\eex

\newpage
\Subsubsection{Cardinality of finite sets}\label{sss.Repr}
We recall that a \ns~$A$ is called \ti{finite} if \te s an \os\ of \nn s
$(E,\le,e,S)$, an \el~$n$ of~$E$ and a bi\jn\ from~$A$ onto the \il\ $[e,n]$,
i.e.\ if $A$ is \ep\ to $[e,n]$. The empty set~$\vn$ is also a finite set by
\df\ (see \E\df s \rfa1{d4.16} and \rfa1{d4.17}). In order to treat nonempty
and empty finite sets on the same footing we first extend the notion of
equipotence introduced in \E\df\ \rfa1{d4.16}.

\bdf3.1
Let $A$ and $B$ be sets (not necessarily distinct). The sets $A$ and~$B$ are
said to be \ti{\ep\/} or \ti{equipollent\/} or \ti{to have the same \ca} if
\ti{either} they are both empty \ti{or} both $A$ and~$B$ are not empty and
\te s a~bi\jn\ from $A$ onto~$B$. We write $A\ax B$ (resp.~$A\not\ax B$) if
$A$ and~$B$ are (resp.\ are not) \ep.
\edf

Clearly \E\df\ \rf{d3.1} extends \rfa1{d4.16}.

\blm3.2
Let $A,B,C$ be sets. Then the \fw\ holds\dw
\bea3.1
&A\ax B \hbox{ and } B\ax C \hbox{ implies } A\ax C,\\
&A\ax A, \lb{3.2} \\
&A\ax B \hbox{ implies } B\ax A. \lb{3.3}
\e
\elm

\bex3.3
Prove Lemma \rf{l3.2}.
\eex

We now give an \ev t \df\ of a \nf. It is based on the \fw\ observation.

\blm3.4
Let $(E,\le,e,S)$ be an \os\ of \nn s. Then
\beq3.4
[e,n] = [e,S(n)) \qh{\fa} n\in E.
\e
\elm

\proof
Let $n\in E$. We apply Lemma \rfa1{l3.43} with $a:=e$, $b:=S(n)$ and we find $\zo e,S(n)
=[e,c]$ \fs \el\ $c\in\zo e,S(n) $ \sf ying $S(c)=S(n)$. From the in\ji\
of~$S$ we infer $c=n$. Hence $[e,n]=\zo e,S(n) $.
\endproof

The next lemma allows us to treat both empty and \nf s ``on the same
footing''.

\blm3.5
A set $A$ is finite $($in the sense of \E\df\ \rfa1{d4.17}$)$ iff \te\ an \os\
of \nn s $(E,\le,e,S)$ and $m\in E$ \st $A$~is \ep\ $($in the sense of \E\df\
\rf{d3.1}$)$ to $[e,m)$.
\elm

\proof \

\ti{If\/}: Suppose $A\ax\zo e,m $ \fs $m\in E$. If $m\ne e$, then $m\in R(S)$,
hence there is $n\in E$ \st $m=S(n)$, hence $\zo e,m \ax[e,n]$ by Lemma
\rf{l3.4}. From \er{3.1} we obtain $A\ax[e,n]$, thus $A$~is finite in the
sense of \E\df\ \rfa1{d4.17}. If $m=e$, then $\zo e,m =\vn$ by \era1{3.29},
since there is no $x\in E$ \st $e\le x<e$, otherwise $e<e$, by \era1{3.11},
\cd ing the \df\ of~$<$. By \E\df\ \rf{d3.1}, $A$~is empty since $A\ax
\zo e,m $ and $\zo e,m =\vn$. Since $A=\vn$, it is finite in the sense of
\E\df\ \rfa1{d4.17}.

\ti{Only if\/}: Suppose $A$ is finite in the sense of \E\df\ \rfa1{d4.17}.
If $A=\vn$, and $(E,\le,e,S)$ is an \os\ of \nn s, then $A\ax\zo e,e $ in
the sense of \E\df\ \rf{d3.1}, since $\zo e,e =\vn$ (by the ``if'' part of
the proof). If $A\ne\vn$, then \te\ a~set of \nn s $(E,\le,e,S)$ and $n\in E$
\st $A\ax[e,n]$. But $[e,n]\ax\zo e,S(n) $ by Lemma \rf{l3.4}. Hence $A\ax
\zo e,S(n) $ by \er{3.1}.
\endproof

The next theorem shows in particular that if $A$ and~$B$ are finite sets, then
one can find a~common \os\ of \nn s in Lemma \rf{l3.5}.

\bth3.6
Let $(E,\le,e,S)$ be an \os\ of \nn s and let $A$ be a finite set. Then \te s
\ooo $m\in E$ \st $A\ax\zo e,m $ and one writes\glossary{$\#_E(A)$}
\bga3.5
\#_E(A)=m \qh{$($not standard notation$)$.}\\
\hbox{\E\Ip}\#_E(A)=e \hbox{ iff }A=\vn. \lb{3.6}
\e
If moreover $(E',\le',e',S')$ is an \os\ of \nn s and if $\vf:E\to E'$ is the
bi\jn\ introduced in Theorem \rfa1{t1.4}, then
\beq3.7
\#_{E'}(A)=\vf(m).
\e
\eth

\proof
We first consider the case $A=\vn$. Since $\zo e,e =\vn$, we have $A=\vn$ iff
$A\ax\zo e,e $ by \E\df\ \rf{d3.1}. \Mo if $\zo e,e \ax\zo e,m $ \fs $m\in E$,
we have $m=e$. Indeed, if $m\ne e$, then $e<m$ by \era1{3.13} and $e\in
\zo e,m $. Then $\zo e,m \ne\vn$ and $\zo e,m \not\ax\zo e,e $. We now assume
$A\ne\vn$. In view of Lemma \rf{l3.5}, \te\ a~set of \nn s $(\wt E,\wt \le,
\wt e,\wt S)$ and an \el\ $\wt m\in \wt E$ \st $A\ax\zo {\wt e,\wt m} \wt{\ }$,
where $\zo {\wt e,\wt m} \wt{\ }:=\{\wt x\in \wt E: \wt e\mathrel{\wt \le}\wt x\mathrel{\wt<}\wt m
\}$. Let $\wt\vf:\wt E\to E$ be the bi\jn\ introduced in Theorem \rfa1{t1.4}.
Note that $\zo{\wt E,\wt m} \wt{\ }\ne\vn$ since $A\ne\vn$. We recall that
$\wt\vf$ is an \ois sm by Theorem \rfa1{t3.36}. \If that \fe $\wt x\in
\zo{\wt e,\wt m} \wt{\ }$, we have $\wt\vf(\wt e)\le \wt\vf(\wt x)<\wt\vf(\wt m)$,
hence $\wt\vf$ maps $\zo{\wt e,\wt m} \wt{\ }$ \ti{into} $\zo{\wt\vf(\wt e),
\wt\vf(\wt m)} $. Conversely, \fe $x\in\zo{\wt\vf(\wt e),\wt\vf(\wt m)} $, i.e.\
$\wt\vf(\wt e)\le x<\wt \vf(\wt m)$, we have $\wt e\mathrel{\wt\le} \wt\vf{}\Inv(x)\mathrel{\wt<}
\wt m$ since $\wt\vf{}\Inv$ is strictly in\cre. \E\Tf setting $\wt x:=
\wt\vf{}\Inv(x)\in\zo{\wt e,\wt m} \wt{\ }$, we find that $x$~belongs to the
image of $\zo{\wt e,\wt m} \wt{\ }$ under~$\wt\vf$. Hence $\wt\vf$~maps $\zo{
\wt e,\wt m} \wt{\ }$ \ti{onto} $\zo{\wt\vf(\wt e),\wt\vf(\wt m)} $. Since
$\wt\vf$ is in\jc, its \rt ion to $\zo{\wt e,\wt m} \wt{\ }$ is also in\jc,
hence $\zo{\wt e,\wt m} \wt{\ }\ax\zo{\wt\vf(\wt e),\wt\vf(\wt m)} $.

From $A\ax\zo{\wt e,\wt m} \wt{\ }$ and $\zo{\wt e,\wt m} \wt{\ } \ax \zo
{\wt\vf(\wt e),\wt\vf(\wt m)} $ we infer $A\ax\zo{\wt\vf(\wt e),\wt\vf(\wt m)}
$ by \er{3.1}. Recalling that $\vf(\wt e)=e$ and setting $m:=\wt\vf(\wt m)
\in E$, we obtain $A\ax\zo e,m $. Note that since $A\ne\vn$, $m\in E\sms e$.
We now prove the \uq\ of such an \el\ $m\in E\sms e$. Suppose \te s $p\in E
\sms e$ \st $A\ax\zo e,p $. Then $\zo e,m \ax\zo e,p $ in view of \er{3.1},
\er{3.3}. Since $m,p\in E\sms e$, \te\ $n,q\in E$ \st $m=S(n)$ and $n=S(q)$
by (1.N1). Using Lemma \rf{l3.4}, \er{3.1}, \er{3.3} we obtain $[e,n]\ax
[e,q]$. From Corollary \rfa1{c4.13}\,(ii) we infer $n=q$. Hence $m=S(n)=S(q)
=p$, which implies \uq.

Now let $\vf:E\to E'$ be as in Theorem \rfa1{t1.4}. Proceeding as above, we
find $\zo e,m \ax\zo \vf(e),\vf(m) '$. Since $\vf(e)=e'$, we have $A\ax
\zo e',\vf(m) '$ by \er{3.3}. Hence \er{3.7} holds in view of the \uq\ of
$\vf(m)$.
\endproof

\bnt3.7 \

\hph i,i, In the remaining part of this section, $(E,\le,e,S)$ will denote
a~\ti{fixed\/} but arbitrary \os\ of \nn s. We shall use the notation $\#(A)$
instead of $\#_E(A)$ for a finite set~$A$ unless explicitly mentioned.\glossary{$\#(A)$}

\hph ii,, We recall that in view of (1.N1) the map $S:E\to E\sms e$ is bi\jc.
We shall denote by~$P$ the inverse of~$S$, i.e.\ the map from $E\sms e$
onto $E$ \sf ying $P\circ S=\id_E$ and $S\circ P=\id_{E\sms e}$.
\ent

\blm3.8
Let $A$ be a \nf. Then $\#(A)\in E\sms e$ and we have
\beq3.8
A\ax[e,P(\#(A))].
\e
\elm

\proof
We have $A\ax\zo e,\#(A) $ and $\#(A)\in E\sms e$ by Theorem \rf{t3.6}. But
$\#(A)=S(P(\#(A)))$. Hence $\zo e,\#(A) \ax[e,P(\#(A))]$ by Lemma \rf{l3.4}.
Thus $A\ax[e,P(\#(A))]$ by \er{3.1}.
\endproof

\Wanp strengthen \E\Pr\ \rfa1{p4.27}\,(A)(iii) and (B)(i).

\bpr3.9
Let $X,Y$ be \ns s and let $f$ be a~map from~$X$ into $Y$.

\hph i,i, If $f$ is in\jc\ and $Y$ is finite, then $X$ is finite and $\#(X)
\le\#(Y)$.

\hph ii,, If $X$ is finite, then $f(X)$ is finite and we have
\beq3.9
S(e)\le \#(f(X)) \le \#(X).
\e
\epr

\proof \

(i): $X$ is finite by \E\Pr\ \rfa1{p4.27}\,(A)(iii). Since both $X$ and~$Y$ are
\nf s, we have $X\ax[e,P(\#(X))]$ and $Y\ax[e,P(\#(Y))]$ by \er{3.8}. Let
$\vf:X\to[e,P(\#(X))]$ and $\psi:Y\to[e,P(\#(Y))]$ be bi\jn s. Then
$\wt f:=\psi\circ f\circ\vf\Inv : [e,P(\#(X))]\to[e,P(\#(Y))]$ is in\jc\ as
\cm\ of in\jc\ maps. We obtain $P(\#(X))\le P(\#(Y))$ by Lemma \rfa1{l4.11}.
Since $S$ is in\cre, we arrive at $\#(X)=S(P(\#(X))) \le S(P(\#(Y))) =\#(Y)$.

(ii): We can view $f$ as a sur\jc\ map from~$X$ onto $Y:=f(X)$ which is
nonempty since $X$~is nonempty. Then $Y$ is finite by \E\Pr\
\rfa1{p4.27}\,(B)(i). \E\Tf $\#(Y)\ne e$ by Theorem \rf{t3.6}, hence $S(e)\le
\#(Y)$ by \era1{3.13}. It remains to show $\#(Y)\le\#(X)$. Let $\wt f$~be as
in part~(i) of the proof. Since $f$~is sur\jc, so is~$\wt f$, as \cm\ of
sur\jc\ maps. We obtain $P(\#(X))\le P(\#(Y))$ by Corollary \rfa1{c4.13}\,(i).
\E\Tf as above we arrive at $\#(Y)\le\#(X)$.
\endproof

\brm3.10
\E\Pr\ \rf{p3.9}\,(i) is often used in the \fw\ \ev t form. If $X$ and~$Y$ are
\nf s with $\#(Y)<\#(X)$, and $f$~is a map from $X$ into~$Y$, then \te\
$x,\wh x\in X$, $x\ne\wh x$, \st $f(x)=f(\wh x)$.
\erm

\bco3.11
Let $X,Y$ be \ep\ sets. If $X$ $($resp.~$Y)$ is finite, then $Y$ $($resp.~$X)$
is finite and $\#(X)=\#(Y)$.
\eco

\proof
By interchanging $X$ and~$Y$ we may restrict ourselves to the case
$X$~is finite. If $X=\vn$, then $Y=\vn$ by \E\df\ \rf{d3.1}, and $\#(X)=e=
\#(Y)$ by \er{3.6}. If $X\ne\vn$, then $Y\ne\vn$ by \E\df\ \rf{d3.1}. Since
$X\ax Y$, \te s a~bi\jn\ $f:Y\to X$, hence $Y$~is finite and $\#(Y)\le
\#(X)$ by \E\Pr\ \rf{p3.9}\,(i). Then $f\Inv :X\to Y$ is also bi\jc, hence
in\jc. \If from the same \Pr\ that $\#(X)\le\#(Y)$, hence $\#(X)=\#(Y)$.
\endproof

We now consider a converse of Corollary \rf{c3.11}.

\blm3.12
Let $X,Y$ be finite sets \sf ying $\#(X)=\#(Y)$. Then $X\ax Y$.
\elm

\proof
We have $X\ax\zo e,\#(X) $ and $Y\ax\zo e,\#(Y) $. Since $\#(X)=\#(Y)$, we have
$\zo e,\#(X) =\zo e,\#(Y) $. Hence $X\ax \zo e,\#(Y) $ and $\zo e,\#(Y) \ax Y$
by \er{3.3}. Then $X\ax Y$ follows by \er{3.1}.
\endproof

\Wanp generalize Lemma \rfa1{l4.14} and Theorem \rfa1{t4.18}\,(iii).

\bpr3.13
Let $X,Y$ be \nf s \sf ying $\#(X)=\#(Y)$.
\bea3.10
&\hbox{If $X\sbs Y$, then $X=Y$.}\\
&\hbox{If $f:X\to Y$ is in\jc, then $f$ is bi\jc.}\lb {3.11} \\
&\hbox{If $f:X\to Y$ is sur\jc, then $f$ is bi\jc.}\lb {3.12}
\e
\epr

\proof \

\er{3.10}: In view of Theorem \rf{t3.6} and Lemma \rf{l3.8} \te\ bi\jn s
$\vf:X \to [e,n]$, $\psi:Y\to[e,n]$ where $n:=P(\#(X))=P(\#(Y))$. Since $X
\sbs Y$, we have $\psi(X)\sbs \psi(Y)=[e,n]$. By \as\ $\#(X)=\#(Y)$, hence
$X\ax Y$ by Lemma \rf{l3.12}. Since $\psi$~is a bi\jn\ we have $\psi(X)\ax X$
and $Y\ax \psi(Y)$ by \E\df\ \rf{d3.1}. Using \er{3.1} we obtain $\psi(X)\ax
Y$, $Y\ax \psi(Y)$, and $\psi(X)\ax \psi(Y)$. By the \df\ of~$\psi$, $\psi(Y)
=[e,n]$, hence $\psi(X)\ax[e,n]$ by \er{3.2}, \er{3.1}. Since $\psi(X)\sbs
[e,n]$, we infer from Lemma \rfa1{l4.14} that $\psi(X)=[e,n]$. Note that
$[e,n]=\psi(Y)$, hence $\psi(X)=\psi(Y)$, and $X=Y$ follows since $\psi$~is
bi\jc.

\er{3.11}, \er{3.12}: Let $\vf$, $\psi$ and $n$ be as in the proof of
\er{3.10}. Define $\wt f:[e,n]\to[e,n]$ by setting $\wt f:=\psi\circ f\circ
\vf\Inv$. Then $\wt f$ is in\jc\ (resp.\ sur\jc) in the case of \er{3.11}
(resp.\ \er{3.12}). Since $[e,n]$ is finite, $\wt f$~is bi\jc\ by Theorem
\rfa1{t4.18}\,(iii) both in case \er{3.11} and in case \er{3.12}. \If that
$f=\psi\circ \wt f\circ\vf$ is bi\jc\ in both cases.
\endproof

We now consider a converse of \E\Pr\ \rf{p3.9}.

\bpr3.14
Let $X,Y$ be \nf s.

\hph i,i, If $\#(X)\le \#(Y)$, then \te s an in\jc\ map $f:X\to Y$.

\hph ii,, If $\#(Y)\le \#(X)$, then \te s a sur\jc\ map $g:Y\to X$.
\epr

\bex3.15
Prove \E\Pr\ \rf{p3.14}.
\eex

We now state a strengthening of Theorem \rfa1{t4.18}\,(i).

\blm3.16
Let $A,B$ be sets \st $B\sbs A$ with $A$ finite. Then $B$ is finite and
\beq3.13
\#(B)\le\#(A) \qh{with \et y iff} B=A.
\e
\elm

\proof
If $B=\vn$, then $B$ is finite and $\#(B)=e$ by \er{3.6}. Then $e\le\#(A)$
follows from \era1{3.13}, with \et y iff $A=\vn$ by \er{3.6} again, i.e.\ iff
$B=A$.

If $B\ne\vn$ then $A\ne\vn$ since $B\sbs A$. Thus both $A$ and $B$ are not
empty and the inclusion map $f:B\to A$ defined by $f(x):=x$, $x\in B$, is
in\jc. Then $B$ is finite and $\#(B)\le\#(A)$ by \E\Pr\ \rf{p3.9}\,(i). If
\et y holds, then $A=B$ in view of \er{3.10}. Finally, if $A=B$, then $\#(A)
=\#(B)$.
\endproof

We now \es\ some \rl\ between the \ca\ of $A\cap B$ and the \ca\ of $A,B$
where $A,B$ are finite sets. Since $A\cap B\sbs A$, we have $A\cap B$ finite
and $\#(A\cap B)\le\#(A)$ with \et y iff $A\cap B=A$ by \er{3.13}. But
$A\cap B=A$ iff $A\sbs B$. Similarly $\#(A\cap B)\le\#(B)$ with \et y iff
$A\cap B=B$ iff $B\sbs A$. Hence we have both $\#(A\cap B)\le\#(A)$ and
$\#(A\cap B)\le\#(B)$.

We shall reformulate this statement using the \fw\ \df.

\bdf3.17
Let $(X,\le)$ denote a \tos\ and let $x,y\in X$.\glossary{$\max(x,y)$}\glossary{$\min(x,y)$}
\bea3.14
\max(x,y)&:=\bca y &\hbox{if }x\le y,\\ x &\hbox{if }y\le x,\eca\\
\min(x,y)&:=\bca x &\hbox{if }x\le y,\\ y &\hbox{if }y\le x.\eca \lb{3.15}
\e
\edf

\blm3.18
Let $(X,\le)$ be a \tos\ and let $a,b,c\in X$. Then we have
\bea3.16
{}&\hbox{if }a\le c \hbox{ and }b\le c, \hbox{ then }\max(a,b)\le c,\\
{}&\hbox{if }c\le a \hbox{ and }c\le b, \hbox{ then }c\le \min(a,b), \lb{3.17} \\
&a\le\max(a,b), \quad b\le\max(a,b), \lb{3.18} \\
\noalign{\eject}
&\min(a,b)\le a,\q \min(a,b)\le b, \lb{3.19} \\
&\max(a,a)=a, \q \min(a,a)=a, \lb{3.19a} \\
&\max(a,b)=\max(b,a), \q \min(a,b)=\min(b,a), \lb{3.20} \\
&\max(a,\max(b,c)) = \max(\max(a,b),c), \lb{3.21} \\
&\min(a,\min(b,c)) = \min(\min(a,b),c), \lb{3.22} \\
&\max(a,\min(a,b)) = a, \q \min(a,\max(a,b))=a, \lb{3.23} \\
&b=\max(a,b) \hbox{ iff }a\le b, \q b=\min(a,b)\hbox{ iff }b\le a,\lb{3.24} \\
&a\le b \hbox{ implies } \max(a,c)\le \max(b,c) \hbox{ and }
\min(a,c)\le\min(b,c). \lb{3.25}
\e
\elm

\bex3.19
Prove Lemma \rf{l3.18}.
\eex

\blm3.20
Let $A,B$ be finite sets. Then $A\cap B$ is finite and
\beq3.25a
\#(A\cap B)\le \min(\#(A),\#(B)) \qh{with \et y iff} A\sbs B \hbox{ or }
B\sbs A.
\e
\elm

\proof
In view of the discussion preceding Lemma \rf{l3.18} it suffices to prove the
in\et y. \If from \er{3.17} and the fact that $\min(\#(A),\#(B))$ is equal to
$\#(A)$ or $\#(B)$.
\endproof

Our next goal is to show that the union of two finite sets is finite and to
investigate is \ca. The main ingredient of the proof is the next lemma.

\bnt3.21
In what follows $+$ denotes the \ad\ $+_E$ defined in \E\Pr\ \rf{p1.2} with
$(\wt E,\wt e,\wt S):=(E,e,S)$.
\ent

\blm3.22
Let $a,b,c\in E$ with $a\le b$. Then
\beq3.26
\zo a,b \ax \zo a+c,b+c .
\e
\elm

\proof
If $a=b$, then $a+c=b+c$, hence $\zo a,b $ and $\zo a+c,b+c $ are both empty.
\E\Tf $\zo a,b $ and $\zo a+c,b+c $ are \ep.

If $a<b$, then $\zo a,b $ and $\zo a+c,b+c $ are both nonempty.
We define a~map $f:\zo a,b \to E$ by setting $f(x):=x+c$, $x\in\zo a,b $. In
view of \er{1.63} we have $a+c\le f(x)<b+c$ whenever $a\le x<b$. Hence $f$
maps $\zo a,b $ \ti{into} $\zo a+c,b+c $. The map~$f$ is \ti{in\jc} in view
of \er{1.8}. We now show that $f$~is sur\jc. Let $y\in\zo a+c,b+c $. In view
of \er{1.46} \te\ $p\in E$ and $q\in E\sms e$ \st
\beq3.27
y=(a+c)+p \qh{and} b+c=y+q.
\e
We have $(a+c)+p\nde1.5 = a+(c+p)\nde 1.6 = a+(p+c)\nde 1.5 = (a+p)+c$. Set
\beq3.28
x:=a+p.
\e
Clearly $a\le x$ by \er{1.46}. We now show that $x<b$. From \er{3.27} we
obtain $b+c=((a+c)+p)+q$. But $((a+c)+p)+q=((a+p)+c)+q \nde1.5 = (a+p)
+(c+q) \nde1.6 = (a+p)+(q+c)=((a+p)+q)+c$. \E\Tf $b+c=((a+p)+q)+c$, hence
by \er{1.8} $b=(a+p)+q$. Since $a+p=x$ by \er{3.28}, we infer $x<b$ by
\er{1.46}. \If that $x\in\zo a,b $, hence $f(x):=x+c=(a+p)+c\nde3.27 = y$.
\csq\ $f$~is sur\jc, hence bi\jc\ and \er{3.26} holds.
\endproof

In the next theorem we are using Notation \rf{n3.7} and \rf{n3.21}.

\bth3.23
Let $A,B$ be finite sets. Then $A\cup B$ is finite and the \fw\ holds.
\bea3.29
{}&\max(\#(A),\#(B)) \le \#(A\cup B) \le \#(A)+\#(B), \\
&\max(\#(A),\#(B)) = \#(A\cup B) \qh{iff} B\sbs A \hbox{ or }A\sbs B,
\lb{3.30} \\
&\#(A\cup B)=\#(A)+\#(B) \qh{iff} A\cap B=\vn, \lb{3.31} \\
&\#(A\cup B)+\#(A\cap B) = \#(A)+\#(B). \lb{3.32}
\e
\eth

\proof \

``\ti{$\max(\#(A),\#(B))\le \#(A\cup B)$ and \er{3.30}}'': $A\sbs A\cup B$
implies $\#(A)\le \#(A\cup B)$ and $\#(A)=\#(A\cup B)$ iff $A=A\cup B$ by
\er{3.13}. Similarly, $B\sbs A\cup B$
implies $\#(B)\le \#(A\cup B)$ and $\#(B)=\#(A\cup B)$ iff $B=A\cup B$. Thus
$\max(\#(A),\#(B))\le \#(A\cup B)$ by \er{3.16}. \E\et y holds in \er{3.30} iff
$A=A\cup B$ or $B=A\cup B$ by \er{3.13} iff $B\sbs A$ or $A\sbs B$.

``\ti{$\#(A\cup B)=\#(A)+\#(B)$ if $A\cap B=\vn$}'':  If $A=\vn$, then $A\cup
B=B$, $\#(A)\nde3.6 = e$ and $\#(A\cup B)=\#(B)\nde1.7 = e+\#(B) = \#(A)+
\#(B)$. The case $B=\vn$ is similar. We now suppose both $A$ and~$B$ nonempty
and finite, \te\ $m,n\in E\sms e$ \st $A\ax\zo e,m $ and $B\ax \zo e,n $ by
\er{3.5}, \er{3.6}. Hence \te\ bi\jc\ maps $\vf:A\to\zo e,m $ and $\psi:B\to
\zo e,n $. By \er{1.7} and Lemma \rf{l3.22} with $a:=e$, $b:=n$ and $c:=m$, we
have $\zo e,n \ax\zo m,m+n $. Since $n\ne e$, \te s a bi\jn\ $\th:\zo e,n \to
\zo m,m+n $. Note that $A\cap B=\vn$ by \as\ and $[e,m)\cap \zo m,m+n =\vn$
by \er{3.27}. Indeed, if \te s $x\in E$ \st $x<m$ and $m\le x$, then $x<x$ by
\era1{3.12}, a~\cd ion. \Mo $\zo e,m \cup \zo m,m+n = \zo e,m+n $. Indeed, if
$x\in E$ \sf ies $x<m$, then $m\le m+n$ by \er{1.46}, hence $x<m+n$ by
\era1{3.12}. Thus $\zo e,m \cup \zo m,m+n \sbs \zo e,m+n $. Conversely, if
$x\in \zo e,m+n $, then $e\le x$ and either $x<m$ or $m\le x$ since the order
is total. In the former case $x\in\zo e,m $ and in the latter case $x\in\zo m,m+n
$. \If that $\zo e,m \cup \zo m,m+n = \zo e,m+n $. We now define a map $f:
A\cup B\to \zo e,m+n $ by setting
\beq3.35
f(x):=\bca
\vf(x) & \hbox{if }x\in A,\\
(\th\circ\psi)(x) & \hbox{if }x\in B.
\eca
\e
Note that $\th\circ\psi: B\to\zo m,m+n $ is a bi\jn\ as a \cm\ of two bi\jn s.
Let $\g\in\zo e,m+n $. Then either $\g\in\zo e,m $ and $f(x)=\g$ with $x=\vf
\Inv(\g)$, or $\g\in\zo m,m+n $ and $f(x)=\g$ with $x=(\th\circ\psi)\Inv(\g)$.
Thus $f$~is \ti{sur\jc}. Now let $x,y\in A\cup B$ be \st $f(x)=f(y)$. Set $\g:=
f(x)$. If $\g\in\zo e,m $, then $x,y\in A$ and $x=y$ since $\vf:A\to\zo e,m $
is in\jc. If $\g\in\zo m,m+n $, then $x,y\in B$ and $x=y$ since $\th\circ\psi:
B\to\zo m,m+n $ is in\jc. Thus $f$~is \ti{in\jc}, hence \ti{bi\jc}. \If that
$A\cup B\ax\zo e,m+n $. \E\Tf $\#(A\cup B)\nde3.5 = m+n\nde3.5 = \#(A)+\#(B)$
since $A\ax\zo e,m $ and $B\ax \zo e,n $.

``\er{3.32}'': Note $A=(A\sm B)\cup (A\cap B)$, $(A\sm B)\cap(A\cap B)=\vn$,
and $B=(B\sm A)\cup(A\cap B)$, $(B\sm A)\cap(A\cap B)=\vn$. Hence by \er{3.31}
$\#(A)=\#(A\sm B)+\#(A\cap B)$, $\#(B)=\#(B\sm A)+\#(A\cap B)$. \E\Tf $\#(A)
+\#(B)=\bigl(\#(A\sm B)+\#(A\cap B)\bigr)+\bigl(\#(B\sm A)+\#(A\cap B)\bigr)
=\#(A)+\bigl(\#(B\sm A)+\#(A\cap B)\bigr) \nde1.5 = \bigl(\#(A)+\#(B\sm A)
\bigr)+\#(A\cap B) \nad*= \#(A\cup B)+\#(A\cap B)$. In $\nad*=$ we used
$A\cup (B\sm A)= A\cup B$, $A\cap(B\sm A)=\vn$, and \er{3.31}. Thus $\#(A)+
\#(B)= \#(A\cup B)+\#(A\cap B)$.

``\ti{If $\#(A)+\#(B)=\#(A\cup B)$, then $A\cap B=\vn$}'': By what precedes
$\#(A\cup B)+\#(A\cap B)=
\#(A)+\#(B) = \#(A\cup B)+e$ by \as\ and by \er{1.5}. Then $\#(A\cap B)=e$ by
\er{1.6}, \er{1.8}, and $A\cap B=\vn$ by \er{3.7}.

``$\#(A\cup B)\le\#(A)+\#(B)$'': $\#(A)+\#(B)\nde3.31 = \#(A\cup B)+\#(A\cap
B) \ge \#(A\cup B)$ by \er{1.46}.

This completes the proof of Theorem \rf{t3.23}.
\endproof

\bnt3.24
It is customary to denote the \ca\ of the empty set by~$0$, and the \ca\ of
a~singleton $A=\{a\}$ with $a\in A$ by~$1$. \csq\ in the remaining part of
this section we shall replace $e$~by~$0$, $S(e)$~by~$1$ and $E\sms 0$ by~$E^*$
(see \er{2.17}) as in Section \ref{sss.Mult}.
\ent

We now establish the \rl\ between the \ca\ of the Cartesian product of two
finite sets and the \mlc\ introduced in \E\Pr~\rf{p2.4}.

\bth3.25
Let $A$ and $B$ be finite sets. Then the set $A\t B$ is finite and
\beq6.28b
\#(A\t B)=\#(A)\cdot \#(B).
\e
\eth

\proof
If $A=\vn$ or $B=\vn$, then $A\t B=\vn$ and $\#(A\t B)=0=\#(A)\cdot\#(B)$, by
\er{3.6} and \er{2.13}. So we may suppose that $A$ and
$B$ are not empty. Let $m:=\#(A)$ and $n:=\#(B)$. Then $m,n\in E^*$ and by
\er{3.5} \te\ bi\jn s $\vf:A\to[0,m)$ and $\phi:B\to [0,n)$. Let $f:A\t B\to
[0,m)\t[0,n)$ be defined by $f(x,y):=(\vf(x),\phi(y))$, \fe $(x,y)\in A\t B$.
$f$~is also a bi\jn. Indeed, set $g(i,j):=(\vf\Inv(i),\phi\Inv(j))$ \fe
$(i,j) \in [0,m)\t[0,n)$. Then we have $g\circ f=\id_{A\t B}$ and
$f\circ g=\id_{[0,m)\t[0,n)}$.

\Mo set $h:[0,m)\t[0,n)\to E$ defined by
\beq6.29a
h(i,j):=i\cdot n+j \qh{\fe} i\in[0,m) \hbox{ and }j\in[0,n).
\e
We claim that
\beq6.30
h([0,m)\t[0,n)) \sbs [0,m\cdot n),
\e
and
\beq6.31
h:[0,m)\t[0,n)\to[0,m\cdot n)
\e
is bi\jc. We first prove \er{6.30}. Since $i<m$, we have $i+1=S(i)\le m$
by~\er{2.1} and \er{3.17}. Since $j<n$, we have $i\cdot n+j<i\cdot n+n=
(i+1)\cdot n\le m\cdot n$ by \er{1.64}, \er{2.12}, \er{2.10} and \er{1.63}. Hence $0\le
h(i,j)<m\cdot n$ by \er{3.11}.

{\it Surjectivity}. Let $b\in[0,m\cdot n)$. Then by Theorem \rf{t1.38}
(recall that $n\in E\sms 0$) \te\ $i,j\in E$ \st $h(i,j)=b$ and
$j\in[0,n)$. It remains to show that $i<m$. We have $i\cdot
n\nad{\er{1.46}}\le i\cdot n+j=b<m\cdot n$. Hence $i<m$ by \er{3.11}, \er{2.19}
since $n\in E^*$.

{\it Injectivity}. Let $(i,j),(k,l)\in[0,m)\t[0,n)$ be \st
$h(i,j)=h(k,l)$. Then by the \uq\ part of the conclusion of
Theorem~\rf{t1.38}, we obtain $i=k$ and $j=l$.

We now combine the results obtained above. Since $f:A\t B\to[0,m)\t[0,n)$
is bi\jc\ and $h:[0,m)\t[0,n)\to[0,m\cdot n)$ is also bi\jc, then $h\circ f
:A\t B\to[0,m\cdot n)$ is bi\jc, hence $A\t B$ is finite by \er{3.8} and
\E\df\ \rfa1{d4.17}, $\#(A\t B)\nad{\er{3.5}}= m\cdot n =\#(A)\cdot\#(B)$.
\endproof

\bex6.13
Let $A$ and $B$ be \ns s. Show that
\beq6.32
A\t B\ax B\t A.
\e
\eex

We conclude this section by showing that the set of all maps ``between''
finite sets is finite.

\bth3.26
Let $A$ and $B$ be \emph{nonempty} finite sets and let $B^A$ be the set of all
maps from $A$ to~$B$. Then
\beq3.38
B^A \hbox{ is finite and }\#(B^A)=(\#(B))^{\#(A)}.
\e
\eth

\proof
In view of Theorem \rf{t3.6} \te\ $m,n\in E^*$, a~bi\jn\ $\vf:A\to\zo 0,n $
and a~bi\jn\ $\psi:B\to\zo 0,m $.
We want to construct a bi\jn\ from $B^A$ onto $\zo0,m ^{\zo0,n }$. To this end we
set
\bea6.48
{}&Rf:=\psi\circ f\circ\vf\Inv, \q f\in B^A,\\
&Sg:=\psi\Inv\circ g\circ\vf, \q g\in\zo0,m ^{\zo0,n }, \lb{6.49}\\
&\xymatrix{
A \ar[r]^f  & B \ar[d] ^\psi\\
\zo0,n  \ar[r]^{Rf} \ar[u]^{\vf\Inv} & \zo0,m
}\qquad
\xymatrix{
\zo0,n  \ar[r]^{g}  & \zo0,m  \ar[d]^{\psi\Inv}\\
A \ar[r]^{Sg}\ar[u] ^\vf  & B } \non
\e
So $f\mt Rf$ is a map from $B^A$ into $\zo0,m ^{\zo0,n }$ and $g\mt Sg$ is a map
from $\zo0,m ^{\zo 0,n }$ into~$B^A$.

We have
\[
\bal
S\circ Rf&=\psi\Inv\circ Rf\circ\vf = \psi\Inv\circ \psi\circ f\circ \vf\Inv
\circ\vf=f, \\
\hbox{and }R\circ Sg&=\psi\circ\psi\Inv\circ g\circ\vf\circ\vf\Inv=g.
\eal
\]
hence
\[
S\circ R=\id_{B^A}, \q R\circ S=\id_{\zo0,m ^{\zo0,n }}.
\]
It follows that $R$ is bi\jc\ and that $S$ is its inverse.

As a con\sq\ it suffices to prove that $\zo0,m ^{\zo0,n }$ is finite and $\#\bigl(
\zo0,m ^{\zo0,n }\bigr)=m^n$. Indeed, in this case since $B^A\ax\zo0,m ^{\zo0,n }$,
$\zo0,m ^{\zo0,n }$ is finite, and $\#(\zo0,m ^{\zo0,n })=m ^n$, we infer by Corollary \rf{c3.11}
that $B^A$ is finite and $\#(B^A)=\#(\zo0,m ^{\zo0,n })$. Hence $\#(B^A)=m^n
=(\#(B))^{\#(A)}$.

We proceed by induction on $n\in E^*$. Let $m\in E^*$ be fixed. Set
\[
M:=\bigl\{n\in E^*: \zo0,m ^{\zo0,n } \hbox{ is finite and }\#\bigl(
\zo0,m ^{\zo0,n }\bigr)=m^n\}.
\]
We have

$1\in M$ since every map $g:\{0\}\to\zo0,m $ is characterized by its value
$g(0)\in \zo0,m $. \E\Tf there is a bi\jn\ between $\zo0,m ^{\{0\}}$ and $\zo0,m $.
Hence $\zo0,m ^{\{0\}}$ is finite and $\#(\zo0,m ^{\{0\}})\nde3.5 = m=m^1$ (by~\er{2.23}).

\ti{$n\in M$ implies $n+1\in M$}:
Suppose $n\in M$. We first claim that $\zo0,m ^{\zo0,n }\t \zo0,m $ is finite
and that $\#(\zo0,m ^{\zo0,n }\t\zo0,m )=m^{n+1}$. Since $n\in M$, $\zo0,m ^
{\zo0,n }$ is finite and $\#(\zo0,m ^{\zo0,n })=m^n$. Note that $\zo0,m $ is
finite, $\id_{\zo0,m }$ maps $\zo0,m $ onto itself and $\#(\zo0,m )\nde3.5 = m$.
Then the claim follows from Theorem \rf{t3.25} with $A:=\zo0,m ^{\zo0,n }$
and $B:=\zo0,m $.

We now construct a bi\jn\ from $\zo0,m ^{\zo0,n+1 }$ onto
$\zo0,m ^{\zo0,n }\t \zo0,m $. Given $g\in\zo0,m ^{\zo0,n+1 }$ we set
$p(g):=$ the \rt ion of $g$ to $\zo0,n $, i.e.\ $p(g)(k):=g(k)$, $k\in\zo0,n $, and
set $Vg:=(p(g),g(n+1))\in\zo0,m ^{\zo0,n }\t\zo0,m $.

We claim that the map $V$ defined above is a bi\jn\ from $\zo0,m ^{\zo0,n+1 }$
onto $\zo0,m ^{\zo0,n }\t\zo0,m $. Indeed, let $W:=\zo0,m ^{\zo0,n }\t\zo0,m
\to\zo0,m ^{\zo0,n+1 }$ be defined by
\[
W((h,r))(i):=\bca
h(i) &\hbox{if }0\le i< n,\\
r    &\hbox{if }i=n,
\eca
\]
where $h\in\zo0,m ^{\zo0,n }$ and $r\in\zo0,m $. Then
\[
W\circ V(g)(i)=\bca
p(g)(i)=g(i), & \hbox{if } 0\le i< n\\
g(n+1), & \hbox{if } i=n
\eca
\ {}=g(i),\ i\in\zo0,n+1 ,
\]
\fe $g\in\zo0,m ^{\zo0,n+1 }$.
\Mo $V\circ W(h,r)=(h,r)$ \fe $h\in\zo0,m ^{\zo0,n }$ and $r\in\zo0,m $. \E\Tf\
$V$ is bi\jc\ and $V\Inv=W$. \If that $\zo0,m ^{\zo0,n }\t \zo0,m \ax
\zo0,m ^{\zo0,n+1 }$. From Corollary \rf{c3.11} with $X:=\zo0,m ^{\zo0,n }
\t\zo0,m $ and $Y:=\zo0,m ^{\zo0,n+1 }$ we infer that $\zo0,m ^{\zo0,n+1 }$ is
finite and  $\#(\zo0,m ^{\zo0,n+1 })=m^{n+1}$.

\E\Tf we have shown that $n+1\in M$, hence $M=E^*$ and \er{3.38} holds. This
completes the proof of Theorem~\rf{t3.26}.
\endproof

\bex3.27
Let $A$ be a \nf\ and let $n:=\#(A)\in E^*$. Let $\cP(A)$ denote the set of all
subsets of~$A$ (including the empty set). Show that $\cP(A)$ is finite and that
\[
\#(\cP(A))=2^n
\]
where $2:=S(1)$.
\eex

\ti{Hint\/}: Find a bi\jn\ between $\cP(A)$ and the set of all maps from~$A$
into $\{0,1\}$.

\bex3.29
Let $(E,0,S)$ be a set of \nn s with $1:=S(0)$. \E\fe $m\in E$, let $S^m$
denote the $m$-fold \IT\ of the \su\ \f\ $S:E\to E$ \wrt $E$ introduced in
\E\df\ \rfa1{d2.1}, and let $I(S)$ denote the set of \IT s of~$S$. Prove:

\hph i,ii, $(I(S),\circ,\id_E)$, where $\circ$ denotes the \cm\ of self-maps
of~$E$, is an \am.

\hph ii,i, The map $\Phi:E\to I(S)$ defined by $\Phi(m):=S^m$, $m\in E$, is a
bi\jn.

\hph iii,, The binary \op\ $+'$ on $E$ defined by
\beq3.44
m +' n:=\Phi\Inv (S^m\circ S^n),\q m,n\in E,
\e
is equal to the binary \op\ $+$ introduced in \E\Pr\ \rf{p1.2}.

\hph iv,, $\Phi:(E,+,0)\to (I(S),\circ,\id_E)$ is a monoid-\is sm.
\eex

\bex3.30
Let $(E,0,S)$ be a set of \nn s with $1:=S(0)$. Let $\Psi:E\to \cP(E)$ be
defined by setting:
\beq3.45
\Psi(n):= \zo0,n , \q n\in E.
\e
Show:
\bea0.34
&\Psi(e)=\vn \qh{and } \Psi(n)\ne\vn \qh{\fa $n\in E\sms e$.}\\
&\Psi(S(e))=\{e\}. \lb{0.35} \\
&\Psi(S(n)) = \Psi(n)\cup \{n\} \qh{and } \Psi(n)\cap\{n\}=\vn \qh{\fa $n\in E$}.
\lb{0.36} \\
&\Psi:(E,\le) \to (\cP(E),\subset) \qh{is in\jc\ and in\cre.} \lb{0.37} \\
&\bigcup_{n\in E} \Psi(n)=E. \lb{0.38} \\
&\hbox{If $A$ is a finite set, then \te s \ooo $n\in E$ \st}  \lb{0.39} \\
&\quad A\ax\Psi(n), \hbox{ and} \nonumber\\
&\#_E(A)= n. \lb{0.40}
\e

\eex

\newpage

\Subsubsection{Composite sums and products}\label{sss.Comp}

In this section, as in Section \ref{sss.Mult}, we assume that $(E,0,S)$ is
a~set of \nn s with distinguished \el~$0$ whose \su\ is denoted by~$1$. \Mo
we use the notation $E^*:=E\sms0$, and the same notation as in Section
\ref{sss.Mult} for \ad, \mlc\ and \epc. Finally, we shall use the notation
$\#(A)$ instead of~$\#_E(A)$ for the \ca\ \wrt $(E,0,S)$ of a finite set~$A$.

As a motivation for the introduction of the notion of \cme sums and products,
we consider the \fw\ problem where these notions naturally occur. Let $A$
and~$B$ denote \nf s.
In view of Theorem \rf{t3.26}, the set of all maps from~$A$ into~$B$, $B^A$,
is finite and $\#(B^A)=\#(B)^{\#(A)}$. We denote by $\Bij(A,B)$ ($\Bij(A)$
if ${A=B}$) the set of all bi\jn s from $A$ onto~$B$. As a~subset of $B^A$,\glossary{$\Bij(A,B)$}\glossary{$\Bij(A)$}
the set $\Bij(A,B)$ is finite and
\beq4.1
\#(\Bij(A,B)) \le \#(B^A)
\e
by Lemma \rf{l3.16}. \Mo $\Bij(A,B)\ne\vn$ iff $A\ax B$ in view of \E\df\
\rf{d3.1}. Hence, if $A\not\ax B$, then $\#(\Bij(A,B))=0$ by \er{3.6}. We are
interested in $\#(\Bij(A,B))$ when $A\ax B$.

\blm4.1
Let $A,B$ be \nf s.

\hph i,ii, If $A\ax B$, then
\bea4.2
{}&\Bij(A)\ax \Bij(A,B),\\
&\Bij(A,B)\ax \Bij(B,A), \lb{4.3} \\
&\Bij(A) \ax \Bij(B), \lb{4.4} \\
&\#(\Bij(A)) = \#(\Bij([1,\#(A)])). \lb{4.5}
\e

\hph ii,i, If $A\sbs B$ and $A\ne B$, then
\beq4.6
\#(\Bij(A)) < \#(\Bij(B)).
\e

\hph iii,, If $\Bij(A)\ax \Bij(B)$, then $A\ax B$.
\elm

\proof \

(i) Since $A\ax B$, \te s $\vf\in\Bij(A,B)$.

\er{4.2}: If $f\in\Bij(A)$, then $\vf\circ f\in\Bij(A,B)$ since the \cm\ of
bi\jn s is a~bi\jn. Thus we may define a~map $\Phi:\Bij(A)\to\Bij(A,B)$ by
setting $\Phi(f):=\vf\circ f$.

\ti{In\ji\ of $\Phi$}: Let $f,g\in\Bij(A)$ be \st $\Phi(f)=\Phi(g)$.
Then $\vf(f(x))=\vf(g(x))$, $x\in A$, hence $f(x)=\vf\Inv(\vf(f(x)))=\vf\Inv
(\vf(g(x))=g(x)$, $x\in A$. Thus $f=g$, and $\Phi$ is in\jc.

\ti{Sur\ji\ of $\Phi$}: Let $h\in \Bij(A,B)$ and set $f:=\vf\Inv\circ h$. Then
$f\in\Bij(A)$ and $\Phi(f)=\vf\circ(\vf\Inv\circ h) = (\vf\circ\vf\Inv)\circ h
={\id_B}\circ h = h$.

\er{4.3}: Let $f\in\Bij(A,B)$, then $f\Inv\in\Bij(B,A)$. Define $\check\Phi:
\Bij(A,B) \to \Bij(B,A)$ by setting $\check\Phi(f):=f\Inv$. Similarly, if
$g\in\Bij(B,A)$, then $g\Inv\in\Bij(A,B)$. Then $\hat\Phi:\Bij(B,A)\to
\Bij(A,B)$ defined by $\hat\Phi(g):=g\Inv$, \sf ies $(\hat\Phi\circ
\check\Phi)(f)=\hat\Phi(\check\Phi(f))=\hat\Phi(f\Inv)=(f\Inv)\Inv=f$. Hence
$\hat\Phi\circ\check\Phi = \id_{\Bij(A,B)}$. Similarly, $\check\Phi\circ
\hat\Phi=\id_{\Bij(B,A)}$, hence $\check\Phi$ is a~bi\jn, and $\hat\Phi
=(\check\Phi)\Inv$.

\er{4.4}: By \er{3.2} we may interchange $A$ and~$B$. Hence $\Bij(A)\ax \Bij(B,A)$.
\er{4.2} and \er{3.3} imply $\Bij(B)\ax\Bij(B,A)$. Then \er{4.4} follows
from~\er{3.3} and \er{4.3}.


\er{4.5}:  From \er{3.26}
with $a:=0$, $b:=P(\#(A))$, and $c:=1$, we infer $[0,P(\#(A))]\ax
[1,P(\#(A))+1]$. But $P(\#(A))+1 = S(P(\#(A)))= (S\circ P)(\#(A))=\#(A)$.
Hence $[0,P(\#(A))]\ax[1,\#(A)]$. \E\Tf $A\nde3.8 \ax [0,P(\#(A))]
\ax [1,\#(A)]$, hence $A\ax[1,\#(A)]$ by \er{3.1}. Thus $\Bij(A)\ax \Bij([1,\#(A)])$
by \er{4.4}. Then \er{4.5} follows from Corollary \rf{c3.11}.

\ssk
(ii) \er{4.6}: We first construct an \ti{in\jc} map $\Phi:\Bij(A)\to\Bij(B)$.
Given $f\in\Bij(A)$, set $j(f)(x):=f(x)$ for $x\in A$ and
$j(f(x)):=x$ for $x\in B\sm A$. We claim that $j(f)\in \Bij(B)$. Indeed,
define $h:B\to B$ by setting $h(x):=f\Inv(x)$ for $x\in A$ and $h(x):=x$ for
$x\in B\sm A$. Then $(h\circ j(f))(x)= f\Inv\circ(f(x)) = (f\Inv\circ f)(x)
=\id_A(x)=x$ for $x\in A$, and $h\circ f(x)=x$ for $x\in B\sm A$. Hence
$h\circ j(f)=\id_B$. Similarly $j(f)\circ h=\id_B$, hence $j(f)\in \Bij(B)$
and $h=j(f)\Inv$. Define $\Phi:\Bij(A)\to\Bij(B)$ by setting $\Phi(f):=
j(f)$, $f\in \Bij(A)$. We now show that $\Phi$~is \ti{in\jc}. Indeed, let
$f,g\in\Bij(A)$ be \st $\Phi(f)=\Phi(g)$. Then $f(x)=\Phi(f)(x)=\Phi(g)(x)=
g(x)$, $x\in A$. Thus $f=g$. Hence $\Phi$~is in\jc. \If that $\Phi$~is a~bi\jn\
from $\Bij(A)$ onto $\Phi(\Bij(A))$. By \df\ $\Phi(\Bij(A))\sbs \Bij(B)$,
hence by \er{3.13} $\#(\Phi(\Bij(A)))\le \#(\Bij(B))$. We now show that
$\Phi(\Bij(A))\ne \Bij(B)$. Indeed, since $A\sbs B$,
$A\ne B$, \te\ $a\in A$ and $b\in B\sm A$. Define a map $\si_{a,b}\in\Bij(B)$
by setting $\si_{a,b}(a):=b$, $\si_{a,b}(b):=a$ and $\si_{a,b}(x)=x$ otherwise.
Since $\si_{a,b}\circ\si_{a,b}=\id_B$, $\si_{a,b}\in\Bij(B)$. But $\si_{a,b}
(b)=a\ne b$, hence $\si_{a,b}\notin\Phi(\Bij(A))$. Hence by \er{3.13} again
$\#(\Phi(\Bij(A))) < \#(\Bij(B))$. Finally, $\Phi(\Bij(A))\nad*= \#(\Phi(\Bij
A)))< \#(\Phi(\Bij(B)))$. Then \er{4.6} follows. In $\nad*=$ we used \er{4.2}
since $\Phi:\Bij(A)\to \Phi(\Bij(A))$ is a bi\jn.

\ssk
(iii) \If from $\Bij(A)\ax \Bij(B)$ that $\#(\Bij(A))=\#(\Bij(B))$ by Corollary
\rf{c3.11}. From \er{4.5} we obtain $\#(\Bij([1,\#(A)]))=\#(\Bij([1,\#(B)])$.
If $\#(A)<\#(B)$, then $[1,\#(A)]\sbs[1,\#(B)]$ and $[1,\#(A)]\ne[1,\#(B)]$
since $\#(B)\not\in[1,\#(A)]$. We get a \cd ion from \er{4.6}, \er{4.5}. \E\Tf $\#(A)
\not<\#(B)$. Exchanging $A$ and~$B$, we obtain $\#(B)\not<\#(A)$. Hence
$\#(A)=\#(B)$.
\endproof

\brm4.2
\If from the proof that \er{4.2}--\er{4.4} hold whenever $A,B$ are
\ti{infinite}.
\erm

\bex4.3
Does (iii) in Lemma \rf{l4.1} hold if $A,B$ are \ti{infinite}?
\eex

\E\et y \er{4.5} motivates

\bdn4.4
We shall denote by $\Pi_E$ (or simply by~$\Pi$, when no confusion arises)
the map $\Pi:E^*\to E^*$ defined by\glossary{$\Pi(n)$}
\beq4.7
\Pi(n):=\#(\Bij([1,n])), \q n\in E^*.
\e
\edn

\bpr4.5
The map $\Pi$ \sf ies
\bea4.8
{}&\Pi(1)=1, \\
&\Pi \hbox{ is strictly in\cre,} \lb{4.9} \\
&\Pi(n+1)=(n+1)\cdot \Pi(n), \q n\in E^*. \lb{4.10}
\e
\epr

\proof \

\er{4.8}: If $n=1$, then $[1,n]=\{1\}$ and $\Bij(\{1\})=\{\id_{\{1\}}\}$. Hence
$\#(\{\id_{\{1\}}\})=1$.

\er{4.9}: Let $n,m\in E^*$ with $n<m$. Then $[1,n]\sbs[1,m]$ and $[1,n]\ne
[1,m]$ since $m\notin[1,n]$. Then \er{4.9} follows from \er{4.6}.

\er{4.10}: As in the proof of \er{4.6}, we use the in\jc\ map $\Phi:\Bij([1,n])
\to\Bij([1,n+1])$ defined by $\Phi(f)(x):=f(x)$, $x\in[1,n]$, $\Phi(f)(n+1):=n+1$.
Note that the range of the map~$\Phi$ is $X_{n+1}:=\{h\in\Bij([1,n+1]):
h(n+1)=n+1\}$. Since $\Phi$ is in\jc, and $X_{n+1}\ax \Bij([1,n])$, we have
by Corollary \rf{c3.11} and \er{4.7}
\beq4.11
\#(X_{n+1})=\Pi(n).
\e
We define the \fw\ subsets of $\Bij([1,n+1])$:
\beq4.12
X_k:=\{h\in\Bij([1,n+1]): h(n+1)=k\}, \q k\in[1,n].
\e
We claim that $\{X_k\}_{k\in[1,n+1]}$ is a \ti{\pt} of $\Bij([1,n+1])$ (see
\E\df\ \rfa1{d4.5}). Indeed, we have $X_k\cap X_l=\vn$ whenever $k\ne l$,
$k,l\in[1,n+1]$, since there is no selfmap~$h$ of $[1,n+1]$ \sf ying $h(n+1)
=k$ and $h(n+1)=l$ for $k\ne l$, $k,l\in[1,n+1]$. \Mo if $h\in\Bij([0,n+1])$
then $h(n+1)=l$ \fs $l\in[1,n+1]$. \E\Tf $h\in X_l$, hence $\Bij([0,n+1])
\sbs \bcl_{k\in[0,n+1]}X_k\sbs \Bij([0,n+1])$. Finally we show that
$X_l\ne\vn$ \fe $l\in[1,n+1]$. Clearly $\id_{[1,n+1]}\in X_{n+1}$. Given
$u,v\in[1,n+1]$, we define a~map $\si_{u,v}:[1,n+1]\to[1,n+1]$ by setting
\beq4.13
\si_{u,v}(x):=\bca
v & \hbox{for }x=u,\\
u & \hbox{for }x=v,\\
x & \hbox{for }x\in [1,n+1]\sms{u,v} \hbox{ (if $\ne\vn$)}.
\eca
\e
One verifies that \fa $u,v\in[1,n+1]$:
\bga4.14
\si_{u,u}=\id_{[0,n+1]},\q \si_{u,v}=\si_{v,u},\\
\si_{u,v}\circ\si_{u,v}=\id_{[1,n+1]}. \lb{4.15}
\e
\If from \er{4.15} that $\si_{u,v}\in\Bij([1,n+1])$ \fa $u,v\in[1,n+1]$. The
map $\si_{u,v}$ is called a \ti{\tp} when $u\ne v$. Clearly $\si_{n+1,k}\in
X_k$ \fa $k\in[1,n+1]$, hence, \Ip $X_k\ne\vn$ \fa $k\in[1,n+1]$. This
completes the proof of the claim.

Next we show
\beq4.16
X_k \ax X_{n+1} \qh{\fa} k\in[1,n].
\e
To this end we construct bi\jn s from $X_{n+1}$ onto~$X_k$, $k\in[1,n]$.
Observe that if ${h\in X_{n+1}}$, then $\si_{n+1,k}\circ h\in X_k$ \fe
$k\in[1,n]$. Let $k\in[1,n]$. Then $\si_{n+1,k}\circ h\in\Bij([1,n+1])$ and
$(\si_{n+1,k}\circ h)(n+1) = \si_{n+1,k}(h(n+1)) = \si_{n+1,k}(n+1)=k$. Thus
$\si_{n+1,k}\circ h\in X_k$. We define $\Phi_k:X_{n+1}\to X_k$ by setting
$\Phi_k(h):=\si_{n+1,k}\circ h$, $h\in X_{n+1}$. We claim that $\Phi_k$ is
\ti{in\jc}. Indeed, if $h,\wt h\in X_{n+1}$ and $h\ne \wt h$, \te s $x\in
[1,n+1]$ \st $h(x)\ne\wt h(x)$ by in\ji\ of~$h$. Then $(\Phi_h(k))(x) =
(\si_{n+1,k}\circ h)(x)= \si_{n+1,k}(h(x))\ne \si_{n+1,k}(\wt h(x)) =
(\Phi_k(\wt h))(x)$ by in\ji\ of $\si_{n+1,k}$. \If that $\Phi_k(h)\ne
\Phi_k(\wt h)$. Hence $\Phi_k$ is in\jc. We now show that $\Phi_k$ is
\ti{sur\jc}. Let $g\in X_k$. Then $\si_{k,n+1}\circ g\in\Bij([1,n+1])$ and
$\si_{k,n+1}\circ g(n+1)=\si_{k,n+1}(k)=n+1$. Hence $\si_{k,n+1}\circ g\in
X_{n+1}$. \E\Tf $\Phi_k(\si_{k,n+1}\circ g) = \si_{n+1,k}\circ (\si_{k,n+1}
\circ g) = (\si_{n+1,k}\circ\si_{k,n+1})\circ g \nad{\er{4.14},\er{4.15}}=
{\id_{[1,n+1]}}\circ g=g$. \If that $\Phi_k$ is \ti{sur\jc}. \csq\ $\Phi_k$ is
bi\jc, hence $X_{n+1}\ax X_k$. From Corollary \rf{c3.11}, \er{4.11} and
\er{4.16}, we obtain
\beq4.17
\#(X_k)=\Pi(n) \qh{\fe} k\in[1,n+1].
\e
We complete the proof of \er{4.10} by means of the \fw\ lemma.

\blm4.6
Let $\O$ be a \nf, let $m\in E^*$ and let $\{B_k\}_{k\in[1,m]}$ be a~\pt\
of~$\O$ \st $\#(B_k)=n\in E$ \fs $n\in E^*$ and \fa $k\in[1,m]$. Then
\beq4.18
\#(\O)=mn.
\e
\elm

\proof[Proof of Lemma \rf{l4.6}]
We define a map $S:[1,m]\to \cP(\O)$, where $\cP(\O)$ denotes the power set
of~$\O$ (see Exercise \rfa1{ex3.4}\,(ii)) by setting
\beq4.19
S_k:=\bigcup_{l\in[1,k]}B_l,\qh{with} k\in[1,m].
\e
If $m=1$, then $\O=B_1=S_1$ and $\#(\O)=\#(B_1)=n\nde2.13 = 1\cdot n$, hence
\er{4.18} holds. We suppose $m\in E^*\sms1$. We claim that $S_k\cap B_{k+1}
=\vn$ whenever $1\le k<k+1\le m$. Indeed, suppose for \cd ion that \te\
$k\in[1,m]$ \st $k<k+1\le m$, and $\o\in S_k\cap B_{k+1}$. Then $\o\in B_l$
\fs $l\in[1,k]$ by \df\ of~$S_k$ and $\o\in B_{k+1}$. We have $l\ne k+1$,
since $l\le k<k+1$ by \er{3.11}. Hence $B_l\cap B_{k+1}=\vn$, since
$\{B_i\}_{i\in[1,m]}$ is a~\pt\ of~$\O$. This \cd s $\o\in S_k\cap B_{l+1}$.
\E\Tf the claim holds. \If that $S_{k+1}=S_k\cup B_{k+1}$ by \er{4.19} and
$S_k\cap B_{k+1}=\vn$, whenever $1\le k<k+1\le m$. In view of \er{3.31}, we obtain
\beq4.20
\#(S_{k+1})=\#(S_k)+\#(B_{k+1}) \qh{for} k\in[1,m).
\e
We now define a map $\vf:[0,m]\to E$ by setting $\vf(0):=0$ and $\vf(k):=
\#(S_k)$, $k\in[1,m]$. In view of \er{4.20}, we have $\vf(k+1)=\vf(k)+n=
n+\vf(k)$, $k\in[1,m)$. Since $\vf(0)=0$ and $\vf(1)=\#(S_1)=\#(B_1)=n$, we
also have $\vf(k+1)=\vf(k)+n=n+\vf(k)$ for $k=0$. \csq, we obtain
\beq4.21
\vf(k+1)=\vf(k)+n=n+\vf(k) \qh{for} k\in[0,m).
\e
Note that the map $\psi:[0,m]\to E$, defined by $\psi(k):=kn$, $k\in[0,m]$,
\sf ies
\beq4.22
\psi(k+1)=\psi(k)+n=n+\psi(k), \q k\in\zo 0,m ,
\e
in view of \er{2.15}, \er{2.13}.
It remains to show that $\vf=\psi$, since, in this case, $\#(\O)=\#(S_m)=
\vf(m)=\psi(m)=mn$. To this end we invoke Lemma \rfa1{l3.46} with $a:=0$,
$b:=m$ and $M:=\{k\in[0,m]: \vf(k)=\psi(k)\}$. We have $0\in M$, since
$\vf(0)=0\nad{\er{2.13}}= 0\cdot m=\psi(0)$. We suppose $k\in M\cap
\zo0,m $. Then $\vf(k+1)\nde4.21 = n+\vf(k)\nad{k\in M\cap\zo0,m }=n+\psi(k)
\nde4.22 = \psi(k+1)$. Hence $k+1\in M$.
\If from Lemma \rfa1{l3.46} that $M=[0,m]$. This completes the proof of
Lemma \rf{l4.6}.
\endproof

\Wanp complete the proof of \E\Pr\ \rf{p4.5}. We use Lemma \rf{l4.6} with
$\O:=\Bij([1,n+1])$, $m:=n+1$, $B_k:=X_k$, $k\in[1,n+1]$, and $n:=\Pi(n)$.
We obtain $\#(\Bij([1,n+1])=(n+1)\cdot\Pi(n)$. Hence $\Pi(n+1)=(n+1)
\cdot\Pi(n)$ by \er{4.7}.
\endproof

Our next goal is to show that \pp ies \er{4.8} together with \er{4.10}
\ch ize the map $\Pi:E^*\to E^*$. In other words, there is at most one
selfmap of~$E^*$ \sf ying \er{4.8} and \er{4.10}. We cannot use the \uq\
part of the \DRT\ (Theorem \rfa1{t1.6}) with $(E,e,S):=(E^*,1,S)$ where $S(x):=
x+1$, $x\in E^*$, since $\Pi(n+1)$ is not of the form $f(\Pi(n))$, $n\in E^*$,
where $f$~is a~selfmap of~$E^*$. Indeed, for each $n\in E^*$, we have a map
$f_{n+1}:E^*\to E^*$ \st $\Pi(n+1)=f_{n+1}(\Pi(n))$. The map $f_{n+1}$ is
defined by $f_{n+1}(x):=(n+1)\cdot x$, $x\in E^*$. However, we can invoke
the \fw\ \Pr\ to prove \uq.

\bpr7.1
Let $F$ be a \ns, let $n\in E$ and let $\{f_i\}_{i\in[1,n]}$ be maps
from $F$ into itself. Then \fe $a\in F$ \te s one and only one map
$\psi:[0,n]\to F$ \sf ying
\beq7.4
\psi(k+1)=f_{k+1}(\psi(k)) \qh{\fe $k\in[0,n)$}
\e
and
\beq7.5
\psi(0)=a.
\e
\epr

\proof
\ti{\E\uq}. Suppose that $\psi,\wt\psi:[0,n]\to  E$ \sf y \er{7.4}
and~\er{7.5}. Set $M:=\{p\in[0,n]:\psi(p)=\wt\psi(p)\}$. We want to show that
$M=[0,n]$. To this end we invoke Lemma~\rfa1{l3.46}. We have

$0\in M$, since $\psi(0)=a=\wt\psi(0)$. \Mo

$p\in M$ \ti{implies $p+1\in M$ whenever $p\in[0,n)$}: Let $p\in M\cap[0,n)$.
Then $\psi(p+1)\nad{\er{7.4}}=f_{p+1}(\psi(p))\nad{p \in M}=
f_{p+1}(\wt\psi(p))\nad{\er{7.4}}= \wt\psi(p+1)$.
Hence $p+1\in M$ and $M=[0,n]$.

\ti{\E\ex}. Set $M:=\{p\in[0,n]: \hbox{\te s a map }\psi_p:[0,p]\to  E$ \st
$\psi_p(0)=a$ and $\psi_p(k+1)=f_{k+1}(\psi_p(k))$ \fe $k\in[0,p)\}$. We want
to show that $M=[0,n]$. As above, we invoke Lemma~\rfa1{l3.46}. We have

$0\in M$ since $\psi_0:\{0\}\to E$ defined by $\psi_0(0)=a$ \sf ies
\cn s \er{7.4}, \er{7.5}, since $[0,0)$ is empty.

\ti{$p\in M\cap[0,n)$ implies $p+1\in M$}: Let $p\in M\cap [0,n)$ and let $\psi_p
:[0,p]\to E$ \sf y \cn s \er{7.5}, \er{7.4} with $n:=p$. Set
\beq7.6
\psi_{p+1}(k):=\bca
\psi_p(k) &\hbox{for }k\in[0,p],\\
f_{p+1}(\psi_p(p)) &\hbox{for }k=p+1.
\eca
\e
Then $\psi_{p+1}:[0,p+1]\to E$ \sf ies
\[
\bal
\psi_{p+1}(0)&\nad{\er{7.6}}=\psi_p(0)\nad{p\in M}=a,\\
\psi_{p+1}(k+1)&\nad{\er{7.6}}=\psi_p(k+1)\nad{p\in M}=f_{k+1}(\psi_p(k))
\nad{\er{7.6}}=f_{k+1}(\psi_{p+1}(k)) \qh{for }k\in[0,p),\\
\psi_{p+1}(p+1)&\nad{\er{7.6}}=f_{p+1}(\psi_p(p))\nad{\er{7.6}}=
f_{p+1}(\psi_{p+1}(p)).
\eal
\]
It follows that $p+1\in M$ and $M=[0,n]$. Then $\psi:=\psi_n$ \sf ies
\er{7.4}, \er{7.5}.
\endproof

The next theorem extends both \DRT\ and \E\Pr~\rf{p7.1}.

\bth7.2
Let $F$ be a \ns\ and let $\{f_i\}_{i\in E^*}$ be a \sq\ of maps from $F$ into
itself. Then \fe $a\in F$ \te s one and only one map $\psi: E\to F$ \sf ying
\er{7.5} and
\beq7.7
\psi(k+1)=f_{k+1}(\psi(k)) \qh{\fe }k\in E.
\e
\Mo if \te\ $n\in E$ and $\wt\psi:[0,n]\to E$ \st $\wt\psi(0)=a$ and
$\wt\psi(k+1)=f_{k+1}(\wt\psi(k))$ \fe $k\in[0,n)$, then $\psi(k)=\wt
\psi(k)$ \fe $k\in[0,n]$.
\eth

\brm7.3
It is obvious that Theorem \rf{t7.2} extends Theorem \rfa1{t1.6}. Concerning
\E\Pr~\rf{p7.1}, note that given $f_i:F\to F$, $i\in[1,n]$, one can always
define $f_i:=\id_F$ for $i>n$. Let $\psi: E\to F$ be the map whose \ex\ and
\uq\ is guaranteed by Theorem~\rf{t7.2}. Then the \rt ion of~$\psi$
to~$[0,n]$ \sf ies \er{7.4}--\er{7.5}, which proves the \ex\ part of
\E\Pr~\rf{p7.1}. The \uq\ part obviously follows from the \uq\ part of
\E\Pr~\rf{p7.1}.
\erm

\proof[Proof of Theorem \rf{t7.2}]
\E\fe $n\in E$ let $\psi_n$ denote the unique map from $[0,n]$ to~$F$ \sf
ying \er{7.4}--\er{7.5} whose \ex\ and \uq\ is guaranteed by \E\Pr\ \rf{p7.1}.
Let $n\in E^*$ and $m\in[0,n)$. Then $\psi_n\bigr|_{[0,m]}$, the \rt ion
of~$\psi_n$ to $[0,m]$, \sf ies \er{7.5} and \er{7.4} where $n$ is replaced
by~$m$. In view of the \uq\ part of \E\Pr~\rf{p7.1},
$\psi_n\bigr|_{[0,m]}=\psi_m$. Define a~map $\psi:E\to F$ by setting $\psi(p)
:=\psi_p(p)$, $p\in E$. Note that $\psi_0(0)=a$, hence $\psi$ \sf ies
\er{7.5}. We show that it \sf ies \er{7.7}. Let $k\in E$, then $\psi(k+1)
:=\psi_{k+1}(k+1)\nde7.4 = f_{k+1}(\psi_{k+1}(k))\nad*= f_{k+1}(\psi_k(k))
=f_{k+1}(\psi(k))$. In $\nad*=$ we used $\psi_{k+1}|_{[0,k]}=\psi_k$.
\E\Tf $\psi$ \sf ies \er{7.5} and \er{7.7}.

\Mo on every \il\ $[0,n]$, $n\in E$, the \rt ion of~$\psi$ to~$[0,n]$ is
uniquely defined in view of \E\Pr~\rf{p7.1}. Hence $\psi: E\to F$ is uniquely
defined.
\endproof

We are now in a position to define \ti{\cme sums} in $( E,+,0)$. Since we shall
also need \ti{\cme products}, we consider a more general situation.

Let $X$ be a \ns\ and let ${\qu}:X\t X\to X$ be a binary \op.

We want to define
expressions of the form $((a_0\qu a_1)\qu a_{1+1})\qu \dots$ and of the form
$\ldots \qu(a_{1+1}\qu (a_1\qu a_0))$. Note that we can reduce the second
case to the first one by defining ${\tqu}:X\t X\to X$ as follows:
$x\tqu y:=y\qu x$, $x,y\in X$. So we \rt\ ourselves to the first case. Given
$n\in E$ and a finite \sq\ of \el s of~$X$, $a:[0,n]\to X$, we define \rc
vely $\psi:[0,n]\to X$ by setting
\bea7.8
\psi(0)&:=a_0\\
\psi(k+1)&:=f_{k+1}(\psi(k)),\ k\in[0,n),
\hbox{ where }f_i(x):=x\qu a_i, \ i\in[1,n],\ x\in X. \lb{7.9}
\e
Clearly $f_i$ is a selfmap of $X$, $i\in[1,n]$.

\def\namespec{Definition and Notation}
\begin{dspc}\lb{d7.4}
Let $X$ be a \ns\ and let ${\qu}:X\t X\to X$ be a binary \op. Let $n\in E$
and $a:[0,n]\to X$. Let $\psi:[0,n]\to X$ be the unique map \sf ying
\er{7.8}--\er{7.9}, whose \ex\ and \uq\ is guaranteed by \E\Pr~\rf{p7.1}. Set\index{composite!operation}\glossary{$\QU_{k=0}^p a_k$}
\beq7.10
\QU_{k=0}^p a_k:=\bca
a_0 &\hbox{for }p=0,\\
\psi(p) &\hbox{for }p\in[1,n].
\eca
\e
Note that if $p+1\le n$, then
$$
\QU_{k=0}^{p+1} a_k\nad{\er{7.10},\er{7.9}}= \Bigl(\QU_{k=0}^p a_k\Bigr)\qu a_{p+1}.
$$
The letter $k$ can be replaced by any other letter distinct from $p$ and~$a$.
It is called a \ti{dummy index}. So we have
\beq7.11
\QU_{k=0}^p a_k=\QU_{l=0}^p a_l, \q \hbox{\fe} p\in[0,n].
\e
We also introduce the \fw\ \gn. Let $b:[m,m+n]\to X$ be given where
$m,n\in  E$. Then
\beq7.12
\QU_{l=m}^{m+p} b_l:=\QU_{k=0}^p a_k \hbox{ where }a_k:=b_{k+m} \q
\hbox{\fa $k\in[0,p]$ and \fe} p\in[0,n].
\e
\end{dspc}
Observe that if $m:=0$ in \er{7.12}, then $\QU_{l=m}^{m+p} b_l \nde7.12 =
\QU_{k=0}^p b_k \nde7.11 = \QU_{l=0}^p b_l$. Hence \df\ \er{7.12} is
consistent with \df\ \er{7.10}.

In the next \Pr\ we prove \pp ies of the \cme \op\ $\QU_{k=0}^p$ which are
inherited from \pp ies of the binary \op~$\qu$.

\begin{prp}[see \cite{Alg}]\lb{p7.5}
Let $X$ be a \ns\ and let ${\qu}:X\t X\to X$ be a binary \op\ and let
$\QU_{k=0}^p$, $\QU_{l=m}^{m+p}$ be the \cme \op s defined in \er{7.10}, \er{7.12}.

\hph i,ii, Suppose the law $\qu$ is \ti{\asc e}. Let $a:[0,m+n]\to X$ where
$m\in E$, $n\in E^*$. Then
\beq7.13
\Bigl(\QU_{i=0}^m a_i\Bigr) \qu \Bigl(\QU_{j=m+1}^{m+n} a_j\Bigr)
=\QU_{k=0}^{m+n} a_k.
\e

\hph ii,i, Suppose \ti{in \ad} the law $\qu$ is \ti{\cmt e}. Let
$a:[0,m]\to X$ and let $\pi:[0,m]\to [0,m]$ be a bi\jn. Then
\beq7.14
\QU_{k=0}^m a_{\pi(k)}=\QU_{k=0}^m a_k.
\e

\hph iii,, Suppose in \ad\ \te s a \ti{\nel} $e\in X$ $($i.e.\ $e\qu x=x\qu e
=x$ \fa $x\in X)$ \st $x\qu y=e$ implies $x=y=e$
\fe $x,y\in X$. Let $a:[0,m]\to X$, where $m\in E$. Then
\beq7.15
\QU_{k=0}^m a_k=e \qhq{implies} a_k=e \hbox{ \fe }k\in[0,m].
\e
\epr

\proof
(i) By induction on $n\in E^*$.

\er{7.13} \ti{holds for} $n=1$: Note that $\QU_{j=m+1}^{m+1}a_j
\nad{\er{7.12}}=\QU_{k=0}^0 a_{k+m+1}\nde7.10 = a_{m+1}$. Hence
$\Bigl(\QU_{i=0}^m a_i\Bigr)\qu \Bigl(\QU_{j=m+1}^{m+1}a_j\Bigr)=
\Bigl(\QU_{i=0}^m a_i\Bigr)\qu a_{m+1}\nde7.10 = \QU_{i=0}^{m+1} a_i$.

\ti{If \er{7.13} holds for $n$, then \er{7.13} holds for $n+1$}: Suppose \er{7.13}
holds for~$n$. Then
\bmlg
\Bigl(\QU_{i=0}^m a_i\Bigr)\qu \Bigl(\QU_{j=m+1}^{m+n+1} a_j\Bigr)
\nde7.10 =\Bigl(\QU_{i=0}^m a_i\Bigr)\qu \Bigl(\Bigl(\QU_{j=m+1}^{m+n}a_j\Bigr)
\qu a_{m+n+1}\Bigr)\\
{}\nad\ast=\Bigl(\Bigl(\QU_{i=0}^m a_i\Bigr)\qu
\Bigl(\QU_{j=m+1}^{m+n}a_j\Bigr)\Bigr)\qu a_{m+n+1}
\nad{**} = \Bigl(\QU_{j=0}^{m+n}a_j\Bigr)\qu a_{m+n+1}
\nde7.10 = \QU_{j=0}^{m+n+1}a_j.
\e
In $\nad\ast=$ the \asc ity of $\qu$ is used, and in $\nad{**}=$ the \In\
hypothesis is used.

\medskip
(ii) By induction on $m\in E$.

\er{7.14} \ti{holds for} $m=0$: Note that $[0,0]=\{0\}$, hence $\pi(0)=0$.

\ti{If \er{7.14} holds for $m$, then \er{7.14} holds for $m+1$}: Suppose \er{7.14}
holds for~$m$. Let $\pi:[0,m+1]\to[0,m+1]$ be bi\jc. Since $\pi$ is bi\jc,
\te s exactly one $l\in[0,m+1]$ \st $\pi(l)=m+1$.

\ti{Case 1}. $\pi(m+1)=m+1$. Then
\bmlg
\QU_{k=0}^{m+1} a_{\pi(k)}\nde7.10 = \Bigl(\QU_{k=0}^{m} a_{\pi(k)}\Bigr)\qu
a_{\pi(m+1)}\\ \nad*= \Bigl(\QU_{k=0}^m a_k\Bigr)\qu a_{\pi(m+1)} =
\Bigl(\QU_{k=0}^{m} a_{k}\Bigr)\qu a_{m+1}
\nde7.10 = \QU_{k=0}^{m+1} a_{k}.
\e
In $\nad*=$ the \In\ hypothesis is used. Thus \er{7.14} holds for $m:=m+1$.

\ti{Case 2}. $\pi(0)=m+1$. Then $\QU_{k=0}^{m+1} a_{\pi(k)}\nde7.13 =
a_{\pi(0)} \qu \Bigl(\QU_{k=1}^{m+1} a_{\pi(k)}\Bigr)=a_{m+1}\qu
\Bigl(\QU_{k=1}^{m+1}a_{\pi(k)}\Bigr) \nad*= \Bg(\QU_{k=1}^{m+1} a_{\pi(k)})
\qu a_{m+1}\nde7.12 = \Bg(\QU_{i=0}^m a_{\pi(i+1)})\qu a_{m+1}$.
Define $\rho:[0,m]\to[0,m+1]$ by setting $\rho(i):=\pi(i+1)$, $i\in[0,m]$.
As a \cm\ of the in\jc\ map $i\mt i+1=S(i)$, $i\in[0,m]$, and of the in\jc\
map~$\pi$, $\rho$~is in\jc. \Mo there is no $i\in[0,m]$ \st $\rho(i)=m+1$,
since $\pi(i+1)=m+1=\pi(0)$ implies $i+1=0$ by the in\ji\ of~$\pi$, and
$i+1=0$ leads to a \cd ion. \If that $\rho$~is an in\jc\ selfmap of $[0,m]$,
hence a bi\jn\ by \er{3.11}. From the \In\ hypothesis we infer $\QU_{i=0}^m
a_{\rho(i)} = \QU_{i=0}^m a_i$. Hence $\QU_{k=0}^{m+1}a_{\pi(k)} =\Bg(
\QU_{i=0}^m a_{\pi(i+1)}) \qu a_{m+1} = \Bg(\QU_{i=0}^m a_{\rho(i)})\qu
a_{m+1} = \Bg(\QU_{i=0}^m a_i)\qu a_{m+1} \nde7.12 = \QU_{i=0}^{m+1}a_i$. Thus
\er{7.14} holds for $m:=m+1$.

\ti{Case 3}. $\pi(l)=m+1$ where $0<l<m+1$.
\bmlg
\QU_{k=0}^{m+1} a_{\pi(k)}\nde7.13 = \Bigl(\QU_{k=0}^{l} a_{\pi(k)}\Bigr)
\qu \Bigl(\QU_{j=l+1}^{m+1} a_{\pi(j)}\Bigr)
\nde7.10 = \Bigl(\Bigl(\QU_{k=0}^{l-1} a_{\pi(k)}\Bigr)\qu a_{\pi(l)}\Bigr)
\qu\Bigl(\QU_{j=l+1}^{m+1} a_{\pi(j)}\Bigr)\\
{}\nad\ast=\Bigl(\QU_{k=0}^{l-1} a_{\pi(k)}\Bigr)\qu
\Bigl(a_{\pi(l)}\qu \QU_{j=l+1}^{m+1} a_{\pi(j)}\Bigr)
\nad{\ast\ast}= \Bigl(\QU_{k=0}^{l-1} a_{\pi(k)}\Bigr)\qu
\Bigl(\Bigl(\QU_{j=l+1}^{m+1} a_{\pi(j)}\Bigr)\qu a_{\pi(l)}\Bigr)\\
{}=\Bigl(\QU_{k=0}^{l-1} a_{\pi(k)}\Bigr) \qu
\Bigl(\Bigl(\QU_{j=l+1}^{m+1} a_{\pi(j)}\Bigr)\qu a_{m+1}\Bigr)
\nad\ast= \Bigl(\Bigl(\QU_{k=0}^{l-1} a_{\pi(k)}\Bigr)\qu
\Bigl(\QU_{j=l+1}^{m+1} a_{\pi(j)}\Bigr)\Bigr)\qu a_{m+1}\\
{}\nde7.12 =
\Bigl(\Bigl(\QU_{k=0}^{l-1} a_{\pi(k)}\Bigr)\qu \Bigl(\QU_{j=l}^{m} a_{\pi(j+1)}
\Bigr)\Bigr)\qu a_{m+1}.
\e
In $\nad\ast=$ (resp.\ $\nad{\ast\ast}=$) the \asc ity (resp.\ \cmt ity)
of~$\qu$ is used.

Set
\[
\wt\pi(k):=\bca
\pi(k) &\hbox{for }k\in [0,l-1],\\
\pi(k+1)&\hbox{for }k\in[l,m].
\eca
\]
We claim that $\wt\pi:[0,m]\to[0,m+1]$ is a sur\jc\ selfmap of $[0,m]$.
Indeed, there is no $k\in[0,m]$ \st $\wt\pi(k)=m+1$, since either $k\in[0,l-1]$
and $\wt\pi(k)=\pi(k)\ne\pi(l)=m+1$, or $k\in[l,m]$ and $\wt\pi(k)=\pi(k+1)
\ne\pi(l)=m+1$, in view of the in\ji\ of~$\pi$. \csq, $\wt\pi$ is a selfmap
of $[0,m]$. We now show that $\wt\pi([0,m])=[0,m]$. Let $i\in[0,m]$, then \te
s $k\in[0,m+1]$ \st $\pi(k)=i$, since $\pi$~is sur\jc. Since $i\ne m+1=\pi(l)$,
we have $k\in[0,l-1]\cup[l+1,m+1]$ by the in\ji\ of~$\pi$. If $k\in[0,l-1]$,
then $\wt\pi(k)=\pi(k)=i$. If $k\in[l+1,m+1]$, then $\wt\pi(k-1)=\pi(k)=i$,
where $k-1\in[l,m]\sbs[0,m]$. This completes the proof of the claim. Then
$\wt\pi$ is bi\jc\ by \er{3.12}. Thus
\bmlg
\Bigl(\Bigl(\QU_{k=0}^{l-1} a_{\pi(k)}\Bigr)\qu\Bigl(\QU_{j=l}^{m}
a_{\pi(j+1)}\Bigr)\Bigr) \qu a_{m+1}
= \Bigl(\Bigl(\QU_{k=0}^{l-1} a_{\tilde\pi(k)}\Bigr)\qu
\Bigl(\QU_{j=l}^{m} a_{\tilde\pi(j)}\Bigr)\Bigr)\qu a_{m+1}\\
{}\nde7.13 = \Bigl(\QU_{k=0}^{m} a_{\tilde\pi(k)}\Bigr)\qu a_{m+1}
\nad*= \Bigl(\QU_{k=0}^{m} a_{k}\Bigr)\qu a_{m+1}
\nde7.10 = \QU_{k=0}^{m+1} a_{k}.
\e
In $\nad*=$ the \In\ hypothesis is used. Thus \er{7.14} holds for $m:=m+1$.

\smallskip
(iii) We proceed by induction on $m\in E$.

$m=0$: $a_0\nde7.10 = \QU_{k=0}^0 a_k=e$.

$m$ \ti{implies} $m+1$: $e=\QU_{k=0}^{m+1} a_k\nde7.10 = \Bigl(\QU_{k=0}^m
a_k\Bigr) \qu a_{m+1}$. Hence $\QU_{k=0}^m a_k=e$ and $a_{m+1}=e$. The
induction hypothesis implies $a_k=e$ for $k\in[0,m]$.
\endproof

\brm7.6
In \er{7.10}, \er{7.12}, the dummy index takes values in a nonempty finite \il\
of~$ E$. When the law $\qu$ is \asc e and \cmt e the \cme \op\
$\QU_{k=0}^n a_k$ does not depend on the order of the operands
$a_0,a_1,\dots$. \E\Pr\ \rf{p7.5}\,(ii) allows us to consider a more general
situation where the \il\ $[0,n]$ is replaced by an arbitrary \ti{nonempty finite}
set. Indeed, if $A$ is a nonempty finite set and $a$ is a map from $A$
into~$X$, we may consider the \cme \op\ $\QU_{j=0}^{\#(A)-1} a_{\vf(j)}$
where $\vf$ is a bi\jn\ from $[0,\#(A))$ onto~$A$ and $\#(A)-1:=P(\#(A))$
(see Notation \rf{n3.7}\,(ii)). Note that
$[0,\#(A))=[0,S(P(\#(A)))]\nde3.4 = [0,P(\#(A)] =[0,\#(A)-1]$.
It follows from \E\Pr\ \rf{p7.5}\,(ii) that this \cme
\op\ is \ti{independent\/} of the bi\jn~$\vf$, that is, if
$\psi:[0,\#(A)-1]\to X$ is a bi\jn\ then
\beq7.16
\QU_{j=0}^{\#(A)-1} a_{\vf(j)}=\QU_{j=0}^{\#(A)-1} a_{\psi(j)}.
\e
Indeed, if we set $b_j:=a_{\vf(j)}$, $j\in[0,\#(A)-1]$ and $\phi:=\vf\Inv
\circ \psi:[0,\#(A)-1]\to [0,\#(A)-1]$, we have $a_{\psi(l)}=a_{\vf\circ
\phi(l)} =b_{\phi(l)}$, $l\in[0,\#(A)-1]$. \E\Tf since $\phi$ is a bi\jn\
from $[0,\#(A)-1]$ into itself we obtain
\[
\QU_{j=0}^{\#(A)-1} a_{\vf(j)}=\QU_{j=0}^{\#(A)-1} b_j\nde7.14 =
\QU_{j=0}^{\#(A)-1} b_{\phi(j)}=
\QU_{j=0}^{\#(A)-1} a_{\psi(j)},
\]
which proves \er{7.16}.
\erm

This observation justifies the \fw\ \df.

\bdf7.7
Let $(X,\qu,e)$ be an abelian monoid, let $\O$ be a \ns, let $a$ be a map
from $\O$ into~$X$ and let $A$ be a \ti{finite} subset of~$\O$. Set\glossary{$\QU_{i\in A} a_i$}
\beq7.17
\QU_{i\in A}a_i:=\bca
e &\hbox{if }A=\vn,\\
\QU_{j=0}^{\#(A)-1}a_{\vf(j)} & \hbox{if }A\ne\vn \hbox{ and $\vf$ is a
bi\jn\ from $[0,\#(A){-}1]$ onto }A.
\eca \kern-16pt
\e
\edf

We now consider the case where $X$ is the abelian monoid $(E,+,0)$.
In that case the symbol $\QU$ is usually replaced by~$\sum$:\glossary{$\sum_{i=0}^m a_i$}
\beq7.25
\sum_{i=0}^m a_i :=\QU_{i=0}^m a_i, \q
\sum_{i\in A} a_i := \QU_{i\in A} a_i.
\e

In the next \Pr\ we collect some results which will be used in the next
section.

\bpr7.9
Let $m,n,p\in E$, $a,b:[0,n]\to E$ and let $A$ be a nonempty finite set.
Then the \fw\ holds\dw
\bea7.26
{}&\sum_{k=0}^n (a_k+b_k) = \sum_{k=0}^n a_k+ \sum_{k=0}^n b_k,\\
&\sum_{k=0}^n m\cdot a_k = m\cdot \sum_{k=0}^n a_k, \lb{7.27}\\
&\sum_{k=0}^n 1= n+1, \q \sum_{k=1}^n 1 =n \qh{and} \sum_{k\in A}1=\#(A). \lb{7.28}
\e
\Mo if $a_k\le b_k$ \fe $k\in[0,n]$, then
\beq7.29
\sum_{k=0}^n a_k \le \sum_{k=0}^n b_k.
\e
If in \ad\ $a_k<b_k$ for some $k\in[0,n]$ then
\beq7.30
\sum_{k=0}^n a_k< \sum_{k=0}^n b_k.
\e
Finally, we have
\beq7.31
(1+p)^{n+1}=1+p\sum_{k=0}^n(1+p)^k.
\e
\epr

\proof

\er{7.26}: By induction on $n$.

$n=0$: $\suml_{k=0}^0(a_k+b_k)\nde7.10 = a_0+b_0\nde7.10 =\suml_{k=0}^0 a_k
+\suml_{k=0}^0 b_k$.

\ti{$n$ implies $n+1$}: Suppose \er{7.26} holds. Then
\bmlg
\sum_{k=0}^{n+1}(a_k+b_k)\nde7.10 = \Bigl(\sum_{k=0}^n(a_k+b_k)\Bigr)+(a_{n+1}+b_{n+1})
\nde7.26 = \Bigl(\sum_{k=0}^n a_k+\sum_{k=0}^n b_k\Bigr)+(a_{n+1}+b_{n+1})\\
{}\nde1.36 = \sum_{k=0}^n a_k+\Bigl(\sum_{k=0}^n b_k+a_{n+1}\Bigr)+b_{n+1}
\nde 1.6 = \sum_{k=0}^n a_k+\Bigl(a_{n+1}+\sum_{k=0}^n b_k\Bigr)+b_{n+1}\\
{}\nde{1.36} = \Bigl(\Bigl(\sum_{k=0}^n a_k\Bigr)+a_{n+1}\Bigr)
+\Bigl(\Bigl(\sum_{k=0}^n b_k\Bigr)+b_{n+1}\Bigr)
\nde7.10 = \sum_{k=0}^{n+1}a_k+\sum_{k=0}^{n+1}b_k.
\e

\er{7.27}: 
By \In\ on $n$.

$n=0$: $\suml_{k=0}^0 m\cdot a_k\nde7.10 = m\cdot a_0 = m\cdot \suml_{k=0}^0
a_k$.

\ti{$n$ implies $n+1$}: Suppose \er{7.27} holds. Then $\suml_{k=0}^{n+1}
m\cdot a_k \nde7.10 = \Bigl(\suml_{k=0}^n m\cdot a_k\Bigr)+m\cdot a_{n+1}
\nde7.27 = \Bigl(m\cdot \suml_{k=0}^n a_k\Bigr)+m\cdot a_{n+1}\nde2.14 =
m\cdot\Bigl(\Bigl(\suml_{k=0}^n a_k\Bigr)+a_{n+1}\Bigr)\nde7.10 =
m\cdot\suml_{k=0}^{n+1}a_k$.

\er{7.28}: Observe that $\suml_{k=0}^n 1:=\suml_{k=0}^n a_k$ \fa $n\in E$ where $a:E\to E$
\sf ies $a_k:=1$ \fa $k\in E$. Set $\vf(n):=\suml_{k=0}^n a_k$, $n\in E$.
Clearly $\vf(0)=1$ since $\suml_{k=0}^0 a_k\nde7.10 = a_0=1$. \Mo $\vf(n+1)=
\vf(n)+1$, $n\in E$. Indeed, $\vf(n+1)=\suml_{k=0}^{n+1}a_k\nde7.10 = \bigl(
\suml_{k=0}^n a_k\bigr)+a_{n+1} = \vf(n)+1$, $n\in E$. Define $f:E\to E$ by
setting $f(x)=x+1$, $x\in E$. Then $\vf(n+1)=f(\vf(n))$, $n\in E$. \E\oh
$\psi:E\to E$ defined by $\psi(n):=n+1$, $n\in E$, \sf ies $\psi(0)=1$ and
$\psi(n+1)=(n+1)+1=f(\psi(n))$, $n\in E$. \If from the \uq\ part of Theorem
\rfa1{t1.6} that $\vf=\psi$. Hence $\vf(n)=n+1$, $n\in E$. \Mo
$n+1=\suml_{k=0}^n1\nde7.13 = \Bigl(\suml_{k=0}^0 1\Bigr)+\suml_{k=1}^n 1
\nde7.10 = 1+\suml_{k=1}^n 1\nde1.6 = \Bigl(\suml_{k=1}^n 1\Bigr)+1$, $n\in E$.
Hence, $n=\suml_{k=1}^n1$, $n\in E$, by \er{1.8}.

Finally, $\suml_{k\in A}1 :=\suml_{k\in A}a_k \nde7.17 = \suml_{k=0}^{\#(A)-1}
a_{\vf(k)}$ where $\vf$ is a bi\jn\ from $[0,\#(A))$ onto~$A$ since $A\ne
\vn$. But $a_{\vf(k)}=1$ for $k\in[0,\#(A)-1]$, hence $\suml_{k\in A}1=
\suml_{k=0}^{\#(A)-1}1=(\#(A)-1)+1$ by what precedes. Since $\#(A)-1=P(\#
(A))$, we have $(\#(A)-1)+1 = S(P(\#(A)))=\#(A)$. \csq\ $\suml_{k\in A}1
=\#(A)$.

\er{7.29}, \er{7.30}: In view of \er{1.46} \te\ $c_k\in E$ \st $b_k=a_k+c_k$
\fe $k\in[0,n]$. Thus $\suml_{k=0}^n b_k = \suml_{k=0}^n (a_k+c_k)\nde7.26 =
\suml_{k=0}^n a_k+\suml_{k=0}^n c_k$, and \er{7.29} follows from \er{1.46}.
If \te s $\ov k\in[0,n]$ \st $c_{\bar k}\in E^*$, then $\suml_{k=0}^n c_k
\in E^*$, and \er{7.30} follows from \er{1.46}. Suppose for \cd ion that
$\suml_{k=0}^n c_k=0$, then $c_k=0$ \fa $k\in[0,n]$ by \er{7.15}. A~\cd ion.

\er{7.31}: By induction on $n\in E$.

$n=0$: $(1+p)\nde7.10 = 1+p(1+p)^0 \nde2.23 = 1+p\cdot1\nde2.13 = 1+p$.

\ti{$p$ implies $p+1$}: Suppose \er{7.31} holds. Then
\bmlg
(1+p)^{(n+1)+1}\nad{\er{2.25},\er{2.24}}= (1+p)^{n+1}(1+p)\nde 7.31 =
\Bigl(1+p\sum_{k=0}^n (1+p)^k\Bigr)(1+p)\\
{}\nad{\er{2.15},\er{2.13},\er{2.12},\er{2.11}}= (1+p)+p(1+p)\sum_{k=0}^n (1+p)^k
\nad{\er{7.27},\er{2.11}} = (1+p)+p\sum_{k=0}^n \bigl((1+p)(1+p)^k\bigr)\\
{}\nad{\er{2.25},\er{2.24}}=(1+p)+p\sum_{k=0}^n (1+p)^{k+1}
\nde7.12 = (1+p)+p\sum_{k=1}^{n+1}(1+p)^k \\
{}\nad{\er{2.23},\er{2.12},\er{2.13}} = (1+p(1+p)^0)
+p\sum_{k=1}^{n+1}(1+p)^k \\
{}\nad{\er{2.11},\er{2.14}} = 1+p\Bigl((1+p)^0+\sum_{k=1}^{n+1}(1+p)^k\Bigr)
\nad{\er{7.13},\er{7.10}} = 1+p\sum_{k=0}^{n+1}(1+p)^k.\hbox
to8pt{\hskip10pt $\Box$\hss}
\e

We now give a \gn\ of both Theorem \rf{t3.23} \er{3.31} and Lemma \rf{l4.6}.

\bpr4.15
Let $\O$ be a \nf. Let $m\in E^*$ and let $\{B_k\}_{k\in[1,m]}$ be a family of
subsets of~$\O$ \sf ying $\bcl_{k\in[1,m]}B_k=\O$, and $B_k\cap B_l=\vn$
whenever $k\ne l$. Then
\beq4.44
\#(\O)=\sum_{k=1}^m \#(B_k).
\e
\epr

\proof
We define $\{S_k\}_{k\in[1,m]}$ as in the proof of Lemma \rf{l4.6} and obtain
\er{4.20}. We also define $\psi:[0,m]\to E$ by setting $\psi(0)=0$ and $\psi(k)
:=\#(S_k)$, $k\in[1,m]$. Observe that $\#(\O)=\#(S_m)=\psi(m)$. \Mo $\psi(k+1)
\nde4.20 = \psi(k)+\#(B_{k+1})$ for $k\in\zo1,m $. For $k=0$, we have $\psi(1)
=\#(S_1)=\#(B_1)$. Set $b_k:=\#(S_k)$, $k\in[1,m]$, $b_0:=0$. Then
\beq4.45
\psi(0)=b_0 \qh{and } \psi(k+1)=\psi(k)+b_{k+1} \hbox{ for }k\in\zo0,m .
\e
Thus $\psi$ \sf ies \er{7.8}, \er{7.9} with $n:=m$, $a_i:=b_i$, $i\in[1,m]$,
and $f_i(x):=x+b_i$, $i\in[1,m]$. \E\Tf $\psi(m)=\suml_{k=0}^m b_k$ in view
of~\er{7.10}. Using the \asc ity of the \ad\ we obtain $\suml_{k=0}^m b_k\nde7.13 =
\suml_{k=0}^0 b_k+\suml_{k=1}^m b_k\nde7.10 = 0+\suml_{k=1}^m b_k\nde1.7 =
\suml_{k=1}^m b_k$. Hence $\#(\O)=\#(S_m)=\psi(m)=
\suml_{k=1}^m b_k = \suml_{k=1}^m \#(B_k)$.
\endproof

\bex4.16
Prove \er{4.44} using the \f\ $\wt\psi:[0,m-1]\to E$ defined by $\wt\psi(k)
:=\#(S_{k+1})$, $k\in[0,m-1]$, and \er{7.12}.
\eex

We now return to the map $\Pi:E^*\to E^*$ defined in \er{4.7}. \If from \E\Pr\
\rf{p4.5} \er{4.8}, \er{4.10}, that $\Pi(n)$, $n\in E^*$, is a \cme product\index{composite!product}
in the monoid $(E,\cdot,1)$. Such a \cme product is usually denoted by
$\prodl_{k=0}^p a_k$ instead of $\QU_{k=0}^p a_k$. Let $k\in E^*$. Then\glossary{$\prodl_{k=0}^p a_k$}
$\Pi(k+1)=(k+1)\Pi(k)$ by \er{4.10}. \Mo $\Pi(1)=1$ by \er{4.8}. Hence
$1=\Pi(0+1)=(0+1)\Pi(0)$ provided $\Pi(0):=1$.

\def\namespec{Proposition and Definition}
\begin{tspc}
\E\te s one and only one map $\psi:E\to E^*$ \sf ying
\bga4.46
\psi(0)=1,\\
\psi(k+1)=(k+1)\psi(k), \q k\in E. \lb{4.47}
\e
The \nm\ $\psi(n)$, $n\in E$, is usually denoted by $n!$ and is called\glossary{$n!$}
\emph{$n$~factorial}. The map $($\f$)$ $n\mt n!$ from $E$ into $E^*$ is called
the \emph{factorial map}. It plays an important role in mathematics.
\bga4.48
n!=\prod_{k=1}^n k, \q n\in E^*, \\
(n+p)!=n!\cdot \prod_{k=n+1}^{n+p} k, \q n,p\in E^*. \lb{4.49}
\e
If $A$ is a \nf, and $\Bij(A)$ denotes the set of bi\jn s of~$A$ onto~$A$,
then
\beq4.50
\#(\Bij(A))=\#(A)!
\e
\end{tspc}

\proof
\E\ex\ and \uq\ of the map $\psi$ is a direct con\sq\ of Theorem \rf{t7.2}
with $F:=E^*$, $f_i(x):=i\cdot x$, $i,x\in E^*$, where $\cdot$ is the \mlc\
in~$E^*$, which is \asc e. Setting $a_0:=1$ and $a_i:=i$, $i\in[1,n]$, $n\in
E^*$, we find that $\psi(k+1)=f_{k+1}(\psi(k))$, $k\in[0,n]$, by \er{4.47}.
\E\Tf we have $\psi(n)=\prodl_{k=0}^n a_k$ for $n\in E$. For $n\in E^*$, we
have $\psi(n)=\Bigl(\prodl_{k=0}^0 a_k\Bigr)\cdot \Bigl(\prodl_{k=1}^n a_k
\Bigr)$ by \er{7.13}. Since $\prodl_{k=0}^0 a_k\nde7.10 = a_0:=1$, \er{4.48}
follows from \er{2.13}.

\Mo for $n\in E^*$ and $p\in E^*$, we have $n+p\in E^*$, and $(n+p)! =
\prodl_{k=0}^{n+p} a_k \nde7.13 = \Bigl(\prodl_{k=0}^n a_k\Bigr)\cdot
\Bigl(\prodl_{k=n+1}^{n+p} a_k\Bigr) = n!\prodl_{k=n+1}^{n+p} k$, which is
\er{4.49}. In $\nde7.13 = $ the \asc ity of the \mlc\ on~$E^*$ is used.

Finally, we prove \er{4.50}. Define $\wt\Pi:E\to E^*$ by setting $\wt\Pi(0):=1$
and $\wt\Pi(k):=\Pi(k)$, $k\in E^*$. Then $\wt\Pi(0+1)=\wt\Pi(1)=\Pi(1)
\nde4.8 = 1=(0+1)\cdot 1=(0+1)\cdot\wt\Pi(1)$. Since $\wt\Pi(n)=\Pi(n)$ for
$n\in E^*$, we have $\wt\Pi(k+1)=(k+1)\wt\Pi(k)$, $k\in E^*$, by \er{4.10}.
Hence $\wt\Pi(k+1)=(k+1)\wt\Pi(k)$ \fa $k\in E$. \If from the
``\uq\ part'' of Theorem \rf{t7.2} that $\wt\Pi(n)=\psi(n)$ \fa $n\in E^*$%
. \E\Tf $\Pi(n)=\psi(n)$, $n\in E^*$. Hence $\#(\Bij([1,n]))\nde
4.7 = \Pi(n)=\psi(n)=n!$ for $n\in E^*$. Since $\#(\Bij(A))\nde4.5 = \#(\Bij
([1,\#(A)]))$, we infer \er{4.50}.
\endproof

In the remaining part of this section we investigate \pp ies of a \cme \op\
of an abelian monoid whose index is a finite set.

\bpr4.18
Let $(X,\qu,e)$ be an abelian monoid.

\smallskip \hph i,ii,
Let $\O$ be a \emph{nonempty} set, let $a$ be a map from~$\O$ into~$X$
and let $A$ be a \emph{finite} subset of~$\O$. Let $\QU_{i\in A}a_i$ denote
the \cme \op\ defined in \er{7.17}. We recall that if the binary \op\ is
denoted by~$+$ $($resp.~$\cdot)$ then $\QU_{i\in A}$ is denoted by
$\suml_{i\in A}$ $($resp.~$\prodl_{i\in A})$ and is called \cme sum\index{composite!sum} $($resp.\
product\/$)$.\index{composite!product} We also recall that if $A$ is \emph{empty}, then, by \df\
\er{7.17},
\beq1.125
\QU_{i\in A}a_i=e.
\e
\Mo if $a:A\to X$ is a \emph{constant} map, i.e.\ \te s $x\in X$ \st
$a_i:=x$ \fa $i\in A$, then
\beq1.126
\QU_{i\in A}a_i = \#(A)\dqu x.
\e
In this case $\QU_{i\in A}a_i$ is also written $\QU_{i\in A}x$.

\smallskip \hph ii,i,
We assume that $\O$ is a \emph{\nf} which is a \emph{finite disjoint}
union of subsets $\zb AiI$, where $I$~is a \emph{\nf}, i.e.\ $\O=\bigcup
\limits_{i\in I}A_i$ and $A_k\cap A_l=\vn$ whenever $k\ne l$, $k,l\in I$. In
this case the \fw\ holds.
\beq1.127
\QU_{\o\in\O}a_\o = \QU_{i\in I}\Bigl(\QU_{\o\in A_i}a_\o\Bigr)=: \QU_{i\in I}\QU_{\o\in A_i}a_\o.
\e
\E\Ip if $a(\O)$ denotes the \emph{range} of the map~$a$, i.e.\
$a(\O):=\{x\in X: \hbox{\te s}$ $\o\in\O$ \st $a_\o=x\}$, and if $A_x$,
$x\in a(\O)$, denotes the subset of~$\O$ consisting of \el s $\o\in\O$ \sf
ying $a_\o=x$, then we have
\beq1.128
\QU_{\o\in\O}a_\o = \QU_{x\in a(\O)}\#(A_x) \dqu x.
\e

\smallskip \hph iii,,
We now suppose that $I$ is a \emph{\nf} and that $a,b$ are maps from~$I$
into~$X$. The \fw\ holds.
\beq1.129
\QU_{i\in I}(a_i\qu b_i)=\Bigl(\QU_{i\in I}a_i\Bigr)\qu
\Bigl(\QU_{i\in I}b_i\Bigr).
\e

\smallskip \hph iv,,
Let $I$ be a \emph{\nf}. Let $a$ be a map from~$I$ into~$X$. Let $(\wt X,
\tqu,\wt e)$ be an abelian monoid and let $\vf:X\to \wt X$ \sf y
\beq1.130
\vf(x\qu y)=\vf(x)\tqu \vf(y) \qh{\fa} x,y\in X.
\e
Then the \fw\ holds.
\beq1.131
\vf\Bigl(\QU_{i\in I}a_i\Bigr)=\tQU_{i\in I}(\vf(a_i)),
\e
where $\tQU_{i\in I}$ denotes the \cme \op\ in $(\wt X,\tqu,\wt e)$.
\epr

\proof \

\er{1.126}: Set $c:=\QU_{\o\in A}a_\o$. If $A=\vn$, then $\#(A)\nde{3.6} =
0$ and $c=e$ by \er{1.125}. Hence $c=\#(A)\dqu x$ by \er{2.3}. If $A\ne\vn$,
then $c\nde{7.17} = \QU_{j=0}^p a_{\vf(j)}$ where $p:=\#(A)-1$ and $\vf$~is
a bi\jn\ from $[0,\#(A)-1]$ onto~$A$. Note that $a_{\vf(j)}=x$ \fa
$j\in[0,\#(A)-1]$. In view of \er{7.17}, $\QU_{j=0}^0 a_{\vf(j)}=x$. Hence
$c=x$ whenever $\#(A)=1$, which proves \er{1.126} for $\#(A)=1$, since
$1\dqu x=x$ by \er{2.3}. Next we assume $\#(A)>1$. Then ${c=\QU_{j=0}
^{\#(A)-1}}x\nde{7.12} = \QU_{j=1}^{\#(A)}x$, since $(\#(A)-1)+1 =
S(P(\#(A)))=\#(A)$. Define $\vf:E\to X$ by setting $\vf(0):=e$ and $\vf(n)
:=\QU_{j=1}^nx$, $n\in E^*$. Then $\vf$~\sf ies $\vf(1)=\suml_{j=1}^1x
\nde7.12 = \suml_{j=0}^0x\nde7.10 = x$. Hence $\vf(0+1)=\vf(1)=x=x\qu e=
x\qu\vf(0)$. \Mo $\vf(k+1) = \suml_{j=1}^{k+1}x \nde7.12 = \suml_{j=0}^k x
\nde7.13 = \Bigl(\suml_{j=0}^0 x\Bigr)\qu \Bigl(\suml_{j=1}^kx\Bigr)
\nde7.10 = x\qu\vf(k)$. We used the \asc ity of~$\qu$ in~\er{7.13}. Thus
$\vf$~\sf ies $\vf(0)=e$ and $\vf(k+1)=x\qu \vf(k)$ \fa $k\in E$. \If from
\E\Pr\ \rf{p1.13} and \E\df s and Notation \rf{d1.14} that $\vf(n)=n\dqu x$.
Hence $\vf(\#(A))=\#(A)\dqu x$. \E\Tf $\QU_{i\in A}a_i\nde7.11 = c=\vf(\#(A))=
\#(A)\dqu x$.

\smallskip
\er{1.127}: We first assume $\#(I)=1$. Let
$\a$~denote the unique \el\ of~$I$. Then \te s a unique bi\jn\ $\vf:[0,1)
\to I:=\{\a\}$. We have $\O=A_\a$ and $\QU_{i\in I}\Bigl(\QU_{\o\in A_i}
a_\o\Bigr)\nde{7.17} = \QU_{k=0}^0 \Bigl(\QU_{j\in A_{\vf(k)}}a_j\Bigr)
\nde{7.10} = \QU_{j\in A_{\vf(0)}}a_j = \QU_{j\in A_\a}a_j \nde7.11 = \QU_{\o\in A_\a}
a_\o=\QU_{\o\in\O}a_\o$.

$\#(I):=S(1)$: Let $\a$ and $\b$ denote the \el s of $I$. Then $\vf:[0,S(1))\to I$
defined by $\vf(0):=\a$ and $\vf(1):=\b$ is a bi\jn. We have $\O=A_\a\cup
A_\b$, $A_\a\cap A_\b=\vn$ and $\O\ne\vn$. We distinguish three cases. \ti{Case
$A_\a=\vn$}. Then $\O=A_\b\ne\vn$, hence $\QU_{\o\in\O}a_\o = \QU_{\o\in A_\b}
a_\o = e\qu \QU_{\o\in A_\b}a_\o \nde1.125 = \QU_{\o\in A_\a}a_\o \qu
\QU_{\o\in A_\b}a_\o = \QU_{\o\in A_{\vf(0)}}a_\o \qu \QU_{\o\in A_{\vf(1)}}
a_\o \nde7.10 = \QU_{l=0}^1 \Bigl(\QU_{\o\in A_{\vf(l)}}a_\o \Bigr)
\nde7.17 = \QU_{i\in I}\Bigl(\QU_{\o\in A_i}a_\o\Bigr)$.
The case $A_\b=\vn$ is similar.
\ti{Case $A_\a$ and $A_\b$ nonempty}. In view of Theorem \rf{t3.23}, $A_\a
\cup A_\b$ is finite and $\#(\O)=\#(A_\a\cup A_\b)= \#(A_\a)+\#(A_\b)$, where
$\#(\O)$, $\#(A_\a)$ and $\#(A_\b)$ belong to~$E^\ast$. Proceeding as in the
proof of Theorem \rf{t3.23}, we construct a~bi\jn\
$\vf:[0,\#(A_\a)+\#(A_\b)) \to\O$ \st its \rt ion to $[0,\#(A_\a))$ is
a~bi\jn\ onto~$A_\a$ and its \rt ion to $[\#(A_\a),\#(A_\a)+\#(A_\b))$ is
a~bi\jn\ onto~$A_\b$. \E\Tf $\QU_{\o\in\O}a_\o \nde7.17 =
\QU_{l=0}^{\#(A_\a)+\#(A_\b)-1}a_{\vf(l)} \nde7.13 =
\QU_{l=0}^{\#(A_\a)-1}a_{\vf(l)}\qu \QU_{l=\#(A_\a)}^{\#(A_\a)+\#(A_\b)-1}
a_{\vf(l)} \nad{\er{7.17},\er{7.12}}= \QU_{\o\in A_\a}a_\o \qu
\QU_{l=0}^{\#(A_\b)-1}a_{\vf(l)}\nde7.17 = \QU_{\o\in A_\a}a_\o
\qu \QU_{\o\in A_\b}a_\o$. As in the case $A_\a=\vn$, we obtain
$\Bigl(\QU_{\o\in A_\a}a_\o\Bigr)\qu\Bigl(\QU_{\o\in A_\b}a_\o\Bigr)=
\QU_{i\in I}\Bigl(\QU_{\o\in A_i}a_\o\Bigr)$. This completes the proof of the
case $\#(I)=S(1)$.

Note that we proved that if $\O=A_\a\cup A_\b$ with $A_\a\cap
A_\b=\vn$, $A_\a,A_\b$ finite and $\O\ne\vn$, the \fw\ holds:
\beq1.138
\QU_{\o\in\O}a_\o=\Bigl(\QU_{\o\in A_\a}a_\o\Bigr)\qu
\Bigl(\QU_{\o\in A_\b}a_\o\Bigr) = \QU_{i\in\{\a,\b\}}\Bigl(
\QU_{\o\in A_i}a_\o\Bigr).
\e

Observe that if $(X,\qu)$ is an abelian \sg\ and $\QU_{i\in A}$ is defined
by \er{7.17} only for $A\ne\vn$, then \er{1.138} holds whenever both
$A_\a$ and~$A_\b$ are not empty.

We now set $M:=\bigl\{m\in E^*\sms1: \hbox{\er{1.127} holds for }\#(I)=m\}$.
We want to prove that $M=E^*\sms1$. Since $S(1)\in M$, it is \sft\
(see Remark \rf{r4.19}) to show that
$n+1\in M$ whenever $n\in M$ (induction hypothesis). We assume that $n\in M$
and we will prove that $n+1\in M$. We have $\#(I)=n+1$, and $\O$~is a \nf\
which is the disjoint union of (finite) sets $\zb AiI$. Let $\g\in I$ be \st
$A_\g\ne\vn$ and set
$J:=I\sms\g$. Set $\O':=\bigcup\limits_{j\in J}A_j$. We have $\O=\O'\cup
A_\g$ and $\O'\cap A_\g=\vn$. In view of \er{1.138} we have
\beq1.139
\QU_{\o\in\O}a_\o=\QU_{\o\in\O'}a_\o \qu \QU_{\o\in A_\g}a_\o.
\e
We claim that
\beq1.140
\QU_{\o\in\O'}a_\o=\QU_{j\in J}\Bigl(\QU_{\o\in A_j}a_\o\Bigr).
\e
Indeed, if $\O'=\vn$, then $\bcl_{j\in J} A_j=\vn$, hence $A_j=\vn$ \fa $j\in J$.
Hence $\QU_{\o\in A_j}a_\o=e$ \fa $j\in J$ and $\QU_{j\in J}\Bigl(
\QU_{\o\in A_j}a_\o\Bigr)=\QU_{j\in J}e\nde1.126 = \#(J)\dqu e\nde2.3 =e$.
But $e\nde1.125 = \QU_{\o\in\O'}A_j$, since $\O'=\vn$. Thus \er{1.140} holds.
If $\O'\ne\vn$, we have $\#(J)=n$ by~\er{1.8}, since $\#(J)+1\nde3.5 = \#(J)+ \#(\{\g\})
\nde3.31 = \#(I)=n+1$. Since $\O'=\bcl_{j\in J}A_j$,
with $A_l\cap A_k=\vn$ whenever $l\ne
k$, $l,k\in J$, \er{1.140} follows from the induction hypothesis $n\in M$.
Combining \er{1.139} and \er{1.140} we arrive at
\beq1.141
\QU_{\o\in\O}a_\o=\Bigl(\QU_{j\in J}\Bigl(\QU_{\o\in A_j}a_\o\Bigr)\Bigr) \qu
\QU_{\o\in A_\g}a_\o.
\e
We now define a map $b:I\to X$ by setting $b_i:=\QU_{\o\in A_i}a_\o$ for $i\in
J$ and $b_\g:=\QU_{\o\in A_\g}a_\o$. The \RHS\ of \er{1.141} becomes
$\QU_{i\in J}b_i \qu \QU_{i\in\{\g\}}b_i$ since $\QU_{i\in\{\g\}}b_i\nde1.126
= \#(\{\g\})\dqu b_\g \nde2.3 = b_\g$. \If from \er{1.138} with $\O:=I$,
$A_\a:=J$, $A_\b:=\{\g\}$, $a_\o:=b_\o$, $\o\in\O$, that
$\QU_{i\in J}b_i \qu \QU_{i\in \{\g\}}b_i = \QU_{i\in J\cup\{\g\}}b_i
=\QU_{i\in I}b_i =\QU_{i\in I}\Bigl(\QU_{\o\in A_i}a_\o\Bigr)$, since
$I=J\cup\{\g\}$ and $J\cap\{\g\}=\vn$.  Consequently
the \RHS\ of \er{1.141} is equal to $\QU_{i\in I}\Bigl(\QU_{\o\in
A_i}a_\o\Bigr)$, which implies \er{1.127} for $\#(I)=n+1$. Hence $n+1\in M$.
\If that \er{1.127} holds \fa $n\in E^*\sms1$. Since \er{1.127} holds for $n=1$,
the proof is complete.

\smallskip
\er{1.128}: We set $I:=a(\O)$ which is nonempty since $\O\ne\vn$ and which is
finite by Theorem \rfa1{t4.18}\,(ii). Given $x\in I$, set $A_x:=\{\o\in\O:
a_\o=x\}$. Then $A_x\ne\vn$ \fe $x\in I$ and $A_x\cap A_y=\vn$ \fa $x,y\in
I$, $x\ne y$. We may thus apply \er{1.127} and obtain $\QU_{\o\in\O}a_\o
\nde1.127 = \QU_{x\in a(\O)}\Bigl(\QU_{\o\in A_x}a_\o\Bigr) =
\QU_{x\in a(\O)}\Bigl(\QU_{\o\in A_x}x\Bigr) \nde1.126 = \QU_{x\in a(\O)}
\#(A_x)\dqu x$. This completes the proof of \er{1.128}.

\smallskip
\er{1.129}: We use induction on $n:=\#(I)\in E^*$. Set $M:=\{m\in E^*:
\hbox{\er{1.129} holds for}$ every set $I$ with $\#(I)=m\}$.
We have: $1\in M$:
Indeed, let $I$ be a set containing only one \el\ which we denote by~$\a$,
and let $a,b:I\to X$. Then $\suml_{i\in I}a_i \nde1.126 = \#(I)\dqu a_\a =\break
1\dqu a_\a \nde2.3 = a_\a$. Similarly $\QU_{i\in I}b_i=b_\a$ and
$\QU_{i\in I}(a_i\qu b_i)=a_\a \qu b_\a$. Hence $\QU_{i\in I}(a_i\qu
b_i)=\break a_\a \qu b_\a = \QU_{i\in I}a_i \qu \QU_{i\in I}b_i$, and $1\in M$.
\ti{$n\in M$ implies $n+1\in M$}: We suppose that $n\in M$ (induction
hypothesis) and let $I$~be a~set of $n+1$ \el s. Let $a,b$ be maps from~$I$
into~$X$. Define $c:I\to X$ by setting $c_i:=a_i\qu b_i$, $i\in I$. Let
$\g\in I$ and set $J:=I\sms\g$. Then $\#(J)=n$ (see the proof of \er{1.127})
and set $A_\a:=J$, $A_\b:=\{\g\}$. In view of \er{1.138} we obtain $\QU_{i\in
I}c_i = \QU_{i\in J}c_i \qu \QU_{i\in\{\g\}}c_i$, i.e.\ $\QU_{i\in I}
(a_i\qu b_i) = \QU_{i\in J}(a_i\qu b_i)\qu (a_\g\qu b_\g)$. Since $n\in M$,
we have $\QU_{i\in J}(a_i\qu b_i) = \Bigl(\QU_{i\in J}a_i\Bigr) \qu
\Bigl(\QU_{i\in J}b_i\Bigr)$. Then $\QU_{i\in I}(a_i\qu b_i)=\Bigl(\QU_{i\in
J}a_i \qu \QU_{i\in J}b_i\Bigr)\qu{(a_\g\qu b_\g)} \nde1.36 = \QU_{i\in J}a_i \qu
\Bigl(\QU_{i\in J} b_i\qu a_\g\Bigr) \qu b_\g \nde1.6 = \QU_{i\in J}a_i
\qu\Bigl(a_\g \qu \QU_{i\in J}b_i\Bigr) \qu b_\g \nde1.36 = \Bigl(\QU_{i\in J}a_i
\qu a_\g\Bigr) \qu {\Bigl(\QU_{i\in J}b_i\qu b_\g\Bigr)} \nde1.138 =
\Bigl(\QU_{i\in I}a_i\Bigr)\qu \Bigl(\QU_{i\in I}b_i\Bigr)$. Hence $n+1\in
M$, $M=E^*$ and \er{1.129} is proved.

\smallskip
\er{1.131}: Since the proof is quite similar to the proof of \er{1.129}, we
only sketch it. We set $M:=\{m\in E^*: \hbox{\er{1.131} holds for sets $I$
with }\#(I)=m\}$. $1\in M$: $\vf\Bigl(\QU_{i\in I}a_i\Bigr) = \vf(a_\a)$
where $I:=\{\a\}$. But $\vf(a_\a)=\tQU_{i\in I}\vf(a_i)$. \ti{$n\in M$
implies $n+1\in M$}: Let $\g\in I$, where $\#(I)=n+1$. Then by \er{1.138}
$\tQU_{i\in I}\vf(a_i)= \Bigl(\tQU_{i\in J}\vf(a_i)\Bigr) \tqu \Bigl(
\tQU_{i\in\{\g\}}\vf(a_i)\Bigr) \nad{n\in M}= \vf\Bigl(\QU_{i\in J}a_i\Bigr)
\tqu \vf(a_\g)$. Using \er{1.130} we obtain $\tQU_{i\in I}
\vf(a_i)=\vf\Bigl(\Bigl(\QU_{i\in J}a_i\Bigr)\qu a_\g\Bigr)$. But $a_\g
=\QU_{i\in\{\g\}}a_i$ and $\Bigl(\QU_{i\in J}a_i\Bigr)\qu
\Bigl(\QU_{i\in\{\g\}}a_i\Bigr) \nde1.138 = \QU_{i\in I}a_i$. Note that
$I$, $J$ and~$\{\g\}$ are not empty (we do not use $\vf(e)=\wt e$). \E\Tf we
obtain $\vf\Bigl(\QU_{i\in I}a_i\Bigr) = \tQU_{i\in I}\vf(a_i)$. Hence
$n+1\in M$, $M=E^*$ and \er{1.131} is proved.
\endproof

\brm4.19
We used in the proof of \E\Pr\ \rf{p4.18} the \fw\ argument. If $M$ is a
subset of $E\sms{0,1}$ and \sf ies (i)~$S(1)\in M$ and (ii)~if $x\in M$, then
$x+1\in M$, then $M=E\sms{0,1}$. This is a con\sq\ of~\er{2.2}. Indeed, set
$\wt M:=\{0,1\}\cup M$. Then $\{0,1\}\cap M=\vn$ since $M\subset E\sms{0,1}$.
\Mo if $x=0$, then $x+1=1\in\{0,1\}\subset\wt M$; if $x=1$, then $x+1=S(1)
\nad{\rm(i)}\in M\sbs\wt M$ by (ii). \E\Tf $x\in\wt M$ implies $x+1\in\wt M$.
\If from \er{2.2} that $\wt M=E$, hence $M=\wt M\sms{0,1}=E\sms{0,1}$.
\erm

A similar proof holds for $E\sms0$ instead of $E\sms{0,1}$. See also
Section~\ref{sss.Sets}.

\bex4.20
Let $(X,\qu,e)$, $\O$, $a$ and $A$ be as in \E\Pr\ \rf{p4.18}, and let $n\in
E$. Prove
\beq4.62
n\dqu \QU_{i\in A}a_i = \QU_{i\in A}(n\dqu a_i).
\e
\eex

\bex4.21
Let $(X,\qu,e):=(E,+,0)$, let $\O,A$ be as in \E\Pr\ \rf{p4.18}, let $a,b:\O
\to E$ be \st $a_i\le b_i$, \fa $i\in\O$. Prove
\beq4.63
\QU_{i\in A}a_i \le \QU_{i\in A}b_i.
\e
If, in \ad, \te s $j\in A$ \st $a_j<b_j$, prove
\beq4.64
\QU_{i\in A}a_i < \QU_{i\in A}b_i.
\e
\eex

\begin{exe}\lb{ex7.8}
Let $(X,\qu,e)$ be an abelian monoid and let $a$ be a map from a \nf\
$\O$ into~$X$.

\hph i,ii, Let $A,B$ be two subsets of $\O$. Show that
\beq7.18
\Bigl(\QU_{i\in A\cup B}a_i\Bigr) \qu \Bigl(\QU_{i\in A\cap B}a_i\Bigr)=
\Bigl(\QU_{i\in A}a_i\Bigr) \qu \Bigl(\QU_{i\in B}a_i\Bigr).
\e

\hph ii,i, Let $\O:=I\t J$ where $I:=[1,m]$, $J:=[1,n]$, $m,n\in E^*$. Let
$a:\O\to X$. Show that
\bea7.19
\QU_{(i,j)\in\O}a_{ij}&=\QU_{i\in I}\Bigl(\QU_{j\in J}a_{ij}\Bigr),\\
\QU_{j\in J}\Bigl(\QU_{i\in I}a_{ij}\Bigr)
&=\QU_{i\in I}\Bigl(\QU_{j\in J}a_{ij}\Bigr), \lb{7.20}
\e
where $a_{ij}$ denotes the value of $a$ at $(i,j)$.

\hph iii,, Let $x$ be an \el\ of $X$. In (ii) set $a_{ij}:=x$ \fe $(i,j)\in\O$.
Show that
\beq7.21
\QU_{(i,j)\in\O}a_{ij}=(m\cdot n)\dqu x.
\e

\hph iv,, Let $m\in E$ and let $\O:=\{(i,j)\in[0,m]\t[0,m]: j\le i\}$. Let $a:\O\to
X$. Show that
\beq7.22
\QU_{i=0}^m\Bigl(\QU_{j=0}^i a_{ij}\Bigr)=
\QU_{j=0}^m\Bigl(\QU_{i=j}^m a_{ij}\Bigr).
\e

\hph v,i, Let $m\in E$. Show (by induction) that \te s a unique $p\in E$ \st
$(m+1)(m+2)=2p$. Show that if in~(v) $a_{ij}=x$, $x\in X$, \fe $(i,j)\in\O$,
then
\beq7.23
\QU_{i=0}^m\Bigl(\QU_{j=0}^i a_{ij}\Bigr)=p\dqu x.
\e

\hph vi,, Let $\O,\O'$ be nonempty \ep\ finite sets and let $f:\O'\to \O$ be
a bi\jn. Show that
\beq4.71
\QU_{\o'\in\O'}a_{f(\o')} = \QU_{\o\in\O}a_\o.
\e
Note that if $\O'=\O$, then \er{4.71} is \gn\ of \er{7.14}.
\eex

\newpage
\Subsubsection{Decimal expansion of \nn s}\label{sss.Decimal}\index{decimal expansion of a natural number}
The \ds\ of \rp ation of \nn s is the most widely used system. It consists of
representing \nn s as a sum of \ml es of powers of the \nm\ ten. The main
purpose of this section is to present a theorem justifying this \rp ation.
In this section $(E,0,S)$ denotes a set of \nn s with distinguished \el~$0$.
\Mo $1:=S(0)$, $E^*:=E\sms0$ and $\le$ denotes the \nog\ on~$E$. Given $m,n
\in E$, $m+n$, $m\cdot n$ (or~$mn$) and $m^n$ denote the \art\ \op s
introduced in Sections \ref{sss.Add} and \ref{sss.Mult}, i.e.\ $m+n:=m+_En$,
$m\cdot n:=m\cdot_E n$ and $m^n$ denotes the $n$-fold \IT\ of~$m$ \wrt
$(E,0,S)$ in the abelian monoid $(E,\cdot,1)$.

The simplest system for \rp ing \el s
of~$E^*$ is the \ti{unary system} which uses only one symbol. This symbol is
simply repeated. For example, the \nm\ one is \rp ed by \6, $S(\6)$ by \2,
$S(\2)$ by~\3, and so on. It is commonly used for ``small'' \nm s and in
theoretical computer science. A~modification of this system consists of using
abbreviations for some \nm s. For example \5\ is replaced by
\kr\ by some card (e.g.~jass) players. The Romans used the symbol V for \5,
and consistently IV for \4\ and VI, VII, VIII for \kr\6, \kr\2, \kr\3. They
used the symbol X for \kr\kr, and IX for \kr\4, XI, XII, XIII, XIV instead of
\kr\kr\6, \kr\kr\2, \kr\kr\3, \kr\kr\4. They used the symbol L for \rp
ing \Xk, C~for LL, D for \Ck\ and M for DD.

It is not our purpose to present a history of the main systems of numeration
(see e.g.~\cite{Nrs}).
Instead, we \rt\ ourselves to the \ds\ and some of its \gn s.\index{decimal system} In this system
every \nn\ is \rp ed by a ``finite \sq'' of symbols called
\ti{digits}. These symbols are 0, 1, 2, 3, 4, 5, 6, 7, 8, 9. They are called
Hindu-arabic numerals. In the \ti{binary} system of numeration,\index{binary system of numeration} only two of
these digits (\ti{bits})\index{bits} are used, namely the symbols 0 and~1. In both systems the
symbol~0 \rp s the \nm~0 which is the smallest \el\ of~$(E,\le)$.
Similarly, the symbol~1 \rp s the \nm\ one which is the \im \su\ of~0,
i.e.\ $1:=S(0)$. In the \ds\ the symbols 2, 3, 4, 5, 6, 7, 8, 9 \rp\ in
this order the \nm s $S(1)$, $S(2)$, $S(3)$, $S(4)$, $S(5)$, $S(6)$, $S(7)$,
$S(8)$. These \nm s are called \ti{two}, \ti{three}, \ti{four}, \ti{five},
\ti{six}, \ti{seven}, \ti{eight}, \ti{nine} in English. The \im \su\ of
the \nm\ nine, $S(9)$, is called \ti{ten} in English and \ti{decem} in Latin.
In the \ds\ it is \rp ed by 10 and in the binary system by 1010. (The \nm\
ten is famous in England because Downing Street \nm~10 is the office of the
Prime Minister.) The \nm\ ten and its powers play an essential role in the
\ds. In the binary system the \nm\ two plays a
fundamental role. In this system the powers of two play the role of the
powers of ten in the \ds.  In both systems the goal is to \rp\ every \el\
of~$E$ as a sum of some \ml es of a power of a
distinguished \nm~$B$ greater than one, called the \ti{base} of the system.
The \nm\ zero is denoted by~0 in all bases.

For convenience we shall use the notation
\beq5.1
\bal
D_{10}&:=\hbox{the set of digits}:=\{0,1,2,3,4,5,6,7,8,9\},\\
\dot D_{10}&:=D_{10}\sm\{0\}.
\eal
\e
We first consider the case of sums of special \ml es of $10^0$ and $10^1$ of the form
$m\cdot 10^1+n\cdot 10^0$, where $m$ and~$n$ are \el s of~$D_{10}$. If
$m=n=0$, then the corresponding \nm\ is zero. We could \rp\ it by~00, but it
is customary to use only the \rp ation~0. Similarly if $m:=0$ and $n\in \dot
D_{10}$, then the corresponding \nm\ is~$n$ and it is \rp ed by~$n$. In
these cases the \el s of~$D_{10}$ have two meanings. Either a digit or a \nm\
less than ten. If $m:=1$ and $n:=0$, then 10 is the \rp ation of $S(9)$.
Similarly the symbol $1n$ with $n\in\dot D_{10}$ \rp s $1\cdot 10+n$. Observe
that $S(19)=S(10+9)= (10+9)+1=10+(9+1)=10+10=2\cdot 10
=2\cdot 10^1=2\cdot 10^1+0\cdot 10^0$. So it is \rp ed by 20.

We shall prove that every \nm\ $k\in[10^1,10^2)$ can be \rp ed in a unique
fashion as $m\cdot10+n$ where $m\in \dot D_{10}$ and $n\in D_{10}$. It will
be denoted by $mn$. A~similar result holds for $k\in[10^2,10^3)$, but with
$m\cdot10+n$ replaced by $m\cdot10^2+n\cdot10^1+p$. We can repeat this
process by replacing $[10^2,10^3)$ by $[10^N,10^{N+1})$, where $N\in E^*$. In
that case we are led to consider sums of several terms which are not
identical in general, i.e.\ \ti{\cme sums}.

In this section, Theorem \rf{t1.38} (division algorithm) will play an
important role. For convenience we reformulate it here with a somewhat different
notation.

\bth5.1
\E\fe $a\in E^*$, \te\ maps $Q_a:E\to E$ and $\F_a:E\to\zo0,a $ \st
\beq5.2
m=Q_a(m)\cdot a+\F_a(m)\qh{\fa} m\in E.
\e
\Mo if $m\in E$ \sf ies
\beq5.3
m=q\cdot a+r \qh{\fs} q\in E \hbox{ and }r\in\zo0,a ,
\e
then $q=Q_a(m)$ and $r=\F_a(m)$.
\eth

\proof
We apply Theorem \rf{t1.38} with $\ees:=(E,0,S)$ and observe that $n\cdot a$
in \er{1.59}, the $n$-fold \IT\ of~$a$ \wrt $(E,0,S)$ in the monoid $(E,
+,0)$, is equal by \er{2.8} to $n\cdot_E a$, which we also denote
by~$n\cdot a$. In view of Theorem \rf{t1.38}, if $m\in E$, then \te s a unique
pair $(q,r)\in E\t \zo0,a $ \st \er{5.3} holds. Hence we may define $Q_a$
(resp.~$\F_a$) by setting $Q_a(m):=q$ and $\F_a(m):=r$. The ``\uq'' part of
Theorem \rf{t5.1} is the same as the ``\uq'' part of Theorem \rf{t1.38}.
\endproof

Since $Q_a$ is a selfmap of $E$, we can define the
\sq\ of its \ti{\IT s}, which we denote by $\{Q_a^i\}_{i\in E}$ (see~\er{2.4}).

\blm5.2
Let $a\in E^*$ and $m\in E$. Then \fe $n\in E$,
\beq5.4
m=Q_a^{n+1}(m)\cdot a^{n+1} + \sum_{k=0}^n \F_a(Q_a^k(m))\cdot a^k.
\e
\elm

\proof
By \In\ on $n\in E$. Set $M:=\{n\in E: \er{5.4} \hbox{ holds \fa} m\in E\}$.

$0\in M$: Observe that $Q_a^1(m)\nde2.7 = Q_a(m)$ and $Q_a^0(m)\nde2.7 = \id_E(m)=m$. Then
$0\in M$ in view of \er{5.2}.

\ti{$n\in M$ implies $n+1\in M$}: Let $n\in M$. In view of \er{5.2} we may
replace $Q_a^{n+1}(m)$ in \er{5.4} by $Q_a((Q_a^{n+1}(m)))\cdot a + \F_a
(Q_a^{n+1}(m))$. We obtain $m\nad{\er{2.6},\er{1.6}}=\bigl(Q_a^{(n+1)+1}(m)\cdot a+\break\F_a(Q_a^{n+1}(m))\bigr)
\cdot a^{n+1} + \suml_{k=0}^n \F_a(Q_a^k(m))\cdot a^k \nad{\er{2.15},\er{1.5}}= (Q_a
^{(n+1)+1}(m)\cdot a)\cdot a^{n+1} + \F_a(Q_a^{n+1}(m))\cdot a^{n+1}+
\suml_{k=0}^n\F_a(Q_a^k(m))\cdot a^k \nad{\er{2.11},\er{2.24},\er{2.25}} = Q_a^{(n+1)+1}(m)\cdot a^{(n+1)+1}
+\F_a(Q_a^{n+1}(m))\cdot a^{n+1}+\suml_{k=0}^n\F_a(Q_a^k(m))\cdot a_k$.
But $\F_a(Q_a^{n+1}(m))\cdot a^{n+1}+\suml_{k=0}^n\F_a(Q_a^k(m))\cdot a^k
=\suml_{k=0}^{n+1} \F_a(Q_a^{k}(m))\cdot a^k$ by \er{7.10}. Hence
$m=Q_a^{(n+1)+1}(m)\cdot a^{(n+1)+1} + \suml_{k=0}^{n+1} \F(Q_a^k(m))\cdot
a^k$. \E\Tf $n+1\in M$ and $M=E$.
\endproof

\blm5.3
Let $a\in E^*\sms1$ and $m\in E^*$. Then \te\ $N\in E$ and a map
$\a:[0,N]\to \zo0,a $ \st
\beq5.5
m=\sum_{k=0}^N \a_ka^k  \hbox{ and } \a_N\ne0.
\e
\elm

\proof
Suppose for \cd ion that $Q_a^{n+1}(m)\ne0$ \fa $n\in E$. Then\break $m\ge Q_a^{n+1}
(m)\cdot a^{n+1}$ by \er{5.4} and \er{1.46}. Since $Q_a^{n+1}(m)\ne0$ \fa
$n\in E$, we have \fa $n\in E$, $Q_a^{n+1}(m)\cdot a^{n+1}\ge1\cdot a^{n+1}=
a^{n+1}$ by \er{2.19}, \er{2.13} since $a^{n+1}\in E^*$ in view of \er{2.30}. Therefore,
by the \tr ity of~$\ge$, we obtain $m\ge a^{n+1}$ \fa $n\in E$. Since $m\in E\sms
{0,1}$, the \sq\ $\{a^{n+1}\}_{n\in E}$ is strictly in\cre\ by \er{2.36}. Set
$c_0:=0$ and $c_k:=a^{k+1}$ for $k\in E^*$. Then the \sq\ $\zb ckE$ \sf ies
the \as s of Lemma \rf{l1.32} with $\ees:=(E,0,S)$. Thus \te s \ooo $p\in E$
\st $m\in\zo c_p,c_{p+1} $. \E\Tf $m<a^{(p+1)+1}$ \cd ing $m\ge a^{n+1}$ \fa
$n\in E$. \If that \te s $\wt N\in E$ \st $Q_a^{\td N+1}(m)=0$, hence $m=\suml
_{k=0}^{\td N} \F_a(Q_a^k(m))\cdot a^k$. Set $\a_k:=\F_a(Q_a^k(m))$, $k\in
[0,\wt N]$. Then by the \df\ of~$\F_a$, we have $\a_k\in[0,a)$ \fa $k\in[0,\wt N]$.
If $\a_k=0$ \fa $k\in[0,\wt N]$, then $\a_k\cdot a^k\nad{\er{2.13}}=
0$ \fa $k\in[0,\wt N]$, hence $\suml_{k=0}^{\wt N}\a_k\cdot a^k=\suml_{k=0}
^{\wt N}0 \nde7.17 = \suml_{k\in[0,\td N]}0 \nde1.126 = \#([0,\td N])\cdot0
\nde2.9 = 0$.
A~\cd ion, since $m\ne 0$ by \as. \E\Tf the set $K:=\{k\in[0,\wt N]:\break\a_k\ne0
\}$ is not empty. As a~subset of $[0,\wt N]$ $(=\zo0,S({\wt N})) $, the set~$K$
is finite. In view of Theorem \rfa1{t3.39},
$K$~possesses a~greatest \el, which we denote by~$N$. Clearly $N\le \wt N$.
If $N:=\wt N$, then $\a_N\ne0$ and \er{5.5} holds. If $N<\wt N$, then $\suml
_{k=0}^{\td N}\a_k\cdot a^k\nde7.13 = \suml_{k=0}^N \a_k\cdot a^k+\suml_{k=
N+1}^{\td N} 0\cdot a^k$, since $N+1\le \wt N$ by \era1{3.17}. Since $\suml
_{k=N+1}^{\td N}0\cdot a_k \nad{\er{2.13}} = \suml_{k=N+1}^{\td N} 0=0$
as above, we have $m=\suml_{k=0}^N \a_k\cdot a^k$ with $\a_N\ne0$, which is
\er{5.5}.
\endproof

In the next lemma we \es\ the \uq\ of the \nm\ $N$ and of the map~$\a$ in
formula \er{5.5}.

\blm5.4
Let $a\in E^*\sms1$, let $N,M\in E$ and let $\a:[0,N]\to\zo0,a $, $\b:[0,M]\to
\zo0,a $ be \st
\beq5.6
\sum_{k=0}^N \a_k\cdot a^k = \sum_{l=0}^M \b_l\cdot a^l, \q \a_N\ne0,\
\b_M\ne0.
\e
Then $N=M$ and $\a=\b$.
\elm

\proof \

(i) $N=M$. Suppose for \cd ion that $N\ne M$. \E\wlg we may suppose $N<M$. We
claim that
\beq5.7
\sum_{k=0}^N \a_k\cdot a^k< a^{N+1}.
\e
Indeed, let $k\in[0,N]$. Then $\a_k<a$, hence $\a_k+1\le a$ by \era1{3.17}.
Since $a\in E^*\sms1$, i.e.\ $a>1$, \te s $p\in E^*$ \st $a=p+1$. From
$\a_k+1\le p+1$, we infer $\a_k\le p$ by \er{1.65}. Hence $\a_k\cdot a^k\le
p\cdot a^k$ by \er{2.19} since $a^k\in E^*$. From \er{7.29} we find
$\suml_{k=0}^N\a_k\cdot a^k \le \suml_{k=0}^N p\cdot a^k$. But $\suml_{k=0}^N
p\cdot a^k = p\cdot\suml_{k=0}^N a^k$ by \er{7.27}. Hence $1+\suml_{k=0}^N
\a_k\cdot a^k \nde1.64 {\le} 1+ p\cdot \suml_{k=0}^N a^k = 1+p\cdot \suml_{k=0}^N(1+p)^k
\nde7.31 = (1+p)^{N+1}=a^{N+1}$. \E\Tf $1+\suml_{k=0}^N \a_k\cdot a^k\le
a^{N+1}$. Hence $\suml_{k=0}^N\a_k\cdot a^k = 0+\suml_{k=0}^N \a_k\cdot a^k
\nde{1.63} < 1+\suml_{k=0}^N\a_k\cdot a^k \le a^{N+1}$. Then the claim
\er{5.7} follows from \era1{3.12}.

We now show that
\beq5.8
a^M\le \sum_{l=0}^M \b_l\cdot a^l
\e
under the \as s $\b_l<a$, $l\in\zo0,M $ and $\b_M\ne0$. If $M:=0$, then $a^M=
a^0 \nde2.23 = 1 \le\b_0\nde2.13 = \b_0\cdot1 = \b_0\cdot a^0\nde7.10 = \suml
_{l=0}^0 \b_l\cdot a^l = \suml_{l=0}^M \b_l\cdot a^l$. Hence \er{5.8} follows
from the \tr ity of~$\le$. We now suppose $M\ge1$. We have $\suml_{l=0}^M
\b_l\cdot a^l \nde7.13 = \suml_{l=0}^{M-1}\b_l\cdot a^l + \suml_{l=M}^M\b_l
\cdot a^l \nde7.12 = \suml_{l=0}^{M-1}\b_l\cdot a^l + \b_M\cdot a^M
\nde1.46 {\ge} \b_M\cdot \a^M$. Since $\b_M\in E^*$ and $a^M\in E^*$
we have $\b_M\ge1$
and $\b_M\cdot a^M \nde2.19 \ge 1\cdot a^M \nde2.13 = a^M$. \E\Tf \er{5.8}
holds by the \tr ity of~$\le$.

Since $N<M$, we have $N+1\le M$ by \era1{3.17}, hence $a^{N+1}\le a^M$ by
\er{2.33} and $a\in E\sms{0,1}$. From \er{5.8} we infer $a^{N+1}\le
\suml_{l=0}^M \b_l\cdot a^l$ and from \er{5.7} we infer $\suml_{k=0}^N
\a_k\cdot a^k < a^{N+1}$,
hence $\suml_{k=0}^N \a_k\cdot a^k < \suml_{l=0}^M \b_k\cdot a^k$, by
\era1{3.12} a~\cd ion. Hence $M=N$, and part (i) is proved.

(ii) We now suppose that \er{5.6} holds with $M=N$. We do not need the \as s
$\a_N\ne0$ and $\b_N\ne0$ in this part of the proof. We proceed by \In\ on
$N\in E$. Set $A:=\{N\in E: \suml_{k=0}^N \a_k\cdot a^k=\suml_{k=0}^N\b_k\cdot a^k
\hbox{ implies }\a=\b\}$.

$0\in A$: From $\suml_{k=0}^0\a_k\cdot a^k=\suml_{k=0}^0\b_k\cdot a^k$ we infer
$\a_0\nad{\er{2.13}}= \a_0\cdot 1\nde2.23 = \a_0\cdot a^0 \nde7.10 =
\suml_{k=0}^0 \a_k\cdot a^k = \suml_{k=0}^0\b_k\cdot a^k = \b_0\cdot a^0 =\b_0
\cdot1=\b_0$. Hence $\a_0=\b_0$, and $0\in A$.

\ti{$N\in A$ implies $N+1\in A$}: We assume $\suml_{k=0}^{N+1}\a_k\cdot a^k =
\suml_{k=0}^{N+1}\b_k\cdot a^k$. We have to show $\a=\b$. $\suml_{k=0}^{N+1}
\a_k\cdot a^k\nde7.13 = \suml_{k=0}^0 \a_k\cdot a^k + \suml_{k=1}^{N+1}\a_k
\cdot a^k =\a_0+\suml_{k=1}^{N+1}\a_k\cdot a^k \nde1.6 = \bigl(\suml_{k=1}^{N+1}
\a_k\cdot a^k\bigr) +\a_0 \nde7.12 = \bigl(\suml_{l=0}^N\a_{l+1}\cdot a^{l+1}
\bigr) +\a_0 \nde2.25 = \bigl(\suml_{l=0}^N \a_{l+1}\cdot (a^1\cdot a^l)\bigr)
+\a_0$.

Note that $\a_{l+1}\cdot(a^1\cdot a^l)\nde2.11 = (\a_{l+1}\cdot a^1)
\cdot a^l \nde2.12 = (a^1\cdot\a_{l+1})\cdot a^l \nde2.11 = \a^1\cdot
(\a_{l+1}\cdot a^l)\nde2.23 = a\cdot(\a_{l+1}\cdot a^l)$. Hence $\suml_{l=0}^N\a_{l+1}
\cdot(a^1\cdot a^l) = \suml_{l=0}^Na\cdot(\a_{l+1}\cdot a^l) \nde7.27 = a\cdot
\suml_{l=0}^N\a_{l+1}\cdot a^l \nde2.12 = \bigl(\suml_{l=0}^N \a_{l+1}\cdot
a^l\bigr)\cdot a$. \E\Tf $\suml_{k=0}^{N+1}\a_k\cdot a^k =\bigl(\suml_{l=0}^N
\a_{l+1}\cdot a^l\bigr)\cdot a+\a_0$. Similarly $\suml_{k=0}^{N+1}\b_k\cdot a^k
=\bigl(\suml_{l=0}^N \b_{l+1}\cdot a^l\bigr)\cdot a+\b_0$. Equating the two
\LHS s, we obtain $\bigl(\suml_{l=0}^N \a_{l+1}\cdot a^l\bigr)\cdot a+\a_0 =
\bigl(\suml_{l=0}^N\b_{l+1}\cdot a^l\bigr)\cdot a+\b_0$. We can apply Theorem
\rf{t5.1} and we obtain $\suml_{l=0}^N \a_{l+1}\cdot a^l = \suml_{l=0}^N
\b_{l+1}\cdot a^l$ and $\a_0=\b_0$. Since $N\in A$, we obtain from the first
\et y $\a_{l+1}=\b_{l+1}$ for $l\in[0,N]$, which together with $\a_0=\b_0$
implies $\a=\b$. Hence $N+1\in A$. \E\Tf $A$~is an \iv\ subset of~$E$, so
$A=E$, by \er{2.2}. This completes the proof of part (ii) and of Lemma \rf{l5.4}.
\endproof

\brm5.5
\If from the proof of \er{5.7} that $\suml_{k=0}^N \a_k\cdot a^k < a^{N+1}$
whenever $\a_k\in\zo0,a $, $k\in[0,N]$. From the proof of \er{5.8} we obtain
$\a_N\cdot a^N\le \suml_{k=0}^N\a_k\cdot a^k$ whenever $\a_k\in E$,
$k\in[0,N]$. Hence if $\a_k\in\zo0,a $, $k\in[0,M]$ and $\a_N\ge1$, then
$\suml_{k=0}^N\a_k\cdot a^k\in\zo a^N,a^{N+1} $.
\erm

Motivated by Lemmata \rf{l5.3} and \rf{l5.4}, we introduce the \fw\ \df\ and
notation.

Given $a\in E^*\sms1$, and $N\in E$, we define sets $E\sbn a{,N}$, $E\sbn a{}$,
by setting
\bea5.9
E\sbn a{,N}&:= \{\a:[0,N]\to\zo0,a \hbox{ \st}\a(N)\ne0\},\\
E\sbn a{}&:=\{0\}\cup \bigcup_{N\in E}E\sbn a{,N}. \lb{5.10}
\e
\Mo we define a map $\vf\sbn a{}:E\sbn a{}\to E$ by setting
\beq5.11
\vf\sbn a{}(0):=0 \qh{and }\vf\sbn a{}(\a):=\suml_{k=0}^N\a_k\cdot a^k
\qh{\fe $\a\in E\sbn a{,N}$ and every $N\in E$}.
\e

\bth5.6
Let $(E,0,S)$ be a set of \nn s, let $1:=S(0)$, $E^*:=E\sms0$, and let
$a\in E^*\sms1$. Then the map $\vf\sbn a{} : E\sbn a{}\to E$ defined in
\er{5.11} is bi\jc.

\Mo if $S\sbn a{}$ denotes the selfmap of $E\sbn a{}$ defined by
\beq5.12
S\sbn a{}:=\vf\sbn a{}\Inv \circ S\circ\vf\sbn a{},
\e
then $(E\sbn a{},0,S\sbn a{})$ is a set of \nn s.\glossary{$(E\sbn a{},0,S\sbn a{})$}
\eth

\brm5.7
We made an abuse of notation when using the same symbol $0$ for the
distinguished \el\ of the set~$E$ and of the set $E\sbn a{}$. One could
replace $(E,0,S)$ by $(E,e,S)$ or $\ees$, as we did until Section~\ref{sss.Mult}
in order to avoid using the same symbol for the \nel\ of a monoid. Similarly,
$1$~could be replaced by $S(e)$, or $\wt S(\wt e)$. The main motivation for
using~$0$ as the distinguished \el\ of a set of \nn s is the fact that the
symbol~$0$ is commonly used for the \ca\ of the empty set (and the symbol~$1$
is used for the \ca\  of a~set consisting of one \el).
\erm

\proof[Proof of Theorem \rf{t5.6}]
The \ti{sur\ji} of $\vf\sbn a{}$ is a direct con\sq\ of Lemma \rf{l5.3} and
the fact that $\vf\sbn a{}(0)=0$.

\ti{In\ji}: If $\a\in E\sbn a{}\sms0$, then there is $N\in E$ \st $\a$~is
a~map from $[0,N]$ into $\zo0,a $ with $\a(N)\ne0$. \E\Tf $\vf\sbn a{}(\a):=
\suml_{k=0}^N \a_k\cdot a^k\ge a^N\ne0$ by \er{5.8}, \er{2.31}. Hence
$\vf\sbn a{}(\a)\ne\vf\sbn a{}(0)$. If $\a,\b\in E\sbn a{}\sms 0$ and
$\vf\sbn a{}(\a)=\vf\sbn a{}(\b)$, then $\a=\b$ by Lemma \rf{l5.4}. \If that
$\vf\sbn a{}$ is in\jc, hence bi\jc.

We now show that $S\sbn a{}$ is a bi\jn\ from $E\sbn a{}$ onto $E\sbn a{}
\sms0$. In \er{5.12} $\vf\sbn a{}:E\sbn a{} \to E$ is a bi\jn\ and $S:E\to E\sms0$
is a bi\jn. Thus $S\circ\vf\sbn a{}:E\sbn a{}\to E\sms0$ is a bi\jn. Since $\vf_a(0)=0$
by \er{5.11}, the \rt ion of~$\vf\sbn a{}$ to $E\sbn a{}\sms0$ is a bi\jn\ from
$E\sbn a{}\sms0$ onto $E\sms0$. \E\Tf the \rt ion of $\vf\sbn a{}\Inv$ to $E\sms0$
is a bi\jn\ from $E\sms0$ onto $E\sbn a{}\sms0$. Thus $S\sbn a{} = \vf\sbn a{}\Inv
\circ S\circ\vf\sbn a{}$ is a bi\jn\ from $E\sbn a{}$ onto $E\sbn a{}\sms0$, which
proves (1.N1). We now
show that (1.N2) is \sf ied. Let $A\sbs E\sbn a{}$ \sf y $0\in A$ and
$S\sbn a{}(A)\sbs A$. We claim that $\vf\sbn a{}(A)$ \sf ies $0\in\vf\sbn a{}
(A)$ and $S(\vf\sbn a{}(A))\sbs \vf\sbn a{}(A)$. Indeed, $0\nde5.11 =
\vf\sbn a{}(0)$, hence $0\in\vf\sbn a{}(A)$, since $0\in A$. From \er{5.12}
we obtain $\vf\sbn a{}\circ S\sbn a{} = S\circ \vf\sbn a{}$. Hence
$S(\vf\sbn a{}(A)) = S\circ \vf\sbn a{}(A) = (\vf\sbn a{}\circ S\sbn a{})(A) =
\vf\sbn a{}(S\sbn a{}(A))\sbs\vf\sbn a{}(A)$, since $S\sbn a{}(A)\sbs A$. From (1.N2) we infer
$\vf\sbn a{}(A)=E$, hence $A=\vf\sbn a{}\Inv(E)=E\sbn a{}$. \csq\
(1.N2)  holds and $(E\sbn a{},0,S\sbn a{})$ is a~set of \nn s.
\endproof

\bnt5.8
Let $(E\sbn a{},0,S\sbn a{})$ be the set of \nn s introduced in Theorem
\rf{t5.6}, and let $\vf\sbn a{}:E\sbn a{}\to E$ be as in \er{5.11}. We shall
denote by $\psi\sbn a{}$ the inverse of~$\vf\sbn a{}$, by $\le\sbn a{}$ the
\nog\ of $(E\sbn a{},0,S\sbn a{})$ defined in \E\df\ \rfa1{d3.24}, by
$+\sbn a{}$ the \ad\ on $E\sbn a{}$ defined in \E\Pr\ \rf{p1.2}, and
by~$\cdot\sbn a{}$ the \mlc\ on $E\sbn a{}$ defined in \E\df\ and Notation
\rf{d2.3}.
\ent

\bth5.9
Under the \as s of Theorem \rf{t5.6}, the \fw\ holds\dw
\bga5.13
\vf\sbn a{}(0)=0, \q \vf\sbn a{}(1)=1,\\
\vf\sbn a{}(\a+\sbn a{} \b)=\vf\sbn a{}(\a)+\vf\sbn a{}(\b),\q \a,\b\in
E\sbn a{}, \lb{5.14} \\
\vf\sbn a{}(\a\cdot\sbn a{}\b)=\vf\sbn a{}(\a)\cdot \vf\sbn a{}(\b), \q
\a,\b\in E\sbn a{}. \lb{5.15}
\e
The map $\vf\sbn a{}$ is an \ois sm from $(E\sbn a{},\le\sbn a{})$ onto
$(E,\le)$.
\eth

\proof \

\er{5.13}: $\vf\sbn a{}(0)=0$ by \er{5.11}; observe that $1\in E\sbn a{,0}$,
hence $\vf\sbn a{}(1)\nde5.11 = \suml_{k=0}^0 1\cdot a^k \nde7.10 = 1\cdot a^0
\nad{\er{2.13},\er{2.23}}=1$.

\er{5.14}: We define a binary \op\ $\qu\sbn a{}$ on $E\sbn a{}$ by setting for
$\a,\b\in E\sbn a{}$:
\beq5.16
\a\qu\sbn a{}\b:= \psi\sbn a{}(\vf\sbn a{}(\a)+\vf\sbn a{}(\b)).
\e
We claim that $\qu\sbn a{}$ \sf ies \er{1.2}, \er{1.3} with $\ees:=(E\sbn a{},
0,S\sbn a{})$. Indeed, given $0,\a,\b\in E\sbn a{}$, we have: $\a\qu\sbn a{}0
\nde5.16 = \psi\sbn a{}(\vf\sbn a{}(\a)+\vf\sbn a{}(0)) \nde5.13 =
\psi\sbn a{}(\vf\sbn a{}(\a)+0)\nde1.2 = \psi\sbn a{}(\vf\sbn a{}(\a))=
\id_{E\sbn a{}}(\a)=\a$.

\Mo $\a\qu\sbn a{}S\sbn a{}(\b) = \psi\sbn a{}\bigl(\vf\sbn a{}(\a)+\vf
\sbn a{}(S\sbn a{}(\b))\bigr) \nde5.12 = \psi\sbn a{}
\bigl(\vf\sbn a{}(\a) +S(\vf\sbn a{}
(\b))\bigr)\break \nde1.3 = \psi\sbn a{}\bigl(S(\vf\sbn a{}(\a)+\vf\sbn a{}(\b))
\bigr) \nde5.12 = S\sbn a{}(\psi\sbn a{}(\vf\sbn a{}(\a)+\vf\sbn a{}(\b)))
\nde5.16 = S\sbn a{}(\a\qu\sbn a{}\b)$. Thus the claim is proved. As a~con\sq\
of the \uq\ part of \E\Pr\ \rf{p1.2}, and the \df\ of $+\sbn a{}$, we obtain
$\a+\sbn a{}\b = \a\qu\sbn a{}\b$, $\a,\b\in E\sbn a{}$. Hence $\vf\sbn a{}
(\a+\sbn a{}\b)\nde5.16 = \vf\sbn a{}\bigl(\psi\sbn a{}(\vf\sbn a{}(\a)
+\vf\sbn a{}(\b))\bigr) =\vf\sbn a{}(\a)+\vf\sbn a{}(\b)$, which implies
\er{5.14}.

\er{5.15}: As in the proof of \er{5.14}, we define $\qu\sbn a{}: E\sbn a{}
\t E\sbn a{} \to E\sbn a{}$ by setting
\beq5.17
\a\qu\sbn a{}\b := \psi\sbn a{}(\vf\sbn a{}(\a)\cdot \vf\sbn a{}(\b)), \q
\a,\b\in E\sbn a{}.
\e
We claim that $\a\qu\sbn a{}\b$ \sf ies \er{2.9}, \er{2.10}, with $(E,0,S):=
(E\sbn a{},0,S\sbn a{})$. Indeed, given $0,1,\a,\b\in E\sbn a{}$, we have:
$\a\qu\sbn a{}0 \nde5.17 = \psi\sbn a{}(\vf\sbn a{}(\a)\cdot \vf\sbn a{}(0))
\nde5.13 = \psi\sbn a{}(\vf\sbn a{}(\a)\cdot0)\nde2.9 = \psi\sbn a{}(0)$. Since
$\psi\sbn a{}$ is the inverse of $\vf\sbn a{}$ which is bi\jc, and $\vf\sbn a{}
(0)=0$, we have $\psi\sbn a{}(0)=0$. \Mo $\a\qu\sbn a{}(\b+\sbn a{}1)\nde5.17
= \psi\sbn a{}(\vf\sbn a{}(\a)\cdot \vf\sbn a{}(\b+\sbn a{}1)) \nde5.14 =
\psi\sbn a{}\bigl(\vf\sbn a{}(\a)\cdot (\vf\sbn a{}(\b)+\vf\sbn a{}(1))\bigr)
\nde5.13 = \psi\sbn a{}\bigl(\vf\sbn a{}(\a)\cdot(\vf\sbn a{}(\b)+1)\bigr)
\nde2.10 = \psi\sbn a{}\bigl((\vf\sbn a{}(\a)\cdot\vf\sbn a{}(\b)\bigr)+
\vf\sbn a{}(\a)\bigr)$.

Note that in view of \er{5.13}, \er{5.14}, the bi\ji\ of~$\vf\sbn a{}$, and
Lemma \rf{l1.8}, $\vf\sbn a{}$~is an \is sm from $(E\sbn a{},+\sbn a{},0)$
onto $(E,+,0)$. \E\Tf $\psi\sbn a{}:=\vf\sbn a{}\Inv$ is an \is sm from
$(E,+,0)$ onto $(E\sbn a{},+\sbn a{},0)$. \If that $\psi\sbn a{}\bigl((
\vf\sbn a{}(\a)\cdot \vf\sbn a{}(\b)) + \vf\sbn a{}(\a)\bigr) = \psi\sbn a{}
(\vf\sbn a{}(\a)\cdot \vf\sbn a{}(\b)) +\sbn a{} \psi\sbn a{}(\vf\sbn a{}(\a))
\nde5.17 = \a\qu\sbn a{}\b +\sbn a{} \a$. Thus the claim is proved.

As a con\sq\ of the \uq\ part of \E\Pr\ \rf{p2.4}, and of the \df\
of~$\cdot\sbn a{}$, we obtain $\a\qu\sbn a{}\b = \a\cdot\sbn a{}\b$ \fa $\a,\b
\in E\sbn a{}$. Hence \er{5.15} follows as in the proof of \er{5.14}.
Concerning the last assertion of the theorem, we observe that since
$(E\sbn a{},{\le}\sbn a{})$ is totally ordered and $\vf\sbn a{}$ is bi\jc, it
is \sft, in view of Lemma \rfa1{l3.34} and \E\df\ \rfa1{d3.33}, to show that
$\vf\sbn a{}$ is in\cre. Let $\a,\b\in E\sbn a{}$ with $\a\le\sbn a{}\b$. In
view of \er{1.46}, \te s $\g\in E\sbn a{}$ \st $\b=\a+\sbn a{}\g$. Then
$\vf\sbn a{}(\b) \nde5.14 = \vf\sbn a{}(\a)+\vf\sbn a{}(\g)$. Hence by
\er{1.46} again, we have $\vf\sbn a{}(\a)\le\vf\sbn a{}(\b)$. This
completes the proof of Theorem \rf{t5.9}.
\endproof

\brm5.10 \

\hph i,i, \If from Theorem \rf{t5.6}, \er{5.13}, \er{5.14} and Lemma \rf{l1.8}
that $\vf\sbn a{}$ is an \is sm from the monoid $(E\sbn a{},+\sbn a{},0)$ onto
the monoid $(E,+,0)$.

\hph ii,, \Mo it follows from \er{5.13}, \er{5.15} that $\vf\sbn a{}$
is an \is sm from the monoid $(E\sbn a{},\cdot\sbn a{},1)$ onto the monoid
$(E,\cdot,1)$.
\erm

\bex5.11
Let $\ees$ be a set of \nn s. Construct a~set of \nn s $(E,0,S)$ from the set
$\ees$ \st $1:=S(0)$, $2:=S(1)$, $3:=S(2)$, $4:=S(3)$, $5:=S(4)$, $6:=S(5)$,
$7:=S(6)$, $8:=S(7)$, $9:=S(8)$, $10:=S(9)$. Proceed as in Lemma \rf{l0.27}.

Let $(E_{(10)},0,S_{(10)})$ be the set of \nn s introduced in Theorem \rf{t5.6}
with $a:=10$. Show:

\hph i,ii, $E_{(10)}$ is what is usually denoted by $\N=\{0,1,2,\dots\}$.\glossary{$\N$}

\hph ii,i, $S_{(10)}n=n+_{(10)}1$ \fa $n\in\N$.

\hph iii,, Let $m:=\suml_{k=0}^N \a_k\cdot10^k$ where $N\in\N$, $\a:[0,N]\to
D_{10}$, $\a(N)\in\dot D_{10}$. Let $n:=m+1$ and let $M\in\N$, $\b:[0,M]\to
D_{10}$, $\b(M)\in\dot D_{10}$, be \st $n=\suml_{l=0}^M \b_l\cdot10^l$.

{\setbox0=\hbox{(iii)(1)}\leftskip\wd0 \setbox1=\hbox{(1) }\hangindent\wd1
\leavevmode\llap{(1) }Show that if $\a_k:=9$ \fa $k\in[0,N]$, then $M=S(N)$,
$\b(M)=1$ and $\b_l=0$ \fa $l\in[0,N]$. (Hint: use \er{7.31}, \er{7.27}.)

\leavevmode\llap{(2) }Show that if \te s $\ov k\in[0,N)$ \st $\a_k:=9$ \fa
$k\in[0,\ov k]$,\break $\a(S(\ov k))\ne9$, then $M=N$, $\b_l=0$ for $l\in[0,\ov k]$,
$\b_{S(\bar k)}=S(\a_{S(\bar k)})$, $\b_l=\a_l$ for $l>S(\ov k)$.

\leavevmode\llap{(3) }Show that if $\a_0:\ne9$, then $\b_0:=S(\a_0)$ and
$\b_l=\a_l$ \fa $l\in\oz0,N $.

}
\hph iv,, Let $m:=1943= 1\cdot10^3 + 9\cdot10^2 + 4\cdot10^1 + 3\cdot10^0$ and
let $n:=75=7\cdot10^1 + 5\cdot10^0$. Show that $m+n=2018$ and $m\cdot n=
145725$.

\hph v,i, $n+_{(10)}m$ (resp.\ $n\cdot_{(10)}m$) is the usual \ad\ (resp.\
\mlc) of $m$~and~$n$ in~$\N$, $m,n\in\N$ (use \er{1.2}, \er{1.3} and \er{2.9},
\er{2.10}).
\eex

\bnt5.12
In the sequel we shall use the notation $(\N,0,S)$ or simply~$\N$ instead of
$(E_{(10)},0,S_{(10)})$ introduced in Exercise \rf{ex5.11}.

The \ad\ (resp.\ \mlc) $+_{(10)}$ (resp.\ $\cdot_{(10)}$) will be simply
denoted by~$+$ (resp.~$\cdot$).
\ent

\brm5.13
The set $(E_{(2)},0,S_{(2)})$ introduced in Theorem \rf{t5.6} with $(E,0,S):=
(\N,0,S)$, $a:=2$, is usually called the \ti{set of binary \rp ations of \el
s of~$\N$}.\index{binary representation of elements of $\N$}
\erm

We conclude this section by collecting some \pp ies of~$\N$ obtained in
Chapters 1~and~2.

\bpr5.14 \

\hph i,vii, $(\N,+,0)$ is an infinite \pn\ monoid with $1$ as \Gn. \E\Ip $(\N,+,0)$ is
a \PM.

\hph ii,vi, Let $m,n\in \N$ and let $m\dpl n$ denote the $m$-fold \IT\ of~$n$
in the monoid $(\N,+,0)$ \wrt the set of \nn s $(\N,0,S)$. Let $\cdot$ denote
the binary \op\ in~$\N$ defined by $m\cdot n:=m\dpl n$, $m,n\in \N$. Then
$(\N,\cdot,1)$  is an \am.

\hph iii,v, Let $\Na:=\N\sms0$. Then $\Na$ is a \sbm\ of $(\N,\cdot,1)$, and
the monoid $(\Na,\cdot,1)$ is a \PM.

\hph iv,ii, Let $m\in\N$ and let $\d_m:\N\to\N$ denote the map defined by $\d_m
(n):=m\cdot n$, $n\in\N$. Then $\d_m$~is an endo\mf\ of the monoid $(\N,+,0)$.
\Mo every endo\mf~$\vf$ of $(\N,+,0)$ is equal to $\d_m$ with $m:=\vf(1)$.

\hph v,iii, The \nog\ $\le$ of the \PM\ $(\N,+,0)$ is a well-\og, \Ip a total
\og.

\hph vi,ii, If $A$ is a subset of $\N$ \sf ying $0\in A$, and $n+1\in A$ whenever $n\in A$,
then $A=\N$.

\hph vii,i, The map $\th_m:(\N,\le) \to(\N,\le)$, $m\in\N$, defined by $\th_m
(n):=m+n$, $n\in\N$, is in\jc\ and strictly in\cre. $\th_m$~is sur\jc\ iff
$m=0$.

\hph viii,, The map $\d_m$ defined in {\rm(iv)} is in\jc\ and strictly in\cre\ iff
$m\ge1$. $\d_m$~is sur\jc\ iff $m=1$.
\epr

%% file: DETOUR3.TEX
\def\h{\hskip2pt minus0.3pt}
\Section{Divisibility}[Divisibility]\label{s.3}
\Subsubsection{Least common multiple and greatest common divisor}\label{sss.lcm}

We first observe that the \mlc\ defined in \era2{2.8} with $E:=\N$ is \asc e
by \era2{2.11} and \cmt e by \era2{2.12}. \Mo $1$~is the \nel\ by \era2{2.13}.
Thus $(\N,\cdot,1)$ is an abelian monoid. It \sf ies \era2{1.9}
with ${\qu}:={\cdot}$ by \E\Pr\ \rfa2{p2.7}\,(iv). It does not \sf y \era2{1.8},
since $2\cdot0 \nda2{2.9} = 0 =3\cdot 0$ does not imply $2=3$. Thus
$(\N,\cdot,1)$ is \ti{not\/} a~\PM. The subset $\Na(:=\N\sms0)$ of~$\N$ is
a~\sbm\ of $(\N,\cdot,1)$, since $1\in\Na$ and $a,b\in\Na$ implies $a\cdot b
\in\Na$ by \E\Pr\ \rfa2{p2.7}\,(i).

In view of \E\Pr\ \rfa2{p2.7}\,(iii), the abelian monoid $(\Na,\cdot,1)$ \sf
ies not only \era2{1.9}, but also \era2{1.8}. Hence $(\Na,\cdot,1)$ is a~\PM.

In what follows we shall also refer to \era2{1.5}--\era2{1.9} for \pp ies of the
\mlc~$\cdot$.

\smallskip

The set $(\N,0,S)$ is a set of \nn s with \nog~$\le$ (see \E\df\
\rfa1{d3.24}). We shall also denote by~$\le$ the \rt ion to~$\Na$ of the
\nog\ of~$\N$. Thus $(\Na,\le)$ becomes an \os\ (see Remark \rfa1{r3.3}).

\bex1.1
Let $(\Na,\le)$ be as above. Let $S'$ denote the \rt ion to~$\Na$ of the
\su\ \f~$S$, which is well-defined since $R(S)=\Na$. Then $(\Na,1,S')$ is
a~set of \nn s (see the discussion preceding \E\df\ \rfa1{d1.3}). Denote
by~$\le'$ the \nog\ of $(\Na,1,S')$. Show that $\le$ and $\le'$ are identical.
\eex

\begin{dfn}[{\cite[p.~319]{Lattice}}]\lb{d1.2}
An \ti{\oam\/} $(X,\qu,e,\le)$ is an abelian monoid $(X,\qu,e)$ with an\index{ordered abelian monoid}
ordering $\le$ \sf ying the \cn
\beq8.1
a\le b \qhq{implies}a\qu c \le b\qu c \q\hbox{for all }a,b,c \in X.
\e
Two \oam s are called \ti{\is c} if they are both monoid- and order-\is c.
\edf

\If from \era2{2.19} that both $(\N,\cdot,1,\le)$ and $(\Na,\cdot,1,\le)$
are \ti{ordered\/} abelian monoids (notice $0\cdot x=0$, $x\in\N$).

\bex1.3 \

\hph i,i, Let $X:=\{0,1\}$ and let $\lor$ (resp.~$\land$) be defined by\glossary{$\lor$}\glossary{$\land$}
\bea8.2
&0\lor0=0 \q\q 0\lor1=1\lor 0=1 \q\q 1\lor1=1,\\
&0\land0=0 \q\q 0\land1=1\land 0=0 \q\q 1\land1=1, \lb{8.3}
\e
and let $\le$ be the only ordering on $X$ \sf ying $0< 1$. Show that both
$(X,\lor,0)$ and $(X,\land,1)$ are \oam s.

\hph ii,, Let $X:=\{0,1\}$ and let ${\qu}:X\t X\to X$ be defined by $0\qu0=1\qu
1=0$, $0\qu1=1\qu0=1$. Let $\le$ be the same ordering as in (i). Show that
$(X,\qu,0)$ is an abelian monoid but not an \oam\ \wrt $\le$.
\eex

We now introduce \ti{another} ordering on $(\Na,\cdot,1)$. We first consider
a more general situation.

\advance\abovedisplayskip by-2pt
\advance\belowdisplayskip by-2pt
\bpr1.4
Let $(X,\qu,e)$ be a \PM.
Then the \rl\ on $X$ defined by\glossary{$\lequ$}
\bml8.7
x\lequ y,\ x,y\in X, \hbox{ if \te s } z\in X \\ \hbox{\st} y=x\qu z,\
\hbox{or \ev tly }y=z\qu x,
\e
is an order \rl\ which makes $(X,\qu,e,\lequ)$ an \oam. \Mo if
\beq8.9
a \stackrel\qu< b \ \hbox{means } a\lequ b \hbox{ and }a\ne b,
\e
then
\beq8.8
e \stackrel\qu< x \qh{\fe}x\in X\sm\{e\},
\e
and
\beq8.10
a \stackrel\qu< b \qh{implies} a\qu c \stackrel\qu< b\qu c \qh{\fa} a,b,c\in X.
\e
The order \rl\ \smash{$\lequ$} is sometimes called the \emph{natural\/} \og\ of the \PM\ $(X,\qu,e)$. 
\epr

\proof

(i) \smash{$\lequ$} is an \ti{order } on $X$. For simplicity we just write $\le$
instead of \smash{$\lequ$}.

\ti{\E\tr ity}: Let $x,y,z\in X$ be \st $x\le y$ and $y\le z$. Then, by
\er{8.7}, \te\ $p,q\in X$ \st $y=x\qu p$, $z=y\qu q$. Hence $z=(x\qu p)\qu q
=x\qu(p\qu q)$ in view of the \asc ity of~$\qu$.

\ti{Reflexivity}: $x=x\qu e$ \fe $x\in X$ by \df\ of~$e$. Hence $x\le x$.

\ti{Antisymmetry}: Let $x,y\in X$ be \st $x\le y$ and $y\le x$. Then \te\
$p,q\in X$ \st $y=x\qu p$ and $x=y\qu q$. Hence
$e\qu x\nda2{1.7} = x= (x\qu p)\qu q \nda21.5 = x\qu(p\qu q) \nda21.6 =
(p\qu q)\qu x$.
It follows from \era2{1.8} that $e=p\qu q$ and from \era2{1.9}
that $p=q=e$. Hence $y=x\qu e=x$.

(ii) $(X,\le)$ is an \ti{\oam\/} and \er{8.10} holds. It suffices to prove
that \er{8.10} holds. Let $a,b\in X$ be \st $a<b$. Then \te s $p\in X$ \st
$b=a\qu p$ and $b\ne a$. In view of the \df\ of~$e$, and of~\era2{1.8}, $b\ne a$
iff $p\ne e$. Then $b\qu c=(a\qu p)\qu c=a\qu (p\qu c)=a\qu (c\qu p)$ in view
of the \cmt ity of~$\qu$. \Mo $a\qu(c\qu p)=(a\qu c)\qu p$. Since $p\ne e$,
$a\qu c<b\qu c$.

(iii) \er{8.8} follows from $x= e\qu x$ \fe $x\in X\sm\{e\}$.
\endproof

\blm8.13
Let $(X,\qu,e)$ and $(\wt X,\tqu,\wt e)$ be P-monoids and let $\vf:X\to \wt X$.

\hph i,i, If $\vf$ \sf ies
\beq8.16
\vf(x\qu y)=\vf(x)\tqu \vf(y) \qh{\fe}x,y\in X
\e
then $\vf$ is a \emph{monoid-homo\mf} and an \emph{order-homo\mf} \wrt
the \og s $\lequ$ and $\stackrel{\tqu}\le$.

\hph ii,, If in \ad\ $\vf$ is \emph{bi\jc}, then $\vf$ is a \emph{monoid-\is sm}
as well as an \emph{\ois sm}.
\elm

\proof \

(i) We have to show that $\vf(e)=\wt e$. We have $\wt e \tqu\vf(e) = \vf(e)
=\vf(e\qu e) \nde8.16 =\break\vf(e)\tqu \vf(e)$. Then $\wt e=\vf(e)$ by~\era2{1.8}.
We now prove that $\vf$ is in\cre. Let $x,y\in X$ be \st $x\lequ y$, that is,
\te s $z\in X$ \st $y=x\qu z$. Then $\vf(y)=\break \vf(x\qu z)\nde8.16 =
\vf(x)\tqu\vf(z)$. Hence $\vf (x)\stackrel{\tqu}\le \vf(y)$.

(ii) follows from Lemma~\rfa2{l1.8}\,(i) and part (i).
\endproof

\brm8.9
The terminology P-monoid which stands for positive monoid is not standard.
Since the term \ti{positive monoid\/} appears in the literature with a
slightly different meaning, we have chosen to call such a monoid a P-monoid.

The notion of P-monoid is a special case of the notion of divisibility monoid
treated in \cite{Lattice}.
\erm

\advance\abovedisplayskip by2pt
\advance\belowdisplayskip by2pt
\brm1.7
\If from \era2{1.46} in Lemma \rfa2{l1.23} that the \og\ $\stackrel+\le$ in
$(\N,+,0)$ is identical to the \nog\ of the set of \nn s $(\N,0,S)$.
\erm

\def\namespec{Definitions and Notations}
\begin{dspc}\lb{d1.8}
The \og\ $\le$ on the \PM\ $(\Na,\cdot,1)$ defined in~\er{8.7}
is usually denoted by $|$ (or~$/$), i.e., given $a,b\in\Na$,\glossary{$|$}
\bml8.52
a|b \qh{if \te s $c\in\Na$ \st $b=c\cdot a$
or $b=a\cdot c$,}\\
\hbox{and $a\nmid b$ if $a|b$ does not hold.}
\e

If $a|b$ (pronounce: \ti{$a$ divides $b$}), then
\beq8.53
\hbox{the \ti{quotient\/} $\frac ba$ ($b$ \ti{over} $a$) is the only \el\
of $\Na$ \st $b=\frac ba\cdot a$}.
\e
Then $b$ is called the \ti{numerator} and $a$ the \ti{denominator} of the
quotient $\frac ba$.\glossary{$\frac ab$}

Note that if $a|b$ and if $a,b$ are viewed as \el s of~$\N$, then $b$~is
the $\frac ba$-th \IT\ of~$a$ \wrt the \ad, hence a \ti{\ml e} of~$a$ (see
Theorem \rfa2{t1.38}). If $a|b$, then $a$~is usually called a \ti{divisor
of}~$b$. We set
\bea8.54
M(a)&:=\{b\in\Na: a|b\} \qh{for} a\in\Na,\\
D(b)&:=\{a\in\Na: a|b\} \qh{for} b\in\Na. \lb{8.55}
\e
Hence, obviously
\beq8.56
a|b \iff a\in D(b) \iff b\in M(a), \ a,b\in\Na.
\e
\end{dspc}

We have defined two \og s in $\Na$, namely $\le$ and $|$.
Recall that $(\Na,\cdot,1)$ is an \oam\ for $\le$. The same holds for~$|$
by \E\Pr\ \rf{p1.4}. Note that $1$ is the least \el\
of~$\Na$ for both \og s (indeed $1|a$ \fe $a\in\Na$ since $a=1\cdot a=
a\cdot 1$). However, $(\Na,\le)$ and $(\Na,|)$ are \ti{not\/} \ois c,
but one \og\ is included in the other one. More precisely, we have
\beq8.57
a|b \qhq{implies} a\le b \qh{\fa } a,b\in\Na.
\e
Indeed, $a|b$ implies $b=c\cdot a$ for some $c\ge1$. Hence $b=c\cdot a
\ge 1\cdot a=a$ by \era2{2.19}, \era2{2.13}.

The converse is not true. Indeed, $2<3$ but $3\notin M(2)$. Indeed, if
$3=2\cdot p$ for some $p\in\Na$, then either $p=1$ or $p\ge 2$. But
$2\cdot 1=2<3$, and $p\ge2$ implies $2\cdot p\ge 2\cdot 2=2+2=(2+1)+1=3+1
>3$. \E\Tf the identity map from $(\Na,|)$ onto $(\Na,\le)$ is a nice
example of a bi\jn\ between two \os s, the latter being even totally
ordered, which is in\cre\ and \st its inverse is \ti{not\/} in\cre.
(See Lemma \rfa1{l3.34} and the example following it.)

\smallskip
We recall that the \og\ $\le$ of the \PM\ $(\N,+,0)$ is a well-\og\ (every
\nss\ of~$\N$ has a least \el), and that every nonempty bounded (\ev tly
finite) subset of~$\N$ has a~greatest \el. None of these \pp ies holds
in~$(\N,\,|\,)$. Indeed, $2\nmid3$ by the \uq\ part of Theorem \rfa2{t1.38} with
$\ees:=(\N,0,S)$ since $3=1\cdot2+1$, and $3\nmid2$ by \er{8.57} since $2<3$.
\E\Tf the finite set $\{2,3\}$ has neither a~least nor a~greatest \el.
However, the \fw\ holds in the \os\ $(\Na,\,|\,)$. Since $1|n$ \fa $n\in\Na$ by
\era2{2.13}, the \nm~$1$ is the least \el\ of $(\Na,\,|\,)$, hence every \nss\
of~$\Na$ is bounded below. \E\Tf in $(\Na,\,|\,)$ as in $(\N,\le)$ the notions
``bounded above'' and ``bounded'' are \ev t.

\bnt1.9
Let $A$ be a nonempty finite subset of~$\Na$. In view of \era2{7.17} with
$(X,{\qu},e):=(\Na,\cdot,1)$, $\O:=A$, and $a:\O\to X$ being equal to the
\rt ion to~$A$ of the identity in~$X$, we shall write $\prodl_{n\in A}n$ for
the \cme product in $(\Na,\cdot,1)$ of all \el s of~$A$. Of course the dummy
index~$n$ can be replaced by any not already used letter.
\ent

\blm1.10
Let $A$ be a \nss\ of~$\Na$. The \fw\ assertions are \ev t\/{\rm:}

\hph i,ii, $A$ is bounded \wrt the \og~$|$,

\hph ii,i, $A$ is bounded \wrt the \og~$\le$,

\hph iii,, $A$ is finite.
\elm

\proof \

(i) \ti{implies} (ii): Follows from \er{8.57}.

(ii) \ti{implies} (iii): Let $m$ be an \ub\ of~$A$. Then $A\subset[0,m]$. As
a~subset of the finite set $[0,m]$, $A$~is finite by Theorem \rfa1{t4.18}\,(i).

(iii) \ti{implies} (i): Follows from the next lemma.
\endproof

\blm9.10
Let $A$ be a finite nonempty subset of\/ $\Na$. Then \fe $m\in A$,
\beq9.15
m\Bigm|\prod_{n\in A}n.
\e
\elm

\proof
By \In\ on $\#(A)\in\Na$.
If $A=\{k\}$, $k\in\Na$, then $\prodl_{n\in A}n=k$, so \er{9.15} holds.

If $\#(A)>1$ and $m\in A$, then
\[
\prod_{n\in A}n = \Bigl(\prod_{n\in A\sms m}n\Bigr)\cdot m =
m\cdot \Bigl(\prod_{n\in A\sms m}n\Bigr)
\]
by \era2{1.138}. Note that $\#(A)=\#((A\sms m)\cup\{m\})\nda23.31 = \#(A\sms m)
+\#(\{m\})= \#(A\sms m)+1$.
\endproof

Let $A$ be a \nss\ of an \os\ $(X,\le)$. If $A$ possesses a greatest
\el~$\hat m$ (resp.\ a~least \el~$\check m$), then by Lemma \rfa1{l3.38}
$\hat m$~(resp.~$\check m$)  is the least \ub\ (resp.\ greatest \ub) of~$A$.
Conversely, if $\hat m$ (resp.~$\check m)\in X$ is an \ub\ (resp. \lo) of~$A$,
and if $\hat m$~(resp.~$\check m$) belongs to~$A$, then clearly $\hat
m$~(resp.~$\check m$) is the greatest (resp.\ least) \el\ of~$A$.

\hbox to\textwidth{\vbox{\advance\hsize by-120pt
A bounded subset of an \os\ does not need to have a least \ub. For example,
if $X:=\{a,b,c,d,e\}$ and $a\le b$, $a\le c$, $a\le d$, $a\le e$, $b\le d$,
$b\le e$, $c\le d$, $c\le e$, then $(X,\le)$ is an \os, the subset $\{b,c\}$
has neither a greatest nor a~least \el, the \el~$a$ is the greatest \lo\ of
$\{b,c\}$, $\{d,e\}$ is the set of \ub s of $\{b,c\}$, and this set has no
least \el. Hence $\{b,c\}$ is bounded and does not possess a least \ub.
}\hfil \vbox{\hsize100pt
\[
\xymatrix{d&&e\\
b\ar[u]\ar[urr]&&c\ar[ull]\ar[u]\\
&a\ar[ur]\ar[uur]\ar[uul]\ar[ul]}
\]}}

\bex1.12 \

\hph i,i, Prove that the set $(X,\le)$ defined above is an \os.

\hph ii,, Prove analogous statements for the reverse order.

\eex

\If from Lemma \rfa1{l3.38} that every bounded subset of $(\N,\le)$ possesses
a~least \ub\ and a greatest \ub. It turns out that the same statement holds
for $(\Na,\,|\,)$.

\def\namespec{Definitions and Notations}
\begin{dspc}\lb{d1.12}
Let $A$ be a \nss\ of an \os\ $(X,\le)$. An~\el\ $x\in X$ is called the
\ti{least \ub\/}\index{least upper bound} or the \ti{supremum}\index{supremum} of~$A$, denoted by $\lub A$ or $\sup A$,
if the set of \ti{\ub s} of~$A$, denoted by $\UB(A)$, is not empty and $x$~is
the least \el\ of~$\UB(A)$. Similarly, $x\in X$ is called the \ti{greatest
\lo\/}\index{greatest lower bound} or the \ti{infimum}\index{infimum} of~$A$, denoted by $\glb A$ or $\inf A$, if $x$~is
the least \ub\ of~$A$ for the reverse \og\ on~$X$. We shall denote the set
of \lo s of~$A$ by $\LB(A)$.
\end{dspc}

We first prove that $\sup\{a,b\}$ exists in $(\Na,\,|\,)$ \fa $a,b\in\Na$,
$a\ne b$.  Let $a,b\in\Na$ with $a\ne b$.
Since the set $\{a,b\}$ is equal to the set $\{b,a\}$, we may suppose $a<b$
without loss of generality. Observe that
\beq1.16
\UB(\{a\})=M(a), \q \UB(\{a,b\})=M(a) \cap M(b).
\e
The first \et y is obvious and the second  holds since $M(a)\cap M(b)\ne\vn$.
Indeed, $ab\in M(a)\cap M(b)$. If $a|b$, for example if $a=2$ and $b=4$, then
any \ml e of~$4$ in~$\Na$ is a \ml e of~$2$. Hence $M(4)\subset M(2)$, and
$\UB(\{2,4\})=M(4)$.

\Mo $4|c$ \fa $c\in M(4)$, hence $4$ is the least \el\ of $M(4)$ in $(\Na,\,|\,)$,
and $4=\sup\{2,4\}$. More generally, if $a|b$ then $\UB(\{a,b\})=M(a)\cap
M(b)$ and $b=\sup\{a,b\}$. Note that $b=\sup\{a,b\}$ also in $(\N,\le)$. We
now consider the general case. Let $a,b\in\Na$, $a\ne b$. We first observe
that \fe $x\in\Na$, we have
\beq1.17
M(x)\cup\{0\} = I(x) \qh{and } M(x)\cap\{0\}=\vn,
\e
where $I(x)$ is the \pn\ \sbm\ of $(\N,+,0)$ generated by~$x$ (see Lemma
\rfa2{l1.21}). \If that $I(a)\cap I(b) = (M(a)\cup\{0\}) \cap (M(b)\cup\{0\})
= (M(a)\cap M(b)) \cup (\{0\}\cap M(b)) \cup (M(a)\cap\{0\}) \cup (\{0\}\cap
\{0\}) = (M(a)\cap M(b))\cup \{0\}$, since $M(a)\cap\{0\}$ and $M(b)\cap\{0\}$
are empty by \er{1.17}. Since $(M(a)\cap M(b))\cap\{0\} = M(a)\cap (M(b)\cap
\{0\}) \nde1.17 = M(a)\cap\vn =\vn$, we have shown that
\beq1.18
M(a)\cap M(b) = (I(a)\cap I(b))\sms0.
\e
One easily verifies that the intersection of a nonempty family of \sbm s of
a~monoid $(X,{\qu},e)$ is a possibly trivial \sbm\ of~$X$. \E\Ip $I(a)\cap
I(b)$ is a \sbm\ of $(\N,+,0)$. It is not trivial since $ab\in I(a)\cap I(b)$.
Not all nontrivial \sbm s of $(\N,+,0)$ are \pn\ (see Exercise \rfa2{ex1.44}
and Remark \rf{r1.15}). It turns out that the intersection of two \pn\ \sbm s
of~$(\N,+,0)$ is \pn. In the next lemma we give a~\ch ization of \pn\ \sbm s
of $(\N,+,0)$. We shall need the \fw

\def\namespec{Definition and Notations}
\begin{dspc}\lb{1.13}
Let $x,y\in\N$ be \st $x\le y$. Then\glossary{$y-x$}
\beq1.19
y-x:=p
\e
where $p$ is the only \el\ of~$\N$ \sf ying
\[
y=x+p\,(=p+x).
\]
\end{dspc}

\blm8.20
Let $A$ be a nontrivial \sbm\ of $(\N,+,0)$. Then $A=I(a)$ for some
$a\in\N\sms0$ iff \fa $x,y\in A$ \st $x\le y$, we have $y-x\in A$ where $y-x$ is
defined in \er{1.19}. In that case $a$ is the least \el\ of
$A\sms0$ \wrt the \og~$\le$, hence there is \emph{at most one} \el~$a$ of
$\N\sms0$ \st $A=I(a)$.
\elm

\proof
Let $A$ be a nontrivial \sbm\ of $(\N,+,0)$.

``\ti{Only if\/}'': Suppose \te s $a\in\Na$ \st $A=I(a)$.
Let $x,y\in A$ be \st $x\le y$. By \df\ \te\ $m,n\in\N$ \st
$x=n\dpl a\nda22.28 =n\cdot a$ and $y=m\dpl a\nda22.28 =m\cdot a$. It follows
from \era2{2.20} that $n\le m$. Then \te s $p\in\N$ \st $m=n+p$. In view of
\era2{2.3} we have $m\dpl a=n\dpl a+p\dpl a$, hence $y-x=m\dpl a-n\dpl a=p\dpl
a$. Thus $y-x\in I(a)=A$.

``\ti{If\/}'': Since $(\N,\le)$ is well-ordered (Theorem~\rfa1{t3.22}) and
$A\sms0\ne \vn$, \te s a least \el\ of $A\sms0$ which we denote by~$a$. Since
$A$ is a nontrivial monoid we have $I(a)\sbs A$ by Lemma~\rfa2{l1.21}. We claim
that $I(a)=A$. Let $b\in A$. Then by Theorem \rfa2{t1.38} \te s a unique pair
$(n,r)\in\N\t\N$ \st $b=n\cdot a+r$ with $r\in[0,a)$. If $r>0$, then
$b-n\cdot a\in A$ by \as. \E\Tf $0<r=b-n\cdot a\in A$ and $r<a$, which is
impossible since $a$ is the least \el\ of $A\sms0$. \E\Tf $r=0$ and $b=n\cdot
a\in I(a)$, which implies $A\subset I(a)$, hence $A=I(a)$. Since the map $\vf:(\N,\le) \to
(I(a),\le)$ defined in \E\Pr\ \rfa2{p1.13} is strictly in\cre\ by
Lemma~\rfa2{l1.21} and Lemma~\rfa2{l1.37}, and since $0=0\cdot a$, $a=1\cdot a$, we have $a<x$ \fe
$x\in I(a)\sms{0,a}=A\sms{0,a}$. It follows that $a$ is the least \el\
of~$A\sms 0$.
\endproof

\brm1.15
The subset $M$ of $\N$ defined by $M:=\{2\cdot m+3\cdot n: m,n\in\N\}$ is
a~\sbm\ of $(\N,+,0)$. Then $2=2\cdot1+3\cdot0$ belongs to~$M$ as well as
$3=2\cdot0+3\cdot1$. But $3-2=1\notin M$ since $2\cdot m+3\cdot n$, $m,n\in\N$,
is either equal to~$0$ or greater than~$1$. \E\Tf it is not \pn\ in view of
Lemma \rf{l8.20}.
\erm

\Wanp show that the intersection of two \pn\ \sbm s of $(\N,+,0)$ is \pn.

\bpr1.16
If $x,y\in\Na$, $x\ne y$, then $xy \in M(x)\cap M(y)$,
\beq1.20
I(x)\cap I(y)=I(\lcm(x,y)),
\e
where\glossary{$\lcm$}
\beq1.21
\lcm(x,y) \hbox{ is the least \el\  of $M(x)\cap M(y)$ in }(\N,\le).
\e
\epr

\proof
In view of Lemma \rfa2{l1.21}\,(ii), both $I(x)$ and $I(y)$ are \sbm s of
$(\N,+,0)$. Thus $0\nda22.3 = 0\cdot x= 0\cdot y\in I(x)\cap I(y)$. If
$u,v\in I(x)\cap I(y)$, then \te\ $k,l,m,n\in\N$ \st $u=k\cdot x=l\cdot y$,
$v=m\cdot x=n\cdot y$. Thus $u+v\nda22.3 = (k+m)\cdot x=(l+n)\cdot y\in
I(x)\cap I(y)$. \Mo $x\cdot y\ne0$ by \E\Pr\ \rfa2{p2.7}\,(i). \E\Tf $I(x)\cap
I(y)$ is a nontrivial \sbm\ of $(\N,+,0)$. We now show that it is \pn. By
Lemma \rf{l8.20} it is \sft\ to show that if $u,v\in I(x)\cap I(y)$ with $u\le v$,
then $v-u\in I(x)\cap I(y)$. If $k,l,m,n\in\N$ are as above, then $k\cdot x\le
m\cdot x$ and $l\cdot y\le n\cdot y$. Since $x,y\in\Na$, we infer from
\era2{2.20} that $k\le m$ and $l\le n$. By \era2{1.46} \te\ $p,q\in\N$ \st
$m=k+p$, $l=n+q$. Hence $v=m\cdot x=(k+p)\cdot x\nda22.3 = k\cdot x +
p\cdot x=u+p\cdot x$. Similarly $v=u+q\cdot y$. \If from \er{1.19} that
$v-u=p\cdot x\in I(x)$ and $v-u=q(y)\in I(y)$, hence
$v-u\in I(x)\cap I(y)$. From Lemma \rf{l8.20} we
infer that \te s $d\in\Na$ \st $I(x)\cap I(y)=I(d)$, where $d$ is the least
\el\ \wrt $\le$ of $(I(x)\cap I(y))\sms 0$. But $(I(x)\cap I(y))\sms0 =
M(x)\cap M(y)$ by \er{1.18}. Hence $d\nde1.21 = \lcm(x,y)$, which implies
\er{1.20}.
\endproof

\bco1.17
If $x,y\in\Na$, $x\ne y$, then $\lcm(x,y)$ is also the least \el\ of
$M(x)\cap M(y)$ \wrt the \og~$|$.
\eco

\proof
$M(x)\cap M(y)\nde 1.18 = (I(x)\cap I(y))\sms0 \nde1.20 = I(\lcm(x,y))\sms0
=\{n\cdot\lcm(x,y):{n\in\Na}\}$. Hence $\lcm(x,y)|n\cdot \lcm(x,y)$ \fa
$n\in\Na$.
\endproof

\bds1.18
If $x,y\in\Na$, $x\ne y$, then $M(x)\cap M(y)$ is called the set of \ti{common
\ml es} of $x,y$ (or of the set $\{x,y\}$) and $\lcm(x,y)$ defined in \er{1.21}
is called the \ti{least common \ml e} of $x,y$ (or of the set $\{x,y\}$).
\eds

We summarize in

\bpr1.19
Let $a,b\in\Na$, $a\ne b$. Then the set of common \ml es of $\{a,b\}$ contains
$ab$ and is equal to the set of \ub s of $\{a,b\}$ in $(\Na,\,|\,)$, i.e.
\beq1.22
M(a)\cap M(b) = \UB(\{a,b\}).
\e
\Mo
\beq1.23
\lcm(a,b)=\sup\{a,b\} \hbox{ in }(\Na,\,|\,).
\e
\epr

We next show that $\inf\{a,b\}$ exists in $(\Na,\,|\,)$ \fa $a,b\in\Na$,
$a\ne b$. In view of \er{8.55} the set of \lo s of $\{a,b\}$, $a,b\in\Na$,
$a\ne b$, is equal to $D(a)\cap D(b)$. Clearly $1\in D(a)\cap D(b)$. \Mo if
$c\in D(a)\cap D(b)$, then $c\le a$ by \er{8.57}, hence $D(a)\cap D(b)$ is
bounded in $(\N,\le)$. \E\Tf $D(a)\cap D(b)$ has a greatest \el\ by
Theorem \rfa1{t3.39}.

\bdn1.21
Let $x,y\in \Na$, $x\ne y$. The set $D(x)\cap D(y)$ is called the set of
\ti{common divisors} of $\{x,y\}$. Its greatest \el\ \wrt the \og~$\le$ is
called the \ti{\GCD} of $\{x,y\}$ and is denoted by $\gcd(x,y)$.\glossary{$\gcd$}
\edn

\bpr1.22
If  $x,y,a\in \Na$, $x\ne y$, then
\bea1.24
{}&\gcd(x,y)=\inf\{x,y\} \hbox{ in }(\Na,|), \\
&x\cdot y= \gcd(x,y)\cdot \lcm(x,y), \lb{1.25} \\
&\lcm(ax,ay)=a\lcm(x,y), \lb{1.26} \\
&\gcd(ax,ay)=a\gcd(x,y). \lb{1.27}
\e
\epr

\proof \

\er{1.24}, \er{1.25}: Since $xy\in M(x)\cap M(y)$, we infer from Corollary
\rf{c1.17} that $\lcm(x,y)|xy$. The idea of the proof is to show that
$d:=\frac{xy}{\lcm(x,y)}$ (see \er{8.53}) is the greatest \el\ of $D(x)\cap
D(y)$ \wrt the \og~$|$. Then $d=\gcd(x,y)$ by \er{8.57}.

(i) $d\in D(x)\cap D(y)$: Since $\lcm(x,y)\in M(x)\cap M(y)$, \te\ $b,c\in\Na$
\st $\lcm(x,y)=bx=cy$. From $xy=d\lcm(x,y)=dbx$, we obtain $y=db$ by
\era2{1.6}, \era2{1.8}. Hence $d|y$. Similarly from $xy=dcy$ we obtain $d|x$,
hence $d\in D(x)\cap D(y)$.

(ii) \ti{If $z\in D(x)\cap D(y)$, then $z|d$ and $d=\gcd(x,y)$}: Let $z\in
D(x)\cap D(y)$. Thus \te\ $b,c\in\Na$ \st $x=bz$, $y=cz$. Hence $xy=(bz)y
\nda21.5 = b(zy)\nda21.6 = b(yz)\nda21.5 = (by)z$. Similarly $xy=x(cz)=(xc)z
=(cx)z$. \If that $\frac{xy}z =by\in M(y)$, $\frac{xy}z=cx\in M(x)$, hence
$\frac{xy}z\in M(x)\cap M(y)$. By Corollary \rf{c1.17}, $\lcm(x,y)\big|
\frac{xy}z$, i.e.\ \te s $k\in\Na$ \st $\frac{xy}z=k\lcm(x,y)$. But $xy
\nde8.53 = \frac{xy}z\cdot z=(k\lcm(x,y))z \nda21.5 = k(\lcm(x,y)\cdot z)
\nda21.6 = k(z\lcm(x,y)) = (kz)\lcm(x,y)$. Thus $d=\frac{xy}{\lcm(x,y)}=kz$,
hence $z|d$, and $z\le d$ by \er{8.57}.

(iii) In view of (i) and (ii), $d=\gcd(x,y)$, every \el\ of $D(x)\cap D(y)$
divides $\gcd(x,y)$, hence \er{1.24} holds.

\er{1.25}: $xy\nde8.53 = \frac{xy}{\lcm(x,y)}\cdot \lcm(x,y)=d\cdot \lcm(x,y)
=\gcd(x,y)\cdot\lcm(x,y)$.

\er{1.26}: (iv) $ax|a\lcm(x,y)$; $ay|a\lcm(x,y)$: \E\te\ $b,c\in\Na$ \st
$\lcm(x,y)=bx=cy$. Hence $a\lcm(x,y)=abx=axb=acy=ayc$. \If that $ax|a\lcm(x,y)$
and $ay|a\lcm(x,y)$.

(v) \ti{$ax|z$ and $ay|z$ implies $\lcm(ax,ay)|z$}: We first observe that if $\a,
\b,\g\in\Na$ then
\beq1.28
\b|\g \qh{iff }\a\b|\a\g.
\e
Indeed, if $\b|\g$, then $\g=\b\d$ \fs $\d\in\Na$, hence $\a\g=\a\b\d$, and
$\a\b|\a\g$. If $\a\b|\a\g$, then $\a\g=\a\b\d$ \fs $\d\in\Na$, hence $\g=\b
\d$ by \era2{1.6} and \era2{1.8}, and $\b|\g$.

We suppose $ax|z$ and $ay|z$. Then $a|ax$ and $ax|z$ implies $a|z$. Set
$c=\frac za$, then $z=ac$. Thus $ax|ac$ and $ay|ac$. \E\Tf $x|c$ and $y|c$
by \er{1.28}. Hence $c\in M(x)\cap M(y)$, \Tf $\lcm(x,y)|c$ by Corollary
\rf{c1.17}. \csq\ $\lcm(x,y)|\frac za$, hence $a\lcm(x,y)|a\cdot\frac za$ by
\er{1.28}. Using \era2{1.6} and \er{8.53}, we obtain $a\lcm(x,y)|z$.

(vi) \If from (iv), (v) and \er{8.57} that $a\lcm(x,y)\in M(ax)\cap M(ay)$ and
that $a\lcm(x,y)=\lcm(ax,ay)$.

\er{1.27}: $(ax)(ay)\nde1.25 = \gcd(ax,ay)\cdot \lcm(ax,ay) \nde1.26 =
\gcd(ax,ay)\cdot a\cdot \lcm(x,y)$. From $axay=\gcd(ax,ay)\cdot a\cdot
\lcm(x,y)$ we obtain $a(xy)a=(\gcd(ax,ay)\cdot\lcm(x,y))\cdot a$ using
\era2{1.5} and \era2{1.6}. Hence by \era2{1.8}, we get $a(xy)=\gcd(ax,ay)
\cdot \lcm(x,y)$. Since $\lcm(x,y)|xy$, we have $\gcd(ax,ay)\cdot \lcm(x,y)
=a(xy)\nde{8.53} = a\bigl(\frac{xy}{\lcm(x,y)}\cdot \lcm(x,y)\bigr)\nde1.25 =
a(\gcd(x,y)\cdot\lcm(x,y)) \nda2{1.5} = (a\cdot\gcd(x,y))\cdot\lcm(x,y)$. Then
\er{1.27} follows from \era2{1.8}.
\endproof

We now present an efficient algorithm for computing the {\bf g}reatest
{\bf c}ommon {\bf d}ivisor (gcd) of two \nm s, \el s of~$\Na$.
This algorithm is described in Euclid's \ti{Elements} and is usually
called \ti{Euclid's algorithm}.

Let $(a_0,a_1)\in\Na\t \Na$ be \st $a_0>a_1$. By Theorem \rfa2{t1.38} with
$\ees:=(\N,0,S)$ and \era2{2.8}, \te s a unique pair $(q_0,a_2)\in\N\t\N$
\sf ying
\bga8.66
a_0=q_0a_1+a_2,\\
0\le a_2<a_1. \lb{8.67}
\e
Then, the \fw\ holds
\bea8.68
\gcd(a_0,a_1)&=a_1 &&\hbox{if }a_2=0,\\
\kern100pt \gcd(a_0,a_1)&=\gcd(a_1,a_2) &&\hbox{if }a_2>0.\kern100pt \lb{8.69}
\e
So if $a_2=0$, $\gcd(a_0,a_1)$ is simply equal to $a_1$. If $a_2>0$, then
the pair $(a_1,a_2)\in\Na\t\Na$ \sf ies the \as\ of Theorem~\rfa2{t1.38},
in view of \er{8.67}. Applying the ``division algorithm'' again we obtain
$(q_1,a_3)\in\N\t\N$ \st the \fw\ holds
\bga8.70
a_1=q_1a_2+a_3,\\
0\le a_3<a_2. \lb{8.71}
\e
Then by \er{8.69} and \er{8.68} we obtain $\gcd(a_0,a_1)=\gcd(a_1,a_2)=a_1$
if $a_3=0$ and $\gcd(a_0,a_1)=\gcd(a_1,a_2)=\gcd(a_2,a_3)$ if $a_3>0$.
Observe that in both cases $a_0>a_1>a_2>a_3$. This suggests that if we
repeat this process, then eventually $a_k$ will be equal to~$0$ for some
$k\ge2$, and $a_{k-1}=\gcd(a_0,a_1)$. This will be proved in
Theorem~\rf{t8.49}. We establish \er{8.68} and \er{8.69} in
Lemma~\rf{l8.48}. For convenience we first collect some useful identities.

\blm8.46
Let $a,b,c,d\in\N$. Then
\beq8.72
b=(b-a)+a \qhq{if}a\le b,
\e
If $a\le b$ and $c\le d$, then
\beq8.73
b-a=d-c \qhq{implies} a+d=b+c.
\e
If $a+d=b+c$ and $a\le b$, then
\beq8.74
c\le d \qhq{and} b-a=d-c.
\e
\E\Ip if $a+d=b$ then
\beq8.75
b-a=d.
\e
If $b\le c$, then $ab\le ac$ and
\beq8.76
a(c-b)=ac-ab.
\e
Suppose now that $a,b,c>0$. If $a|b$ and $a|c$, then
\beq8.77
a|(b+c) \qhq{and} \frac{b+c}a=\frac ba+\frac ca\,.
\e
If $b< c$ and if $a|b$, $a|c$, then
\bga8.78
\frac ba< \frac ca, \q a|(c-b) \qhq{and}
\frac{c-b}a = \frac ca-\frac ba\,.\\
\hbox{If }a|b, \hbox{ then }a|bc \hbox{ and }\frac {bc}a=\frac{b}{a}\cdot c, \lb{8.79} \\
ca|cb \qh{iff } a|b,\ \hbox{ and } \frac ba=\frac{cb}{ca}. \lb{1.42}
\e
\elm

\bex8.47
Prove Lemma \rf{l8.46}.
\eex

\blm8.48
Let $a_0,a_1\in\Na$, $q_0,a_2\in\N$ be \st \er{8.66} and \er{8.67}
hold. Then \er{8.68} and \er{8.69} hold.
\elm

\proof

\ti{Case $a_2=0$.} Then $a_1|a_0$ by \er{8.66}. Since $a_1|a_1$, $a_1$ is
a common divisor of $a_1$ and~$a_0$. Then \er{8.68} follows from $a_1<a_0$
and \E\df\ and Notation \rf{d1.21}.

\ti{Case $a_2>0$.} Let $b>0$. If $b|a_1$ and $b|a_2$, then $b|q_0a_1$ by
\era2{1.6}, \er{8.79}, hence $b|a_0$ by \er{8.77}. Using notation \er{8.55} we have
$D(a_1)\cap D(a_2)\sbs D(a_0)$. From $D(a_1)\cap D(a_2)=D(a_1)\cap D(a_1)
\cap D(a_2)\sbs D(a_1)\cap D(a_0)=D(a_0)\cap D(a_1)$, we get
\beq8.80
D(a_1)\cap D(a_2)\sbs D(a_0)\cap D(a_1).
\e
If $b|a_0$ and $b|a_1$, then $b|a_0$, $b|q_0a_1$ by \era2{1.6}, \er{8.79} and $q_0a_1
<a_0$ by \er{8.66}. By \er{8.78}, $b|(a_0-q_0a_1)$, hence by \er{8.75}, $b|a_2$. So we
have $D(a_0)\cap D(a_1)\sbs D(a_2)$, hence $D(a_0)\cap D(a_1)\sbs D(a_1)
\cap D(a_2)$. From \er{8.80} we obtain $D(a_1)\cap D(a_2)=D(a_0)\cap
D(a_1)$ which implies \er{8.69}.
\endproof

We now formulate Euclid's algorithm as an iteration of some map. Set
\beq8.81
A:=\{(x,y)\in\Na\t\N: x>y\}.
\e
We define a map $f:A\to\N\t\N$ by setting
\bea1.45
f(x,0)&:=(x,0) \qh{if} (x,0)\in A,\\
f(x,y)&:=(y,z) \qh{if} (x,y)\in A,\ y>0, \lb{1.46}
\e
where $(q,z)$ is the unique pair in $\N\t\N$ \st
\beq1.47
x=qy+z, \q 0\le z<y.
\e
We recall that the \ex\ and \uq\ of such a pair is a con\sq\ of Theorem
\rfa2{t1.38} with $\ees:=(\N,0,S)$ and \era2{2.8}. We claim that $f$~maps $A$
into itself. Indeed, if $y=0$, then $f(x,y)=(x,y)\in A$ since $(x,y)\in A$.
If $y>0$, then $(y,z)\in A$ by \er{1.47}. \E\Ip if $x:=a_0$, $y:=a_1$, $a_0>a_1
>0$, then $f(a_0,a_1)=(a_1,a_2)$ in \er{8.66}, \er{8.67}. If $a_2=0$, then
$a_1=\gcd(a_0,a_1)$ and $f\circ f(a_0,a_1)=f(a_1,a_2)=(a_1,a_2)$. If $a_2>0$,
then $f\circ f(a_0,a_1)=(a_2,a_3)$ in \er{8.70}, \er{8.71}.
Since $f$ maps~$A$ into itself (we still denote by~$f$ the map $f:A\to A$),
we may apply Theorem \rfa1{t1.6} and obtain a unique \sq\ of self-maps of~$A$,
denoted by $(f^n)_{n\ge0}$ \sf ying
\beq8.84
f^0:=\id_A \qhq{and} f^{n+1}:=f\circ f^n \hbox{ \fe }n\in\N.
\e
Set
\beq8.85
(f_1^n(x,y),f_2^n(x,y)):=f^n(x,y) \qh{\fe }n\in\N,\ (x,y)\in A.
\e
Note that if $y=f^0_2(x,y)=0$, $(x,y)\in A$, then $f_2^1(x,y)=0$, hence by
induction on~$n$ we have $f_2^n(x,y)=0$ \fe $n\in\N$. The case $y=f^0(x,y)
>0$ is treated in the \fw\ theorem.

\begin{thm}[Euclid's algorithm]\lb{t8.49}
Let $(f^n)_{n\ge0}$ be as in \er{8.84} and let $(x,y)\in A$ with $y>0$. Then
\te s one and only one $\ov n>0$ \st $f^{\bar n}_2(x,y)=0$ and $f_2^{\bar n-1}
(x,y)>0$. \Mo
\bea8.86
f_2^{\bar n-1}(x,y)&=\gcd(x,y),\\
f_2^{n+1}(x,y)&<f_2^n(x,y) \qh{\fe} n\in[0,\ov n-1], \lb{8.87}\\
f_2^n(x,y)&=0 \qh{\fe}n\ge \ov n. \lb{8.88}
\e
\eth

\proof
Define \sq s $a,b:\N\to \N$ by setting
\bea8.89
a_n&:=f_2^n(x,y), \qh{where} (x,y)\in A,\ y>0,\ n\in\N,\\
b_n&:=\bca
y-n & \hbox{for }0\le n\le y,\\
0   & \hbox{for }n>y.
\eca \lb{8.90}
\e
We claim that $a_n\le b_n$ \fe $n\in\N$. Set $M:=\{m\in\N: a_m\le b_m\}$.
Then

$0\in M$ since $a_0:=f_2^0(x,y)=y=b_0$.

\ti{$n\in M$ implies $n+1\in M$}: Suppose $n\in M$. If $a_n=0$, then $a_{n+1}
=0$. Indeed, if $a_n=0$, then $f^n_2(x,y)=0$ by \er{8.89}. Then $f^{n+1}(x,y)
\nde8.84 = f(f^n_1(x,y),f^n_2(x,y)) = f(f^n_1(x,y),0)\nde1.45 = (f^n_1(x,y),
0)$. Hence $a_{n+1}=f_2^{n+1}(x,y)=0$.
\E\Tf $a_{n+1}\le b_{n+1}$, by \era1{3.13}. Then $n+1\in M$ when $a_n=0$.

If $a_n>0$, then $f^n_2(x,y)>0$, hence $f^{n+1}(x,y)=f(f^n_1(x,y),
f^n_2(x,y)) \nde1.46 = (f^n_2(x,y),z)$ where $z\in\zo0,f^n_2(x,y) $ \sf ies
$f^n_1(x,y)=q\cdot f^n_2(x,y)+z$ \fs $q\in\N$. \E\Tf $a_{n+1}\nde8.89 =
f^{n+1}_2(x,y)=z < f^n_2(x,y)=a_n$. Thus $a_{n+1}<a_n$, hence $a_{n+1}+1
\le a_n$ by \era1{3.17}. Since $n\in M$, $a_n\le b_n$, hence $0<a_n\le b_n$.
Thus $0<b_n$ by \era1{3.12}. From \er{8.90} we infer $n<y$, hence $n+1\le y$
by \era1{3.17}. From \er{8.90}, \er{8.72}, we obtain $b_{n+1}+(n+1)=y
=b_n+n$. Using \er{8.9}, \er{8.8}, \er{8.16}, we get $b_{n+1}+1=b_n$.
\csq\ $a_{n+1}+1 \le a_n\le b_n=b_{n+1}+1$, hence $a_{n+1}+1\le b_{n+1}+1$,
i.e.\ $b_{n+1}+1 = (a_{n+1}+1)+c$ \fs $c\in\N$. Using \er{8.9}, \er{8.8} and
\er{8.16} again, we arrive at $a_{n+1}\le b_{n+1}$. Then $n+1\in M$ when
$a_n>0$.

\If that $M$ is \iv, hence $M=\N$, which proves the claim.

Since $b_y=0$ and $a_y\le b_y$, we have $a_y=0$. Set $B:=\{m\in\N:a_m=0\}$.
Since $(\N,\le)$ is well-ordered and $B\ne\vn$,
$B$~has a least \el\ which we denote by~$\ov n$. Since $a_0=y>0$, $\ov n>0$.
Then $f_2^{\bar n}(x,y)=a_{\bar n}=0$ and $f_2^{\bar n-1}(x,y)=a_{\bar n-1}
>0$. Since $\ov n$ is the least \el\ of~$B$, $a_n>0$ for $n<\ov n$. So if
$n<\ov n$, then $f_2^{n+1}(x,y)=a_{n+1}<a_n=f_2^n(x,y)$ by what precedes,
since $a_n>0$. This proves \er{8.87}. Next we prove
\er{8.88}, equivalently $f_2^{\bar n+p}(x,y)=0$ \fe $p\in\N$.
We have $f_2^{\bar n+p}(x,y)\nde8.89 = a_{\bar n+p}$, $p\in\N$. Set $\wt M
:=\{p\in\N: a_{\bar n+p}=0\}$. Clearly $0\in \wt M$, since $a_{\bar n}=0$.
If $p\in \wt M$, then $f_2^{\bar n+p}(x,y)=a_{\bar n+p}=0$, hence
$f^{\bar n+p+1}(x,y)= f(f_1^{\bar n+p}(x,y),f_2^{\bar n+p}(x,y))\nde1.45 =
(f_1^{\bar n+p}(x,y),0)$, and $f_2^{\bar n+p+1}(x,y)=0$. Thus $a_{\bar
n+p+1}=0$, and $p+1\in\wt M$. \If that $\wt M=\N$, which proves \er{8.88}.

Finally we prove \er{8.86}. Since $f_2^{\bar n-1}(x,y)>0$
and $f_2^{\bar n}(x,y)=0$, it follows from \er{1.46}, \er{1.47} and
Lemma~\rf{l8.48} that $\gcd\bigl(f_1^{\bar n-1}(x,y),f_2^{\bar
n-1}(x,y)\bigr) =f_2^{\bar n-1}(x,y)$. \Mo $(x,y)\nde8.84 =
(f_1^0(x,y),f_2^0(x,y))$. Hence we have to show
\beq8.91
\gcd\bigl(f_1^0(x,y),f_2^0(x,y)\bigr)=\gcd\bigl(f_1^{\bar n-1}(x,y),
f_2^{\bar n-1}(x,y)\bigr).
\e
We shall prove that
\beq8.92
\gcd\bigl(f_1^0(x,y),f_2^0(x,y)\bigr)=\gcd\bigl(f_1^{n}(x,y),
f_2^{n}(x,y)\bigr)
\e
holds \fe $n\in[0,\ov n-1]$. If $\ov n=1$, there is nothing to prove.
So we may suppose $\ov n>1$. We invoke Lemma~\rfa1{l3.46} with $E:=\N$,
$a:=0$, $b:=\ov n-1>0$, and $M:=\{n\in[0,\ov n-1]: \er{8.92}$ holds$\}$.
Clearly $0\in M$. \Mo if $n\in[0,\ov n-2] \cap M$, then $n+1\in M$. Indeed,
since $\bar n$ is the least \el\ of~$B$, $f_2^n(x,y)\nde8.89 = a_n>0$
for $n\in[0,\bar n-2]$. Thus if $n\in[0,\bar n-2]\cap M$, then
$\gcd\bigl(f_1^{n}(x,y),f_2^{n}(x,y)\bigr)=\gcd\bigl(
f_1^{n+1}(x,y),f_2^{n+1}(x,y)\bigr)$ by \er{1.46}--\er{8.85} and
Lemma~\rf{l8.48}, hence \er{8.92} also holds for $n:=n+1$.
\E\Tf $n+1\in M$. From Lemma~\rfa1{l3.46} we infer that $M=[0,\ov n-1]$, hence \er{8.91} holds. This
completes the proof of Theorem~\rf{t8.49}.
\endproof

We next establish a useful \pp y of the gcd of two \nn s.

\bpr8.50
Let $x,y\in \Na$, $x\ne y$, and let
\beq8.94
V(x,y):=\{\ell x-my\in\Na: \ell,m\in\N \hbox{ \st} my<\ell x\}.
\e
Then
\bea1.58
&x \in V(x,y),\ y\in V(y,x),\\
&V(x,y)=V(y,x), \lb{1.59} \\
&V(x,y)\cup\{0\} = I(\gcd(x,y)), \lb{1.60}
\e
where $I(\gcd(x,y))$ is the \pn\ \sbm\ of $(\N,+,0)$ generated by $\gcd(x,y)$.
\E\Ip \te s $(\ell,m)\in\Na\t\N$ \st
\beq1.61
\ell x=my+ \gcd(x,y).
\e
\epr

\proof
\er{1.58}:
$x=1\cdot x-0\cdot y\in V(x,y)$ and $y=1\cdot y-0\cdot x\in V(y,x)$.

\er{1.59}:
Let $u:=\ell x-my\in V(x,y)$. Set $a:=\max(\ell,m)+1$,
where $\max(\ell,m)$ is defined in \era2{3.14}. Then $a\ge\ell+1$ and $a\ge
m+1$, hence $a>\ell$, $a>m$ and $a\ge1$. Since $x\ge1$ we obtain $ax=xa
\nda2{2.19} \ge 1\cdot a=a$. Hence $ax>m$. Similarly $ay\ge a>\ell$, hence $ay>\ell$.
Set $c:=ax-m$ and $d:=ay-\ell$, so $ay=d+\ell$, $ax=c+m$ by \er{8.72}. From
$(ay)x=(ax)y$, we obtain $(d+\ell)x=(c+m)y$, hence $dx+\ell x=cy+my$ by
\era2{2.15}. Since $my\le \ell x$, we obtain $dx\le cy$ and
$cy-dx=\ell x-my$ by \er{8.74}. Since $\ell x-my=u>0$, we have $cy-dx>0$,
hence $u\in V(y,x)$. Interchanging $x$ and~$y$, we obtain
$V(y,x)\sbs V(x,y)$, hence \er{1.59} holds.

\er{1.60}: Since $(\N,\ge)$ is well-ordered
and $V(x,y)=V(y,x)\ne\vn$, $V(x,y)$ possesses a least \el\
which we denote by~$\a$. Thus \te\ $\ov \ell,\ov m,\ov c,\ov d\in\N$ \st
\bea8.96
\a&=\ov\ell x-\ov m y=\ov cy-\ov dx>0,\\
\a&\le u \qh{\fe} u\in V(x,y)=V(y,x). \lb{8.97}
\e
We claim that
\beq8.98
\gcd(x,y)=\a.
\e

$\gcd(x,y)\le\a$:
We have $\gcd(x,y)>0$, $\gcd(x,y)|x$ and $\gcd(x,y)|y$.
By \er{8.79}, $\gcd(x,y)|\ov\ell x$ and $\gcd(x,y)|\ov my$. Since $\ov my
<\ov\ell x$, we deduce from \er{8.78} that\break $\gcd(x,y)|(\ov\ell x-\ov my)$,
hence $\gcd(x,y)|\a$, and $\gcd(x,y)\le\a$ follows from \er{8.57}.

$\a\le\gcd(x,y)$:
We show that $\a|x$ and $\a|y$. From \E\df\ and Notation \rf{d1.21} we infer
$\a\le\gcd(x,y)$.

$\a|x$: Since $\a\in\Na$, we obtain from Theorem \rfa2{t1.38} and \era2{2.8}
that $x=q_1\a+r_1$ for some $q_1\in\N$ and $r_1\in[0,\a)$. We claim that
$r_1=0$, hence $\a|x$. Suppose for \cd ion that $r_1>0$. From \er{8.96},
$\a=\ov cy-\ov dx>0$, hence $q_1\ov dx\le q_1\ov cy$ and $q_1\a=q_1\ov cy
-q_1\ov dx$ by \er{8.76}. \E\Tf $x=r_1+(q_1\ov cy-q_1\ov dx)$, hence by
\er{8.72} $(1+q_1\ov d)x=x+q_1\ov dx=r_1+q_1\ov cy$. Thus $r_1=(1+q_1\ov d)x
-q_1\ov cy$. Since $r_1>0$, we have $r_1\in V(x,y)$, hence $\a\le r_1$ by
\er{8.97}, \cd ing $r_1\in[0,\a)$.

$\a|y$: Similarly we have $y=q_2\a+r_2$ with $q_2\in\N$ and $r_2\in[0,\a)$.
From $\a=\ov\ell x-\ov my$ we deduce that $r_2=(1+q_2\ov m)y -q_2\ov\ell x
>0$, if $r_2>0$. Hence $r_2\in V(y,x)=V(x,y)$, and $r_2\ge\a$. A~\cd ion. \E\Tf
$r_2=0$ and $\a|y$.
\If that \er{8.98} holds. From \er{8.96} and \er{8.98} we infer \er{1.61}.

$I(\gcd(x,y))\sbs V(x,y)\cup\{0\}$: $0\dpl \gcd(x,y)\nda22.3 = 0\in\{0\}$.
Let $n\in\Na$. Then $n\dpl \gcd(x,y)\nda22.8 = n\cdot(\bar\ell x-\bar my)
\nde8.76 = n(\bar\ell x)-n(\bar my)\nda22.11 = (n\ell)x-(n\bar m)y\in V(x,y)$.

$\{0\}\sbs I(\gcd(x,y))$: Follows from $0\dpl \gcd(x,y)$.

$V(x,y)\sbs I(\gcd(x,y))$: Let $\ell,m\in\N$ be \st $\ell x-my>0$. By \df,
$\gcd(x,y)|x$ and $\gcd(x,y)|y$. By \er{8.79}, $\gcd(x,y)|\ell x$ and $\gcd(x,y)|my$.
From $\ell x-my>0$ and \er{8.78} we infer $\gcd(x,y)|(\ell x-my)$. Hence \te s
$d\in\Na$ \st $\ell x-my=d\gcd(x,y)$. Thus $\ell x-my\in I(\gcd(x,y))$.

\brm1.27
Interchanging $x$ and $y$ in \E\Pr\ \rf{p8.50} we find that \te\ $(\wt m,
\wt\ell)\in\Na\t\N$ \st
\beq1.62
\wt\ell y=\wt m x+\gcd(x,y),
\e
since $\gcd(y,x)=\gcd\{y,x\}=\gcd\{x,y\}=\gcd(x,y)$.
\erm

Let $a,b,c\in\Na\sms1$. If $ab|c$ then $a|c$ and $b|c$. Indeed, since $a|ab$
(resp.\ $b|ab$) we have $a|c$ (resp.\ $b|c$) by \tr ity of~$|$. The converse
is not true since $aa\nmid a$. However, the \fw\ holds.

\blm1.28
Let $a,b,c\in\Na\sms1$ be \st $a|c$ and $b|c$. If $\gcd(a,b)=1$, then $ab|c$.
\elm

\proof
We have $c\in M(a)$ and $c\in M(b)$ by \er{8.54}, hence $c\in M(a)\cap M(b)$.
By Corollary \rf{c1.17}, $\lcm(a,b)$ is the least \el\ of $M(a)\cap M(b)$ in
$(\Na,|)$ and $\lcm(a,b)|c$ by \er{1.23}. By \er{1.25}, $ab=\gcd(a,b) \cdot
\lcm(a,b)=1 \cdot\lcm(a,b)=\lcm(a,b)$. Hence $ab|c$.
\endproof

A \gn\ of Lemma \rf{l1.28} is given in Lemma \rfa4{l2.30}.

We conclude this section by proving a lemma (sometimes called Euclid's lemma),
which will play an important role in the next section.

\blm1.29
Let $a,b,q\in\Na$. If $\gcd(a,q)=1$, and $q|ab$, then $q|b$.
\elm

\proof
Since $\gcd(a,q)=1$, it follows from \E\Pr\ \rf{p8.50}, \er{1.61} that \te s
$(\ell,m)\in\Na\t \N$ \st $\ell a=mq+1$. \E\Tf $\ell ab=(\ell a)b=(mq+1)b
=mqb+b$. Since $q|ab$, \te s $d\in\Na$ \st $ab=dq$. Hence $(\ell d)q=(mb)q+b$.
Since $b\in\Na$, $(mb)q< (\ell d)q$, hence by \era2{2.20} and $q\in\Na$,
$mb< \ell d$. Using \er{1.19},
\er{8.76} we deduce that $b=(\ell d)q-(mb)q=(\ell d-mb)q$. Thus $q|b$.
\endproof

\bdf1.30
Two \el s $x$ and $y$ of $\Na\sms1$ are called \ti{coprime}\index{coprime} or \ti{relatively
prime}\index{relatively prime} if $\gcd(x,y)=1$, that is, if $1$~is the only common divisor of
$x$ and~$y$.
\edf

\newpage

\Subsubsection{The fundamental theorem of \art}\label{sss.Fundament}
We recall that $0$ (resp.~$1$) is the least \el\ of the P-monoid $(\N,+,0)$
(resp.\ $(\Na,\cdot,1)$) \wrt the \og\ defined in \er{8.7}. \Mo $1$~is the
least \el\ of $\N\sms0$ in $(N,\le)$ and $1$~is the \Gn\ of the \pn\ monoid
$(\N,+,0)$, i.e.\ $\N=I(1)$. In contrast to $(\N,+,0)$, $\Na\sms1$ has no
least \el\ in~$(\Na,|)$. However, it contains \ti{minimal\/} \el s.

\bdf9.1
Let $(X,\le)$ be an \os. An \el\ $a\in X$ \st $a>x$ (resp.\ $a<x$) for no
$x\in X$ is called a \ti{minimal\/}\index{element!minimal} (resp.\ \ti{maximal\/})\index{element!maximal} \el\ of
$(X,\le)$.
\edf

\brm9.2
The least \el\ (if it exists) of an \os\ $(X,\le)$ is the unique minimal \el\
of~$X$. Even if an \os\ possesses a unique minimal \el~$a$, this \el~$a$ is
not necessarily the least \el\ of~$X$. Similarly for maximal \el s. One shows
using \In\ on the \nme\ that a \ti{finite\/} \os\ possesses at least
a minimal and a maximal \el.
\erm

\bdf9.6
A \ti{\Pn} (or \ti{prime})\index{number!prime} is a minimal \el\ of $\Na\sms1$ in $(\Na,|)$. In
other words, a \Pn\ is an \el\ of $\Na\sms1$ that has no divisors other
than~$1$ and itself.
\edf

Clearly the \nm\ $2$  is a \Pn\ since  any divisor of~$2$ is less than or
equal to~$2$ by \er{8.57}.

\bnt9.7
We denote by $P$ the set of all \Pn s.
\ent

A \ml e of $2$ in $(\N,+,0)$ is called an \ti{even} \nm\ and a \nm\ in~$\N$
which is not even is called an \ti{odd\/} \nm.
So $2$ is the only even \Pn\ and $3$ is the least (\wrt $\le$) odd \Pn.

A \nm\ $n\in \Na\sms1$ which is not a \Pn\ is by \E\df\ \rf{d9.6} the product
$a\cdot b$ of two (not necessarily distinct) factors $a,b\in\Na\sms{1,n}$.
This motivates the \fw\ \df.

\bdf9.8
An \el\ of $\Na\sms1$ which is not a \Pn\ is called a \ti{\cme}\nm.\index{number!composite}
\edf

We now show that every \cme \nm\ $n$ possesses a factorization $a\cdot b$
where $a$ is a prime.

\blm9.9
\E\fe $n\in\Na\sms1$ \te s $p\in P$ \st $p|n$.
\elm

\proof
Set $A:=\{n\in\Na\sms1: \hbox{\te s no $p\in P$ \st}p|n\}$. Suppose for \cd
ion that $A\ne\vn$. Since $(\N,\le)$ is well-ordered \te s a least \el\ of
the subset $A$ of~$\N$, which we denote by~$\bar n$. Since $\bar n\in A$,
$\bar n$~is a \cme \nm\ since every \Pn\ divides itself. \E\Tf $\bar n=a\cdot
b$ for some $a,b\in\Na\sms{1,\bar n}$.

Since $a|\bar n$ and $b|\bar n$, we have $a\le\bar n$ and $b\le\bar n$ by
\er{8.57}. Since $a\ne\bar n$ and $b\ne\bar n$, we conclude that $a<\bar n$
and $b<\bar n$. Hence neither $a$ nor~$b$ belongs to~$A$ since $\bar n$ is
the least \el\ of~$A$. \E\Tf there is $p\in P$ \st $p|a$. Since $p|a$ and
$a|\bar n$, we infer by the \tr ity of~$|$ that $p|\bar n$. A~\cd ion, since
$\bar n\in A$. Hence $A=\vn$.
\endproof

From Lemma \rf{l9.10} we deduce that a \cme product of a finite \ns\ of factors is a
\ml e of each of its factors.

We are now in a position to prove that the set $P$ of \Pn s is
infinite, a result already contained in Euclid's \ti{\E\el s}.

\bpr9.11
The set of \Pn s $P$ is \ct y infinite.
\epr

\proof
\ti{Step 1}. \ti{Every nonempty finite subset of $P$ has a nonempty
complement in~$P$} (Euclid).

Let $A$ be a finite nonempty subset of $P$. Set
\beq9.16
m:=\Bigl(\prod_{p\in A} p\Bigr)+1.
\e
Then clearly $m\in\Na\sms1$. In view of Lemma~\rf{l9.10}, we have \fe $p\in A$:
\[
m=\Bigl(\prod_{q\in A\sms p} q\Bigr)\cdot p+1.
\]
From the \uq\ part of Theorem \rfa2{t1.38} we infer that \fe $p\in A$ there is
no $n\in\Na$ \st $m=n\cdot p+0=n\cdot p$. \E\Tf none of the primes belonging
to~$A$ divides~$m$. By Lemma~\rf{l9.9} there is $\bar p\in P$ \st $\bar p|m$.
Hence $\bar p\in P\sm A$.

\smallskip
\ti{Step 2}. \ti{$P$ is \ct y infinite}.

As a \nss\ of~$\N$, $P$ is either bounded or unbounded. If $P$ is bounded, \te\
$n\in\N$ and a bi\jn\ $b:[0,n]\to P$ by Lemma~\rfa1{l4.10}. This is impossible since by
Step~1, $P$~should contain an \el\ not belonging to the range of~$b$. \E\Tf
$P$~is unbounded in~$\N$, hence $P$~is \ep\ to~$\N$ by Theorem~\rfa1{t3.49}. Thus
$P$~is \ct y infinite in view of \E\df~\rfa1{d4.21}.
\endproof

Our next goal is to show that every \nn\ greater than one is either a prime
or a product of primes. To this end we shall use a form of mathematical
induction different from the one used previously.

\begin{lem}[second form of mathematical induction, see \cite{Alg}]\lb{l9.12}
Let $(W,\le)$ be a well-ordered set and let $e$ be its least \el. Let $M$ be
a subset of~$W$. Then $M=W$ provided it \sf ies the \fw\ two \cn s{\rm:}

\hph i,i, $e\in M$,

\hph ii,, for all $x\in W$, $x\in M$ whenever $\{y\in W:e\le y<x\}\subset M$.
\elm

\proof
Suppose for \cd ion that $M\ne W$. Then the \ns\ $\{x\in W:x\notin M\}$ has a
least \el\ $\bar x$ since $(W,\le)$ is well-ordered. By \df\ of~$\bar x$,
${\{x\in W: e\le y<\bar x\}}$ is a subset of~$M$ and is nonempty since
$e\in M$. From \cn~(ii) we infer that $\bar x\in M$. A~\cd ion.
\endproof

We are now in a position to prove that every \nn\ greater than one is either
prime or a product of primes.

\blm9.13
Let $m\in\Na\sms1$. Then \te\ $N\in\N$ and a map $f:[0,N] \to P$ not
necessarily in\jc\ \st
\beq9.17
m=\prod_{j=0}^N f_j
\e
where $[0,N]$ denotes an \il\ of $(\N,\le)$ $($see \era1{3.27}$)$.
\elm

\proof
Let $M:=\{m\in\Na\sms1 \hbox{ \st \er{9.17} holds}\}$. We claim (i)~$2\in M$,
(ii)~if $[2,n)\sbs M$ for some $n\in\Na$, $n>2$, then $n\in M$.

$2\in M$: $2$ is a \Pn. \E\Tf \er{9.17} holds with $N=0$ and $f_0:=2$.

\ti{$[2,n)\sbs M$ implies $n\in M$}: Suppose $[2,n)\sbs M$ for some
$n\in\Na$, $n>2$. Then either $n\in P$ or $n\notin P$. If $n\in P$, then
\er{9.17} holds with $N=0$ and $f_0:=n$. \E\Tf $n\in M$. If $n\notin P$, \te\
$a,b\in\Na\sms1$ \st $n=ab$. Note that $a\ne n$ since otherwise $nb=ab=n
\nda21.13 = n1$, hence $b=1$ by \era2{2.16}\,(iii), \cd ing $b\in\Na\sms1$.
Similarly one proves $b\ne n$. From $a|n$ and $b|n$, we
infer $a\le n$ and $b\le n$ by \er{8.57}. Since $a\ne n$ and $b\ne n$, we have
$a<n$ and $b<n$. Hence
$a,b\in[2,n)\sbs M$. Consequently, \te\ $s,t\in\N$ and maps
$f:[0,s]\to P$, $g:[0,t]\to P$ \st $a=\prodl_{i=0}^s f_i$, $b=\prodl_{j=0}^t
g_j$. Define a map $h:[0,s+1+t]\to P$ by setting
\[
\bal
h_i:=f_i &\hbox{ for }i\in[0,s],\\
h_{s+1+k}:=g_k &\hbox{ for }k\in [0,t].
\eal
\]
Then
\[
n=a\cdot b=\prod_{i=0}^s f_i\cdot \prod_{j=0}^t g_j =
\prod_{i=0}^s h_i \cdot\prod_{j=0}^{t} h_{s+1+j} \nad *=
\prod_{i=0}^s h_i \cdot\prod_{i=s+1}^{s+1+t} h_i\nad{**}=
\prod_{i=0}^{s+1+t}h_i.
\]
In $\nad*=$ we used \era2{7.12} and in $\nad{**}=$ \era2{7.13} which is a con\sq\
of the \asc ity of the \mlc. \E\Tf \er{9.17} holds with $N:=s+1+t$ and
$f:=h$. In view of Lemma \rf{l9.12} where $W=\Na\sms1$, $e:=2$ and $\le$ is
the \rt ion of the \og\ of~$\N$ to~$\Na\sms1$, we have $M=\Na\sms1$.
\endproof

\bex9.14
Find for each $n\in\{1,2,3,4,5\}$ a finite subset of $P$ denoted by~$A$
\st $\#(A)=n$ and $\bigl(\prodl_{p\in A }p\bigr)+1$ is not a \Pn.
\eex

We now consider the problem of ``\uq'' of the
factorization \er{9.17}. The \fw\ \pp y  of \Pn s is \fd\ in the proof of
\uq. This \pp y was already known to Euclid (and proved by him in \ti{\E\el
s}, Book~VII). We shall call it \ti{\pp y}~(P).
$$
\ga
\hbox{If a \Pn\ $p$ divides a product of two \nn s greater than zero,}\\
\hbox{then $p$ divides one of its factors.}
\ega \leqno\rm(P)
$$


Note that a \cme \nm\ $n$ cannot \sf y \pp y (P). Indeed, by \df\ \te\
$a,b\in\Na\sms1$ \st $n=ab$. Hence $n|ab$ but neither $n|a$ nor~$n|b$ holds,
since $a,b<n$.

\E\pp y (P) for \Pn s is a con\sq\ of Lemma \rf{l1.29}.

\begin{lem}[also sometimes called Euclid's lemma]\lb{l9.15} \
\beq9.18
\hbox{For all }a,b\in\Na \hbox{ and all $p\in P$, }
p|ab \hbox{ implies $p|a$ or }p|b.
\e
\elm

\proof[Proof of {\rm(P)}]
Let $a,b\in\Na$ and $p\in P$ be \st $p|ab$. By \df\ of a \Pn\ we have either
$\gcd(a,p)=1$ or $\gcd(a,p)=p$. If $\gcd(a,p)=1$, then $p|b$ by
Lemma~\rf{l1.29}. If $\gcd(a,p)=p$, then $p|a$.
\endproof

We now prove a \gn\ of \pp y (P).

\blm9.17
Let $N\in\N$, $f:[0,N]\to\Na$ and let $p\in P$. If $p$ divides
$\prodl_{k=0}^N f_k$, then \te s $\ell\in[0,N]$ \st $p$ divides $f_\ell$.
\elm

\proof
By induction on $N$. Set
\[
\cM:=\{N\in\N: \hbox{the assertion of Lemma \rf{l9.17} holds for }N\}.
\]

$0\in\cM$ since $f_0=\prodl_{k=0}^0 f_k$.

\ti{$M\in\cM$ implies $M+1\in\cM$}: Suppose $M\in\cM$ and let $p\in P$ \st
$p|\prodl_{j=0}^{M+1} f_j$ where $f_j\in\Na$, $j\in[0,M+1]$. Then
$\prodl_{j=0}^{M+1}f_j= \bigl(\prodl_{j=0}^M f_j\bigr)\cdot f_{M+1}$. If
$p|f_{M+1}$ then $M+1\in\cM$. If not, we have $p|\prodl_{j=0}^M f_j$ by \pp
y~(P) since $p\in P$. From the \iv\ hypothesis $M\in\cM$ we infer that \te s
$\ell\in[0,M]$ \st $p|f_\ell$. Hence $M+1\in\cM$. \E\Tf $\cM=\N$.
\endproof

We are now in a position to prove the ``\uq'' of the factorization \er{9.17}.

\blm9.18
Let $f:[0,M]\to P$ and $g:[0,N]\to P$ with $M,N\in\N$. Set
\[
a:=\prod_{i=0}^M f_i, \q  b:=\prod_{j=0}^N g_j
\]
and
\beag
A&:=\{p\in P: p=f_i \hbox{ for some }i\in[0,M]\}\\
B&:=\{p\in P: p=g_j \hbox{ for some }j\in[0,N]\}
\e
Then $a=b$ implies $A=B$.
\elm

\proof\

$A\sbs B$: Let $p\in A$, then $p\in P$ and $p=f_i$ for some $i\in[0,M]$. It
follows from Lemma~\rf{l9.10} that $p|a$. If $a=b$, then $p|b$. If $N=0$,
then $b=g_0$ by \era2{7.10}, hence $p\in B$. If $N>0$, then $p$ divides some
factor of~$b$ by Lemma~\rf{l9.17}. Thus $p|g_j$ for some $j\in[0,N]$. Since
$g_j\in P$ and $p\ne1$, we have $p=g_j$, hence $p\in B$.

$B\sbs A$ follows from ``$A\sbs B$'' by interchanging $f$ and~$g$.
\E\Tf $A=B$.
\endproof

We have shown that if two \cme products of primes are equal, then the primes
occurring in these products are the same. We are now interested in the
\ti{\nm} of occurrences of these primes in these products. This \nm\ will be
called the \ml icity of the prime.\index{multiplicity of the prime}

\bdf9.19
Let $f:[0,M]\to P$ for some $M\in\N$ and let $p\in R(f)$ ($:=$ the range
of~$f$ in~$P$). Then $\#\bigl(\{i\in[0,M]: f_i=p\}\bigr)$ is called the
\ti{\ml icity} of~$p$ in the \cme product $\prodl_{j=0}^M f_j$.
\edf

\begin{thm}[\E\uq\ part of the \fd\ theorem of \art]\lb{t9.20}
Let $M,N,f$, $g,a,b,A,B$ be as in Lemma \rf{l9.18}. Suppose $a=b$, then $A=B$,
and \fa $p\in A$ the \ml icity of $p$ in $\prodl_{i=0}^M f_i$ is equal to the \ml
icity of~$p$ in $\prodl_{j=0}^N g_j$. \Mo $M=N$ and the sum of the \ml
icities of the primes $p$ belonging to~$A$ is equal to $M+1$.
\eth

\proof
$A=B$ by Lemma \rf{l9.18}.
Let $p\in A=B$. Set $M_p:=\{i\in[0,M]: f_i=p\}$ and $N_p:=\{j\in[0,N]:
g_j=p\}$. If $A=B=\{p\}$, $\prodl_{i=0}^M f_i=p^{M+1}$ and
$\prodl_{j=0}^N g_j=p^{N+1}$ by \era2{1.126} in $(\Na,\cdot,1)$.
Since $a=b$, we have $p^{M+1}=p^{N+1}$, hence
$M+1=N+1$ by \era2{2.34} since $p\ge2$. Then $M=N$ by \era2{1.8}.
\Mo $M+1=\#(M_p)$ and
$N+1=\#(N_p)$ hence $M_p=N_p$. Suppose now that $A=B$ contains more than one
prime, and let $p\in A=B$. Let $m_p$ (resp.~$n_p$) denote the \ml icity
of~$p$ in $\prodl_{i=0}^M f_i$ (resp.\ $\prodl_{j=0}^N g_j$). \E\te s (see
Exercise \rf{ex9.21}) a~bi\jn~$\vf:[0,M] \to [0,M]$ (resp.\
$\psi:[0,N]\to[0,N]$) \st $f_{\vf(i)}$ (resp.~$g_{\psi(i)}$) equals~$p$ \fe $i\in
[0,m_p)$ (resp.~$[0,n_p)$). By \era2{7.13} $\prodl_{i=0}^M f_i=p^{m_p}\cdot c$ and
$\prodl_{j=0}^N f_i=p^{n_p}\cdot d$ where $c,d\in\Na$ are \cme products of
primes not equal to~$p$. If $a=b$, we claim that $m_p=n_p$. Indeed, suppose
for \cd ion that $m_p\ne n_p$. Suppose $m_p>n_p$, i.e.\ $m_p=n_p+r$,
$r\in\Na$. We have $p^{n_p+r}\cdot c=p^{n_p}\cdot d$, hence
$p^{n_p}\cdot(p^r\cdot c)=p^{n_p}\cdot d$, by~\era2{2.11}, \era2{2.25}.
From \era2{2.16}\,(iii) we
infer $p^r\cdot c=d$. Since $r\ge1$, $p|p^r$, hence $p|(p^r\cdot c)$
by \er{8.79}, thus $p|d$.
By Lemma~\rf{l9.17} $p$~divides some $g_j$ in the \cme
product~$d$. \E\Tf $p=g_j$, a \cd ion. The case $n^p>m^p$ follows by
interchanging $f$~and~$g$.

\E\Tf $m_p=n_p$. Finally, we prove that $M+1=\suml_{p\in A}m_p
=\suml_{p\in B}n_p=N+1$. The equality $\suml_{p\in A}m_p=\suml_{p\in B}n_p$
is clear since $A=B$ and $m_p=n_p$. We now prove $M+1=\suml_{p\in A}m_p$.
Observe that $[0,M]\nda23.5 =\zo0,M+1 $ is the disjoint union of $\{M_p\}_{p\in A}$.
Hence $M+1=\#(\zo0,M+1 )=\#([0,M])\nda22.13 = \#([0,M])\cdot1 \nda21.126 = \suml
_{i\in[0,M]}1\nda21.127 = \suml_{p\in A}\bigl(\suml_{i\in M_p}1\bigr)
\nda21.126 = \suml_{p\in A}(\#(M_p)\dpl 1)\nad{\era2{2.8},\era2{2.13}} = \suml_{p\in
A}\#(M_p)=\suml_{p\in A}m_p$ by \E\df\ \rf{d9.19}.
 The proof of $N+1=\suml_{p\in B}n_p$ is similar.
\endproof

\bex9.21
Let $M\in\Na$, $C,D$ be \nss s of~$[0,M]$ \st $C\cup D=[0,M]$ and $C\cap D=\vn$.
Construct a bi\jn\ $\vf:[0,M]\to [0,M]$ \st $\vf(C)=[0,\#(C)-1]$ and
$\vf(D)=[\#(C),M]$. (Use \In\ on $M\in\Na$.)
\eex

\blm9.22
Let $M,f,a$ and $A$ be as in Lemma \rf{l9.18} and let $m_p$ be the \ml icity
of $p\in A$. Then $m_p$ is the only \el\ of the set
\[
\{m\in\Na: p^m|a \hbox{ and }p^{m+1}\nmid a\}.
\]
\elm

\proof
As in the proof of Theorem \rf{t9.20} we find that $a=p^{m_p}\cdot c$ where
$c$~is either equal to~$1$ or is a \cme product of primes different from~$p$.
In the first case clearly $p^{m_p}|p^{m_p}$ and $p^{m_p+1}\nmid p^{m_p}$ since
$p^{m_p}<p^{m_p+1}$. In the second case clearly $p^{m_p}|a$ but $p^{m_p+1}
\nmid a$ since otherwise as in the proof of Theorem \rf{t9.20} $p$~would
divide~$c$, hence by Lemma \rf{l9.17} $p$~would divide a prime different
from~$p$, which is impossible.
\endproof

The \fw\ theorem is a reformulation of Lemma \rf{l9.13}.

\begin{thm}[\E\ex\ part of the \fd\ theorem of \art]\lb{t9.22}
Let $n$ be a \nn\ greater than one. Then
\beq9.19
n=\prod_{\sbk{p\in P\\p|n}}p^{m_p}
\e
where $m_p$ is the only positive \nn\ \st $p^{m_p}$ divides~$n$ and
$p^{m_p+1}$ does not divide~$n$.
\eth


\proof
By Lemma \rf{l9.13} $n=\prodl_{j=0}^N f_j$ where $f:[0,N]\to P$. Then
$\prodl_{j=0}^N f_j= \prodl_{j\in[0,N]}f_j$ by \era2{7.17} where $A:=[0,N]$
and $\vf:=\id_A$. We now apply \era2{1.128} where $(X,\qu,e):=(\Na,\cdot,1)$,
$\O:=[0,N]$, $a:\O\to X$ is defined by $a_\o:=f_\o$, $a(\O)=f([0,N])$, and
$A_p:=\{p\in P: f_\o=p\}$. In view of \E\df\ \rf{d9.19} $\#(A_p)=m_p$, the
\ml icity of~$p$ in the \cme product $\prodl_{j=0}^N f_j$. Note that $f([0,N])
=\{p\in P: p|n\}$. Indeed, if $p\in f([0,N])$, then $p|n$ by Lemma \rf{l9.10}.
Conversely, if $p|n$, then by Lemma \rf{l9.17} \te s $j\in[0,N]$ \st $p|f_j$
hence $p=f_j$ since $f_j$ is a prime. \E\Tf by \era2{1.128} we obtain
$\suml_{j\in[0,N]} f_j = \prodl_{\sbk{p\in P\\p|n}} m_p\ddt p=
\prodl_{\sbk{p\in P\\p|n}}p^{m_p}$.
\endproof

\newpage
\Subsubsection{The lattice $(\Na,\cdot,1)$}\label{sss.latt}

\begin{dfn}[\cite{Lattice}]\lb{d3.1}
A \ti{\lt} is an \os\ \st every subset consisting of two \el s possesses\index{lattice}
a~supremum and an~infimum.
\edf

\bxs3.2 \

\hph i,ii, $(\N,\le)$ is a lattice: $\sup\{x,y\}=\max(x,y)$, $\inf\{x,y\}=
\min(x,y)$ by \E\df\ \rfa2{d3.17} and Definitions and Notations \rf{d1.12},
where $x,y\in\N$, $x\ne y$.

\hph ii,i, $(\Na,|)$ is a \lt: $\sup\{x,y\}=\lcm(x,y)$, $\inf\{x,y\}=
\gcd(x,y)$, by \er{1.21}, \er{1.24}, where $x,y\in\Na$, $x\ne y$.

\hph iii,, The \rt ion of the \og~$|$ of $(\Na,|)$ to the subset $\{2,3,5\}$
is the \et y \rl. Hence $(\{2,3,5\},|)$ is \ti{not\/} a~\lt.
\exs

Observe that if $(X,\le)$ is an \os\ and $\ge$ denotes the reverse \og\ (see
Exercise \rfa1{ex3.4}), if $A$~denotes a \nss\ of~$X$, and if $\sup A$
(resp.\ $\inf A$) exists, then $\sup A$ (resp.\ $\inf A$) is equal to
$\inf A$ (resp.\ $\sup A$) \wrt
the reverse \og~$\ge$. \If that if $(X,\le)$ is a~\lt, then so is $(X,\ge)$.

\smallskip
Our first goal is to show that if $(X,\le)$ is an \os\ \st $\sup\{x,y\}$
exists \fa $x,y\in X$, $x\ne y$, then $\sup A$ exists \fe \ti{nonempty finite}
subset of~$X$. To this end we introduce the \fw\ binary \op\ on~$X$.
We set
\beq3.1
x\lor y:=\sup{\{x\}\cup\{y\}}, \quad x,y\in X.
\e
If $x=y$, then $\{x\}\cup\{y\}=\{x\}$, $\UB(\{x\})$ is not empty since $x\in
\UB(\{x\})$ and $x$~is the least \el\ of $\UB(\{x\})$. Hence $\sup\{x\}=x$,
thus
\beq3.2
x\lor x=x \qh{\fa} x\in X.
\e
We now suppose $x\ne y$. Since $\{x\}\cup\{y\}=\{y\}\cup\{x\}$, we have
\beq3.3
x\lor y=y\lor x \qh{\fa}x,y\in X.
\e
\Mo by \er{3.1}, we have
\beq3.4
x\le x\lor y,\ y\le x\lor y \qh{\fa}x,y\in X,
\e
and
\beq3.5
\hbox{if }x\le a,\ y\le a \hbox{\ \fs}a\in X,\qh{then} x\lor y\le a.
\e
If $x\le y$, then $y=\sup\{x\}\cup\{y\}$ by \er{3.2} and by Lemma \rfa1{l3.38},
hence $y=x\lor y$. Conversely, if $y=x\lor y$, then $x\nde3.4 \le x\lor
y=y$, hence $x\le y$. Thus we have
\beq3.6
x\le y \qh{iff } y=x\lor y \hbox{ \fa}x,y\in X.
\e

We next show that the \op~$\lor$ is \ti{\asc e}. We first prove
\beq3.7
(x\lor y)\lor z = \sup{\{x\}\cup\{y\}\cup\{z\}} \qh{\fa}x,y,z\in X.
\e
Set $\a:=(x\lor y)\lor z$, $x,y,z\in X$.

$\a\in UB(\{x\}\cup\{y\}\cup\{z\})$: $x\lor y\nde3.4 \le \a$, $z\nde3.4 \le
\a$, $x\nde3.4 \le x\lor y$, $y\nde3.4 \le x\lor y$, hence by \tr ity
of~$\le$, $x\le\a$, $y\le \a$.

\ti{$\a$ is the least \ub\ of $\{x\}\cup\{y\}\cup\{z\}$}: Let $a\in X$ \sf y
$x\le a$, $y\le a$ and $z\le a$. Then $x\lor y\nde3.5 \le a$, hence
$\a=(x\lor y)\lor z\nde3.5 \le a$. Thus \er{3.7} holds.

Since $\{x\}\cup\{y\}\cup\{z\}=\{y\}\cup\{z\}\cup\{x\}$ \fa $x,y,z\in X$,
$(x\lor y)\lor z\nde3.7 = \sup\{x\}\cup\{y\}\cup\{z\}=\sup\{y\}\cup\{z\}\cup
\{x\}\nde3.7 = (y\lor z)\lor x\nde3.3 = x\lor (y\lor z)$. Thus
\beq3.8
(x\lor y)\lor z=x\lor (y\lor z) \qh{\fa} x,y,z\in X.
\e
Now let $A$ be a \ti{nonempty finite subset\/} of~$X$, and let $\vf:[0,
\#(A)-1]\to A$ be a bi\jn\ (see the discussion preceding \era2{7.16}). Set
$\LOR_{j=0}^{\#(A)-1}\vf(j)$, where $\LOR_{j=0}^{\#(A)-1}$ is the \cme \op\
defined in \era2{7.10} with $f_i(x):=x\lor \vf(i)$, $x\in X$,
$i\in[0,\#(A)-1]$ in \era2{7.9}. In view of \er{3.8}, \er{3.3} and
\era2{7.14}, $\LOR_{j=0}^{\#(A)-1}\vf(j)=\LOR_{j=0}^{\#(A)-1}\psi(j)$, where
$\psi$~is any bi\jn\ from $[0,\#(A)-1]$ onto~$A$. Note that $(X,\lor)$ is
an abelian \sg, but it has no \nel\ in general.

\bpr3.3
Let $(X,\le)$ be an \os\ \st every subset consisting of two \el s has a
supremum. Let $A$~be a \nfs\ of~$X$. Then $\sup A$ exists and
\beq3.9
\sup A = \LOR_{i=0}^{\#(A)-1}\vf(i),
\e
where $\vf$ is any bi\jn\ from $[0,\#(A)-1]$ onto~$A$ and
$\LOR_{i=0}^{\#(A)-1}$ is the \cme \op\ on~$X$ induced by the binary
\op~$\lor$ defined in \er{3.1}.

\E\Ip if $\#(A)>1$,
\beq3.10
\sup A=(\sup A\sms a)\lor a \qh{\fa}a\in A.
\e
\Mo if $B$ is a \nf,
\beq3.11
\sup A\cup B = (\sup A)\lor (\sup B).
\e
\E\Ip if $A\sbs B$, then
\beq3.12
\sup A\le \sup B.
\e
\epr

\proof \

\ti{$\sup A$ exists and \er{3.9} holds}: We proceed by \In\
on~$n:=\#(A)\in\Na$. Set $M:=\{n\in\Na: \hbox{\fe \nfs\ $A$ of $X$ with
}\#(A)=n,\ \sup A \hbox{ exists and}\break \hbox{\er{3.9} holds}\}$.

$1\in M$: Let $A$ be a singleton of $X$ and let $\vf:\{0\}\to A$ be \st
$\{\vf(0)\}=A$. Then $\sup A=\sup\{\vf(0)\}=\vf(0)\nda27.10 = \LOR_{i=0}^0
\vf(i)=\LOR_{i=0}^{\#(A)-1}\vf(i)$. Thus $1\in M$.

\ti{$n\in M$ implies $n+1\in M$}: Let $n\in M$ and let $A$ be a finite subset
of~$X$ with $\#(A)=n+1$. Let $\vf:[0,n]\to A$ be a bi\jn. Then $A$~is the
disjoint union of $A\sms{\vf(n)}$ and $\{\vf(n)\}$. From \era2{3.31} we infer
$\#(A\sms{\vf(n)})=n$. Since $n\in M$, $\sup A\sms{\vf(n)}$ exists and
\beq3.13
\sup A\sms{\vf(n)}= \LOR_{i=0}^{n-1}\vf(i).
\e
\E\Tf $\LOR_{i=0}^{\#(A)-1}\vf(i) = \LOR_{i=0}^n \vf(i)\nad{\era2{7.9},
\era2{7.10}}=\Bigl(\LOR_{i=0}^{n-1}\vf(i)\Bigr)\lor \vf(n)\nde3.13 =
(\sup A\sms{\vf(n)})\lor\vf(n)$. Set $\a:=(\sup A\sms{\vf(n)})\lor\vf(n)$. We
claim:
\beq3.14
\sup A \hbox{ exists and }\sup A=\a.
\e

$\a\in\UB(A)$: From \er{3.4} we infer $\sup A\sms{\vf(n)}\le \a$ and $\vf(n)
\le\a$. From the \df\ of a supremum and the \df\ of~$\vf$, we obtain
$\vf(i)\le \sup A\sms{\vf(n)}$ \fa $i\in[0,n-1]$. By \tr ity of~$\le$, we have
$\vf(i)\le\a$, $i\in[0,n-1]$, hence, also $\vf(i)\le\a$ for $i\in[0,n]$. Thus
$\a\in \UB(A)$.

\ti{$\a$ is the least \ub\ of $A$}: Let $a\in \UB(A)$. Then $\vf(i)\le a$,
$i\in[0,n]$. \E\Ip $\vf(i)\le a$, $i\in[0,n-1]$, hence $\sup A\sms{\vf(n)}
\le a$ as the least \ub\ of $A\sms{\vf(n)}$. Since $\vf(n)\le a$, we have
$(\sup A\sms{\vf(n)})\cup \vf(n)\le a$ by \er{3.4}. Hence $\sup A$ exists and
$\sup A=\a$, which proves \er{3.14}. \If that $\LOR_{i=0}^{\#(A)-1}\vf(i)=\a
=\sup A$. \E\Tf $n+1\in M$. Since $M$~is \iv\ in~$\Na$, $M=\Na$, and \er{3.9}
holds.

\er{3.11}: From Theorem \rfa2{t3.23}, $A\cup B$ is finite and nonempty, so
$\sup A$, $\sup B$ and $\sup A\cup B$ exist by~\er{3.9}. Set $\a:=\sup(A\cup
B)$. Then $x\le\sup A\nde3.4 \le \a$ \fa $x\in A$. Hence $x\le\a$ \fa $x\in A$.
Similarly $x\le\a$ \fa $x\in B$, hence $x\le\a$ \fa $x\in A\cup B$, and $\a\in
\UB(A\cup B)$. Let $a\in\UB(A\cup B)$. Then $x\le a$ \fa $x\in A$ and \fa $x
\in B$. Hence $\sup A\le a$, $\sup B\le a$ and $\a:=(\sup A)\lor(\sup B)\le a$
by~\er{3.4}. Thus $\a=\sup(A\cup B)$.

\er{3.10} follows from \er{3.11} with $A:=A\sms a$, $B:=\{a\}$, and $\sup\{a\}
=a$.

\er{3.12}: If $A\sbs B$, then $B=A\cup B$. Hence $\sup A\nde3.4 \le (\sup A)
\lor(\sup B)\nde3.11 = \sup(A\cup B)=\sup B$.
\endproof

\bex3.4
Let $(X,\le)$ be as in \E\Pr\ \rf{p3.3}. Let $n\in\Na\sms1$ and let $\{A_k\}
_{k\in[0,n]}$ be a family of \nfs s of~$X$. Show that $\sup \bcl_{k\in[0,n]}
A_k$ exists and that
\beq3.15
\sup\bcl_{k\in[0,n]}A_k = \LOR_{k=0}^n \sup A_k.
\e
\eex

\brm3.5
If $(X,\le)$ is an \os\ and if $\ge$ denotes the reverse \og\ of~$\le$ (see
Exercise \rfa1{ex3.4}), then:

\noindent $x$ is an \ub\ of $A$ \wrt $\le$ iff $x$ is a \lo\ of $A$ \wrt $\ge$;

\noindent $x$ is the greatest \el\ of $A$ \wrt $\le$ iff $x$ is the least \el\
of $A$ \wrt $\ge$;

\noindent $x$ is the supremum of $A$ \wrt $\le$ iff $x$ is the infimum of $A$ \wrt
$\ge$.
\erm

\Mo the reverse \og\ of $\ge$ is $\le$. \E\Tf if an \os\ $(X,\le)$ \sf ies the
\cn: $\inf\{x,y\}$ exists \fa $x,y\in X$, $x\ne y$, then it \sf ies the \as\
of \E\Pr\ \rf{p3.3} for the reverse \og~$\ge$. In this case one uses the
notation
\beq3.16
x\land y:=\inf\{x\}\cup\{y\}, \q x,y\in X.
\e
Then the binary \op\ $\land$ on $X$ \sf ies
\bea3.17
& x\land x=x, \q x\in X,\\
& x\land y=y\land x, \q x,y\in X, \lb{3.18} \\
& (x\land y)\land z= x\land(y\land z), \q x,y,z\in X, \lb{3.19} \\
& x\ge x\land y,\ y\ge x\land y, \q x,y\in X, \lb{3.20} \\
& \hbox{if }x\ge a \hbox{ and }y\ge a \hbox{ \fs} a\in X, \hbox{ then }
x\land y\ge a, \hbox{ \fa} x,y\in X, \lb{3.21} \\
& x\ge y \hbox{ iff }y=x\land y, \ x,y\in X. \lb{3.22}
\e

\bex3.6
Find and prove the analogue of \E\Pr\ \rf{p3.3} when $\inf\{x,y\}$ exists \fa
$x,y\in X$, $x\ne y$.
\eex

\Wanp extend \er{1.23}, \er{1.24} to the case where the set $\{x,y\}$ is
replaced by a \nfs\ of~$\Na$. We recall that if $A$ is a~\nfs\ of~$\Na$,
$\prodl_{a\in A}a$ is a common \ml e, and $1$~is a~common divisor of~$A$. \Mo
the set $\Bca_{a\in A} M(a)$ is bounded below by~$1$ \wrt the \og~$\le$ and
$\Bca_{a\in A}D(a)$ is bounded above by every $b\in A$ \wrt the \og~$\le$. Therefore,
by the well-\og\ of $(\N,\le)$ and by Theorem \rfa1{t3.39},
$A$~has a~least common \ml e, denoted by $\lcm(A)$, and a \GCD, denoted
by $\gcd(A)$, where \ti{least and greatest are
understood \hbox{\wrt $\le$}}. Since $(\Na,|)$ is a \lt, it follows from \E\Pr~\rf{p3.3}
and Remark \rf{r3.5} that $\sup A$ and $\inf A$ exist in $(\Na,|)$. We claim
\beq3.23
\lcm(A) = \sup A, \q \gcd(A) = \inf A
\e
\fe \nfs\ of~$\Na$, where $\sup$ and $\inf$ are defined in~$(\Na,|)$.

We prove the first \et y in \er{3.23}. We have $\lcm(A) \le \sup A$, since
$\sup A$~is a common \ml e of~$A$ and $\lcm(A)$ is the least common \ml e of~$A$
\wrt $\le$. \E\oh $\sup A \le \lcm(A)$. Indeed, by \df\ $\sup A | x$ \fe
common \ml e~$x$ of~$A$. \E\Ip $\sup A|\lcm(A)$. From \er{8.57} we infer
$\sup A\le \lcm(A)$, which completes the proof of the first \et y of~\er{3.23}.
The proof of the second \et y is similar.

Clearly, if $(X,\le)$ is a \lt, then $\sup A$ and $\inf A$ exist \fe \nfs\
of~$X$ by \E\Pr\ \rf{p3.3} and Remark \rf{r3.5}.

\smallskip
We now show that if $X$ is a \ns, and if $\di:X\t X\to X$ is a binary \op\
\st \fa $x,y,z\in X$:
\bea3.24
x\di x&= x \qh{(idempotency),}\\
x\di y&= y\di x \qh{(\cmt ity),} \lb{3.25} \\
(x\di y)\di z&= x\di(y\di z) \qh{(\asc ity),} \lb{3.26}
\e
then the \rl\ $\dile$ defined by
\beq3.27
x\dile y \qh{if }y=x\di y, \q x,y\in X,
\e
is an \og\ on $X$ \sf ying the \as\ of \E\Pr\ \rf{p3.3}. \Mo
\beq3.28
x\di y = \sup \{x\} \cup \{y\}, \q x,y\in X,
\e
where the supremum is defined in $(X,\dile)$.

Let $x,y,z\in X$.

\ti{\E\tr ity}: From $x\dile y$, $y\dile z$, we obtain $z\nde3.27 = y\di z
\nde3.27 = (x\di y)\di z \nde3.26 = \break x\di (y\di z) \nde3.27 = x\di z$, hence
$x\dile z$.

\ti{Reflexivity}: $x\nde3.24 = x\di x$ implies $x\dile x$ by \er{3.27}.

\ti{Anti\sy y}: From $x\dile y$ and $y\dile x$, we obtain $x\nde3.27 = x\di
y\nde3.25 = y\di x\nde3.27 = y$. Hence $x=y$.

\ti{$x\di y$ is an \ub\ of $\{x\}\cup\{y\}$}: $x\di(x\di y)\nde3.26 = (x\di
x)\di y\nde3.24 = x\di y$, hence $x\dile x\di y$. In view of \er{3.25} and
$\{x\}\cup\{y\} =\{y\}\cup\{x\}$, we may interchange $x$ and~$y$, and we find
$y\dile y\di x=x\di y$.

\ti{$x\di y$ is the least \ub\ of $\{x\}\cup\{y\}$}: Let $a$ be an \ub\ of
$\{x\}\cup\{y\}$, i.e.\ $x\dile a$ and $y\dile a$, or by
\er{3.27} $a=x\di a$, $a=y\di a$. Then $a=y\di a\nde3.25 = a\di y = (x\di a)
\di y\nde3.25 = (a\di x)\di y\nde3.26 = a\di(x\di y) \nde3.25 = (x\di y)\di
a$. Hence $x\di y\dile a$.

\bds3.7 \

\hph i,i, An \os\ $(X,\le)$ \st $\sup\{x\}\cup\{y\}$ (resp.\
$\inf\{x\}\cup\{y\}$) exists \fa $x,y\in X$ is called a \ti{join-semi\lt}\index{join-semilattice}
(resp.\ \ti{meet-semi\lt}),\index{meet-semilattice} and the binary \op~$\lor$ (resp.~$\land$) on~$X$
defined by \er{3.1} (resp.~\er{3.16}) is called \ti{join} (resp.\
\ti{meet\/}).

\hph ii,, A \ns\ $X$ together with a binary \op~$\di$ $(X,\di)$ \sf ying
\er{3.24}, \er{3.25} and \er{3.26}, i.e.\ an \idp\ \cmt e \sg, is called a
\ti{semi\lt}.\index{semilattice}
\eds

We summarize the above results in

\bpr3.8 \

\hph i,ii, If $(X,\le)$ is a join-\slt, then $(X,\lor)$ is a \slt, and
\beq3.29
x\le y \hbox{ \ iff \  }y=x\lor y \qh{\fa} x,y\in X.
\e

\hph ii,i, If $(X,\le)$ is a meet-\slt, then $(X,\land)$ is a \slt, and
\beq3.30
x\le y \hbox{ \ iff \ }x=x\land y \qh{\fa} x,y\in X.
\e

\hph iii,, If $(X,\di)$ is a \slt, then the \rl\ $\dile$ on~$X$ defined  by
\beq3.31
x\dile y \hbox{ \ if \ } y=x\di y, \q x,y\in X,
\e
is an \og\ on $X$. \Mo $(X,\dile)$ is a join-\slt\ with join~$\lor$, and
\beq3.32
x\di y=x\lor y \qh{\fa} x,y\in X.
\e

\hph iv,, If $(X,\di)$ is a \slt, then the \rl\ $\dile$ on~$X$
defined  by
\beq3.33
x\dile y \hbox{ \ if \ } x=x\di y, \q x,y\in X,
\e
is an \og\ on $X$. \Mo $(X,\dile)$ is a meet-\slt\ with meet~$\land$, and
\beq3.34
x\di y=x\land y \qh{\fa} x,y\in X.
\e
\epr

\brm3.9
If $(X,\le)$ is a \lt\ with join $\lor$ and meet $\land$, then, in view of
\er{3.29} and \er{3.30}, $u\le v$ iff $v=u\lor v$ iff $u=u\land v$ \fa
$u,v\in X$. Since $x\le x\lor y$ \fa $x,y\in X$,
\beq3.35
x=x\land(x\lor y) \qh{\fa} x,y\in X.
\e
Since $x\land y\le x$ \fa $x,y\in X$, $x=(x\land y)\lor x \nde3.25 = x\lor
(x\land y)$ \fa $x,y\in X$, hence
\beq3.36
x=x\lor (x\land y) \qh{\fa} x,y\in X.
\e
\erm

\E\rl s \er{3.35} and  \er{3.36} are usually called ``absorption'' \rl s. It
turns out that these \rl s guarantee that a \ns~$X$ equipped with two binary
\op s $\lor$~and~$\land$ \sf ying \er{3.24}, \er{3.25} and \er{3.26} be a
\lt\ \wrt the \og~$\stackrel\lor\le$ (or~$\stackrel\land\ge$).

\begin{prp}[\cite{Lattice}]
Let $X$ be a \ns\ equipped with two binary \op s $\lor$ and~$\land$ \sf ying
\er{3.24}, \er{3.25} and \er{3.26}. Then
\beq3.37
y=x\lor y \hbox{ \ iff \ } x=x\land y \qh{\fa}x,y\in X,
\e
iff the ``absorption'' \rl s \er{3.35} and \er{3.36} hold. In this case,
$(X,\le)$ is a \lt\ where $x\le y$ if $y=x\lor y$ \fa $x,y\in X$ and
$\lor$ $($resp.~$\land)$ is the join- $($resp.\ meet-$)$ \op\ in $(X,\le)$.
\epr

\proof \

\ti{Only if\/}: follows from Remark \rf{r3.9}.

\ti{If\/}: We suppose that \er{3.35} and \er{3.36} hold. Let $x,y\in X$. If
$y=x\lor y$, then $x\nde3.35 = x\land (x\lor y) = x\land y$. Conversely, if
$x=x\land y$, then $y\nde3.36 = y\lor (y\land x)\nde3.25 = y\lor(x\land y)
=y\lor x\nde3.25 = x\lor y$.

Suppose \er{3.37} holds, then the \rl\ $\stackrel\lor\le$ defined by $x\stackrel
\lor\le y$ if $y=x\lor y$, $x,y\in X$, is an \og\ on~$X$ and $(X,\stackrel\lor\le)$
is a join-semi\lt\ with join~$\lor$ by \E\Pr\ \rf{p3.8}\,(iii). Note that
the \rl~$\stackrel\land\le$ defined on~$X$ by $x\stackrel\land\le y$ if $x=x\land y$,
$x,y\in X$, is an \og\ on~$X$ and $(X,\stackrel\land\le)$ is a meet-semi\lt\
with meet~$\land$ by \E\Pr\ \rf{p3.8}\,(iv). In view of \er{3.37} $\stackrel
\lor\le$ and $\stackrel\land\le$ are identical. If we denote $\stackrel\lor\le$
and $\stackrel\land\le$ by~$\le$, then $(X,\le)$ is a \lt\ with join~$\lor$
and meet~$\land$.
\endproof

\brm3.11
Observe that \cn s \er{3.35} and \er{3.36} imply $x\lor x=x$ and $x\land x=x$
\fa $x\in X$. Indeed, from \er{3.35} where $y$~is replaced by $x\land y$, we
obtain $x=x\land(x\lor\break (x\land y)) \nde3.36 = x\land x$, $x\in X$. Similarly,
$x\nde3.36 = x\lor (x\land(x\lor y))\nde3.35 = x\lor x$, $x\in X$.
\erm

Recall that in the case of the \os\ $(\Na,|)$ we proved the \ex\ of
$\inf\{x,y\}$, $x,y\in\Na$, $x\ne y$, by showing that $\frac{x\cdot
y}{\sup\{x,y\}}$ is the infimum of~$\{x,y\}$ (see the proof of \E\Pr\
\rf{p1.22}. It turns out that this can be done in an arbitrary \PM\
$(X,\qu,e)$ with \og~$\lequ$ defined in \er{8.7}. Thus if $(X,\lequ)$ is a
join-semi\lt, then $(X,\lequ)$ is a \lt. To this end, we introduce some \df\
and prove some preparatory lemma.

\bdn3.12
Let $(X,\qu,e)$ be a \PM\ with \ti{natural\/} \og\ $\lequ$ defined in
\er{8.7}. Let $x,y\in X$ be \st $x\lequ y$. Then
\beq3.38
y \mqu x:=p
\e
where $p$ is the unique \el\ of~$X$ \sf ying $y=x\qu p$ ($=p\qu x)$.
\edn

\bxs3.13 \

\hph i,i, If $(X,\qu,e):=(\N,+,0)$, then
$y\stackrel{\smash{\scriptstyle+}}{{\smash-}\vrule height3pt depth0pt width0pt} x= y-x$, defined in
\er{1.19}.

\hph ii,i, If $(X,\qu,e):=(\Na,\cdot,1)$, then $y
\stackrel{\smash{\scriptstyle\cdot}}{{\smash-}\vrule height3pt depth0pt width0pt} x= \frac
yx$, defined in \er{8.53}.
\exs

\blm8.36
Let $(X,\qu,e)$ be a P-monoid with \nog\ $\lequ$ simply denoted by~$\le$.
Let $\th_a$ and $\d_a$ be the maps from~$X$ into itself defined by:
\bea8.35
\theta_a(x)&:=a\qu x\ ({}=x\qu a) \qh{\fa}a,x\in X, \\
\d_n(x)&:=n\cdot x \qh{\fa}n\in\N \hbox{ and }x\in X, \lb{8.36}
\e
where $n\cdot x$ denotes the $n$-th \IT\ of~$x$ $($for
simplicity we write $n\cdot x$ instead of $n\dqu x)$. Then the \fw\
holds{\rm:}
\bea8.37
{}&x\le y \hbox{ implies }\th_a(x)\le \th_a(y), \ a,x,y\in X,\\
&x\le y \hbox{ implies }\d_n(x)\le \d_n(y), \ n\in\N,\ x,y\in X, \lb{8.38}\\
&x<y \hbox{ iff }\th_a(x)<\th_a(y), \ a,x,y\in X, \lb{8.39}\\
&x< y \hbox{ implies }\d_n(x)< \d_n(y), \ n\in\Na,\ x,y\in X, \lb{8.40}\\
&\th_a \hbox{ is in\jc\ \fe }a\in X, \lb{8.41}\\
&z\mqu y\le \ (\hbox{resp.\ }<)\ z\mqu x \qh{\fa}x,y,z\in X \lb{8.42}\\
&\qquad \hbox{ \st $x,y\le z$ and $x\le{}$ $($resp.\ $<{})\ y$ where $\mqu$
is defined in \er{3.38}.} \non
\e
\elm

\proof \

\er{8.37} follows from \er{8.39} and \er{8.38} follows from \er{8.40}.

\er{8.39}: $x<y$ iff $p:=y\mqu x\in X\sms e$ iff $y=x\qu p$. Since
$\th_a(y)=\th_a(x\qu p)=a\qu x\qu p=\th_a(x)\qu p$, $x<y$ implies $\th_a(x)
<\th_a(y)$. Conversely, if $\th_a(x)<\th_a(y)$, then $\th_a(y)=\th_a(x)\qu q
$ where $q:=\th_a(y)\mqu\th_a(x)\in X\sms e$. Then $y\qu a=a\qu y=a\qu x\qu q=x\qu q
\qu a$, which implies $y=x\qu q$ by \era2{1.8}. Hence $x<y$.

\er{8.40}: Let $x<y$, $p:=y\mqu x\in X\sms e$. Then $y=x\qu p$ and
$n\cdot (x\qu p)\nda22.3 =\break (n\cdot x)\qu (n\cdot p)$. We claim that
$n\cdot p\in X\sms e$ when $n\in\Na$ and $p\in X\sms e$. Indeed, $n\cdot p
=\QU_{k\in\zo0,n } a_k$ with $a_k:=p$, $0\le k< n$, by \era2{1.126} and \era2{3.5}.
Suppose for \cd ion that $\QU_{k\in\zo0,n } a_k=e$, then by
\era2{7.17}, \era2{7.15} $p=a_1=e$,
a~\cd ion. Since $n\cdot p\in X\sms e$, we infer $\d_n(x)<\d_n(y)$.

\er{8.41} follows from \era2{1.8}.

\er{8.42}: Let $x,y,z\in X$ be \st $x,y\le z$ and $x<y$. Then $z=x\qu (z
\mqu x)$, $z=y\qu (z\mqu y)$ and $y=x\qu p$ where $p=y\mqu x\in X\sms e$. Then
$z=x\qu (z\mqu x)=\break y\qu (z\mqu y)=x\qu p\qu (z\mqu y)$. By \era2{1.8}, $z\mqu x
=p\qu (z\mqu y)$, hence $z\mqu y<z\mqu x$. If $x=y$, then $p=e$ and
$z\mqu x=z\mqu y$.
\endproof

We are now in a position to show that if a P-monoid $X$ possesses a join binary
\op, then $X$ is a \lt\ \wrt its \nog.

\blm8.37
Let $(X,\qu,e)$ be a nontrivial P-monoid. If $\sup\{x,y\}$ exists \fa ${x,y
\in X}$ with $x\ne y$, then $(X,\qu,e)$ is a~\lt\ \wrt its \nog\ and
\beq8.43
x\qu y=(x\lor y)\qu (x\land y) \qh{\fa}x,y\in X.
\e
\elm

\proof
If $x=y$, then $x=x\lor x=x\land x$ and \er{8.43} holds. Suppose $x\ne y$.
Set $x\lor y:=\sup\{x,y\}$. By \df\
of~$\le$, we have $x\le x\qu y$ and $y\le x\qu y$, hence $x\qu y$ is an
\ti{\ub\/} of $\{x,y\}$. Hence $x\lor y \le x\qu y$,
since $x\lor y$ is the \ti{least \ub\/} of
$\{x,y\}$. Set $p:=(x\qu y)\mqu (x\lor y)$. We claim that
\beq8.44
p=\inf\{x,y\}.
\e

(i) \ti{$p$ is a \lo\ of $\{x,y\}$}: From \er{8.42} with $z:=x\qu y$, $x:=x$,
$y:=x\lor y$, we obtain $p=(x\qu y)\mqu (x\lor y)\le (x\qu y)\mqu x=y$.
Similarly with $z:=x\qu y$, $x:=y$ and $y:=x\lor y$, we obtain
$p=(x\qu y)\mqu (x\lor y)\le (x\qu y)\mqu y=x$.

(ii) \ti{$p$ is the greatest \lo\ of $\{x,y\}$}: Let $b\in X$ be \st $b\le x$
and $b\le y$. We have to show that $b\le p$. Since $b\le x\le x\qu y$ and
$b\le y\le x\qu y$, we obtain from \er{8.42} that
\[
y=(x\qu y)\mqu x\le (x\qu y)\mqu b \qhq{and}
x=(x\qu y)\mqu y\le (x\qu y)\mqu b.
\]
Hence $(x\qu y)\mqu b$ is an \ub\ of $\{x,y\}$. By \df\ $x\lor y$ is the least
\ub\ of $\{x,y\}$, then $(x\lor y)\le (x\qu y)\mqu b$. From \er{8.37} with
$a:=b\qu p$, $x:=x\lor y$ and $y:=(x\qu y)\mqu b$, we obtain
\[
(x\lor y) \qu b\qu p \le ((x\qu y)\mqu b)\qu b\qu p.
\]
Note that $(x\lor y)\qu b\qu p=(x\lor y)\qu p\qu b=((x\lor y)\qu p)\qu b
=(x\qu y)\qu b$, and\break $((x\qu y)\mqu b)\qu b\qu p=(((x\qu y)\mqu b)\qu b)\qu p
=(x\qu y)\qu p$. Hence $(x\qu y)\qu b\le\break (x\qu y)\qu p$.
From \er{8.39} with $a:=x\qu y$, $x:=b$  and $y:=p$ we conclude that
$b\le p$. So \er{8.44} holds.

(iii) $p=x\land y$ \ti{and} \er{8.43}: Since $p=\inf
\{x,y\}$ if $x\ne y$, we have $x\land y=p$, by \er{3.16}. \E\Tf $x\qu y
=(x\lor y)\qu p=(x\lor y)\qu (x\land y)$ if $x\ne y$.
Hence \er{8.43} holds and $(X,\qu,e)$ is a~\lt.
\endproof

\bex3.16
Let $(X,\qu,e)$ be a nontrivial \PM\ \st $\sup\{x,y\}$ exists \fa $x,y\in X$
with $x\ne y$. Let $A$ be a \nfs\ of~$X$ such that
$\inf \{x,y\} =  e$  for all  $x,y  \in A$  with $x \ne y$. Show
\beq3.49
\sup A = \prod_{a\in A} a.
\e
\eex

\bdf8.34
A P-monoid $(X,\qu,e)$ which is a \lt\ \wrt its \nog\ $\lequ$ is called an
\ti{L-monoid\/}\index{L-monoid}. The join and meet \op s are denoted by $\stackrel\qu\lor$
and $\stackrel\qu\land$ (or simply by $\lor$ and $\land$ if no confusion
arises).
\edf

\brm8.35
The term \ti{L-monoid\/} is not standard. Note that the term $\ell$-monoid
is used for ordered monoids which are \lt s (see~\cite[p.~323]{Lattice}).
An~L-monoid is an $\ell$-monoid but the converse is not always true. For
example, the monoid $(X,\lor,0)$ of Exercise \rf{ex1.3} where $X:=\{0,1\}$
and $\lor$ \sf ies \er{8.2} is an $\ell$-monoid but not an L-monoid (since
it is not \cnc e).
\erm

We now prove that in an L-monoid the map $\th_a$
defined in \er{8.35} is join and meet preserving and that an L-monoid is
a~\dsb e \lt.\index{distributive lattice}

\bth3.19
Let $(X,\qu,e)$ be an L-monoid and let $\le$ denote its \nog. Then \er{8.43}
holds as well as
\bea8.45
\th_a(x\lor y)&=\th_a(x) \lor \th_a(y) \qh{for all} a,x,y\in X, \\
\th_a(x\land y)&=\th_a(x) \land \th_a(y) \qh{for all} a,x,y\in X. \lb{8.46}
\e
\Mo $(X,\le)$ is a \emph{\dsb e \lt}, i.e.\
\bea3.52
x\lor(y\land z)&=(x\lor y)\land(x\lor z)\qh{\fa}x,y,z\in X, \\
x\land(y\lor z)&=(x\land y)\lor(x\land z)\qh{\fa}x,y,z\in X. \lb{3.53}
\e
\eth

\proof \er{8.43} is a direct con\sq\ of Lemma \rf{l8.37}.

\er{8.45} --- \ti{$\th_a(x\lor y)$ is an \ub\ of $\{\th_a(x),\th_a(y)\}$}:
Since $x\le x\lor y$ and $y\le x\lor y$, we have in view of \er{8.37}
$\th_a(x)\le\th_a(x\lor y)$ and $\th_a(y)\le\th_a(x\lor y)$.

\ti{$\th_a(x\lor y)$ is the least \ub\ of $\{\th_a(x),\th_a(y)\}$}:
Let $b\in X$ be \st $\th_a(x)\le b$ and $\th_a(y)\le b$. We have to show
that $\th_a(x\lor y)\le b$. Set $p:=b\mqu \th_a(x)$. Then
$b=\th_a(x)\qu p=a\qu(x\qu p)$. Hence $a\le b$. Set $c:=b\mqu a$.
Then $b=a\qu c$, $\th_a(x),\th_a(y)\le a\qu c$. So $a\qu x\le a\qu c$
and $a\qu y\le a\qu c$. It follows from \er{8.39} that $x\le c$ and
$y\le c$. Hence $c$~is an \ub\ of $\{x,y\}$. \E\Tf $x\lor y\le c$.
Then $\th_a(x\lor y)=a\qu (x\lor y)\nde8.37 \le a\qu c=b$. This completes
the proof of \er{8.45}.

\er{8.46}: From $x\qu y\nde8.43 = (x\lor y)\qu(x\land y)$ we obtain
$(a\qu a)\qu(x\qu y)=\break a\qu a\qu(x\lor y)\qu(x\land y)$. Hence
$(a\qu x)\qu(a\qu y)=(a\qu(x\lor y))\qu(a\qu(x\land y))$, i.e.\
$\th_a(x)\qu\th_a(y)=\th_a(x\lor y)\qu\th_a(x\land y)$. Using \er{8.45}
we obtain $(\th_a(x)\lor\th_a(y))\qu\break (\th_a(x)\land\th_a(y))=
(\th_a(x)\lor\th_a(y))\qu\th_a(x\land y)$. Using \era2{1.8} we obtain
\er{8.46}.

\er{3.52}: To prove $x\lor(y\land z)=\inf \{x\lor y,x\lor z\}$:

\ti{$x\lor(y\land z)$ is a \lo\ of $\{x\lor y,x\lor z\}$}.
Note that in every \lt\ $(L,\le)$ we have
\bea8.47
a\le b &\qhq{implies} a\lor c\le b\lor c\\
a\le b &\qhq{implies} a\land c\le b\land c \lb{8.48}
\e
for all $a,b,c\in L$.

\bex8.39
Prove \er{8.47} and \er{8.48}.
\eex

Since $y\land z\le y$ and $y\land z\le z$, we get from \er{8.47}:
$x\lor(y\land z)$ is a \lo\ of $\{x\lor y,x\lor z\}$.

\ti{$x\lor(y\land z)$ is the greatest \lo\ of $\{x\lor y,x\lor z\}$}.
To prove that $a\le x\lor y$ and $a\le x\lor z$ implies $a\le x\lor
(y\land z)$. We have by \er{8.37}, \er{8.43}: $a\qu(x\land y)\le
(x\lor y)\qu(x\land y)=x\qu y$, similarly $a\qu(x\land z)\le
(x\lor z)\qu(x\land z)=x\qu z$. Hence $[a\qu(x\land y)]\land
[a\qu(x\land z)]\le (x\qu y),(x\qu z)$, so
\bmlg
a\qu[(x\land y)\land(x\land z)]\nde8.46 = [a\qu (x\land y)]\land
[a\qu(x\land z)] \le (x\qu y)\land(x\qu z)\\ {}\nde8.46 =
x\qu(y\land z)=(x\lor (y\land z))\qu(x\land(y\land z)).
\e
Since $(x\land y)\land(x\land z)=x\land(y\land z)$ by
\er{3.24}, \er{3.25}, \er{3.26}, we obtain
$a\qu\break (x\land y\land z)\le [x\lor(y\land z)]\qu(x\land y\land z)$.
From \er{8.39} with $\th_{x\land y\land z}$ we obtain $a\le x\lor(y\land z)$.

This completes the proof of \er{3.52}.

\er{3.53}: Imitate the proof of \er{3.52} or use the \fw:
$$
\dsl{\indent
(x\land y)\lor(x\land z)\nde3.52 = [(x\land y)\lor x]\land [(x\land y)\lor
 z]\nad{\er{3.25},\er{3.36}}= x\land [z\lor(x\land y)]\hfill\cr
{}\nde3.52 = x\land [(z\lor x)\land(z\lor y)]\nde3.26 =
[x\land(z\lor x)]\land (z\lor y) \cr
\hfill \nad{\er{3.25},\er{3.35}} =
x\land(z\lor y)\nde3.25 = x\land(y\lor z).\indent\Box}
$$

\goodbreak
\brm3.21
Not all \lt s are \dsb e. As an example consider the set
$X=\{a,b,c,d,e\}$ endowed with the \og\ \sf ying $a\le b\le e$, $a\le c\le
e$, and $a\le d\le e$.
\[
\xymatrix{&e\\ b\ar[ur]&c\ar[u]&d\ar[ul]\\ &a\ar[ul]\ar[u]\ar[ur]}
\]
One verifies that $(X,\le)$ is a \lt. However,
 $b\lor (c\land d)=b\lor a=b$, $b\lor c=e$, $b\lor d=e$ and
$(b\lor c)\land(b\lor d)=e\ne b$, $b\land (c\lor d)=b\land e=b$,
$b\land c=a$, $b\land d=a$ and
$(b\land c)\lor(b\land d)=a\ne b$.
\erm

\begin{exe}[\cite{Lattice}] \lb{ex8.33}
Let $(X,\le)$ be a \lt.

\hph i,i, Show that the \fw\ statement holds
\beq8.32
\bal
x\land(z\lor y)&=x\land [(z\lor x)\land(z\lor y)],\\
x\lor(z\land y)&=x\lor [(z\land x)\lor(z\land y)],
\eal
\e
\fa $x,y,z\in X$.

\hph ii,, Show that \er{3.52} holds iff \er{3.53} holds.
\eex

\bex3.22
Let $Y$ be a \ns\ and let $(\cP(Y),\sbs)$ denote the set of all subsets of~$Y$
ordered by inclusion (see Exercise \rfa1{ex3.4}\,(ii)). Show that
$(\cP(Y),\sbs)$ is a \ti{\dsb e} \lt\ with union as join and intersection as
meet.
\eex

\bdf8.40
A Boolean algebra $B$ is a \ti{\dsb e \lt} having a least \el~$0$ and\index{Boolean algebra}
a greatest \el~I \st the \fw\ condition is \sf ied:
\beq8.49
\hbox{\E\fe $a\in B$ \te s a unique \el\ $a^c\in B$}
\e
(called the complement of $a$) \st\index{complement}
\beq8.50
a\lor a^c={\rm I} \qhq{and} a\land a^c=0.
\e
\edf

Note that
\beq8.51
\bal
&(a^c)^c=a \qh{\fe} a\in B;\\
&(a\lor b)^c=a^c\land b^c;\ (a\land b)^c=a^c\lor b^c
\qh{\fe} a,b\in B.
\eal
\e

\bex8.41
Show that if $Y$ is a \ns\ and $\cP(Y)$ is the power set ordered by
inclusion, then $(\cP(Y),\sbs)$ is a Boolean algebra where $0:=\vn$ and
$I:=Y$.
\eex

\begin{dfn}[\cite{Lattice}] \lb{d8.42}
A \ti{sub\lt}\index{sublattice} of a \lt\ $(L,\le)$ is a nonempty subset $M$ of~$L$ \st
$x\lor y\in M$ and $x\land M$ whenever $x,y\in M$.
\edf

\begin{cau}
A nonempty subset $N$ of a \lt\ $(L,\le)$ may be a \lt\ \wrt the \rt ion
of the \og\ of~$L$ without being a sub\lt\ of~$L$.
\end{cau}

\bxa8.43
The \os\ $(X,\le)$ defined by $X=\{a,b,c,d,e\}$, $a\le b$, $a\le c$,
$a \le d$, $a\le e$, $b\le d$, $b\le e$, $c\le d$, $c\le e$, $d\le e$,
whose diagram is
\[
\xymatrix{&e\\ &d\ar[u]\\ b\ar[ur]\ar[uur]&&c\ar[ul]\ar[uul]\\
&a\ar[ul]\ar[uu]\ar[ur]}
\]
is a \lt\ with $b\lor c=d$. The subset $Y=X\sms d$ is \ti{not\/} a sub\lt\
of~$X$ but it is a \lt\ \wrt the \rt ion of the \og~$\le$.
\exa

\bex3.27
Let $(X,\qu,e)$ be a nontrivial L-monoid and let $a,b\in X$ with $a<b$.
Set $B:=\{x\in X: a\le x\le b\}$. Show that $B$ is a \dsb e sub\lt\ with
greatest \el~$b$ and least \el~$a$. Can $B$ be made a Boolean algebra?
\eex

We conclude this section by giving some applications and some reformulation of the \fd\ theorem of
\art. As a first application we show that $\N\t\N$ and~$\N$ are \ep. Since
$\Na\ax\N$, it is \sft\ by \era1{4.1} to show $\N\t\N\ax\Na$. In view of
Schr\"oder--Bernstein's Theorem \rfa1{t4.30}, it is also \sft\ to find in\jc\
maps from~$\Na$ to~$\N\t\N$ and from $\N\t\N$ to~$\Na$. The map $\psi:\Na
\to\N\t\N$ defined by $\psi(m):=(m,0)$, $m\in\Na$, is clearly in\jc. \E\oh
the map $\vf:\N\t\N \to \Na$ defined by $\vf(m,n):=2^{m+1}3^{n+1}$, is in\jc.
Indeed, let $m,n,m',n'\in\N$ and suppose $2^{m+1}3^{m+1}=2^{m'+1}3^{n'+1}$. We
want to show that $m=m'$ and $n=n'$. To this end we apply Theorem \rf{t9.20}.

Define $f:[0,m+n+1]\to\Na$ and $g:[0,m'+n'+1]\to\Na$ by setting $f(i):=2$,
$i\in[0,m]$, $f(i):=3$, $i\in[m+1,m+n+1]$; $g(j):=2$, $j\in[0,m']$, $g(j):=3$,
$j\in[m'+1,m'+n'+1]$. Let $M:=m+n+1$, $N:=m'+n'+1$, $A:=\{2,3\}$ and $B:=\{2,
3\}$. Note that $a:=\prodl_{i=0}^M f(i)=2^{m+1}3^{n+1}$ and $b:=\prodl_{j=0}^N
g(j)=2^{m'+1}3^{n'+1}$, and that $m+1$ (resp.\ $m'+1$) is the \ml icity of~$2$
in the \cme product $\prodl_{i=0}^M f(i)$ (resp.\ $\prodl_{j=0}^N g(j)$).
Similarly for $n+1$ (resp.\ $n'+1$). \If from Theorem \rf{t9.20} that $m+1=
m'+1$ and $n+1=n'+1$, hence $m=m'$, $n=n'$ by \era2{1.8}. Thus $\vf$~is in\jc\
and $\N\t\N$ is \ep\ to~$\Na$, hence also to~$\N$.

\bex3.28 \

\hph i,ii, Construct an explicit bi\jn\ between $\N\t\N$ and $\N$.

\hph ii,i, Let $X,Y,X',Y'$ be \ns s \st $X\ax X'$ and $Y\ax Y'$. Show that
$X\t Y\ax X'\t Y'$.

\hph iii,, Prove by \In\ on $n\in\N$ that
\beq3.60
\N^{[0,n]} \hbox{ is \ct y infinite \fa }n\in\N.
\e
\eex

We now show that $\N^{\N}$ is \ti{un\ct e}. By Cantor's Theorem \rfa1{t4.26},
$\cP(\N)$ is un\ct e. We claim that $\cP(\N)$ and $\{0,1\}^\N$ are \ep. \E\fe
subset $A$ of~$\N$, we define the \ti{indicator \f} $1_A:\N\to\{0,1\}$ by
setting $1_A(x):=1$, $x\in A$, and $1_A(x):=0$, ${x\notin A}$. The map $A\mt
1_A$ from $\cP(\N)$ into $\{0,1\}^\N$ is bi\jc. Indeed, \fe $g:\N\to\{0,1\}$
set $A:=\{x\in\N: g(x)=1\}$. We have $g=1_A$, hence $A\mt 1_A$ is sur\jc.
\Mo if $A,B\in\cP(\N)$ and $1_A=1_B$, then $A=\{x\in\N: 1_A(x)=1\}=\{{x\in\N}:
1_B(x)=1\}=B$. Hence, the map $A\mt1_A$ is in\jc. Thus the claim is proved.
Since $A\mt1_A$ is in\jc\ and $\cP(\N)$ is un\ct e, $\{0,1\}^\N$ is un\ct e
by \E\Pr\ \rfa1{p4.27}\,(A)(ii). In view of the same \Pr, it suffices to find
an in\jc\ map from $\{0,1\}^\N$ into $\N^\N$. Given $f:\N\to\{0,1\}$, we
define $j(f):\N\to\N$ by setting $j(f)(n):=f(n)$ \fa $n\in\N$.
One verifies that the map $j:\{0,1\}^\N \to\N^\N$ is in\jc.
Thus
\beq3.61
\N^\N \hbox{ is un\ct e.}
\e

We next consider a subset of $\N^\N$ \ep\ to~$\N$ which contains in an
appropriate sense $\N^{[0,n]}$ \fa $n\in\N$. Given $n\in\N$, let $j_n:\N
^{[0,n]}\to\N^\N$ be defined by
\beq3.62n
(j_n(a))(m):= \bca
a_m &\hbox{for }m\le n,\\
0 &\hbox{for }m>n,
\eca \qquad a\in\N^{[0,n]}.
\e
One easily verifies that $j_n$ is in\jc. We denote by $\N^{[0,n]}_0$ the range
of~$j_n$ in~$\N^\N$, and we define ${\bf0}\in\N^\N$ by setting
\beq3.63n
{\bf0}(k) =0 \qh{\fa }k\in\N.
\e
Observe that \fa $n\in\N$
\beq3.64n
\N^{[0,n]}_0 = \{a\in\N^\N: a(k)=0 \hbox{ \fa}k>n\}.
\e
\E\Ip
\beq3.65n
{\bf0}\in\N^{[0,n]}_0 \hbox{ \fa}n\in\N.
\e
Since $j_n:\N^{[0,n]}\to\N^{[0,n]}_0$ is bi\jc, $\N^{[0,n]}$ and $\N^{[0,n]}_0$
are \ep\ \fa $n\in\N$.

Note that
\beq3.66n
\N^{[0,n]}_0 \sbs \N^{[0,m]}_0, \q n\le m,\ n,m\in\N,
\e
and let
\beq3.67n
\N^\N_0 := \bigcup_{n\in\N} \N^{[0,n]}_0
\e
(not standard notation).

Observe that if $a\in\N^\N_0$, then either $a=\bf0$ or \te s \ooo $k\in\N$
denoted by~$d(a)$ \st $a_{d(a)}\ne0$ and $a_k=0$ \fa $k>d(a)$. Indeed, if
$a\in\N^\N_0\sms{\bf0}$, then by \er{3.67n} \te s $n\in\N$ \st $a\in\N^{[0,n]}_0$.
\E\Tf $a_i=0$ \fa $i>n$ and \te s $\ov k\in[0,n]$ \st $a_{\bar k}\ne0$. Set
$B:=\{k\in[0,n]:a_k\ne0\}$. Then $B$~is nonempty and bounded, hence by Theorem
\rfa1{t3.39} $B$~possesses a greatest \el\ denoted by~$d(a)$. \csq, $a_{d(a)}
\ne0$ and $a_i=0$ \fa $i>d(a)$.

Our goal is to find a bi\jn\ $\vf:\N_0^\N \to\Na$. To this end we recall that
the set of primes $P$ is infinite by \E\Pr\ \rf{p9.11}, hence unbounded in
$(\N,\le)$ by Lemma \rfa1{l4.20}. \If from Theorem \rfa1{t3.49} that $(P,\le)$
is \ois c to $(\N,\le)$. \E\Tf \te s a strictly in\cre\ \sq\ $\zb pk\N$ \sf
ying $P=\bcl_{k\in\N}\{p_k\}$. \E\Ip $p_0:=2$, $p_1:=3$, $p_2:=5$. We now
define $\vf:\N_0^\N\to \Na$ by setting
\beq3.68n
\bca
\vf({\bf0}):=1,\\
\vf(a):= \prodl_{k=0}^{d(a)} p_k^{a_k},\q a\in\N_0^\N\sms{\bf0},
\eca
\e
where $d(a)$ is  defined above and $\prodl_{k=0}^{d(a)}$ is the \cme product
in the monoid $(\Na,\cdot,1)$.

\bpr3.30n
The map $\vf$ defined in \er{3.68n} is bi\jc.
\epr

\proof
\ti{Sur\ji}: Let $n\in\Na$. If $n:=1$ then $\vf({\bf0})=1$. If $n>1$, then
$n\nde9.19 = \prodl_{\sbk{l\in P\\l|n}} l^{m_l}$ where $m_l$ is the \el\
of~$\Na$ \st $l^{m_l}|n$ and $l^{m_l+1}\nmid n$. Set $A:=\{l\in P: l|n\}$.
Then $A$ is nonempty and bounded since $l|n$ implies $l\le n$ by \er{8.57}.
By Theorem \rfa1{t3.39} $A$~possesses a greatest \el\ in $(\N,\le)$ which we
denote by~$\ov l$. In view of the \df\ of the \sq\ $\zb pk\N$, \te s \ooo $\ov k
\in\N$ \st $\ov l=p_{\bar k}$. Thus if a \Pn~$l$ divides~$n$, then $l=p_k$
\fs $k\in[0,\ov k]$. Set $m:=\prodl_{k=0}^{\bar k} p_k^{a_k}$ where $a_k:=0$
if $p_k\nmid n$ and $a_k:=m_{p_k}$ if $p_k|n$. Set $I:=\{k\in[0,\ov k]:
p_k|n\}$ and $J:=\{k\in[0,\ov k]: p_k\nmid n\}$. Then $[0,\ov k]=I \cup J$
and $I\cap J=\vn$.

Thus $m=\prodl_{k=0}^{\bar k} p_k^{a_k} \nda27.17 = \prodl_{k\in[0,\bar k]}
p_k^{a_k} \nda21.138 = \prodl_{k\in I}p_k^{a_k} \cdot \prodl_{k\in J}
p_k^{a_k}\nda22.23 = \prodl_{k\in I}p_k^{a_k}\cdot \prodl_{k\in J} 1\nad* =
{\bg(\prodl_{k\in I}p_k^{a_k})\cdot1} \break\nda22.13 = \prodl_{k\in I}p_k^{a_k}
= \prodl_{k\in I}p_k^{m_{p_k}}= \prodl_{\sbk{p\in P\\p|n}}p^{m_p}=n$. In $\nad*=$
we used $\prodl_{k\in J}1 \nda21.128 = \#(J)\ddt 1 \nad{\era2{2.3}\,\rm I4} =
1$. Thus $n=\vf(a)$. Hence $\vf$~is sur\jc.

\ti{In\ji}: Note that if $a\in\N_0^\N\sms{\bf0}$, then $\vf(a)>1$. Indeed, if
$a\in\N_0^{[0,n]}$ \fs $n\in\N$ and $\vf(a)=1$, then we have $\prodl_{k=0}^n
p_k^{a_k}=1$. Since $(\Na,\cdot,1)$ is a \PM, we obtain $p_k^{a_k}=1$ \fa
$k\in[0,n]$ by \era2{7.15}. Since $p_k\ne1$, $k\in[0,n]$, we have $a_k=0$,
$k\in[0,n]$, by \era2{2.32}. \E\Tf $a=\bf0$. Thus if $\vf(a)=\vf(b)=1$,
$a,b\in\N_0^\N$, then $a=b=\bf0$. Suppose now $\vf(a)=\vf(b)>1$. Then $a,b
\in\N_0^\N\sms{\bf0}$.

Let $d(a),d(b)\in\N$ be as above. We have $\prodl_{i=0}^{d(a)} p_i^{a_i}=
\prodl_{j=0}^{d(b)} p_j^{b_j}$, and we want to show that $a=b$. We first prove
that $d(a)\not< d(b)$. Suppose for \cd ion that \te s $r\in\Na$ \st $d(b)=
d(a)+r$. Then $\prodl_{i=0}^{d(a)} p_i^{a_i}= \bigl(\prodl_{i=0}^{d(a)}p_i
^{a_i}\bigr)\cdot1 \nda27.13 = \bigl(\prodl_{i=0}^{d(a)}p_i^{a_i}\bigr)
\cdot\bigl(\prodl_{j=d(a)+1}^{d(a)+r} p_j^{a_j}\bigr)$, hence $\suml_{j=d(a)+1}
^{d(a)+r} p_j^{a_j}=1$ by \era2{2.16}\,(iii) since $\prodl_{i=0}^{d(a)}
p_i^{a_i}\ne0$. Thus by \era2{7.12} $\suml_{j=1}^r p_{d(a)+j}^{a_{d(a)+j}}=1$
which implies $p_{d(b)}^{d(b)}=1$ by \era2{7.15}, hence $d(b)=0$ by \era
2{2.32}, a \cd ion. Interchanging $a$ and~$b$ we find $d(b)\not< d(a)$, hence
$d(a)=d(b)$. We now show that $a_k=b_k$ \fa $k\in[0,d(a)]$. Suppose for \cd
ion that \te s $\ov k\in[0,d(a)]$ \st $a_{\bar k}\ne b_{\bar k}$. If $d(a)=
d(b)=0$, then $\ov k=0$ and $p_0^{a_0}= \prodl_{i=0}^{d(a)}p_i^{a_i}=
\prodl_{j=0}^{d(b)}p_j^{b_j}=p_0^{b_0}$. Then $a_0=b_0$ by \era2{2.34},
a~\cd ion. Thus $a=b$. We now assume $d(a)=d(b)>0$. If $\ov k=0$, then
$\prodl_{i=0}^{d(a)}p_i^{a_i}=p_0^{a_0}\cdot c$ and $\prodl_{j=0}^{d(b)}p_j
^{b_j} = p_0^{b_0}\cdot d$ where $c,d$ are (composite) products of primes
different from~$p_0$. From $p_0^{a_0}\cdot c =p_0^{b_0}\cdot d$ and $a_0\ne
b_0$ we infer as in the proof of Theorem \rf{t9.20} that $p_0$ divides
$c$ or~$d$ \cd ing Lemma \rf{l9.17}. Hence $a=b$. Finally, if $\ov k\in\oz
0,d(a) $, then $\prodl_{i=0}^{d(b)}p_{\si(i)}^{a_{\si(i)}} = \prodl_{j=0}^{d(b)}
p_{\si(j)}^{b_{\si(j)}}$ where $\si$~is the \Pm\ of $[0,d(a)]$ \sf ying
$\si(0):=\ov k$, $\si(\ov k):=0$, and $\si(k):=k$ otherwise. From $\prodl_{i=0}
^{d(b)}p_{\si(i)}^{a_{\si(i)}} = \prodl_{j=0}^{d(a)}p_{\si(j)}^{b_{\si(j)}}$
we obtain as above $p_{\bar k}^{a_{\bar k}}\cdot c =
p_{\bar k}^{b_{\bar k}}\cdot d$ \fs $c,d$ (composite) products of primes
different from~$p_{\bar k}$. As above we obtain a \cd ion. Hence $a=b$. This
concludes the proof of the in\ji\ of~$\vf$.
\endproof

\brm3.31n
\If from \E\Pr\ \rf{p3.30n} that $\N_0^\N$ is \ep\ to~$\Na$ hence $\N_0^\N$~is
\ct y infinite. Since $\N^{[0,n]}\ax \N^{[0,n]}_0$ and $\N^{[0,n]}_0\sbs\N_0
^\N$ \fa $n\in\N$, it follows at once that $\N^{[0,n]}$ is \ct y infinite \fa
$n\in\N$. Indeed, $\N^{[0,0]}\sbs\N^{[0,n]}$ implies $\N^{[0,n]}$~is infinite
and $\N^{[0,n]}$ is \ct e \fa $n\in\N$ since $\N^{[0,n]}\sbs\N_0^\N$ \fa
$n\in\N$ by \E\Pr\ \rfa1{p4.27}\,(A)(iv).
\erm

We now mention \ad al \pp ies of the map $\vf$ defined in \er{3.68n}. An \ad,
denoted by~$+$, can be defined on $\N_0^\N$, which makes $(\N_0^\N,+,{\bf0})$
an L-monoid. Then $\vf$ becomes a~monoid- and an \ois sm between
$(\N_0^\N,+,{\bf0})$ and the L-monoid $(\Na,\cdot,1)$. \E\Ip $\vf$~is join
and meet preserving.

\bpr3.29
Let $a,b\in\N^\N$. Let $+$ denote the binary \eq\ on $\N^\N$ defined by
\beq3.65a
(a+b)(m):= a(m)+b(m) \qh{\fa}m\in\N.
\e
Then $(\N^\N,+,{\bf0})$ is an L-monoid.

\E\Ip given $a,b\in\N^\N$,
\bea3.66
a \stackrel+\le b \hbox{ iff }a(m)\le b(m) &\quad \hbox{\fa}m\in \N,\\
(a\lor b)(m)=\max(a(m),b(m)) &\q \hbox{\fa}m\in \N, \lb{3.67} \\
(a\land b)(m)=\min(a(m),b(m)) &\q \hbox{\fa}m\in \N. \lb{3.68}
\e
\Mo $\N_0^\N$ is a submonoid of $(\N^\N,+,{\bf0})$, and $(\N_0^\N,+,{\bf0})$
is a \PM\ whose \nog\ is identical to the \rt ion to $\N_0^\N$ of the \og~$
\stackrel+\le$. The ordered set $(\N_0^\N,\stackrel+\le)$ is a sub\lt\ of the
\lt\ $(\N^\N,\stackrel+\le)$, hence if $a,b\in\N_0^\N$, $a\ne b$, then
$\sup\{a,b\}$ $($resp.\ $\inf\{a,b\})$ \wrt $\stackrel+\le$ in~$\N_0^\N$ is
given by \er{3.67} $($resp.~\er{3.68}$)$.
\epr

\bex3.30
Prove \E\Pr\ \rf{p3.29}.
\eex

\brm3.31 \

\hph i,i, $(a\lor b)+(a\land b)=a+b \qh{\fa} a,b\in\N^\N$.

\hph ii,, Let $M$ be the \sbm\ of $(\N,+,0)$ defined in Remark \rf{r1.15}.
One verifies that
$M=\N\sms1$. Clearly $(M,+,0)$ is a P-monoid. However, its \nog\ is not
identical to the \rt ion of the \nog\ $\le$ of $(\N,+,0)$. Let us denote by
$\stackrel M\le$ the \nog\ of $(M,+,0)$. Indeed, given $x,y\in M$, $x\stackrel
M\le y$ holds iff \te s $p\in M$ \st $y=x+p$. Since $1\notin M$, we have
$2\stackrel M{\not<}3$ but $2<3$.
\erm

\blm3.32
Let $(X,\qu,e)$ be a \PM, and let $M$ be a \sbm\ of~$X$. Then $(M,{\qu},e)$
is a \PM\ whose \nog\ is denoted by $\stackrel M\le$. The \og\ $\stackrel M\le$
on~$M$ is identical to the \rt ion to~$M$ of the \og\ $\lequ$ on~$X$ iff \fa
$x,y\in M$ \st $x\lequ y$ we have $y\mqu x\in M$.
\elm

\bex3.33
Prove Lemma \rf{l3.32}.
\eex

\bth3.34
Let $\vf:\N_0^\N \to (\Na,\cdot,1)$ be as in \er{3.68n}. Then $\vf$~is bi\jc\
and moreover, \fa $a,b\in\N_0^\N${\rm:}
\bea3.69
\vf(a+b) &= \vf(a)\cdot \vf(b), \\
a\stackrel+\le b\ &\hbox{iff }\vf(a)|\vf(b), \lb{3.70} \\
\vf(a\lor b) &= \lcm(\vf(a),\vf(b)), \lb{3.71} \\
\vf(a\land b) &= \gcd(\vf(a),\vf(b)). \lb{3.72}
\e
\eth

\blm3.35
Let $(X,\le,\lor,\land)$ and $(X',\le',\lor',\land')$ be \lt s. If $f:X\to Y$
is an \ois sm, then
\beq3.73
\vf(x\lor y)=\vf(x)\lor' \vf(y); \quad \vf(x\land y)=\vf(x)\land' \vf(y),
\e
\fa $x,y\in X$.
\elm

\bex3.36
Prove Lemma \rf{l3.35} and Theorem \rf{t3.34}.
\eex

We conclude this section by considering endo- and auto\mf s of
$(\Na,\cdot,1)$.

\bex9.26 \

\hph i,ii, Let $\phi$ be an endo\mf\ of $(\Na,\cdot,1)$ and let $n$ be
equal to $\prodl_{k=0}^m p_k^{\a_k}$, $m\in\N$, $\a_k\in\N$, $k\in[0,m]$.
Show that
\[
\phi(n)=\prod_{k=0}^{m} \phi(p_k)^{\a_k}.
\]

\hph ii,i, Find such a $\phi$ which is \ti{not\/} in\cre\ \wrt the \og~$\le$.

\hph iii,, Show that if $\phi_m(p_k)=(p_k)^m$, $m\in\Na$, $k\in\N$, then $\phi_m(u)
=u^m$, $u\in\Na$. Show that $\phi_m$ is in\cre\ \wrt the \og~$\le$.

\hph iv,, Are there other \hm sms which are strictly in\cre\ \wrt $\le$?

\hph v,i, Is the identity the only auto\mf?
\eex

%% file: DETOUR4.TEX
\Section{Prime fields}[Prime fields]\label{s.4}
\Subsubsection{Principal semirings}\label{sss.pr.smr}

The set $\N$ equipped with the \ad~$+$ and the \mlc~$\cdot$ is an abelian
monoid with \nel~$0$ \wrt the \ad, and is a \sg\ (see \E\df\ \rfa1{d2.2})
\wrt the \mlc. \Mo it \sf ies the \fw\ \cn s:
\bea1.1n
{}&0\cdot x=x\cdot 0=0 &&\kern-60pt\hbox{\fa} x\in X, \\
&x\cdot (y+z) = (x\cdot y)+(x\cdot z) &&\kern-60pt\hbox{\fa} x,y,z\in X, \lb{1.2n} \\
&(x+y)\cdot z=(x\cdot z)+(y\cdot z) &&\kern-60pt\hbox{\fa} x,y,z\in X. \lb{1.3n}
\e

\bdf1.1n
A set $X$ containing a distinguished \el~$e_0$ (also denoted by~$0$ when no
confusion arises) and equipped with two binary \op s, denoted by $\qu_0$ and
$\qu_1$ (also denoted by $+$ and~$\cdot$ when no confusion arises), will be
called (nonstandard terminology) a \ti{\sr}\index{semiring} if $(X,\qu_0,e_0)$ is an
\ti{abelian monoid\/}, $(X,\qu_1)$ is a \ti{\sg} and \er{1.1n}--\er{1.3n} hold.

A \sr\ $(X,\qu_0,\qu_1,e_0)$ is called \ti{\cmt e} if the \sg\ $(X,\qu_1)$ is\glossary{$(X,\qu_0,\qu_1,e_0)$}
abelian, and is called a \ti{\cmt e \sr\ with unity\/} (or \ti{identity\/})
if \te s ${e_1\in
X\sms{e_0}}$ \st $(X,\qu_1,e_1)$ is an abelian monoid. In this case we shall
use the notation $(X,{\qu_0},{\qu_1},e_0,e_1)$.
\edf

A more general \df\ of \sr\ can be found in \cite{Wiki}.

\bex1.2n
Show that in a \ti{\cmt e} \sr\ \cn s \er{1.2n} and \er{1.3n} are \ev t.
\eex

\bxs1.3n \

\hph i,i, $(\N,+,\cdot,0,1)$ is a \cmt e \sr\ with unity.\glossary{$(\N,+,\cdot,0,1)$}

\hph ii,, A \ti{\dsb e} \lt\ $(X,\lor,\land)$ with a \ti{least \el\/}~$e_0$
is a \cmt e \sr. \E\Ip the \dsb e \lt\ $(\Na,|)$ is a \cmt e \sr\ (without
unity).
\exs

\bex1.4n
Prove the assertions of Examples \rf{xa1.3n}.
\eex

\bcn1.5n
In this book we shall use the term \ti{\sr} for a \cmt e \sr\ with unity.
\ecn

The ``additive'' monoid $(\N,+,0)$ of the \sr~$\N$ is \pn\ with $1$~as \Gn,
since $n\nad{\era2{2.13}}= n\cdot 1 \nda22.8 = n\dpl1$ \fa
$n\in\N$.

\bdf1.6n
A \sr\ $(X,\qu_0,\qu_1,e_0,e_1)$ is called \ti{\pn\/}\index{principal semiring} (non-standard terminology)
if the monoid $(X,\qu_0,e_0)$ is \pn\ with \Gn~$e_1$.
\edf

\bdf1.7n
Let $(X,\qu_0,\qu_1,e_0,e_1)$ and $(\wt X,\tqu_0,\tqu_1,\wt e_0,\wt e_1)$ be
\sr s. A~map\break ${\vf:X\to \wt X}$ is called a \ti{\sr}- (or simply a~\ti{ring}-)\index{semiring-homomorphism}
\ti{\hm sm}, if the \fw\ holds:\index{ring-homomorphism}
\beq1.4n
\vf:(X,\qu_i,e_i)\to(\wt X,\tqu_i,\wt e_i), \ i=0,1,\ \hbox{are monoid-\hm
sms.}
\e
If, in \ad, $\vf$ is bi\jc, then $\vf$ is called a (\ti{semi})\ti{ring-\is sm}.
\edf

\bex1.8n
Show that if $\vf:X\to\wt X$ is a ring-\is sm, then so is $\vf\Inv$. Formulate
and prove an analogue of Lemma \rfa2{l1.8} where ``monoid'' is replaced by
``\sr''.
\eex

It turns out that all \ti{infinite} \pn\ \sr s are ring-\is c.

\bth1.9n
Let $(X,\qu_0,\qu_1,e_0,e_1)$ be an \emph{infinite} \pn\ \sr. Then

\hph i,ii, The \sq\ of \IT s of $e_1$ in $(X,\qu_0,e_0)$
\beq1.5n
\vf(n):=n \mathrel{\lower1pt\hbox{$\stackrel{\lower7pt\hbox{\smash{\hbox to0pt{$
\scriptstyle\sqcap$\hss}$\scriptstyle\sqcup$}$_0$}}\cdot$}}e_1, \q n\in\N,
\e
$($see \E\df s and Notation \rfa2{d1.14}$)$, is a ring-\is sm from $\N$ onto $X$.

\hph ii,i, The map $\vf$ is the only ring-\is sm from $\N$ onto $X$.

\hph iii,, The identity is the only ring-\is sm $($ring-auto\mf$)$ of~$X$.

\hph iv,, If $\wt X$ is an infinite \pn\ \sr\ and $\wt\vf:\N\to\wt X$ is
defined by \er{1.5n} where $\qu_0$~$($resp.~$e_1)$ is replaced by~$\tqu_0$
$($resp.~$\wt e_1)$, then $\wt\vf\circ \vf\Inv:X\to \wt X$ is a ring-\is sm,
the only one from~$X$ onto~$\wt X$.
\eth

\proof \

(i) In view of Theorem \rfa2{t1.42}, $(X,\qu_0,e_0)$ is a \pn\ infinite \PM\ with
\Gn~$e_1$, hence by Theorem \rfa2{t1.25}\,(i) with $\wt E:=\N$, $(M,\qu,e):=
(X,\qu_0,e_0)$ and $a:=e_1$, the map $\vf$ defined in~\er{1.5n} is a monoid-\is
sm from $(\N,+,0)$ onto $(X,\qu_0,e_0)$. Since $\vf$~is bi\jc, it suffices by
Lemma \rfa2{l1.8} to show that $\vf:(\N,\cdot,1) \to (X,\qu_1,e_1)$ is a \hm
sm. We have $\vf(1)\nde1.5n = 1\dquz e_1\nda22.3 = e_1$. \Mo if $m,n\in\N$,
then $\vf(m\cdot n)\nde1.5n = (m\cdot n)\dquz e_1\nda22.8 = (m\dpl n)\dquz
e_1$. Since the map $\vf:(\N,+,0)\to{(X,\qu_0,e_0)}$ is a monoid-\hm sm,
we obtain $(m\dpl n)\dquz
e_1\nad{\era2{2.3}\,\rm I3} = m\dquz (n\dquz e_1)=\break m\dquz ((n\dquz e_1)\qu_1e_1)$
by \E\df\ {\rf{d1.1n}}.
The map $\d:(X,\qu_0,e_0)\to (X,\qu_0,e_0)$ defined by $\d(x):=(n\dquz e_1)
\qu_1x$, $x\in X$, \sf ies $\d(e_0)=(n\dquz e_1)\qu_1e_0 \nde1.1n = e_0$ and
$\d(x\qu_0 y)=(n\dquz e_1)\qu_1(x\qu_0 y)\nde1.2n = ((n\dquz e_1)\qu_1 x)
\qu_0 ((n\dquz e_1)\qu_1 y)$, $x,y\in X$. \If that the map~$\d$ is an endo\mf\
of $(X,\qu_0,e_0)$. Hence by Lemma \rfa2{l1.22}, we obtain $m\dquz ((n\dquz
e_1)\qu_1e_1) = (n\dquz e_1)\qu_1(m\dquz e_1)$. Since $(X,\qu_1,e_1)$ is
abelian, $(n\dquz e_1)\qu_1(m\dquz e_1) = (m\dquz e_1)\qu_1(n\dquz e_1)
=\vf(m)\qu_1 \vf(n)$. Thus $\vf(m\cdot n)=\vf(m)\qu_1 \vf(n)$, and $\vf$~is
a~\hm sm from $(\N,\cdot,1)$ onto $(X,\qu_1,e_1)$.

(ii) follows from Theorem \rfa2{t1.25}\,(iii) with $\wt E:=\N$, $(M,\qu,e)
:=(X,\qu_0,e_0)$, $a:=e_1$.

(iii) If $\psi$ is an auto\mf\ of $X$, then $\psi\circ\vf:\N\to X$ is an \is
sm. Then $\psi\circ\vf=\vf$ by (ii). Hence $\psi=\psi\circ(\vf\circ\vf\Inv)=(\psi\circ\vf)
\circ\vf\Inv= \vf\circ\vf\Inv=\id_X$.

(iv) follows from Lemma \rfa2{l1.8} and (iii).
\endproof

\brm1.10n
In the proof of Theorem \rf{t1.9n} we used the fact that the selfmaps $\d_a$
of $(X,\qu_0,e_0)$, $a\in X$, defined by $\d_a(x):=a\qu_1x$, $x\in X$, are
endo\mf s. \If from Lemma \rfa2{l1.37}\,(i) and \era2{2.8}, \era2{2.12}
that every endo\mf\ of $(N,+,0)$
is a \mlc\ by an \el\ of~$\N$. Theorem~\rf{t1.9n} implies that the same holds
for all infinite \pn\ \sr s.
\erm

By \df\ a \sr\ contains at least two \el s. We claim that \te\ exactly two
\sr s consisting of $0$~and~$1$ where $0$ (resp.~$1$) is the \nel\ of the
``additive'' (resp.\ ``\mlv'') monoid. Set $X:=\{0,1\}$. We showed in
Examples \rfa2{xa1.5}\,(iii) that there are two binary \op s $\qu^{(1)}$
and~$\qu^{(2)}$ for which $(X,\qu^{(i)},0)$, $i=1,2$, are abelian monoids.

These binary \op s are defined by
\beq1.6n
0\qn10=0,\q 0\qn11=1,\q 1\qn10=1, \q 1\qn11=1,
\e
and
\beq1.7n
0\qn20=0,\q 0\qn21=1,\q 1\qn20=1, \q 1\qn21=0.
\e
If $(X,\qn i,\qu_1,0,1)$, $i=1,2$, are \sr s, then the binary \op\ $\qu_1$
\sf ies
\beq1.8n
0\qu_1 0=0,\q 0\qu_11=0, \q 1\qu_10=0, \q 1\qu_11=1.
\e
Indeed, the first three \et ies follow from \er{1.1n} and the last \et y (in
fact the last three \et ies) follows from the fact that $1$~is a \nel. In
view of Example \rfa2{xa1.5}\,(iii) with $\a:=1$, $\b:=0$, ${\qu_1}:={\qn1}$, we
find that $(X,\qu_1,1)$ is an abelian monoid. We now show that
$(X,\qu_0,\qu_1,0,1)$ is a~\sr\ whenever ${\qu_0}:={\qn1}$ or~$\qn2$ and $\qu_1$
is defined by \er{1.8n}. It suffices to prove \er{1.1n} and \er{1.2n} (by
Exercise \rf{ex1.2n}). \E\cn\ \er{1.1n} follows from \er{1.8n}.

\er{1.2n}: If $x:=0$ in \er{1.2n}, then $0\qu_1(y\qu_0z)\nde1.1n =0
= 0\qu_00\nde1.1n = (0\qu_1x)\qu_0(0\qu_1y)$ \fa $y,z\in X$, since $0$~is the
\nel\ of~$\qu_0$.

If $x:=1$ in \er{1.2n}, then $1\qu_1(y\qu_0z) = y\qu_0z = (1\qu_1y)\qu_0(1
\qu_1z)$ \fa $y,z\in X$, since $1$~is the \nel\ of~$\qu_1$.

We summarize in

\bpr1.11n
Let $X:=\{0,1\}$. Then $(X,\qu_0,\qu_1,0,1)$ is a \sr\ with~$0$ $($resp.~$1)$
iff\/ $\qu_1$ \sf ies \er{1.8n} and $\qu_0$ is equal either to $\qn1$ defined in
\er{1.6n} or to $\qn2$ defined in~\er{1.7n}.
\epr

If $(X,\qu_0,\qu_1,e_0,e_1)$ is a \sr, if $\O$ is a \ns\ and $X^\O$ denotes
the set of all maps from~$\O$ into~$X$, then we can define binary \op s
$\tqu_i$, $i=0,1$, on~$X^\O$ and \el s $\wt e_0,\wt e_1$ of~$X^\O$ \st
$(X^\O, \tqu_0,\tqu_1,\wt e_0,\wt e_1)$ is a~\sr.

\bpr1.12n
Let $\O$ be a \ns.

\hph i,i, Let $(X,\qu,e)$ be a \emph{monoid}. Given $\vf,\psi\in  X^\O$
define ${\tqu}:X^\O\t X^\O\to X^\O$ by setting
\beq1.9n
(\vf \tqu \psi)(\o):=\vf(\o)\qu \psi(\o), \q \o\in\O.
\e
Define $\wt e \in X^\O$ by setting
\beq1.10n
\wt e(\o):=e, \q \o\in\O.
\e
Then $(X^\O,\tqu,\wt e)$ is a monoid. \Mo if $(X,\qu,e)$ is an abelian
$($resp.\ C-, P-$)$ monoid, then so is $(X^\O,\tqu,\wt e)$.

\hph ii,, Let $(X,\qu_0,\qu_1,e_0,e_1)$ be a \emph{\sr}. Given $\vf,\psi\in
X^\O$ define ${\tqu_i}:X^\O\t X^\O\to X^\O$, $i=0,1$, by setting
\beq1.11n
(\vf\tqu_i \psi)(\o):= \vf(\o)\qu_i \psi(\o), \q \o\in\O,\ i=0,1.
\e
Define $\wt e_i\in X^\O$, $i=0,1$, by setting
\beq1.12n
\wt e_i(\o):=e_i, \q \o\in\O,\ i=0,1.
\e
Then $(X^\O,\tqu_0,\tqu_1,\wt e_0,\wt e_1)$ is a \sr.
\epr

\proof
(i) The proof is similar to the proof of Example \rfa2{xa1.5}\,(iv) with
$(X_1,\qu_1,e_1) = (X_2,{\qu_2},e_2)$, i.e., when $\O:=\{1,2\}$, \Tf it is
omitted.

(ii) $(X^\O,\tqu_i,\wt e_i)$, $i=0,1$, are abelian monoids by~(i). \Mo given
$\vf\in X^\O$, we have $(\wt e_0\tqu_1\vf)(\o)\nde1.11n = \wt e_0(\o)\qu_1
\vf(\o) \nde1.12n = e_0\qu_1\vf(\o) \nde1.1n = e_0$, $\o\in\O$. Hence $\wt
e_0(\o) \tqu_1 \vf=\wt e_0$. Since $\wt e_1\tqu_1\vf = \vf\tqu_1\wt e_1$ by~(i),
\er{1.1n} holds. We now prove \er{1.2n}. Given $\xi,\vf,\psi\in X^\O$, we
have $\bigl(\xi\tqu_1 (\vf\tqu_0\psi)\bigr)(\o) \nde1.11n = \xi(\o) \qu_1
(\vf\tqu_0\psi)(\o) \nde1.11n = \xi(\o)\qu_1(\vf(\o)\qu_0\psi(\o)) \nde1.2n =
(\xi(\o)\qu_1\vf(\o)) \qu_0 (\xi(\o)\qu_1\psi(\o)) = \bigl((\xi\tqu_1\vf)
(\o)\bigr)\qu_0 \bigl((\xi\tqu_1\psi)(\o)\bigr) = \bigl((\xi\tqu_1\vf)\tqu_0
(\xi\tqu_1\psi)\bigr)(\o)$, $\o\in\O$. Hence $\xi\tqu_1(\vf\tqu_0\psi)  =
(\xi\tqu_1\vf)\tqu_0(\xi\tqu_1\psi)$, which proves \er{1.2n}.
\endproof

\bnt1.13n
It is customary to write $a\qu_1 b\qu_0 c$ (resp.\ $a\qu_0 b\qu_1 c$) instead
of $(a\qu_1b)\qu_0c$ (resp.\ $a\qu_0(b\qu_1c)$). The \op\ $\qu_1$ is
performed first.
\ent

\bex1.14n \

\hph i,ii, Let $X:=\{0,1\}$. As in Exercise \rfa3{ex1.3} we denote by~$\lor$
the binary \op\ defined by \era3{8.2} and by~$\land$
the binary \op\ defined by \era3{8.3}. Observe that $\lor$~is equal to $\qn1$
defined in \er{1.6n} and that $\land$ is equal to $\qu_1$ defined in
\er{1.8n}. If $X$~is equipped with the \rt ion of the \nog\ of~$\N$ to
$\{0,1\}$, then $x\lor y\nda23.14 = \max(x,y)= \sup\{x,y\}$, $x\land
y\nda23.15 = \min(x,y)=\inf\{x,y\}$, $x,y\in X$. Show that $(X,\le)$ is a
\dsb e \lt\ with least \el~$0$ and greatest \el~$1$. Define $0^c=1$ and
$1^c=0$. Show that $X$~is a Boolean algebra (see \E\df\ \rfa3{d8.40}). By
\E\Pr\ \rf{p1.11n}, $(X,\lor,\land,0,1)$ is also a \sr. Using the symbols
$\wt\lor, \wt\land,\wt0,\wt1$ instead of $\tqu_0,\tqu_1,\wt e_0,\wt e_1$ in
\E\Pr\ \rf{p1.12n}, we find by the same \Pr\ that $(X^\O,\wt\lor,\wt\land,
\wt0, \wt1)$ is a \sr. \Mo if $f\in X^\O$ and $\supp(f):=\{\o\in\O:
f(\o)=1\}$, then $f$~is the \ti{indicator \f} $1_{\supp(f)}:\O\to X$ (see
Section~\ref{sss.latt})  where $1_A(\o):=1$ for $\o\in A$, $1_A(\o):=0$ for
$\o\notin A$, $A$~subset of~$\O$. \E\Ip
\beq1.13n
\wt0 = 1_\vn \qh{and }\wt1 = 1_\O.
\e
Show that if $A,B$ are subsets of~$\O$, then
\bea1.14n
1_A \mathrel{\wt\lor} 1_B &= 1_{A\cup B}, \\
1_A \mathrel{\wt\land} 1_B &= 1_{A\cap B}. \lb{1.15n}
\e

Given $f\in X^\O$, set $f^c(\o):=(f(\o))^c$, $\o\in\O$. Show that $(X^\O,
\wt\lor,\wt\land,\wt0,\wt1)$ equipped with the \op\ $\wt f\mt \wt f{}^c$ is
a~Boolean algebra.

\hph ii,i, Show that if in part (i) we replace $\lor$ by $\qn2$ defined in
\er{1.7n}, then
\beq1.16n
1_A \mathrel{\wt{\qn2}} 1_B= 1_{A\dot\cup B},\qh{ $A,B$ subsets of~$\O$},
\e
where $A\mathrel{\dot\cup} B$ is defined by
\beq1.17n
A\mathrel{\dot\cup} B:=(A\cup B)\sm(A\cap B).
\e
\E\Ip $A\mathrel{\dot\cup}A=\vn$ and $A\mathrel{\dot\cup}A^c=\O$.

\hph iii,, Show that $(X^\O,\wt\lor,\wt 0)$ is \ti{not\,} a \Cm, whilst
$(X^\O,\wt{\qn2},\wt 0)$ is a \Cm.

\hph iv,, Show that if $\O$ contains at least two \el s, then both
$(X^\O,\wt\lor,\wt 0)$ and $(X^\O,\wt{\qn2},\wt 0)$ are not \pn\ monoids.
\eex

We now construct finite \pn\ \sr s with \ti{\cnc e} ``additive'' monoids of
arbitrary \ca~$n$ larger than one. The \sr\ $(\{0,1\},\qn2,\qu_1,0,1)$ with
$\qn2$ defined in \er{1.7n} and $\qu_1$ defined in \er{1.8n} is such an
example  ($n=2$).

Let $n\in\Na\sms1$ and set
\beq1.1
\N_n:=\{k\in\N: 0\le k<n\}.
\e
By Theorem \rfa2{t1.38} and \era2{2.8}, given $x\in\N$ \te
s one and only one pair $(q,r)\in\N\t\N_n$ \st
\beq1.2
x=q\cdot n+r.
\e
Let $\phi_n:\N\to\N_n$ be defined by
\beq1.3
\phi_n(x):=r.
\e
Observe that $\phi_{n}(x)$ is the last digit (the only one\glossary{$\phi_n$}
if $x\in[0,n)$) of the \rp ation of~$x$ with base~$n$. For example,
$\phi_{10}(12)=2$, $\phi_{10}(7)=7$.

The map $\phi_n$ \sf ies
\beq1.4
\phi_n(x+mn)=\phi_n(x) \qh{for all} x,m\in\N.
\e
Indeed, we have $x+mn\nde{1.2} = (qn+r)+mn=mn+(qn+r)=(mn+qn)+r=(m+q)n+r$. Hence
$\phi_n(x+mn) =r=\phi_n(x)$.

\E\Ip if $m=0$:
\bea1.5
\phi_n(x)=x &\q\hbox{\fe} x\in\N_n,\\
\phi_n(0)=0 &\q\hbox{and }\phi_n(1)=1. \lb{1.6}
\e

Observe that if $x,y\in\N_n$ with $x+y\in\N_n$ then $\phi_n(x+y)=x+y=\phi_n(x)
+\phi_n(y)$. If $x+y\notin\N_n$, of course \et y cannot hold, however, the
\fw\ holds:
\beq1.7
\phi_n(x+y)=\phi_n(\phi_n(x)+\phi_n(y)) \qh{for all} x,y\in\N_n.
\e
Indeed, from the \ex\ part of Theorem \rfa2{t1.38} we have $x=q\cdot
n+\phi_n(x)$, $y=\wt q\cdot n+\phi_n(y)$ for some $q,\wt q\in\N$. Hence
$x+y=(q\cdot n+\phi_n(x))+(\wt q\cdot n+\phi_n(x))=
q\cdot n+(\phi_n(x)+\wt q\cdot n)+\phi_n(y)
=q\cdot n+(\wt q\cdot n+\phi_n(x))+\phi_n(y)
=(q\cdot n+\wt q\cdot n)+(\phi_n(x)+\phi_n(y))
= (q+\wt q)\cdot n+(\phi_n(x)+\phi_n(y))$. Then \er{1.7} follows
from \er{1.3}.

Identity \er{1.7} says that the last digit (\wrt a base~$n$) of the sum $x+y$
is the last digit of the sum of the last digit of~$x$ and the last digit
of~$y$. This motivates the \df\ of a binary \op\ on~$\N_n$ which we shall
denote by~$+_n$.\glossary{$+_n$}

Set
\beq1.8
x+_n y:=\phi_n(x+y) \qh{for all} x,y\in\N_n.
\e
We shall call $x+_ny$ the \ti{$n$-sum} of $x$ and $y$. It turns out that
$(\N_n,+_n,0)$ is a \pn\ C-monoid and that the map $\phi_n:(\N,+,0) \to
(\N_n,+_n,0)$ is a sur\jc\ \hm sm.

\bpr1.8
Let $n\in\Na\sms1$ and let $\N_n:=[0,n)$. Given $x\in\N$, let $\phi_n(x)$
denote the unique \el\ of~$\N_n$ \sf ying
\beq1.9
x=q\cdot n+\phi_n(x) \qh{for some} q\in\N.
\e
Let $+_n$ be the $n$-sum of two \el s of~$\N_n$ defined in \er{1.8}. Then

\hph i,ii, $(\N_n,+_n,0)$ is a \pn\ C-monoid with $1$ as a \Gn.

\hph ii,i, \E\fe $a\in\N_n$ \te s one and only one \el\ $a\Inv$ in $\N_n$, called\index{inverse}
the inverse of~$a$, \st
\bga1.10
a\Inv +_n a=0,\\
a\Inv=0 \hbox{ if }a=0 \qh{and } a\Inv=n-a \hbox{ if }a\in(0,n), \lb{1.28}
\e
where $a,n$ are viewed as \el s of $(\N,+,0)$.

\hph iii,, The map $\phi_n:(\N,+,0)\to (\N_n,+_n,0)$ is a sur\jc\ monoid-\hm sm.

\hph iv,, If $\vf_1:\N\to \N_n$ denotes the \sq\ of \IT s of~$1$ in
$(\N_n,+_n,0)$, then
\beq1.11
\vf_1(m)=\phi_n(m) \qh{\fe} m\in\N.
\e
\epr

\proof
$(\N_n,+_n,0)$ is a C-monoid.

\ti{\E\asc ity}: Let $a,b,c\in\N_n$.
\bmlg
(a+_n b)+_n c\nde1.8 = \phi_n((a+_n b) + c)\nde1.8 =
\phi_n(\phi_n(a+b)+c)\\
{}\nde1.5 = \phi_n(\phi_n(a+b) +\phi_n(c))
\nde1.7 = \phi_n((a+b)+c) \nda2{1.5} = \phi_n(a+(b+c))\\
{} \nde 1.7 = \phi_n(\phi_n(a)+\phi_n(b+c))
{}\nde1.5 = \phi_n(a+\phi_n(b+c))
\nde 1.8 =\phi_n(a+(b+_nc)) \nde1.8 = a+_n(b+_nc).
\e

\ti{\E\cmt ity}: $a+_nb\nde1.8 = \phi_n(a+b) \nda21.6 = \phi_n(b+a)\nde1.8 =
b+_na$.

\ti{\E\nel\/}: $0+_n a\nde1.8 = \phi_n(0+a)\nda21.7 = \phi_n(a)\nde1.5 = a$.

\ti{\E\cnc ity}: Suppose $a+_nb=a+_nc$. Then $\phi_n(a+b) =
\phi_n(a+c)=:r \in\N_n$. \E\te\ $q,\wt q\in\N$ \st
\[
\bal a+b&= q\cdot n+r,\\
a+c&= \wt q\cdot n+r.
\eal
\]
Suppose $q\le \wt q$. Then \te s $p\in\N$ \st $\wt q=p+q$. So we obtain
$a+ b=q\cdot n+r$ and $a+c=(p+q)\cdot n+r \nda22.15 = p\cdot n+q\cdot n+r$.
Hence $a+c=a+b+p\cdot n$. By the \cnc ity of~$+$, we obtain $c=b+p\cdot n$.
Since $c\in\N_n$, we have $b+p\cdot n<n$. If $p\ge1$, then $p\cdot n\ge
1\cdot n=n$ and $n\le p\cdot n\le b+p\cdot n<n$, a~\cd ion. \E\Tf $p=0$. It
follows that $a+c=a+b$, hence $c=b$. The case $\wt q\le q$ is similar.

\ti{$\phi_n:(\N,+,0)\to(\N_n,+_n,0)$ is a sur\jc\ \hm sm}: $\phi_n(0)\nde1.6
=0$ and\break $\phi_n(x+y)\nde 1.7 = \phi_n(\phi_n(x)+\phi_n(y))\nde1.8 = \phi_n
(x)+_n \phi_n(y)$, $x,y\in\N_n$. The surjectivity of~$\phi_n$ follows
from~\er{1.5}.

\ti{$(\N_n,+_n,0)$ is \pn\/}: Let $\vf_1$ denote the \sq\ of \IT s
of~1 in $(\N_n,+_n,0)$ (see \E\Pr\ \rfa2{p1.13}). We prove by induction on~$m$
that \er{1.11} holds. Indeed, $\vf_1(0)=0\nda21.17 = \F_n(0)$. Suppose $\vf_1(m)
=\phi_n(m)$ for some $m\in\N$. Then $\vf_1(m+1)\nda21.18 = 1+_n\vf_1(m)=
1+_n\phi_n(m) \nde1.8 = \phi_n(1+\phi_n(m))\nde1.6 = \phi_n(\phi_n(1)+
\phi_n(m)) \nde1.7 = \phi_n(1+m) \nda21.6 = \phi_n(m+1)$. It follows that
\er{1.11} holds.

By \df\ of $\vf_1$, $\{\vf_1(m):m\in\N\}\subset\N_n$. Now let $m\in\N_n$.
Then $m\nde1.5 = \phi_n(m)\nde1.11 = \vf_1(m)\in\{\vf_1(m):m\in\N\}$. Hence
$\N_n$ is equal to the range of~$\vf_1$.

\ti{\er{1.10}{\rm:} \uq}. Let $a\in\N_n$. Suppose \te\ $b,c\in\N_n$ \st $b+_na=0$
and $c+_na=0$. Then $b+_na=c+_na$ and by \cmt ity and \cnc ity we obtain
$b=c$.

\ti{\E\ex}: If $a=0$, then clearly $a\Inv=0$ since $0+_n0=0$. If $a\in[1,n)$,
then by \era2{1.46} \te s $p\in\N$ \st $n=p+a$. Notice that $p\in\N_n$. So
$\phi_n(n)=\phi_n(p+a)=p+_na$ and $\phi_n(n)=0$ by \er{1.4} where $x:=0$,
$m:=1$, and by~\er{1.6}. \E\Tf $p$ is~$a\Inv$.
\endproof


\bex1.17
Let $n\in\Na\sms1$, $a,b\in\N_n$. Show:
\bea1.30
{}&(a\Inv)\Inv=a, \q 0\Inv=0,\\
&(a+_nb)\Inv = a\Inv+_nb\Inv. \lb{1.31}
\e
\eex

We want to
show that \te s a binary \op\ $\cdot_n$ on~$\N_n$ which makes\break $(\N_n,+_n,
\cdot_n,0,1)$ a \sr. Suppose that it is the case, then the map $\d_a:\N_n \to\N_n$,\glossary{$\cdot_n$}
$a\in\N_n$, defined by $\d_a(x):=a\cdot_n x$, $x\in \N_n$, is an endo\mf\
of $(\N_n,+_n,0)$, which \sf ies $\d_a(1)=a$ since $(\N_n,\cdot_n,1)$ is an
abelian monoid. Now, if $m\in\N_n$, then $m=\F_n(m)$ by \er{1.5}, hence
$m$~is the $m$-th \IT\ of~1 in $(\N_n,+_n,0)$ in view of \er{1.11}. Using
the notation introduced in \era2{1.25}, we obtain $m=m\dpln1$. In view of
Lemma~\rfa2{l1.22} we have
\beq3.22
\d_a(m\dpln1)=m\dpln\d_a(1).
\e
Since $\d_a(1)=a$, $\d_a(m)=\d_a(m\dpln1)\nde3.22 =
m\dpln\d_a(1)=m\dpln a$. Since $a=a\dpln1$, $m\dpln a=
m\dpln (a\dpln1)\nda22.3 = (m\dpl a)\dpln1$ where $m\dpl a$ is the $m$-th
\IT\ of~$a$ in $(\N,+,0)$. Notice that \era2{2.3} holds since $(\N_n,+_n,0)$
is abelian. From \era2{2.8} it follows that $m\dpl a=m\cdot a$ where $\cdot$
is the \mlc\ in~$\N$. \E\Tf $a\cdot_n m=\d_a(m)=(m\cdot a)\dpln1
\nde1.11 = \F_n(m\cdot a)\nda22.12 = \F_n(a\cdot m)$. We have thus proved that
if the \op\ ${\cdot_n}:\N_n\t\N_n \to \N_n$ makes $(\N_n,+_n,\cdot_n,0,1)$
a~\sr, then
\beq1.33
x\cdot_n y=\F_n(x\cdot y), \q x,y\in\N_n,
\e
where $\cdot$ denotes the \mlc\ in~$\N$.

\bex3.14
Let $n\in\Na\sms1$ and let $\vf:\N_n\to\N_n$. Show that $\vf$ is an endo\mf\
of $(\N_n,+_n,0)$ iff it \sf ies
\beq3.23
\vf(m)=m\dpln \vf(1) \qh{\fe}m\in\N_n.
\e
\eex

\bpr3.15
Let $n\in\Na\sms1$, let $\F_n:\N\to\N_n$ $(:=[0,n))$ be defined as in
\er{1.3}, let $+_n$ be as in \er{1.8}, and let $\cdot_n$ be as in \er{1.33}.
Then

\hph i,i, $(\N_n,+_n,\cdot_n,0,1)$ is a \ti{\pn\ \sr} and $(\N_n,+_n,0)$ is a \Cm,

\hph ii,, $\F_n:(\N,+,\cdot,0,1) \to (\N_n,+_n,\cdot_n,0,1)$ is a \ti{sur\jc\
\sr-\hm sm}.
\epr

\proof \

1) $(\N_n,+_n,0)$: $(\N_n,+_n,0)$
is a finite \pn\ \Cm\ by \E\Pr\ \rf{p1.8}\,(i). \Mo $\F_n:
(\N,+,0)\to(\N_n,+_n,0)$ is a sur\jc\ monoid-\hm sm by \E\Pr\ \rf{p1.8}\,(iii).

2) $(\N_n,\cdot_n,1)$: We first prove an analogue of \er{1.7} where $+$ is
replaced by~$\cdot$. The \fw\ holds:
\beq3.26
\F_n(x\cdot y)=\F_n(\F_n(x)\cdot \F_n(y)) \qh{\fa} x,y\in\N_n.
\e
Indeed, from the \ex\ part of Theorem \rfa2{t1.38} we have $x=q\cdot n+
\F_n(x)$, $y=\wt q\cdot n+\break \F_n(y)$ for some $q,\wt q\in\N$. Hence $x\cdot y=
(q\cdot n+\F_n(x))\cdot (\wt q\cdot n+\F_n(y)) = [q\cdot \wt q\cdot n +\break
q\cdot\F_n(y)+\wt q\cdot \F_n(x)]\cdot n + \F_n(x)\cdot \F_n(y)$. Using
part~1) of the proof we obtain $\F_n(x\cdot y)=\F_n\bigl(
[q\cdot \wt q\cdot n + q\cdot\F_n(y)+\wt q\cdot \F_n(x)]\cdot n\bigr)
+_n \F_n(\F_n(x)\cdot \F_n(y)) = \F_n(r\cdot n)+_n \F_n(\F_n(x)\cdot
\F_n(y))$ \fs $r\in\N$. We have $\F_n(r\cdot n)
=0$ by \er{1.4}, \er{1.7} with $x:=0$. \E\Tf \er{3.26} holds.

Since $\F_n$ maps $\N$ into $\N_n$, $x\cdot_n y\in\N_n$ \fa $x,y\in\N_n$.

\ti{\E\asc ity of $\cdot_n$}: $x,y,z\in\N_n$.
\bgg
(x\cdot_ny)\cdot_nz \nde1.33 = \F_n((x\cdot_ny)\cdot_nz) \nde1.33 =
\F_n((\F_n(x\cdot y))\cdot z) \nde1.5 = \F_n(\F_n(x\cdot y)\cdot \F_n(z))\\
\nde3.26 = \F_n((x\cdot y)\cdot z)
\nda22.11 = \F_n(x\cdot (y\cdot_n z)) \nde3.26 = \F_n(\F_n(x)\cdot
\F_n(y\cdot z))\\
\nde1.5 = \F_n(x\cdot \F_n(y\cdot z))
\nde1.33 = \F_n(x\cdot (y\cdot_n z)) \nde1.33 = x\cdot_n(y\cdot_n z).
\e

\ti{\E\cmt ity of $\cdot_n$}: $x,y\in\N_n$. $x\cdot_n y\nde1.33 =
\F_n(x\cdot y)\nda22.12 = \F_n(y\cdot x)\nde1.33 = y\cdot_n x$.

\ti{$1$ is the \nel\/}: $x\in\N_n$. $1\cdot_nx\nde1.33 = \F_n(1\cdot x)
\nda22.13 = \F_n(x)\nde 1.5 = x$.

It follows that $(\N_n,\cdot_n,1)$ is an abelian monoid. \Mo $\F_n:
(\N,\cdot,1)\to (\N_n,\cdot_n,1)$ \sf ies $\F_n(1)\nde1.5 = 1$ and
$\F_n(x\cdot y)\nde3.26 = \F_n(\F_n(x)\cdot \F_n(y)) \nde1.33 = \F_n(x)\cdot_n
\F_n(y)$, $x,y\in\N_n$. Hence $\F_n$ is a monoid-\hm sm from $(\N,\cdot,1)$
onto $(\N_n,\cdot_n,1)$. It is sur\jc\ by part~1) of the proof.

3) \ti{Proof of \er{1.1n}}: $x\in\N_n$. $0\cdot_nx \nde1.33 = \F_n(0\cdot x)
\nad{\era2{2.12},\era2{2.9}} = \F_n(0)\nde1.5 = 0$.

\er{1.2n}: $x,y,z\in\N_n$.
\bgg
x\cdot_n y+_n x\cdot_nz \nde1.33 = \F_n(x\cdot y)+_n \F_n(x\cdot z)
\nde1.7 = \F_n(\F_n(x\cdot y)+\F_n(x\cdot z)) \\ \nde1.7 = \F_n(x\cdot y
+x\cdot z)
\nda22.14 = \F_n(x\cdot(y+z)) \nde3.26 = \F_n(\F_n(x)\cdot \F_n(y+z))\\
\nde1.8 = \F_n(\F_n(x)\cdot (y+_nz)) \nde1.5 =
\F_n(x\cdot(y+_nz)) \nde1.33 = x\cdot_n(y+_nz).
\e

\er{1.3n}: $x,y,z\in\N_n$. $(x+_ny)\cdot_nz \nad*= z\cdot_n(x+_ny)\nde1.2n =
z\cdot_n x+_n z\cdot_n y\nad*= x\cdot_nz+_n y\cdot_nz$, where in $*$ the
\cmt ity of $\cdot_n$ is used.
\endproof

\bxs1.19 \

\hph i,i, The \fw\ tables are ``\ad\ tables'' of $\N_4$ and $\N_5$:
\[
\begin{array}{c|cccc}
\ +_4\ &\ 0\ &\ 1\ &\ 2\ &\ 3\ \\ \hline 0&0&1&2&3 \\ 1&1&2&3&0 \\ 2&2&3&0&1 \\ 3&3&0&1&2
\end{array}
\qquad
\begin{array}{c|ccccc}
\ +_5\ &\ 0\ &\ 1\ &\ 2\ &\ 3\ &\ 4\ \\ \hline 0&0&1&2&3&4 \\ 1&1&2&3&4&0 \\ 2&2&3&4&0&1 \\ 3&3&4&0&1&2 \\ 4&4&0&1&2&3
\end{array}
\]

\hph ii,, The \fw\ tables are ``\mlc\ tables'' of $\N_4$ and $\N_5$:
\[
\begin{array}{c|cccc}
\ \cdot_4\ &\ 0\ &\ 1\ &\ 2\ &\ 3\ \\ \hline 0&0&0&0&0 \\
1&0&1&2&3 \\ 2&0&2&0&2 \\ 3&0&3&2&1
\end{array}
\qquad
\begin{array}{c|ccccc}
\ \cdot_5\ &\ 0\ &\ 1\ &\ 2\ &\ 3\ &\ 4\ \\ \hline 0&0&0&0&0&0 \\
1&0&1&2&3&4 \\ 2&0&2&4&1&3 \\ 3&0&3&1&4&2 \\ 4&0&4&3&2&1
\end{array}
\]
\exs

\bex1.20 \

\hph i,i, Prove the assertions of Examples \rf{xa1.19}.

\hph ii,, Find the tables of $+_n$, $\cdot_n$ for $n=6,7,10$.
\eex

In what follows we will only consider \sr s with \cnc e ``additive'' monoid.
For convenience we introduce the \fw\ (nonstandard) \df.

\bdf1.21
A \ti{\cnc e} \sr\ $(X,\qu_0,\qu_1,e_0,e_1)$ is a set $X$ containing two
distinguished \el s $e_0,e_1$, equipped with two binary \op s $\qu_0,\qu_1$,
and \sf ying the \fw\ \cn s:
\bea1.36
{}&(X,\qu_0,e_0) \hbox{ is a \Cm.} \\
&(X,\qu_1,e_1) \hbox{ is an \am.} \lb{1.37} \\
&\hbox{If }\d_a:X\to X,\ a\in X, \hbox{ denotes the map defined by } \non\\
&\hskip30pt \d_a(x):=a\qu_1 x \ ({}=x\qu_1 a),\ x\in X, \lb{1.38} \\
&\hbox{then $\d_a$ is an endo\mf\ of the monoid $(X,\qu_0,e_0)$ \fe} a\in X. \non
\e
If, in \ad, the monoid $(X,\qu_0,e_0)$ is \pn, then $X$ is called a \ti{\pn\
\cnc e \sr}.\index{principal c-semiring}

The monoid $(X,\qu_0,e_0)$ (resp.\ $(X,\qu_1,e_1)$) is called the
\ti{additive} (resp.\ \ti{\mlv}) monoid of the \sr~$X$.
\edf

We proved in Theorem \rf{1.9n} that an \ti{infinite} \pn\ \cnc e \sr\ is \is
c to $(\N,+,\cdot,0,1)$. It turns out that a \ti{finite} \pn\ \cnc e \sr\ is
\is c to $(\N_n,+_n,\cdot_n,0,1)$ where $n:=\#(X)$. We first show that
a~finite \pn\ \Cm\ is \is c to $(\N_n,+_n,0)$ where $n:=\#(X)$. We recall
that if  $(X,\qu,e)$ is a nontrivial monoid and if $a\in X$, then $\vf_a:\N
\to X$, the \sq\ of \IT s of~$a$ in $(X,\qu,e)$, also denoted by $\vf_a:=
x\dqu a$, $x\in \N$, \sf ies the \fw\ \pp ies:
\bea1.40
{}&\vf_a(0)=e, \q \vf_a(1)=a, \\
&\vf_a(m+n)=\vf_a(m)\qu \vf_a(n), \q m,n\in\N, \lb{1.41} \\
&\vf_e(m)=e, \q m\in\N. \lb{1.42}
\e
Indeed, \er{1.40} is \era2{2.3}\,I0,I1, \er{1.41} is \era2{2.3}\,I2 and
\er{1.42} is \era2{2.3}\,I4 with $E:=\N$ and $M:=X$.

We also recall that $I(a)$ is the range of $\vf_a$ in~$X$, \Ip $I(e)=\{e\}$
by~\er{1.42}. If $a\ne e$, then $I(a)$ is a~\sbm\ of~$X$, $(I(a),\qu,e)$ is a
\pn\ \am\ with \Gn~$a$ by Lemma \rfa2{l1.21}\,(ii) and $\vf_k:\N\to
(I(a),\qu,e)$ is a sur\jc\ \hm sm by Lemma \rfa2{l1.21}\,(iii). In the next lemma, we
show among other things that if $(X,\qu,e)$ is finite and is a \Cm, then \te
s \ooo $m\in\Na\sms1$ \st
\beq1.43
\vf_a(m)=e \qh{and } \vf_a(k)\ne e \qh{\fa} k\in(0,m).
\e
\Mo the \rt ion of $\vf_a$ to $[0,m)$, which we denote by $\ov\vf_a$, is an
\is sm from $(N_m,+_m,0)$ onto $(I(a),\qu,e)$. As a con\sq, if $(X,\qu,e)$ is
a finite \pn\break \Cm\ with \Gn~$a$, then by Lemma \rfa2{l1.8} $(X,\qu,e)$ is \is
c to $(\N_n,+_n,0)$ with $n:=\#(X)$.

\blm1.22
Let $(X,\qu,e)$ be a finite \Cm\ with $n:=\#(X)\ge2$, let $a\in X\sms e$ and
let $\vf_a:\N\to X$ be the \sq\ of \IT s of~$a$ in $(X,\qu,e)$. Then the \fw\
holds{\rm:}

\hph i,ii, \E\te s \ooo $m\in\Na\sms1$ \st \er{1.43} holds.

\hph ii,i, If $k,l\in[0,m)$ and  $\vf_a(k)=\vf_a(l)$, then $k=l$.

\hph iii,, $\vf_a(qm)=e$ \fa $q\in\N$.

\hph iv,, $\vf_a(k)=\vf_a(\F_m(k))$ \fa $k\in \N$.

\hph v,i, $\vf_a(k)=e$ with $k\in\N$ implies $k=qm$ \fs $q\in\N$.

\hph vi,, If $\ov\vf_a$ is the \rt ion of $\vf_a$ to $[0,m)$, then $\ov\vf_a:
(\N_m,+_m,0) \to (I(a),\qu,e)$ is an \is sm.
\elm

\proof \

(i) Since $\N$ is infinite and $X$ is finite, $\vf_a$ is \ti{not\/} in\jc\ by
\E\Pr\ \rfa1{p4.27}\,(A)(i). Thus \te\ $k,l\in\N$, $k<l$, \st $\vf_a(k)=
\vf_a(l)$ (the \og\ of~$\N$ being total). Since $k<l$, $l=p+k$ \fs $p\in\Na$.
Then $e\qu\vf_a(k)=\vf_a(k)=\vf_a(l)\nde1.41 = \vf_a(p)\qu\vf(k)$. By \cnc
ity \era2{1.8}, $\vf_a(p)=e$. \csq, $A:=\{l\in\Na: \vf_a(l)=e\}$ is not
empty. Since $(N,\le)$ is well-ordered, the nonempty subset~$A$ of~$\N$
possesses a least \el\ $m\in\Na$. Observe that $m\ne1$, since $\vf_a(1)
\nde1.40 = a\ne e$.

(ii) Since the \og\ of $\N$ is total, interchanging $k$ and~$l$ if necessary,
we may suppose $k\le l$. Let $p\in\N$ be \st $l=p+k$. Then $p\in[0,m)$ by
\era1{3.12} since $p\le l$ and $l<m$. \Mo $e\qu\vf_a(k)=\vf_a(k)=\vf_a(l)
\nde1.41 = \vf_a(p)\qu\vf_a(k)$. By \cnc ity again, $\vf_a(p)=e$. Since
$p\in[0,m)$ and $\vf_a(p)=e$, we infer $p=0$ from~(i). Hence $k=l$.

(iii) $\vf_a(qm)=(qm)\dqu a=(q\dpl m)\dqu a \nad{\era2{2.3}\,\rm I3} =
q\dqu(m\dqu a) = q\dqu\vf_a(m) \nad{\rm(i)}= q\dqu e=\vf_e(q) \nde1.42 = e$, \fa
$q\in\N$.

(iv) Let $k\in\N$. Then $k=qm+\F_m(k)$ \fs $q\in\N$ by \er{1.2}, \er{1.3}
since $m\in\Na$. Then $\vf_a(k)\nde1.41 = \vf_a(qm)\qu\vf_a(\F_m(k))
\nad{\rm(iii)} = e\qu\vf_a(\F_m(k))=\vf_a(\F_m(k))$.

(v) Let $k\in\N$ be \st $\vf_a(k)=e$. Then $e\nad{\rm(iv)}= \vf_a(\F_m(k))$.
Since $\F_m(k)\in[0,m)$ we have $\F_m(k)=0$ by~(i). Hence \te s $q\in\N$ \st
$k=qm+\F_m(k) =qm$ by \er{1.2}, \er{1.3} since $m\in\Na$.

(vi) \ti{$\vf_a$ is in\jc}: follows from (ii).

\ti{$\vf_a$ is sur\jc}: Let $b\in I(a)$. By \df\ of $I(a)$, \te s $k\in\N$
\st $\vf_a(k)=b$. So $b=\vf_a(k) \nad{\rm(iv)}= \vf_a(\F_m(k)) = \ov\vf_a(\F_m
(k)) = \ov\vf_a(\F_m(k))$, since $\F_m(k)\in[0,m)$.

\ti{$\ov\vf_a$ is a \hm sm}: $\ov\vf_a(0)=\vf_a(0)\nde1.40 = e$. Let $k,l\in
[0,m)$. Then $\ov\vf_a(k+_ml)\nde1.8 = \ov\vf_a(\F_m(k+l)) \nad{\rm(iv)}=
\vf_a(\F_m(k+l)) = \vf_a(k+l) \nde1.41 = \vf_a(k)
\qu\vf_a(l) = \ov\vf_a(k) \qu \ov\vf_a(l)$.

\If from Lemma \rfa2{l1.8}\,(i) that $\ov\vf_a$ is an \is sm.
\endproof

Observe that the \sq\ $\vf_a:\N\to X$ of Lemma \rf{l1.22} \sf ies $\vf_a(k+m)
=\vf(k)$ \fa $k\in\N$. Indeed, $\vf_a(k+m)\nde1.41 = \vf_a(k) \qu \vf_a(m)
\nad{\rm(i)}= \vf_a(k)\qu e =\vf_a(k)$, $k\in\N$. Such a~\sq\ is called
\ti{$m$-periodic}.

\bds1.11
Let $F$ be a nonempty set, let $T$ be an \el\ of~$\Na$ and let $f$~be a map
(\sq) from~$\N$ into~$F$. The map $f$ is called \ti{$T$-periodic}\index{T-periodic@$T$-periodic} if it \sf ies
\beq1.14
f(x+T)=f(x) \qh{\fe}x\in\N.
\e
Notice that $f$ is a constant map iff it is $1$-periodic.
If $f$ is $T$-periodic, then $T$ is called the \ti{period\/} of the map~$f$.
If there is no $T'\in\Na$ with $T'<T$ \st $f$~is $T'$-periodic, then $T$ is
called the \ti{minimal period\/} of~$f$.\index{period}\index{period;minimal}
\eds

\blm1.12
Let $F,T$ and $f$ be as in \E\df s \rf{d1.11}.

\hph i,ii, If $f$ is $T$-periodic, then $f$ is $kT$-periodic for all $k\in\Na$.

\hph ii,i, If $f$ is $T$-periodic, it has a minimal period.

\hph iii,, If $f$ is $T$-periodic, then the image of $\N$ under~$f$ in~$F$ is
equal to the image of $[0,T)$ under~$f$ in~$F$.

\hph iv,, Let $g:[0,T)\to F$. Then \te s exactly one $T$-periodic map $f:\N\to F$
\st its \rt ion to $[0,T)$ is equal to~$g$. The \f~$f$ is called the
\emph{$T$-periodic extension} of~$g$.

\elm

\bex1.13
Prove Lemma \rf{l1.12}. (Hint: (i)~use induction, (ii)~use the well-\og\
of~$\N$, (iii)~use (i) and the division algorithm.)
\eex

\bex1.26
Show that the map $\vf_a$ of Lemma \rf{l1.22} is $m$-periodic with
\ti{minimal\/} period~$m$. Show that $\vf_a$ is the $m$-periodic \ext\
of~$\ov\vf_a$.
\eex

\brm1.27
In contrast with infinite \pn\ C-monoids, \ti{finite} \pn\ C-monoids may
possess more than one \Gn. For example in $(\N_4,+_4,0)$ the \IT s of $1,2,3$
are:
$$
\bal
&0,1,2,3,0,1,2,3,0,\dots\\
&0,2,0,2,0,2,0,2,0,\dots\\
&0,3,2,1,0,3,2,1,0,\dots
\eal
$$
We see that 1 and 3 are \Gn s of $(\N_4,+_4,0)$.
\erm

\Wanp prove an analogue of Theorem \rf{t1.9n}.

\bth1.28
Let $(X,\qu_0,\qu_1,e_0,e_1)$ be a \ti{finite} \pn\ \cnc e \sr. Then

\hph i,ii, $\ov\vf$, the \rt ion to $[0,n)$, $n:=\#(X)$, of the \sq\ of \IT s
of~$e_1$\break in $(X,{\qu_0},e_0)$, defined in \er{1.5n}, is a ring-\is sm from
$(\N_n,+_n,\cdot_n,0,1)$ onto\break ${(X,\qu_0,\qu_1,e_0,e_1)}$.

\hph ii,i, The map $\ov\vf$ is the only ring-\is sm from $\N_n$ onto~$X$.

\hph iii,, The identity is the only ring-auto\mf\ of~$X$.

\hph iv,, If $\wt X$ is a finite \pn\ \cnc e \sr\ and $\wt\vf:\N_n\to \wt X$
is defined in \er{1.5n}, where $e_1$ is replaced by~$\wt e_1$, then
$\ov{\wt\vf} \circ \ov\vf\Inv : X\to \wt X$ is the only ring-\is sm from~$X$
onto~$\wt X$.
\eth

\proof \

(i): \ti{$\vf$ is a monoid-\is sm from $(\N_n,+_n,0)$ onto $(X,\qu_0,e_0)$ \sf
ying $\ov\vf(1)=e_1$}: By \as\ $(X,\qu_0,e_0)$ is a finite \pn\ \Cm\ with
\Gn~$e_1$, i.e. $X=I(e_1)$. In view of Lemma \rf{l1.22}\,(vi) with $a:=e_1$,
we have $m=\#(I(e_1))$ since $\ov\vf$~is a~bi\jn\ from $\zo0,m $ onto
$I(e_1)$. Since $X=I(e_1)$, and $\#(X)=n$, we have $m=n\ge2$. \E\Tf
$\ov\vf:(\N_n,+_n,0) \to (X,\qu_0,e_0)$ is a monoid-\is sm by Lemma
\rf{l1.22}\,(vi). \Mo $\ov\vf(1)=\vf(1)=\vf_{e_1}(1)\nde1.40 = e_1$.

\ti{$\ov\vf$ is a monoid-\is sm from $(\N_n,\cdot_n,1)$ onto $(X,\qu_1,e_1)$}:
In view of Lemma \rfa2{l1.8}\,(i) it is \sft\ to show that $\ov\vf$~is a \hm sm.
We have $\ov\vf(1)=\vf(1)=e_1$. We have to show that $\ov\vf(k\cdot_n l)=
\ov\vf(k)\qu_1 \ov\vf(l)$ \fa $k,l\in\zo0,n $. Let $k,l\in\zo0,n $. Then
$\ov\vf(k\cdot_n l)= \vf(k\cdot_n l)\nde1.33 = \vf(\F_n(kl)) = \vf_{e_1}(\F_n
(kl)) = \vf_{e_1}(kl)$ by Lemma \rf{l1.22}\,(iv). But $\vf_{e_1}(kl)=(kl)\dquz
e_1 \nda22.8 = (k\dpl l)\dquz e_1 \nad{\era2{2.3}\,\rm I3} = k\dquz (l\dquz e_1)
= k\dquz((l\dquz e_1)\qu_1 e_1)\nad{\er{1.38}}= k\dquz
(\d_{l\dqumz e_1}(e_1))$. Since $\d_{l\dqumz e_1}$ is an endo\mf\ of
$(X,\qu_0,e_0)$, we obtain $k\dquz(\d_{l\dqumz e_1}(e_1))
= \d_{l\dqumz e_1}(k\dquz e_1)$ by Lemma \rfa2{l1.22}. \Mo
$\d_{l\dqumz e_1}(k\dquz e_1) \nde1.38 = (l\dquz e_1)\qu_1(k\dquz
e_1)\nde1.37 = (k\dquz e_1)\qu_1(l\dquz e_1) = \vf_{e_1}(k) \qu_1
\vf_{e_1}(l) = \vf(k)\qu_1\vf(l) = \ov\vf(k)\qu_1 \ov\vf(l)$. Hence $\ov\vf
(k\cdot_n l)=\ov\vf(k)\qu_1\ov\vf(l)$.

Thus we have proved that $\ov\vf:\N_n\to X$ is a ring-\is sm (see \E\df\
\rf{d1.7n}).

(ii) Let $\ov\psi:\N_n\to X$ be a ring-\is sm. We want to show that $\ov\vf=
\ov\psi$. To this end we define $\psi:\N\to X$ by setting $\psi(k):=\ov\psi
(\F_n(k))$, $k\in\N$. Observe that $\psi(k)=\ov\psi(\F_n(k))=\ov\psi(k)$ \fa
$k\in\zo0,n $. Hence $\psi$~is an \ext\ of~$\ov\psi$. We now show that
$\vf(k)=\psi(k)$, $k\in\N$, using \In\ on $k\in\N$. Set $A:=\{l\in\N:
\vf(l)=\psi(l)\}$. We have $\psi(0)=\vf(0)$. Hence $0\in A$.

\ti{$k\in A$ implies $k+1\in A$}: Suppose $k\in A$, then $\psi(k+1)=\ov\psi
(\F_n(k+1))\nde1.7 = \ov\psi(\F_n(\F_n(k)+\F_n(1)))\nde1.6 = \ov\psi(\F_n
(\F_n(k)+1)) \nde1.8 = \ov\psi(\F_n(k)+_n1) = \ov\psi(\F_n(k))\qu\ov\psi(1)$
since $\ov\psi:(\N_n,+_n,0)\to(X,\qu_0,e_0)$ is a \hm sm. Since $\ov\psi:
(\N_n,\cdot_n,1)\to (X,\qu_1,e_1)$ is a \hm sm, we have $\ov\psi(1)=e_1$.
\E\Tf $\psi(k+1)=\ov\psi(\F_n(k))\qu e_1 = \psi(k)\qu e_1 \nad{k\in A}=
\vf(k)\qu_0 e_1 = \vf(k)\qu_0\vf(1) = \vf_{e_1}(k)\qu_0 \vf_{e_1}(1)\nde1.41
= \vf_{e_1}(k+1)=\vf(k+1)$. Hence $k+1\in A$, and $A$~is \iv\ in $(\N,\le)$.
\If that $A=\N$ and $\vf=\psi$. The \rt ion of~$\vf$ to $\zo0,n $ is~$\ov\vf$
by \df\ of~$\ov\vf$, and the \rt ion of~$\psi$ to $\zo0,n $ is $\ov\psi$ since
$\psi$~is an \ext\ of~$\ov\psi$. \If that $\ov\vf=\ov\psi$.

(iii) Let $\ov\psi$ be a ring-auto\mf\ of the \sr~$X$. Then $\ov\rho:= \ov\psi
\circ \ov\vf:\N\to X$ is a ring-\is sm, by Lemma \rfa2{l1.8}. Thus $\ov\rho=
\ov\vf$ by~(ii). Hence $\id_X=\ov\vf\circ\ov\vf{}\Inv= \ov\rho\circ\ov\vf{}
\Inv = (\ov\psi \circ \ov\vf)\circ \ov\vf{}\Inv=
\ov\psi\circ(\ov\vf\circ\ov\vf{}\Inv)=\ov\psi\circ\id_X=\ov\psi$.

(iv) By Lemma \rfa2{l1.8}, $\ov{\wt\vf}\circ\ov\vf{}\Inv$ is a ring-\is sm
from $X$ onto~$\wt X$. Let $\ov\psi:X\to \wt X$ be a ring-\is sm. Then
$\ov\psi{}\Inv\circ(\ov{\wt\vf}\circ\ov\vf{}\Inv)$ is a ring-\is sm of~$X$.
By (iii) it is equal to $\id_X$, thus $\ov\psi{}\Inv\circ(\ov{\wt\vf}\circ\ov
\vf{}\Inv)=\id_X$. Hence $\ov\psi=\ov\psi\circ\id_X=\ov\psi\circ(\ov\psi{}\Inv
\circ(\ov{\wt\vf}\circ\ov\vf{}\Inv)) = (\ov\psi\circ\ov\psi{}\Inv)\circ
(\ov{\wt\vf}\circ\ov\vf{}\Inv)=\id_X\circ(\ov{\wt\vf}\circ\ov\vf{}\Inv)=
\ov{\wt\vf}\circ\ov\vf{}\Inv$.
\endproof

We conclude this section by considering some divisibility criteria. Let $m\in
\Na\sms1$ and let
\beq1.44
m=\sum_{k=0}^N \a_k 10^k, \q N\in\N,\ \a_k\in[0,9] \hbox{ for }k\in[0,N],\
\a_N\ne0,
\e
be its \ti{\dr} introduced in Lemma \rfa2{l5.3}. We recall that the \rp ation
\er{1.44} is \ti{unique} by Lemma \rfa2{l5.4}. The \nm~$\a_k$, $k\in[0,N]$, is
called the $k$-th \ti{digit\/}\index{digit} of~$m$ and $\a_0$~is called the \ti{last\/}
digit of~$m$ even if $N=0$. Let $d\in \Na\sms1$. We are interested in \cn s
on the digits of~$m$ insuring the divisibility of~$m$ by~$d$. In view of
\era3{8.57} we assume $d\le m$, thus
\beq1.45
2 \le d \le m.
\e
We first observe that
\beq1.46
d|m \qh{iff }\phi_d(m)=0,
\e
where $\phi_d(m)$ is defined in \er{1.3}, \er{1.2}. Indeed, if $d|m$, then
\te s $q\in\Na$ \st $m=qd$. Hence $\phi_d(m)=\phi_d(qd) = \phi_d(0+qd)
\nde1.4 = \phi_d(0)\nde1.6 = 0$.

Conversely, if $\phi_d(m)=0$, then $m=qd$ \fs $q\in\N$ by \er{1.2}, \er{1.3}.
Since $m\ne0$, we have $q\ne0$ by \era2{2.13}. Hence $d|m$.

\bpr1.29
Let $d,m\in\N$ be \st \er{1.45} holds and let \er{1.44} be the \dr\ of~$m$.
Then
\beq1.47
\phi_d(m) = \phi_d\Bg(\sum_{k=0}^N \a_k\rho_{d,k}),
\e
where
\beq1.48
\rho_{d,k} := \phi_d(10^k), \q k\in[0,N].
\e
\epr

\brm1.30
$\rho_{d,0} = \phi_d(10^0) \nda22.23 = \phi_d(1)\nde1.5 = 1$.
\erm

The proof of \E\Pr\ \rf{p1.29} is a con\sq\ of the \fw\ lemmata.

\blm1.31
Let $n\in\Na\sms1$, $N\in\N$ and $b:[0,N]\to\N$. Then
\beq1.49
\phi_n\Bg(\sum_{k=0}^N b_k) = \F_n\Bg(\sum_{k=0}^N {\F_n(b_k)}).
\e
\elm

\proof
Set $A:=\{M\in\N: \hbox{\er{1.49} holds for }N:=M\}$.

``$0\in A$'': $\sum_{k=0}^0 b_k=b_0$ by \df. Hence $\F_n\bg(\suml_{k=0}^0 b_k) =
\F_n(b_0) \nde1.5 = \F_n(\F_n(b_0)) = \F_n\bg(\suml_{k=0}^0 {\F_n(b_k)})$.

``$1\in A$'': $\F_n\bg(\suml_{k=0}^1 b_k) = \F_n(b_0+b_1) \nde1.7 =
\F_n(\F_n(b_0)+\F_n(b_1))= \F_n\bg(\suml_{k=0}^1 {\F_n(b_k)})$.

``\ti{$M\,{\in}\, A$ implies $M{+}1\,{\in}\, A$}'': Let $M\in A$, then $\F_n\bg(\suml_{k=0}
^{M+1} b_k)\nda27.13  =\F_n\bigl(\bg(\suml_{k=0}^M b_k) +b_{M+1}\bigr)\break
\nde1.7 = \F_n\bigl(\F_n\bg(\suml_{k=0}^M b_k)+\F_n(b_{M+1})\bigr) \nad{M\in
A} = \F_n\bigl(\F_n\bigl(\suml_{k=0}^M \F_n(b_k)\bigr)+\F_n(b_{M+1})\bigr)
\nde1.5 = \F_n\bigl(\F_n\bigl(\suml_{k=0}^M \F_n(b_k)\bigr)\break+ \F_n\bigl(
\F_n(b_{M+1})\bigr)\bigr) \nde1.7  =\F_n\bigl(\bigl(\suml_{k=0}^M \F_n(b_k)
+\F_n(b_{M+1})\bigr) \nda27.13 = \F_n\bigl(\suml_{k=0}^{M+1}\F_n(b_k)\bigr)$.
Thus $M+1\in A$. Since $A$ is \iv\ in $(\N,0,S)$, $A=\N$.
\endproof

\blm1.32
Let $n\ge2$, $b,c\in\N$. Then
\beq1.50
\F_n(bc) = \F_n(b\cdot \F_n(c)).
\e
\elm

\proof
By \er{1.2}, \er{1.3}, $c=q\cdot n+\F_n(c)$ \fs $q\in\N$. Hence $bc\nda22.14
= bqn+b\F_n(c) = (bq)n +b\F_n(c)$. Thus $\F_n(bc)=\F_n((bq)n + b\F_n(c))
\nde1.4 = \F_n(b\F_n(c))$.
\endproof

\proof[Proof of \E\Pr\ \rf{p1.29}]
$\F_n(m)\nde1.44 = \F_d\bg(\suml_{k=0}^N {(\a_k10^k)}) \nde1.49 = \F_d\bigl(
\suml_{k=0}^N \F_d(\a_k10^k)\bigr) \nde1.50 = \\ \F_d\bg(\suml_{k=0}^N {\F_d(\a_k\F_d
(10^k))}) \nde1.48 = \F_d\bg(\suml_{k=0}^N\a_k\rho_{d,k})$.
\endproof

\If from \er{1.46} and \E\Pr\ \rf{p1.29} that $d|m$ iff
$\F_d\bg(\suml_{k=0}^N \a_k\rho_{d,k})=0$. This does not mean that $d\big|
\suml_{k=0}^N\a_k\rho_{d,k}$ since it may happen that $\suml_{k=0}^N\a_k
\rho_{d,k}=0$. In \era3{8.52}, $a|b$ is only defined for $a,b\in\Na$. For example, if
$d=m=10$, then $d|m$ and $N=2$, $\a_0=0$, $\a_1=1$. Hence $\suml_{k=0}^1
\a_k\rho_{d,k}=\rho_{d,1}=\F_{10}(10)=0$. \E\oh if $\suml_{k=0}^N \a_k\rho_{d,k}
\ne0$, then the proof of the ``if'' part of \er{1.46} shows that $\F_d\bg(
\suml_{k=0}^N \a_k\rho_{d,k})=0$ implies $d|\suml_{k=0}^N \a_k\rho_{d,k}$.

\bth1.33
Let $d,m\in\N$ \sf y \er{1.45} and let \er{1.44} be the \dr\ of~$m$ with $N\ge1$.
Then $d|m$ iff either
$\suml_{k=0}^N \a_k\rho_{d,k}=0$ or $d|\suml_{k=0}^N \a_k\rho_{d,k}$.

\E\Ip if $d\nmid 10^N$ then
\beq1.51
d|m \qh{iff }d|\suml_{k=0}^N \a_k\rho_{d,k}.
\e
\eth

\proof
The first assertion follows from \er{1.46}, \er{1.47}, \er{1.48} and
\beq1.52
\F_n(l)=0 \qh{iff either} l=0 \hbox{ or }n|l, \q n,l\in\N,\ n\ge2.
\e
The second assertion follows from the first one provided we can show that if
$d\nmid10^N$, then $\suml_{k=0}^N\a_k\rho_{d,k}\ne0$. Observe that
if $d\nmid10^N$, then $\F_d(10^k)\ne0$ for $k\in[1,N]$.
Since $N\ge1$, note that $d\nmid10^N$ implies $d\nmid 10^k$ \fa
$k\in[1,N]$. If $N=1$, there is nothing to prove. If $N>1$ and $k\in[1,N-1]$,
then $10^N=10^{N-k}\cdot10^k$ with $k,N-k\ge1$. Thus if $d|10^k$, then
$d|10^N$ hence by contraposition $d\nmid 10^N$ implies $d\nmid 10^k$. \If from
\er{1.52} that if $d\nmid 10^N$, then $\F_d(10^k)\ne0$ \fa $k\in[1,N]$. \E\Tf
if $\suml_{k=0}^N \a_k\F_d(10^k)=0$, then $a_k\F_d(10^k)=0$ \fa $k\in[0,N]$,
by \era2{7.15} with $(X,\qu,e):=(\N,+,0)$. Since $\F_d(10^k)\ne0$ \fa
$k\in[1,N]$, then $\a_k=0$, $k\in[1,N]$, by \E\Pr\ \rfa2{p2.7}\,(ii). \If
that $\suml_{k=0}^N \a_k\F_d(10^k) = \a_0\F_d(1)+ \suml_{k=1}^N 0\F_d(10^k)
\nda27.27 = \a_o1 + 0\suml_{k=1}^N \F_d(10^k) =\a_0+0 = \a_0$. Thus if
$\suml_{k=0}^N\a_k\F_d(10^k)=0$, then $\a_k=0$ \fa $k\in[0,N]$, hence $m=
\suml_{k=0}^N 0\cdot10^k =  0\suml_{k=0}^N 10^k=0$, a~\cd ion. Hence if
$N\ge1$ and $d\nmid 10^N$, then $\suml_{k=0}^N \a_k\rho_{d,k}\ne0$.
\endproof

We next give some \pp ies of the set $\{\rho_{d,k}\}_{k\in\N}$.

\bpr1.34
Let $n,a\in\Na\sms1$. Then
\beq1.53
\F_n(a^k) = k\ddtn \F_n(a), \q k\in\N,
\e
where $k\ddtn \F_n(a)$ denotes the $k$-fold \IT\ of $\F_n(a)$ in the monoid
$(\N_n,\cdot_n,1)$. \E\Ip if $\F_n(a^l)=0$ \fs $l\in\Na$, then
\beq1.54
\F_n(a^k)=0 \qh{\fa} k\ge l.
\e
\Mo \te\ $k_0\in\N$ and $T\in\Na$ \st
\beq1.55
\F_n(a^{k+T}) = \F_n(a^k) \qh{\fa} k\ge k_0.
\e
If $k_0=0$, then the \sq\ $\{\F_d(a^k)\}_{k\in\N}$ is called
\emph{$T$-periodic} $($see \E\df\ \rf{d1.11}$)$, otherwise it is called
\emph{eventually ($T$-)periodic}.\index{periodic!eventually}

Finally, if $k_0$ is as in \er{1.55}, then
\beq1.56
\F_n(a^k)\ne\F_n(a^{k'}) \qh{if} k,k'\in[0,k_0) \hbox{ and }k\ne k'.
\e
\epr

\proof \

\er{1.53}: Note that $\F_n(l)\in\N_n$, $l\in\N$, by \er{1.3}. We use \In\ on
$k\in\N$. Set $A:=\{l\in\N:$ \er{1.53} holds for $k:=l\}$. Then

``$0\in A$'': $\F_n(a^0)\nda22.23 = \F_n(1)\nde1.6 = 1 = 0\ddtn \F_n(a)$.

``\ti{$k\in A$ implies $k+1\in A$}'': Let $k\in A$. Then $\F_n(a^{k+1})\nda22.25 =
\F_n(a^k\cdot a^1) \nde1.33 = \F_n(a^k)\cdot_n \F_n(a^1) \nad{k\in A,\era2{2.24},
\era2{2.3}\,\rm I1} = (k\ddtn \F_n(a))\cdot_n (1\ddtn\F_n(a)) \nad{\era2{2.3}\,
\rm I2} = (k+1)\ddtn \F_n(a)$.
Since $A$ is \iv\ in $(\N,0,S)$, we have $A=\N$.

\er{1.54}: Suppose $\F_n(a^l)=0$, $l\ge1$, and let $k\ge l$. Let $p\in\N$ be
\st $l=p+k$. Then $\F_n(a^k)=k\ddtn \F_n(a) = (p+k)\ddtn \F_n(a)
\nad{\era2{2.3}\,\rm I2} = (p\ddtn \F_n(a))\cdot_n (k\ddtn \F_n(a)) =\break
(p\ddtn \F_n(a))\cdot_n0 \nde1.33 = \F_n((p\ddtn\F_n(a))\cdot0) \nad{
\era2{2.13}} = \F_n(0)\nde1.6 = 0$.

\er{1.55}: In view of \era2{2.3}\,I0 and I2, the map $k\mt k\ddtn\F_n(a)$ is a
\hm sm from $(\N,+,0)$ into $(\N_n,\cdot_n,0)$. This \hm sm is not in\jc,
since otherwise its range would be infinite by \E\Pr\ \rfa1{p4.27}\,(A)(i).
Indeed, $\N$~is infinite and the range of this \hm sm is a subset of~$\N_n$,
which is finite by Theorem \rfa1{t4.18}. Hence we may apply Lemma
\rfa2{l1.41}, and from its proof (see~\era2{1.67}) we find $k_0:=m\in\N$,
$T:=p\in \Na$, $s\in\N$, \st
\beq1.57
(k_0+l+sT)\ddtn \F_n(a) = (k_0+l)\ddtn \F_n(a), \q l\in[0,T).
\e
Let $k\ge k_0$. Then by \er{1.2} \te\ $q\in\N$, $r\in[0,T)$ \st $k-k_0=qT+r$,
i.e.\ $k=k_0+qT+r$. Hence $k+T = k_0+(q+1)T+r$. By \er{1.57} with $l:=r$,
$s:=q$ and $s:=q+1$, we obtain $k\ddtn \F_n(a) = (k_0+qT+r)\ddtn \F_n(a) =
(k_0+r)\ddtn \F_n(a) \nde1.57 = (k_0+(q+1)T+r)\ddtn\F_n(a) = (k+T)\ddtn
\F_n(a)$.

\er{1.56}: directly follows from Lemma \rfa2{l1.41}\,(i).
\endproof

Recall that by Theorem \rfa3{t9.22} every \nn\ $d>1$ is a product of prime
powers, i.e.\ $d=\prodl_{i=1}^M p_i^{m_i}$ where $M\in\Na$, $p_i$~prime,
$m_i\in\Na$, $i\in [1,M]$, and $p_i\ne p_j$ whenever $i\ne j$, $i,j\in[1,M]$.
\If from Lemma \rf{l2.16} and Lemma \rfa3{l9.17} that $\gcd(p_i^{m_i},p_j^{m_j})=1$ if $i\ne j$,
$i,j\in[1,M]$. \Mo in view of Lemma \rfa3{l1.28} and its \gn, Lemma
\rf{l2.30}, if $m\in\Na$, $d\le m$, then $d|m$ iff $p_i^{m_i}|m$ \fa
$i\in[1,M]$. \E\Tf we first consider divisibility criteria for $d=p^l$,
$p$~prime and $l\in\Na$. Since $2$~and~$5$ are the only prime divisors
of~$10$, we distinguish two cases: case~A where $p\notin\{2,5\}$ and case~B
where $p\in\{2,5\}$.

\ti{Case} A: $d=p^l$, $p\notin\{2,5\}$, $l\in\Na$.

We claim $\F_d(10^k)\ne0$ \fa $k\in\N$. Suppose for \cd ion that
$\F_d(10^k)=0$ \fs $k\in\N$. Then $k\ne0$ since $\F_d(10^0)=\F_d(1)=1\ne0$.

If $k\ge1$ then $10^k\nda22.28 = 2^k\cdot 5^k$, and $p^l|(2^k\cdot5^k)$ by
\er{1.52}. As observed above, ${\gcd(p^l,2^k)=1}$ and $\gcd(p^l,5^k)=1$ since
$p\notin \{2,5\}$. Hence by Lemma \rfa3{l1.29} we have $p^l|2^k$ and
$p^l|5^k$. A~\cd ion, since $p^l|2^k$ (resp.~$p^l|5^k$) implies $\gcd(p^l,2^k)
=p^l$ (resp.\ $\gcd(p^l,5^k)=p^l$) and $p^l\ne1$. This completes the proof of
the claim. Now let $m\in\Na$, $m\ge p^l$, be given by \er{1.44} with $N\ge1$, then by
\er{1.51} we obtain
\beq1.58
p^l|m \qh{iff }p^l\big| \sum_{k=0}^N\a_k\F_{p^l}(10^k).
\e

We now consider some simple examples.

$d:=3$: Note that $\F_3(10)=1$, hence $\F_3(10^k)\nde1.53 = k\ddtt3 1\nad
{\era2{2.3}\,\rm I4} = 1$. \csq,
\beq1.59
3|m \qh{iff }3\big|\sum_{k=0}^N\a_k.
\e
Observe that if $N=10$ then $m\ge10^{10}$, but $\suml_{k=0}^{10} \a_k\le99$,
thus $\suml_{k=0}^{10} \a_k$ is equal to $\b_0\ne0$ or $\b_0+\b_1\cdot10$,
$\b_1\ne0$, $\b_0,\b_1\in[0,9]$. Hence $3|\suml_{k=0}^{10} \a_k$ iff $3|\b_0$
or $3|(\b_0+\b_1\cdot10)$.

$d:=9$: $\F_9(10)=1$, hence as above $\a_k=1$ \fa $k\in[0,N]$. Thus
\beq1.60
9|m \qh{iff } 9\big|\sum_{k=0}^N \a_k.
\e

$d:=27$: $\F_{27}(10)=10$, $\F_{27}(10^2)=19$, $\F_{27}(10^3) = 10\cdot_{27}
\F_{27}(10^2)=\F_{27}(190)=1$, hence $\F_{27}(10^3)=1$, $\F_{27}(10^4) =
10\cdot_{27}\F_{27}(10^3)=10$. \If that $\{\F_{27}(10^k)\}_{k\in\N}$ is
3-periodic, with $\F_{27}(10^0)=1$, $\F_{27}(10^1)=10$, $\F_{27}(10^2)=19$.
Hence if $k=i+3l$, $l\in\N$, $i\in\zo0,3 $, then $\F_{27}(10^k)=
\F_{27}(10^i)= \F_{27}(10^{\F_3(k)})$. Thus
\beq1.61
27|m \qh{iff }27\big|\sum_{k=0}^N\a_k\F_{27}(10^{\F_3(k)}).
\e
\E\Ip if $N=7$, $27|m$ iff $27|(\a_0+\a_3+\a_6)+10(\a_1+\a_4+\a_7)
+19(\a_2+\a_5)$.

$d:=81$: $\F_{81}(10^0)=1$, $\F_{81}(10^1)=10$, $\F_{81}(10^2)=19$,
$\F_{81}(10^3)=28$, $\F_{81}(10^4)=37$, $\F_{81}(10^5)=46$,
$\F_{81}(10^6)=55$, $\F_{81}(10^7)=64$, $\F_{81}(10^8)=73$,
$\F_{81}(10^9)=1$. \If that $\{\F_{81}(10^k)\}$ is 9-periodic. For example
if $N=9$, then $81|m$ iff $81|(\a_0+\a_9)+10\a_1+19\a_2+28\a_3+37\a_4
+46\a_5+55\a_6+64\a_7+73\a_8$. \E\Ip if $m=1111111020$, then $\suml_{k=0}^N
\a_k\F_{81}(10^k) = 0+20+0+28+37+46+55+64+73+1=324$. Since $81|324$, $81|m$.

$d:=7$: $\F_7(10^0)=1$, $\F_7(10^1)=3$, $\F_7(10^2)=\F_7(3\cdot3)=2$,
$\F_7(10^3) =\F_7(3\cdot2)=6$, $\F_7(10^4)=\F_7(3\cdot6)=4$,
$\F_7(10^5)=\F_7(3\cdot4)=5$, $\F_7(10^6)=\F_7(3\cdot5)=1$. Hence
\beq1.62
\{\F_7(10^k)\}_{k\in\N} \qh{is 6-periodic:} 1,3,2,6,4,5,1 \dots
\e

$d:=11$: $\F_{11}(10^0)=1$, $\F_{11}(10^1)=10$,
$\F_{11}(10^2)=\F_{11}(10\cdot 10)=1$.
\beq1.63
\{\F_{11}(10^k)\}_{k\in\N} \qh{is 2-periodic:}1,10,1\dots
\e

$d:=13$: $\F_{13}(10^0)=1$, $\F_{13}(10^1)=10$, $\F_{13}(10^2)=9$, $\F_{13}(10^3)
=\F_{13}(10\cdot9)=12$, $\F_{13}(10^4)=\F_{13}(10\cdot12)=3$, $\F_{13}(10^5)=
\F_{13}(10\cdot3)=4$, $\F_{13}(10^6)=\F_{13}(10\cdot4)=1$. Hence
\beq1.64
\{\F_{13}(10^k)\}_{k\in\N} \qh{is 6-periodic: }1,10,9,12,3,4,1,\dots.
\e

\ti{Case} B:

$d:=2$: $\F_2(10^0)=1$, $\F_2(10^1)=0$, $\F_2(10^k)=0$, $k>1$. Thus $\{\F_2(2^k)
\}_{k\in\N}$ is 1-periodic from $k=1$: $1,0,0,\dots$. From Theorem \rf{t1.33} we infer:
\beq1.65
\hbox{if $N\ge1$, then $2|m$ iff }\a_0\in\{0,2,4,6,8\}.
\e

$d:=2^l$, $l\ge1$: We have $2^l|10^k$ iff $k\ge l$. If $l\le k$, $k=l+p$ \fs
$p\in\N$. Then $2^k\nda22.25 = 2^l\cdot 2^p$ and $10^k=(2\cdot5)^k \nda22.28 =
2^k\cdot 5^k = (2^l\cdot 2^p)\cdot 5^k \nda21.5 = 2^l\cdot(2^p\cdot 5^k)$.
Hence $2^l|10^k$.

If $2^l|10^k$, then $2^l|(2^k\cdot5^k)$ and by Lemma \rfa3{l1.29}, $2^l|2^k$ since
$\gcd(2^l,5^k)=1$. Thus $l\le k$ by Lemma \rf{l2.16}.

In view of \er{1.52}, $\F_{2^l}(10^k)=0$ iff $k\ge l$. Hence by Theorem \rf{t1.33}
\beq1.66
2^l|m \qh{iff \ either $\N\ge l$ and $\a_k=0$, $k\in[0,l-1]$, or \ }
d\big|\sum_{k=0}^{\min(N,l-1)} \a_k \F_{2^l}(10^k).
\e

$d=4$: $\F_4(10^0)=1$, $\F_4(10^1)=2$, $\F_4(10^2)=0$.

$d=8$: $\F_8(10^0)=1$, $\F_8(10^1)=2$, $\F_8(10^2)=4$, $\F_8(10^3)=0$.

$d=16$: $\F_{16}(10^0)=1$, $\F_{16}(10^1)=10$, $\F_{16}(10^2)=4$, $\F_{16}(10^3)
=8$, $\F_{16}(10^4)=0$.

For example, if $N\ge4$, then $16|m$ iff $\a_0=\a_1=\a_2=\a_3=0$, or $16|(\a_0
+10\a_1+4\a_2+8\a_3)$. \E\Ip $16|1104$ since $16|(4+0+4+8)$.

\bex1.35 \

\hph i,i, Compute $\{\F_{17}(10^k)\}_{k\in\N}$, $\{\F_{19}(10^k)\}_{k\in\N}$.
(Hint: use \er{1.50}, $\F_d(10^{k+1})= \F_d(\F_d(10)\cdot\F_d(10^k))$, $k\in\N$).

\hph ii,, Find an analogue of \er{1.66} for $5^l|m$, $l\in\Na$.
\eex

\brm1.36
Using Lemma \rfa3{l1.28}, one can find divisibility criteria for $d=d_1\cdot
d_2$ where $\gcd(d_1,d_2)=1$. For example, if $d=6$, then
\beq1.67
6|d \qh{iff }2|d \hbox{ and }3|d \qh{iff }\a_0\in\{0,2,4,6,8\} \hbox{ and }
3\big|\sum_{k=0}^N \a_k \hbox{ whenever }N\ge1.
\e
\erm

More on ``divisibility criteria'' can be found at the end of Section
\ref{sss.indecomp}, and in Section~\ref{sss.ipr.fld} Exercise \rf{ex5.25}.

\newpage

\Subsubsection{Indecomposable \pn\ c-\sr s}\label{sss.indecomp}

\bnt2.1
In what follows we shall use the shorthand notation \ti{c-semiring}\index{c-semiring} for \cnc
e \sr\ (which we recall, means \cmt e \sr\ with unity and with \cnc e
additive monoid!).
\ent

In Remark \rfa2{r1.6} we introduced the notion of direct product of two monoids.
Similarly we define the direct product of two \sr s.

\bpr2.2
Let $(X_i,\qn i_0,\qn i_1,e_0\en i,e_1\en i)$, $i=1,2$, be \sr s. Let $\qu_0$
and $\qu_1$ be  the binary \op s on $X_1\t X_2$ defined by\dw
\bea2.1
(x_1,x_2)\qu_0(y_1,y_2)&:= (x_1\qn 1_0 y_1,x_2\qn 2_0 y_2),\\
(x_1,x_2)\qu_1(y_1,y_2)&:= (x_1\qn 1_1 y_1,x_2\qn 2_1 y_2), \lb{2.2}
\e
where $x_i,y_i\in X_i$, $i=1,2$.

Let $e_0$ and $e_1$ be the \el s of $X_1\t X_2$ defined by
\beq2.3
e_0:=(e\en 1_0,e\en 2_0), \q e_0:=(e\en 1_1,e\en 2_1).
\e
Then

\hph i,ii, $(X_1\t X_2,\qu_0,\qu_1,e_0,e_1)$ is a \sr\ called the
\emph{direct product\/} of the \sr s $(X_i,\qn i_0,\qn i_1,e\en i_0,e\en
i_1)$, $i=1,2$.

\hph ii,i, If $\pi_i:X_1\t X_2\to X_i$, $i=1,2$, are defined by
\beq2.4
\pi(x_1,x_2):=x_i, \q i=1,2,
\e
then $\pi_1$ and $\pi_2$ are sur\jc\ ring-\hm sms. The maps $\pi_1$ and
$\pi_2$ are called \emph{projection}\index{projection} maps of $X_1\t X_2$ onto $X_1$
and~$X_2$.

\hph iii,, If $X_i'$, $i=1,2$, are \sr s and $j_i:X_i\to X_i'$, $i=1,2$, are
ring-\is sms, then $j:X_1\t X_2 \to X_1'\t X_2'$ defined by
\beq2.5
j(x_1,x_2)=(j_1(x_1),j_2(x_2)), \q x_i\in X_i,\ i=1,2,
\e
is a ring-\is sm.

\hph iv,, $X_1\t X_2$ is a c-\sr\ iff $X_1$ and $X_2$ are c-\sr s.

\hph v,i, If $X_1\t X_2$ is a \pn\ \sr, then so are $X_1$ and $X_2$.
\epr

\proof\

(i) $(X_1\t X_2,\qu_0,e_0)$ and $(X_1\t X_2,\qu_1,e_1)$ are abelian monoids
by Remark~\rfa2{r1.6}. We next prove \er{1.1n} and \er{1.2n}.

\er{1.1n}: $(e_0\en 1,e_0\en 2)\qu_1 (x_1,x_2) \nde2.2 = (e_0\en1 \qn1_1 x_1,
e_0\en2 \qn2_1 x_2) = (x_1\qn1_1 e_0\en1, x_2\qn2_1 e_0\en2) \nde1.1n =
(e_0\en1,e_0\en2)$, $x_i\in X_i$, $i=1,2$.

\er{1.2n}: $(x_1,x_2)\qu_1((y_1,y_2)\qu_0(z_1,z_2)) \nde2.1 = (x_1,x_2)
\qu_1(y_1\qn1_0 z_1,y_2\qn2_0z_2) \nde2.2 =\break (x_1\qn1_1 (y_1\qn1_0 z_1),
x_2\qn2_1(y_2\qn2_0z_2)) \nde1.2n = ((x_1\qn1_1 y_1)\qn1_0 (x_1\qn1_1z_1),
(x_2\qn2_1y_2) \qn2_0 (x_2\qn2_1 z_2)) \nde2.1 = (x_1\qn1_1y_1,
x_2\qn2_1 y_2) \qu_0 (x_2\qn1_1 z_1,x_2\qn2_1 z_2)\nde2.2 = ((x_1,x_2)
\qu_1(y_1,y_2))\qu_0 ((x_1,x_2)\qu_1(z_1,z_2))$, $x_i,y_i,z_i\in X_i$,
$i=1,2$.

\ssk
(ii) $\pi_1$ is \ti{sur\jc} since $x_1\nde2.4 = \pi_1(x_1,e\en2_0)$ \fa
$x_1\in X_1$.

\ti{$\pi_1:(X_1,\qn1_0,e\en1_0)\t (X_2,\qn2_0,e\en2_0) \to
(X_1,\qn1_0,e\en1_0)$ is a monoid-\hm sm}:\break $\pi_1(e\en1_0,e\en2_0)\nde2.4 =
e\en1_0$ and $\pi_1((x_1,x_2) \qu_0 (y_1,y_2)) \nde2.1 = \pi_1(x_1\qn1_0
y_1, x_2\qn2_0y_2) \nde2.4 = x_1\qn1_0 y_1\nde2.4 =
 (\pi_1(x_1,x_2))\qn1_0(\pi_1(y_1,y_2))$.

\ti{$\pi_1:(X_1,\qn1_1,e\en1_1)\t (X_2,\qn2_1,e\en2_1) \to
(X_1,\qn1_1,e\en1_1)$ is a monoid-\hm sm}:\break $\pi_1(e\en1_1,e\en2_1)=e\en1_1$
and $\pi_1((x_1,x_2)\qu_1(y_1,y_2))\nde2.2 = \pi_1(x_1\qn1_1y_1,x_2\qn2_1y_2)
={x_1\qn1_1y_1} = (\pi_1(x_1,x_2))\qn1_1(\pi_1(y_1,y_2))$. Hence $\pi_1$ is
a~sur\jc\ ring-\hm sm. The proof for~$\pi_2$ is similar.

\ssk
(iii) and (iv) are left as exercises.

(v) Let $X_1\t X_2$ be \pn. In view of \E\df\ \rf{d1.6n}, \E\df\ \rfa2{d1.18}
and \er{2.3}, $X_1\t X_2=\bcl_{n\in\N}n\dquz (e\en1_1,e\en2_1)$, i.e.\ given
$(x_1,x_2)\in X_1\t X_2$, \te s $n\in\N$ \st $(x_1,x_2)=n\dquz (e\en1_1,
e\en2_1)$. We claim
\beq2.6
m\dquz(y_1,y_2) = (m\dqun1 y_1,m\dqun2 y_2), \q m\in\N, \ y_i\in X_i,\ i=1,2.
\e
Indeed, since $\pi_1$ and $\pi_2$ are ring-\hm sms, $\pi_1$ and $\pi_2$ are monoid-\hm sms
from $(X_1\t X_2,\qu_0,e_0)$ into $(X_1,\qu_0\en1,e_0\en1)$ and
$(X_2,\qu_0\en2,e_0\en2)$, \rsp. Using Lemma \rfa2{l1.22} we obtain for $m\in\N$,
$i=1,2$, $y_i\in X_i$, $\pi_i(m\dquz(y_1,y_2))\nda2{1.45} = m\dqun i \pi_i
(y_1,y_2)=m\dqun i y_i$.
Then $m\dquz(y_1,y_2)= \bigl(\pi_1(m\dquz(y_1,y_2)), \pi_2(m\dquz(y_1,y_2))\bigr)
=\bigl(m\dqun1 \pi_1(y_1,y_2),\break m\dqun2 \pi_2(y_1,y_2)\bigr) =
(m\dqun1 y_1,m\dqun2 y_2)$, which implies \er{2.6}. \If that if $x_1\in X_1$,
$x_2\in X_2$, then $x_i$ is an \IT\ of $e\en i_1$ in $(X_i,\qu\en i_0,
e\en i_0)$, $i=1,2$. Thus both $X_1$ and $X_2$ are \pn\ \sr s.
\endproof

\bex2.3 \

\hph i,i, Prove (iii) and (iv) of \E\Pr\ \rf{2.2}.

\hph ii,, Let $X_1$ and $X_2$ be \sr s. Show that $X_1\t X_2$ and $X_2\t X_1$
are ring-\is c even if $X_1$ and~$X_2$ are not ring-\is c.
\eex

We now introduce a \pp y of \sr s which is preserved under ring-\is sms.

\bdf2.4
A \sr\ $X$ is called \ti{\dcm} (not standard terminology)\index{decomposable semiring} if \te\ \sr s $X_1$
and~$X_2$ \st $X$~is ring-\is c to $X_1\t X_2$. Otherwise $X$~is called
\ti{in\dcm}.
\edf

\blm2.5
Let $X$ and $X'$ be \is c \sr s. Then

\hph i,ii, if $X$ is \cnc e, then so is $X'$,

\hph ii,i, if $X$ is \pn, then so is $X'$,

\hph iii,, if $X$ is \dcm\ $($resp.\ in\dcm$)$, then so is $X'$.
\elm

\proof
Let $f:X\to X'$ be a ring-\is sm.

(i) Let $a',b,'c'\in X'$ be \st $a'\qu_1' c'= b'\qu_1 c'$. Let $a,b,c\in X$
be \st $f(a)=a'$, $f(b)=b'$, $f(c)=c'$. Then $f(a\qu_1 c)=f(a)\qu_1' f(c)
=a'\qu'_1 c' = b'\qu'_1 c'=f(b\qu_1c)$. Hence $a\qu_1 c=f\Inv(f(a\qu_1c))=
f\Inv(f(b\qu_1c))
= b\qu_1c$. If $X$~is \cnc e then $a=b$, hence $a'=f(a)=f(b)=b'$.

(ii) Let $a'\in X'$ and $a\in X$ be \st $a'=f(a)$. Since $X$~is \pn, \te s
$n\in\N$ \st $a=n\dquz e_0$. Then $a'=f(n\dquz e_0)\nda2{1.45} =n
\mathrel{\lower1pt\hbox{$\stackrel{\lower7pt\hbox{\smash{\hbox to0pt{$\scriptstyle\sqcap'$\hss}$\scriptstyle\sqcup_0$}}}\cdot$}}
f(e_0)=n
\mathrel{\lower1pt\hbox{$\stackrel{\lower7pt\hbox{\smash{\hbox to0pt{$\scriptstyle\sqcap'$\hss}$\scriptstyle\sqcup_0$}}}\cdot$}}
e_0'$. Hence $X=I(e_0')$.

(iii) Suppose $X$ is \dcm\ and let $X_1,X_2$ be \sr s \st $g:X\to X_1\t X_2$
is an \is sm. Then $g\circ f\Inv:X'\to X_1\t X_2$ is an \is sm. Hence $X'$ is
\dcm. The case $X$~\dcm\ follows by contraposition.
\endproof

\brm2.6
Lemma \rf{l2.5}\,{\rm(ii)} holds under the weaker \as{\rm:} \te s $f:X\to X'$,
$f$~sur\jc\ \hm sm (see \E\Pr\ \rf{p2.2}\,{\rm(v)}).
\erm

\bpr2.7
An infinite \pn\ \sr\ is in\dcm.
\epr

\proof
Suppose for \cd ion that the infinite \pn\ \sr~$X$ is \dcm. Then \te\ \sr s
$(X_i,\qu\en i_0,\qu\en i_1, e\en i_0,e\en i_1)$, $i=1,2$, \st the direct
product $X_1\t X_2$ is \is c to~$X$. Since $X$~is principal, $X_1\t X_2$ is
\pn\ with \Gn\ $(e\en 1_1,e\en2_1)$ by Lemma \rf{l2.5}\,(ii) and \er{2.3}.
In view of \E\Pr\ \rf{p2.2}\,(v), $X_1$~(resp.~$X_2$) is \pn\ with
\Gn~$e\en1_1$ (resp.~$e\en2_1)$. Since $X$~is infinite, $X_1\t X_2$ is
infinite by \E\Pr\ \rfa1{p4.27}\,A\,(i). \Mo by Theorem \rfa2{t3.25},
$X_1$~or~$X_2$ is infinite. Suppose $X_1$ infinite. \If that $(X_1,\qu\en1_0,
e\en1_0)$ is an infinite \pn\ monoid with \Gn~$e\en1_1$. From Theorem
\rfa2{t1.25} we infer that this monoid is a \pn\ \PM. Since $X_1\t X_2$ is
\pn\ with \Gn\ $(e\en 1_1,e\en2_1)$, \te s $n\in\N$ \st $(e\en 1_0,e\en2_1)
=\break n\dquz(e\en 1_1,e\en2_1)\nde2.6 = (n\dqun1 e\en1_1,n\dqun2 e\en2_1)$, hence
$n\dqun1 e\en1_1=e\en1_0$. By Theorem \rfa2{t1.25} the \sq\ of \IT s of the
\Gn~$e\en1_1$ of~$X_1$ is in\jc. Hence $n=0$ by \era2{2.3}\,I0. \csq,
$(e\en 1_0,e\en2_1)=0\dquz (e\en 1_1,e\en2_1)=(e\en 1_0,e\en2_0)$.
A~\cd ion, since $e\en2_1 \ne e\en2_0$. If $X_2$~is infinite, we replace
$(e\en 1_1,e\en2_0)$ by $(e\en 1_0,e\en2_1)$ and we find $e\en1_1 \ne e\en1_2$.
A~\cd ion.
\endproof

We next consider the problem of the ``indecomposability'' of a \ti{finite}
\pn\ c-\sr.

\bdf2.8
The \nme\ of a \ti{finite} monoid (resp.\ \sr)~$X$ is said to be the
\ti{order}\index{order of a finite monoid} of the monoid (resp.\ \sr)~$X$ and is denoted by $\#(X)$ or $|X|$.
\edf

We first look for \sft\ \cn s for a finite \pn\ c-\sr~$X$ to be in\dcm, \ev
tly, for necessary \cn s for~$X$ to be \dcm. We observe that if $X$~is a \dcm\
finite \sr\ and $X_1,X_2$ are as in \E\df\ \rf{d2.4}, then $X$ and~$X_1\t X_2$
are \ep. By Corollary \rfa2{c3.11}, $X_1\t X_2$ is finite and $|X_1\t X_2|=
|X|$. Since $X_i=\pi_1(X_1,X_2)$, $i=1,2$, $\pi_1,\pi_2$ are \ti{sur\jc} maps
and $X_1,X_2$ are finite by \E\Pr\ \rfa1{p4.27}\,(B)\,(i). \Mo $|X_1\t X_2|
=|X_1|\cdot|X_2|$ by Theorem \rfa2{t3.25}. Notice that $|X_i|\ge2$, $i=1,2$.
\If that $|X|=|X_1|\cdot|X_2|$ is a \ti{\cme}\nm\ (see \E\df\ \rfa3{d9.8}).
Thus we have proved

\blm2.9
A finite \sr\ of prime order is in\dcm.
\elm

We now consider the case where the order of $X$ is a \cme \nm.

\blm2.10
Let $X$ be a finite \dcm\ \pn\ c-\sr\ of order~$n$. Then \te\ $k,l\ge2$ \st
$kl=n$ and $\N_k\t\N_l$ is a \pn\ c-\sr.
\elm

\proof
From the discussion preceding Lemma \rf{l2.9} we know that if $X_1$ and~$X_2$
are as in \E\df\ \rf{d2.4}, then $X_1,X_2$ are finite and $n:=|X|=|X_1|\cdot
|X_2|$, $|X_1|,|X_2|\ge2$. \Mo by Lemma \rf{l2.5}\,(ii), $X_1\t X_2$ is a finite
\pn\ c-\sr, and by \E\Pr\ \rf{p2.2}\,(iv),\,(v), $X_1$~and~$X_2$ are also
finite \pn\ c-\sr s. In view of Theorem \rf{t1.28}, $X_1$ (resp.~$X_2$) is
\is c to the \sr\ $\N_k$ (resp.~$\N_l$) with $k:=|X_1|$ and $l:=|X_2|$.
Observe that $k,l\ge2$ and $kl=|X_1|\cdot|X_2|=n$.

Since $X_1$ (resp.~$X_2$) and $\N_k$ (resp.\ $\N_l$) are \is c, $X_1\t X_2$
and $\N_k\t \N_l$ are also \is c by \E\Pr\ \rf{p2.2}\,(iii). \If by Lemma
\rf{l2.5}\,(ii) that $\N_k\t \N_l$ is \pn.
\endproof

By \df\ the \sr\ $\N_k\t\N_l$, $k,l\ge2$, is \pn\ if the additive monoid
of $\N_k\t\N_l$ is \pn\ with \Gn~$e_1$, where $e_1$ is the \nel\ of the \mlv\
monoid. Since $1$~is the \nel\ of the \mlc\ in $N_k,\N_l$, $e_1=(1,1)$
by~\er{2.3}. We now investigate the set of \IT s of $(1,1)$ in $\N_k\t\N_l$.
We denote by~$\qu_0$ the \ad\ in $\N_k\t\N_l$ defined in \er{2.1} where
$\qu^{(1)}_0:=+_k$ and $\qu^{(2)}_0:=+_l$, by $e_0:=(0,0)$ the \nel\ of the \ad\
in $\N_k\t\N_l$, and by $I((1,1))$ the set of \IT s of~$(1,1)$.

By Lemma \rfa2{l1.21}, $I((1,1))$ is a \sbm\ of $(\N_k\t\N_l,\qu_0,(0,0))$.
Hence\break $(I((1,1)),{\qu_0},(0,0))$ is a finite \Cm. By Lemma \rf{l1.22}, the order
of this monoid is equal to the least \el\ of the set $\{n\in\Na: n\dquz(1,1)
=(0,0)\}$.

\bdf2.11
If $(Y,\qu,e)$ is a finite \Cm\ and $a\in Y$, then the order of the monoid
$(I(a),\qu,e)$ is called the \ti{order} of~$a$. \E\Ip the order of~$e$ is~$1$.
\edf

\blm2.12
If $(Y,\qu,e)$ is a finite \Cm, $a\in Y\sms e$, and $m\in\Na$ is the order
of~$a$, then
\bea2.7
&n\dqu a\ne e,\q n\in(0,m),\\
&m\dqu a=e, \lb{2.8} \\
&(n+qm)\dqu a=n\dqu a \q\hbox{\fa}n,q\in\N. \lb{2.9}
\e
Conversely, if $m\in \Na$ \sf ies \er{2.7}, \er{2.8}, then $m$ is the order
of~$a$.
\elm

\proof
We apply Lemma \rf{l1.22} with $(X,\qu,e):=(Y,\qu,e)$, $a:=a$. We find by~(i)
the \ex\ of a unique \el\ of $\Na\sms1$ \sf ying \er{1.43}. For convenience we
denote this \el\ by~$\wt m$. From~(vi) we infer that $\zo0,{\wt m} $ and $I(a)$
are \ep. By Corollary \rfa2{c3.11} we find $\#(I(a)) = \#([0,\wt m)) \nda23.5 =
\wt m$. Hence \er{1.43} and Lemma \rf{l1.22}\,(iii) hold with $m$ replaced by
$\wt m=\#(I(a))$. But $\#(A)=m$ by \E\df\ \rf{d2.11}. \E\Tf \er{1.43} and
Lemma \rf{l1.22}\,(iii) hold with $m:=\#(A)$. Then \er{2.7}, \er{2.8} follow
from \er{1.43}. We now prove \er{2.9}. Let $n,q\in\N$, then $(n+qm)\qu a
\nad{\era2{2.3}\,\rm I2} = (n\dqu a) \qu (qm\dqu a) \nad{\era2{2.3}\,\rm I3} =
(n\dqu a)\qu\break (qm\dqu a) \nda22.8 = (n\dqu a)\qu ((q\dpl m)\qu a)
\nad{\era2{2.3}\,\rm I3} = (n\dqu a)\qu(q\dqu(m\qu a))
= (n\dqu a) \qu (q\qu e)\break \nad{\era2{2.3}\,\rm I4} = (n\dqu a)\qu e \nda21.7 = n\dqu a$.

Let $m\in\Na$ \sf y \er{2.7}, \er{2.8}. We apply Lemma \rf{l1.22} with
$(X,\qu,e)=(Y,\qu,e)$, $a:=a$, and use the \uq\ part of~(i) to conclude by~(v)
that $I(a)$ is \ep\ to~$\zo0,m $, hence $m\nda23.5 = \#(\zo0,m )=\#(I(a))$ by
Corollary \rfa2{c3.11}. Hence $m$ is the order of~$a$.
\endproof

\bxs2.13
The \IT s of $(1,1)$ in $\N_k\t\N_l$ are:

\hph i,ii, $k=l=2$:
\[
\substack{\scriptstyle n=0\\\textstyle (0,0)}\
\substack{\scriptstyle 1\\\textstyle (1,1)}\
\substack{\scriptstyle 2\\\textstyle (0,0)}\
\substack{\scriptstyle 3\\\textstyle (1,1)} \q m=2, \ kl=4.
\]

\hph ii,i, $k=2$, $l=3$:
\[
(0,0) (1,1) (0,2) (1,0) (0,1) (1,2) (0,0) \q m=6,\ kl=6.
\]

\hph iii,, $k=4$, $l=6$:
\[
(0,0) (1,1) (2,2) (3,3) (0,4) (1,5) (2,0) (3,1) (0,2) (1,3) (2,4) (3,5) (0,0)
\q m=12,\ kl=2.
\]
In cases (i) and (iii) $(1,1)$ is \ti{not\/} a \Gn\ of $\N_k\t \N_l$.
\exs

\blm2.14
Let $k,l\in\Na\sms1$. Then the order of $(1,1)$ in $(\N_k\t\N_l,\qu_0,
(0,0))$ is equal to $\lcm(k,l)=\frac{kl}{\gcd(k,l)}$, where
$\lcm(k,k):=k$ and $\gcd(k,k):=k$.
\elm

\proof
\E\fa $n\in\N$, $n\dquz (1,1)\nde2.6 = (n\dpls k1,n\dpls l 1)$. Then
$n\dpls k1=0$ iff $n=qk$ \fs $q\in\N$ by Lemma \rf{l2.12}. Similarly
$n\dpls l1=0$ iff $n=\wt q l$ \fs $\wt q\in\N$. \E\Tf $n\dquz(1,1)=(0,0)$ iff
$n$~is either equal to~$0$ or is a common \ml e of~$k$ and~$l$. Thus the order
of $(1,1)$ is the least common \ml e of $k$ and~$l$ (see
\era3{1.21}). Recalling \era3{1.25} we have \fa $k,l\ge2$, $k\ne l$, $\lcm(k,l)
\cdot\gcd(k,l)=kl$.
\endproof

\bco2.15
Let $k,l\in\Na\sms1$. Then the monoid $(\N_k\t\N_l,\qu_0,(0,0))$ is \pn\ with
\Gn\ $(1,1)$ iff
\beq2.10
k\ne l \qh{and }\gcd(k,l)=1.
\e
\eco

\proof
Since $I(1,1)\sbs \N_k\t\N_l$ and $|\N_k\t\N_l|\nda26.28b = |N_k|\,|\N_l|
\nda23.5 =kl$,  we have $I(1,1)=\N_k\t\N_l$ iff $I(1,1)=kl$ iff $\gcd(k,l)=1$
by Lemma \rf{l2.14}.
\endproof

Combining Corollary \rf{c2.15} with Lemma \rf{l2.10}, we find that a \sft\
\cn\ for a finite \pn\ c-\sr\ of order~$n$ to be in\dcm\ is: ``if $n=kl$
with $k,l\ge2$, then $\gcd(k,l)>1$''. In the next lemma we show that this
\cn\ is \sf ied whenever $n$ is equal to $p^m$, $p$~prime and $m\ge2$.

\blm2.16
Let $q$ be a \Pn\ and let $m\in\N$.  Then the set of divisors of~$q^m$ is
the set $\{q^i\in\Na: 0\le i\le m\}$.
\elm

\proof
Set $A_m:=\{q^i\in\Na: 0\le i\le m\}$ and let $B_m$ denote the set of
divisors of~$q^m$, \fe $m\in\N$.

$A_m\sbs B_m$: Let $m\in\N$ and let $i\in[0,m]$. Let $j\in\N$ be \st $m=i+j$.
Then $q^m=q^{i+j}\nda22.25 = q^i\cdot q^j$. Since $q^0\nda22.23 = 1$,
$q^i$ divides $q^m$, and $q^i\in B_m$.

$B_m\sbs A_m$: By induction on $m$. The case $m=0$ is clear. Suppose the
statement true for $m\in\N$. Let $x\in\Na$ be a divisor of $q^{m+1}=q\cdot
q^m$. Then either $\gcd(x,q)=1$ or $\gcd(x,q)=q$, since $q$~is prime. If
$\gcd(x,q)=1$, then $x$ divides $q^m$ by Euclid's lemma \rfa3{l9.15}. Hence
by the induction hypothesis $x\in A_m\sbs A_{m+1}$. If $\gcd(x,q)=q$, then
$q$~divides~$x$, that is, \te s $r\in\Na$ \st $x=qr$. Since
$x$~divides~$q^{m+1}$, \te s $s\in\Na$ \st $q^{m+1}=xs$, hence $q\cdot q^m
=q^{m+1}=xs=qrs$.

Since $q\ne0$, we obtain $q^m=r\cdot s$ by \cnc ity. \E\Tf $r$~divides $q^m$
and by the induction hypothesis $r\in A_m$. Hence $x=qr\in A_{m+1}$. Thus
$B_{m+1}\sbs A_{m+1}$ for all $m\in\N$.

\If that the set $\{m\in\N: B_m\sbs A_m\}$ is \iv\ in~$\N$. \E\Tf $B_m\sbs A_m$
\fa $m\in\N$. \csq, $A_m=B_m$ \fa $m\in\N$.
\endproof

\Wanp present the first main result of this section.

\bth2.17
Let $X$ be a finite \pn\ c-\sr\ of order $n\ge2$. Then $X$ is in\dcm\ iff
$n=p^m$ \fs $p$~prime and $m\in\Na$.
\eth

\proof\

\ti{Sufficiency}: Let $p$ be prime and $m\in\Na$. If $m=1$, then $X$~is
in\dcm\ by Lemma \rf{l2.9}. Let $m\ge2$. Suppose for \cd ion that $X$~is
\dcm. Then by Lemma \rf{l2.10}, \te\ $k,l\ge2$ \st $kl=p^m$ and $\N_k\t \N_l$
is a~\pn\ c-\sr. From Corollary \rf{c2.15} we infer that $\gcd(k,l)=1$. \E\oh
since $k$~is a~divisor of~$p^m$ and $k\ge2$, it follows from Lemma \rf{l2.16}
that $k=p^i$ \fs $i\in\oz1,m $. Hence $p$~divides~$k$. Similarly,
$p$~divides~$l$. \E\Tf $p$ is a~common divisor of $k$ and~$l$, hence
$\gcd(k,l)\ge p>1$. A~\cd ion.

\ti{Necessity}: Suppose that $n$ is not a power of a \Pn. By Lemma
\rfa3{l9.9} there is a \Pn\ $p$ dividing~$n$. Let $A:=\{m\in\Na: p^m| n\}$. Then
$1\in A$ since $p|n$. \Mo $A$~is bounded above. Indeed, define a~map $c:\N\to
\N$ by setting $c_m:=p^m-1$, $m\in\N$. Then $c_0=0$ and $c$~is strictly
in\cre\ (use \era2{2.33}). By Lemma \rfa2{l1.32} \te s $\ov m\in\N$ \st
$n\in\zo c_{\bar m},c_{\bar m+1} $, hence $n< c_{\bar m+1}$. Thus if $m\in A$,
then $p^m|n$. Hence $p^m\le n$ by \era3{8.57}, $c_m=p^m-1<p^m\le n<c_{\bar m
+1}$. \E\Tf $m<\ov m+1$, otherwise $c_{\bar m+1}\le c_m< c_{\bar m+1}$,
a~\cd ion, hence $A$ is bounded above. By Theorem \rfa1{t3.39},
$A$~has a greatest \el\ which we denote by~$m_p$. Thus $p^{m_p}|n$ and
$p^{m_p+1}\nmid n$. Set $k:=p^{m_p}$ and $l:=\frac nk$. By \as\ $l\ne1$. \Mo
$p\nmid l$, otherwise $l=pr$, $r\in\Na\sms1$, and $n=p^{m_p}\cdot pr =
p^{m_p+1}\cdot r$, \cd ing $p^{m_p+1}\nmid n$. \If from Lemma \rf{l2.16} that
$1$~is the only common divisor of $k$ and~$l$, i.e.\ $\gcd(k,l)=1$. In view of
Corollary \rf{c2.15}, the monoid $(\N_k\t \N_l,\qu_0,(0,0))$ is \pn\ with
\Gn~$(1,1)$, hence the finite c-\sr\ $\N_k\t\N_l$ is \pn. Since the order of
$\N_k\t\N_l$ is equal to~$kl$ by \era2{6.28b}, $\N_k\t\N_l$ is ring-\is c to $\N_{kl}$ by
Theorem \rf{t1.28}. By the same theorem $X$~is ring-\is c to~$\N_{kl}$ since
$n=kl$, hence by \er{1.4n} and Lemma \rfa2{l1.8} $X$~is \is c to $\N_k\t\N_l$.
Thus $X$~is \dcm.
\endproof

We saw that if $k,l\ge2$, then $\N_k\t\N_l$, the product of the \sr s $\N_k$
and~$\N_l$, is a finite c-\sr\ of order~$kl$ by \E\Pr\ \rf{p2.2} and Theorem
\rfa2{t3.25}. If in addition $\gcd(k,l)=1$, then $\N_k\t\N_l$ is \pn\ with
\Gn~$(1,1)$ by Corollary \rf{c2.15}. Then by Theorem \rf{t1.28}, $\N_k\t\N_l$
and $\N_{kl}$ are \is c and the \sq\ of \IT s of~$(1,1)$ in the additive
monoid, \rt ed to $\zo0,kl $ is the \ti{only} \is sm from $\N_{kl}$ onto
$\N_k\t\N_l$. The aim of the remaining part of this section is to present
an \ext\ of this result to a finite product of \sr s
$\{\N_{n_i}\}_{i\in[1,N]}$, $n_i\in\Na\sms1$, $i\in[1,N]$, $N\ge2$.

We first introduce the notion of
\ti{finite product\/} of sets (resp.\ monoids, \sr s). Let $N\in\Na\sms1$ and\index{finite product of sets, monoids, semirings}
let $I:=[1,N]$ ($:=\{i\in\N: 1\le i\le N\}$), and let $\zb XiI$ be a finite
\sq\ of nonempty sets (resp.\ monoids, \sr s). Set\glossary{$\prodl_{i\in I}X_i$}
\beq2.11
\prod_{i\in I} X_i:= \Bigl\{f:I\to\bigcup_{i\in I}X_i\mid f(i)\in X_i \hbox{
\fa} i\in I\Bigr\}.
\e

If $(X_i,\qn i,e\en i)_{i\in I}$ are monoids, set
\bga2.12
(f\qu g)(i):= f(i)\qn i g(i) \qh{\fa} i\in I,\ f,g\in \prod_{i\in I}X_i,\\
e(i):=e\en i \qh{\fa} i\in I. \lb{2.13}
\e

If $(X_i,\qn i_0, \qn i_1,e\en i_0, e\en i_1)_{i\in I}$ are \sr s, set
\bga2.14
(f \qu_j g)(i):= f(i)\qn i_j g(i) \qh{\fa} i\in I,\ j\in\{0,1\}, \hbox{ and }
f,g\in \prod_{i\in I}X_i, \\
e_j(i):=e\en i_j \qh{\fa} i\in I,\ j\in\{0,1\}. \lb{2.15}
\e

\bpr2.18
Let $\prodl_{i\in I}X_i$ be as in \er{2.11}. Then $\prodl_{i\in I}$ is a \nf
\ and
\beq2.16
\#\Bigl(\prod_{i\in I} X_i\Bigr) = \prod_{i\in I} \#(X_i).
\e
\epr

\bex2.19 \

\hph i,ii, Prove \E\Pr\ \rf{p2.18}.

\ti{Hint\/}:
Use induction on $N\in\Na\sms1$. Observe that $\bigl(\prodl_{i\in[1,N]}X_i
\bigr)\t X_{N+1}$ and $\prodl_{i\in[1,N+1]}X_i$ are \ep.

\hph ii,i, Let $X_1,X_2,X_3$ be sets. Show $X_1\t X_2 \ax X_2\t X_1$ (see
\E\df\ \rfa1{d4.1}), $X_1\t X_2\t X_3\ax (X_1\t X_2)\t X_3 \ax X_1\t (X_2\t
X_3)$, $(X_1\t X_2)\t(X_3\t X_4)\ax X_1\t(X_2\t X_3)\t X_4$.

\hph iii,, Let $\zb Xi{[1,n]}$, $n\in\Na\sms1$, be sets. Let $\pi:[1,n]\to
[1,n]$ be bi\jc. Show that $\prodl_{i=1}^n X_i \ax \prodl_{i=1}^n X_{\pi(i)}$.
\eex

\bpr2.20
Let $N\in\Na\sms1$ and let $I:=[1,N]$.

\hph i,i, Let $\zb XiI$ be a family of monoids $($resp.\ \Cm s$)$. Let
$\prodl_{i\in I}X_i$, $\qu$, $e$ be as in \er{2.11}, \er{2.12} and \er{2.13},
\rsp. Then $\bigl(\prodl_{i\in I}X_i,{\qu},e\bigr)$ is a monoid $($resp.\ \Cm$)$, and
the \pj s $\pi_j: \prodl_{i\in I}X_i\to X_j$, $j\in I$, defined by
\beq2.17
\pi_j(x) := x(j), \q j\in I,\ x\in\prod_{i\in I}X_i
\e
are sur\jc\ monoid-\hm sms.

\hph ii,, Let $\zb XiI$ be a finite family of \sr s $($resp.\ c-\sr s$)$. Let
$\prodl_{i\in I}X_i,\qu_j,e_j$, $j\in\{0,1\}$ be as in \er{2.11}, \er{2.14},
\er{2.15}, \rsp. Then $(\prodl_{i\in I}X_i,\qu_0,\qu_1,e_0,e_1)$ is a \sr\
$($resp.\ c-\sr$)$, and the \pj s $\pi_j$, $j\in I$, defined in \er{2.17} are
sur\jc\ ring-\hm sms.
\epr

\bex2.21
Prove \E\Pr\ \rf{p2.20}.
\eex

\bpr2.22
Let $N\in\Na\sms1$, $I:=[1,N]$ and let $\zb XiI$ be a finite family  of
monoids $($resp.\ \sr s$)$. Let $Y$ be a monoid $($resp.\ \sr$)$ and let
$\vf:Y \to \prodl_{i\in I}X_i$. Then $\vf$~is a monoid- $($resp.\ \sr-$)$
\hm sm iff the maps $\pi_j\circ\vf : Y\to X_j$, $j\in I$, are monoid-
$($resp.\ \sr-$)$ \hm sms.
\[
\xymatrix{Y\ar[rr]^\vf \ar[drr]_{\pi_j\circ\vf} && \prodl_{i\in I}X_i\ar[d]^{\pi_j} \\
&&X_j} \quad j\in I.
\]
\epr

\proof \ti{Only if\/}: directly follows from \E\Pr\ \rf{p2.20}\,(ii), Lemma \rfa2{l1.8}
and Exercise \rf{ex1.8n}.

\ti{If\/}: it is \sft\ to consider the case of \sr s. Let $+$ and~$\cdot$
denote the \ad\ and \mlc\ in~$Y$ with \crs\ \nel s $0$~and~$1$. \E\fa $j\in I$
we have $\pi_j\circ \vf(0)=e^{(j)}_0$ by \as. Hence $\vf(0)=e_0$ by \er{2.15}.
Similarly $\vf(1)=e_1$.

Let $y_1,y_2\in Y$. Then, \fa $j\in I$ we have $\pi_j\circ \vf(y_1+y_2) =
(\pi_j\circ \vf)(y_1) \qu^{(j)}_0 (\pi_j\circ\vf)(y_2)$ by \as. Hence
$\vf(y_1+y_2) = \vf(y_1)\qu_0 \vf(y_2)$ by \er{2.14}. Similarly
$\vf(y_1\cdot y_2) = \vf(y_1)\qu_1 \vf(y_2)$.
\endproof

We recall that the maps $\F_n:\N \to\N_n$, $n\in\Na\sms1$, introduced in
\er{1.3} are sur\jc\ by \er{1.4} and \sf y \fa $x,y\in\N$:
\beq2.18
\left\{\ga
\F_n(x+y) \nde1.7 = \F_n(\F_n(x)+\F_n(y)) \nde1.8 = \F_n(x)+_n \F_n(y),\\
\F_n(x\cdot y) \nde3.26 = \F_n(\F_n(x)\cdot \F_n(y)) \nde1.33 =
\F_n(x)\cdot_n \F_n(y).
\ega\right.
\e

\Mo $\F_n(0)=0$, $\F_n(1)=1$ by \er{1.6}. Hence
\[
\F_n : (\N,+,\cdot,0,1) \to (\N_n,+_n,\cdot_n,0,1), \q n\in\Na\sms1,
\]
are sur\jc\ ring-\hm sms.

Now let $k,m\in\Na\sms1$ with $k<m$.
\beq2.19
\gathered
\hbox{Let $\F_{k\lea m}$ denote the \rt ion of $\F_k$}\\
\hbox{to $\N_m$ with range contained in $\N_k$.}
\endgathered
\e
Then $\F_{k\lea m} : \N_m\to \N_k$ is sur\jc\ by \er{1.5}, $\F_{k\lea
m}(0)=0$ and $\F_{k\lea m}(1)=1$. In general, $\F_{k\lea m}$ is \ti{not\/}
a~ring-\hm sm. Indeed, if $m:=5$ and $k:=3$, then $\F_{3\lea 5}(4+_51)
=\F_{3\lea 5}(0)=0$ and $\F_{3\lea 5}(4)+_3 \F_{3\lea 5}(1) = 1+_31=2 \ne0$.
However, if $k$ divides~$m$, then $\F_{k\lea m}$ is a ring-\hm sm.

\bpr2.23
Let $k,m\in\Na\sms1$ be \st $k$ divides~$m$. Then the map
\[
\F_{k\lea m}:(N_m,+_m,\cdot_m,0,1) \to (N_k,+_k,\cdot_k,0,1)
\]
defined  in \er{2.19} is a sur\jc\ ring-\hm sm.
\epr

\proof
Let $l:=\frac mk$, i.e., $kl=m$. Let $x,y\in\N_m$.
\beq2.20
\F_{k\lea m}(x+_my) = \F_{k\lea m}(x) +_k \F_{k\lea m}(y):
\e
Since $x,y$ and $x+_my$ belong to $\N_m$, we have $\F_{k\lea m}(x)=\F_k(x)$,
$\F_{k\lea m}(y)=\F_k(y)$ and $\F_{k\lea m}(x+_my)=\F_k(x+_my)$. In view of
\er{1.8}, \er{1.2}, \er{1.3}, \te s $p\in\N$ \st $x+y=pm+(x+_my)$. But $pm=
p(kl)=p(lk)=(pl)k$. Hence $x+y=(pl)k+(x+_my)$. Thus $\F_{k\lea m}(x+_my)=
\F_k((x+_my)+(pl)k)\nde1.4  = \F_k(x+y) \nde2.18 = \F_k(x)+_k\F_k(y)$, which implies
\er{2.20}.
\beq2.21
\F_{k\lea m}(x\cdot_my) = \F_{k\lea m}(x)\cdot_k\F_{k\lea m}(y):
\e
Since $x,y,x\cdot_m y$ belong to $\N_m$, we have $\F_{k\lea m}(x\cdot_m y)
= \F_k(x\cdot _m y)$. In view of \er{1.33}, \er{1.2}, \er{1.3}, \te s $p\in\N$
\st $xy=pm+(x\cdot_m y)$. But as above $pm=(pl)k$. Hence $\F_{k\lea m}
(x\cdot_m y)\nde1.4 = \F_k((x\cdot_my)+(pl)k) = \F_k
(x y)\nde2.18 = \F_k(x)\cdot_k \F_k(y)$, which implies \er{2.21}.
\Mo $\F_{k\lea m}(0)=0$ and $\F_{k\lea m}(1)=1$.
\endproof

\Wanp state the second main result of this section.

\bth2.24
Let $N\in \Na\sms1$, $n_j\in\Na\sms1$, $j\in[1,N]$, and
\beq2.22
M:= \prod_{j\in[1,N]}n_j.
\e
Let $\Phi:\zo0,M \to \prodl_{j\in[1,N]}\zo0,n_j $ be the map defined by
\beq2.23
\Phi(x)(j) := \F_{n_j}(x), \q x\in\zo0,M ,\ j\in[1,N],
\e
where $\F_{n_j}$ is defined in \er{1.3}. Then $\Phi$ is a ring-\hm sm from
$(\N_M,+_M,\cdot_M,0,1)$ into the \sr\
$\bigl(\prodl_{j\in[1,N]}\N_{n_j},\qu_0,\qu_1,e_0,e_1\bigr)$ where
$\qu_0,\qu_1,e_0,e_1$ are defined in \er{2.14}, \er{2.15}. \Mo $\Phi$~is a
ring-\is sm iff
\beq2.24
\gcd(n_j,n_k)=1\qh{\fa $j,k\in[1,N]$ \st }j\ne k.
\e
\eth

\proof
\ti{The map $\Phi$ is a ring-\hm sm}:

We apply \E\Pr\ \rf{2.22} with $Y$ being the \sr\ $(\N_M,+_M,\cdot_M,0,1)$,
$X_i$~being the \sr s $(\N_{n_i}, +_{n_i},\cdot_{n_i},0,1)$, $i\in I$, and
$\vf:=\Phi$ defined in \er{2.22}. Since \fa $j\in I$, $n_j$~divides
$M=\prodl_{i\in I}n_i$ by Lemma \rfa3{l9.10}, the maps $\vf\circ \pi_j =
\F_{n_j\lea M}$, $j\in I$, are ring-\hm sm by \E\Pr\ \rf{p2.23}. \If from
\E\Pr\ \rf{p2.22} that $\Phi$~is a ring-\hm sm.

\smallskip
\ti{$\Phi$ is a ring-\hm sm iff \er{2.24} holds}: Since $\Phi$~is a ring-\hm
sm it is \sft\ to prove that $\Phi$ is a bi\jn. We denote by~$+$
(resp.~$\cdot$) the \ad\ (resp.\ \mlc) in the \sr\ $\prodl_{i\in I}\N_{n_i}$
and by $\bf0$ (resp.~$\bf1$) the \crs\ \nel s. We claim
\beq2.25
\Phi(n)= n\dpl{\bf1}  \qh{\fa} n\in[0,M],
\e
where the \RHS\ of \er{2.25} is the $n$-th \IT\ of~$\bf1$ in the monoid
$\bigl(\prodl_{i\in I}\N_{n_i},+,{\bf0}\bigr)$.  Indeed, \er{2.25} is \ev t to
\beq2.26
\pi_j(\Phi(n)) = \pi_j(n\dpl{\bf1}) \qh{\fa} j\in[0,N].
\e
\Mo \fa $j\in\N$ and $n\in\zo0,M $, we have $\pi_j(\Phi(n))=\F_{n_j}(n) \nde1.11 =
n\dpls{n_j} {1}$ the  $n$-th \IT\ of~$1$ in $(\N_{n_j},+_{n_j},0)$. But
$1=\pi_j ({\bf1})$, hence $n\dpls{n_j}1 = n\dpls{n_j}\pi_j({\bf1})$. Since
the \pj\ $\pi_j:\bigl(\prodl_{i\in[1,N]}\N_{n_i},+,{\bf0}\bigr) \to (\N_{n_j},
+_{n_j},0)$ is a monoid-\hm sm by \E\Pr\ \rf{p2.20}(ii), we have $n\dpls
{n_j}\pi_j({\bf1})=\pi_j(n\dpl{\bf1})$ by Lemma \rfa2{l1.22}. \E\Tf
\er{2.26}, hence \er{2.25} hold.

We denote by $I({\bf1})$ the set of \IT s of~$\bf1$ in the monoid
$\bigl(\prodl_{i\in[1,N]}\N_{n_i},+,{\bf0}\bigr)$. Since this monoid is a finite
\Cm, we can apply Lemma \rf{l1.22} with $n:=\#\bigl(\prodl_{i\in[1,N]}\N_{n_i}
\bigr) \nde2.16 = \prodl_{i\in[1,N]}n_i \nde2.22 = M$, and with $a:={\bf1}$.
Let $\vf_{\bf1}$ denote the \sq\ of \IT s of~$\bf1$ in this monoid. \If from
\er{2.25} that \ti{the map $\Phi$ is the \rt ion to $\zo0,M $ of~$\vf_{\bf1}$}.
In view of Lemma \rf{l1.22}\,(i) \te s exactly one $m\in\Na$ \st $\vf_{\bf1}
(m)=0$ and $\vf_{\bf1}(k)\ne0$, $k\in(0,m)$.

We next show that $m$ is the \ti{least common multiple of} $\cM:= \bcl_{j\in
[1,N]}\{n_j\}$ in~$\Na$. We observe that if $l\in\Na$, then $l\dpl {\bf1}
={\bf0}$ iff $0=\pi_j({\bf0})= \pi_j(l\dpl{\bf1}) \nda21.45 = l\dpls{n_j}
\pi_j({\bf1})= l\dpls{n_j}1$ in the monoid $(\N_{n_j},+_{n_j},0)$ \fa $j\in
[1,N]$. Given $j\in[1,N]$ and $l\in\Na$, we claim that $l\dpls{n_j}1=0$ iff
$l$~is a \ml e of~$n_j$. Indeed, by \er{1.11} with $n:=n_j$ and $m:=l$, we
find $l\dpls{n_j}1 = \phi_{n_j}(l)$. But in view of the \df\ of~$\phi_{n_j}$
in \er{1.3}, we find that $\phi_{n_j}(l)=0$ iff $l$~is a \ml e of~$n_j$. Thus
the claim is proved. \If that if $l\in\Na$ then $l\dpl{\bf1}={\bf0}$ iff
$l$~is a common \ml e of~$\cM$. Since $m\dpl{\bf1}={\bf0}$, we infer that
$m$~is a~common \ml e of~$\cM$, and since $k\dpl{\bf1}\ne{\bf0}$ for
$k\in(0,m)$, we infer that $m=\lcm(\cM)$.

\csq, we find by Lemma \rf{l1.22}\,(vi) with $(X,\qu,e):=\bigl(\prodl_{i\in[1,N]}
\N_{n_j},+,{\bf0}\bigr)$, $n:=M$, $a:=\bf1$, that $\vf_{\bf1}:\zo0,\lcm(\cM) \to
I({\bf1})$ is \ti{bi\jc}. \E\Ip we have
\beq2.27
\lcm(\cM) = \#(I({\bf1})).
\e
Since $\Phi$ is the \rt ion to $\zo0,M $ of $\vf_{\bf1}$ we have
\beq2.28
R(\Phi) \sbs I({\bf1}),
\e
where $R(\Phi)$ denotes the range of $\Phi$.

We now suppose that \cn\ \er{2.24} holds. \If from Lemma \rf{l2.26} below that
$\lcm(\cM) = \prodl_{i\in[1,N]}n_i$, hence $\lcm(\cM)=M$ by \er{2.22}. \E\Tf
$\vf_{\bf1} : \zo0,M \to I({\bf1})$ is bi\jc, hence in\jc. Since $\Phi$ is the
\rt ion of~$\vf_{\bf1}$ to $\zo0,M $, we find that $\Phi$ is in\jc. Finally,
since $\Phi:\zo0,M \to \prodl_{j\in[1,N]}\zo0,n_j $ is in\jc\ and $\#(\zo0,M
)=M= \prodl_{j\in[1,N]}n_j = \prodl_{j\in[1,N]}\#(\zo0,n_j ) \nde2.16
= \#\Bigl(\prodl_{j\in[1,N]}\zo0,n_j \Bigr) \nde2.22 = M$, we infer from \E\Pr\
\rfa2{p3.13} that $\Phi$ is \ti{bi\jc}.

We now consider the case where \cn\ \er{2.24} is \ti{not\/} \sf ied, and
first treat the case $N:=2$. We claim that $\lcm(\cM)<n_1n_2$. Indeed, if
$n_1=n_2$, then $\lcm(\{n_1\}\cup\{n_2\})=\lcm\{n_1\}=n_1 < n_1n_2$ by
\era2{2.12}, \era2{2.19} since $n_1\in\Na$, $n_1=n_11$ and $1<n_2$. If
$n_1\ne n_2$, then $n_1n_2\nde1.8 = \lcm(n_1,n_2)\cdot \gcd(n_1,n_2)
\nda22.19 > \lcm(n_1,n_2)\cdot1 \nda22.18 = \lcm(n_1,n_2) = \lcm(\cM)$, since
$\gcd(n_1,n_2)>1$, and $\lcm(\cM)\in \Na$. Hence $\lcm(\cM)<n_1n_2$, which
completes the proof of the claim. \If from \er{2.27} that
$\#(I({\bf1}))<n_1n_2$. From \er{2.28} and \era2{3.13}, we infer $\#(R(\Phi))\le
\#({\bf1})$. Hence $\#(R(\Phi))<n_1n_2$ by \era1{3.11}. Since $R(\Phi)\sbs
\zo0,n_1 \t\zo0,n_2 $ and $\#(\zo0,n_1 \t\zo0,n_2 )=n_1n_2$ by \era2{6.28b},
we conclude from \era2{3.13}  that $R(\Phi)\ne\zo0,n_1 \t\zo0,n_2 $. Thus
$\Phi$~is not sur\jc, hence \ti{not\/} bi\jc.

Finally, we consider the  case $N>2$. Let $k,l\in[1,N]$ be \st $k\ne l$ and
$\gcd(n_k,n_l)>1$. We define $\wt\Phi:\zo0,n_kn_l \to \N_{n_k}\t\N_{n_l}$ by
setting $\wt\Phi(n):=(\F_{n_k}(n),\F_{n_l}(n))$, $n\in\zo0,n_kn_l $. From the
first part of the proof we have $\wt\Phi(n) = n
\mathrel{\lower2pt\hbox{$\stackrel{\lower1pt\hbox{$\scriptstyle\wt
+$}}\cdot$}} \wt{\bf1}$, $n\in\N$, where $\wt +$ denotes the \ad\ in the
additive monoid of $\N_{n_k}\t\N_{n_l}$ and $\wt{\bf1}:= (1,1) \in
\N_{n_k}\t\N_{n_l}$. Let $I(\wt{\bf1})$ denote the set of \IT s of~$\wt{\bf1}$
in this monoid. From \er{2.27} we obtain $\#(I(\wt{\bf1}))=\lcm(\{n_k\}\cup
\{n_l\})$. As above we have $\lcm(\{n_k\}\cup\{n_l\})<n_kn_l$ since
$\gcd (n_k,n_l)>1$. Hence
\beq2.29
\#(I(\wt{\bf1})) < n_kn_l.
\e
We now suppose for \cd ion that $\Phi$~is bi\jc. Let $(\a,\b)\in\N_{n_k}
\t\N_{n_l}$. Define $y\in\prodl_{j\in[1,N]}\N_{n_j}$ by setting $y(k):=\a$,
$y(l):=\b$ and $y(i):=0$ \fa $i\in[1,N]\sms{k,l}$. Since $\Phi$~is bi\jc\ \te
s one (and only one) $n\in\zo0,M $ \st $\Phi(n)=y$. \E\Ip
$\phi_{n_k}(n)=\a$ and $\phi_{n_l}(n)=\b$. Since $(\a,\b)$ is arbitrary in
$\N_{n_k}\t\N_{n_l}$, every \el\ of $\N_{n_k} \t\N_{n_l}$ is an \IT\
of~$\wt{\bf1}$ in the additive monoid of $\N_{n_k}\t\N_{n_l}$. Thus $\N_{n_k}
\t\N_{n_l}\sbs I(\wt{\bf1})$ hence by \era2{6.28b} we have $n_kn_l =
\#(\N_{n_k}\t\N_{n_l}) \le \#(I(\wt{\bf1}))$, \cd ing \er{2.29}.
\endproof

We recall that in the proof of Theorem \rf{t2.24} we used Lemma \rf{l2.26}.
This lemma requires some result which we establish in the next lemma.

\blm2.25
Let $N\in\Na\sms1$, let $a:[1,N]\to\Na$ be a finite \sq\ in~$\Na$ and let
$b\in\Na$. Then
\beq2.30
\Bigl(\LOR_{i\in[1,N]}a_i\Bigr) \land b = \LOR_{i\in[1,N]}(a_i \land b)
\e
where $\LOR_{i\in[1,N]}$ is the \cme \op\ on $(\Na,|)$ induced by the
\op~$\lor$ defined by $x\lor x:=x$, $x\lor y:=\lcm(x,y)$, $x,y\in\Na$, $x\ne
y$ and $x\land x:=x$, $x\land y:=\gcd(x,y)$, $x,y\in\Na$, $x\ne y$.
\elm

\proof
By \In\ on $N\in\Na\sms1$.

$N=2$: Let $a_1,a_2,b\in\Na$ be \st $a_1\land a_2=1$. Then $(a_1\lor a_2)\land b \nda33.18 = b\land
(a_1\lor a_2) \nda33.53 = (b\land a_1)\lor (b\land a_2) \nda33.18 = (a_1\land
b)\lor (a_2\land b)$, which is \er{2.30} for $N:=2$.

\ti{$N$ ``implies'' $N+1$}: We suppose that \er{2.30} holds for~$N$. Let
$a:[1,N+1]\to \Na$ and $b\in\Na$. Then $\bigl(\LOR_{i\in[1,N+1]}a_i\bigr)
\land b \nad*=\bigl(\bigl(\LOR_{i\in[1,N]}a_i\bigr)\lor a_{N+1}\bigr) \land b
\nad{**}= \bigl(\bigl(\LOR_{i\in[1,N]}a_i\bigr)\land b\bigr) \lor \break(a_{N+1}
\land b) \nad{*{*}*}=\bigl(\LOR_{i\in[1,N]}(a_i\land b)\bigr)\lor
(a_{N+1}\land b) \nad*=
\LOR_{i\in[1,N+1]}(a_i\land b)$. In $\nad*=$ the \df\ of $\LOR_{i\in[1,N+1]}$
is used, in $\nad{**}=$ the case $N=2$ is used and in $\nad{*{*}*}=$ the \iv\
\as\ is used.
\endproof

\blm2.26
Let $N\in\Na\sms1$ and let $a:[1,N]\to\Na\sms1$ be a finite \sq\ of pairwise
coprime \el s of~$\Na$, i.e.
\beq2.31
\gcd(a_i,a_j)=1 \qh{whenever} i\ne j,\ i,j\in[1,N].
\e
Set $\cM :=\bcl_{i\in[1,N]}\{a_i\}$. Then
\beq2.32
\lcm(\cM) = \prod_{i\in[1,N]}a_i.
\e
\elm

\proof
Observe that in view of \cn\ \er{2.31} the map $i\mt a_i$ from $[1,N]$ into
$\Na\sms1$ is in\jc. \Mo $\lcm(\cM)\nda33.23 = \sup\cM \nda33.9 =
\LOR_{i\in[1,N]} a_i$. We proceed by \In\ on $N\in\Na\sms1$.

$N=2$: $\LOR_{i\in[1,2]}a_i = a_1\lor a_2 = (a_1\lor a_2)\cdot 1 \nde2.31 =
(a_1\lor a_2)\cdot (a_1\land a_2) \nda31.25 = a_1a_2 =\prodl_{i\in[1,2]}a_i$.

\ti{$N$ implies $N+1$}: Suppose \er{2.32} holds. $\LOR_{i\in[1,N+1]}a_i \nad*= \bigl(
\LOR_{i\in[1,N]}a_i\bigr)\lor a_{N+1} \nad{**}= \bigl(\prodl_{i\in[1,N]}
a_i\bigr)\lor a_{N+1}$. In $\nad*=$ the \df\ of $\LOR_{i\in[1,N+1]}$ is used
and in $\nad{**}=$ the \In\ hypothesis is used (observe that if \er{2.31}
holds for $i,j\in[1,N+1]$, it also holds for $i,j\in[1,N]$). We claim that
$\bigl(\prodl_{i\in[1,N]}a_i\bigr) \land a_{N+1}=1$. Indeed, since
$\prodl_{i\in[1,N]}a_i = \LOR_{i\in[1,N]}a_i$, we have $\bigl(\prodl_{i\in[1,N]}
a_i\bigr)\land a_{N+1} = \bigl(\LOR_{i\in[1,N]}a_i\bigr)\land a_{N+1} \nde2.30 =
\LOR_{i\in[1,N]}(a_i \land a_{N+1}) \nde2.31 = \LOR_{i\in[1,N]}1$. Observe that
$1\lor x=x\lor 1 =x$ \fa $x\in\Na$, hence $(\Na,\lor,1)$ is an abelian monoid.
\E\Tf $\LOR_{i\in[1,N]}1\nda21.126 = \#([1,N] \mathrel{\lower4pt
\hbox{$\stackrel\lor\cdot$}}1 = N\mathrel{\lower4pt\hbox{$\stackrel\lor\cdot$}}1
\nad{\era2{2.3}\,\rm I4}=1$ in the \am\ $(\Na,\lor,1)$, which proves the claim. \If from \era3{1.25} that $\bigl(
\prodl_{i\in[1,N]}a_i\bigr)\lor a_{N+1} = \bigl(\prodl_{i\in[1,N]}a_i\bigr)
\cdot a_{N+1} = \prodl_{i\in[1,N+1]}a_i$. Thus $\LOR_{i\in[1,N+1]}a_i =
\prodl_{i\in[1,N+1]}a_i$.
\endproof

\brm2.27
If $x\in\N$ and $n\in\Na$, then $\F_n(x)$ is the unique \el\ of~$\zo0,n $
\sf ying $x=q\cdot n+\phi_n(x)$ \fs $q\in\N$ (see \er{1.2}, \er{1.3}). The
\nm\ $\F_n(x)$ is usually called the \ti{remainder of the Euclidean division
of~$x$ by~$n$}. Let $\zb mi{[1,N]}$, $N>1$, be a family of \el s of $\Na\sms1$
which are pairwise coprime, and let $\zb ai{[1,N]}$ be a family of \el s
of~$\N$ \sf ying $a_i<m_i$ \fa $i\in[1,N]$. Theorem \rf{t2.24} tells us that
\te s exactly one \nm\ $x\in\bigl[0,\prodl_{i\in[1,N]}m_i)$ \st the remainder
of the Euclidean division of~$x$ by~$m_i$ is $a_i$ \fa $i\in[1,N]$. This
result is known as the \ti{Chinese remainder theorem}.\index{Chinese remainder theorem} The proof of this
theorem given here is \ti{not\/} constructive. We invite the reader
to find in the literature historic and constructive aspects as well as
applications of this theorem.
\erm

Let $n\in\Na\sms1$. Then by the \fd\ theorem of \art, Theorem \rfa3{t9.22},
$n$~is either a power of a prime \nm\ or a \cme product of powers of distinct
prime \nm s, that is, $n=\prodl_{i\in[1,N]}n_i$ with $n_i=p_i^{m_i}$,
$p_i$~prime, $p_i\ne p_j$ whenever $i\ne j$, $m_i\in\Na$ \fa $i,j\in[1,N]$ \fs
$N\in\Na\sms1$. \If from Lemma \rf{l2.16} that $\gcd(m_i,m_j)=1$ whenever
$i\ne j$, $i,j\in[1,N]$. \E\Tf by Theorem \rf{t2.24} and Theorem \rf{2.17}
the \sr\ $(\N_n,+_n,\cdot_n,0,1)$ is \is c to a product of indecomposable
finite \pn\ c-\sr s.

More generally we have:

\bth2.27
Every finite \pn\ c-\sr\ is either indecomposable or \is c to a finite product
of indecomposable finite c-\sr s.
\eth

\bex2.28
Prove Theorem \rf{t2.27}.
\eex

\Wanp prove an \ext\ of Lemma \rfa3{l1.28}.

\blm2.30
Let $N\in\Na\sms1$, let $a:[1,N]\to \Na\sms1$ be \st \er{2.31} holds, and let
$c\in\Na\sms1$. If $a_i|c$ \fa $i\in[1,N]$, then $\prodl_{i\in[1,N]}a_i|c$.
\elm

\bex2.31
Prove Lemma \rf{l2.30}.
\eex

\bex2.32
Let $m\in\Na\sms1$ with \dr\ \er{1.44}. Find \cn s on the digits of~$m$
for~$m$ to be divisible by 132, 135, 315 and 2772 (see Remark \rf{r1.36}).
\eex

\newpage
\Subsubsection{Cyclic groups}\label{sss.fcg}

The additive monoid of a \pn\ c-\sr\ is by \df\ a \pn\ \Cm. In the first part
of this section we investigate \pp ies of \ti{finite} \pn\ \Cm s which are not
shared by infinite \pn\ ones.

Let $(X,\qu,e)$ be a \ti{finite} \pn\ \Cm, and let $a\in X\sms e$ be a \Gn\
of~$K$, i.e. $X=I(a)$. Then the \fw\ \pp ies hold.

\bit
\item[(A)] \E\fe $x\in X$, \te s \ooo \el\ $y\in X$ \st $x\qu y=e$. The \el\
$y\in X$ is called the \ti{inverse} of~$x$, and is sometimes denoted by $x\Inv$.
\eit

We recall that an infinite \pn\ \Cm\ is a \PM\ by Theorem \rfa2{t1.42}. Hence,
if $x\qu y=e$, then $x=y=e$.

\bit
\item[(B)] Every nontrivial \sbm\ of $X$ is \pn.
\eit

An example of a \sbm\ of a \pn\ infinite \Cm\ which is \ti{not\/} \pn\
is given in Exercise \rfa2{ex1.44} and Remark \rfa3{r1.15}.

\bit
\item[(C)] Unless $\#(X)=2$, $X$~possesses more than one \Gn.
\eit

Recall that by Theorem \rfa2{t1.42}, an infinite \pn\ \Cm\ possesses exactly
one \Gn.

One of the goals of this section is to \es\ a formula for the \nm\ of \Gn s of
a finite \pn\ \Cm.

We now prove \pp y (A).

If $(X,\qu,e):=(\N_n,+_n,0)$, $n\ge2$, then $X$ is a finite \pn\ \Cm\ and
\pp y (A) holds by \E\Pr\ \rf{p1.8}\,(i) and~(ii). More generally, every finite
\Cm\ (not necessarily \pn) possesses \pp y~(A).

\bpr3.1
Let $(X,\qu,e)$ be a \emph{finite} \Cm. Then, \fa $x,z\in X$ \te s \ooo\
$y\in X$ \st
\beq3.1
x\qu y=z.
\e
\epr

\proof
Let  $x\in X$ and let $\th_x:X\to X$ be the map defined by
\beq3.2
\th_x(y):=x\qu y\,(=y\qu x) \qh{\fa}y\in X.
\e
The map $\th_x$ is \ti{in\jc}. Indeed, if $u,v\in X$ \sf y $\th_x(u)=\th_x(v)$,
then by \cmt ity $u\qu x=v\qu x$ and by \cnc ity \era2{1.8} $u=v$. Since
$X$~is finite, $\th_x$~is also \ti{sur\jc} by \E\Pr\ \rfa2{p3.13} \era2{3.11}.
Hence \fe $x\in X$ \te s one (and only one) $y\in X$ \st $\th_x(y)=e$.
\endproof

\bdf3.2
Let $(M,\qu,e)$ be a monoid. An \el\ $x\in M$ is called \ti{invertible} if\index{invertible element}
\te s $y\in M$ \st
\beq3.3
x\qu y=y\qu x=e.
\e
The \el\ $y$ is the only \el\ of~$M$ \sf ying \er{3.3}. Indeed, if $x\qu z=
z\qu x=e$ with $z\in M$, then $z=z\qu e=z\qu(x\qu y) = (z\qu x)\qu y = e\qu y
=y$. The \el\ $y\in X$ \sf ying \er{3.3} is called the \ti{inverse} of~$x$
in $(M,\qu,e)$, and is denoted by $x\Inv$.

\ssk
A \ti{group} $(G,\qu,e)$ is a monoid \st \ti{all\/} \el s are invertible\index{group}
\cite{Groups}. See Remark \rf{r3.57}.

\ssk
A \sbm\ $H$ of a monoid $(M,\qu,e)$ is called a \ti{subgroup} of~$M$ if the\index{subgroup}
monoid $(H,\qu,e)$ is a group. Clearly $(\{e\},\qu,e)$ is a group which is
called the \ti{trivial\/} subgroup of~$M$.
\edf

If the binary \op\ of a group is denoted by~$+$, then $a\Inv$ is usually
called the \ti{\ng} of~$a$ and is denoted by $-a$. If the \op\ is denoted
by~$\cdot$, then $a\Inv$ is called the \ti{reciprocal\/} of~$a$ and is
denoted by~$a^{-1}$.

\blm3.3
Let $(G,\qu,e)$ be a group and let $H$ be a \sbm\ of~$G$. Then $H$~is a subgroup
of~$G$ iff \fa $a\in H$, $a\Inv$, the inverse of~$a$ in $(G,\qu,e)$, belongs
to~$H$. In this case $a\Inv$ is also the inverse of~$a$ in $(H,\qu,e)$.
\elm

\proof \

\ti{If\/}: Let $a\in H$ with inverse $a\Inv$ in $(G,\qu,e)$. Then $a\qu a\Inv
=e=a\Inv\qu a$. Since $a\Inv\in H$, $a\Inv$ is also the inverse of~$a$ in
$(H,\qu,e)$. Since $a$ is arbitrary in~$H$, $(H,\qu,e)$ is a group, hence
$H$~is a~subgroup of~$G$.

\ti{Only if\/}: If $H$ is a subgroup of~$G$, $H$~is a \sbm\ of~$G$ and the monoid
$(H,\qu,e)$ is a group. \E\Tf \fa $a\in H$, \te s $b\in H$ \st $a\qu b= e=
b\qu a$ (that is, $b$~is the inverse of~$a$ in the monoid $(H,\qu,e)$). Since
$H\sbs G$, $a\in G$. Hence \te s $a\Inv\in G$ \st $a\qu a\Inv=e =
a\Inv \qu a$. We show
that $a\Inv=b$. Indeed, $a\Inv = a\Inv\qu e=a\Inv\qu(a\qu b)= (a\Inv\qu a)
\qu b=e\qu b=b$.
\endproof

\bpr3.4
Let $(X,\qu,e)$ be a group. Then \fa $a,b\in X$, and $c_k\in X$, $k\in[1,n]$,
$n\in\Na${\rm:}
\bga3.4
(a\Inv)\Inv = a,\\
(a\qu b)\Inv = b\Inv \qu a\Inv, \lb{3.5}\\
\Bigl(\QU_{k=1}^n c_k\Bigr)\Inv = \QU_{j=1}^n (c_{n+1-j})\Inv, \lb{3.6} \\
(n\dqu a)\Inv = n\dqu a\Inv.\lb{3.7}
\e
\epr

\bex3.5
Prove \E\Pr\ \rf{p3.4}. (\ti{Hint\/}: use \In\ on $n\in\Na$ in \er{3.6},
\er{3.7} and \E\Pr\ \rfa2{p1.13} in \er{3.7}).
\eex

\blm3.6
Let $(M,\qu,e)$ be a monoid and let $M^\t$ denote the set of all invertible\glossary{$M^\t$}
\el s of~$M$. Then $M^\t$ is a \emph{subgroup} of~$M$, \Ip $(M^\t,\qu,e)$ is
a group.
\elm

\proof
\ti{$M^\t$ is a \sbm\ of $M$}: Clearly $e\in M^\t$. Let $a,b\in M^\t$ and let
$u,v\in M$ be \st\ $a\qu u=u\qu a=e$, $b\qu v=v\qu b=e$. Then $v\qu u\in M$
and $(a\qu b)\qu\break (v\qu u) \nda2{1.36} = a\qu(b\qu v)\qu u=a\qu e\qu u =a\qu u
=e$. Hence $a\qu b\in M^\t$ and \er{3.5} holds.

\ti{$(M^\t,\qu,e)$ is a subgroup of $M$}: Since $M^\t$ is a submonoid of~$M$,
$(M^\t,\qu,e)$ is a monoid (we still denote the
\rt ion of~$\qu$ to $M^\t \t M^\t$ by~$\qu$). Let $a\in M^\t$ and $u\in M$
be \st $a\qu u=u\qu a=e$. \If that $u\in M^\t$ and that $a$ is invertible
in~$M^\t$ with inverse~$u$. Hence $u=a\Inv$ and \er{3.4} holds.
\endproof

\bxs3.7
Let $F$ be a \ns\ and let $(F^F,\circ,\id_F)$ be the monoid introduced in
Example \rfa1{xa2.3}\,(v). Then a self-map $f$ of~$A$ is invertible in the
sense of \E\df\ \rf{d3.2} if \te s $g\in F^F$ \st $f\circ g=g\circ f=\id_A$,
that is, iff $f$~is bi\jc\ and its inverse~$g$ is the inverse of~$f$ in the
set-theoretical sense. Hence the notation $f\Inv$ for~$g$ is appropriate.
If we denote by $\Bij(F)$ or~$S_F$ the set of bi\jn s of~$F$, then
$S_F=(F^F)^\t$, \Ip $(S_F,\circ,\id_F)$ is a group by Lemma \rf{l3.6}.
Cayley's theorem states that every group~$G$
is monoid-\is c to a subgroup of the group $(S_G,\circ,\id_G)$.
\exs

\bdf3.8
A group $(G,\qu,e)$ is called \ti{finite} (resp.\ \ti{infinite}) if $G$~is
finite (resp.\ infinite). If $G$ is finite, $\#(G)$ is called the \ti{order}
of~$G$ and is also denoted by~$|G|$. If the monoid $(G,\qu,e)$ is abelian
(resp.\ \pn), then $G$~is called \ti{abelian} (resp.\ \ti{\pn}).\index{group!principal}
\edf

\brm3.9
The notion of \pn\ group is not standard but, as we shall see, is very
convenient.
\erm

\bpr3.10
\

\hph i,ii, A finite \Cm\ is a finite \ag\ and conversely.

\hph ii,i, A finite \pn\ \Cm\ is a \pn\ group and conversely.

\hph iii,, \E\fa $n\in\Na\sms1$, $(\N_n,+_n,0)$ is a \pn\ group with~$1$
as a \Gn.

\hph iv,, If $(G,\qu,e)$ is a \pn\ group with \Gn~$a$, then $\#(G)\ge2$
and the map $\ov\vf_a:
\zo0,\#(G) \to G$ defined by $\ov\vf_a(k):=k\dqu a$, $k\in\zo0,\#(G) $ is a
monoid-\is sm from $(\N_{\#(G)},+_{\#(G)},0)$ onto~$G$.

\hph v,i, If $(G_i,\qn i,e\en i)$, $i=1,2$, are \ep\ \pn\ groups with \Gn s
$a_1,a_2$ respectively, then the map $\vf:G_1\to G_2$ defined by
\beq3.8
\vf(k\qn 1 a_1):= k\qn2 a_2, \q k\in\zo0,\#(G_1)
\e
is a monoid-\is sm.

\hph vi,, Let $(G_i,\qu_i,e_i)$, $i=1,2$, be groups. Then a \hm sm $f:G_1\to
G_2$ is \emph{in\jc} iff $\ker f:=\{x\in G_1: f(x)=e_2\}=\{e_1\}$. $\ker f$~is
usually called the \emph{kernel} of~$f$.\index{kernel}
\epr

\proof \

(i) A finite \Cm\ $(X,\qu,e)$ is a group by \E\Pr\ \rf{p3.1}, which is finite
since $X$~is finite. A~\Cm\ is abelian by \df, thus $X$~is an \ag.

Let $(X,\qu,e)$ be a finite \ag. Then $(X,\qu,e)$ is a finite abelian monoid.
We show that it is \cnc e. Let $x,y,z\in X$ be \st $x\qu z=y\qu z$. Then
$x=x\qu e=x\qu(z\qu z\Inv) = (x\qu z)\qu z\Inv = (y\qu z)\qu z\Inv =y\qu
(z\qu z\Inv) = y\qu e=y$.

(ii) Let $(X,\qu,e)$ be a finite \pn\ \Cm. Then $(X,\qu,e)$ is a finite group
by~(i), and is \pn\ since the monoid $(X,\qu,e)$ is \pn.

Let $(X,\qu,e)$ be a \pn\ group and let $a\in X\sms e$ be \st $X=I(a)$.
Suppose for \cd ion that $X$~is infinite. Then by Theorem \rfa2{t1.42},
$(X,\qu,e)$ is a \PM. Since $a\qu a\Inv=e$, we have $a=a\Inv=e$. A~\cd ion,
since $a\ne e$, Thus $X$~is finite. \Mo $X$~is abelian by Lemma \rfa2{l1.21}\,
(ii). By~(i) it is a \Cm.

(iii) Follows from \E\Pr\ \rf{p1.8}\,(i), (ii).

(iv) First observe that $\#(G)\ge2$. Indeed, $\#(\{e,a\}) \nda2{3.31} =
\#(\{e\}) + \#(\{a\})$ since $e\ne a$, $\#(\{e\})=\#(\{a\})=1$, since
$\{e\},\{a\}$ are \ep\ to $\zo0,1 $, and $2=1+1=\#(\{e,a\}) \nda2{3.13} \le
\#(G)$ since $\{e,a\}\sbs G$. By~(ii) $(G,\qu,e)$ is a finite \Cm, hence by
Lemma \rf{l1.22} \te s $m\ge2$ \st $\ov\vf_a:(\N_m,+_m,0) \to (I(a),\qu,e)$
is an \is sm. \If that $m=\#(\N_m)=\#(I(a))=\#(G)$, which proves~(iv) since
$G=I(a)$.

(v) The map $\vf$ is well-defined, since \fe $x\in G_1$ there is \ooo $k\in
\zo0,\#(G_1) $ \st $x=k \mathrel{\lower1pt\hbox{$\stackrel{\lower7pt\hbox{\smash{\hbox to0pt{$
\scriptstyle\sqcap$\hss}$\scriptstyle\sqcup$}$^{(1)}$}}\cdot$}} a_1$ in view of~(iv). Hence by \er{3.8},
$\vf$~is a well-defined map from $G_1$ into $\{ k\mathrel{\lower1pt\hbox{$\stackrel{\lower7pt\hbox{\smash{\hbox to0pt{$
\scriptstyle\sqcap$\hss}$\scriptstyle\sqcup$}$^{(2)}$}}\cdot$}} a_2: k\in
\zo0,\#(G_1) \}$. Since $G_1$ and~$G_2$ are \ep,\break
$\#(G_1)=\#(G_2)$ by Corollary \rfa2{c3.11}, and
$\{ k\mathrel{\lower1pt\hbox{$\stackrel{\lower7pt\hbox{\smash{\hbox to0pt{$
\scriptstyle\sqcap$\hss}$\scriptstyle\sqcup$}$^{(2)}$}}\cdot$}} a_2: k\in
\zo0,\#(G_1) \}=G_2$ by (iv). Set $n:=\#(G_1)$ ($=\#(G_2)$), and let
$\ov\vf_{a_i} :\zo0,n \to G_i$ denote the \rt ion of $\vf_{a_i}$ to $\zo0,n $,
$i=1,2$. Then $\vf=\vf_{a_2}\circ \vf\Inv_{a_1}$. Since both $\vf_{a_1}$ and
$\vf_{a_2}$ are \is sms, $\vf$~is an \is sm by Lemma \rfa2{l1.8}.

(vi) \ti{Only if\/}: Since $f$ is a \hm sm, $f(e_1)=e_2$. If $f(x)=e_2$,
$x\in G_1$, then $x=e_1$ by the in\ji\ of~$f$.

\ti{If\/}: Let $x,y\in G_1$ be \st $f(x)=f(y)$. Let $z\in G_1$ be \st
$y\qu_1 z=e_1$. Then $f(x\qu_1 z)=f(x)\qu_2 f(z)= f(y)\qu_2 f(z) = f(y\qu_1 z)
= f(e_1) = e_2$. Hence $x\qu_1 z=e_1$. Thus $x= x\qu_1 e_1 = x\qu_1 (z\qu_1
z\Inv) = (x\qu_1 z)\qu_1 z\Inv = e_1\qu_1 z\Inv = (y\qu_1 z)\qu_1 z\Inv =
y\qu_1(z\qu_1 z\Inv) = y\qu_1 e_1 = y$.
\endproof

\bdf3.11 \

\hph i,i, A (\cmt e) \sr\ (with unity) is usually called a (\cmt e) \ti{ring}\index{ring}
(with unity) if its additive monoid is a group~\cite{Alg}. In this chapter,
unless explicitly mentioned, ``ring'' will always mean ``\cmt e ring with
unity''.

\hph ii,, A \ti{\pn\ ring} is a \pn\ \sr\ (see \E\df\ \rf{d1.6n}) \st its
additive monoid is a group.
\edf

\brm3.12
The notion of \pn\ ring is not standard but convenient, and it should not be
confused with the notion of \pn\ ideal ring~\cite{Alg}.
\erm

\bpr3.13 \

\hph i,ii, The additive monoid of a \pn\ ring is a \pn\ group. \E\Ip
a~\pn\ ring is finite.

\hph ii,i, A finite \pn\ c-\sr\ is a \pn\ ring and conversely.

\hph iii,, \E\fe $n\in\Na\sms1$, $(\N_n,+_n,\cdot_n,0,1)$ is a \pn\ ring.
\epr

\proof
(i) The additive monoid of a ring is a group, and of a \pn\ \sr\ is \pn, hence
it is a \pn\ group. It is finite by \E\Pr\ \rf{p3.10}.

(ii) The additive monoid of a finite \pn\ c-\sr\ is finite and \cnc e by \as.
Since it is \pn, it is abelian by Lemma \rfa2{l1.21}\,(ii), hence it is a finite
\pn\ \Cm. By \E\Pr\ \rf{p3.10}\,(ii), it is a \pn\ group. The \am\ of a \pn\
ring is a \pn\ group by~(i), hence it is a finite \pn\ \Cm\ by \E\Pr\
\rf{p3.10}\,(ii).

(iii) $(\N_n,+_n,\cdot_n,0,1)$ is a \pn\ \sr\ by \E\Pr\ \rf{p1.8}\,(i) and
$(\N_n,+_n,0)$ is a \pn\ monoid with \Gn~$1$ by \E\Pr\ \rf{p3.10}\,(iii).
\endproof

\brm3.14
The full analogue of \E\Pr\ \rf{p3.10} for \pn\ rings can be found in
Theorem \rf{t1.28} where finite \pn\ \cnc e \sr\ has to be replaced by
\pn\ ring in view of (ii) of the preceding \Pr.
\erm

\begin{ntt}[\cite{Groups}]\lb{n3.15a}
Let $G$ be a \pn\ group. Set
\beq3.10
{\rm Gen}(G):=\{x\in G: I(x)=G\}.
\e
\ent

\bpr3.15n \

\hph i,i, Let $(R,\qu_0,\qu_1,e_0,e_1)$ be a \pn\ ring and $n:=\#(R)\ge2$. Let
$\vf_{e_1}$ be the unique ring-\is sm from $(\N_n,+_n,\cdot_n,0,1)$ onto~$R$.
Then
\beq3.9
a\qu_1 b = (\vf_{e_1}\Inv(a) \cdot_n \vf_{e_1}\Inv(b))\dquz e_1, \q a,b\in R.
\e

\hph ii,, Let $G$ be a \pn\ group of order~$n$. Given $c\in{\rm Gen}(G)$, set
\beq3.11
a\qu_1^{(c)}b := (\vf_{c}\Inv(a) \cdot_n \vf_{c}\Inv(b))\dquz c, \q a,b\in G.
\e
Then $(G,\qu,\qu_1^{(c)},e_0,c)$ is a \pn\ ring. \Mo $\qu_1^{(c)}$ is the only
binary \op\ on~$G$ \st $(G,\qu,\qu_1^{(c)},e_0,c)$ is a \pn\ ring.
\epr

\bex3.16
Prove \E\Pr\ \rf{p3.15n}.
\eex

\bxs3.18
Observe that if in \E\Pr\ \rf{p3.15n} $(R,\qu_0,e_0):=(\N_n,+_n,0)$ with
$n\ge2$, $c\in\zo0,n $ is a \Gn\ of $(\N_n,+_n,0)$ and $a:=k
\mathrel{\lower1pt\hbox{$\stackrel{\lower7pt\hbox{\smash{\hbox to0pt{$
\scriptstyle\sqcap$\hss}$\scriptstyle\sqcup$}$_n$}}\cdot$}}c$, $b:=l
\mathrel{\lower1pt\hbox{$\stackrel{\lower7pt\hbox{\smash{\hbox to0pt{$
\scriptstyle\sqcap$\hss}$\scriptstyle\sqcup$}$_n$}}\cdot$}}c$, $k,l\in
[0,n)$, then we have from \er{3.11}:
\beq3.12
a \qu_1^{(c)} b
=(k\cdot_n l)\dpls n c.
\e

\ti{Case} $n:=3$: $c:=2$ is a \Gn\ since $0\dpls 32=0$, $1\dpls32=2$,
$2\dpls32=1$. Hence, the \mlc\ table of $\qu_1^{(1)}$ is
\[
\begin{array}{c|ccc}
\ \cdot_3\ &\ 0\ &\ 1\ &\ 2 \ \\
\hline
0&0&0&0\\
1&0&1&2\\
2&0&2&1
\end{array}
\]
and the \mlc\ table of $\qu_1^{(2)}$ is
\[
\begin{array}{c|ccc}
\ \qu_1^{(2)}\ &\ 0\ &\ 1\ &\ 2 \ \\
\hline
0&0&0&0\\
1&0&2&1\\
2&0&1&2
\end{array}
\]

\ti{Case} $n=4$. $c=3$ is a \Gn\ since $0\dpls43=0$, $1\dpls43=3$,
$2\dpls43=2$, $3\dpls43=1$.
\[
\begin{array}{c|cccc}
\ \cdot_4\ &\ 0\ &\ 1\ &\ 2\ &\ 3\ \\
\hline
0&0&0&0&0\\
1&0&1&2&3\\
2&0&2&0&2\\
3&0&3&2&1
\end{array}
\qquad
\begin{array}{c|cccc}
\ \qu_1^{(3)}\ &\ 0\ &\ 1\ &\ 2\ &\ 3\ \\
\hline
0&0&0&0&0\\
1&0&3&2&1\\
2&0&2&0&2\\
3&0&1&2&3
\end{array}
\]
\exs

\bex3.19
Give the \mlc\ tables for the cases $n:=5$, $c:=2,3$; $n:=6$, $c:=5$;
$n:=10$, $c:=3,7,9$.
\eex

We now \es\ \pp y (B) and some other \pp ies of \sbm s of \pn\ groups.

\bpr3.20
Every nontrivial \sbm\ of a \pn\ group is a \pn\ subgroup.
\epr

In the proof we shall use the \fw\ lemma.

\blm3.21
A \emph{finite} \sbm\ of an \ag\ is a subgroup.
\elm

\proof
Let $(G,\qu,e)$ be an \ag\ and let $M$ be a nontrivial \sbm. Then $(M,\qu,e)$
is a finite abelian monoid. We show that $(M,\qu,e)$ is \cnc e. Let $x,y,z\in M$
be \st $x\qu z=y\qu z$. Then $x=x\qu e = x\qu(z\qu z\Inv) = (x\qu z)\qu
z\Inv = (y\qu z)\qu z\Inv = y\qu (z\qu z\Inv)=y\qu e=y$. Hence $(M,\qu,z)$
is a finite \Cm, and by \E\Pr\ \rf{p3.1} it is a~group. Thus $M$~is a~subgroup
of~$G$.
\endproof

\brm3.22
We shall see later that the \as\ ``abelian'' can be dropped.
\erm

\proof[Proof of \E\Pr\ \rf{p3.20}]
Let $(X,\qu,e)$ be a \pn\ group and let $M$ be a nontrivial \sbm.
Let $a\in X\sms e$ be a~\Gn\ of~$X$, and
let $\vf_n:\N\to X$ be the \sq\ of \IT s of~$a$. Set $\wt M:=\{n\in \N:
\vf_a(n)\in M\}$. We show that $\wt M$ is a \sbm\ of $(\N,+,0)$. Indeed,
$\vf_a(0)=e\in M$, hence $0\in\wt M$. Let $m,n\in\wt M$, then \te\ $b,c\in M$ \st $\vf_a(m)=b$ and
$\vf_a(n)=c$. But $\vf_a(m+n) \nad{\era2{2.3}\,{\rm I2}} = \vf_a(m)\qu
\vf_a(n) = b\qu c\in M$. Hence $m+n\in \wt M$, and $\wt M$ is  a~\sbm\ of
$(\N,+,0)$. We next show that $\wt M$ is a \ti{\pn\/} \sbm\ of $(\N,+,0)$.
To this end we apply Lemma \rfa3{l8.20}. Let $m,n\in \wt M$ be \st $m\le n$,
and let $p\in\N$ be \st $n=p+m$. Then $\vf_a(n)=\vf_a(p+m) \nad{\era2{2.3}
\,{\rm I2}} = \vf_a(p)\qu\vf_a(m)$. Hence $\vf_a(p)=\vf_a(p)\qu (\vf_a(m)
\qu(\vf_a(m))\Inv) = (\vf_a(p)\qu\vf_a(m))\qu(\vf_a(m))\Inv =\vf_a(n)\qu
(\vf_a(m))\Inv$, where $(\vf_a(m))\Inv$ is the inverse of $\vf_a(m)$ in the
group $(X,\qu,e)$. Since $(X,\qu,e)$ is a \pn\ group, $X$~is finite by \E\Pr\
\rf{p3.10}\,(ii), hence $M$ as a subset of~$X$ is finite by Theorem
\rfa1{t4.18}\,(ii). By Lemma \rf{l3.21}, $M$~is a \ti{subgroup} of the
group~$X$. Thus $\vf_a(m)\Inv\in M$ since $\vf_a(m)\in M$. \E\Tf $\vf_a(p)
= \vf_a(n)\qu(\vf_a(m))\Inv \in M$. Hence $p\in\wt M$,
i.e.\ $n-m\in \wt M$. Thus in view of Lemma \rfa3{l8.20}, $\wt M$~is a~\pn\
\sbm\ of $(\N,+,0)$, i.e.\ there is $\ov n\in\Na$ \st $\wt M=I(\ov n)$. We now
show that $M\sbs I(\vf_a(\ov n))$. Indeed, if $c\in M$, then \te s $k\in \N$
\st $c=\vf_a(k)$ since $X$~is a \pn\ group, hence $k\in\wt M$. \Mo \te s
$l\in \N$ \st $k=l\ov n$ since $\wt M$ is a \pn\ \sbm\ of $(\N,+,0)$ with
\Gn~$\ov n$. Hence $c=\vf_a(k)=\vf_a(
l\ov n) = (l\ov n)\dqu a \nad{\era2{2.3}\,\rm I3}= l\dqu(\ov n\dqu a)\in
I(\vf_a(\ov n))$. \E\Tf $M\sbs I(\vf_a(\ov n))$. \E\oh $\ov n\in\wt M$, hence
$\vf_a(\ov n)\in M$ by \df\ of~$\wt M$. By Lemma \rfa2{l1.21}\,(i) with
$(M,\qu,e):=(X,\qu,e)$, $M_0:=M$, $a:=\vf_a(\ov n)$, we find that $I(\vf_a(
\ov n))$, the range of $\vf_{\vf_a(\bar n)}$ in~$X$, is contained in~$M$ since
$\vf_a(\ov n)\in M$. \csq, $M=I(\vf_a(\ov n))$.
\endproof

It turns out that the order of a \sbm\ of a \pn\ group~$G$ divides the order
of~$G$ (see \E\df\ \rf{d2.8}), and that \fe divisor~$d$ of $\#(G)$, \te s \ooo
\sbm\ of order~$d$.

\blm3.23
Let $(G,\qu,e)$ be a \pn\ group, let $c\in G\sms e$ and let $k\in\zo1,\#(I(c)) $.
Then
\beq3.13
\#(I(c)) = \gcd(k,\#(I(c)))\cdot \#(I(k\dqu c)).
\e
\elm

\proof
$I(c)$ is a \sbm\ of~$G$ by Lemma \rfa2{l1.21}\,(ii) with $(M,\qu,e):=
(G,\qu,e)$, $M_0:=M$, $a:=c$ and $\wt E:=\N$. Since $G$ is a finite \Cm\ by
\E\Pr\ \rf{p3.10}\,(ii), we may apply Lemma \rf{l1.22} with $(X,\qu,e):=
(G,\qu,e)$, $a:=c$. We find that $m\in\Na\sms1$ defined in~(i) \sf ies
$m\nda23.5 = \#(\zo0,m )=\#(I(c))$ by Corollary \rfa2{c3.11} since $\zo0,m $
and $I(c)$ are \ep\ by~(vi). \If from~(v) that
\beq3.14
\hbox{if $l\dqu c=e$ \fs $l\in\N$, then $l=q\cdot \#(I(c))$ \fs $q\in\N$.}
\e

We now apply Lemma \rf{l1.22} again with $(X,\qu,e):=(G,\qu,e)$, but with
$a:=k\dqu c$. We find $m'\in\Na\sms1$ \st $m'=\#(I(k\dqu c))$. \Mo by~(iii) with $q:=1$, we get
$m'\dqu (k\dqu c)=e$. Using \era2{2.3}\,I3 we find $(m'k)\dqu c=e$. \If from
\er{3.14} that
\beq3.15
m'k = q\cdot \#(I(c)) \qh{\fs}q\in\N.
\e
Since both $m'$ and $k$ belong to $\Na$, the \LHS\ of \er{3.15} belongs
to~$\Na$ by \E\Pr\ \rfa2{p2.7}\,(i). \E\Tf $q\cdot\#(I(c))\in\Na$, hence
$q\in\Na$ by \era2{2.13}. \If that $m'k$ is a common \ml e of
$k$ and~$\#(I(c))$. We now show that $m'k$ is the \ti{least\/} common \ml e
of $k$ and~$\#(I(c))$. Suppose, for \cd ion, that $m'k$ is not the least common
\ml e of $k$ and~$\#(I(c))$. Then \te s $\bar l\in\zo1,m' $ \st $\bar lk=q
\#(I(c))$ \fs $q\in\Na$. This implies $\bar l\dqu(k\dqu c)\nad{\era2{2.3}\,\rm
I3}= (\bar lk)\dqu c= (q\#(I(c)))\dqu c = q\dqu((\#(I(c)))\dqu c)\nde2.8 =
q\dqu e\nad{\era2{2.3}\,\rm I4}= e$.

\E\Tf $\bar l\dqu(k\dqu c)=e$ with $\bar l\in(0,m')$, \cd ing \er{1.43}. Thus
\beq3.16
m'k = \lcm(k,\#(I(c))).
\e
\E\ml ying both sides of \er{3.16} by $\gcd(k,\#(I(c)))$ and using \era3{1.25}
we arrive at\break $\gcd(k,\#(I(c)))\cdot m'k = k\cdot \#(I(c))$. By \era2{2.11},
\era2{2.12} and \E\Pr\ \rfa2{p2.7}\,(iii), we obtain $\gcd(k,\#(I(c)))\cdot m'
= \#(I(c))$, which is \er{3.13} since $m'=\#(I(k\dqu c))$.
\endproof

\bpr3.24
Let $(G,\qu,e)$ be a \pn\ group.

\hph i,i, If $H$ is a \sbm\ of $G$, then $\#(H)$ divides $\#(G)$.

\hph ii,, \E\fe divisor $d$ of $\#(G)$, $d\ge 2$, \te s \ooo \sbm\ of~$G$
of order~$d$.
\epr

\proof
Let $a\in G\sms e$ be a \Gn\ of~$G$.

(i) Let $H$ be a \sbm\ of $G$. If $H$ is trivial, then $\#(H)=1$ hence $1 | \#(G)$. If
$H$~is nontrivial, then $(H,\qu,e)$ is a \pn\ group by \E\Pr\ \rf{p3.20}. Let
$b\in H\sms e$ be a \Gn\ of~$H$. Then there is $k\in\zo1,\#(G) $ \st $b=k
\dqu a$ by Lemma \rf{l1.22}. From Lemma \rf{l3.23} with $c:=a$ we obtain
$\#(G) = \#(I(a)) = \gcd(k,\#(G)) \cdot \#(I(k\dqu a)) = \gcd(k,\#(G)) \cdot
\#(I(b)) = \gcd(k,\#(G))\cdot \#(H)$. Hence $\#(H)|\#(G)$.

(ii) {\it Existence}. Set $n:=\#(G)\ge2$ and let $d\ge2$ be a divisor of~$n$.
Let $k$ be the unique \el\ of~$\Na$ \st $n=kd$. Set $b:=k\dqu a$. We claim that
$I(b)$ is a \sbm\ of~$G$ with $\#(I(b))=d$. In view of Lemma \rfa2{l1.21}\,(ii),
with $(M,\qu,e):=(G,\qu,e)$, $M_0:=M$ and $a:=b$, $I(b)$~is a \sbm\ of~$G$.
We now show that $\#(I(b))=d$. By Lemma \rf{l1.22}\,(i) with $(X,\qu,e):=
(G,\qu,e)$ and $a:=b$, \te s a unique \el\ $m\in\Na\sms1$ \sf ying \er{1.43}.
\Mo by (vi) $\zo0,m $ is \ep\ to~$I(b)$. Thus $m\nda23.5 = \#(\zo0,m )=
\#(I(b))$ by Corollary \rfa2{c3.11}. \E\Tf it suffices to show that $d$~\sf ies
\er{1.43}, that is, $d\dqu b=e$ and $d'\dqu b\ne e$ for $d'\in\zo1,d $. But
$d\dqu b = d\dqu(k\dqu a)\nad{\era2{2.3}\,\rm I3} = dk\dqu a\nda22.12 =
kd\dqu a = n\dqu a\nde2.8 = e$. \E\oh we have $l\dqu a\ne e$ whenever
$l\in\zo1,n $ by \er{2.7}. Let $d'\in\zo1,d $. Then $d'k\in\Na$ by \E\Pr\
\rfa2{p2.7}\,(iii) and $d'k\nda22.19 < dk\nda22.12 = kd=n$. Hence $d'k\in
\zo1,n $. Thus $d'\dqu b= d'\dqu (k\qu a)=d'k\dqu a\ne e$.

(iii) \ti{\E\uq}. Let $k$, $b$ and $I(b)$ be as in part~(ii). Let $Y$~be a
\sbm\ of~$G$ of order $d\ge2$. We want to show that $Y=I(b)$. By \E\Pr\
\rf{p3.20}, $(Y,\qu,e)$ is a \pn\ group, hence \te s $c\in G\sms e$ \st
$Y=I(c)$. \E\Ip $\#(I(c))=\#(Y)=d$. By~\er{2.8}, we have $d\dqu c=e$. \E\oh
since $c\in G=I(a)$, \te s $l\in\Na$ \st $c=l\dqu a$. \E\Tf $e=d\dqu c = d\dqu
(l\dqu a) = dl\dqu a$. By Lemma \rf{l1.22}\,(iii) with $(X,\qu,e):=
(G,\qu,e)$, $a:=a$, and $m=I(a)=n$, we obtain $dl=qn$ \fs $q\in\Na$, since
$dl\in\Na$. But $n=kd$  and $dl=ld$, thus by \cnc ity ($d\in\Na$), we obtain
$l=qk$. \E\Tf $c=qk\dqu a=q\dqu(k\dqu a) = q\dqu b\in I(b)$. Thus $c\in
I(b)$. In view of Lemma \rfa2{l1.21}\,(i) with $(M,\qu,e):=(G,\qu,e)$,
$M_0:=I(b)$, $a:=c$, we have $I(c)\sbs I(b)$. Since $\#(I(c))=\#(Y)=d=
\#(I(b))$, we have $I(c)=I(b)$ by \era2{3.13}. Hence $Y=I(c)=I(b)$.
\endproof

We now \es\ \pp y (C) and a formula for  the \nm\ of \Gn s of a \pn\ group.
We recall that a finite \pn\ \Cm\ is a \pn\ group.

\blm3.25
Let $(X,\qu,e)$ and $(X',\qu',e')$ be monoids and let $f:X\to X'$ be an \is
sm. Then

\hph i,ii, If $X$ is a group, then so is $X'$ and
\beq3.17
f(x\Inv) = f(x)\Inv \qh{\fa} x\in X.
\e

\hph ii,i, If $(X,\qu,e)$ is a \pn\ group with $a\in X\sms e$ as a \Gn, then
$(X',\qu',e')$ is a \pn\ group with $f(a)$ as a~\Gn.

\hph iii,, If $X$ and $X'$ are groups and $f$ is a \hm sm, then \er{3.17} holds.
\elm

\proof
(i) We first show that every \el\ $x'\in X'$ is invertible (see \E\df\
\rf{d3.2}). Since $f$~is bi\jc, \te s one and only one $x\in X$ \st $f(x)=x'$. Since
$X$~is a group we have $x\qu x\Inv=e$ and $x\Inv\qu x=e$. Hence $e'\nad*=
f(e)= f(x\qu x\Inv)\nad*= f(x)\qu' f(x\Inv)$, and $e'\nad*=
f(e)= f(x\Inv\qu x)\nad*= f(x\Inv)\qu' f(x)$. In $\nad*=$ we used the fact
that $f$~is a \hm sm. \If that $x'\qu'f(x\Inv)=e'$ and $f(x\Inv)\qu' x'=e'$,
hence $x'$~is invertible with inverse equal to $f(x\Inv)$. \E\Tf
$(X',\qu',e')$ is a group and $(x')\Inv=f(x\Inv)$. Thus \er{3.17} holds.

(ii) By (i), $(X',\qu',e')$ is a group. It suffices to show that
$X'=I(f(a))$.

$f(I(a))\sbs I(f(a))$: Let $x'\in f(I(a))$ and let $x:=f\Inv(x')$. Then
$x\in I(a)$. \E\te s $k\in\N$ \st $x=k\dqu a$. Then $x'=f(x)=f(k\dqu
a)\nda21.45 = k
\mathrel{\lower1pt\hbox{$\stackrel{\lower7pt\hbox{\smash{\hbox to0pt{$\scriptstyle\sqcap'$\hss}$\scriptstyle\sqcup$}}}\cdot$}}
f(a) \in I(f(a))$.

$I(f(a))\sbs f(I(a))$: Let $x'\in I(f(a))$. Then \te s $k\in\N$ \st $x'=k
\mathrel{\lower1pt\hbox{$\stackrel{\lower7pt\hbox{\smash{\hbox to0pt{$\scriptstyle\sqcap'$\hss}$\scriptstyle\sqcup$}}}\cdot$}}
f(a)$, and $f\Inv(x')=f\Inv(k
\mathrel{\lower1pt\hbox{$\stackrel{\lower7pt\hbox{\smash{\hbox to0pt{$\scriptstyle\sqcap'$\hss}$\scriptstyle\sqcup$}}}\cdot$}}
f(a))\nda21.45 = k\dqu f\Inv(f(a))=k\dqu a$. Hence $x'=f\circ f\Inv(x') =\break
f(k\dqu a)\in f(I(a))$.

\If that $f(I(a))=I(f(a))$. But $I(a)=X$ by \as, thus $f(X)=I(f(a))$. Since
$f$~is sur\jc, $X'=f(X)$, hence $X'=I(f(a))$.

(iii) The proof is left as an exercise.
\endproof

\Wanp \es\ \pp y (C).

\bpr3.26
Let $(X,\qu,e)$ be an \emph{abelian} group. Then

\hph i,i, the map ${\rm Inv}:X\to X$ defined by ${\rm Inv(x)}:=x\Inv$,
$x\in X$, \sf ies ${\rm Inv}\circ {\rm Inv}=\id_X$ and is an auto\mf\ of
$(X,\qu,e)${\rm;}

\hph ii,, if in  \ad, $X$ is a \pn\ group with $a\in X\sms e$ as a \Gn, then
$a\Inv$ is also a \Gn\ of~$X$. \Mo $a=a\Inv$ iff $\#(X)=2$.
\epr

\proof
(i) Let $x\in X$. Then $x\qu x\Inv = x\Inv\qu x=e$. \Mo $(x\Inv)\Inv\qu
x\Inv=e$, hence $(x\Inv)\Inv = (x\Inv)\Inv\qu e = (x\Inv)\Inv \qu (x\Inv\qu
x) = ((x\Inv)\Inv \qu x\Inv)\qu x = e\qu x=x$. Hence ${\rm Inv}\circ{\rm Inv}
=\id_X$ and Inv is a bi\jn\ of~$X$. We have ${\rm Inv}(e)=e$ since $e\qu
e=e$.

Given $a,b\in X$, we have $(a\qu b)\qu (b\Inv\qu a\Inv) \nda21.36 = a\qu
(b\qu b\Inv)\qu a\Inv = a\qu e\qu a\Inv = a\qu a\Inv =e$. Hence $(a\qu b)\Inv
=b\Inv\qu a\Inv \nad*= a\Inv\qu b\Inv$, where in $\nad*=$ we used the \cmt
ity of~$\qu$. \E\Tf ${\rm Inv}(a\qu b)={\rm Inv}(a)\qu{\rm Inv}(b)$, hence
Inv is an endo\mf\ of~$X$. Since it is bi\jc, it is an auto\mf\ by Lemma
\rfa2{l1.8}.

(ii) follows from Lemma \rf{l3.25}\,(ii) with $(X',\qu',e'):=(X,\qu,e)$,
$f:=\rm Inv$, and part~(i). If $\#(X)=2$, then $X=\{e,a\}$. We have $a\qu a=e$,
since $a\qu a=a$ implies $a\qu a=e\qu a$, hence $a=e$ by \cnc ity, a \cd ion.
Thus $a=a\Inv$. Now we suppose $a\Inv=a$, hence $2\dqu a=a\qu a= a\qu a\Inv=
e$. Since $1\dqu a=a\ne e$, $2\dqu a=e$, we have $m=2$ in Lemma \rf{l1.22}
with $(X,\qu,e):=(X,\qu,e)$ and $a:=a$. \E\Tf $I(a)$ is \ep\ to $\zo0,2 $ by
Lemma \rf{l1.22}\,(vi) and Corollary \rfa2{c3.11}. \If that $\#(I(a))=2$,
hence $\#(X)=2$ since $X=I(a)$.
\endproof

We are now interested in the \nm\ of \Gn s of a \pn\ group. Recall that we
denoted by ${\rm Gen}(G)$ (see \er{3.10}) the set of \Gn s of a \pn\ group~$G$.
We shall show that if $G$ and~$G'$ are \is c \pn\ groups, then $\#({\rm Gen}(G)
)=\#({\rm Gen}(G'))$. Since $G$~is \is c to $(\N_n,+_n,0)$ with $n=\#(G)$,
$\#({\rm Gen}(G))=\#({\rm Gen}(\N_n,+_n,0))$. It turns out that $k\in\zo1,n $
is a \Gn\ of $(\N_n,+_n,0)$ iff $\gcd(k,n)=1$. Thus $\#({\rm Gen}(G))=
\#(\{k\in\zo1,n : \gcd(k,n)=1\})$.

The map $\phi:\Na\to\Na$ defined by
\beq3.18
\phi(n):=\#\bigl(\{k\in[1,n]: \gcd(k,n)=1\}\bigr), \q n\in\Na,
\e
was introduced by Euler and is called the \ti{Euler totient \f} or the\index{totient function}
\ti{Euler $\phi$-\f}.

\blm3.27 \

\hph i,i, Let $G$ and $G'$ be \is c \pn\ groups. Then
\beq3.19
\#({\rm Gen}(G)) = \#({\rm Gen}(G')).
\e

\hph ii,, Let $n\in\Na\sms1$. Then
\beq3.20
{\rm Gen}((\N_n,+_n,0)) = \{k\in[1,n]: \gcd(k,n)=1\}.
\e
\elm

\proof
(i) The sets Gen$(G)$ and Gen$(G')$ are not empty by \as. Let $f:G\to G'$ be
an \is sm. By Lemma \rf{l3.25}\,(ii), if $x\in{\rm Gen}(G)$, then $f(x)\in
{\rm Gen}(G')$, hence $f({\rm Gen}(G)) \sbs {\rm Gen}(G')$. Similarly, since
$f\Inv:G' \to G$ is an \is sm by Lemma \rfa2{l1.8}, $f\Inv({\rm Gen}(G')
\sbs{\rm Gen}(G)$. Then ${\rm Gen}(G')=f(f\Inv({\rm Gen}(G')))\sbs f({\rm Gen}
(G))$. \E\Tf ${\rm Gen}(G')=f({\rm Gen}(G))$, hence ${\rm Gen}(G')$ and
${\rm Gen}(G)$ are \ep. By \E\Pr\ \rf{p3.10}, $G$~and~$G'$ are finite, and so
are ${\rm Gen}(G)$, ${\rm Gen}(G')$ as subsets of finite sets by Lemma
\rfa2{l3.16}. Then $\#({\rm Gen}(G))=\#({\rm Gen}(G'))$ by Corollary
\rfa2{c3.11}.

(ii) Let $n\ge2$ and let $k\in\zo0,n $. We denote $(\N_n,+_n,0)$ by~$\N_n$.
Clearly $0\notin{\rm Gen}(\N_n)$ since $I(0)=\{0\}$ by \era2{2.3}\,I4. \Mo
$1\in{\rm Gen}(\N_n)$ by \E\Pr\ \rf{p1.8}\,(i). Then $k\in\zo1,n $ belongs
to ${\rm Gen}(\N_n)$ iff $\#(I(k))=n$. Indeed, if $I(k)=\N_n$, then $\#(I(k))
=\#(\N_n)=n$. Conversely, if $\#(k)=n$, then $I(k)=\N_n$ by \era2{3.13}. We
have $\#(I(1))=n$ and by Lemma \rf{l3.23} with $(G,\qu,e):=(\N_n,+_n,0)$ and
$c:=1$, $\#(I(1))\nde3.13 = \gcd(k,\#(I(1)))\cdot\#(I(k\dpln 1)$. But
$k\dpln1 \nde1.11 = \F_n(k)\nde1.5 = k$. Hence $n=\gcd(k,n)\cdot \#(I(k))$.
Thus, if $\gcd(k,n)=1$, then $n=1\cdot \#(I(k))\nda22.13 = \#(I(k))$, and
$k\in{\rm Gen}(\N_n)$. \E\oh if $k\in{\rm Gen}(\N_n)$, then $\#(I(k))=n$, and
$n\cdot 1\nda22.18 = n=\gcd(k,n)\cdot n\nda22.12 = n\cdot\gcd(k,n)$.
\If that $1=\gcd(k,n)$ by \E\Pr\ \rfa2{p2.7}. \csq, ${\rm Gen}(\N_n)=
\{k\in\zo1,n : \gcd(k,n)=1\}$. Since $\gcd(n,n)=n\ge2>1$, $\gcd(n,n)\ne1$ by
\era1{3.12}. Hence ${\rm Gen}(\N_n)=\{k\in[1,n]: \gcd(k,n)=1\}$.
\endproof

\bth3.28
Let $G$ be a \pn\ group of order $n\ge2$. Then
\beq3.21
\#({\rm Gen}(G))=\F(n)
\e
where $\F(n)$ is defined in \er{3.18}.
\eth

\proof
By \E\Pr\ \rf{p3.10}\,(iv), $(G,\qu,e)$ is monoid-\is c to $(\N_n,+_n,0)$
since $n=\#(G)$. Thus $({\rm Gen}(G)) \nde3.19 = \#({\rm Gen}(\N_n)) \nad
{\er{3.20},\er{3.18}}= \F(n)$ since $\#(\N_n)=n$.
\endproof

We now investigate \pp ies of the Euler $\F$-\f. Clearly
\beq3.22n
\F(1)=1,
\e
since $\gcd(1,1)=1$. Since $\gcd(1,2)=1$ and $\gcd(2,2)=2$, we have
\beq3.23n
\F(2)=1.
\e
Using \E\Pr\ \rf{p3.26}\,(ii) and Theorem \rf{t3.28}, we obtain
\beq3.24
\F(n) \hbox{ is even (i.e.\ $2|\F(n)$) whenever }n>2.
\e

\bex3.29
Prove \er{3.24}.
\eex

If $n:=p$ prime and $k\in[1,p]$, then $\gcd(p,p)=p$ and $\gcd(k,p)=1$ \fa
$k\in[1,p)$, since the set of common divisors of $k$ and~$p$ is~$\{1\}$.
Thus $\F(p)=\#(\zo1,p )\nda23.26 = {\#(\zo0,p-1 )}\nda23.5 = p-1$. Hence
\beq3.25
\F(p)=p-1 \qh{if $p$ is prime.}
\e

\begin{prp}[Gauss]\label{4.p3.30}
Let $n\in\Na$, then
\beq3.26n
n=\sum_{\sbk{1\le d\le n\\d|n}}\F(d).
\e
\epr

\proof[First proof]
We first give a proof based on the results above  about \pn\ groups.
Then we shall give a direct proof. The case
$n=1$ is trivial. Let $n\ge2$ and let $(G,\qu,e)$ (for example
$(\N_n,+_n,0)$) be a \pn\ group of order~$n$. Let $D(n):=\{d\in\Na: d|n\}$
(see \era3{8.55}). In view of \E\Pr\ \rf{p3.24}\,(ii), \te s \ooo \sbm\
of~$G$ of order~$d$, which we denote by~$H_d$, \fe $d\in D(n)\sms 1$. By
\E\Pr\ \rf{p3.20}, since $\#(H_d)\ge2$, $(H_d,\qu,e)$ is a \pn\ group. Hence
\fe $d\in D(n)\sms1$, \te s $a\in H_d\sms e$ \st $I(a)=H_d$. Set
\beq3.27
A_d:={\rm Gen}(H_d), \q d\in D(n)\sms 1,
\e
and
\beq3.28
A_1:=\{e\}.
\e
We  claim that $\zb Ad{D(n)}$ is a \pt\ of~$G$.

We have
\beq3.29
A_d\ne\vn \qh{\fa} d\in D(n).
\e
Indeed, $e\in A_1$. Let $d\in D(n)\sms1$ and let $H_d$ denote the \pn\ \sbm\
of order~$d$. By \df\ $H_d$ possesses a \Gn, hence $A_d\ne\vn$.
\beq3.30
\hbox{If $d,d'\in D(n)$ and $d\ne d'$, then $A_d\cap A_{d'}=\vn$.}
\e
Suppose for \cd ion that \te s $x\in A_d\cap A_{d'}$. Then $I(x)$ the set of
\IT s of~$x$ in~$G$ \sf ies $I(x)=H_d$ and $I(x)=H_{d'}$ since $x$~is a \Gn\
of both $H_d$ and $H_{d'}$. However $H_d\ne H_{d'}$ since $\#(H_d)=d\ne d'=
\#(H_{d'})$. A~\cd ion.
\beq3.31
G=\bigcup_{d\in D(n)}A_d.
\e
It suffices to show that $G\sbs\bcl_{d\in D(n)}A_d$. Clearly $e\in A_1\sbs
\bcl_{d\in D(n)}A_d$. Let $x\in G\sms e$. By Lemma \rfa2{l1.21}\,(ii) with
$M=M_0:=G$ and $a:=x$, $I(x)$ is a \sbm\ of~$G$. Set $d':=I(x)$, then $d'\in
D(n)$ and $I(x)=H_{d'}$ by \E\Pr\ \rf{p3.24}\,(i) and (ii). Since $H_{d'}=
I(x)$, $x\in{\rm Gen}(H_{d'})=A_{d'}\sbs \bcl_{d\in D(n)}A_d$. In view of
\E\df\ \rfa1{d4.5}, $\zb Ad{D(n)}$ is a \pt\ of~$G$. From \era2{7.28} and
\era2{1.127} with $\O:=G$, $I:=D(n)$, $a_\o:=1$ \fa $\o\in\O$, we obtain
\beq3.32
\#(G)=\sum_{d\in D(n)} \#(A_d).
\e
Note that $A_1:=\{e\}\ax\zo0,1 $, hence $\#(A_1)\nda23.5 = 1\nde3.22n = \F(1)$,
and $\#(A_d)\nde3.27 = \#({\rm Gen}(H_d))\break\nde3.21 = \F(d)$ since $\#(H_d)=d$,
$d\in D(n)\sms1$. From $\#(G)=n$, we infer \er{3.26n}. Note that, in view of
\er{3.29}--\er{3.31}, $\zb Ad{D(n)}$ is  a \pt\ of~$G$ induced by the \ev ce
\rl~$\sim$ on~$G$ defined by $x\sim y$ if $\#(I(x))=\#(I(y))$, $x,y\in G$.
\endproof

We now give a proof of \er{3.26n}, based on \pp ies of ``$\gcd$''. Let $n\in\Na
\sms1$ and $i\in [1,n]$. By \df\ $\gcd(i,n)|n$, thus $\gcd(i,n)\in D(n)$. We
define a \rl~$\sim'$ on $[1,n]$ by setting
\beq3.33
i \sim' j \qh{if } \gcd(i,n)=\gcd(j,n),\q i,j\in[1,n].
\e
Clearly the \rl\ $\sim'$ is an \ev ce \rl. Let $[x]$ denote the $\sim'$-\ev ce
class containing $x\in[1,n]$. Note that
\beq3.34
\gcd(d,n)=d \qh{\fa} d\in D(n).
\e
Thus we have
\beq3.35
[d] = \{i\in[1,n]: \gcd(i,n)=d\} \qh{\fa}d\in D(n).
\e
We claim that the map $d\mt [d]$ from $D(n)$ into~$\cA$, the set of $\sim'$-\ev
ce classes, is \ti{bi\jc}. Indeed, if $A\in\cA$, then $A\ne\vn$ by \E\df\
\rfa1{d4.5}. Let $j\in A\sbs[1,n]$, then $A=[j]=\{i\in[1,n]: g(i,n)=g(j,n)\}
\nde3.35 = [d']$ with $d':=\gcd(j,n)\in D(n)$. \E\oh if $[d]=[d']$, $d,d'\in
D(n)$, then $d\sim' d'$, hence $d=\gcd(d,n)=\gcd(d',n)=d'$, which proves the
claim. \If that $[1,n]=\bcl_{d\in D(n)}[d]$ and $[d]\cap[d']=\vn$ whenever
$d\ne d'$, $d,d'\in D(n)$. Using \era2{7.28} and \era2{1.127} as above we obtain
\beq3.36
\#([1,n])=\sum_{d\in D(n)}\#([d]).
\e
Recall that if $n\in\Na$, then $[1,n]=\zo1,n+1 \nda26.30 \ax \zo0,n $, hence
$\#([1,n])=\#(\zo0,n )$ by \era1{4.1} and by Corollary \rfa2{c3.11}. Thus
$\#([1,n])=\#_\N ([1,n])\nda23.5 = n$.

We next prove
\beq3.37
\#([d]) = \F(\tfrac nd) \qh{\fa}d\in D(n),
\e
which implies by \er{3.36}
\beq3.38
n = \sum_{d\in D(n)}\F\Bg(\frac nd).
\e
Let $d\in D(n)$, then \fe $i\in[d]$ we have $d|i$, hence $i\le n$ by \era3{8.57}.
Thus $[d]=\{i\in M(d)\cap [1,n]: \gcd(i,d)=d\}$, where $M(d)$ is defined in
\era3{8.54}. Set $k:=\frac nd$. If $j\in[1,k]$, then $jd\in M(d)\cap [1,n]$
by \era2{2.16}\,(i) and \era2{2.19}. Conversely, since $d\in\Na$, if
$i\in M(d)\cap[1,n]$, then $d\in D(i)$ by \era2{1.13} and $j:=\frac id\in
[1,k]$ by \era2{2.20}. \If that $[d]=\{jd\in[1,n]: j\in[1,k] \hbox{ and }
\gcd(jd,kd)=d\}$. By \era3{1.27} and \era2{1.12} we have $\gcd(jd,kd)=d\gcd
(j,k)$, $j\in[1,k]$. By \era2{1.13} $d=d1$ and by \era2{2.16}\,(iii), we
obtain $[d]=\{jd\in[1,n]: j\in[1,k] \hbox{ and }\gcd(j,k)=1\}$. The map $j\mt
jd$ from $[1,k]$ into $[d,n]$ is bi\jc\ by \era2{2.19} and \era2{2.20}. Hence
$[d]=\{jd\in[d,n]: j\in[1,k] \hbox{ and }\gcd(j,k)=1\}$ and $[d]\ax \{j\in
[1,k]: \gcd(j,k)=1\}$. By Corollary \rfa2{c3.11} and \er{3.18} we obtain
$\#([d])=\F(k)$. Since $k=\frac nd$, the proof of \er{3.37} is complete.

Finally, we observe that the map $f:D(n)\to D(n)$ defined by $f(i):=\frac di$,
$i\in D(n)$, is a bi\jn. Indeed, $f\circ f=\id_{D(n)}$ since $f(f(i))=f(\frac
di)=\frac d{(\frac di)}=i$ since $\frac di i =d$ \fa $i\in D(n)$. \E\Tf by
\era2{4.71} we obtain $\suml_{d\in D(n)}\F(\frac nd) = \suml_{d\in D(n)}
\F(f(\frac nd)) = \suml_{d\in D(n)}\F(d)$. Thus we proved
\beq3.39n
n=\sum_{d\in D(n)} \F\Bg(\frac nd)= \sum_{d\in D(n)}\F(d)\qh{\fa}n\in\Na.
\e

Since the proof of \era2{4.71} was left as an exercise, we give it here. Let
$(X,\qu,e)$, $\O,\O'$, $f:\O'\to\O$ and $a:\O\to X$ be as in \era2{4.71}, that
is, $(X,\qu,e)$ is an \am, $\O$~and~$\O'$ are nonempty finite \ep\ sets,
$f:\O'\to\O$ is a bi\jn\ and $a$~is a map from~$\O$ into~$X$. We have to show
$\QU_{\o\in\O}a_\o = \QU_{\o'\in\O'} a_{f(\o')}$. Let $\vf:[0,\#(\O)-1]\to\O$
be a bi\jn. Such a bi\jn\ exists since $\O$ is finite and nonempty by Lemma
\rfa2{l3.5}. \Mo $\#(\O')=\#(\O)$ by Corollary \rfa2{c3.11}. Note that
$\o'\mt a_{f(\o')}$ is a map from~$\O'$  into~$X$. We have
$$
\xymatrix{\O'\ar[r]^f &\ \O\  \ar[r]^a &X\\
&\sbk{[0,\#(\O)-1]\\ = \\ [0,\#(\O')-1].} \ar[u]_\vf}
$$
By \E\df\ \rfa2{d7.7}, $\QU_{\o'\in\O'}a_{f(\o')} = \suml_{j=0}^{\#(\O')-1}
a_{f((f\Inv\circ \vf)(j))}$ since $f\Inv\circ\vf : [0,\#(\O')-1] \to\O'$ is a
bi\jn. Since $[0,\#(\O')-1]=[0,\#(\O)-1]$, we obtain $\suml_{j=0}^{\#(\O')-1}
a_{f((f\Inv\circ \vf)(j))} \nda12.1 = \suml_{j=0}^{\#(\O)-1}
a_{(f\circ f\Inv)\circ \vf(j)} = \suml_{j=0}^{\#(\O)-1} a_{\id_\O\circ\vf(j)}
\nda12.2 = \suml_{j=0}^{\#(\O)-1}a_{\vf(j)} \nda27.17 = \QU_{\o\in\O}a_\o$.
This completes the proof of \era2{4.71}.
\ssk

A \gn\ of \er{3.26n} is given in \er{3.72}.

\bex3.31
Let $(G,\qu,e):=(\N_{12},+_{12},0)$, and let $H_d$, $d\in D(12)\sms1$ be the
\sbm\ of~$\N_{12}$ of~$G$ of order~$d$. Show that $D(12)=\{1,2,3,4,6,12\}$,
and that $H_2=\{0,6\}$, $H_3=\{0,4,8\}$, $H_4=\{0,3,6,9\}$, $H_6=\{0,2,4,6,8,
10\}$ and $H_{12}=\zo0,12 $. Let $A_d:=\{x\in H_d: \#(I(x))=d\}$. Show that
$A_2=\{6\}$, $A_3=\{4,8\}$, $A_4=\{3,9\}$, $A_6=\{2,10\}$, $A_{12}=\{1,5,7,
11\}$. Show that $A_d\cap A_{\td d}=\vn$ whenever $d\ne \wt d$, $d,\wt d\in
D(12)$, and $\bcl_{d\in D(n)\sms1}A_d=\zo1,12 $, hence $\zo0,12 =\{0\}
\cup \bcl_{d\in D(n)\sms1}A_d$.

Let $B_d:=\{i\in\zo1,12 : \gcd(i,12)=d\}$, $d\in D(12)$. Show that $B_2=
\{2,10\}$, $B_3=\{3,9\}$, $B_4=\{4,8\}$, $B_6=\{6\}$. Show that $B_d =
A_{\frac{12}d}$. Compute $\F(k)$, $k\in\zo1,12 $, and show that $\suml_{d|12}
\F(d)=12$.
\eex

\bpr3.32
Let $p$ be a \Pn\ and let $k\in\Na$. Then
\beq3.39
\F(p^k)=p^{k-1}(p-1),
\e
\Ip $\F(p)=p-1$.
\epr

\proof
In view of Lemma \rf{l2.16}, we have
\beq3.40
D(p^m)=\{p^i\in\Na: 0\le i\le m\},\q m\in\N.
\e
Hence, by \E\Pr\ \rf{p3.30},
\[
p^k\nde3.26n = \sum_{d\in D(p^k)}\F(d) \nde3.40 = \sum_{0\le i\le k}\F(p^i)
\nad{\era2{7.9},\era2{7.10}}= \sum_{0\le i\le k-1}\F(p^i)+\F(p^k) \nde3.26n =
p^{k-1}+\F(p^k).
\]
Since $p\in\N\sms{0,1}$, $p-1,p\in\N$, $p-1<p$, we have $p^{k-1}<p^k$ by
\era2{2.33}, and $p^k\nda38.72 = (p^k-p^{k-1})+p^{k-1}$. Hence $\F(p^k)+
p^{k-1} \nda21.6 = p^{k-1}+\F(p^k)=p^k= (p^k-p^{k-1})+p^{k-1}$. From
\era2{1.8} we obtain $\F(p^k)=p^k-p^{k-1}$. But $p^k=p^{k-1}\cdot p$, hence
$p^k-p^{k-1}=p^{k-1}\cdot p- p^{k-1}\cdot1 \nda38.76 = p^{k-1}(p-1)$ since
$p>1$. Thus \er{3.39} holds.
\endproof

\brm3.33
\If from \er{3.39} that we have a formula for $\F(n)$, $n\in[2,29]\sms
{6,10,12,14,15,18,20,21,22,24,26,28}$. Notice that $6=2\cdot3$, $10=2\cdot5$,
$12=2^2\cdot3$, $14=2\cdot7$, $18=2\cdot3^2$, $20=2^2\cdot5$, $21=3\cdot7$,
$22=2\cdot11$, $24=2^3\cdot3$, $26=2\cdot13$, $28=2^2\cdot7$. Observe that in
all these cases we have $m=k\cdot l$ where $\gcd(k,l)=1$. In general
$\F(k\cdot l)\ne \F(k)\cdot\F(l)$ since $\F(p^2)\nde3.39 = p(p-1)\ne (p-1)^2
=(\F(p))^2$ when $p$ is prime. It turns out however that
\beq3.41
\F(k\cdot l)=\F(k)\cdot\F(l), \q k,l\in\Na \qh{whenever}\gcd(k,l)=1.
\e
\erm

We now give a proof of \er{3.41}, and of some \gn, based on the \fw\
observations.

We first recall that if $(M,\qu,e)$ is a monoid, then $(M^\t,\qu,e)$ is a
group by Lemma \rf{l3.6}.

\bdf3.34
Let $(M,\qu,e)$ be a monoid. The set of invertible \el s of~$M$ (see \E\df\
\rf{d3.2}) is denoted by~$M^\t$ and the group $(M^\t,\qu,e)$
(also denoted by $(M,\qu,e)^\t$) is called the
\ti{group of invertible \el s of\,}~$M$. If $(R,\qu_0,\qu_1,e_0,e_1)$ is
a~ring, then the invertible \el s of $(R,\qu_1,e_1)$ are usually called the
\ti{units}\index{unit} of the ring~$R$, $(R^\t,\qu_1,e_1)$ is also denoted by~$U(R)$ or
${\rm Unit}(R)$, and $(U(R),\qu_1,e_1)$ is called the \ti{group of units}
of~$R$.
\edf

\blm3.35
Let $(M',\qu',e')$ and $(M'',\qu'',e'')$ be monoids, and let $f:M'\to M''$ be
an \is sm. Then $f$~maps $(M')^\t$ into $(M'')^\t$ and $f|_{(M')^\t} : (M')^\t\to
(M'')^\t$ is an \is sm.
\elm

\proof
Let $a\in(M')^\t$ and let $b\in M'$ be \st $a\qu'b=e'$. Then $f(a)\qu''f(b)=\break
f(a\qu'b)=f(e')=e''$, hence $f(a)\in(M'')^\t$. Thus $f((M')^\t)\sbs(M'')^\t$.
Since $f\Inv$ is an \is sm by
Lemma \rfa2{l1.8}, $f\Inv((M'')^\t)\sbs (M')^\t$, hence $(M'')^\t = f\circ
f\Inv((M'')^\t) \sbs f((M')^\t)$. \E\Tf $f((M')^\t)=(M'')^\t$. Clearly
$f|_{(M')^\t}:(M')^\t \to (M'')^\t$ is a \hm sm and is in\jc. From $f((M')^\t)
=(M'')^\t$ we infer that  $f\Inv((M'')^\t)=(M')^\t$ since $f$~is bi\jc. Let
$y\in(M'')^\t$. Then $x:=f\Inv(y)\in(M')^\t$, and $f|_{(M')^\t}(x)=f(x)=y$.
Hence $f|_{(M')^\t}:(M')^\t \to (M'')^\t$ is sur\jc. Finally, $f|_{(M')^\t}$
is an \is sm by Lemma \rfa2{l1.8}.
\endproof

\blm3.36
Let $(M_i,\qn i,e\en i)$, $i\in[1,n]$, $n\ge2$, be a family of monoids, and
let $\bigl(\prodl_{i\in[1,n]}M_i,\qu,e\bigr)$ be their product defined in
\er{2.11}, \er{2.12}, \er{2.13}. If $a_i\in M_i^\t$, $i\in[1,n]$, and $b_i
\in M_i$ are \st $a_i\qn ib_i=e\en i$, $i\in[1,n]$, then $a,b\in\prodl_{i\in
[1,n]}M_i$ defined by
\beq3.42
a(i):=a_i, \q b(i):=b_i, \q i\in[1,n],
\e
\sf y $a\qu b=e$. Hence $a\in \bigl(\prodl_{i\in[1,n]}M_i\bigr)^\t$. Conversely,
if $a\in\bigl(\prodl_{i\in[1,n]}M_i\bigr)^\t$, then $\pi_i(a)\in M_i^\t$ \fa
$i\in[1,n]$, where $\pi_i$ is defined in \er{2.17}. \E\Ip
\beq3.43
\Bigl(\prod_{i\in[1,n]}M_i\Bigr)^\t = \prod_{i\in[1,n]} M_i^\t.
\e
\elm

\proof
If \er{3.42} holds, then $a\qu b=e$ is a direct con\sq\ of \er{2.12}, \er{2.13},
hence $a\in\bigl(\prodl_{i\in[1,n]}M_i\bigr)^\t$ by \E\df\ \rf{d3.34}. The
converse follows from \E\Pr\ \rf{p2.20}, and \er{3.43} directly follows from
what precedes.
\endproof

\bco3.37 \

\hph i,i, The direct product of a finite family of groups is a group.

\hph ii,, If $\{R_i\}_{i\in[1,n]}$, $n\ge2$, is a family of rings, then
\beq3.44
{\rm Unit}\Bigl(\prod_{i\in[1,n]}R_i\Bigr) = \prod_{i\in[1,n]}{\rm Unit}
(R_i).
\e
\eco

\proof
(i) is a direct con\sq\ of Lemma \rf{l3.36} and the fact that a monoid $M$ is
a group iff $M=M^\t$.

(ii) is a direct con\sq\ of Lemma \rf{l3.36}.
\endproof

The next \Pr\ is a key of the proof of \er{3.41}.

\bpr3.38
In the ring $(\N_n,+_n,\cdot_n,0,1)$, $n\ge2$, the \fw\ holds{\rm:}
\beq3.45
{\rm Gen}((\N_n,+_n,0)) = (\N_n,\cdot_n,1)^\t.
\e
\epr

\proof
We first observe
\beq3.46
k\dpln l = \F_n(k)\cdot_n l, \q k\in\N,\ l\in\zo0,n .
\e
Indeed, $k\dpln l \nad{\er{1.11},\er{1.5}} = k\dpln(l\dpln1) \nad{\rm\era2{2.3}
\,I3,\er{2.8}} = (k\dpln l)\dpln1 \nde1.11 = \F_n(kl)\nde3.26 = \F_n(\F_n(k)
\F_n(l))\nde1.5 = \F_n(\F_n(k)\cdot l)\nde1.33 = \F_n(k)\cdot_n l$.

``$\sbs$'': Let $l\in{\rm Gen}(\N_n,+_n,0)$. Then $I(l)=\zo0,n $. Thus \te s
$k\in\N$ \st $k\dpln l=1$. By \er{3.46} $\F_n(k)\cdot_n l=1$, hence
$l\cdot_n \F_n(k)=1$ in view of \E\Pr\ \rf{p3.15}\,(i). Since $\F_n(k)\in
\N_n$, $l\in(\N_n,\cdot_n,1)^\t$.

``$\supset$'': Let $l\in(\N_n,\cdot_n,1)^\t$. Then  \te s $k\in\N_n$ \st
$l\cdot_n k=1$. Let $m\in\N_n$. Then $m \nde1.5 = \F_n(m)\nda22.13 =
\F_n(1m)\nde1.33 = 1\cdot_n m= (l\cdot_n k)\cdot_n m = m\cdot_n(l\cdot_n k)
= m\cdot_n(k\cdot_nl)
\nad{\rm\era2{2.3}\,I3}= (mk)\cdot_n l\nde3.46 = mk\dpln l$. \If that $\N_n
\sbs I(l)$, thus $I(l)=\N_n$ since $I(l)\sbs \N_n$. \E\Tf $l\in{\rm Gen}
((\N_n,\cdot_n,0))$.
\endproof

\Wanp prove \er{3.41}. Let $k,l\in\N$ be \st $\gcd(k,l)=1$. Then $\F(kl)
\nde3.21 = \#({\rm Gen}(\N_{kl},+_{kl},0)) \nde3.45 = \#((\N_{kl},\cdot_{kl},
1)^\t)$. In view of Theorem \rf{t2.24} with $N:=2$, $n_1:=k$, $n_2:=l$, and
$\gcd(k,l)=1$, the ring $(\N_{kl},+_{kl},\cdot_{kl},0,1)$ is \is c to the
product of the rings $(\N_k,+_k,\cdot_k,0,1)$ and $(\N_l,+_l,\cdot_l,0,1)$.
Hence $\#((\N_{kl},\cdot_{kl},1)^\t)\nde3.43 = \#\bigl((\N_k,\cdot_k,1)^\t
\t(\N_l,\cdot_l,1)^\t\bigr)\nda26.28b = \#((\N_k,\cdot_k,1)^\t)\cdot
\#((\N_l,\cdot_l,1)^\t)\nde3.45 = \#\bigl({\rm Gen}((\N_k,+_k,0))\bigr)\cdot
\#\bigl({\rm Gen}((\N_l,+_l,0))\bigr)\nde3.21 = \F(k)\cdot\F(l)$.
\ssk

Observe that we now have a formula \fa \nm s less than~30. Since $30=2\cdot
3\cdot5$ and $\gcd(2,3)=1$, $\gcd(2,5)=1$ and $\gcd(3,5)=1$, proceding as above
where Theorem \rf{t2.24} is used with $N:=3$, $n_1:=2$, $n_2:=3$ and $n_3:=5$,
we obtain $\F(30)=\F(2\cdot3\cdot5)=\F(2)\F(3)\F(5)=1\cdot2\cdot4=8$. We could
also write $30$ as $6\cdot5$ and observe that $\gcd(6,5)=1$, hence $\F(30)=
\F(6)\cdot\F(5)=(\F(2)\cdot\F(3))\cdot\F(5)=(1\cdot2)\cdot4=8$. Using the full
strength of Theorem \rf{t2.24} we obtain

\bpr3.39
Let $N\ge2$, $n_j\ge2$, $j\in[1,N]$, be \st $\gcd(n_j,n_k)=1$ \fa $j,k\in[1,N]$
\st $j\ne k$, and let $M:=\prodl_{j\in[1,N]}n_j$. Then
\beq3.47
\F(M)=\prod_{j\in[1,N]}\F(n_j).
\e
\epr

\bex3.40
Prove \E\Pr\ \rf{p3.39}.
\eex

If in \E\Pr\ \rf{p3.39} the $n_j$'s are of the form $n_j:=p_j^{m_j}$ with
$p_j$~prime, $m_j\in\Na$, $j\in[1,N]$, and $p_j\ne p_k$ whenever $j\ne k$,
then $\gcd(p_j^{m_j},p_k^{m_k})=1$ in view of Lemma \rf{l2.16}. \E\Tf by
\er{3.47}
\[
\F(M)=\prod_{j\in[1,N]}\F_j(p_j^{m_j})\nde3.39 = \prod_{j\in[1,N]}
p_j^{m_j-1}(p_j-1).
\]

\smallskip
In view of Theorem \rfa3{t9.22}, every \nn\ greater than one is a (finite)
product of powers of primes. \E\Tf we have

\bth3.41
Let $n\in\Na$. Then
\beq3.48
\F(n)= \prod_{\sbk{p\in \cP\\p|n}} p^{m_p-1}(p-1),
\e
where $\cP$ is the set of \Pn s, $m_p$ is the only \po\ \nn\ \st $p^{m_p}$
divides~$n$ and $p^{m_p+1}$ does not divide~$n$.
\eth

\bex3.42
Prove Theorem \rf{t3.41}.
(Hint: Show that the set $A:=\{p\in\cP: p|n \hbox{ and }{p{\ne} n}\}$ is finite
and find a~bi\jn\ from $[1,\#(A)]$ onto~$A$.)
\eex

\bex3.43
Compute $\F(100)$, $\F(120)$ and $\F(155001)$.
\eex

We conclude this part of Section \ref{sss.fcg} dedicated to the study of
``property~C'' by stating \ev t \ch izations of a \Gn\ of $(\N_n,+_n,0)$,
$n\ge2$.

\bpr3.44
Let $n\ge2$ and let $(\N_n,+_n,\cdot_n,0,1)$ be the ring introduced in \E\Pr\
\rf{p3.15}\,{\rm(i)}. Let $k\in(0,n)$, then the \fw\ assertions are \ev t\/{\rm:}

\hph i,ii, $k$ is a \Gn\ of the group $(\N_n,+_n,,0)$,

\hph ii,i, $\#(I(k))$, the order of $k$, is equal to $n$,

\hph iii,, $\gcd(k,n)=1$,

\hph iv,, $k$ is a unit of the ring $(\N_n,+_n,\cdot_n,0,1)$,

\hph v,i, the map $\d_k:(\N_n,+_n,0) \to (\N_n,+_n,0)$ defined by
\beq3.49
\d_k(l):= k\cdot_n l, \q l\in\N_n,
\e
is an auto\mf.
\epr

\bex3.45
Prove \E\Pr\ \rf{p3.44}.
\eex

We now consider \ev t \df s of a group.

\bdf3.46
Let $(M,\qu)$ be a \sg, and let $a\in M$. The maps $\la_a:X\to X$ (resp.\
$\rho_a:X\to X$) defined by
\beq3.50
\la_a(x):=a\qu x, \q \rho_a(x):=x\qu a, \q x\in X,
\e
are called \ti{left-} (resp.\ \ti{right-}) \ti{\tl s} by~$a$. If the
\op~$\qu$ is \cmt e, we use the notation $\th_a$ instead of
$\la_a\,(=\rho_a)$.
\edf

\blm3.47
Let $(X,\qu)$ be a \sg. Set
\beq3.51
x\tqu y=y\qu x, \q x,y\in X.
\e
Then

\hph i,i, $(X,\tqu)$ is a \sg, called the \emph{converse} \sg\ of $(X,\qu)$.
\E\fe $a\in X$ $\la_a$ $($resp.\ $\rho_a)$ is the right- $($resp.\ left-$)$
\tl\ by~$a$ in $(X,\tqu)$.

\hph ii,, If $(X,\qu)$ is a monoid $($resp.\ group$)$, so is $(X,\tqu)$.
\elm

\bex3.48
Prove Lemma {\rm\rf{l3.47}}.
\eex

\blm3.49
Let $(X,\qu,e)$ be a monoid, and let $a,b\in X$ be \st $a\qu b=e$. Then
$\rho_a,\la_b$ are in\jc\ and $\la_a,\rho_b$ are sur\jc.
\elm

\proof \

\ti{In\ji\ of $\rho_a$}: Let $x,y\in X$ be \st $x\qu a\nad*= y\qu a$. Then
$x=x\qu e={x\qu(a\qu b)}=(x\qu a)\qu b \nad*=(y\qu a)\qu b=y\qu(a\qu b)=
y\qu e=y$.

\ti{Sur\ji\ of $\la_a$}: Let $z\in X$. Then $z=e\qu z=(a\qu b)\qu z=a\qu
(b\qu z)={\la_a(b\qu z)}$.
\endproof

\bex3.50
Prove the in\ji\ of $\la_b$ and the sur\ji\ of~$\rho_b$ in Lemma \rf{l3.49}.
\eex

\bpr3.51
Let $(X,\qu,e)$ be a group. Then \emph{all} left- and right-\tl s are bi\jc.
\epr

\proof
Let $x\in X$. Then, by Lemma \rf{l3.49} and $x\qu y\nde3.3 = e$ \fs $y\in X$,
$\rho_x$~is in\jc\ and $\la_x$~is sur\jc. Similarly, from $y\qu x\nde3.3 = e$,
we infer that $\la_x$ is in\jc\ and $\rho_x$ is sur\jc.
\endproof

We now consider a stronger version of \E\Pr\ \rf{p3.51}.

\bth3.52
Let $(X,\qu)$ be a \sg. Suppose that all left- and right-\tl s are
\emph{sur\jc}. Then \te s a \nel\ $e\in X$ \st $(X,\qu,e)$ is a group.
\E\Ip all left- and right-\tl s are bi\jc, by \E\Pr\ \rf{p3.51}.
\eth

\proof
Let $c\in X$. Since $\rho_c$ is sur\jc, \te s $e_l\in X$ \st
\beq3.52
e_l\qu c=c.
\e

We now prove that
\beq3.53
e_l\qu a=a \hbox{\ \it \fa} a\in X.
\e
Let $a\in X$. \E\te s $b\in X$ \st $c\qu b=a$ since $\la_c$ is sur\jc. Then
$e_l\qu a=e_l\qu (c\qu b) = (e_l\qu c)\qu b \nde3.52 = c\qu b=a$.

\ti{$a\qu e_l=a$ \fa $a\in X$}: Let $a\in X$. \E\te s $b\in X$ \st $b\qu a
\nad*=e_l$ since $\rho_a$ is sur\jc. \E\te s $c\in X$ \st $c\qu b=e_l$ since
$\rho_b$~is sur\jc. Then $a\qu b=e_l\qu(a\qu b)$ by \er{3.53}. \Mo $e_l\qu
(a\qu b)= (c\qu b)\qu(a\qu b)\nda21.36  = \break c\qu(b\qu a)\qu b=c\qu e_l\qu b=
c\qu(e_l\qu b)\nde3.53 = c\qu b=e_l$. Hence $a\qu b=e_l$. \E\Tf $a\qu e_l\nad*=
a\qu(b\qu a)=(a\qu b)\qu a= e_l\qu a=a$.

\If that $e_l$ is a (the) \nel\ of $(X,\qu)$, hence $(X,\qu,e_l)$ is a monoid.
We set $e:=e_l$.

Let $x\in X$. Since $\la_x$ is sur\jc, \te s $y\in X$ \st $x\qu y=e$.
Similarly, since $\rho_x$ is sur\jc\ \te s $z\in X$ \st $z\qu x=e$. But
$y=e\qu y =(z\qu x)\qu y = z\qu(x\qu y)=z\qu e=z$. Hence $y=z$ and $(X,\qu,e)$
is a group since all \el s of~$X$ are invertible.
\endproof

\brm3.53
If $(X,\qu,e)$ is a monoid \st all \tl s are sur\jc, then by Theorem \rf{t3.52}
\te s $\wt e\in X$ \st $\wt e$~is a \nel\ of $(X,\qu)$. But $e=\wt e\qu e=
\wt e$ since both $e$~and~$\wt e$ are \nel s. \E\Tf $(X,\qu,e)$ is a group.
In what follows we show that if $(X,\qu,e)$ is a monoid and all left-\tl  s
are sur\jc, then $(X,\qu,e)$ is a group. Similarly for the right-\tl s.
\erm

\bpr3.54
Let $(X,\qu,e)$ be a monoid. Then the \fw\ assertions are \ev t.

\hph i,ii, $(X,\qu,e)$ is a group $($see \E\df\ \rf{d3.2}$)$.

\hph ii,i, All \tl s of $(X,\qu,e)$ are bi\jc.

\hph iii,, All \tl s of $(X,\qu,e)$ are sur\jc.

\hph iv,, All left-\tl s are sur\jc.

\hph v,i, All right-\tl s are sur\jc.

\epr

\proof \

(i) \ti{implies} (ii): See \E\Pr\ \rf{p3.51}.

(ii) \ti{implies} (iii) and (iii) \ti{implies} (iv): Trivial.

(iv) \ti{implies} (i): It suffices to prove that every \el\ of~$X$ is
invertible. Let $x\in X$. \E\te s $y\in X$ \st $x\qu y=e$ since $\la_x$ is
sur\jc. \Mo \te s $z\in X$ \st $y\qu z=e$ by the sur\ji\ of~$\la_y$. Then
$y\qu x= (y\qu x)\qu e = (y\qu x)\qu(y\qu z) \nda21.36 = y\qu(x\qu y)\qu z =
y\qu e\qu z= y\qu z= e$. Hence \er{3.3} holds.

\If that (i), (ii), (iii) and (iv) are \ev t. We now show that (i) is \ev t
to~(v).

(i) \ti{implies} (v): By what precedes (i) implies (iii) which trivially
implies~(v).

(v) \ti{implies} (i): Let $\tqu$ be the binary \op\ of~$X$ defined in
\er{3.51}. Then $(X,\tqu,e)$ is a monoid by Lemma \rf{l3.47}. Let $\wt\la_x$
(resp.~$\wt\rho_x$), $x\in X$, denote the left- (resp.\ right-) \tl\ by~$x$
in $(X,\tqu,e)$. By Lemma \rf{l3.47}, $\wt\la_x=\rho_x$ and $\wt\rho_x=\la_x$,
$x\in X$. \E\Tf if $\rho_x$ is sur\jc\ \fa $x\in X$, so is $\wt\la_x$ \fa
$x\in X$. Since (iv) implies (i), we infer that $(X,\tqu,e)$ is a~group.
Notice that
$\raise1pt\hbox{\hbox
to0pt{$\tilde{\tilde{\scriptstyle\sqcap}}$\hss}$\scriptstyle\sqcup$}={\qu}$,
hence the converse of $(X,\tqu,e)$ is $(X,\qu,e)$ which is a~group by Lemma
\rf{l3.47}.
\endproof

In the next \Pr\ the \cn\ of the \ex\ of a \nel\ in the semigroup $(X,\qu)$ is
weakened.

\begin{prp}[\cite{Alg}]\lb{p3.55}
Let $(X,\qu)$ be a semigroup.

\hph i,i, Suppose \te s $e_l\in X$ \st
\beq3.54
e_l\qu x=x \qh{\fa} x\in X.
\e
If \fe $x\in X$ \te s $y\in X$ \st
\beq3.55
y\qu x=e_l,
\e
then $e_l$ is a \nel\ for the \op~$\qu$ and $(X,\qu,e_l)$ is a group.

\hph ii,, Similarly, if \te s $e_r\in X$ \st $x\qu e_r=x$ \fa $x\in X$, and
\fa $x\in X$ \te s $y\in X$ \st
\beq3.56
x\qu y=e_r,
\e
then $(X,\qu,e_r)$ is a group.
\epr

\proof[Proof of {\rm(i)}]
We first show that if $x\in X$ and $y\in X$ \sf y \er{3.55}, then $x\qu
y=e_l$. By \as\ \er{3.55} \te s $z\in X$ \st $z\qu y=e_l$. Hence $x\qu y
\nde3.54 = e_l\qu{(x\qu y)}= (z\qu y)\qu(x\qu y)\nda21.36 = z\qu(y\qu x)\qu y
\nde3.55 = z\qu e_l\qu y =z\qu y =e_l$. \If that \fe $x\in X$,
\te s $y\in X$ \st $x\qu y = y\qu x=e_l$. \E\Tf it remains to show that $e_l$
is a \nel\ of $(X,\qu)$. Let $x\in X$ and $y\in X$ \sf y \er{3.55} and $x\qu
y =e_l$. Then $x\qu e_l = x\qu(y\qu x)=(x\qu y)\qu x=e_l\qu x\nde3.54 = x$.
Since $e_l\qu x=x$ \fa $x\in X$, $e_l$~is a (hence the) \nel\ of $(X,\qu)$.
Since $(X,\qu,e_l)$ is a monoid \st all \el s are invertible, $X$~is a~group.
However, if for \ti{some} $x\in X$ \te s $y\in X$ \st $x\qu y=e$, then $y\qu x$
is not necessarily equal to~$e$.
\endproof

\bex3.56
Prove part (ii) of \E\Pr\ \rf{p3.55}.
\eex

\brm3.57
\E\Pr\ \rf{p3.55} shows that in \E\df\ \rf{d3.2} \cn s $x\qu y=e$ and $y\qu
x=e$  of \er{3.3} can be replaced by one of them.
\erm

Our next goal is to show that every group $(G,\qu,e)$ is monoid-\is c to
a~\ti{subgroup} of the group $(S_G,\circ,\id_G)$ (see Examples \rf{xa3.7}).

Let $(X,\qu,e)$ be a group and set
\beq3.57
L_X := \{\la_x\in S_X: x\in X\}.
\e
We have $\id_X=\la_e\in L_X$. \Mo if $a,b\in X$, then $\la_a\circ \la_b=
\la_{a\qum b}$, where $\circ$ denotes the \cm\ of $\la_a$ and~$\la_b$. Indeed,
\fe $x\in X$ we have $(\la_a\circ\la_b)(x) = \la_a(\la_b(x)) =\break a\qu(b\qu x) =
(a\qu b)\qu x=\la_{a\qum b}(x)$. \If that $L_X$ is a \sbm\ of~$S_X$, hence
$(L_X,\qu,e)$ is a monoid. \E\fe $a\in X$, $\la_a\circ \la_{a\Inv}= \la_{a
\qum a\Inv}=\la_e=\id_X$. Similarly $\la_{a\Inv}\circ\la_a=\id_X$, hence
$\la_{a\Inv}$ is the inverse of $\la_a$ in $(S_X,\circ,\id_X)$. \E\Tf
$(L_X,\circ,\id_X)$ is a group, hence $L_X$ is a subgroup of $(S_X,\circ,\id_X)$.

Define a map $L:X\to L_X$ by setting
\beq3.58
L(x):= \la_x, \q x\in X.
\e
Then we have $L(e)=\la_e=\id_X$ and for $a,b\in X$, $L(a\qu b)=\la_a\circ\la_b
= L(a)\circ L(b)$. Hence $L(X,\qu,e)\to (L_X,\circ,\id_X)$ is a monoid-\hm sm.
\Mo $L$~is in\jc. Indeed, if $L(a)=L(b)$, $a,b\in X$, then $\la_a=\la_b$,
hence $a=\la_a(e)=\la_b(e)=b$. By the \df\ of~$L_X$, $L$~is sur\jc. Thus $L$~is
a monoid-\is sm by Lemma \rfa2{l1.8}\,(i).

We summarize the results above in

\begin{prp}[Cayley \cite{Groups}]\lb{p3.58}
Every group $X$ is \is c to a subgroup of~$S_X$. More precisely, the map $L:X\to
L_X$ defined in \er{3.58} is a monoid-\is sm.
\epr

\bex3.59 \

\hph i,vii, Let $X$ be a finite set with $\#(X)=2$. Show that the group
$(S_X,\circ,\id_X)$ is \is c to the \ag\ $(\N_2,+_2,0)$.

\hph ii,vi, Show that if a set $X$ contains at least three \el s, then $S_X$
is \ti{not\/} abelian.

\hph iii,v, Show that if $n\ge3$ there is \ti{no} set $X$ \st $(\N_n,+_n,0)$
is \is c to~$S_X$.

\hph iv,ii, Let $X,Y$ be \ns s and let $\vf:X\to Y$ be a bi\jn. Let $\Phi:S_X
\to S_Y$ be the map defined by
\beq3.59
\Phi(f):=\vf\circ f\circ\vf\Inv, \q f\in S_X.
\e
Show that $\Phi$ is an \is sm from $(S_X,\circ,\id_X)$ onto
$(S_Y,\circ,\id_Y)$.

\hph v,iii, Let $(X,\qu,e)$ be a group, let $(X',\qu',e')$ be a monoid and
let $\vf:X\to X'$ be a (monoid-)\is sm. Show that $(X',\qu',e')$ is a group.

\hph vi,ii, Let $X$ be a \nf\ with $\#(X)\ge3$. Show that $\#(X)<\#(S_X)$.

\hph vii,i, If $X:=\N$, then $S_X$ is \ti{infinite} since it contains the
infinite set of \ti{\tp s} $\si_{a,b}$, $a,b\in\N$, $a\ne b$, defined
by $\si_{a,b}(x):=b$ for $x:=a$, $\si_{a,b}(x):=a$ for $x:=b$ and
$\si_{a,b}(x):=x$ \fa $x\in\N\sms{a,b}$. Clearly $\si_{a,b}\circ \si_{b,a}=
\id_\N$, hence $\si_{a,b}\in S_X$ and $(\si_{a,b})\Inv=\si_{b,a}=\si_{a,b}$.
Is $S_\N$ uncountable?

\hph viii,, Let $(X,\qu,e)$ be a group and let $(S_X,\wt\circ,\id_X)$ be the
converse group of $(S_X,\circ,\id_X)$. Let $R_X:=\{\rho_x\in S_X: x\in X\}$.
Show that $R_X$ is a subgroup of $(S_X,\wt\circ,\id_X)$. Define $R:X\to R_X$ by
setting $R(x):=\rho_x$, $x\in X$. Show that $R$~is a monoid-\is sm from
$(X,\qu,e)$ onto $(R_X,\wt\circ,\id_X)$.
\eex

We proved in Lemma \rf{l3.21} that every finite \sbm\ of an \ag\ is a~subgroup.
We now show that the \as\ ``abelian'' can be removed (see Remark \rf{r3.22}).

\blm3.60
Let $(G,\qu,e)$ be a group and let $M$ be a finite \sbm\ of~$G$. Then
$(M,\qu,e)$ is a group.
\elm

\proof
Clearly $(M,\qu,e)$ is a finite monoid. Let $a\in M$ and let $\la_a$
(resp.~$\rho_a$) be the left- (resp.\ right-) \tl s of~$G$ defined in \er{3.50}.
Since $M$~is a monoid, $\la_a|_M$ (resp.\ $\rho_a|_M$), the \rt ion of $\la_a$
(resp.~$\rho_a$) to~$M$ maps $M$ into~$M$. Since $G$ is a~group,
$\la_a$~and~$\rho_a$ are in\jc\ by \E\Pr\ \rf{p3.51}. So are $\la_a|M$
and~$\rho_a|_M$. Since $M$~is finite, $\la_a|M$
and~$\rho_a|_M$ are bi\jc\ by \E\Pr\ \rfa2{p3.13}. Since $\la_a|_M$ (resp.\
$\rho_a|_M$) are left- (resp.\ right-) \tl s of $(M,\qu,e)$, $(M,\qu,e)$
is a group by \E\Pr\ \rf{p3.54}.
\endproof

In \E\Pr\ \rf{p3.24} we proved that if $H$ is a nontrivial \sbm\ of a \pn\
group~$G$, then $\#(H)$ divides $\#(G)$.  The next theorem, usually called
\ti{Lagrange's theorem}, shows that this result holds for every finite
group~$G$.

\bth3.61
If $G$ is a finite group and $H$~is a nontrivial \sbm\ of~$G$, then $\#(H)$
divides $\#(G)$.
\eth

\proof
Let $H$ be a nontrivial \sbm\ of~$G$. By Lemma \rf{l3.60}, $H$~is a subgroup
of~$G$, and by Lemma \rfa2{l3.16} $H$~is finite with $\#(H)\le \#(G)$.
By \E\Pr\ \rf{p3.51}, $\la_x(H)$ is \ep\ to~$H$ \fa $x\in G$ and by
Corollary \rfa2{c3.11}, $\#(\la_x(H))=\#(H)\ne 0$. Set $\cA:=\{A\sbs G:
\hbox{\te s $x\in G$ \st}A=\la_x(H)\}$. Since $\cA\sbs \cP(G)$, the power
set of~$G$, and $G$~is finite, $\cA$~is finite by Exercise \rfa2{ex3.27} and
Theorem \rfa1{t4.18}\,(i). $\cA\ne\vn$, since $H=\la_e(H)\in \cA$. \Mo $\#(A)
=\#(H)$ \fa $A\in\cA$ by what precedes. We claim that $\cA$~is a \pt\ of~$G$
(see \E\df\ \rfa1{d4.5}). Indeed, \fe $A\in\cA$, \te s $x\in G$ \st $A=\la_x
(H)$, hence by Corollary \rfa2{c3.11} $\#(A)=\#(\la_x(A))=\#(H)\ne0$. Thus $A\ne\vn$ by \era2{3.6}.
\Mo \fe $x\in G$, $x\in\la_x(H)$, since $x=x\qu e$ and $e\in H$. Hence \fe
$x\in G$, \te s $A\in\cA$ \st $x\in A$. Thus $G\sbs \bcl_{A\in \cA}A$. Clearly
$\bcl_{A\in \cA}A\sbs G$, hence $G=\bcl_{A\in \cA}A$. Finally, we show that \fa
$A,B\in\cA$ \st $A\ne B$, we have $A\cap B=\vn$, \ev tly, \fa $A,B\in \cA$ \st
$A\cap B\ne\vn$, we have $A=B$. Let $A,B\in\cA$ and $z\in A\cap B$. Then \te\
$x,y\in G$ \st $A=\la_x(H)$, $B=\la_y(H)$ and $z\in G$ \st $z\in \la_x(H)\cap
\la_y(H)$. \csq, $z=x\qu h$, $z=y\qu\wt h$ \fs $h,\wt h\in H$. Since
$(H,\qu,e)$ is a group, \te s $h\Inv\in H$ \sf ying $h\qu h\Inv=h\Inv\qu h
=e$. Thus $x\nda21.7 = x\qu e = x\qu(h\qu h\Inv)\nda21.5 = (x\qu h)\qu
h\Inv = (y\qu \wt h)\qu h\Inv \nda21.5 = y\qu(\wt h\qu h\Inv)$. \Mo $x\qu h'
= (y\qu(\wt h\qu h\Inv))\qu h'\nda21.5 = y\qu((\wt h\qu h\Inv)\qu h')$ \fe $h'\in H$.
Since $H$~is a \sbm\ of~$G$, $(\wt h\qu h\Inv)\qu h'\in H$, we infer that
$\la_x(H)\sbs \la_y(H)$. Hence $A\sbs B$. Interchanging $A$ and~$B$, we find
$B\sbs A$, hence $A=B$. This completes the proof of the claim that $\cA$~is
a~\pt\ of~$G$. In view of \era2{7.28}, \era2{1.127} with $\O:=G$, $I:=\cA$,
$i:=A\in\cA$, $a_\o:=1$ \fa $\o\in\O$, and $(X,\qu,e):=(\N,+,0)$, we obtain
$\#(G) = \suml_{A\in\cA}\#(A) =  \suml_{A\in\cA}\#(H) \nda21.126 =
\#(\cA)\dpl \#(H)\nda22.8 = \#(\cA)\#(H)$.
Hence
\beq3.60
\#(G)= \#(\cA)\#(H).
\e
This completes the proof of the theorem.
\endproof

\brm3.62 \

\hph i,ii, Subsets of $G$ of the form $\{x\qu h: h\in H\}$, $x\in G$, are called
\ti{left-cosets} of~$H$ in~$G$. \If from the proof of Theorem \rf{t3.61} that
the \nm\ of left-cosets of~$H$ in~$G$ is equal to $\frac{\#(G)}{\#(H)}$, called
the \ti{index} of~$H$.

\hph ii,i, Similarly sets of the form $\{h\qu x: h\in H\}$, $x\in G$, are called
\ti{right-cosets} of~$H$ in~$G$. Using Lemma \rf{l3.47} we find that
$\frac{\#(G)}{\#(H)}$ is also the \nm\ of right-cosets of~$H$ in~$G$. Left-
(resp.\ right-) cosets of~$H$ in~$G$ are usually denoted by $xH$ (resp.~$Hx$),
$x\in G$.

\hph iii,, The set of left- (resp.\ right-) cosets of~$H$ is a \pt\ of~$G$
\crs\ to the \ev ce \rl\ $x\sim y$ if $y\Inv \qu x\in H$ or $\la_x(H) =
\la_y(H)$, $x,y\in H$ (resp.\ $x\sim y$ if $x\qu y\Inv\in H$ or $\rho_x(H) =
\rho_y(H)$, $x,y\in H$).
\erm

We consider a first application of Theorem \rf{t3.61} in the next exercise.

\bex3.63
Let $\zb ni{[1,m]}$, $m\ge2$, be a finite \sq\ of \el s of~$\Na$. Then
\beq3.61
\prod_{i=1}^m (n_i!) \qh{divides } \Bigl(\sum_{i=1}^m n_i\Bigr)!.
\e
\eex

The aim of this exercise is to give a proof of \er{3.61} using Lagrange's
theorem. We proceed in several steps.

(i) Give a \pt\ $\cA:=(A_i)_{i\in[1,m]}$ of $\bigl[1,\suml_{i=1}^m n_i\bigr]$
\st $\#(A_i)=n_i$ \fa $i\in[1,m]$.

(ii) Set $X:=\bigl[1,\suml_{i=1}^m n_i\bigr]$ and let $(S_X,\circ,\id_X)$ be
the group of \Pm s of~$X$. Let $S_{X,\cA}$ denote the subset of~$S_X$
consisting of bi\jn s of~$X$ leaving $A_i$, $i\in[1,m]$, invariant, i.e.
\beq3.62a
S_{X,\cA}:=\{f\in S_X: f(A_i)=A_i,\ i\in[1,m]\}.
\e
Show that $S_{X,\cA}$ is a \sbm\ of $S_X$.

(iii) Show that $S_{X,\cA}$ is \ep\ to $\prodl_{i=1}^m S_{A_i}$.

(iv) Use \era2{4.50}, \er{2.16} and Lagrange's theorem to prove \er{3.61}.

(v) Show that \er{3.61} holds if $n_i\in\N$, $i\in[1,m]$, $m\in\Na$, with
$0!=1$.

\begin{ntt}[multiindex notation]\lb{n3.64}
Let $\a:[1,n]\to\N$, $n\in\Na$, and let
\beq3.62
|\a|:=\sum_{i=1}^n \a_i \qh{and }\a! = \prod_{i=1}^n \a_i!.
\e
Then the \fw\ notation is used:
\beq3.63
\binom{|\a|}{\a!}:=\frac{|\a|!}{\a!} = \frac{\Bigl(\sum_{i=1}^n \a_i\Bigr)!}
{\prod_{i=1}^n \a_i!}\,.
\e
The map $\a$ is usually called a \ti{multiindex} and $\binom{|\a|}{\a!}$ is
called a \ti{multinomial \cf\/} (for reasons which will be given in the next
section).

In the case $n=2$, another notation is used:
\beq3.64
\binom{\a_1+\a_2}{\a_1}:=\binom{|\a|}{\a!}=\frac{(\a_1+\a_2)!}{\a_1!\a_2!}\,.
\e
\E\Ip if $n:=\a_1+\a_2$ and $k=\a_1$, then\glossary{$\binom nk$}
\beq3.65
\binom nk :=\frac{n!}{k!(n-k)!}\,.
\e
The \nm s $\binom nk$ are called binomial \cf s.
\ent

We now consider some other con\sq s of Theorem \rf{t3.61}.

\bpr3.65
Let $(G,\qu,e)$ be a finite group with $|G|:=\#(G)>1$.

\hph i,i, \E\fe $a\in G$, the $|G|$-th \IT\ of~$a$ is equal to~$e$.

\hph ii,, If $|G|$ is prime, then $G$ is \is c to $(\N_{|G|},+_{|G|},0)$.
\epr

\proof \

(i) Since $|G|>1$, \te s $a\in G\sms e$. Let $I(a)$ denote the set of \IT s
of~$a$ in $(G,\qu,e)$. Since $G$~is finite, $I(a)$~is finite by Theorem
\rfa1{t4.18}\,(i). By Lemma \rfa2{l1.21}\,(ii) with $M_0=M=G$, $I(a)$~is an abelian
\sbm\ of~$G$. \E\Tf $\#(I(a))$ divides~$|G|$ by Theorem \rf{t3.61}. Let
$k\in\Na$ be \st $|G|=k\cdot\#(I(a))$. We claim that $\#(I(a))\dqu a=e$.
Observe that \fe $l\in\N$, the $l$-th \IT\ of~$a$ in $(I(a),\qu,e)$ is equal
to the $l$-th \IT\ of~$a$ in $(G,\qu,e)$. Conversely, if $l\dqu a$ denotes
the $l$-th \IT\ of~$a$ in $(G,\qu,e)$, then $l\dqu a\in I(a)$ by Lemma
\rfa2{l1.21} with $M:=G$ and $M_0:=I(a)$, thus $l\dqu a$ is also the $l$-th
\IT\ of~$a$ in $(I(a),\qu,e)$. \E\Tf $(I(a),\qu,e)$ is a \pn\ monoid (see
\E\df\ \rfa2{d1.18}). \Mo as an abelian \sbm\ of a group it is a \Cm. Indeed,
if $a,b,c\in I(a)$ \sf y $a\qu c=b\qu c$, then $a=b$, since all \tl s are
in\jc. \If that $(I(a),\qu,e)$ is a finite \pn\ \Cm. By Lemma \rf{l2.12} we
find that $\#(I(a))\dqu a=e$ since $\#(I(a))$ is the order of~$a$. This \es es
the claim. Finally, we have $|G|\dqu a= (k\cdot \#(I(a)))\dqu a
\nad{\rm\era2{2.3}\,I3} = k\dqu (\#(I(a))\dqu a) = k\dqu e \nad{\rm
\era2{2.3}\,I4} = e$.

(ii) As in part (i), let $a\in G\sms e$ and let $I(a)$ denote the set of \IT
s of~$a$ in $(G,\qu,e)$. Since $I(a)$ is a nontrivial \sbm\ of~$G$,
$\#(I(a))$ divides~$|a|$ by Theorem \rf{t3.61}. Since $|G|$ is prime, either
$\#(I(a))=1$ or $\#(I(a))=|G|$. But $\#(I(a))>1$ since $e,a\in I(a)$, $e\ne
a$. \E\Tf $\#(I(a))=|G|$. Since $I(a)\sbs G$, we infer $I(a)=G$ by
\era2{3.13}. In view of part~(i), $G$~is a \pn\ group hence (ii) follows from
part~(i) of the proof of Theorem \rf{t1.28} with $e_1$ replaced by~$a$.
\endproof

The Little Fermat Theorem states that if $p$ is a \Pn\ and if $a$ is an \el\
of $\Na$ which is not a \ml e of~$p$, then
\beq3.66
a^{p-1}=kp+1 \qh{\fs}k\in\N,
\e
in other words, $p$ divides $a^{p-1}-1$.

\bex3.66
Show that if $a$ is a \ml e of an \el\ $q\in\Na$, then $a^m$ is also a \ml e
of~$q$, and that \er{3.66} does not hold if $a$ is a \ml e of~$p$.
\eex

The \cn\ ``\ti{$a$ is not a \ml e of~$p$}'' can be rewritten as $\F_p(a)\ne0$
in view of \er{1.2}, \er{1.3}. By \er{1.2} and \er{1.3}, \er{3.66}
can be rewritten as
\beq3.67
\F_p(a^{p-1})=1.
\e

Euler proved the \fw\ \gn\ of the Little Fermat Theorem: if $n\in\Na\sms1$
and if $a\in\Na$ \sf ies $\gcd(a,n)=1$, then
\beq3.68
\F_n(a^{\F(n)})=1,
\e
where $\F(n)$ is the $\F$-Euler \f\ defined in \er{3.18}.

Note that if $n=p$ is prime, then $\F(p)=p-1$ by \er{3.25}. \Mo if
$a\in\Na$ and $\gcd(a,p)=1$, then $a$ is not a \ml e of~$p$, since if $a=kp$
\fs $k\in\Na$, then $\gcd(kp,p)= p>1$. \E\Tf Euler's theorem is a \gn\ of
Fermat's theorem. We next show that \er{3.68} is a con\sq\ of \E\Pr\
\rf{p3.65}.

\proof[Proof of Fermat--Euler's theorem]
We first show
\beq3.69
\F_n(a^k) = \F_n((\F_n(a))^k), \q n\ge2,\ k,a\in\Na.
\e
By \er{1.53} we have $\F_n(a^k)=k\ddtn \F_n(a)$. Since $\F_n(\F_n(a))=\F_n(a)$
by \er{1.3}, \er{1.5}, we obtain $k\ddtn \F_a(a) = k\ddtn \F_n(\F_n(a))\nde
1.53 = \F_n((\F_n(a))^k)$, which proves \er{3.69}.

We now suppose that $n\in\Na\sms1$, $a\in\Na$ and $\gcd(a,n)=1$, and  show
\beq3.70
\gcd(\F_n(a),n)=1.
\e
Indeed, by \er{1.2}, \er{1.3}, \te s $q\in\N$ \st $a=qn+\F_n(a)$. Thus
\era3{8.66}, \era3{8.67} hold with $a_0:=a$, $q_0:=q$, $a_1:=n$ and
$a_2:=\F_n(a)$. We show that $\F_n(a)\ne0$. Indeed, if $\F_n(a)=0$, then
$a=qn$. Since $a\ne0$, we infer $q\ne0$ from \era2{2.29}, \era2{2.12}. Since
$n|qn$, we obtain $\gcd(a,n)=\gcd(qn,n)=n>1$ \cd ing $\gcd(a,n)=1$. Since
$a_2\ne0$ in \era3{1.19}, we deduce from \era3{8.69} that $\gcd(a,n)=
\gcd(n,\F_n(a))$. But $\gcd(n,\F_n(a))=\gcd(\F_n(a),n)$ by \era3{1.24}. \E\Tf
\er{3.70} holds.

In view of \er{1.53} it remains to prove
\beq3.71
\F(n)\ddtn \F_n(a)=1 \qh{if}\gcd(a,n)=1.
\e
We recall $\{k\in[1,n]:\gcd(k,l)=n\}\nde3.20 = {\rm Gen}((\N_n,+_n,0))
\nde3.45 = (\N_n,\cdot_n,1)^\t$, and $(\N_n,\cdot_n,1)^\t$ is a group by
Lemma \rf{l3.6}. Thus if $\gcd(a,n)=1$, then $\gcd(\F(a),n)=1$ by what precedes,\glossary{$(\N_n,\cdot_n,1)$}
hence $\F(a)\in(\N_n,\cdot_n,1)^\t$. \Mo by \er{3.18} $\F(n)=\#((\N_n,\cdot_n,
1)^\t)$ and \er{3.71} holds by \E\Pr\ \rf{p3.65}\,(i).
\endproof

\bxa3.66
Let $a,m\in\Na\sms1$. We want to compute the last digit of~$a^m$, i.e.\
$\F_{10}(a^m)$. Set $\a:=\F_{10}(a)$. In view  of \er{3.69} it is \sft\ to
compute $\F_{10}(\a^m)$, where $\a\in[0,9]$. If $\a:=0$, then $\a^m=0$ by
\era2{2.29}, and $\a^m=1$ if $\a:=1$ by \era2{2.27}. If $\a:=5$, then
$\a^2=25$, hence $\F_{10}(5^2)=5$. More generally, if $\F_{10}(a^m)=5$ \fs
$m\ge 2$, then $\F_{10}(\a^{m+1}) \nad{\era2{2.25},\era2{2.24}}=
\F_{10}(\a^m\cdot\a) \nde3.26 = \F_{10}(\F_{10}(\a^m)\cdot\F_{10}(\a))
\nde1.5 = \F_{10}(5\cdot\a) = \F_{10}(25) = 5$. Using \In\ we find
$\F_{10}(5^m)= 5$ \fa $m\ge2$. Similarly one shows that $\F_{10}(6^m)=6$ \fa
$m\ge2$, since $\F_{10}(6\cdot 6)=\F_{10}(36)=6$. It remains to compute
$\F_{10}(\a^m)$ for $\a\in\{2,3,4,7,8,9\}$.

If $\a\in\{3,7,9\}$, then $\gcd(\a,10)=1$. Since $10=2\cdot5$ and
$\gcd(2,5)=1$, $\F(10)\nde3.47 = \F(2)\cdot\F(5)\nde3.25 = 1\cdot4=4$. Hence
$\F_{10}(\a^4)=1$ by \er{3.68}. From \er{1.2}, \er{1.10n} we infer\break that
$m=4q+\F_4(m)$ \fs $q\in\N$. Hence $\F_{10}(\a^m) = \F_{10}(\a^{4q+\F_4(m)})\break
\nda22.25 = \F_{10}(\a^{4q}\cdot\a^{\F_4(m)}) \nde3.26 =
\F_{10}(\F_{10}(\a^{4q}) \cdot \F_{10}(\a^{\F_4(m)}))$. But $\F_{10}(\a^{4q})
\nde2.26 = \F_{10}((\a^4)^q) \nde3.69 = \F_{10}((\F_{10}(\a^4))^q) =
\F_{10}((1)^q) \nde2.27 = \F_{10}(1)\nde1.6 = 1$. Hence $\F_{10}(\a^m) =
\F_{10}(1\cdot\F_{10}(\a^{\F_4(m)})) \nda22.13 =
\F_{10}(\F_{10}(\a^{\F_4(m)})) = \F_{10}(\a^{\F_4(m)})$.
\exa

\bex3.67 \

\hph i,ii, Compute $\F_{10}(\a^{\F_4(m)})$ for $\a\in\{3,7,9\}$ and
$\F_4(m)\in [0,3)$.

\hph ii,i, Compute $\F_{10}(2^m)$, $m\ge2$, \Ip for $m=2n,3n$, $n\in\N$.

\hph iii,, Compute $\F_{10}(3^m)$, $m\ge2$, \Ip for $m=2n$, $n\in\N$.

In (ii), (iii) make use of \er{3.71}.
\eex

We now give a \gn\ of formula \er{3.26n}.

\bpr3.68
Let $(G,\qu,e)$ be a finite group of order $n>1$, and let $c_d\in\N$ denote
the \nm\ of \pn\ \sbm s of~$G$ of order $d\in\oz1,n $. Then
\beq3.72
n=1+\sum_{\sbk{1<d\le n\\d|n}} c_d\F(d),
\e
where $\F$ is the Euler-phi \f.
\epr

\brm3.69
If $G$ is a \pn\ group then $c_d=1$ \fa $d\in\oz1,n $ \st $d$~divides~$n$ by
\E\Pr s \rf{p3.24}\,(ii) and \rf{p3.20}. Hence in this case formula \er{3.72}
reduces to formula \er{3.26n} since $\F(1)=1$.
\erm

For the proof of \E\Pr\ \rf{p3.68} it is convenient to introduce the notion
of order of an \el\ of a finite group. In the case of a finite \ag\ $(X,\qu,e)$
(\ev tly finite \Cm) we defined in \E\df\ \rf{d2.11} the order of an \el~$a$
of~$X$ as the order of the monoid $(I(a),\qu,e)$. In the case of a finite
group $(G,\qu,e)$, if $a\in G$ then $I(a)$ is a \sbm\ of~$G$ by Lemma
\rfa2{l1.21} (observe that $I(e)=\{e\}$ by \era2{2.3}\,I4). \Mo $I(a)$ is
finite as a subset of a finite set by Theorem \rfa1{t4.18}. \E\Tf $\#(I(a))$
is well-defined.

\bdf3.70
Let $(G,\qu,e)$ be a finite group and let $a\in G$. The order of $I(a)$ is
called the \ti{order of~$a$} (in~$G$) and is denoted by~$|a|$ (not standard
notation). \E\Ip $|e|=1$.
\edf

\blm3.71
Let $(G,\qu,e)$ be a finite group and let $a\in G$. Then
\bea3.73
&|a| \hbox{ divides }|G|, \\
&n\dqu a\ne e,\q n\in(0,|a|), \lb{3.74}\\
&|a| \dqu a =e, \lb{3.75} \\
&(n+qm)\dqu a = n\dqu a \qh{\fa $n,q\in\N$} \lb{3.76}
\e
\elm

\proof
If $a=e$, then clearly $1\bigm| |G|$, \er{3.73} holds, and \er{3.75}, \er{3.76}
follow from \era2{2.3}\,I4. If $a\ne e$, then $I(a)$ is a nontrivial finite
abelian \sbm\ of~$G$. \csq, $|a|:=\#(I(a))$ divides~$|G|$ by Theorem \rf{t3.61},
and the other assertions hold by Lemma \rf{l2.12} since $(I(a),\qu,e)$ is a
finite \Cm\ (see the proof of \E\Pr\ \rf{p3.65}\,(i)).
\endproof

\proof[Proof of \E\Pr\ \rf{p3.68}]
Let $G$ be a nontrivial finite group, let
\bga3.77
n:=|G|,\\
A_d:=\{a\in G: |a|=d\}, \q d\in[1,|G|], \lb{3.78}
\e
Then in view of Lemma \rf{l3.71}, $G=\bcl_{\sbk{1\le d\le n\\d|n}}A_d$. Clearly
$A_d\cap A_{d'}=\vn$ whenever $d\ne d'$, $d,d'\in[1,n]$. \Mo $A_1=\{e\}$.
\If from \era2{1.126}, \era2{1.127} that
\beq3.79
n=1+\sum_{\sbk{1<d\le n\\d|n}}\#(A_d).
\e
Let $d\in\oz1,n $ be \st $d|n$. If $A_d\ne\vn$, let $a\in A_d$. Then, since
$I(a)$ is a \pn\ \sbm\ of order~$d$, the set of \pn\ \sbm s of~$G$ of order~$d$
is not empty.

Since $G$ is finite, its power set $\cP(G)$ is finite by Exercise \rfa2{ex3.27},
hence the set of \pn\ \sbm s of~$G$ of order~$d$, being a subset of $\cP(G)$,
is finite by Theorem \rfa1{t4.18}\,(i). \E\Tf \te\ a \nf~$I$ and an in\jc\ map
$i\mt H_{d,i}$ from $I$ into the set of \pn\ \sbm s of~$G$ of order~$d$. Set
\beq3.80
c_d:=\#(I).
\e
We claim that $\{{\rm Gen}(H_{d,i}\}_{i\in I}$ is a \pt\ of~$A_d$. Indeed,
\fe $i\in I$, ${\rm Gen}(H_{d,i})\ne\vn$ since $H_{d,i}$ is a \ti{\pn\/}
\sbm\ of~$G$. \Mo \fe $a\in A_d$, $a\in I(a)$ by \era2{2.3}\,I1. Since $I(a)$
is a \pn\ \sbm\ of~$G$, \te s $i\in I$ \st $I(a)=H_{d,i}$. \If that $a\in
{\rm Gen}(H_{d,i})$, hence $A_d\sbs \bcl_{i\in I}{\rm Gen}(H_{d,i})$. \E\oh
if $a\in {\rm Gen}(H_{d,i})$ \fs $i\in I$, then $I(a)=H_{d,i}$, hence
$\#(I(a)) =\#(H_{d,i})=d$. Thus $a\in A_d$. \E\Tf $\bcl_{i\in I}{\rm
Gen}(H_{d,i}) \sbs A_d$, and $A_d=\bcl_{i\in I}{\rm Gen}(H_{d,i})$. Finally,
${\rm Gen}(H_{d,i})\cap {\rm Gen}(H_{d,j})=\vn$ whenever $i,j\in I$, $i\ne
j$. Indeed, suppose for \cd ion that \te s $a\in{\rm Gen}(H_{d,i})\cap {\rm
Gen}(H_{d,j})$ \fs $i,j\in I$, $i\ne j$. By Lemma \rfa2{l1.21}\,(i) with
$(M,\qu,e):= (G,\qu,e)$, $M_0:=H_{d,i}$, we have $I(a)\sbs H_{d,i}$. Observe
that $\#(I(a))=d$ since $a\in H_{d,i}\sbs A_d$. By \df, $\#(H_{d,i})=d$. \E\Tf
$\#(I(a))=\#(H_{d,i})$, hence $I(a)=H_{d,i}$ by \era2{3.13}. Similarly
$I(a)=H_{d,j}$, hence $H_{d,i}=H_{d,j}$. A~\cd ion, since $i\ne j$, and the
map $i\mt H_{d,i}$ is in\jc. This completes the proof that $\{{\rm
Gen}(H_{d,i})\}_{i\in I}$ is a \pt\ of~$A_d$. Observe that
\beq3.81
\#({\rm Gen}(H_{d,i})) = \F(d), \q i\in I,
\e
in view of Theorem \rf{3.28} since $\#(H_{d,i})=d$ and $(H_{d,i},\qu,e)$ is
a \pn\ group by Lemma \rf{l3.60}. \If from \er{3.80}, \er{3.81} and
\era2{1.126} and \era2{1.127} that
\beq3.82
\#(A_d)=c_d\F(d).
\e

We now consider the case where $A_d$ is empty. Then $\#(A_d)=0$ by \era2{3.6}
and \er{3.82} holds with $c_d=0$ by \era2{2.13}. \If from
\er{3.79} that \er{3.72} holds.
\endproof

\Wanp give a \ch ization of \pn\  groups in the class of finite groups.

\begin{cor}[\cite{Groups}, Theorem 2.17]\lb{c3.72}
A nontrivial finite group of order~$n$ is \pn\ iff there is at most one \pn\
\sbm\ of order~$d$ \fe divisor $d$ of~$n$.
\eco

\proof
\ti{Only if\/}: follows from \E\Pr s \rf{p3.24}\,(ii) and \rf{p3.20}.

\ti{If\/}: by \as\ $c_d$ equals zero or one in \er{3.72}. Set $d_d:=1-c_d$,
hence $1=c_d+d_d$ (by \era3{8.72}, \era2{1.6}). Then $\suml_{d|n,1< d\le n}
\F(d)=\suml_{d|n,d>1}(c_d+d_d)\F(d) \nda21.129 = \suml_{d|n,d>1}c_d\F(d)
+\suml_{d|n,d>1}d_d\F(d)$. Hence
\bmlg
n+0\nda21.7 = n\nde3.26n = \suml_{d|n,d\ge 1}
\F(d) = 1+\suml_{d|n,1<d\le n}\F(d) = 1+\Bigl(\suml_{d|n,1<d\le n}c_d\F(d) +
\suml_{d|n,1<d\le n}d_d\F(d)\Bigr)\\ \nda21.5 = \Bigl(1+\suml_{d|n,1<d\le n}c_d\F(d)\Bigr)
+\suml_{d|n,1<d\le n}d_d\F(d) = n+\suml_{d|n,1<d\le n}d_d\F(d).
\e
Thus by \era2{1.6},
\era2{1.8} we obtain $\suml_{1<d\le n}d_d\F(d)=0$. Hence $d_d\F(d)=0$ \fa
$d\in\oz1,n $ \st $d|n$ by \era2{7.15}, \era2{7.17}. \E\Ip $d_n\F(n)=0$. Since
$\F(n)\ne0$ by \er{3.18}, we infer that $d_n=0$ in view of \E\Pr\
\rfa2{p2.7}\,(ii). Since $d_n=0$, we have $c_n=1$. \E\Tf \te s a \pn\ \sbm~$H_n$
of~$G$ of order~$n$. Since $\#(H_n)=n=\#(G)$, we have $G=H_n$ by \era2{3.13}.
Thus $G$ is \pn.
\endproof

Our next goal is to give an example of a group possessing an
infinite \pn\ \sbm. The group $(S_\N,\circ,\id_N)$ is an infinite
group. A \pn\ \sbm\ of~$S_\N$ is
of the form $I(f)$, the set of \IT s of an \el\ $f\in
S_\N\sms{\id_\N}$. Indeed, by Lemma \rfa2{l1.21}, $(I(f),\circ,\id_\N)$ is an
\am. Suppose that $I(f)$ is finite. Let $m:=\#(I(f))\ge2$. Since $I(f)$ is a
\sbm\ of the group $(S_\N,\circ,\id_\N)$, the \pn\ monoid $(I(f),\circ,\id_\N)$
is a \Cm\ hence a \pn\ group by \E\Pr\ \rf{p3.10}\,(ii). \E\Tf $f^m:= m\dci f$,
the $m$-th \IT\ of~$f$ in $(I(f),\circ,\id_\N)$ as well as in $(S_\N,\circ,
\id_\N)$ by Lemma \rfa2{l1.21}, \sf ies $f^m=\id_\N$ by \E\Pr\ \rf{p3.65}\,(i).
Conversely, if \te s $m'\in\N$ \st $f^{m'}=\id_\N$, then $I(f)$ is finite.
Indeed, note that the map $m\mt f^m$ from $(\N,+,0)$ to $(S_\N,\circ,\id_\N)$
is a \hm sm by \era2{2.3}\,I4, I5. If $f^{m'}=\id_\N$ \fs $m'\in\Na$, then
this \hm sm is \ti{not\/} in\jc, hence by Lemma \rfa2{l1.41}\,(iii) its range
$I(f)$ is finite. \csq, $I(f)$ is infinite iff $f^{m'}\ne\id_{S_\N}$ \fa
$m'\in\Na$.

An~\el\ $f\in S_\N$ is said to have
infinite order if there is no $m\in\Na$ \st $f^m=\id_\N$. We now construct
such a~\f~$f$. Let $\zb ck\N$ be the \sq\ of \el s of~$\N$ defined by
\beq3.83
c_0:=0, \q c_n:=\sum_{k=1}^n (k+1), \q n\in\Na.
\e
We have $c_n<c_m$ whenever $n<m$, $n,m\in\N$. Set
$A_n:=\zo c_n,c_{n+1} $, $n\in\N$. \If from Lemma \rfa2{l1.32} that $\zb
An\N$ is a~\pt\ of~$\N$. Since $\zo0,c_{n+1} $ is the disjoint union of
$\zo0,c_n $ and~$A_n$, we have $c_{n+1}\nda23.5 = \#(\zo0,c_{n+1} )\nda23.31
= \#(\zo0,c_n )+\#(A_n) = c_n+\#(A_n)$, hence $\#(A_n)=c_{n+1}-c_n$,
$n\in\N$. If $n=0$, then $\#(A_0)=c_1-c_0=c_1-0 ={1+1=2}$. If $n\ge1$, then
$\#(A_n)=\suml_{k=1}^{n+1}(k+1) - \suml_{k=1}^n(k+1) = \suml_{k=n+1}^{n+1}
(k+1)=(n+1)+1=n+2$. Hence
\beq3.84
\#(A_n)=n+2 \qh{\fa }n\in\N.
\e
We want to construct a \ti{bi\jc} map $f:\N\to\N$ \sf ying $f^n(c_n)=c_n+n$,
$n\in\N$.

Thus \Ip
\beq3.87
f^n\ne \id_\N \qh{\fa}n\in\Na.
\e
We require
\beq3.88
f(A_n)\sbs A_n \qh{\fa}n\in\N,
\e
and $f|_{A_n}$, the \rt ion of $f$ to $A_n$, to be a ``circular \Pm'' of the
form
\beq3.89
\bca
f|_{A_n}(c_n+k):= (c_n+k)+1 \ \hbox{ for }k\in\zo0,n+1 ,\\
f|_{A_n}(c_n+(n+1)):=c_n  & \hbox{\fa} n\in\N.
\eca
\e
For example, if $n:=1$, then $c_n=2$, $c_{n+1}=5$, $A_2=\zo2,5 $ and
$2 \nad f\to 3\nad f\to 4\nad f\to 2$. It turns out that $f^n(c_n)=(f|_{A_n})^n
(c_n)=c_n+n$ \fa $n\in\N$ and that $f$ is bi\jc. We now give a proof.

Since $\zb An\N$ is a \ti{\pt} of $\N$, $f$ is well-defined by \er{3.89}.
\Mo since $c_n$ and $c_n+(k+1)$ belong to $\zo c_n,c_n+(n+2) $ \fa $n\in\N$,
\er{3.88} holds. Observe that \fa $n\in\N$ $f|_{A_n}:A_n\to A_n$ is sur\jc,
hence bi\jc\ by Theorem \rfa1{t4.18} since $A_n$ is finite. Define $g:\N\to \N$
by setting
\beq3.90
g|_{A_n}(l):= (f|_{A_n})\Inv(l) \qh{\fa} l\in A_n,\ n\in\N.
\e
Then $g\circ f|_{A_n}=\id|_{A_n}$ and $f\circ g|_{A_n}=\id|_{A_n}$, $n\in\N$.
\E\Tf\ $g\circ f=\id_\N$ and $f\circ g=\id_\N$. Thus $f$~is bi\jc\ and $g=f
\Inv$.

Let $\vf_n:\N_{n+2}\to A_n$ and $\psi_n:A_n\to \N_{n+2}$, $n\in\N$, be defined
by
\beq3.91
\bca
\vf_n(k):=k+c_n, & k\in\N_{n+2},\\
\psi_n(l):= l-c_n, & l\in A_n,\ n\in\N.
\eca
\e
Then $\psi_n\circ\vf_n= \id_{\N_{n+2}}$, $\vf_n\circ \psi_n=\id_{A_n}$, hence
$\vf_n$ is bi\jc\ and $\psi_n=\vf_n\Inv$, $n\in\N$. Observe that
\beq3.92n
f|_{A_n}= \vf_n\circ \la_1\circ \vf_n\Inv, \q n\in\N,
\e
where $\la_1:\N_{n+2}\to \N_{n+2}$, defined by $\la_1(x):=x+_{n+2}1$, $x\in
(\N_{n+2},+_{n+2},0)$, belongs to $S_{\N_{n+2}}$ by \er{3.50} and \E\Pr\
\rf{p3.54}. By Exercise \rf{ex3.59}\,(iv) the map $\si\mt \vf_n\circ\si\circ
\vf_n\Inv$ from $S_{\N_{n+2}}$ into $S_{A_n}$
is an \is sm from $(S_{\N_{n+2}},\circ,\id_{\N_{n+2}})$
onto $(S_{A_n},\circ,\id_{A_n})$. \E\Tf by \era2{1.45}
\beq3.93n
(f|_{A_n})^l = \vf_n\circ\la_1^l \circ\vf_n\Inv, \q n,l\in\Na,
\e
where $\la_1^l$ is the $l$-th \IT\ of $\la_1$ in the group $(S_{\N_{n+2}},
\circ,\id_{\N_{n+2}})$.

\blm3.73
Let $m\ge2$ and $l\in\N$. Let $\la_1(x):=x+_m1$, $x\in\zo0,m $, and let
$\la_1^l$ denote the $l$-th \IT\ of~$\la_1$ in $(S_{\zo0,m },\circ,
\id_{\zo0,m })$. Then
\beq3.92
(\la_1^l)(x)=\F_m(x+l), \q x\in\zo0,m ,\ l\in\N.
\e
\elm

\proof
By \In\ on $l\in\N$.

$l:=0$: $(\la_1)^0(x)=\id_{\zo0,m }(x)=x\nde1.5 = \F_m(x)$ \fa $x\in\zo0,m $.

$l$ \ti{implies} $l+1$: $\la_1^{l+1}(x) \nda12.26a = \la_1^l \circ \la_1(x) = \la_1^{l}
(x+_m1) = \F_m((x+_m1)+l) \nde1.8 =\break \F_m(\F_m(x+1)+l) \nde1.8 = \F_m(x+1)
+_m l\nde1.7 = \F_m(\F_m(\F_m(x+1))+\F_m(l)) \nde1.5 = \F_m({\F_m(x+1)}+\F_m(l))
\nde1.7 = \F_m((x+1)+l)=\F_m(x+(1+l)) = \F_m(x+(l+1))$ \fa $x\in\zo0,m $.

Hence \er{3.92} holds.
\endproof

We now claim that $f^l|_{A_n}=(f|_{A_n})^l$ for all $l,n\in\N$. Let $n\in\N$.
We have $f^l(A_n)\sbs A_n$ \fa $l,n\in\N$ by \era1{2.31} and \er{3.88}. Let
$n\in\N$ and define $h:A_n\to A_n$ by setting $h:=f|_{A_n}$. Set $M:=\{k\in\N:
f^k|_{A_n}=g^k\}$ where $g^k$ is the $k$-th \IT\ of~$g$ viewed as a selfmap
of~$A_n$. We have $0\in M$, since $f^0|_{A_n}=\id|_{A_n}= \id_{A_n}=g^0$. We
next show that $k\in M$ implies $k+1\in M$. Let $k\in M$ and $x\in A_n$. Then
$f^{k+1}(x) = (f\circ f^k)(x)= f(f^k(x)) \nad{k\in M}= f(g^k(x))= f|_{A_n}
(g^k(x))= g(g^k(x))=g^{k+1}(x)$. Since $x$~is arbitrary in $A_n$ we obtain
$f^{k+1}|_{A_n}= g^{k+1}$, hence $k+1\in M$. \If that $M$~is \iv\ in $(\N,\le)$,
hence $M=\N$, which proves the claim. \E\Tf \fa $n\in\N$ we have $f^n(c_n)=
f^n|_{A_n}(c_n) = (f|_{A_n})^n(c_n) \nde3.93n = (\vf_n\circ \la_1^n\circ
\vf_n\Inv)(c_n)\nde3.91 = \vf_n(\la_1^n(0))\nde3.92 = \vf_n(\F_{n+2}(n))
\nde1.5 = \vf_n(n)\nde3.91 = c_n+n$.

We have proved

\blm3.74
The map $f:\N\to\N$ defined in \er{3.89} is bi\jc\ and \sf ies
\beq3.93
f^m\ne\id_\N \qh{\fa}m\in\Na
\e
where $m$ denotes the $m$-th \IT\ of~$f$ in $(S_\N,\circ,\id_\N)$.
\elm

The set $I(f)$ is not a subgroup of $(S_\N,\circ,\id_\N)$ but $I(f)\cup
I(f\Inv)$ is such a subgroup.

\bpr3.75
Let $f$ be as in \er{3.89}. Then $I(f)\cup I(f\Inv)$ is a countably infinite
abelian subgroup of $(S_\N,\circ,\id_\N)$.
\epr

\proof
We set
\beq3.94
\<f> := I(f)\cup I(f\Inv).
\e

\ti{$\<f>$ is a \sbm\ of $S_\N$}: $\id_\N\nad{\era2{2.3}\,\rm I0} = 0^0\circ f
\in I(f)\sbs\<f>$. Let $a,b\in\<f>$. If $a,b\in I(f)$ then $a\circ b\in I(f)$
since $I(f)$ is a \sbm\ of~$S_\N$ by Lemma \rfa2{l1.21}\,(ii). Hence $a\circ b
\in\<f>$. Similarly $a\circ b\in\<f>$ if $a,b\in I(f\Inv)$. If $a\in I(f)$
and $b\in I(f\Inv)$, then \te\ $m,n\in\N$ \st $a\nda22.4 = f^m$ and $b=
(f\Inv)^n$.

In case $m=n$, $a\circ b=f^m\circ(f\Inv)^m\nde3.7 = f^m\circ(f^m)\Inv= \id_\N
\in I(f)\sbs\<f>$.

In case $m>n$, \te s $p\in\N$ \st $m=p+n$. Hence $a\circ b=f^{p+n}\circ
(f\Inv)^n \nde3.7 = f^{p+n}\circ (f^n)\Inv \nda22.6 = (f^p\circ f^n)\circ
(f^n)\Inv \nda12.1 = f^p \circ(f^n\circ(f^n)\Inv) = f^p\circ \id_\N\nda12.2
= f^p \in I(f)\sbs\<f>$.

In case $m<n$, \te s $p\in\N$ \st $n=p+m$. Hence $a\circ b=f^m\circ (f\Inv)
^{p+m} \nde3.7 = f^m\circ (f^{p+m})\Inv \nda22.6 = f^m\circ((f^p\circ f^m))\Inv
\nde3.5 = f^m\circ((f^m)\Inv \circ(f^p)\Inv) \nde1.4 = (f^m\circ(f^m)\Inv)
\circ(f^p)\Inv = \id_\N\circ(f^p)\Inv \nda12.2 = (f^p)\Inv \nde3.7 = (f\Inv)^p
\in I(f\Inv)\sbs \<f>$.

This completes the proof of the case $a\in I(f)$, $b\in I(f\Inv)$. The proof
of the case $a\in I(f\Inv)$ and $b\in I(f)$ is similar and is left as an
exercise.

\ssk
\ti{$\<f>$ is a subgroup of $S_\N$}: If $a\in I(f)$, then $a=f^m$ \fs $m\in\N$.
Then $a\Inv=(f^m)\Inv = (f\Inv)^m \in I(f\Inv)\sbs\<f>$. Similarly, if $a\in
I(f\Inv)$, then $a=(f\Inv)^m$ \fs $m\in\N$, and $a\Inv=((f\Inv)^m )\Inv
\nde3.7 = ((f\Inv)\Inv)^m = f^m\in I(f)\sbs \<f>$.

\ti{$(\<f>,\circ,\id_\N)$ is an \ag}: Let $a,b\in\<f>$. If $a,b\in I(f)$
(resp.\ $I(f\Inv)$) then $a\circ b=b\circ a$ by Lemma \rfa2{l1.21}\,(ii).
If $a\in I(f)$ and $b\in I(f\Inv)$, then \te\ $m,n\in\N$ \st $a=f^m$ and
$b=(f\Inv)^n$. Set $g:=f\Inv$. We have $f\circ g=g\circ f$. Then $f^m\circ
g^n = g^n\circ f^m$ by \era2{1.23} in \E\Pr\ \rfa2{p1.13} with $(\wt E,\wt e,
\wt S):=(\N,0,S)$ and $(M,\qu,e):=(S_\N,\circ,\id_\N)$. The case $a\in
I(f\Inv)$ and $b\in I(f)$ is of course similar.

\ssk
\ti{$\<f>$ is countably infinite}: We first show

\ti{$I(f)$ is \ep\ to $\N$}. Let $\vf_f:\N\to S_\N$ denote the \sq\ of \IT s
of~$f$ in $(S_\N,\circ,\id_\N)$, that is, $\vf_f(m)=f^m$, $m\in\N$. Since
$I(f)$ is the range of~$\vf_f$ in~$S_\N$, $\vf_f:\N\to I(f)$ is sur\jc. We now
show that this map is also in\jc. Suppose, for \cd ion, that \te\ $m,n\in\N$,
$m\ne n$, \st $\vf_f(m)=\vf_f(n)$. \E\wlg we may suppose $m>n$, that is,
$m=p+n$ \fs $p\in\Na$. Then $f^p\circ f^n\nda22.6 = f^{p+n}=f^m=f^n\nda12.2 =
\id_\N\circ f^n$. By the \cnc ity of~$\circ$ in $(S_\N,\circ,\id_\N)$, we
obtain $f^p=\id_\N$, \cd ing \er{3.93} since $p\in\Na$. \E\Tf $\vf_f:\N\to
I(f)$ is bi\jc.

\ti{$I(f\Inv)$ is \ep\ to $\N$}. Since $(\<f>,\circ,\id_\N)$ is a group, the
map $\INV:\<f>\to \<f>$ defined by $\INV(a):=a\Inv$, $a\in\<f>$, \sf ies
$\INV \circ \INV=\id_{\<f>}$ in view of \er{3.4}. \E\Tf $\INV$ is bi\jc. \If
that the map $\INV\circ\vf_f :\N\to I(f\Inv)$ is bi\jc, since the inverse
of~$f$ as an \el\ of the group $\<f>$ is~$f\Inv$.

\ti{$\<f>$ is \ep\ to $\N$}. We first write $\N$ as the disjoint union of
even and odd \nm s, that is, $\N=A\cup B$, where $A:=\{n\in\N: \F_2(n)=0\}$,
$B:=\{n\in\N: \F_2(n)=1\}$, and $\F_2:\N\to\zo0,2 $ is defined in \er{1.3}.
In view of \er{1.2}, $A=\{n\in\N: n=2m$ \fs $m\in\N\}$ and $B=\{n\in\N: n=2m+1$
\fs $m\in\N\}$. We define $\psi:\N\to\<f>$ by setting
\beq3.95
\bca
\psi(2m):= \vf_f(m), &m\in\N\\
\psi(2m+1):=(\INV\circ\vf_f)(m+1), &m\in\N.
\eca
\e
In view of what precedes, $\psi$ is well defined.

\ti{Sur\ji\ of $\psi$}. If $a\in I(f)$, then there is $m\in\N$ \st
$a=f^m$. Hence $a=\vf_f(m)=\psi(2m)$. If $a\in I(f\Inv)\sms{\id_N}$,
then there is $m\in\N$ \st $a=(f\Inv)^{m+1}$. But $(f\Inv)^{m+1} \nde3.7 =
(f^{m+1})\Inv$, hence $a=(\INV\circ\vf_f)(m+1)=\psi(2m+1)$.

\ti{In\ji\ of $\psi$}. Let $i,j\in\N$, $i\ne j$. We show $\psi(i)\ne\psi(j)$.

\ti{Case $i,j\in A$}. We have $i=2m$, $j=2n$ \fs $m,n\in \N$. Since $i\ne j$,
$m\ne n$. Note that $\psi(i)=\psi(2m)=\vf_f(m)$,
$\psi(j) =\psi(2n)=\vf_f(n)$. Since $\vf_f$ is in\jc\ (see above), $\vf_f(m)
\ne\vf_f(n)$, hence $\psi(i)\ne\psi(j)$.

\ti{Case $i,j\in B$}. We have $i=2m+1$, $j=2n+1$ \fs $m,n\in\N$. Since $i\ne j$,
$2m+1\ne 2n+1$, hence we
obtain $m\ne n$. Note that $\psi(i)=\psi(2m+1)=(\INV\circ\vf_f)(m+1)$, $\psi(j)=
\psi(2n+1)=(\INV\circ\vf_f)(n+1)$. Hence $(\INV\circ\vf_f)(m+1)\ne(\INV\circ\vf_f)(n+1)$
(in view of the in\ji\ of $\INV\circ\vf_f$), which implies $\psi(i)\ne\psi(j)$.

\ti{Case $i\in A$, $j\in B$}. We have $i=2m$, $j=2n+1$ \fs $m,n\in\N$ (observe
that $i\ne j$). Suppose for \cd ion that $\psi(i)=\psi(j)$. Then $f^m=\vf_f
(m)=\psi(2m)=\psi(i)=\psi(j)=\psi(2n+1)=(\INV\circ\vf_f)(n+1)=(f^{n+1})\Inv$.  From $f^m=
(f^{n+1})\Inv$ we infer $f^m\circ f^{n+1}=\id_N$. However, $f^m\circ f^{n+1}
\nda21.42 = f^{m+(n+1)}\nda21.5 = f^{(m+n)+1}=\id_\N$. Since $(m+n)+1\in\Na$,
we obtain a \cd ion as a con\sq\ of \er{3.93}. \E\Tf $\psi(i)\ne\psi(j)$.

\ti{Case $i\in B$, $j\in A$}. Set $i:=j$, $j:=i$.

This completes the proof of the bi\ji\ of~$\psi$. Hence $\<f>$ is countably
infinite in view of \E\df\ \rfa1{d4.21}.
\endproof

\bdf3.76
A nontrivial group $G$ is said to be \ti{cyclic}\index{cyclic group} if there is $a\in G\sms e$
\st $G=I(a)\cup I(a\Inv)\nde3.94 = \<a>\nde3.4 = \<a\Inv>$.
\edf

\bpr3.77 \

\hph i,ii, A cyclic group is finite iff it is \pn.

\hph ii,i, All in\fcg s are \is c.

\hph iii,, The only auto\mf s of an in\fcg\ are the identity and the map
$\INV$ defined in \E\Pr\ {\rm\rf{p3.26}\,(i)}.
\epr

\bex3.78
Prove \E\Pr\ \rf{p3.77}.
\eex

\bex3.79
Let $n\ge2$ and let $A$ be a set containing at least two \el s. Let $G:=\{
f:A\to\N_n\}$, and set $(f+g)(i):=f(i)+_ng(i)$, $i\in A$, $f,g\in G$;
$0(i):=0$, $i\in A$. Show that $(G,+,0)$ is an \ag\ which is not cyclic. Show
that $G$ is un\ct e whenever $A:=\N$.
\eex

\bex3.80
Show that every nontrivial subgroup of a cyclic group is cyclic.
\eex

Recall that a finite \Cm\ is an \ag\ by \E\Pr\ \rf{p3.10}\,(i). We now show
that an \ti{infinite} \Cm\ can be embedded in an \ag. To this end we need the
\fw\ \df\ and lemma.

\bdf3.82
Let $(M,\qu,e)$ be a monoid and let $R$ be an \ev ce \rl\ on~$M$. Then the
binary \op~$\qu$ and the \rl~$R$ are called \ti{compatible} if the \fw\ holds:
\beq3.98
x\mathrel R x' \hbox{ and } y\mathrel R y' \qh{implies } (x\qu y)\mathrel R
(x'\qu y') \hbox{ \fa} x,y,x',y'\in M.
\e
\edf

\blm3.83
Let $(M,\qu,e)$ be a monoid and let $R$ be a compatible \ev ce \rl\ on~$M$.
Let $M/R$ denote the set of \ev ce classes $($see \E\df\ \rfa1{d4.6}$)$. Then

\hph i,ii, $[x]\hqu [y]:=[x\qu y]$, $x,y\in M$ is a well-defined binary \op\
on $M/R$.

\hph ii,i, $(M/R,\hqu,[e])$ is a monoid called the \emph{\qt\ monoid\/} $M/R$.
If $M$ is abelian, then so is $M/R$.

\hph iii,, The map $\pi:M\to M/R$ defined by $\pi(x):=[x]$, $x\in M$, is a
\emph{sur\jc\ \hm sm}.
\elm

\proof \

(i) It is \sft\ to show that if $x\mathrel R x'$, $y\mathrel R y'$, $x,y,x',y'
\in M$, then $[x\qu y]=[x'\qu y']$. By \er{3.98} we have $(x\qu y)\mathrel R
(x'\qu y')$, hence $x'\qu y'\in [x\qu y]$ by \era1{4.4} and $[x\qu y] =\break
[x'\qu y']$ by Exercise \rfa1{ex4.4}\,(i).

(ii) Let $x,y,z\in M$.

\ti{\E\asc ity}: $([x]\hqu [y])\hqu[z] = ([x\qu y]\hqu[z]) = [(x\qu y)\qu z]
= [x\qu (y\qu z)]\break = [x]\hqu [y\qu z] = [x]\hqu ([y]\hqu [z])$.

\ti{\E\cmt ity}: If $M$ is abelian, then $[x]\hqu[y] = [x\qu y] = [y\qu x] = [y]\hqu[x]$.

\ti{\E\nel\/}: $[x]\hqu[e] = [x\qu e] = [x] = [e\qu x] = [e]\hqu[x]$.

(iii) $\pi(e)=[e]$. Let $x,y\in M$, then $\pi(x\qu y)=[x\qu y] = [x]\hqu[y]
= \pi(x)\hqu \pi(y)$. Let $A\in M/R$, then $A\ne\vn$ by \E\df\ \rfa1{d4.5}\,(iii).
Let $x\in A$. Then $[x]=A$ by \era1{4.4} since $x\mathrel Rx$. Hence $\pi(x)
=A$.
\endproof

\Wanp prove the promised embedding theorem.

\bth3.84
Let $(X,\qu,e)$ be an infinite \Cm\ and let $(X\t X,\tqu,(e,e))$ denote
the $($direct\/$)$ product of~$X$ with itself, where
\beq3.99
(a,b)\tqu(c,d) := (a\qu c,b\qu d), \q a,b,c,d\in X.
\e

\hph i,ii, Let $\sim$ denote the  \rl\ on $X\t X$ defined by
\beq3.100
(a,b) \sim (c,d) \qh{if }a\qu d=b\qu c,\ a,b,c,d\in X.
\e
Then $\sim$ is a compatible \ev ce \rl\ on $(X\t X,\tqu,(e,e))$.

Let $[(a,b)]$ denote the \ev ce class containing $(a,b)$, $a,b\in X$.

\hph ii,i, Let $X\t X/{\sim}$ denote the set of $\sim$-\ev ce classes of $X\t X$.
Set
\bga3.101
[(a,b)] \hqu [(c,d)]:=[(a,b)\tqu(c,d)], \q a,b,c,d\in X,\\
\wh X:=X\t X/{\sim}\qh{and }\wh e:=[(e,e)]. \lb{3.102}
\e
Then the binary \op~$\hqu$ on $X\t X/{\sim}$ is well-defined and
$(\wh X,\hqu,\wh e)$ is an abelian group.

\hph iii,, Let $j:X\to \wh X$ be defined by $j(x):=[(x,e)]$, $x\in X$. Then $j$~is
an in\jc\ \hm sm.

\hph iv,, \E\fe $\wh x\in \wh X$, \te\ $a,b\in X$ \st $\wh x=j(a)\hqu
j(b)\Inv$, where $j(b)\Inv$ denotes the inverse of $j(b)$ in the
group~$\wh X$. \Mo $j(b)\Inv=[(e,b)]$.

\hph v,i, If $\wh Y$ is a subgroup of~$\wh X$ \st $j(X)\subset \wh Y$, then
$\wh Y=\wh X$.
\eth

\proof \

(i) \ti{$\sim$ is an \ev ce \rl\ on $X\t X$}:

\ti{Reflexivity}: Let $(a,b)\in X\t X$. Then $a\qu b \nda21.6 = b\qu a$, hence
$(a,b)\sim(a,b)$ by \er{3.100}.

\ti{\E\sy y}: Let $(a,b),(c,d)\in (X\t X)$ be \st $(a,b)\sim(c,d)$. Then
$a\qu d\nde3.100 = b\qu c \nda21.6 = c\qu b$, hence $c\qu b = a\qu d\nda21.6 = d\qu a$. Thus
$(c,d)\sim(a,b)$ by \er{3.100}.

\ti{\E\tr ity}: Let $(a,b),(c,d),(u,v)\in X\t X$ be \st $(a,b)\sim(c,d)$ and
$(c,d)\sim(u,v)$. Then $a\qu d\nde3.100 = b\qu c$ and $c\qu v=d\qu u$. By
\er{3.100} it is \sft\ to show $a\qu v=b\qu u$. We have $(a\qu v)\qu c
\nda21.5 = a\qu(v\qu c) \nda21.6 = a\qu(c\qu v) = a\qu(d\qu u) \nda21.5 =\break
(a\qu d)\qu u = (b\qu c)\qu u\nda21.6 = u\qu(b\qu c)\nda21.5 = (u\qu b)\qu c$.
Hence $(a\qu v)\qu c = (u\qu b)\qu c$. Then $a\qu v = u\qu b$ by \era2{1.8}.

\ti{$\qu$ and $\sim$ are compatible}: Suppose $(a_i,b_i)\in X\t X$, $i=1,2,3,4$.
We first prove the case where $(a_2,b_2)=(a_4,b_4)$, i.e.
\beq3.103
(a_1,b_1)\sim(a_3,b_3) \qh{implies }(a_1,b_1)\tqu(a_2,b_2)\sim(a_3,b_3)\tqu
(a_2,b_2).
\e
We have $a_1\qu b_3 = a_3\qu b_1$ by \er{3.100} and \era2{1.6}. \Mo $(a_1,b_1)\tqu(a_2,b_2)\nde3.99 = \break
(a_1\qu a_2,b_1\qu b_2)$ and $(a_3,b_3)\tqu(a_2,b_2)=(a_3\qu a_2,b_3\qu b_2)$.
Then $(a_1\qu a_2)\qu(b_3\qu b_2) \nda21.36 = a_1 \qu{(a_2\qu b_3)}\qu b_2
\nda21.6 = a_1\qu(b_3\qu a_2)\qu b_2 \nda21.36 = (a_1\qu b_3)\qu(a_2\qu b_2)
= (a_3\qu b_1)\qu\break(a_2\qu b_2) \nda21.36 = a_3\qu(b_1\qu a_2)\qu b_2
\nda21.6 = a_3\qu(a_2\qu b_1)\qu b_2 \nda21.36 = (a_3\qu a_2)\qu(b_1\qu b_2)$.
Hence $(a_1\qu a_2)\qu(b_3\qu b_2)= (a_3\qu a_2)\qu(b_1\qu b_2) \nda21.6 =
(b_1\qu b_2)\qu(a_3\qu a_2)$. \E\Tf by \er{3.100}, $(a_1,b_1)\tqu(a_2,b_2)
\sim (a_3,b_3)\tqu(a_2,b_2)$.

We now prove the general case, i.e.
\bml3.104
(a_1,b_1)\sim(a_3,b_3) \hbox{ and } (a_2,b_2)\sim(a_4,b_4)\\
\qh{implies } (a_1,b_1)\tqu(a_2,b_2) \sim(a_3,b_3)\tqu(a_4,b_4).
\e
By \er{3.103} we have $(a_1,b_1)\tqu(a_2,b_2) \sim (a_3,b_3)\tqu(a_2,b_2)$.
Since $\qu$~is \cmt e, so is $\tqu$ by Example \rfa2{xa1.5}\,(iv). Hence
$(a_3,b_3)\tqu(a_2,b_2) = (a_2,b_2)\tqu(a_3,b_3)\nde3.103 {\sim}
(a_4,b_4)\tqu (a_3,b_3) \nda21.6 = (a_3,b_3)\tqu (a_4,b_4)$.
From $(a_1,b_1)\tqu(a_2,b_2)\sim(a_3,b_3)\tqu(a_2,b_2)$ and
$(a_3,b_3)\tqu(a_2,b_2)\sim(a_3,b_3)\tqu(a_4,b_4)$, we infer $(a_1,b_1)\tqu
(a_2,b_2) \sim (a_3,b_3)\tqu(a_4,b_4)$ by the \tr ity of~$\sim$ proved above.

(ii) \If from Lemma \rf{l3.83} that $(\wh X,\hqu,\wh e)$ is an \ti{\am}. Let
$\wh x\in \wh X$. We show that \te s $\wh y\in\wh X$ \st $\wh x\hqu \wh y=
\wh e$. Let $(a,b)\in X\t X$ be such that $\wh x=[(a,b)]$. Set $\wh y:=
[(b,a)]$. Then $\wh x\hqu\wh y = [(a,b)]\hqu[(b,a)] \nad{\er{3.101},\er{3.99}}=\break
[(a\qu b,b\qu a)]\nda21.6 = [(a\qu b,a\qu b)]$. From \er{3.100} we have
$(a\qu b,a\qu b)\sim (e,e)$ hence $[(a\qu b,a\qu b)]=[(e,e)]\nde3.102 = \wh e$.
Hence $(\wh X,\hqu,\wh e)$ is a group by \era2{1.6} and \E\df\ \rf{d3.2}.

(iii) $j(e)=[(e,e)]=\wh e$. Let $x,y\in X$. Then $j(x\qu y) =[(x\qu y,e)]
\nda21.7 = [(x\qu y,e\qu e)] \break\nde3.99 = [(x,e)\tqu(y,e)] \nde3.101 = [(x,e)]
\hqu[(y,e)] = j(x)\hqu j(y)$. Let $x,y\in X$ be \st $j(x)=j(y)$. Then $[(x,e)]
=[(y,e)]$. Hence $(x,e)\in [(y,e)]$ by Exercise \rfa1{ex4.4}\,(i), and $(x,e)
\sim (y,e)$ by \era1{4.4}. Thus $x\nda21.7 = x\qu e \nde3.100 = e\qu y
\nda21.7 = y$.

(iv) Let $\wh x\in \wh X$ and let $(a,b)\in X\t X$ be \st $\wh x=[(a,b)]$
(recall that $\pi:X\t X\to \wh X$ is sur\jc). We have $(a,b)\nad{\era2{1.7},
\er{3.99}} = (a,e)\tqu(e,b)$, hence $\wh x=[(a,b)] = [(a,e)\tqu(e,b)] \nde3.101
= [(a,e)]\hqu[(e,b)]$. Since $[(e,b)]\hqu[(b,e)]\nde3.101 = [(e,b)\tqu(b,e)]
\nda21.7 = [(b,b)]\nde3.100 = [(e,e)]\nde3.102 = \wh e$, $[(e,b)]$ is the
inverse of $[(b,e)]$ in the group $(\wh X,\hqu,\wh e)$, which we denote by
$[(b,e)]\Inv$. \E\Tf $\wh x=[(a,e)]\hqu[(b,e)]\Inv = j(a)\hqu(j(b))\Inv$.

(v) Let $\wh Y$ be a subgroup of $\wh X$ \st $j(X)\sbs\wh Y$. Let $\wh x\in
\wh X$. Then by (iv) \te s $(a,b)\in X\t X$ \st $\wh x=j(a)\hqu(j(b))\Inv$,
where $j(b)\Inv$ is the inverse of~$j(b)$ in the group~$\wh X$. By \as\
$j(a),j(b)\in \wh Y$, hence $(j(b))\Inv\in\wh Y$ by Lemma \rf{l3.3} since
$\wh Y$ is a subgroup of~$\wh X$. \E\Tf $\wh x=j(a)\hqu (j(b))\Inv\in \wh Y$
since $\wh Y$ is a subgroup of~$\wh X$. Since $\wh x$ is arbitrary in~$\wh X$,
it follows that $\wh X\sbs\wh Y$, which together with $\wh Y\sbs \wh X$
implies $\wh Y=\wh X$. This completes the proof of Theorem \rf{t3.84}.
\endproof

\bex3.85
Let $X,\wh X,\wh e,j$ be as in Theorem \rf{t3.84}.

\hph i,i, Show that $j(X)=\wh X$ iff $X$ is a group. Set
\beq3.105
j(X)\Inv:=\{\wh y\in\wh X: \hbox{\te s $x\in X$ \st} \wh y=(j(x))\Inv \}.
\e
Show that $j(X)\cap j(X)\Inv=\wh X$ iff $(X,\qu,e)$ is a group.

\hph ii,, Show that $j(X)\cap j(X)\Inv =\{\wh e\}$ iff $X$ is a P-monoid.
\eex

\ssk
We recall that an \ti{infinite} \pn\ \Cm\ $(M,\qu,e)$ with \Gn~$a$ is a \PM\
(see \E\df\ \rfa2{d1.4}) by Theorem \rfa2{t1.42} and that the map $\vf_a:
\N\to M$ defined by $\vf_a(n)
:=n\dqu a$ is a monoid-\is sm by Theorem \rfa2{t1.25} with $\wt E:=\N$. \Mo
by Lemma \rfa3{l8.13}, Remark \rfa3{r1.7} $\vf_a$~is an order-\is sm (see
\E\df\ \rfa1{d3.33}) where both $\N$ and~$M$ are equipped with their \nog s
(see \E\df\ \rfa1{d3.24} with $E:=\N$ and \E\df\ and Notation \rfa3{d3.12}).
\E\Ip $(M,\lequ)$ is totally ordered.

We now consider Theorem \rf{t3.84} in special cases, namely when the monoid~$X$
is a\break  \PM\ whose \nog\ is a total \og\ and in the case where $M$~is a \pn\
monoid.

\bco3.85
Let $X,\wh X$ and $j$ be as in Theorem \rf{t3.84}.

\hph i,i, Suppose in \ad\ that $X$ is a \emph{\PM} whose \nog\ is a \emph{total
\og} $($see \E\df\ \rfa1{d3.19}$)$. Then \fe $\wh x\in\wh X\sms{\wh e}$, \te s
\ooo $c\in X$, called the \emph{modulus} of~$\wh x$, \st either $\wh x=j(c)$
or $\wh x=(j(c))\Inv$. We set
\beq3.106
|\wh x|:=c, \q |\wh e|:=e,
\e
and\glossary{$\sgn(z)$}
\beq3.107
\sgn(\wh x):=\bca
+ &\hbox{if }\wh x=j(|\wh x|),\\
- &\hbox{if }\wh x=(j(|\wh x|))\Inv.
\eca
\e
Note that $\sgn(\wh e)$ is \emph{not} defined. \Mo setting
\bga3.108
\wh X_+:=\{\wh x \in \wh X: \sgn(\wh x)=+\}, \q
\wh X_-:=\{\wh x \in \wh X: \sgn(\wh x)=-\},\\
\wh X \hbox{ is the \emph{disjoint} union of }\wh X_+,\ \{\wh e\}
\hbox{ and }\wh X_-. \lb{3.109}
\e
The sets $\wh X_+\cup\{\wh e\}$ and $\wh X_-\cup\{\wh e\}$ are \sbm s of the
group $(\wh X,\hqu,\wh e)$. The \rt ion to $\wh X_+\cup\{\wh e\}$ $($resp.\
$\wh X_-\cup\{\wh e\})$ of the map $\INV:\wh X\to\wh X$ defined in \E\Pr\
\rf{p3.26}\,{\rm(i)} is a monoid-\is sm from $\wh X_+\cup\{\wh e\}$ $($resp.\ $\wh
X_-\cup\{\wh e\})$ onto $\wh X_-\cup\{\wh e\}$ $($resp.\ $\wh X_+\cup\{\wh
e\})$.

Let $\wh x\in\wh X_+$ with $\wh x=j(x)$, $x\in X\sms e$, and $\wh y\in\wh
X_-$ with $\wh y=(j(y))\Inv$, $y\in X\sms e$, then
\beq3.110
\wh x\hqu \wh y = \bca
j(p) & \hbox{if }x=y\qu p,\ p\in X\sms e,\\
(j(p))\Inv & \hbox{if }y=x\qu p,\ p\in X\sms e.
\eca
\e

Finally, set for $\wh x,\wh y\in\wh X${\rm:}
\beq3.111
\wh x\le \wh y \qh{if } \wh y\hqu\wh x{}\Inv \in \wh X_+\cup\{\wh e\}.
\e
Then $\le$ is a \emph{total \og\ on~$\wh X$}.

\Mo
\bga3.112
j:(X,\nad{\qum}\le )\to (\wh X,\le) \hbox{ is strictly in\cre},\\
\hbox{and if $\wh x,\wh y,\wh z\in\wh X$, then}\non\\
\wh x\le \wh y \qh{implies } \wh x\hqu\wh z \le \wh y{}\hqu \wh z.\lb{3.113}
\e
Thus $\wh X$ is an \emph{infinite totally ordered \ag}.

\hph ii,, Suppose $X$ is a \emph{\pn\ monoid} with \Gn~$a$, i.e.\ $X=I(a)$.
\E\Ip $X$ is a totally ordered \PM. Then
\bea3.114
\wh X_+\cup\{\wh e\} &= I(j(a)),\\
\wh X_-\cup\{\wh e\} &= I((j(a))\Inv),\lb{3.115}
\e
where $I(j(a))$ denotes the set of \IT s of $j(a)$ in $(\wh X_+\cup\{\wh e\},\hqu,  \wh
e)$ and $I((j(a))\Inv)$ denotes the set of \IT s of $(j(a))\Inv$ in $(\wh X_-\cup\{\wh e\},
\hqu,\wh e)$. \E\Ip $(\wh X,\hqu,\wh e)$ is an \emph{\icg} and $\wh X= \langle
j(a)\rangle$ $($see \er{3.94}$)$.
\eco

\bex3.86
Let $X,\wh X,j$ be as in Corollary \rf{c3.85}\,(i). Prove
\bga3.116
\wh x \le\wh y \hbox{ implies }\wh y{}\Inv \le \wh x{}\Inv,
 \q \wh x,\wh y\in \wh X,\\
|\wh x{}\Inv|=|\wh x| \hbox{ and }
|\wh x\hqu\wh y| \le |\wh x|\qu |\wh y|, \q \wh x,\wh y\in \wh X,
\e
where $\le$ denotes the \nog\ of $X$.
\eex

\proof[Proof of Corollary \rf{c3.85}]\

(i) ``\ti{$(j(X),\hqu,\wh e)$ is an infinite \PM}'': Since $X$ is infinite
and $j$~is in\jc, $j(X)$ is infinite by \E\Pr\ \rfa1{p4.27}\,(A)(i). \Mo
$\wh e=j(e)\in j(X)$, and if $\wh u,\wh v\in j(X)$, then \te\ $x,y\in X$ \st
$\wh u=j(x)$, $\wh v=j(y)$, hence $\wh u\hqu\wh v=j(x)\hqu j(y) = {j(x\qu y)}
\in j(X)$ by Theorem \rf{t3.84}\,(iii). Thus $j(X)$ is a \sbm\ of $(\wh X,
\hqu,\wh e)$. Since $(\wh X,\hqu,\wh e)$ is a \Cm, $(j(X),\hqu,\wh e)$ is a
\Cm. Let $\wh u,\wh v\in j(X)$ be \st $\wh u\hqu\wh v=\wh e$. Then $j(x\qu y)
=j(x)\hqu j(y)=\wh u\hqu \wh v=\wh e=j(e)$. Since $j$~is in\jc, $x\qu y=e$ and
since $(X,\qu,e)$ is a \PM, we have $x=y=e$. Hence $\wh u=j(x)=j(e)=\wh e$,
and $\wh v=\wh e$.

``\ti{$\le$ defined in \er{3.111} is an \og\ on $X$}'': Note that by \er{3.111}
\fa $\wh x,\wh y\in\wh X$ $\wh x\le \wh y$ iff \te s $\wh p\in j(X)$ \st
$\wh y\hqu\wh x{}\Inv=\wh p$ iff \te s $\wh p\in j(X)$ \st $\wh y=\wh p\hqu
\wh x = \wh x\hqu\wh p$.

\ti{\E\tr ity} follows from $\wh y=\wh x\hqu\wh p$ and $\wh z=\wh y\hqu\wh q$,
$\wh x,\wh y\in\wh X$, $\wh p,\wh q\in j(X)$ implies $\wh z=(\wh x\hqu\wh p)\qu
\wh q=\wh x\hqu(\wh p\hqu\wh q)$ where $\wh p\hqu\wh q\in j(X)$.

\ti{Reflexivity} follows from $\wh x=\wh x\hqu \wh e$ where $\wh e\in j(X)$,
$\wh x\in\wh X$.

\ti{Anti\sy y} follows from $\wh y=\wh x\hqu\wh p$, $\wh x=\wh y\hqu \wh q$,
$\wh x,\wh y\in \wh X$, $\wh p,\wh q\in j(X)$, $\wh y\hqu\wh e = \wh y=(\wh y
\hqu\wh q)\hqu\wh p= \wh y\hqu (\wh q\qu \wh p)$, hence $\wh q\hqu\wh p=\wh e$
by \cnc ity. Then $\wh q=\wh p=\wh e$ since $j(X)$ is a \PM, and $\wh y=\wh x
\hqu \wh e=\wh x$.

``\er{3.112}, \er{3.113}'': Recall that for $x,y\in X$ $x\nad{\qum}\le y$ if
$y=x\qu p$ \fs $p\in X$ by \era3{8.7}. Since $j:X\to \wh X$ is an in\jc\ \hm
sm, we obtain $j(y)=j(x)\hqu j(p)$, $j(p)\in j(X)$. Then $j(x)\le j(y)$. If
$j(x)=j(y)$, then $x=y$. Thus $j$~is strictly in\cre. Let $\wh x,\wh y,\wh z
\in\wh X$ be \st $\wh y=\wh x\hqu\wh p$ \fs $\wh p\in j(X)$. Then $\wh y\hqu
\wh z=(\wh x\hqu\wh p)\hqu\wh z={\wh x\hqu(\wh p\qu\wh z)}= \wh x\hqu \wh(\wh p
\qu\wh z)=\wh x\hqu(\wh z\hqu\wh p)=(\wh x\hqu \wh z)\hqu \wh p$. Thus \er{3.113}
holds.

``\ti{$(\wh X,\le)$ is totally ordered\/}'': Let $\wh x,\wh y\in \wh X$ and
$a,b,c,d\in X$ be \st $\wh x:=[(a,b)]$, $\wh y:=[(c,d)]$. We claim that the
\fw\ holds:
\beq3.118
\aligned
{\rm (i)}\hskip30pt &\wh x<\wh y \qh{iff }a\qu d\nad{\qum}< b\qu c,\hskip100pt \\
{\rm (ii)}\hskip30pt &\wh y<\wh x \qh{iff }b\qu c\nad{\qum}< a\qu d.\hskip100pt
\endaligned
\e
Indeed, $\wh x<\wh y$ iff \te s $\wh p\in j(X)\sms{\wh e}$ \st $\wh y=\wh x
\hqu \wh p$. Since $j$~is in\jc\ and since $j(e)=\wh e$, \te s \ooo $p\in X
\sms e$ \st $\wh p=j(p)=[(p,e)]$ by Theorem \rf{t3.84}\,(iii). By \er{3.101},
\er{3.99}, $[(c,d)]=[(a,b)]\hqu [(p,e)]= {[(a\qu p,b\qu e)]}\nda21.7 = [(a\qu p,
b)]$. \If from Exercise \rfa1{ex4.4} that $(a\qu p,b)\sim (c,d)$. Hence $b\qu
c \nde3.100 = (a\qu p)\qu d = a \qu(p\qu d) = a\qu(d\qu p)= (a\qu d)\qu p$.
Thus by \df\ of $\nad{\qum}<$, $a\qu d \nad{\qum}< b\qu c$. Conversely, if
$a\qu d\nad{\qum}< b\qu c$, we have $b\qu c=(a\qu d)\qu p$ \fs $p\in X\sms e$,
hence $(a\qu p,b)\sim(c,d)$ by \er{3.100}. Thus $\bigl[((a\qu p),b)\bigr] =
[(c,d)]$. Since $\bigl[((a\qu p),b)\bigr]=[(a\qu p,b\qu e)] = [(a,b)\hqu
[(p,e)]$, we obtain $\wh y=\wh x\hqu j(p)$, where $j(p)\in\wh X\sms{\wh e}$.
Thus $\wh x<\wh y$, which completes the proof of \er{3.118}\,(i).
\er{3.118}\,(ii) follows by interchanging $\wh x$ and~$\wh y$. Since $(X,
\nad{\qum}\le)$ is totally ordered, we have either $a\qu d=b\qu c$, or
$a\qu d<b\qu c$, or $b\qu c<a\qu d$, \ev tly by \er{3.100} and \er{3.118},
either $\wh x=\wh y$, or $\wh x<\wh y$, or $\wh y<\wh x$. Hence $(\wh X,\le)$
is totally ordered.

We claim that the \fw\ holds.
\bea3.119
\{\wh x\in\wh X: \wh x>\wh e\} &= j(X)\sms{\wh e},\\
\{\wh x\in\wh X: \wh x<\wh e\} &= j(X)\Inv\sms{\wh e},\lb{3.120}
\e
where $j(X)\Inv$ is defined in \er{3.105}. Indeed, given $\wh x\in\wh X$,
$\wh x>\wh e$ iff $\wh x\ge \wh e$ and $\wh x\ne \wh e$ by \era3{8.9}. \Mo
$\wh x\ge\wh e$ iff $\wh x=\wh e\hqu \wh p=\wh p$ \fs $\wh p\in j(X)$. Thus
\er{3.119} holds. Note that $\wh x\le\wh e$ iff $\wh e=\wh x\hqu\wh p$ \fs
$\wh p\in j(X)$ iff $\wh x=\wh p{}\Inv$ \fs $\wh p\in j(X)$ iff $\wh x\in
j(X)\Inv$. Hence $\wh x<\wh e$ iff $\wh x\in j(X)\Inv\sms{\wh e}$. Thus
\er{3.120} holds. \csq, we have
\beq3.121
j(X) \cap j(X)\Inv=\{\wh e\}.
\e
Indeed, $j(X)\cap j(X)\Inv = ((j(X)\sms{\wh e}) \cup \{\wh e\}) \cap
((j(X)\Inv\sms{\wh e})\cup\{\wh e\}) \nad{\er{3.119},\er{3.120}}= \{{\wh x\in
\wh X}: \wh x\ge\wh e\}\cap\{\wh x\in\wh X: \wh x\le \wh e\}=\{\wh e\}$ since
the \og~$\le$ is total.

``\ti{Modulus of $\wh x$}'': Let $\wh x\in\wh X\sms{\wh e}$. Since $(\wh X,\le)$
is totally ordered, we have either $\wh x>\wh e$ thus $\wh x\in j(X)\sms{\wh e}$
by \er{3.119} or $\wh x<\wh e$ thus $\wh x\in j(X)\Inv\sms{\wh e}$ by
\er{3.120}. Since $j$~is in\jc, $j(e)=\wh e$ and $\wh e{}\Inv=\wh e$, we have
$j(X)\sms{\wh e}=j(X\sms e)$ and $j(X)\Inv\sms{\wh e}=j(X\sms e)\Inv$. \E\Tf
if $\wh x\in X\sms{\wh e}$, then either $\wh x=j(c)$ \fs $c\in X\sms e$ or
$\wh x=j(d)\Inv$ \fs $d\in X\sms e$. Note that such a~$c$ (resp.~$d$) is unique
if it exists in view of the in\ji\ of~$j$. \Mo if $c$ \ti{and\/}~$d$ exist
then $c\ne d$. Indeed, if $c=d$, then $j(c)=j(d)\Inv$, hence $j(c)\qu j(d)=
\wh e$. From \er{3.121} we infer $j(c)=j(d)=\wh e$, hence $c=d=e$, and
$\wh x=j(c)=j(e)=\wh e$, a~\cd ion since $\wh x\ne\wh e$. \E\Tf $|\wh x|$ the
modulus of~$\wh x$ is well-defined, as well as $\sgn(\wh x)$ by \er{3.107}.
\If from \er{3.119} (resp.\ \er{3.120}) that $\wh X_+$ (resp.\ $\wh X_-$)
defined in \er{3.108} \sf ies
\beq3.122n
\wh X_+=\{\wh x\in\wh X: \wh x>\wh e\},\q
\wh X_-=\{\wh x\in\wh X: \wh x<\wh e\}.
\e
Then \er{3.109} follows from ``$(\wh X,\le)$ is totally ordered''.

``\ti{$(j(X)\Inv,\hqu,\wh e)$ is an infinite \PM\ and $\INV|_{j(X)}: j(X)\to
j(X)\Inv$ is a monoid-\is sm}'': We have $\wh e=\wh e{}\Inv=(j(e))\Inv \in
j(X)\Inv$ by \er{3.105}. Let $\wh u,\wh v\in j(X)\Inv$ and let $x,y\in X$ be
\st $\wh u=(j(x))\Inv$, $\wh v=(j(y))\Inv$. Then $\wh u\hqu\wh v= (j(x))\Inv
\hqu(j(y))\Inv \nde3.5 = (j(y)\hqu j(x))\Inv = (j(x)\hqu j(y))\Inv = (j\qu(x
\qu y))\Inv\in j(X)\Inv$. Hence $j(X)\Inv$ is a \sbm\ of $(\wh X,\hqu,\wh e)$.
Since $(\wh X,\hqu,\wh e)$ is a \Cm, so is the monoid $(j(X)\Inv,\hqu,\wh e)$.
$\INV$, the selfmap of~$X$, is bi\jc, hence $\INV|_{j(X)}$ its \rt ion to
$j(X)$ is in\jc. By the \df\ of $j(X)\Inv$, $\INV|_{j(X)}$ maps $j(X)$ onto
$(j(X))\Inv$ hence is bi\jc. Since $\INV|_{j(X)}(\wh x)=\wh x{}\Inv$, $\wh
x\in \wh X$, and $\wh e{}\Inv=\wh e$, $(\wh x\hqu\wh y)\Inv= \wh x{}\Inv\hqu
\wh y{}\Inv$, $\wh x,\wh y\in\wh X$, $\INV|_{j(X)}$ is a monoid-\is sm. \csq,
$j(X)\Inv$ is an infinite \PM\ in view of \E\Pr\ \rfa1{p4.27}\,(A)(i) and
Exercise \rfa1{ex1.9}. We recall that $\wh X_+\cup\{\wh e\}$ (resp.\ $\wh X_-\cup
\{\wh e\}$) are equal to $j(X)$ (resp.\ $j(X)\Inv$) by \er{3.119}, \er{3.120}
and \er{3.122n}.

It remains to prove \er{3.110}. If $x=y\qu p$, then
$\wh x\hqu\wh y = j(x)\hqu (j(y))\Inv = (j(y)\hqu j(p))\hqu (j(y))\Inv =
(j(p)\hqu(j(y)))\hqu (j(y))\Inv =j(p)\hqu(j(y)\hqu(j(y))\Inv) = j(p)\hqu\wh e
= j(p)$. The proof of the case $y=x\qu p$ is similar, and is left to the
reader. This completes the proof of part~(i).

(ii) \er{3.114}: $\wh X_+\cup\{\wh e\} \nad{\er{3.122n},\er{3.119}} = j(X)$. By \as\ $X=\bcl_{n
\in\N} n\dqu a$. Since ${j:X\to\wh X}$ is in\jc, $j$~is a bi\jn\ from~$X$ onto
$j(X)$. Hence $j(X)=j\bigl(\bcl_{n\in\N}n\dqu a\bigr) = \bcl_{n\in\N} j(n\dqu
a)$. Since $j:X\to\wh X$ is a monoid-\hm sm from $(X,\qu,e)$ into $(\wh X,
\hqu,\wh e)$, it follows from \era2{1.45} that $j(n\dqu a)=n\dhqu a$. Hence
$j(X)=\bcl_{n\in\N}n\dhqu j(a) \nda21.26 = I(j(a))$.

\er{3.115}: $\wh X_-\cup\{\wh e\}\nad{\er{3.122n},\er{3.120}} = j(X)\Inv \nde3.105 = \INV(j(X))
=\INV\bigl(\bcl_{n\in\N}n\dhqu j(a)\bigr) \nad*= \break \bcl_{n\in\N}\INV(n\dhqu
j(a)) \nad{**}= \bcl_{n\in\N}n\dhqu(\INV(j(a))) = \bcl_{n\in\N} n\dhqu (j(a))
\Inv$. In $\nad*=$ we used the bi\ji\ of $\INV$ and in $\nad{**}=$ we used
the fact that $\INV$ is a monoid-auto\mf\ of $(\wh X,\hqu,\wh e)$. Thus
$\wh X_-\cup\{\wh e\} = I((j(a))\Inv)$ by \era2{1.26}.

Since $X\nde3.109 = X_+\cup\{e\}\cup X_- \nad{\er{3.122n},\er{3.119},\er{3.120},
\er{3.114},\er{3.115}}= I(j(a))\cup I((j(a))\Inv)\nde3.94 = \<j(a)>$, $(\wh X,
\hqu,\wh e)$ is a cyclic group in view of \E\df\ \rf{d3.76}. \Mo since $j:X\to
\wh X$ is in\jc\ and $X$~is infinite, $\wh X$~is infinite by \E\Pr\
\rfa1{p4.27}\,(A)(i).

This completes the proof of Corollary \rf{c3.85}.
\endproof

\bex3.87
Let $X,\wh X$ and $j$ be as in Corollary \rf{c3.85}.
Let $a,b\in X$ be \st $a\lequ b$ and let $a',b'\in X$ \sf y
$a\qu b'=b\qu a'$. Show: $a'\lequ b'$, $b'\mqu a' = b\mqu a$, $(a,b) \sim
(e,b\mqu a)$ and $[(a,b)] \cap (\{e\}\t X) = \{(e,b\mqu a)\}$.
\eex

The ``prototype'' of \icg s is the group of \ig s.

\bdf3.87
The set of \ti{\ig s}, usually denoted by~$\Z$, is the disjoint union of~$\{0\}$\glossary{$\Z$}
$(0\in\N)$ and of the product $\{+,-\}\t\Na$. The \el s of $\{+\}\t\Na$ are
usually denoted by $+n$ or simply by~$n$, $n\in\Na$, and are called \ti{\po\
\ig s}. The \el s of $\{-\}\t\Na$ are usually denoted by~$-n$, $n\in\Na$, and
are called \ti{\ng\ \ig s}. The set of \po\ (resp.\ \ng) \ig s is usually denoted\glossary{$\Z_{\ge0}$}
by $\Z_{>0}$ (resp.\ $\Z_{<0}$). The \el s of the set $\Z_{\ge0}:=\Z_{>0}\cup
\{0\}$ (resp.\ $\Z_{\le0}:=\Z_{<0}\cup\{0\}$) are called \ti{non\ng} \ig s
(resp.\ \ti{non\po} \ig s).
\edf

\bex3.88
Show that the map $h:\N\to\Z$ defined by
\beq3.122
h(0):=0, \q h(2n):=n, \q h(2n-1):=-n, \q\ n\in\Na,
\e
is a bi\jn. Thus $\Z$ is a countably infinite set (see \E\df\ \rfa1{d4.21}).
\eex

\bdf3.89
The binary \op\ on $\Z$, called \ti{\ad} and denoted by~$+$, is defined as
follows. Let $+_\N$ (resp.\ $-_\N$, see \era2{1.48}) be the \ad\ (resp.\ \sbt)
in~$\N$, then:
\bea3.123
z+0&:=z \hbox{ and } 0+z:=z &&\hbox{\fa}z\in\Z,\\
m+n&:=m+_\N n &&\hbox{\fa}m,n\in\Na, \lb{3.124}\\
(-m)+(-n)&:= -(m+_\N n) &&\hbox{\fa}m,n\in\Na, \lb{3.125}\\
m+(-n)&:= m-_\N n &&\hbox{\fa}m,n\in\Na \hbox{ \st} n\le m, \lb{3.126}\\
m+(-n)&:= -(n-_\N m) &&\hbox{\fa}m,n\in\Na \hbox{ \st} m<n, \lb{3.127}\\
(-m)+n&:= n-_\N m &&\hbox{\fa}m,n\in\Na \hbox{ \st} m\le n, \lb{3.128}\\
(-m)+n&:= -(m-_\N n) &&\hbox{\fa}m,n\in\Na \hbox{ \st}n<m. \lb{3.129}
\e
\edf

\bpr3.90
Let $\Z$, $+$ be as in \E\df s \rf{d3.87} and \rf{d3.89}. Then\/{\rm:}

\hph i,ii, $(\Z,+,0)$ is an \ag\ containing $\N$ as a \sbm.

\hph ii,i, The \ad\ in $\Z$ \rt ed to~$\N$ is identical to the \ad\ on~$\N$
introduced in \E\Pr\ \rfa2{p1.2} where $\ees:=(\N,0,S)$.

\hph iii,, \E\fe $n\in\Na$, $-n$ is the inverse of $n$, and $n$~is the inverse
of~$(-n)$.

\hph iv,, $(\Z,+,0)$ is an \icg\ generated by~$1$, i.e.\ $\Z=I(1)\cup I(-1)$
where $I(1)$ $($resp.\ $I(-1))$ denotes the set of \IT s of~$1$ $($resp.\ $(-1))$
in $(\Z,+,0)$.

\Mo if we denote by $-z$ the inverse of~$z$, \fe $z\in\Z$, then $w+(-z)$,
$w,z\in\Z$, is usually denoted by $w-z$, and is called the \emph{difference}
between $w$ and~$z$. If $w,z\in\N$ and $w\ge z$, then $w-z$ coincides with
$w-z$ defined in \era2{1.48}. The binary \op\ $(w,z)\mt w-z$ is called
\emph{\sbt}.

Finally, if we denote by $\le$ the \rl\ defined on $\Z$ by
\beq3.130
z\le w \qh{if } w-z\in \Z_{\ge0},\q z,w\in\Z,
\e
then $\le$ is a total \og\ on~$\Z$. Then \era2{1.63}--\era2{1.66} hold with
$x,y,z\in\Z$ where $<$~is defined in \era1{3.7}.
\epr

\brm3.91 \

\hph i,i, If $m,n\in\N$, $m<n$, then $m-n$ has no meaning in~$\N$.
However, in~$\Z$ $m-n$ is just $n-m$ with the \ng\ sign (see \er{3.127}).

\hph ii,, We prove \E\Pr\ \rf{p3.90} using Corollary \rf{c3.85}. One could
of course prove it directly. The verification of the \asc ity of the \ad\ is
tedious.

\erm

\proof[Proof of \E\Pr\ \rf{p3.90}]
We apply Corollary \rf{c3.85} with $(X,\qu,e):=(\N,+,0)$. Since $\N=I(1)$ by
\era2{2.18}, \era2{2.8}, we may apply Corollary \rf{c3.85}\,(ii).
We use the same notation as in Corollary \rf{c3.85}, namely $\wh X,\hqu, \wh
e, \wh X_+,\wh X_-,j$ and $a:=1$. The idea of the proof is to show that \te s
a bi\jc\ map $f:\wh X\to \Z$ \st $f(\wh e)=0$ and
\beq3.131
z+w = f(f\Inv(z)\hqu f\Inv(w)), \q z,w\in\Z.
\e
By Example \rfa2{xa1.5}\,(v), Remark \rfa2{r1.6}\,(iii) and \E\df\ \rfa2{d1.7}
we find that $(\Z,+,0)$ is an \am\ and that $f:\wh
X\to \Z$ is a monoid-\is sm. Then the remaining \pp ies of~$\Z$ easily follow.

\er{3.131}: We first define $g:\Z\to\wh X$ by setting
\beq3.132
g(0):=\wh e, \q g(n):=n\dhqu j(1), \q g(-n):= n\dhqu (j(1))\Inv,\q\ n\in\Na.
\e
Since $\wh e\nda22.3 = 0\dhqu j(1)$ and $\wh e\nda22.3 = 0\dhqu (j(1))\Inv$,
we have $g(\Z_{\ge0})=I(j(1))\nde3.114 = \break\wh X_+\cup\{\wh e\}$, and
$g(\Z_{\le0}) = I(j(1))\Inv)\nde3.115 = \wh X_-\cup\{\wh e\}$, hence $g(\Z)=
g(\Z_{\ge0}\cup \Z_{\le0}) = g(\Z_{\ge0})\cup g(\Z_{\le0})=\wh X_+\cup
\{\wh e\}\cup \wh X_-\nde3.109 = \wh X$. \If that $g$~is \ti{sur\jc}.
We now show that $g$~is \ti{in\jc}.
Since $j(\N,+,0)\to(\wh X,\hqu,\wh e)$ is a \hm sm, we have\break $n\dhqu j(1)\nda
21.45 = j(n\dpl1)$, $n\in\N$. Since $\INV:\wh X\to\wh X$ is an auto\mf, we
obtain $n\dhqu (j(1))\Inv \nda21.45 = (n\dhqu j(1))\Inv = (j(n\dpl 1))\Inv$.
But $n\dpl1 \nda22.8 = n\cdot1\nda22.18 = n$, $n\in\N$. \If that if $n\in\N$,
$g(n)=j(n)$ and if $n\in\Na$, $g(-n)=j(n)\Inv$. \E\Tf $g(\Z_{\ge0})\sbs \wh X
_+\cup\{\wh e\}$ and $g(\Z_{<0})\sbs \wh X_-$. Now let $z,w\in\Z$ be \st $g(z)=
g(w)$. If $g(z)\in \wh X_+$, then $g(w)\in \wh X_+$ and $z,w\in\Z_{\ge0}$, since
$\wh X_+\cup\{\wh e\}$ and $\wh X_-$ are disjoint. Thus $j(z)=g(z)=g(w)=j(w)$.
From the in\ji\ of~$j$ we infer $z=w$. Similarly if $g(z)\in\wh X_-$, then
$g(w)\in\wh X_-$, hence $z,w\in\Z_{<0}$. Let $m,n\in\Na$ be \st $z=-m$ and
$w=-n$. Then $g(-m)=j(m)\Inv$, $g(-n)=j(n)\Inv$, hence $j(m)\Inv = g(z) =
g(w) = j(n)\Inv$. \E\Tf $j(m)=g(z)\Inv =g(w)\Inv = j(n)$, hence $m=n$ by the
in\ji\ of~$j$. Thus $z=w$. Finally, if $g(z)=\wh e$, $z\in\Z$, then $z=0$.
Indeed, if $z=n\in\N$, then $g(z)=j(n)$. Since $j(0)=\wh e$ and $j$~is in\jc,
we have $z=n=0$. If $z=-n$, $n\in\Na$, then $g(z)=j(n)\Inv=\wh e$. Hence
$j(n) =((j(n))\Inv)\Inv =\wh e{}\Inv =\wh e$. A~\cd ion since $n\ne 0$ and
$j$~is in\jc. Thus $g(z)=\wh e$, $z\in\Z$, implies $z=0$. \E\Tf $g(z)=g(w)
=\wh e$, $z,w\in\Z$, implies $z=0=w$. This concludes the proof of the in\ji\
of~$g$, hence also of the \ti{bi\ji} of~$g$. We now define $f:\wh X\to \Z$ by
setting
\beq3.133
f:=g\Inv.
\e
Clearly $f$ is bi\jc\ and $f\Inv=g$. \Mo since $g(0)=\wh e$, $f(\wh e)
=f(g(0))=0$.

We now verify \er{3.131}. We set
\beq3.134
z\mathrel{\tilde+}w:=f(g(z)\hqu g(w)), \q z,w\in\Z.
\e
We verify $z\mathrel{\tilde+}w = z+w$ \fa $z,w\in\Z$.

$z\in\Z$, $w:=0$: $z\mathrel{\tilde+}0 \nde3.134 = f(g(z)\hqu g(0))\nde3.132
= f(g(z)\hqu\wh e) \nda21.7 = f(g(z))\nde3.133 =z\nde3.123 = z+0$. The case
$z:=0$, $w\in\Z$, follows from the \cmt ity of~$\hqu$.

$z:=m$, $w:=n$, $m,n\in\Na$: $z\mathrel{\tilde+}w = f(g(m)\hqu g(n))\nde3.132
= f((m\dhqu j(1))\hqu\break(n\dhqu j(1))\nda22.3 = f((m+_\N n)\hqu j(1))\nde 3.132
= f(g(m+_\N n))\nde3.133 = m+_\N n\nde3.124 = m+n=z+w$.

$z:=-m$, $w:=-n$, $m,n\in\Na$: $f(g(-m)\hqu g(-n))\nde3.132 = f(m\dhqu (j(1))
\Inv \hqu\break n\dhqu(j(1))\Inv) \nda22.3 = f((m+_\N n)\dhqu (j(1))\Inv) \nde3.132 =
f(g(-(m+_\N n))) = -(m+_\N n)\nde3.125 = (-m)+(-n) = z+w$.

$z:=m$, $w:=-n$, $m,n\in\Na$, $n\le m$: $f(g(m)\hqu g(-n))\nde3.132 =
f(m\dhqu j(1)\hqu n\dhqu (j(1))\Inv)$. Let $p\in\N$ be \st $m=p+n$. We have
$(p+n)\dhqu j(1)\hqu\break n\dhqu(j(1))\Inv) \nad{\era2{2.3},\era2{1.5}}= (p\dhqu
j(1)) \hqu((n\dhqu j(1))\hqu (j(1))\Inv))\nad{\era2{2.3}\,\rm I5}=
(p\dhqu j(1))\hqu(n\dhqu(j(1)\hqu (j(1))\Inv)) = (p\dhqu j(1))\hqu (n\dhqu
\wh e)\nad{\era2{2.3}\,\rm I4} = p\dhqu j(1))\hqu\wh e\nda21.7 = p\dhqu j(1)
\nda21.48 = (m-_\N n)\dhqu j(1)$. Thus $ m\mathrel{\tilde+}(-n) =
f((m-_\N n)\dhqu j(1))\nde3.132 = f(g(m-_\N n))= m-_\N n\nde3.126 = z+w$.

$z:=m$, $w:=-n$, $m,n\in\Na$, $m<n$: Let $p\in\Na$ be \st $n=m+p$. Then
$f(g(m) \hqu (g(-n)))\nad{\er{3.132},\era2{2.3}}= f(m\dhqu j(1)\hqu (m\dhqu
(j(1))\Inv \hqu p\dhqu(j(1))\Inv))\nad{\era2{1.5},\era2{2.3},\era2{1.7}}=
f(m\dhqu \wh e\hqu p\dhqu((j(1))\Inv))\nda2{2.3} = f(\wh e\hqu
p\dhqu((j(1))\Inv)) \nda21.7 = f(p\dhqu(j(1))\Inv)\nda21.48 =\break f((n-_\N m)\dhqu
(j(1))\Inv)\nde3.132 = f(g(-(n-_\N m))) = -(n-_\N m)\nde3.127 = n+(-m)=z+w$.
Observe that in view of the \cmt ity of~$\hqu$ we have
\beq3.135
z\mathrel{\tilde+}w = w\mathrel{\tilde+}z, \q z,w\in\Z.
\e

$z:=-m$, $w:=n$, $m,n\in\Na$, $m\le n$: $z\mathrel{\tilde+}w =
(-m)\mathrel{\tilde+} n\nde3.135 = n\mathrel{\tilde+}(-m)\nde3.126 = n-_\N m=
(-m)+n=z+w$.

$z:=-m$, $w:=n$, $m,n\in\Na$, $n<m$: $z\mathrel{\tilde+}w =
(-m)\mathrel{\tilde+}n \nde3.135 = n\mathrel{\tilde+}(-m)\nde3.127 = -(m-_\N
n) =(-m)+n=z+w$.

This completes the proof of \er{3.131}.

As mentioned above, it follows from Example \rfa2{xa1.5}\,(v), Remark
\rfa2{r1.6} and Lemma \rfa2{l1.8} that $(\Z,+,0)$ is
an \am\ and that $f:\wh X\to\Z$ is an iso\mf. Then $(\Z,+,0)$ is a group by
Lemma \rf{l3.25}, which proves (i),~(ii) since $0\in\N\sbs\Z$, and $x,y\in
\Z_{\ge0}$ implies $x+y\in\Z_{\ge0}$ by \er{3.123}, \er{3.124}. Then (iii)
follows from \er{3.126}, \er{3.128}.

(iv) $\Z=\Z_{\ge0}\cup \Z_{\le0}$ by \E\df\ \rf{d3.87}. \Mo $\Z_{\ge0}=\N
=I(1)$ and $\Z_{\le0}=\{-n\in\Z:n\in\N\} \nad{\rm(iii)} = \{n\Inv\in\Z:
n\in\N\} = \{(n\dpl1)\Inv: n\in\N\} \nda21.45 = \{n\dpl1\Inv: n\in\N\}
\nad{\rm(iii)}= I(1\Inv)$. Thus $\Z=I(1)\cup I(1\Inv)$ is a cyclic group with
\Gn~1 (or~$-1$) in view of \E\df\ \rf{d3.87}. $\Z$~is \ct y infinite by
Exercise \rf{ex3.88}.
\endproof

\bex3.92
Show that the \rl\ $\le$ on $\Z$ defined by \er{3.130} is a total \og\
on~$\Z$. Prove \era2{1.63}--\era2{1.66}, with $x,y,z\in\Z$.
\eex

\bex3.93 \

\hph i,ii, Show that the \sbt\ on~$\Z$ is neither \cmt e nor \asc e.

\hph ii,i, Let $a,b,c\in\Z$. Prove
\bga3.136
(a-b)-c = a-(b+c),\\
-(a+b)=(-a)+(-b), \q -(a-b)=b-a, \lb{3.137}\\
(a-b)+(b-c)=a-c. \lb{3.138}
\e

\hph iii,, Let $\sum$ denote the \cme \ad\ in $(\Z,+,0)$. Let $n\in\Na$ and
${a:[0,n]\to\Z}$. Prove
\bea3.139
\sum_{k=1}^n (a_k-a_{k-1}) &= a_n-a_0,\\
\sum_{k=0}^{n-1} (a_{k+1}-a_k) &= a_n-a_0. \lb{3.140}
\e
\eex

We now introduce the notion of \ti{\ng\ \IT s} of an \el\ of a~group. Recall
that if $(M,\qu,e)$ is a monoid and $a\in M$, then the map $\vf:\N\to M$ defined by
$\vf(n):=\break n \dqu a$ is a \hm sm from $(\N,+,0)$ into $(M,\qu,e)$, since
$\vf(0) \nda22.3 = e$ and $\vf(m+n)\break \nda22.3 = (m\dqu a)\qu (n\dqu a)
$, $m,n\in\N$. If $M$ is a group, we want to extend $\vf$ to~$\Z$ in such
a way that the extended \f~$\ov\vf$ be a \hm sm from $(\Z,+,0)$ into
$(M,\qu,e)$. We claim that a necessary \cn\ for~$\ov\vf$ to be a \hm sm is
$\ov\vf(-n)=n\dqu a\Inv$ \fe $n\in\Na$. Indeed, we have $a+\ov\vf(1\Inv)\nda
22.3 =\ov\vf(1) + \ov\vf(1\Inv)= \ov\vf(1 + 1\Inv)= \ov\vf(0)
\nda22.3 = e$. Hence
$\ov\vf(1\Inv) = a\Inv$. Then we apply Lemma \rfa2{l1.22} to prove the claim.

\bdf3.94
Let $(G,\qu,e)$ be a group (not necessarily abelian), and let $a\in G$. Given
$z\in\Z$, we set
\beq3.141
z\dqu a:=\bca n\dqu a &\hbox{if } z=n\in\N,\\
n\dqu a\Inv & \hbox{if }z=-n,\ n\in\Na. \eca
\e
\edf

The set of $z$-fold \IT s of~$a$ will be called the set of \ti{signed \IT s}\index{signed iterates}
of~$a$ in $(G,\qu,e)$. This terminology is not standard.

\setbox9=\hbox{SI0}
\bpr3.142
Let $(G,\qu,e)$ be an \ag. Given $z\in\Z$ and $x\in G$, let $z\dqu x$ $($or
simply~$z\cdot x)$ denote the $z$-fold \IT\ of~$x$ in $(G,\qu,e)$. Then \fa
$z,w\in\Z$ and $x,y\in M$ we have
\begin{gather}
\bca
\lca SI0 {0\cdot x = e,}\\
\lca SI1 {1\cdot x = x,}
\eca \non\\
\bca
\lca SI2 {(z+w)\cdot x = z\cdot x\qu w\cdot x,}\lb{3.142}
\eca \\
\bca
\lca SI4 {z\cdot e=e,}\\
\lca SI5 {z\cdot(x\qu y)=z\cdot x \qu z\cdot y,}\\
\lca SI6 {z\cdot x\Inv = (z\cdot x)\Inv,}\\
\lca SI7 {-z\cdot x=z\cdot x\Inv.}
\eca \non
\end{gather}
\epr

\bex3.143
Prove \E\Pr\ \rf{p3.142}.
\eex

\brm3.144
The \ext\ to \ng\ \IT s of \era2{2.3}\,I3 will be considered in Section 4.5.
\erm

\bex3.145
Let $(G,\qu,e)$ be a group and let $a\in G$. Prove:

\hph i,ii, $\{z\dqu a\in G: z\in \Z\} = I(a)\cup I(a\Inv)=:\<a>$.

\hph ii,i, The map $\vf:\Z\to\<a>$ defined by $\vf(z):=z\dqu a$, $z\in\Z$, is
a sur\jc\ \hm sm.

\hph iii,, If $\<a>$ is infinite, then $\vf$ is an iso\mf.

\eex

\brm3.101
Observe that $g(z)$ defined in \er{3.132} is the $z$-fold \IT\ of $j(1)$ in
the group $(\wh \N,\hpl,\wh 0)$ defined in Theorem \rf{t3.84}.
\erm

%% file: DETOUR44.TEX
\newpage
\Subsubsection{Finite prime fields}\label{sss.pr.fld}

The rings $(\N_p,+_p,\cdot_p,0,1)$, $p$ prime, possess a \pp y not shared by
the rings $\N_n$, $n$~composite \nm s. Indeed, every nonzero \el~$k$ of the
\mlv\ monoid $(\N_p,\cdot_p,1)$ \sf ies $\gcd(k,p)=1$, hence by \E\Pr\
\rf{p3.44} $k$~is a~unit of the ring $\N_p$ (see \E\df\ \rf{d3.34}). Since
a~unit of the ring $\N_p$ is an invertible \el\ of the monoid $(\N_p,\cdot_p,
1)$, it follows from Lemma \rf{l3.6} that the \rt ion to $(\N_p\sms0)\t
(\N_p\sms0)$ of the binary \op~$\cdot_p$ is a binary \op\
on $\N_p\sms0$ (still denoted by~$\cdot_p$)
with~$1$ as \nel. \Mo if $\N_p^\t:=\N_p\sms0$, then $(\N_p^\t,\cdot_p,1)$ is
an \ag. In case $n\in \Na\sms1$ is not a prime, that is, \te\ $a,b\in(1,n)$
\st $ab=n$, we have $a\cdot_n b\nde1.33 = \F_n(ab) \nad{\er{1.4}\er{1.6}} =
\F_n(n)=0$. Thus the \rt ion of~$\cdot_n$ to $(\N_n\sms0)\t(\N_n\sms0)$ is
not a binary \op\ on $\N_n\sms0$.

\begin{dfn}[\cite{Groups}]\lb{d4.1}
A (\cmt e) ring (with unity) $(F,+,\cdot,0,1)$ is called a \ti{field\/} if\index{field}
every nonzero \el\ is a unit. The set $F\sms0$ is usually denoted by~$F^\t$
and $(F^\t,\cdot,1)$, the group of units of~$F$, is usually called the
\ti{\mlv\ group} of~$F$.
\edf

\bpr4.2 \

\hph i,i, A \sr\ $($ring-$)$\is c to a field is a field.

\hph ii,, Let $(F,+,\cdot,0,1)$ be a field, let $F'$ be a set \ep\ to~$F$ and
let $f:F\to F'$ be a bi\jn. Set for $a',b'\in F'$
\bea4.1
a'+'b'&:= f\bigl(f\Inv(a')+f\Inv(b')\bigr),\\
a'\cdot' b'&:= f\bigl(f\Inv(a')\cdot f\Inv(b')\bigr),\lb {4.2} \\
0':=f(0),& \qquad 1':=f(1).\lb{4.3}
\e
Then $(F',+',\cdot',0',1')$ is a field and $f$ is a $($ring-$)$ \is sm.
\epr

\bex4.3
Prove \E\Pr\ \rf{p4.2}.
\eex

\brm4.4
Let $(X_i,\qu_0^{(i)},\qu_1^{(i)},e_0^{(i)},e\en i_1)$, $i\in[1,n]$,
$n\in\Na\sms1$, be a finite set of \sr s and let $X:=\prodl_{i\in[1,n]}X_i$ be
the product of \sr s defined in \E\Pr\ \rf{p2.20}. If $X_i$, $i\in[1,n]$, are
rings, then $X$~is a ring since the product of groups is a group by
Corollary \rf{c3.37}\,(i), and the group of units of~$X$ is the product of
the groups of units $U(X_i)$, $i\in[1,n]$. However, if $X_i$, $i\in[1,n]$,
are fields, then $X$~is \ti{not\/} a~field. Indeed, let $\ve_1\in X$ be
defined by $\ve_1(1):=e\en1_1$, $\ve_1(i):=e\en i_0$, $i\in[1,n]$, $\ve_2\in X$
be defined by $\ve_2(1):=e\en2_0$, $\ve_2(2):=e\en2_1$ and $\ve_2(i):=e\en i_0$,
$i\in[1,n]\sms{1,2}$. Then $\ve_1 \qu_1 \ve_2=e_0$. Since both $\ve_1$
and~$\ve_2$ belong to $X\sms{e_0}$, $\ve_1,\ve_2$ cannot be units of~$X$.
In the special but important case where all $X_i$'s, $i\in[1,n]$, are equal,
the ring~$X$ is also a \ti{\vs\ over the field $X_1$}. The notion of \vs\ plays an important
role in mathematics and in its applications. It will be introduced in Section \ref{sss.fdvs}.
\erm

We now give another \ch ization of a field. Let $(X,+,\cdot,0,1)$ be a~\ti{\sr},
let $a\in X\sms0$, $b\in X$, and consider the problem of finding $x\in X$ \sf
ying the \ti{\eq}
\beq4.4
ax+b=0.
\e
(Note that if $a=0$, then \er{4.4} reduces to $b=0$.) The symbol $x$ is
usually called the \ti{\uk} of the \eq. An \el\ $x\in X$ \sf ying \er{4.4}
is called a \ti{\so} to \eq~\er{4.4}. For example, in case $X:=\N$, in view of
\era2{1.9}, \eq~\er{4.4} has no \so\ if $b\ne0$. However, if $X$~is a field,
then \eq~\er{4.4} possesses \ooo \so
\beq4.5
\ov x:=a\mo\cdot(-b).
\e
In \er{4.5} $a\mo$ denotes the inverse of~$a$ in the group $(X^\t,\cdot,1)$,
and $-b$~denotes the inverse of~$b$ in the group $(X,+,0)$ (see \E\df\
\rf{d3.2}). Indeed, $a\ov x+b = a(a\mo\cdot(-b))+b = (aa\mo)(-b)+b =
1(-b)+b=(-b)+b=0$. \Mo if $\ov x,\ov y\in X$ are \so s to \er{4.4}, then
$a\ov x+b=0=a\ov y+b$, hence $a\ov x+b=a\ov y+b$. By \era2{1.8} we have
$a\ov x=a\ov y$, thus $\ov x=1\,\ov x=(a\mo a)\ov x=a\mo(a\ov x) = a\mo(a
\ov y) = (a\mo a)\ov y = 1\,\ov y=\ov y$. \E\Tf \eq\ \er{4.4} possesses \ooo
\so\ \fa $a\in X\sms0$ and $b\in X$.

Conversely, if $X$ is a \sr\ and \eq\ \er{4.4} has \ooo \so\ \fa $a\in X\sms0$,
$b\in X$, then $X$~is a field. Indeed, if $a:=1$ and $b\in X$, then \er{4.4}
reduces to $x+b=0$. Hence the \so~$x$ is the \ng\ of~$b$. Since $b$ is
arbitrary in~$X$, $(X,+,0)$ is a group, hence $X$~is a ring. If $a\in X\sms0$
and $b=-1$, then \er{4.4} becomes $ax+(-1)=0$, and the \so~$x$ \sf ies
$ax=ax+0 = ax+((-1)+1) = (ax+(-1))+1 = 0+1=1$. Hence $a$~is a unit of the
ring~$X$. Since $a\in X\sms0$ is arbitrary, $X$~is a field. Observe that only
\ti{the \ex\ of a \so} to \er{4.4} was used in the proof. We summarize in

\bpr4.5
Let $X$ be a \sr. Then the \fw\ are \ev t\/{\rm:}

\hph i,ii, $X$ is a field.

\hph ii,i, \E\eq\ \er{4.4} has \ooo \so\ \fa $a\in X\sms0$ and $b\in X$.

\hph iii,, \E\eq\ \er{4.4} has a \so\ \fa $a\in X\sms0$ and $b\in X$.
\epr

From now on we assume that $(F,+,\cdot,0,1)$ is a \ti{field\/}. We
investigate some \gn s of \eq~\er{4.4}. Let $a,b,c,d,f,g\in F$. We consider
the problem of finding $x,y\in F$ \st the \ti{system} of \eq s
\beq4.6
\bca
ax+by = f,\\
cx+dy = g
\eca
\e
is \sf ied. A pair $(x,y)\in F\t F$ \sf ying \er{4.6} is called a \ti{\so} to
system \er{4.6}. We are primarily interested in finding necessary and \sft\
\cn s on the \cf s $a,b,c,d$ of the \uk s $x,y$ guaranteeing \ex\ and \uq\ of
\so s to~\er{4.6} for \ti{all\/} \RHS s $f,g\in F$. Observe that in the case
of \eq\ $ax=f$, $a,f\in F$, such a \cn\ is $a\ne0$. We begin by looking for
a \ti{\sft\/} \cn\ for \uq\ of \so s. System \er{4.6} is called \ti{\hg} if
$f=g=0$ and \ti{in\hg} otherwise. If \er{4.6} is \hg, then the pair $(0,0)\in
F\t F$ is a \so\ in view of \er{1.1n} and \era2{1.7}. Such a \so\ is called
the \ti{trivial\/} \so. We are looking for a \cn\ on the \cf s implying that
$(0,0)$ is the only \so. Suppose $x,y\in F$ \sf y \er{4.6} with $f=g=0$.
\E\ml ying both sides of the first \eq\ by~$d$, and of the second by~$b$, we
obtain
\[
d(ax+b)=d0\nde1.1n = 0=b0=b(cx+dy).
\]
Hence by \er{1.2n}, $d(ax)+d(by)=b(cx)+b(dy)$. By the \asc ity and \cmt ity
of~$\cdot$, we have $(ad)x+(bd)y=(bc)x+(bd)y$. Hence by the \cnc ity of~$+$,
we arrive at $(ad)x=(bc)x$. If $x\ne0$, that is, if $x\in F^\t$, then
$xx\mo=1$ and $ad=(ad)1=(ad)(xx\mo)=((ad)x)x\mo = ((bc)x)x\mo =(bc)(xx\mo)=
bc1=bc$. Hence if $x\ne0$, then $ad=bc$. By contraposition, $ad\ne bc$ implies
$x=0$.

We now assume \cn
\beq4.7
ad\ne bc.
\e
Since $x=0$, we obtain $by =0$ from the first \eq\ in view of \er{1.1n} and
\era2{1.7}. Similarly $dy=0$ follows from the second one. By \er{4.7} and
\er{1.1n}, $b$~and~$d$ cannot be both equal to zero. If $b\ne0$, then $y=0$,
otherwise $by\in F^\t$. Similarly, if $d\ne0$, $y=0$. \If that $y=0$. \E\Tf
we proved that if \cn\ \er{4.7} is \sf ied, and $(x,y)$ is a \so\ of
\er{4.6} with $f=g=0$, then $x=y=0$.

It turns out that \cn\ \er{4.7} is also \ti{\sft\/} for the unique
solvability of system \er{4.6} for \ti{all\/} \RHS s $f,g\in F$. The \fw\
notation is standard.

\bnt4.6
Let $\a,\b,\g\in F$ and $\d\in F^\t$. Then
\bea4.8
{}&\a-\b := \a+(-\b),\\
&\frac\g\d :=\g\d\mo = \d\mo\g. \lb{4.9}
\e
\ent

\bpr4.7
Let $a,b,c,d\in F$ \sf y \er{4.7} and let $f,g\in F$. Set
\beq4.10
\ov x:=\frac{fd-gb}{ad-cb}\,, \qquad \ov y:=\frac{ag-cf}{ad-cb}\,.
\e
Then $(\ov x,\ov y)$ is the unique \so\ to system \er{4.6}.
\epr

In the proof we shall use the \fw\ identities:
\bea4.11
-(\a+\b)= (-\a)+(-\b), \q &\a,\b\in F,\\
-(\a\b)= (-\a)\b = \a(-\b),\q &\a,\b\in F. \lb{4.12}
\e

\bex4.8
Prove \er{4.11} and \er{4.12}.
\eex

\proof[Proof of \E\Pr\ \rf{p4.7}]
Set
\beq4.13
\D:=ad-bc.
\e

\ti{$(\ov x,\ov y)$ is a \so}: Since $\D\in F^\t$, it is \sft\ in view of
\er{1.3n} to show
\beag
ax\D + by\D &= f\D \\
\hbox{and }\ cx\D + dy\D &= g\D.
\e
We have $ax\D+by\D = (a(fd-gb))+(b(ag-cf)) = \nad{\er{1.2n},\er{4.12}} =
(afd+(-agb))+(bag+(-bcf)) \nad{\era2{1.6},\era2{1.36},\er{4.12}} = fad +
((-agb)+(agb)) + f(-(bc)) \nda21.7 = fad + f(-(bc))\nde1.2n = f\D$.
Similarly $cx\D+dy\D = (c(fd-gb))+(d(ag-cf)) = c(-gb)+dag = g(ad-bc)=g\D$.

\ti{\E\uq}: The idea is to show that if $(x',y'),(x'',y'')\in F\t F$ are \so
s to \er{4.6} with \RHS\ $f,g$, then $u:=x'-x''$ and $v:=y'-y''$ \sf y the
\hg\ system. Hence by the discussion preceding \E\Pr\ \rf{p4.7} we find
$u=v=0$. Hence $x''=0+x''=(x'+(-x''))+x'' = x'+((-x'')+x'') = x'+0 =x'$.
Similarly, one finds $y''=y'$.

It remains to show that $u,v$ \sf y the \ti{\hg} system \er{4.6}. We have
$ax'+by' = f =ax''+by''$. Hence $0=(ax'+by')-(ax'+by')=ax'+by'-(ax'' +by'')
\nde4.11 = ax'+by' +(-(ax''))+(-(by'')) = ax'+(-(ax''))+by'+(-(by''))
\nde4.12 = ax'+a(-x'')+by'+(-by'')
\nde1.2n = a(x'-x'') + b(y'-y'')=au+bv$. The proof of $cu+dv=0$ is similar.
\endproof

Next, we show that \cn\ \er{4.7} is necessary for \uq\ of \so s to system
\er{4.6} with $f=g=0$.

\bpr4.9
The \hg\ system \er{4.6} possesses a nontrivial \so\ if \cn\ \er{4.7} is
violated.
\epr

\proof
We assume
\beq4.14
ad=bc \qh{and } f=g=0.
\e

\ti{Case $b=c=0$}: System \er{4.6} reduces to two uncoupled \eq s $ax=0$ and
$dy=0$. Since $bc=0$, we must have $a=0$ or $d=0$ by \er{4.14}.
If $a=0$, then $x=1$ and $y=0$ is a \so, and if $b=0$ then $x=0$ and $y=1$ is
a~\so.

\ti{Case $b=0$, $c\ne0$}: Since $bc=0$, $ad=0$ by \er{4.14}, hence $a=0$ or
$d=0$. If $a=0$, then \er{4.6} reduces to one \eq\ $cx+dy=0$. Hence $x=
-c\mo dy$ by \er{4.5}. Thus $y=1$ and $x=-c\mo d$ is a \so. If $d=0$, $x=0$
and $y=1$ is a \so.

\ti{Case $b\ne0$, $c=0$}: As above $ad=0$, hence $a=0$ or $d=0$. If $a=0$,
$x=1$ and $y=0$ is a \so. If $d=0$, $x=1$ and $y=-b\mo a$ is a \so.

\ti{Case $b\ne0$, $c\ne0$}: Since $bc\ne0$, we have $ad\ne0$ by \er{4.14},
hence $a\ne0$ and $d\ne0$ by \er{1.1n}. From \er{4.14} we obtain $\frac dc=
\frac ba$. Set $\la:= \frac dc\in F^\t$. Then $d=\la c$ and $b=\la a$. Thus
\er{4.6} reads $a(x+\la y)=0$ and $c(x+\la y)=0$. Since $a,c\in F^\t$, we
obtain $x+\la y=0$. \E\Tf $y=1$ and $x=-\la$ is a \so.
\endproof

\bex4.10 \

\hph i,ii, Let $f,g\in F$ and assume $ad=bc$. Find necessary and \sft\ \cn s
on $f,g$ for \er{4.6} to have a \so. (\ti{Hint\/}: distinguish cases as in
the proof of \E\Pr\ \rf{p4.9}).

\hph ii,i, Show that if \er{4.6} possesses a \so\ \fa $f,g\in F$, then
\er{4.7} holds.

\hph iii,, Show that if $(\ov x,\ov y)$ is a \so\ to \er{4.6} and $(u,v)$
is a \so\ to the \hg\ system \er{4.6}, then $(\ov x+u,\ov y+v)$ is a \so\ to
\er{4.6}.
\eex

\brm4.11
\E\Pr s \rf{p4.9} and \rf{p4.7} show that \cn\ \er{4.7} is necessary and
\sft\ for system \er{4.6} to have at most one \so\
\fa $f,g\in F$. \E\Tf if the \hg\ system \er{4.6} possesses only the trivial
\so, then it possesses at least one \so\ \fa $f,g\in F$. Exercise
\rf{ex4.10}\,(ii) shows that if system \er{4.6} possesses at least a \so\ \fa
$f,g\in F$, then it possesses at most one \so, \Ip the \hg\ system \er{4.6}
possesses only the trivial \so.
\erm

We consider a more general situation. Let $I,J$ be \nf s, and let $m:=\#(I)$,
$n:=\#(J)$. Let $a:I\t J\to F$, $f:I\to F$ and $x:J\to F$. The \fw\
system of $m$~\eq s in $n$~\uk s
\beq4.15
\sum_{j\in J} a_{ij}x_j = f_i, \q i\in I,
\e
is called \ti{\hg} if $f_i=0$, $i\in I$, and \ti{in\hg} otherwise. Here
$\suml_{j\in J}$ denotes the \cme sum \era2{7.17} in the \am\ $(F,+,0)$. In
case \er{4.15} is \hg\ the \so\ $x_j=0$, $j\in J$, is called the \ti{trivial
\so}. Any other \so\ is called \ti{nontrivial}. \If from \E\Pr\ \rf{p4.9}
that if $m:=1$, $n:=2$, then the \hg\ system \er{4.15} possesses a nontrivial
\so. Indeed, if $c=d=0$ in \er{4.6}, then $ad=0=bc$, and the second \eq\ can
be omitted.

\blm4.12
The \hg\ system \er{4.15} possesses a nontrivial \so\ if $m<n$, $m,n\in\Na$.
\elm

The \fw\ identities will be used in the proof. Let $I$ denote a \nf, let
$\a:I\to F$ and let $\la\in F$. Then
\bea4.16
\la \sum_{i\in I} \a_i &= \sum_{i\in I}\la\a_i,\\
- \sum_{i\in I} \a_i &= \sum_{i\in I}(-\a_i).\lb{4.17}
\e
Observe that \er{4.16} is a con\sq\ of \era2{1.130} with $(X,\qu,e)=(\wt X,
\tqu,\wt e):=(F,+,0)$ and $\vf(x)=\la x$, $\la,x\in F$, since $\la(x+y)=\la x
+\la y$ by \er{1.2n}. \Mo \er{4.17} follows from \er{4.16} with $\la:=-1$.
Indeed, $(-1)\a\nde4.12 = -(1\a)=-\a$, $\a\in F$.

\proof[Proof of Lemma \rf{l4.12}] \

(i) We first consider the case $n:=m+1$, and proceed by \In\ on $m\in\Na$. Set
\bmlg
M:=\{k\in\Na: \hbox{every \hg\ system of the form \er{4.15}}\\
\hbox{with $\#(I)=k$, $\#(J)=k+1$ has a nontrivial \so}\}.
\e
We show that $M$ is \iv\ in $\Na$, hence $M=\Na$.

$1\in M$: follows from \E\Pr\ \rf{p4.9} with $c=d=0$.

\ti{$k\in M$ implies $k+1\in M$}: Suppose $k\in M$. Let $I,J$ be \nf s with
$\#(I)=k+1$ and $\#(J)=(k+1)+1$, let $a:I\t J\to F$. We are looking for
$\ov x:J\to F$ \st $\zb{\ov x{}}jJ$ is a nontrivial \so\ to the system
\beq4.18
\sum_{j\in J}a_{ij}x_j=0, \q i\in I.
\e
If all \cf s $a_{ij}$, $i\in I$, $j\in J$, are equal to zero, then $\ov x_j:=1$,
$j\in J$, is a nontrivial \so\ to \er{4.18}, by \er{1.1n} and $\suml_{j\in J}
0\nda21.126 = \#(J)\dpl 0\nad{\era2{2.3}\,\rm I4} = 0$.

We now suppose that \te\ $r\in I$, $l\in J$ \st $a_{rl}\ne0$. We claim that
$\zb xiJ$ is a \so\ to \er{4.18} iff the \fw\ system holds:
\bga4.19
x_l = a_{rl}\mo \Bigl(-\Bigl(\sum_{j\in J\sms l}a_{rj}x_j\Bigr)\Bigr),\\
\sum_{j\in J} a_{ij}x_j=0, \q i\in I\sms r. \lb{4.20}
\e

\ti{Only if\/}: If $\zb xiJ$ \sf ies \er{4.18}, then clearly \er{4.20} holds.
By \era2{1.138}, \eq
\beq4.21
\sum_{j\in J}a_{rj}x_j=0
\e
can be rewritten as
\beq4.22
a_{rl}x_l + \sum_{j\in J\sms l} a_{rj}x_j = 0.
\e
Since $a_{rl}\ne0$, \er{4.19} follows from \er{4.4}, \er{4.5}.

\ti{If\/}: Let $\zb xiJ$ \sf y \er{4.19}, \er{4.20}. \E\ml ying both sides of
\er{4.19} by $a_{rl}$, we find
\[
a_{rl}x_l = 1\Bigl(-\Bigl(\sum_{j\in J\sms l}a_{rj}x_j\Bigr)\Bigr) =
-\sum_{j\in J\sms l}a_{rj}x_j.
\]
Hence
\[
a_{rl}x_l + \sum_{j\in J\sms l}a_{rj}x_j = \Bigl(-\sum_{j\in J\sms l}
a_{rj}x_j\Bigr) + \sum_{j\in J\sms l}a_{rj}x_j=0.
\]
Thus $\zb xiJ$ \sf ies \er{4.22}, hence also \er{4.21}, which together with
\er{4.20} is \er{4.18}.

\medskip
The next step of the proof consists of showing that $\zb xiJ$ \sf ies \er{4.19},
\er{4.20} iff it \sf ies \er{4.19} and
\beq4.23
\sum_{j\in J\sms l} b_{ij}x_j =0, \q i\in I\sms r,
\e
where
\beq4.24
b_{ij}:= a_{ij} - a_{il}a_{rl}\mo a_{rj}, \q i\in I\sms r,\ j\in J\sms l.
\e
It is \sft\ to prove that if $\zb xiI$ \sf ies \er{4.19} then
\beq4.25
\sum_{j\in J}a_{ij}x_j  = \sum_{j\in J\sms l}b_{ij}x_j, \q i\in I\sms r.
\e
Let $\zb xiI$ \sf y \er{4.19}, and let $i\in I\sms r$. Then
\beq4.26
\sum_{j\in J} a_{ij}x_j \nda21.138 = \Bigl(\sum_{j\in J\sms l}a_{ij}x_j\Bigr)
+a_{il}x_l.
\e
\Mo $a_{il}x_l \nde4.19 = a_{il}a_{rl}\mo \bigl(-\bigl(\suml_{j\in J\sms l}
a_{rj}x_j\bigr)\bigr) \nde4.12 = -\bigl((a_{il}a_{rl}\mo) \suml_{j\in J
\sms l}a_{rj}x_j\bigr) \nde4.16 = \\
-\bigl(\suml_{j\in J\sms l}a_{il}a_{rl}\mo a_{rj}x_j\bigr) \nde4.17 =
\suml_{j\in J\sms l}(-(a_{il}a_{rl}\mo a_{rj}x_j)) \nde4.12 =
\suml_{j\in J\sms l}(-(a_{il}a_{rl}\mo a_{rj}))x_j$.

Hence \er{4.26} becomes $\suml_{j\in J}a_{ij}x_j = \suml_{j\in J\sms l}a_{ij}x_j
+ \suml_{j\in J\sms l}(-(a_{il}a_{rl}\mo a_{rj})x_j) \nda21.129 =\\
\suml_{j\in J\sms l}\bigl(a_{ij}x_j + (-(a_{il}a_{rl}\mo a_{rj}))x_j\bigr)
\nde1.3n = \suml_{j\in J\sms l}\bigl(a_{ij} + (-(a_{il}a_{rl}\mo a_{rj}))\bigr)
x_j \nad{\er{4.8},\er{4.24}} = \suml_{j\in J\sms l}b_{ij}x_j$. Hence
\er{4.25} holds.

Note that $k+1 = \#(I)\nda23.31 = \#(I\sms r)+\#(\{r\}) = \#(I\sms r)+1$.
Hence $\#(I\sms l)=k$ by \era2{1.8}. Similarly $\#(J\sms l)=k+1$. Since $k\in
M$, system \er{4.23} has  a nontrivial \so\ $\zb{\ov y}j{J\sms l}$. Let
$s\in J\sms l$ be \st $\ov y_s\ne0$. Set $\ov x_j:=\ov y_j$, $j\in J\sms l$,
and $\ov x_l:=a_{rl}\mo \bigl(-\suml_{j\in J\sms l}a_{rj}\ov x_j\bigr)$.
Then $\zb{\ov x}jJ$ \sf ies \er{4.19} and \er{4.23}. In view of \er{4.25},
it \sf ies \er{4.19}, \er{4.20}. Hence $\zb{\ov x}jJ$ is a \so\ to \er{4.18}.
Since $\ov x_s=\ov y_s\ne 0$, $\zb{\ov x}jJ$ is a nontrivial \so, and $k+1\in M$. This
completes the proof of part~(i).

\ssk
(ii) We now consider the case $n>m+1$. Set $p:=n-(m+1)$. We have $p\in\Na$,
and let $K$ be a finite set \sf ying $I\cap K=\vn$ and $\#(K)=p$. We define
$\wt a:(I\cup K)\t J\to F$ by setting
\beq4.27
\bca
\wt a_{ij}:=a_{ij}, &i\in I,\ j\in J,\\
\wt a_{ij}:=0,      &i\in K,\ j\in J.
\eca
\e
Since $\#(J)=n=(m+1)+p=(m+p)+1= \#(I\cup K)+1$, \te s $\ov x:J\to F$ \st
$\suml_{j\in J}\wt a_{ij}\ov x_j=0$, $i\in I\cup K$, in view of part~(i).
Hence $\suml_{j\in J}a_{ij}\ov x_j = \suml_{j\in J} \wt a_{ij}\ov x_j=0$,
$i\in I$. This completes the proof of part~(ii) and of the lemma.
\endproof

\brm4.13
System \er{4.15} is usually called a system of \ti{linear} \eq s for reasons\index{linear equations}
which will be clear later (see Section \ref{sss.fdvs}). As we shall see, Lemma \rf{l4.12} plays a \fd\
role in the theory of systems of linear \eq s. A~more general (non\cmt e)
version of Lemma \rf{l4.12} can be found for example in \cite[Theorem~1]{Artin}.
\erm

Observe that the proof given here is \ti{not\/} constructive, i.e.\ it does
not provide us with an algorithm for computing a nontrivial \so.

We now concentrate our attention to systems of the form \er{4.15} when the
number of \uk s is \ti{equal\/} to the \nm\ of \eq s. A~first important
con\sq\ of Lemma \rf{l4.12} is

\bpr4.13
Suppose that the \hg\ system \er{4.15} with $m=n$ has only the trivial \so. Then
the \crs\ in\hg\ system possesses a \so\ \fa \RHS s $\zb fiI$.
\epr

\proof
Set $\wt I:=I$, $\wt J:=J\cup\{l\}$ where $l\notin J$. Define $c:\wt I\t \wt
J\to F$ by setting
\beq4.28
\bca
c_{ij}:=a_{ij}, & i\in I,\ j\in J,\\
c_{il}:=f_i, &i\in I,
\eca
\e
where $f:I\to F$ is a given \RHS\ of \er{4.15}.

In view of Lemma \rf{l4.12} the \hg\ system $\suml_{j\in \td J}c_{ij}x_j=0$,
$i\in \wt I$, possesses a \ti{nontrivial\/} \so\ $\zb{\ov x}j{\td J}$.
Suppose for \cd ion that $\ov x_l=0$. Then $0=\suml_{j\in \td J} c_{ij}x_j
\nda21.138 = \suml_{j\in J}c_{ij}\ov x_j + c_{il}\ov x_l
\nad{\er{4.28},\er{1.1n}}= \suml_{j\in J}a_{ij}\ov x_j=0$, $i\in I$. By \as\
$\suml_{j\in J}a_{ij}\ov x_j=0$, $i\in I$, implies $\ov x_j=0$, $j\in J$. Thus
$\ov x_j=0$ \fa $j\in \wt J$, which is impossible since $\zb{\ov x}j{\td J}$
is a nontrivial \so. \E\Tf $\ov x_l\ne0$, and we obtain by \era2{1.138}
$c_{il}\ov x_l+\suml_{j\in J}c_{ij}\ov x_j=0$ \fa $i\in I$. \E\ml ying both sides
by $\ov x{}_l\mo$ and using \er{4.28}, \er{4.3}, \er{4.16} we arrive at $f_i+
\suml_{j\in J}a_{ij}\ov x{}_l\mo \ov x_j=0$, $i\in I$. Setting
$\ov y_j:=-\ov x_l\mo  \ov x_j$, $j\in J$, and using \er{4.8}, \er{4.12},
\er{4.17} we obtain $f_i-\suml_{j\in J}a_{ij}\ov y_j=0$, $i\in I$. \E\Tf
$\suml_{j\in J}a_{ij}\ov y_j =f_i$, $i\in I$.
\endproof

We now introduce another formulation of system \er{4.15}. Let $I,J,a,x$
and~$f$ be as in \er{4.15}. We introduce a map $L_a:F^J \to F^I$ by setting
\beq4.29
(L_ax)(i) := \sum_{j\in J} a_{ij}x_j, \q i\in I.
\e
Then \er{4.15} becomes
\beq4.30
L_ax=f
\e
where the \et y takes place in $F^I$.

Instead of having $m$ \eq s in $n$ \uk s in~$F$, we have one \eq\ in one \uk\
in~$F^I$.

\bnt4.14
Let $K$ be a \nf\ and let $(F,+,0)$ denote the additive group of the
field~$F$. We denote by $+_K$ (resp.~$0_K$) the binary \op\ in~$F^K$ (resp.\
the \el\ of~$F^K$) defined by
\bga4.31
(x+_K y)(k):= x(k)+y(k), \q k\in K,\\
0_K(k):=0, \q k\in K. \lb{4.32}
\e
If no confusion arises one writes $+$ (resp.~$0$) instead of  $+_K$
(resp.~$0_K$).
\ent

\blm4.15
Let $(F^K,+_K,0_K)$ be as in Notation \rf{n4.14}. Then

\hph i,i, $(F^K,+_K,0_K)$ is an \ag.

\hph ii,, Let $L_a:F^J\to F^I$ be as in \er{4.29}. Then $L_a$ is a group-\hm sm,
i.e.
\bea4.33
L_a0_J &= 0_I, \\
L_a(x+_Jy)&=L_ax+_I L_ay, \q x,y\in F^J, \lb{4.34}\\
L_a(-x)&=-L_a(x), \q x\in F^J, \lb{4.35}
\e
where $-x$ denotes the \ng\ of~$x$ in $(F^J,+_J,0_J)$ and $-L_a(x)$ denotes
the \ng\ of $L_a(x)$ in $(F^I,+_I,0_I)$.
\elm

\proof \

(i) The proof that $(F^K,+_K,0_K)$ is an \am\ is similar to the proof of
Example \rfa2{xa1.5}\,(iv). \E\Tf it is omitted. Then (i) follows from
Corollary \rf{c3.37}\,(i).

(ii) \er{4.33}: Let $i\in I$. Then $(L_a0_J)(i)\nad{\er{4.29},\er{4.32}}=
\suml_{j\in J}a_{ij}0 \nde1.1n = \suml_{j\in J} 0\nda21.126 = {\#(J)\dpls K
0_K}\break \nad{\era2{2.3}\,\rm I4}=0$. Thus \er{4.33} holds.

\er{4.34}: Let $i\in I$. Then $(L_a(x+_Jy))(i)\nde4.29 = \suml_{j\in J}
a_{ij}(x+_Jy)(j)\nde4.31 = \suml_{j\in J}a_{ij}(x(j) + y(j)) \nde1.2n =
\suml_{j\in J}(a_{ij}x(j)+a_{ij}y(j)) \nda21.129 = \suml_{j\in J}a_{ij}x(j)
+\suml_{j\in J}a_{ij}y(j) = (L_ax)(i)+(L_ay)(i) \nde4.31 = (L_ax+_I L_ay)(i)$.

\er{4.35}: Since $x+_J (-x)=0_J$, we have $L_a(x)+_I L_a(-x) \nde4.34  =
L_a(x+_J(-x)) = L_a(0_J)\nde4.33 = 0_I$. Hence $L_a(-x)=-L_a(x)$.
\endproof

\brm4.16
\E\pp y \er{4.34} of the map $L_a$ is usually called \ti{additivity}. It
implies that if~$x$ (resp.~$y$) is a \so\ to \er{4.28} \ev tly \er{4.30}
for the \RHS~$f$ (resp.~$g$), then $x+_Jy$ is a \so\ for the \RHS\ $f+_Jg$
(sometimes called \ti{principle of superposition}). Note that $L_a$ is
\ti{sur\jc} iff by \df\ system \er{4.28} possesses at least one \so\ \fa
\RHS s $f\in F^I$. \Mo $L_a$~is \ti{in\jc} iff $\ker L_a=0_I$ by \E\Pr\
\rf{p3.10}\,(vi), since by \er{4.33}, \er{4.34} $L_a$~is a \hm sm. Note that
$\ker L_a=0_I$ means that system \er{4.28} possesses only the trivial \so.
\erm

\bpr4.17
Let $I,J,a$ be as in \er{4.15} and let $L_a:F^J \to F^I$ be as in \er{4.29}.
Then

\hph i,i, If $L_a$ is in\jc, then $\#(J)\le \#(I)$.

\hph ii,, If $\#(J)=\#(I)$ and $L_a$ is in\jc, then $L_a$ is also sur\jc.
\epr

\proof \

(i) By contraposition. Suppose $\#(J)>\#(I)$, then by Lemma \rf{l4.12}, system
\er{4.15} with $f=0_I$ possesses a nontrivial \so. Then $L_a$ is \ti{not\/} in\jc\
by Remark \rf{r4.16}.

(ii) Since $L_a$ is in\jc, the \hg\ system \er{4.15} has only the trivial \so\
by Remark \rf{r4.16}.  Since $\#(J)=\#(I)$, $L_a$~is sur\jc\ by \E\Pr\ \rf{p4.13}.
\endproof

\brm4.18
If the field $F$ is \ti{finite} as are the fields $(\N_p,+_p,\cdot_p,0,1)$,
$p$~prime, then $F^J$ and~$F^I$ are finite by Theorem \rfa2{t3.26} and
$\#(F^J)=\#(F^I)$ are equal by \era2{3.38}. Hence, part~(ii) of \E\Pr\
\rf{p4.17} follows from \E\Pr\ \rfa2{p3.13}. Even more, if $L_a$~is sur\jc,
then it is also in\jc\ by the same \Pr. Since \te\ infinite fields we shall
give another proof of this fact even if $F$~is infinite.
\erm

\bex4.19
Show that if $\ov x\in F^J$ is a \so\ to \er{4.30}, then \ti{all\/} \so s to
\er{4.30} are of the form $\ov x+_J h$, where $h\in \ker L_a$.
\eex

\blm4.20
Let $I,J$ be \nf s, let $a:I\t J\to F$ and $b:J\t I\to F$. Let $c:I\t I\to F$
be defined by
\beq4.36
c_{ik}:=\sum_{j\in J} a_{ij}b_{jk}, \q i,k\in I.
\e
If $L_a: F^J\to F^I$, $L_b:F^I\to F_J$ and $L_c:F^I\to F^I$ are defined as in
\er{4.29}, then
\beq4.37
(L_a)\circ (L_b)=L_c,
\e
where $\circ$ denotes the \cm\ of two maps.
\elm

\proof
$(L_a)\circ(L_b)=L_c$ means $(L_a\circ L_b)(x)=L_cx$ \fa $x\in F^I$, which in
turn means
\[
((L_a\circ L_b)(x))(i) = (L_cx)(i) \qh{\fa}i\in I.
\]
Let $x\in F^I$ and let $i\in I$. Then $((L_a\circ L_b)(x))(i)=(L_a(L_bx))(i)
\nde4.29 = \suml_{j\in J}a_{ij}(L_bx)(j) \nde4.29 = \suml_{j\in J}a_{ij}
\bigl(\suml_{k\in I}b_{jk}x_k\bigr)\nde4.16 = \suml_{j\in J}\bigl(\suml_{k\in I}
a_{ij}b_{jk}x_k\bigr)$. We apply \E\Pr\ \rfa2{p4.18}\,(ii) with $\O:=J\t I$,
$\O=\bcl_{j\in J}B_j$, $B_j:=I\t \{j\}$. Then $\O$~is the disjoint union of the
subsets~$B_j$, hence by \era2{1.127} with $(X,\qu,e):=(F,+,0)$, we obtain
$\suml_{j\in J}\bigl(\suml_{k\in I}a_{ij}b_{jk}x_k\bigr) =
\suml_{(j,k)\in J\t I}a_{ij}b_{jk}x_k$. Since $J\t I$ is the disjoint union
of the subsets~$A_i$, defined by $A_i:=\{i\}\t J$, we find again by \era2{1.127},
$\suml_{(j,k)\in J\t I}a_{ij}b_{jk}x_k = \suml_{k\in I}\bigl(\suml_{j\in J}
a_{ij}b_{jk}x_k\bigr)$. But $\suml_{k\in I}\bigl(\suml_{j\in J}a_{ij}b_{jk}x_k
\bigr) = \suml_{k\in I}\bigl(\suml_{j\in J}x_k(a_{ij}b_{jk})\bigr)
\nde4.16 = \break \suml_{k\in I}x_k\bigl(\suml_{j\in J}a_{ij}b_{jk}\bigr) \nde4.36 =
\suml_{k\in I}x_kc_{ik} = \suml_{k\in I}c_{ik}x_k \nde4.29 = (L_cx)(i)$.
\endproof

\Wanp prove the result announced in Remark \rf{r4.18}.

\bth4.21
Let $I,J$ be \nf s with $\#(I)=\#(J)$. Let $F$ be a field, let $a:I\t J\to F$,
and let $L_a:F^J\to F^I$ be the map defined in \er{4.29}. Then $L_a$ is in\jc\
iff it is sur\jc.
\eth

\proof
\ti{Only if\/}: follows from \E\Pr\ \rf{p4.17}\,(ii).

\ti{If\/}: Let $\ve\en k\in F^I$ be defined by $\ve\en k(i) = e_{ki}$, \fa
$k,i\in I$, where $e_{ki}=1$ if $k=i$ and $e_{ki}=0$ if $k\ne i$. If $L_a$
is sur\jc, \te s $b\en k\in F^J$ \st
\beq4.38
L_ab\en k = \ve\en k \qh{\fe}k\in I.
\e
Set
\beq4.39
b_{jk}:=b\en k(j), \q k\in I,\ j\in J.
\e
Then $(L_ab\en k)(i)=\ve \en k(i)=e_{ki}=e_{ik}$, $i,k\in I$. We obtain
$\suml_{j\in J}a_{ij}b_{jk} = \suml_{j\in J}a_{ij}b\en k(j) =(L_ab\en k)(i)
= e_{ik}$, $i,k\in I$.

Observe that $e\in F^{I\t I}$. By Lemma \rf{l4.20}, $L_a\circ L_b=L_e$. But
$L_ex=x$ \fa $x\in F^I$. Indeed, $(L_ex)(i)=\suml_{j\in J}e_{ij}x_j =
\suml_{j\in J\sms i}0x_j + 1x_i= \suml_{j\in J\sms i}0 + x_i =
0+x_i = x_i$, $i\in I$. Hence $L_e=\id_{F^I}$. \E\Tf $L_a\circ L_b=\id_{F^I}$.
\If that $L_b$~is in\jc. Indeed, if $L_bx=L_by$, $x,y\in F^J$, then $x=\id_
{F^I}=(L_a\circ L_b)x = L_a(L_bx)= L_a(L_by)= (L_a\circ L_b)y =\id_{F^J}y=y$.
In view of \E\Pr\ \rf{p4.17}\,(ii), $L_b$~is also sur\jc, hence bi\jc. Thus
$L_a =L_a\circ \id_{F^I}= L_a\circ(L_b\circ L_b\Inv) = (L_a\circ L_b)\circ
L_b\Inv =\id_{F^I}\circ L_b\Inv=L_b\Inv$. Since $L_b\Inv$ is sur\jc, so
is~$L_a$.
\endproof

\bco4.22
If $\#(I)=\#(J)$, then system \er{4.5} possesses \ooo\so\ for \ti{all\/} \RHS
s $f\in F^I$ iff the trivial \so\ is the only \so\ to the \crs\ \hg\ system.
\eco

After having considered results valid in an arbitrary field~$F$, we now turn
to the study of a special class of fields called \ti{\Pf s}.

\bdf4.23
A subset $E$ of a field $(F,+,\cdot,0,1)$ is called a \ti{subfield\/}\index{subfield} of~$F$ if

\hph a,, $E$ is a subgroup (see \E\df\ \rf{d3.2}) of the \ag\ $(F,+,0)$,

\hph b,, $E$ is a \sbm\ (see \E\df\ \rfa1{d2.2}) of the \am\ $(F,\cdot,1)$,

\hph c,, Every \el\ of $E\sms0$ has an inverse in the monoid $(E,\cdot,1)$
(see \E\df\ \rf{d3.2}).
\edf

\blm4.24
If $E$ is a subfield of a field $(F,+,\cdot,0,1)$ and if $+$~$($resp.~$\cdot)$
denotes the \rt ion of~$+$ $($resp.~$\cdot)$ to $E\t E$, then $(E,+,\cdot,0,1)$
is a~field.
\elm

\proof
Since $E$ is a subgroup of $(F,+,0)$, we infer that $(E,+,0)$ is an \ag\ by
\E\df\ \rf{d3.2}.
Since $E$ is a \sbm\ of $(F,\cdot,1)$, we conclude that $(F,\cdot,1)$ is an
\am\ by \E\df\ \rfa1{d2.2}.
As a subset of~$F$, $E$~\sf ies \er{1.1n}, \er{1.2n} and \er{1.3n}. \E\Tf
$(E,+,\cdot,0,1)$ is a (\cmt e) ring (with $1$~as identity) (see \E\df\
\rf{d3.11}\,(i)). \Mo in view of~(c), every \el\ of~$E\sms0$ is a unit of the
ring~$E$. \E\Tf is a field by \E\df\ \rf{d4.1}.
\endproof

\bdf4.25
A subfield $E$ of a field $F$ is called \ti{proper} if $E\ne F$. A~field\index{subfield!proper}
which has \ti{no} proper subfield is called a \ti{\Pf\/}.\index{field!prime}
\edf

\bxa4.26
The finite fields $(\N_p,+_p,\cdot_p,0,1)$ with $p$~prime are \Pf s. Indeed,
if $E$~is a subfield of~$\N_p$, then $(E,+_p,0)$ is a subgroup of
$(\N_p,+_p,0)$ by \E\df\ \rf{d4.23}. Since $1\in E$, the group $(E,+_p,0)$ is
not trivial, hence by Theorem \rf{t3.61}, $\#(E)$ divides $p=\#(N_p)$. Since
$p$~is prime, $\#(E)\in\{1,p\}$. Since $(E,+_p,0)$ is not trivial, $\#(E)=p$.
Hence $E=\N_p$ by \era2{3.10}.
\exa

\blm4.27 \

\hph i,ii, The intersection of a nonempty family of \sbm s $($resp.\
subgroups, subfields\/$)$ of a monoid $($resp.\ group, field\/$)$~$X$ is
a~\sbm\ $($resp.\ subgroup, subfield\/$)$ of~$X$.

\hph ii,i, If\/ $Y$  is a \sbm\ $($resp.\ subgroup, subfield\/$)$ of a monoid
$($resp.\ group, field\/$)$~$X$, and $Z$~is a \sbm\ $($resp.\
subgroup, subfield\/$)$ of\/~$Y$, as monoid $($resp.\ group, field\/$)$, then
$Z$~is a \sbm\ $($resp.\ subgroup, subfield\/$)$ of the monoid $($resp.\
group, field\/$)$~$X$.

\hph iii,, Let $X$ be a monoid $($resp.\ group, field\/$)$, and let $Y,Z$
be \sbm s $($resp.\ subgroups, subfields\/$)$ of~$X$. If $Z\sbs Y$, then
$Z$~is a \sbm\ $($resp.\ subgroup, subfield\/$)$ of~$Y$ as monoid $($resp.\
group, field\/$)$.
\elm

\proof \

(i) ($\a$) Let $(X,\qu,e)$ be a \ti{monoid\/} and let $\cA$ be a nonempty
family of \ti{\sbm s} of~$X$. Set
\beq4.40
Y:=\bigcap_{M\in\cA} M.
\e
Then $Y$ is a subset of $X$. Since $e\in M$ \fa $M\in\cA$, $e\in Y$. Let
$a,b\in Y$, then $a,b\in M$ \fa $M\in\cA$. Since $M$~is a monoid \fa
$M\in\cA$, $a\qu b\in M$ \fa $M\in\cA$. Hence $a\qu b\in Y$. Since $a,b$ are
arbitrary in~$Y$, $Y$~is a \ti{\sbm} of~$X$.

($\b$) If $X$ is a \ti{group} and if $M$ is a \ti{subgroup} of~$X$ \fa
$M\in\cA$, then $Y$ defined in \er{4.40} is a \sbm\ of~$X$ since a group
(resp.\ subgroup) is a monoid (resp.\ \sbm). We now show that the monoid
$(Y,\qu,e)$ is a group, that is, \fe $a\in Y$ \te s $b\in Y$ \st $a\qu b=e=
b\qu a$. Let $a\in Y$. By \er{4.40} $a\in M$ \fa $M\in\cA$. Let $M\in\cA$.
Since $M$ is a subgroup of~$X$, $(M,\qu,e)$ is a group hence \te s $b\in M$
\st $a\qu b=e=b\qu a$. Since $Y\sbs X$, $a\in X$ and since $(X,\qu,e)$ is
a~group, $a$~possesses an inverse $a\Inv$ in $(X,\qu,e)$. In view of Lemma
\rf{l4.33}, $a\Inv=b$. Hence $a\Inv\in M$. Since $M$~is arbitrary in~$\cA$,
$a\Inv\in Y$, and $a\qu a\Inv=e=a\Inv\qu a$. Thus every \el\ of the monoid
$(Y,\qu,e)$ possesses an inverse, hence $(Y,\qu,e)$ is a group. Since $Y\sbs
X$ and $e\in Y$, $Y$~is a \ti{subgroup} of~$X$.

($\g$) Let $(X,\qu_0,\qu_1,e_0,e_1)$ be a \ti{field\/} and let $\cA$ be a
nonempty family of \ti{subfields} of~$X$. Let $Y$ be defined as in \er{4.40}.
We claim that $Y$ is a \ti{subfield\/} of~$X$.

\lhp$\g$,a, \E\fa $M\in\cA$, $M$~is a subfield of~$X$, hence $(M,\qu_0,e_0)$
is a subgroup of the group $(X,\qu_0,e_0)$. As an intersection of a nonempty
family of subgroups of $(X,\qu_0,e_0)$, $Y$~is a \ti{subgroup} of $(X,\qu_0,
e_0)$ by part~$(\b)$ of the proof.

\lhp$\g$,b, Similarly, $Y$ is a \ti{\sbm\/} of $(X,\qu_1,e_1)$.

\lhp$\g$,c, We first observe that \cn~(c) in \E\df\ \rf{d4.23} is \ev t to
the \fw\ \cn: $E\sms0$ is a \sbm\ of the monoid $(E,\cdot,1)$ and the monoid
$(E\sms0,\cdot,1)$ is a group. Indeed, suppose that \cn~(c) holds in \E\df\
\rf{d4.23}, then by Lemma \rf{l3.6}, $E\sms0$ is a subgroup of the monoid
$(E,\cdot,1)$, that is, $E\sms0$ is a \sbm\ of the monoid $(E,\cdot,1)$ and
$(E\sms0,\cdot,1)$ is a group. Conversely, if $E\sms0$ is a \sbm\ of the
monoid $(E,\cdot,1)$ and $(E\sms0,\cdot,1)$ is a group, then \fe $a\in E\sms0$,
\te s $b\in E\sms0$ \st $a\cdot b=1=b\cdot a$.

\leavevmode
\hphantom{$(\g)$} \csq, it is \sft\ to show that $Y\sms{e_0}$ is a \sbm\ of
the monoid $(X,\qu_1,e_1)$ and that the monoid $(Y\sms0,\qu_1,e_1)$ is a group.
We first claim that $Y\sms{e_0}$ is a \sbm\ of the monoid $(X,\qu_1,e_1)$. Let
$M\in\cA$. Then $M$ is a subfield of~$X$, hence in view of what precedes and
\E\df\ \rf{d4.23}\,(c), $M\sms{e_0}$ is a \sbm\ of $(M,\qu_1,e_1)$ and
$(M\sms{e_0},\qu_1,e_1)$ is a group. We claim that $M\sms{e_0}$ is a \sbm\ of
the group $(X\sms{e_0},\qu_1,e_1)$. Clearly $M\sms{e_0}\sbs X\sms{e_0}$ and
$e_1\in M\sms{e_0}$ since $e_1\in M$ and $e_1\ne e_0$. Let $a,b\in M\sms{e_0}$.
Then $a\qu_1 b\in M\sms{e_0}$ since $(M\sms{e_0},\qu_1,e_1)$ is a group, \Tf
since $a,b$ are arbitrary in $M\sms{e_0}$, $M\sms{e_0}$ is a \sbm\ of the group
$(X\sms{e_0},\qu_1,e_1)$. By part~$(\a)$ of the proof, $\Bca_{M\in \cA}
M\sms{e_0}$ is a \sbm\ of the monoid $(X\sms{e_0},\qu_1,e_1)$. Observe that
$\Bca_{M\in\cA}M\sms{e_0} = \bigl(\Bca_{M\in\cA} M)\sms{e_0} = Y\sms{e_0}$.
Hence $Y\sms{e_0}$ is a \sbm\ of $(X\sms{e_0},\qu_1,e_1)$, and the claim is
proved. We now show that $(Y\sms{e_0},\qu_1,e_1)$ is a group. Let $M\in\cA$.
Then $M\sms{e_0}$ is a \sbm\ of the group $(X\sms{e_0},\qu_1,e_1)$ and
$(M\sms{e_0},\qu_1,e_1)$ is a group. Hence $M\sms{e_0}$ is a subgroup of the
group $(X\sms{e_0},\qu_1,e_1)$. Since $Y\sms{e_0} = \Bca_{M\in\cA}M\sms{e_0}$
and $(X\sms{e_0},\qu_1,e_1)$ is a group, $Y\sms{e_0}$ is a subgroup of
$(X\sms{e_0},\qu_1,e_1)$ by part~($\b$) of the proof. \If that
$(X\sms{e_0},\qu_1,e_1)$ is a group, hence \cn~(c) of \E\df\ \rf{d4.23} is
\sf ied by~$Y$. \csq, $Y$~is a subfield of the field~$X$.

\ssk
(ii) ($\a$): Let $(X,\qu,e)$ be a monoid, $Y$~be a \sbm\ of the monoid
$(X,\qu,e)$ and $Z$ be a \sbm\ of the monoid $(Y,\qu,e)$. We have $Z\sbs Y
\sbs X$, and $e\in Z$ since $Z$~is a \sbm\ of $(Y,\qu,e)$. Let $a,b\in Z$.
Since $a,b\in Y$ and $Z$~is a \sbm\ of $(Y,\qu,e)$, $a\qu b\in Z$. Since
$a,b\in Z$ are arbitrary and $Z\sbs X$, $Z$~is a \sbm\ of the monoid $(X,\qu,e)$.

$(\b)$: \If from $(\a)$ that $Z$ is a \sbm\ of the \ti{group} $(X,\qu,e)$. It
suffices to show that $(Z,\qu,e)$ is a group. As a subgroup of the group~$X$,
$(Y,\qu,e)$ is a group, and as a subgroup of the group~$(Y,\qu,e)$, $(Z,\qu,e)$
is a group. Hence $Z$~is a \ti{subgroup} of~$X$.

$(\g)$: Let $(X,\qu_0,\qu_1,e_0,e_1)$ be a \ti{field\/}, and let $Y,Z$ be as
in~(iii).

\leavevmode\hphantom{$(\g)$}
(a) Since $(X,\qu_0,e_0)$ is a group, it follows from $(\b)$ that
$(Z,\qu_0,e_0)$ is a subgroup of $(X,\qu_0,e_0)$.

\leavevmode\hphantom{$(\g)$}
(b) Similarly $(Z,\qu_1,e_1)$ is a \sbm\ of the monoid $(X,\qu_1,e_1)$.

\leavevmode\hphantom{$(\g)$}
(c) By Lemma \rf{l4.24}, $(Y,\qu_0,\qu_1,e_0,e_1)$ is a field, since $Y$~is
a subfield of the field~$X$. Since $Z$ is a subfield of the field~$Y$,
$Z\sms{e_0}$ is a \sbm\ of $Y\sms{e_0}$ and $(Z\sms{e_0},\qu_1,e_1)$ is
a~group by \E\df\ \rf{d4.23}\,(c) and Lemma \rf{l3.6}. Clearly $e_1\in
Z\sms{e_0}$ since $e_1\in Z\sms{e_0}$ and $e_1\ne e_0$. \Mo $Z\sms{e_0}
\sbs X\sms{e_0}$. If $a,b\in Z\sms{e_0}$, then $a\qu_1 b \in Z\sms{e_0}$
since $(Z\sms{e_0},\qu_1,e_1)$ is a group. Hence $Z\sms{e_0}$ is a \sbm\ of
$(X\sms{e_0},\qu_1,e_1)$. Since $(Z\sms{e_0},\qu_1,e_1)$ is a group, every
\el\ of $Z\sms{e_0}$ has an inverse. \E\Tf \cn~(c) is \sf ied by~$Z$, hence
$Z$~is a \ti{subfield\/} of the field $X$ in view of parts (a), (b) of the
proof.

\ssk
(iii) ($\a$) As a \sbm\ of the monoid $(X,\qu,e)$, $Z$~\sf ies $e\in Z$ and
$a\qu b\in Z$. Since $Z\sbs Y$, $Z$~is a \sbm\ of the monoid $(Y,\qu,e)$.

($\b$) Since a group (resp.\ subgroup) is a monoid (resp.\ \sbm), $Z$~is
a~\sbm\ of the group~$Y$ by~($\a$). Since $Z$~is a subgroup of the group
$(X,\qu,e)$, \fa $a\in Z$ \te s $b\in Z$ \st $a\qu b=b\qu a=e$. Hence $Z$~is
a subgroup of the group $(Y,\qu,e)$.

($\g$) Let $(X,\qu_0,\qu_1,e_0,e_1)$ be a field. Then $Z$~is a subgroup of the
group $(Y,\qu_0,e_0)$ by~($\b$). Similarly $Z$~is a \sbm\ of the monoid
$(Y,\qu_1,e_1)$ by~($\a$). Since $Z$~is a subfield of the field~$X$, \fa
$a\in Z\sms{e_0}$ \te s $b\in Z\sms{e_0}$ \st $a\qu b=b\qu a=e$. Hence $Z$~is
a~subfield of the field $(Y,\qu_0,\qu_1,e_0,e_1)$.
\endproof

\bdf4.28
Let $F$ be a field. A subfield $E$ of the field~$F$ is called a \ti{\pf\/}\index{subfield!prime}
of~$F$ if, as a field, $E$~is a \Pf.
\edf

\bpr4.29 \

\hph i,ii, The \isc\ of all subfields of a field~$F$ is a \pf\ of~$F$.

\hph ii,i, A field possesses at most one \pf.

\hph iii,, Let $E$ be a subfield of a field~$F$. Then the \pf\ of the field~$E$
is equal to the \pf\ of~$F$.
\epr

\proof \

(i) Since the field $F$ is a subfield of~$F$, the family $\cA$ of all subfields
is not empty, hence $E:=\Bca_{M\in\cA}M$ is a subfield of~$F$ by Lemma
\rf{l4.27}\,(i). By Lemma \rf{l4.24}, $E$~can be viewed as a field. Let $K$~be
a subfield of the field~$E$, $K$~is also a subfield of~$F$ by Lemma
\rf{l4.27}\,(ii). \If that $K\in\cA$, hence $E\sbs K$, and $K=E$. Thus
$E$~has no proper subfield.

(ii) Suppose $E_1$ and $E_2$ are \pf s of a field~$F$. Then $E_1\cap E_2\sbs
E_1$. But $E_1\cap E_2$ is also a subfield of~$F$ by Lemma \rf{l4.27}\,(i).
Hence $E_1\sbs E_1\cap E_2$ by the \df\ of~$E_1$. Thus $E_1=E_1\cap E_2$.
Similarly $E_2=E_1\cap E_2$, hence $E_1=E_2$.

(iii) Let $\check E$ (resp.\ $\check F$) denote the \pf\ of the field~$E$
(resp.~$F$). Then $\check E$ is a subfield of the field~$E$, and $E$~is a
subfield of the field~$F$. Hence $\check E$ is a subfield of~$F$ by Lemma
\rf{l4.27}\,(ii). \E\Tf $\check F\sbs \check E$ by \df\ of~$\check F$. \E\oh
$\check F$ and~$E$ are both subfields of~$F$, and $\check F\sbs E$ since
$\check F\sbs \check E\sbs E$. By Lemma \rf{l4.27}\,(iii) we infer that
$\check F$ is a subfield of the field~$E$. Hence by \df\ of~$\check E$,
$\check E\sbs\check F$. \E\Tf $\check F=\check E$.
\endproof

\blm4.30
Let $F,F'$ be fields and let $f:F\to F'$ be a ring-\is sm.

\hph i,i, If $E$ is a subfield of $F$, then $f(E)$ is a subfield of~$F'$.

\hph ii,, If $\check F$ is the \pf\ of~$F$, then $f(\check F)$ is the \pf\
of~$F'$.
\elm

\bex4.31
Prove Lemma \rf{l4.30}.
\eex

We showed in Example \rf{xa4.26} that the finite fields $(\N_p,+_p,\cdot_p,0,
1)$, $p$~prime, are \Pf s hence also finite \pf s. Our next goal is to show
that the \ca\ of a \ti{finite} \Pf\ is a \Pn~$p$ and that this field is \is c to
$(\N_p,+_p,\cdot_p,0,1)$. To this end we shall use the \fw\ \df\ and lemma.

\bdf4.32
A subset $Y$ of a \sr\ $(X,\qu_0,\qu_1,\e_0,e_1)$ is called a \ti{sub\sr}\index{subsemiring}
of~$X$ if

\hph a,, $Y$ is a \sbm\ of the monoid $(X,\qu_0,e_0)$,

\hph b,, $Y$ is a \sbm\ of the monoid $(X,\qu_1,e_1)$.
\edf

One easily verifies that if $Y$ is a sub\sr\ of the ring~$X$, then
$(Y,\qu_0,\qu_1,e_0,e_1)$ is a \sr\ (where $\qu_i$ is the \rt ion to $Y\t Y$
of the binary \op\ $\qu_i$ in~$X$, $i=1,2$).

\blm4.33
Let $(X,\qu_0,\qu_1,e_0,e_1)$ be a \sr\ $($not necessarily \pn$)$, and let
${\vf:\N\to X}$ be defined as in \er{1.5n}. Let $I(e_1)$ denote the range of
the map~$\vf$ in~$X$, that is, the set of \IT s of~$e_1$ in the monoid
$(X,\qu_0,e_0)$. Then $I(e_1)$ is a \emph{sub\sr} of the \sr~$X$ and the map
$\vf:(\N,+,\cdot,0,1)\to (I(e_1),\qu_0,\qu_1,e_0,e_1)$ is a \ti{sur\jc} ring-\hm
sm.
\elm

\proof
(a) $\vf(0):=0\mathrel{\lower1pt\hbox{$\stackrel{\lower7pt\hbox{\smash{\hbox
to0pt{$\scriptstyle\sqcap$\hss}$\scriptstyle\sqcup_0$}}}\cdot$}}
e_1=e_0$ by \era2{2.3}\,I0. Let $m,n\in\N$. Then $\vf(m+n):=(m+n)
\mathrel{\lower1pt\hbox{$\stackrel{\lower7pt\hbox{\smash{\hbox
to0pt{$\scriptstyle\sqcap$\hss}$\scriptstyle\sqcup_0$}}}\cdot$}}e_1
\nad{\era2{2.3}\,{\rm I2}}= (m\mathrel{\lower1pt\hbox{$\stackrel{\lower7pt\hbox{\smash{\hbox
to0pt{$\scriptstyle\sqcap$\hss}$\scriptstyle\sqcup_0$}}}\cdot$}}e_1)
\qu_0(n\mathrel{\lower1pt\hbox{$\stackrel{\lower7pt\hbox{\smash{\hbox
to0pt{$\scriptstyle\sqcap$\hss}$\scriptstyle\sqcup_0$}}}\cdot$}}e_1)=
\vf(m)\qu_0\vf(n)$. Let $x,y\in I(e_1)$. Then \te\ ${m,n\in\N}$ \st $x=\vf(m)$,
$y=\vf(n)$. We have $x\qu_0y=\vf(m)\qu_0\vf(n)=\vf(m+n)\in I(e_1)$. Since
$e_0\in I(e_1)$, we infer that $I(e_1)$ is a \sbm\ of $(X,\qu_0,e_0)$.

(b) $\vf(1):=1\mathrel{\lower1pt\hbox{$\stackrel{\lower7pt\hbox{\smash{\hbox
to0pt{$\scriptstyle\sqcap$\hss}$\scriptstyle\sqcup_0$}}}\cdot$}}e_1=e_1$ by
\era2{2.3}\,I1. Let $m,n\in\N$, then $\vf(mn)\nda22.8 = {\vf(m\dpl n)}=\vf(m)
\qu_1 \vf(n)$. The last \et y follows from part~(i) of the proof of Theorem
\rf{t1.9n}. Observe that the proof of this \et y does not require the \sr~$X$
to be infinite or \pn. As in part~(a) one shows that $I(e_1)$ is a \sbm\ of the monoid
$(X,\qu_1,e_1)$.

(c) We conclude from (a), (b) that $I(e_1)$ is a sub\sr\ of the \sr~$X$ and
that $\vf:\N\to I(1)$ is a ring-\hm sm. Since $I(e_1)$ is the range of~$\vf$,
the map~$\vf$ is sur\jc.
\endproof

\blm4.34
A finite \sr\ of a field  $F$ is a subfield of~$F$.
\elm

\proof
Let $(F,\qu_0,\qu_1,e_0,e_1)$ be a field, hence a \sr, and let $E$~be a
sub\sr\ of~$F$.

(a) \ti{$E$ is a subgroup of $(F,\qu_0,e_0)$}: By \as\ $E$ is a finite \sbm\
of the \ag\ $(F,\qu_0,e_0)$. Then $E$ is a subgroup of~$F$ by Lemma
\rf{l3.21}.

(b) \ti{$E$ is a \sbm\ of $(F,\qu_1,e_1)$}: By \E\df\ \rf{d4.32}\,(b).

(c) \ti{Every \el\ of $E\sms{e_0}$ has an inverse in the monoid
$(E,\qu_1,e_1)$}: We first show that $E\sms{e_0}$ is a \sbm\ of the group
$(F\sms{e_0},\qu_1,e_1)$. Clearly $E\sms{e_0}\sbs F\sms{e_0}$ and $e_1\in F\sms{e_0}$ since
$e_0\ne e_1$. Let $a,b\in E\sms{e_0}$. Then $a\qu_1b\in
F\sms{e_0}$. Otherwise $a\qu_1 b=e_0$, and $a=a\qu_1e_1 = a\qu_1(b\qu_1
b\Inv) = (a\qu_1b)\qu_1 b\Inv = e_0\qu_1 b\Inv \nde1.1n = e_0$, a~\cd ion.
Here $b\Inv$ is the inverse of~$b$ in $(F,\qu_1,e_1)$. \If that $E\sms{e_0}$
is a \sbm\ of the group $(F\sms{e_0},\qu_1,e_1)$. Since $F$~is finite, $E\sms
{e_0}$ is finite by \era2{3.13}, and as a finite \sbm\ of an \ag, $E\sms{e_0}$
is a \ti{subgroup} of the group $(F\sms{e_0},\qu_1,e_1)$ by Lemma \rf{l3.21}.
\E\Tf \fe $a\in E\sms{e_0}$ \te s $b\in E\sms{e_0}$ \st $a\qu_1 b=e_1$, which
proves~(c).

(d) In view of \E\df\ \rf{d4.23}, $E$ is a subfield of $F$.
\endproof

\Wanp \ch ize finite \pf s.

\bth4.35
Let $(F,\qu_0,\qu_1,e_0,e_1)$ be a field and let $\check F$ denote its \pf.
Let $\vf:\N\to F$ be the map defined in \er{1.5n} and let $I(e_1)$ denote its
range in~$F$.

\hph i,i, If $\check F$ is finite and $n:=\#(\check F)$, then the \fw\
holds{\rm:}
\beq4.41
\vf(n)=e_0 \qh{and } \vf(k)\ne e_0 \qh{\fa}k\in\zo1,n .
\e
\hph ii,, If \te s $n\in\Na\sms1$ \st \er{4.41} holds, then $n$~is a
\emph{\Pn} which we denote by~$p$. \Mo we have
\beq4.42
\check F = I(e_1),
\e
and $\ov\vf$, the \rt ion to $\zo0,p $ of $\vf$, is a ring-\is sm from the
field $(\N_p,+_p,\cdot_p,0,1)$ onto the field $(\check F,\qu_0,\qu_1,e_0,e_1)$.
\E\Ip $\check F$ is finite and $\#(\check F)=p$.
\eth

\bdf4.36
Let $F$ be a field. If \te s $n\in\Na\sms1$ \st \er{4.41} holds, then $n$~is
called the \ti{\ch istic} of the field~$F$\index{characteristic of a field} and is denoted by $\chr(F)$.
If $\vf(n)\ne e_0$ \fa $n\in\Na$, then $\chr(F):=0$.
\edf

\proof[Proof of Theorem \rf{t4.35}] \

(i) By \as\ the abelian additive group $(\check F,\qu_0,e_0)$ of~$\check F$ is
finite. Hence it is a finite \Cm\ by \E\Pr\ \rf{p3.10}\,(i). Then (i)~follows
from Lemma \rf{l1.22} (i) and~(vi).

(ii) By Lemma \rfa2{l1.21}\,(i)--(ii) with $\wt E:=\N$, $(M,\qu,e):=(F,\qu_0,
e_0)$, $M_0:=\check F$ and $a:=e_1\in \check F\sms{e_0}$, we find that $I(e_1)$
is a \ti{\sbm} of the group $(\check F,\qu_0,e_0)$. \Mo $\vf:(\N,+,0)\to
(I(e_1),\qu_0,e_0)$ is a \ti{sur\jc\ \hm sm} by~(iii) of the same lemma. Since
$\vf(0)=e_0=\vf(n)$, $0\ne n$, the map $\vf$ is \ti{not in\jc}. \If from Lemma
\rfa2{l1.41} with $(\wt E,+_{\td E},\wt e):=(\N,+,0)$ that $I(e_1)$, the range
of~$\vf$ is \ti{finite}. \If from Lemma \rf{l4.33} with $X:=\check F$ that
$I(e_1)$ is a \ti{sub\sr} of the field~$\check F$. \E\Tf $I(e_1)$ is a
\ti{subfield\/} of~$\check F$ by Lemma \rf{l4.34}. Since $\check F$ is a \Pf,
$I(e_1)=\check F$. Observe that $(\check F,\qu_0,e_0)$, as a finite \ag, is a
finite \Cm\ by \E\Pr\ \rf{p3.10}\,(i). Hence by Lemma \rf{l2.12} with
$(Y,\qu,e):=(\check F,\qu_0,e_0)$, $a:=e_1$, and by \er{4.41} we find
$\#(\check F)=n$. We claim that $n$~is a \Pn. Otherwise \te\ $k,l\in\Na\sms
{1,n}$ \st $kl=n$. Since $k$~divides~$n$, we obtain $k\le n$ by \era3{8.57},
hence $k\in(1,n)$. Similarly $l\in(1,n)$. From \er{4.41} we have $\vf(k),
\vf(l)\in I(e_1)\sms{e_0}=\check F\sms{e_0}$. Since $\check F\sms{e_0}$ is
a \sbm\ of $(\check F,\qu_1,e_1)$ and $\check F=I(e_1)$, $\vf(k)\qu_1\vf(l)
\ne e_0$. However by Lemma \rf{l4.33} and \er{4.42}, the map~$\vf$ is a \hm sm
from $(\N,\cdot,0)$ into $(\check F,\qu_1,e_1)$, hence $\vf(k)\qu_1\vf(l)=
\vf(kl)=\vf(n)=e_0$. A~\cd ion. \E\Tf $n$~is a \ti{\Pn} which we denote
by~$p$.

\Mo by Lemma \rf{l1.22}\,(i) and (v) with $(X,\qu,e):=(\check F,\qu_0,e_0)$,
$n:=p$ and $a:=e_1$, we find that $\ov\vf$, the \rt ion of~$\vf$ to $\zo0,p $,
is a \ti{monoid-\is sm} from $(\N_p,+_p,0)$ onto $(I(e_1),\qu_0,e_0)$. By~(iv)
of the same lemma, we have
\beq4.43
\vf(k)=\vf(\F_p(k)) \qh{\fa} k\in\zo0,p .
\e
\Wanp prove that $\ov \vf$ is a \ti{\hm sm from} $(\N_p,\cdot_p,1)$ into
$(\check F,\qu_1,e_1)$. Let $k,l\in\zo0,p $. Then $\ov\vf(k\cdot_p l)=\vf
(k\cdot_p l)\nde1.33 = \vf(\F_p(kl))\nde3.26 = \vf(\F_p(\F_p(k)\cdot \F_p
(l)))\nde4.43 = \vf(\F_p(k)\cdot\F_p(l))\nde1.5 = \vf(k\cdot l)\nad*= \vf(k)
\qu_1\vf(l) = \ov\vf(k)\qu_1\ov\vf(l)$, where in $\nad*=$ we used Lemma
\rf{l4.33}. Since $\ov\vf:\N_p \to \check F$ is bi\jc, it is also an \is sm
from $(\N_p,\cdot_p,1)$ onto $(\check F,\qu_1,e_1)$. \csq, it is a ring-\is sm
from $\N_p$ onto~$\check F$.
\endproof

\bpr4.37
Let $E$ be a field of \ch istic $p$, $p$~\Pn. Let $F$ be a field \st $E$~is
a subfield of~$F$ or, more generally, \st \te s an in\jc\ ring-\hm sm
$j:E\to F$. Then $\chr(F)=p$.
\epr

\proof
Let $e_0$ (resp.\ $e_1$) denote the \nel\ of the additive (resp.\ \mlv) monoid
of~$E$. Let $\vf:\N\to E$ denote the \sq\ of \IT s of~$e_1$ in the additive
monoid of~$E$. By \as\ we have $\vf(0)=e_0$, $\vf(1)=e_1$, $\vf(p)=e_0$ and
$\vf(k)\ne e_0$ \fa $k\in\zo1,p $. Since $j$~is a ring-\hm sm from~$E$ into~$F$,
$j(e_0)$ (resp.\ $j(e_1)$) is the \nel\ of the additive (resp.\ \mlv) monoid
of~$F$. Let $\psi:\N\to F$ denote the \sq\ of \IT s of $j(e_1)$ in the
additive monoid of~$F$. By Lemma \rfa2{l1.22} we obtain
\beq4.44
\psi(n)=j(\vf(n)) \qh{\fa} n\in\N.
\e
Since $\vf(p)=e_0$, $\psi(p)=j(\vf(p))=j(e_0)=\psi(0)$. In view of the in\ji\
of the map~$j$, we have $j(\vf(k))\ne j(e_0)=\psi(0)$ \fa $k\in\zo1,p $. \If
from Theorem \rf{t4.35} that $\chr(F)=p$. Note that if $E$~is a subfield
of~$F$, then the inclusion map $j:E\to F$ defined by $j(x)=x$, $x\in E$, is an
in\jc\ \hm sm.
\endproof

\bex4.38
Let $E$ be a subfield of the field~$F$. Let $j:E\to F$ be defined by
$j(x):=x$, $x\in E$. Show that $j$ is an in\jc\ ring-\hm sm from~$E$ as a
field into the field~$F$.
\eex

\brm4.39
The \pf\ of a finite field is finite by \E\Pr\ \rfa2{p4.18}\,(i), hence every
finite field has \ch istic~$p$ \fs $p$~prime. \E\te\
\ti{infinite} fields with \ch istic~$p$ (see Remark \rfa5{r2.32}).
\erm

We next show that the \mlv\ group of a finite \Pf\ of
order larger than two is \ti{\pn} (= \ti{cyclic}) and we give some application of
this result.

We first give some \df s.

\def\namespec{Definitions and Notation}
\begin{dspc}\lb{d4.40}\

\hph i,i, Let $(R,+,\cdot,0,1)$ be a \ti{ring}, let $a\in R$ and let $k\in\N$.
Then the $k$-th \IT\ of~$a$ in the \mlv\ monoid $(R,\cdot,1)$ is denoted
by~$a^k$ and is called the $k$-th \ti{power} of~$a$.

\hph ii,, Let $(F,+,\cdot,0,1)$ be a \ti{field\/}. An \el\ $\o\in F$ \sf ying
$\o^n=1$ \fs $n\in\Na$ is called a $n$-th \ti{\ru{}}, or simply, a \ru{}.\index{root of unity} It is
called a \ti{primitive} $n$-th \ru{}\index{primitive $n$-th root of unity} if $\o^n=1$ and $\o^m\ne 1$ for $1\le m<n$.
\end{dspc}

Since $1^n\nad{\era2{2.3}\,{\rm I4}}= 1$ \fa $n\in\Na$, $1$~is a $n$-th \ru{} \fa
$n\in\Na$. Since $\o^1\nad{\era2{2.3}\,{\rm I1}} = \o$, $1$~is the only $1$-th \ru{}.

\ssk
We now prove that $(\N_p\sms0,\cdot_p,1)$, $p$ odd prime, is cyclic. Observe that
$(\N_2\sms0,\cdot_2,1)$ is the trivial group, hence not cyclic by \df. As a
con\sq\ of Corollary \rf{c3.72} it suffices to show that \fe divisor~$d$
of~$p-1$ with $d\in(1,p-1)$, \te s at most one \pn\ \sbm\ of order~$d$. Let
$H$ be such a subgroup, then the $d$-th \IT\ of every \el\ in $(H,\cdot_p,1)$
is equal to~$1$ by \E\Pr\ \rf{p3.65}\,(i). By Lemma \rfa2{l1.21}\,(i) if $a\in
H$ and $d\ddtt pa$ is the $d$-th \IT\ of~$a$ in the group $(H,\cdot_p,1)$
then $d\ddtt pa$ is also the $d$-th \IT\ of~$a$ in the group $(\N_p\sms0,
\cdot_p,1)$ and in the \mlv\ monoid $(\N_p,\cdot_p,1)$. In view of Definitions
and Notation \rf{d4.40} $d\ddtt p a=a^d$. Since $a^d=1$, $a$~is a $d$-th \ru{}
of the \ti{field\/} $(\N_p,+_p,\cdot_p,0,1)$. Thus every \el\ of~$H$ is a $d$-th \ru{} of the
field~$\N_p$. Note that $\#(H)=d$, hence there are at least~$d$ $d$-th \ru s
in~$\N_p$. Suppose now for \cd ion that \te s \ti{another} \sbm~$H'$ of
order~$d$, then \te s at least one \el\ $b\in H'\setminus H$, hence $b$~is also
a $d$-th \ru{} in the field~$\N_p$. It turns out that in a field~$F$ there are
\ti{at most\/}~$d$ $d$-th \ru s \fe $d\in\Na$ (see \E\Pr\ \rf{p4.41} below).
This leads to a \cd ion. \E\Tf \fe divisor $d\in(1,p-1)$ of $p-1$ there is
at most one \sbm\ of $(\N_p\sms0,\cdot_p,1)$ of order~$d$. Hence $(\N_p\sms0,
\cdot_p,1)$ is a \fcg\ of order $p-1$, by Corollary \rf{c3.72}. This result
is due to Gauss.

\bpr4.41
Let $F$ be a field. Then there are at most $n$ $n$-th roots of unity in~$F$.
\epr

See the proof below the proof of Theorem \rf{t4.48}.

\begin{cor}[\cite{Groups}]\lb{c4.42}\

\hph i,i, If $(F,+,\cdot,0,1)$ is a field and $G$ is a finite nontrivial
subgroup of the monoid $(F,\cdot,1)$, then $G$ is \emph{cyclic}.

\hph ii,, If $F$ is finite with $\#(F)>2$, then its \mlv\ group is \emph{cyclic}.
\E\Ip $(\N_p^\t,\cdot_p,1)$, $p$~odd prime, is cyclic.
\eco

\bex4.43
Prove Corollary \rf{c4.42}. (Imitate the proof of the case $F:=\N_p$ given
above \E\Pr\ \rf{p4.41}.)
\eex

\bxs4.44 \

\hph i,i, Let $p$ be an odd \Pn. Then the field $(\N_p,+_p,\cdot_p,0,1)$
possesses \ti{exactly} $p-1$ $(p-1)$-th \ru s. Indeed, since $([1,p-1],\cdot_p,
1)$ is a group, we have $\o^{p-1}=1$ \fa $\o\in[1,p-1]$ by \E\Pr\
\rf{p3.65}\,(i).

We give an example of four \el s of a \ti{ring}~$R$ \sf ying the \eq\ $x^2=1$.

\hph ii,, Let $(F,+,\cdot,0,1)$ be a \ti{field\/} and let $(F\t F,\wt+,\td
\cdot,(0,0),(1,1))$ denote the direct product of~$F$ with itself defined in
\E\Pr\ \rf{p2.2}. Then $F\t F$ is a \ti{ring} since the direct product
$(F\t F,\wt+,(0,0))$ is a group by Corollary \rf{c3.37}\,(i). However, the
ring $F\t F$ is \ti{not\/} a field, since $(1,0)\mathrel{\wt\cdot} (0,1)=(0,0)$.
\ti{Suppose that\/} $\chr(F)\ne2$ (see \E\df\ \rf{d4.36}), then $1\ne-1$ where
$-1$ denotes the inverse of~$1$ in the additive group $(F,+,0)$. Indeed, if
$1=-1$, then $0=1+(-1)=1+1$, hence $\chr(F)=2$. \E\Ip the fields $\N_p$,
$p$~odd prime, have \ch istic equal to $p\ne2$. We claim
\beq4.45
(-1)^2=1.
\e
Indeed, $(-1)^2=(-1)\cdot(-1)\nde4.12 = -(1\cdot(-1))\nda21.7 = -(-1)=1$
since $(-1)+1=0$. \csq, if $\a_1:=(1,1)$, $\a_2:=(1,-1)$, $\a_3:=
(-1,1)$ and $\a_4:=(-1,-1)$, we obtain $(\a_i)^2=(1,1)$, $i\in[1,4]$. Thus in
the ring $F\t F$ the \eq\ $\a^2=(1,1)$ possesses at least 4~\so s. Conversely,
observe that if $a,b\in F$, then $(a,b)^2=(1,1)$ implies $a^2=b^2=1$. By what
precedes $1$~and~$-1$ \sf y the \eq\ $x^2=1$ in~$F$. One verifies that $x^2-1
=(x+1)(x-1)$ \fe $x\in F$. Note that if $u,v\in F\sms0$, then $u\cdot v\in
F\sms0$, then $u\cdot v\in F\sms0$ since otherwise $v=1\cdot v=(u\mo\cdot u)
\cdot v= u\mo\cdot(u\cdot v)=u\mo\cdot0\nde4.1 =0$, a~\cd ion. Thus if $x^2-1=
(x+1)(x-1)=0$, then either $x=-1$ or $x=1$. \csq, $1$~and~$-1$ are the only
$2$-th \ru s in the field~$F$, and $\a_1,\a_2,\a_3,\a_4$ are the only \so s
to the \eq\ $(a,b)^2=(1,1)$ in $F\t F$.
\exs


\blm4.45
Let $(R,+,\cdot,0,1)$ be a ring, $\a,\b\in R$, $\a\ne\b$ and let $n\in\N$.
Then
\beq4.46
\a^{n+1}-\b^{n+1}= (\a-\b)\cdot \sum_{k=0}^n \a^{n-k}\b^k.
\e
\elm

\proof
$n\in\Na$: Set $A:=\a\suml_{k=1}^n \a^{n-k}\b^k$ and $B:= \b\suml_{k=0}^{n-1}
\a^{n-k}\b^k$. Then\\
$(\a-\b)\cdot\suml_{k=0}^n \a^{n-k}\b^k \nad{\er{4.8},\er{1.3n}} =
\a\suml_{k=0}^n \a^{n-k}\b^k + (-\b)\suml_{k=0}^n \a^{n-k}\b^k \nad{\era2{7.13},\er{4.12},\er{1.2n}} =
(\a\a^n+A) + \bigl(-\bigl(\b\suml_{k=0}^n \a^{n-k}\b^k\bigr)\bigr) \nad{\er{4.8},\era2{7.13},\er{1.2n}} =
(\a\a^n+A) - (B+\b\b^n) \nad{\era2{2.3}\,{\rm I1,I2},\er{4.11}}=
(\a^{n+1}+A)+((-B)+(-\b^{n+1})) \nda21.36 =
\a^{n+1}+(A+(-B))+(-\b^{n+1})$. Thus it suffices to show that ${A=B}$. We have
$A\nde4.16 = \suml_{k=1}^n \a(\a^{n-k}\b^k) = \suml_{k=1}^n (\a\a^{n-k})\b^k
\nad{\era2{2.3}\,{\rm I1,I2}} = \suml_{k=1}^n \a^{1+(n-k)}\b^k$.
\E\oh $B\nde4.16 = \suml_{k=0}^{n-1}\b(\a^{n-k}\b^k) = \suml_{k=0}^{n-1}
(\b\a^{n-k})\b^k = \suml_{k=0}^{n-1}(\a^{n-k}\b)\b^k = \suml_{k=0}^{n-1}
\a^{n-k}(\b\b^k) = \suml_{k=0}^{n-1}\a^{n-k}\b^{k+1}\break \nda27.12 = \suml_{l=1}^{n}
\a^{n-(l-1)}\b^l \nda27.11 = \suml_{k=1}^{n}\a^{n-(k-1)}\b^k$. Since $k-1\le n$
whenever $k\in[1,n]$, \te s $p\in\N$ \st $n=p+(k-1)$. Hence $n+1 = (p+(k-1))
+1 = p+((k-1)+1)=p+k$. \E\Tf $p=(n+1)-k$ and $n-(k-1)=p=(n+1)-k=(1+n)-k$. \If
that $B=\suml_{k=1}^n \a^{(1+n)-k}\b^k=A$.

$n=0$: $\a^{0+1}-\b^{0+1}= \a-\b = (\a-\b)\cdot1 = (\a-\b)\cdot \suml_{k=0}^0
\a^{0-k}\b^k$.

Note that \er{4.8} holds in ring.
\endproof

\bxs4.46
Let $R$ be a ring and let $\a,\b\in R$, $\a\ne\b$. Then
\bea4.47
\a^2-\b^2 &= (\a-\b)(\a+\b),\\
\a^3-\b^3 &= (\a-\b)(\a^2+\a\b+\b^2),\lb{4.48}\\
\a^4-\b^4 &= (\a-\b)(\a^3+\a^2\b+\a\b^2+\b^3)\lb{4.49},\\
\a^5-\b^5 &= (\a-\b)(\a^4+\a^3\b+\a^2\b^2+\a\b^3+\b^4).\lb{4.50}
\e
Observe that \er{4.46}--\er{4.50} trivially hold if $\a=\b$.
\exs

\blm4.47
Let $R$ be a ring, let $n\in\N$, $\ov x\in R$ and $a:[0,n+1]\to R$ with
$a_{n+1}=1$. Suppose that $\suml_{k=0}^{n+1} a_k\ov x{}^k=0$. Then \te s
$b:[0,n]\to R$ with $b_n=1$ \st
\beq4.51
\sum_{k=0}^{n+1} a_kx^k = (x-\ov x)\cdot \sum_{l=0}^n b_lx^l \qh{\fa}
x\in R\sms{\ov x}.
\e
\elm

\proof
Since $0+0\nde1.1n = 0$, $-0=0$. Let $x\in R\sms{\ov x}$. Then $\suml_{k=0}^{n+1}a_kx^k =
{\suml_{k=0}^{n+1}a_kx^k +0} = \suml_{k=0}^{n+1}a_kx^k -
\suml_{k=0}^{n+1}a_k\ov x{}^k = \suml_{k=0}^{n+1}a_kx^k
+\Bigl(-\suml_{k=0}^{n+1}a_k \ov x{}^k\Bigr) \nde4.17 = \suml_{k=0}^{n+1}a_k
x^k +\suml_{k=0}^{n+1}(-(a_k\ov x{}^k)) \nda21.129 =\break \suml_{k=0}^{n+1}\bigl((a_k
x^k) + (-(a_k\ov x{}^k))\bigr) \nde4.12 = \suml_{k=0}^{n+1}\bigl((a_kx^k)
+(a_k(-\ov x{}^k))\bigr) \nde1.2n = \suml_{k=0}^{n+1}a_k(x^k +(-\ov x{}^k))
\nad{\era2{2.3}\,{\rm I0},\era2{7.13}} = a_0(1+(-1)) + \suml_{k=1}^{n+1}a_k
(x^k+(-\ov x{}^k)) \nad{\er{1.1n},\er{4.8}}= 0+\suml_{k=1}^{n+1}a_k
(x^k-\ov x{}^k) \nde4.46 =\\ \suml_{k=1}^{n+1}\Bigl(a_k(x-\ov x)\suml_{l=0}
^{k-1} x^{(k-1)-l}\ov x{}^l\Bigr) \nad{\era2{1.6},\er{4.16}}=
(x-\ov x)\cdot\suml_{k=1}^{n+1}a_k\suml_{l=0}^{k-1} x^{(k-1)-l}\ov x{}^l$.

We now claim $\suml_{k=1}^{n+1}a_k\suml_{l=0}^{k-1}x^{(k-1)-l}\ov x{}^l  =
\suml_{l=0}^n b_lx^l$ where
\beq4.52
b_l:=\sum_{j=l}^n a_{j+1}\ov x{}^{j-l}, \q l\in[0,n].
\e
Indeed, $\suml_{k=1}^{n+1}a_k\suml_{l=0}^{k-1}x^{(k-1)-l}\ov x{}^l
\nde4.31 = \suml_{j=0}^{n}a_{j+1} \suml_{l=0}^{j} x^{j-l}\ov x{}^l
\nda27.14 = \suml_{j=0}^{n}a_{j+1}\suml_{l=0}^{j}x^l\ov x{}^{j-l}
=\suml_{j=0}^{n}a_{j+1}\break \suml_{l=0}^{j}\ov x{}^{j-l}x^l \nde4.16 =
\suml_{j=0}^{n}\suml_{l=0}^{j}a_{j+1}\ov x{}^{j-l}x^l \nad*=
\suml_{l=0}^{n}\suml_{j=l}^{n}a_{j+1}\ov x{}^{j-l}x^l \nde4.16 =
\suml_{l=0}^{n}x^l\suml_{j=l}^{n}a_{j+1}\ov x{}^{j-l} \nde4.52 =
\suml_{l=0}^n x^lb_l = \suml_{l=0}b_lx^l$.

In $\nad*=$ we used \era2{7.22} with $(X,\qu,e):=(R,+,0)$, $E:=\N$, $m:=n$,
$i:=j$, $j:=l$, $\lower2pt\hbox{$\QU$}:=\sum$, $a_{jl}:=a_{j+1}\ov x{}^{j-l}x^l$. Since the
proof of \era2{7.22} was left as an exercise, we give a proof now.

\proof[Proof of \era2{7.22}]
We use \era2{1.127} twice with $\O:=\bigl\{(k,l)\in[0,m]\t[0,m]: l\le
k\bigr\}$. In the first case, $A_i:=\{(i,j)\in\O: j\le i\}$, $i\in[0,m]$.
Then $\QU_{\o\in\O}a_\o \nda2{1.127} = \QU_{i\in[0,m]}\bigl(\QU_{(i,j)\in
A_i} a_{ij}\bigr) = \QU_{i=0}^m \bigl(\QU_{j=0}^i a_{ij}\bigr)$.
In the second case, $A_j:=\{(i,j)\in\O: j\le i\}$, $j\in[0,m]$. Then
$\QU_{\o\in\O}a_\o \nda2{1.127} = \QU_{j\in[0,m]}\bigl(\QU_{(i,j)\in A_j}
a_{ij}\bigr) = \QU_{j=0}^m \bigl(\QU_{i=j}^m a_{ij}\bigr)$.
\endproof

We conclude the proof of the lemma by showing that $b_n=1$. We have $b_n
\nde4.52 = a_{n+1}\ov x{}^0 \nad{\era2{2.3}\,{\rm I0}}=1\cdot 1=1$.
\endproof

\bth4.48
Let $F$ be a field, let $d\in\Na$ and let $a:[0,d]\to F$ with $a_d=1$.
Suppose that \te\ $\zb xl{[1,d]}\sbs F$ \st $x_i\ne x_j$ whenever $i\ne j$,
$i,j\in[1,d]$, and \st $\suml_{k=0}^d
a_k(x_l)^k=0$ \fe $l\in[1,d]$. Then $\suml_{k=0}^d a_kx^k\ne 0$ \fa $x\in
F\sm\zb xi{[1,d]}$.
\eth

\proof
By \In\ on $d\in\Na$. Set $M:=\{l\in\Na: \hbox{Theorem \rf{t4.48} holds for }d:=l\}$.

$1\in M$: Let $a_0\in F$ and $x_1\in F$ be \st $a_0+x_1=0$, i.e.\ $x_1=-a_0$.
If $x\in F\sms{x_1}$, we have $a_0+x=a_0+(x-x_1)+x_1 = a_0+x_1+(x-x_1) = 0
+(x-x_1)= x-x_1\ne0$.

\ti{$d\in M$ implies $d+1\in M$}: Let $d\in M$, let $a\in[0,d+1]$ with
$a_{d+1}=1$ and let $x_i$, $i\in[1,d+1]$ be \st $\suml_{k=0}^{d+1} a_k(x_i)^k
=0$, and $x_i\ne x_j$, $i\ne j$, $i,j\in[1,d+1]$. By Lemma \rf{l4.47}, \te s
$b:[0,d]\to F$ with $b_d=1$ \st
\beq4.53
\sum_{k=0}^{d+1} a_kx^k = (x-x_{d+1})\cdot \sum_{l=0}^d b_lx^l, \q
x\in F\sms{x_{d+1}}.
\e
Then $0=\suml_{k=0}^{d+1} a_k(x_i)^k = (x_i-x_{d+1})\cdot \suml_{l=0}^d
b_l(x_i)^l$ \fa $x_i$, $i\in[1,d]$. Since $x_i\ne x_{d+1}$ \fa $i\in[1,d]$,
we have $\suml_{l=0}^d b_l(x_i)^l=0$ \fa $i\in[1,d]$. Indeed, if $y,z\in F\sms0
$, then $y\cdot z\in F\sms0$, since $F\sms0$ is a subgroup of the monoid
$(F,\cdot,1)$. Since $d\in M$ and since $x_i\ne x_j$, $i\ne j$, $i,j\in[0,d]$,
$b_d=1$, we have $\suml_{l=0}^d b_lx^l\ne0$ \fa $x\in F\sm\zb xi{[1,d]}$. From
\er{4.53} we infer $\suml_{k=0}^{d+1}a_kx^k\ne0$ \fa $x\notin\zb xi{[1,d+1]}$.
Hence $d+1\in M$. Thus $M=\Na$.
\endproof

\proof[Proof of \E\Pr\ \rf{p4.41}]
We apply Theorem \rf{t4.48} with $d\in\Na\sms1$, and $a_0:=-1$, $a_i:=0$,
$i\in[1,d-1]$. The case $d:=1$ is trivial.
\endproof

We now consider an application and extension of Corollary \rf{c4.42}\,(ii).
We denote by $(\N_n,+_n,\cdot_n,0,1)$, $n\ge3$, the \pn\ rings introduced in\glossary{$(\N_n,+_n,\cdot_n,0,1)$}
\E\Pr\ \rf{p3.15}. Let $\N_n^\t$ denote the set of units of the ring~$\N_n$
(see \E\df\ \rf{d3.34}). Then $\N_n^\t = \{k\in[1,n]: \gcd(k,n)=1\}$ (see
\E\Pr\ \rf{p3.44}), $\N_n^\t$ is a subgroup of the monoid $(\N_n,\cdot_n,1)$
and $\#(\N_n^\t)=\F(n)$ (see \er{3.18}). If $n:=p$ odd prime, then $\N_p$~is
a field and $\N_p^\t=\N_p\sms0$. Corollary \rf{c4.42}\,(ii) tells us that the
group $(\N_p^\t,\cdot_p,1)$ is cyclic. \E\Gn s of such groups are called
\ti{primitive roots} by Gauss, and their \nm\ is given by $\F(p-1)$. For
example, if $p=3$, $\F(2)=2-1=1$ \er{3.23n}, and $2$ is the primitive root since
$2\cdot_32=\F_3(2\cdot2)=\F_3(4)=1$. If $p=5$, $\F(4)=\F(2^2)\nde3.39 = 2^1
\cdot1= 2$. One verifies that $\N_5^\t = \{1,2,3,4\}$, $2^1=2$, $2^2=\F_5(4)
=4$, $2^3=\F_5(8)=3$, $2^4=\F_5(16)=1$ and $3^1=\F_5(3)=3$, $3^2=\F_5(9)=4$,
$3^3=\F_5(27)=\F_5(3\cdot4)=\F_5(12)=2$, $3^4=\F_5(2\cdot3)=\F_5(6)=1$. Hence
the two primitive roots are $2$~and~$3$. Observe that $4^1=4$, $4^2=\F_5(16)
=1$, hence the order of~$4$ is~$2$ and $4$~is \ti{not\/} a primitive root
of~$\N_5$.

\bex4.49
Find the primitive roots of $\N_p$ for $p=7,11,13,17,19$.
\eex

It turns out that the group of units of the ring $\N_n$ with $n=p^k$, $p$~odd
prime, $k\in\Na$, is \ti{cyclic}. The proof of this result for $k>1$ is based
on the case $k=1$ and uses the \fw\ lemma.

\blm4.50
Let $a,b,n\in\Na$ with $n\ge2$. Then \te s $c_n\in\N$ \st
\beq4.54
(a+b)^n = a^n + na^{n-1}b+\frac{n(n-1)}{2}\,a^{n-2}b^2+c_nb^3.
\e
\elm

\proof
By \In\ on $n\in\Na\sms1$. Set $M:=\{k\in\Na\sms1: \hbox{\er{4.54} holds with }
n:=k\}$.

$2\in M$: $(a+b)^2 = (a+b)(a+b) = a(a+b) +b(a+b) = (a^2+ab)+(ba+b^2) =(a^2+ab)
+(ab+b^2) = a^2+(ab+ab)+b^2 = a^2+2ab+b^2$. Hence $2\in M$ with $c_2:=0$.

\ti{$k\in M$ implies $k+1\in M$}: $(a+b)^{k+1}=(a+b)(a+b)^k = a(a+b)^k
+b(a+b)^k \nad{k\in M}= (a^{k+1}+ka^kb + \frac{k(k-1)}2 a^{k-1}b^2 + ac_kb^3)
+(a^kb+ka^{k-1}b^2 + \frac{k(k-1)}2 a^{k-2}b^3+c_kb^4 = a^{k+1}+(k+1)a^kb
+(\frac{k(k-1)}2+k) a^{k-1}b^2 + (ac_k +\frac{k(k-1)}2 a^{k-2}+c_kb)b^3 =
a^{k+1}+(k+1)a^{(k+1)-1}b + \frac{(k+1)((k+1)-1)}2 a^{(k+1)-2}b^2
+c_{k+1}b^3$, where $c_{k+1}:=(a+b)c_k+\frac{k(k-1)}2 a^{k-2}$. Hence $k+1\in
M$.

\If that $M$ is \iv\ in $\Na\sms1$, hence $M=\Na\sms1$.
\endproof

\bpr4.51
Let $p$ be an odd \Pn\ and let $k\in\Na$. Then the group of units of the ring
$(\N_{p^k},+_{p^k},\cdot_{p^k},0,1)$ $($see \E\df\ \rf{d3.34}$)$ is
\emph{cyclic}.
\epr

The case $k:=1$ follows from Corollary \rf{c4.42}\,(ii). We shall use the
\fw\ notation:

\hph i,i, Given $a,l\in\N$, $a^l$ denotes the $l$-th \IT\ (power) of~$a$ in
the monoid $(\N,\cdot,1)$.

\hph ii,, Given $n\in\Na\sms1$, $a\in\zo0,n $ and $l\in\N$, $l\ddtn a$
denotes the $l$-th \IT\ of~$a$ in the monoid $(\N_n,\cdot_n,1)$.

We recall \er{3.71}:
\beq4.55
l\ddtn a=\F_n(a^l),\q l\in\N,\ a\in\zo0,n ,\ n\ge2.
\e
Note
\beq4.56
\F_n(a^{l+1})=\F_n(\F_n(a^l)\cdot a), \q l\in\N,\ a\in\zo0,n ,\ n\ge2,
\e
where $\cdot$ denotes the \mlc\ in~$\N$. Indeed, $\F_n(a^{l+1})=
\F_n(a^l\cdot a)\nde3.26 = \F_n(\F_n(a^l)\cdot\F_n(a))\nde1.5 = \F_n(\F_n
(a^l)\cdot a)$. Observe that if $l\in\Na$, $n>2$ and $a\in\zo2,n $, then
\beq4.57
l\ddtn a=1 \qh{iff \te s $q\in\Na$ \st }a^l-1=q\cdot n.
\e
Indeed, note that $a^l>1$ by \era2{2.23}, \era2{2.33} since $a\ge2$, and that
$l\ddtn a=1$ iff $\F_n(a^l)=1$ by \er{4.55} iff $a^l=q\cdot n+1$ \fs $q\in\N$
iff $a^l-1=q\cdot n$ \fs $q\in\N$ iff $a^l-1=q\cdot n$ \fs $q\in\Na$ by
\era2{2.13} since $a^l-1\in\Na$. \E\Tf we have
\beq4.58
l\ddtn a=1 \qh{iff $n$ divides} a^l-1.
\e

\ssk
Before giving the proof of \E\Pr\ \rf{p4.51}, we give an example.

\bxa4.52
Let $p:=3$ and $k:=2$. By \E\Pr\ \rf{p3.44}
$\N^\t_{3^2}:=\{l\in\zo1,9 :\break \gcd(l,9)=1\} = \{1,2,4,5,7,8\} = {\rm Gen}
((\N_9,+_9,0))$. We have $6=\#({\rm Gen}(\N_9,+_9,0))\nde3.21 = \F(3^2)
\nde3.39 = 3^{2-1}\cdot(3-1)$. In order to show that $(\N^\t_9,\cdot_9,1)$ is
cyclic, it is \sft\ in view of \E\df\ \rf{d3.70} and Lemma \rf{l3.71} to find an
\el\ $a\in\N_9^\t$ of order $6=\#(\N_9^\t)$. For example if $a=2$, we have
$1\ddtt9 2=\F_9(2^1)=2$, $2\ddtt9 2=\F_9(2^2)=4$, $3\ddtt9 2=\F_9(2^3)=8$,
$4\ddtt9 2=\F_9(16)=7$, $5\ddtt9 2\nde4.56 = \F_9(7\cdot2)=5$, $6\ddtt9 2=
\F_9(5\cdot2)=1$. \If that $\N_9^\t=I(2)$, the set of \IT s of~$2$ in
$(\N_9^\t,\cdot_9,1)$, hence the group $(\N_9^\t,\cdot_9,1)$ is cyclic (\pn).

\Mo by Lemma \rf{l3.27}\,(i) the \nm\ of \Gn s of $(\N^\t_9,\cdot_9,1)$ is
equal to the \nm\ of \Gn s of $(\N_6,+_6,0)$, which is
$\F(6)=\F(2\cdot 3)\nde3.48 = 1\cdot2=2$. The \Gn s of $(\N_6,+_6,0)$ are
$1$~and its \ng~$5$. By \E\Pr\ \rf{p3.10} the map $\ov\vf_2:(\N_6,+_6,0)
\to (\N_9^\t,\cdot_9,1)$ defined by $\ov\vf_2(l):=l\ddtt9 2$ is an \is sm.
Hence $\ov\vf_2(5)=5\ddtt9 2=\F_9(2^5)=5$ is also a \Gn\ by \E\Pr\ \rf{p3.26}\,%
(ii). Indeed, $\F_9(5^1)=5$, $\F_9(5^2)=7$, $\F_9(5^3)=\F_9(7\cdot5)=8$,
$\F_9(5^4)=\F_9(8\cdot5)=4$, $\F_9(5^5)=\F_9(4\cdot5)=2$, $\F_9(2\cdot5)=1$.
Summarizing,
\beq4.59
\hbox{2 and 5 are the \Gn s of the cyclic group }(\N_9^\t,\cdot_9,1).
\e
\exa

\bex4.53 \

\hph i,ii, Show that $2$ and $3$ are \Gn s of $(\N_5^\t,\cdot_5,1)$ and
$(\N_{25}^\t,\cdot_{25},1)$.

\hph ii,i, Are $7$ and $8$ \Gn s of $(\N_{25}^\t,\cdot_{25},1)$\,?

\hph iii,, What is the \nm\ of \Gn s of $(\N_{25}^\t,\cdot_{25},1)$\,?

\hph iv,, What are the \Gn s of $(\N^\t_{25},\cdot_{25},1)$\,?
\eex

We now present a proof of \E\Pr\ \rf{p4.51} for the case $k=2$. We look for
an \el\ $a\in(\N^\t_{p^2},\cdot_{p^2},1)$ whose order is equal to
$\#(\N^\t_{p^2})\nde3.21 = \F(p^2)\nde3.39 = p^{2-1}(p-1) = p(p-1)\ge6$. Let
$a\in\N_{p^2}^\t\sms1$. Then $\#(I(a))$, the order of~$a$, which we denote by~$m$,
divides $\#(N_{p^2}^\t)=p(p-1)$ by \er{3.73}. Let $l\in\Na$
be \st $l\cdot m=p(p-1)$. We would like to have $l=1$. To this end we choose for
the \el~$a$, a primitive root $r$ of $\N_p^\t$, that is, a \Gn\ of the group
$(\N_p^\t,\cdot_p,1)$, which exists by Corollary \rf{c4.42}\,(ii). We have
\beq4.59a
(p-1)\ddtt p r=1  \qh{in} (\N_p^\t,\cdot_p,1)
\e
by \er{3.75}. Since $\gcd(r,p)=1$, we also have $\gcd(r,p^2)=1$. Indeed, by
Lemma \rf{l2.16}, the divisors of~$p^2$ are $1$,~$p$ and~$p^2$, \Tf
$\gcd(r,p^2)\in\{1,p,p^2\}$. If $\gcd(r,p^2)=p$, then $p|r$ which is impossible
since $\gcd(r,p)=1$. If $\gcd(r,p^2)=p^2$, then $p^2|r$, hence $p|r$, which is
impossible. \E\Tf $r\in\N^\t_{p^2}$. By the \df\ of~$m$ we have $m\ddtt{p^2}r=1$.
Since $r\ge2$, we have $m\ge2$, hence $r^m>1$ by \era2{2.33}. \If from
\er{4.58} that $p^2|r^m-1$. \csq, $p|r^m-1$, hence $m\ddtt p r=1$ by \er{4.58}.
Since $(\N_p^\t,\cdot_p,1)$ is a finite \Cm\ with $\#(\N_p^\t)=p-1\ge2$, and
since $r$ is a \Gn\ of $(\N_p^\t,\cdot_p,1)$, we infer from Lemma
\rf{l1.22}\,(v) and from $m\ddtt p r=1$ that $m=q\cdot(p-1)$ \fs $q\in\N$.
Since $m\ne0$, we have $q\ne0$ by \era2{2.13}. Since $p(p-1)=l\cdot m$,
we have $p(p-1)=l\cdot q(p-1)$. Since $p-1\ne0$, we obtain $p=lq$ by \era2{2.12}
and \E\Pr\ \rfa2{p2.7}\,(iii). \If that either $l=1$, $q=p$ or $l=p$, $q=1$
since $p$ is prime. In the first case, $m=p(p-1)$, hence $\N^\t_{p^2}=I(r)$,
and $(\N^\t_{p^2},\cdot_{p^2},1)$ is cyclic. If $l\ne1$, then $q=1$ and
\beq4.60
m=p-1 \qh{and }p^2|r^{p-1}-1.
\e
In this case, the idea is to replace $r$ by $r+p$. Set
\beq4.61
r':=r+p.
\e
We shall show that $r'\in\N^\t_{p^2}$, that its order in $(\N^\t_{p^2},
\cdot_{p^2},1)$ denoted by~$m'$ \sf ies $l'\cdot m'=p(p-1)$ \fs $l'\in\Na$
and that \te s $q'\in\Na$ \st $m'=q'\cdot(p-1)$. As above we obtain $l'q'=p$.
With this choice of~$r'$, we show that $p^2\nmid (r')^{p-1}-1$, hence
$q'\ne1$, which implies $l'=1$. Thus in this case the order of~$r'$ is
$p(p-1)$, and $(\N^\t_{p^2},\cdot_{p^2},1)$ is \ti{cyclic}.

$r'\in\N^\t_{p^2}$: Since $r\in[2,p-1]$, $2<p$, we have $2p<p^2$, hence
$r'\in[2,2p-1]\sbs\zo2,p^2-1 $. \Mo since $p\nmid r$, we infer $p\nmid r+p$,
and a~fortiori $p^2\nmid r+p$. Therefore, by Lemma \rf{l2.16},
$\gcd(r',p^2)=1$. Thus $r'\in\N^\t_{p^2}$ by \E\Pr\ \rf{p3.44}.
As above we show that \te s $l'\in\Na$ \st $l'\cdot m'=p(p-1)$, and $m'\ddtt
{p^2}r'=1$. \E\Tf $p^2|(r')^{m'}-1$ by \er{4.58}, hence $p|(r')^{m'}-1$. We
claim that $p|r^{m'}-1$. Indeed, $(r')^{m'}=(r+p)^{m'}\nde4.54 = r^{m'}
+pm'r^{m'-1}+dp^2$ \fs $d\in\N$. Thus, since $(r')^{m'}\ge1$ and $r^{m'}\ge1$,
we have $(r')^{m'}-1=(r^{m'}-1)+pm'r^{m'-1}+dp^2$. Since $p|(r')^{m'}-1$, we
obtain $0\nad{\er{1.4},\er{1.6}} = {\F_p((r')^{m'}-1)} = \F_p\bigl((r^{m'}-1)
+p(m'r^{m'-1}+dp)\bigr)\nde1.4 = \F_p(r^{m'}-1)$. \If that $p|r^{m'}-1$, thus
$m'\ddtt p r=1$. As above we conclude from Lemma \rf{l1.22}\,(v) that \te s $q'\in\Na$ \st $m'=q'\cdot
(p-1)$. Hence $l'\cdot q'=p$. If $q'=1$, then $m'=p-1$ and $p^2|(r')^{p-1}-1$.
We now show that $p^2\nmid (r')^{p-1}-1$. To this end we invoke the \fw\ lemma.

\blm4.53
Let $p$ be an odd \Pn\ and let $r\in[2,p-1]$. Suppose $p|r^{p-1}-1$ and
$p^2| r^{p-1}-1$. Then $p|(r+p)^{p-1}-1$ and $p^2\nmid (r+p)^{p-1}-1$.
\elm

We postpone the proof of Lemma \rf{l4.53}, and apply Lemma \rf{l4.53}.

Since $p^2|r^{p-1}-1$ by \er{4.60}, we obtain $p^2\nmid (r')^{p-1}-1$
($r':=r+p$). \E\Tf $q'\ne1$ and $l'=1$, which implies that $r'$~is a \Gn\ of
$(\N^\t_{p^2},\cdot_{p^2},1)$. This concludes the proof that
$(\N^\t_{p^2},\cdot_{p^2},1)$ is cyclic.
\endproof

\proof[Proof of Lemma \rf{l4.53}]
Since $p^2|r^{p-1}-1$, \te s $q\in\Na$ \st $r^{p-1}-1=qp^2$. By \er{4.54}
$(r+p)^{p-1}=r^{p-1}+(p-1)r^{p-2}p+dp^2$ where $d:=\bigl(\frac{(p-1)(p-2)}2
r^{p-3}+c_{p-1}p\bigr)$. Note that $p\ge3$, hence $d\in\Na$ is well-defined.
Since $r\ge2$ and $p\ge3$, we have $(r+p)^{p-1},r^{p-1}\ge2$, hence $(r+p)
^{p-1}-1 = (r^{p-1}-1) + p((p-1)r^{p-2}+dp) = p(qp+dp+(p-1)r^{p-2})$. \If that
$p|(r+p)^{p-1}-1$. \Mo $\F_{p^2}((r+p)^{p-1}-1) = \F_{p^2}(p(p-1)r^{p-2}
+p^2(q+d)) \nde1.4 = \F_{p^2}(p(p-1)r^{p-2})$. Since $(r+p)^{p-1}-1 \ge2$ and
$p(p-1)r^{p-2}\ge2$, it follows from \er{1.46} that $p^2|(r+p)^{p-1}-1$ iff
$\F_2((r+p)^{p-1}-1)=0$ iff $\F_{p^2}(p(p-1)r^{p-2})=0$ iff $p^2|p(p-1)
r^{p-2}$ iff $p|(p-1)r^{p-2}$. It remains to show that $p\nmid (p-1)r^{p-2}$.
In view of Euclid's lemma \era3{9.18} $p\nmid(p-1)r^{p-2}$ iff $p\nmid (p-1)$
and $p\nmid r^{p-2}$. Since $p-1<p$, $p\nmid p-1$ by \era3{8.57}. It remains
to show that $p\nmid r^{p-2}$. Since $\gcd(r,p)=1$ by the \df\ of~$r$, we have
$p\nmid r$. Set $A:=\{l\ge2: p|r^l\}$. We claim that $A=\vn$. Suppose for \cd
ion that $A\ne\vn$. Since $(\N,\le)$ is well-ordered, $A$~possesses a least \el
\ which we denote by~$\ov l$. We have $\ov l\ge2$, and $p|r\cdot r^{\bar l-1}$.
Since $p\nmid r$, we obtain $p|r^{\bar l-1}$ by \era3{9.18}. A~\cd ion, since
$\ov l$~is the least \el\ of~$A$. \E\Tf $A=\vn$, hence $p\nmid r^{p-2}$,
$p\nmid (p-1)r^{p-2}$ and $p\nmid (r+p)^{p-1}-1$.
\endproof

In the proof of \E\Pr\ \rf{p4.51}, $k\ge3$, we shall use the \fw\ lemma.

\blm4.54
Let $a\in\Na\sms1$, $k\in\Na$ and $p$ be an odd prime. Suppose $p^k|a-1$ and
$p^{k+1} \nmid a-1$, then
\beq4.62
p^{k+1}|a^p-1 \qh{and }p^{k+2}\nmid a^p-1.
\e
\elm

\proof
From \E\Pr\ \rfa2{p7.9} \era2{7.31} with $E:=\N$, $n,p\in\N$ and $b:=1+p$ we
obtain
\beq4.63
b^{n+1}=1+(b-1)\suml_{l=0}^n b^l, \q n\in\N,\ b\in\Na.
\e
Let $a,k,p$ be as in Lemma \rf{l4.54}. In \er{4.63} set $b:=a$, $n:=p-1$,
which is allowed since $p\ge3$, and observe that $a^p>1$ by \era2{2.33} since
$a\in\N\sms{0,1}$, $p>0$ and $a^0\nda22.33 = 1$. \E\Tf \er{4.63} can be
rewritten as
\beq4.64
a^p-1 = (a-1)\suml_{l=0}^{p-1} a^l.
\e
By \as\ \te s $c\in\Na$ \st
\beq4.65
a-1 = cp^k \qh{and }p\nmid c.
\e
We now show that \te s $d\in\N$ \st
\beq4.66
\sum_{l=0}^{p-1} a^l=p+dp^2.
\e
We claim that \te\ $\zb el{[0,p-1]}\sbs\N$ \st
\beq4.67
a^l = 1+lcp^k + e_lp^{2k}, \q l\in[0,p-1].
\e
Indeed, for $l:=0$, $a^0\nda22.33 = 1= 1+0cp^k+0p^{2k}$. Hence \er{4.67}
holds with $e_0:=0$. For $l:=1$, $a^1\nda22.24 = a \nde4.65
=1+1cp^k+0p^{2k}$. Hence \er{4.67} holds with $e_1:=0$.

Let $l\in[2,p-1]$. Since $cp^k\in\Na$, we have $a^l=(1+cp^k)^l \nde4.54 =
1+lcp^k +\frac{l(l-1)}2 (cp^k)^2 +c_l(cp^k)^3
\nad{\era2{2.28},\era2{2.26},\era2{2.25},\era2{2.15}}=
1+lcp^k +\bigl[\frac{l(l-1)}2c^2 + c_lc^3p^k\bigr]p^{2k}$. Hence \er{4.67}
holds with $e_l:=\frac{l(l-1)}2c^2+c_lc^3p^k$ for $l\in[2,p-1]$.

\If from \er{4.67}, \era2{7.26} and \era2{7.27} that
\beq4.68
\sum_{l=0}^{p-1}a^l = \sum_{l=0}^{p-1}1 + \Bigl(\sum_{l=0}^{p-1}l\Bigr)cp^k
+\Bigl(\sum_{l=0}^{p-1}e_l\Bigr)p^{2k}.
\e
We have $\suml_{l=0}^{p-1}1 \nda27.28 = p$, and $\suml_{l=0}^{p-1}l =
0+\suml_{l=1}^{p-1}l= \frac{p(p-1)}2$. The last \et y follows from
\beq4.69
\suml_{l=1}^n l=\frac{n(n+1)}2\qh{\fa} n\in\Na.
\e
Formula \er{4.69} can be proved for example by \In\ on $n\in\Na$.
\E\Tf \er{4.68} becomes by using \era2{2.11}, \era2{2.12}, \era2{2.25}:
\beq4.69a
\suml_{l=0}^{p-1}a^l = p+\biggl(\frac{p-1}2\,c\biggr)p^{k+1} +
\Bigl(\suml_{l=0}^{p-1}e_l\Bigr)p^{2k}.
\e
If $k=1$, then \er{4.66} holds with $d:=\frac{p-1}2c+\suml_{l=0}^{p-1}e_l$.
If $k\ge2$, then \er{4.66} holds with
\[
d:=\frac{p-1}2\,cp^{k-1} +\Bigl(\sum_{l=0}^{p-1}e_l\Bigr)p^{2(k-1)},
\]
which proves the claim. Combining \er{4.64}--\er{4.66} we obtain $a^p-1 =
cp^{k+1}+dcp^{k+2} = c(1+dp)p^{k+1}$. Since $c,1+dp\in\Na$, we have $c(1+dp)\in\Na$,
hence $p^{k+1}|a^p-1$ and $\frac{a^p-1}{p^{k+1}}=c(1+dp)$. Suppose for \cd
ion that $p^{k+2}|a^p-1$. Set $q:=\frac{a^p-1}{p^{k+2}}\in\Na$. Then $p^{k+1}
c(1+dp) = p^{k+2}\cdot q=p^{k+1}pq$. From \E\Pr\ \rfa2{p2.7}\,(iii) we infer
$c(1+dp) =pq$, hence $p|c(1+dp)$. Since $p\nmid c$, we infer by Euclid's
lemma \era3{9.18} that $p|1+dp$. But $\F_p(1+dp)\nde1.4 = \F_p(1)\nde1.6 =
1$. Hence $p\nmid c(1+dp)=pq$, a~\cd ion. \E\Tf $p^{k+2}\nmid a^p-1$ and the
proof of the lemma is complete.
\endproof

\Wanp prove \E\Pr\ \rf{p4.51}.

\proof[Proof of \E\Pr\ \rf{p4.51}]
Case $k\ge2$.
Let $r\in[2,p-1]$ be a \Gn\ of the cyclic group $(\N_p^\t,\cdot_p,1)$. As in
the proof of the case $k=2$, $r$~\sf ies \er{4.59a}. Hence $p|r^{p-1}-1$ by
\er{4.58} since in view of \era2{2.23}, \era2{2.33} $r^{p-1}>1$. If $p^2\nmid
r^{p-1}-1$, we set $r':=r$. If not, we set $r':=r+p$. By Lemma \rf{l4.53},
\beq4.70
p|(r')^{p-1}-1 \qh{and } p^2\nmid (r')^{p-1}-1
\e
in both cases.

Since in both cases $r'\in [2,2p-1]$ and $2p-1\le p^k-1$, we infer $r'\in\N_{p^k}$. We claim
that $\gcd(r',p^k)=1$. If $r'=r$, then $\gcd(r',p)=1$ since $r$~is a unit of
the ring $(\N_p,+_p,\cdot_p,0,1)$ (see \E\Pr\ \rf{p3.44}). By Lemma \rf{l2.16}
the divisors of~$p^k$ are $\{p^i\}_{0\le i\le k}$. \E\Tf if $\gcd(r',p^k)=p^i$,
$i\in[1,k]$, then $p|r'$, which is impossible by \era3{8.57} since $r'=r<p$.
\If that $\gcd(r',p^k)=1$ if $r'=r$, hence in this case $r'$~is a unit of the
ring $(\N_{p^k},+_{p^k},\cdot_{p^k},0,1)$. If $r'=r+p$, then $\gcd(r',p^k)=1$.
Otherwise, as in the first case, $p|r+p$, hence $p|r$, a~\cd ion. Thus in both
cases $r'$ belongs to the group $(\N^\t_{p^k},\cdot_{p^k},1)$. Let $m\in
[2,\F(p^k)]$ denote the order of~$r'$ in the group $\N^\t_{p^k}$. By Lemma
\rf{l3.71} \er{3.73} $m$~divides $\#(N^\t_{p^k}) = \F(p^k)\nde3.39 = p^{k-1}
(p-1)$. Let $l\in\Na$ be \st
\beq4.71
l\cdot m = p^{k-1}(p-1).
\e
We want to show that $l=1$, that is, $\#(I(r'))=m=p^{k-1}(p-1) = \#(\N^\t
_{p^k})$, \ev tly, $\N^\t_{p^k}=I(r')$. We claim that $m$ is a \ml e of~$p-1$.
If $r'=r$, then $m\ddtt{p^k}r=1$ by \df\ of~$m$. Since $r\ge2$, we have $m\ge2$,
hence $r^m>1$ by \era2{2.33}. From \er{4.58} we infer $p^k|r^m-1$.
Since $p|p^k$, $p|r^m-1$, hence $m\ddtt pr=1$.
\If from Lemma \rf{l1.22}\,(v) as in the proof of the case $k=2$ that $m$~is
a \ml e of $p-1$, the order of~$r$ in~$\N^\t_p$. We now consider the case
$r':=r+p$. We have $m\ddtt{p^k}r'=1$, hence $p^k|(r')^m -1$ by \er{4.58}.
Thus $p|(r')^m-1$ and proceeding as in the proof of the case $k=2$, we find
$p|r^m-1$. Thus $m\ddtt p r=1$, hence, as above, $m$~is a \ml e of $p-1$,
which proves the claim. \E\Tf \te s $q\in\Na$ \st
\beq4.72
m=q\cdot(p-1).
\e
Combining \er{4.71} and \er{4.72}, we find $lq(p-1)=p^{k-1}(p-1)$. By
\era2{2.12} and \E\Pr\ \rfa2{p2.7}\,(iii) we obtain $lq=p^{k-1}$ since
$p-1\ne0$. \E\Tf $q|p^{k-1}$, hence by Lemma \rf{l2.16}, $q\in\{p^i\in\Na:
0\le i\le k-1\}$.  Thus by \er{4.72}
\beq4.73
m=p^i\cdot(p-1) \qh{\fs}i\in[0,k-1].
\e
We claim that \fa $n\in\N$:
\beq4.74
p^{n+1}|(r')^{p^n\cdot(p-1)}-1\qh{and } p^{n+2}\nmid(r')^{p^n\cdot(p-1)}-1.
\e
We proceed by \In\ on $n\in\N$. Set $M:=\{n\in\N: \hbox{\er{4.74} holds}\}$.
$0\in M$ follows from \er{4.70}. \ti{$n\in M$ implies $n+1\in M$}: assume
$n\in M$, that is, \er{4.74} holds for~$n$. We apply Lemma \rf{l4.54} with
$a:=(r')^{p^n\cdot(p-1)}$, $k:=n+1$, and we obtain
$p^{n+2}|\bigl((r')^{p^n\cdot(p-1)}\bigr)^p-1$ and $p^{n+3}\nmid
\bigl((r')^{p^n\cdot(p-1)}\bigr)^p-1$. By \era2{2.26}, \era2{2.11},
\era2{2.12}, \era2{2.24}, \era2{2.25}, we obtain
$p^{n+2}|(r')^{p^{n+1}\cdot(p-1)}-1 $ and
$p^{n+3}\nmid(r')^{p^{n+1}\cdot(p-1)}-1$. Hence $n+1\in M$, and $M=\N$, which
proves the claim.

Recall that $p^k|(r')^m-1$, that is by \er{4.73},
\beq4.75
p^k|(r')^{p^i\cdot(p-1)}-1, \qh{\fs}  i\in[0,k-1].
\e
Suppose for \cd ion that $0\le i<k-1$. Then $i+1<k$, hence $p^{i+1}|p^k$
and by the \tr ity of~$|$ we obtain $p^{i+1}|(r')^{p^i(p-1)}-1$. By \er{4.74}
$p^{i+2}\nmid(r')^{p^i\cdot(p-1)}-1$, hence a~fortiori
$p^{l}\nmid(r')^{p^i\cdot(p-1)}$ \fa $l\ge i+2$. Since $i+2<k+1$, we have
$(i+2)+1 \le k+1$ by \era1{3.17}, hence $i+2\le k$ by \era2{1.65}. \E\Tf
$p^k\nmid (r')^m-1$, a~\cd ion. \If that $i=k-1$, hence $m=p^{k-1}\cdot
(p-1)=\F(p^k)$. \csq, $\N^\t_{p^k}=I(r')$, hence $(\N^\t_{p^k},\cdot
_{p^k},1)$ is cyclic. \Mo $r$~or~$r+p$ is a \Gn\ of $\N^\t_{p^k}$.
\endproof

\bco4.55
Let $p$ be an odd \Pn\ and let $k\in\Na$ and let $a$~be a \Gn\ of the cyclic
group $(\N^\t_{p^k},\cdot_{p^k},1)$. Then the map $\ov\vf_a:
(\N_{p^{k-1}(p-1)},+_{p^{k-1}(p-1)},0) \to (\N^\t_{p^k},\cdot_{p^k},1)$
defined by $\ov\vf_a(l):=l\ddtt{p^k}a$, $l\in[0,p^{k-1}(p-1)-1]$ is an \is
sm.
\eco

\proof
Direct con\sq\ of \E\Pr\ \rf{p3.10}\,(iv) and $\#(\N^\t_{p^k})=\F(p^k)
=p^{k-1}(p-1)$.
\endproof

We now consider the \sc\ of the group of units of the rings $(\N_{2^k},
+_{2^k},\cdot_{2^k},0,1)$, $k\in\Na$.

\bpr4.56 \

\hph i,ii, $(\N^\t_2,\cdot_2,1)=(\{1\},\cdot_{2^k},1)$ the trivial group.

\hph ii,i, $(\N^\t_4,\cdot_4,1)=(\{1,3\},\cdot_4,1)$ and $\vf_{(4)}:
(\N_2,+_2,0)\to(\N^\t_4,\cdot_4,1)$ defined by
\beq4.76
\vf_{(4)}(0):=1,\ \vf_{(4)}(1):=3
\e
is an \is sm.

\hph iii,, Let $k\in\N$, $k\ge3$, and let
\[
\vf_{(2^k)}:(\N_2,+_2,0)\t(\N_{2^{k-2}},+_{2^{k-2}},0)\to (\N^\t_{2^k},
\cdot_{2^k},1)
\]
be defined by
\beq4.77
\bca
\vf_{(2^k)}(0,l):= l\ddtt{2^k}3,\\
\vf_{(2^k)}(1,l):=(2^k-1)\cdot_{2^k}\vf_{(2^k)}(0,l),\q l\in[0,2^{k-2}).
\eca
\e
Then $\vf_{(2^k)}$ is an \is sm.
\epr

\proof
We first observe that $2|2^k$ since $k\in\Na$ and $2^k=2\cdot2^{k-1}$. \Mo
\beq4.78
\N^\t_{2^k} = \{l\in[1,2^k-1]: l \hbox{ is odd (i.e.\ }\F_2(l)=1)\}.
\e
Indeed $\N^\t_{2^k}\nad{\er{3.45},\er{3.20}} = \{l\in[1,2^k]:
\gcd(l,2^k)=1\}$ and for $q\in\Na$ $q|2^k$ iff $q=2^i$ \fs $i\in[0,k]$ by
Lemma \rf{l2.16}. Therefore, for $l\in[1,2^k]$, $\gcd(l,2^k)=1$ iff $2\nmid
l$ iff $\F_2(l)=1$.

(i) If $k=1$, then $\N^\t_{2^k}=\{1\}$ by \er{4.78}.

(ii) If $k=2$, then $\N^\t_{2^k}=\{1,3\}$ by \er{4.78}.
Since $2\ddtt4 3=3\cdot_43 =\F_4(3^2)=1$, $\N^\t_{2^2}=I(3)$, hence
$(\N^\t_{2^2}, \cdot_4,1)$ is cyclic of order~$2$, hence \er{4.76} follows
from \E\Pr\ \rf{p3.10}\,(iv).

(iii) Let $k\ge3$. Set
\beq4.79
\ve_k:=2^k-1.
\e
Since $\ve_k\in[1,2^k -1]$ and $\ve_k$ is odd, $\ve_k\in\N^\t_{2^k}$.

In the first part of the proof we show:
\beq4.80
\#(I(3)) = 2^{k-2}, \q I(\ve_k)=\{1,\ve_k\}.
\e

In the second part of the proof, we show that the map $f:[0,1]\t [0,2^{k-2}-1]
\to\N^\t_{2^k}$ defined by
\beq4.81
f((\a,\b)):=(\a\ddtt{2^k}\ve_k)\cdot_{2^k}(\b\ddtt{2^k}3)
\e
is a bi\jn.

In the last part of the proof, we show that $f$~is a \hm sm from $(\N_2,+_2,0)\t
(\N_{2^{k-2}},+_{2^{k-2}},0)$ onto $(\N^\t_{2^k},\cdot_{2^k},1)$, and that
$f$~is equal to the map $\vf_{(2^k)}$ defined in \er{4.77}.

\er{4.80}: \E\fa $l\ge3$, $3^{(2^{l-2})}>1$ by \era2{2.23}, \era2{2.33}. We
claim \fa $l\ge3$:
\beq4.82
2^l| 3^{(2^{l-2})}-1.
\e
We proceed by \In\ on $l\ge3$. Set $M:=\{l\in\N,\ l\ge3: \hbox{\er{4.82} holds}\}$.
The \nm~3 belongs to~$M$, since $2^3=8$ and $3^{(2^{3-2})}-1=8$. Suppose $l\in M$.
Then $l\ge3$ and \er{4.82} holds. \E\Tf $l-2\in\Na$, $2^{l-2}\in\Na$, and
$3^{(2^{l-2})}\in\Na$ by \era2{2.30} with $E:=\N$. We have $\F_2(3^1)=\F_2(2+1)
\nde1.4 = \F_2(1)\nde1.6 = 1$. If $n\in\Na$ and $\F_2(3^n)=1$, then
$\F_2(3^{n+1})=\F_2(3^n\cdot3^1)\nde3.26 = \F_2\bigl(\F_2(3^n)\cdot\F_2(3^1)
\bigr)=\F_2(1\cdot1)=1$. Hence by \In\ $\F_2(3^n)=1$ \fa $n\in\Na$, \Ip
$\F_2\bigl(3^{(2^{l-2})}\bigr)=1$, hence \te s $q_l\in\N$ \st $3^{(2^{l-2})}
=2q_l+1$. Since $3^{(2^{l-2})}>1$, we have $q_l\in\Na$ by
\era2{2.13}. Hence $3^{(2^{l-2})}+1= (2q_l+1)+1=2q_l+2=2(q_l+1)$.
\If that $2|3^{(2^{l-2})}+1$. Since $2^l|3^{(2^{l-2})}-1$, we conclude that
$2^{l+1}|(3^{(2^{l-2})}-1)(3^{(2^{l-2})}+1)$. But \fa $a\in\Na$, we have $a^2
\ge a\ge1$ by \era2{2.19} and
\beq4.83
(a-1)(a+1)=a^2-1.
\e
Indeed, setting $b:=a-1$, we have $a=b+1$, hence $(a-1)(a+1)=b(b+2)= b^2+2b=
(b^2+2b+1)-1 =(b+1)(b+1)-1 =a^2-1$. \E\Tf $\bigl(3^{(2^{l-2})}-1\bigr)
\bigl(3^{(2^{l-2})}+1\bigr) \nde4.83 = \bigl(3^{(2^{l-2})}\bigr)^2-1 \nda22.26 =
3^{(2^{l-2}\cdot2)}-1 = 3^{(2^{l-1})}-1 = 3^{(2^{(l+1)-2})}-1$ (since
$l+1\ge2$). \E\Tf $2^{l+1}|3^{(2^{(l+1)-2})}$, hence $l+1\in M$, and
\er{4.82} holds \fa $l\ge3$.

By \er{4.82} with $l:=k$, and by \er{4.58}, since $3^{(2^{k-2})}>1$, we have
\beq4.84
2^{k-2} \ddtt{2^k}3=1.
\e
By Lemma \rf{l1.22} with $(X,\qu,e):=(\N^\t_{2^k},\cdot_{2^k},1)$ and $a:=3$,
\te s $m\ge2$ \st \er{1.43} holds. Then $m$~is the order of~$a$ by \E\df\
\rf{d2.11} and $m|2^{k-2}$ by Lemma \rf{l1.22}\,(v). By Lemma \rf{l2.16} and
$m\ge2$, we have
\beq4.85
m=2^i \qh{with }i\in[1,k-2].
\e
\E\Ip if $k=3$, then $m=2$. Indeed, $3\cdot_8 3=\F_8(3\cdot3)=1$. We now show
that
\beq4.86
m=2^{k-2}  \qh{for \ti{all}} k\ge3.
\e
We suppose $k>3$. We claim that \fa $n\in\Na$:
\beq4.87
2^{n+2}|3^{(2^n)}-1 \qh{and } 2^{n+3}\nmid 3^{(2^n)}-1.
\e
(Note that $3^{(2^n)}>1$ \fa $n\in\N$.)
We proceed by \In\ on $n\in\Na$.  Set $M:=\break\{n\in\Na: \hbox{\er{4.87}
holds}\}$.

$1\in M$, since $2^3=8|3^2-1$ and $16\nmid 3^2-1$. Suppose $n\in M$. Then
\er{4.87} holds. \Mo $3^{(2^{n+1})}-1 = (3^{(2^n)})^2 -1 \nde4.83 =
\bigl(3^{(2^n)}-1\bigr)\bigl(3^{(2^n)}+1\bigr)$. Proceeding as above (since
$3^{(2^n)}$ is odd and $3^{(2^n)}+1$ is even), we find
that $2|3^{(2^n)}+1$. Since $n\in M$, $2^{n+2}|3^{(2^n)}-1$, hence $2^{n+3}
|3^{(2^{n+1})}-1$. Suppose \ti{for \cd ion} that $2^{n+4}|3^{(2^{n+1})}-1$.
Then $2^{n+4}|{\bigl(3^{(2^n)}-1\bigr)}\bigl(3^{(2^n)}+1\bigr)$. As above
we find that $2|3^{(2^n)}+1$, hence, since $2^{n+2}|3^{(2^n)}-1$, we obtain $2|
\dfrac{3^{(2^n)}-1}{2^{n+2}}\cdot\dfrac{3^{(2^n)}+1}2$. Since $2^{n+3}\nmid
3^{(2^n)}-1$, $2\nmid\dfrac{3^{(2^n)}-1}{2^{n+2}}$, hence by Euclid's lemma
\era3{9.18}, $2|\dfrac{3^{(2^n)}+1}2$, that is, $4|3^{(2^n)}+1$. But
$3^{(2^n)} =3^{(2\cdot2^{n-1})}\nda22.26 = (3^2)^{(2^{n-1})}=
9^{(2^{n-1})} = (1+8)^{(2^{n-1})}$. If $n:=1$, then $3^{(2^n)}+1= 9+1=10$,
which is not divisible by~$4$. If $n>1$, then $3^{(2^n)}+1=(1+8)^{(2^{n-1})}
\nde4.54 = (1+8d)+1 = 2+8d$ \fs $d\in \N$ since $2^{n-1}\ge2$. But $4\nmid
2+8d$, hence $4\nmid 3^{(2^n)}+1$ whenever $n\in\Na$. A~\cd ion. \E\Tf
$2^{n+4}\nmid\bigl(3^{(2^{n-1})}-1\bigr)$, hence $n+1\in M$. \csq, $M=\Na$,
hence \er{4.87} holds \fa $n\in\Na$.

Suppose now for \cd ion $1\le i<k-2$ in \er{4.85}. Then $2^{i+3}\nmid
3^{(2^i)} -1$ by \er{4.87}. Hence, a~fortiori, $2^l\nmid 3^{(2^i)}-1$ \fa
$l\ge i+3$. Since $i+3<k-2+3=k+1$, we have $i+2<k$, hence $i+3\le k$ by
\era2{1.65}. Thus we may choose $l:=k$, and obtain $2^k\nmid 3^{(2^i)}-1$, \cd ing
$m\ddtt{2^k} 3=1$ with $m=2^i$, in view of \er{4.58} and $3^{(2^i)}>1$. Thus
\er{4.86} holds, and the first \et y in \er{4.80} holds.

We next show $I(\ve_k)=\{1,\e_k\}$.

We have
\beq4.88
2\ddtt{2^k}\ve_k=1 \qh{\fe}k\ge3.
\e
Indeed, note that $2^k>2>1$ by \era2{2.33}, \era2{2.23}. Hence $(\ve_k)^2-1 =
(2^k-1)^2-1 \nde4.83 = ((2^k-1)-1)((2^k-1)+1)=(2^k-2)\cdot2^k$. \E\Tf
$2^k|(\ve_k)^2-1$ and $2\ddtt{2^k}\ve_k=1$ by \er{4.58}. We have
$0\ddtt{2^k}\ve_k=1$ and $1\ddtt{2^k}\ve_k=\ve_k$ by \era2{2.3}\,(I0,I1). If
$n=2l$, $l\in\Na$, then $n\ddtt{2^k}\ve_k \nad{\era2{2.3}\,\rm I3}=
n\ddtt{2^k}(2\ddtt{2^k}\ve_k) \nde4.88 = n\ddtt{2^k}1 \nad{\era2{2.23}\,\rm
I4} =1$. If $n=2l+1$, $l\in\Na$, then $n\ddtt{2^k}\ve_k \nad{\era2{2.3}\,\rm
I2}=  (1\ddtt{2^k}\ve_k)\cdot_{2^k}(2l\ddtt{2^k}\ve_k)=
\ve_k\cdot_{2^k}1=\ve_k$. \If that $I(\ve_k)=\{1,\ve_k\}$. This completes the
proof of \er{4.80}.

\ssk
\er{4.81}: We first observe that $\#([0,1]\t[0,2^{k-2})) \nda26.28b =
\#([0,1])\cdot\#([0,2^{k-2}))\break = 2\cdot2^{k-2}=2^{k-1}$. \E\oh $\#(\N^\t
_{2^k})=\F(2^{k})=2^{k-1}$. In view of \era2{3.12} it suffices to show that
$f$~is sur\jc. Set $A:=\{y\in\N^\t_{2^k}: \hbox{\te s } l\in[0,2^{k-2})\break
\hbox{\st} f(0,l)=y\}$, and $B:=\{y\in\N^\t_{2^k}: \hbox{\te s }{l\in[0,2^{k-2})}\hbox{
such that}\break f(1,l)=y\}$. We shall show that $\#(A)=\#(B)=2^{k-2}$ and that $A\cap
B =\vn$. Hence by \era2{3.31}, $\#(A\cup B)=2\cdot2^{k-2}=2^{k-1}$. Since
$A\cup B\sbs \N^\t_{2^k}$ and $\#(A\cup B)=\#(\N^\t_{2^k})$ we obtain that
$A\cup B$, the range of~$f$, is equal to~$\N^\t_{2^k}$.

Since $f(0,l):=(0\ddtt{2^k}\ve_k)\cdot_{2^k}(l\ddtt{2^k}3) = 1\cdot_{2^k}
(l\ddtt{2^k}3) = l\ddtt{2^k}3$, $l\in\zo0,2^{k-2} $ and $2^{k-2}$ is the order of~$3$ in
the group $\{\N^\t_{2^k},\cdot_{2^k},1)\}$, $\zo0,2^{k-2} $ and~$A$ are \ep\ by
Lemma \rf{l1.22}, and $\#(A)=2^{k-2}$ by Corollary \rfa2{c3.11}. \Mo $f(1,l)
:= (1\ddtt{2^k}\ve_k)\cdot_{2^k}(l\ddtt{2^k}3) = \ve_k\cdot_{2^k}f(0,l)$.
Since the \mlc\ by~$\ve_k$ in the group $(\N^\t_{2^k},\cdot_{2^k},1)$ is
bi\jc\ by \E\Pr\ \rf{p3.51}, $A$~and~$B$ are \ep, hence
$\#(B)=\#(A)=2^{k-2}$. It remains to show that $A\cap B=\vn$.
Suppose for \cd ion that $A\cap B\ne\vn$. Then \te\ $l,l'\in\zo0,2^{k-2} $ \st
$$
l\ddtt{2^k}3 = \ve_k\cdot_{2^k}(l'\ddtt{2^k}3).
$$
If $l=l'$, we obtain $1=\ve_k$ by \cnc ity, which is
impossible. If $l>l'$, and $l=l''+l'$ with $l''\in(0,l-l')$, then by \era2{2.3}\,I2 we
obtain $(l''\ddtt{2^k}3)\cdot_{2^k}(l'\ddtt{2^k}3) = \ve_k\cdot_{2^k}
(l'\ddtt{2^k}3)$. By \cnc ity in
$(\N^\t_{2^k}, \cdot_{2^k},1)$ we obtain $l''\ddtt{2^k}3=\ve_k$,
which implies $\F_{2^k}(3^{l''})
= \F_{2^k}(\ve_k)$ by \er{3.71} since $\e_k=2^k-1\in\zo0,2^k $.

Note that $\F_{2^k}(\ve_k)=2^k-1$ by \er{1.5}. \Mo $3^{l''}=q2^k+
\F_{2^k}(3^{l''})=q2^k+
\F_{2^k}(\ve_k) = q2^k+2^k-1$ \fs $q\in\N$. \If that $3^{l''}+1 = (q+1)2^k =
(q+1)2^{k-3}\,8$ since $k>3$. \E\Tf $8|3^{l''}+1$. We claim that $\F_8(3^m)
\in\{1,3\}$ \fa $m\in\N$. Indeed, if $m:=0$, then $\F_8(3^m)=\F_8(1)\nde1.5 =
1$. If $m:=1$, then $\F_8(3^m)=\F_8(3)\nde1.5 = 3$. If $m:=2$, then $\F_8(3^m)
=\F_8(1+8)\nad{\er{1.4},\er{1.5}}= 1$. If $m:=2n$, $n\in\Na$, then
$\F_8(3^{2n})\nde3.71 = 2n\ddtt8 3 = (n\dpl2)\ddtt8 3 \nad{\era2{2.3}\,\rm I3}
= n\ddtt8(2\ddtt8 3)= n\ddtt8 1\nad{\era2{2.3}\,\rm I4}= 1$. If $m:=2n+1$,
$n\in\Na$, then $\F_8(3^{2n+1})= (2n+1)\ddtt8 3 \nad{\era2{2.3}\,\rm I2}=
(2n\ddtt8 3)\cdot_8(1\ddtt8 3)\nad{\era2{2.3}\,\rm I1} = 1\cdot_8 3=3$, which
proves the claim. \If that $3^{l''}= q8+r$ \fs $q\in\N$ and $r\in\{1,3\}$. Thus
$8\nmid 3^{l''}+1$. A~\cd ion. Hence $l\not>
l'$. We now consider the case $l<l'$, that is, $l'=l''+l$ with $l''\in(0,l-l')$.
Thus $l\ddtt{2^k}3 = \ve_k\cdot(l+l'')\ddtt{2^k}3$. By \er{4.88} $\ve_k
\ddtt{2^k}\ve_k=1$, hence $\ve_k\cdot_{2^k}(l\ddtt{2^k}3)\nda22.3 = (l\ddtt{2^k}
3)\cdot_{2^k}(l''\ddtt{2^k}3)$. By \cmt ity and \cnc ity in $(\N^\t_{2^k},
\cdot_{2^k},1)$ we obtain $\ve_k=l''\ddtt{2^k}3$ in $(\N^\t_{2^k},\cdot_{2^k},
1)$, hence $\F_{2^k}(\ve_k)=\F_{2^k}(l''\ddtt{2^k}3)$. But we just showed
above that this is impossible. \E\Tf $A\cap B=\vn$ and the proof of \er{4.81}
is complete.

\ti{$\vf_{(2^k)}$ is an \is sm}: We first show that the maps $\vf_{(2^k)}$
and~$f$ are identical. Indeed, $f(0,l)\nde4.81 = (0\ddtt{2^k}\ve_k)\cdot_{2^k}
(l\ddtt{2^k}3)\nda22.3 = 1\cdot_{2^k}(l\ddtt{2^k}3) = l\ddtt{2^k}3\nde4.77 =
\vf_{(2^k)}(0,l)$ and $f(1,l)\nde4.81 = (1\ddtt{2^k}\ve_k)\cdot_{2^k}(l
\ddtt{2^k}3)\nda22.3 = \ve_k\cdot_{2^k}\vf_{(2^k)}(0,l)\nad{\er{4.77},
\er{4.79}}= \vf_{(2^k)}(1,l)$ \fa $l\in\zo0,2^{k-2} $.

We now show that $\vf_{(2^k)}$ is a \hm sm. We have $\vf_{(2^k)}(0,0)\nda22.3
= {1\cdot_{2^k}1=1}$. Let $(x',y'),(x'',y'')\in(\N_2,+_2,0)\t (\N_{2^{k-2}},
+_{2^{k-2}},0)$. Then $\vf_{(2^k)}((x',y')\mathrel{\wt+}(x'',y''))\nad*=\break
\vf_{(2^k)}(x'+_2x'',y'+_{2^{k-2}}y'')\nad{**}= ((x'\ddtt2 \ve_k)\cdot_{2^k}
(x''\ddtt2\ve_k))\cdot_{2^k}((y'\ddtt{2^{k-2}}3)\cdot_{2^k}(y''\ddtt{2^{k-2}}3))
\nda21.36 =\break  (x'\ddtt2 \ve_k)\cdot_{2^k}((x''\ddtt2 \ve_k)\cdot_{2^k}(y'
\ddtt{2^{k-2}}3))\cdot(y''\ddtt{2^{k-2}}3) \nda21.6 = (x'\ddtt2 \ve_k)
\cdot_{2^k}((y'\ddtt{2^{k-2}}3)\cdot_{2^k}({x''\ddtt2 \ve_k}))\cdot(y''
\ddtt{2^{k-2}}3) \nda21.36 =((x'\ddtt2 \ve_k)\cdot_{2^k}(y'\ddtt{2^{k-2}}
3))\cdot_{2^k}((x''\ddtt2 \ve_k)\cdot(y''\ddtt{2^{k-2}}3)) \nad{*{*}*}=
\vf_{(2^k)}(x',y')\cdot_{2^k}\vf_{(2^k)}(x'',y'')$. In $\nad*=$ we denoted
by~$\wt+$ the binary \op\ of the direct product of $(\N_2,+_2,0)$ and
$(\N_{2^{k-2}},+_{2^{k-2}},0)$ (see Example \rfa2{xa1.5}\,(iv)). In
$\nad{**}=$ we used Lemma \rf{l1.22}\,(vi). In $\nad*=$ and $\nad{*{*}*}=$ we
used the fact that $\vf_{(2^k)}$ and~$f$ are identical. \If that $\vf_{(2^k)}$
is a \hm sm. Since $\vf_{(2^k)}=f$ and $f$~is a bi\jn, $\vf_{(2^k)}$ is an
\is sm by Lemma \rfa2{l1.8}. This completes the proof of \E\Pr\ \rf{p4.56}.
\endproof

We are ``almost'' in a position to prove one of the main results of this
section. We shall need the \fw\ lemma.

\blm4.57
Let $(\N_k,+_k,0)$, $(\N_l,+_l,0)$, $k,l\ge2$, be the groups introduced in
\E\Pr\ \rf{p1.8}. Then the finite group $(\N_k,+_k,0)\t(\N_l,+_l,0)$ is
cyclic iff $k\ne l$ and $\gcd(k,l)=1$.
\elm

The proofs of \E\Pr\ \rf{p4.56}, Lemma \rf{l4.57} and Theorem \rf{t4.58} are
based on the proofs given in~\cite{Rham}.

\proof \

\ti{If\/}: follows from Corollary \rf{c2.15}.

\ti{Only if\/}: We show that if either $k=l$ or $\gcd(k,l)>1$, then
$\N_k\t\N_l$ is not cyclic. Suppose $k=l$ and suppose for \cd ion that
$\N_k\t\N_k$ is cyclic hence \pn. The order of $\N_k\t\N_k$ is $k^2$ by
\era2{6.28b}. Then $\N_k\t\{0\}$ and $\{0\}\t\N_k$ are distinct subgroups of
order~$k$. A~\cd ion in view of \E\Pr\ \rf{p3.24}\,(ii).

Suppose $k\ne l$, $d:=\gcd(k,l)>1$, and suppose for \cd ion that $\N_k\t\N_l$
is cyclic. By \E\Pr\ \rf{p3.24}\,(ii) \te s a \sbm\ $M$ of $\N_k$ of order~$d$,
hence $M\t\{0\}$ is a \sbm\ of $\N_k\t\N_l$ of
order~$d$, hence $M\t\{0\}$ is  a \sbm\ of $\N_k\t\N_l$ of order~$d$.
Similarly \te s a \sbm~$N$ of $\N_l$ of order~$d$. Hence both $M\t\{0\}$
and $\{0\}\t N$ are \sbm s of $\N_k\t\N_l$ of order~$d$. A~\cd ion as above
since $\N_k\t\N_l$ is cyclic.
\endproof

\bth4.58
Let $n\in\N$, $n\ge3$. Then $(\N^\t_n,\cdot_n,1)$ is cyclic iff $n=4$, $p^k$
or $2p^k$ \fs odd \Pn~$p$ and some $k\in\Na$.
\eth

\proof
We shall use the symbol $\cong$ for ``(monoid)-\is c to''. We recall that
\fa $m\ge2$ $(\N_m,+_m,0)$ is cyclic by \E\Pr s \rf{p1.8} and
\rf{p3.10}\,(ii).

\ti{If\/}: $n=4$. By \E\Pr\ \er{p4.56}\,(ii) $(N^\t_4, \cdot_4,1)\cong
(\N_2,+_2,0)$, hence $(\N^\t_4, \cdot_4,1)$ is cyclic by Lemma
\rf{l3.25}\,(ii).

$n=p^k$. By Corollary \rf{c4.55}, $\N^\t_{p^k} \cong(\N_{p^{k-1}(p-1)},
+_{p^{k-1}(p-1)},0)$, hence $(\N^\t_{p^k},\cdot_{p^k},1)$ is cyclic by Lemma
\rf{l3.25}\,(ii).

$n=2p^k$. Since $2<p^k$ and $\gcd(2,p^k)=1$ by Lemma \rf{l2.16}, the ring
$(\N_{2p^k},+_{2p^k},\cdot_{2p^k},0,1)$ is ring-\is c to the ring
$(N_2,+_2,\cdot_2,0,1)\t (\N_{p^k},+_{p^k},\cdot_{p^k},0,1)$ by Theorem
\rf{t2.24}. \E\Ip the monoid $(\N_{2p^k},\cdot_{2p^k},1)$ is \is c to the
monoid $(\N_2,\cdot_2,1)\t (\N_{p^k},\cdot_{p^k},1)$. By Lemma \rf{l3.36}
$(\N^\t_{2p^k},\cdot_{2p^k},1) \cong(\N^\t_2,\cdot_2,1)\t(\N^\t_{p^k},
\cdot_{p^k},1)$. But $(\N^\t_2,\cdot_2,1)=(\{1\},\cdot_2,1)$. Hence
$(\N^\t_{2p^k},\cdot_{2p^k},1)\cong (\N^\t_{p^k},\cdot_{p^k},1)\cong
(\N_{p^{k-1}(p-1)},\cdot_{p^{k-1}(p-1)},1)$. Hence by Lemma \rfa2{l1.8}\,(ii)
and Lemma \rf{l3.25}\,(ii), $(\N^\t_{2p^k},\cdot_{2p^k},1)$ is cyclic.

\def\labelenumi{{\rm(\arabic{enumi})}}
\ti{Only if\/}: In view of the ``if'' part it suffices to show that if $n\ge3$
and $n\ne4$, $n\ne p^k$, $n\ne 2p^k$, $p$~odd prime, $k\in\N^\t$, then
$(\N^\t_n,\cdot_n,1)$ is \ti{not\/} cyclic. According to Theorem \rfa3{t9.22},
we need only to consider four cases:
\ben
\item $n:=2^l$, $l\ge3$;
\item $n:=m$ with $m:=\prodl_{i=1}^N p_i^{m_i}$, $N\ge2$, $p_i$ odd prime,
$i\in[1,N]$, $p_i\ne p_j$ if $i\ne j$, and $m_i\in\Na$;
\item $n:=2^lm$, $l\ge1$, $N\ge2$;
\item $n:=2^lm$, $l\ge2$, $N=1$.
\een

\ti{Case} (1): By \E\Pr\ \rf{p4.56}, $(\N^\t_{2^l},\cdot_{2^l},1)\cong
(\N_2,+_2,0)\t (\N_{2^{l-2}},+_{2^{l-2}},0)$. If $l=3$, then $2=2^{l-2}$
and if $l>3$, then $2|2^{l-2}$ and $\gcd(2,2^{l-2})=2$. Therefore, by Lemma
\rf{l4.57} and Lemma \rf{l3.25}\,(ii), $(\N_{2^{l}},+_{2^{l}},1)$ is \ti{not\/}
cyclic.

\ti{Case} (2): By Lemma \rf{l2.16}, $\gcd(p_i^{m_i},p_j^{m_j})=1$ \fa $i,j\in
[1,N]$, $i\ne j$. \E\Tf by Theorem \rf{t2.24} and Lemma \rf{l3.36} we obtain
\beq4.89
(\N^\t_m,\cdot_m,1) \cong \prod_{i=1}^N (\N^\t_{p_i^{m_i}},\cdot_{p_i^{m_i}},1).
\e
We claim
\beq4.90
\prod_{i=1}^N (\N^\t_{p_i^{m_i}},\cdot_{p_i^{m_i}},1)\quad (N\ge2)\quad
\hbox{is \ti{not\/} cyclic},
\e
hence by Lemma \rf{l3.25}\,(ii), $(\N^\t_m,\cdot_m,1)$ is not cyclic.

Indeed, by Corollary \rf{c4.55} and Lemma \rfa2{l1.8} \te\ \is sms
$j_i:(\N^\t_{p_i^{m_i}},\cdot_{p_i^{m_i}},1)\break\to (\N_{p_i^{m_i-1}(p_i-1)},
\cdot_{p_i^{m_i-1}(p_i-1)},0)$, $i=1,2$. One verifies that if $(X_i,\qu_i,e_i)$,
$(X_i',\qu_i',e_i')$, $i=1,2$, are monoids and $j_i:X_i\to X_i'$, $i=1,2$, are
\is sms, then $j:X_1\t X_2 \to X_1'\t X_2'$ defined by $j(x_1,x_2)=(j_1(x_1),
j_2(x_2))$ is an \is sm. Since $2|p_i^{m_i-1}(p_i-1)$, $i=1,2$, then either
$p_1^{m_1-1}(p_1-1)$ and $p_2^{m_2-1}(p_2-1)$ are equal or \sf y $\gcd(
p_1^{m_1-1}(p_1-1),\break p_2^{m_2-1}({p_2-1}))\ge2$. \E\Tf by Lemma \rf{l4.57},
$\prodl_{i=1}^2 (\N_{p_i^{m_i-1}(p_i-1)},+_{p_i^{m_i-1}(p_i-1)},0)$ is
\ti{not\/} cyclic and by what precedes and Lemma \rf{l3.25}\,(ii),
$\prodl_{i=1}^2 (\N^\t_{p_i^{m_i}},\cdot_{p_i^{m_i}},1)$ is also \ti{not\/}
cyclic. If $N=2$, the proof of case~(2) is complete. If $N>2$, then the subset
consisting of $x\in\prodl_{i=1}^N (\N^\t_{p_i^{m_i}},\cdot_{p_i^{m_i}},1)$
\sf ying $\pi_i(x)=1$, $i\in[3,N]$ (see \er{2.17}), is a nontrivial \sbm\ of
$\prodl_{i=1}^N (N^\t_{p_i^{m_i}},\cdot_{p_i^{m_i}},1)$, which is \is c to
$\prodl_{i=1}^2 (\N_{p_i^{m_i}},\cdot_{p_i^{m_i}},1)$, hence \ti{not\/}
cyclic by case~(2) and Lemma \rf{l3.25}\,(ii). \csq, $\prodl_{i=1}^N (N^\t_{p_i^{m_i}},
\cdot_{p_i^{m_i}},1)$ cannot be cyclic by \E\Pr\ \rf{p3.20} and \E\Pr\
\rf{p3.10}\,(ii). The proof of case (2) is complete.

\ti{Case} (3): By Lemma \rf{l2.16}, $2^l\ne m$ and $\gcd(2^l,m)=1$. By Theorem
\rf{t2.24} and Lemma \rf{l3.36},
\[
(N^\t_{2^lm},\cdot_{2^lm},1) \cong (\N^\t_{2^l},\cdot_{2^l},1)\t(\N^\t_m,
\cdot_m,1).
\]
But $(\{1\},\cdot_{2^l},1)\t(\N^\t_m,\cdot_m,1)$ is a nontrivial \sbm\ of
$(\N^\t_{2^l},\cdot_{2^l},1)\t(\N^\t_m, \cdot_m,1)$, which is \is c to
$(\N^\t_m,\cdot_m,1)$, hence \ti{not\/} cyclic by case~(2) and Lemma \rf{l3.25}\,(ii).
As in case~(2), we conclude that $(\N^\t_{2^l},\cdot_{2^l},1)\t(\N^\t_m,
\cdot_m,1)$ is not cyclic by invoking \E\Pr\ \rf{p3.20}. \E\Tf $(\N^\t_
{2^lm},\cdot_{2^lm},1)$ is \ti{not\/} cyclic by Lemma \rf{l3.25}\,(ii).

\ti{Case} (4): Set $p:=p_1$ and $k:=m_1$. By Lemma \rf{l2.16} $2^l\ne p^k$
and $\gcd(2^l,p^k)=1$. As above
\[
(\N^t_{2^lp^k},\cdot_{2^lp^k},1) \cong (\N^\t_{2^l},\cdot_{2^l},1)
\t(\N^\t_{p^k},\cdot_{p^k},1).
\]
If $l\ge3$ then $(\N^\t_{2^l},\cdot_{2^l},1)$ is \ti{not\/} cyclic by case~(1).
Hence $(\N^\t_{2^l},\cdot_{2^l},1)\t(\{1\},\cdot_{p^k},1)$ is a nontrivial
\sbm\ of $(\N^\t_{2^l},\cdot_{2^l},1)\t(\N^\t_{p^k},\cdot_{p^k},1)$, which is
\ti{not\/} cyclic. By \E\Pr\ \rf{p3.20} and Lemma \rf{l3.25}\,(ii), we
conclude as above that $(\N^\t_{2^lp^k},\cdot_{2^lp^k},1)$ is not cyclic.

It remains to consider the case $l=2$, that is, $n:=4p^k$. By \E\Pr\
\rf{p4.56}\,(ii), $(\N^\t_4,\cdot_4,1)\cong(\N_2,+_2,0)$ and by Corollary
\rf{c4.55}, $(\N^\t_{p^k},\cdot_{p^k},1)\cong(\N_{p^{k-1}(p-1)},\cdot
_{p^{k-1}(p-1)},1)$. Therefore, as in the proof of case~(2), we infer
$(\N^\t_4,\cdot_4,1)\t(\N^\t_{p^k},\cdot_{p^k},1) \cong (\N_2,+_2,0)\t
(\N_{p^{k-1}(p-1)},+_{p^{k-1}(p-1)},0)$. Since either $2$ and $p^{k-1}(p-1)$
are equal or $\gcd(2,p^{k-1}(p-1))=2$, we infer by Lemma \rf{4.57} that
$(\N^\t_4,\cdot_4,1)\t(\N^\t_{p^k},\cdot_{p^k},1)$ is \ti{not\/} cyclic. Hence
by Lemma \rf{l3.25}\,(ii) $(\N^\t_{2^lp^k},\cdot_{2^lp^k},1)$ is not cyclic.
This completes the proof of case~(4), hence also of Theorem \rf{t4.58}.
\endproof

In the proof of \E\Pr\ \rf{p4.51} we used formula \er{4.54}. Observe that if
$n=2$, then $(a+b)^2=(a+b)(a+b) \nda22.15 = a(a+b)+b(a+b)\nda22.14 = (a^2+ab)
+(ba+b^2) \nda21.36 = a^2+(ab+ba)+b^2 \nda22.12 = a^2+(ab+ab)+b^2 \nda22.8 =
a^2+2ab+b^2$. \E\Tf $c_2=0$ in formula \er{4.54} if $n=2$.

In the remaining part of this section we consider the so-called ``\ti{binomial
formula\/}'',\index{binomial formula} which expresses a power of a sum of terms as a sum of
products of factors. \E\Ip it allows us to give an explicit expression for
the factor~$c_n$ in formula \er{4.54} when $n\ge3$. We first consider a somewhat
more general situation where a (composite) product of a sum of two terms is
\rp ed as a \cme sum of products of two factors. In what follows
$(X,+,\cdot,0,1)$ denotes a (\cmt e) \sr\ (with unity). Let $a,b,c,d\in X$.
Then
\beq4.91
(a+b)(c+d) = ac+ad+bc+bd.
\e
Indeed, $(a+b)(c+d) \nde1.3n = a(c+d)+b(c+d) \nde1.2n = (ac+ad)+(bc+bd)$.

Now let $a:[1,2]\t [1,n]\to X$.

We want to express $\prodl_{j=1}^n\bigl(\suml_{i=1}^2a_{ij}\bigr)$ as a \cme
sum of \cme products of $n$~factors. If $n=3$, we obtain from \er{4.91}:
$\prodl_{j=1}^3 \bigl(\suml_{i=1}^2 a_{ij}\bigr) = (a_{11}+a_{21})
(a_{12}+a_{22})(a_{13}+a_{23})=\break(a_{11}+a_{21})(a_{12}a_{13} + a_{12}a_{23} +
a_{22}a_{13} + a_{22}a_{23}) = a_{11}a_{12}a_{13} +a_{11}a_{12}a_{23} +
a_{11}a_{22}a_{13} + a_{11}a_{22}a_{23} + a_{21}a_{12}a_{13}+a_{21}a_{12}a_{23}
+a_{21}a_{22}a_{13}+ a_{21}a_{22}a_{23}$.

We observe that if we write
$a_{ij}$, $i\in[1,2]$, $j\in[1,3]$, as
\beq4.92
\bmn a_{11}\\a_{21}\emn ,  \q \bmn a_{12}\\a_{22}\emn, \q \bmn a_{13}\\a_{23}\emn,
\e
each term is a product of the form $a_{r1}a_{s2}a_{t3}$, $r,s,t\in[1,2]$,
which is a product of \el s of each column. The index~$i$ indicates which \el\
of the column is chosen. Note that all possible choices of ``\sq s''
$(r,s,t) \in [1,2]\t[1,2]\t[1,2]$ occur. \E\Tf we can rewrite $\prodl_{j=1}
^3 \bigl(\suml_{i=1}^2 a_{ij}\bigr)$ as
\[
\sum_{g\in[1,2]^{[1,3]}}\Bigl(\prod_{j=1}^3 a_{g(j)j}\Bigr),
\]
where $[1,2]^{[1,3]}$ denotes the set of all maps from $[1,3]$ into $[1,2]$.
More generally, we have

\bpr4.59
Let $(X,+,\cdot,0,1)$ be a $($\cmt e$)$ \sr\ $($with unity$)$, in particular
$(\N,+,\cdot,0,1)$. Let $n\ge2$ and let $a:[1,2]\t[1,n]\to X$. Then
\beq4.93
\prod_{j=1}^n \Bigl(\sum_{g\in[1,2]^{[1,n]}}a_{ij}\Bigr)=
\sum_{g\in[1,2]^{[1,n]}}\Bigl(\prod_{j=1}^n a_{g(j)j}\Bigr).
\e
\epr

\proof
By \In\ on $n\in \Na\sms1$. Set $M:=\{n\ge2: \hbox{\er{4.93} holds}\}$. Then
$2\in M$ by \er{4.91}.

{\allowdisplaybreaks
\ti{$n\in M$ implies $n+1\in M$}: Let $a:[1,2]\t[1,n+1]$, and let
$n\in M$. Then
\bgg
\sum_{g\in[1,2]^{[1,n+1]}}\Bigl(\prod_{j=1}^{n+1} a_{g(j)\,j}\Bigr)
\nda27.13 = \sum_{g\in[1,2]^{[1,n+1]}}\Bigl(\Bigl(\prod_{i=1}^n
a_{g(j)\,j}\Bigr)\cdot a_{g(n+1)\,n+1}\Bigr)\\
\nda21.138 =\sum_{\sbk{g\in[1,2]^{[1,n]}\\g(n+1)=1}}\Bigl(\Bigl(\prod_{i=1}^n
a_{g(j)\,j}\Bigr)\cdot a_{g(n+1)\,n+1}\Bigr)+
\sum_{\sbk{g\in[1,2]^{[1,n]}\\g(n+1)=2}}\Bigl(\Bigl(\prod_{i=1}^n
a_{g(j)\,j}\Bigr)\cdot a_{g(n+1)\,n+1}\Bigr)\\
{}= \sum_{g\in[1,2]^{[1,n]}}\Bigl(\Bigl(\prod_{i=1}^n
a_{g(j)\,j}\Bigr)\cdot a_{1\,n+1}\Bigr)+
\sum_{g\in[1,2]^{[1,n]}}\Bigl(\Bigl(\prod_{i=1}^n
a_{g(j)\,j}\Bigr)\cdot a_{2\,n+1}\Bigr)\\
\nda21.6 = \sum_{g\in[1,2]^{[1,n]}}\Bigl(a_{1\,n+1}\cdot\Bigl(\prod_{i=1}^n
a_{g(j)\,j}\Bigr) \Bigr)+
\sum_{g\in[1,2]^{[1,n]}}\Bigl(a_{2\,n+1}\cdot\Bigl(\prod_{i=1}^n
a_{g(j)\,j}\Bigr) \Bigr) \\
\nad*= a_{1\,n+1}\cdot\sum_{g\in[1,2]^{[1,n]}}\Bigl(\prod_{i=1}^n
a_{g(j)\,j}\Bigr) +
a_{2\,n+1}\cdot\sum_{g\in[1,2]^{[1,n]}}\Bigl(\prod_{i=1}^n
a_{g(j)\,j}\Bigr)  \\
\nad{\er{1.3n},n\in M} = (a_{1\,n+1}+a_{2\,n+1})\cdot\prod_{j=1}^n \Bigl(\sum_{i=1}^2
a_{ij}\Bigr) \\
\nda21.6 = \Bigl(\prod_{j=1}^n \Bigl(\sum_{i=1}^2 a_{ij}\Bigr)
\Bigr)\cdot\Bigl(\sum_{i=1}^2 a_{i\,n+1}\Bigr) \nda27.13 =
\prod_{j=1}^{n+1}\Bigl(\sum_{i=1}^2 a_{ij}\Bigr).
\e}%
In $\nad*=$ we used \era2{1.131} with $(X,\qu,e):=(X,+,0)$, $\wt X:=X$,
$\vf(x):=a_{1\,n+1}\cdot x$ (resp.\ $a_{2\,n+1}\cdot x$), $I:=[1,2]^{[1,n]}$,
$i:=g$, $a_i:= \prodl_{j=1}^n a_{g(j)\,j}$, and where \era2{1.130} follows
from~\er{1.2n}.
\endproof

We now consider the special case $n=3$ where all columns in \er{4.92} are
identical and, more generally, the case where $a_{ij}=a_{i1}$ \fa $i\in[1,2]$
and all $j\in[1,n]$. Thus we can omit the second index in $a_{ij}$ and we
simply set $a_i:=a_{ij}$, $(i,j)\in[1,2]\t[1,n]$. \Mo we replace the
index~$i$ by~$l$ setting $l=1$ whenever $i=1$ and $l=0$ whenever $i=2$.
Finally, we set $a:=a_1$ and $b:=a_0$. Then, using \era2{1.126} where
$(X,\qu,e):=(X,\cdot,1)$, $\prodl_{j=1}^n \bigl(\suml_{i=1}^2 a_{ij}\bigr)$
becomes $n\ddt (a+b)$, the $n$-th \IT\ of $a+b$ in $(X,\cdot,1)$, since
$\#([1,n])=n$.

\bnt4.60
Let $(X,+,\cdot,0,1)$ be a (\cmt e) \sr\ (with unity), and let $x\in X$,
$n\in\N$. We use the \fw\ notation (which coincides with the standard notation
in the \sr\ $(\N,+,\cdot,0,1)$):
\bea4.94
nx&:=n\dpl x,\\
x^n&:=n\ddt x. \lb{4.95}
\e
\ent

Accordingly, formula \er{4.93} becomes\index{multinomial formula}
\beq4.96
(a+b)^n = \sum_{g\in[0,1]^{[1,n]}}\Bigl(\prod_{j=1}^n c_{g(j)}\Bigr),
\e
where $c_0:=b$, $c_1:=a$.

Recall that if $E$ is a \ns, then there is a \ti{bi\jn} between $\cP(E)$, the
set of all subsets of~$E$, and $[0,1]^E$, the set of all maps from~$E$ into
$[0,1]$. Indeed, if $A\in \cP(E)$, then $1_A:E\to[0,1]$ defined by $1_A(x):=1$
if $x\in A$ and $1_A(x):=0$ if $x\notin A$. The \f\ $1_A$ is called the
\ti{indicator \f} of~$A$ (see Exercise \rf{ex1.14n}) or the \ti{\ch istic \f}
of~$A$. \E\oh given $f:E\to[0,1]$, the subset of~$E$ defined by $\{x\in E:
f(x)=1\}$ is called the \ti{support\/} of~$f$ and is sometimes denoted by
$\supp f$.

We have
\bea4.97
\supp(1_A)=A \q&\hbox{\fa}A\in \cP(E),\\
1_{\supp f}=f \q &\hbox{\fa}f\in[0,1]^E. \lb{4.98}
\e
In view of \er{4.97}, \er{4.98} and the \df\ of~$c_i$, we obtain
\beq4.99
\sum_{g\in[0,1]^{[1,n]}} \Bigl(\prod_{j=1}^n c_{g(j)}\Bigr) =
\sum_{A\in \cP([1,n])}\Bigl(\prod_{j=1}^n d_j(A)\Bigr)
\e
where
\beq4.100
d_j(A):=a \hbox{ if }j\in A, \q\ d_j(A):=b \hbox{ if }j\notin A.
\e
Given $A\sbs[1,n]$, we have
\bmlg
\prod_{j=1}^n d_j(A) \nda27.17 = \prod_{j\in[1,n]}d_j(A) \nda21.138 =
\Bigl(\prod_{\sbk{j\in[1,n]\\j\in A}}d_j(A)\Bigr)
\Bigl(\prod_{\sbk{j\in[1,n]\\j\notin A}}d_j(A)\Bigr)\\
{}= \Bigl(\prod_{\sbk{j\in[1,n]\\j\in A}}a\Bigr)
\Bigl(\prod_{\sbk{j\in[1,n]\\j\notin A}}b\Bigr)\nad{\er{4.54},\er{4.95}}=
a^{\#(A)}\cdot b^{\#(A^c)} \nda23.31 = a^{\#(A)}\cdot b^{n-\#(A)}.
\e
\csq, we obtain
\beq4.101
(a+b)^n = \sum_{A\in\cP([1,n])} a^{\#(A)}b^{n-\#(A)}.
\e
Since $a+b=b+a$, we obtain $\suml_{A\in\cP([1,n])} a^{\#(A)}b^{n-\#(A)}=
(a+b)^n=(b+a)^n=\break\suml_{A\in\cP([1,n])} b^{\#(A)}a^{n-\#(A)}$. Hence
by \era2{1.6}
\beq4.102
(a+b)^n=\sum_{A\in\cP([1,n])}a^{\#(A)}b^{n-\#(A)}
=\sum_{A\in\cP([1,n])} a^{n-\#(A)}b^{\#(A)}.
\e
Finally, since $[1,n]=\bcl_{k\in[0,n]}\{A\in \cP{[1,n]}: \#(A)=k\}$ and
$\{A\in \cP{[1,n]}: \#(A)=k\}\cap\break\{A\in \cP{[1,n]}: \#(A)=l\}=\vn$ whenever
$k,l\in[0,n]$, $k\ne l$, we obtain in view of \era2{1.127}
\[
\sum_{A\in\cP([1,n])} a^{n-\#(A)}b^{\#(A)}\nad{\era2{1.126},\er{4.94},\er{4.95}}=
\sum_{k=0}^n \sum_{\sbk{A\in\cP([1,n])\\ \#(A)=k}} a^{n-\#(A)}b^{\#(A)}.
\]

\bnt4.61 \

\hph i,i, Let $n\in\Na$ and let $k\in[0,n]$.
\beq4.103
\cP_k([1,n]):= \{A\in\cP([1,n]): \#(A)=k\}.
\e

\indent \hph ii,, Let $n\in\Na$ and let $k\in[1,n]$.
\beq4.104
\Inj([1,k],[1,n]):= \{f:[1,k]\to[1,n] \mid f \hbox{ is in\jc}\}.
\e
Using \er{4.103} and \era2{1.126} we obtain
\beq4.105
(a+b)^n = \sum_{k=0}^n \#(\cP_k([1,n]))a^{n-k}\cdot b^k.
\e
\ent

\bpr4.62
Let $n\in\Na$. Then
\bea4.106
\#(\Inj([1,k],[1,n])) =\frac{n!}{(n-k)!}, &\q k\in[1,n],\\
\#(\cP_k([1,n]))\cdot k! = \frac{n!}{(n-k)!}, &\q k\in[0,n], \lb{4.107}
\e
where
\beq4.108
0!:=1 \qh{and } m!:=\prod_{l=1}^m l,\ m\in\Na.
\e
\epr

\brs4.63 \

\hph i,ii, If $k=n$ in \er{4.106}, then $\Inj([1,k],[1,n])\nda23.11 =
\Bij([1,n])$, hence
\[
\#(\Inj([1,n],[1,n]))\nda24.50 = n! \nde4.108 = \frac{n!}{(n-n)!}\,.
\]

\indent \hph ii,i, If $k=n$ in \er{4.107}, then $\cP_n([1,n])\nda23.13 = \{[1,n]\}$.
Hence
\[
\#(\cP_n([1,n])\cdot n! = 1\cdot n! = n! \nde4.108 = \frac{n!}{(n-n)!}\,.
\]

\hph iii,, If $k=0$ in \er{4.107}, $\cP_0([1,n])=\{\vn\}$, hence $\#(\cP_0([1,n]
))\cdot 0! = 1\cdot 1=1 = \frac{n!}{(n-0)!}$.
\ers

\proof[Proof of \er{4.106}]
In view of Remark \rf{r4.63}\,(i) we may assume $k<n$. We define a \rl\ $\sim$
on~$\Bij([1,n])$ by setting $f\sim g$ if $f|_{[1,k]}=g|_{[1,k]}$, \fa $f,g
\in \Bij([1,n])$. One easily verifies that $\sim$~is an \ev ce \rl. We denote
by~$\cA$ the set of \ev ce classes. Since the \ev ce classes form a \pt\ of
$\Bij([1,n])$, we have $\Bij([1,n])=\bcl_{\a\in\cA}\a$. By \era2{4.50}
$\Bij([1,n])$ is finite and $(\#(\Bij([1,n]))=n!$. We claim that $\cA$ is \ep\
to $\Inj([1,k],[1,n])$. Since all \el s of an \ev ce class have the same \rt
ion to $[1,k]$, we may define a map $\Phi:\cA \to \Inj([1,k],[1,n])$ by setting:
Given $\a\in\cA$, $\Phi(\a):=f|_{[1,k]}$ \fe $f\in\a$. Observe that $\Phi(\a)$ is
an in\jc\ map from $[1,k]$ into $[1,n]$.

\ti{$\Phi$ is in\jc}: Let $\a_1,\a_2\in\cA$ be \st $\Phi(\a_1)=\Phi(\a_2)$. Let
$f\in\a_1$, $g\in\a_2$. Then $f|_{[1,k]}= \Phi(\a_1)=\Phi(\a_2)=g|_{[1,k]}$. Hence
$f\sim g$ and $\a_1=[f]=[g]=\a_2$ (see Exercise \rfa1{ex4.4}). \E\Tf $\Phi$~is
in\jc.

\ti{$\Phi$ is sur\jc}: Let $\vf\in\Inj([1,k],[1,n])$. Then $\vf$~is a bi\jn\
from $[1,k]$ onto $\vf([1,k])$. Hence $\#(\vf([1,k]))=\#([1,k])$ by Corollary
\rfa2{c3.11}. Note that $\#([k+1,n])=\#([1,n]\sm[1,k])\nda23.31 = n-k>0$.
Similarly $\#([1,n]\sm\vf([1,k])=n-k$. Hence by Lemma \rfa2{l3.12} \te s a
bi\jn\ $\psi:[k+1,n]\to[1,n]\sm\vf([1,k])$. We define a map $f:[1,n]\to
[1,n]$ by setting $f(x):=\vf(x)$, $x\in[1,k]$, and $f(x):=\psi(x)$,
$x\in[k+1,n]$. Then one verifies that $f\in\Bij([1,n])$ and $\Phi([f])=\vf$.
Since $\vf$~is arbitrary in $\Inj([1,k],[1,n])$, $\Phi$~is sur\jc.

\csq, $\cA$ and $\Bij([1,k],[1,n])$ are \ep\ hence $\cA$~is finite and
$\#(\cA)=\#(\Bij([1,k],[1,n]))$ by Corollary \rfa2{c3.11}. Let $m:=\#(\cA)$.

We claim that $\#(\a)=(n-k)!$ \fa $\a\in\cA$. Recall that if $X,Y$ are \ns s,\break
$f:X\to Y$ is a bi\jn\ and $A,B$ are \nss s of~$X$ \st $A\cup B=X$, $A\cap B
=\vn$, then $f(A)\cup f(B)=Y$ and $f(A)\cap f(B)=\vn$. Let $\a\in\cA$ and let
$f\in\a$. Then \fa $g\in\a$ we have $g|_{[1,k]}=f|_{[1,k]}$ and $g([1,k])=
f([1,k])$. \E\Tf\ $g\in\a$ is completely \ch ized by $f|_{[1,k]}$ and $g|_{[k+1,
n]}$. Since $g\in\a$ is bi\jc, $g|_{[k+1,n]}$ maps $[k+1,n]$ onto $A:=[1,n]\sm
f([1,k])$. Note that $g|_{[k+1,n]}: [k+1,n]\to A$ is in\jc, hence bi\jc. \If
that \te s a bi\jn\ between $\a$ and $\Bij([k+1,n],A)$ defined by $g\mt g|_{[k+1,
n]}$. Thus by Corollary \rfa2{c3.11} $\#(\a)=\#(\Bij([k+1,n],A))\nda24.2 =\break \#
(\Bij([k+1,n])) \nda24.5 = \#(\Bij([1,n-k]))\nda24.50 = (n-k)!$ since $n-k=
\#([1,n]) - \#([1,k])\nda23.31 = \#([k+1,n])$. This proves the claim. \If from
Lemma \rfa2{l4.6} that $n!=\#(\Bij([1,n]))=\#(\cA)\cdot\#(\a)=m\cdot (n-k)!$.
Thus $m=\frac{n!}{(n-k)!}$. Note that $n!=\bg(\prodl_{i=1}^k i)\cdot\bg(\prodl_
{i=k+1}^n i)$ by \era2{7.13}. Hence $m=\bg(\prodl_{i=k+1}^n i)$. Since $m=\#
(\Bij([1,k],[1,n])$, \er{4.106} holds.
\endproof

\proof[Proof of \er{4.107}]
In view of Remarks \er{r4.63}\,(ii)(iii) we may assume $n\ge2$ and
$k\in{[1,n-1]}$. \E\fe $A\in\cP_k([1,n])$ we set $B_A:=\{f:[1,k]\to A \mid
f\hbox{ is in\jc}\}$. Clearly $B_A\sbs \Inj([1,k],[1,n])$ \fe $A\in\cP_k
([1,n])$, hence $\bcl_{A\in\cP_k([1,n])}B_A \sbs \Inj([1,k],[1,n])$.
Conversely, if $f\in\Inj([1,k],[1,n])$, then $f:[1,k]\to f([1,k])$ is bi\jc,
hence\break $\#(f([1,k]))=\#([1,k])=k$ by Corollary \rfa2{c3.11}. Thus $f([1,k])$
belongs to $\cP_k([1,n])$, and $f\in\bcl_{A\in\cP_k([1,n])}B_A$. \If that
$\Inj([1,k],[1,n])=\bcl_{A\in\cP_k([1,n])}B_A$. \Mo $B_A=\{
f:[1,k]\to A\mid f \hbox{ bi\jc}\}$ \fe $A\in \cP_k([1,n])$ by \era2{3.11}.
\E\Tf $\#(B_A)=k!$ by \era2{4.2}, \era2{4.50}. \E\Ip $B_A\ne\vn$ \fe
$A\in\cP_k ([1,n])$. Finally, we show that if $A,A'\in\cP_k([1,n])$ with $A\ne
A'$, then $B_A\cap B_{A'}=\vn$. Suppose for \cd ion that \te s $f\in B_A\cap
B_{A'}$. Since $f\in \Bij([1,k],A)\cap \Bij([1,k],A')$, we have $A=f([1,k])
=A'$, a~\cd ion. \If that $\{B_A\}_{A\in\cP_k([1,k])}$ is a \pt\ of
$\Inj([1,k],[1,n])$. Since $\cP_k([1,n])\sbs \cP([1,n])$, which is finite by
Exercise \rfa2{ex3.27}, $\cP_k([1,n])$ is finite by \E\Pr\ \rfa2{p4.18}\,(i).
\Mo $\cP_k([1,n])\ne\vn$ since $[1,k]\in \cP_k([1,n])$. Set $m:=\cP_k([1,n])$.
Then $m\in \Na$ and \te s a bi\jn\ $\vf:[1,m]\to\cP_k([1,n])$. Thus
$\{B_{\vf(i)}\}_{i\in[1,m]}$ is a \pt\ of $\Inj([1,k],[1,n])$ \st
$\#(B_{\vf(i)})=k!$ \fa $i\in[1,m]$. Then, by Lemma \rfa2{l4.6} we obtain
\[
\frac{n!}{(n-k)!} \nde4.106 = \#(\Inj([1,k],[1,n])) = m\cdot k!
=\#(\cP_k([1,n]))\cdot k!,
\]
which implies \er{4.107} for $n\ge2$, $k\in[1,n-1]$. In view of Remarks
\rf{r4.63}\,(ii)(iii) the proof of \er{4.107} is complete.
\endproof

\E\ml ying both sides of \er{4.107} by $(n-k)!$ we obtain by \er{3.10},
$n!=\cP([1,n])\cdot k!\cdot(n-k)!$. Hence
\beq4.109
k!(n-k)! | n!, \q\ 0\le k\le n,\ n\in \N,
\e
and
\beq4.110
\cP_k([1,n]) = \frac{n!}{k!(n-k)!}, \q\ 0\le k\le n,\ n\in \N.
\e
Combining \er{4.105} and \er{4.110} we obtain

\def\namespec{The binomial formula}
\begin{tspc} \lb{t4.64}
Let $(X,+,\cdot,0,1)$ be a \sr. Then
\beq4.111
(a+b)^n = \sum_{k=0}^n \frac{n!}{k!(n-k)!}\,a^{n-k}b^k
\e
\fa $a,b\in X$ and $n\in\Na$.
\end{tspc}

\brs4.65 \

\hph i,i, We deduce from \er{4.111} that in formula \er{4.54}:
\beq4.112
c_n = \sum_{k=3}^n \frac{n!}{k!(n-k)!}\,a^{n-k}b^{k-3},\q\ n\ge3.
\e

\hph ii,, For $0\le k\le n$:
\beq4.113
\binom nk = \binom n{n-k},
\e
since $n-(n-k)=k$.
\ers

In what follows we will use the notation $a^{b^c}$ instead of $a^{(b^c)}$
for $a,b,c\in\N$.

\bpr4.65
Let $(X,+,\cdot,0,1)$ be a field of \ch istic $p$~prime. Let\break ${a:[1,m]\to X}$,
$m\ge2$, and $l\in\Na$. Then
\bga4.114
\Bigl(\sum_{i=1}^m a_i\Bigr)^{p^l} = \sum_{i=1}^m a_i^{p^l},\\
(a_1-a_2)^{p^l} = a_1^{p^l} - a_2^{p^l}. \lb{4.115}
\e
\epr

\proof \

\ti{\er{4.114} with $l:=1$, $m:=2$}: We claim that $p|\binom pk$ \fa $k\in
[1,p-1]$. If $k:=1$, then $\binom p1= \frac{p!}{1!(p-1)!} = \frac{(p-1)!p}
{(p-1)!1} = \frac p1=p$ by \er{4.108}, \era2{7.13}, \era3{1.42} and \era3{8.53},
hence $p|\binom p1$. \E\Tf we may assume $p\ge3$ and $k\in[2,p-2]$. We have
$p!=\prodl_{i=1}^p i = \Bigl(\prodl_{i=1}^{p-k} i\Bigr)\Bigl(\prodl_{i=
p-(k-1)}^p i\Bigr) = (p-k)! \cdot \Bigl(\prodl_{i=p-(k-1)}^{p-1} i\Bigr)\cdot p$
by \er{4.108} and \era2{7.13}, since $p-(k-1)\le p-1$ (recall $k\ge2$).
\E\Tf $\frac{p!}{(p-k)!} = \Bigl(\prodl_{i=p-(k-1)}^{p-1}i\Bigr)\cdot p$. But
$k!|\frac{p!}{(p-k)!}$ by \er{4.107}. We claim that $k!\nmid p$. Suppose, for
\cd ion, that $k!|p$. Since $k\ge2$, $2|k!$ by \er{4.108} and \era2{7.13},
hence $2|p$. A~\cd ion since $p$~is prime, $p\ge3$ and $2<3$. \E\Tf $k!\nmid p$.
\E\oh $\gcd(k!,p)\in[1,p]$ since $p$ is prime. Suppose, for \cd ion, that
$p|k!$. By \er{4.108} and Lemma \rfa3{l9.17} (Euclid's lemma), $p|i$ \fs $i\in
[1,k]$, which is impossible by \era3{8.57} since $i\le k< p$. Therefore,
since $\gcd(k!,p)=1$, and $k!\nmid p$, we obtain $k!|\prodl_{i=p-(k-1)}^{p-1}i$
by Lemma \rfa3{l9.17}.

\E\Tf \te s $d\in\Na$ \st $\prodl_{i=p-(k-1)}^{p-1} i=d\cdot k!$. Hence
$\frac{p!}{(p-k)!}=d\cdot k!\cdot p$. \If that $\binom pk = \frac{p!}{k!
(n-k)!}= \frac{d\cdot k!\cdot p}{k!}=dp$. Thus $p|\binom pk$, which completes
the proof of the claim.

We now claim that $\binom pk x=0$ \fe $x\in X$ and every $k\in[1,p-1]$. We
recall that $\binom pk x\nde 4.94 {:=} \binom pk\dpl x$ in the abelian group
$(X,+,0)$. We define $\vf:X\to X$ by setting $\vf(y):=x\cdot y$ ($=y\cdot x$)
\fa $y\in X$. Then $\vf$~is an endo\mf\ of $(X,+,0)$ by \er{1.1n}, \er{1.2n}.
Recall that $x=1\cdot x=x\cdot1$, hence $\binom pk\dpl x=\binom pk\dpl\vf(1)
\nda21.45 = \vf(\binom pk\dpl 1) = \vf((d\cdot p)\dpl 1)
\nad{\era2{2.3}\,{\rm I3}} = \vf(d\dpl(p\dpl1)) = \vf(d\dpl0)$ since $X$~is
a field of \ch istic~$p$. But $d\dpl0  = 0$ by \era2{2.3}\,I4, hence
$\vf(d\dpl0)=\vf(0)=0$, which proves the claim.

We now prove \er{4.114} for $m=2$. We have $\bigl(\suml_{i=1}^2 a_i\bigr)^p
\nde4.111 = \suml_{k=0}^p \binom pk a_1^{p-k}a_2^k \nda2{7.13} = \binom p0
a_1^pa_2^0 + \suml_{k=1}^{p-1} \binom pka_1^{p-k}a_2^k + \binom pp a_1^0a_2^p
= a_1^p + \suml_{k=1}^{p-1}0 +a_2^p$ by Remarks \rf{r4.63}\,(ii)(iii), and
what precedes. But $\suml_{i=1}^{p-1}0\nda21.126 = (p-1)\dpl0 \nad{\era2{2.3}
\,{\rm I4}}=0$. \E\Tf \er{4.114} holds for $m=2$, $n=p$.

\ti{Case $m:=2$, $n=p^l$, $l\ge1$}: We proceed by \In\  on $l\in\Na$. Let
$M=\{l\in\Na:\break \hbox{\er{4.114} holds with }m:=2\}$. We just proved that
$1\in M$. Suppose that $l\in M$. Then $(a_1+a_2)^{p^{l+1}}=(a_1+a_2)^{p^l
\cdot p^1}=((a_1+a_2)^{p^l})^p$ by \er{4.95}, \era2{2.25}, \era2{2.24},
\era2{2.23}\,I3. But $(a_1+a_2)^{p^l}\nad{l\in M}= a_1^{p^l}+a_2^{p^l}$,
hence $(a_1+a_2)^{p^{l+1}}= (a_1^{p^l}+a_2^{p^l})^p \nad{1\in M}= (a_1^{p^l})^p
+(a_2^{p^l})^p = a_1^{p^{l+1}}+a_2^{p^{l+1}}$. Hence $l+1\in M$. Thus $M=\Na$.

\ti{Case $m\ge2$, $n:=p^l$, $l\ge1$}: By \In\ on $m\ge2$. \er{4.114} holds for
$m=2$. Suppose it holds for~$m$ $(\ge2)$. Then $\Bigl(\suml_{i=1}^{m+1}a_i
\Bigr)^{p^l}= \Bigl(\Bigl(\suml_{i=1}^{m}a_i\Bigr)+a_{m+1}\Bigr)^{p^l}=
\Bigl(\suml_{i=1}^{m}a_i\Bigr)^{p^l}+(a_{m+1})^{p^l} = \suml_{i=1}^{m}(a_i)
^{p^l}+(a_{m+1})^{p^l} = \suml_{i=1}^{m+1}(a_i)^{p^l}$.

\er{4.115}: $a_1-a_2:=a_1+(-a_2)$, hence $(a_1-a_2)+a_2= (a_1+(-a_2))+a_2 =
a_1+\break((-a_2)+a_2) = a_1+0=a_1$. Hence $((a_1-a_2)+a_2)^{p^l}=a_1^{p^l}$ and by
\er{4.114} $((a_1-a_2)+a_2)^{p^l} = (a_1-a_2)^{p^l}+a_2^{p^l}$. Thus
$(a_1-a_2)^{p^l}=a_1^{p^l}-a_2^{p^l}$.
\endproof

\bex4.66
Let $k,l\in\Na\sms1$. Show

\hph i,ii, $kl\ddtt{kl}(k+_{kl} l)= (kl\ddtt{kl}k)+_{kl}(kl\ddtt{kl}l) $.

\hph ii,i, $2^n\ddtt{2^n}(1+_{2^n}1)\ne (2^n\ddtt{2^n}1)+_{2^n}(2^n\ddtt{2^n}1)$
\fa $n\in\Na\sms1$.
Does this contradict \E\Pr\ \rf{p4.65}?

\hph iii,, $6\ddtt6(1+_61) = (6\ddtt61)+_6(6\ddtt61)$.

\eex

\bex 4.67 \

\hph i,i, Prove Pascal's identity:
\beq4.116
\binom nk + \binom n{k-1} = \binom {n+1}k, \q\ 1\le k\le n.
\e

\hph ii,, Prove \er{4.116} by \In\ on $n\ge2$.
\eex

One can ``visualize'' \er{4.116} in the \fw\ triangle:

\setbox0=\hbox{11} \sz=\wd0
\centerline{Pascal's triangle}
\[
\begin{matrix}
&&&&\binom 00&&&&\\
&&&\binom 10&&\binom 11&&&\\
&&\binom 20&&\binom 21&&\binom 22&&\\
&\binom 30&&\binom 31&&\binom 32&&\binom 33&\\
\binom 40&&\binom 41&&\binom 42&&\binom 43&&\binom 44
\end{matrix} \qquad \qquad
\begin{array}{ccccccccccc}
\hphantom{11}&\hphantom{11}&\hphantom{11}&\hphantom{11}&\hphantom{11}&\hbox
to\sz{\hfil1\hfil}&\hphantom{11}&\hphantom{11}&\hphantom{11}&\hphantom{11}&\hphantom{11}\\
&&&&1&&1&&&&\\
&&&1&&2&&1&&&\\
&&1&&3&&3&&1&&\\
&1&&4&&6&&4&&1&\\
1&&5&&10&&10&&5&&1
\end{array}
\]

We next present a \gn\ of \E\Pr\ \rf{p4.59} and the binomial formula \rf{t4.64}.

\bpr4.68
Let $(X,+,\cdot,0,1)$ be a \sr, let $a:[1,m]\t[1,n] \to X$, $m,n\ge2$. Then
\beq4.117
\prod_{j=1}^n \Bigl(\sum_{i=1}^m a_{ij}\Bigr) = \sum_{g\in[1,m]^{[1,n]}}
\Bigl(\prod_{j=1}^n a_{g(j)j}\Bigr).
\e
\epr

\bex4.69
Prove \er{4.117} by \In\ on $m\ge2$.
\eex

\begin{prp}[Multinomial formula]\lb{p4.70}
Let $(X,+,\cdot,0,1)$ be a \sr. Let $m,n\ge2$, $a:[1,m]\to X$. Then
\beq4.118
\Bigl(\sum_{i=1}^m a_i\Bigr)^n = \sum_{\sbk{\a:[1,m]\to \N\\|\a|=n}}
\binom{|\a|}{\a!}\ddt a^\a
\e
where $|\a|$ and $\a!$ are defined in \er{3.62}, $\binom{|\a|}{\a!}$ is
defined in \er{3.63}, and
\beq4.119
a^\a := \prod_{i=1}^m a_i^{\a_i}.
\e
\epr

\bex4.71 \

\hph i,i, Prove \er{4.118} by \In\ on $m\ge2$.

\hph ii,, Derive \er{4.118} from \er{4.117}.

\eex

\newpage
\Subsubsection{Infinite prime fields}\label{sss.ipr.fld}

One of the main goals of this section is to prove the \ex\ of a \ct y infinite\glossary{$\Q$}
\Pf~$\Q$ called the field of \ti{\ra\ \nm s},\index{field!of rational numbers} to study some of its \pp ies and
to show that every infinite \Pf\ is ring-\is c to~$\Q$.
We give two different proofs of the \ex\ of the field~$\Q$.

In the first one, we embed $(\Na,\cdot,1)$ in an \ag\glossary{$\Q_{>0}$} $(\Q_{>0},\cdot,
1)$, the \mlv\ group of \po\ \ra\ \nm s. Then we introduce on $\Q_{\ge0}:=
\Q_{>0}\cup\{0\}$ an \ad~$+$ and a \mlc~$\cdot$ which make\glossary{$\Q_{\ge0}$}
$(\Q_{\ge0},+,\cdot,0,1)$ a \sr\ \st
the monoid $(\Q_{\ge0},+,0)$ is a \PM\ totally ordered by its \nog. We finally
embed $\Q_{\ge0}$ into a field by means of Corollary \rf{c3.85}\,(i). This
field is denoted by~$\Q$.

In the second approach we first introduce a \mlc~$\cdot$ on the group $(\Z,+,0)$
which makes $(\Z,+,\cdot,0,1)$ an \ido. We then prove the classical theorem
allowing to embed an \ido\ in a field (\cite{Nrs,Alg}).

\ssk
Let $(\N,\cdot,1)$ denote the \mlv\ \am\ of the \sr\ $(\N,+,\cdot,0,1)$. We
recall that $\Na$~is a \sbm\ of $(\N,\cdot,1)$ by \E\Pr\ \rfa2{p2.7}\,(i),
and by the same \Pr\ $(\Na,\cdot,1)$ is a \PM. Since the map $n\mt n+1$ from
$\N$ into~$\Na$ is bi\jc, $\Na$~is \ct y infinite. Thus we may apply Theorem
\rf{t3.84} with $(X,\qu,e):=(\Na,\cdot,1)$. Then $\wh X$ is usually called
the set of \ti{\po\ \ra\ \nm s} and is sometimes denoted by $\Q_{>0}$. The
binary \op~$\hqu$ is usually denoted by~$\cdot$ or is simply omitted and is
called \ti{\mlc}. It is customary to denote an \el\ of~$\Q_{>0}$ by a
\ti{fraction}. If $[(a,b)]$ denotes the $\sim$-\ev ce class containing
$(a,b)\in\Na\t\Na$, then one writes
\beq5.1
\frac ab:=[(a,b)].
\e
The symbol $\frac ab$ is called a \ti{fraction}\index{fraction} with \ti{numerator}~$a$ and\index{numerator}
\ti{denominator}~$b$.\index{denominator}

If $a,b,c,d\in\Na$, then $\frac ab=\frac cd$ means $[(a,b)]=[(c,d)]$, \ev tly
$(a,b)\sim(c,d)$ or $ad=bc$ by \er{3.100}. Hence, given $a,b,c,d\in\Na$:
\beq5.2
\frac ab=\frac cd \qh{iff }ad=bc,
\e
where the first \et y sign denotes the \et y in $\Q_{>0}$.

\Mo $[(a,b)]\hqu[(c,d)] \nde3.101 = [(a,b)\tqu(c,d)]
\nde3.99 = [(a\qu c,b\qu d)]\nde5.1 = \frac{ac}{bd}$, where $ac$ (resp.~$bd$)
is the product of $a$ and~$b$ (resp.\ $c$ and~$d$) in $(\Na,\cdot,1)$.

Thus we may define the \mlc\ in $\Q_{>0}$ by
\beq5.3
\frac ab\cdot \frac cd:= \frac{ac}{bd}\,, \q a,b,c,d\in\Na.
\e
The \nel\ is $\wh e=[(e,e)]=[(1,1)]=\frac11$. Hence by Theorem \rf{t3.84},
$(\Q_{>0},\cdot,\frac11)$ is an \ti{\ag}, and $i:\Na\to\Q_{>0}$ defined by
$i(a):=[(a,1)] \nde5.1 = \frac a1$, $a\in\Na$, is an \ti{in\jc\ \hm sm}.
Note that for $a,b\in\Na$
$$
\frac ab\cdot \frac ba\nde5.3 = \frac{ab}{ba}\nad{\er{5.2},\era2{1.6}} = \frac 11\,.
$$
\If that $\frac ba$ is the inverse of $\frac ab$ in $(\Na,\cdot,\frac11)$.

It is customary to write $a$ instead of~$i(a)$ (abuse of notation), i.e.\ to
write $a$~instead of $[(a,1)]$. \E\Ip one writes $1$ instead of~$\frac11$. So
we have
\beq5.4
a=\frac a1, \q a\in\Na.
\e
By doing so, $\Na$ can be viewed as  a \sbm\ of~$\Q_{>0}$.

One may also define another binary \op\ on~$\Q_{>0}$ called \ti{division}.
Since $\Q_{>0}$ is a group, given $q,r \in \Q_{>0}$, there is a unique
$x\in\Q_{>0}$ \st $q\cdot x=r$. Such an $x$~is called the \ti{quotient\/}
of~$r$ by~$q$ and the binary \op\ on~$\Q_{>0}$ defined by $(r,q)\mt x$ is
called \ti{division}. Note that if $r,q\in\Na$ then $x=\frac rq$. Indeed,
we have $q\cdot\frac rq\nde5.4 = \frac q1\cdot\frac rq \nde5.3 =
\frac{q\cdot r}q \nde5.2 = \frac r1 \nde5.4 = r$. Hence the fraction
$\frac rq$ can be interpreted (\ti{a~posteriori}) as the quotient of~$r$
by~$q$. Now for arbitrary \el s $q,r\in\Q_{>0}$, it is natural to use the
notation
\beq5.5
\frac rq := \hbox{the quotient of $r$ by }q.
\e

\bex5.1 \

\hph i,ii, Show that if $a,b,c,d\in\Na$, then
\beq5.6
\frac{\,\frac ab\,}{\,\frac cd\,} = \frac{ad}{bc}\,.
\e

\hph ii,i, Establish formula \er{5.3} for $a,b,c,d\in\Q_{>0}$.

\hph iii,, Show that $\frac ba$ is the inverse of $\frac ab$ in $(\Q_{>0},\cdot,1)$
\fa $a,b\in\Q_{>0}$.

\eex

If in Exercise \rf{ex3.145} $(G,\qu,e):=(\Q_{>0},\cdot,1)$ and $q\in\Q_{>0}$,
then the $z$-\IT\ of~$q$ with $z\in\Z$ is denoted by~$q^z$ ($q$~to the
power~$z$). \E\Ip if $q:=\frac ab$, $a,b\in\Na$ and $z:=-n$, $n\in\N$, then
$$
\biggl(\frac ab\biggr)^{-n}= \biggl(\frac ba\biggr)^{n}\,.
$$

\bex5.2 \

\hph i,i, Show that the map $z\mt q^z$ from~$\Z$ into~$\Q_{>0}$,
$q\in\Q_{>0}$, is an in\jc\ \hm sm iff $q\ne1$.

\hph ii,, Show:
\bga5.7
\biggl(\frac ab\biggr)^{-z} = \Bgg(\frac ba)^z, \q a,b\in\Na,\ z\in\Z,\\
(q\cdot r)^z = q^z\cdot r^z, \q q,r\in\Q_{>0},\ z\in\Z.\lb{5.8}
\e
\eex

We now set
\beq5.9
\Q_{\ge0}:=\Q_{>0} \cup\{0\} \qh{with }0\in\N.
\e
Then $\Q_{\ge0}$ is called the set of \ti{non\ng\ \ra\ \nm s}.
We now introduce a \mlc~$\cdot$ on~$\Q_{\ge0}$
which extends the \mlc\ on~$\Q_{>0}$. We set
\beq5.10
q\cdot r:=\bca
0 &\hbox{if $q$ or $r$ is equal to 0},\\
q\cdot r, &\hbox{the product of $q$ and $r$ in $\Q_{>0}$}.
\eca
\e
Let $q,r,s\in\Q_{\ge0}$. Then $(q\cdot r)\cdot s = q\cdot(r\cdot s)$. Indeed,
if $q$ or $r$ or~$s$ equals~$0$, then the above \et y is a direct con\sq\ of
\er{5.10}. Otherwise, it follows from the \asc ity of~$\cdot$ in $\Q_{>0}$.
Clearly $q\cdot r=r\cdot q$ follows from \er{5.10} and the \cmt ity of~$\cdot$
in~$\Q_{>0}$. \E\Tf $(\Q_{\ge0},\cdot,1)$ is an \ti{\am\/} (which is not \cnc
e since $0\cdot0=0=0\cdot1$ and $0\ne1$. However, $\Q_{>0}$ is a
\ti{subgroup} of $(\Q_{\ge0},\cdot,1)$.

We want to introduce an \ad~$+$ on~$\Q_{\ge0}$ which makes $(\Q_{\ge0},+,\cdot,
0,1)$ a (\cmt e) \sr\ (with unity) containing~$\N$ (by \er{5.4}, \er{5.9}), \st
$\N$~is a \sbm\ of $(\Q_{\ge0},+,0)$. In this case if
$q\in\Q_{\ge0}$, then $q+0=q=0+q$. \Mo if $q,r\in\Q_{>0}$,
then \te\ $a,b,c,d\in\Na$ \st $q=\frac ab$, $r=\frac cd$. Then
$q+r= \frac ab+\frac cd \nde5.2 = \frac{ad}{bd}+\frac{bc}{bd} \nde5.3 =
\frac{ad}1 \cdot \frac1{bd} + \frac{bc}1\cdot \frac1{bd} \nde1.3n =
\bg(\frac{ad}1+\frac{bc}1)\cdot\frac1{bd} \nad*= (ad+bc)\cdot\frac1{bd} =
\frac{ad+bc}1 \cdot \frac1{bd}\nde5.3 = \frac{ad+bc}{bd}$. In $\nad*=$ we used
\er{5.4} and $\N$~being \sbm\ of $(\Q_{\ge0},+,0)$. \csq, we define the \ad~$+$
on~$\Q_{\ge0}$ by setting for $q,r\in\Q_{\ge0}$:
\bga5.11
q+0:=q, \q 0+q:=q, \\
\frac ab+\frac cd :=\frac{ad+bc}{bd} \qh{for} q:=\frac ab\,,\ r:=\frac cd\,,\
a,b,c,d\in\Na.\lb{5.12}
\e
However, this \df\ only makes sense if it does not depend on the choice of the
\rp atives $\frac ab$ of~$q$, and $\frac cd$ of~$r$. We need to show that
$\frac ab+\frac cd = \frac{a'}{b'}+\frac{c'}{d'}$ if $\frac ab=\frac{a'}{b'}$
and $\frac cd=\frac{c'}{d'}$, $a,b,c,d,a',b',c',d'\in\Na$. This is done in
a more general situation in the \fw\ lemma.

\blm5.3
Let $(X,+,\cdot,0,1)$ be a $($\cmt e$)$ ring $($with unity$)$, and let
$a,b,c,d$, $a',b',c',d'\in X$. If
\beq5.13
ab'=ba' \qh{and } cd'=dc',
\e
then
\beq5.14
(ad+bc)b'd' = bd(a'd'+b'c').
\e
\elm

\proof
$(ad+bc)b'd' \nde1.3n = (ad)(b'd')+(bc)(b'd') \nda21.36 = a(db')d' +b(cb')d'
\nda21.6 = a(b'd)d' + b(b'c)d' \nda21.36 = (ab')(dd') + (bb')(cd') \nde5.13 =
(ba')(dd')+(bb')(dc') \nda21.36 = b(a'd)d'+b(b'd)c' \nda21.6 = b(da')d' +
b(db')c' \nda21.36 = (bd)(a'd')+(bd)(b'c')\nde1.2n = bd(a'd'+b'c')$.
\endproof

Thus the \ad\ in $\Q_{\ge0}$ is well defined. We now show that $(\Q_{\ge0},+,0)$ is a \PM.

\era2{1.7}: Let $q\in\Q_{\ge0}$. Then $q+0=q=0+q$ by \er{5.11}.

\era2{1.5}: Let $q,r,s\in\Q_{\ge0}$. If $q$ or $r$ or $s$ equals zero then \era2{1.5}
easily follows from \era2{1.7}. If $q:=\frac ab$, $r:=\frac cd$, $s:=\frac ef$,
$a,b,c,d,e,f\in \Na$, then:
\bgg
\frac ab+\Bigl(\frac cd+\frac ef\Bigr)=\frac ab+\frac{cf+de}{df}=
\frac{adf+b(cf+de)}{b(df)}\nad{\er{1.2n},\era2{1.5}} = \frac{adf+bcf+bde}{bdf}\\
\Bigl(\frac ab+\frac cd\Bigr)+\frac ef = \frac{ad+bc}{bd} +\frac ef =
\frac{(ad+bc)f+bde}{(bd)f}\nad{\er{1.3n},\era2{1.5}} = \frac{adf+bcf+bde}{bdf}\,.
\e

\era2{1.6}: Let $q:=\frac ab$ and $r:=\frac cd$. Then
$$
\frac ab+\frac cd\nde5.12 = \frac{ad+bc}{bd} = \frac{da+cb}{db} =\frac{cb+da}
{db} \nde5.12 = \frac cd+\frac ab\,.
$$
If $q$ or $r$ equals $0$, then \era2{1.6} follows from \era2{1.7}.

\era2{1.8}: Let $q,r,s\in\Q_{\ge0}$ be \st $q+s=r+s$. We have to show $q=r$.
If $s:=0$, then $q=r$ by \era2{1.7}. We suppose $s:=\frac ef$, $e,f\in\Na$. If
$r:=0$, then $q+\frac ef=\frac ef$. We have to show $q=0$. Suppose for \cd ion
that $q=\frac ab$, $a,b\in\Na$. Then $\frac{af+be}{bf} = \frac ef = \frac{be}
{bf}$. From $\frac{af+be}{bf} = \frac{be}{bf}$ we obtain $af+be \nde5.4 =
\frac{af+be}{1}\nde5.2 = \frac{(af+be)bf}{1(bf)} \nde5.3 = \frac{(af+be)}{bf}
\cdot \frac{bf}1 = \frac{be}{bf}\cdot\frac{bf}1 \nde5.3 = \frac{(be)(bf)}{(bf)1}
\nde5.2 = \frac{be}1 \nde5.4 = be$. Hence $af+be = be=0+be$. By \cnc ity in
$(\N,+,0)$ we have $af=0$. Then $a$ or~$f$ equals zero, by \E\Pr\
\rfa2{p2.7}\,(ii). A~\cd ion. Finally, if $q:=\frac ab$, $r:=\frac cd$, $a,b,c,d
\in\Na$, then we have $\frac{af+be}{bf} = \frac{cf+de}{df}$. We want to infer
$ad=bc$. By \er{5.2} we obtain $(af+be)(df)=(bf)(cf+de)$. Hence by \er{1.3n},
\er{1.2n} $(af)(df)+(be)(df) = (bf)(cf)+(bf)(de)$. \Mo $(be)(df) = (be)(fd) =
b(ef)d = b(fe)d = (bf)(ed) = (bf)(de)$. Hence by \era2{1.8} in $(\N,+,0)$ we
have $(af)(df) = (bf)(cf)$, thus $a(fd)f=b(fc)f$. By \E\Pr\ \rfa2{p2.7}\,(iii)
we obtain $a(fd)=b(fc)$. Hence $(ad)f = a(df) = a(fd) = b(fc) = b(cf) = (bc)f$.
By the same \Pr\ we arrive at $ad=bc$, which completes the proof of \era2{1.8}.

\era2{1.9}: Let $q,r\in\Q_{\ge0}$ be \st $q+r=0$. We have to show $q=r=0$. If
$q,r\in\Q_{>0}$, then $q+r\in\Q_{>0}$ by \er{5.12}. Hence if $q+r=0$, $q,r\in
\Q_{\ge0}$, then $q=0$ or $r=0$. If $q=0$ (resp.\ $r=0$), then $r=0$ (resp.\
$q=0$) by \era2{1.7}.

\If that $(\Q_{\ge0},+,0)$ is a \ti{\PM}.
\ssk

We next show that $\Q_{\ge0}$ is \ti{totally ordered} by the \nog\ of
$(\Q_{\ge0},+,0)$. To this end we first prove that if $q,r\in\Q_{>0}$, then
$q<r$ (in the sense of \era3{8.7}) iff $ad<bc$ whenever $q=\frac ab$,
$r=\frac cd$, $a,b,c,d\in\Na$.

\ti{Only if\/}: By \era3{8.7} \te s $s\in\Q_{>0}$ \st $r=q+s$. Let $a,b,c,d,
e,f\in\Na$ be \st $q=\frac ab$, $r=\frac cd$ and $s=\frac ef$. Then by
\er{5.12} $\frac cd=\frac{af+be}{bf}$. Thus $cbf\nde5.2 = d(af+be)\nde1.2n =
daf+dbe$. Hence $daf<cbf$ by \era2{1.46} and $da<cb$ by \era2{2.20}. Thus
$ad<bc$.

\ti{If\/}: Let $q=\frac ab$, $r=\frac cd$, $a,b,c,d\in\Na$ be \st $ad<bc$.
By \era2{1.46} \te s $p\in\Na$ \st $bc=ad+p$. Set $s:=\frac p{bd}\in\Q_{>0}$.
Then $r=\frac cd \nde5.2 = \frac{bc}{bd} = \frac{ad+p}{bd} \nde5.2 =
\frac{(ad+p)(bd)}{(bd)(bd)} \nde1.3n = \frac{(ad)(bd)+p(bd)}{(bd)(bd)}
\nde5.12 = \frac{ad}{bd}+\frac p{bd} \nde5.2 = \frac ab+s=q+s$. \E\Tf $q<s$
by \era3{8.7}.

Note that if $r\in\Q_{\ge0}$, then $r\in\Q_{>0}$ iff $r>0$. Indeed, if $r>0$,
then $r\ne0$ by \era1{3.8}, hence $r\in\Q_{>0}$ by the \df\ of $\Q_{\ge0}$.
Conversely, if $r\in\Q_{>0}$, then $r\nde5.11 = 0+r$, hence $r>0$ by
\era3{8.7}, \era3{8.9}. This ``justifies'' the notation $\Q_{>0}$.

\Wanp show that the natural order $\le$ is total. Let $q,r\in\Q_{\ge0}$ be \st
$q\ne r$. If $q=0$, then $r\in\Q_{>0}$ hence $q< r$. Let now $q,r\in\Q_{>0}$
and $q=\frac ab$, $r=\frac cd$ with $a,b,c,d\in\Na$. Since the order of $(\N,
+,0)$ is total, we have either $ad<bc$ or $bc<ad$ in view of \er{5.2} and
$q\ne r$. If $ad<bc$, then $q<r$ by what precedes. If $bc<ad$, we obtain
$r<q$ by exchanging the roles of $q$ and~$r$.

Summarizing, we proved

\bpr5.4
Let $\Q_{\ge0}$ be defined as in \er{5.9}, let $+$ be defined by \er{5.11},
\er{5.12}, then $(\Q_{\ge0},+,0)$ is a \PM\ and its natural ordering defined in
\era3{8.7} \sf ies
\beq5.15
q<r \qh{iff }ad<bc
\e
whenever $q=\frac ab$, $r=\frac cd$, $a,b,c,d\in\Na$. \E\Ip the \nog~$\le$ is
total.
\epr

\Wanp apply Corollary \rf{c3.85}\,(i) with $(X,\qu,e):=(\Q_{\ge0},+,0)$.
Let $\wh X, \wh X_+, \wh X_-$ and~$j$ be as in Corollary \rf{c3.85}\,(i). Set
\beq5.16
\left\{\bal
+q &:= j(q), \q q\in\Q_{>0},\\
0&:=\wh e,\\
-q&:= (j(q))\Inv.
\eal \right.
\e
Then by \er{3.107}, \er{3.108}, \er{3.109}:
\beag
\wh X_+&= \{+q: q\in\Q_{>0}\},\\
\wh X_-&= \{-q: q\in\Q_{>0}\},\\
\wh X\ &\hbox{is the disjoint union of }\wh X_+,\{0\}\hbox{, and }\wh X_-.
\e
We denote the \ad\ $\hqu$ in $\wh X$ by~$+$. There is no ambiguity with the
\ad\ in $\Q_{\ge0}$ since we have different notations for the \el s of~$\wh
X$ and~$\Q_{\ge0}$ (except~$0$). In view of Corollary \rf{c3.85}\,(i), $(\wh
X,+,0)$ is an \ag\ \sf ying
\beq5.17
\left\{\bal
(+q)\Inv &= -q,\\
0\Inv &= 0,\\
(-q)\Inv &= +q, \quad\ q\in\Q_{>0}.
\eal \right.
\e
We now introduce a \mlc\ on~$\wh X$, which we denote by~$\cdot$, by setting
\beq5.18
0\cdot \wh x:=0 \qh{and } \wh x\cdot 0=0 \qh{\fa } \wh x\in\wh X,
\e
and \fa $q,r\in\Q_{>0}$:
\beq5.19
\bal
(+q)\cdot (+r)&:=+(q\cdot r),\\
(-q)\cdot (-r)&:=+(q\cdot r),\\
(+q)\cdot (-r)&:=-(q\cdot r),\\
(-q)\cdot (+r)&:=-(q\cdot r),
\eal
\e
where $q\cdot r$ is the product of $q$ and $r$ in $\Q_{>0}$.

\bex5.5
Show that the \mlc\ on~$\wh X$ defined by \er{5.18}, \er{5.19} is the only
one which makes $(\wh X,+,\cdot,0,+1)$ a \cmt e \sr\ with unity $+1$ \sf ying
$(+q)\cdot(+r)=+(q\cdot r)$, $q,r\in\Q_{>0}$.
\eex

\bth5.6
$(\wh X,+,\cdot,0,+1)$ is a field.
\eth

\proof
Since $(\wh X,+,0)$ is an \ag, it suffices to show that $(\wh X,\cdot,+1)$ is
an \am, that $\wh X\sms0$ is a subgroup of $(\wh X,\cdot,+1)$ and \er{1.1n},
\er{1.2n} or \er{1.3n} hold.

We first observe
\beq5.20
\wh x,\wh y\in\wh X\sms0 \qh{implies }\wh x\cdot\wh y\in\wh X\sms0,
\e
and we next show
\beq5.21
(\wh X\sms0,\cdot,+1) \hbox{ is an \ag.}
\e
It is convenient to write $\wh x\in\wh X\sms0$ as an \el\ of the product
$\{+,-\}\t \Q_{>0}$, that is, to write $\wh x=(\sgn(\wh x),|\wh x|)$ defined
in \er{3.107} and \er{3.106}. We define a binary \op\ on $\{+,-\}$, which we
denote by~$\cdot$, by setting
\beq5.22
\left\{\bal
+ \cdot +&:=+,\\
- \cdot - &:=+,\\
+\cdot -&:=-,\\
-\cdot +&:=-.
\eal \right.
\e
\If from \er{5.19} that
\beq5.23
\sgn(\wh x\cdot\wh y) = \sgn(\wh x)\cdot\sgn(\wh y), \q \wh x, \wh y\in\wh
X\sms0.
\e
Note that $+$ is a \nel\ for $\cdot$. From Examples \rfa2{xa1.5}\,(iii) we
infer that $(\{+,-\},\cdot,+)$ is \asc e and \cmt e. Since $-\cdot- = +$,
$(\{+,-\},\cdot,+)$ is an \ag\ (\is c to $(\N_2,+_2,0)$). For a more explicit
argument see the proof of Theorem \rf{t5.7} below \er{5.35}--\er{5.45}.

We know that $(\Q_{>0},\cdot,+1)$ is a group (see the line below \er{5.3}).

Observe that the map $|\cdot| : \wh X\sms0\to\Q_{>0}$ \sf ies
\beq5.24
|\wh x\cdot\wh y| = |\wh x|\cdot|\wh y|.
\e
\E\Tf $(\wh X\sms0,\cdot,+1)$ is the product of the groups $\{+,-\}$ and
$\Q_{>0}$ by \E\Pr\ \rf{p2.20}. Then $(\wh X\sms0,\cdot,+1)$ is an \ag\ by
Corollary \rf{c3.37}.

We now show that $(\wh X,\cdot,+1)$ is an \am. The \asc ity and \cmt ity of the
product~$\cdot$ follows from \er{5.18} when one of the factors is equal
to~$0$. Otherwise the \asc ity and \cmt ity of~$\cdot$ follows from the
\asc ity and \cmt ity of~$\cdot$ in $\wh X\sms0$. Finally
$(+1)\cdot \wh x=\wh x$ \fa $\wh x\in\wh X$ by \er{5.18} and $\sgn((+1)\cdot
\wh x)= {+}\cdot\sgn(\wh x)=\sgn(\wh x)$, $|(+1)\cdot\wh x|=|{+1}|\cdot|\wh x| =
1\cdot|\wh x| = |\wh x|$.

In view of \er{5.20}, \er{5.21} it remains to show that \er{1.1n},
\er{1.2n}, \er{1.3n} hold. Clearly \er{1.1n} follows from \er{5.18}.

Since the \mlc~$\cdot$ is \cmt e, \er{1.2n} implies \er{1.3n} (see Exercise
\rf{ex1.2n}). Observe that if $\wh x$ or $\wh y$ or~$\wh z$ equals~0, then
\er{1.2n} follows from \er{5.18} and \era2{1.7}. \E\Tf it is \sft\ to prove
\er{1.2n} for $\wh x,\wh y,\wh z\in\wh X\sms0$. We first consider the case
where $\wh x,\wh y,\wh z\in \wh X_+$, which we call case~(i),

\ti{Case} (i): Let $\wh x:=+\frac ab$, $\wh y:=+\frac cd$, $\wh z:=+\frac ef$
with $a,b,c,d,e,f\in\Na$. We have $\bg(+\frac cd)+\bg(+\frac ef) \nde5.16 =
j\bg(\frac cd)+j\bg(\frac ef)\nad*= j\bg(\frac cd+\frac ef)\nde5.12 =
j\bg(\frac{cf+de}{df}) \nde5.16 = +\frac{cf+de}{df}$, where in $\nad*=$ we
used the fact that $j$~is a \hm sm. \E\Tf $\wh x\cdot(\wh y+\wh z) =
\bg(+\frac ab)\cdot\bg({\bg(+\frac cd)+\bg(+\frac ef)}) = \bg(+\frac ab)\cdot
\bg(+\frac{cf+de}{df})\nde5.19 = +\frac ab\cdot \frac{cf+de}{df} \nde5.12 =
+\frac{a(cf+de)}{b(df)} \nad{\era2{2.14},\era2{1.5}} = +\frac{acf+ade}{bdf}$.
We next use the \fw\ formula in $\Q_{>0}$:
\beq5.25
\frac km+\frac lm = \frac{k+l}m, \q k,l,m\in\Na.
\e

\proof[Proof of \er{5.25}]
$\frac km+\frac lm\nde5.12 = \frac{km+ml}{mm} \nad{\era2{1.6},\era2{2.15}} =
\frac{(k+l)m}{mm} \nde5.3 = \frac{k+l}m \,\frac mm \nde5.2 =
\frac{k+l}m\,\frac11 = \frac{(k+l)1}{m1}\break \nad{\era2{1.6},\era2{2.13}} =
\frac{k+l}m$.
\endproof

By \er{5.25} we obtain $\frac{acf+ade}{bdf}=\frac{acf}{bdf}+\frac{ade}{bdf}
\nda21.6 = \frac{acf}{bdf}+\frac{aed}{bfd} = \frac{ac}{bd}+\frac{ae}{bf}$.
\csq, $+\frac{acf+ade}{bdf} = j\bg(\frac{acf+ade}{bdf}) = j\bg(\frac{ac}{bd} +
\frac{ae}{bf}) \nad*= j\bg(\frac{ac}{bd})+j\bg(\frac{ae}{bf})  =
j\bg(\frac{a}{b} \,\frac cd) + j\bg(\frac ab\,\frac ef) \nad*=
j\bg(\frac{a}{b})\cdot j\bg(\frac cd) + j\bg(\frac ab)\cdot j\bg(\frac ef)
\nde5.16 = \wh x\cdot \wh y+\wh x\cdot\wh z$. In $\nad*=$ we used the fact
that $j$~is a ring-\hm sm.

For the proof of the remaining cases we use the \fw\ \et ies:
\beq5.26
\ga
(\wh x\cdot\wh y)\Inv= (\wh x{}\Inv)\cdot\wh y = \wh x\cdot \wh y{}\Inv,
\q \wh x,\wh y\in\wh X\sms0,\\
\hbox{where}\q \wh x+\wh x{}\Inv =\wh x{}\Inv+\wh x=0, \q \wh x\in\wh X\sms0.
\ega
\e
Note that \er{5.26} is an \im con\sq\ of \er{5.17}, \er{5.19} and $(\wh
x{}\Inv)\Inv=\wh x$,\break ${\wh x\in\wh X\sms0}$.

\Wanp prove the \fw\ cases.

\ti{Case} (ii): $\wh x,\wh y\in\wh X_+$, $\wh z\in\wh X_-$ with $|\wh z|
< |\wh y|$.

Let $\wh x:=+q$, $\wh y:=+r$, $\wh z:=-s$ with $q,r,s \in\Q_{>0}$ and $s<r$.
Let $p\in\Q_{>0}$ be \st
\beq5.27
r=s+p.
\e
\E\Tf $\wh y=j(r)=j(s+p)=j(s)+j(p)$, and $\wh z=-s\nde5.16 = j(s)\Inv$. Thus
$\wh y+\wh z=j(s)+j(p)+(j(s))\Inv=j(p)$, and $\wh x\cdot(\wh y+\wh z)=
\wh x\cdot j(p)$, \E\oh $\wh x\cdot\wh y= \wh x\cdot(j(s)+j(p)) = \wh x\cdot
j(s)+\wh x\cdot j(p)$ by case~(i), and $\wh x\cdot\wh z=\wh x\cdot j(s)\Inv\nde5.26 =
(\wh x\cdot j(s))\Inv$. Hence $\wh x\cdot \wh y+\wh x\cdot \wh z = \wh x
\cdot j(s)+\wh x\cdot j(p)+(\wh x\cdot j(s))\Inv = \wh x\cdot j(p) =
\wh x\cdot(\wh y+\wh z)$.

\ti{Case} (iii): $\wh x,\wh y\in\wh X_+$, $\wh z\in\wh X_-$ and $|\wh z|
=|\wh y|$. We have $\wh y+\wh z=0$, hence $\wh z=\wh y{}\Inv$. \E\Tf
$\wh x\cdot(\wh y+\wh z)=\wh x\cdot0 \nde5.18 = 0$. \E\oh $\wh x\cdot\wh z =
\wh x\cdot(\wh y{}\Inv)\nde5.26 = (\wh x\cdot\wh y)\Inv$, hence $\wh x\cdot\wh y
+\wh x\cdot\wh z= \wh x\cdot \wh y+(\wh x\cdot \wh y)\Inv =0=\wh x\cdot(\wh y+\wh z)$.

\ti{Case} (iv): $\wh x, \wh y\in\wh X_+$, $\wh z\in\wh X_-$ and $|\wh y|<|\wh z|$.
Let $\wh y=j(r)$, $\wh z=j(s)\Inv$ with $r,s\in\Q_{>0}$ and $r<s$. Let $p\in\Q
_{>0}$ be \st $s=r+p$. Then $\wh y+\wh z=j(r)+j(r+p)\Inv=j(r)+(j(r)+j(p))\Inv
\nad*= j(r)+j(r)\Inv+j(p)\Inv=j(p)\Inv$. In $\nad*=$ we used \E\Pr\
\rf{p3.26}\,(i). Hence $\wh x\cdot(\wh y+\wh z)=\wh x\cdot j(p)\Inv$. \E\oh
$\wh x\cdot \wh y+\wh x\cdot \wh z= \wh x\cdot j(r)+\wh x\cdot(j(r+p))\Inv
\nde5.26 = \wh x\cdot j(r) + (\wh x\cdot j(r+p))\Inv = \wh x\cdot j(r) +
(\wh x\cdot (j(r)+j(p)))\Inv \nad{\rm case\ (i)} = \wh x\cdot j(r)+ (\wh x
\cdot j(r)+\wh x\cdot j(p))\Inv = \wh x\cdot j(r) + (\wh x\cdot j(r))\Inv +
(\wh x\cdot j(p))\Inv = (\wh x\cdot j(p))\Inv \nde5.26 = \wh x\cdot j(p)\Inv =
\wh x\cdot(\wh y+\wh z)$.

\ti{Case} (v): $\wh x\in\wh X_+$, $\wh y\in\wh X_-$, $\wh z\in\wh X_+$.

We have $\wh x\cdot(\wh y+\wh z)=\wh x\cdot(\wh z+\wh y)$. Interchanging $\wh z$
and~$\wh y$, we are either in case (ii), or~(iii), or~(iv). \E\Tf $\wh x\cdot
(\wh z+\wh y)=\wh x\cdot\wh z + \wh x\cdot \wh y$. Hence by \cmt ity $\wh x\cdot
(\wh y+\wh z)=\wh x\cdot\wh y+\wh x\cdot \wh z$.

\ti{Case} (vi): $\wh x\in\wh X_+$, $\wh y,\wh z\in \wh X_-$.

We have $\wh x\cdot(\wh y+\wh z) = (\wh x\cdot(\wh y+\wh z)\Inv)\Inv
= (\wh x\cdot(\wh y{}\Inv+\wh z{}\Inv))\Inv \nad{\rm case\ (i)}= (\wh x\cdot
\wh y{}\Inv + \wh x\cdot\wh z{}\Inv)\Inv = \wh x\cdot \wh y+ \wh x
\cdot\wh z$.

Observe that we proved \er{1.2n} for $\wh x\in\wh X_+$ and $\wh y,\wh z\in
\wh X\sms0$. It remains to show that if $\wh x\in\wh X_-$,
$\wh y,\wh z\in\wh X\sms0$, then \er{1.2n} holds. In this case
$\wh x\cdot(\wh y+\wh z) = (\wh x{}\Inv\cdot(\wh y+\wh z))\Inv =
(\wh x{}\Inv\cdot\wh y+\wh x{}\Inv\cdot\wh z)\Inv = \wh x\cdot\wh y
+\wh x\cdot\wh z$.

This completes the proof of \er{1.2n}, hence also the proof of Theorem \rf{t5.6}.
\endproof

\bdf5.7
The field $(\wh X,+,\cdot,0,+1)$ introduced in Theorem \rf{t5.6} will be called
the field of \ti{\ra\ \nm s}, and will be denoted by~$\Q$. Usually the \el s of
$\Q_{\ge0}$ ($=\wh X_+\cup\{0\}$) are simply denoted by~$q$ instead of~$+q$.
Then the map~$j$ is simply the inclusion map $j(q):=q$, $q\in\Q_{\ge0}$.
\edf

We now consider another approach to the field of \ra\ \nm s. The first step
of this approach consists of introducing on the additive group of \ig s $(\Z,+,0)$
a~\mlc~$\cdot$ which makes $(\Z,+,\cdot,0,1)$ a~ring. Recall that we
introduced a \mlc\ on the monoid $(\N,+,0)$ by setting
\beq5.28
m\cdot n:= m\dpl n, \q m,n\in\N,
\e
(see \E\df\ and Notation \rfa2{d2.3}). We proved that $(\N,\cdot,1)$ is an
\am\ (see \era2{2.11}--\era2{2.13}) and that $\Na:=\N\sms0$ is a \sbm\ of
$(\N,\cdot,1)$ \st $(\Na,\cdot,1)$ is a \PM\ (see \E\Pr\ \rfa2{p2.7}), \Ip
a~\Cm. In a similar fashion we define a \mlc\ on~$\Z$ by setting
\beq5.29
z\cdot w:= z\dpl w, \q z,w\in\Z,
\e
where $z\dpl w$ denotes the $z$-fold \IT\ of~$w$ in $(\Z,+,0)$ introduced in
\er{3.142}.

\bth5.7
Let $(\Z,+,0)$ denote the \ag\ of \ig s introduced in \E\Pr\ \rf{p3.90}, and
let $\cdot$ denote the binary \op\ on~$\Z$ introduced in \er{5.29}. Then

\hph i,ii, $(\Z,\cdot,1)$ is an \am,

\hph ii,i, $\Z\sms0$ is a \sbm\ of $(\Z,\cdot,1)$ \st $(\Z\sms0,\cdot,1)$ is
a \Cm,

\hph iii,, $(\Z,+,\cdot,0,1)$ is a $($\cmt e$)$ ring $($with unity$)$.

\eth

\proof \

(i), (ii): \ti{$0\cdot w=0$ for $w\in\Z$}: follows from \er{3.142}\,SI0 in $(\Z,+,0)$.

\ti{$z\cdot0 =0$ for $z\in\Z$}: follows from \er{3.142}\,SI4 in $(\Z,+,0)$.

In what follows we use the notation
\beq5.30
\Z^\t:= \Z\sms0.
\e
Below $z\Inv$ denotes the inverse of~$z$ in $(\Z,+,0)$, that is, $z+z\Inv=z\Inv+z=0$.
We recall that the map $\INV:\Z\to\Z$ is an auto\mf\ of $(\Z,+,0)$ and $\INV
\circ\INV = I_\Z$ by \E\Pr\ \rf{p3.26}.

We prove
\beq5.31
z\Inv\cdot w = (z\cdot w)\Inv, \q z,w\in\Z^\t.
\e
Note that \er{5.31} is a con\sq\ of \er{3.142}\,SI6 and SI7. Since the proof of
\E\Pr\ \rf{p3.142} was  left as an exercise, we prove \er{5.31} here. The same
remark applies to \er{5.32}.

If $z:=m\in\Na$, then $z\Inv=-m$ and $z\Inv\cdot w\nad{\er{5.29},\er{3.141}}=
m\dpl w\Inv\nda21.45 = (m\dpl w)\Inv\nde5.29 = (z\cdot w)\Inv$. If $z:=-m$,
$m\in\Na$, then $z\Inv\cdot w = m\cdot w$. But $(-m)\cdot w \nde5.29 =
(-m)\dpl w\nde3.141 = m\dpl w\Inv \nde5.29 = m\cdot w\Inv$. We have $m\cdot w
+(-m)\cdot w\nde3.141 = m\cdot w+m\cdot w\Inv \nad{\era2{2.3}\,\rm I5}=
m\cdot (w+w\Inv)=m\cdot 0=0$. Hence $m\cdot w= ((-m)\cdot w)\Inv$, that is,
$z\Inv \cdot w=(z\cdot w)\Inv$.

Now we want to show
\beq5.32
z\cdot w\Inv= (z\cdot w)\Inv, \q z,w\in\Z^\t.
\e
If $z:=m\in\Na$, then $z\cdot w\Inv\nde5.29 = m\dpl w\Inv$. But $m\dpl w\Inv
+ m\dpl w\nad{\era2{2.3}\,\rm I5} = m\dpl (w\Inv+w) = m\dpl 0=0$. Hence $m\dpl
w\Inv =(m\dpl w)\Inv \nde5.29 = (z\cdot w)\Inv$. \csq, $z\cdot w\Inv=(z\cdot
w)\Inv$.

If $z:=-m$, $m\in\Na$, then $(-m)\cdot w\Inv\nde5.29 = (-m)\dpl w\Inv
\nde3.142 = m\dpl(w\Inv)\Inv = m\dpl w = z\Inv\dpl w\nde5.29 =
z\Inv\cdot w \nde5.31 = (z\cdot w)\Inv$.

Next we prove
\beq5.33
z\cdot w = w\cdot z, \q z,w\in\Z.
\e
Since $0\cdot w=0=z\cdot0$, it is \sft\ to prove \er{5.33} for $z,w\in\Z^\t$.
Let $m,n\in\Na$. Then $m\cdot n\nda22.12 = n\cdot m$. \Mo $(-m)\cdot
n\nde5.31 = (m\cdot n)\Inv = (n\cdot m)\Inv \nde5.32 = n\cdot(-m)$.
Similarly, $m\cdot(-n)\nde5.32 = (m\cdot n)\Inv =(n\cdot m)\Inv \nde5.31 =
(-n)\cdot m$. Finally, $(-m)\cdot(-n)\nde5.31 = (m\cdot(-n))\Inv \nde5.32 =
((m\cdot n)\Inv)\Inv = m\cdot n=n\cdot m=(-n)\cdot(-m)$.

As a con\sq\ we have
\beq5.34
z\cdot 1=1\cdot z=z, \q z\in\Z.
\e
Indeed, $z\cdot 1\nde5.33 = 1\cdot z\nde5.29  = 1\dpl z\nad{\er{3.142}\rm\,
SI1}=z$, $z\in\Z$. In order to show that $(\Z,\cdot,1)$ is an \am, it remains
to show that the \mlc~$\cdot$ is \ti{\asc e}. To this end we use the approach
used in the proof of Theorem~\rf{t5.6}. We define a map $\sgn:\Z^\t\to\{+,-\}$
by setting\glossary{$\sgn(z)$}
\beq5.35
\sgn(z):=\bca
+ &\hbox{if }z\in\Z_{>0},\\
- &\hbox{if }z\in\Z_{<0},\eca
\e
and a map $|\ |:\Z^\t \to\Na$ by setting\glossary{$|z|$}
\beq5.36
|z|:=\bca
z &\hbox{if }z\in\Z_{>0},\\
-z &\hbox{if }z\in\Z_{<0}.\eca
\e
\If from \E\Pr\ \rf{p3.90} that the map $f:\{+,-\}\t\Na \to\Z^\t$ defined by
\beq5.37
f((\a,a)):=\a a, \q \a\in\{+,-\},\ \a\in\Na,
\e
is a bi\jn\ and that
\beq5.38
f\Inv(z)=(\sgn(z),|z|), \q z\in\Z^\t.
\e
If $\cdot$ denotes the binary \op\ defined in \er{5.22}, then
$(\{+,-\},\cdot,+)$ is an \ag\ as observed in the proof of Theorem \rf{t5.6}.
\E\oh $(\Na,\cdot,1)$ is a \Cm. \If from \er{5.35}, \er{5.36}, \er{3.141},
\er{5.29}, \er{5.31}, \er{5.32}, \er{3.4} that
\bga5.39
\sgn(z\cdot w) = \sgn(z)\cdot\sgn(w),\\
|z\cdot w|=|z|\cdot|w|, \lb{5.40}\\
\sgn(1)={+},\q |1|=1. \lb{5.41}
\e
If we define the binary \op\ $\qu$ on $\{+,-\}\t\Na$ by setting
\beq5.42
(\a,a)\qu(\b,b) := (\a\cdot\b,ab), \q \a,\b\in\{+,-\},\ a,b\in\Na,
\e
then we find by \er{5.39}, \er{5.40}:
\beq5.43
z\cdot w= f\bg({f\Inv(z)\qu f\Inv(w)}), \q z,w\in\Z^\t.
\e
\Mo by \er{5.41}
\beq5.44
1=f((+,1)).
\e
\If from Example \rfa2{xa1.5}\,(v) that $(\Z^\t,\cdot,1)$ is a monoid, and
from \er{5.43}, \er{5.44} that $f$~is a monoid-\is sm from the monoid
$(\{+,-\}, \cdot,+)\t(\Na,\cdot,1)$ onto the monoid $(\Z^\t,\cdot,1)$.

\E\csq,
\beq5.45
(x\cdot y)\cdot z = x\cdot(y\cdot z) \qh{\fa }x,y,z\in\Z^\t.
\e
Since $0\cdot z=z\cdot 0=0$ \fa $z\in\Z$, \er{5.45} holds \fa $x,y,z\in\Z$.
\If that $(\Z,\cdot,1)$ is an \am. Since $(\Z^\t,\cdot,1)$ is a monoid,
$\Z^\t$~is a \sbm\ of $(\Z,\cdot,1)$. \Mo as a group, $(\{+,-\},\cdot,+)$ is
a \Cm\ as well as $(\Na,\cdot,1)$. Hence the direct product $\{+,-\}\t \Na$
also is a \Cm\ by Example \rfa2{xa1.5}\,(iv). Since $f$~is a monoid-\is sm,
$(\Z^\t,\cdot,1)$ is a \Cm. This completes the proof of (i) and~(ii).

(iii): Clearly \er{1.1n} holds. In view of Exercise \rf{ex1.2n} it is \sft\
to prove \er{1.2n}. By~\er{1.1n} and \era2{1.7}, it suffices to prove
\er{1.2n} for $x,y,z\in\Z^\t$. If $x\in\Z_{>0}$ then \er{1.2n} follows from
\era2{2.3}\,I5 with $(M,\qu,e):=(\Z,+,0)$. If $x\in\Z_{<0}$, then
$x\Inv\in\Z_{>0}$. Hence $x\Inv\cdot(y+z)=x\Inv\cdot y+x\Inv\cdot
z$, for $y,z\in\Z$. \E\Tf $x\cdot(y+z)\nad{\er{5.31},\er{3.4}} = (x\Inv\cdot(y+z))\Inv = (x\Inv\cdot y+x\Inv \cdot
z)\Inv =(x\Inv \cdot y)\Inv +(x\Inv \cdot z)\Inv \nad{\er{5.31},\er{3.4}} = x\cdot y+x\cdot
z$. This completes the proof of~\er{1.2n}.
Since $(\Z,+,0)$ is an \ag, $(\Z,+,\cdot,0,1)$ is a \ti{ring} in view of
\E\df\ \rf{d3.11}\,(i).
\endproof

\bex5.8
Let $n\in\Na\sms1$ and let $(M,\circ,(0,1))$ denote the direct product of the
monoids $(\N_n,+_n,0)$ and $(\Na,\cdot,1)$ (see Examples \rfa2{xa1.5}\,(iv),
Remark \rfa2{r1.6}\,(i)). Set $\wh M:=M\cup\wh0$. Define a binary
\op~$\mhc$ on~$\wh M$ by setting $a\mhc\wh0:=\wh0$, $\wh0\mhc
\wh a:=\wh0$ \fa $a\in\wh M$, and $a\mhc b:=a\circ b$ \fa $a,b\in M$.

\hph i,ii, Show that $(\wh M,\mhc,(0,1))$ is an \am.

\hph ii,i, Show that if $n:=2$, then \te s a binary \op~$\hpl$ on~$\wh M$
\st $(\wh M,\hpl,\mhc,\wh0,(0,1))$ is a \cmt e ring with unity.

\hph iii,, Is there an $n\in\Na\sms{1,2}$  \st a binary \op~$\hpl$ on~$\wh M$
exists which makes $(\wh M,\hpl,\mhc,\wh0,(0,1))$ a \cmt e ring?
\eex

Before embedding the ring $\Z$ into a field we consider some results
involving the ring~$\Z$. We first \es\ the analogue of \era2{2.3}\,I3 for
$z$-fold \IT s of an \el\ of a group.

\setbox9=\hbox{SI0}
\bpr5.9
Let $(G,\qu,e)$ be an \ag, let $z,w\in\Z$ and let $x\in G$. Then
\beq5.46
\lca SI3 {(z\cdot w)\dqu x=z\dqu(w\dqu x),}
\e
where $a\dqu y$ with $a\in\Z$ and $y\in G$ denotes the $a$-fold \IT\ of~$y$
in $(G,\qu,e)$ defined in \er{3.141}.
\epr

\proof \

\ti{$z$ or $w$ equals zero}: If $z=0$ or $w=0$ then $z\cdot w\nde1.1n = 0$ and
$(z\cdot w)\dqu x=e$ by \er{3.142}\,SI0. If $z=0$ then $z\dqu (w\dqu x)=e$. If
$w=0$ then $w\dqu x=e$ and $z\dqu(w\cdot x)=e$ by \er{3.142}\,SI4. Let
$m,n\in\Na$.

\ti{$z:=m$, $w:=n$}: follows from \era2{2.3}\,I3.

\ti{$z:=m$, $w:=-n$}: $z\cdot w\nde5.32 = (m\cdot n)\Inv$, hence $(z\cdot
w)\dqu x\nde3.141 = (m\cdot n)\dqu x\Inv$, \E\oh $z\dqu (w\dqu x)\nde3.141 =
m\dqu (n\dqu x\Inv) \nad{\rm I3}= (m\cdot n)\dqu x\Inv$.

\ti{$z:=-m$, $w:=n$}: $z\cdot w \nde5.31 = (m\cdot n)\Inv$, hence $(z\cdot w)
\dqu x\nde3.141 = (m\cdot n)\dqu x\Inv$. \E\oh $z\dqu (w\dqu x) = z\dqu (n
\dqu x) = m\dqu (n\dqu x)\Inv\nda21.45 = m\cdot (n\cdot x\Inv)\nad{\rm I3}=
(m\cdot n)\dqu x\Inv$.

\ti{$z:=-m$, $w:=-n$}: $z\cdot w\nad{\er{5.31},\er{5.32}}= m\cdot n$, hence
$(z\cdot w)\dqu x=(m\cdot n)\dqu x$. \E\oh $z\dqu (w\dqu x)= z\dqu
((-n)\dqu x) \nde3.141 = z\dqu (n\dqu x\Inv) = (-m)\dqu (n\dqu
x\Inv)\nde3.141 = \break m\dqu ((n\dqu x\Inv)\Inv) \nda21.45 =
m\dqu(n\dqu(x\Inv)\Inv) =  m\dqu (n\dqu x)\nad{\rm I3}= (m\cdot n)\dqu x$.
\endproof

We next consider \rl s between the ring $(\Z,+,\cdot,0,1)$ and the rings
$(\N_n,+_n,\cdot_n,0,1)$, $n\ge2$. To this end we shall use the \fw\ lemma
which is an \ext\ to~$\Z$ of the ``division algorithm'' introduced in Theorem
\rfa2{t1.38} with $(\wt E,\wt e,\wt S):=(\N,0,S)$.

\blm5.10
Let $a\in\Z\sms0$ and $b\in\Z$. Then \te s one and only one pair $(q,r)\in
\Z\t\N$ \sf ying
\beq5.47
b=q\cdot a+r
\e
and
\beq5.48
0\le r<|a|.
\e
\elm

\brm5.11
We claim that if $a\in\N\sms0$ and $b\in\N$ in Lemma \rf{l5.10}, then $q\in\N$,
which shows that the \ex\ part of Lemma \rf{l5.10} implies the \ex\ part of
Theorem \rfa2{t1.38}.

Indeed, suppose for \cd ion that $q\in\Z_{<0}$, then $|q|\ge1$, $q=-|q|$ and
$q\cdot a=(-|q|)\cdot a\nde5.31 = -(|q|\cdot a)$. Hence $-(|q|\cdot a)+r=b$.
Adding $|q|\cdot a$ to both terms of the \et y, we find $r=|q|\cdot a+b$. Note that
$r,|q|,a$ and~$b$ belong to~$\N$. Since $b\ge0$, we have $|q|\cdot a=
|q|\cdot a+0\le |q|\cdot a+b$, by~\era2{1.64}. Hence $|q|\cdot a\le r$.
Since $a\in\Na$ and $1\le|q|$, we have $a=1\cdot a\le |q|\cdot a$,
by~\era2{2.19}. Hence $a\le|q|\cdot a\le r$, so $a\le r$. A~\cd ion,
since $r<|a|=a$. This proves the claim.
\erm

We now show that the \uq\ part of Lemma \rf{l5.10} implies the \uq\ part
of Theorem \rfa2{t1.38}. Indeed, let $a\in\N\sms0$, $b\in\N$, $(q_i,r_i)\in
\N\t\N$, $i=1,2$, \sf y \era2{1.59}, \era2{1.60}. Then they clearly \sf y
\er{5.47}, \er{5.48}. Hence $(q_1,r_1)=(q_2,r_2)$ by Lemma \rf{l5.10}. It
follows that Lemma \rf{l5.10} is a \gn\ of Theorem \rfa2{t1.38}.

\proof[Proof of Lemma \rf{l5.10}]\

\ti{\E\ex}: The case $a\in\Z_{>0}$ and $b\in\Z_{\ge0}$ is a direct con\sq\
of Theorem \rfa2{t1.38}.

Suppose $a\in\Z_{>0}$ and $b\in\Z_{<0}$. Set $b':=-b\in\Z_{>0}$. Then by
Theorem \rfa2{t1.38}, \te s $(q',r')\in\N\t\N$ \st $b'=q'\cdot a+r'$ and
$0\le r'<a$. \E\Tf $b\nad*= -b'=-(q'\cdot a+r')\nad*= -(q'\cdot a)
+(-r')\nde5.31 = (-q')\cdot a+(-r')$. In $\nad*=$ we used \E\Pr\ \rf{p3.26}.
If $r'=0$, then set $q:=-q'$ and
$r:=-r'=0$. Then $b=q\cdot a+0$, $0=r<a$ and we are done. If $r'>0$, then
$0<r'<a$. \E\te s $p\in\Na$ \st $a=r'+p$. Since $r',p\in\Na$, we have
$0<a-r'=p<a$. So $r:=a-r'$ \sf ies $0\le r<a$. Since $r'\nda21.5 =a-r$, we obtain
$b'=q'\cdot a+a+(-r)\nad{\era2{2.13},\er{1.3n}}=(q'+1)\cdot a+(-r)$.
Consequently, $b=-b' \nde3.137 = -((q'+1)\cdot a)+(-(-r))\nde5.31 = (-(q'+1))\cdot a+r$.
Setting $q:=-(q'+1)$ we obtain \er{5.47}, \er{5.48}.

Finally, we consider the case $a\in\Z_{<0}$. Set $a'=-a$. By what precedes,
\te s $(q',r')\in\Z\t\N$ \st $b=q'\cdot a'+r'$ with $r'\in[0,a')$. Setting
$q:=-q'$, and using $q'=-q$, $(-q)\cdot(-a)\nad{\er{5.31},\er{5.32}} = q\cdot a$,
$a'=-a =|a|$, we arrive at $b=q\cdot a+r'$ with $0\le r'<|a|$, which
is \er{5.47}, \er{5.48} with $r:=r'$.

\ti{\E\uq}: Let $a\in\Z\sms0$, $b\in\Z$, $(q_i,r_i)\in\Z\t\N$, $i=1,2$, be
\st $b=q_i\cdot a+r_i$, $0\le r_i<|a|$, $i=1,2$. Then $q_1\cdot a+r_1=q_2
\cdot a+r_2$, which implies $q_1\cdot a-q_2\cdot a=r_2-r_1$. Note that
$q_1\cdot a-q_2\cdot a=q_1\cdot a+(-(q_2\cdot a))\nde5.31 = q_1\cdot a
+(-q_2)\cdot a\nde 1.3n = (q_1+(-q_2))\cdot a=(q_1-q_2)\cdot a$. Taking
absolute values we obtain $|(q_1-q_2)\cdot a|=|r_2-r_1|$. \Mo $|(q_1-q_2)
\cdot a|\nde5.40 = |q_1-q_2|\cdot |a|$ and $r_1-r_2\le r_1<|a|$,
$r_2-r_1\le r_2<|a|$, since $r_1,r_2\ge0$. \E\Tf $|r_1-r_2|<|a|$. It follows
that $|q_1-q_2|\cdot |a|<|a|=1\cdot |a|$. Since $|a|\in\Na$, we infer from
\era2{2.20} that $|q_1-q_2|<1$. The only non\ng\ \ig\ less than~1 is~0, \Tf
$|q_1-q_2|=0$, hence $q_1-q_2=0$ and $q_1=q_2$. Finally, from $q_1\cdot a
+r_1=q_1\cdot a+r_2$ we obtain $r_1=r_2$ by \cnc ity of the \ad.
\endproof

In what follows Lemma \rf{l3.83} will play an important role as it did in
particular in the proof of Theorem \rf{t3.84}. Observe that if $(M,\qu,e)$ is
a group, then so is the \qt\ monoid $M/R$ introduced in Lemma \rf{l3.83}.
Indeed, it suffices to show that if $\wh x\in M/R$,
then \te s $\wh y\in M/R$ \st $\wh x \hqu \wh y=[e]$. Let $x\in M$ be \st
$\wh x=[x]$. Since $M$~is a group, \te s $y\in M$ \st $x\qu y=e$. Set $\wh
y:=[y]$. Then $\wh x\hqu \wh y=[x\qu y]=[e]$.

Now let $(G,\qu,e)$ be a group not necessarily abelian, and let $H$~be a
subgroup of~$G$. Then the \rl~$\sim$ on~$G$ defined by
\beq5.50
x\sim y \qh{if } x\qu y\Inv\in H \hbox{ \fa}x,y\in G,
\e
is an \ti{\ev ce} \rl.

Recall that if $H$ is a subgroup of $G$, then $e\in H$, $x\Inv\in H$ \fa
$x\in H$, and $x\qu y\in H$ \fa $x,y\in H$. Thus \fa $x,y,z\in G$, $x\sim x$
since $x\qu x\Inv=e\in H$, moreover, $x\sim y$ implies $y\sim x$ since $x\qu
y\Inv\in H$ implies $(x\qu y\Inv)\Inv \in H$, thus $y\qu x\Inv\nad{\er{3.4},
\er{3.5}}=\break (x\qu y\Inv)\Inv\in H$; finally, $x\sim y$ and $y\sim z$ implies
$x\qu y\Inv\in H$, $y\qu z\Inv\in H$, hence $x\qu z\Inv=x\qu e\qu z\Inv=
x\qu(y\Inv\qu y)\qu z\Inv\nda21.36 = (x\qu y\Inv)\qu(y\qu z\Inv)\in H$.

It is important to note that the \rl~$\sim$ defined in \er{5.50} and the
binary \op~$\qu$ in~$G$ are not necessarily compatible in the sense of \E\df\
\rf{d3.82}. Indeed, if $x\sim x'$, $y\sim y'$, $x,y,x',y'\in G$, then
$(x\qu y)\qu (x'\qu y')\Inv \nde3.5 = (x\qu y)\qu(y'{}\Inv\qu x'{}\Inv)
\nda21.36 =\break x\qu(y\qu y'{}\Inv)\qu x'{}\Inv = x\qu h\qu x'{}\Inv$ \fs $h\in H$. If $G$~is
abelian, then $x\qu h\qu x'{}\Inv= h\qu x\qu x'{}\Inv=h\qu h'$ \fs\ $h'\in H$.
Since $h\qu h'\in H$, $x\qu h \qu x'{}\Inv\in H$. Thus if $G$ is \ti{abelian},
then $\sim$ and~$\qu$ are \ti{compatible}. More generally, if \fa $x\in G$
and $h\in H$, \te s $h''\in H$ \st $x\qu h=h''\qu x$, then $x\qu h\qu x'{}\Inv
=h''\qu x\qu x'{}\Inv = h''\qu h'\in H$ \fs $h'\in H$.
In this case the subgroup $H$ is called \ti{normal\/}. \E\te\ not normal
subgroups \st $\sim$ and~$\qu$ are not compatible. For the sake of completeness
we give below such an example. Summarizing, we have

\blm5.11
Let $(G,\qu,e)$ be a $($not necessarily abelian\/$)$ group and let $H$~be
a subgroup of~$G$. Then the \rl~$\sim$ defined in \er{4.50} is an \emph{\ev ce
\rl}. If $G$ is abelian or, more generally, if $H$ \sf ies \cn
\beq5.51
\hbox{\fe $x\in G$ and every $h\in H$, \te s $h'\in H$ \st }x\qu h=h'\qu x,
\e
then the \rl\ $\sim$ and the binary \op\ $\qu$ are \emph{compatible} in the
sense of \E\df\ \rf{d3.82}. In that case, the \qt\ group $G/{\sim}$ introduced
in Lemma \rf{l3.83} is also denoted by~$G/H$.
\elm

\bdf5.12
A subgroup $H$ of a group $(G,\qu,e)$ \sf ying \er{5.51} is called a
\ti{normal\/} subgroup.
\edf

\bxa5.13
Let $(S_3,\circ,\id)$ denote the group of bi\jn s (\Pm s) of $\{1,2,3\}$ where
$\circ$~denotes the \cm\ of two bi\jn s. By \era2{4.50} $\#(S_3)=3!=6$. Let
$\si_{1,2}$ denote the \el\ of~$S_3$ defined by $\si_{1,2}(1):=2$, $\si_{1,2}(2):=1$
and $\si_{1,2}(3):=3$. Then $\si_{1,2}\circ \si_{1,2}=\id$, hence $H:=\{\id,\si_{1,2}\}$
is a subgroup of~$S_3$ (\is c to $(\N_2,+_2,0)$). In view of \er{5.50},
$f,g\in S_3$ \sf y $f\sim g$ iff either $f=g$ or $f=\si_{1,2}\circ g$.
Let $f:=\si_{1,3}\in S_3$ where $\si_{1,3}(1):=3$, $\si_{1,3}(2):=2$,
$\si_{1,3}(3):=1$. Note that $\si_{1,3}\circ\si_{1,3}=\id$.

Set $g:=\si_{1,2}\circ f$, $f':=\si_{1,2}\circ f$, $g':=\si_{1,2}\circ g$. Then
$f'\sim f$ and $g'\sim g$. Since $\sim$ is \sy ic, we have $f\sim f'$, $g\sim
g'$. Thus $f\circ g= f\circ(\si_{1,2}\circ f)=(f\circ\si_{1,2})\circ f$. \Mo
$f'\circ g'=(\si_{1,2}\circ f)\circ(\si_{1,2}\circ g)= (\si_{1,2}\circ f)\circ
(\si_{1,2}\circ(\si_{1,2}\circ f))= (\si_{1,2}\circ f)\circ (\si_{1,2}\circ
\si_{1,2})\circ f = (\si_{1,2}\circ f)\circ f$. We claim $f\circ g\not\sim
f'\circ g'$, that is, neither $f\circ g=f'\circ g'$ nor $f\circ g=\si_{1,2}\circ
(f'\circ g')$. We first show that $f\circ g\ne f'\circ g'$. Indeed, suppose
for \cd ion that $f\circ g=f'\circ g'$, then $f\circ\si_{1,2}=\si_{1,2}\circ
f$ by \cnc ity. However, $\si_{1,3}\circ\si_{1,2} \ne \si_{1,2}\circ\si_{1,3}$
since $(\si_{1,3}\circ\si_{1,2})(3) = \si_{1,3}(\si_{1,2}(3)) =
\si_{1,3}(3)=1$, and $(\si_{1,2}\circ\si_{1,3})(3)=\si_{1,2}(1)=2\ne1$.
A~\cd ion. Suppose now for \cd ion that $f\circ g=\si_{1,2}\circ(f'\circ g')$.
Then by \asc ity $f\circ\si_{1,2}\circ f=\si_{1,2}\circ\si_{1,2}\circ f\circ f
=f\circ f$ since $\si_{1,2}\circ\si_{1,2}=\id$. By \cnc ity we obtain
$f\circ \si_{1,2}=f$. Since $f=f\circ\id$, we infer $\si_{1,2}=\id$. A~\cd ion,
since $\si_{1,2}(1)=2\ne1=\id(1)$.

Therefore, the \rl\ $\sim$ and the binary \op~$\circ$ are not compatible.
Hence $H$ is \ti{not\/} normal in view of Lemma \rf{l5.11}.
\exa

We now apply Lemma \rf{l5.11} to the \ag~$(\Z,+,0)$. We have

\blm5.14
Let $H$ be a subgroup of $(\Z,+,0)$. Then \te s \ooo $n\in\N$ \st $H=n\Z$
where $n\Z$ is defined by\glossary{$n\Z$}
\beq5.52
n\Z := \{n\cdot z: z\in \Z\}= \{z\dpl n: z\in\Z\}.
\e
\elm

\bex5.15
Prove Lemma \rf{l5.14}.
\eex

Note that $0\Z=\{0\}$ the trivial subgroup of $(\Z,+,0)$, and $1\Z=\Z$. It
turns out that the \qt\ group $\Z/n\Z$, $n\ge2$, is \is c to the \fcg\glossary{$\Z/n\Z$}
$(\N_n,+_n,0)$. More precisely, we have

\bpr5.16
Let $n\in\N$, $n\ge2$, and let $n\Z$ denote the subgroup of $(\Z,+,0)$
defined in \er{5.52}. Let $(\Z/n\Z,\hpl,\wh0)$ denote the \qt\ group of~$\Z$
by $n\Z$ introduced in Lemma \rf{l3.83}. Let $(\N_n,+_n,0)$ denote the
cyclic group introduced in \E\Pr\ \rf{p1.8} $($see also \E\Pr\
\rf{p3.10}\,{\rm(i))}. Let $f:\N_n\to \Z/n\Z$ be defined by
\beq3.222
f(k):=[k],\ k\in\N_n,\qh{where $[k]$ denotes the \ev ce class
containing $k$.}
\e
Then $f$ is an \emph{\is sm} from $\N_n$ onto $\Z/n\Z$.
\epr

\proof
Since the group $(\Z,+,0)$ is abelian, the subgroup $n\Z$, $n\ge2$, introduced
in Lemma \rf{l5.14} \sf ies \er{5.51}, where $(G,\qu,e):=(\Z,+,0)$ and $H:=
n\Z$. \E\Tf the \ev ce \rl~$\sim$ on~$\Z$ defined by \er{5.50} and the \ad~$+$
in~$\Z$ are compatible in view of Lemma \rf{l5.11}. By Lemma \rf{l3.83} with
$(M,\qu,e):=(\Z,+,0)$ and $R={\sim}$, the \qt\ monoid $(\Z/{\sim},\hpl,[0])$,
also denoted by $(\Z/n\Z,\hpl,\wh0)$, is abelian. It is a group in view of
the discussion preceding \er{5.50}, since $(\Z,+,0)$ is a group. \Mo the
map~$f$ defined in \er{3.222} is the \rt ion to $\N_n:=[0,n)$ of the map
$\pi:\Z\to \Z/n\Z$ introduced in Lemma \rf{l3.83}. \E\Ip $f(0)=[0]=\pi(0)
=\wh0$ by Lemma \rf{l3.83}\,(iii). We now show that $f:(\N_n,+_n,0) \to
(\Z/n\Z,\hpl,\wh0)$ is a monoid-\is sm.

``$f(x+_ny)=f(x)\hpl f(y)$, $x,y\in\N_n$'': Let $x,y\in\N_n$. Then $f(x+_ny)
=\pi(x+_ny)$. From \er{1.8}, \er{1.9} we obtain $x+y=q\cdot n+(x+_ny)$ where
$+$ (resp.~$\cdot$) denotes the \ad\ (resp. \mlc) in~$\N$, hence also in~$\Z$.
Thus $x+_ny = x+y+(q\cdot n)\Inv$ where $(q\cdot n)\Inv$ is the inverse of
$q\cdot n$ in $(\Z,+,0)$. By Lemma \rf{l3.83}\,(ii) we have $\pi(x+y+
(q\cdot n)\Inv) = \pi(x+y)\hpl \pi((q\cdot n)\Inv)$. Note that
$(q\cdot n)\Inv \nde5.31 = q\Inv\cdot n
\nde5.33 = n\cdot q\Inv\nde5.52 \in n\Z$. Hence $\pi((q\cdot n)\Inv) =
[(q\cdot n)\Inv]$ with $(q\cdot n)\Inv \in n\Z$. Observe that $n\Z\nda21.7 =
\{z\in\Z:z+0\Inv\in n\Z\}=
\{z\in\Z: z\sim0\} = \{z\in\Z: 0\sim z\} \nda14.4 = [0]$. From $(q\cdot n)
\Inv\in[0]$, we infer $[(q\cdot n)\Inv]=[0]$ by Exercise
\rfa1{ex4.4}\,(i). Since $[0]=\pi(0)$, we obtain by Lemma \rf{l3.83}\,(ii)
$\pi(x+y+(q\cdot n)\Inv) = \pi(x+y)\hpl \pi((q\cdot n)\Inv) = \pi(x+y)\hpl
\pi(0)= \pi((x+y)+0) = \pi(x+y)= \pi(x)\hpl \pi(y)$. \csq, $f(x+_ny)=\pi(x+_ny)=
\pi(x+y)= \pi(x)\hpl \pi(y)=f(x)\hpl f(y)$.

``\ti{$f$ is in\jc}'': Since $f$ is a \hm sm it suffices to show that $\ker f
=\{0\}$ by \E\Pr\ \rf{p3.10}\,(vi). Let $x\in\N_n$ be \st $f(x)=\wh0$, that
is, $[x]=\wh0$. Since $\wh0=\pi(0)=[0]=n\Z$ by what precedes, $x\in n\Z$ in
view of Exercise \rfa1{ex4.4}\,(i). Thus $x=q\cdot n$ \fs $q\in\Z$, hence
$n|q|\nde5.36 = |n|\,|q|\nde5.40 = |nq|= |x|\nde5.36 = x<n$. By \era2{2.20},
we conclude $|q|\in\zo0,1 $. Since $|q|\in\N$, we obtain $q=0$, hence $x=0$.

``\ti{$f$ is sur\jc}'': Let $A\in\Z/n\Z$. By \df\ of an \ev ce \rl, \te s
$b\in\Z$ \st $A=[b]$. By Lemma \rf{l5.11} \te s $(q,r)\in\Z\t\N$ \st $b=q\cdot
n+r$, and $r\in\N_n$. Then $A=[b]=\pi(b) = \pi(q\cdot n+r)=\pi(q\cdot n)
\hpl \pi(r) = [0]\hpl \pi(r)=\pi(0)\hpl\pi(r) = \pi(0+r)=\pi(r)=f(r)$. Since
$A$~is arbitrary in $\Z/n\Z$, $f$~is sur\jc.
\endproof

Our next goal is to show that the \ev ce \rl\ R on~$\Z$ defined by
$x\mathbin{\rm R}y$ if $x-y\in n\Z$ is compatible with the \mlc\ $\cdot$
in the monoid $(\Z,\cdot,1)$. We consider a more general
situation. We assume that $(X,+,\cdot,0,1)$ is a ring and that $H$~is a
subgroup of the additive group $(X,+,0)$. We denote by~R the \ev ce \rl\ defined by
$x\mathbin{\rm R}y$ if $x+(-y)\in H$, $x,y\in X$. In view of Lemma \rf{l5.11}
with $(G,\qu,e):=(X,+,0)$, R~and~$+$ are compatible since $(X,+,0)$ is
abelian. We want to know under
which \cn s this \ev ce \rl\ \sf ies: $x\mathbin{\rm R}x'$, $y\mathbin{\rm
R}y'$  implies $x\cdot y\mathbin{\rm R}x'\cdot y'$, $x,y,x',y'\in X$. We have
$x'=x+h$, $y=y'+k$ \fs $h,k\in H$, since $(X,+,0)$ is an \ag.
Then $x'\cdot y'=(x+h)\cdot (y+k)
\nde1.2n = (x+h)\cdot y+(x+h)\cdot k\nde1.3n = x\cdot y+h\cdot y+x\cdot k
+h\cdot k$. The question is whether $h\cdot y+x\cdot k+h\cdot k\in H$. In
case $k=0$, this amounts to show that $h\cdot y\in H$. Since $h\in H$,
$y\in X$ are arbitrary, we need the \fw\ \cn:
\beq5.53
h\cdot x\in X\qh{\fa $h\in H$ and all} x\in X.
\e
Note that if \cn\ \er{5.53} is \sf ied, then in the case considered above
we have $h\cdot y,x\cdot k$ and $h\cdot k\in H$. Hence $h\cdot y+x\cdot
k+h\cdot k\in H$, since $H$~is a subgroup.

\If that if the subgroup $H$ of $(X,+,0)$ is nontrivial ($H\ne\{0\}$), proper
($H\ne X$), and \sf ies \cn~\er{5.53}, then the \op~$\cdot$ and the \rl~R
are compatible in the monoid $(X,\cdot,1)$. If we set
\bga5.54
[x]\hcd[y]:=[x\cdot y], \q x,y\in X,\\
\wh1:=[1], \lb{5.55}
\e
where $[x]$ denotes the R-\ev ce class containing $x\in X$, then, by Lemma
\rf{l3.83}, $(X/H,\hcd,\wh1)$ is an \am. Note that $\wh0\cdot[x] = [0\cdot x]
= [0]$, $x\in X$, and $[x]\hcd([y]\hpl[z]) = [x]\hcd[y+z] \nde5.54 = [x\cdot
(y+z)] \nde1.2n = [x\cdot y]\hpl [x\cdot z] \nde5.54 = [x]\hcd[y] \hpl [x]
\hcd [z]$, $x,y,z\in X$. \If that $(X/H,\hpl,\hcd,\wh0,\wh1)$ is a (\cmt e)
ring (with unity).

We have thus proved

\bpr5.17
Let $(X,+,\cdot,0,1)$ be a $($\cmt e\/$)$ ring $($with unity\/$)$ and let
$H$~be a nontrivial $(H\ne\{0\})$, proper $(H\ne X)$ subgroup of the group
$(X,+,0)$ \sf ying \er{5.53}. Then, the \rl~$\sim$ on~$X$ defined by
\beq5.54a
x\sim y \hbox{ if }x-y\in H, \q x,y\in X,
\e
is \emph{compatible} with $+$ in the group $(X,+,0)$ and with~$\cdot$ in the
monoid $(X,\cdot,1)$. \Mo $(X/H,\hpl,\hcd,\wh 0,\wh1)$ where $[x]\hpl[y]:=
[x+y]$, $[x]\hcd[y]:=[x\cdot y]$, $x,y\in X$, $\wh0:=[0]=H$, $\wh1:=[1]$, is
a $($\cmt e\/$)$ \emph{ring} $($with unity\/$)$, called \qt\ ring or factor
ring.
\epr

\bxa5.18
We give an example of a ring $X$ and of a subgroup~$H$ of $(X,+,0)$ \st \cn\
\er{5.53} is violated. Let $X$~be the product ring $\Z\t\Z$ (see \E\Pr\ \rf{p2.2}
and Corollary \rf{c3.37}). Let $H:=\{(z_1,z_2)\in\Z\t\Z: z_1=z_2\}$. Then
$H$~is an abelian subgroup of the additive group of $\Z\t\Z$. Take
$h:=(1,1)\in H$ and $z:=(1,0)\in\Z\t\Z$. Then $h\cdot z=(1,0)\notin H$.
\exa

\bco5.19
If $X$ in \E\Pr\ \rf{p5.17} is the ring $(\Z,+,\cdot,0,1)$, then \emph{every}
nontrivial proper subgroup of $(\Z,+,0)$ is of the form $n\Z$ \fs $n\ge2$ by
Lemma \rf{l5.14}. Clearly $n\Z$ \sf ies \er{5.53} since $(n\cdot w)\cdot
z=n\cdot (w\cdot z)$ \fa $w,z\in\Z$. \E\Tf $(\Z/n\Z,\hpl,\hcd,\wh 0,\wh1)$ is
a $($\cmt e\/$)$ ring $($with unity\/$)$.
\eco

\bex5.20
Show that the \rl\ $z\sim w$ in~$\Z$ defined by $z-w\in n\Z$ \fs $n\ge2$, is
equal to the \rl: \te s $q\in\Z$ \st $z-w=q\cdot n$.
\eex

\bnt5.21
The \rl\ $z-w\in n\Z$ is usually written $z\Eq w\mod n $.
\ent

The compatibility of $\sim$ with $+$ and $\cdot$ is identical to $z\Eq z'\mod
n $, $w\Eq w'\mod n $\break implies $z+ w\Eq z'+ w'\mod n $, $z\cdot w\Eq z'\cdot
w'\mod n $.

\Wanp give an \ext\ of \E\Pr\ \rf{p5.16}.

\bth5.22
Let $n\ge2$ and let $(\Z/n\Z,\hpl,\hcd,\wh0,\wh1)$ denote the ring introduced
in Corollary \rf{c5.19}. Let $(\N_n,+_n,\cdot_n,0,1)$ denote the ring
introduced in \E\Pr\ \rf{p3.15} $($see also \E\Pr\ \rf{p3.10}\,{\rm(i))}. Let
$f:\N_n\to \Z/n\Z$ be as in \E\Pr\ \rf{p5.16}. Then $f$ is a ring-\is sm from
the ring $\N_n$ onto the factor ring $\Z/n\Z$.
\eth

\brm5.23
A subgroup $H$ of the additive group of a ring $R$ \sf ying \cn\ \er{5.53} is
called an \ti{ideal\/} of~$R$ (see \E\df\ \rfa5{d4.38}).
\erm

\proof[Proof of Theorem \rf{t5.22}]
The map $f:(\N_n,{+_n},0) \to (\Z/n\Z, {\hpl},\wh 0)$ is bi\jc\ and a monoid-\hm
sm by \E\Pr\ \rf{p5.16}. In view of Lemma \rfa2{l1.8}, it is \sft\ to show
that $f:(\N_n,\cdot_n ,1)\to (\Z/n\Z,\hcd,\wh1)$ is a monoid-\hm sm. We have
$f(1):= [1]\nde5.55 = \wh1$. Let $x,y\in\N_n$. Then $f(x\cdot_n y)=[x\cdot_n
y]$. By \er{1.33}, \er{1.9}, \te s $q\in\N$ \st $x\cdot y=q\cdot n+x\cdot_n
y$. \E\Tf $[x\cdot y] = [q\cdot n+x\cdot_n y] = [q\cdot n]\hpl[x\cdot_n y] =
n\Z \hpl [x\cdot_n y]= \wh0\hpl[x\cdot_n y]=[x\cdot_n y]$. Hence $[x\cdot_n y]
=[x\cdot y] = [x]\hcd[y] = f(x)\hcd f(y)$. Thus $f(x\cdot_n y) =f(x)\hcd
f(y)$.
\endproof

\bdf5.24
Let $z,w\in\Z^\t$. Then \ti{$z$ divides $w$} (notation $z|w$) if $|d|\big||w|$.
\edf

\bex5.25 \

\hph i,i, Let $m\in\Na$ be \st $m\ge11$, and let \er{1.44} be its \dr. Show
that $11|m$ iff $\suml_{k=0}^N (-1)^k\a_k|m$ where $\suml_{k=0}^N$ is the \cme
sum in $(\Z,+,0)$.

\hph ii,, Let $\rho_0:=1$, $\rho_1:=3$, $\rho_2:=2$, $\rho_3:=-1$, $\rho_4:=-3$,
$\rho_5:=-2$, where $-l$, $l\in[0,7)$, denotes the inverse of~$l$ in the group
$(\N_7,+_7,0)$. Set $\rho_{i+6q}:=\rho_i$, $i\in[0,6)$, $q\in\N$. Show that
$7|m$ iff $7\big|\suml_{k=0}^N \rho_k\a_k$ where $\suml_{k=0}^N$ is the \cme
sum in $(\Z,+,0)$ and $-l$, $l\in[0,7)$, is the inverse of~$l$ in the group
$(\Z,+,0)$.
\eex

We now return to the problem of embedding the ring $(\Z,+,\cdot,0,1)$ in a
field, that is, of finding a field $(\wh F,\hpl,\hcd,\wh0,\wh1)$ and an
\ti{in\jc} ring-\hm sm $j:\Z\to\wh F$. Note that if a ring $(R,+,\cdot,0,1)$
can be embedded in a field~$\wh F$ then the \fw\ \pp y of~$R$ holds:
\beq5.55a
x\cdot y\ne0 \qh{whenever }x,y\in R\sms0.
\e

\bdf5.26
An \el\ $x$ of a ring $(R,+,\cdot,0,1)$ is called a \ti{\zd} if $x\in R\sms0$
and if \te s $y\in R\sms0$ \st $x\cdot y=0$.\index{zero divisor}
\edf

For example, in the ring $(\N_6,+_6,\cdot_6,0,1)$, $2$~is a \zd\ since $2
\cdot_6 3\nde1.33 = \F_6(2\cdot3)=\F_6(0+6)\nde1.4 = \F_6(0)\nde1.6 = 0$. Of
course $3$~is also a \zd\ in~$\N_6$.

\blm5.27
Let $R$ be a ring. Then
\bga5.56
\hbox{a unit is \emph{not} a \zd\rm;}\\
\hbox{if $R$ is \emph{finite} and $x\in R$, then $x$ is a unit iff $x$ is not a \zd.}
\lb{5.57}
\e
\elm

\proof \

\er{5.56}: Let $x$ be a unit of $R$. Then by \df\ \te s $z\in R\sms0$ \st
$x\cdot z=1$. Let $y\in R\sms0$. It suffices to show that if $x\cdot y=0$, then
$y=0$. If $x\cdot y=0$, then $y\cdot x=0$ and $0\nde4.1 = (y\cdot x)\cdot z =
y\cdot(x\cdot z)= y\cdot 1=y$. We used the fact that $(R,\cdot,1)$ is an \am.

\er{5.57}: Let $\d:R\to R$ be the map defined by $\d(z):=x\cdot z$, $x\in R$,
where $x$~is not a zero divisor. We show that $\d$~is \ti{in\jc}. Let $z_1,z_2
\in R$ be \st $\d(z_1)=\d(z_2)$, and let $w\in R$ be the inverse of~$z_2$ in
$(R,+,0)$. Then $0\nde4.1 = x\cdot0 =x\cdot (z_2+w)\nde4.2 = x\cdot z_2 +
x\cdot w = x\cdot z_1 + x\cdot w \nde4.2 = x\cdot(z_1+w)$. Since $x$~is not
a \zd, $z_1+w=0$. Thus $z_2 = (z_1+w)+z_2 = z_1+(w+z_2) = z_1+(z_2+w) = z_1+0
=z_1$. Since $\d$ is in\jc\ and $R$~is finite, $\d$~is sur\jc\ by \era2{3.11}.
Hence \te s $y\in R$ \st $x\cdot y=\d(y)=1$. Clearly $y\ne0$, since otherwise
$x\cdot y=0\ne1$. Thus $x$~is a unit.
\endproof

\bdf5.28
A ring $R$ \st no $x\in R\sms0$ is a \zd\ is called an \ti{\ido}\index{integral domain} or simply
a \ti{domain}.\index{domain}
\edf

\bdf5.29
Some authors (\cite{Alg}, \cite{Is}) call $0$ a \zd. Then a domain is a ring
\st $0$~is the only \zd.
\edf

Note that \ti{no} \el\ of $\Z\sms0$ is a \zd\ since $z,w\in\Z\sms0$ implies
$z\cdot w\in\Z\sms0$. Indeed, if $z,w \in \Z\sms0$, then $|z|,|w|\in\Na$, hence
$|z|\,|w|\ne0$ by \E\Pr\ \rfa2{p2.7}\,(i). Thus $|z\cdot w|\nde5.40 = |z|\,|w|
\ne0$, and $z\cdot w\in \Z\sms0$. However, only $1$~and~$-1$ are units of~$\Z$,
since $z\cdot w=1$ implies $z,w\in\Z\sms0$ (otherwise $z\cdot w=0$) and $|z|\,|w|
\nde5.40 = |z\cdot w|=1$ implies $|z|=|w|=1$ by \E\Pr\ \rfa2{p2.7}\,(iv), thus
$z,w\in\{1,-1\}$.
Thus the ring $\Z$ is an \ido.

\ssk
A finite domain is a field, since by Lemma \rf{l5.27} every nonzero \el\ is
a unit. Note that if a ring~$R$ is embedded in a field~$\wh F$, then $R$ is a
domain. Let $R,\wh F,j$ be \st $j:R\to\wh F$ is an in\jc\ ring-\hm sm. Let
$x,y\in R\sms0$. Then $j(x),j(y)\in \wh F\sms0$ since $j(0)=\wh0$ and $j$~is
in\jc. Since $\wh F$~is a field, we infer $j(x)\hcd j(y)\in\wh F\sms{\wh0}$.
Thus $j(x\cdot y) = j(x)\hcd j(y)\in\wh F\sms{\wh0}$, hence $x\cdot y\in R
\sms0$.

In the next theorem (attributed to Steiniz in \cite[pp.~22--23]{Nrs}) it is
shown that every infinite domain can be embedded in a field.

We first recall that if $(D,+,\cdot,0,1)$ is a domain then
\bea5.58
&(D,+,0) \hbox{ is an \ag.}\\
&(D,\cdot,1) \hbox{ is an \am.} \lb{5.59}\\
&\hbox{\er{1.1n}, \er{1.2n} and \er{1.3n} hold.} \lb{5.60}\\
&\hbox{If }x,y\in D^\t, \hbox{ then }x\cdot y\in D^\t, \hbox{ where }
D^\t :=D\sms0. \lb{5.61}
\e

Since $D^\t\sbs D$, $1\in D^\t$ and \er{5.61} holds, $D^\t$~is a \sbm\ of
$(D,\cdot,1)$ and $(D^\t,\cdot,1)$ is an \am. We claim that $(D^\t,\cdot,1)$
is a \Cm. Defining, as in the proof of \er{5.57}, the map $\d:D\to D$ by
setting $\d(z):=x\cdot z=z\cdot x$, $x\in D^\t$, $z\in D$, we find that $\d$~is
in\jc. Hence $y\cdot x=z\cdot x$, $y,z\in D^\t$, implies $y=z$. Hence
\era2{1.8} holds and the claim is proved.

Thus we have:
\beq5.62
(D^\t,\cdot,1) \hbox{ is a \Cm.}
\e

If $D$ is \ti{finite}, then so is $D^\t$ and $(D^\t,\cdot,1)$ is a group by Lemma \rf{l5.27},
hence $D$ is a \ti{field\/}.
From now on we assume $D$~\ti{infinite}, hence we are in a position to apply
Theorem \rf{t3.84} with $(X,\qu,e):=(D^\t,\cdot,1)$. Set
\beq5.63
\wh D{}^\t := \wh X \qh{and }\wh D:=\wh D{}^\t \cup\{\wh0\},
\e
where $\wh0$ does not belong to $D^\t$,
\beq5.64
{\hcd}:= {\hqu} \qh{and }\wh1:=[(1,1)].
\e
We also define a map $\wh\jmath :D\to \wh D$ by setting
\beq5.65
\wh\jmath (0):=\wh0 \qh{and }\wh\jmath (x):=j(x), \q x\in\wh D{}^\t.
\e
By Theorem \rf{t3.84} $(\wh D{}^\t,\cdot,\wh1)$ is an \ti{\ag}.

Next, we extend the binary \op\ $\hcd$ to~$\wh D$ by setting
\beq5.66
\wh0 \hcd \wh x:=\wh 0,\ \wh x\in\wh D \qh{and }
\wh x\hcd \wh 0:=\wh 0, \ \wh x\in \wh D{}^\t.
\e
As in the proof of Theorem \rf{t5.6} one shows that
\beq5.67
(\wh D,\hcd,\wh 1) \hbox{ is an \am,}
\e
and $\wh D{}^\t$ is a \sbm\ of $\wh D$, since $\wh D{}^\t\sbs \wh D$ and
$\wh 1\in \wh D{}^\t$, $\wh x,\wh y\in\wh D{}^\t$ implies $\wh x\hcd\wh y
\in \wh D{}^\t$.

Recall that the \ev ce \rl\ used in Theorem \rf{t3.84} is defined by
\beq5.68
(x,y)\sim(x',y') \qh{if }x\cdot y'=y\cdot x', \q
x,y,x',y' \in D^\t.
\e
Motivated by Theorem \rf{t5.6} we want to define the \ad\ $\hpl$ on $\wh D$ by
setting
\beq5.69
\wh0 \hpl \wh x:= \wh x, \ \wh x\in \wh D \qh{and } \wh x\hpl \wh 0:=\wh x,
\ \wh x\in\wh D{}^\t,
\e
and
\beq5.70
\frac ab \hpl \frac cd:= \frac{ad+bc}{bd}
\e
for $a,b,c,d\in D^\t$ where
\beq5.71
\frac ab:=[(a,b)], \hbox{ the $\sim$-\ev ce class containing }(a,b).
\e
However, $\frac{ad+bc}{bd}$ is not defined when $ad+bc=0$. We would like to set
\beq5.72
\frac ab+\frac cd:=\wh 0 \qh{if }ad+bc=0,\ a,b,c,d \in D^\t.
\e
However, this \df\ only makes sense if $\frac{a'}{b'}=\frac ab$, $\frac{c'}{d'}
=\frac cd$ and $\frac ab+\frac cd=0$ implies ${a'd'+b'c'=0}$, $a,b,c,d,a',b',c',d'
\in D^\t$. From Lemma \rf{l5.3} we obtain \er{5.14}. Thus $bd(a'd'+b'c')=0$ if
$ad+bc=0$. Since $D$~is a \ti{domain} and $bd\in D^\t$, we obtain $a'd'+b'c'=0$.
Thus $\frac ab+\frac cd$ is well-defined when $ad+bc=0$. The same applies for
the case $ad+bc\ne0$ in \er{5.70}. In this case, $a'd'+b'c'\ne0$ by what precedes, and
from \er{5.14} we infer $(ad+bc,bd)\sim (a'd'+b'c',b'd')$ whenever $(a,b)\sim
(a',b')$ and $(c,d)\sim(c',d')$. Thus, the \fw\ \df\ makes sense.
\beq5.73
\frac ab \hpl \frac cd := \bca
\dfrac{ad+bc}{bd}, & \hbox{if } ad+bc\ne0,\\
\wh0, & \hbox{if }ad+bc=0,
\eca \q\ a,b,c,d\in D^\t.
\e

\ti{\E\cmt ity of $\hpl$}:
\beq5.74
\wh x\hpl \wh y=\wh y\hpl \wh x, \q \wh x,\wh y\in\wh D.
\e
If $\wh x$ or $\wh y$ equals $\wh 0$, \er{5.74} follows from \er{5.69}. If
$\wh x,\wh y\in \wh D{}^\t$ and $\wh x=\frac ab$, $\wh y=\frac cd$, $a,b,c,d
\in D^\t$ then $ad+bc=0$ iff $cb+da=0$ in view of the \cmt ity of $+$
and~$\cdot$ in~$D$. Hence $\wh x\hpl\wh y=\wh0$ implies $\wh y\hpl\wh x=\wh0 =
\wh x\hpl\wh y$. If $ad+bc\ne0$, then $cb+da\ne0$, and $\wh x\hpl \wh y =
\frac ab\hpl \frac cd \nde5.73 = \frac{ad+bc}{bd}\nad*= \frac{cb+da}{db}
\nde5.73 = \frac cd\hpl \frac ab = \wh y\hpl\wh x$. In $\nad*=$ we used the
\cmt ity of~$+$ and~$\cdot$ in~$D$.

\brm5.32
\E\df\ \er{5.73} requires distinguishing two cases in the proof of the \cmt ity
of the \ad~$\hpl$. A~fortiori more cases have to be considered in the proof of
the \asc ity. \E\Tf it would be convenient to have a \df\ of~$\wh0$ as the
``fraction'' $\frac0b$, $b\in D^\t$. It turns out that the \ext\ of \rl\
\er{5.68} to $D\t D^\t$ is an \ev ce \rl, which allows us to define~$\wh0$ as
$[(0,b)]$, $b\in D^\t$.
\erm

\blm5.33
Let $(x,y), (x',y')\in D\t D^\t$ and set
\beq5.77
(x,y)\sim(x',y') \hbox{ if } x\cdot y' = y\cdot x'.
\e
Then $\sim$ is an \emph{\ev ce \rl} on $D\t D^\t$. \Mo
\beq5.78
(x,y)\sim (0,y'),\q x\in D,\ y,y'\in D^\t \hbox{ iff }x=0.
\e
\elm

\bex5.34
Show that the \ext\ to $D\t D\sms{(0,0)}$ of the \rl\ defined in \er{5.77} is
an \ev ce \rl.
\eex

\proof
Let $(x,y),(x'y'),(x'',y'')$ belong to $D\t D^\t$.

\ti{Reflexivity}: $x\cdot y=y\cdot x$ by \era2{1.6}, hence $(x,y)\sim(x,y)$
by \er{5.77}.

\ti{\E\sy y}: Suppose $(x,y)\sim(x',y')$. Then $x\cdot y' \nde5.77 = y\cdot
x'$, hence $x'\cdot y\nda21.6 = y\cdot x'=x\cdot y'\nda21.6 = y'\cdot x$.
Thus $(x',y')\sim(x,y)$ by \er{5.77}.

\ti{\E\tr ity}: Suppose $(x,y)\sim(x',y')$ and $(x',y')\sim(x'',y'')$. We
have $x\cdot y'=y\cdot x'$ and $x'\cdot y''=y'\cdot x''$ by \er{5.77}. Hence
$(x\cdot y'')\cdot y' = x\cdot(y''\cdot y') = x\cdot(y'\cdot y'') = (x\cdot y')
\cdot y'' = (y\cdot x')\cdot y'' = y\cdot(x'\cdot y'') = y\cdot(y'\cdot x'') =
y\cdot(x''\cdot y') = (y\cdot x'')\cdot y'$. Since $y'\in D^\t$, we obtain
$x\cdot y'' = y\cdot x''$ by \er{5.62}. Hence $(x,y)\sim(x'',y'')$
 by \er{5.77}.

``\ti{If\/}'': Let $a=0$, $b,d\in D^\t$. Then $a\cdot d=0\cdot d\nde4.1 =
0\nde4.1 = b\cdot 0$, hence $(a,b)\sim(0,d)$ by \er{5.77}.

``\ti{Only if\/}'': Let $a\in D$, $b,d\in D^\t$ be \st $(a,b)\sim(0,d)$. Then
$a\cdot d\nde5.77 = b\cdot 0\nde4.1 = 0$. Since $D$ is a domain, we have
$a=0$.
\endproof

\bth5.35
Let $(D,+,\cdot,0,1)$ be a domain and let $D^\t:=D\sms0$. Let $\sim$ denote
the \ev ce \rl\ on $D\t D^\t$ introduced in Lemma \rf{l5.33}.

\hph i,ii, If $(a,b)$, $(c,d)$, $(a',b')$, $(c',d')\in D\t D^\t$ \sf y $(a,b)
\sim(a',b')$ and $(c,d)\sim(c',d')$, then $b\cdot d,b'\cdot d'\in D^\t$ and
\bga5.79
(a\cdot d+b\cdot c, b\cdot d)\sim(a'\cdot d'+b'\cdot c',b'\cdot d'),\\
(a\cdot c,b\cdot d)\sim(a'\cdot c',b'\cdot d').\lb{5.80}
\e

\hph ii,i, Let $(a,b)\in D\t D^\t$. Set
\beq5.81
\frac ab:=[(a,b)]
\e
where $[(a,b)]$ denotes the $\sim$-\ev ce class containing $(a,b)$, and set
\beq5.82
\wh D:=D\t D^\t/{\sim}.
\e
Let $(a,b),(c,d)\in D\t D^\t$, then
\beq5.83
\frac ab=\frac cd \qh{iff } a\cdot d=b\cdot c.
\e
\Mo the \fw\ binary \op s $\hpl$ and $\hcd$ on $\wh D$ are well-defined.
\bea 5.84
\frac ab \hpl \frac cd &:= \frac{a\cdot d+b\cdot c}{b\cdot d}\,, \q \frac ab\,,
\frac cd\in\wh D,\\
\frac ab \hcd \frac cd &:= \frac{a\cdot c}{b\cdot d}\,, \q \frac ab\,,
\frac cd\in\wh D.\lb{5.85}
\e
Set
\beq5.86
\wh0:=\frac01\,,\q \wh1:=\frac11\,.
\e
Then $(\wh D,\hpl,\hcd,\wh0,\wh1)$ is a \emph{field}.

\hph iii,, Define $j:D\to\wh D$ by setting
\beq5.87
j(a):=\frac a1\,, \q a\in D.
\e
Then $j$ is an \emph{in\jc\ ring-\hm sm} from $D$ into~$\wh D$. \Mo \fe
$\wh x\in\wh D$ \te s $(a,b)\in D\t D^\t$ \st
\beq5.88
\wh x= j(a)\hcd j(b)\Inv
\e
where $j(b)\Inv:=\frac1b$ the inverse of $j(b)$ in the group $(\wh D\sms{\wh 0},
\hcd,\wh1)$.

\hph iv,, If $\wh D_0$ is a subfield of $\wh D$ $($see \E\df\ \rf{d4.23}$)$ \st
$j(D)\sbs \wh D_0$, then\break $(\wh D_0,\hpl,\hcd,\wh0,\wh1)=(\wh D,\hpl,\hcd,
\wh0,\wh1)$.

\hph v,i, Let $(D',+',\cdot',0',1')$ be a domain, and let $\si:D\to D'$ be
a ring-\is sm. Let $a,c\in D$ and $b,d\in D^\t$. If $\frac ab=\frac cd$, then
$\frac{\si(a)}{\si(b)}=\frac{\si(c)}{\si(d)}$.

Let $\si:\wh D\to\wh D{}'$ be defined by
\beq5.89
\si(\wh x):=\frac{\si(a)}{\si(b)}, \q \wh x:=\frac ab\,, \q
(a,b)\in D\t D^\t,\ \wh x\in \wh D.
\e
Then $\si$ is \emph{well-defined} and $\si$ is a \emph{ring-\is sm}.
\eth

\begin{dfn}[\cite{Alg}]\lb{d5.35}
The field $\wh D$ introduced in Theorem \rf{t5.35} is called the \ti{field of
\qt s} of~$D$.\index{field!of quotients}
\edf

\brm5.36
Note that if in \er{5.84} we allow $(a,b),(c,d)\in (D\t D)\sms{(0,0)}$, then
$\frac{a\cdot d+b\cdot c}{b\cdot d}$ may not belong to $(D\t D)\sms{(0,0)}$.
Indeed, if $a=c=1$ and $b=d=0$, then $(1,0)\in(D\t D)\sms{(0,0)}$ and
$\frac{a\cdot d+b\cdot c}{b\cdot d}=\frac{0+0}{0\cdot0}=\frac 00$.
\erm

\proof[Proof of Theorem \rf{t5.35}]\

(i) Let $x\in D^\t$. If $y\in D^\t$ then $x\cdot y\ne0$ by \er{5.63}. Since
$b,d,b',d'\in D^\t$, we have $b\cdot d,b'\cdot d'\in D^\t$. In view of
\er{5.77}, \er{5.79} means $(a\cdot d+b\cdot c)\cdot(b'\cdot d')=(b\cdot d)
\cdot(a'\cdot d'+b'\cdot c')$, which is a con\sq\ of Lemma \rf{l5.3}.
Similarly, \er{5.79} means $(a\cdot c)\cdot(b'\cdot d') = (b\cdot d)\cdot
(a'\cdot c')$. From $a\cdot b'=b\cdot a'$ and $c\cdot d'=d\cdot c'$, we infer
$(a\cdot c)\cdot(b'\cdot d')\nda21.36 = a\cdot(c\cdot b')\cdot d' =a\cdot
(b'\cdot c)\cdot d' = (a\cdot b')\cdot(c\cdot d')= (b\cdot a')\cdot(c\cdot d')
= (b\cdot a')\cdot(d \cdot c')= b\cdot(a'\cdot d)\cdot c'=b\cdot(d\cdot a')
\cdot c' = (b\cdot d)\cdot(a'\cdot c')$.

(ii) $\frac ab=\frac cd$ iff $[(a,b)]=[(c,d)]$ by \er{5.81}, and $[(a,b)] =
[(c,d)]$ iff $(a,b)\in[(c,d)]$ by Exercise \rfa1{ex4.4}\,(i) iff $(a,b)\sim
(c,d)$ by \era1{4.4}. Thus \er{5.83} holds.

The \op s $\hpl$ and $\hcd$ are well-defined by \er{5.84}, \er{5.85} in view
of~(i).

\bex5.37
Let $x,y\in D$ and $u,z\in D^\t$. Show
\beq5.90
\frac{x\cdot u}{z\cdot u} = \frac{u\cdot x}{u\cdot z}=\frac xz\,; \q
\frac xz\hpl \frac yz = \frac{x+y}z\,.
\e
\eex

Let $\wh x,\wh y,\wh z\in\wh D$, let $a,c,e\in D$ and let $b,d,f\in D^\t$ be
\st $x=\frac ab$, $y=\frac cd$, $z=\frac ef$.

\ssk
\bq 5.91 $(X,\hpl,\wh0)$ is an \ag\e

``\ti{\E\asc ity}'': $\bg(\frac ab\hpl\frac cd)\hpl \frac ef \nde5.84 = \frac
{a\cdot d+b\cdot c}{b\cdot d}\hpl \frac ef\nde5.84 = \frac{((a\cdot d+b\cdot c)
\cdot f)+(b\cdot d)\cdot e}{(b\cdot d)\cdot f}$,

$(a\cdot d+b\cdot c)\cdot f+(b\cdot d)\cdot e\nde4.3 = ((a\cdot d)\cdot f+(b\cdot c)\cdot f)+
(b\cdot d)\cdot e\nad{\er{5.58},\er{5.59}} = (a\cdot d\cdot f+b\cdot c\cdot f)
+b\cdot d\cdot e\nde5.58 = a\cdot d\cdot f+b\cdot c\cdot f+b\cdot d\cdot e$.
Similarly $(b\cdot d)\cdot f=b\cdot d\cdot f$. Hence
\beq5.92
(\wh x+\wh y)\hpl \wh z = \frac{a\cdot d\cdot f+b\cdot c\cdot f+b\cdot d
\cdot e}{b\cdot d\cdot f}\,.
\e
\E\oh we have
\bmlg
\wh x\hpl(\wh y\hpl\wh z) = \frac ab \hpl\Bg(\frac cd\hpl \frac ef) =
\frac ab\hpl\frac{c\cdot f+d\cdot e}{d\cdot f}=\frac{a\cdot (d\cdot f)+(b\cdot
(c\cdot f+d\cdot e))}{b\cdot(d\cdot f)}\\
{}=\frac{a\cdot(d\cdot f)+(b\cdot(c\cdot f)+b\cdot(d\cdot e))}{b\cdot (d\cdot
f)} = \frac{a\cdot d\cdot f+b\cdot c\cdot f+b\cdot d\cdot e}{b\cdot d\cdot f}
\nde5.92 = (\wh x\hpl\wh y)\hpl\wh z.
\e

``\ti{\E\cmt ity}'': $\wh x\hpl\wh y = \frac ab\hpl \frac cd \nde5.84 =
\frac{a\cdot d+b\cdot c}{b\cdot d}\nad{\er{5.58},\er{5.59}} = \frac{c\cdot b
+d\cdot a}{d\cdot b} = \frac cd\hpl\frac ab =\wh y\hpl\wh x$.

``$\wh0\hpl\wh x=\wh x$'': $\wh0\hpl\wh x\nde5.86 = \frac01\hpl\frac ab\nde5.84
= \frac{0\cdot b+1\cdot a}{1\cdot b}\nad{\er{5.58},\er{5.59}} = \frac{0+a}b
\nde5.58 = \frac ab=\wh x$.

``\ti{\E\ex\ of inverse}'': Let $\wh x=\frac ab$. By \er{5.58} \te s $a'\in D$
\st $a+a'=0$. Set $\wh x{}':=\frac{a'}b$. Then $\wh x\hpl\wh x{}'=\frac ab\hpl
\frac{a'}b\nde5.90 = \frac{a+a'}b = \frac 0b \nde5.86 = \wh0$.
\endproof

\bq 5.93 $(\wh X,\hcd,\wh1)$ \ti{is an \am}\e

``\ti{\E\asc ity}'': $(\wh x\hcd \wh y)\hcd\wh z = (\frac ab\hcd \frac cd)\hcd
(\frac ef) \nde5.85 = \frac {a\cdot c}{b\cdot d}\hcd\frac ef\nde5.85 =
\frac{(a\cdot c)\cdot e}{(b\cdot d)\cdot f} \nde5.59 = \frac{a\cdot(c\cdot e)}
{b\cdot(d\cdot f)} \nde5.85 = \frac ab\hcd \frac{c\cdot e}{d\cdot f} =
\frac ab\hcd (\frac cd\hcd \frac ef)= \wh x\hcd(\wh y\hcd \wh z)$.

''\ti{\E\cmt ity}'': $\wh x\hcd\wh y = \frac ab\hcd\frac cd \nde5.85 =
\frac{a\cdot c}{b\cdot d}\nde5.59 = \frac{c\cdot a}{d\cdot b}\nde5.85 =
\frac cd\hcd \frac ab=\wh y\hcd \wh x$.

``$\wh 1\hcd\wh x=\wh x$'': $\wh 1\hcd\wh x\nde5.86 = \frac11\hcd\frac ab
\nde5.85 = \frac{1\cdot a}{1\cdot b}\nde5.59 = \frac ab=\wh x$.

\ssk
\bq 5.94 $\wh0\hcd\wh x:$\e
$\wh0\hcd\wh x \nde5.86 = \frac01\hcd\frac ab\nde5.85 = \frac{0\cdot a}
{1\cdot b} \nad{\er{4.1},\er{5.59}} = \frac0b\nde5.86 = \wh0$.

\ssk
\bq 5.95 $\wh x\hcd(\wh y\hpl\wh z)=\wh x\hcd \wh y\hpl\wh x\hcd\wh z$\e
\bmlg
\wh x\hcd(\wh y+\wh z) = \frac ab\hcd\Bg(\frac cd\hpl\frac ef)\nde5.84 =
\frac ab\hcd \frac{c\cdot f+d\cdot e}{d\cdot f} \nde5.85 = \frac{a\cdot (c
\cdot f+d\cdot e)}{b\cdot(d\cdot f)}\\ {}\nde4.2 = \frac{a\cdot(c\cdot f)+
a\cdot(d\cdot e)}{b\cdot(d\cdot f)}\nad{\er{5.58},\er{5.59}}= \frac{a\cdot c
\cdot f+a\cdot d\cdot e}{b\cdot d\cdot f} \nde5.90 = \frac{a\cdot c\cdot f}
{b\cdot d\cdot f}\hpl\frac{a\cdot d\cdot e}{b\cdot d\cdot f}\\{}\nde5.85 =
\frac ab\hcd \frac{c\cdot f}{d\cdot f}\hpl \frac ab\hcd \frac{d\cdot e}
{d\cdot f}\nde5.90 = \frac ab\hcd\frac cd \hpl \frac ab\hcd\frac ef =
\wh x\hcd\wh y \hpl \wh x\hcd\wh z.
\e

\If from \er{5.91}, \er{5.93}, \er{5.94} and \er{5.95} that $(\wh X,\hpl,\hcd,
\wh 0,\wh 1)$ is a (\cmt e) ring (with unity). It remains to show that if
$x\in \wh D{}^\t$, then $\wh x$ is a \ti{unit\/}. We have $\wh x=\frac ab$ with
$a,b\in D^\t$. Set $\wh x{}'':=\frac ba$. Then $\wh x\hcd \wh x{}'' =\frac ab
\hcd \frac ba\nde5.85 = \frac{a\cdot b}{b\cdot a}\nde5.59 = \frac{a\cdot b}
{a\cdot b}\nde5.83 = \frac 11\nde5.86 =\wh 1$.

This completes the proof of (ii).

\ssk
(iii) ``\ti{In\ji\ of $j$}'': Let $a,b\in D$ be \st $j(a)=j(b)$. Then $\frac a1
\nde5.87 = j(a)=j(b)\nde5.87 = \frac b1$. From $\frac a1=\frac b1$ we infer by
\er{5.83} $a\cdot 1=1\cdot b$. Hence $a\nde5.59 = a\cdot1=1\cdot b\nde5.59 =b$.

``\ti{$j:(D,+,0)\to (\wh D,\hpl,\wh0)$ is a \hm sm}'': $j(0)\nde5.87 = \frac01
\nde5.86 = \wh0$. Let $a,b\in D$. Then $j(a+ b)\nde5.87 = \frac{a+b}1
\nde5.90 = \frac a1\hpl\frac b1 = j(a)\hpl j(b)$.

``\ti{$j:(D,\cdot,1)\to (\wh D,\hcd,\wh1)$ is a \hm sm}'': $j(1)\nde5.87 =
\frac11 \nde5.86 = \wh1$. Let $a,b\in D$. Then $j(a\cdot b) \nde5.87 = \frac
{a\cdot b}1 \nde5.59 = \frac{a\cdot b}{1\cdot 1}\nde5.85 = \frac a1\hcd
\frac b1=j(a)\hcd j(b)$.

Let $\wh x\in\wh D$. Then \te s $(a,b)\in D\t D^\t$ \st $\wh x=\frac ab$. We
have $\frac ab\nde5.59 = \frac{a\cdot1}{1\cdot b}\nde5.85 = \frac a1\hcd
\frac1b \nde5.87 = j(a)\hcd\frac1b$. Note that $\frac b1\cdot \frac1b\nde5.59
= \frac bb \nad{\er{5.59},\er{5.83}} = \frac11\nde5.86 = \wh1$. Hence $\frac1b$
is the inverse of $j(b)$ in the group $(\wh D\sms{\wh0},\hcd,\wh1)$, which we
denote by $j(b)\Inv$. \E\Tf $\wh x=j(a)\hcd j(b)\Inv$.

\ssk
(iv) In view of the \df\ of a subfield it is \sft\ to show that the set $\wh D
_0$ equals the set~$\wh D$, hence that $\wh D$ is contained in~$\wh D_0$. By Lemma
\rf{l4.24} $(\wh D_0,\hpl,\hcd,\wh0,\wh1)$ is a field where $\hpl$ (resp.~$\hcd$)
denotes the \rt ion to~$\wh D_0$ of the \ad\ (resp.\ \mlc) of~$\wh D$. Since
$\wh0\in\wh D_0$, it suffices to show that every $\wh x\in \wh D\sms{\wh0}$ belongs
to~$\wh D_0$. Let $\wh x\in\wh D{}\sms{\wh0}$. By~(iii) \te\ $a,b\in wh D{}\sms{\wh0}$ \st $\wh x
=j(a)\hpl j(b)\Inv$. By \as\ $j(a),j(b)\in\wh D_0$. Since $j(0)=\wh 0$ and
$j$~is in\jc\ by~(iii), $j(a),j(b)\in\wh D_0\sms{\wh0}$. Since $(\wh D_0,\hpl,
\hcd,\wh0,\wh1)$ is a field, $(\wh D_0\sms{\wh0},\hcd,\wh1)$ is a group,
hence a subgroup of $(\wh D\sms{\wh0},\hcd,\wh1)$. By Lemma \rf{l3.3},
$j(b)\Inv$ the inverse of $j(b)$ in the group $(\wh D\sms{\wh0},\hcd,\wh1)$ is
also the inverse of $j(b)$ in the subgroup $(\wh D_0\sms{\wh0},\hcd,\wh1)$. \E\Tf $j(b)\Inv\in
\wh D_0\sms{\wh0}$, and $\wh x=j(a)\hcd j(b)\Inv\in \wh D_0\sms{\wh0}$.

\ssk
(v) By \er{5.83}, $\frac ab=\frac cd$ implies $a\cdot d=b\cdot c$. Since
$\si$~is a ring-\is sm, $\si(b)$ and~$\si(d)$ belong to $D'\sms{0'}$. Indeed,
$\si$~is in\jc, $\si(0)=0'$ and $b\ne0$, $d\ne0$. Thus $\frac{\si(a)}{\si(b)}$,
$\frac{\si(c)}{\si(d)}$ are well-defined. \Mo since $\si:(D,\cdot,1)\to
(D',\cdot',1')$ is a monoid-\hm sm, $\si(a)\cdot'\si(d)=\si(a\cdot d)=
\si(b\cdot c)=\si(b)\cdot'\si(c)$. Hence by \er{5.83}, $\frac{\si(a)}{\si(b)}
=\frac{\si(c)}{\si(d)}$. Thus $\wh\si$ is well-defined by \er{5.89}.

``\ti{$\wh\si:(\wh D,\hpl,\wh0)\to(\wh D{}',\hpl',\wh0)$ is a \hm sm}'':
We have $\wh\si(\wh0)\nde5.86 = \wh\si(\frac01)\nde5.89 = \frac{\si(0)}{\si(1)}$.
Since $\si$~is a ring-\hm sm, we infer $\si(0)=0'$ and $\si(1)=1'$. \E\Tf
$\wh\si(\wh0)=\frac{0'}{1'}\nde5.86 = \wh0{}'$. Let $\wh x,\wh y\in\wh D$ and
$(a,b),(c,d)\in D\t D^\t$ be \st $\wh x=\frac ab$, $\wh y=\frac cd$. Then
$\si(\wh x+\wh y)\nde5.84 = \si\bg(\frac{a\cdot d+b\cdot c}{b\cdot d})$, where
$b\cdot d\in D^\t$. \E\Tf  $\si\bg(\frac{a\cdot d+b\cdot c}{b\cdot d})\nde5.89 =
\frac{\si(a\cdot d+b\cdot c)}{\si(b\cdot d)} \nad*= \frac{\si(a\cdot d)+'
\si(b\cdot c)}{\si(b\cdot d)}\nad*= \frac{\si(a)\cdot' \si(d)+' \si(b)\cdot'
\si(c)}{\si(b)\cdot'\si(d)} \nde5.90 = \frac{\si(a)\cdot' \si(d)
}{\si(b)\cdot'\si(d)} \hpl \frac{\si(b)\cdot'\si(c)}{\si(b)\cdot'\si(d)}
\nde5.89 = \frac{\si(a)}{\si(b)} \hpl \frac{\si(c)}{\si(d)} = \wh\si(\wh x)
\hpl\wh\si(\wh y)$.  In $\nad*=$ we used the fact that $\si$~is a ring-\hm sm.

``\ti{$\wh\si:(\wh D,\hcd,\wh1)\to(\wh D{}',\hcd',\wh1{}')$ is a \hm sm}:
$\si(\wh1)\nde5.86 = \si(\frac11)\nde5.89 = \frac{\si(1)}{\si(1)}=\frac{1'}{1'}
\nde5.86 = \wh1{}'$.

Let $\wh x,\wh y\in\wh D$, $(a,b),(c,d)\in D\t D^\t$ be \st $\wh x=\frac ab$,
$\wh y=\frac cd$. Then $\wh\si(\wh x\hcd\wh y)=\wh\si(\frac ab\hcd \frac cd)
\nde5.85 = \si(\frac{a\cdot c}{b\cdot d})\nde5.88 = \frac{\si(a\cdot c)}
{\si(b\cdot d)}= \frac{\si(a)\cdot'\si(c)}{\si(b)\cdot' \si(d)}\nde5.85 =
\frac{\si(a)}{\si(b)}\hcd' \frac{\si(c)}{\si(d)}\nde5.88 = \si(\frac ab)
\hcd'\si(\frac cd) = \si(\wh x)\hcd'\si(\wh y)$.

\If that $\wh\si$ is a ring-\hm sm. In view of Lemma \rfa2{l1.8} it remains to
show that $\si$~is \ti{bi\jc}. Let $x'\in\wh D{}'$. \E\te\ $a'\in D'$ and
$b'\in D'\sms{0'}$ \st $\wh x{}'=\frac{a'}{b'}$. Let $\si\Inv:D'\to D$ denote
the inverse of~$\si$. Set $\wh x:=\frac{\si\Inv(a')}{\si\Inv(b')}$. Note that
$\si\Inv(b')\ne0$. Thus $\wh x\in\wh D$ and $\wh\si(\wh x)\nde5.88 =
\frac{\si(\si\Inv(a'))}{\si(\si\Inv(b'))}=\frac{a'}{b'}=\wh x{}'$.
Thus $\si$~is \ti{sur\jc}. Let $\wh x,\wh y\in\wh D$ be \st $\wh\si(\wh x)=
\wh\si(\wh y)$. Let $(a,b),(c,d)\in D\t D^\t$ be \st $\wh x=\frac ab$, $\wh y
=\frac cd$. Then $\frac{\si(a)}{\si(b)}\nde5.88 = \si(\frac ab)=\si(\wh x)=
\si(\wh y)=\si(\frac cd)=\frac{\si(c)}{\si(d)}$. From $\frac{\si(a)}{\si(b)}=
\frac{\si(c)}{\si(d)}$ we infer by \er{5.83} that $\si(a)\cdot'\si(d) =
\si(b)\cdot'\si(c)$. Hence $\si(a\cdot d)=\si(a)\cdot'\si(d)= \si(b)\cdot'
\si(c)=\si(b\cdot c)$. Thus $a\cdot d=\si\Inv(\si(a\cdot d)) =\si\Inv(\si(b
\cdot c))=b\cdot c$. Then $\wh x=\frac ab=\frac cd=\wh y$ follows from
\er{5.83}. This completes the proof of~(v), hence the proof of Theorem
\rf{t5.35}.
\endproof

\bex5.36
Let $(X,+,\cdot,0,1)$ be an infinite (\cmt e) \sr\ (with unity) \st $(X,+,0)$
is a \Cm. Prove the \fw:

\hph i,ii, \E\te\ a ring $(\wh X,\hpl,\hcd,\wh0,\wh1)$ and an in\jc\ (\sr)-\hm
sm\break ${j:X\to \wh X}$ \st \fe $\wh x\in\wh X$, $\wh x=j(a)-j(b)$ \fs $a,b\in X$.

\hph ii,i, Suppose in \ad\ that $X$ \sf ies the \fw\ \cn:
$$
\hbox{\E\fa $a,b,c,d\in X$ \st $a\ne b$, $a\cdot c+b\cdot d$ implies $c=d$.}
\leqno(*)
$$
Then $\wh X$ introduced in (i) is a domain.

\hph iii,, If the additive monoid of $X$ is \pn, then $X$ \sf ies \cn~$(*)$.
\eex

We now show that the field $\Q$ introduced in Theorem \rf{t5.6} and the field
of \qt s of the domain $(\Z,+,\cdot,0,1)$ are ring-\is c. In order to avoid a
collision of notation we introduce the \fw\ one: Given $a,b,c,d\in\Na$,
\beq5.96
(a,b)\sim' (c,d)\qh{if } ad=bc.
\e
We denote by $[(a,b)]'$ the $\sim'$-\ev ce class containing $(a,b)$ and by
$\frp ab$ the \crs\ fraction defined in \er{5.1}. \Mo we identify $\Na$ with
$\{(a,1): a\in\Na\}$ by means of \er{5.4}.

\bpr5.37
Let $(\Q,+,\cdot,0,1)$ denote the field introduced in Theorem \rf{t5.6} and
let $(\wh \Z,\hpl,\hcd,\wh0,\wh1)$ denote the field of \qt s of the domain
$(\Z,+,\cdot,0,1)$ introduced in Theorem \rf{t5.35}.

\hph i,i, Let $a,b,c,d\in\Na$ \sf y $(a,b)\sim'(c,d)$, and let $\frac
ab,\frac cd\in\wh \Z$. Then $\frac ab=\frac cd$ in~$\wh \Z$.

\hph ii,, Let $\th : \Q\to \wh \Z$ be defined by
\beq5.97
\bca
\displaystyle{\th\Bg(+\frp ab):= \frac ab\,,} & a,b\in\Na,\\
\th(0):=\wh 0, \\
\displaystyle{\th\Bg(-\frp ab):=\frac{-a}b,} & a,b\in\Na.
\eca
\e
Then $\th$ is a ring-\is sm.
\epr

\proof \

(i) If $(a,b)\sim'(c,d)$, then $ad=bc$ by \er{5.96}. Hence $(a,b)\sim(c,d)$
in $\Z\t \Z^\t$ since $\Na=\Z_{>0}$. The conclusion follows from \er{5.83}.

(ii) In view of (i), $\th$ is well-defined by \er{5.97}. We first show that
$\th$~is a bi\jn. To this end we define a map $\rho:\wh Z\to\Q$ \sf ying
\beq5.98
\rho\circ\th = \id_{\Q} \qh{and }\th\circ\rho=\id_{\wh \Z}.
\e
We set
\beq5.99
\bca
\displaystyle{\rho\Bg(\frac 0b):=0}, & b\in\Z^\t,\\
\displaystyle{\rho\Bg(\frac ab):=\sgn(ab)\frp{|a|}{|b|}\,,} & a,b\in\Z^\t,
\eca
\e
where $\sgn(z)$ and $|z|$ are defined in \er{5.35} and \er{5.36}.

``\ti{$\rho$ is well-defined\/}'': We need to show that $(0,b)\sim(c,d)$,
$c\in\Z$, $b,d\in\Z^\t$, implies $c=0$, and $(a,b)\sim(c,d)$ implies $\sgn(ab)
=\sgn(cd)$ and $\frp{|a|}{|b|}=\frp{|c|}{|d|}$, $a,b,c,d\in\Z^\t$. The first assertion follows
from Lemma \rf{l5.33}. From $(a,b)\sim(c,d)$ we infer $ad=bc$ by \er{5.77}.
Hence $|a|\,|d|\nde5.40 = |ad|=|bc|\nde5.40 = |b|\,|c|$, and $(|a|,|b|)\sim'
(|c|,|d|)$ since $|a|,|b|,|c|,|d|\in\Na$. Thus\break $\frp{|a|}{|b|}=\frp{|c|}{|d|}$.
\Mo $\sgn(a)\cdot\sgn(d)\nde5.39 = \sgn(ad)=\sgn(bc)\nde5.39 = \sgn(b)\cdot
\sgn(c)$. Recall that $(\{+,-\},\cdot,+)$ is an \ag\ by \er{5.22} and that
$\sgn(z)\cdot\sgn(z)=+$, $z\in\Z^\t$. \E\Tf $\sgn(a)=\sgn(a)\cdot{+} =
\sgn(a)\cdot(\sgn(d)\cdot\sgn(d)) = (\sgn(a)\cdot\sgn(d))\cdot\sgn(d) =
(\sgn(b)\cdot\sgn(c))\cdot\sgn(d) = \sgn(b)\cdot(\sgn(c)\cdot\sgn(d))=
(\sgn(c)\cdot\sgn(d))\cdot\sgn(b)$. Hence $\sgn(ab)=\sgn(a)\cdot\sgn(b) =
((\sgn(c)\cdot\sgn(d))\cdot\sgn(b))\cdot\sgn(b)=(\sgn(c)\cdot\sgn(d))\cdot
(\sgn(b)\cdot\sgn(b))=(\sgn(c)\cdot\sgn(d))\cdot{+} = \sgn(c)\cdot\sgn(d)
=\sgn(cd)$.

\medskip
``$\rho\circ\th=\id_{\Q}$'': $(\rho\circ\th)(0)\nde5.97 = \rho(\wh0)\nde5.86 =
\rho(\frac01)\nde5.99 = 0$.

$(\rho\circ\th)(+\frp ab)\nde5.97 = \rho(\frac ab)\nde5.99 = \sgn(ab)\frp
{|a|}{|b|} \nad{\er{5.39},\er{5.36}}= (\sgn(a)\cdot\sgn(b))\frp ab \nad{
\er{5.35},\er{5.22}}= (+\frp ab)$, $a,b\in\Na$.

$(\rho\circ\th)(-\frp ab)\nde5.97 = \rho(\frac{-a}b)= \sgn((-a)\cdot b)
\frp{|-a|}{|b|} \nad{\er{5.35},\er{5.39},\er{5.36}} = ({-}\cdot{+})\frp{|a|}
{|b|} \nde5.22 = -\frp ab$, $a,b\in\Na$.

\medskip
``$\th\circ\rho=\id_{\wh\Z}$'': $(\th\circ\rho)(\wh0)\nde5.96 = (\th\circ\rho)
(\frac01)\nde5.99 = \th(0)\nde5.97 = \wh0$.

$(\th\circ\rho)(\frac ab)\nde5.99 = \th\bg({\sgn(ab)}\frp{|a|}{|b|})$, $a,b\in
\Z^\t$. If $\sgn(ab)={+}$, then $\th\bigl(\sgn(ab)\frp{|a|}{|b|}\bigr)\nde5.97
= \frac{|a|}{|b|}$, and by \er{5.22} either $\sgn(a)=\sgn(b)={+}$, or
$\sgn(a)=\sgn(b)={-}$. In the first case $|a|=a$, $|b|=b$ by \er{5.36}, hence
$\frac{|a|}{|b|}=\frac ab$. In the second case $|a|=-a$, $|b|=-b$ by \er{5.36},
hence $\frac{|a|}{|b|}=\frac{-a}{-b}\nde5.83 = \frac ab$, since $(-a)\cdot b
\nde5.31 = (a\cdot b)\Inv \nde5.32 = a\cdot(-b)\nde5.33 = (-b)\cdot a$,
$a,b\in\Z^\t$.

From \er{5.98} we infer that $\th$ and $\rho$ are bi\jc\ and that $\rho=\th\Inv$,
$\th=\rho\Inv$.

We now show that $\th:(\Q,+,0)\to(\wh \Z,\hpl,\wh0)$ is a (monoid)-\hm sm.
We already showed that $\th(0)=\wh0$. Let $q,r\in\Q$. We need to show
$\th(q+r)=\th(q)\hpl\th(r)$. The case $q=0$ or $r=0$ follows from $\th(0)=\wh0$
and \era2{1.7}. Thus we suppose $q,r\in\Q\sms0$.

``$\sgn(q)=\sgn(r)={+}$'': Let $a,b,c,d\in\Na$ be \st $q:=+\frp ab$, $r:=+\frp
cd$. Then $q+r\nde5.12 = \frp{ad+bc}{bd}$, hence $\th(q+r)\nde5.97 = \frac
{ad+bc}{bd}\nde5.84 = \frac ab\hpl \frac cd\nde5.97 = \th(+\frp ab)\hpl
\th(+\frp cd) = \th(q)\hpl\th(r)$.

``$\sgn(q)=\sgn(r)={-}$'': Let $a,b,c,d\in\Na$ be \st $q:=-\frp ab$, $r:=-\frp
cd$. For $s\in\Q$ we denote by $s\Inv$ the inverse of~$s$ in the group
$(\Q,+,0)$. In view of \er{5.17} $q=(+\frp ab)\Inv$ and $r=(+\frp cd)\Inv$.
Then $q+r=((+\frp ab)+(+\frp cd))\Inv$ by \E\Pr\ \rf{p3.26}\,(i). For $z\in\wh\Z$
we also denote by $z\Inv$ the inverse of~$z$ in the group $(\wh\Z,\hpl,\wh0)$.

We claim
\beq5.100
\th(s\Inv) = (\th(s))\Inv, \q s\in\Q.
\e
Indeed, since $0+0=0$ (resp.\ $\wh0\hpl\wh0=\wh0$), $0\Inv=0$ and $\wh0{}\Inv=
\wh0$. Thus $\th(0\Inv)=\th(0)\nde5.97 = \wh0=\wh0{}\Inv=(\th(0))\Inv$. If
$s:=+\frp ef$, $e,f\in\Na$, then $s\Inv\nde5.17 = -\frp ef$. Hence $\th(s\Inv)
=\th(-\frp ef)\nde5.97 = \frac{-e}f$. Note that $\frac ef\hpl\frac{-e}f
\nde5.90 = \frac{e+(-e)}{f}=\frac 0f \nad{\er{5.77},\er{5.78},\er{5.83}}=
\frac01\nde5.86 = \wh0$. Hence $\frac{-e}f=(\frac ef)\Inv$. Finally, if
$s:=-\frp ef$, then $s\Inv=+\frp ef$. Then $\th(s\Inv)\nde5.97 = \frac ef =
(\frac{-e}f)\Inv \nde5.97 = (\th(s))\Inv$. This completes the proof of
\er{5.100}.

We obtain $\th(q+r)=\bigl(\th((+\frp ab)+(+\frp cd))\bigr)\Inv \nde5.100 = \bigl(\th
(+\frp ab)\hpl \th(+\frp cd)\bigr)\Inv = \bigl(\th
(+\frp ab)\bigr)\Inv \hpl \bigl(\th(+\frp cd)\bigr)\Inv \nde5.100 = \th(
(+\frp ab)\Inv) \hpl \th((+\frp cd)\Inv) = \th(-\frp ab)\hpl\th(-\frp cd) =
\th(q)\hpl\th(r)$.

``$q+r=0$'': We have $r=q\Inv$, $q=r\Inv$.  $\th(q+r)=\th(0)\nde5.97 = \wh0 =
\th(q)\hpl (\th(q))\Inv \nde5.100 = \th(q)\hpl \th(q\Inv)=\th(q)\hpl\th(r)$.

``$\sgn(q)={+}$, $\sgn(r)={-}$, $q>r\Inv$'':  Since $r\in\Q\sms0$ and $\sgn(r)
={-}$, $r\Inv\in \Q\sms0$ and $\sgn(r\Inv)={+}$ by \er{5.17}, \er{5.16}. Thus
$q,r\Inv\in\Q_{\ge0}$ and since $q>r\Inv$, \te s $s\in\Q_{>0}$ \st $q=s+r\Inv$
by \E\Pr\ \rf{p5.4}. Hence $q+r=(s+r\Inv)+r=s+(r\Inv+r)=s$, and $\th(q+r)=
\th(s)$. Since $\sgn(q)=\sgn(r\Inv)={+}$, we infer from what precedes $\th(q)=
\th(s)\hpl\th(r\Inv)$. Hence $\th(s)=\th(q)\hpl(\th(r\Inv))\Inv$. \E\Tf
$\th(q+r)=\th(s)=\th(q)\hpl(\th(r\Inv))\Inv \nde5.100 = \th(q)\hpl
\th((r\Inv)\Inv)\nda43.4 = \th(q)\hpl\th(r)$.

``$\sgn(q)={+}$, $\sgn(r)={-}$, $r\Inv>q$'': Let $s\in\Q_{>0}$ be \st $r\Inv=
q+s$. Then $r\Inv+s\Inv=q$, hence $q+r=s\Inv$, and $\th(q+r)=\th(s\Inv)
\nde5.100 = (\th(s))\Inv$. \E\oh from $r\Inv=q+s$ we infer $\th(r)\Inv
\nde5.100 = \th(r\Inv) = \th(q)\hpl\th(s)$, since $\sgn(q)=\sgn(s)={+}$. \E\Tf
$\wh0= \th(r)\hpl(\th(r))\Inv \nda21.5 = \th(r)\hpl\th(q)\hpl\th(s)$, hence
$\th(s)\Inv = \th(r)\hpl\th(q)$. Finally, $\th(q+r)=\th(s)\Inv= \th(r)\hpl
\th(q)\nda21.6 = \th(q)\hpl\th(r)$.

Note that if $q\in\Q\sms0$ with $\sgn(q)={+}$, and $r\in\Q\sms0$ \sf y
$q+r=0$, then $r\Inv=q$, and $\sgn(r\Inv)={+}$. In view of \E\Pr\ \rf{p5.4}
the \nog\ of $(\Q_{\ge0},+,0)$ is total, hence $\Q_{\ge0}\sms0=\Q_{>0} =
\{q\in\Q_{>0}: q=r\Inv\}\cup\{q\in\Q_{>0}: q>r\Inv\}\cup\{q\in\Q_{>0}: q<r\Inv\}$,
where $r\in\Q\sms0$, $\sgn(r)={-}$. \If from what precedes that if $q,r\in\Q
\sms0$, $\sgn(q)={+}$ and $\sgn(r)={-}$, then $\th(q+r)=\th(q)\hpl\th(r)$.

It remains to consider the case

``$\sgn(q)={-}$, $\sgn(r)={+}$'': Set $q':=q\Inv$ and $r':=r\Inv$, then
$\sgn(q')\nde5.16 = {+}$ and $\sgn(r')\nde5.16 = {-}$. We have $\th(r'+q')=
\th(r')\hpl\th(q')$ from what precedes. Hence $\th(q+r)\nda43.4 = ((\th(q+r))
\Inv)\Inv\nde5.100 = (\th((q+r)\Inv))\Inv\nda43.5 = (\th(r\Inv+q\Inv))\Inv =
(\th(r'+q'))\Inv = (\th(r')\hpl\th(q'))\Inv \nda43.5 = \th(q')\Inv \hpl\th(r')
\Inv \nde5.100 = \th(q'{}\Inv)\hpl\th(r'{}\Inv)\nda43.5 = \th(q)\hpl\th(r)$.

This completes the proof that $\th:(\Q,+,0)\to(\wh\Z,\hpl,\wh0)$ is a \hm sm.
It remains to prove

``\ti{$\th:(\Q,\cdot,1)\to(\wh\Z,\hcd,\wh1)$ is a \hm sm}'': We have $\th(1)
=\th(\frp11)\nde5.97 = \frac11\nde5.86 = \wh1$. Let $q,r\in\Q$ be \st $q=0$
(resp.\ $r=0$), then $\th(q)=\th(0)=\wh0$ (resp.\ $\th(r)=\wh0$), $q\cdot r=0$,
and $\th(q\cdot r)=\th(0)=\wh0=\th(q)\hcd\th(r)$, in view of \er{4.1}. We may
suppose $q,r\in\Q\sms0$. In view of the \cmt ity of the \mlc\ it is \sft\ to
consider three cases: ``$\sgn(q)=\sgn(r)={+}$'', ``$\sgn(q)=\sgn(r)={-}$'' and
``$\sgn(q)={+}$, $\sgn(r)={-}$''.

``$\sgn(q)=\sgn(r)={+}$'': Let $a,b,c,d\in\Na$ be \st $q:=+\frp ab$,
$r:=+\frp cd$. Then $\th(q)\nde5.97 = \frac ab$, $\th(r)=\frac cd$ and
$\th(q)\hcd \th(r)\nde5.85 = \frac{ab}{cd}$.
By \er{5.19}, \er{5.3}, $q\cdot r=+\frp{ac}{bd}$. Hence $\th(q\cdot
r)\nde5.97 = \frac{ac}{bd}= \th(q)\hcd \th(r)$.

``$\sgn(q)=\sgn(r)={-}$'': Let $a,b,c,d\in\Na$ be \st $q:=-\frp ab$,
$r:=-\frp cd$. Then $\th(q)\nde5.97 = \frac {-a}b$, $\th(r)=\frac {-c}d$,
$\th(q)\hcd \th(r)\nde5.85 = \frac{(-a)(-c)}{bd}= \frac{ac}{bd}$ by
\er{5.31}, \er{5.32}, \era4{3.4}. \Mo $q\cdot r=+\frp{ac}{bd}$ by \er{5.19},
\er{5.3}. Hence $\th(q\cdot r)\nde5.97 = \frac{ac}{bd}=\th(q)\hcd\th(r)$.

``$\sgn(q)={+}$, $\sgn(r)={-}$'': Let $a,b,c,d\in\Na$ be \st $q=+\frp ab$,
$r=-\frp cd$. Then $\th(q)\nde5.97 = \frac ab$, $\th(r)\nde5.97 = \frac{-c}d$,
$\th(q)\hcd\th(r) = \frac{a(-c)}{bd} \nde5.32 = \frac{-(ac)}{bd}$. By
\er{5.19}, \er{5.3}, $q\cdot r=-\frp{ac}{bd}$. Hence $\th(q\cdot r)\nde5.97 =
\frac{-(ac)}{bd}=\th(q)\hcd \th(r)$.

This completes the proof of \E\Pr\ \rf{p5.37}.
\endproof

\brm5.38
In \E\df\ \rf{d5.7}, the field $\wh X$ introduced in Theorem \rf{t5.6} was
called the field $\Q$ of \ra\ \nm s. It is customary (see \cite{Nrs},
\cite{Alg}) to call~$\wh \Z$ the field of \ra\ \nm s and to denote it by~$\Q$.
In view of \E\Pr\ \rf{p5.37} we shall call both $\wh X$ and~$\wh \Z$ \rp
ations of the field of \ra\ \nm s and denote them by~$\Q$. The context
indicates which \rp ation is used.
\erm

We now show that the field of \ra\ \nm s is a \Pf.

\blm5.39
Let $\Q$ be the field introduced in Theorem \rf{t5.6}, and let $F$~be a
subfield of~$\Q$ $($see \E\df\ \rf{d4.23}$)$, then $F=\Q$.
\elm

\proof
By \E\df\ \rf{d4.23}\,(a), $F$ is a subgroup of the \ag\ $(\Q,+,0)$. \E\Ip
$F\sbs\Q$ and $0\in F$. We claim that it suffices to show that $\{+q: q\in
\Q_{>0}\}\sbs F$. Indeed, note that for all $q\in\Q_{>0}$, $-q$~is the inverse
of~$+q$ in the group $(\Q,+,0)$ by \er{5.17}, \er{5.16}. Since for all $q\in
\Q_{>0}$, the inverse of~$+q$ in $(\Q,+,0)$ is the same as the inverse of~$+q$
in the group $(F,+,0)$ by Lemma \rf{l3.3}, it follows that $-q\in F$ for all
$q\in\Q_{>0}$. Therefore, if $\{+q:q\in\Q_{>0}\}\sbs F$, then $\Q=\{+q: q\in
\Q_{>0}\}\cup\{0\}\cup \{-q:q\in\Q_{>0}\}\sbs F$. \csq, $F=\Q$, and the claim
is proved. We next show that $+q\in F$ \fa $q\in\Q_{>0}$, that is, $+\frac ab
\in F$ \fa $a,b\in\Na$. We claim that it suffices to show that $+\frac a1\in F$
\fa $a\in\Na$. Observe that if $+\frac a1\in F$, then $+\frac 1a\in F$.
Indeed, $(+\frac a1)\cdot(+\frac 1a)\nde5.19 = +(\frac a1\cdot\frac1a)\nde5.3
= +\frac{a1}{1a}\nad{\er{2.12},\er{5.2}}= +\frac11\nde5.4 = +1$ \fa
$a\in\Na$, hence $+\frac1a$ is the inverse of $+\frac a1$ in the group
$(\Q_{\ge0}\sms0,\cdot,+1)$. By \E\df\ \rf{d4.23}\,(c) and Lemma \rf{l3.3}, we
conclude that $+\frac1a$ is also the inverse of $+\frac a1$ in the group
$(F\sms0,\cdot,+1)$, hence \Ip $+\frac1a\in F$. Now let $a,b\in\Na$ and $q:=
\frac ab$. Then $+q=+\frac ab\nde5.19 = (+\frac a1)\cdot(+\frac 1b)\in F$
since both $+\frac a1$ and $+\frac 1b$ belong to $F\sms0\sbs F$. Thus the
claim is proved.

\csq, it remains to show that $+\frac a1\in F$ \fa $a\in\Na$. We prove it by
\In\ on $a\in\Na$. Set $A:=\{a\in\Na: +\frac a1\in F\}$.

``$1\in A$'': Since $(F,\cdot,+1)$ is a \sbm\ of the field $\Q$ by \E\df\
\rf{d4.23}\,(b).

``\ti{$a\in A$ implies $a+1\in A$}'': Suppose $a\in A$. Then $+\frac{a+1}1
\nde5.25 = +(\frac a1+\frac 11)\nde5.16 = j(\frac a1+\frac 11)\nad*= j(\frac
a1) + j(\frac 11) = (+\frac a1)+(+\frac11)\in F$ since both $+\frac a1$ and
$+\frac11$ belong to~$F$ being a \sbm\ of $(\Q,+,0)$ by \E\df\
\rf{d4.23}\,(a). Hence $a+1\in A$, and $A$ is \iv\ in $(\Na,\le)$. Thus
$A=\Na$. In $\nad*=$ we used the fact that $j:(\Q_{\ge0},+,0) \to (\Q,+,0)$
is a  \hm sm.
\endproof

The next lemma implies that $\wh \Z$ is a \Pf.

\blm5.40
Let $F,F'$ be $($ring-$)$ \is c fields. If $F$ is prime, then so is~$F'$.
\elm

\proof
Let $\si:F\to F'$ be a ring-\is sm, and let $H'$ be a subfield of~$F'$. Then
$\si\Inv(H')$, the image of~$H'$ under $\si\Inv$, is a subfield of~$F$ by
Lemma \rf{l4.30}\,(i) (see also Lemma \rfa5{l4.22}). Since $F$~is a \Pf,
$\si\Inv(H')=F$, and $\si$ being a bi\jn\ we obtain $H'=\si(\si\Inv(H'))=
\si(F)=F'$.
\endproof

\bth5.41
Let $F$ be an \emph{infinite} \Pf, then $F$ is ring-\is c to the field~$\Q$.
If $F'$~is an infinite \Pf, then $F$ and~$F'$ are ring-\is c.
\eth

\proof
Let $(F,\qu_0,\qu_1,e_0,e_1)$ be an infinite \Pf. We first prove that the
field~$F$ is ring-\is c to the field~$\wh\Z$ introduced in Theorem \rf{t5.35}.
It suffices to show the \ex\ of an in\jc\ ring-\hm sm $j:\wh\Z\to F$. Indeed,
in this case $j(\wh\Z)$ is a subfield of~$F$ by Lemma \rfa5{l4.22}. Since
$F$~is a \Pf, we infer that $j(\wh\Z)=F$, hence $j$~is bi\jc, thus a ring-\is sm.
We first suppose $\vf:\wh\Z\to F$ is an in\jc\ ring-\is sm. Then
\bea5.101
\vf(\wh0)=e_0,&\q \vf(\wh x\hpl\wh y)=\vf(\wh x)\qu_0\vf(\wh y),\ \wh x,\wh y\in\wh\Z,\\
\vf(\wh1)=e_1,&\q \vf(\wh x\hcd\wh y)=\vf(\wh x)\qu_1\vf(\wh y),\ \wh x,\wh y\in\wh\Z.
\lb{5.102}
\e
By \era2{2.3}\,I1 and Theorem \rf{t4.35}\,(ii) we obtain
\beq5.103
n{\mathrel{\lower1pt\hbox{$\stackrel{\lower7pt\hbox{\smash{\hbox
to0pt{$\scriptstyle\sqcap$\hss}$\scriptstyle\sqcup_0$}}}\cdot$}}}e_1\ne e_0
\q\hbox{\fa} n\in\Na,
\e
since $F$ is \ti{infinite}. Since $\vf:(\wh\Z,\hpl,\wh0)\to (F,\qu_0,e_0)$ is
a \hm sm, we obtain from \era2{1.45}:
\beq5.104
\vf(n\dhp \wh1) = n\dqum0 e_1, \q n\in\N.
\e

We now prove by \In\ over $n\in\N$ that $n\dhp\wh1=\frac n1$, $n\in\N$. Recall that
$\wh0=\frac01$, $\wh1=\frac11$ by \er{5.86} where $D=\Z$ and $D^\t=\Z^\t$, and
where $\Z_{\ge0}$ is identified with~$\N$. We have $0\dhp1
\nad{\era2{2.3}\,\rm I0}=\wh0=\frac01$. Suppose
$m\dhp 1=\frac m1$ \fs $m\in\N$, then $(m+1)\dhp1 \nad{\era2{2.3}\,\rm I2}=
(m\dhp\wh1)\hpl(1\dhp\wh1)\nad{\era2{2.3}\,\rm I1}= \frac m1\hpl\frac11 \nde5.84 = \frac{m+1}1$.
\E\Tf we have:
\beq5.105
\vf\Bg(\frac n1)=n\dqum0 e_1, \q n\in\N.
\e
For $\wh z\in\wh \Z$ we denote by $\wh z{}\Inv$ the inverse of~$\wh z$ in
the group $(\wh\Z,\hpl,\wh0)$. \E\Ip $(\frac n1)\Inv=\frac{-n}1$, $n\in\N$,
since $\frac n1\hpl\frac{-n}1 \nde5.90 = \frac{n+(-n)}1= \frac01=\wh 0$. Hence
$e_0\nde5.101 = \vf(\wh0) = \vf(\frac n1\hpl\frac{-n}1)\nde5.101 = \vf(\frac n1)
\qu_0\vf(\frac{-n}1)$. \E\Tf $\vf(\frac{-n}1)=\vf(\frac n1)\Inv$, the inverse
of $\vf(\frac n1)$ in the group $(F,\qu_0,e_0)$. By \er{5.105} $\vf(\frac{-n}1)
= (n\dqum0 e_1)\Inv$, $n\in\N$. Note that the map $\wh z\mt n\dqum0 \wh z$,
$n\in\N$, from $\wh\Z$ into~$\wh\Z$ is an endo\mf\ of the group $(\wh\Z,+,\wh0)$ by
\era2{2.3}\,I0, I2. We obtain
\beq5.106
(n\dqum0 \wh z)\Inv \nda21.45 = n\dqum0(\wh z{}\Inv)\nde3.141 = (-n)\dqum0
\wh z, \q n\in\N,\ \wh z\in F.
\e
Observe that $(-n)\dqumz e_1\ne e_0$, $n\in\Na$. Indeed, if $(-n)\dqumz e_1
=e_0$, then $e_0=e_0\Inv = ((-n)\dqum0 e_1)\Inv\nde5.106 = ((n\dqum0 e_1)\Inv
)\Inv = n\dqu e_1=e_0$, \cd ing \er{5.103} when $n\in\Na$.
Hence $\vf(\frac{-n}1)=(-n)\dqum0 e_1$, $n\in\N$. From \er{5.105} we obtain
\beq5.107
\vf\Bg(\frac z1) = z\dqum0 e_0\in F\sms{e_0}, \q z\in\Z^\t.
\e
If $w\in\Z^\t$, then $\frac w1\hcd \frac1w\nde5.85 = \frac{w\cdot1}{1\cdot w}
\nde5.83 =\frac 11\nde5.86 =\wh 1$. Hence $\frac1w$ is the inverse of $\frac w1$
in the group $(\wh\Z\sms{\wh0},\hcd,\wh1)$. By \er{5.102} we have $e_1=\vf(\wh1)
\nad{\er{5.85},\er{5.83}}= \vf(\frac w1\hcd\frac1w)\nde5.102 = \vf(\frac w1)
\qu_1 \vf(\frac1w)$. Hence $\vf(\frac1w)=(\vf(\frac w1))^{-1}$ the inverse of
$\vf(\frac w1)$ in the group $(F\sms{e_0},\qu_1,e_1)$. From \er{5.107} we
obtain $\vf(\frac1w)=(w\dqum0 e_1)^{-1}$, $w\in\wh\Z^\t$. Finally, if $z\in\wh\Z$
and $w\in\Z^\t$, we have $\vf(\frac zw)
\nde5.85 = \vf(\frac z1\cdot\frac1w) \nde5.102 = \vf(\frac z1)\qu_1\vf(\frac1w)
=(z\dqum0 e_1)\qu_1(w\dqum0 e_1)^{-1}$. This motivates the \fw\ \df.

Let $\wh z\in\wh\Z$ and let $(z,w)\in\Z\t\Z^\t$ be \st $\wh z=\frac zw$. Set
\beq5.108
j(\wh z):=(z\dqum0 e_1)\qu_1 (w\dqum0 e_1)^{-1}.
\e
This \df\ makes sense if $\frac zw=\frac{z'}{w'}$ with $(z',w')\in\Z\t\Z^\t$
implies
\beq5.109
(z\dqum0 e_1)\qu_1 (w\dqum0 e_1)^{-1}=(z'\dqum0 e_1)\qu_1 (w'\dqum0 e_1)^{-1}.
\e
Note that \fa $(z,w),(z',w')\in \Z\t\Z^\t$ \er{5.109} is \ev t to
\beq5.110
(z\dqum0 e_1)\qu_1 (w'\dqum0 e_1)=(z'\dqum0 e_1)\qu_1 (w\dqum0 e_1).
\e

For the proof of \er{5.110} we use the \fw\ identities:
\beq5.111
mn\dqum0 e_1 = (m\dqum0 e_1)\cdot(n\dqum0 e_1), \q m,n\in\Na,
\e
and
\beq5.112
x\Inv \qu_1 y=(x\qu_1 y)\Inv,\q x\Inv\qu_1 y\Inv = x\qu_1 y,\q x,y\in F,
\e
where $x\Inv$ denotes the inverse of $x$ in $(F,\qu_0,e_0)$.

\bex5.42
Prove \er{5.111} by \In\ on $n$ and \er{5.112}.
\eex

We claim
\beq5.113
zw\dqum0 e_1 = (z\dqum0 e_1)\qu_1(w\dqum0 e_1),\q z,w\in\Z^\t.
\e
If $\sgn(zw)={+}$, then by \er{5.39}, \er{5.22} either $\sgn(z)=\sgn(w)={+}$
or $\sgn(z)={\sgn(w)={-}}$. In the first case $z=m$, $w=n \in\Na$, hence
\er{5.113} follows from \er{5.111}. If $\sgn(z)=\sgn(w)={-}$, then $z=-m$,
$w=-n$ \fs $m,n\in\Na$. Then $zw=mn$ by \er{5.22}. Hence $(zw\dqum0 e_1)
\nde5.111 = (m\dqum0 e_1)\cdot(n\dqum0 e_1)\nde5.112 = (m\dqum0 e_1)\Inv
\qu_1 (n\dqum0 e_1)\Inv \nde3.17 = (m\dqum0(e_1)\Inv)\cdot (n\dqum0(e_1)\Inv)
\nde3.142 = ((-m)\dqum0 e_1)\cdot((-n)\dqum0 e_1) =(z\dqum0 e_1)\cdot
(w\dqum0 e_1)$. If $\sgn(zw)={-}$, by the \cmt ity of~$\cdot$, it suffices
to consider $z=m$, $w=-n$, $m,n\in\Na$. Then $zw=-(mn)$ by \er{5.22}. Hence
$(zw\dqum0 e_1)=(-(mn)\dqum0 e_1)\nde3.142 =\break (mn\dqum0 (e_1)\Inv)\nde3.17 =
(mn\dqum0 e_1)\Inv \nde5.111 = ((m\dqum0 e_1)\qu_1(n\dqum0 e_1))\Inv \nde5.112
= (m\dqum0 e_1)\Inv \qu_1(n\dqum0 e_1)\nde3.17 = (m\dqum0 (e_1)\Inv)\qu_1
(n\dqum0 e_1)\nde3.142 = ((-m)\dqum0 e_1)\qu_1(n\dqum0 e_1)=(z\dqum0
e_1)\cdot(w\dqum0 e_1)$. The proof of the claim is complete.

\Wanp prove \er{5.110}. Let $z,z',w,w'\in\Z^\t$. Then $(z\dqum0 e_1)\qu_1
(w'\dqum0 e_1)\nde5.113 = (zw'\dqum0 e_1)=(wz'\dqum0 e_1)\nde5.113 = (w\dqum0
e_1)\qu_1(z'\dqum0 e_1) \nda21.6 = (z'\dqum0 e_1)\qu_1(w\dqum0 e_1)$, which
proves \er{5.110} for $z,w,z',w'\in\Z^\t$ \sf ying $zw'=wz'$. We now consider
the case where $z=0$. Note that $(0,w)\sim(z',w')$, $w,w'\in\Z^\t$, $z'\in\Z$,
iff $z'=0$ by Lemma \rf{l5.33} with $D:=\Z$. We have $zw=0=z'w'$. Hence
$(0\dqum0 e_1)\qu_1(w'\dqum0 e_1)\nad{\era2{2.3}\,\rm I0}= e_0\qu_1(w'\dqum0
e_1)\nde4.1 = e_0=(0\dqum0 e_1)\qu_1(w'\dqum0 e_1)$. This completes the proof
of \er{5.110}. \E\Tf the map $j:\wh\Z\to F$ is well-defined by \er{5.108}.

We next show that $j$ is a ring-\hm sm. We first prove \er{5.102}.

``$j(\wh1)=e_1$'': $j(\wh1)=j(\frac11)\nad{\era2{2.3}\,\rm I1}= (1\dqum0e_1)\cdot(1\dqum0e_1)^{-1}=
e_1\cdot(e_1)^{-1}=e_1$.

``$j(\wh x\hcd\wh y)=j(\wh x)\qu_1 j(\wh y)$, $\wh x,\wh y\in\wh\Z$'':
Let $a,c\in\Z$, $b,d\in\Z^\t$ be \st $\wh x=\frac ab$, $\wh y=\frac cd$. Then
$j(\wh x\hcd\wh y)=j(\frac{ac}{bd})\nde5.108 = (ac\dqum0e_1)\qu_1(bd\dqum0e_1)
^{-1}\nde5.113 = \bigl((a\dqum0e_1)\qu_1(c\dqum0e_1)\bigr)\qu_1\bigl((b\dqum0
e_1)\qu_1(d\dqum0e_1)\bigr)^{-1} \nad{\er{3.4},\era2{1.6}}= \bigl((a\dqum0e_1)
\qu_1(c\dqum0e_1)\bigr)\qu_1\bigl((b\dqum0e_1)^{-1}\qu_1(d\dqum0e_1)^{-1}\bigr)
\nad{\era2{1.36},\era2{1.6}}= (a\dqum0e_1)\qu_1\bigl((b\dqum0e_1)^{-1}\qu_1
(c\dqum0e_1)\bigr)\qu_1(d\dqum0e_1)^{-1} \nda21.36 = \bigl((a\dqum0e_1)\qu_1
(b\dqum0e_1)^{-1}\bigr)\qu_1\bigl((c\dqum0e_1)\qu_1(d\dqum0e_1)^{-1}\bigr)=
j(\frac ab)\qu_1j(\frac cd) = j(\wh x)\qu_1(\wh y)$.

We next prove \er{5.101}.

``$j(\wh0)=e_0$'': $j(\wh0)=j(\frac01)\nde5.108 = (0\dqum0e_1)\cdot(1\dqum0
e_1)^{-1}\nad{\era2{2.3}\,\rm I0,I1}= e_0\cdot e_1^{-1}\nde4.1 = e_0$. Note
that $e_1\Inv=e_1$ since $e_1\qu e_1=e_1$.

``$j(\wh x\hpl\wh y)=j(\wh x)\qu_0 j(\wh y)$, $\wh x,\wh y\in\wh\Z$'': Let
$a,c\in\Z$, $b,d\in\Z^\t$ be \st $\wh x=\frac ab$, $\wh y=\frac cd$. Then
$j(\wh x+\wh y)\nde5.84 = j(\frac{ad+bc}{bd})\nde5.85 = j(\frac{ad+bc}1
\cdot\frac1{bd})\nde5.102 = j(\frac{ad+bc}1)\qu_1 j(\frac1{bd})\nde5.108 =
\bigl(((ad+bc)\break\dqum0e_1)\qu_1(1\dqum0e_1)^{-1}\bigr)\qu_1 (1\dqum0e_1)\qu_1
((bd)\dqum0e_1)^{-1} = \bigl(((ad+bc)\dqum0e_1)\qu_0e_1\bigr)\qu_1
(e_1\qu_0(bd\dqum0e_1)^{-1}) = ((ad+bc)\dqum0e_1)\qu_1(bd\dqum0e_1)^{-1}
\nad{\era2{2.3}\,\rm I2}= ((ad\dqum0e_1)\qu_0(bc\dqum0e_1)) \qu_1 ({(b\dqum0e_1)}
\cdot(d\dqum0e_1))^{-1} = ((ad\dqum0e_1)\qu_1((b\dqum0e_1)\qu_1(d\dqum0e_1))
^{-1}) \qu_0 \bigl((bc\dqum0e_1)\qu_1 ({(b\dqum0e_1)}\qu_1(d\dqum0e_1))^{-1}\bigr)$.

Observe $(ad\dqum0e_1)\qu_1((b\dqum0e_1)\qu_1(d\dqum0e_1)^{-1})\nad{\er{5.113},
\er{3.4}} = (a\dqum0e_1)\qu_1(b\dqum0e_1)^{-1} = j(\frac ab)= j(\wh x)$.
Similarly $(bc\dqum0e_1)\qu_1((b\dqum0e_1)\qu_1(d\dqum0e_1))^{-1} = ((b
\dqum0e_1)\qu_1{(c\dqum0e_1)})\qu_1 ((b\dqum0e_1)\qu_1(d\dqum0e_1))^{-1} =
(c\dqum0e_1)\qu_1(d\dqum0e_1)^{-1} = j(\frac cd)=j(\wh y)$. \E\Tf $j(\wh x
+\wh y)=j(\wh x)\qu_0 j(\wh y)$.

We now prove that $j$ is \ti{in\jc}. Since $j:(\wh\Z,\hpl,\wh0)\to(F,+,0)$ is
a \hm sm, it suffices to show that $\{\wh z\in\wh \Z: j(\wh z)=e_0\}=\wh0$, in
view of \E\Pr\ \rf{p3.10}\,(vi). Let $\wh z=\frac ab$, $a\in\Z$, $b\in\Z^\t$
be \st $j(\wh z)=e_0$, that is, $(a\dqum0e_1)\qu_1(b\dqum0e_1)^{-1} =e_0$ by
\er{5.108}. Hence $a\dqum0e_1=a\dqum0e_1 \qu_1e_1 =
a\dqum0e_1 \qu_1((b\dqum0e_1)^{-1}\qu_1(b\dqum0e_1)) = (a\dqum0e_1\qu_1
(b\dqum0e_1)^{-1})\qu_1 b\dqum0e_1 = e_0\qu_1 b\dqum0e_1\nde4.1 = e_0$.
Hence $a\dqum0e_1 = e_0$. In view of \er{5.103} $a=0$ since $0\dqum0e_1=e_0$.
\If that $\wh z=\frac0b= \frac01
=\wh0$. This completes the proof that $j$ is an in\jc\ ring-\hm sm. In view of
the argument given at the beginning of the proof, $j$~is a ring-\is sm.

Let $F'$ be an infinite \Pf\ and let $j':\wh \Z\to F'$ be a ring-\is sm whose
\ex\ is guaranteed by the first part of the proof. Then $j'\circ j\Inv:F\to F'$
is a ring-\is sm by Lemma \rfa2{l1.8}.
\endproof

Our next goal is to show that the field $\Q$ introduced in Theorem \rf{t5.6}
is \ct y infinite. We first show that $\Q_{>0}$ the \mlv\ group of \ra\ \nm s
introduced at the beginning of this section is \ti{\ct y infinite}. We recall
that $\N\t\N$ is \ct y infinite (see Remark \rfa3{r3.31n}).
A proof based on Schr\"oder--Bernstein's theorem is given just before Exercise
\rfa3{ex3.28}. We also recall that $\Q_{>0}$ is the \qt\ monoid $\Na\t\Na/
{\sim}$ where $(a,b)\sim(c,d)$ if $ad=bc$.

\blm5.43
Let $a,b\in\Na$. Then \te s \ooo pair $(c,d)\in\Na\t\Na$ \st
\beq5.114
(a,b)\sim(c,d) \qh{and } \gcd(c,d)=1.
\e
\elm

\proof
``\ti{\E\ex}'': Note that $\gcd(a,b)$ divides $a$ and~$b$. Set $a':=\frac a
{\gcd(a,b)}$, $b':=\frac b{\gcd(a,b)}$. Then $\gcd(a,b)\cdot\gcd(a',b')=
\gcd(a,b)\cdot\gcd(\frac a{\gcd(a,b)},\frac b{\gcd(a,b)})\nda31.27 =
\bigl(\gcd(a,b)\cdot\frac a{\gcd(a,b)}\cdot\gcd(a,b)\cdot\frac b
{\gcd(a,b)}\bigr)\nad{\era3{8.53},\era2{2.18}}=\gcd(a,b)\cdot1$. Hence
$\gcd(a',b')=1$ by \cnc ity. \Mo $ab'=a\frac b{\gcd(a,b)}\nde5.4 = \frac a1
\cdot \frac b{\gcd(a,b)}\nde5.3 = \frac{ab}{\gcd(a,b)}$ and $ba'=b\cdot
\frac a{\gcd(a,b)}=\frac{ba}{\gcd(a,b)}=\frac{ab}{\gcd(a,b)}$. Hence
$ab'=ba'$, that is, $(a,b)\sim(a',b')$.

``\ti{\E\uq}'': Let $a,b,c,d\in\Na$ be \st $(a,b)\sim(c,d)$, $\gcd(a,b)=
\gcd(c,d)=1$. We have $a|ad$, hence $a|bc$ since $ad=bc$. By Lemma
\rfa3{l1.29}, $a|c$ since $\gcd(a,b)=1$. Similarly, $c|bc$ hence $c|ad$.
Since $\gcd(c,d)=1$, we have $c|a$ by the same lemma. From $a|c$ and $c|a$ we
infer $a=c$. Thus $da=ad=ba$. By \cnc ity, $b=d$. Hence $(a,b)=(c,d)$.

\bigskip
Let $V:=\{(c,d)\in\Na\t\Na: \gcd(c,d)=1\}$. Let $z\in\Q_{>0}$. By \df\ \te s
$(a,b)\in \Na\t\Na$ \st $z=\frac ab$. By the \ex\ part of Lemma \rf{l5.43},
\te s $(c,d)\in V$ \st $\frac ab=\frac cd$. \Mo if \te s $(c',d')\in V$ \st
$z=\frac{c'}{d'}$, then $(c,d)\sim(c',d')$ by \er{5.2}, hence $(c,d)=(c',d')$
by the \uq\ part of Lemma \rf{l5.43}. \E\Tf \te s a map $f:\Q_{>0}\to V$ \st
$f(z)$ is the only pair $(c,d)\in V$ \sf ying $z=\frac cd$. Let $g:V\to
\Q_{>0}$ be defined by $g((c,d)):=\frac cd$ where $(c,d)\in V$. Clearly
$g\circ f=\id_{\Q_{>0}}$ and $f\circ g=\id_V$. \If that the maps $f$ and~$g$
are bi\jc. Let $j$ denote the inclusion map from~$V$ into $\N\t\N$, then the
map $j\circ f:\Q_{>0}\to \N\t\N$ is in\jc, as a \cm\ of two in\jc\ maps. From
\E\Pr\ \rfa1{p4.27}\,(iv) we infer that $\Q_{>0}$ is \ti{\ct e}. \E\oh the
map $i:\Na\to\Q_{>0}$ defined by $i(n):=\frac n1$ is in\jc. Indeed, from
$i(n)=i(m)$, $n,m\in\Na$, we infer $n=n\cdot 1=1\cdot m=m$ by \er{5.2}.
\If from \E\Pr\ \rfa1{p4.27}\,(i) that $\Q_{>0}$ is \ti{infinite}. Hence
$\Q_{>0}$~is \ti{\ct y infinite}, by \E\df\ \rfa1{d4.21}. Since $\N\approx\Na$ (see
\E\df\ \rfa1{d4.6}), \te s a bi\jn\ ${h:\Q_{>0}\to\Na}$. We now define a map $F:\Q\to
\Z$ by setting
\beq5.115
F(q):=\bca
+h(|q|) &\hbox{if }\sgn(q)={+},\\
0 &\hbox{if }q=0,\\
-h(|q|) &\hbox{if }\sgn(q)={-}.
\eca
\e
We claim that $F$ is a bi\jn. We define $G:\Z\to\Q$ by setting
\beq5.116
G(z):=\bca
+h\Inv(+m) &\hbox{if }z=+m,\ m\in\Na,\\
0 &\hbox{if }z=0,\\
-h\Inv(m) &\hbox{if }z=-m,\ m\in\Na.
\eca
\e
Then $G\circ F=\id_\Q$ and $F\circ G=\id_{\Z}$. \E\Tf both $F$ and~$G$ are
bi\jc\ maps.

Finally, we define a map $H:\Z\to\N$ by setting
\beq5.117
H(z):=\bca
2m &\hbox{if }z=+m,\ m\in\Na,\\
0 &\hbox{if }z=0,\\
2m-1 &\hbox{if }z=-m,\ m\in\Na.
\eca
\e
One verifies that the map $H$  is bi\jc. Hence the map $H\circ F:\Q\to\N$ is
a bi\jn, thus $\Q$ is \ct y infinite.
\endproof

\advance\abovedisplayskip by-2pt
\advance\belowdisplayskip by-2pt
In view of Theorem \rf{t5.41} we obtain

\bpr5.44
An infinite \Pf\ is \emph{\ct y infinite}.
\epr

We conclude this section by introducing an \og~$\le$ on the field~$\Q$
introduced in Theorem \rf{t5.6} and by investigating some of its \pp ies.

\bpr5.45
Let $\Q$ be as in \E\df\ \rf{d5.7}, and let $\le$ denote the \rl\ on~$\Q$
defined by
\beq5.118
r\le s \q\ \hbox{if}\q s+r\Inv\in\Q_{\ge0}, \q\ s,r\in\Q,
\e
where $r\Inv$ is the inverse of $r$ in the group $(\Q,+,0)$. Then $\le$~is a
\emph{total \og} on~$\Q$ $($see \E\df\ \rfa1{d3.19}$)$. \Mo for $r\in\Q \sms0$,
\beq5.119
r\in\Q_{\ge0} \qh{iff } \sgn r={+}.
\e
\epr

\begin{exe}[see Appendix A]\label{\sss ex5.46} \

\hph i,i, Prove \E\Pr\ \rf{p5.45} (see the discussion \fw\ the proof of Lemma 3.21 of Appendix~A).

\hph ii,, Show that \fa $x,y,z\in\Q$
\bea5.120
x<y & \hbox{ implies } x+z<y+z,\\
x<y & \hbox{ and }z>0 \hbox{ implies }x\cdot z<y\cdot z.\lb{5.121}
\e
\eex

\bdf5.47
A field $(F,+,\cdot,0,1)$ is called an \ti{ordered field\/}\index{field!ordered} if \te s a total
\og~$\le$ on~$F$ \sf ying \er{5.120} and \er{5.121} \fa $x,y,z\in F$.
\edf

\bex5.48
Let $(F,+,\cdot,0,1,\le)$ be an ordered field. Set $P:=\{x\in F: x\ge0]$,
$\dot P:=P\sms0$. For $x\in F$ let $-x$ denote the inverse of~$x$ in
$(F,+,0)$ and $x\mo$ the inverse of~$x$ in $(F\sms0,\cdot,1)$. Let $-P:=
\{-x: x\in P\}$ and $-\dot P:=\{-x:x\in\dot P\}$. Show that $F$ is the
disjoint union of $\dot P$, $\{0\}$ and~$-\dot P$. Prove the \fw: Let $x,y\in F$,
then
\bea5.122
&x\le y \hbox{ (resp.\ }x<y) \hbox{ iff }y+(-x)\in P \hbox{ (resp.\ }\dot P),\\
&x\in\dot P,\ y\in-\dot P \hbox{ implies }x\cdot y\in -\dot P, \lb{5.123}\\
&x\in-\dot P,\ y\in-\dot P \hbox{ implies }x\cdot y\in \dot P, \lb{5.124}\\
&1 = (1\cdot 1) = (-1)\cdot(-1)\in \dot P,\lb{5.125}\\
&n\dpl1\in\dot P \hbox{ \fa} n\in\Na, \hbox{ \Ip the \ch istic}\lb{5.126}\\
&\hskip40pt \hbox{ of an ordered field is zero,}\nonumber\\
&\hbox{If }x\in \dot P, \hbox{ then }x\mo\in\dot P.\lb{5.127}
\e
Show that the \og\ introduced in \E\Pr\ \rf{p5.45} is the only \og\ which
makes $\Q$ an ordered field.
\eex

Show that \fa $a,b\in\dot P$ we have: $a+b\in\dot P$, $a\cdot b\in\dot P$,
$-a\notin\dot P$ and  \fa $a\in F$ we have $a=0$ or $a\in\dot P$ or $-a\in
\dot P$. Compare \E\df\ \rf{d5.47} with the \df\ given in \cite[p.~232]{Alg}.

We now consider some \pp ies of the \tos\ $(\Q_{\ge0},\le)$ where $\le$
denotes the \rt ion to $\Q_{\ge0}$ of the order on~$\Q$ introduced in \E\Pr\
\rf{p5.45}. Observe that this \rt ion is identical to the \nog\ of the \PM\
$(\Q_{\ge0},+,0)$ (see \E\Pr\ \rf{p5.4}). Since $q+0=q$, $q\in\Q_{\ge0}$,
$0$~is the least \el\ of $(\Q_{\ge0},\le)$. In contrast to the \os\ $(\N,\le)$,
$(\Q_{\ge0},\le)$ is not well-ordered. Indeed, the set $A:=\{\frac1n\in
\Q_{\ge0}: n\in\Na\}$ is bounded below by~$0$, and $A$ has no least \el.
Suppose, for \cd ion, that \te s $\ov n\in\Na$ \st $\frac1{\,\ov n\,}\le
\frac1n$ \fa $n\in\Na$, that is, $n\le\ov n$ \fa $n\in\Na$ by \er{5.15}. Then
$\ov n+1\le \ov n$, a~\cd ion.

Recall that if $a<b$, $a,b\in\Na$, then \te s $n\in\Na$ \st $n\cdot a>b$.
Indeed, by Theorem \rfa2{t1.38} and \era2{2.8} \te s a pair $(q,r)\in\N\t[0,a)$
\sf ying $b=q\cdot a+r$. Since $r<a$ we have $q\cdot a+r<q\cdot a+a = (q+1)
\cdot a$. Thus $b<n\cdot a$ with $n:=q+1$.

\blm5.49
Let $x<y$, $x,y\in\Q_{>0}$. Then \te s $n\in\Na$ \st $\frac n1\cdot x>y$.
\elm

\proof
Let $a,b,c,d\in\Na$ be \st $x=\frac ab$, $y=\frac cd$. In view of \er{5.15}
we have $ad<bc$. By what precedes \te s $n\in\Na$ \st $nad>bc$.
Hence $\frac{na}b>\frac cd$ by \er{5.15}. Thus $\frac n1\cdot x= \frac n1\cdot
\frac ab\nde5.3 = \frac{na}{1b}= \frac{na}b > \frac cd=y$ by \er{5.15}. Hence
$\frac n1\cdot x>y$.
\endproof

Note that $\frac n1\cdot x= n\dpl x$ (use \In\ on $n\in\Na$. See also the proof of Theorem 4.30 in Appendix A).

\begin{dfn}[\cite{Alg}]\lb{d5.50}
An ordered field $(F,+,\cdot,0,1,\le)$ is called \ti{archimedean}\index{archimedean ordered field} if \fe \el\
$x\in F$ \sf ying $x>0$, \te s $n\in\Na$ \st $n\dpl 1>x$.
\edf

\blm5.50
Let $(F,+,\cdot,0,1,\ge)$ be an ordered field and let $\dot P:=\{x\in F: x>0\}$.
The \fw\ assertions are \ev t\/{\rm:}

\hph i,ii, \E\fe $x\in \dot P$ \te s $n\in\Na$ \st $n\dpl 1>x$.

\hph ii,i, \E\fe $y\in\dot P$ \te s $m\in\Na$ \st $(m\dpl1)\mo<y$.

\hph iii,, \E\fe pair $(x,y)\in\dot P\t \dot P$ \sf ying $x<y$ \te s $k\in\Na$
\st $(k\dpl x)>y$.
\elm

\bex5.51
Prove Lemma \rf{l5.50}. See the discussion preceding Definition 3.36 in Appendix~A.
\eex

We now return to $A:=\{\frac1n\in\Q_{\ge0}: n\in\Na\}$. Clearly $0$ is a \lo\
of~$A$. We claim that $0$~is the only \lo\ of~$A$ in~$\Q_{\ge0}$. Suppose for
\cd ion that \te s $a\in\Q_{>0}$ \st $a<\frac1n$ \fa $n\in\Na$. Then by
\er{5.121} $n\cdot a< n\cdot\frac1n=1$. \E\Ip $a<1$. However, by Lemma
\rf{l5.49} \te s $\ov n\in\Na$ \st $1<\ov n\cdot a$, hence $\ov n\cdot a<1
<\ov n\cdot a$, a \cd ion.

\E\Tf $0$~is the greatest \lo\ of~$A$ or the infimum of~$A$,
$0=\inf A$ (see \E\df\ and Notations \rf{n1.13n}). Note that $0$~is the first
example so far of an infimum of a subset of a \tos\ which is not the least
\el\ of this subset.

\bex5.53
Let $x\in\Q$, and let $\le$ denote the total \og\ of~$\Q$ introduced in \E\Pr\
\rf{p5.45}. Prove the \fw:
\bea5.138
x&= \inf\{y\in\Q: x<y\},\\
x&= \sup\{y\in\Q: y<x\}. \lb{5.139}
\e
(\ti{Hint\/}: if $a,b\in\Q$, $a<b$, then $a<\frac12(a+b)<b$.)
\eex

In what follows we give an example of a \nss\ of~$\Q_{\ge0}$ which is bounded
above and which has no supremum.

Let $[0,1]:=\{x\in\Q: 0\le x\le 1\}$ and let $f:[0,1]\to[0,1]$ be \st $f(0)=0$,
$f(1)=1$ and $f(x)<f(y)$ whenever $0\le x<y\le 1$. Let  $\a\in[0,1]\sms{0,1}$
and set $A:=\{x\in[0,1]:f(x)<\a\}$, $B:=\{y\in[0,1]: \a<f(y)\}$. Then $0\in A$
and $1\in B$. \Mo every $x\in A$ is a \lo\ of~$B$ and every $y\in B$ is an \ub\
of~$A$. Suppose in \ad\ that \te s $M\in\Na$ \st $f(t)-f(s)\le M(t-s)$
whenever $0\le s<t\le 1$. We claim that $A$~has no greatest \el\ and $B$~has
no least \el. Suppose, for \cd ion, that \te s $\ov x\in A$ \st $x\le \ov x$
\fa $x\in A$. Since $\ov x\in A$, we have $f(\ov x)<\a$. Note that $\a<1=f(1)$.
\Mo $\ov x\ne1$, otherwise $1=f(1)=f(\ov x)<\a<1$. Thus $\ov x<1$. Set $h_1:=
1-\ov x>0$. We want to find $\ov h>0$ \st $\ov h<h_1$ and $f(\ov x+\ov h)<\a$,
so that $\ov x<\ov x+\ov h$ and $\ov x+\ov h\in A$.
If $0<h<h_1$, then $\ov x<\ov x+h<1$ and $f(\ov x+h)-f(\ov x)<Mh$ by \as. We claim that it
suffices to find $\ov h$ \st $M\ov h <\a-f(\ov x)$. Indeed, since $0<\a-f(\ov x)$, we may invoke Lemma
\rf{l5.49} with $x:=\a-f(\ov x)$, $y:=M$. \csq, \te s $\ov n\in\Na$ \st
$\ov n(\a-f(\ov x))>M$. Set $\ov h:=\frac1{\,\ov n\,}$. Then $M\cdot \ov h
=\ov h\cdot M \nde5.121 < \ov h\cdot \ov n(\a-f(\ov x)) = 1\cdot(\a-f(\ov x))
=\a-f(\ov x)$. \E\Tf $f(\ov x+\ov h)-f(\ov x)\le M\ov h<\a-f(\ov x)$, hence
$f(\ov x+\ov h)= (f(x+\ov h)-f(\ov x))+f(\ov x) \nde5.120 < (\a-f(\ov x))+
f(\ov x)=\a$. Thus $f(\ov x+\ov h)<\a$, and we have found $\ov h>0$ \st
$\ov x<\ov x+\ov h<1$ and $\ov x+\ov h\in A$. Hence $A$~has no greatest \el.
The proof that $B$~has no least \el\ is similar and is left as an exercise.

Suppose now that \te s $\xi\in(0,1)$ \st $f(\xi)=\a$. We claim $\xi=\sup A$.
Indeed, since $A$ has no greatest \el, \ti{no} \el\ of~$A$ is an \ub\ of~$A$.
Since $f(x)<\a=f(\xi)$ \fa $x\in A$ and $f$~is strictly in\cre, $\xi$~is an
\ub\ of~$A$ as well as all \el s of~$B$. Note that $[0,1]$ is the disjoint
union of $A$, $\{\xi\}$ and~$B$, and $\xi$~is the least \el\ of $\{\xi\}\cup
B$. \E\Tf $\xi$~is the least \ub\ of~$A$, that is, $\xi=\sup A$. Similarly one
shows that $\xi=\inf B$.

Finally, we suppose that there is \ti{no} $\xi\in(0,1)$ \st $f(\xi)=\a$. In that
case $A$~has \ti{no} supremum and $B$~has \ti{no} infimum, since the set of
\ub s of~$A$ is~$B$ (resp.\ \lo s of~$B$ is~$A$) and $A$ has no greatest \el\
(resp.\ $B$~has no least \el).

We have thus proved

\blm5.54
Let $[0,1]:=\{x\in\Q_{\ge0}: 0\le x\le 1\}$, let $M\in\Na$ and let $f:[0,1]\to
[0,1]$ \sf y
\beq5.140
0=f(0)\le f(x)<f(y)\le f(1)=1,
\e
and
\beq5.141
f(y)-f(x) \le M(y-x),
\e
for $0\le x<y\le 1$.

Let $\a\in(0,1)$, $A:=\{x\in[0,1]:f(x)<\a\}$, and $B:=\{y\in[0,1]: \a<f(y)\}$.
Then $\sup A$ $($resp.\ $\inf B)$ exists iff \te s $\xi\in(0,1)$ \st $f(\xi)=\a$.
In this case $\xi=\sup A$ $($resp.\ $\xi=\inf B)$.
\elm

\bex5.55
Prove that under the \as s of Lemma \rf{l5.54}, the set $B$ has no least \el,
and that $\xi=\inf B$ if $f(\xi)=\a$.
\eex

\bxs5.55
We choose $f(x):=x^n$, $x\in[0,1]$, $n\in\Na$, in Lemma \rf{l5.54}. If $n=1$ then
clearly $f$~\sf ies \er{5.140}, \er{5.141}. Assume $n\ge2$, and let $0\le x<y
\le 1$. Then $f(y)-f(x)=y^n-x^n \nde4.46 = (y-x)\suml_{k=0}^{n-1} y^{n-1-k}
x^k$. If $x=0$, then $x^n=x\cdot x^{n-1}=0\cdot x^{n-1}=0$, hence $f(0)=0$.
If $x=1$, $f(x)=1^n=1$. Since $y>0$ and $x\ge0$, $\suml_{k=0}^{n-1}y^{n-1-k}
x^k = y^{n-1}\cdot1 + \suml_{k=1}^{n-1}y^{n-1-k}x^k \ge y^{n-1} >0$. Hence
$f(y)-f(x)\ge (y-x)\cdot y^{n-1}>0$. This completes the proof of \er{5.140}.
Since $y^{n-1-k}x^k \le 1\cdot x^k\le 1$, we have $\suml_{k=0}^{n-1}y^{n-1-k}
x^k \le \suml_{k=0}^{n-1}1=n\cdot1=n$. Thus \er{5.141} holds with $M:=n$. Let
$\a:=\frac1p$, $p$~prime. Then $0<\a<1$. If $n:=1$ and $\xi:=\frac1p$, then
$f(\xi)=\xi=\frac1p=\a$. Hence $\sup\{x\in[0,1]: x<\frac1p\}=\frac1p$, and
$\inf\{y\in[0,1]:\frac1p<y\}=\frac1p$ as in Exercise \rf{ex5.53}.

We claim that there is no $\xi\in(0,1)$ \st $\xi^n=\frac1p$ when $n\ge2$.
Suppose for \cd ion that \te\ $a,b\in\Na$ \st $\frac ab<1$ and $(\frac ab)^n
=\frac 1p$. By Lemma \rf{l5.43} we may assume $\gcd(a,b)=1$, hence, \Ip
$p\nmid a$ or $p\nmid b$. In view of \er{5.3}, using \In\ on $n\in\Na$, one
finds $(\frac ab)^n=\frac{a^n}{b^n}$. Hence $pa^n=b^n$ by \er{5.2}. Thus
$p|b^n$. We now show that $p|b$. Suppose, for \cd ion, that $p\nmid b$, then
$\phi_p(b)\ne0$. By \E\Pr\ \rf{p3.15}\,(ii) and \era2{1.45} we have $\F_p(b^n)
=\F_p(n\ddt b)=n\ddtt p \F_p(b)$ the $n$-fold \IT\ of $\F_p(b)$ in $(\N_p,
\cdot_p,1)$. Since $\F_p(b)\in\N_p\sms0$ and $\N_p\sms0$ is a subgroup of the
monoid $(\N_p,\cdot_p,1)$ by Corollary \rf{c4.42}\,(ii), $n\ddtt p\F_p(b)$ belongs to $\N_p\sms0$ by Lemma
\rfa2{l1.21}\,(i). \E\Tf $\F_p(b^n)\ne0$, hence $p\nmid b^n$, a~\cd ion. \If
that $p|b$. Let $c\in\Na$ \sf y $b=cp$. Then $b^n=(cp)^n \nda22.28 = c^np^n$.
Since $b^n=pa^n$, we obtain $a^n=p^{n-1}c^n$ by \cnc ity. Since $n\ge2$,
$p|a^n$. As above we infer $p|a$ which implies $\gcd(a,b)\ge p>1$, a~\cd ion.
\E\Tf there is \ti{no} $\xi\in(0,1)$ \st $\xi^n=\frac1p$. From Lemma \rf{l5.54}
we find that the set $\{x\in[0,1]:x^n<\frac1p\}$ (resp.\ $\{y\in[0,1]: \frac1p
<y^n\}$) has \ti{no} supremum (resp.\ \ti{no} infimum) when $n\ge2$.
\exs

\advance\abovedisplayskip by2pt
\advance\belowdisplayskip by2pt

We next consider an important example of a bounded in\cre\ \sq\ of \el s of
$\Q_{\ge0}$ having a supremum.

\bpr5.56
Let $x\in\Q_{\ge0}$ be \st $0<x< 1$, and let
\beq5.142
s_n:=\sum_{k=0}^n x^k, \q n\in\Na.
\e
Then
\bga5.143
1=s_0<s_n<s_{n+1}<\frac1{1-x} \q \hbox{\fa}n\in\Na.\\
\frac1{1-x} = \sup\{s_n\in\Q_{\ge0}: n\in\N\}. \lb{5.144}
\e
\epr

\proof \

``\er{5.143}'': $s_0=\suml_{k=0}^0 x^k \nda27.10 = x^0=1$. Let $n\in\N$.
$s_{n+1}=\suml_{k=0}^{n+1} x^k \nda27.13 = \suml_{k=0}^n x^k+x^{n+1} >
\suml_{k=0}^n x^k = s_n$, since $x^{n+1}\in\Q_{>0}$. Here $\suml_{k=0}^{n+1}$ is the
\cme sum in the monoid $(\Q_{\ge0},+,0)$.

``\er{5.144}'': Set $A:=\{s_n\in\Q_{\ge0}: n\in\N\}$.

``\ti{$s_n<\frac1{1-x}$ \fa $n\in\Na$}'': Note that $1-x\in\Q_{>0}$, hence
$\frac1{1-x}=(1-x)\mo$ is well-defined. Let $n\in\N$.
Then $xs_n = x\suml_{k=0}^n x^k = \suml_{k=0}^n xx^k$ by \er{4.46} in the field
$(\Q,+,\cdot,0,1)$. Note that $x^{k+1}=(k+1)\ddt x$ in the monoid $(\Q_{>0},
\cdot,1)$, hence $(k+1)\ddt x=(1\ddt x)\cdot(k\ddt x)= x\cdot x^k$ by
\era2{2.3}\,I2. Thus $x\suml_{k=0}^n x^k = \suml_{k=0}^n x^{k+1}= \suml_{k=1}^
{n+1}x^k$ by \er{4.31}. Hence $1+\suml_{k=1}^{n+1}x^k = \suml_{k=0}^{n+1}x^k=
s_n+x^{n+1}$, where $x^{n+1}>0$. \E\Tf
\beq5.145
1+xs_n = s_n+x^{n+1}.
\e
Since $x^{n+1}>0$, we have $1+xs_n>s_n$, and since $x<1$ and $s_n>0$, we have
$xs_n<1s_n=s_n$ by \er{5.121}. \E\Tf $1+xs_n > (s_n-xs_n)+xs_n$, hence $1>s_n
-xs_n$ by \er{5.120}. Note that $s_n-xs_n = 1s_n+(-xs_n)\nde4.9 = (1+(-x))s_n
=(1-x)s_n$. Thus $(1-x)s_n<1$, hence $s_n=(1-x)\mo(1-x)s_n < (1-x)\mo$ by
\er{5.121}. This completes the proof of \er{5.144}.

``\ti{$\frac1{1-x}$ is the least \ub\ of $A$}'': We have to show that there is
no $y\in\Q_{\ge0}$ \st $s_n\le y<\frac1{1-x}$ \fa $n\in\N$. Suppose, for \cd
ion, that such a $y$ exists, that is, $(1-x)s_n\le (1-x)y<1$ by \er{5.121}.
From \er{5.145} we have $1-x^{n+1} = s_n-xs_n = (1-x)s_n\le (1-x)y<1$ \fa
$n\ge0$. We claim that \te s $\ov n\in\N$ \st $1-\ov x{}^{n+1}>(1-x)y$. This
would imply the \cd ion $(1-x)y<1-\ov x{}^{n+1}\le(1-x)y$. Note that $z:=1-
(1-x)y>0$. We show $\ov x{}^{n+1}<z$ \fs $\ov n\in\N$. Since $x<1$, $\frac1x
>1$, hence \te s $h>0$ \st $\frac1x=1+h$. Thus $(\frac1x)^{n+1}=(1+h)^{n+1}$,
$n\in\N$.

From \er{4.46} with $F:=\Q$, $\a:=1+h$, $\b:=1$, we obtain $(1+h)^{n+1}-1 =\break
((1+h)-1)\cdot \suml_{k=0}^n (1+h)^{n-k}1^k$. Since $1+h>1$, $(1+h)^{n-k}\ge1$,
$n\in[0,k]$. Hence $\suml_{k=0}^n (1+h)^{n-k}1^k \ge \suml_{k=0}^n 1\nda21.126
= n+1$. \E\Tf $(1+h)^{n+1}-1 \ge h(n+1)$, $n\in\N$, and $x^{n+1}=(\frac1{1+h}
)^{n+1}= \frac1{(1+h)^{n+1}}\le \frac1{1+h(n+1)}$, $n\in\N$. It remains to
find $\ov n\in\N$ \st $\frac1{1+h(\ov n+1)}<z$. It suffices to find $\ov n
\in\N$ \st $\frac1{h\ov n}<z$, since $\frac1{1+h(\ov n+1)} < \frac1
{h(\ov n+1)}<\frac1{h\ov n}$, $\ov n\in\Na$. Such an $\ov n\in\Na$ exists by
Lemma \rf{l5.49} with $y:=\frac1h$ and $x:=z$.
\endproof

We conclude this section by showing how \ra\ \nm s can be ``approximated'' by
decimal \nm s.

We first introduce the notions of integral and fractional parts of a \po\ \ra\
\nm. We use the notation
\beq5.136
\zo0,1 :=\{\a\in\Q_{\ge0}:\a<1\}.
\e

Observe that $\zo0,1 $ as a subset of $\N$ is equal to~$\{0\}$, thus from now
on we use the notation $\zo0,1 $ only for its \df\ in \er{5.136}.

\bpr5.52
\E\fe $x\in\Q_{\ge0}$ \te s \ooo pair $(n,\rho)\in\N\t\zo0,1 $ \sf ying
\beq5.129
x=n+\rho.
\e
Then $n$ is called the \emph{integral part}\index{integral part} of~$x$, denoted by $\Int(x)$, and
$\rho$ is called the \emph{fractional part}\index{fractional part} of~$x$, denoted by $\Fr(x)$.
\epr

\proof \

``$x=0$'': If $x=0$, then $x=0+0$ with $0\in \N\cap\zo0,1 $. If $0=n+\rho$,
$n\in\N$, $\rho\in\zo0,1 $, then $n=0$ and $\rho=0$ by \era2{1.9} since
$(\Q_{\ge0},+,0)$ is a \PM\ in view of \E\Pr\ \rf{p5.4}.

``$x>0$'': Let $k,l,k',l'\in\Na$ be \st $x=\frac kl=\frac{k'}{l'}$. By \er{5.2}
we have
\beq5.130
kl'=lk'.
\e
By Theorem \rfa2{t1.38} with $(\wt E,\wt e,\wt S):=(\N,0,S)$ and \era2{2.8},
\te s a unique pair $(q,r)\in\N\t\N_l$ and a unique pair $(q',r')\in\N\t\N
_{l'}$ \sf ying
\bea5.131
k&=ql+r,\\
k'&=q'l'+r'.\lb{5.132}
\e
After \ml ying both sides of \er{5.131} by $l'$ and both sides of \er{5.132}
by~$l$, and using the \cmt ity and the \asc ity of the \mlc\ we obtain
\beq5.133
qll'+rl' = kl'\nde5.130 = lk'=q'll'+r'l.
\e
We claim $q=q'$. Suppose first $q'\le q$. Let $p\in\N$ be \st $q=q'+p$. Then
from \er{5.133}, \er{4.3} we have $q'll'+pll'+rl' = q'll'+r'l$, hence by \cnc
ity $pll'+rl'=r'l\nda22.19 < ll'$ since $r'<l'$ and $l\in\Na$. Thus $pll'<pll'
+rl'<ll'$, hence $pll'<ll'$ by \era1{3.11}. Since $ll'\in\Na$, $ll'\nda22.18 =
1ll'$, we obtain $p<1$ by \era2{2.20}. Hence $p=0$, and $q=q'$. The case
$q\le q'$ follows from interchanging $k,l$ and $k',l'$, and the claim is
proved. Since $q=q'$ in \er{5.133}, we obtain $rl'=r'l$. Summarizing, we obtain
\beq5.134
q=q' \qh{and } rl'=r'l.
\e

If $r=0$ then $r'=0$ by \er{5.134}, \era2{2.13}, and $l|k$, $l'|k'$. Set $n:=
\frac kl=\frac{k'}{l'}$ and $\rho:=0$. Then $x=n+\rho$ where $n\in\N$ and
$\rho\in\zo 0,1 $. If $q=q'=0$, then $k=r>0$, $k'=r'>0$ by \er{5.131},
\er{5.132}. Set $n:=0$ and $\rho:=\frac rl\nad{\er{5.134},\er{5.2}}=\frac{r'}
{l'}$. Note that $0<\frac rl<1$ by \er{5.15}, since $r<l$. Set $n:=0$, $\rho:=
\frac rl$. Then $x=r+\rho$ where $n\in\N$ and $\rho\in\zo0,1 $, since $\frac rl
=\frac kl$. Finally, if $q=q'\ne0$ and $r,r'\in\Na$, then set
\beq5.143a
n:=q=q', \q \rho:=\frac rl= \frac{r'}{l'}\,.
\e
Then, as above, we have $n\in\N$ and $\rho\in\zo0,1 $. \Mo $n+\rho = q+\frac rl
= \frac{ql}l + \frac rl \nde5.25 = \frac{ql+r}l \nde5.131 = \frac kl=x$.

We now suppose $x=n+\rho=n'+\rho'$ where $x>0$, $n,n'\in\N$ and $\rho,\rho'
\in\zo0,1 $. We have to show $n=n'$ and $\rho=\rho'$. It suffices to show that
the case $n\ne n'$ is impossible. Indeed, if $n=n'$, then by \cnc ity $\rho=
\rho'$. If $n>n'$ then \te s $p\in\Na$ \st $n=p+n'$. Hence by \cnc ity $p+\rho
=\rho'$. Then $\rho'<1\le p\le p+\rho=\rho'$. Thus $\rho'<\rho'$, which is
impossible. If $n'>n$ then \te s $p\in\Na$ \st $n'=p+n$. By \cnc ity $\rho=
p+\rho'$. Then $\rho<1\le p\le p+\rho'=\rho$. Thus $\rho<\rho$, which is
impossible.
\endproof

\bdf5.57
A \ra\ \nm\ $x\in\Q$ is called a \ti{decimal \nm}\index{decimal number} if either $x=0$, or $x=\frac
m{10^n}$, or $-x=\frac m{10^n}$ \fs $m\in\Na$ and some $n\in\N$. We shall
denote the set of decimal \nm s by $\Q_{10}$ (not standard notation).
\edf

Clearly $\Z\sbs \Q_{10}$ (take $n=0$). \Mo the sum and the product of decimal
\nm s is a decimal \nm.

\bex5.58
Show:

\hph i,ii, $\Q_{10}$ is a sub\sr\ of $\Q$.

\hph ii,i, $\Q_{10}$ as a \sr\ is a ring.

\hph iii,, $\Q_{10}$ is a domain.

\hph iv,, The field of \qt s of $\Q_{10}$ is ring-\is c to $\Q$.
\eex

Let $x=\frac m{10^n}$, $m\in\Na$ and $n\in\Na$ be \st $m>10^n$, and let
$m=\suml_{k=0}^N \a_k10^{k}$ be its \dr, $N\ge n$, $\a_N\ne0$, $\a_k\in
[0,9]$ for $k\in [0,N]$. Then $x=\suml_{k=0}^N \a_k10^{k-n}$ with $N\ge n$,
and $x=\suml_{k=n}^N \a_k10^{k-n}+\suml_{k=0}^{n-1} \a_k10^{k-n}$, where
\beq5.147
10^{-l} := (10^l)\mo, \q l\in\Na.
\e
Since $\a_N\ne0$, we have $\suml_{k=n}^N \a_k10^{k-n}=\suml_{l=0}^{N-n}
\a_{l+n}10^l \in\Na$. \E\oh $\suml_{k=0}^{n-1} \a_k10^{k-n}= 10^{-n}
\suml_{k=0}^{n-1} \a_k10^{k}<10^{-n}\cdot10^n=1$ by \era2{5.7}. Thus
$\suml_{k=0}^{n-1} \a_k10^{k-n}\in\zo0,1 $. \If that $\suml_{k=n}^N
\a_k10^{k-n}=\Int(x)$ and $\suml_{k=0}^{n-1} \a_k10^{k-n}=\Fr(x)$ by \E\Pr\
\rf{p5.52}.

\bdf5.59
A decimal \nm\ $x\in(0,1)$ will be called a (positive) \ti{decimal fraction}.\index{decimal fraction}
\edf

If $x$ is a decimal fraction and $x=\frac m{10^n}$, then $m\in\Na$, $n\in\Na$
with $m<10^n$. \If that $x=\suml_{l=1}^M \b_l10^{-l}$ \fs $M\in\Na$, $\b_l\in
[0,9]$ with $\b_{\bar l}\ne0$ \fs $\ov l\in[1,M]$. Clearly if $\ov l<M$, then
$x=\suml_{l=1}^{\bar l} \b_l10^{-l}$, and it is usually written as
$0{,}\b_1\b_2\ldots\b_l$ (or $0.\b_1\b_2\ldots\b_l$) or
$0{,}\b_1\b_2\ldots\b_{\bar l}$ (or $0.\b_1\b_2\ldots\b_{\bar l}$).
For example, if $x=0.381=\suml_{l=1}^3 \b_l10^{-l}$ with $\b_1=3$, $\b_2=8$
and $\b_3=1$, then $10x = \suml_{l=1}^3 \b_l 10^{1-l}$, since $10\cdot 10\mo
=10\cdot \frac1{10}=1=1^0$, $10\cdot10^{-2}= 10\cdot\frac1{100}\nde5.2 =
\frac1{10}= 10\mo$, and $10\cdot10^{-k}=10\cdot\frac1{10^k}\nde5.2 = \frac1
{10^{k-1}}=10^{-(k-1)}=10^{1-k}$, $k\in\Na$. Thus $10x = \b_1+\suml_{k=2}^3
\b_k10^{1-k}$. Note that $\b_1=3\in\Na$ and $8\cdot10\mo+1\cdot10^{-2}\in
(0,1)$. Thus $3=\b_1=\Int(10x)$, and $8\cdot10\mo+1\cdot10^{-2}=0.81=\Fr(10x)$.
Then $10\cdot0.81=8.1$, $\b_2=8=\Int(10\cdot0.81)$, $\b_3=1=\Fr(8.1)$.

\bex5.60
Let $x=\suml_{l=1}^M \b_l10^{-l}$, $\b_l\in[0,9]$ for all $l\in[1,M]$, $M\in\Na$. Set $x_0:=x$, $y_k:=\Int(10
x_{k-1})$, $x_k:=\Fr(10x_{k-1})$ for $k\in[1,M]$. Show that $y_l=\b_l$,
$l\in[1,M]$.
\eex

\bex5.61
Let $x\in\Q_{>0} \cap (0,1)$. Define \sq s $\zb yk\Na$, $\zb xk\N$ recursively
by $x_0:=x$, $y_k=\Int(10 x_{k-1})$, $x_k=\Fr(10x_{k-1})$ \fa $k\in\Na$. Show

\hph i,vii, $y_k\in[0,9]$, $k\in\Na$.

\hph ii,vi, $x=\suml_{k=1}^n y_k10^{-k}+x_n10^{-n}$, $n\in\Na$.

\noindent Set $s_n:=\suml_{k=1}^n y_k10^{-k}$, $n\in\Na$. Show

\hph iii,v, $s_n\in\Q_{10}\cap(0,1)$.

\hph iv,ii, $s_n\le s_{n+1}\le x$, $n\in\Na$.

\hph v,iii, $x-s_n\le 10^{-n}$, $n\in\Na$.

\hph vi,ii, $x=\sup\{s_n:n\in\Na\}$.

\hph vii,i, $x=\max\{s_n:n\in\Na\}$ iff $x$ is a decimal \nm.

\hph viii,, $\sup\{\suml_{k=1}^n 9\cdot 10^{-k}:n\in\Na\}=1$.
\eex

\brm5.62
It can be shown that the \sq\ $\zb yk\Na$ introduced in Exercise \rf{ex5.61}
is \ti{eventually periodic},\index{eventually periodic sequence} that is, \te\ $L,T\in\Na$ \st $y_{k+T}=y_k$ \fa
$k\ge L$.

Conversely, if a \sq\ $\zb zk\Na$, $z_k\in[0,9]$ \fa $k\in\Na$ with $z_{\bar k}
\ne0$ \fs $\ov k\in\Na$ is eventually periodic, then \te s $x\in\Q\cap(0,1)$
\st $x=\sup\{s_n\in\Q\cap (0,1): n\in\Na\}$ where $s_n:=\suml_{k=1}^n z_k
10^{-k}$.

\Mo if $x=\frac ab$, $a,b\in\Na$, $a<b$ with $\gcd(a,b)=1$, then $L=1$ iff
$\gcd(b,10)=1$.
\erm

%% file: DETOUR5.TEX
\Section{Finite fields}[Finite fields]\label{s.5}
\Subsubsection{Finite-dimensional vector spaces}\label{sss.fdvs}

Let $F$ be a finite field and let $K$ denote its \Pf. We assume $F\ne K$.
Our first goal is to show that $(F,+,0)$, the additive group of~$F$, is
monoid-\is c to the group $(K,+,0)^{[1,n]}$ \fs $n\in\Na$. As a subset of
a finite set, $K$~is finite, hence the field~$K$ is ring-\is c to the field
$(\N_p,+_p,\cdot_p,0,1)$ \fs \Pn~$p$ in view of Theorem \rfa4{t4.35}. \If from
\era2{3.38} that $\#(F)=p^n$. The goal of this chapter is to show that \fe
prime~$p$ and every $n\in\Na$ \te s a field of order~$p^n$ which is unique up
to ring-\is sms.

We first introduce the important algebraic notion of \vsf0(not necessarily
a~\Pf).\glossary{$(K,V)$}

\bdf1.1
A pair $(K,V)$ is called a \ti{\vsf0}$K$, if\index{vector space over a field}

\hph i,ii, $V$ is an \ag\ $(V,+,0)$, whose \el s are called \ti{vectors},\index{vector}

\hph ii,i, $K$ is a field $(K,+,\cdot,0,1)$, whose \el s are called \ti{scalars},\index{scalar}

\hph iii,, \te s a map, called \ti{\mlc\ by scalars}\index{multiplication by scalars} $(\la,x)\mt \la x$, from $(K,V)$ into~$V$
\sf ying \fa $\la,\mu\in K$ and all $x,y\in V$:
\beq1.1
\aligned
&\bca
{\rm I0} & 0x=0,\\ {\rm I1} & 1x=x,
\eca \\
&\bca
{\rm I2} & (\la+\mu)x = \la x+\mu x, \\ {\rm I3} & (\la\mu)x = \la(\mu x),
\eca \\
&\bca
{\rm I4} & \la0 = 0, \\ {\rm I5} & \la(x+y) = \la x+\la y.
\eca
\endaligned
\e
\edf

On the \RHS\ of I0, I4 and on the \LHS\ of~I4, $0$~denotes the \nel\ of the
group $(V,+,0)$. Note the analogy between axioms \er{1.1} and \pp ies
\era2{2.3} of the \IT-map in a monoid $(M.\qu,e)$. If $V:=\{0\}$ then
$(K,V)$ is called a \ti{trivial\/} \vs.\index{trivial vector space}

\bxs1.2 \

\hph i,i, Let $I$ be a \ns, and let $K$ be a field. Set $V:=K^I$, define
$0\in V$ by setting $0(i):=0$, $i\in I$, where $0\in K$, and given $a,b\in V$
define $a+b\in V$ by setting $(a+b)(i):=a(i)+b(i)$, $i\in I$, where the \ad\
on the \RHS\ is the \ad\ in~$K$. One verifies that $(V,+,0)$ is an \ag. \Mo
given $\la\in K$ and $a\in V$, define $\la a$ by setting $(\la a)(i)=\la
a(i)$, $i\in I$, where the \mlc\ on the \RHS\ is the \mlc\ in~$K$. Then
\er{1.1} directly follows from the axioms of a field. In case $I:=[1,n]$,
$n\in\Na$, $K^{[1,n]}$ is usually denoted by~$K^n$. If $a\in K^I$ and $i\in I$,\glossary{$K^n$}
then $a(i)$ is called the \ti{$i$-th \co}\index{i-th coordinate@$i$-th coordinate} of the vector~$a$.

\hph ii,, Let $E$ be a subfield (not necessarily a \pf) of a field~$F$.
Let $(V,+,0):=(F,+,0)$ and define the scalar \mlc\ $\la a$, $\la\in E$,
$a\in F$, by setting $\la a:=\la\cdot a$ where the \mlc\ on the \RHS\ is the
\mlc\ in $(F,+,\cdot,0,1)$. Then \er{1.1} directly follows from the axioms of
a~field.
%
\exs

\bex1.3
Give the details of the proofs in Examples (i)--(ii).
\eex

\bdf1.4
Let $(K,V_1)$ and $(K,V_2)$ be \vs s over a field~$K$. A~map $L:V_1\to V_2$
is  called a \ti{linear map}\index{linear map} (or \ti{\Op}) if it \sf ies{\rm:}
\bea1.2
&L(x+y) = Lx+Ly \q \hbox{\fa} x,y\in V_1,\\
&L(\la x) = \la(Lx)\q \hbox{\fa} \la \in K \hbox{ and all } x\in V_1. \lb{1.3}
\e
\edf

\E\Ip $L0 \nad{\rm I0}= L(00) \nde1.3 = 0L0 \nad{\rm I0}= 0$.

A \ti{bi\jc} \lop\ is called a \ti{linear \is sm}.\index{linear isomorphism}

\blm1.5 \

\hph i,ii, The inverse of a bi\jc\ linear map is a linear map.

\hph ii,i, The \cm\ of linear maps $($resp.\ \is sms$)$ is a linear map $($resp.\
\is sm$)$.

\hph iii,, The identity in a \vs\ is a linear \is sm.

\elm

\bex1.6
Prove Lemma \rf{l1.5}.
\eex

\bxs1.7 \

\hph i,i, Let $K$ be a field, let $m,n\in \Na$ and let $c:[1,m]\t[1,n]\to K$.
Define $L:K^n\to K^m$ by setting
\beq1.4
(Lx)(i) := \sum_{j=1}^n c_{ij}x_j, \q x\in K^n,\ i\in[1,m],
\e
where $\sum$ is the \cme sum in the \ag\ $(K^n,+,0)$.
Then $L$ is a \ti{\lop}. Indeed, let $\la\in K$, $x,y\in K^n$, $i\in[1,m]$.
Then $(L(\la x+y))(i) = \suml_{j=1}^n {c_{ij}(\la x+y)(j)} = \suml_{j=1}^n
c_{ij}((\la x)(j)+y(j)) = \suml_{j=1}^n c_{ij}((\la x(j)+y(j)) = \suml_{j=1}^n
(c_{ij}\la x(j)+c_{ij}y(j)) \nda21.129 = \suml_{j=1}^n c_{ij}\la x(j) +
\suml_{j=1}^n c_{ij}y(j) = \suml_{j=1}^n \la c_{ij}x(j) + \suml_{j=1}^n
c_{ij}y(j) \nda21.131 = \la\suml_{j=1}^n c_{ij}x(j) + \suml_{j=1}^n c_{ij}
y(j) = \la((Lx)(i)) + (Ly)(i) = (\la (L(x)(i)))+(Ly)(i) = (\la (Lx)+Ly)(i)$. Hence
$L(\la x+y)=\la (Lx)+Ly$. If $\la:=1$, then $L(x+y)\nad{\rm I1}= L(1x+y) =
1(Lx)+Ly \nad{\rm I1}= Lx+Ly$, which proves \er{1.2}. If $y:=0$, then
$L(\la x)=L(\la x+0)= \la (Lx)+L0$. But \fa $i\in[1,m]:(L0)(i) = \suml_{j=1}^n
c_{ij}0(j) = \suml_{j=1}^n c_{ij}0 =\suml_{j=1}^n 0\nda21.126 = \#([1,n])\dpl0
\nad{\era2{2.3}\,\rm I4} = 0$, where $\#([1,n])\dpl0$ is the $\#([1,n])$-th
\IT\ of~$0$ in $(K^n,+,0)$. Thus $L0=0$ and $L(\la x)=\la (Lx)$, which proves
\er{1.3}.

\hph ii,, Let $K$ be a field, $l\in\Na$ and
\beq1.5
\psi_k (z):=z(k), \q z\in K^l \qh{and} k\in[1,l].
\e
The maps $\psi_k$ are called the \ti{\co-maps}\index{co-ordinate maps} of~$K^l$. Then $\psi_k:K^l
\to K$ is \ti{linear} \fe $k\in[1,l]$. Indeed, if
$x,y\in K^l$, $k\in[1,l]$, then $\psi_k(x+y):=(x+y)(k):=x(k)+y(k) =
\psi_k(x)+\psi_k(y)$. \Mo if $\la\in K$, $x\in K^l$ and $k\in[1,l]$, then
$\psi_k(\la x):=(\la x)(k):=\la x(k)=\la\psi_k(x)$.
\exs

Observe that if $l\in\Na$, $x,y\in K^l$, then
\beq1.6
x=y \qh{iff } \psi_k(x) = \psi_k(y)
\e
\fa $k\in[1,l]$ in view of the \df\ of $K^l$.

Linear maps from a \vs\ $(K,V)$ into~$K$ are usually called \ti{\lf s}.\index{linear functional}

Our first goal is to show that if $m,n\in\Na$, then \te s a linear \is sm
between the \vs\ $K^m$ and $K^n$ iff $m=n$.

\bnt1.8
Let $(K,V)$ be a \vs\ over~$K$ and let $x,y\in V$. Then $\sum$ denotes the
\cme sum in $(V,+,0)$, and
\bga1.7
-x \hbox{ denotes the inverse of $x$ in the group }(V,+,0),\\
x-y := x+(-y). \lb{1.8}
\e
\ent

\blm1.9
Let $(K,V)$ be a \vs\ over~$K$, let $\la\in K$, $x\in V$. Let $J$ be a
nonempty \emph{finite} set, $\a_j\in K$, and $x_j,y_j\in V$ \fa $j\in J$. Then
\bea1.9
& -x = (-1)x, \\
& \sum_{j\in J}\a_jx = \Bg(\sum_{j\in J}\a_j)x, \q
\sum_{j\in J} \la x_j = \la \sum_{j\in J}x_j,\q \sum_{j\in J}0 = 0,\lb{1.10}\\
& \sum_{j\in J}(-x_j) = -\sum_{j\in J}x_j, \lb{1.11}\\
& \sum_{j\in J}(x_j\pm y_j) = \Bigl(\sum_{j\in J}x_j\Bigr) \pm \Bigl(\sum_{j\in J}
y_j\Bigr). \lb{1.12}
\e
\elm

\proof \

\er{1.9}: $x+(-1)x \nad{\rm I1} = (1x)+(-1)x \nad{\rm I2} = (1+(-1))x=0x
\nad{\rm I0}= 0$.

\er{1.10}: The first \et y follows from \er{1.1}\,I2, \era2{1.127} and
\In\ on $\#(J)\in\Na$. The second \et y
follows from \era2{1.131} with $I:=J$, $(X,\qu,e):=(V,+,0)$,
$\wt X:=X$, $a:=x$, $\vf:V\to V$ defined by $\vf(z):=\la z$, $\la\in K$,
$z\in V$. Observe $\vf$ \sf ies \era2{1.130} by~I5. \Mo $\suml_{j\in J}0
\nad{\rm I0}= \suml_{j\in J}0\, 0= 0\suml_{j\in J}0\nad{\rm I0}=0$.

\er{1.11}: $\suml_{j\in J}(-x_j) \nde1.9 = \suml_{j\in J}(-1)x_j \nde1.10
= (-1)\suml_{j\in J}x_j \nde1.9 = -\suml_{j\in J}x_j$.

\er{1.12}: $\suml_{j\in J}(x_j-y_j)\nde1.8 = \suml_{j\in J}(x_j+(-y_j))
\nda21.129 = \suml_{j\in J}x_j +\suml_{j\in J}(-y_j) \nde1.9 = \suml_{j\in J}
x_j+\bigl(-\suml_{j\in J}y_j\bigr) \nde1.8 = \suml_{j\in J}x_j- \suml_{j\in J}
y_j$.
\endproof

\blm1.10
Let $V_1,V_2$ be \vs s over $K$ and let $L:V_1\to V_2$  be a \lop.
Let $J$ be a \nf\ and let $x_j\in V_1$ \fa $j\in J$. Then
\bga1.13
L\Bigl(\suml_{j\in J}x_j\Bigr) = \suml_{j\in J}(Lx_j),\\
\lb{1.14}
L(-x)=-Lx, \qh{\fa}x\in V.
\e
\elm

\proof
\er{1.13} follows from \er{1.2} and \era2{1.131}.

\er{1.14}: $L(-x)\nde1.9 = L((-1)x)\nde1.3 = (-1)Lx \nde1.9 = -Lx$.
\endproof

\bpr1.11
Let $K$ be a field and let $m,n\in\Na$. If \te s an in\jc\ linear map from
$K^n$ into $K^m$, then $n\le m$.
\epr

\proof
Let $L$ be an in\jc\ linear map from $K^n$ into $K^m$. We denote by $\ve_j$,\glossary{$\ve_j$}
$j\in[1,n]$, the vectors of $K^n$ whose $i$-th \co s are equal to\glossary{$\d_{ij}$}
\begin{gather}
\d_{ij}:=\bca 1 & \hbox{if }i=j,\\0 & \hbox{if }i\ne j,\eca \q
i,j\in[1,n]:\nonumber\\
\ve_j(i):= \d_{ij}, \q i,j\in[1,n].  \lb{1.15}
\end{gather}
We claim that if $x\in K^n$ and
\beq1.16
\sum_{j=1}^n x(j)\ve_j=0
\e
(where the \et y takes place in $K^n$), then $x=0$.
Indeed, let $\psi_k:K^n\to K$, $k\in[1,n]$, be the \co\ maps defined in
\er{1.5} with $l:=n$, then \er{1.16} becomes: $\psi_i\bigl(\suml_{j=1}^n
x(j)\ve_j\bigr) = \psi_i(0)\nde1.3 = 0$, $i\in[1,n]$.

By linearity we obtain $\psi_i\bigl(\suml_{j=1}^n x(j)\ve_j\bigr) \nde1.13 =
\suml_{j=1}^n \psi_i(x(j)\ve_j) \nde1.3 = \suml_{j=1}^nx(j)\psi_i(\ve_j)
\nad{\er{1.5},\er{1.15}} = \suml_{j=1}^nx(j)\d_{ij} \nda21.127 =
\bigl(\suml_{j\in[1,n]\sms i}x(j)0\bigr) + x(i)1 =\bigl(\suml_{j\in[1,n]\sms i}0) +x(i)
\nde1.10 = 0 +x(i)= x(i)$, $i\in[1,n]$. Hence $\suml_{j=1}^nx(j)\ve_j=0$
implies $x(i)=\psi_i\bigl(\suml_{j=1}^nx(j)\ve_j\bigr)\nde1.3 = 0$, $i\in
[1,n]$. Since $x(i)=0=0(i)$, $i\in[1,n]$, we infer $x=0$, which proves the
claim.

Now set $e_j:=L\ve_j\in K^m$, $j\in[1,n]$, and let $x\in K^n$. Then
$\suml_{j=1}^n x(j)e_j = \suml_{j=1}^n x(j)L\ve_j \nde1.9 = \suml_{j=1}^n
L(x(j)\ve_j)\nde1.13 = L\bigl(\suml_{j=1}^nx(j)\ve_j\bigr)$. Suppose
$\suml_{j=1}^n x(j)e_j=0$ in~$K^m$ where $x\in K^n$. Then $L\bigl(
\suml_{j=1}^nx(j)\ve_j\bigr)=0$ in~$K^m$. Since $L(0)=0$ by \er{1.3}, and
since $L$~is \ti{in\jc}, we infer $\suml_{j=1}^nx(j)\ve_j=0$ by \E\Pr\
\rfa4{p3.10}\,(vi), hence by what
precedes $x=0$.

Finally, set $a_{ij}:=\psi_i(e_j)$, $j\in[1,n]$, where $\psi_i$ is defined in
\er{1.5} with $l:=m$, $j\in[1,m]$. We obtain \fe $x\in K^n$ \sf ying
$\suml_{j=1}^n x(j)e_j=0$, \ev tly $\psi_i\bigl(\suml_{j=1}^nx(j)e_j\bigr)
=0$ \fa $i\in[1,m]$, that $x=0$. But $\psi_i\bigl(\suml_{j=1}^nx(j)e_j\bigr)
\nde1.13 = \suml_{j=1}^n\psi_i(x(j)e_j)\nde1.3 = \suml_{j=1}^n x(j)\psi_i
(e_j) = \suml_{j=1}^nx(j)a_{ij} = \suml_{j=1}^n a_{ij}x(j)$ \fa $i\in[1,m]$.
\If that the system $\suml_{j=1}^n
a_{ij}x(j)=0$, $i\in[1,m]$, possesses only the trivial \so\ $x=0$, hence by
Lemma \rfa4{l4.12} we have $m\not<n$, i.e.\ $n\le m$.
\endproof

\Wanp \es\ the first basic result of this section.

\bth1.12
Let $K$ be a field and let $m,n\in\Na$. Then the \vs s $K^m$ and~$K^n$ are
\emph{\is c} iff $m=n$.
\eth

\proof \

\ti{If}: If $m=n$, then the identity in $K^m$ is an \is sm by Lemma \rf{l1.5}\,(iii).

\ti{Only if}: Let $L$ be an \is sm from $K^n$ onto $K^m$. Since $L:K^n\to K^m$
is in\jc\ and linear, $n\le m$ by \E\Pr\ \rf{p1.11}. Since $L\Inv:K^m\to K^n$
is an \is sm by Lemma \rf{l1.5}\,(i), $L\Inv:K^m\to K^n$ is in\jc\ and linear,
hence $m\le n$ by the same \Pr. Thus $m=n$.
\endproof

We now investigate when a \vs\ $(K,V)$ is \is c to~$K^n$ \fs $n\in\Na$. Note
that in view of Theorem \rf{t1.12} and Lemma \rf{l1.5}\,(ii) a \vs\ $(K,V)$
which is \is c to~$K^n$ \fs $n\in\Na$ is \ti{not\/} \is c to~$K^m$ if $m\ne n$.

\bdf1.13
A \vs\ $(K,V)$ over a field~$K$ is called \ti{$n$-\dm al\/}\index{n-dimensional vector space@$n$-dimensional vector space} or is said to have\index{dimension of a vector space}
\ti{\dm}~$n$ if \te\ $n\in\Na$ and an \is sm $L:K^n\to V$. A~nontrivial \vs\
$(K,V)$ is called \ti{finite-\dm al\/} if \te s $n\in\Na$ \st $(K,V)$ is
$n$-\dm al, otherwise $(K,V)$ is called \ti{infinite-\dm al}.
\edf

\bxa1.14
$K^n$ is a $n$-\dm al \vs\ over the field~$K$ in view of Lemma \rf{l1.5}\,(iii).
\exa

Observe that every vector $x\in K^n$ \sf ies
\beq1.17
x=\sum_{j=1}^n \psi_j(x)\ve_j
\e
where $\ve_j$ is defined in \er{1.15} and $\psi_j$ in \er{1.5} with $l:=n$ for
$j\in[1,n]$. Indeed, set $y:=\suml_{j=1}^n \psi_j(x)\ve_j$. Then \fa $i\in
[1,n]$, $\psi_i(y)=\psi_i\bigl(\suml_{j=1}^n\psi_j(x)\ve_j\bigr) \nad{\er{1.13},\er{1.3}}=
\suml_{j=1}^n
\psi_j(x)\psi_i(\ve_j) \nde1.15 = \suml_{j=1}^n \psi_j(x)\d_{ij} \nde1.17 = \psi_i(x)$.
Thus $y=x$, by \er{1.6}.

\bdf1.15
Let $(K,V)$ be a nontrivial \vs\ over~$K$.

\hph i,ii,  Let $A$ be a nonempty \ti{finite} subset of $V\sms 0$, let $I$~be
a finite set with $\#(I)=\#(A)$, let $\vf:I\to A$ be a bi\jn\ and let $\a:I\to
K$. The vector $\suml_{i\in I}\a_i\vf_i\in V$ is called a \ti{\lc} of the
vectors $\zb \vf iI$.

\hph ii,i, The nontrivial \vs\ $V$ is called \ti{\fg} if \te s a finite subset~$A$ of
$V\sms0$ \st every vector of~$V$ is a \lc\ of the \el s of~$A$.

\hph iii,, Let $B$ be a \nss\ of $V$. The set of all \lc s of \el s of~$B$ is
called the \ti{linear span}\index{linear span} (or simply the \ti{span}) of~$B$ and is denoted by\glossary{$\spn(A)$}
$\spn(B)$.
\edf

\bxa1.16 \

\hph i,i, The \vs\ $(K,K^n)$, $n\in\Na$, is \fg. Indeed, if $n=1$, set $I:=
\{1\}$, $\wh A:=\{\ve_1\}$, then $\#(I)=\#(\wh A)=1$ and $\{\la_1\ve_1: \la_1
\in K\}=K^1$, since \fe $x\in K^1$, $x=\psi_1(x)\ve_1$ by \er{1.17}. Hence
$K^1=\spn(\{\ve_1\})$. Now let $n\in\Na\sms1$. We claim that $K^n=\spn(\wh A)$
where $\wh A:=\{\ve_j\in K^n: j\in[1,n]\}$. Clearly $\wh A\sbs V\sms0$ since
$\ve_j(j)=1$, $j\in[1,n]$. Set $I:=[1,n]$. The map $j\mt \ve_j$ from~$I$
into~$\wh A$ is sur\jc\ by \df\ of~$\wh A$ and in\jc\ since for $j,k\in I$,
$j\ne k$, we have $\e_j(j)=1\ne0=\ve_k(j)$. \E\Tf $\#(I)=\#(\wh A)$ by
Corollary \rfa2{c3.11}, and $K^n=\spn(\wh A)$ by \er{1.17}.

\hph ii,, Every $n$-\dm al \vs\ $(K,V)$, $n\in\Na$, is \fg. Set $I:=[1,n]$,
$A:=\{L\ve_j: j\in[1,n]\}$ where $L:K^n\to V$ is an \is sm. The map $j\mt L\ve_j$
from~$I$ into~$A$ is sur\jc\ by \df\ of~$A$. \Mo this map is the \cm\ of the
in\jc\ map $j\mt \ve_j$ from~$I$ into~$\wh A:=\{\ve_j\in K^n: j\in[1,n]\}$
and $\ve_j\mt L\ve_j$ which is in\jc.
\E\Tf $\#(I)=\#(A)$. Let $x\in V$ and let $\wh x:=L\Inv x$. Then $x=L\wh x\nde
1.17 = L\bigl(\suml_{j\in I}\psi_j(\wh x)\ve_j\bigr) \nad{\er{1.13},\er{1.3}}=
\suml_{j\in I}\psi_j(L\Inv x)L\ve_j$. Hence $V=\spn(A)$.
\exa

Our next goal is to show that a \fg\ \vs\ over a field $K$ is finite-\dm al. To this end
we first introduce the notion of \ti{basis}\index{basis of a vector space} of a \fg\ \vs. Let $(K,V)$ be a
nontrivial \vs\ over a field~$K$ and let $A$ be a \nfs\ of $V\sms0$. We suppose that
$V=\spn(A)$. We claim that \fe $x\in V$ \te s a map $\a\en x:A\to K$ \st
\beq1.18
x= \sum_{a\in A}\a\en x(a)a.
\e
Indeed, by \E\df\ \rf{d1.15}\,(i), \te\ a \nfs\ $A$ of $V\sms0$, a finite set~$I$
\ep\ to~$A$ and a bi\jn\ $\wt\vf:I\to A$ \st \fe $x\in V$, \te s a~map
$\wt\a{}\en x:I\to K$ \sf ying
\beq1.19
x = \sum_{i\in I}\wt\a{}\en x(i)\wt\vf(i).
\e
Set $l:=\#(I)= \#(A)$ (by Corollary \rfa2{c3.11}). Let $\wh\vf$ be a bi\jn\
from $[0,l-1]$ onto~$I$. Then by \era2{7.17} we obtain
\beq1.20
\sum_{i\in I} \wt\a{}\en x(i)\wt\vf(i) = \sum_{j=0}^{l-1} \wt\a{}\en x
(\wh\vf(j))\wt\vf(\wh\vf(j)).
\e
Set $\vf:=\wt\vf\circ\wh\vf : [0,l-1]\to A$. Observe that $\vf$ is a bi\jn\
as \cm\ of two bi\jn s, and set $\a\en x:=\wt\a{}\en x\circ\wh\vf\circ
\vf\Inv :A\to K$. Then $\suml_{j=0}^{l-1}\wt\a{}\en x(\wh\vf(j))\wt\vf
(\wh\vf(j)) = \suml_{j=0}^{l-1}\a\en x(\vf(j))\vf(j) \nda27.17 =
\suml_{a\in A} \a\en x(a)a$. Thus by \er{1.19}, \er{1.20}, we arrive at
\er{1.18}, which proves the claim.

\bdf1.17
Let $(K,V)$ be a \fg\ \vsf0$K$. A~non-empty finite subset $A$ of $V\sms0$ is called a \ti{basis}
of the \vs\ $V$ if \fe $x\in V$, \te s one and \ti{only one} map $\a\en x:A
\to K$ \st \er{1.18} holds.
\edf

\brm1.18
The \df\ of a  basis given here is not the usual one. On the one hand, it is
not usual to ``index'' the \el s of the basis~$A$ by $a\in A$, and \oh the \uq\
\cn\ is usually replaced by the \lic\ of the \el s of the basis~$A$. We shall
introduce this notion in a moment.
\erm

In what follows we shall use as ``index set'' the set $[1,\#(A)]$ and denote
by~$\vf$ a bi\jn\ from $[1,\#(A)]$ onto~$A$. Thus by \era2{7.17} and
\era2{7.12}, $\suml_{a\in A}\a\en x(a)a = \suml_{i=1}^{\#(A)}\a\en x(\vf(i))
\vf(i)$. Finally, we define $\la\en x:=\a\en x\circ\vf:[1,\#(A)]\to K$. Observe
that if $\a\en x,\b\en x:A\to K$ are distinct, so are $\a\en x\circ\vf$,
$\b\en x\circ\vf$ and conversely.

These considerations justify that the \fw\ \df\ of a basis is \ev t to \E\df\
\rf{d1.17}.

\def\namespec{Equivalent definition}
\begin{dspc} \lb{d1.19}
Let $(K,V)$ be a \fg\ \vsf0$K$. Let $A$ be a \nfs\ of $V\sms0$ and let $\vf:
[1,\#(A)]\to A$ be a bi\jn. The set $A$ is called a \ti{basis} of the \vs~$V$
if \fe $x\in V$, \te s \ti{\ooo}map $\la\en x:[1,\#(A)]\to K$ that the \fw\
holds:
\beq1.21
x=\sum_{i=1}^{\#(A)} \la\en x(i)\vf(i).
\e
\end{dspc}

\bxa1.20
Let $K$ be a field, let $I$ be a nonempty \ti{finite} set and let $K^I$ be the
\vs\ over~$K$ introduced in Example \rf{xa1.2}\,(i). \E\fe $l\in I$ we denote
by~$e_l$ the \el\ of~$K^I$ defined by
\beq1.22
e_l(l'):=\bca
1 &\hbox{if }l'=l,\\
0 &\hbox{otherwise.}
\eca
\e
We claim that $A:=\{e_l:l\in I\}$ is a \ti{basis} of~$K^I$. We first observe
that the map $\wt\vf:I\to A$ defined by $\wt\vf(l):=e_l$ is \ti{bi\jc}. It is
clearly sur\jc. Let $l,k\in I$ be \st $l\ne k$. Then $e_l(l)\nde1.22 = 1\ne 0
\nde1.22 = e_k(l)$. Thus $\wt\vf$ is in\jc, hence bi\jc. Since $I$~is a \nf,
\te s a bi\jn\ $\wh\vf:[1,\#(I)]\to I$. Indeed, by Theorem \rfa2{t3.6} with
$(E,\le,e,S) := (\N,\le,0,S)$ and $A:=I$, $\#(I)\in\Na$, \te s a~bi\jn\
$\vf':[0,\#(I))\to I$. Hence the map $\wh\vf:[1,\#(I)]\to I$ defined by
$\wh\vf(n) := \vf'(n-1)$, $n\in[1,\#(I)]$, is a bi\jn\ from $[1,\#(I)]$
onto~$I$. Define $\vf:[1,\#(I)]\to A$ by setting $\vf:=\wt\vf\circ\wh\vf$.
Then $\vf$~is bi\jc\ as the \cm\ of two bi\jn s and $\#(A)=\#(I)$ by
Corollary \rfa2{c3.11}.

Let $\psi_l:K^I \to K$ denote the $l$-th \co\ map defined by $\psi_l(z):=
z(l)$, $z\in K^I$, $l\in I$. Proceeding as in Example \rf{xa1.14} we find
\beq1.23
x=\sum_{l\in I} \psi_l(x)e_l \qh{\fa} x\in K^I.
\e
By \era2{7.17}, \era2{7.12} we obtain $x=\suml_{i=1}^{\#(I)}\psi_{\wh\vf(i)}
(x)e_{\wh\vf(i)}$, $x\in K^I$. Since $e_{\wh\vf(i)}=\wt\vf(\wh\vf(i))
=\vf(i)$, $i\in[1,\#(I)]$, we infer $x=\suml_{i=1}^{\#(I)}\psi_{\wh\vf(i)}(x)
\vf(i)$, $x\in K^I$. Setting $\la\en x(i):=\psi_{\wh\vf(i)}(x)$, $i\in
[1,\#(I)]$, $x\in K^I$, we find
\beq1.24
x=\sum_{i=1}^{\#(I)}\la\en x(i)\vf(i) \qh{\fa} x\in K^I,
\e
where $\la\en x$ is a map from $[1,\#(I)]$ into $K$, $x\in K^I$.

It remains to prove the ``\uq'' \cn. We have
\[
\xymatrix{[1,\#(I)] \ar[r]^{\wh\vf} \ar@/_/[rr]_\vf
& I \ar[r]^{\td\vf} &A}
\]
where $\wh\vf$, $\wt\vf$ and $\vf$ are bi\jn s. Let $x\in X$ and let
$\mu:[1,\#(I)] \to K$ be \st $x=\suml_{i=1}^{\#(I)} \mu(i)\vf(i)$. We want to
show that $\mu=\la\en x$. Let $i\in[1,\#(I)]$ and let $l\in I$. Then we have
\beq1.25
\psi_l(\vf(i))=\bca
1 & \hbox{if }\wh\vf(i)=l,\\
0 & \hbox{otherwise.}
\eca
\e
Indeed, $\psi_l(\vf(i))= \psi_l(\wt\vf(\wh\vf(i)))=\psi_l(e_{\wh\vf(i)})$,
hence \er{1.25} follows from \er{1.22}. Using the linearity of $\psi_l$'s, we
obtain \fe $l\in[1,\#(I)]$:
\bgg
x(l)=\psi_l(x)\nde1.24 = \psi_l\Bigl(\sum_{i=1}^{\#(I)}\la\en x(i)\vf(i)\Bigr)
\nad{\er{1.13},\er{1.3}}= \sum_{i=1}^{\#(I)}\la\en x(i)\psi_l(\vf(i)) \hskip40pt\\
\nde1.25 = \Bigl(\sum_{i\in[1,\#(I)]\sms{\wh\vf{}\Inv(l)}} \la\en x(i)0\Bigr)
+\la\en x(\wh\vf{}\Inv(l))1 =0+\la\en x(\wh\vf{}\Inv(l))
= \la\en x(\wh\vf{}\Inv(l))\\
\hbox{and similarly } x(l)=\psi_l(x)=\psi_l\Bigl(\sum_{i=1}^{\#(I)}\mu(i)\vf(i)\Bigr)
= \mu(\wh\vf{}\Inv(l)).
\e
Hence $\la\en x(\wh\vf{}\Inv(l))=\mu(\wh\vf{}\Inv(l))$ \fe $l\in I$.
Since $\wh\vf:[1,\#(I)]\to I$ is bi\jc, we conclude that $\la\en x=\mu$. \If
that $A:=\{e_j:j\in I\}$ is a basis and that the \nme\ of~$A$ is equal to~$\#(I)$.
\exa

\bex1.21
Let $K^I$ be the \vs\ and let $A$ be the basis introduced in Example \rf{xa1.20}.
Suppose $\#(I)\ge2$. Let $l\in I$ and set $b_{l'}:=e_{l'}$ \fa $l'\in
I\sms l$ and set $b_l:=\la e_l+e_k$, $\la\in K\sms0$, where $k\in I\sms l$.
Show that $\{b_{l'}\in K^I: l'\in I\}$ is a basis of~$K^I$.
\eex

As promised we now introduce the notion of linear (in)dependence for sets of
vectors in a nontrivial \vs\ $(K,V)$.

Let $(K,V)$ be a nontrivial \vs\ over $K$, let $A$ be a \nfs\ of $V\sms0$,
and let $\vf:[1,\#(A)]\to A$ be a bi\jn. We set $a_i:=\vf(i)$,
$i\in[1,\#(A)]$, and introduce a map $\cR: K^{\#(A)}\to V$
(also denoted by $\cR_A$) defined by
\beq1.27
\cR(\la):=\sum_{i=1}^{\#(A)}\la_ia_i, \q \la\in K^{\#(A)}.
\e
We denote by $N(\cR)$ the \ti{null-space} of~$\cR$, i.e.
\beq1.29
N(\cR):=\{\la\in K^{\#(A)}: \cR(\la)=0\}.
\e

\blm1.23
Let $\cR$ be as in \er{1.27}. Then

\hph i,ii, $\cR$ is linear.

\hph ii,i, $\cR$ is in\jc\ iff $N(\cR)=\{0\}$.

\hph iii,, The range of $\cR$ is equal to $\spn(A)$ $($see \E\df\
\rf{d1.15}\,{\rm (iii))}.

\hph iv,, $\cR$ is bi\jc\ iff $A$ is a basis of $V$.

\hph v,i, If $A$ and $B$ are bases of $V$, then
\beq1.28
\#(A)=\#(B).
\e
\elm

\proof\

(i) Let $\la,\mu\in K^d$. Then $\cR(\la+\mu)= \suml_{i=1}^d
(\la+\mu)_ia_i = \suml_{i=1}^d(\la_i+\mu_i)a_i \nad{\er{1.1}\,\rm I2} = \suml_{i=1}^d(\la_ia_i +
\mu_ia_i)\nde1.12 = \suml_{i=1}^d\la_ia_i + \suml_{i=1}^d\mu_ia_i =
\cR(\la)+\cR(\mu)$.

Let $\a\in K$ and $\la\in K^d$. Then $\cR(\a\la)= \suml_{i=1}^d(\a\la)_i
a_i=\suml_{i=1}^d\a\la_ia_i \nde1.10 = \a\suml_{i=1}^d\la_ia_i =\a\cR(\la)$.

(ii) Since $\cR$ is linear, $\cR$ is a monoid-\hm sm from $(K^{\#(A)},+,0)$
into $(V,+,0)$ by \er{1.2}. Then (ii) follows from \E\Pr\ \rfa4{p3.10}\,(vi)
with $\ker\cR=N(\cR)$.

(iii) directly follows from the \df\ of $\spn(A)$.

(iv) By \E\df\ \rf{d1.19}, $A$ is a basis of~$V$ iff $\spn(A)=V$ and $\cR$~is
in\jc\ iff $\cR$ is bi\jc.

(v) If $A$ is a basis of $V$, then $\cR_A:K^{\#(A)}\to V$ is a linear \is sm
by (i), (iv) and \E\df\ \rf{d1.4}. If $B$~is a basis of~$V$, then $\cR_B:
K^{\#(B)}\to V$ is a linear \is sm. Then $\cR_B\Inv \circ \cR_A : K^{\#(A)}
\to K^{\#(B)}$ is a linear \is sm by Lemma \rf{l1.5}. Hence $\#(A)=\#(B)$
 by Theorem \rf{t1.12}.
\endproof

\blm1.25
Let $\cR$ be as in \er{1.27}. Then $N(\cR)\ne\{0\}$ iff \te s $x\in A$ \st $x\in
\spn(A\sms x)$.
\elm

\proof \

\ti{Only if\/}: \E\te s $\la\in K^{\#(A)}$ and $j\in [1,\#(A)]$ \st
$\la_j\ne0$ and $\suml_{i=1}^{\#(A)}\la_ia_i=0$. Hence $\la_ja_j=
\suml_{i\in\{j\}}\la_ia_i=-\suml_{i\in[1,\#(A)]\sms j}\la_ia_i$ by
\era2{1.138}, thus
\[
a_j = \la_j\mo\Bigl(
-\sum_{i\in[1,\#(A)]\sms j}\la_ia_i\Bigr) \nad{\er{1.10},\er{1.11}}=
\sum_{i\in[1,\#(A)]\sms j}(-\la_j\mo \la_i)a_i,
\]
where $\la_j\mo$ denotes
the inverse of $\la_j$ in $(K\sms0,\cdot,1)$. Then $a_j\in\spn(A\sms{a_j})$.

\ti{If\/}: \E\te\ $j\in[1,\#(A)]$ and $\mu:[1,\#(A)]\sms j \to K$ \st $a_j
=\!\!\suml_{i\in[1,\#(A)]\sms j}\!\mu_ia_i$. Set $\mu_j:=-1$. Then
$\suml_{i\in[1,\#(A)]}\mu_ia_i=0$, where $\mu_j\ne0$. Hence $\mu:[1,
\#(A)]\to K$ \sf ies $\mu\ne0$ and $\cR(\mu)=0$.
\endproof

\bdf1.26
Let $(K,V)$ be a \fg\ \vs\ over~$K$. A~\nfs\ $A$ of $V\sms0$ is called
\ti{linearly dependent\/}\index{linearly dependent subset} if \te s $x\in A$ \st $x\in \spn(A\sms x)$, and
\ti{\li\/}\index{linearly independent subset} otherwise.
\edf

\bpr1.27
Let $(K,V)$ be a \fg\ \vs\ over~$K$. Let $A$~be a \nfs\ of $V\sms0$. Then
$A$~is a basis iff

\hph i,i, $V=\spn(A)$,

\hph ii,, $A$ is \li.
\epr

\bex1.28
Prove \E\Pr\ \rf{p1.27}.
\eex

\bth1.29
Let $(K,V)$  be a nontrivial \vs\ over~$K$. Suppose \te s a \nfs\ $A$ \st
$V=\spn(A)$. Then \te s a \nss\ $B$ of~$A$ \st $B$ is a basis of~$V$ and
$V$~is $\#(B)$-\dm al with $\#(B)\le\#(A)$.
\eth

\proof
We proceed by \In\ on $\#(A)\in\Na$. Set $M:=\{n\in\Na: $ Theorem \rf{t1.29}
holds with $\#(A)=n\}$.

$1\in M$: If $\#(A)=1$, then $A=\{a\}$ \fs $a\in V\sms0$ and $V=\{\la a\in V:
\la\in K\}$. We show that $\cR(\la)=0$, $\la\in K$, implies $\la=0$. Let $\la
\in K$ be \st $\cR(\la)\nde1.27 = \la a=0$. Suppose for \cd ion that $\la\ne0$.
Then $a\nad{\er{1.1}\,\rm I1}= 1a = (\la\mo\la)a \nad{\er{1.1}\,\rm I3}=
\la\mo(\la a) = \la\mo 0\nad{\er{1.1}\,\rm I4}=0$, a~\cd ion. Since $\spn(\{
a\})=V$ and $N(\cR)=\{0\}$, $\{a\}$~is a basis of~$V$, by Lemma \rf{l1.23}\,(ii).

\ti{$n\in M$ implies $n+1\in M$}: Let $n\in M$ and let $A$ be \st
$\#(A)=n+1$, and $V=\spn(A)$. Then either $A$ is a basis or in view of \E\Pr\ \rf{p1.27}
$A$~is linearly dependent. Let $\vf:[1,n+1]\to A$ be a bi\jn\ and set
$a_i:=\vf(i)$, ${i\in[1,n+1]}$. In view of \E\df\ \rf{d1.26} \te s
$j\in[1,n+1]$ \st ${a_j\in\spn(A\sms{a_j})}$. Set $\wh A:=A\sms{a_j}$, then
$\#(\wh A)=n$. We claim that $\spn(\wh A)=V$. Clearly $\spn(\wh A)\sbs V$. We
now show that every \el\ of~$V$ is a \lc\ of \el s of~$\wh A$. Since
$a_j\in\spn(\wh A)$, \te s $\wh\la:[1,\#(A)]\sms j\to K$ \st $a_j=
\suml_{i\in[1,\#(A)]\sms j}\wh\la_ia_i$. Let $x\in V$. Since $V=\spn(A)$, \te
s $\mu:[1,\#(A)]\to K$ \st $x=\suml_{i\in[1,\#(A)]}\mu_ia_i$. But
$\suml_{i\in[1,\#(A)]}\mu_ia_i= \bigl(\suml_{i\in[1,\#(A)]\sms j}\mu_ia_i\bigr)
+\mu_ja_j = \bigl(\suml_{i\in[1,\#(A)]\sms j}\mu_ia_i\bigr) +
\mu_j\suml_{i\in[1,\#(A)]\sms j}\wh\la_ia_i\nad{\er{1.10},\er{1.12},\er{1.1}
\,\rm I2} = \break \suml_{i\in[1,\#(A)]\sms j}
(\mu_i+\mu_j\wh\la_i)a_i$. \E\Tf $x\in\spn(\wh A)$, and $V=\spn(\wh A)$.
Since $\#(\wh A)=n$ and since $n\in M$, we infer that \te s~$B$, a~subset
of~$\wh A$, \st $B$ is a basis of~$V$. \If from Lemma \rf{l1.23}\,(i) and (iv)
that $V$~is \is c to $K^{\#(B)}$, hence
$\#(B)$-\dm al. Since $B\sbs A$, we have $\#(B)\le\#(A)$ by \era2{3.13}.
\endproof

\brm1.28
We showed in Theorem \rf{t1.29} that a finitely generated \vs~$V$ over a
field~$K$ possesses a basis~$B$, hence $V$~is finite-\dm al. Conversely, if
$(K,V)$ is a $d$-\dm al \vs\ over~$K$, then $V$ possesses a basis~$B$ with
$\#(B)=d$.
\erm

\bex1.29
Let $(K,V)$ be a $d$-\dm al \vs\ over~$K$, and let $S:K^d\to V$ be a linear
\is sm. Let $\zb \ve j {[1,d]}$ be the basis of~$K^d$ introduced in \er{1.15}.
Show that $\zb{S\ve}j{[1,n]}$ is a basis of~$V$.
\eex

We now assume that $(K,V)$ is a $d$-\dm al \vs\ over~$K$ and $\zb bi{[1,d]}$
is a basis of~$V$. Let $A$ be a \nfs\ of~$V$ \st $\#(A)=d$. It turns out that
$A$~is a basis iff $\spn(A)=V$ iff $A$~is \li.

\bpr1.31 \

\hph i,ii, Let $(K,V)$ be a $d$-\dm al \vs\ over a field~$K$, and let $A$~be
a finite subset of $V\sms0$ \st $\#(A)=d$. If $\spn(A)=V$, then $A$
is a basis of~$V$. If $A$ is \li, then $A$ is a basis of~$V$.

\hph ii,i, Let $(K,V_1)$, $(K,V_2)$ be \vs s over~$K$, let
$L:V_1\to V_2$ be linear and let $A\sbs V_1\sms0$ be nonempty and finite.
If $L$ is sur\jc\ and $\spn(A)=V_1$, then $V_2=\spn(L(A))$.
If $L$ is in\jc\ and $A$ is \li, then $L(A)$ is \li.

\hph iii,, If in {\rm(ii)} $V_1$ and~$V_2$ are $d$-\dm al, then $L$~is in\jc\ iff $L$ is sur\jc.
\epr

\proof \

(i) Suppose $\spn(A)=V$. By Theorem \rf{t1.29} \te s a \nss\ $B$ of~$A$ \st
$B$~is a basis of~$V$, $V$~is $\#(B)$-\dm al and $\#(B)\le\#(A)$. Since $V$~is
$d$-\dm al, $\#(B)=d$ by Lemma \rf{l1.23}\,(i), (iv). Since $\#(A)=d$, $\#(B)=\#(A)$. By \era2{3.13}, $B=A$,
hence $A$~is a basis of~$V$.

Suppose $A$ is \li. We have to show that $\spn(A)=V$. Suppose, for \cd ion,
that $\spn(A)\ne V$. Let $x\in V\sm\spn(A)$, and set $B:=A\cup\{x\}$.
Note that $x\ne0$. We claim that
$B$~is \li. Let $\la:A\to K$ and $\a\in K$ be \st $\suml_{a\in A}\la_a a+\a x
=0$. Then $\a=0$, otherwise $x=-\a\Inv \suml_{a\in A}\la_aa = \suml_{a\in A}
\mu_aa$ with $\mu_a:=-\a\Inv \la_a$, $a\in A$, by \er{1.10}. Thus
$x\in\spn(A)$, a~\cd ion. Since $\a=0$, $\suml_{\a \in A}\la_aa=0$, hence $\la_a=0$ \fa
$a\in A$ since $A$~is \li. \If that $\la_a=0$ \fa $a\in A$ and $\a=0$, hence
$B$~is \li. Observe that $\#(B)\nda23.31 = \#(A)+\#(\{x\})=d+1$. \E\Tf $\cR_B:
K^{d+1}\to V$ defined in \er{1.27} is linear and in\jc\ by Lemma \rf{l1.23}\,(i),
(ii). Let $C$ be a basis of~$V$. Then $\cR_C:K^d\to V$ is a linear \is sm by
Lemma \rf{l1.23}\,(i), (iv). \E\Tf the map $\cR_C\Inv \circ \cR_B: K^{d+1}\to
K^d$ is linear and in\jc, \cd ing \E\Pr\ \rf{p1.11} since $d<d+1$.

(ii) We suppose $\spn(A)=V_1$, $L$ sur\jc, and show that $V_2\sbs \spn(L(A))$. Let
$y\in V_2$. Since $L$~is sur\jc, \te s $x\in V_1$ \st $Lx=y$. Since
$\spn(A)=V_1$, \te s $\la:A\to K$ \st $x=\suml_{a\in A}\la_aa$. \E\Tf\ $Lx =
L\bigl(\suml_{a\in A}\la_aa\bigr) \nde1.13 = \suml_{a\in A}(L(\la_aa))\nde1.3
= \suml_{a\in A}\la_a L(a)\in\spn(L(A))$. Thus $V_2\sbs \spn(L(A))$.
Since $\spn(L(A))\sbs V_2$, $\spn (L(A)) =V_2$.

Suppose $A$ \li\ and $L$ in\jc. Let $\la:A\to K$ and suppose $\suml_{a\in A}
\la_aL(a)=0$. Then by linearity $L\bigl(\suml_{a\in A}\la_aa\bigr)=0$. By the
in\ji\ of~$L$, $\suml_{a\in A}\la_aa=0$. Since $A$ is \li, $\la_a=0$ \fa $a\in A$.
Hence $L(A)$ is \li.

(iii) Suppose $L$ in\jc, and let $\zb ei{[1,d]}$ be a basis of~$V_1$. Then
$\{L(e_i)\}_{i\in[1,d]}$ is \li\ in~$V_2$ by~(ii). Since $L$~is in\jc,
$L|_{\zb ei{[1,d]}}:\zb ei{[1,d]} \to \zb{Le}i{[1,d]}$ is a bi\jn. Thus
$\#(\zb ei{[1,d]})=\#(\zb {Le}i{[1,d]})=d$. Since $V_2$~is $d$-\dm al, $\zb{Le}i
{[1,d]}$ is a basis of~$V_2$ by~(i). \E\Tf given $y\in V_2$, \te s $\la:[1,d]
\to K$ \st $y=\suml_{i=1}^d \la_iL(e_i)$. By the linearity of~$L$, we obtain
$y=\suml_{i=1}^d L(\la_ie_i) = L\bigl(\suml_{i=1}^d \la_ie_i\bigr)$. Set
$x:=\suml_{i=1}^d \la_ie_i$, then $Lx=y$. Since $y$~is arbitrary in~$V_2$,
$L$~is sur\jc.

Suppose $L$ is sur\jc\ and let $\zb ei{[1,d]}$ be a basis of~$V_1$. Since
$\spn(\zb{e}i{[1,d]})=V_1$, $\spn(\zb{Le}i{[1,d]})=V_2$ by~(ii). By Theorem
\rf{t1.29} \te s a subset $B$ of $\zb{Le}i{[1,d]}$ \st $B$~is a basis of~$V_2$.
Since $\dim V_2=d$, $\#(B)=d$ by Lemma \rf{l1.23}. \Mo $\#(B)\le \#(\zb{Le}i
{[1,d]})$ by \era2{1.33}, and $\#(\zb{Le}i{[1,d]})\le\#(\zb{e}i{[1,d]})$
by~\era2{3.9}. Finally, $\#(\zb{e}i{[1,d]})=d$ since $\zb{e}i{[1,d]}$ is
a basis of~$V_1$. Thus $\dim V_2=d\le\#(\zb{Le}i{[1,d]}\le d=\dim V_1=\dim
V_2$. \If that $\#(\zb{Le}i{[1,d]})=d$. Since $\spn(\zb{Le}i{[1,d]})=V_2$,
$\zb{Le}i{[1,d]}$ is a basis of~$V_2$ by~(i). Let $u,v\in V_1$ be \st
$Lu=Lv$. Let $\la,\mu:[1,d]\to K$ be \st $u=\suml_{i=1}^d \la_ie_i$ and
$v=\suml_{i=1}^d \mu_ie_i$. Then, by linearity, $\suml_{i=1}^d \la_iLe_i=
\suml_{i=1}^d L(\la_ie_i)=Lu=Lv=\suml_{i=1}^d \mu_iLe_i$. Hence $\suml_{i=1}^d
\la_i(Le_i) = \suml_{i=1}^d \mu_i(Le_i)$. Since $\zb{Le}i{[1,d]}$ is a basis
of~$V_2$, we have $\la_i=\mu_i$, $i\in[1,d]$ by \E\df\ \rf{d1.17}. \csq, $u=\suml_{i=1}^d \la_ie_i
=\suml_{i=1}^d \mu_ie_i=v$, hence $u=v$. \If that $L$~is in\jc.
\endproof

\Wanp prove the theorem about the \ca\ of a finite field mentioned at the
beginning of this section.

\bth1.30
Let $F$ be a finite field. Then \te\ $p$~prime and $n\in\Na$ \st $\#(F)=p^n$.
\Mo $\chr(F)=p$, and the \pf\ of~$F$ is \is c to the field $(\N_p,+_p,\cdot_p,
0,1)$.
\eth

\proof
Let $K$ denote the prime subfield of~$F$ introduced in Theorem \rfa4{t4.35}.
Since $K\sbs F$ and $F$ is finite, $K$~is finite and $\#(K)\le\#(F)$ by
\era2{3.13}. \Mo $\#(K)=p$ \fs prime~$p$ by the same theorem. Set
$A:=F\sms0$. Then $F=\spn(A)$ since $0\in\spn(A)$ and $A\sbs\spn(A)$. Since
$0,1\in F$, $A$~is nonempty. Since $A\sbs F$, $A$~is finite. \E\Tf the \vs\
$(K,F)$ introduced in Example \rf{xa1.2}\,(ii) is \fg\ in view of \E\df\
\rf{d1.15}\,(ii). By Theorem \rf{t1.29}, $F$~is finite-\dm al. Let $n\in\Na$
denote the \dm\ of the \vs~$F$ over~$K$. \If from \E\df s \rf{d1.13} and
\rf{d1.4} that $F$ and~$K^n$ are \ep, hence $\#(F)=\#(K^n)=\#(K^{[1,n]})$ by
Corollary \rfa2{c3.11}. Since $\#(K^{[1,n]})=\#(K)^{\#([1,n])}$ by
\era2{3.38}, we obtain $\#(F)=p^{\#([1,n])}=p^n$.
\endproof

We conclude this section by giving the simplest example of a finite field which
is not a \Pf. We look for a field $(F,+,\cdot,{\bf0},{\bf1})$ with $\#(F)=2^2$.
\If from the proof of Theorem \rf{t1.30} that if such field exists, then its
additive group $(F,+,{\bf0})$ is \is c to the product $(\N_2,+_2,0)\t(N_2,+_2,0)$,
whose \el s are $(0,0)$, $(1,0)$, $(0,1)$ and $(1,1)$. The ``table''
of $(\N_2,+_2,0)\t(N_2,+_2,0)$ is:

\def\tab#1,#2,#3,#4,#5,{\begin{array}{c|cccc}#1&#2&#3&#4&#5\\\hline
#2&#2&#3&#4&#5\\#3&#3&#2&#5&#4\\#4&#4&#5&#2&#3\\#5&#5&#4&#3&#2\end{array}}
\btb1.32 \arraycolsep8pt
$$
\tab+,{(0,0)},{(1,0)},{(0,1)},{(1,1)},
$$
\etb

We denote by ${\bf0},{\bf1},{\boldsymbol u},\bv$ the \el s of~$F$, and by $\vf$ the
monoid-\is sm from $(F,+,{\bf0})$ onto $(\N_2,+_2,0)\t(N_2,+_2,0)$. We have
$\vf({\bf0})=(0,0)$ and choose $\vf({\bf1}):=(1,0)$, $\vf({\boldsymbol u}):=(0,1)$,
$\vf(\bv):=(1,1)$. It remains to determine the ``\mlc\ table'' of the monoid $(F,\cdot,{\bf1})$. From Corollary
\rfa4{c4.42}\,(i) we know that the \mlv\ group $(F\sms{\bf0},\cdot,{\bf1})$ is cyclic, hence
\is c to $(\N_3,+_3,0)$. By Theorem \rfa4{t3.28} and \era4{3.25} it has $\F(3)=2$ \Gn s. Hence ${\boldsymbol u},\bv$
are \Gn s. Thus ${\boldsymbol u}^0={\bf1}$, ${\boldsymbol u}^1={\boldsymbol u}$, ${\boldsymbol u}^2=\bv$ and ${\boldsymbol u}^3={\bf1}$,
similarly $\bv^0={\bf1}$, $\bv^1=\bv$, $\bv^2={\boldsymbol u}$ and $\bv^3={\bf1}$. This implies
that the ``\mlc\ table'' of $(F,\cdot,{\bf1})$ is

\btb1.33 \arraycolsep8pt
$$
\begin{array}{c|cccc}
\cdot&{\bf0}&{\bf1}&{\boldsymbol u}&\bv\\\hline
{\bf0}&{\bf0}&{\bf0}&{\bf0}&{\bf0}\\
{\bf1}&{\bf0}&{\bf1}&{\boldsymbol u}&\bv\\
{\boldsymbol u}&{\bf0}&{\boldsymbol u}&\bv&{\bf1}\\
\bv&{\bf0}&\bv&{\bf1}&{\boldsymbol u}
\end{array}
$$
\etb

\bex1.34
Show that $(F,+,\cdot,{\bf0},{\bf1})$ is a field.
\eex

\newpage
\Subsubsection{Formal \pl s}\label{sss.formp}

The \ex\ of a finite field of order~$p^n$ \fe prime~$p$ and every \nm\
$n\in\Na\sms1$ is usually \es ed using the notion of \fp s. In Algebra \fp s
are called \pl s. Since the \df\ of \pl s in Analysis is not the same as in
Algebra, we prefer to use the term \ti{\fp}.

In this section $K$ denotes a fixed arbitrary field, finite or infinite, not
necessarily prime. We denote by $K^\N$ the \vs\ over~$K$, introduced in\glossary{$K^\N$}
Example \rf{xa1.2}\,(i) with $I:=\N$, of all \sq s of \el s of~$K$. We use the
notation
\beq2.1
0(i):=0,\ i\in\N,\qh{and }\1(i):=1, \ i\in\N.
\e
\E\fe $i\in\N$ we introduce the $i$-th \co\ map $\psi_i:K^\N\to K$ defined by
\beq2.2
\psi_i(a):=a(i), \q\ a\in K^\N,\ i\in\N.
\e
Clearly, if $a,b\in K^\N$,
\beq2.3
\psi_i(a)=\psi_i(b) \hbox{ \fa}i\in\N \hbox{ iff }a=b.
\e
It directly follows from the \df\ of the \ad\ and the \mlc\ by scalars that
$\psi_i$ is a \ti{linear} map \fa $i\in\N$. We also introduce the notation
$\ve_k$, $k\in\N$, for the vectors defined by
\beq2.4
\ve_k(i):=\d_{ki}, \q k,i\in\N.
\e
Thus we have
\beq2.5
\psi_i(\ve_k)=\d_{ki}, \q k,i\in\N.
\e

If $J$ is a \nfs\ of~$\N$, $\la:J\to K$ and $\suml_{j\in J}\la_j\ve_j=0$, then
$\la_k=0$ \fa $k\in J$. Indeed, let $k\in J$, then $0= \psi_k(0) = \psi_k\bigl(
\suml_{j\in J}\la_j\ve_j\bigr)= \suml_{j\in J}\la_j\psi_k(\ve_j)
= \suml_{j\in J}\la_j\d_{jk} = \la_k$.

We now introduce the notion of \ti{support\/} of an \el\ $a\in K^\N$, denoted by\index{support}
$\supp(a)$.\glossary{$\supp(a)$}

\bdf2.1
Let $a\in K^N$. Then
\beq2.6
\supp(a):= \{i\in\N: \psi_i(a)\ne0\}.
\e
The \el\ $a$ is called \ti{finitely supported\/} if $\supp(a)$ is finite
(possibly empty).
\edf

\bnt2.2
The set of finitely supported \el s of~$K^\N$ is denoted by $K^\N_0$
(nonstandard notation).\glossary{$K^\N_0$}
\ent

\bxs2.3
$\supp(0)=\vn$, $\supp({\1})=\N$, $\supp(\ve_k)=\{k\}$, $\supp(\ve_0
+\ve_1)=\{0,1\}$.
\exs

\blm2.4
Let $a,b\in K^\N_0$ and $\la\in K$. Then
\bea2.7
\supp(a+b) &\sbs \supp(a)\cup\supp(b),\\
\supp(\la a) &= \supp(a) \qh{if} \la\ne0. \lb{2.8}
\e
\E\Ip $a+b\in K^\N_0$ and $\la a\in K^\N_0$.
\elm

\bex2.5
Prove Lemma \rf{l2.4}.
\eex

\bdf2.6
Let $(K,V)$ be a \vs\ over~$K$. A~subset $W$ of~$V$ is called a
(\ti{vector})-\ti{subspace} (or \ti{\lss}) of~$V$ if \fa $x,y\in W$ and
$\la\in K$, $x+y\in W$ and $\la x\in W$.\index{linear subspace}
\edf

Note that if $W$ is a \lss\ of~$V$, then the \rt ion of the \ad\ and of the
\mlc\ by scalars to~$W$ is well defined. It is usual to denote these \rt ions
by the same letters.

\blm2.7
Let $(K,V)$ be a \vs\ over~$K$ and let $W$ be a \lss\ of~$V$. Then $(K,W)$ is
a \vs. \E\Ip $K_0^\N$ is a \vs\ over~$K$.
\elm

\bex2.8 \

\hph i,i, Prove Lemma \rf{l2.7}.

\hph ii,, Let $V_1,V_2$ be \vs s over $K$, and let $L:V_1\to V_2$ be a \lop.
Show that the range of~$L$ ($R(L):=\{y\in V_2: \hbox{\te s }x\in V_1 \hbox{
\st} Lx=y\}$) is a \lss\ of~$V_2$.

Show that the null-space of $L$ ($N(L):=\{x\in V_1: Lx=0\}$) is a \lss\
of~$V_1$.
\eex

We now investigate some \pp ies of the \vs\ $K_0^\N$. We first show that
$(K,K_0^\N)$ is \ti{infinite-\dm al\/}. Suppose, for \cd ion, that
$(K,K_0^\N)$ is $d$-\dm al with $d\in\Na$. Since the inverse of an \is sm is
an \is sm, \te s a linear in\jc\ map $f:K_0^\N \to K^d$. \E\oh \fe $n\in\Na$
the map $g:K^n \to K_0^\N$ defined by $g(\la):=\suml_{j=1}^n \la_j\ve_j$ is
linear (see the proof of Lemma \rf{l1.23}). We already observed that
$\suml_{j=1}^n \la_j\ve_j=0$ implies $\la_k=0$ \fa $k\in[1,n]$. Hence the
map~$g$ is in\jc\ by Lemma \rf{l1.23}. \If that the map $f\circ g:K^n\to K^d$
is also linear and in\jc. Choosing $n>d$ we obtain a \cd ion by \E\Pr\
\rf{p1.11}.

\bdf2.9
Let $(K,V)$ be a nontrivial \vs\ over~$K$ not necessarily \fg. A~\nss\ $A$ of
$V\sms 0$ is called a \ti{basis} of~$V$ if the \fw\ holds:

\hph i,i, Every $x\in V$ is a \lc\ (see \E\df\ \rf{d1.15}\,(i)) of \el s
of~$A$.

\hph ii,, Every \nfs\ $B$ of~$A$ is \li\ (see \E\df\ \rf{d1.26}).
\edf

\brm2.10
In view of \E\Pr\ \rf{p1.27}, \E\df\ \rf{d2.9} is consistent with \E\df s\
\rf{d1.17} and \rf{d1.19} in case $V$~is \fg.
\erm

\bpr2.11
The set $\zb\ve k\N$ is a basis of $K_0^\N$.
\epr

\proof \

(i) Clearly $0=0\ve_0$, and if $\la\ve_0=0$, $\la\in K$, then $\la=0$. Let
$x\in K_0^\N\sms 0$. We claim
\beq2.9
x=\sum_{j\in\supp(x)}\psi_j(x)\ve_j.
\e
Indeed, set $y:=\suml_{j\in\supp(x)}\psi_j(x)\ve_j$. Note that $\supp(x)
\ne\vn$. Then
\[
\psi_x(y)\nad{\er{1.13},\er{1.3}}=\sum_{j\in\supp(x)} \psi_j(x)\psi_k(e_j)\nde2.5 =
\sum_{j\in\supp(x)}\psi_j(x)\d_{kj}=\psi_k(x)
\]
whenever $k\in\supp(x)$ and
is equal to~$0$ whenever $k\notin\supp(x)$. Since $\psi_k(x)=0$ whenever
$k\notin\supp(x)$, we have $\psi_k(y)=\psi_k(x)$ \fa $k\in\N$, hence $y=x$.

(ii) Observe that the map $i\mt\ve_i$ from $\N$ into $\zb\ve j\N$ is bi\jc.
We denote this map by~$\vf$. Let $B$ be a \nfs\ of $\zb \ve j\N$, and let
$J:=\vf\Inv(B)$. We already showed that if $\la:J\to K$ is \st $\suml_{j\in J}
\la_j\ve_j=0$, then $\la_j=0$ \fa $j\in J$. Thus $B$~is \li. Since $B$~is an
arbitrary \nfs\ of $\zb\ve j\N$, we conclude that $\zb\ve j\N$ is \li.
\endproof

Since $\zb\ve j\N$ and $\N$ are \ep, the basis $\zb\ve j\N$ of~$K_0^\N$ is \ti{\ct y
infinite}.

\bnt2.12
\beq2.10
E:=\zb\ve j\N.
\e
\ent

\bex2.13
Let $X,Y$ be \vs s over a field $K$. Suppose that $X$~is $d$-\dm al and that
$\zb ei{[1,d]}$ is a basis of~$X$. Let $l:[1,d]\to Y$ be given. Define $L:X\to
Y$ by setting $Lx:=\suml_{i=1}^d x_il_i$ where $x=\suml_{i=1}^d x_ie_i$,
$x_i\in K$, $i\in[1,d]$. Show
that $L$~is the only linear map from $X$ into~$Y$ \sf ying $Le_i=l_i$,
$i\in[1,d]$.
\eex

We now consider an analogue of Exercise \rf{ex2.13} where $X:=K_0^\N$ with
basis $E$ defined in \er{2.10}. To this end we introduce the notion of degree
of an \el\ of $K_0^\N\sms0$.

If $a\in K_0^\N\sms0$, then its support is nonempty and finite. Hence by
Theorem \rfa1{t3.39} and Lemma \rfa1{l4.20}, $\supp(a)$ has a greatest \el,\glossary{$\deg(a)$}
which we denote by~$\deg(a)$:\index{degree of a formal polynomial}
\beq2.11
\deg(a):= \hbox{the greatest \el\ of }\supp(a), \q a\in K_0^\N \sms0.
\e
Let $a, b\in K_0^\N\sms0$ and $\la\in K\sms0$. By Lemma \rf{l2.4}, $a+b,
a-b, \la a \in K_0^\N$. From \er{2.6} we infer:
\bea2.12
\deg(a\pm b) &\le \max(\deg(a),\deg(b)) \q \hbox{if }a\pm b\ne0,\\
\deg(\la a) &= \deg(a). \lb{2.13}
\e
\Mo if $a\in K_0^\N\sms0$, then $a \nde2.9 = \suml_{i\in\supp(a)} a_i\ve_i
\nde2.2 = \suml_{i\in\supp(a)}\psi_i(a)\ve_i$. Since $\supp(a)\sbs [0,\deg(a)]$
and $a_j=0$ for $j\in[0,\deg(a)]\setminus\supp(a)$, we have $\suml_{i=0}
^{\deg(a)}a_i\ve_i = \suml_{i\in\supp(a)}a_i\ve_i + \suml_{i\in[0,\deg(a)]\sm
\supp(a)}a_i\ve_i = \suml_{i\in\supp(a)}a_i\ve_i + 0= \suml_{i\in\supp(a)}
a_i\ve_i$.

Similarly, if $n>\deg(a)$, we have $\suml_{i=0}^n a_i\ve_i = \suml_{i=0}^
{\deg(a)}a_i\ve_i + \suml_{i=\deg(a)+1}^n a_i\ve_i = \suml_{i=0}^{\deg(a)}
a_i\ve_i +0 = \suml_{i=0}^{\deg(a)}a_i\ve_i$.

Summarizing, let $a\in K_0^\N\sms0$ and $n\ge\deg(a)$:
\beq2.14
a = \sum_{i\in \supp(a)}a_i\ve_i= \sum_{i=0}^{\deg(a)}a_i\ve_i= \sum_{i=0}^{n}
a_i\ve_i.
\e

In the same way one proves the \fw: let $a\in K_0^\N\sms0$, let $Y$~be a \vs\
over~$K$ and let $f:\N\to Y$. Then, if $n\ge\deg(a)$:
\beq2.15
\sum_{i\in \supp(a)}a_i f_i =\sum_{i=0}^{\deg(a)}a_i f_i =\sum_{i=0}^{n}a_i f_i.
\e

\bpr2.14
Let $Y$ be a \vs\ over a field~$K$, let $E$ denote the basis of $K_0^\N$
defined in \er{2.10} and let $f:\N\to Y$. Define $L:K_0^\N\to Y$ by setting
\beq2.16
L(a):= \bca
0 &\hbox{if }a=0 ,\\
\suml_{i=0}^{\deg(a)} a_if_i &\hbox{whenever }a\ne 0.
\eca
\e
Then $L$ is \emph{linear} and is the only linear map from $K_0^\N$ into~$Y$
which \sf ies
\beq2.17
L(\ve_j)=f_j \qh{\fa} j\in\N.
\e
\epr

\proof
\ti{$L$ is linear}: \er{1.2}.

\ti{Case $a=0$ or $b=0$}: Suppose $a=0$. Then $L(a+b)=L(0+b)=L(b)=0+L(b)
\nde2.16 = L(0)+L(b)=L(a)+L(b)$. The case $b=0$ is similar.

\ti{Case $a\ne0$, $b\ne0$ and $a+b=0$}: Then $b=-a=(-1)a$. Then $\deg(b)
\nde2.13 = \deg(a)$ and $L(a+b)=L(0)\nde2.16 = 0 = \suml_{i=0}^{\deg(a)}0 =
\suml_{i=0}^{\deg(a)}(a_if_i-a_if_i) \nde1.13 = \suml_{i=0}^{\deg(a)}a_if_i
+\suml_{i=0}^{\deg(a)}(-a_if_i) \nde2.16 = L(a) + \suml_{i=0}^{\deg(a)}
b_if_i = L(a) + \suml_{i=0}^{\deg(b)}b_if_i = L(a)+L(b)$.

\ti{Case $a\ne0$, $b\ne0$ and $a+b\ne0$}: $L(a+b)\nde2.16 = \suml_{i=0}^
{\deg(a+b)}(a+b)_if_i = \suml_{i=0}^{\max(\deg(a),\deg(b))}(a+b)_if_i$
by \er{2.12} and \er{2.15}.
\bmlg
\sum_{i=0}^{\max(\deg(a),\deg(b))}(a+b)_if_i
= \sum_{i=0}^{\max(\deg(a),\deg(b))}(a_if_i+b_if_i)\\ {}\nde1.13 = \suml_{i=0}^
{\max(\deg(a),\deg(b))}a_if_i +\suml_{i=0}^{\max(\deg(a),\deg(b))}b_if_i
\nad\ast= \sum_{i=0}^{\deg(a)}a_if_i +\sum_{i=0}^{\deg(b)}b_if_i
=L(a)+L(b).
\e
In $\nad\ast=$ we used $\deg(a),\deg(b) \le \max(\deg(a),\deg(b))$ and
\er{2.15}.

\medskip

\er{1.3}, $\la=0$: $L(\la a)=L(0)=0=0 L(a)=\la L(a)$.

$a=0$: $L(\la a)=L(0)=0=\la 0 = \la L(0)=\la L(a)$.

$\la\ne0, a\ne0$: $L(\la a)=\suml_{i=0}^{\deg(\la a)}(\la a)_if_i =
\suml_{i=0}^{\deg(\la a)}\la(a_if_i)
\nde1.10 = \la \suml_{i=0}^{\deg(\la a)}a_if_i =\la L(a)$.

\er{2.17}: $L(\ve_j) = \suml_{i=0}^{\deg(\ve_j)}\ve_j(i)f_i \nde2.15 =
\suml_{i\in\supp(\ve_j)}\ve_j(i)f_i = \ve_j(j)f_j = f_j$, $j\in\N$.

\ti{\E\uq}: Let $\wt L:K_0^\N\to Y$ be linear and \sf y $\wt L(\ve_j)=f_j$,
$j\in\N$. Then $L(\ve_j)=\wt L(\ve_j)$, $j\in\N$. $\wt L(0)=0$ since $L$~is
linear and $L(0)=0$. Let $a\in K_0^\N\sms0$. Then $\wt L(a)\nde2.14 =
\wt L\bigl(\suml_{i\in\supp(a)} a_i\ve_i\bigr)\nde1.13 = \suml_{i\in\supp(a)}
\wt L(a_i\ve_i) \nde1.3 = \suml_{i\in\supp(a)}a_i \wt L(\ve_i) =\suml_{i\in\supp(a)}
a_i L(\ve_i) \nde1.3 = \suml_{i\in\supp(a)}L(a_i\ve_i) \nde1.13 =
L\bigl(\suml_{i\in\supp(a)} a_i\ve_i\bigr)=L(a)$.
\endproof

As a first application of \E\Pr\ \rf{p2.14} we show that if $X,Y$ are \vs s over
a~field~$K$ which have a \ct y infinite basis, then they are \ti{\is c}. It suffices to
show that if $X$~is a \vs\ over~$K$ with a \ct y infinite basis $\zb ei\N$,
then \te s an \is sm $\vf_X:K_0^\N\to X$. Indeed, if $\vf_Y:K_0^\N \to Y$ is an
\is sm, then $\vf_Y\circ\vf_X\Inv$ is an
\is sm from $X$ onto~$Y$.

\bco2.15
Let $X$ be a \vs\ over a field~$K$ with a \ct y infinite basis $\zb ei\N$.
Then the linear map $L:K_0^\N\to X$ \sf ying $L\ve_i= e_i$, $i\in\N$,
introduced in \E\Pr\ \rf{p2.14}, is an \emph{\is sm}.
\eco

\proof \

\ti{Sur\ji}: Let $x\in X$. If $x=0$, then $L(0)=0$. If $x\ne0$, since $\zb
ei\N$ is a basis of~$X$, \te\ a \nfs~$J$ of~$\N$ and $\la:J\to X$ \st
$x=\suml_{j\in J}\la_je_j$. Set $a:=\suml_{j\in J}\la_j\ve_j$. Then by
linearity of~$L$, $L(a)=L\bigl(\suml_{j\in J}\la_j\ve_j\bigr) \nde2.13 =
\suml_{j\in J} L(\la_j\ve_j) \nde1.3 = \suml_{j\in J}\la_jL(\ve_j) =
\suml_{j\in J} \la_je_j=x$. Thus $x$ belongs to the range of~$L$.

\ti{In\ji}: Since $L$ is linear, it suffices to show that $\{a\in K_0^\N:
L(a)=0\}=\{0\}$ by Lemma \rf{l1.23}.
We show that if $L(a)=0$ then $a=0$. If, for \cd ion, $a\ne0$, then
$\suml_{j=0}^{\deg(a)} a_j\ve_j=L(a)=0$. But the set $\zb \ve j{[0,\deg(a)]}$
is \li, hence $a_j=0$ \fa $j\in[0,\deg(a)]$, and $a=0$, a~\cd ion.
\endproof

Summarizing, we have proved

\bth2.16
If $(K,V_i)$, $i=1,2$, are \vs s over a field $K$ with a \ct y infinite basis,
then $V_1$ and~$V_2$ are linear \is c.
\eth

\E\fp s over a field $K$ are \el s of the \vs\ $K_0^\N$ endowed with an \ad al
binary \op~$\cdot$ which makes among other things $(K_0^\N,+,\cdot)$
a $K$-algebra.

\bdf2.17
Let $X,Y,Z$ be \vsf1$K$. A map $f:X\t Y\to Z$ is called \ti{bilinear} if the\index{bilinear map}
\fw\ holds:
\bea2.17a
x&\mt f(x,y) \hbox{ is linear \fe}y\in Y,\\
y&\mt f(x,y) \hbox{ is linear \fe}x\in X. \lb{2.18}
\e
\edf

\begin{dfn}[\cite{Nrs}]\lb{d2.18}
A $K$-algebra $(K,V,\cdot)$ is a \vsf0$K$ together with a \ti{bilinear} binary\index{K-algebra}
\op\ $\cdot:V\t V\to V$. The algebra is called \ti{\asc e} (resp.\ \ti{\cmt e})
if the binary \op~$\cdot$ is \asc e (resp.\ \cmt e). An \el\ $e\in V$ is called\index{unity}
a \ti{unity} if $e\cdot x=x\cdot e=x$ \fa $x\in V$. The algebra is called
\ti{finite}- (resp.\ \ti{infinite})-\ti{\dm al\/} if the \vs\ is finite-
(resp.\ infinite)-\dm al.\index{K-algebra!finite-dimensional}
\edf

One easily verifies that an algebra possesses at most one unity.

In the next lemma we collect some \pp ies of bilinear maps defined on
$K_0^\N\t K_0^\N$.

\blm2.19
Let $K$ be a field and let $V$ be a \vs\ over~$K$. Let $B:K_0^\N \t K_0^\N
\to V$ be a bilinear map.

\hph i,ii, Let $a,b\in K_0^\N$. Then
\beq2.19
B(a,b) = \bca
0 &\hbox{if }a=0 \hbox{ or }b=0,\\
\suml_{(i,j)\in[0,\deg(a)]\t[0,\deg(b)]} a_ib_jB(\ve_i,\ve_j)
&\hbox{if }a,b\in K_0^\N\sms0.
\eca
\e

\hph ii,i, $B$ is \ti{\sy ic}, that is, $B(a,b)=B(b,a)$ \fa $a,b\in K_0^\N$,
iff $B(\ve_i,\ve_j)=B(\ve_j,\ve_i)$, $i,j\in\N$.

\hph iii,, Let $V:=K_0^\N$ and let $\cdot$ be a bilinear binary \op\
on~$K_0^\N$. Then $\cdot$ is \asc e iff
\beq2.20
(\ve_i\cdot \ve_j)\cdot \ve_k = \ve_i\cdot (\ve_j\cdot \ve_k),
\q i,j,k\in\N.
\e

\hph iv,, \Mo $e\in K_0^\N$ is a unity of the $K$-algebra $K_0^\N$ iff
\beq2.21
(e\cdot\ve_i) = (\ve_i\cdot e)= \ve_i\in\N.
\e
\elm

\proof \

(i) The case $a=0$ or $b=0$ directly follows from \er{1.3}, \er{2.17a},
\er{2.18}. Let $a,b\in K_0^\N \sms 0$. Then
\bmlg
B(a,b)= B\Bigl(\suml_{i=0}^{\deg(a)}a_i\ve_i,
\suml_{j=0}^{\deg(b)}b_j\ve_j\Bigr)\nad{\er{2.17a},\er{1.13},\er{1.3}}
=\suml_{i=0}^{\deg(a)}a_iB \Bigl(\ve_i,\suml_{j=0}^{\deg(b)}b_jB(\ve_i,\ve_j) \Bigr)\\
{}\nad{\er{2.18},\er{1.13},\er{1.3}}
=\suml_{i=0}^{\deg(a)}a_i \Bigl(\suml_{j=0}^{\deg(b)}b_jB(\ve_i,\ve_j) \Bigr)
\nde1.10 = \suml_{i=0}^{\deg(a)}\Bigl(\suml_{j=0}^{\deg(b)}a_ib_j
B(\ve_i,\ve_j)\Bigr).
\e
 Then \er{2.19} follows from \era2{7.19}.

(ii) \ti{Only if\/} is trivial.

\ti{If\/}: $B(b,a) =
\suml_{(j,i)\in[0,\deg(b)]\t[0,\deg(a)]}b_ja_iB(\ve_j,\ve_i)=
\suml_{(j,i)\in[0,\deg(b)]\t[0,\deg(a)]}b_ja_iB(\ve_i,\ve_j)=\break
\suml_{(j,i)\in[0,\deg(b)]\t[0,\deg(a)]}a_ib_jB(\ve_i,\ve_j)=
\suml_{(i,j)\in[0,\deg(a)]\t[0,\deg(b)]}a_ib_jB(\ve_i,\ve_j)=B(a,b)$.

(iii) \ti{Only if\/} is trivial.

\ti{If\/}: Let $a,b,c\in K_0^\N$. If $a=0$, then $((a\cdot b)\cdot c) =
((0\cdot b)\cdot c) =0\cdot c =0 = 0\cdot(b\cdot c)= a\cdot(b\cdot c)$
by~\er{1.3}. The proofs of the cases $b=0$ and $c=0$ are similar.

We suppose $a,b,c\in K_0^\N\sms0$. Let $y,z\in K_0^\N$. Then the map $x\mt
((x\cdot y)\cdot z)$ from $K_0^\N$ into $K_0^\N$ is the \cm\ of the linear
map $x\mt (x\cdot y)$ and the linear map $\wh x\mt (\wh x\cdot z)$ from
$K_0^\N$ into~$K_0^\N$. By Lemma \rf{l1.5}\,(ii), the map $x\mt ((x\cdot y)
\cdot z)$ is linear. \E\Tf by part~(i), $((a\cdot b)\cdot c) =
\suml_{i=0}^{\deg(a)}a_i((\ve_i\cdot b)\cdot c)$. For the same reason the map
$x\mt ((y\cdot x)\cdot z)$ is linear. Thus $((a\cdot b)\cdot c)=
\suml_{i=0}^{\deg(a)}a_i \suml_{j=0}^{\deg(b)}b_j((\ve_i\cdot\ve_j),c)
\nde1.10 = \suml_{i=0}^{\deg(a)}\suml_{j=0}^{\deg(b)}a_ib_j((\ve_i\cdot\ve_j)
\cdot c)$. The map $x\mt ((y\cdot z)\cdot x)$ is linear since $\cdot$~is
bilinear. Thus $((a\cdot b)\cdot c)=\suml_{i=0}^{\deg(a)}\suml_{j=0}^{\deg(b)}
\suml_{k=0}^{\deg(c)}a_ib_jc_k((\ve_i\cdot\ve_j)\cdot \ve_k)$. Using \er{2.20}
we find $(a\cdot b)\cdot c = \suml_{i=0}^{\deg(a)}\suml_{j=0}^{\deg(b)}
\suml_{k=0}^{\deg(c)}(\ve_i\cdot(\ve_j\cdot\ve_k))$ which is equal to $(a\cdot
(b\cdot c))$ by the same arguments as above.

(iv) The \ti{only if\/} part is trivial.

\ti{If\/}: We have to show that $(a\cdot e)=(e\cdot a)=a$ \fa $a\in K_0^\N$.
The case $a=0$ is clear. Let $a\in K_0^\N\sms0$. Then $(a\cdot e)=\bigl(\bigl(
\suml_{i=0}^{\deg(a)}a_i\ve_i\bigr)\cdot e\bigr) \nad{\er{1.13},\er{1.3}}=
\suml_{i=0}^{\deg(a)}a_i(\ve_i\cdot e) \nde2.21 = \suml_{i=0}^{\deg(a)}a_i
\ve_i = a = \suml_{i=0}^{\deg(a)}a_i\ve_i \nde2.21 = \suml_{i=0}^{\deg(a)}
a_i(e\cdot \ve_i) \nad{\er{1.13},\er{1.3}} = (e\cdot a)$.
\endproof

\brm2.20
A linear map $L:K_0^\N\to V$, where $V$ is a \vs\ over~$K$, is \ti{completely
determined\/} by its \rt ion to the basis~$E$. Indeed, if $a=0$, then $La=0$,
and if $a\ne0$, then $La=\suml_{i=0}^{\deg(a)}a_iL(\ve_i)$. Part~(i) of
Lemma \rf{l2.19} tells us that the same holds for a bilinear map. A~map
$f:K_0^\N\t K_0^\N\t K_0^\N\to V$ is called \ti{trilinear} if \fa $y,z\in
K_0^\N$, the maps $x\mt f(x,y,z)$, $x\mt f(y,x,z)$ and $x\mt f(y,z,x)$ are
linear. One shows as in~(i) that if $f$~is trilinear and $a,b,c\in K_0^\N
\sms0$, then
\beq2.22
f(a,b,c) = \suml_{i=0}^{\deg(a)}\suml_{j=0}^{\deg(b)}\suml_{k=0}^{\deg(c)}
a_ib_jc_k f(\ve_i,\ve_j,\ve_k).
\e
Thus, another way to prove (iii) is to show that the maps $(a,b,c)\mt
((a\cdot b)\cdot c)$ and $(a,b,c)\mt (a\cdot(b\cdot c))$ are trilinear, and
that they coincide on $E\t E\t E$.
\erm

In \E\Pr\ \rf{p2.14} we showed that given a map $\ell$ from~$E$ into a \vs\
$Y$ over~$K$, \te s \ooo \ti{linear} map $L:K_0^\N\to Y$ \st $L(a)=\ell(a)$
\fa $a\in E$. Such a map is usually called an \ti{\ext\ by linearity}\index{extension by linearity}
to~$K_0^\N$ of the map~$\ell$ defined on the basis~$E$. It turns out that
given a map $B_0:E\t E\to V$ \te s \ooo bilinear map $B:K_0^\N\t K_0^\N\to V$
\st $B\big|_{E\t E}=B_0$.

One can reduce the proof of this \ext\ ``by linearity'' by means of the \fw\
lemma.

\blm2.21
Let $X,Y$ be \vsf1$K$, let $I$ be a \nfs\ of~$\N$, let $\a:I\to K$ and let
$L_i$, $i\in I$, be a set of linear maps from~$X$ into~$Y$. Define a map
$L:X\to Y$ by setting
\beq2.23
L(x):= \sum_{i\in I}\a_iL_i(x), \q x\in X.
\e
Then $L$ is \emph{linear}.
\elm

\bex2.22
Prove Lemma \rf{l2.21}.
\eex

\bpr2.23
Let $K$ be a field, let $V$ be a \vs\ over~$K$ and let $B_0:E\t E\to V$ be
given, where $E$ is defined in \er{2.10}. Then \te s \ooo bilinear map
$B:K_0^\N \t K_0^\N \to V$ \sf ying
\beq2.24
B(\ve_i,\ve_j) = B_0(\ve_i,\ve_j), \q i,j\in\N.
\e
\epr

\proof
In view of \E\Pr\ \rf{p2.14}, given $j\in\N$, \te s one (and only one)
\ti{linear} map $L_j:K_0^\N \to V$ \sf ying $L_j(\ve_i)=B_0(\ve_i,\ve_j)$ \fa
$i\in\N$. Set $B(x,\ve_j):=L_j(x)$ \fa $x\in K_0^\N$ and all $j\in\N$. Again
by \E\Pr\ \rf{p2.14}, given $x\in K_0^\N$, \te s one (and only one) \ti{linear}
map $\wt L_x:=K_0^\N \to V$ \sf ying $\wt L_x(\ve_j)=B(x,\ve_j)$ \fa $j\in\N$.
Define $B:K_0^\N \t K_0^\N \to V$ by setting
\beq2.25
B(x,y):= \wt L_x(y) \qh{\fa $x,y\in K_0^\N$.}
\e
Since $B(x,\ve_j)=L_j(x)$, $x\in K_0^\N$, we have $B(\ve_i,\ve_j) = B_0(\ve_i,
\ve_j)$, $i,j\in\N$. By \er{2.25} $y\mt B(x,y)$ is \ti{linear} \fa $x\in
K_0^\N$. It remains to show that $x\mt B(x,y)$ is linear \fa $y\in K_0^\N$. If
$y=0$, then $B(x,y)=0$ \fa $x\in K_0^\N$, thus $x\mt B(x,0)$ is linear.
Finally, if $y\ne0$, then $y=\suml_{j=0}^{\deg (y)} y_j\ve_j$. Hence $B(x,y) =
\suml_{j=0}^{\deg (y)}y_j B(x,\ve_j) =
\suml_{j=0}^{\deg (y)}y_j L_j(x)$, $x\in K_0^\N$. Since $L_j$'s are linear \fa
$j\in[0,\deg(y)]$, the maps $x\mt B(x,y)$ are \ti{linear} \fa $y\in K_0^\N$
by Lemma \rf{l2.21}.

The \uq\ follows from Lemma \rf{l2.19}\,(i).
\endproof

\bth2.24
Let $K$ be a field, let $K_0^\N$ be the \vs\ over~$K$ of all finitely
supported \sq s of \el s of~$K$, and let $E$ denote the basis of $K_0^\N$
defined in \er{2.10}. Then

\hph i,i, \E\te s \ooo bilinear binary \op\ on~$K_0^\N$ denoted by~$\cdot$,
\sf ying
\beq2.26
\ve_i\cdot \ve_j=\ve_{i+j}, \q i,j\in\N.
\e
The \crs\ $K$-algebra is \asc e, \cmt e with $\ve_0$ as unity.

\hph ii,, If $a,b\in K_0^\N$ and $a=0$ or $b=0$, then $a\cdot b=0$. If $a,b\in
K_0^\N\sms0$, then for $k\in\N$
\beq2.27
(a\cdot b)(k) = \bca
\suml_{\sbk{i+j=k\\(i,j)\in[0,\deg(a)]\t[0,\deg(b)]}} a_ib_j &
\hbox{whenever }k\in[0,\deg(a)+\deg(b)],\\
0 &\hbox{whenever }k >\deg(a)+\deg(b),
\eca
\e
where $\deg(a),\deg(b)$ are defined in \er{2.11},
\bga2.28
(a\cdot b)(\deg(a)+\deg(b))=a(\deg(a))\cdot b(\deg(b))\ne0,\\
a\cdot b\ne0 \hbox{ and } \deg(a\cdot b)=\deg(a)+\deg(b). \lb{2.29}
\e
\eth

\proof \

(i) \ti{\E\ex} and \ti{\uq} follow from \E\Pr\ \rf{p2.23}.

\ti{\E\asc ity} follows from Lemma \rf{l2.19}\,(iii) since $((\ve_i\cdot
\ve_j)\cdot\ve_k)\nde2.26 = (\ve_{i+j}\cdot \ve_k) \nde2.26 = \ve_{(i+j)+k}
\nda2{1.5} = \ve_{i+(j+k)} \nde2.26 = \ve_i\cdot \ve_{j+k} \nde2.26 =
\ve_i\cdot(\ve_j\cdot\ve_k)$ \fa $i,j,k\in\N$.

\ti{\E\cmt ity} follows from Lemma \rf{l2.19}\,(ii) since $\ve_i\cdot\ve_j
\nde2.26 = \ve_{i+j} \nda21.6 = \ve_{j+i} \nde2.26 = \ve_j\cdot\ve_i$ \fa
$i,j\in\N$.

\ti{$\ve_0$ is the unity\/}: follows from Lemma \rf{l2.19}\,(iv) since
$\ve_0\cdot\ve_i \nde2.26 = \ve_{0+i} \nda21.7 = \ve_i \nda21.7 =\ve_{i+0}
\nde2.26 = \ve_i\cdot\ve_0$ \fa $i\in\N$.

(ii) Let $a,b\in K_0^\N$. If $a=0$ or $b=0$, then $a\cdot b=0$ by \er{1.3}.
Suppose $a,b\in K_0^\N\sms0$. Then by \er{2.19} and \er{2.26} we have
\beq2.30
a\cdot b = \sum_{(i,j)\in [0,\deg(a)]\t[0,\deg(b)]} a_ib_j\ve_{i+j}.
\e
We now use \era2{1.127} with $\O:=[0,\deg(a)]\t[0,\deg(b)]$, $I:=[0,\deg(a)
+\deg(b)]$ and $A_k:=\{(i,j)\in\O: i+j=k\}$, $k\in I$. This is possible since
$A_k\cap A_l=\vn$ if $k\ne l$, $k,l\in I$ and $\bcl_{k\in I}A_k=\O$. The first
assertion is clear as well as $\bcl_{k\in I}A_k\sbs\O$. We now show that $\O
\sbs\bcl_{k\in I}A_k$. Let $(i,j)\in\O$. Then $i+j\ge0$, and $i+j\le \deg(a)
+\deg(b)$ since $i+j\le \deg(a)+j\le \deg(a)+\deg(b)$ by \era2{1.63},
\era2{1.64}. \E\Tf $(i,j)\in A_{i+j}$, $i+j\in I$. As an example let $\deg(a)
:=2$, $\deg(b):=3$. Then
\bgg
A_0:=\{(0,0)\},\q A_1:=\{(1,0),(0,1)\},\q A_2:=\{(2,0),(1,1),(0,2)\}, \\
A_3:=\{(2,1),(1,2),(0,3)\},\q A_4:=\{(2,2),(1,3)\},\q A_5:=\{(2,3)\},\\
\begin{picture}(130,42)(0,0)\put(-29,30){$\O:=\{$}\put
(0,30){(0,0),}\put(34,30){(0,1),}\put(68,30){(0,2),}\put(102,30){(0,3),}\put
(0,15){(1,0),}\put(34,15){(1,1),}\put(68,15){(1,2),}\put(102,15){(1,3),}\put
(0,0){(2,0),}\put(34,0){(2,1),}\put(68,0){(2,2),}\put(102,0){(2,3).\}}\put
(22,21){\line(3,2){12}}\put(56,21){\line(3,2){12}}\put(90,21){\line(3,2){12}}\put
(22,6){\line(3,2){12}}\put(56,6){\line(3,2){12}}\put(90,6){\line
(3,2){12}}\end{picture}
\e
If we set $\O:=[0,\deg(a)]\t[0,\deg(b)]$ and $I:=[0,\deg(a)+\deg(b)]$, and use
\era2{1.127}, \er{2.30} becomes
\[
a\cdot b = \sum_{k\in I}\sum_{\sbk{i+j=k\\(i,j)\in\O}}a_ib_j\ve_{i+j} =
\sum_{k\in I}\sum_{\sbk{i+j=k\\(i,j)\in\O}}a_ib_j\ve_k.
\]
Note that
\[
\sum_{\sbk{i+j=k\\(i,j)\in\O}}a_ib_j\ve_k \nde1.10 =
\Bigl(\sum_{\sbk{i+j=k\\(i,j)\in\O}}a_ib_j\Bigr)\ve_k.
\]
Using $l$ instead of~$k$ in $\suml_{k\in I}$, we find for $k\in\N$:
\[
(a\cdot b)(k) = \psi_k(a\cdot b)=\psi_k\Bigl(\sum_{l\in I}\Bigl(
\sum_{\sbk{i+j=l\\(i,j)\in\O}}a_ib_j\Bigr)\ve_l\Bigr) =
\sum_{l\in I}\Bigl(\sum_{\sbk{i+j=l\\(i,j)\in\O}}a_ib_j\Bigr)\psi_k(\ve_l).
\]
If $k\notin I$, then $\psi_k(\ve_l)=0$ \fa $l\in I$, hence $(a\cdot b)(k)=0$.
Observe $k\notin I$ iff\break $k>\deg(a)+\deg(b)$.

If $k\in I$, then $\psi_k(\ve_l)=0$ iff $k\ne l$, and $\psi_k(e_k)=1$. Hence
$(a\cdot b)(k)=\suml_{\sbk{i+j=k\\(i,j)\in\O}}a_ib_j$. This completes the
proof of \er{2.27}.

\er{2.28}: Observe that if $(i,j)\in\O$ and $i+j=\deg(a)+\deg(b)$, then
$i=\deg(a)$ and $j=\deg(b)$. Otherwise $i+j<\deg(a)+\deg(b)$ by \era2{1.63},
\era2{1.64} with strict in\et ies, and by \era2{3.11}, \era2{3.12}. \E\Tf
$\suml_{\sbk{i+j=\deg(a)+\deg(b)\\(i,j)\in\O}}a_ib_j = a(\deg(a))\cdot
b(\deg(b))$. By \df\ $a(\deg(a)),b(\deg(b))\in K\sms0$, hence their product
belongs to $K\sms0$ since $K$ is a field. This proves \er{2.28}.

\er{2.29}: Since $a\cdot b(k)=0$ \fa $k>\deg(a)+\deg(b)$ and $a\cdot b(k)\ne0$
for $k:=\deg(a)+\deg(b)$, we have $\deg(a\cdot b)=\deg(a)+\deg(b)$.

This completes the proof of Theorem \rf{t2.24}.
\endproof

Since the  binary \op\ $\cdot$ on $K_0^\N$ is \asc e, \cmt e with \nel~$\ve_0$,
$(K_0^\N,\cdot,\ve_0)$ is an \am. \Mo since $\ve_{i+1}=\ve_i\cdot\ve_1
=\ve_1\cdot\ve_i$ \fa
$i\in\N$,  $\ve_i$~is the $i$-th \ti{\IT} of~$\ve_1$ in the monoid $(K_0^\N,\cdot,
\ve_0)$ (see \E\df\ \rfa2{d1.14}). If $a,b\in K_0^\N\sms0$, then $a\cdot b\in
K_0^\N\sms0$ by \er{2.28}, \Tf\ $K_0^\N\sms0$ is a \sbm\ of $(K_0^\N,\cdot,
\ve_0)$. It turns out that the monoid $(K_0^\N\sms0,\cdot,\ve_0)$ is a \Cm.
Indeed, if $a,b,c\in K_0^\N\sms0$ \sf y $a\cdot c=b\cdot c$, then $(a- b\cdot
c)=a\cdot c-b\cdot c=0$. Since $c\ne0$, $a-b=0$, hence $a=b$. Observe that
$E$~is a \sbm\ of $(K_0^\N\sms0,\cdot,\ve_0)$. The monoid $(E,\cdot,\ve_0)$
is \ti{\pn\/} with \Gn~$\ve_1$, since $E=\bcl_{i\in\N}\{\hbox{$i$-th \IT\ of }
\ve_1\}$. The map $n\mt \ve_n$ is a monoid-\is sm of $(\N,+,0)$ onto $(E,\cdot,
\ve_0)$. If we denote $\ve_1$ by~$X$ and $\ve_i$ by~$X^i$, $i\in\N$, then for\glossary{$X^k$}
$a\in K_0^\N\sms0$, $a=\suml_{i=0}^{\deg(a)} a_iX^i$. An expression of the form
$\suml_{i=0}^{\deg(a)} a_ix^i$ where $a\in K_0^\N\sms0$, $x\in K$, is usually
called a \ti{\pl\/}. Since $X$ is not an \el\ of~$K$, $\suml_{i=0}^{\deg(a)}
a_iX^i$ is called a \ti{\fp\/} over~$K$.\index{formal polynomial}

\bdf2.25
The $K$-algebra introduced in Theorem \rf{t2.24} will be called the
\ti{$K$-algebra of \fp s over~$K$} and will be denoted by $(K,K_0^\N,\cdot)$.
\edf

Formula \er{2.27} can be rewritten in a somewhat simpler way.

\blm2.26
Let $a,b\in K_0^\N\sms0$. Then
\beq2.31
(a\cdot b)(k) = \sum_{l=0}^k a_{k-l}b_l = \sum_{l=0}^k a_lb_{k-l}.
\e
\elm

\proof
Set $m:=\max\{\deg(a),\deg(b)\}$. Observe that if $(i,j)\in([0,m]\t[0,m])
\sm\break ([0,\deg(a)]\t[0,\deg(b)])$, then $i\notin [0,\deg(a)]$ or $j\notin[0,\deg(b)]$,
hence $a_ib_j=0$. \E\Tf\
\[
\suml_{\sbk{(i,j)\in[0,\deg(a)]\t[0,\deg(b)]\\
i+j=k}}a_ib_j = \suml_{\sbk{(i,j)\in[0,m]\t[0,m]\\i+j=k}}a_ib_j.
\]
\Mo
$\{(i,j)\in[0,m]\t[0,m]: i+j=k\} = \{(i,j)\in[0,k]\t[0,k]: i+j=k\}$ since
$i+j\le k$ implies $i,j\le k$. \Mo $\{(i,j)\in[0,k]\t[0,k]: i+j=k\}=
\bcl_{l=0}^k A_l$ with $A_l:=\{(l,k-l)\}$, $l\in[0,k]$. Similarly
$\{(i,j)\in[0,k]\t[0,k]: i+j=k\}=\bcl_{l=0}^k B_l$ with $B_l:=\{(k-l,l)\}$,
$l\in[0,k]$.

Clearly the sets $A_l$ and $B_l$, $l\in[0,k]$, are pairwise disjoint. Then
\er{2.31} follows from \era2{1.127}.
\endproof

\brm2.27
The $K$-algebra of \fp s can be viewed as a (\cmt e) ring (with unity).
Indeed, $(K_0^\N,+,0)$ is an \ag, $(K_0^\N,\cdot,\ve_0)$ is an \am, and
\era4{1.1n}, \era4{1.2n}, \era4{1.3n} are con\sq s of the bilinearity
of~$\cdot$.
\erm

Usually, in Algebra, one defines the ring of (formal) \pl s over a~field~$K$,
by introducing on $K_0^\N$ the usual \ti{\ad}~$+$ by setting $(a+b)(k):=a(k)
+b(k)$, $k\in\N$, and the ``\ti{\mlc}'', which we denote for a moment by~$*$,
by setting $(a*b)(k):=\suml_{l=0}^{k}a_{k-l}b_l$, $k\in\N$. This \mlc\ is in an
appropriate sense a~``\cv''. Then $(K_0^\N,+,0)$ is easily shown to be an \ag.
One shows directly that $*$~is \asc e, \cmt e, that $\ve_0*a=a$ \fa $a\in
K_0^\N$, and that \era4{1.1n}, \era4{1.2n}, \era4{1.3n} hold. By Lemma
\rf{l2.26}, the ring of (formal) \pl s $(K_0^\N,+,*,0,\ve_0)$ coincides with
$(K_0^\N,+,\cdot,0,\ve_0)$.

\bex2.28
Show that if $(K_0^\N,+,0)$ is viewed as a \vs\ over $K$ with basis~$E$, then
the binary \op~$*$ defined in Remark \rf{r2.27} is bilinear and \sf ies
$\ve_i*\ve_j=\ve_{i+j}$ \fa $i,j\in\N$. \csq, by Lemma \rf{l2.19} the binary
\op s $*$ and~$\cdot$ coincide.
\eex

Introducing the \vs\ \sc\ on $K_0^\N$ allows us to use the basis~$E$ and to
interpret $a$ as $\suml_{k=0}^{\deg(a)}a_kX^k$, where $X^k$ has a precise
meaning. In \cite[p.~223]{Is} $X$ is called an \ti{indeterminate}.\index{indeterminate}

\brm2.29
It is interesting to observe that the ``\cv''~$*$ can also be introduced on
$(K^\N,+,0)$. It is \asc e, \cmt e and $a*\ve_0=\ve_0*a=a$ \fa $a\in K^\N$.
Then $(K^\N,+,*,0,\ve_0)$ is a (\cmt e) ring (with unity), which can also be
viewed as a \cmt e, \asc e $K$-algebra with unity. However, $E$~is not a basis
of the \vs\ $(K,K^\N)$. It can be shown that a basis of $(K,K^\N)$ has to be
un\ct e.
\erm

\bex2.30
Show that in the ring $(K_0^\N,+,\cdot,0,\ve_0)$ an \el\ $a\in K_0^\N$ is
a unity iff $a=\la \ve_0$ \fs $\la\in K\sms0$.
\eex

It turns out (see \cite{Analysis}) that in the ring $(K^\N,+,\cdot,0,\ve_0)$ an
\el\ $a\in K^\N$ is a unity iff $a(0)\ne0$.

\bnt2.31
The ring of \fp s over a field~$K$ is usually denoted by $K[x]$ or $K[X]$. We\glossary{$K[X]$}
prefer the notation $K[X]$. The ring $(K^\N,+,*,0,\ve_0)$ is usually called
the ring of ``\ti{formal power series}''\index{formal power series} \cite{Analysis} and is sometimes denoted by
$K\{x\}$ or $K\{X\}$. We shall also use the notation $K[X]$ (resp.\ $K\{X\}$)
for the \crs\ $K$-algebras.\glossary{$K\{X\}$}
\ent

\brm2.32
Both rings $K[X]$ and $K\{X\}$ have the \pp y that $a*b\ne0$ whenever $a\ne0$
and $b\ne0$. \If from \E\df\ \rfa4{d5.26} and Theorem \rfa4{t5.35} that the ring $R:=K[X]$ (resp.\ $K\{X\}$)
can be embedded in a field~$F$, i.e.\ \te s an in\jc\ ring-\hm sm $j:R\to F$.
\E\Ip if the prime field of~$K$ is finite, the \Pf\ of~$F$ is also finite, but
$F$ is \ti{infinite}. In the next chapter we are rather interested in proving
the \ex\ of \ti{finite} fields with finite \pf s.
\erm

\bex2.33
Show that a nontrivial linear subspace of a finite-\dm al \vs\ is a finite-\dm
al \vs.
\eex

\newpage
\Subsubsection{Kronecker's theorem}\label{sss.Kronecker}

In Section \ref{sss.pr.fld} we introduced finite rings $(\N_n,+_n,\cdot_n,0,1)$,
$n\ge2$, by means of the division algorithm (Theorem \rfa2{t1.38}) in the \sr\
$(\N,+,\cdot,0,1)$. For special values of~$n$ these finite rings are finite
fields. In this section we \es\ a division algorithm for \fp s on a field~$K$.
This allows us to construct finite $K$-algebras which in special cases are
fields extending the field~$K$.

\begin{thm}[Division algorithm]\lb{t3.1}
Let $K$ be a field, let $K[X]$ be the $K$-algebra of \fp s over~$K$, and let
$a\in K[X]\sms0$. Then, \fe $b\in K[X]$, \te s \ooo pair $(q,r)\in K[X]\t
K[X]$ \sf ying\/{\rm:}
\bea3.1
{}&b=q\cdot a+r, \\
&\hbox{either } r=0 \hbox{ or } \deg(r)<\deg(a). \lb{3.2}
\e
\eth

\proof \

\ti{\E\uq}: Let $a,b\in K[X]$ with $a\ne0$, and let $(q_i,r_i)\in K[X]\t
K[X]$, $i=1,2$, \sf y \er{3.1}, \er{3.2}. Set $\wt q:=q_1-q_2$, $\wt
r:=r_1-r_2$. We show $\wt q=0$ and $\wt r=0$. From $q_1\cdot a+r_1= q_2\cdot
a+r_2$ and the bilinearity of~$\cdot$ we obtain
\beq3.3
\wt q\cdot a+\wt r=0.
\e
If $\wt q=0$, then $a\cdot \wt q=0$, hence $\wt r=0$ from \er{3.3} and we are
done. Suppose, for \cd ion, that $\wt q\ne0$. Since $a\ne0$, we infer from
\er{2.29} that $\wt q\cdot a\ne 0$ and $\deg(\wt q\cdot a)=\deg(\wt q)
+\deg(a)$. Hence
\beq3.4
\deg(\wt q\cdot a)\ge \deg(a).
\e
On the other hand, $\wt r\ne0$ by \er{3.3} since $\wt q\cdot a\ne0$. Then,
either $r_1\ne0$ and ${r_2=0}$, or $r_1=0$ and $r_2\ne0$, or both $r_1$
and~$r_2$ are different from~$0$. In the first case $\wt r=r_1$, hence
$\deg(\wt r)=\deg(r_1)\nde3.2 < a$. In the second case $\wt r=-r_2$, hence $\deg(\wt
r)= \deg(-r_2)\nde2.13 = \deg(r_2)\nde3.2 < \deg(a)$. Finally, if $r_1\ne0$,
$r_2\ne0$, then $\deg(\wt r)={\deg(r_1-r_2)}\nde2.12 \le
\max(\deg(r_1),\deg(r_2))\break  \nde3.2 < \deg(a)$. Hence we have $\deg(\wt r)<a$.
\csq, $\deg(\wt q\cdot a)\nde3.1 = \deg(-\wt r)\nde2.13 = {\deg(\wt r)<a}$,
hence $\deg(\wt q\cdot a)<a$, \cd ing \er{3.4}.

\ti{\E\ex}: If $b=0$, then $(q,r):=(0.0)$ \sf ies \er{3.1}, \er{3.2}, since
$0\cdot a=0$ and $0+0=0$. If $b\ne0$, we use Lemma \rfa3{l9.12} with $(W,\le)
:=(\N,\le)$, $e:=0$ and $M:=\{k\in\N: \hbox{\te s } (q,r)\in K[X]\t K[X]$
\sf ying \er{3.1}, \er{3.2} when $\deg(b)=k\}$. Then $M=\N$ provided the \fw\
two \cn s hold:
\bea3.5
&0\in M,\\
&\hbox{\fe} k\in\Na, \hbox{ if }m\in M \hbox{ \fa}m\in[0,k), \hbox{ then
}k\in M. \lb{3.6}
\e
We first prove \er{3.5}. Let $b\in K[X]\sms0$ with $\deg(b)=0$. Then $b=b(0)
\ve_0$ with $b(0)\in K\sms0$.

If $\deg(a)=0$, then $a=a(0)\ve_0$ with $a(0)\in K\sms0$. Then $q:=b(0)a(0)
\mo \ve_0$, $r:=0$ \sf y $(b(0)a(0)\mo\ve_0)\cdot(a(0)\ve_0) \nad*= (b(0)a(0)
\mo a(0))(\ve_0\cdot \ve_0)\nde2.26 = b(0)\ve_0=b$. Here $a(0)\mo$ denotes the
inverse of $a(0)$ in the group $(K\sms0,\cdot,1)$. In $\nad*=$ we used the
bilinearity of~$\cdot$. Thus \er{3.1} holds. Since $r=0$, \er{3.2} holds.

If $\deg(a)>0$, then $q:=0$ and $r:=b$ \sf y \er{3.1} since $b=0\cdot a+b$
and \er{3.2} since $\deg(r)=\deg(b)=0<\deg(a)$. This completes the proof
of~\er{3.5}.

We now prove \er{3.6}. Let $b\in K[X]\sms0$ be \st $k:=\deg(b)>0$. We suppose
$m\in M$ \fa $m\in\zo0,k $, and we show $k\in M$.

We distinguish two cases: $\deg(b)<\deg(a)$ and $\deg(b)\ge \deg(a)$.

If $\deg(b)<\deg(a)$, then $q:=0$ and $r:=b$ \sf y \er{3.1} since $0\cdot
a+r=b$, and \er{3.2} since $r\ne0$ and $\deg(r)=\deg(b)<\deg(a)$. Then $k\in M$.

If $\deg(b)\ge \deg(a)$, we denote by $l$ the degree of~$a$, hence $a(l)\in
K\sms0$ and we define $c\in K[X]$ by setting
\beq3.7
c:=(\la\ve_{k-l}\cdot a), \q \la:= b(k)a(l)\mo \in K\sms0,
\e
which is well-defined since $l\le k$.

If $c=b$, then $q:=\la\ve_{k-l}$ and $r:=0$ \sf y \er{3.1}, \er{3.2}. Hence
$k\in M$.

If $c\ne b$, set
\beq3.8
\wh b:=b-c.
\e
Since $\wh b\ne0$, $\deg(\wh b)$ is well-defined. Since $\la\ve_{k-l}\ne0$,
$a\ne0$, we have $(\la\ve_{k-l}\cdot a)\nde2.29 \ne 0$, $\deg(\la\ve_{k-l}\cdot
a)\nde2.29 = \deg(\la\ve_{k-l})+\deg(a) \nde2.13 = \deg(\ve_{k-l})+l \nde2.11 =
(k-l)+l = k+((-l)+l) = k+0 =k$. Thus
\beq3.9
\deg(c) = k = \deg(b).
\e
We next show that the choice of~$\la$ implies that $b(k)=c(k)$.
In this case $\deg(\wh b)<k$, so $\deg(\wh b)\in M$.
We have $c(k)=(\la\ve_{k-l}\cdot a)(k) = (\la\ve_{k-l}\cdot a)(\deg(\la\ve_{k-l})
+\deg(a))\nde2.28 = (\la\ve_{k-l})(\deg(\la\ve_{k-l})\cdot a(\deg(a)))=
\la\ve_{k-l}(k- l)\cdot a(l) = \la\cdot a(l) = (b(k)a(l)\mo)a(l)=
b(k)(a(l)\mo a(l))= b(k)$.
Thus $c(k)=b(k)$. By \er{3.9} $c(n)=0=b(n)$ \fa $n>k$. Since $\wh b\ne0$
and $\wh b=b-c$, we obtain $\deg(\wh b)<k$, hence $\deg(\wh b)\in M$. \E\Tf
\te s $(\wh q,\wh r)\in K[X]\t K[X]$ \sf ying $\wh b=\wh q\cdot a+\wh r$ and
either $\wh r=0$ or $\deg(\wh r)<\deg(a)$. Thus $b=\wh b+c$ \sf ies $b=\wh q
\cdot a+\wh r + b(k)a(l)\mo \ve_{k-l}\cdot a = \bigl(\wh q +b(k)a(l)\mo
\ve_{k-l}\bigr)\cdot a +\wh r$, in view of the bilinearity of~$\cdot$. \Mo
either $\wh r=0$ or $\deg(\wh r)<\deg(a)$. Hence $(\wh q+b(k)a(l)\mo \ve_{k-l},
\wh r) \in K[X]\t K[X]$ \sf ies \er{3.1}, \er{3.2}. Thus $k\in M$. This
completes the proof of Theorem \rf{t3.1}.
\endproof

In \er{3.1} the \fp\ $b$ is usually called the \ti{dividend\/}, $a$~($\ne0$)
the \ti{divisor}, $q$~the \ti{quotient\/} and $r$ the \ti{remainder}. If\index{remainder}
$r=0$, then $b$ is called a \ti{\ml e} of~$a$ and one says that $a$
\ti{divides}~$b$.

Given a formal \pl\ $a\in K[X]\sms0$, Theorem \rf{t3.1} allows us to define\glossary{$P_a$}
a~map $P_a:K[X] \to K[X]$ by setting\glossary{$N(P_a)$}
\beq3.10
P_ab:=\hbox{the remainder of $b$ in \er{3.1}--\er{3.2}}.
\e
Observe that the map $P_a$ is the analogue of the map $\F:\N\to\N$ defined
in \era4{1.3}, where $n\in\Na$. Here $a$~plays the role of~$n$ in \era4{1.3}.
In view of \er{3.2}, the range of~$P_a$, which we denote by $R(P_a)$, lies
in the set consisting besides the zero formal \pl\ of formal \pl s of degree
less than the degree of~$a$. We shall denote this set by~$K_a$:\glossary{$K_a$}
\beq3.11
K_a:=\bigl\{r\in K[X]: r=0 \hbox{ or } r\ne 0 \hbox{ and }\deg(r)<\deg(a)\bigr\}.
\e
If $\deg(a)=0$, then $K_a=\{0\}$, \Tf we shall assume
\beq3.12
a\in K[X]\sms0 \qh{with} \deg(a)\ge1.
\e
As already observed in the proof of Theorem \rf{t3.1},
\beq3.13
P_ab=b \qh{whenever }b\in K_a.
\e
This follows from $b=0\cdot a+b$ and from the \uq\ part of Theorem \rf{t3.1}. As
a~con\sq\ we have
\beq3.14
R(P_a)=K_a \qh{and } P_a\circ P_a=P_a,
\e
where $P_a\circ P_a$ means the \cm\ of~$P_a$ and $P_a$ in the set of all
maps from $K[X]$ into itself. Indeed, $P_a(P_au)\nde3.13 = P_au$ since
$P_au\in K_a$, $u\in K[X]$. $P_a$~is said to be \ti{idempotent\/}.\index{idempotent}

We claim that $P_a:K[X]\to K[X]$ is a linear map. We first prove \er{1.2}.
Let $b_1,b_2\in K[X]$. Then $b_i=q_i\cdot a+P_ab_i$, $i=1,2$, \fs $q_1,q_2\in
K[X]$. Adding both \eq s we obtain $b_1+b_2=(q_1+q_2)\cdot a+(P_ab_1+P_ab_2)$
where we used the bilinearity of~$\cdot$. If $P_ab_1+P_ab_2=0$, then $P_ab_1
+P_ab_2=P_a(b_1+b_2)$. If $P_ab_1+P_ab_2\ne0$ and $P_ab_1=0$ then $P_ab_2\ne0$
and $\deg(P_ab_1+P_ab_2)=\deg(P_ab_2)<\deg(a)$, hence $P_ab_1+P_ab_2=P_a(b_1+b_2)$.
Similarly if $P_ab_1+P_ab_2\ne0$ and $P_ab_2=0$. Finally if $P_ab_1,
P_ab_2, P_a(b_1+b_2)\in K[X]\sms0$, then $\deg(P_ab_1+P_ab_2)\nde2.12 {\le}
\max(P_ab_1,P_ab_2)\nde3.2 < \deg(a)$. Thus $P_a(b_1+b_2)=P_ab_1+P_ab_2$,
which proves \er{1.2}.

Next we prove \er{1.3}. Let $\la\in K$ and $b\in K[X]$. Then $b=q\cdot a+P_ab$
\fs $q\in K[X]$. \E\ml ying by~$\la$ we obtain $\la b=(\la q)\cdot a+\la P_ab$,
using the bilinearity of~$\cdot$. If $\la=0$ or $b=0$, then $P_a(\la b)=P_a0
\nde3.13 = 0=\la b$. If $\la\ne0$ and $b\ne0$, we have $\deg(\la P_ab)
\nde2.13 = \deg(P_ab)\nde3.2 < \deg(a)$. Hence $\la P_ab=P_a(\la b)$, which
completes the proof of the \ti{linearity} of~$P_a$.

\bex3.2
Let $V$ be a \vsf0$K$ and let $P:V\to V$ be linear and \sf y $P\circ
P=P$. Show that \fe $x\in V$ \te s \ooo $y\in R(P)$, the range of~$P$,
\ooo $z\in N(P)$, the null-space of~$P$, \st $x=y+z$.
\eex

\bdf3.3
Let $V$ be a \vsf0$K$. A~linear map $P:V\to V$ \sf ying $P\circ P=P$ is called
a \ti{\pj} (on~$R(P)$ along $N(P)$).
\edf

\bex3.4
Let $V$ be a \vsf0$K$ and let $P:V\to V$ be a \pj. Define $Q:V\to V$ by setting
$Qx:=x-Px$, $x\in V$. Show that $Q$~is a \pj\ with $R(Q)=N(P)$, $N(Q)=R(P)$.
\eex

We now return to the \vs\ $K[X]$ and the \pj\ $P_a$ introduced in \er{3.10}.
\If from Theorem \rf{t3.1} that
\beq3.15
N(P_a)=\{b\in K[X]: \hbox{\te s $q\in K[X]$ \st} b=q\cdot a\}.
\e

\blm3.5
The range of $P_a$, $R(P_a)$, is a $\deg(a)$-\dm al \vs\ over~$K$ and $A:=\{\ve_i: i\in\zo
0,\deg(a) \}$ is a basis of~$R(P_a)$.
\elm

\proof
By Exercise \rf{ex3.2}, $R(P_a)$ is a \lss\ of $K[X]$ since $P_a$~is linear.
By Lemma \rf{l2.7}, $R(P_a)$ is a \vs\ over~$K$. By Theorem \rf{t3.1} $R(P_a)=
\{b\in K[X]: \hbox{either $b=0$ or }\deg(b)<\deg(a)\}$, hence $A\sbs R(P_a)$,
since $\deg(\ve_i)=i$, $i\in\N$. As a \nss\ of $\zb\ve i\N$, $A$~is \li\ since
$\zb\ve i\N$ is a basis of $K_0^\N$. Let $b\in R(P_a)$. If $b=0$, then $b=0\cdot
\ve_0$. If not, $\deg(b)<\deg(a)$, hence $b=\suml_{l=0}^{\deg(b)}\psi_l(b)\ve
_l$. Hence $R(P_a)=\spn(A)$.
Since $\#(A)=\deg(a)$, $R(P_a)$ is a $\deg(a)$-\dm al \vs\ over~$K$.
\endproof

\bnt3.6
We denote by $K_a$ the \vs\ $R(P_a)$.
\ent

As in Section \ref{sss.pr.fld} we introduce an \ad\ $+_a$ and a \mlc\ $\cdot_a$\glossary{$+_a$}\glossary{$\cdot_a$}
on~$K_a$ by setting
\bea3.16
u+_av &:= P_a(u+v), \q u,v\in K_a,\\
u\cdot_av &:= P_a(u\cdot v), \q u,v\in K_a.\lb{3.17}
\e
Observe that $u,v\in K_a$ implies $u+v\in K_a$, hence
\beq3.18
u+_a v = u+v, \q u,v\in K_a,
\e
in contrast to the situation in Section \ref{sss.pr.fld}.

The binary \op\ $\cdot_a$ in $K_a$ is the analogue of~$\cdot_n$ defined in \era4{1.33}.
We now show that the \op~$\cdot_a$ is bilinear. Indeed, let $y\in K_a$, and
consider the map $f(x):=x\cdot_a y$, $x\in K_a$. The map~$f$ is the \cm\ of
three linear maps. The first one is the inclusion $i:K_a \to K[X]$
($i(x)=x$, $x\in K_a$) which is
($K$-)linear. The second map is the map $x\mt x\cdot y$ from $K[X]$ into
itself which is ($K$-)linear by Theorem \rf{t2.24}. The third one is $P_a:K[X] \to
K_a$ which is ($K$-)linear. Indeed, if $u,v\in K_a$, then $P_a(u+v) \nde3.16 =
u+_a v\nde3.13 = P_a(u)+_a P_a(v)$. Hence \er{1.2} holds. If $\la\in K$ and
$u\in K_a$, then $\la u\in K_a$ since $K_a\nde3.14 = R(P_a)$ is a \lss\ of
$K[X]$. Thus $P_a(\la u)\nde3.13 = \la u\nde3.13 = \la P_a(u)$. Hence \er{1.3}
holds. Then $x\mt x\cdot_a y=P_a((ix)
\cdot y)$, $y\in K_a$, is linear, hence the binary \op\ $\cdot_a:K_a\t K_a
\to K_a$ is \ti{bilinear} (see \E\df\ \rf{d2.17}), and $(K_a,\cdot_a)$ is a
$\deg(a)$-\dm al $K$-algebra (see \E\df\ \rf{d2.18}). We now show that $\cdot_a$
is \ti{\asc e}. As in the case of~$\cdot_n$ we use the \fw:
\beq3.19
P_a(b\en1 \cdot b\en2) = P_a(P_ab\en1 \cdot P_ab\en2), \q b\en1,b\en2\in K[X],
\e
which we now prove. We have by Theorem \rf{t3.1}
\[
b\en i = q\en i\cdot a +P_ab\en i, \q i=1,2.
\]
Using the bilinearity of $\cdot$, we obtain $b\en1\cdot b\en2 = (q\en1\cdot
a) \cdot (q\en2\cdot a) + (q\en1\cdot a)\cdot(P_ab\en2) + (P_ab\en1)\cdot
(q\en2\cdot a) + (P_ab\en1)\cdot (P_ab\en2) \nad*= [(q\en1\cdot q\en2\cdot
a)+ (q\en1\cdot(P_ab\en2)) + (q\en2\cdot P_a(b\en1))]\cdot a + (P_a(b\en1)
\cdot P_a(b\en2)$. In $\nad*=$ we used the \asc ity and \cmt ity of~$\cdot$.
Setting
\[
q:= (q\en1\cdot q\en2\cdot a) + [(q\en1\cdot (P_a b\en2)) + (q\en2 \cdot P_a
(b\en1))],
\]
we obtain $b\en1\cdot b\en2 = q\cdot a+ P_a(b\en1)\cdot P_a(b\en2)$. Since
$q\cdot a=q\cdot a+0$, we have $P_a(q\cdot a)=0$ by Theorem \rf{t3.1}. Hence,
by the linearity of~$P_a$, we arrive at \er{3.19}.

\ti{\E\asc ity of $\cdot_a$}: Let $b\en1,b\en2,b\en3 \in K[X]$. Then $(b\en1
\cdot_a b\en2)\cdot_a b\en3 \nde3.17 =\break P_a(P_a(b\en1\cdot b\en2)\cdot b\en3)
\nde3.19 = P_a(P_a(P_a(b\en1\cdot b\en2))\cdot P_ab\en3) \nde3.14 = P_a(P_a(
b\en1 \cdot b\en2)\cdot P_a(b\en3)) \nde3.19 = P_a((b\en1\cdot b\en2)\cdot
b\en3) \nad*= P_a(b\en1\cdot(b\en2 \cdot b\en3)) \nde3.19 = P_a(P_ab\en1
\cdot P_a(b\en2\cdot b\en3)) \nde3.14 =\break P_a(P_ab\en1\cdot P_a(P_a(b\en2 \cdot
b\en3))) \nde3.19 = P_a(b\en1\cdot P_a(b\en2\cdot b\en3)) \nde3.17 =
b\en1\cdot_a  (b\en2\cdot_a b\en3)$. In $\nad*=$ we used the \asc ity
of~$\cdot$.

The \ti{\cmt ity} of $\cdot_a$ follows from $b\en1\cdot_a b\en2 = P_a(b\en1
\cdot b\en2) \nad*= P_a(b\en2\cdot b\en1)= b\en2\cdot_a b\en1$. In $\nad*=$
we used the \cmt ity of~$\cdot$.

Finally, we show that $\ve_0$ is the unity \el\ of the $K$-algebra
$(K_a,\cdot_a)$. Indeed,\break $b\cdot_a\ve_0 \nde3.17 = P_a(b\cdot\ve_0) \nad*=
P_a(b) \nad*= P_a(\ve_0\cdot b) = \ve_0\cdot_a b$, moreover, $b=P_a(b)=P_a(b\cdot
\ve_0)= b\cdot_a \ve_0$, $b\in K_a$. In $\nad*=$ we
used the fact that $\ve_0$ is the unity of $K[X]$.

\If that $(K_a,\cdot_a)$ is a $\deg(a)$-\dm al \asc e \cmt e $K$-algebra with
unity~$\ve_0$. \E\Ip $(K_a,\cdot_a,\ve_0)$ is an \am. In view of \er{3.19}\break
$P_a: (K[X]),\cdot,\ve_0) \to (K_a,\cdot_a,\ve_0)$ is a (monoid)-\hm sm.
Indeed, $P_a\ve_0=\ve_0$ and\break $P_a(b\en1\cdot b\en2) \nde3.19 = P_a(P_a b\en1
\cdot P_ab\en2) \nde3.17 = P_ab\en1 \cdot_a P_ab\en2$, $b\en1,b\en2 \in K_a$.

Summarizing, we have proved

\bpr3.7
Let $K$ be a field, $K[X]$ be the $K$-algebra of \fp s over~$K$, let $a\in K[X]\sms0$
with $\deg(a)\ge1$ and let $P_a:K[X]\to K[X]$ be as in \er{3.10}. Then $P_a$
is a \pj\ $($see \E\df\ \rf{d3.3}$)$ \st
\[
\aligned
R(P_a) &= \spn(\zb \ve i{\zo0,\deg(a) }),\\
N(P_a) &= \{b\in K[X]: \hbox{\te s }q\in K[X] \hbox{ \st}b=q\cdot a\}.
\endaligned
\]
Let $(K,K_a)$ denote $R(P_a)$ viewed as a $\deg(a)$-\dm al \vs\ over~$K$, and
let $\cdot_a$ denote the binary \op\ in~$K_a$ defined in \er{3.17}. Then
$(K_a,\cdot_a)$ is an \asc e, \cmt e $K$-algebra with unity $\ve_0$.

\Mo the sur\jc\ \lop\ $P_a:K[X]\to K_a$ is a monoid-\hm sm from $(K[X],\cdot,
\ve_0)$ onto $(K_a,\cdot_a,\ve_0)$.
\epr

\brm3.8
If we ``forget'' the \mlc\ by scalars in the $K$-algebra $K_a$ of \E\Pr\
\rf{p3.7}, we find that $(K_a,+,\cdot_a,0,\ve_0)$ is a (\cmt e) ring (with
unity).
\erm

\bdf3.9
An algebra $(A,\circ)$ is said to be \ti{without zero divisors} if $u\circ v
\ne0$ whenever $u,v\in A\sms0$.
\edf

\blm3.10
Let $K$ be a field and let $(A,\circ)$ be a \emph{finite-\dm al} \asc e
\cmt e $K$-algebra with \ue~$e$. Then the ring $(A,+,\circ,0,e)$ is
a~\emph{field} iff $(A,\circ)$ is without zero divisors.
\elm

\brm3.11
$K[X]$ is an \ti{infinite-\dm al\/} \asc e \cmt e $K$-algebra with \ue~$\ve_0$
and without zero divisors. But $(K[X],+,\cdot ,0,\ve_0)$ is not a~field. Indeed,
there is no $b\in K[X]\sms0$ \st $\ve_1\cdot b=\ve_0$, since in this case $\deg(\ve_1\cdot b)
\nde2.29 = \deg(\ve_1)+\deg(b)\ge \deg(\ve_1)=1$ and $\deg(\ve_0)=0$.
\erm

\proof[Proof of Lemma \rf{l3.10}]\

``\ti{If}'': Let $u\in A\sms0$. Define a map $L:A\to A$ by setting $Ly:=u\circ
y$, $y\in A$. Then $L$~is a linear map on the finite-\dm
al \vs~$A$. In view of \E\Pr\ \rf{p1.31}, $L$~is sur\jc\ iff it is in\jc,
since $A$~is finite-\dm al. Since $L$~is linear, $L$~is in\jc\ iff its
null-space $N(L)=\{0\}$, by Lemma \rf{l1.23}. But $N(L):=\{z\in A: u\circ z
=0\}$. Since $(A,\circ)$ is without zero divisors, $N(L)=\{0\}$. Hence $L$~is
sur\jc\ and in particular \te s $y\in A$ \st $u\circ y=e$. \E\Tf $u$~is
a~unit (see \E\df\ \rfa4{d3.34}), and $(A,+,\circ,0,e)$ is a field.

``\ti{Only if}'': If $(A,+,\circ,0,e)$ is a field and $u,v\in A\sms0$, then
$u\circ v\in A\sms0$ since $u,v$ are units of the ring $(A,+,\circ,0,e)$.
\endproof

\begin{dfn}[\cite{Alg}, \cite{Artin}]\lb{d3.11}
A formal \pl\ $a\in K[X]\sms0$ is called \ti{\rd\ in}~$K$ if \te\ $b,c\in
K[X]\sms0$ with $\deg(b)\ge1$ and $\deg(c)\ge1$ \st $a=b\cdot c$. A~formal \pl\ $a\in
K[X]\sms0$ with $\deg(a)>0$ is called \ti{ir\rd\ in}~$K$ if it is not \rd\
in~$K$.\index{formal polynomial!reducible}\index{formal polynomial!irreducible}
\edf

\bxs3.12
Let $K$ be a field. Then $a\in K[X]$ defined by $a=\ve_2-\ve_0$ is \rd\ since
$(\ve_1+\ve_0)\cdot(\ve_1-\ve_0) = \ve_1\cdot\ve_1 - \ve_1\cdot\ve_0 + \ve_0
\cdot \ve_1 - \ve_0\cdot\ve_0=\ve_2-\ve_0$, and $\deg(\ve_1+\ve_0) = \deg(\ve_1-\ve_0)=1$.

If $a\in K[X]$ is \rd, then $a=b\cdot c$ with $b,c\in K[X]\sms0$ and $\deg(b),
\deg(c)\ge1$. \E\Tf $\deg(a)\nde2.29 = \deg(b)+\deg(c)\ge 2$. Thus $a:=\a\ve_1
+\b\ve_0$, $\a\in K\sms0$, $\b\in K$, is \ti{ir\rd} since $\deg(a)=1$.
\exs

It turns out that the irreducibility of~$a$  in \E\Pr\ \rf{p3.7} guarantees the
absence of \zd s in $(K_a,\cdot_a)$.

\bth3.13
Let $(K_a,\cdot_a)$ be the $K$-algebra introduced in \E\Pr\ \rf{p3.7}. Then
the ring $(K_a,+,\cdot_a,0,\ve_0)$ is a field iff $a$~is \emph{ir\rd}.
\eth

\proof \

\ti{Only if\/}: Suppose that the ring $(K_a,+,\cdot_a,0,\ve_0)$ is a field and
suppose for \cd ion that $a$~is \rd. Let $b,c\in K[X]\sms0$ be \st $\deg(b),
\deg(c)\ge1$ and $a=b\cdot c$. Then $b\cdot_a c=P_a(b\cdot c)=P_a(a)=0$
by~\er{3.15} with $q:=\ve_0$. Since $K_a$ is a field and $b\ne0$, $b$~is a
unity, hence \te s $\wt b\in K_a$ \st $b\cdot_a\wt b=\wt b\cdot_a b=\ve_0$.
\E\Tf $0=b\cdot_ac = \ve_0\cdot_a(b\cdot_ac) = (b\cdot_a\wt b)\cdot_a(b\cdot_a
c) \nda21.36 = b\cdot_a(\wt b\cdot_ab)\cdot c = b\cdot_a\ve_0\cdot c = b\cdot
c =a$, a~\cd ion.

\ti{If\/} (\cite[Lemma, p.~24]{Artin}): In view of \E\Pr\ \rf{p3.7} and Lemma
\rf{l3.10} it suffices to show that if $a\in K[X]\sms0$ with $\deg(a)\ge1$ is
ir\rd, then $(K_a,\cdot_a)$ is without \zd s. If $\deg(a)=1$, then $K_a =
\{b\in K[X]: b=b(0)\ve_0,\ b(0)\in K\}$ by Lemma~\rf{l3.5}. Let $b,c\in K_a
\sms0$, we show that $b\cdot_ac\ne0$. Indeed, $b=b(0)\ve_0$, $c=c(0)\ve_0$
with $b(0),c(0)\in K\sms0$ \sf y $b\cdot c=b(0)\ve_0\cdot c(0)\ve_0 =b(0)c(0)
(\ve_0\cdot\ve_0) = b(0)c(0)\ve_0 = P_a(b(0)c(0)\ve_0) = P_a(b\cdot c)=
b\cdot_ac$. Since $b(0)c(0)\ne0$, $b\cdot_ac\ne0$. We now consider the case
$\deg(a)\ge2$. Suppose for \cd ion
that \te\ $b,c\in K_a\sms0$ \st $b\cdot _ac=0$. Then by \er{3.17} we have
$P_a(b\cdot c)=0$. In view of \er{3.10} and Theorem \rf{t3.1} \te s $q\in K[X]$
\st $b\cdot c=q\cdot a$. Since $b,c\in K_a\sms0$, $b\cdot c\ne0$ by \er{2.29}.
Thus $q\ne0$. Set
\beq3.20
I:=\bigl\{u\in K_a\sms0: \hbox{\te s $\wt q\in K[X]\sms0$ \st }u\cdot c
=\wt q\cdot a\bigr\}.
\e
Then $I\ne\vn$ since $b\in I$. Define a map $\vf:I\to\N$ by setting $\vf(u):=
\deg(u)$, $u\in I$. It is well-defined since $u\in I$ implies $u\ne0$. The image
of~$\vf$, $\vf(I)$ is a \nss\ of~$\N$ since $I\ne\vn$. By Theorem \rfa1{t3.22},
$\vf(I)$~has a~least \el~$\ov n\in\N$, i.e.\ \te s $\ov b\in I$ \st
\beq3.21
\deg(\ov b)\le \deg(u) \qh{\fe}u\in I.
\e
Since $\ov b\in K_a\sms0$, we have $\ov b\in K[X]\sms0$ and using Theorem
\rf{t3.1} we obtain the \ex\ of a pair $(\ov q,\ov r)\in K[X]\t K[X]$
\sf ying
\bga3.22
a=\ov q\cdot \ov b + \ov r,\\
\hbox{either }\ov r=0 \hbox{ or }\deg(\ov r)<\deg(\ov b). \lb{3.23}
\e
We first show that $\ov r\ne0$ leads to a \cd ion. \E\ml ying both sides of
\er{3.22} by~$c$, we obtain $c\cdot a=c\cdot \ov q\cdot \ov b+ c\cdot \ov r$. Since $\ov b\in I$,
\te s $\wh q\in K[X]\sms0$ \st $\ov b\cdot c=\wh q\cdot a$. \E\Tf we have $c\cdot \ov r=c\cdot a - \ov q\cdot \ov b
\cdot c=c\cdot a - \ov q\cdot \wh q\cdot a = (c-\ov q\cdot \wh q)\cdot a$. Hence $\ov r\cdot c=(c-\ov q\cdot \wh q)\cdot a$. Since
$\ov r\ne0$ and $c-\ov q\cdot\wh q\in K[X]$, we have $c-\ov q\cdot\wh q\ne0$,
and $\ov r\in I$ by \er{3.20}.
Hence $\deg(\ov r)\ge \deg(\ov b)$ by \er{3.21}, \cd ing \er{3.23}.

We now show that $\ov r=0$ also leads to a \cd ion. We first show that $\deg(\ov q)\ge 1$.
Since $a\ne0$ and $\ov r=0$, we have $\ov q\ne0$ and we know that $\ov b\ne0$.
Since $b\in I$, $\deg(\ov b)\le \deg(b)$ by~\er{3.21}. \Mo since $b\in K_a\sms0$,
$\deg(b)<\deg(a)$ by Lemma \rf{l3.5}, Notation \rf{n3.6}. Hence $\deg(\ov b)
<\deg(a)$ by \era1{3.11}. \Mo $\deg(a)=
\deg(\ov q)+\deg(\ov b)$ by \er{2.29}. \E\Tf $\deg(\ov q)=\deg(a)-\deg(\ov b)
\ge1$. We now show that $\deg(\ov b)\ge1$. Suppose for \cd ion that $\deg(\ov b)
=0$. Then $\ov b=\ov b(0)\ve_0$ with $\ov b(0)\in K\sms0$. Since $\ov b\in I$,
we have $\ov b\cdot c = \wt q\cdot a$ \fs $\wt q\in K[X]\sms0$. Hence
$\deg(\wt q\cdot a)\nde2.29 = \deg(\wt q)+\deg(a)\ge \deg(a)$. However, $\deg(\wt q\cdot a)=\deg(\ov b\cdot c)
=\deg(\ov b_0\ve_0\cdot c) = \deg(\ov b_0c)\nde2.13 = \deg(c)< \deg(a)$ since
$c\in K_a\sms0$. A~\cd ion. \E\Tf $a=\ov q\cdot \ov b$ by \er{3.22}, with $\deg(\ov q),
\deg(\ov b)\ge1$. Hence $a$~is \rd\ \cd ing the \as\ on~$a$. \If that $b\cdot_a
c\ne0$.
\endproof

We now consider in Theorem \rf{t3.1} the important case $\deg(a)=1$ with $a:=\ve_1-\a\ve_0$,
$\a\in K$. Then either $r=0$ or $\deg(r)=0$. \E\ev tly $r=r_0\ve_0$ with
$r_0\in K$. \If that $a$~divides~$b$ iff $r_0=0$. Observe that $r_0$ only
depends on~$b$ and~$\a$.

\blm3.14
Let $K$ be as in Theorem \rf{t3.1},
let $a:=\ve_1-\a\ve_0$, $\a\in K$, and let $b\in K[X]\sms0$. Then
\beq3.24
P_a(b)= \Bigl(\sum_{i=0}^{\deg(b)} b_i\a^i\Bigr) \ve_0,
\e
where $\a^i$ is the $i$-th \IT\ of~$\a$ in $(K,\cdot,1)$, the \mlv\ monoid of
the field~$K$, and $b_i\a^i$ is the product of~$b_i$ and $\a^i$ in the field~$K$.
\elm

\proof
Since $P_a:K[X]\to K[X]$ is linear and since $b=\suml_{i=0}^{\deg(b)} b_i
\ve_i$, it suffices in view of \er{1.3}, \er{1.13} to show \fa $j\in\N$
\beq3.25
P_a(\ve_j)=\a^j \ve_0.
\e
We now prove \er{3.25} by \In\ on $j\in\N$. Let $A:=\{n\in\N: \hbox{\er{3.25}
holds with }j:=n\}$.

$0\in A$: Note that $\deg(\ve_0)=0<\deg(a)$, hence $\ve_0\in K_a$ and
$P_a(\ve_0)=\ve_0$. \Mo $\ve_0\nad{\er{1.1}\,{\rm I1}}= 1\ve_0 = \a^0\ve_0$,
since $\a^0$ is the $0$-th \IT\ of~$\a$ in the monoid $(K,\cdot,1)$.

\ti{$n\in A$ implies $n+1\in A$}: We first show $P_a(\ve_1)=\a\ve_0$.

\ti{Case $\a:=0$}: If $\a:=0$, then $a=\ve_1$. Hence by \er{2.26} $\ve_1=\ve_0
\cdot\ve_1+0$ and $P_a(\ve_1)\nde3.10 = 0 \nad{\er{1.1}\,{\rm I0}}=0\ve_0$.

\ti{Case $\a\ne0$}: We have $\ve_1=(\ve_1-\a\ve_0)+\a\ve_0 =a+\a\ve_0 =\ve_0
\cdot a+\a\ve_0$. Hence $P_a\ve_1=\a\ve_0$ since $\deg(\a\ve_0)=0 <\deg(a)$.

We now suppose $n\in A$. Then $P_a(\ve_{n+1})\nde2.26 = P_a(\ve_1\cdot\ve_n)
= P_a(\ve_1)\cdot_a P_a(\ve_n)$ by \E\Pr\ \rf{p3.7}. Thus $P_a(\ve_1\cdot\ve_
n) \nad{n\in A} = \a\ve_0\cdot_a \a^n\ve_0 = (\a\cdot\a^n)\ve_0\cdot\ve_0
\nde2.26 = \a^{n+1}\ve_0$, hence $n+1\in A$. \If that $A=\N$.
\endproof

\bco3.14
Let $K$ be a field, let $\a\in K$ and let $b\in K[X]\sms0$. Then \te s
$q\in K[X]$ \st the \fw\ holds{\rm:}
\beq3.26
b = q\cdot (\ve_1-\a\ve_0)
\e
iff
\beq3.27
\sum_{i=0}^{\deg(b)} b_i\a^i =0.
\e
\eco

\proof
By Theorem \rf{t3.1} with $a:=\ve_1-\a\ve_0$, \er{3.26} holds iff $P_a(b)=0$.
\E\oh $P_a(b)=0$ iff \er{3.27} holds by \er{3.24} (see the proof of $1\in M$
in the proof of Theorem \rf{t1.29}).
\endproof

\bdf3.15
Let $K$ be a field, let $b\in K[X]$ and let $\a\in K$. Let $\Phi_\a :K[X]\to K$
be the map defined by\glossary{$\Phi_\a$}
\beq3.28
\Phi_\a (b):=\bca
0 &\hbox{for }b=0,\\
\suml_{i=0}^{\deg(b)}b_i\a^i & \hbox{for }b\in K[X]\sms0.
\eca
\e
The map $\Phi_\a$ is usually called the \ti{evaluation at $\a$} \cite{Is}.\index{evaluation map}
\edf

\If from Lemma \rf{l3.14} that
\beq3.29
\Phi_\a(b) = \psi_0(P_{\ve_1-\a\ve_0}(b)),
\e
where $\psi_0$ is defined in \er{2.2} and $P_{\ve_1-\a\ve_0}$ in \er{3.10},
since $P_{\ve_1-\a\ve_0}(0)=0$.

\blm3.16
The evaluation map $\Phi_\a:K[X] \to K$ defined in \er{3.28} \sf ies
\bea3.30
&\Phi_\a(b+c) = \Phi_\a(b)+\Phi_\a(c), &&\kern-20pt b,c \in K[X],\\
&\Phi_\a(\la b)= \la\Phi_\a(b), &&\kern-20pt \la\in K,\ b\in K[X], \lb{3.31} \\
&\Phi_\a(b\cdot c) = \Phi_\a(b) \Phi_\a(c), &&\kern-20pt b,c\in K[X], \lb{3.32} \\
&\Phi_\a(\ve_i) = \a^i, &&\kern-20pt i\in \N, \lb{3.33} \\
&\Phi_\a(\ve_1-\a\ve_0) =0. \lb{3.34}
\e
\elm

\bex3.17
Prove Lemma \rf{l3.16} (observe
\beq3.35
\psi_0(b\cdot c) = \psi_0(b)\cdot \psi_0(c), \q b,c\in K[X],
\e
by \er{2.27}).
\eex

We now consider a simple example. Let $K:=(\N_2,+_2,\cdot_2,0,1)$ and let
$a\in K[X]\sms0$ with $\deg(a)=2$:
\beq3.36
a = \ve_2+\la\ve_1+\mu\ve_0 \qh{\fs}\la,\mu\in K.
\e
In view of \E\df\ \rf{d3.11}, $a$ is \ti{\rd} iff \te\ $b,c\in K[X]\sms0$ \st
$a=b\cdot c$ with $\deg(b),\deg(c)\ge1$. By \er{2.29} we must have $\deg(b)=
\deg(c)=1$. Hence $b=\ve_1+\rho\ve_0$, $c=\ve_1+\si \ve_0$, with $\rho,\si
\in K$. \If that $a=b\cdot c= \ve_2 + (\rho+_2\si)\ve_1 + (\rho\cdot_2\si)
\ve_0$. If $\rho=\si=1$, then $\rho+_2\si=0$ and $\rho\cdot_2\si=1$, thus
$\la=0$ and $\mu=1$. If $\rho=\si=0$, then $\la=0=\mu$. If $\rho=1$, $\si=0$
or $\rho=0$, $\si=1$, then $\la=1$, $\mu=0$. \E\Tf $a$ is \ti{ir\rd}
iff $\la=\mu=1$, i.e.
\beq3.37
a = \ve_2+\ve_1+\ve_0.
\e

Note that $\suml_{i=0}^2a_i0^i = \suml_{i=0}^20^i =1$ and
$\suml_{i=0}^2a_i1^i = \suml_{i=0}^21^i = 1+_21+_21 = 1$, where
$\suml_{i=0}^2$ is the \cme sum in $(\N_2,+_2,0)$. By Lemma \rf{l3.5},
$R(P_a)$, the range of~$P_a$, is a two-\dm al \vs\  over
$(\N_2,+_2,\cdot_2,0,1)$ with basis $\{\ve_0,\ve_1\}$. \E\Tf
\[
R(P_a) = \{0\ve_0+_a0\ve_1, 1\ve_0+_a0\ve_1, 0\ve_0+_a1\ve_1, 1\ve_0+_a\ve_1).
\]
Set ${\bf0}:=0\ve_0+_a0\ve_1$, ${\bf1}:=1\ve_0+_a0\ve_1$,
$\bu:=0\ve_0+_a1\ve_1$, $\bv:=1\ve_0+_a1\ve_1$, and
\beq3.38
F:=\{{\bf0},{\bf1},\bu,\bv\}.
\e
The ``\ad\ table'' is
\btb3.18 \arraycolsep8pt
\[
\begin{array}{c|cccc}
+_a&{\bf0}&{\bf1}&{\boldsymbol u}&\bv\\\hline
{\bf0}&{\bf0}&{\bf1}&{\bu}&{\bv}\\
{\bf1}&{\bf1}&{\bf0}&{\boldsymbol v}&\bu\\
{\boldsymbol u}&{\bu}&{\boldsymbol v}&{\bf0}&{\bf1}\\
\bv&{\bv}&\bu&{\bf1}&{\bf0}
\end{array}
\]
\etb

We now compute the ``\mlc\ table''.

Since $(F,+,{\bf0})$ is a \vs\ over~$K$ and the binary \op\ $\cdot_a :F\t F\to
F$ is bilinear, it is completely determined by its \rt ion to the basis
$\{{\bf1},\bu\}$. We have ${\bf1}\cdot_a{\bf1} \nde3.17 =
P_a(\ve_0\cdot\ve_0) \nde2.26 = P_a(\ve_0)=\ve_0={\bf1}$ since the \rt ion to
$R(P_a)$ of~$P_a$ is the identity on $R(P_a)$. Similarly, ${\bf1}\cdot_a\bu =
P_a(\ve_0\cdot\ve_1) \nde2.26 = P_a(\ve_1)=\ve_1=\bu$ and $\bu\cdot_a{\bf1} =
P_a(\ve_1\cdot\ve_0)=P_a(\ve_1)=\bu$.

The case $\bu\cdot_a\bu$ is different since $\bu\cdot_a\bu = P_a(\ve_1\cdot
\ve_1) \nde2.26 = P_a(\ve_2)$ and $\ve_2\notin R(P_a)$. However, observe that
$P_a(a)={\bf0}$ since $a=\ve_0\cdot a+0$ (see \E\Pr\ \rf{p3.7}). Hence
${\bf0} = P_a(a) = P(a)\bigl(\suml_{i=0}^2 \ve_i\bigr) = \suml_{i=0}^2
P_a(\ve_i) = P_a(\ve_0)+_a P_a(\ve_1)+_a P_a(\ve_2) = {\bf1}+_a\bu+_a
P_a(\ve_2)$. \E\Tf
\beq3.39
P_a(\ve_2) = -({\bf1}+_a\bu).
\e
Observe that ${\bf1}+_a\bu=\bv$. \Mo $\bv=1\bv$ by \er{1.1}\,I1, and
$\bv+_a\bv =1\bv+_a1\bv = (1+_21)\bv$ by \er{1.1}\,I2. Since $1+_21=0$, we
have $\bv+_a\bv = 0\bv={\bf0}$ by \er{1.1}\,I0. \If that
$-\bv=\bv$. \csq, $\bu\cdot_a \bu=P_a(\ve_2)=-\bv=\bv$. We
obtain

\btb3.19 \arraycolsep8pt
$$
\begin{array}{c|cc}
\cdot_a&{\bf1}&{\boldsymbol u}\\\hline
{\bf1}&{\bf1}&{\boldsymbol u}\\
{\boldsymbol u}&{\boldsymbol u}&\bv
\end{array}
$$
\etb

Since ${\bf0}\nad{\er{1.1}\,{\rm I0}}= 0{\bf1}$, using the bilinearity of
$\cdot_a$ we have ${\bf0}\cdot_a\bw = (0{\bf1})\cdot_a\bw = 0({\bf1}\cdot_a
\bw) \nad{\er{1.1}\,{\rm I0}}={\bf0}$ \fa $\bw \in F$. Similarly $\bw\cdot_a
{\bf0}={\bf0}$ \fa $\bw\in F$. \Mo since $\bv={\bf1}+_a\bu$, we obtain
${\bf1} \cdot_a\bv = {\bf1}\cdot_a({\bf1}+_a\bu) \nad*=
({\bf1}\cdot_a{\bf1}) + ({\bf1}\cdot_a\bu)= {\bf1}+_a\bu=\bv$. In
$\nad*=$ we used the bilinearity of~$\cdot_a$. Note that since $\cdot$ is
\cmt e, so is~$\cdot_a$. Hence $\bv\cdot_a{\bf1}=\bv$.

We next compute $\bu\cdot_a\bv$. We have $\bu\cdot_a\bv = \bu\cdot_a({\bf1}
+_a\bu) = \bu\cdot_a{\bf1} + \bu\cdot_a\bu = \bu+_a\bv = {\bf1}$.
Hence $\bu\cdot_a\bv=\bv\cdot_a\bu={\bf1}$.

It remains to compute $\bv\cdot_a\bv$. We have $\bv\cdot_a\bv = ({\bf1}+_a
\bu)\cdot_a({\bf1}+_a\bu) = ({\bf1}\cdot_a(1+_a\bu)) +_a (\bu\cdot_a(1+_a\bu))
={\bf1}\cdot_a{\bf1} + {\bf1}\cdot_a\bu + \bu\cdot_a{\bf1} + \bu\cdot_a\bu =
{\bf1}+_a (\bu+_a\bu) + \bu\cdot_a\bu = {\bf1}+_a{\bf0}+\bv = {\bf1}+\bv=\bu$.
We obtain
\btb3.20 \arraycolsep8pt
$$
\begin{array}{c|cccc}
\cdot_a&{\bf0}&{\bf1}&{\boldsymbol u}&\bv\\\hline
{\bf0}&{\bf0}&{\bf0}&{\bf0}&{\bf0}\\
{\bf1}&{\bf0}&{\bf1}&{\boldsymbol u}&\bv\\
{\boldsymbol u}&{\bf0}&{\boldsymbol u}&\bv&{\bf1}\\
\bv&{\bf0}&\bv&{\bf1}&{\boldsymbol u}
\end{array}
$$
\etb

By Theorem \rf{t3.13}, the ring $(F,+_a,\cdot_a,{\bf0},{\bf1})$ is a
\ti{field\/} since $a$~is \ti{ir\rd}. Note that Tables \rf{tb3.18} and
\rf{tb1.32} are identical, as well as Tables \rf{tb3.20} and \rf{tb1.33} (see
Exercise \rf{ex1.34}).

Set $E:=\{{\bf0},{\bf1}\}$. \If from Table \rf{tb3.18} that $E$~is a \sbm\ of
$(F,+_a,{\bf0})$. Since ${\bf0}+_a{\bf0}={\bf0}$ and ${\bf1}+_a{\bf1}={\bf0}$,
$(E,+_a,{\bf0})$ is a subgroup of $(F,+_a,{\bf0})$. \Mo it follows from Table
\rf{tb3.20} that $E$ is a \sbm\ of $(F,\cdot_a,{\bf1})$. Since ${\bf1}\cdot_a
{\bf1}={\bf1}$, $E$~is a \ti{subfield\/} of $(F,+_a,\cdot_a,{\bf0},{\bf1})$ in
view of \E\df\ \rfa4{d4.23}. By Lemma \rfa4{l4.24}, $(E,+_a,\cdot_a,{\bf0},{\bf1})$
is a~field. As a field, $E$~has obviously no proper subfield hence $E$~is
\ti{the} \pf\ of~$F$ (see \E\df\ \rfa4{d4.25} and \E\Pr\ \rfa4{p4.29}). The field
$(\N_2,+_2,\cdot_2,0,1)$ is \ti{not\/} the same as the field $(E,+_a,\cdot_a,
{\bf0},{\bf1})$ since $0,1\in\N$ and ${\bf0},{\bf1}\in K[X]$. However, they
are \ti{\is c}. Indeed, the map $j:\N_2\to F$ defined by $j(0):={\bf0}$ and
$j(1):={\bf1}$ is bi\jc\ and is a ring-\hm sm, since $j(x+_2y)= j(x)+_a j(y)$,
$j(x\cdot_2y)=j(x)\cdot_aj(y)$, $x,y\in\N_2$. If we view the map $j$ as a map
from $\N_2$ into~$F$, then $j$~is an in\jc\ ring-\hm sm.

\bdf3.21
Let $K_1$ and $K_2$ be fields and let $j:K_1\to K_2$ be an in\jc\ ring-\hm sm. Then
$(K_2,j)$ is said to be a \ti{field \ext} of~$K_1$, and $K_1$ is said to be\index{field extension}
\ti{embedded\/} in~$K_2$.
\edf

\blm3.22
Every ring-\hm sm from a field into a ring is in\jc.
\elm

\proof
Let $(K,+,\cdot,0,1)$ be a field, let $(R,+',\cdot',0',1')$ be a ring, and let
$j:K\to R$ be a ring-\hm sm. Let $x,y\in K$ be \st $j(x)=j(y)$. Suppose for
\cd ion that $x-y\ne0$. Since $K$~is a field \te s $v\in K$ \st $(x-y)\cdot v
=1$. Since $j$~is a ring-\hm sm, $1'=j(1)=j((x-y)\cdot v) = j(x+(-y))\cdot'
j(v) = (j(x)+j(-y))\cdot' j(v)$. Note that $0' = j(0) = j(y+(-y)) = j(y)
+j(-y) = j(x)+j(-y)$. Hence $1'=0'\cdot' j(v)=0'$, a \cd ion. Hence $x-y=0$,
and $y=0+y= (x+(-y))+y = x+((-y)+y) = x+0 =x$.
\endproof

Let $K_1,K_2$ and $j$ be as in \E\df\ \rf{d3.21}. We define a map $\wh\jmath
:K_1^\N \to K_2^\N$ by setting\glossary{$\jh$}
\beq3.40
\wh\jmath(b)(l):= j(b(l)), \q l\in\N,\ b\in K_1^\N.
\e
We shall show that the map $\wh\jmath$ maps $K_1[X]$ into $K_2[X]$ and we
shall still denote the \rt ion of~$\wh\jmath$ to $K_1[X]$ by~$\wh\jmath$. \Mo
we shall show that $\wh\jmath:K_1[X]\to K_2[X]$ is an in\jc\ ring-\hm sm. To
this end we introduce some notation:
\bea3.41
K_i&=(K_i,+\en i,\cdot\en i,0\en i,1\en i), \q i=1,2, \\
\ve_k\en i(l)&:=\bca
0\en i & \hbox{if }k\ne l,\\
1\en i & \hbox{if }k=l,
\eca \q k,l\in\N,\ i=1,2. \lb{3.42} \\
{\bf0}\en i(l)&:= 0\en i,\q l\in\N,\ i=1,2. \lb{3.43} \\
K_i[X]&=(K_i[X],+\en i,\cdot\en i,{\bf0}\en i,\ve_0\en i), \q i=1,2. \lb{3.44}
\e

\allowdisplaybreaks
\bpr3.23
Let $K_1,K_2,j$ be as in \E\df\ \rf{d3.21}, and let $\wh\jmath$ be as in
\er{3.40}. Then
\bea3.45
&\supp(\wh\jmath(b))=\supp(b) \qh{\fa}b\in K_1^\N.\\
&\hbox{\E\fa} b\in K_1[X],\ \wh\jmath(b)\in K_2[X] \hbox{ and } \deg(\wh\jmath
(b)) = \deg(b). \lb{3.46} \\
&\wh\jmath(b+\en1 c)=\wh\jmath(b) +\en2 \wh\jmath(c) \qh{\fa}b,c\in K_1[X].
\lb{3.47} \\
&\wh\jmath(\a b)=j(\a)\wh\jmath(b) \qh{\fa}\a\in K_1 \hbox{ and }b\in K_1[X].
\lb{3.48} \\
&\wh\jmath(b\cdot\en1 c) = \wh\jmath(b)\cdot\en2 \wh\jmath(c)\qh{\fa} b,c\in
K_1[X]. \lb{3.49} \\
&\wh\jmath:K_1[X] \to K_2[X] \qh{is in\jc.} \lb{3.50} \\
&\wh\jmath(\ve_k\en1) = \ve_k\en2 \qh{\fa}k\in\N. \lb{3.51} \\
&\wh\jmath({\bf0}\en1)={\bf0}\en2. \lb{3.52}
\e
\epr

\brm3.24
We recall that $K_i[X]$ is a ring but also a \vs\ over~$K_i$, $i=1,2$. In
\er{3.48} $\a b$ is the \mlc\ of $b\in K_1[X]$ by the scalar~$\a\in K_1$.
Similarly $j(\a)\wh\jmath(b)$ is the \mlc\ of $\wh\jmath(b)\in K_2[X]$ by the
scalar $j(\a)\in K_2$. The map $\wh\jmath$ cannot be linear if $K_1\ne K_2$.
\erm

\proof[Proof of \E\Pr\ \rf{p3.23}]\

\er{3.45}, \er{3.52}: The map $j\in K_1\to K_2$ is a ring-\hm sm, hence in
particular $j(0\en1)=0\en2$. Since it is in\jc, $j(x)=0\en2$ iff $x=0\en1$.
Let $b\in K_1^\N$. If $b(l)=0\en1$, $l\in\N$, then $(\wh\jmath(b))(l)
\nde3.40 = j(b(l)) = j(0\en1)=0\en2$, hence $\N\sm\supp(b) \sbs \N\sm
\supp(\wh\jmath(b))$. Conversely, if $(\wh\jmath(b))(l)=0\en2$, $l\in\N$,
then by \er{3.40} $j(b(l))=0\en1$. Hence by what precedes $b(l)=0\en1$. Hence
$\N\sm\supp(\wh\jmath(b)) \sbs \N\sm\supp(b)$. Thus $\N\sm\supp(b)=
\N\sm\supp(\wh\jmath(b))$, which implies \er{3.45}. Since
$\supp({\bf0}\en1)=\vn$, $\supp(\wh\jmath({\bf0}\en1))=\vn$, we have
$\wh\jmath({\bf0}\en1)={\bf0}\en2$.

\er{3.46}, \er{3.51}: If $b\in K_1[X]$, $\supp(b)$ is finite (possibly
empty), hence, by \er{3.45}, $\supp(\wh\jmath(b))$ is finite and
$\wh\jmath(b) \in K_2[X]$. In this case $\deg(b)=\deg(\wh\jmath(b))$ as a
direct con\sq\ of the \df\ \er{2.11}. Let $k\in\N$, then $\supp(\ve_k\en1) =
\{k\}$. \Mo $(\wh\jmath(\ve_k))(k) \nde3.40 = j(\ve_k(k)) \nde3.42 = j(1\en1)
= 1\en2$ since $j$ is a ring-\hm sm. \If from \er{3.45} that
$\supp(j(\ve_k\en1)) =\{k\}$. Since $\wh\jmath(\ve_k\en1)=1\en2$, \er{3.51}
holds in view of \er{3.42}.

\er{3.47}: Let $b,c\in K_1[X]$ and $l\in\N$. Then $(\wh\jmath(b+\en1 c))(l)
\nde3.40 = j((b+\en1 c)(l)) \nad*= j(b(l)+\en1 c(l)) \nad{**}= j(b(l))+\en2
j(c(l)) \nde3.40 = (\wh\jmath(b))(l)+\en2 (\wh\jmath(c))(l) = (\wh\jmath(b)
+\en2 \wh\jmath(c))(l)$. In $\nad*=$ we used the \df\ of the \ad\ in $K_1[X]$
(resp.\ $K_2[X]$), and in $\nad{**}=$ we used the fact that $j$~is a ring-\hm
sm. Since $l$~is arbitrary in~$\N$, \er{3.47} holds.

\er{3.50}: Let $b,c\in K_1[X]$ be \st $\wh\jmath (b)=\wh\jmath(c)$. We show
that $b=c$. \E\fa $l\in\N$ $j(b(l))=(\wh\jmath(b))(l) = (\wh\jmath(c))(l) =
j(c(l))$. Since $j$ is in\jc, $b(l)=c(l)$ \fa $l\in\N$, hence $b=c$.

\er{3.48}: Let $\a\in K_1$, $b\in K_1[X]$, and $l\in\N$. Then $(\wh\jmath(\a b))
(l) \nde3.40 = j((\a b)(l)) = j(\a\cdot\en1 b(l)) = j(\a)\cdot\en2
j(b(l)) \nde3.40 = j(\a)\cdot\en2 (\wh\jmath(b))(l) = (j(\a)\wh\jmath
(b))(l)$. Since $l\in\N$ is arbitrary, \er{3.48} holds.

\er{3.49}: If $b$ or $c$ is equal to ${\bf0}\en1$, then $b\cdot\en1 c={\bf0}
\en1$ and $\wh\jmath({\bf0}\en1)={\bf0}\en2$. We suppose $b,c\in
K_1[X]\sms{{\bf0}\en1}$. Then by \er{2.19} and \er{2.26}, we have
\[
b\cdot\en1 c = \sum_{(k,l)\in[0,\deg(b)]\t[0,\deg(c)]}b(k)c(l)\ve\en1_{k+l}.
\]
Thus by \era2{1.131} with $(X,\qu,e):=(K_1[X],+\en1,{\bf0}\en1)$, $(\wt X,\tqu,
\wt e):=(K_2[X],+\en2, {\bf0}\en2)$ and $\vf:=\wh\jmath$, we obtain
\[
\wh\jmath\Bigl(\sum_{(k,l)\in[0,\deg(b)]\t[0,\deg(c)]}(b(k)\cdot\en1 c(l))
\ve\en1_{k+l}\Bigr) = \sum_{(k,l)\in[0,\deg(b)]\t[0,\deg(c)]}\wh\jmath((b(k)\cdot\en1
c(l))\ve\en1_{k+l})
\]
since \era2{1.130} holds as a con\sq\ of \er{3.47}. \Mo we have $\wh\jmath((b(k)
\cdot\en1 c(l))\ve\en1_{k+l}) \nde3.48 = j(b(k)\cdot\en1 c(l)) \wh\jmath
(\ve\en1_{k+l}) \nde3.51 = j(b(k)\cdot\en1 c(l))\ve\en2_{k+l}$. Since $j$~is
a ring-\hm sm, we obtain $j(b(k)\cdot\en1 c(l))\ve\en2_{k+l} = (j(b(k))\cdot\en2
j(c(l)))\ve\en2_{k+l} \nde2.25 = (j(b(k))\cdot\en2 j(c(l)))\ve\en2_k \cdot\en2
\ve\en2_l \nad{\er{1.1}\,{\rm I3}} = j(b(k))(j(c(l))(\ve\en2_k\cdot\en2\ve\en2_l))
\nad*= j(b(k))(\ve\en2_k\cdot\en2 j(c(l))\ve\en2_l) \nad*= (j(b(k))\ve\en2_k
\cdot\en2 j(c(l))\ve\en2_l)$. In $\nad*=$ we used the bilinearity of $\cdot\en2$.
\E\Tf
\bmlg
\wh\jmath(b\cdot\en1c) = \sum_{(k,l)\in[0,\deg(b)]\t[0,\deg(c)]}
(j(b(k))\ve\en2_k \cdot\en2 j(c(l))\ve\en2_l) \\
{}=\sum_{k\in[0,\deg(b)]}\sum_{l\in[0,\deg(c)]} (j(b(k))\ve\en2_k j(c(l))
\ve\en2_l) \nad*= \sum_{k\in[0,\deg(b)]}\Bigl(j(b(k))\ve\en2_k \cdot\en2
\sum_{l\in[0,\deg(c)]} c(l)\ve\en2_l\Bigr)\\
\nad*= \Bigl(\sum_{k\in[0,\deg(b)]}j(b(k))\ve_k\en2 \cdot\en2
\sum_{l\in[0,\deg(c)]}j(c(l))\ve_l\en2 \Bigr).
\e
Finally,
\[
\sum_{k\in[0,\deg(b)]}j(b(k))\ve_k\en2 \nde3.48 = \sum_{k\in[0,\deg(b)]}
\wh\jmath (b(k)\ve_k\en1) \nad{\era2{1.131},\er{3.47}} =
\wh\jmath\Bigl(\sum_{k\in[0,\deg(b)]}b(k))\ve_k\en1\Bigr) = \wh\jmath(b).
\]
Similarly, $\suml_{l\in[0,\deg(c)]}j(c(l))\ve_l\en2 = \wh\jmath(c)$. Hence
$\wh\jmath(b\cdot\en1 c)=\wh\jmath(b)\cdot\en2 \wh\jmath(c)$, which completes
the proof of \er{3.49}.
\endproof

We now return to our example where
\[
K_1:=(\N_2,+_2,\cdot_2,0,1), \q K_2:=\bigl(\{{\bf0},{\bf1},\bu,\bv\},+_a,\cdot_a,
{\bf0},{\bf1}\bigr) \hbox{ and } j(0):={\bf0},\ j(1):={\bf1},
\]
with $a:=\ve_2+\ve_1+\ve_0$. Recall that $\Phi_0(a)=1$ and $\Phi_1(a)=1$.

By \er{3.47} $\wh\jmath(a)=\wh\jmath(\ve_2)+_a \wh\jmath(\ve_1)+_a
\wh\jmath(\ve_0)$, and by \er{3.51}
\beq3.53
\eta_k:=\wh\jmath(\ve_k), \q k\in\N,
\e
\sf ies
\beq3.54
\eta_k(l)= \bca
{\bf0} &\hbox{if }k\ne l,\\
{\bf1} &\hbox{if }k= l.
\eca
\e
Thus $\wh\jmath(a)=\eta_2+\eta_1+\eta_0$, and
\beag
\Phi_{\bf0}(\wh\jmath(a))&= {\bf0}\cdot {\bf0}+_a {\bf0}+_a{\bf1}={\bf1},\\
\Phi_{\bf1}(\wh\jmath(a))&= {\bf1}\cdot {\bf1}+_a{\bf1}+_a{\bf1}={\bf1}.
\e
\Mo
\beag
\Phi_\bu(\wh\jmath(a)) &= \bu\cdot_a\bu +_a\bu +_a{\bf1}=
(\bv+_a \bu)+_a{\bf1} = {\bf1}+_a{\bf1}={\bf0},\\
\Phi_\bv(\wh\jmath(a)) &= \bv\cdot_a \bv +_a\bv+_a{\bf1}=
{\bf1}+_a{\bf1}={\bf0}.
\e
In view of Corollary \rf{c3.14} \te\ $q_1,q_2\in K_2[X]$ not equal to zero \st
\beq3.55
\wh\jmath(a) = q_1\cdot_a (\eta_1-\bu \eta_0)
\e
and
\beq3.56
\wh\jmath(a) = q_2\cdot_a (\eta_1-\bv \eta_0).
\e
One verifies that $\wh\jmath(a)=(\eta_1-\bu \eta_0)\cdot_a (\eta_1-\bv
\eta_0)$.

We now revisit \E\Pr\ \rf{p3.7} and Theorem \rf{t3.13}. To this end we need some
\df s.

\bdf3.25
Let $(K,+,\cdot,0,1)$ be a field and let $b$ be a nonzero \ti{formal\/} \pl\
over~$K$ with $\deg(b)\in\Na$. An element $\a\in K$ is called a \ti{root\/}\index{root of a formal polynomial}
of~$b$ in~$K$ if $\suml_{i=0}^{\deg(b)} b(i)\a^i=0$.
\edf

\bpr3.26
Let $K_1,K_2,j$ be as in \E\df\ \rf{d3.21}, let $b\in K_1[X]\sms0$, and let
$\wh\jmath(b)$  be as in \er{3.40}.

\hph i,ii, If $\a\in K_1$ is a root of~$b$ in $K_1$, then $j(\a)\in K_2$ is a root of
$\wh\jmath(b)$ in~$K_2$.

\hph ii,i, If $j$ is a bi\jn, then $\a\in K_1$ is a root of~$b$ in $K_1$ iff $j(\a)\in
K_2$ is a root of $\wh\jmath(b)$ in~$K_2$.

\hph iii,, Let $K$ be a field and let $c\in K[X]\sm \{0\}$ of degree $m\in\Na$.
Then $c$ has at most $m$~roots in~$K$.
\epr

\proof[Proof of \rm(iii)]
Suppose, for \cd ion, that $c$ possesses $m+1$ roots in $K$:\break ${\zb xi{[1,m+1]}
\sbs K}$ where $i\ne j$, $i,j\in[1,m]$, implies $x_i\ne x_j$. Then $\suml_{k=0}
^m b_k(x_i)^k=0$. By Theorem \rfa4{t4.48} we have $\suml_{k=0}^m b_k(x_{m+1})^k
\ne0$. A~\cd ion.
\endproof

\bex3.27
Prove \E\Pr\ \rf{p3.26} (i) and (ii).
\eex

\bex3.28
Let $K$ be a field and $b$ be a \fp\ of degree two. Show that $b$ is ir\rd\
iff it has no roots in~$K$. Give an example of a \fp\ of degree four which is \rd\
and has no roots in~$K$.
\eex

In Theorem \rf{3.14} we showed that the ring $(K_a,+,\cdot_a,0,\ve_0)$
introduced in \E\Pr\ \rf{p3.7} is a field if (and only if) the \fp\ $a\in K[X]$
is ir\rd. In view of Corollary \rf{c3.14}, $a$~has no root in~$K$ if it is
ir\rd. In the next theorem we show that $(K_a,j)$, where $j:K\to K_a$ is
defined in \er{3.57} below, is a field \ext\ of~$K$ and that $\wh\jmath(a)\in
K_a[X]$ defined in \er{3.40} has a root in~$K_a$.

\bth3.29
Let $K$ be a field, $K[X]$ be the $K$-algebra of \fp s over~$K$ and let $a\in
K[X]$ be ir\rd\ with $\deg(a)\ge2$. Then the $K$-algebra $(K_a,+_a,\cdot_a,
0,\ve_0)$ introduced in \E\Pr\ \rf{p3.7} is a \emph{field}, and the map
$j:K\to K_a$ defined by
\beq3.57
j(\a):=\a \ve_0, \q \a\in K,
\e
is an \emph{in\jc\ ring-\hm sm}.

\Mo if $\wh\jmath:K[X]\to K_a[X]$ is the map defined in \er{3.40}, then
$\ve_1\in K_a$ is a \emph{root} of $\wh\jmath(a)$ in $K_a$.
\eth

\proof
The ring $K_a$ is a field by Theorem \rf{t3.13}. We first show that $j$~is a
monoid-\hm sm from $(K,+,0)$ into $(K_a,+_a,0)$: $j(0)\nde3.57 = 0\ve_0
\nad{\er{1.1}\,{\rm I0}}=0$, $j(\a+\b)=\break(\a+\b)\ve_0 \nad{\er{1.1}\,{\rm I2}}
= \a\ve_0+\b\ve_0 = j(\a)+j(\b)\nde3.18 = j(\a)+_a j(\b)$, $\a,\b\in K$.

\ti{$j$ is a monoid-\hm sm from $(K,\cdot,1)$ into $(K_a,\cdot_a,\ve_0)$}:
$j(1) =1\ve_0\nad{\er{1.1}\,{\rm I1}}= \ve_0$. \Mo let $\a,\b\in K$, then
$j(\a\b) \nde3.57 = (\a\b)\ve_0 \nad{\er{1.1}\,{\rm I3}}= \a(\b\ve_0)
\nde2.26 = \a(\b(\ve_0\cdot_a\ve_0))\nad*= \a(\ve_0\cdot_a\b\ve_0) \nad*=
(\a\ve_0)\cdot_a(\b\ve_0)= j(\a)\cdot_a j(\b)$. In $\nad*=$ we used the
bilinearity of~$\cdot_a$.

\ti{$j$ is in\jc}: Let $\a,\b\in K$ be \st $j(\a)=j(\b)$, i.e.\ $\a\ve_0
=\b\ve_0$. We claim that $\a=\b$. Since $j:(K,+,0)\to (K_a,+_a,0)$ is a
\hm sm, it suffices by \E\Pr\ \rfa4{p3.10}\,(vi) to show that $j(\a)=0$,
$\a\in K$, implies $\a=0$. Suppose for \cd ion that $j(\a)=0$ and $\a\in
K\sms0$. Then $e_0\nad{\er{1.1}\,\rm I1} = 1e_0 = (\a\mo\a)e_0 \nad{\er{1.1}
\,\rm I3} = \a\mo(\a e_0)\nde3.57 = {\a\mo 0\nad{\er{1.1}\,\rm I4} = 0}$.
Hence $e_0=0$, a~\cd ion.

\ti{$\ve_1\in K_a$ is a root of $\wh\jmath(a)$}: In this part of the proof of
Theorem \rf{t3.29} we shall use the \fw\ notations. The \ad\ in $K,K[X],K_a$ and
$K_a[X]$ will be denoted by~$+$ and the \crs\ \cme sum by~$\sum$. The \nel\
of the \ad\ will be denoted by~$0$ in~$K$, by~$\bf0$ in $K[X]$, $K_a$, and by
$\wh{\bf0}$ in $K_a[X]$. The \nel\ for the \mlc\ in~$K$ will be denoted by~$1$.
\E\fe $k\in\N$, the \el\ $\ve_k\in K[X]$ is defined by $\ve_k(l):=1$ if $k=l$
and $0$ if $k\ne l$ as in \er{2.4}. The \mlc\ in~$K$ will be denoted by
juxtaposition, i.e., if $\a,\b\in K$, then $\a\b$ denotes their product. Thus
$(K,+,(\ ),0,1)$ denotes the field~$K$. The \mlc\ by scalars in the \vs s
$(K,K[X])$ (resp.\ $(K_a, K_a[X])$) will also be denoted by juxtaposition. The
\mlc\ in $K[X]$ (resp.\ $K_a$, $K_a[X])$ will be denoted by~$\cdot$ (resp.\
$\cdot_a$, $\hat\cdot$). The \nel\ of the \mlc\ in $K[X]$ as well as in $K_a$
is~$\ve_0$. Thus the ring $K[X]$ is $(K[X],+,{\bf0},\ve_0)$ and the field
$K_a$ is $(K_a,+,\cdot_a,{\bf0},\ve_0)$. \E\fe $k\in\N$ we define $\xi_k\in
K_a[X]$ by
\beq3.58
\xi_k(l) := \bca
\ve_0 &\hbox{if }k=l,\\
{\bf0} &\hbox{if }k\ne l,\ l\in\N.
\eca
\e
The \mlc\ $\hat\cdot$ on $K_a[X]$ is by \df\ the unique \ti{bilinear} map from
$K_a[X]\t K_a[X]\to K_a[X]$ \sf ying
\beq3.59
\xi_k \mathbin{\hat\cdot} \xi_l = \xi_{k+l}, \q k,l\in\N.
\e
The \nel\ of the \mlc\ $\hat\cdot$ in $K_a[X]$ is~$\xi_0$. Thus the ring
$K_a[X]$ is $(K_a[X],+,\hat\cdot,\wh{\bf0},\xi_0)$, and $\zb\xi k\N$ is a basis
of the \vs\ $(K_a,K_a[X])$. \E\Ip if $B\in K_a[X]\sms{\hat{\bf0}}$, then
\beq3.60
B = \sum_{k=0}^{\deg(B)} B_k\xi_k.
\e
\Mo if $\Phi_z(B)$ denotes the evaluation of~$B$ at $z\in K_a$, then in view of
what precedes and \er{3.28}:
\beq3.61
\Phi_z(B)=\sum_{k=0}^{\deg(B)} B_k\cdot_a z^k,
\e
where $z^k$ is the $k$-th \IT\ of $z$ in the monoid $(K_a,\cdot_a,\ve_0)$.

After these preparations we are in a position to explicit $\wh\jmath (a)$.
Since $a=\suml_{k=0}^{\deg(a)}a_k\ve_k$, we have $\wh\jmath(a)=\suml_{k=0}
^{\deg(a)}\wh\jmath(a_k\ve_k)$ by \era2{1.131} with $(X,\qu,e):=(K[X],+,{\bf0})$,
$(\wt X,\tqu,\wt e):= (K_a[X],+,\wh{\bf0})$, $\vf:=\wh\jmath$ \sf ying
\era2{1.130} in view of \er{3.47}. Next we use \er{3.48} and obtain
$\wh\jmath(a_k\ve_k)=j(a_k)\wh\jmath(\ve_k)$, $k\in[0,\deg(a)]$. By \er{3.42}
and \er{3.51}, $\wh\jmath(\ve_k)=\xi_k$, $k\in[0,\deg(a)]$. Thus $\wh\jmath(a)
=\suml_{k=0}^{\deg(a)}j(a_k)\xi_k$. Then by \er{3.60}, \er{3.61}, $\Phi
_{\ve_1}(a)= \suml_{k=0}^{\deg(a)}j(a_k)\cdot_a k\ddtt a \ve_1$ where
$k\ddtt a \ve_1$ is the $k$-th \IT\ of~$\ve_1$ in $(K_a,\cdot_a,\ve_0)$.
Note that $\ve_1\in K_a$ since $\deg(a)\ge2$.

We recall that by \E\Pr\ \rf{p3.7} $P_a:(K[X],\cdot,\ve_0) \to (K_a,\cdot_a,
\ve_0)$ is a sur\jc\ monoid-\hm sm and $P_a(\ve_1)=\ve_1$ (since $\deg(a)\ge2$).
We also recall that $\ve_k$, $k\in\N$, is the $k$-th \IT\ of~$\ve_1$ in the
monoid $(K[X],\cdot,\ve_0)$ since $\ve_{i+1}\nde2.26 = \ve_i\cdot\ve_1
\nde2.26 = \ve_1\cdot\ve_i$, $i\in\N$, and $\ve_0$~is the \nel\ of $(K[X],
\cdot,\ve_0)$. \E\Tf by Lemma \rfa2{l1.22}, $P_a\ve_k$ is equal to the $k$-th
\IT\ of~$\ve_1$ in $(K_a,\cdot_a,\ve_0)$. Hence $k\ddtt a \ve_1=P_a \ve_k$,
$k\in\N$. \If that $\Phi_{\ve_1}(a) = \suml_{k=0}^{\deg(a)}j(a_k)\cdot_a P_a(\ve_k)$.
We recall that the binary \op~$\cdot_a$ on the \vs~$K_a$ over~$K$ is bilinear.

\E\Tf $j(a_k)\cdot_a P_a(\ve_k)\nde3.57 = (a_k\ve_0)\cdot_a P_a(\ve_k)\nad*=
a_k(\ve_0\cdot_a P_a(\ve_k)) \nad{**}= a_kP_a(\ve_k)$ where in $\nad*=$ we used
the bilinearity of $\cdot_a$ and in $\nad{**}=$ the fact that $\ve_0$ is the
\nel\ of the \mlc~$\cdot_a$. Observe that $a_kP_a(\ve_k)$ is the scalar \mlc\
of $P_a(\ve_k)$ by~$a_k$. Since $P_a$ is linear (over the field~$K$), $a_kP_a
(\ve_k) = P_a(a_k\ve_k)$. \If that $\Phi_{\ve_1}(a)= \suml_{k=0}^{\deg(a)}P_a
(a_k\ve_k)\nde1.13 = P_a\bigl(\suml_{k=0}^{\deg(a)}a_k\ve_k\bigr) \nde2.9 =
P_a(a)$. But $P_a(a)=0$ since $a=\ve_0\cdot_aa+{\bf0}$ (see \E\Pr\ \rf{p3.7}).
\If that $\ve_1$ is a \ti{root in} $K_a$ of the \fp\ $\wh\jmath(a)\in K_a[X]$
in view of \E\df\ \rf{d3.25}.
\endproof

\bco3.30
Let $K$ be a field and let $a$ be an \ti{ir\rd} \fp\ over~$K$. Then \te s
a field \ext\ $(F,j)$ of~$K$ \st $\wh\jmath(a)$ has a root in~$F$ where $\wh\jmath$
is defined in \er{3.40}.
\eco

We recall that if $a,b\in K[X]$ with $a\ne\bf0$, then we say that
\ti{$a$~divides~$b$} if \te s $q\in K[X]$ \st $b=q\cdot a$. In this case $a$~is
called a \ti{divisor} of~$b$ (see the lines below the proof of Theorem
\rf{t3.1}). Note that if $b\ne\bf0$ then $a\ne\bf0$.

\blm4.21
Let $K$ be a field and let $b\in K[X]\sms{\bf0}$ with $\deg(b)\ge2$. Then \te s $a\in
K[X]\sms{\bf0}$ \emph{ir\rd} \st $a$ divides~$b$.
\elm

\proof
(Compare with the proof of Lemma \rfa3{l9.9}.)

Let $A:=\{b\in K[X]\sms{{\bf0}}: \deg(b)\ge2$ and $b$ has no ir\rd\
divisor$\}$. We have to show $A=\vn$. Suppose, for \cd ion, that $A\ne\vn$.
Since $(\N,\le)$ is well-ordered, the set $d(A):=\{\deg(b)\in\N: b\in A\}$
has a least \el. Then \te s
$u\in A$ \st $2\le \deg(u)\le \deg(b)$ \fa $b\in A$. Since an ir\rd\ \el\ of
$K[X]$ divides itself and since $u\in A$, we infer that $u$~is \rd. \E\Tf
\te\ $v,w\in K[X]\sms{\bf0}$ with $\deg(v),\deg(w)\ge1$
\st $u=v\cdot w$. By \er{2.28} $\deg(u)=\deg(v)
+\deg(w) > \deg(v),\deg(w)$. Hence neither $v$ nor~$w$ belongs to~$A$, since
$\deg(u)\le \deg(b)$ \fa $b\in A$. \E\Tf \te s an ir\rd\ \el\ $a\in K[X]$ \st
$a$ divides~$v$. Since $v$ divides~$u$, we infer that $a$~divides~$u$. A~\cd
ion since $u\in A$. Hence $A=\vn$.
\endproof

As a consequence of Corollary \rf{c3.30} and Lemma \rf{l4.21} we obtain

\bth3.34
Let $a$ be a nonzero formal polynomial over a field $K$ of degree greater than one.
Then there exists a field extension $(F,j)$ such that the formal
polynomial $\wh\jmath(a)\in F[X]$ defined in \er{3.40} has a root in~$F$.
\eth

\proof
If $a$ is irreducible, then use Corollary \rf{c3.30}, otherwise use Lemma
\rf{l4.21}, Corollary \rf{c3.30} and \er{3.32}. Indeed, if $a$~is not ir\rd,
then by Lemma \rf{l4.21} with $b:=a$ \te\ $v,w\in K[X]\sms{{\bf0}}$ \st $v$~is
ir\rd\ and $a=v\cdot w$. Note that $\deg(w)\ne0$, otherwise $w=\la \ve_0$ \fs
$\la\in K\sms0$, hence $a=v\cdot \la\ve_0= \la v\cdot\ve_0= \la v$. One
verifies that $\la v$ is ir\rd\ \cd ing the \as\ on~$a$. We now distinguish
two cases: $\deg(v)=1$ and $\deg(v)\ge2$. If $\deg(v)=1$, then $v=v(1)(\ve_1
-(-v(1)\mo v(0))\ve_0)$ which is well-defined since $v(1)\in K\sms0$. Setting
$\a:=-v(1)\mo v(0)\in K$ we find $\Phi_\a(v)\nde3.31 = {v(1)\Phi_\a (\ve_1-\a
\ve_0)}\nde3.34 = v(1)0\nda41.1n = 0$. Hence $\Phi_\a(a) = \Phi_\a(v\cdot w)
\nde3.32 = \Phi_\a(v)\Phi_\a(w) =0\Phi_\a(w)=0$. \If that $\a$~is a~root
of~$a$ in~$K$, hence the conclusion of Theorem \rf{t3.34} holds with the field
\ext\ $(F,j)$ where $F:=K$ and $j:=\id_K$. If $\deg(v)\ge2$, we may apply
Corollary \rf{c3.30} with $b:=v$. Then \te s a field \ext\ $(F,j)$ of~$K$ and
$\b\in F$ \st $\Phi_\b(\wh\jmath(v))=0$ in~$F$ where $\Phi_\b:F[X]\to F$.
Hence $\Phi_\b(\wh\jmath(a))=\Phi_b(\wh\jmath(v\cdot w))\nde 3.32 =
\Phi_\b(\wh\jmath(v))\Phi_\b(\wh\jmath(w))= 0\Phi_\b(\wh\jmath(w))=0$ in~$F$.
Thus $\wh\jmath(a)$ has a root in~$F$.
\endproof

\brm3.31
In \cite{Artin} Theorem \rf{t3.34} is attributed to Kronecker.
\erm

In the remaining part of this section we consider \rl s between \fp s and
\pl\ maps over a field of \ch istic zero. Let $K$ be a field of \ch istic zero
and let $b\in K[X]$. The map $x\mt \Phi_x(b)$ from~$K$ into~$K$ is called a
\ti{\pl\ map}\index{polynomial map} (or mapping or \f), and is an \el\ of~$K^K$ the set of all
self-maps over~$K$. Recall that $K^K$ can be viewed as a \vs\ over~$K$ if the
\ad\ and the \mlc\ by scalars are defined as in Example \rf{xa1.2}\,(i).
We denote by $P(K)$ the set of all \pl\ maps over~$K$.

\bex3.36
Show that $P(K)$ is a \lss\ of the \vs\ $(K,K^K)$. Thus $P(K)$ is a \vs\
over~$K$.
\eex

We set:
\bea3.62
\wt0(x)&:=\Phi_x({\bf0})=0, \ \wt1(x)=\Phi_x(\ve_0)=1, \q x\in K,\\
\wt\ve_k(x)&:=\Phi_x(\ve_k), \q k\in\Na,\ x\in K, \lb{3.63}
\e
where $\bf0$, $\ve_k$ are defined in \er{2.1}, \er{2.2}.

\bpr3.37 \

\hph i,i, $\zb {\wt\ve} k \N$ is a basis of the \vs\ $(K,P(K))$.

\hph ii,, The map $\Phi:K[X]\to P(K)$ defined by
\beq3.64
(\Phi(b))(x):=\Phi_x(b), \q b\in K[X],\ x\in K,
\e
is a bi\jc\ linear map.
\epr

In the proof we shall use the \fw\ \df\ and lemma.

\bdf3.38
Let $(G,\qu,e)$ be a group and let $\pz0F$ be a field. Let $F^\t:=F\sms0$
and let $\pz2{F^\t}$ denote the \mlv\ group of the field~$F$. A~map
$\rho:G\to F^\t$ \sf ying
\beq3.65
\rho(e)=1 \qh{and }\rho(a\qu b)=\rho(a)\cdot\rho(b), \q a,b\in G,
\e
is said to be a \ti{\ch} of $G$ in $F$. The \ch\ $\rho(x)=1$ \fa $x\in G$ is\index{character}
called the \ti{trivial\/} \ch.
\edf

\begin{lem}[Independence of distinct \ch s) (see \cite{Artin}, Th.~12 and
\cite{Alg}, p.~156] \lb{l3.39}
Let $\zb\rho i {[0,n]}$, $n\in\N$, be a set of \emph{distinct\/} \ch s of
a~group~$(G,\qu,e)$ in a~field~$F$ $($here \emph{distinct\/} means that the map $i\mt
\rho_i$ is in\jc$)$. Let $\zb\la i{[0,n]}$ be a set of \el s of~$F$.
\beq1.117
\hbox{If\/ $\suml_{i=0}^n \la_i\rho_i(y)=0$ $($in $F)$ \fa $y\in G$,
then $\la_i=0$ \fa $i\in[0,n]$.}
\e
\elm

\proof
We use induction on $n\in\N$. Set
\[
M:=\{n\in\N: \hbox{\er{1.117} holds}\}.
\]
We have $0\in M$. Indeed, suppose $\la_0\rho_0(y)=0$ \fa $y\in F$. Choose
$y:=e$, then by \er{3.65} $0=\la_0\rho_0(e)=\la_01=\la_0$. We now suppose
$m\in M$ and we show that $m+1\in M$. Assume that
\beq1.118
\sum_{i=0}^{m+1}\la_i\rho_i(y)=0 \qh{\fa} y\in G,
\e
where $\la_i\in F$, $i\in[0,m+1]$. Replacing $y$ by $ay$ in \er{1.118}, where
$a\in G$ is arbitrary, we obtain using \er{3.65}
\beq1.119
\sum_{i=0}^{m+1} \la_i\rho_i(a)\rho_i(y)=0 \qh{\fa} y,a\in G.
\e
\E\ml ying both sides of \er{1.118} by $\rho_{m+1}(a)$ we get
\beq1.120
\sum_{i=0}^{m+1}\la_i\rho_{m+1}(a)\rho_i(y)=0 \qh{\fa}y,a\in G.
\e
Subtracting \er{1.120} from \er{1.119} and using \era4{4.17} with $\a:=-1$ we
obtain
\beq1.121
\sum_{i=0}^m \la_i(\rho_i(a)-\rho_{m+1}(a))\rho_i(y)=0 \qh{\fa}y,a\in G.
\e
From the induction hypothesis $m\in M$ we infer
\beq1.122
\la_i(\rho_i(a)-\rho_{m+1}(a))=0 \qh{\fa} a\in G \hbox{ and }i\in[0,m].
\e
Since the \ch s $\zb\rho i {[0,m+1]}$ are distinct, given $i\in[0,m]$ \te s
$a\in G$ \st
\beq1.123
\rho_i(a)\ne\rho_{m+1}(a).
\e
\E\Tf \fa $i\in[0,m]$ we have $\la_i=0$. Finally, substituting $\la_i=0$,
$i\in[0,m]$, in~\er{1.118} we obtain $\la_{m+1}\rho_{m+1}(y)=0$ \fa $y\in G$.
Choosing $y:=e$ we have $\la_{m+1}=0$ by \er{3.65}. \E\Tf $m+1\in\N$, hence
$M=\N$.
\endproof

\proof[Proof of \E\Pr\ \rf{p3.37}]\

(i) In view of \E\df\ \rf{d2.9}\,(i) we have to show that every \pl\ map
over~$K$ is a \lc\ of \el s of $\zb{\wt\ve}k\N$. Let $f\in P(K)$. By \df\ \te
s $b\in K[X]$ \st $f(x)=\Phi_x(b)$, $x\in K$. If $b=\bf0$, then $f(x)=0$,
$x\in K$, hence $f\nde3.62 = \wt0$. Thus $f=0\wt\ve_0$. If $b\ne\bf0$, then
$f(x)=\suml_{i=0}^{\deg(b)}b_ix^i \nde3.33 = \suml_{i=0}^{\deg(b)}b_i\Phi_x
(\ve_i) \nde3.63 = \suml_{i=0}^{\deg(b)}b_i\wt\ve_i(x) = \bigl(\suml_{i=0}^{\deg(b)}
b_i \wt\ve_i\bigr)(x)$, $x\in K$. Hence $f=\suml_{i=0}^{\deg(b)}b_i\wt \ve_i$.

In view of \E\df\ \rf{d2.9}\,(ii) we have to show that every \nfs\ of
$\zb{\wt \ve}i\N$ is linearly independent. Let $B$ be a \nfs\ of~$\N$, and let
$\la:B\to K$ be \st $\suml_{k\in B}\la_k\wt\ve_k=\wt0$. Then
$\suml_{k\in B}\la_kx^k\nde3.33 = \suml_{k\in B}\la_k
\Phi_x(\ve_k)\nad{\er{3.62},\er{3.63}} = \suml_{k\in B}\la_k\wt\ve_k(x) = \bigl(
\suml_{k\in B}\la_k\wt\ve_k\bigr)(x)=0$ \fa $x\in K$. \E\Ip $\suml_{k\in B}\la_k
x^k=0$ \fa $x\in K^\t:=K\sms0$. Let $\si_k:K\sms0\to K$, $k\in B$, be defined
by $\si_k(x):=x^k$, $x\in K^\t$, $k\in B$. Note that $(K^\t,\cdot,1)$ is a
group, and $x^k=k\ddt x$, the $k$-fold \IT\ of~$x$ in the group $(K^\t,\cdot,1)$
by \df. By \era2{2.3}\,I5 we have $\si_k(x\cdot y)=k\ddt(x\cdot y)= (k\ddt x)
\cdot(k\ddt y)= \si_k(x)\cdot \si_k(y)$, $x,y\in K^\t$, $k\in B$. \Mo $\si_k(1)
=k\ddt1=1$ by \era2{2.3}\,I0, $k\in B$. \If that $\si_k$, $k\in B$, is a \ti{\ch}
of $(K^\t,\cdot,1)$ in $(K,+,\cdot,0,1)$. We now show the map $k\mt \si_k$
from~$B$ into $\{\si_k\in K^{(K^\t)}: k\in B\}$ is in\jc. Let $k,l\in B$,
$k\ne l$. Then $2^k\ne 2^l$ by \era2{2.33}. Since the field~$K$ has \ch istic
zero, its \pf\ denoted by~$\check K$ is infinite. \E\Tf \te s an in\jc\ ring-\is
sm $j:\Q\to\check K$ by Theorem \rfa4{t5.41}. Hence $j:\Q\to K$ is an in\jc\
ring-\hm sm. \E\Ip $j:(\Q,\cdot,1)\to (K,\cdot,1)$ is an in\jc\ monoid-\hm sm.
By \era2{1.45}, we have $\si_k(j(z))= k\ddt j(z) = j(k\ddt z)=
j(z^k)$. Similarly $\si_l(j(z))=j(z^l)$. Since $z^k\ne z^l$ and $j$ is in\jc,
$\si_k(j(z))\ne\si_l(j(z))$. \E\Tf $\si_k\ne\si_l$, and the map $k\mt \si_k$
is in\jc. By Lemma \rf{l3.39}, $\la_k=0$ \fa $k\in B$. Hence $\zb{\wt\ve}kB$
is linearly independent in $(K,P(K))$. This completes the proof of~(i).

(ii) ``\ti{Sur\ji\ of $\Phi$}'': If $f=\wt0$, then $f=\Phi({\bf0})$ by \er{3.62},
\er{3.64}. If $f\ne\wt0$, then \te s $b\in K[X]\sms{\bf0}$ \st $f(x)=\Phi_x
(b)$, $x\in K$. Hence $f=\Phi(b)$ by \er{3.64}.

``\ti{Linearity of $\Phi$}'': This is a direct con\sq\ of \er{3.30}, \er{3.31}.

``\ti{In\ji\ of $\Phi$}'': Since $\Phi$ is a monoid-\hm sm from $(K[X],+,{\bf0})$
into\break $(P(K),+,\wt0)$, it is \sft\ in view of \E\Pr\ \rfa4{p3.10}\,(vi) to show
that $\{b\in K[X]: \Phi(b)=\wt0\}={\bf0}$. We have $\Phi({\bf0})=\wt0$ by \er{3.62}.
Suppose, for \cd ion, that \te s $b\in K[X]\sms{\bf0}$ \st $\Phi(b)=\wt0$.
Then $b=\suml_{k=0}^{\deg(b)}b_k\ve_k$ and $\Phi_x(b)=\suml_{k=0}^{\deg(b)}
b_kx^k=0$ \fa $x\in K$. Thus $\suml_{k=0}^{\deg(b)}b_k\wt\ve_k=\wt0$. Since
$\{\wt\ve_k\in P(K): k\in[0,N]\}$ is linearly independent in $(K,P(K))$, we
have $b_k=0$, $k\in[0,\deg(b)]$. Thus $b={\bf0}$, a~\cd ion. This completes
the proof of~(ii), hence also of \E\Pr\ \rf{p3.37}.
\endproof

Now we equip $P(K)$ not only with the (pointwise) \ad~$+$ introduced in
Examples \rf{xa1.2}\,(i), but also with the (pointwise) \mlc\ introduced in
the same example.

\bth3.40
Let $(K,+,\cdot,0,1)$ be a field of \ch istic zero, let $(K,P(K))$ be the \vs\ over~$K$ of
\pl\ maps over~$K$ introduced in \E\Pr\ \rf{p3.37} and let $\cdot$~be the
pointwise \mlc\ on~$K^K$ introduced in Examples \rf{xa1.2}\,{\rm(i)}. Then

\hph i,ii, $(K^K,\cdot,\wt1)$ is an \am\ and $P(K)$ is a \sbm\ of $(K^K,\cdot,
\wt1)$, where $\wt1$ is defined in \er{3.62}.

\hph ii,i, $(P(K),+,\cdot,\wt0,\wt1)$ is a $($\cmt e\/$)$ ring with unity.

\hph iii,, The map $\Phi$ defined in \er{3.64} is a ring-\is sm.

\hph iv,, $(P(K),+,\cdot,\wt0,\wt1)$ is a domain, and the field of \qt s of
$P(K)$ is ring-\is c to the field of \qt s of~$K[X]$.
\eth

\proof \

(i) Let $f,g,h\in K^K$ and $x\in K$. Then $((f\cdot g)\cdot h)(x)= (f\cdot g)
(x)\cdot h(x) = (f(x)\cdot g(x))\cdot h(x) \nda21.5 = f(x)\cdot (g(x)\cdot h(x))
= f(x)\cdot ((g\cdot h)(x)) = (f\cdot(g\cdot h))(x)$. Hence $\cdot$ in $P(K)$ is \asc e.
Similarly $(f\cdot g)(x)=f(x)\cdot g(x) \nda21.6 = g(x)\cdot f(x)=(g\cdot f)
(x)$. Hence $\cdot$ in $P(K)$ is \cmt e. $(\wt1\cdot f)(x)=1\cdot f(x)\nda21.7 = f(x)$.
Hence $(K^K,\cdot,\wt1)$ is an \am.

$\wt1(x)=1=\Phi_x(\ve_0)$, $x\in K$, hence $\wt 1\in P(K)$. Let $f,g\in P(K)$.
By \df\ \te\ $b,c\in K[X]$ \st $f(x)=\Phi_x(b)$, $g(x)=\Phi_x(c)$. Hence $(f\cdot g)
(x)=f(x)\cdot g(x)=\Phi_x(b)\cdot \Phi_x(c)\nde3.32 = \Phi_x(b\cdot c)$, $x\in K$.
Thus $f\cdot g\in P(K)$.

(ii) By \E\Pr\ \rf{p3.37} we know that $(P(K),+,\wt0)$ is an \ag. We just
showed that $(P(K),\cdot,\wt1)$ is an \am. It remains to prove \era4{1.1n},
\era4{1.2n}. Let $f,g,h \in P(K)$. Then $(\wt0\cdot f)(x)=\wt0(x)\cdot f(x)
=0\cdot f(x)
\nda41.1n = 0$, $x\in K$, hence \era4{1.1n} holds for $P(K)$. Let $x\in K$,
then $(f\cdot(g+h))(x) = f(x)\cdot((g+h)(x)) = f(x)\cdot (g(x)+h(x))\nda41.2n
= f(x)\cdot g(x)+f(x)\cdot h(x)= (f\cdot g)(x)+(f\cdot h)(x) = ((f\cdot g)
+(f\cdot h))(x)$. Hence \era4{1.2n} holds for $P(K)$.

(iii) In view of \E\Pr\ \rf{p3.37} and Lemma \rfa2{l1.8} it is \sft\ to show
that\break $\Phi:(K[X],\cdot,\ve_0) \to (P(K),\cdot,\wt1)$ is a monoid-\hm sm. We have
$(\Phi(\ve_0))(x)\nde3.64 {:=} \Phi_x(\ve_0)=x^0=1$ \fa $x\in K$, since $x^0=0\ddt x
=1$ in $(K,\cdot,1)$ by \era2{2.3}\,I0. Hence $\Phi(\ve_0)=\wt1$ by \er{3.62}.
Let $b,c\in K[X]$, $x\in K$. Then $(\Phi(b\cdot c))(x)\nde3.64 = \Phi_x(b\cdot c)
\nde3.32 = \Phi_x(b)\cdot \Phi_x(c) \nde3.64 = (\Phi(b))(x)\cdot(\Phi(c))(x) =
(\Phi(b)\cdot \Phi(c))(x)$. Hence $\Phi(b\cdot c)=\Phi(b)\cdot \Phi(c)$.

(iv) Let $f,g\in P(K)\sms{\wt0}$. Then \te\ $b,c\in K[X]\sms{\bf0}$ \st
$f=\Phi(b)$, $g=\Phi(c)$. Thus $f\cdot g=\Phi(b)\cdot \Phi(c)=\Phi(b\cdot c)$. Since
$K[X]$ is a domain, we have $b\cdot c\in K[X]\sms{\bf0}$. Hence $f\cdot g =
\Phi(b\cdot c)\ne\bf0$ since $\Phi({\bf0})=\wt0$ and $\Phi$~is in\jc. \E\Tf $P(K)$
is a domain. The second conclusion of (iv) is a direct con\sq\ of (iii) and
Theorem \rfa4{t5.35}\,(v).
\endproof

\bex3.41
Let $f\in P(K)$ and $b\in K[X]$ be the unique \el\ of~$K[X]$ \sf ying $f=\Phi
(b)$. Then $f$~is said to be \ti{ir\rd} in $P(K)$ if $b$~is ir\rd\ in $K[X]$.
Give a necessary and \sft\ \cn\ on~$f$ in terms of the ring $(P(K),+,\cdot,
\wt0,\wt1)$ for~$f$ to be ir\rd\ in $P(K)$.
\eex

\brm3.42
It can be shown (e.g.~\cite[pp.\ 95, 96]{Alg} Eisenstein's Theorem, or \cite[pp.\ 245, 246]{Is}) that $f\in P(\Q)$ defined by $f(x):=x^n-p$,
$p$~\Pn\ and $n\in\Na$, is ir\rd\ in~$P(\Q)$. Show that \fa $n\ge2$ \te s a
field \ext\ $(F_n,j_n)$ of the field~$\Q$ \st the \dm\ of the \vs\ $(j_n(\Q),
F_n)$ is equal to~$n$.
\erm

\newpage
\Subsubsection{Galois fields}\label{sss.Galois}

The aim of this section is to prove the \ex\ and \uq\ (up to \is sms) of a
finite field of order~$p^n$, $p$~prime and $n\in\Na$. We first observe that
in case $n=1$ the \fw\ holds. If $p$~is prime and $K$~denotes the field
$(\N_p,+_p,\cdot_p,0,1)$, then the \pl\ map
$f:K\to K$ defined by $f(x):=x^p-x$ where $x^p$ denotes the $p$-th \IT\ of~$x$
in the monoid $(\N_p,\cdot_p,1)$ is identically equal to zero. Indeed, since
$(\N_p\sms0,\cdot_p,1)$ is a finite group of order $p-1$ (see the introduction
of Section \ref{sss.pr.fld}), we have $k^{p-1}=1$
\fe $k\in[1,p-1]$ by \E\Pr\ \rfa4{p3.65}\,(i). \E\ml ying both sides of the
\et y by~$k$ we obtain $k^p=k$ or $k^p-k=0$ (i.e.\ $k^p+(-k)=0$ where $-k$ is
the inverse of~$k$ in $(\N_p,+_p,0)$). Since $p>1$ we have $0^p=0$ ($0^p:=
p\ddtt p0 = ((p-1)+1)\ddtt p 0\nad{\era2{2.3}\,\rm I2}= ((p-1)\ddtt p 0)\cdot_p
(1\cdot_p0) \nad{\era2{2.3}\,\rm I1}= ((p-1)\ddtt p 0)\cdot_p0\nda41.1n = 0$
where $p\ddtt p0$ is the $p$-th \IT\ of~$0$ in the monoid $(\N_p,\cdot_p,1)$),
hence $f(k)=0$ \fa $p\in\N_p$. In view of \E\df\ \rf{d3.15} every \el\ of~$K$ is
a~\ti{root\/} of the \fp\ $\ve_p-\ve_1\in K[X]$, hence by
Corollary \rf{c3.14} $\ve_1-k\ve_0$ \ti{divides} $\ve_p-\ve_1$ \fe
$k\in[0,p-1]$. We claim
\beq4.1
\ve_p-\ve_1 = \prod_{k\in\zo0,p }(\ve_1-k\ve_0),
\e
where $\prod$ denotes the \cme product associated with the \mlv\ monoid of the
ring of \fp s $K[X]$.

\E\et y \er{4.1} is a con\sq\ of the \fw\ lemma.

\blm4.1
Let $F$ be a field, let $b\in F[X]\sms0$ be \st $\deg(b)\ge1$ and
$b(\deg(b))=1$. Let $n:=\deg(b)$
and let $a:[1,n]\to F$ be \emph{in\jc}. If $a(i)$, $i\in[1,n]$, are roots of~$b$,
then
\beq4.2
b=\prod_{i\in[1,n]} (\ve_1-a(i)\ve_0).
\e
\elm

\proof
We use \In\ on $n\in\Na$. Set $M:=\{m\in\Na: \text{Lemma \rf{l4.1} holds for }
n:=m\}$.

$1\in M$: By \as\ $b=\ve_1+\b\ve_0$ for some $\b\in F$ and $a(1)+\b=0$ since
$a(1)$ is a root of~$b$. Thus $\b=-a(1)$ and $b=\ve_1+((-a(1))\ve_0)$. Recall
$-a(1)=-(1a(1))\nda44.12 = (-1)a(1)$, hence $\ve_1+((-a(1))\ve_0) = \ve_1+
((-1)a(1))\ve_0 \nad{\er{1.1}\,\rm I3}= \ve_1 +((-1)(a(1)\ve_0))
\nde1.9 = \ve_1 +(-(a(1)\ve_0))
\nde1.8 = \ve_1-(a(1)\ve_0)$. Hence $b=\ve_1-(a(1)\ve_0)$ and $1\in M$.

\ti{$m\in M$ implies $m+1\in M$}: We suppose $m\in M$, let $b\in F[X]\sms0$
with $\deg(b)=m+1$, $b(m+1)=1$, let $a:[1,m+1]\to F$
be in\jc\ and let $a(i)\in F$, $i\in[1,m]$, be roots of~$b$.
By Corollary \rf{c3.14} \te s $q\in F[X]$ \st $b=q\cdot (\ve_1-
a(m+1)\ve_0)$. Clearly $q\ne0$ since $b\ne0$. \Mo $m+1=\deg(b)\nde2.29 =
\deg(q)+1$, hence $\deg(q)=m$. By \er{2.28}, we have $1=b(\deg(b))= q(\deg(q))
\cdot1$, hence $q(\deg(q))=1$. We claim that $a(i)$, $i\in[1,m]$, are roots
of~$q$. We have $0=\Phi_{a_i}(b) \nde3.32 = \Phi_{a_i}(q)\Phi_{a_i}(\ve_1-
a(m+1)\ve_0)\nde3.28 = \Phi_{a_i}(q)(a(i)-a(m+1))$ \fa $i\in[1,m]$. By \as\
$a(i)-a(m+1)\ne0$ \fa $i\in[1,m]$, hence $\Phi_{a_i}(q)=0$ \fa $i\in[1,m]$,
since $F$~is a field. The claim is proved. \If from Lemma \rf{l4.1} with $n:=m$ that
$q=\prodl_{i\in[1,m]}(\ve_1-a(i)\ve_0)$. Thus, in view of the \df\ of~$\prod$,
$m+1\in M$. \E\Tf $M=\Na$.
\endproof

In the next \Pr\ we \ti{suppose} that $F$ is a finite field of order $p^n$,
$n\in\Na$, and show that a \gn\ of \er{4.1} holds.

\bpr4.2
Let $F$ be a finite field of order $p^n$, $p$ prime and $n\in\Na$. Then the
\fw\ \fc\ of the \fp\ $\ve_{p^n}-\ve_1 \in F[X]$ holds{\rm:}
\beq4.3
\ve_{p^n}-\ve_1 = \prod_{\a\in F} (\ve_1-\a\ve_0).
\e
\epr

\proof
The \mlv\ group of the field~$F$ has order $p^n-1$, hence by \E\Pr\
\rfa4{p3.65}\,(i) every \el~$\a$ of $F\sms0$ \sf ies $\a^{p^n-1}=1$ where
$\a^{p^n-1}$ is the $(p^n-1)$-th \IT\ of~$\a$ in the monoid $(F,\cdot,1)$.
Thus we also have as above $\a^{p^n}=\a$ \fa $\a\in F\sms0$. Since $p^n\ge1$,
we have $0^{p^n}=0$, and every \el\ of~$F$ is a root of $\ve_{p^n}-\ve_1$. We
conclude using Lemma \rf{l4.1}.
\endproof

\bdf4.3
Let $F$ be a field and let $b\in F[X]\sms0$ be a nonzero \fp\ of degree
$n\in\Na$. Then the \fp~$b$ is said to \ti{split\/} over~$F$ if \te s\index{split over}
$a:[0,n]\to F$ not necessarily in\jc\ with $a(0)\ne0$ \st
\beq4.4
b=a(0)\prod_{i\in[1,n]}(\ve_1-a(i)\ve_0) \hbox{ if } n\in\Na.
\e
\edf

\bex4.4
Show that the $a_i$'s, $i\in[1,n]$, in \er{4.4} are roots of~$b$ in~$F$.
\eex

\brm4.5
\If from \E\Pr\ \rf{p4.2} that if $F$ is a finite field of order~$p^n$, then
the \fp\ $\ve_{p^n}-\ve_1$ splits over~$F$.
\erm

The next \Pr\ reduces the problem of proving the \ex\ of a~field of order $p^n$
to the problem of finding a field \ext\ $(F,j)$ of the field $\N_p$ \st
$\wh\jmath(\ve_{p^n}-\ve_1)$ splits over~$F$.

\bpr4.6
Let $p$ be a \Pn\ and let $n\in\Na$. Suppose \te s a field \ext\ $(F,j)$ of
the field $(\N_p,+_p,\cdot_p,0,1)$ \st $\wh\jmath(\ve_{p^n}-\ve_1)$ splits
over~$F$ where $\wh\jmath$ is defined in \er{3.40}. Then $\{\a\in F: \a^{p^n}=
\a\}$ is a \emph{subfield} of~$F$ of order~$p^n$.
\epr

The \ex\ of such \ext\ will be guaranteed by a \gn\ of Kronecker's theorem
(see Theorem \rf{t4.20}).
For the proof of \E\Pr\ \rf{p4.6} we need several lemmata.

\blm4.7
Let $F$ be a field of \ch istic $p$ and let $n\in\Na$. Then the set $\{x\in F:
a^{p^n}=\a\}$ is a subfield of~$F$.
\elm

\proof
Set
\beq4.5
G:=\{\a\in F: \a^{p^n}=\a\}.
\e
We verify (a), (b), (c) of \E\df\ \rfa4{d4.23}.

(a) Since $p^n\ge1$, $0\in G$. Let $\a,\b\in G$. Then $(\a+\b)^{p^n}
\nda44.114 = \a^{p^n}+\b^{p^n} = \a+\b$. Hence $\a+\b\in G$. Thus $G$~is
a~\sbm\ of $(F,+,0)$ and $(G,+,0)$ is an (abelian) monoid. Let $\a\in G$.
We show that $-\a\in G$, where $-\a$ is the inverse of~$\a$ in $(F,+,0)$.
Indeed, $(-\a)^{p^n}=(0-\a)^{p^n} \nda44.115 = (0^{p^n})-(\a^{p^n}) = (0)-
(\a^{p^n}) = (0)+(-(\a^{p^n}))= -(\a^{p^n})=-\a$. \If that $(G,+,0)$ is
a subgroup of $(F,+,0)$.

(b) We have $1\in G$ by \era2{2.3}\,I4. Let $\a,\b\in G$. Then $(\a\b)
^{p^n} \nad{\era2{2.3}\,\rm I5}= \a^{p^n}\b^{p^n}=\a\b$. Hence $G$ is
a \sbm\ of~$(F,\cdot,1)$.

(c) Let $\a\in G\sms0$, and let $\b\in F\sms0$ be the inverse of~$\a$
in $(F,\cdot,1)$, i.e.\
$\a\b=1$. As above $(\a\b)^{p^n}=\a^{p^n}\b^{p^n}$. Hence $\a\b^{p^n} =
\a^{p^n}\b^{p^n} = (\a\b)^{p^n}=1^{p^n}\nad{\era2{2.3}\,\rm I4}=1$.
By the \uq\ of the inverse of~$\a$
in $(F\sms0,\cdot,1)$, $\b^{p^n}$ is the inverse of~$\a$, hence $\b^{p^n}=\b$.
\If that $\b\in G$, hence $\b$ is the inverse of~$\a$ in~$G$.
\endproof

We now investigate \pp ies of \fp s of the form \er{4.4}.

\blm4.8
Let $F$ be a field, let $n\in\Na$ and let $a:[0,n]\to F$ with $a(0)\ne0$. Then
the \fp~$b$ defined by the \RHS\ of \er{4.4} \sf ies{\rm:}

\hph i,ii, $b\in F[X]\sms0$,

\hph ii,i, $\deg(b)=n$,

\hph iii,, $b(n)=a(0)$,

\hph iv,,  $\a\in F$ is a root of $b$ in $F$ iff $\a=a(i)$ \fs $i\in[1,n]$.
\elm

\bex4.9
Prove Lemma \rf{l4.8}.
\eex

\blm4.10
Let $F$ be a field and let $b\in F[X]$ be as in Lemma \rf{l4.8}. Then
\beq4.6
a(i)\ne a(j) \hbox{ \fa} i,j\in [1,n],\ i\ne j,
\e
iff
\beq4.7
\sum_{k=1}^{\deg(b)} (k\dpl b(k)) a(i)^{k-1}\ne0 \qh{\fa} i\in[1,n],
\e
where $k\dpl b(k)$ is the $k$-th \IT\ of $b(k)$ in $(F,+,{\bf0})$.
\elm

Before giving the proof of Lemma \rf{l4.10}, we show how \E\Pr\ \rf{p4.6}
follows from Lemma \rf{l4.7} and Lemma \rf{l4.10}.

\proof[Proof of \E\Pr\ \rf{p4.6}]
For convenience set $K:=(\N_p,+_p,\cdot_p,0,1)$ and $b\in K[X]$ defined by
$b:=\ve^{p^n}-\ve_1$. By \as\ \te\ a field $(F,+,\cdot,{\bf0},{\bf1})$ and
an in\jc\ ring-\hm sm $j:K\to F$ \st $\wh\jmath(b)$ splits over~$F$. Observe
$j(0)={\bf0}$ and $j(1)={\bf1}$. We define $\eta_i\in F[X]$, $i\in \N$, by
setting
\beq4.8
\eta_i(j):=\bca
{\bf1} & \hbox{for }j:=i,\\
{\bf0} & \hbox{otherwise.}
\eca
\e
From \er{3.40} we infer
\beq4.9
\eta_i=\wh\jmath(\ve_i), \q i\in\N.
\e
\E\Tf $\wh\jmath(b)(l)\nde3.40 = j((\ve_{p^n}-\ve_1)(l)) = j(\ve_{p^n}(l) -
\ve_1(l)) = j(\ve_{p^n}(l))-j(\ve_1(l)) =\break
(\wh\jmath(\ve_{p^n}))(l) - (\wh\jmath
(\ve_1))(l) \nde4.9 = \eta_{p^n}(l)-\eta_1(l) = (\eta_{p^n}-\eta_1)(l)$ \fa
$l\in\N$. Hence
\beq4.10
\wh\jmath(b) = \eta_{p^n}-\eta_1.
\e
Note that $\deg(\wh\jmath(b))=p^n$. By \as\ \te s $a:[0,p^n]\to F$ with $a(0)
\ne0$ \st
\beq4.11
\wh\jmath(b)=a(0)\prod_{i=1}^{p^n} (\eta_1-a(i)\eta_0).
\e
By Lemma \rf{l4.8}\,(i)(ii)(iii) we have $a(0)={\bf1}$. We now apply Lemma
\rf{l4.10} for proving that the map $a:[1,p^n]\to F$ is \ti{in\jc}. It
suffices to verify \cn\ \er{4.7}, i.e.
\beq4.12
\sum_{k=1}^{p^n}(k\dpl (\wh\jmath(b)(k))) a(i)^{k-1}\ne{\bf0}, \q i\in[1,p^n].
\e
Since $p^n\ne1$, we have $\wh\jmath(b)(p^n)={\bf1}$, moreover, $(\wh\jmath(b))
(1)=-{\bf1}$ and $(\wh\jmath(b))(k)={\bf0}$ otherwise. Hence, using
\era2{2.3}\,I4, \er{4.7} reduces to
\beq4.13
(p^n\dpl{\bf1}) a(i)^{p^n-1} +1\dpl(-{\bf1})a(i)^0 \ne0, \q i\in[1,p^n].
\e
Note that $p^n\dpl{\bf1} = p^n\dpl j(1)\nda21.45 = j(p^n\dpl 1)$, where
$p^n\dpl1$ is the $p^n$-th \IT\ of~$1$ in $(\N_p,+_p,0)$. Hence $p^n\dpl1
\nda41.11 = \F_p(p^n)=0$ since $p$ divides~$p^n$. Thus $p^n\dpl{\bf1} =j(0)
={\bf0}$, and \er{4.13} becomes $1\dpl(-{\bf1})a(i)^0\ne0$, $i\in[1,p^n]$.
By \era2{2.3}\,I1 in $(F,+,{\bf0})$, $1\dpl(-{\bf1})a(i)^0 = (-{\bf1})a(i)^0$.
By \era2{2.3}\,I0 in $(F,\cdot,{\bf1})$, $a(i)^0={\bf1}$, hence $(-{\bf1})
a(i)^0=-{\bf1}\ne{\bf0}$. \If that \cn\ \er{4.7} is \sf ied, thus $\wh\jmath
(b)$ has $p^n$ distinct roots $a(i)$, $i\in[1,p^n]$, in~$F$. In view of
\E\Pr\ \rf{p3.26}\,(iii) $\wh\jmath(b)$ has \ti{exactly} $p^n$~roots in~$F$.
\If that $G$~defined in \er{4.5} is the set of roots of~$\jh(b)$ in~$F$ since
$\a\in G$ iff $\a=a(i)$ \fs $i\in[1,p^n]$. \E\Tf\ $\#(G)=p^n$.
\endproof

We now turn to the proof of Lemma \rf{l4.10}. To this end we first introduce
the important notion of \ti{\dv} of a \fp.\index{derivative of a formal polynomial}

\bdf4.11
Let $K$ be a field and $K[X]$ be the $K$-algebra of \fp s over~$K$. The
\ti{\dv} of an \el\ $b\in K[X]$, denoted by $D(b)$ (or~$b'$), is defined by
\beq4.14
D(b)(l):=(l+1)\dpl b(l+1), \q l\in\N,
\e
where $(l+1)\dpl b(l+1)$ is the $(l+1)$-th \IT\ of $b(l+1)$ in the monoid
$(K,+,0)$. Since $k\dpl0 = 0$, $k\in\N$, by \era2{2.3}\,I4, we have $\supp
(D(b))\sbs \supp(b)$, hence $D(b)\in K[X]$. The map $b\mt D(b)$ from $K[X]$
into $K[X]$ will be denoted by~$D$.
\edf

\brm4.12
Formula \er{4.14} is also used to define the \dv\ of a formal power series
(see for example \cite{Analysis}).
\erm

From \er{4.14} we infer:
\beq4.15
D(\a\ve_0) = {\bf0} \qh{\fa}\a\in K,
\e
where $\bf0$ denotes the \nel\ of $K[X]$.
\bga4.16
D(\ve_i) = i\dpl \ve_{i-1} \qh{\fa} i\in\Na, \\
D(b+c) = D(b)+D(c), \q b,c\in K[X]. \lb{4.17}
\e
Indeed, $D(b+c)(l)=(l+1)\dpl
(b+c)(l+1) = (l+1)\dpl(b(l+1)+c(l+1)) \nad{\era2{2.3}\,\rm I5}= (l+1)\dpl
b(l+1) + (l+1)\dpl c(l+1) = D(b)(l)+D(c)(l)= (D(b)+D(c))(l)$, $l\in\N$.
\beq4.18
D(\la b)=\la D(b), \q \la\in K,\ b\in K[X].
\e
Indeed, $D(\la b)(l) = (l+1)\dpl(\la b)(l+1) = (l+1)\dpl (\la b(l+1))
\nad*= \la(l+1)\dpl b(l+1) = \la(D(b)(l))=(\la D(b))(l)$, $l\in\Na$. In
$\nad*=$ we used the \fw\ lemma.

\blm4.13 \

\hph i,i, Let $(X,+,\cdot,0,1)$ be a \sr\ and let $a,b\in X$. Then
\beq4.19
n\dpl (a\cdot b)=(n\dpl a)\cdot b \qh{\fa} n\in\N.
\e

\hph ii,, Let $F$ be a field, let $(F,V)$ be a \vs\ over~$F$ and
let $\a\in F$, $a\in V$. Then
\beq4.20
n\dpl (\a a) = (n\dpl \a)a = \a(n\dpl a) \qh{\fa} \a\in F,\ a\in V.
\e
\elm

\proof \

(i) Use Lemma \rfa2{l1.22} with $(M_i,\qu_i,e_i):= (X,+,0)$, $i=1,2$, and
$\vf(a) :=a\cdot b$, $a,b\in X$.

(ii) Use Lemma \rfa2{l1.22}, firstly with $(M_i,\qu_i,e_i):=(F,+,0)$,
$i=1,2$, and $\vf(\la):=\la a$, $\la\in F$, $a\in V$, secondly, with
$(M_i,\qu_i,e_i):= (V,+,0)$, $i=1,2$, and $\vf(a):=\la a$, $\la \in F$, $a\in
V$.
\endproof

In view of \er{4.17}, \er{4.18} the map $D:K[X]\to K[X]$ is linear. Hence if
$b\ne0$ with $\deg(b)\ge1$, we have $b=\suml_{i=0}^{\deg(b)} b(i)\ve_i$ and by
\er{4.13}, \er{4.3}
\bmlg
D(b)=\suml_{i=0}^{\deg(b)}b(i)D(\ve_i) = b(0)D(\ve_0) +
\suml_{i=1}^{\deg(b)} b(i)D(\ve_i) \nde4.15 = \suml_{i=1}^{\deg(b)}b(i)D(\ve_i)\\
{}\nde4.16 = \suml_{i=1}^{\deg(b)}b(i)(i\dpl \ve_{i-1}) \nde4.20 =
\suml_{i=1}^{\deg(b)} (i\dpl b(i)\ve_{i-1}).
\e
\E\Tf we proved
\beq4.21
D(b) = \sum_{i=1}^{\deg(b)}(i\dpl b(i))\ve_{i-1} \qh{\fa $b$ with }\deg(b)\ge1.
\e

\If from formula \er{4.21} that if $\deg(b)\ge1$, then the degree of $D(b)$
is at most equal to $\deg(b)-1$ whenever $D(b)\ne{\bf0}$. If the field~$K$
has \ch istic $p$~(prime), then it turns out that $D(\ve_p)={\bf0}$. Indeed,
$K:=(K,+,\cdot, 0,1)$ and $p\dpl1=0$, then $D(\ve_p)=p\dpl \ve_{p-1}=
p\dpl(1\cdot \ve_{p-1}) \nde4.19 = (p\dpl 1)\cdot \ve_{p-1} = 0\cdot\ve_{p-1}
=0$. Thus in such a field the \dv\ of a \fp\ of degree larger than zero may
be equal to~zero. Similarly, $D(\ve_p+\ve_1) = D(\ve_p)+D(\ve_1) = D(\ve_1)
=\ve_0$. In this case $\deg(D(\ve_p+\ve_1))=0 < p-1 =\deg(\ve_p+\ve_1)-1$.

\bex4.14
Let $(K,+,\cdot,0,1)$ be a field \st $n\dpl1\ne0$ \fa $n\in\Na$. Show that if
$b\in K[X]$ and $D(b)={\bf0}$, then \te s $\a\in K$ \st $b=\a\ve_0$. Show
that if $b\in K[X]\sms{\bf0}$ and $\deg(b)\ge1$, then $\deg(D(b)) =
\deg(b)-1$.
\eex

\brm4.15
Recall that a field $(K,+,\cdot,0,1)$ \st $n\dpl 1\ne0$ \fa $n\in\Na$ is said to be of
\ti{\ch istic zero} (see \E\df\ \rf{d4.38}). Most of the fields used in Analysis are of \ch istic
zero.
\erm

We now investigate the \dv\ of a \ti{product\/} of \fp s. We first compute the \dv\
of $\ve_i\cdot\ve_j$, $i,j\in\Na$. Using \er{2.26} we obtain $D(\e_i\cdot\e_j) =
D(\ve_{i+j})\nde4.16 = (i+j)\dpl \ve_{i+j-1} \nad{\era2{2.3}\,\rm I2} = i\dpl
\ve_{i+j-1} + j\dpl \ve_{i+j-1} \nde2.26 = i\dpl(\ve_{i-1}\cdot\ve_j)
+j\dpl (\ve_{j-1}\cdot\ve_i) \nde4.19 = \break (i\dpl\ve_{i-1})\cdot\ve_j + (j\dpl
\ve_{j-1})\cdot\ve_i \nde4.16 = D(\ve_i)\cdot\ve_j + D(\ve_j)\cdot\ve_i =
D(\ve_i)\cdot \ve_j + \ve_i\cdot D(\ve_j)$. If $i=0$, $j\in\N$, $D(\ve_0\cdot
\ve_j)\nde2.26 = D(\ve_j) = {\bf0}\cdot\ve_j + \ve_0\cdot D(\ve_j) =
D(\ve_0)\cdot \ve_j + \ve_0\cdot D(\ve_j)$. Similarly, if $i\in\N$ and $j=0$,
$D(\ve_i\cdot\ve_j)=D(\ve_i)\cdot \ve_j + \ve_i\cdot D(\ve_j)$. Thus we obtain
\beq4.22
D(\ve_i\cdot\ve_j)=D(\ve_i)\cdot \ve_j + \ve_i\cdot D(\ve_j), \q i,j\in\N.
\e
Observe that the maps $B_i:K[X]\t K[X] \to K[X]$, $i=1,2$, defined by $B_1(u,v)
:=D(u\cdot v)$, and $B_2(u,v):=D(u)\cdot v+u\cdot D(v)$, $u,v\in K[X]$, are \sy ic,
i.e.\ $B_i(u,v)=B_i(v,u)$, $u,v\in K[X]$, $i=1,2$. Indeed, $D(u\cdot v) =
D(v\cdot u)$, and $D(u)\cdot v+u\cdot D(v) = u\cdot D(v)+D(u)\cdot v= D(v)\cdot
u + v\cdot D(u)$, $u,v\in K[X]$. \Mo $B_i$, $i=1,2$, are \ti{bilinear} (see
\E\df\ \rf{d2.17}). In view of the \sy y it is sufficient to show that the maps
$u\mt B_i(u,v)$, $v\in K[X]$, $i=1,2$, are linear. But $u\mt B_1(u,v)$ is the
\cm\ of the linear map $u\mt u\cdot v$ and the linear map~$D$, which is also
linear by \er{4.17}, \er{4.18}. Hence $u\mt B_1(u,v)$ is linear by Lemma
\rf{l1.5}\,(ii). Let $u,w\in K[X]$. Then $D(u+w)\cdot v + (u+w)\cdot D(v) =
(D(u)+D(w))\cdot v + (u+w)\cdot D(v) =
D(u)\cdot v + (D(w)\cdot v+u\cdot D(v)) + w\cdot D(v) = D(u)\cdot v +(u\cdot
D(v)+D(w)\cdot v) + w\cdot D(v) =  B_2(u,v)+B_2(w,v)$. Similarly, let $u\in K[X]$ and $\la\in K$,
then $D(\la u)\cdot v +(\la u)\cdot v = (\la D(u))\cdot v + \la(u\cdot v) =
\la(D(u)\cdot v+(u\cdot v))$. Hence $u\mt B_2(u,v)$ is linear. Since $B_1
(\ve_i,\ve_j)=B_2(\ve_i,\ve_j)$, $i,j\in\N$, by \er{4.22}, we infer that
$B_1(u,v)=B_2(u,v)$ \fa $u,v\in K[X]$ in view of \E\Pr\ \rf{p2.23}. Thus we
have shown
\beq4.23
D(u\cdot v)= D(u)\cdot v+u\cdot D(v) \qh{\fa} u,v\in K[X].
\e
Observe that $D(\ve_0)={\bf0}$ follows from \er{4.23}. Indeed, $D(\ve_0)
\nde2.26 = D(\ve_0\cdot \ve_0)\nde4.23 = D(\ve_0)\cdot \ve_0+ \ve_0\cdot
D(\ve_0) = D(\ve_0)+D(\ve_0)$. From ${\bf0}+D(\ve_0) = D(\ve_0)+D(\ve_0)$ we
obtain $D(\ve_0)={\bf0}$. Another con\sq\ of \er{4.23} and \er{2.26} is
$D(\ve_{i+1})=D(\ve_i)\cdot \ve_1 + \ve_i\cdot D(\ve_1) = D(\ve_i) + \ve_i
\cdot D(\ve_1)$, $i\in\N$.

\bex4.16
Let $\wh D:K[X]\to K[X]$ \sf y \er{1.2}, \er{1.3}, \er{4.22} and $\wh D(\ve_1)=\ve_0$.
Show that $\wh D=D$.
\eex

Summarizing, we have

\bpr4.17
Let $D:K[X]\to K[X]$ be the map introduced in \E\df\ \rf{d4.11}. Then $D$~is
linear and \sf ies \er{4.15}, \er{4.16}, \er{4.21} and \er{4.23}.
\epr

\bex4.18
Let $D$ be as in \E\Pr\ \rf{p4.17}.

\hph i,i, Let $u,v,w\in K[X]$, show that
\beq4.24
D(u\cdot v\cdot w) = D(u)\cdot v\cdot w + u\cdot D(v)\cdot w + u\cdot v\cdot D(w).
\e
Find a formula for $D\bigl(\prodl_{i=1}^n u_i\bigr)$ where $u_i\in K[X]$,
$i\in[1,n]$.

\hph ii,, Let $u,v\in K[X]$ and let $D^k$, $k\in\N$, denote the $k$-th \IT\
of~$D$ in the monoid $(K[X]^{K[X]},\circ,\id_{K[X]})$, where $\circ$ denotes
the \cm\ of two selfmaps of~$K[X]$. Prove the Leibniz formula:
\beq4.25
D^n(u\cdot v)= \sum_{k=0}^n \binom nk \ddt \bigl(D^k(u)\cdot D^{n-k}(v)),
\e
where $\binom nk \ddt a$ denotes the $\binom nk$-th \IT\ of~$a$ in $(K[X],
\cdot,\ve_0)$, and $\binom nk=\frac{n!}{k!(n-k)!}$.
\eex

\medskip
Recall that in Corollary \rf{c3.14} we gave a necessary and \sft\ \cn\ on a \fp\
$b\ne0$ for~$b$ to be divisible by $\ve_1-\a\ve_0$ \fs $\a\in K$. \Wanp \ch
ize \fp s $b\ne0$ that are divisible by $(\ve_1-\a\ve_0)^2$ \fs $\a\in K$.

\bpr4.19
Let $K$ be a field and let $b\in K[X]\sms{\bf0}$. Then \te s $\a\in K$ \st
$(\ve_1-\a\ve_0)^2$ divides~$b$ iff
\beq4.26
\Phi_\a(b)=0 \qh{and } \Phi_\a(D(b))=0.
\e
\epr

\proof \

\ti{Only if\/}: Suppose $b=q\cdot (\ve_1-\a\ve_0)^2$ \fs $q\in K[X]\sms{{\bf0}}$ and
$\a\in K$. By \er{3.32}, \er{3.34}, \era4{1.1n}, $b=(q\cdot (\ve_1-\a\ve_0))\cdot(\ve_1
-\a\ve_0)$ \sf ies $\Phi_\a(b)=0$. By \er{4.23} $D(b)=\break D(q\cdot
(\ve_1-\a\ve_0)^2) = D(q)\cdot(\ve_1-\a\ve_0)^2 + q\cdot
D((\ve_1-\a\ve_0)^2)$. Again by \er{4.23} $D((\ve_1-\a\ve_0)\cdot\break
(\ve_1-\a\ve_0)) = D(\ve_1-\a\ve_0)\cdot(\ve_1-\a\ve_0) + D(\ve_1-\a\ve_0)
(\ve_1-\a\ve_0)$. Hence $\Phi_\a(D(b))\nde3.30 = \Phi_\a(D(q)\cdot
(\ve_1-\a\ve_0)^2) +\Phi_\a(D(\ve_1-\a\ve_0)\cdot(\ve_1-\a\ve_0)) + \Phi_\a
(D(\ve_1-\a\ve_0) \cdot(\ve_1-\a\ve_0)) \nde3.32 = \Phi_\a(D(q)\cdot\break
(\ve_1-\a\ve_0)) \Phi_\a(\ve_1-\a\ve_0) + \Phi_\a(D(\ve_1-\a\ve_0))
\Phi_\a(\ve_1-\a\ve_0) +\Phi_\a(D(\ve_1-\a\ve_0))\Phi_\a(\ve_1-\a\ve_0)\break
\nad{\er{3.34},\era4{1.1n}} = 0+0+0 = 0$.

\ti{If\/}: Since $\Phi_a(b)=0$, \te s $q_0\in K[X]$ \st
$b=q_0\cdot(\ve_1-\a\ve_0)$ by Corollary \rf{c3.14}. We have $q_0\ne0$ since $b\ne0$. Then $D(b) \nde4.23 =
D(q_0) \cdot (\ve_1-\a\ve_0)+q_0\cdot D(\ve_1-\a\ve_0)$. Thus $\Phi_\a(D(b))
\nad{\er{3.30}\er{3.32}} = \Phi_\a(D(q_0))\Phi_\a(\ve_1-\a\ve_0) +
\Phi_\a(q_0) \Phi_\a(D(\ve_1-\a\ve_0))$. By \er{3.34} the first term is equal
to zero. By linearity $D(\ve_1-\a\ve_0)= D(\ve_1)-D(\a\ve_0)
\nad{\er{4.16}\er{4.15}} = \ve_0$. Hence $\Phi_\a(D(b)) = 0+\Phi_\a(q_0)
\cdot \Phi_\a(\ve_0)\nde3.33 = \Phi_\a(q_0)$. Since $\Phi_\a(D(b))=0$, we
have $\Phi_\a(q_0)=0$, hence \te s $q_1\in K[X]$ \st $q_0= q_1\cdot
(\ve_1-\a\ve_0)$. \If that $b=q_0\cdot(\ve_1-\a\ve_0)=q_1\cdot
(\ve_1-\a\ve_0)^2$.
\endproof

\proof[Proof of Lemma \rf{l4.10}]
Let $i\in[1,n]$. By Corollary \rf{c3.14}, $\Phi_{a(i)}(b)=0$ since $\ve_1-a(i)
\ve_0$ divides~$b$. Then $a(j)=a(i)$ \fs $j\in[1,n]\sms i$ iff $(\ve_1-a(i)
\ve_0)^2$ divides~$b$. Thus by \E\Pr\ \rf{p4.19} $a(j)=a(i)$ \fs $j\in[1,n]
\sms i$ iff $\Phi_{a(i)}(D(b))=0$. But \cn\ \er{4.7} says in view of \er{4.21}
that $\Phi_{a(j)}(D(b))\ne0$ \fa $j\in[1,n]$. Hence $a(i)\ne a(j)$ \fa
$i,j\in[1,n]$, $i\ne j$.
\endproof

We now return to the proof of the \ex\ of a finite field of order~$p^n$. In view of \E\Pr\
\rf{p4.6} and Lemma \rfa4{l4.24} it suffices to find a field \ext\ $(F,j)$ of~$\N_p$ \st $\wh\jmath
(\ve_{p^n}-\ve_1)$ splits over~$F$. The \ex\ of such a field \ext\ is guaranteed
by the \fw\ theorem.

\bth4.20
Let $K$ be a field and let $b$ be a nonzero \fp\ over~$K$ with $\deg(b)\ge1$.
Then \te s a field \ext\ $(F,j)$ \st $\wh\jmath(b)$ $($see \er{3.40}$)$ splits
over~$F$.
\eth

We now give a proof of Theorem \rf{t4.20} inspired by the nice proof of Theorem 17.8
given in~\cite{Is}.

\proof
Let $b\in K[X]\sms{\bf0}$ be of the form
\beq4.27
b=c  \qh{or } b=c\cdot \prod_{i=1}^r (\ve_1-a(i)\ve_0)
\e
where $c\in K[X]\sms{\bf0}$, $\deg(c)\ge1$, $r\in\Na$ and $a:[1,r] \to K$.

We use \In\ on $\deg(c)\in\Na$. Set
\beq4.28
M:=\{n\in\Na: \hbox{Theorem \rf{t4.20} holds for $\deg(c)=n$}\}.
\e
We have

$1\in M$:  If $\deg(c)=1$, then $c=c(1)(\ve_1-\a\ve_0)$ where $\a:=-c(1)\mo
c(0)$. Hence $b$ splits over~$K$ in both cases in \er{4.27}.

\ti{$n\in M$ implies $n+1\in M$}: We assume that $n\in M$ and suppose that
$\deg(c)=n+1$. By Lemma \rf{l4.21} since $\deg(c)\ge2$ \te\ $a_1\in K[X]$
\ti{ir\rd} and $a_2\in K[X]\sms{{\bf0}}$ \st $c=a_1\cdot a_2$. By Corollary \rf{c3.30}
\te\ an \ext\ field $(K_1,j_1)$ and $\b\in K_1$ \st $\Phi_\b(\wh\jmath_1 (a_1))
=0$ in~$K_1$. Since $\wh\jmath_1(c)=\wh\jmath_1(a_1\cdot a_2)\nde3.49 = \wh\jmath_1(a_1)
\cdot\wh\jmath_1(a_2)$ in~$K_1[X]$, $\Phi_\b(\wh\jmath_1(c))\nde3.32 = \Phi_\b(\wh\jmath_1
(a_1))\Phi_\b(\wh\jmath(a_2))=0$. Setting $\eta_i:= \wh\jmath_1(\ve_i)$,
$i\in\N$, we obtain by Corollary \rf{c3.14} $\wh\jmath_1(c)=\wt c\cdot (\eta_1
-\b\eta_0)$ \fs $\wt c\in K_1[X]\sms{{\bf0}}$. Note that $\deg(c)\nde3.46 = \deg(\wh\jmath(c))
\nde2.29 = \deg(\wt c)+1$. Thus $\deg(\wt c)=n$. Hence if $b=c$ in \er{4.27},
then $\wh\jmath_1(b)=\jh_1(c)=\wt c\cdot(\eta_1-\b\eta_0)$ in $K_1[X]$. If
$b=c\cdot\prodl_{i=1}^n (\ve_1-a(i)\ve_0)$ in \er{4.27}, then $\jh_1(b)
\nde3.49 = \jh_1(c)\cdot \jh\bigl(\prodl_{i=1}^n (\ve_1-a(i)\ve_0)\bigr)$ in $K_1[X]$.
Recall that $\jh_1:K[X]\to K_1[X]$ is a ring-\hm sm by \E\Pr\ \rf{p3.23}.
\E\Ip $\jh_1$~is a  monoid-\hm sm from the \mlv\ monoid of the ring $K[X]$
into the \mlv\ monoid of the ring $K_1[X]$. \E\Tf by \era2{1.131} $\jh_1
\bigl(\prodl_{i=1}^r(\ve_1-a(i)\ve_0)\bigr) = \prodl_{i=1}^r
(\jh_1(\ve_1-a(i)\ve_0))$ where $\prodl_{i=1}^r$ on the \RHS\ of the \et y
is the \cme product in $K_1[X]$. \Mo $\jh_1(\ve_1-a(i)\ve_0)\nde3.47 =
\jh_1(\ve_1)-\jh_1(a(i)\ve_0)\nde3.48 = \eta_1- j_1(a(i))\eta_0$. \csq,
$\jh_1(b) = \wt c\cdot(\eta_1-\b\eta_0)\cdot \prodl_{i=1}^r (\eta_i - j_1
(a(i))\eta_0)$. Thus in both cases of \er{4.27}, $\jh(b)$~is of the form
$\wt c\cdot \prod_{k=1}^t (\eta_k-d(k)\eta_0)$ in $K_1[X]$ where $d:[1,t]
\to K_1$, $t\in\Na$, and $\wt c\in K_1[X]\sms{{\bf0}}$ with $\deg(\wt c)=n$.
Indeed, if $b=c$, then $t:=1$, $d(1):=\b$, and in the other case $t:=r+1$,
$d(k):=j(a(k))$, $k\in[1,r]$, $d(r+1):=\b$.

Since $\deg(\wt c)=n\in M$, \te s a field \ext\ $(F,j_2)$ with $j_2:K_1\to F$
\st $\jh_2(\jh_1(b))$ splits over~$F$. Define $j:K\to F$ by setting $j:=j_2\circ
j_1$. As a \cm\ of two in\jc\ ring-\hm sms, $j$~is an in\jc\ ring-\hm sm. \Mo
$\jh(b)$ splits over~$F$. Indeed, $\jh(b)(l) \nde3.40 = j(b(l)) = (j_2\circ
j_1)(b(l)) = j_2(j_1(b(l))) \nde3.40 = j_2(\jh_1(b)(l)) \nde3.40 = \jh_2(
\jh_1(b))(l)$, $l\in\N$. \If that $n+1\in M$. Thus $M=\N$.
\endproof

Summarizing, we have

\bth4.21
\E\te s a finite field of order~$p^n$ \fe prime~$p$ and every $n\in\Na$.
\eth

\proof
Let $p$ be prime, $n\in\Na$, let $K$ denote the field $(\N_p,+_p,\cdot_p,
0,1)$ and let $\ve_{p^n}-\ve_1\in K[X]$. By Theorem \rf{t4.20} \te s a field \ext\
$(F,j)$ of~$K$ \st $\wh\jmath(\ve_{p^n}-\ve_1)$ splits over~$F$. In view of
\E\Pr\ \rf{p4.6} the subset of~$F$, $G:=\{\a\in F: \a^{p^n}=\a\}$ is a
subfield of $(F,+,\cdot,0,1)$ of order $p^n$ and $(G,+,\cdot,0,1)$ is a
field of order~$p^n$ by Lemma \rfa4{l4.24}.
\endproof

We next turn to the proof of \uq, up to \is sms, of a finite field of
order~$p^n$. Note that in view of \E\Pr\ \rfa4{p4.2}\,(ii) one may expect
\uq\ only up to \is sms.

Let $(E,+,\cdot,0,1)$ be a finite field with $\#(E)=p^n$. Set
\beq4.29
\ve_i(l):= \bca
1 & \hbox{for }l=i, \\
0 & \hbox{for }l\ne i,
\eca \q\ i,l\in \N.
\e
Clearly $\ve_i\in E[X]$ \fa $i\in\N$. Set
\beq4.30
b:=\ve_{p^n}-\ve_1.
\e
Then $b\in E[X]$, and by \E\Pr\ \rf{p4.2} we have
\beq4.31
b=\prod_{\a\in E}(\ve_1-\a\ve_0).
\e

By \E\Pr\ \rfa4{p4.29} the field $E$ possesses \ooo \pf\ which we denote
by~$\check E$. By Lemma \rfa4{l4.24}, $(\check E,+,\cdot,0,1)$ is a field
where $+,\cdot$ denote the \rt ion to~$\check E$ of the \ad\ and \mlc\
of~$E$. \Mo we define the inclusion map $i:\check E \to E$ by setting
\beq4.32
i(x):=x, \q x\in\check E.
\e
Clearly $i$ is an in\jc\ ring-\hm sm, hence $(E,i)$ is a field \ext\
of~$\check E$ (see \E\df\ \rf{d3.21}).

\blm4.22
Let $K_1,K_2$ and $j$ be as in \E\df\ \rf{d3.21}. Then $j(K_1)$  is a
subfield of~$K_2$.
\elm

\proof
Let $K_i:=(K_i,+_i,\cdot_i,0_i,1_i)$, $i=1,2$. We prove (a), (b), (c) of
\E\df\ \rfa4{d4.23}.

(a) We have $0_2=j(0_1)\in j(K_1)$.

Let $x,y\in j(K_1)$ and $\check x,\check y\in K_1$ be \st $x=j(\check x)$,
$y=j(\check y)$. Then $x+_2y= j(\check x)+_2j(\check y) = j(\check x+_1\check y) \in
j(K_1)$. Hence $j(K_1)$ is a \sbm\ of $(K_2,+_2,0_2)$. Let $x,\check x$ be as
above. \E\te s $\check z\in K_1$ \st $\check x+_1\check z=0_1$. Set
$z:=j(\check z)$. Then $x+_2z= j(\check x)+_2 j(\check z)= j(\check x +_1
\check z) = j(0_1)=0_2$. Hence $z\in j(K_1)$ is the  inverse of~$x$ in the
monoid $(j(K_1),+_2,0_2)$. \If that $j(K_1)$ is  a subgroup of
$(K_2,+_2,0_2)$.

(b) We have $1_2=j(1_1)\in j(K_1)$. Let $x,y\in j(K_1)$ and $\check x,\check
y\in K_1$ be \st $x=j(\check  x)$, $y=j(\check y)$. Then $x\cdot_2 y=
j(\check x)\cdot_2j(\check y) = j(\check x\cdot_1\check y) \in j(K_1)$. Hence
$j(K_1)$ is a \sbm\ of $(K_2,\cdot_2,1_2)$.

(c) Let $x\in j(K_1)\sms{0_2}$. Then \te s $\check x\in K_1\sms{0_1}$ \st
$x=j(\check x)$ since $j$ is in\jc. Since $K_1$ is a field \te s $\check y\in
K_1\sms{0_1}$ \st $\check x\cdot_1\check y=1_1$. Set $y:=j(\check y)\in
j(K_1)$. Then $x\cdot_2y = j(\check x)\cdot_2 j(\check y) = j(\check x\cdot_1
\check y) = j(1_1)=1_2$. \If that $x$~has an inverse in $(j(K_1),\cdot_2,
1_2)$. \E\Tf\ $j(K_1)$ is a subfield of~$K_2$.
\endproof

We now return to the field \ext\ $(E,i)$ of $\check E$, and define $\wh\imath:
\check E[X] \to E[X]$ by setting as in \er{3.40},
\beq4.33
(\wh\imath(c))(l):=i(c(l)), \q l\in\N.
\e
Thus $(\wh\imath(c))(l)=c(l)$, $l\in\N$, hence
\beq4.34
\wh\imath(c)=c,\q c\in \check E[X].
\e
Note that $\ve_{p^n}-\e_1\in \check E[X]$ and thus by \er{4.31}, \er{4.34} $\wh\imath(\ve_{p^n}
-\ve_1)\in E[X]$ splits over~$E$ (see \E\df\ \rf{d4.3}). More is true.

\blm4.23
If $H$ is a subfield of~$E$ \st $\ve_{p^n}-\ve_1\in \check E[X]$ splits
over~$H[X]$, then $H=E$.
\elm

\proof
Suppose, for \cd ion, that \te s a
\ti{proper} subfield $H$ of~$E$ \st $\ve_{p^n}-\ve_1$ splits over~$H$. Then,
$\ve_{p^n}-\ve_1\in H[X]$ and
by \er{4.4}, $\ve_{p^n}-\ve_1 = \b(0)\prodl_{i=1}^{p^n}(\ve_1-\b(i)\ve_0)$ \fs
$\b:[0,p^n]\to H\sbs E$. We have $\b(0)=(\ve_{p^n}-\ve_1)(p^n)$ by Lemma
\rf{l4.8}\,(iii), hence $\b(0)=1$. Let $\Phi_\g : E[X] \to E$, $\g\in E$,
denote the evaluation map in~$E$ defined in \er{3.28}. By \er{3.32}, \er{3.33}
with $i:=0$, we find that $\Phi_\g:(E[X],\cdot,\ve_0) \to (E,\cdot,1)$ is a
monoid-\hm sm. Then by \era2{1.130} we have
\beq4.35
\Phi_\g\Bigl(\prod_{i\in I} c_i\Bigr) = \prod_{i\in I} (\Phi_\g(c_i))
\qh{\fa}\g\in E,
\e
where $I$ is a \nf, $c_i\in E[X]$ \fa $i\in I$, $\prodl_{i\in I}c_i$ is the
\cme product in $(E[X],\cdot,\ve_0)$ and $\prodl_{i\in I}(\Phi_\g(c_i))$ is
the \cme product in $(E,\cdot,1)$. Then, \fe $\g\in E$, we obtain on the one
hand $\Phi_\g (\ve_{p^n}-\ve_1) = \Phi_\g \bigl(\prodl_{i=1}^{p^n}(\ve_1-\b(i)
\ve_0)\bigr) \nde4.35 = \break \prodl_{i=1}^{p^n}\bigl(\Phi_\g(\ve_1-\b(i)\ve_0)\bigr)
\nad{\er{3.30},\er{3.31},\er{3.33}} = \prodl_{i=1}^{p^n}(\g-\b(i))$. \E\oh
$\Phi_\g(\ve_{p^n}-\ve_1) \nde4.3 = \Phi_\g\bigl(\prodl_{\a\in E} (\ve_1-\a
\ve_0)\bigr) \nde4.35 = \prodl_{\a\in E}(\Phi_\g(\ve_1-\a\ve_0)) =
\prodl_{\a\in E}(\g-\a) = (\g-\g)\cdot \prodl_{\a\in E\sms\g}(\g-\a) =0\cdot
\prodl_{\a\in E\sms\g}(\g-\a)=0$.

Hence $\prodl_{i=1}^{p^n}(\g-\b(i))=0$ \fe $\g\in E$. Note that if $\g-\b(i)
\ne0$ \fa $i\in[1,p^n]$, then $\prodl_{i=1}^{p^n}(\g-\b(i))\ne0$ since
$(E\sms0,\cdot,1)$ is a monoid. \E\Tf \fe $\g\in E$, \te s $i\in[1,p^n]$ \st
$\g=\b(i)\in H$. Hence $E\sbs H$. Thus $H=E$, a \cd ion.
\endproof

\bdf4.24
Let $(F,j)$ be a field \ext\ of a field~$K$ and let $c\in K[X]$ be \st
$\wh\jmath(c)$ splits over $F[X]$. The field $F$ is called a \ti{splitting
field\/} for~$c$ (see \cite{Artin}) if there is no field \ext\ $(F',j')$ of\index{splitting field}
the field~$K$ \st $j'(K)$ is a \ti{proper} subfield of~$F$ and $j'(c)$ splits
over~$F'$.
\edf

Thus the finite field $E$ is a splitting field for the \fp\ $\ve_{p^n}-\ve_1
\in \check E[X]$.

We now show that if $K$ is a field and $c\in K[X]\sms{\bf0}$ with
$\deg(c)\ge1$, then \te s a field \ext\ $(F,j)$ \st $F$~is a splitting field
for $\wh\jmath(c)$.

\bth4.25
Let $K$ be a field and let $c\in K[X]\sms{\bf0}$ of degree larger than zero.
Then \te s a splitting field for~$c$.
\eth

\proof
By Theorem \rf{t4.20} \te s an \ext\ field $(F,j)$ \st $c$~splits over~$F$. Let $a:[0,\deg(c)] \to F$
be as in \E\df\ \rf{d4.3}, let $A:=\{a(i): i\in[1,\deg(c)]\}$ and $\cA:=
\{H\sbs F: H$ is a subfield of~$F$ \st $A\sbs H\}$. Note that $c$ splits over
a subfield~$H$ of~$F$ iff $H\in \cA$. \E\Ip $F\in\cA$. Set $\ov H:=
\Bca_{H\in\cA}H$. By Lemma \rfa4{l4.27}\,(i) $\ov H$~is a subfield of~$F$,
hence also a field by Lemma \rfa4{l4.24}. \Mo $A\sbs \ov H$. Suppose that
\te s a proper subfield $\wt H$ of~$\ov H$, then \te s $i\in[1,\deg(c)]$ \st
\st $a(i)\notin\wt H$, otherwise $\wt H\in\cA$, hence $\ov H\sbs \wt H$, which
is impossible since $\wt H$ is a \ti{proper} subfield of~$\ov H$. Suppose, for
\cd ion, that $c$~splits over~$\wt H$. Then \te s $\g: [0,\deg(c)]\to\wt H$,
\st $\g(0)\ne0$ and $c = \g(0)\prodl_{i=1}^{\deg(c)}(\ve_1-\g(i)\ve_0)$.
Proceeding as in the proof of Lemma \rf{l4.23} we find
\[
\g(0)\prod_{l\in[1,\deg(c)]}(a(i)-\g(l)) = 0.
\]
Since $\g(0)\nde4.4 = c(\deg(c))\ne0$, \te s $\ov l\in[1,\deg(c)]$ \st $a(i)
-\g(\ov l)=0$. Then $a(i)=\g(\ov l)\in \wt H$, a~\cd ion. Hence $\ov H$ is
a splitting field for~$c$.
\endproof

Now let $(E',+',\cdot',0',1')$ be a finite field of order~$p^n$. By Theorem
\rf{t1.30} $\chr(E)=\chr(E')=p$ and both $\check E$ and ${\check E}'$ are \is c
to the field $\N_p$. Hence \te s an \is sm $\vf:\check E\to {\check E}'$. As
in \er{3.40} we define a map $\wh\vf$ from $E[X]$ into $E'[X]$ by setting
\beq4.36
\wh\vf(c)(l) := \vf(c(l)), \q l\in\N,\ c\in E[X].
\e
By \E\Pr\ \rf{p3.23} we find that $\wh\vf$ is a ring-\is sm from $E[X]$ onto
$E'[X]$. We set $\eta_i:= \wt\vf(\ve_i)$, $i\in\N$, and observe that
$\eta_{p^n}-\eta_1 = \wh\vf(\ve_{p^n}-\ve_1)$. Using the same proof as
for~$E$ we find that $E'$~is a splitting field for $\eta_{p^n}-\eta_1$. It
turns out (see \cite[Theorem~10]{Artin}, \cite[Theorem 17.22]{Is}) that if $K_1,K_2$
are fields, $\psi:K_1\to K_2$ is a (ring-)\is sm, $b_i\in K_i[X]$, $i=1,2$,
\st $b_2=\wt\psi(b_1)$ where $\wt\psi$ is defined as in \er{4.36}, and $F_1$
(resp.~$F_2$) are splitting fields of~$b_1$ (resp.~$b_2$), then \te s a
ring-\is sm extending~$\psi$. \csq, $E$ and~$E'$ are \is c, thus solving
the problem of \uq\ up to \is sms of a finite field. However, we shall
give a different proof using the fact that the \mlv\
group of a finite field is \ti{cyclic} (see Corollary \rfa4{c4.42}\,(ii)).

We first investigate some \pp ies of the set of \ti{monic} \fp s over a
field~$K$.

\bdf4.27
Let $K$ be a field and let $K[X]$ be the ring (resp.\ $K$-algebra) of \fp s
over~$K$. A~nonzero \fp\ $a\in K[X]\sms{\bf0}$ is called \ti{monic} if its\index{formal polynomial!monic}
leading \cf\ is equal to~$1$. We use the (nonstandard) notation $MK[X]$ for\glossary{$MK[X]$}
the set of monic \pl s,
\beq4.37
MK[X]:=\{a\in K[X]\sms{{\bf0}}: a(\deg(a))=1\}.
\e
\edf

Clearly, the basis vectors $\ve_k$, $k\in\N$, defined in \er{2.4}, are monic,
as well as $b$~defined in \er{4.30}. If $(K[X],\cdot,\ve_0)$ denotes the \mlv\
monoid of the ring $(K[X],+,\cdot,{\bf0},\ve_0)$, then $MK[X]\sbs K[X]$,
$\ve_0\in MK[X]$, and $u\cdot v\in MK[X]$ whenever $u,v\in MK[X]$ by \er{2.28},
\er{2.29}. Hence $MK[X]$ is a \sbm\ of $(K[X],\cdot,\ve_0)$. Recall that
$K[X]\sms0$ is a \sbm\ of $(K[X],\cdot,\ve_0)$ by \er{2.29}, and
that the monoid $(K[X]\sms{{\bf0}},\cdot,\ve_0)$ is \cnc e (see \era2{1.8}).
Indeed, if $u,v,w\in K[X]\sms{{\bf0}}$ and $u\cdot w=v\cdot w$, then $(u-v)
\cdot w=\bf0$, hence $u-v=\bf0$ by \er{2.29} since $w\ne\bf0$. Thus $u=v$. \If
that the monoid $(MK[X],\cdot,\ve_0)$ is a \Cm\ (see \E\df\ \rfa2{d1.4}). \Mo
if $u,v\in MK[X]$ \sf y $u\cdot v=\ve_0$, then $\deg(u)+\deg(v)\nde2.29 =
\deg(u\cdot v) = \deg(\ve_0)\nde2.11 = 0$, hence $\deg(u)=\deg(v)=0$. Since
$\ve_0$ is the only monic \pl\ of degree zero, we have $u=v=\ve_0$. \If that
the monoid $(MK[X],\cdot,\ve_0)$ is a \PM\ (see \E\df\ \rfa2{d1.4}). As in the
\PM\ $(\Na,\cdot,1)$ we introduce the \og~$|$ in $MK[X]$. Let $u,v\in MK[X]$.
\beq4.38
u|v \qh{if \te s $w\in MK[X]$ \st $v=u\cdot w(= w\cdot u)$.}
\e
Note that
\beq4.39
\ve_0|u\hbox{ \fa} u\in MK[X]
\e
and that
\beq4.40
u|v \qh{implies} \deg(u)\le\deg(v), \q u,v\in MK[X],
\e
which is an analogue of \era3{8.57} in the case of $(\Na,\cdot,1)$.

Clearly, $\ve_0|u$ and $u|u$ \fa $u\in MK[X]$. In the case of $(\Na,\cdot,1)$,
a \nm\ $n\in\Na\sms1$ is called prime if $1$~and~$n$ are the only divisors
of~$n$. Similarly, a monic \pl\ $a\in MK[X]\sms{\ve_0}$ (i.e.\ of degree larger
than zero) is called \ti{ir\rd} if $\ve_0$ and~$a$ are the only divisors of~$a$.

\bex4.28
Verify that this notion of irreducibility is consistent with \E\df\ \rf{d3.11}.
\eex

The \fw\ analogue of Lemma \rfa3{l9.13} holds.

\blm4.29
Every \el\ $a$ of $MK[X]\sms{\ve_0}$ is either ir\rd\ or a \cme product of
ir\rd\ \el s of $MK[X]\sms{\ve_0}$.
\elm

\proof
We proceed by \In\ using Lemma \rfa3{l9.12} with $(W,\le):=(\Na,\le)$,
$e:=1$. Set $M:=\{n\in\Na: \hbox{Lemma \rf{l4.29} holds for }\deg(a)=n\}$.

$1\in M$: If $\deg(a)=1$, then $a=\ve_1+a(0)\ve_0$, $a(0)\in K$. Suppose
$u\in MK[X]$ divides~$a$. Then \te s $w\in MK[X]$ \st $a=u\cdot w$. Hence
$1=\deg(a)=\deg(u\cdot w)\nde2.29 = \deg(u)+\deg(w)\ge \deg(u)$. If $\deg(u)
=0$, then $u=\ve_0$. If $\deg(u)=1$, then $\deg(w)=0$, hence $w=\ve_0$ and
$u=a$. Thus $a$~is ir\rd.

Suppose $\deg(a)=n>1$, and suppose $\zo1,n \sbs M$. We show that $n\in M$. By
Lemma \rf{l4.21} \te\ $\wt c\in K[X]\sms{\bf0}$ ir\rd\ and $\wt d\in K[X]\sms
{\bf0}$ \st $a=\wt c\cdot\wt d$. Set $c:=(\wt c(\deg(\wt c)))\Inv \wt c$ and
$d:=(\wt d(\deg(\wt d)))\Inv \wt d$. Then, by Exercise \rf{ex4.30} below,
$a=c\cdot d$, where $c,d\in MK[X]$ and $c$~ir\rd. Since $c$ is ir\rd, $\deg(c)
\ge1$. If $\deg(d)=0$, then $d=\ve_0$ and $a=c\cdot \ve_0=c$ is ir\rd. If
$\deg(d)>0$, then $\deg(a)\nde2.29 = \deg(c)+\deg(d)>\deg(d)$, hence $\deg(d)
<n$. Thus $\deg(d)\in\zo1,n \sbs M$. \If that \te\ $N\in\Na$, $\zb di{[1,N]}$,
$d_i\in MK[X]$ ir\rd\ \fa $i\in[1,N]$, \st $d=\prodl_{i=1}^N d_i$, where
$\prodl_{i=1}^N$ is the \cme product in $(MK[X],\cdot,\ve_0)$. Setting $d_0
:=c$, we have $a=c\cdot d=d_0\cdot d = d_0\cdot \prodl_{i=1}^N d_i =
\prodl_{i=0}^N d_i$, where $d_i$~is ir\rd\ \fa $i\in[0,N]$. \E\Tf\ $n\in M$.
By Lemma \rfa3{l9.12} $M=\Na$.
\endproof

\bex4.30
Let $K$ be a field.

\hph i,ii, Let $a\in MK[X]$, $\wt c,\wt d\in K[X]\sms{\bf0}$ be \st $a=\wt c
\cdot\wt d$. Show that \te\ $\la,\mu\in K\sms0$ \st $a=c\cdot d$ with
$c=\la\wt c$, $d=\mu\wt d$, $c,d\in MK[X]$, and $\deg(c)=\deg(\wt c)$, $\deg(d)
=\deg(\wt d)$.

\hph ii,i, Show that if $a\in K[X]$ is ir\rd, then so is $\la a$ \fa $\la\in
K\sms0$.

\hph iii,, Show that if $a\in K[X]\sms{\bf0}$ and $P_a:K[X]\to K[X]$ is defined as in
\er{3.10}, then $P_{\la a}=P_a$ \fa $\la\in K\sms0$. \E\Ip $K_{\la a}=K_a$,
$\cdot_{\la a}=\cdot_a$ \fa $\la\in K\sms0$ in \er{3.17}.
\eex

We now return to the proof of the ``\uq''. As a first step we show that if
$E$~is a finite field of order~$p^n$, $\check E$~is its \pf\ and $\a\in E\sm
\check E$ is a \Gn\ of the cyclic group $E^\t=E\sms0$, then the set $\{\a^i:
i\in\zo0,n \}$ of \IT s of~$\a$ in the monoid $(E^\t,\cdot,1)$ is a basis of
the \vs~$E$ over~$\check E$.

\blm4.31
Let $K$ be a proper subfield of a field~$F$ and let $\a\in F\sm K$. Set
\beq4.41
A_k:=\{\a^i\in F: i\in\zo0,k \}, \q k\in\Na.
\e
Then either $\a^l\notin\spn(A_l)$ \fa $l\in\Na$ or \te s $\ov l\ge2$ \st
\beq4.42
\{\a^i\in F: i\in\zo0,{\ov l} \} \qh{is \li}
\e
as a subset of the \vs\ $(K,F)$, and
\beq4.43
\a^m \in \spn(A_{\bar l}) \q\hbox{\fa} m\ge \ov l.
\e
\elm

\proof
It suffices to consider the case where \te s $\wt l\in\Na$ \st $\a^{\td l}\in
\spn(A_{\td l})$. Since $\spn(A_1)=\{\la\a^0\in F:\la\in K\}= \{\la\in F: \la
\in K\}=K$, and $\a\notin K$, we have $\wt l\ge2$. Set $N:=\{l\in\Na\sms1:
\a^l\in\spn(A_l)\}$. Then $\wt l\in N$. Since $N$~is nonempty, $N$~possesses
a least \el\ by Theorem \rfa1{t3.22}. We denote this least \el\ by~$\ov l$.
Note that $\ov l\ge2$. Thus
\bea4.44
\a^{\bar l}&\in \spn(A_{\bar l}),\\
\a^l &\notin\spn(A_l) \qh{\fa} l\in\zo1,{\ov l} .\lb{4.45}
\e
We first show \er{4.42}. Let $M:=\{k\in\zo1,{\ov l} : A_k \hbox{ is \li}\}$.
We have $1\in M$ since if $\la\a^0=0$, $\la\in K$, then $\la=\la 1=\la\a^0=0$.
In view of Lemma \rfa1{l3.46} it is \sft\ to show that if $A_k$ is \li\ for
$k\in[1,\ov l-1]$, then $A_{k+1}$ is \li. Let $\zb\la i{[0,k]}\sbs K$ be \st
$\suml_{i=0}^k \la_i\a^i=0$. We have $\la_k=0$, otherwise $\a^k=-(\la_k)\Inv
\suml_{i=0}^{k-1}\la_i\a^i$, hence $\a^k\in\spn(A_k)$ \cd ing \er{4.45} since
$k<\ov l$. \E\Tf\ $\suml_{i=0}^{k-1} \la_i\a^i=0$, hence $\la_i=0$ \fa
$i\in[0,k-1]$ since $A_k$ is \li. \If that $\la_i=0$ \fa $i\in[0,k]$, hence
$A_{k+1}$ is \li, and by Lemma \rfa1{l3.46}, $M:=[1,\ov l]$, which proves
\er{4.42}.

We now prove \er{4.43}. We proceed by \In\ on $m\ge\ov l$. Set $M:=\{m\in\N,
m\ge \ov l: \a^m\in\spn(A_{\bar l})\}$.

$\ov l \in M$: by \er{4.44}.

\ti{$m\in M$ implies $m+1\in M$}: $\a^{m+1}=\a(\a^m)\nad{m\in M} = \a\bigl(
\suml_{i=0}^{\bar l-1} \la_i\a^i\bigr) = \suml_{i=0}^{\bar l-1}\a\la_i\a^i =
\suml_{i=0}^{\bar l-1}\la_i\a^{i+1} = \suml_{i=0}^{\bar l-2}\la_i\a^{i+1}
+\la_{\bar l-1}\a^{\bar l} = \suml_{j=1}^{\bar l-1}\la_{j-1}\a^j +
\la_{\bar l-1}\a^{\bar l}$, \fs $\la_i\in K$, $i\in[0,\ov l-1]$. Since
$\a^{\bar l}\in\spn(A_{\bar l})$, we have $\a^{\bar l} = \suml_{j=0}^{\bar l-1}
\mu_j\a^j$ \fs $\mu_j\in K$, $j\in[0,\ov l-1]$. Hence $\a^{m+1}=\la_{\bar l-1}
\mu_0\a^0 + \suml_{j=1}^{\bar l-1}(\la_{j-1}+\la_{\bar l-1}\mu_j)\a^j$. \E\Tf
$\a^{m+1}\in\spn(A_{\bar l})$, hence $M=\{m\in\N: m\ge\ov l\}$.
\endproof

\bex4.32
Show that if in Lemma \rf{l4.31} $a^l\notin\spn(A_l)$ \fa $l\in\Na$, then the
\vs\ $(K,F)$ is infinite-\dm al.
\eex

\blm4.33
Let $(E,+,\cdot,0,1)$ be a finite field of order~$p^n$, $p$~prime, $n\ge2$.
Let $\a\in E^\t$ be a \Gn\ of the cyclic $($\pn\/$)$ group $(E^\t,\cdot,1)$.
Then

\hph i,i, $\{\a^i\in E^\t: i\in\zo0,n \}$ is a basis of the \vs\ $(\check
E,E)$, where $\check E$ is the \pf\ of~$E$.

\hph ii,, Let $\la_i\in \check E$, $i\in[0,n)$, \sf y
\beq4.46
\a^n = \sum_{i=0}^{n-1}\la_i\a^i
\e
and let $a\in \check E[X]$ be defined by
\beq4.47
a(l):=\bca
-\la_l, & l\in\zo0,n ,\\
1, & l=n,\\
0, & l>n,\ l\in\Na,
\eca
\e
then $a$ is monic, ir\rd\ in $\check E[X]$ and $\deg(a)=n$.
\elm

\proof \

(i) Since $(E^\t,\cdot,1)$ is a finite cyclic (\pn) group of order $p^n-1$
and since $\a$ is a \Gn, the map $\ov\vf_\a: \zo0,p^n \to E^\t$ defined by
$\ov\vf_\a (k):=k\ddt \a$ is a monoid-\is sm from $(N_{p^n-1},+_{p^n-1},0)$
onto~$E^\t$ by \E\Pr\ \rfa4{p3.10}\,(iv). \E\Ip
\beq4.48
\a^i\ne \a^j \qh{\fa} i,j\in\zo0,p^n , \q i\ne j.
\e
By \E\Pr\ \rfa4{p3.65}
\beq4.49
\a^{p^n-1}=1.
\e
Thus $\a^{p^n-1}\in\spn(A_{p^n-1})$ (defined in \er{4.41}), hence by Lemma
\rf{l4.31} with $F:=E$, $K:=\check E$, both \er{4.42} and \er{4.43} hold \fs
$\ov l\ge2$. We first show that $\ov l<n+1$. Indeed, if $\ov l\ge n+1$, then
by \er{4.43} $A_{\bar l}$ is a basis of the \vs\ $(\check E,\spn(A_{\bar
l}))$. By Lemma \rf{l1.23}\,(iv) this \vs\ is \is c, thus \ep, to the \vs\ $\check
E^{[1,\bar l]}$ over the field~$\check E$. Note that $\#(\check E^{[1,\bar l]})
\nda23.38 = \#(\check E)^{\#([1,\bar l])} = p^{\bar l}\ge p^{n+1}$ by \era2
{2.33}. Thus $\#(\spn A_{\bar l})\ge p^{n+1}$. This is impossible since
$\spn(A_{\bar l})\sbs E$, hence $\#(\spn(A_{\bar l})) \nda23.13 \le \#(E) =p^n
\nda22.33 < p^{n+1}$. \If that $\ov l\le n$. We now show that $\ov l=n$.
Suppose for \cd ion that $\ov l<n$, then by \er{4.43} $\a^l\in\spn(A_{l})$
\fa $l\ge \ov l$, and by \df\ $\a^l\in \spn(A_{\bar l})$ \fa $l<\ov l$. Hence
$E^\t\sbs \spn(A_{\bar l})$. Clearly, $0\in\spn(A_{\bar l})$, thus $E\sbs
\spn(A_{\bar l})$. Since $A_{\bar l}$ is \li\ we obtain as above
$\#(\spn(A_{\bar l}))=p^{\bar l}$. A~\cd ion since $p^n=\#(E)\le \#(
\spn(A_{\bar l})) = p^{\bar l}<p^n$. \E\Tf $\ov l=n$. \If that $A_n$~is
\li\ in $(\check E,E)$. Since $\dim_{\check E}(E)=n$ and $\#(A_n)=n$, $A_n$~is
a basis of~$E$ by \E\Pr\ \rf{p1.31}\,(i).

(ii) In view of (i) \te s \ooo $\la_i\in\check E$ \fe $i\in[0,n-1]$ \st
\er{4.46} holds. Clearly, $a$~is monic and  $\deg(a)=n\ge2$. We show that
$a$~is ir\rd\ in $\check E[X]$. Suppose \te\ $c,d\in M\check E[X]$ with
$\deg(c),\deg(d)\ge1$ \st $a=c\cdot d$ (see Exercise \rf{ex4.30}\,(i)). Let
$\Phi_\a : E[X]\to E$ denote the evaluation map defined in \er{3.28}.
Observe that $\Phi_\a(a)=0$, hence $0=\Phi_\a(a)= \Phi_\a(c\cdot d)\nde3.32 =
\Phi_\a(c)\cdot \Phi_\a(d)$. Hence $\Phi_\a(c)=0$ or $\Phi_\a(d)=0$ since
$E$ is a field. Suppose $\Phi_\a(c)=0$. Set $m:= \deg(c)$. We have
$n=\deg(a)\nde2.29 = \deg(c)+\deg(d) = m+\deg(d) > m$. \If that $\a^m =
\suml_{i=0}^{m-1}(-(c(i))\a^i) \in \spn(A_m)$ with $m<n=\ov l$, \cd ing
\er{4.45} since $-(c(i))\in\check E$, $i\in[0,m-1]$. The case $\Phi_\a(d)=0$
is similar. Thus $a$~is ir\rd\ in $\check E[X]$.
\endproof

\brm4.34 \

\hph i,ii, Observe that the \fp\ $a\in M\check E[X]$ defined by \er{4.47} has
$\a$ as a root in $E[X]$. As a con\sq\ of (i), this is the only \el\ of
$M\check E[X]$ of degree~$n$ having $\a$ as a root in $E[X]$.

\hph ii,i, If $q\in M\check E[X]$ with $\deg(q)\ge1$, then $q\cdot a\in
\check E[X]$ since $\check E[X]$ is a \sbm\ of $(E[X],\cdot,\ve_0)$ by \er{3.49}
with $K_1:=\check E$, $K_2:=E$ and $j(x):=x$, $x\in K_1$. \Mo $q\cdot a\in
M\check E[X]$ by \er{2.28}, \er{2.29}. Then $\suml_{i=0}^{\deg(a)+\deg(q)}
(q\cdot a)(i)\a^i \nde3.28 = \Phi_\a(q\cdot a) \nde3.32 = \Phi_\a(q)\cdot
\Phi_\a(a)=0$ since $\Phi_\a(a)=0$.

\hph iii,, It turns out that if $b\in M\check E[X]$ \sf ies $\Phi_\a(b)=0$,
then \te s $q\in M\check E[X]$ \st $b=q\cdot a$. This is a con\sq\ of the \fw\
\Pr.
\erm

\bpr4.35
Let $K$ be a proper subfield of a field~$F$ and let $\a\in F\sm K$. Let\glossary{$\wt\Phi_\a$}
$\wt\Phi_\a:K[X]\to F$ denote the \rt ion of the evaluation map $\Phi_\a$
defined in \er{3.28} to $K[X]$. Set
\beq4.50
N(\wt\Phi_\a):= \{c\in K[X]: \wt\Phi_\a(c)=0\}.
\e
If $N(\wt\Phi_\a)\ne\{{\bf0}\}$, then \te s \ooo \el\ $a\in MK[X]$ \st
\beq4.51
N(\wt\Phi_\a) = \{q\cdot a\in K[X]: q\in K[X]\}.
\e
\E\Ip $a\in N(\wt\Phi_\a)$ $($choose $q:=\ve_0)$.

\Mo $\deg(a)\ge2$ and $a$ is ir\rd,
\beq4.52
\deg(a)\le \deg(c) \qh{\fa} c\in N(\wt\Phi_\a)\sms{\bf0}.
\e
\epr

\proof
The map $\deg:N(\wt\Phi_\a)\sms{{\bf0}} \to \N$ is well-defined. The set
$\deg(N(\wt\Phi_\a)\sms{{\bf0}})$ has a least \el\ which we denote by~$\wt a$
since $(\N,\le)$ is well-ordered. \E\Tf \te s $\wt a\in N(\wt\Phi_\a)\sms
{\bf0}$ \st $\deg(\wt a)\le \deg(c)$ \fa $c\in N(\wt\Phi_\a)\sms{\bf0}$.
Set $a:=(\wt a(\deg(\wt a)))\mo \wt a$, then
$\deg(a)\nde2.13 = \deg(\wt a)$, $a\in MK[X]$ and \er{4.52} holds. We now
prove \er{4.51}. Let $c\in N(\wt\Phi_\a)\sms{{\bf0}}$. By Theorem \rf{t3.1}
\te\ $q,r\in K[X]$ \sf ying \er{3.1}, \er{3.2} with $b:=c$. If $r\ne\bf0$, then
$\deg(r)<\deg(a)$, and $r=c-q\cdot a$. Note that $a,q,q\cdot a, c, c-q\cdot a
\in K[X]$. Hence $\wt\Phi_\a(r)=\wt\Phi_\a(c-q\cdot a)$. Since $(c-q\cdot a)+
q\cdot a= (c+(-q\cdot a))+q\cdot a = c$, we have $\wt\Phi_\a(c-q\cdot a)+
\wt\Phi_\a(q\cdot a)= \wt\Phi_\a(c)$. Thus $\wt\Phi_\a(c-q\cdot a)\nde3.30 =
\wt\Phi_\a(c) - \wt\Phi_\a(q\cdot a) = 0-\wt\Phi_\a(q)\cdot \wt\Phi_\a(a) =
-\wt\Phi_\a(q)\cdot0 = 0$. Hence $r\in N(\wt\Phi_\a)\sms{{\bf0}}$ with
$\deg(r)<\deg(a)$, \cd ing \er{4.52}. \If that $r=0$, hence $c=q\cdot a$
\fs $q\in K[X]$, which together with Remark \rf{r4.34}\,(ii) proves
the ``\ex\ part'' of \er{4.51}. We now prove the \uq\ part of the \Pr.
Suppose \te s $a'\in N(\wt\Phi_\a)\sms{\bf 0}$ \sf ying \er{4.51} with $a:=a'$,
and $a'\in MK[X]$. Then \te\ $q,q'\in K[X]$ \st $a'=q\cdot a$ and $a=q'\cdot a'$.
Then $a'=q\cdot (q'\cdot a')= (q\cdot q')\cdot a'$. Since $a'=\ve_0\cdot a'$
we have $(\ve_0-(q\cdot q'))\cdot a'={\bf0}$. From \er{2.29} and $a'\ne\bf0$
we infer $\ve_0-(q\cdot q')=0$, \ev tly $q\cdot q'=\ve_0$. From \er{2.28},
\er{2.29} we deduce $\deg(q)+\deg(q')=0$, hence $\deg(q)=\deg(q')=0$ by
\era2{1.9}. Thus $q=q(0)\ve_0$, $q'=q'(0)\ve_0$, with $q(0)\ne0$, $q'(0)\ne0$.
Since $a'=q\cdot a$, we obtain $a'=q(0)\ve_0\cdot a= q(0)a$. Note that
$\deg(a')=\deg(q(0)a) \nde2.13 = \deg(a)$. Finally, $1=a'(\deg(a')) = (q(0)a)
(\deg(a)) = q(0)a(\deg(a)) =q(0)1= q(0)$, hence $a'=a$. This completes the
proof of the \uq\ part.

Next we show that $\deg(a)\ge2$. If $\deg(a)=0$, then $a=\ve_0$ since $a\in
MK[X]$. Then $N(\wt\Phi_\a)\nde4.51 = \{q\cdot \ve_0: q\in K[X]\} =K[X]$,
i.e.\ $\wt\Phi_\a(c)=0$ \fa $c\in K[X]$. But $\wt\Phi_\a(\ve_0)=\a^0=1\ne 0$.
Hence $\deg(a)\ne0$. If $\deg(a)=1$, then $a=\ve_1+a(0)\ve_0$, $a(0)\in K$.
Since $\wt\Phi_\a(a)=0$, we have $\a+a(0)\a^0 = 0$, hence $\a=-a(0)\in K$,
\cd ing $\a\in F\sm K$. Thus $\deg(a)\ne1$.

Finally, we show that $a$ is ir\rd. Suppose for \cd ion that $a$ is \rd. Then
\te\ $\wt c,\wt d\in K[X]$ with $\deg(\wt c)\ge1$, $\deg(\wt d)\ge1$ \st
$a=\wt c \cdot\wt d$. In view of Exercise \rf{ex4.30}\,(i) \te\ $c,d\in MK[X]$
with $\deg(c),\deg(d)\ge1$ \st $a=c\cdot d$. Then $\wt\Phi_\a(c)\cdot\wt\Phi_\a
(d)\nde3.32 = \wt\Phi_\a(c\cdot d) = \wt\Phi_\a(a)=0$. Since $K$ is a field
$\wt\Phi_\a(c)=0$ or $\wt\Phi_\a(d)=0$. If $\wt\Phi_\a(c)=0$, then $c\in
N(\wt\Phi_\a\sms{{\bf0}})$, hence \te s $q\in K[X]$ \st $c=q\cdot a$ by \er{4.51}.
Clearly, $q\ne\bf0$, since $c\ne\bf0$. Then $a=c\cdot d= q\cdot a\cdot d$, hence
$q\cdot d=\ve_0$ by the \asc ity, \cmt ity and \cnc ity of~$\cdot$ in $(K[X]\sms
{{\bf0}},\cdot,\ve_0)$. By \er{2.28}, \er{2.29} $\deg(q)+\deg(d)=0$, hence
$\deg(q)=\deg(d)=0$ by \era2{1.9}, \cd ing $\deg(d)\ge1$. The case $\wt\Phi_\a
(d)=0$ is similar.

Thus $a$ is ir\rd.
\endproof

\bco4.36
Let $E,\check E,\a$ and $a$ be as in Lemma \rf{l4.33}. If $\a$~is a root in~$E$
of an \el\ $c\in\check E[X]\sms{{\bf0}}$, then \te s $q\in\check E[X]\sms{{\bf0}}$
\st $c=q\cdot a$. \E\Ip $a$~divides $\ve_{p^n}-\ve_1$.
\eco

\proof
Let $\wt\Phi_\a:\check E[X] \to E$ denote the \rt ion of the evaluation map
$\Phi_\a:E[X]\to E$ to $\check E[X]$. Since $c\in N(\wt\Phi_a)\sms{\bf0}$, we may apply
\E\Pr\ \rf{p4.35} with $F:=E$, $K:=\check E$. Then \te s \ooo $\wt a \in
M\check E[X]$ \st $a=q\cdot\wt a$ \fs $q\in\check E[X]$. Clearly, $q\ne\bf0$
since $a\ne\bf0$. \Mo $\deg(\wt a)\ge2$ and $\wt a$~is ir\rd\ in $\check E[X]$.
Since $a$ is ir\rd\ in $\check E[X]$, $\deg(q)$ must be equal to zero. Hence
$q=q(0)\ve_0$ with $q(0)\in\check E\sms0$. Then $a=q(0)\wt a$. Since both
$a$ and~$\wt a$ belong to $M\check E[X]$, we have $q(0)=1$. Indeed, $\deg(a)
=\deg(\wt a)$ by \er{2.13} and $1=a(\deg(a)) = q(0)\wt a(\deg(\wt a))=q(0)$.
Thus $a=\wt a$. Since $c\in N(\wt\Phi_\a)\sms{\bf0}$, $c=q\cdot a$ \fs $q\in
\check E[X]\sms{\bf0}$. Finally, since $\wt\Phi_\a(\ve_{p^n}-\ve_1)=\a^{p^n}
-\a=0$ (see the proof of \E\Pr\ \rf{p4.2}), $\ve_{p^n}-\ve_1\in N(\wt\Phi_\a)
\sms{\bf0}$, hence as above with $c:=\ve_{p^n}-\ve_1$, $a$~divides $\ve_{p^n}
-\ve_1$.
\endproof

\brm4.37
The first part of the proof of \E\Pr\ \rf{p4.35} is based on the fact that the
map $\wt\Phi_\a$ from the ring $K[X]$ into the field~$F$ is a ring-\hm sm by
\er{3.30}, \er{3.32}, \er{3.33} with $i=0$. Note that \er{3.30} implies
$\wt\Phi_\a({\bf0})=0$. The nullspace $N(\wt\Phi_\a)$ defined in \er{4.50} of
such a ring-\hm sm is not only a subgroup of $(K[X],+,{\bf0})$ and a \sbm\
of $(K[X],\cdot,\ve_0)$ but it also \sf ies the stronger \cn\ $x\cdot y\in
N(\wt\Phi_\a)$ \fa $x\in N(\wt\Phi_\a)$ and all $y\in K[X]$.
\erm

\begin{dfn}[\cite{Alg}] \lb{d4.38}
Let $(R,+,\cdot,0,1)$ be a (\cmt e) ring (with unity), and let $I$~be a
subgroup of $(R,+,0)$. Then $I$~is called an \ti{ideal\/} of~$R$ if\index{ideal}
\beq4.53
x\cdot y\in I \qh{\fa $x\in I$ and all }y \in R.
\e
An ideal $I$ of $R$ is called \ti{\pn\/} if \te s $a\in I$ \st\index{ideal!principal}
\beq4.54
I= \{q\cdot a\in R: q\in R\} =: (a).
\e
\edf

Note that if a subset $A$ of $R$ is of the form $\{q\cdot a\in R: q\in R\}$
\fs $a\in R$, then $A$ is an ideal of~$R$.

Clearly, $R=(1)$ and $\{0\}=(0)$ are \pn\ ideals of~$R$.

The proof of \E\Pr\ \rf{p4.35} shows that every ideal of the ring $K[X]$ is \pn.

\bex4.39 \

\hph i,i, Let $F$ be a field and let $a\in F$. Show that $(a)=\{0\}$ if $a=0$
and $(a)=F$ if $a\ne0$.

\hph ii,, Let $R_i$, $i=1,2$, be (\cmt e) rings (with unity), and let $\vf:R_1
\to R_2$ be a ring-\hm sm. Show that $N(\vf):=\{x\in R_1: \vf(x)=0_2\}$ is an
ideal of~$R_1$. Show that

{\advance\parindent by15pt
(1) $\vf$ is in\jc\ if $R_1$ is a field.

(2) If $R_2$ is a field and $R_1$ is not a field then $N(\vf)\ne
\{0_1\}$.

}\eex

\bex4.40
Let $K$ be a proper subfield of a field $F$. Let $\a\in F\sm K$. Let
$\wt\Phi_\a:K[X]\to F$ be the \rt ion to $K[X]$ of the evaluation map $\Phi_a$
defined in \er{3.28}. Prove

\hph i,i, if $c\in N(\wt\Phi_\a)\sms{\bf0}$ is ir\rd\ in $K[X]$, then $N(\wt\Phi_\a)
=(c)$;

\hph ii,, if $c,d\in N(\wt\Phi_\a)\sms{\bf0}$ are monic and ir\rd, then $c=d$.
\eex

We now return to the situation of Lemma \rf{l4.33}. Let $E,\check E,E^\t,\a$
and~$a$ be as in Lemma \rf{l4.33}, and let $\wt\Phi_\a$ be as in \E\Pr\
\rf{p4.35} with $K:=\check E$ and $F:=E$. Since $a\in MK[X]$ is ir\rd\ with
$\deg(a)=n$, we can invoke \E\Pr\ \rf{p3.7} and Theorem \rf{t3.13}, and
consider the \pj~$P_a$ introduced in \er{3.10}. We recall that $R(P_a)$, the
range of~$P_a$, also denoted by~$K_a$, is a $n$-\dm al \vs\ over~$\check E$
with basis $\{\ve_i: i\in\zo0,n \}$. Since $\{\a^i\in E: i\in\zo0,n \}$ is a
basis of the \vs\ $E$ over~$\check E$ by Lemma \rf{l4.33}\,(i), \fe $x\in E$,
\te s \ooo $\wh x_i\in\check E$ \fe $i\in\zo0,n $, \st $x=\suml_{i=0}^n
\wh x_i\a^i$. Define $\wh x\in\check E[X]$ by setting
\beq4.55
\wh x(i):= \bca
\wh x_i, & i\in\zo0,n ,\\
0, & i\ge n,\ i\in\N.
\eca
\e
Then $x=\wt\Phi_\a(\wh x)$. Since $\wh x\in K_a$, we may write $x=\ov\Phi_\a
(\wh x)$ where $\ov\Phi_\a$ is the \rt ion of $\wt\Phi$ to~$K_a$. \If that
the map $\ov\Phi_\a:K_a\to E$ is \ti{sur\jc}. Since both $E$ and~$K_a$ are
$n$-\dm al \vs s over the finite field~$\check E$ with $\#(\check E)=p$, the
\vs s $E$ and~$K_a$ are linear \is c, hence \ep, $\#(E)=\#(K_a)=p^n$. \E\Tf
the map $\ov\Phi_\a$ is also \ti{in\jc} by \era2{3.12}. By Theorem \rf{t3.13},
$(K_a,+_a,\cdot_a,{\bf0},\ve_0)$ is a field, since $a$ is ir\rd. We now show
that $\ov\Phi_\a$ is a \ti{ring-\hm sm}. We have $\ov\Phi_\a({\bf0})=\Phi_\a({\bf0})
=\Phi_\a({\bf0}+{\bf0})\nde3.30 = \Phi_\a({\bf0}) + \Phi_\a({\bf0})$, hence
$\ov\Phi_\a({\bf0})={\bf0}$ by \cnc ity. \Mo $\ov\Phi_\a(\ve_0)=\Phi_a(\ve_0)
=\a^0=1$. Let $\wh x,\wh y\in K_a$. Then $\ov\Phi_\a(\wh x+_a \wh y)\nde3.18 =
\ov\Phi_\a(\wh x+\wh y)=\break \Phi_\a(\wh x+\wh y)\nde3.30 = \Phi_\a(\wh x)
+\Phi_\a(\wh y) =\ov\Phi_\a(\wh x)+\ov\Phi_\a(\wh y)$. Recall $\wh x\cdot_a
\wh y\nde3.17 = P_a(\wh x\cdot\wh y)$. By Theorem \rf{t3.1}, and the \df\
of~$P_a$, $\wh x\cdot \wh y=q\cdot a+P_a(\wh x\cdot\wh y)$ \fs $q\in K[X]$.
Thus $\Phi_\a(\wh x\cdot\wh y)= \Phi_\a(q\cdot a + P_a(\wh x\cdot \wh y))\nde3.30
= \Phi_\a(q\cdot a)+\Phi_\a(P_a(\wh x\cdot\wh y)) \nde3.32 = \Phi_\a(q) \cdot
\Phi_\a(a) + \Phi_a(\wh x\cdot_a\wh y)\nad{\er{4.46},\er{4.47}} = \Phi_\a(q)
\cdot 0 + \ov\Phi_\a(\wh x\cdot_a\wh y)$. Hence $\ov\Phi_\a(\wh x\cdot_a\wh y)
= \Phi_\a(\wh x\cdot\wh y)\nde3.32 = \Phi_\a(\wh x)\cdot \Phi_\a(\wh y) =
\ov\Phi_\a(\wh x)\cdot\ov\Phi_\a(\wh y)$.

We summarize in

\blm4.41
Let $E,\a$ and $a$ be as in Lemma \rf{l4.33}, let $(K_a,+_a,\cdot_a,{\bf0},
\ve_0)$ be the field introduced in Theorem \rf{t3.13}, and let $\ov\Phi_\a$
denote the \rt ion to~$K_a$ of the evaluation map $\Phi_\a$ defined in \er{3.28}.
Then $\ov\Phi_\a:K_a\to E$ is a ring-\is sm. We shall denote its \emph{inverse}\glossary{$\Psi_\a$}
by~$\Psi_\a$.
\elm

We consider a more general situation.

\blm4.42
Let $K$ be a proper subfield of a field~$F$, let $\a\in F\sm K$ and let $\ov a
\in MK[X]$ be ir\rd\ in~$K[X]$. Let $K_{\bar a}$ be the field introduced in
Theorem \rf{t3.13} and let $\ov\Phi_\a:K_{\bar a}\to F$ denote the \rt ion\glossary{$\ov\Phi_\a$}
to~$K_{\bar a}$ of the evaluation map~$\Phi_\a$ defined in \er{3.28}. If
$\wt\Phi_\a(\ov a)=0$, then the map $\ov\Phi_\a$ is an \emph{in\jc} ring-\hm
sm from $K_{\bar a}$ into~$F$.
\elm

\proof
Let $\wt\Phi_\a:K[X]\to F$ denote the \rt ion to $K[X]$ of $\Phi_\a$ defined
in \er{3.28} and let $P_{\bar a}$ be the \pj\ in $K[X]$ defined in \er{3.10}.
We claim
\beq4.56
N(\wt\Phi_\a)=N(P_{\bar a}),
\e
where $N(\wt\Phi_\a)$ is defined in \er{4.50} and $N(P_{\bar a})=
\{c\in K[X]: P_{\bar a}(x)={\bf0}\}$. Note that $N(P_{\bar a})\sbs
N(\wt\Phi_a)$. Indeed, $N(P_{\bar a})\nde3.15 = \{c\in K[X]: \hbox{\te s
$q\in K[X]$ \st} c=q\cdot \ov a\}$. Thus if $c\in N(P_{\bar a})$, we have
$\wt\Phi_\a(c) = \wt\Phi_\a(q\cdot \ov a) \nde3.32 = \wt\Phi_\a(q)\cdot
\wt\Phi_\a(\ov a) = \wt\Phi_\a(q)\cdot 0=0$. Since $\wt\Phi_\a(\ov a)=0$, we
have $N(\wt\Phi_\a)\ne \{{\bf0}\}$. Hence we may invoke \E\Pr\ \rf{p4.35}.
From \er{4.51} $\ov a\in N(\wt\Phi_\a)$, we infer that $\ov a=q\cdot a$ \fs
$q\in K[X]$ and some $a$ ir\rd\ with $\deg(a)\ge2$. Then $q\ne\bf0$ since
$\ov a\in MK[X]$ is not equal to~$\bf0$, and $\deg(\ov a)\nde2.29 = \deg(q) +
\deg(a)\ge2$. Since $\ov a$ is ir\rd, we must have $\deg(q)=0$, i,e.\
$q=q(0)\ve_0$ with $q(0)\in K\sms0$. \E\Tf $\ov a=q(0)a$. Since both $\ov a$
and~$a$ belong to $MK[X]$, we find $q(0)=1$ by \er{2.28}, hence $a=\ov a$.
Now let $c\in N(\wt\Phi_\a)$. By \er{4.51} and $a=\ov a$ we find $c\in
N(P_{\bar a})$. \E\Tf \er{4.56} holds.

\ti{$\ov\Phi_\a$ is a ring-\hm sm}: We have $\ov\Phi_\a({\bf0}) =
\Phi_\a({\bf0}) =0$ and $\ov\Phi_\a(\ve_0)=\Phi_\a(\ve_0)=1$. \Mo if
$u,v\in K_{\bar a}$, then $\ov\Phi_\a(u+_av) \nde3.18 = \ov\Phi_\a(u+v) =
\wt\Phi_\a(u+v) \nde3.30 = \wt\Phi_\a(u)+\wt\Phi_\a(v) = \ov\Phi_\a(u)+
\ov\Phi_\a(v)$.

$\ov\Phi_\a(u\cdot_{\bar a}v)\nde3.17 = \ov\Phi_\a(P_{\bar a}(u\cdot v))=
\wt\Phi_\a(P_{\bar a}(u\cdot v))$. Note that $u\cdot v = P_{\bar a}(u\cdot v)
+((u\cdot v)-P_{\bar a}(u\cdot v))$, and that $P_{\bar a}((u\cdot v) -
P_{\bar a}(u\cdot v)) = P_{\bar a}(u\cdot v) - P_{\bar a}(P_{\bar a}(u\cdot
y)) \nde3.14 = {P_{\bar a}(u\cdot v)-P_{\bar a}(u\cdot v) = \bf0}$. Thus
$(u\cdot v)-P_{\bar a}(u\cdot v) \in N(P_{\bar a}) \sbs N(\wt\Phi_\a)$, by
\er{4.56}. \E\Tf $\wt\Phi_\a(P_{\bar a}(u\cdot v))=\break \wt\Phi_\a(P_{\bar a}
(u\cdot v))+\wt\Phi_\a((u\cdot v)-P_{\bar a}(u\cdot v)) \nde3.30 =
\wt\Phi_\a(u\cdot v) \nde3.32 = \wt\Phi_\a(u)\cdot \wt\Phi_\a(v)
=\ov\Phi_\a(u)\cdot \ov\Phi_\a(v)$.

\ti{In\ji\ of $\ov\Phi_\a$}: Follows from Exercise \rf{ex4.39}\,(ii) since
both $K_a$ and~$F$ are fields and $\ov\Phi_\a$~is a ring-\hm sm.
\endproof

In the next lemma we consider an \ext\ of Theorem \rf{t3.13}.

\blm4.43
Let $K,K'$ be fields, let $\vf:K\to K'$ be a ring-\is sm and let $a\in MK[X]$
be ir\rd\ with $\deg(a)\ge2$. Set
\beq4.57
(\wh\vf(c))(l) := \vf(c(l)),\q l\in \N,\ c\in K[X].
\e
Then

\hph i,ii,
\bea4.58
\wh\vf(c)&={\bf0}' \qh{iff } c={\bf0},\\
\wh\vf(c)&\in K'[X]\sms{{\bf0}'} \qh{and } \deg(\wh\vf(c))=\deg(c) \hbox{ \fa}
c\in K[X]\sms{\bf0}. \lb{4.59}
\e

\hph ii,i, $\wh\vf: K[X]\to K'[X]$ is a ring-\is sm. \Mo
\beq4.60
\wh\vf(\la c) = \vf(\la) \wh\vf(c) \qh{\fa $\la\in K$ and all} c\in K[X].
\e

\hph iii,, Set
\beq4.61
a':=\wh\vf(a).
\e
Then $a'\in MK'[X]$, $a'$ is ir\rd\ in $K'[X]$ and $\deg(a')=\deg(a)$.

\hph iv,, Let $P_a$ $($resp.\ $P_{a'})$ be defined as in \er{3.10}, then
\beq4.62
\wh\vf \circ P_a = P_{a'}\circ \wh\vf.
\e
\elm

\proof \

(i), (ii): Clearly, $\wh\vf(c)\in {K'}^\N$. Since for $\la\in K$ $\vf(\la)=0'$
iff $\la=0$, we have $\supp \wh\vf(c)=\supp(c)$, $c\in K[X]$. \E\Ip $\wh\vf(c)
\in K'[X]$, $c\in K[X]$, $\wh\vf(c)={\bf0}'$ iff $c=\bf0$ and $\deg\wh\vf(c)=
\deg(c)$, $c\in K[X]\sms{\bf0}$. Note that \er{4.57} is \er{3.40} where $K_1
:=K$, $K_2:=K'$ and $j:=\vf$. Thus
\beq4.63
\wh\vf(\ve_k) = \ve'_k, \q k\in\N,
\e
where
\[
\ve'_k(l):= \bca
1' &\hbox{ if }k=l,\\
0' &\hbox{ otherwise.}
\eca
\]
From \er{3.47} and \er{3.49} we deduce that $\wh\vf$~is a ring-\hm sm. \Mo
\er{4.60} follows from \er{3.48}.

Set $\psi:=\vf\Inv$, $\wh\psi(c')(l):=\psi(c'(l))$, $l\in\N$, and $c'\in K'[X]$.
Then $\wh\psi$ maps $K'[X]$ into $K[X]$,
\beq4.64
\wh\psi \circ \wh\vf = \id_{K[X]}, \q \wh\vf\circ \wh\psi= \id_{K'[X]}.
\e
Indeed, $((\wh\psi\circ\wh\vf)(c))(l) = (\wh\psi(\wh\vf(c)))(l) = \psi(\vf
(c(l))) = (\psi\circ\vf)(c(l)) = c(l)$, $l\in\N$.

The proof of the other equality is similar. \If that $\wh\vf$ is bi\jc\ with
inverse~$\wh\psi$ and that $\wh\vf$ is a ring-\is sm.

(iii) By \er{4.59} $\deg(a')=\deg(a)$ and $a(\deg(a))=1$ by \as. Thus
$a'(\deg(a'))=\wh\vf(a)(\deg(a)) = \vf(a(\deg(a))) = \vf(1) = 1'$, hence
$a'\in MK'[X]$. Note that $\deg(a')=\deg(a)\ge2$. Suppose for \cd ion that
$a'$~is \rd. Thus \te\ $c',d'\in K'[X]$ with $\deg(c'),\deg(d')\ge1$ \st
$a'=c'\cdot' d'$. Then $\wh\vf(a) = \wh\vf(c)\cdot' \wh\vf(d) = \wh\vf(c\cdot d)$
by~(ii). Thus $a=\wh\vf\Inv(\wh\vf(a)) = \wh\vf\Inv(\wh\vf(c\cdot d)) = c\cdot
d$, \cd ing the fact that $a$~is ir\rd.

(iv) Let $c\in K[X]$. Then by Theorem \rf{t3.1} and \er{3.10} \te s $q\in K[X]$ \st
$c=q\cdot a+P_a(c)$. Since $\wh\vf$ is a ring-\is sm, we have $\wh\vf(c) =
\wh\vf(q)\cdot' \wh\vf(a) + \wh\vf(P_a(c)) \nde4.61 = \wh\vf(q)\cdot' a' +
\wh\vf(P_ac)$. From the \uq\ part of Theorem \rf{t3.1} with $K:=K'$, $a:=a'$
and \er{3.10}, we obtain $\wh\vf(P_a(c))= P_{a'}(\wh\vf(c))$. Thus $(\wh\vf
\circ P_a)(c) = (P_{a'}\circ\wh\vf)(c)$. Since $c\in K[X]$ is arbitrary,
\er{4.62} follows.
\endproof

\Wanp prove

\bth4.44
Let $E$ and $E'$ be finite fields of order $p^n$, $p$~prime, $n\ge2$. Then $E$
and~$E'$ are ring-\is c.
\eth

\proof
Let $\check E{}$ (resp.\ $\check E{}'$) denote the \pf\ of~$E$ (resp.~$
E'$). Since $n\ge2$, we have $\check E\ne E$ and $\check E{}'\ne E'$. Let
$a\in ME[X]$ ir\rd\ with $\deg(a)=n\ge2$ be defined by \er{4.47}, and let $K_a$
be the field introduced in Theorem \rf{t3.13} with $K:=\check E$. Let $\Psi_\a
:E \to K_a$ be the inverse of the ring-\is sm $\ov\Phi_\a$ introduced in
Lemma \rf{l4.41}.

Since the \ch istic of $E,\check E,E',\check E{}'$ is~$p$, \te s a map
$\ov\vf: (\N_p,+_p,\cdot_p,0,1) \to\check E$ which is a ring-\is sm by
Theorem \rfa4{t4.35}. Similarly, \te s a ring-\is sm $\ov\vf{}':(\N_p,+_p,
\cdot_p,0,1)\to\check E{}'$. Set
\beq4.65
\vf:=\ov\vf{}'\circ \ov\vf{}\Inv : \check E \to \check E{}'.
\e
Then $\vf$ is a ring-\is sm from $\check E$ onto~$\check E{}'$. Define
$\wh\vf: \check E[X]\to\check E{}'[X]$ by \er{4.57}. By Lemma
\rf{l4.43}\,(ii) $\wh\vf$~is a ring-\is sm. We claim that $\wh\vf$ maps $K_a$
into~$K_{a'}$ where $a':=\wh\vf(a)\in \check E{}'[X]$. Let $c\in K_a$.
Then $c=P_ac$ by \er{3.13}. Thus $\wh\vf(c)=\wh\vf(P_ac)\nde4.62 =
P_{a'}(\wh\vf(c))\in R(P_{a'}) \nde3.14 = K_{a'}$. We denote by $\wh\vf\rhu$
this \rt ion with range in~$K_{a'}$. Thus $\wh\vf\rhu:K_a\to K_{a'}$. Since
$\wh\vf$ is in\jc, so is $\wh\vf\rhu$. We now show that $\wh\vf\rhu$ is
sur\jc. Let $d'\in K_{a'}$. Since $\wh\vf$ is sur\jc, \te s $d\in \check
E[X]$ \st $\wh\vf(d)=d'$. Thus $P_{a'}(\wh\vf(d))=P_{a'}(d')=d'=\wh\vf(d)$.
Hence $\wh\vf(d)=P_{a'}(\wh\vf(d))\nde4.62 = \wh\vf(P_a(d))$. By the in\ji\
of~$\wh\vf$, we obtain $d=P_a(d)$, i.e.\ $d\in K_a$. Thus $\wh\vf\rhu$ is
bi\jc. Since both $K_a$ and~$K_{a'}$ are fields and since $\wh\vf$ is a
ring-\hm sm, so is $\wh\vf\rhu$. \If that $\wh\vf\rhu:K_a\to K_{a'}$ is a
ring-\is sm. \csq, $\wh\vf\rhu\circ\Psi_\a$ is a ring-\is sm from~$E$
onto~$K_{a'}$. By Lemma \rf{l4.43}\,(iii) $a'\in E'[X]$ is monic, ir\rd\ and
$\deg(a')=\deg(a)=n\ge2$.

We claim that \te s $\a'\in E'\sm \check E{}'$ \st $\suml_{i=0}^{\deg(a')}
a'(i) (\a')^i=0$. Let $b:=\ve_{p^n}-\ve_1 \in \check E[X]$. By Corollary
\rf{c4.36} \te s $q\in\check E[X]\sms{\bf0}$ \st
\beq4.66
b=q\cdot a.
\e
\If that $\ve'_{p^n}-\ve_1'\nde4.63 = \wh\vf(b)\nde4.66 = \wh\vf(q\cdot a) =
\wh\vf(q) \cdot' \wh\vf(a) = \wh\vf(q)\cdot' a'$.

Let $\Phi'_{\b'}:E'[X]\to E'$ denote the evaluation map defined in \er{3.28}
with $\a:=\b'\in E'$, $K:=E'$. Then $\Phi'_{\b'}(\ve'_{p^n}-\ve_1) \nde4.31 = \Phi'_{\b'}\bigl(
\mathop{\prod'}\limits_{\g'\in E'}(\ve'_1-\g'\ve'_0)\bigr)\nda2{1.131} =
\mathop{\prod'}\limits_{\g'\in E'}\bigl(\Phi'_{\b'}(\ve'_1-\g'\ve'_0)\bigr)
= \mathop{\prod'}\limits_{\g'\in E'\sms{\b'}}\bigl(\Phi'_{\b'}(\ve_1'-\g'\ve_0')\bigr)
\cdot' \Phi'_{\b'}(\ve'_1-\b'\ve_0')$ \fa $\b'\in E'$. Note that
$\Phi'_{\b'}(\ve'_1-\b'\ve_0')=0'$ \fa $\b'\in E'$, hence $\Phi'_{\b'}
(\ve'_{p^n}-\ve'_1)=0'$ \fa $\b'\in E'$. Thus $\wh\vf(q) \cdot' a'=
\ve'_{p^n}-\ve'_1$ has $p^n$~roots in~$E'$. Suppose for \cd ion that $a'$~has
\ti{no} roots in~$E'$, then $\wh\vf(q)$ should have $p^n$~roots in~$E'$ since
$u'\cdot' v'=0'$, $u',v'\in K[X]$ implies $u'=0'$ or $v'=0'$ by \er{2.29}.
However from \er{4.66}, $q\ne\bf0$, $a\ne\bf0$ and \er{2.29} we have $\deg(b)=
\deg(q)+\deg(a)$. Since $\deg(b)=p^n$ and $\deg(a)=n$, we have $\deg(q)<p^n$.
Thus $\deg(\wh\vf(q))\nde4.59 = \deg(q)<p^n$. Since $\wh\vf(q)$ possesses at
most $\deg(\wh\vf(q))$ roots in~$E'$ by \E\Pr\ \rf{p3.26}\,(iii), we obtain
a~\cd ion. \E\Tf $a'$~has at least one root $\a'\in E'$. Since $a'$~is ir\rd\
in $\check E{}'[X]$, $\a'\notin \check E{}'$, otherwise $\ve'_1-\a'\ve'_0$
would divide~$a'$ by Corollary \rf{c3.14}. \Wanp apply Lemma \rf{l4.42} with
$F:=E'$, $K:=\check E{}'$, $\a:=\a'$ and $\wh a:=a'$. We denote by
$\ov\Phi{}'_{\a'}$ the \rt ion to~$K_{a'}$ of the evaluation map
$\Phi'_{\a'}: E'[X]\to E'$ defined in \er{3.28} with $K:=E'$ and $\a:=\a'$.
\If from Lemma \rf{l4.42} that $\ov\Phi{}'_{\a'} :K_{\a'}\to E'$ is an in\jc\
ring-\hm sm. \If that the map $\Lambda:=\ov\Phi{}'_{\a'}\circ \wh\vf\rhu
\circ \Psi_\a:E \to E'$ is an in\jc\ ring-\hm sm. Since $\#(E)=p^n=\#(E')$,
$\La$~is bi\jc\ by \era2{3.11}, hence $\La$~is a ring-\is sm from~$E$\glossary{$\La(\a)$}
onto~$E'$.
\endproof

\brm4.45
Observe that in the proof of Theorem \rf{t4.44}, $\La(\a)=\ov\Phi{}'_{\a'}
\circ \wh\vf\rhu(\Psi_\a(\a)) = \ov\Phi{}'_{\a'}\circ \wh\vf\rhu(\ve_1)
= \ov\Phi{}'_{\a'}(\ve'_1)=\a'$.  Since $E^\t = \bcl_{i\in
[0,p^n-1]}\{\a^i\}$, and since $\La$~is a ring-\is sm, we obtain
$\La\bigl(\bcl_{i\in[0,p^n-1]}\{\a^i\}\bigr) = \bcl_{i\in[0,p^n-1]}
\{\La(\a^i)\} \nda21.45 = \bcl_{i\in[0,p^n-1]}\{(\La(\a))^i\} =
\bcl_{i\in[0,p^n-1]}(\a')^i$. But $\La\bigl(\bcl_{i\in[0,p^n-1]}\{\a^i\}
\bigr) = \La(E^\t)=(E')^\t$. Hence $\a'$ is a \Gn\ of the cyclic group
$((E')^\t, \cdot',1')$.
\erm

\brm4.46
Finite fields are usually called \ti{Galois fields}.\index{Galois fields}
\erm

We conclude this section by some exercises.

\bex4.47
Let $K$ be a field and let $MK[X]$ be the set of monic \fp s over~$K$. Show
that $(MK[X],\cdot,\ve_0)$ is not only a \PM, but also a L-monoid (see \E\df\
\rfa3{d8.34}). (\ti{Hint\/}: Imitate the proof of \E\Pr\ \rfa3{p1.16} and use
Remark \rf{r4.37}).
\eex

The aim of the next exercises is to show that a finite field~$E$ with prime
field~$\check E$ is ring-\is c to a field of matrices with entries in~$\check E$.

\bex4.48
Let $(E,+,\cdot,0,1)$ denote a finite field of order~$p^n$, $p$~prime, $n\ge2$.
Let $\check E$ denote its \pf\ of order~$p$. We recall that $E$~is a \vs\ over
the field $\check E$ of \dm~$n$. \Mo if $\a\in E^\t$ is a \Gn\ of the cyclic group
$(E^\t,\cdot,1)$, then $\{\a^i: i\in\zo0,n \}$ is a basis of $(\check E,E)$
by Lemma \rf{l4.33}. Given $\g\in E$ we define a map $\wh\g:E \to E$ by setting\glossary{$\wh\g(x)$}
\beq4.67
\wh\g(x):= \g\cdot x, \q x\in E.
\e
If we view $E$ as a \vs\ over $\check E$, the map $\wh\g$ is $\check E$-linear
(see \E\df\ \rf{d1.4}). Indeed, $\wh\g(x+y)=\g\cdot(x+y)=\g\cdot x+\g\cdot y
=\wh\g(x)+\wh\g(y)$, $x,y\in E$, and $\wh\g(\la x)=\g\cdot(\la x)=\g\cdot(\la\cdot x)
=(\g\cdot\la)\cdot x =(\la\cdot\g)\cdot x = \la\cdot(\g\cdot x)=\la\cdot\wh\g(x)=
\la\wh\g(x)$, $x\in E$, $\la\in\check E$. Let $(K,V)$ be a
\vsf0$K$, then the set of all linear maps from~$V$ into~$V$ is usually denoted\glossary{$\cL(V)$}
by~$\cL(V)$. Note that the \ti{null \Op} $0$ defined by $0(x):=0$, $x\in V$,
and the \ti{identity \Op} $I$ defined by $I(x):=x$, $x\in V$, are linear maps
from~$V$ into~$V$. The \ns\ $\cL(V)$ can be made into a \vs\ over~$K$ by
defining an \ad~$+$ and a \mlc\ by scalars in the \fw\ way:
\bga4.68
(S+T)(x):= S(x)+T(x), \q x\in V, \ S,T\in\cL(V),\\
(\la T)(x):=\la (T(x)), \q \la\in K,\ x\in V,\ T\in\cL(V).\lb{4.69}
\e

\hph i,xii, Show that $\cL(V)$ equipped with the \op s defined in \er{4.68}
and \er{4.69} is a \ti{\vs\ over~$K$} (see Lemma \rf{l2.21}).

If $S,T\in\cL(V)$, then the \cm\ $S\circ T:V\to V$
$$
V \stackrel S\lea V \stackrel T\lea V
$$
is linear by Lemma \rf{l1.5}\,(ii). \E\Tf we may define a binary \op~$\cdot$
on~$\cL(V)$ by setting
\beq4.70
S\cdot T:=S\circ T, \q S,T\in\cL(V).
\e
It turns out that the map $\cdot$ from $\cL(V)\t \cL(V)$ into $\cL(V)$ is bilinear
(see \E\df\ \rf{d2.17}).

\hph ii,xi, Show that the map $\cdot$ defined in \er{4.70} is \ti{bilinear}
(\Ip $(S_1+S_2)\cdot T = (S_1\cdot T)+(S_2\cdot T)$ and $(\la S)\cdot T=\la
(S\cdot T)$, $S_1,S_2,S,T \in \cL(V)$, $\la\in K$). Since the \cm\ of maps is
\asc e (see \era1{2.1}), so is the binary \op~$\cdot$. \Mo we have (see \era1{2.2})
\beq4.71
S\cdot I = I\cdot S = S, \q S\in\cL(V).
\e
\If that $(\cL(V),+,\cdot,0,I)$ is an \ti{\asc e} (in general non\cmt e)
\ti{$K$-algebra with unity}~$I$ (see \E\df\ \rf{d2.18}). If we ``forget'' the
\mlc\ by scalars, $\cL(V)$~is a \ti{non\cmt e ring with unity}.

If $T\in \cL(V)$ and $i\in\N$, we denote by $T^i$ the $i$-th \IT\ of~$T$ in
the monoid $(\cL(V),\cdot,I)$ (see \E\df\ and Notation \rfa2{d1.14}).
\E\Tf a \lc\ of \IT s of~$T$ is an expression of the form $\suml_{i=0}^N
\la_iT^i$ \fs $N\in\N$, $\la_i\in K$, $i\in[0,N]$, where $\suml_{i=0}^N$ is
the \cme sum in the \am\ $(\cL(V),+,0)$.

We set\glossary{$\cL_T(V)$}
\beq4.72
\cL_T(V):=\{\hbox{the set of all \lc s of \IT s of $T$ in }\cL(V)\}.
\e

\hph iii,x, Show that $\cL_T(V)$ is a \sbm\ of $(\cL(V),+,0)$ and of
$(\cL(V),\cdot,I)$. Show that if $\la\in K$ and $S\in\cL_T(V)$, then $\la S\in
\cL_T(V)$.

\hph iv,ii, Show that $(\cL_T(V),+,\cdot,0,I)$ is an \asc e, \ti{\cmt e}
$K$-algebra (resp.\ ring) with unity (use Lemma \rfa2{l1.21}).

We now return to the case $V:=E$ and $K:=\check E$, and define a map $\La_\a
:E\to \cL_{\hat\a}(E)$ by setting:
\beq4.73
\La_\a(\g):=\wh \g, \q \g\in E.
\e

\hph v,iii, Show that the \fw\ holds:
\bea4.74
&\La_\a(0)=0, \q\ \La_\a(1)=I,\\
&\La_\a(\g+\d)=\La_\a(\g)+\La_\a(\d),\q \g,\d\in E,\lb{4.75}\\
&\La_\a(\g\cdot\d) = \La_\a(\g)\cdot\La_\a(\d),\q \g,\d\in E, \lb{4.76}\\
&\La_\a(\la\g) = \la\La_\a(\g), \q \la\in\check E,\ \g\in E. \lb{4.77}
\e
From \er{4.74}, \er{4.76} and Lemma \rfa2{l1.22} we find
\beq4.78
\La_\a(\a^i)=\wh\a{}^i, \q i\in\N.
\e

\hph vi,ii, Show that $\{\wh\a{}^i\in\cL_{\hat\a}(E): i\in\zo0,n \}$ is a
\ti{basis} of the \vs\ $\cL_{\hat\a}(E)$ over~$\check E$, and that $\La_a$
is a $\check E$-\ti{linear \is sm from $E$ onto $\cL_{\hat\a}(V)$}.

\hph vii,i, Show that $\La_\a$ is an in\jc\ \ti{ring-\hm sm} from~$E$ into
$\cL(E)$. Show that $\cL_{\hat\a}(E)$ \ti{is a field\/} and that $\La_\a:E
\to\cL_{\hat\a}(E)$ is a \ti{ring-\is sm}.

So far we have obtained a \rp ation of the finite field~$E$ by a \ti{field of
linear maps} on the finite-\dm al \vs~$E$ over~$\check E$. We now use the
basis $\{\a^i: i\in\zo0,n \}$ of~$E$ in order to find a \rp ation of~$E$
by a field of \ti{matrices}. We suppose that the reader has some familiarity
with matrix calculus.\index{matrix calculus} However, we recall some basic facts. Let $(K,V)$ be
an~$n$-\dm al \vsf0$K$, let $\zb ej{[1,n]}$, $n\in\Na$, be a basis of~$V$,
and let $T$ be a linear map from $V$ into~$V$. In view of the linearity of~$T$,
the map~$T$ is completely determined by its \rt ion to the basis $\zb ej{[1,n]}$
(see Exercise \rf{ex2.13}). Thus by the \df\ of a basis \te s \fe pair $(i,j)\in
[1,n]\t[1,n]$ \ooo $t_{ij}\in K$ \st the \fw\ holds:
\beq4.79
Te_j = \sum_{i=1}^n t_{ij}e_i, \q j\in[1,n].
\e
Conversely, given $t_{ij}\in K$, $i,j\in[1,n]$, we can define a map $T_0:\zb ej{[1,n]}
\to V$ by means of \er{4.79}. Then \te s \ooo linear map $T:V\to V$ \sf ying
\er{4.79} (see Exercise \rf{ex2.13}). Usually, an object of the form $(t_{ij})$,
$i,j\in[1,n]\t[1,n]$, $n\in\Na$, $t_{ij}\in K$, is called an $n\t n$-matrix\index{square matrix}
(square matrix) with \el s or entries $t_{ij}$. For our purpose we use the \fw\ \df.\index{entries}
\eex

\begin{dfn}[\cite{Is}]\lb{d4.49}
Let $(R,+,\cdot,0,1)$ be a \cmt e ring with unity, and let $n\in\Na$. Let
$M_n(R)$ denote the set of all maps from $[1,n]\t[1,n]$ into~$R$. An \el\
$A\in M_n(R)$ is called an $n\t n$ \ti{matrix} with entries in~$R$, and is\index{matrix}
also denoted by $(a_{ij})$ or~$(A_{ij})$. We use the notation $A_{ij}:=A((i,j
))$, $i,j\in[1,n]$. The null-matrix $0_n$ is defined by $0_{n\,ij}:=0$, $i,j
\in[1,n]$, and the unit (or identity)-matrix $I_n$ by
$$
I_{n\,ij}:=\bca
1 &\hbox{for }i=j,\\
0 &\hbox{otherwise.}
\eca
$$
An \ti{\ad} $+$ on $M_n(R)$ is defined by
\beq4.80
(A+B)_{ij}:=A_{ij}+B_{ij}, \q\ i,j\in[1,n],
\e
and a \ti{\mlc} $\cdot$ on $M_n(R)$ is defined by
\beq4.81
(A\cdot B)_{ij} :=\sum_{k=1}^n A_{ik}B_{kj}, \q\ i,j\in[1,n].
\e
The product $A\cdot B$ is also denoted by juxtaposition $AB$.
\edf

\hph viii,, Show that $(M_n(R),+,\cdot,0_n,I_n)$ is a non\cmt e (unless $n=1$)
ring with unity.

Now let as above $(K,V)$ be an $n$-\dm al, $n\ge2$, \vsf0$K$ (which is of
course also a \cmt e ring with unity), let $\cL(V)$ denote the non\cmt e ring
of linear maps on~$V$, and let $\zb ej{[1,n]}$ be a basis of~$V$. From what
precedes \te s a bi\jn\ between $\cL(V)$ and $M_n(K)$. More precisely, given
$x\in V$ \te s \ooo $\la:[1,n]\to K$ \st $x=\suml_{j=1}^n \la_je_j$, since
$\zb ej{[1,n]}$ is a basis of~$V$. We denote by $\psi_j:V\to K$, $j\in[1,n]$,
the maps defined by
\beq4.82
\psi_j(x):=\la_j, \q j\in[1,n].
\e

\hph ix,ii, Show that $\psi_j:V\to K$ is linear \fe $j\in[1,n]$, and
\bga4.83
\psi_j(e_i)=\d_{ij}:=\bca
1 & \hbox{if }i=j,\\
0 & \hbox{otherwise,}
\eca \\
x=\sum_{j=1}^n \psi_j(x)e_j \q \hbox{\fa }x\in V. \lb{4.84}
\e
\medskip

We now define a map $\Psi:\cL(V)\to M_n(K)$ by setting
\beq4.85
(\Psi(T))_{ij} := \psi_i(Te_j),\q i,j\in[1,n],\ T\in\cL(V),
\e
and a map $\cX:M_n(K)\to \cL(V)$ by setting\glossary{$\cX(A)$}
\beq4.86
\cX(A)x := \sum_{i=1}^n \Bg(\sum_{j=1}^n A_{ij}{\psi_j(x)})e_i, \q x\in V,\
A\in M_n(K).
\e

\hph x,iii, Prove
\beq4.87
\cX \circ \Psi = \id_{\cL(V)},\q \Psi\circ\cX = \id_{M_n(K)}.
\e
Thus $\Psi$ is a bi\jn\ from $\cL(V)$ into $M_n(K)$ and $\cX$~is its inverse.

\hph xi,ii, Show that $\Psi$ is a ring-\is sm.

We now define a \mlc\ by scalars on $M_n(K)$ by setting
\beq4.88
(\la A)_{ij} := \la A_{ij}, \q i,j\in[1,n],\ A\in M_n(K).
\e

\hph xii,i, Show that $(M_n(K),+,0)$ equipped with the scalar \mlc\ defined in
\er{4.88} is a \vs\ over~$K$, and that $\Psi$ is a linear \is sm from $\cL(V)$
onto $M_n(K)$.

We denote the image of $\cL_{\hat \a}(E)$ under~$\Psi$ by $M_n\en\a(\check E)$\glossary{$M_n\en\a(\check E)$}
and of $\wh\a$ by~$\wh A$, thus
\beq4.89
M_n\en\a(\check E):=\Psi(\cL_{\hat\a}(E)); \q \wh A:=\Psi(\wh\a).
\e

\hph xiii,, Prove:

\lhp xiii,1, $M_n\en\a(\check E)$ is a subgroup of $(M_n(\check E),+,0_n)$.

\lhp xiii,2, $M_n\en\a(\check E)$ is a \lss\ of $(\check E,M_n(\check E))$.

\lhp xiii,3, $M_n\en\a(\check E)$ is a \sbm\ of $(M_n(\check E),\cdot,I_n)$.

Set
\beq4.90
\wh A{}^i := \hbox{the $i$-th \IT\ of $\wh A$ in }(M_n(\check E),\cdot,I_n),
\q i\in\N.
\e
Prove

\lhp xiii,4, $\{\wh A{}^i: i\in\zo0,n \}$ is a basis of $(\check E,
M_n\en\a(\check E))$.

\lhp xiii,5, $M_n\en\a(\check E)= \{0_n\} \cup \{\wh A{}^i:i\in[0,p^n-1]\}$.

\lhp xiii,6, If $B,C\in M_n\en\a(\check E)\sm\{0_n\}$, then $BC\in M_n\en\a(\check E)\sms{0_n}$.

\lhp xiii,7, $(M_n\en\a(\check E),+,\cdot,0_n,I_n)$ is a \ti{field}.

\lhp xiii,8, The map $\Psi\circ\La_\a :E\to M_n\en\a(\check E)$ is a ring-\is
sm from $(E,+,\cdot,0,1)$ onto $(M_n\en\a(\check E),+,\cdot,0_n,I_n)$, and is
a linear \is sm from $(\check E,E)$ onto $(\check E,M_n\en\a(\check E))$.

Thus we obtain a \rp ation of the finite field~$E$ by a \ti{field of matrices}.\index{field of matrices}

Let $p$~be a \Pn\ and $n\ge2$. By Theorem
\rf{t4.21} \te s a finite field~$E$ of order $p^n$ and by Corollary
\rfa4{c4.42}\,(ii) \te s a \Gn\ $\a$ of the cyclic group~$E^\t$. Let $\la_i\in
\check E$, $i\in[0,n-1]$, be as in \er{4.46}, and let $a\in\check E[X]$ be the
monic ir\rd\ \fp\ of degree~$n$ defined in \er{4.47}. Set $e_i:=\a^{i-1}$,
$i\in[1,n]$. Then $\zb ei{[1,n]}$ is a basis of $(\check E,E)$ by Lemma
\rf{l4.33}. \Mo $\wh\a\in\cL(E)$, defined in \er{4.67}, \sf ies
\bga4.91
\wh\a e_j = e_{j+1}, \q j\in[1,n-1],\\
\wh\a e_n = \sum_{k=1}^n \la_{k-1}e_k. \lb{4.92}
\e
By \er{4.89}, $\wh A=\Psi(\wh \a)$ and by \er{4.85}
\beq4.93
\wh A_{ij} = \psi_i (\wh\a e_j), \q i,j\in[1,n].
\e
Note that $\psi_i(\wh\a e_j)\nde4.91 = \psi_i(e_{j+1})=\d_{i\,j{+}1}$ for
$i\in[1,n]$, $j\in[1,n-1]$, and $\psi_i(\wh\a e_n)\nde4.92 = \psi_i\bg(
\suml_{k=1}^n \la_{k-1} e_k) = \suml_{k=1}^n \la_{k-1}\psi_i(e_k) =
\suml_{k=1}^n \la_{k-1}\d_{ik}= \la_{i-1}$. \E\Tf
\beq4.94
\wh A_{ij} = \bca
\d_{i\,j{+}1} & \hbox{for }j\in[1,n-1]\\
\la_{i-1} & \hbox{for }j:=n
\eca \q \ \hbox{for }i\in[1,n].
\e

\E\Ip if $n=3$,
$$
\wh A=\bmk \wh A_{11}\ & \wh A_{12}\ & \wh A_{13} \\
\wh A_{21}\ & \wh A_{22}\ & \wh A_{23} \\
\wh A_{31}\ & \wh A_{32}\ & \wh A_{33} \emk =
\bmk 0\ & 0\ & \la_0 \\
1\ & 0\ & \la_1 \\
0\ & 1\ & \la_2 \emk.
$$
We find
$$
\wh A{}^2= \bmk 0\ &\la_0\ &\la_0\cdot\la_2\\
0\ &\la_1\ & \la_0+(\la_1\cdot\la_2)\\
1\ &\la_2\ & \la_1+\la_2^2 \emk, \qquad
\wh A{}^3 = \bmk \la_0\ & \la_0\cdot\la_2 \ & \la_0\cdot\la_1+\la_0\cdot\la_2^2\\
\la_1\ &\la_0+\la_1\cdot\la_2\ &\la_0\cdot\la_2+\la_1^2+\la_1\cdot\la_2^2\\
\la_2\ &\la_1+\la_2^2\ &\la_0+\la_1\cdot\la_2+\la_2\cdot\la_1+\la_2^3\emk,
$$
where $+$ (resp. $\cdot$) denotes the \ad\ (resp. \mlc) in~$\check E$, and
$\la^i$ denotes the $i$-th \IT\ of~$\la$ in $(\check E,\cdot,1)$.

As it should be we find
$$
\wh A{}^3 = \la_0 \wh A{}^0 + \la_1\wh A{}^1 + \la_2\wh A{}^2.
$$

The matrix $\wh A$ defined by \er{4.93} is often called the \ti{companion
matrix} \cite{Galois} of the \fp\ $\suml_{k=0}^n a(k)X^k$ (see Exercise
\rf{ex2.28}), where $a(l)$, $l\in[0,n]$, is defined in \er{4.47}.\index{companion matrix}

We now consider as an example the field of order~$p^n$, $p=2=n$, introduced at
the end of Section~\ref{sss.fdvs}. We have $E:=\{{\bf0},{\bf1},\bu,\bv\}$
equipped with the \op s $+_a$ (resp.~$\cdot_a$) defined in Table \rf{tb3.18}
(resp.\ \rf{tb3.20}), and $\check E=\{0,1\}$. \Mo $\{{\bf1},\bu\}$, $\{{\bf1},
\bv\}$ and $\{\bu,\bv\}$ are bases of the \vs\ $(\check E,E)$ and $\bu$
(resp.~$\bv$) are \Gn s of the cyclic group $(\{{\bf1},\bu,\bv\},\cdot_a,
{\bf1})$. Set
\beq4.95
\a:=\bu.
\e
We have $\a^0:=0\ddtt a \bu={\bf1}$, $\a^1:=1\ddtt a\bu=\bu$, and $\a^2:=2
\ddtt a\bu=\bv$. Then $\{\a^0,\a^1\}$ is a basis of $(\check E,E)$ and $\a^2
=\bv={\bf1}+\bu=\a^0+\a^1$. Thus in \er{4.46} we have $\la_0=1$ and $\la_1=1$.
From \er{4.94} we find that
$$
\wh A = \bmk 0 && \la_0 \\ 1 && \la_1 \emk = \bmk 0 && 1 \\ 1 && 1 \emk.
$$
\E\Tf $M_2\en \a(\check E) = \{0_2,I_2,\wh A,\wh A{}^2\} = \{
\left[\begin{smallmatrix} 0 & 0 \\ 0 & 0 \end{smallmatrix}\right],
\left[\begin{smallmatrix} 1 & 0 \\ 0 & 1 \end{smallmatrix}\right],
\left[\begin{smallmatrix} 0 & 1 \\ 1 & 1 \end{smallmatrix}\right],
\left[\begin{smallmatrix} 1 & 1 \\ 1 & 0 \end{smallmatrix}\right]\}$, and
as a subring of $M_2(\check E)$, $M_2\en\a(\check E)$ is a \ti{field\/}.
\Mo the map $\ov\Psi_\a:E \to M_2\en\a(\check E)$ defined by $\ov\Psi_\a({\bf0})
:=0_2$, $\ov\Psi_\a({\bf1}):=I_2$, $\ov\Psi_\a(\bu):=\wh A$ and $\ov\Psi_\a(\bv)
:=\wh A{}^2$ is a ring-\is sm.

\bex4.49a
Show that if we set $\b:=\bv$ instead of $\a:=\bu$ in \er{4.95} we obtain $M_2
\en\a(E)=M_2\en\b(E)$, hence $\ov\Psi{}_\b\Inv\circ\ov\Psi_\a$ is a ring-auto\mf
\ of the field~$E$. Investigate the more general situation when $\b\ne\a$ is
a root of $a\in\check E[X]$ defined in \er{4.47}.
\eex

We conclude this section by considering matrix \rp ations of the field~$K_a$\index{matrix representation}
introduced in Theorem \rf{t3.13}. We need some \df s. Let $K$~be a field and\glossary{$(\cA,+,\cdot,0,I)$}
let $(\cA,+,\cdot,0,I)$ denote an \asc e $K$-algebra with unity~$I$ which is
not necessarily \cmt e. Let $T\in\cA\sms0$. We denote by~$T^i$, $i\in\N$, the
$i$-th \IT\ of~$T$ in the monoid $(\cA,\cdot,I)$, and by $\<T>$ the set of all\glossary{$\<T>$}
\lc s of \IT s of~$T$.

\bex4.50 \

\hph i,i, Show that $\<T>$ is a \sbm\ of $(\cA,+,0)$ and an abelian \sbm\ of
$(\cA,\cdot,I)$.

\hph ii,, Show that $(\<T>,+,\cdot,0,I)$ is a \cmt e ring with unity.
\eex

Observe that if $K$ is a field, then $M_n(K)$, $n\in\Na$, equipped with the\glossary{$M_n(K)$}
\op s defined in \er{4.80}, \er{4.81} and \er{4.88} is an $n^2$-\dm al \asc e
$K$-algebra with unity and with basis $\{\ve\en{i,j}: i,j\in[1,n]\}$ where
\beq4.96
\ve\en{i,j}_{kl} := \bca
1 & \hbox{if }(k,l)=(i,j),\\
0 & \hbox{otherwise.}
\eca
\e

We showed that if $E$ is a finite field of order~$p^n$, $\check E$~is its \pf,
$\a$~is a \Gn\ of~$E^\t$ and $a\in\check E[X]$ is the monic ir\rd\ \fp\
introduced in \er{4.47}, and if $\wh A\in M_n(\check E)$ is the companion
matrix of~$a$, then $\langle\wh A\rangle$ viewed as a ring is a field \is c
to~$E$.

We now consider the \fw\ \gn. Let $K$ be a field and let $a\in K[X]$ be monic
and ir\rd. Set $n:=\deg(a)$. Let $C(a)$ denote the companion matrix of~$a$\glossary{$C(a)_{ij}$}
defined by
\beq4.97
C(a)_{ij}:=\bca
\d_{i\,j{+}1} & \hbox{for } j\in[1,n-1]\\
-a_{i-1}  & \hbox{for } j:=n
\eca \q \ \hbox{for }i\in[1,n].
\e
Let $(K_a,+_a,\cdot_a,{\bf0},\ve_0)$ denote the field introduced in Theorem
\rf{t3.13}, and let $(K,K_a)$ denote the $n$-\dm al \vs\ over~$K$ with basis
$\{\ve_i:i\in\zo0,n \}$. We define $T\in\cL(K_a)$ by setting
\beq4.98
Tx:= \ve_1\cdot_a x,\q x\in K_a,
\e
and $e_i\in K_a$ by setting
\beq4.99
e_i:=\ve_{i-1}, \q i\in[1,n].
\e
We find the analogue of \er{4.79}:
\beq4.100
Te_j = \sum_{i=1}^n C(a)_{ij}e_i, \q j\in[1,n].
\e

Proceeding as above we obtain

\bpr4.51
The ring $\<C(a)>$ in $M_n(K)$ is a field \is c to the field\break
$(K_a,+_a,\cdot_a,{\bf0},\ve_0)$.\glossary{$\<C(a)>$}
\epr

\bex4.52
Prove \E\Pr\ \rf{p4.51}.
\eex

\brm4.53
If the field $K$ in \E\Pr\ \rf{p4.51} is a field of matrices over a field~$K_0$,
then the field $K_a$ can be \rp ed by matrices made of blocks of matrices
with entries in~$K_0$. 
\erm

\bex4.54
Find a matrix \rp ation for the fields $(K_a,+_a,\cdot_a,{\bf0},\ve_0)$ where
$K:=\Q$, and $a(x):=x^2+1$, $a(x):=x^n-p$, $p$~prime, $n\ge2$, $x\in\Q$ (see
Remark \rf{r3.42}).
\eex

%% file: app1.tex
\Section{Appendix. An introduction to real Analysis}[Appendix]\label{as.1}
\Subsubsection{Introduction}\label{ass.1}

The aim of this Appendix is to introduce the field~$\R$ of real numbers, the cornerstone of Real Analysis. The field~$\R$ is ``maximal'' among \Ar\
\of s since every \Ar\ \of\ can be embedded, in an appropriate sense, in the field~$\R$. We construct the field~$\R$ by means of Dedekind cuts of
$\Q_{>0}$ rather than~$\Q$. The \fd\ \pp y of~$\R$ is its order-\cp ness. Based on this \pp y, $n$-th roots and \ra\ powers of \po\ real \nm s are
\es ed. Using convexity arguments and the ``approximability'' of real \nm s by dyadic ones, we extend the map $r\mt \a^r$ with $\a\in\R_{>0}$,
$r\in\Q$, to~$\R$. It turns out that for $\a\in\R_{>1}$ this \ext\ $x\mt \a^x$, $x\in\R$, is an order- and group-\is sm from the additive group
of~$\R$ onto the \po\ subgroup of its \mlv\ group. Logarithms to base~$\a$ are defined as inverse of these \is sms. In the last section the notion of
\dv\ is \itd\ and the \dfb ility of the map $x\mt \a^x$, $x\in\R$, and of its inverse are \es ed. It turns out that \te s \ooo $\ov\a\in\R_{>0}$ \st
the \dv\ of the map $x\mt\ov\a{}^x$ is equal to the map itself. This \nm\ is usually denoted by~$\ee$, and is sometimes called the Euler \nm. We show
that $\ee-1$ is the supremum of the bounded \sq\ of \ra\ \nm s $n\mt s_n:=\suml_{k=1}^n \frac1{k!}$ with $k!=\prodl_{l=1}^kl$, $k\in\Na$ and
$n\in\Na$. A~\df, sketch and reference of the proof of transcendence of~$\ee$ is given at the end of Section~\ref{ass.9}. This \pp y allows us to
define a field $\ee^\R$, order- and ring-\is c to~$\R$, \st the \mlv\ group of the \pf\ of~$\ee^\R$ consists of \tc\ (\Ip ir\ra) \nm s.

The notions of metric (resp.\ \cp\ metric) spaces are \itd\ in Section~\ref{ass.7} (resp.\ Section~\ref{ass.8}). It turns out that \fe \Ar\ field~$F$
a metric on~$F$ can be defined \st the \ms~$F$ is metrically \cp\ iff the field~$F$ is order-\cp.

In Section \ref{ass.8} some basic topological notions and related theorems are considered. \E\Ip it is shown that an \Ar\ field endowed with its
order-topology is order-\cp\ iff it is \cnt ed.

The author would like to thank Dr.\ Jan Kowalski for his excellent typing and remarks. Without the patience and support of my family, this Appendix
would never have been written.

\ssk
Bourguillon, Eygelshoven and Pavia 2025.

\rightline{Ph.\ Cl\'ement\indent}

\newpage
\Subsubsection{Embeddings and fields of \qt s}\label{ass.2}

We first recall some \df s. For the \df s of monoids (abelian, \Cm s, \hbox{\PM s}, (proper) \sbm s), we refer the reader to \E\df s 1.2.2 and 2.1.4.
In this book, a set~$X$ containing two distinct \el s $0$~and~$1$, and having two binary \op s: \ad\ and \mlc, denoted by~$+$ and by~$\cdot$
(or simply by juxtaposition) \sf ying the \fw\ \as s
\bea 2.1
&(X,+,0) \hbox{ is an \am,} \\
&(X,\cdot,1) \hbox{ is an \am,} \lb{2.2} \\
&0\cdot x = 0 = x\cdot 0 \qh{\fa}x\in X, \lb{2.3} \\
&x(y+z) = (xy) + (xz) \qh{\fa }x,y,z\in X, \lb{2.4} \\
&(x+y)z = (xz)+(yz) \qh{\fa }x,y,z \in X, \lb{2.5}
\e
is called a \ti{\sr}.

See \E\df\ 4.1.1 and Convention 4.1.5.
An important example is the \sr\ $(\N,+,\cdot,0,1)$.
Clearly, there is some redundancy in axioms \er{2.3}, \er{2.4}, \er{2.5}, but we prefer the \sy y of this formulation. Observe that $\Na$~is
a \sbm\ of $(\N,\cdot,1)$ by \E\Pr\ 2.2.10, and that the monoid $(\Na,\cdot,1)$ is a \PM\ (in particular a \hbox{\Cm}). A~\sr~$X$ \st its additive
monoid $(X,+,0)$ is an (abelian) group (see \E\df\ 4.3.2) is called a (\cmt e) \ti{ring} (with unity), see \E\df\ 4.3.11\.(i). \If from Theorem
4.5.8\,(iii) that $\Zf$ is a ring. 

The ring $\Zf$ possesses an important \pp y. By Theorem 4.5.8, $\Z\sms0$ is a \sbm\ of the monoid $(\Z,\cdot,1)$, thus if $x\in\Z\sms0$, there is
no $y\in\Z\sms0$ \st $x\cdot y=0$. A~nonzero \el\ of a ring~$R$ is called a \ti{\zd} if \te s a nonzero \el\ $y$ of~$R$ \st $x\cdot y=0$ (see \E\df\
4.5.28). A~ring~$R$ \st no \el\ of $R\sms0$ is a \zd\ is called an (integral) \ti{domain} (see \E\df\ 4.5.30). It turns out (see page~253) that if
$R$~is a domain then $(R\sms0,\cdot,1)$ is a \Cm. \E\Ip $(\Z\sms0,\cdot,1)$ is a \Cm. An \el~$x$ of a ring~$R$ \st \te s $y\in R$ \sf ying
$x\cdot y=1$ is called a \ti{unit\/} of the ring~$R$ (see \E\df\ 4.3.34). Since $1\cdot1=1$, the unity~$1$ is a unit of the ring~$R$. Clearly, $0$~is
not a unit in view of \er{2.3}. A~unit of~$R$ is an invertible \el\ of the monoid $(R,\cdot,1)$ (see \E\df\ 4.3.2). The set of invertible \el s of
a monoid $(M,\qu,e)$, denoted by~$M^*$ or $(M,\qu,e)^*$ (nonstandard notation), is a \sbm\ of the monoid $(M,\qu,e)$, and $(M^*,\qu,e)$ is a group
(see Lemma 4.3.6 and \E\df\ 4.3.2). A~\ti{field\/} $F$ is a ring \st every nonzero \el\ is a unit (see \E\df\ 4.4.1).

Note that a unit in a ring $R$ is \ti{not\/} a \zd. Indeed, if $x$ is a unit of~$R$, then $x\cdot y=1$ \fs $y\in R$ (then $y$ is also a~unit), and
there is no $z\in R\sms0$ \st $x\cdot z=0$, since we would have $z=1\cdot z=(x\cdot y)\cdot z = (y\cdot x)\cdot z = y\cdot(x\cdot z) = y\cdot 0 =0$.
\If that if $F$~is a field, then $(F\sms0,\cdot,1)$ is a group. The domain $\Zf$ is a domain but not a field. Indeed, if $x,y\in\Z\sms0$ \sf y
$x\cdot y=1$, then $1\nad{(4.5.36)}= |x\cdot y|\nad{(4.5.40)}= |x|\cdot|y|$ where $|x|,|y|$ belong to~$\Na$. Then $|x|=|y|=1$ by \E\Pr\ 2.2.9\,(iv).
Hence $x\in\{1,-1\}\ne \Z\sms0$. A~finite \Cm\ is a group by \E\Pr\ 4.3.10. Hence a \ti{finite} domain is a field.

\ssk
It turns out that an infinite \Cm\ can be embedded in an \ag\ (see Theorem 4.3.84). Mainzer claims in \cite[p.~21]{Nrs}: ``\ti{Historically, it was
also Dedekind who introduced the idea of defining \ig s by pairs from $\N\t\N$}.'' Following \cite[p.~22]{Nrs}, Steinitz~\cite{Stein} showed that
a~domain can be embedded in a field by formation of fractions (see Theorem 4.5.35).

The notion of embedding needs some explanation. We first consider the case of the embedding of a monoid in an \ag. The simplest situation is when a
monoid $M$ is a \ti{subset\/} of an \ag~$G$. Then the binary \op\ of~$M$ has to be the \rt ion to~$M$ of the binary \op\ of~$G$ and the \nel\ of~$M$
has to be the \nel\ of~$G$. Since the group~$G$ is a monoid, the monoid~$M$ has to be a \ti{\sbm\/} of~$G$ (see \E\df s 1.2.2 and the discussion
preceding Example 1.2.3). \Mo the group~$G$ being a \Cm, the monoid~$M$ has to be a~\Cm. \E\Tf a \Cm~$M$, contained in an \ag~$G$, is said to be
embedded in~$G$ if the inclusion map is a monoid-\hm sm.

We now consider the case where the \Cm~$M$ is \ti{not\/} a~subset of a group~$G$. Roughly speaking, a~\Cm~$M$ can be embedded in a~group~$G$, if $M$
can be ``transported'' into a \sbm\ of~$G$ embedded in~$G$. This
``transport'' is obtained by means of an in\jc\ map $\vf$ from~$M$ into~$G$ that \ti{preserves} the monoid \sc\ of~$M$. This is achieved by requiring
that $\vf:(M,\qu,e) \to (G,\hqu,\wh e)$ \sf ies $\vf(e)=\wh e$ and $\vf(x\qu y) = \vf(x)\hqu \vf(y)$ \fa $x,y\in M$. Let $\wh M:= \vf(M)$, $\wh x
:=\vf(x)$, $\wh y:=\vf(y)$. Then $\wh x\hqu \wh y = \vf(x)\hqu\vf(y) =\break \vf(x\qu y)\in \vf(M) = \wh M$. Similarly, we have $\wh x\hqu\wh y =
\wh y\hqu \wh x$, $\wh e\hqu \wh x=\wh x$ \fa $\wh x,\wh y\in \wh M$. \E\Tf $(\wh M,\hqu,\wh e)$ is an \am. Since an \ag\ is a \Cm, we find that
$\wh M$~is a \Cm. Observe that the map $\vf:M\to \wh M$ is a \ti{monoid-\is sm} from $M$ onto~$\wh M$, hence $M$~is \ti{monoid-\is c} (see \E\df\
2.1.7) to the \ti{monoid\/}~$\wh M$.

We now show that the monoid~$\wh M$ is embedded in the group~$G$, that is, the inclusion map from~$\wh M$ into~$G$ denoted by~$i$ is a monoid-\hm sm.
Indeed, $i(\wh e)=\wh e$ and $i(\wh x \hqu \wh y) = \wh x \hqu \wh y = i(\wh x)\hqu i(\wh y)$ \fa $\wh x,\wh y\in\wh M$. Since the \cm\ of two in\jc\
\hm sms is an in\jc\ \hm sm, the map~$\vf$ is an in\jc\ \hm sm from~$M$ into~$G$. Then the monoid~$M$ is said to be embedded in the group~$G$,
In other words, \te s an \ti{in\jc\ monoid-\hm sm into the group~$G$}.

We recall that the \ex\ part of Theorem 4.3.84 states that if $M$~is an
infinite \Cm, there \ti{exists} an \ag~$\wh M$, and an in\jc\ monoid-\hm sm $j:M\to\wh M$. Similarly, the \ex\ part of Theorem 4.5.35 states that
if $D$~is a domain, \te s a \ti{field\/}~$\wh D$ (called the field of \qt s of~$D$), and an in\jc\ ring-\hm sm (see \E\df\ 4.1.7) $j:D \to\wh D$.
An embedding of a field~$K$ in a field~$F$ is also called a \ti{field \ext} of the field~$F$
(see \E\df\ 5.3.23). We recall that a ring-\hm sm from a field~$F_1$ into a field~$F_2$ is in\jc\ (see
Lemma 5.3.24). We conclude this section by mentioning an important example of a field of \qt s. Let $K$~be a field (finite or infinite). Then the
ring of \fp s $K[X]$ (see \E\df\ 5.2.25 and Remark 5.2.27) is a domain (see Remark 5.2.32). The \crs\ field of \qt s is sometimes denoted by $K(X)$
(see \cite[p.~234]{Alg}).

\newpage
\Subsubsection{\E\Ar\ fields}\label{ass.3}
Let $\pz0F$ be a field of \ti{\ch istic zero}. By \E\df\ 4.4.37, $n\dpl1\ne0$ \fa $n\in\Na$, where $n\dpl1$ denotes the $n$-th \IT\ of~$1$ in the
monoid $\pz1F$. Let $\vf:\N\to\pz1F$ be defined by $\vf(n):=n\dpl1$, $n\in\N$, then $\vf$~is \ti{in\jc}. Indeed, suppose for \cd ion that \te\
$m,n\in\N$ with $m\ne n$ \st $\vf(m)=\vf(n)$. \E\wlg we may suppose $n<m$, that is, \te s $k\in\Na$ \st $m=n+k$. Then $\vf(m)=\vf(n+k)
\nad{(2.2.3)\,I2}= \vf(n)+\vf(k) = \vf(m)+\vf(k)$, hence $0+\vf(m) = \vf(m)= \vf(k)+\vf(m)$. By \cnc ity, $\vf(k)=0$. Hence $k\dpl1 = \vf(k)=0$ with
$k\ne0$. A~\cd ion. $\N$~being infinite by Theorem 1.4.18\,(iv), we find that $F$~is \ti{infinite} by \E\Pr\ 1.4.27\,(i) since $\vf$~is in\jc.
We conclude that a field of \ch istic zero is infinite. The converse is not true. Indeed, for every prime~$p$, \te s an infinite field of
\ch istic~$p$ (see Remark 5.2.32). For convenience, we give a proof. Let $n\in\N$, $n\ge2$, and let $(\N_n,+_n,\cdot_n,0,1)$ denote the semiring
introduced in \E\Pr\ 4.1.18\,(i). The monoid $(\N_n,+_n,0)$ is a \pn\ \Cm. Since $\N_n$ is finite, $(\N_n,+_n,0)$ is an \ag\ by \E\Pr\ 4.3.10,
hence $(\N_n,+_n,\cdot_n,0,1)$ is a \ti{ring} (see \E\df\ 4.3.11). Let $p$ be a \Pn\ and let $k\in\N_p\sms0$, that is, $k\in\N$ and $0\le k<p$.
Then $\gcd(k,p)=1$, and $k$~is a \ti{unit\/} of the ring $(\N_p,+_p,\cdot_p,0,1)$ by \E\Pr\ 4.3.44. Hence $\N_p$ is a \ti{field\/} in view of
\E\df\ 4.4.1. We now show that $\N_p$ is a \Pf. Suppose that $H$~is a subfield of~$\N_p$. Then $\#(H)\le \#(\N_p)=p$. By \E\df\ 4.4.24, $(H,+_p,0)$
is a subgroup of the finite group $(\N_p,+_p,0)$, and $1\in H$. \E\Tf $(H,+_p,0)$ is a nontrivial subgroup of $(\N_p,+_p,0)$, hence by Lagrange's
theorem 4.3.61, $\#(H)$~\ti{divides} $\#(\N_p)=p$. Thus $\#(H)=1$ or~$p$. Since $0,1\in H$, we have $\#(H)>1$, hence $\#(H)=p$, and $H=\N_p$ by
Lemma 2.3.16. Thus $\N_p$~is a \ti{\Pf\/}. Now, let $\N_p[X]$ denote the ring of \fp s over $\N_p$ (see \E\df\ 5.2.25 and Remark 5.2.27). Let
$\ve_k\in K[\N_p]$, $k\in\N$, defined in (5.2.4). Clearly, $\ve_k\ne \ve_l$ if $k,l\in\N$, $k\ne l$. Hence the map $i:\N \to K[\N_p]$ is in\jc.
\If that $K[\N_p]$ is infinite. \Mo \fa $a,b\in K[\N_p]\sms0$ we have $a\cdot b\ne0$ by Theorem 5.2.24. Hence the ring $K[\N_p]$ is a \ti{domain}.
By Theorem 4.5.35 \te s a field $\N_p(X)$, the field of fractions of the domain $\N_p[X]$, and an in\jc\ field \hm sm $j:\N_p[X]\to \N_p(X)$.
Hence the field $\N_p(X)$ is \ti{infinite}. We now show that the field of \qt s $\N_p(X)$ has \ch istic~$p$. Let $j$ denote the in\jc\ \hm sm
from $\N_p[X]$ into $\N_p(X)$, and let $\ve_0$ denote the identity of the \mlv\ monoid of $\N_p[X]$. Then $j(\ve_0)$ is the identity of the \mlv\
monoid of $\N_p(X)$ since $j$~is a ring-\hm sm. We have to show that $n\dpl j(\ve_0)\ne0$ for $0\le n<p$ and $p\dpl j(\ve_0)=0$ where $m\dpl
j(\ve_0)$ denotes the $m$-th \IT\ of $j(\ve_0)$ in the additive monoid of $\N_p(X)$. In view of Lemma 2.1.22 we have for $m\in\N$:
$m\dpl j(\ve_0) = j(m\dpl \ve_0)$, where $m\dpl\ve_0$ is the $m$-th \IT\ of~$\ve_0$ in $(\N_p[X],+,0)$. Let $\psi_i$, $i\in\N$, denote the linear map
from $\N_p[X]$ into the field~$\N_p$ introduced in (5.2.2). \E\Ip $\psi_i$~is a monoid-\hm sm from $(N_p[X],+,0)$ into $(\N_p,+_p,0)$. In view of
(5.2.1), (5.2.2), $a=0$ in $\N_p[X]$ iff $\psi_i(a)=0$ \fa $t\in\N$. Note that $\psi_0(\ve_0)=1$ and $\psi_i(\ve_0)=0$ \fa $i\in\Na$. By Lemma
2.1.22 again, $\psi_i(m\dpl\ve_0) = m\dpls p \psi_i(\ve_0)$ \fa $m,i\in \N$. If $i\in\Na$, then $\psi_i(\ve_0)=0$, hence $\psi_i(m\dpl \ve_0) =
m\dpls p 0= 0$ in $(\N_p,+_p,0)$ \fa $m\in\N$ by (2.2.3)\,I4. \E\oh if $i:=0$, we have $\psi_0(m\dpl \ve_0) = m\dpls p \psi_0(\ve_0) =
m\dpls p1$. By Theorem 4.1.28\,(i) with $(X,\qu_0,\qu_1,e_0,e_1):=(\N_p,+_p,\cdot_p,0,1)$ we obtain $m\dpls p1=m$ \fa $m\in\zo0,p $. \E\Ip
$m\dpls p1\ne0$ \fa $m\in(0,p)$, hence
$j(m\dpls p1)\ne0$, $m\in(0,p)$, by in\ji\ of~$j$. Finally, if $m:=p$, then we have $p\dpls p1=0$ by Theorem 4.4.36 with $F=\check F=
(\N_p,+_p,\cdot_p,0,1)$, $\vf(k) := k\dpls p1$, $k\in\N$, $e_0:=0$, $e_1:=1$, $I(1)=[0,p)$ and $\#(I(1))=p$.

\If that $j(p\dpls p e_0)=0$. \E\Tf the field of \qt s $\N_p(X)$ has \ch istic~$p$. Finally, it follows from Theorem 4.4.36\,(iii) that the \Pf\
of a field of \ch istic~$p$ is finite, \ev tly, by contraposition, an infinite \Pf\ has \ch istic zero. \csq, we obtain

\bpr 3.1
Let $F$ be a \ti{\Pf\/}. Then the \fw\ assertions are \ev t\/{\rm:}
\bea 3.1
&\hskip-100pt \hbox{$F$ is \is c to the field $\Q$,\hskip100pt}\\
&\hskip-100pt \hbox{$F$ is infinite,} \lb{3.2} \\
&\hskip-100pt \hbox{$F$ has \ch istic $0$}.\lb{3.3}
\e
\epr

\proof
The \ev ce $\er{3.2}\iff\er{3.3}$ follows from the discussion preceding the \Pr. The implication $\er{3.2}\imp\er{3.1}$ is a con\sq\ of Theorem
4.5.44. The proof of $\er{3.1}\imp\er{3.2}$ goes as follows.
Since $\Na$~is infinite and the map $i:\Na\to\Q$ defined by $i(n):=\frac n1$ is in\jc, $\Q$~is infinite by \E\Pr\ 1.4.27. If $F$~is \is c to~$\Q$,
then \te s a bi\jn\ from $\Q$ onto~$F$, and the conclusion follows from the same \Pr.
\endproof

We recall (see Example 2.1.5\,(v)) that if $(M,\qu,e)$ is a monoid, $M'$~is a set \ep\ to~$M$ (see \E\df\ 1.4.1), and $f:M\to M'$ is a bi\jn, then
\te s a monoid \sc\ on~$M'$ defined by (2.1.11) which makes~$f$ a monoid-\is sm from~$M$ onto~$M'$ (see Remark 2.1.6\,(iii)). Note that the identity
$\id_M$ in the monoid~$M$ is a monoid-\is sm (see \E\df\ 2.1.7 and Lemma 2.1.8). If $\id_M$ is the only auto\mf\ of~$M$, and $f_1,f_2$ are \is sms
from~$M$ onto a monoid $M''$, then $f_1=f_2$.

Indeed, $f\Inv_2\circ f_1$ is an auto\mf\ of~$M$, hence $f_2\Inv\circ f_1 = \id_M$, and $f_2 = f_2\circ\id_M = f_2\circ(f_2\Inv\circ f_1) =
(f_2\circ f_2\Inv)\circ f_1 = \id_M\circ f_1=f_1$. The same result holds for a field instead of a monoid, where (2.1.11) is replaced by
(4.4.1)--(4.4.3) (see \E\Pr\ 4.4.2).

\bpr 3.2
The identity is the only ring-auto\mf\ of the field~$\Q$. \E\Ip the \is sm $\th:\Q\to\wh\Z$ introduced in \E\Pr\ {\rm4.5.40} is the only \is sm
from~$\Q$ onto~$\wh\Z$.
\epr

\proof
Let $\vf$ be an auto\mf\ of~$\Q$. We want to prove
\beq 3.4
\vf(r)= r
\e
\fa $r\in\Q$. By \E\df s 4.1.7 and 2.1.7, $\vf$ is a bi\jc\ endo\mf\ of the monoid $\pz1{\Q}$. By \E\df\ 4.5.7 and Theorem 4.5.6, $\Q:=\wh X$ is the
disjoint union of $\wh X_+$, $\{0\}$ and~$\wh X_-$, where $\wh X_-$ is by (4.5.17) the set of inverses in the group $\pz1{\wh X}$ of the \el s
of~$\wh X$. We have $\vf(0)=0$ by (2.1.15). We claim that \er{3.4} holds \fa $r\in\wh X_+$ if \er{3.4} holds \fa $r\in\wh X_+$. Indeed, let $r\in
\wh X_-$. By \df, \te s $s\in \wh X^+$ \st $r = s\Inv = -s$. We have $s+s\Inv=0$, hence $\vf(s)+\vf(s\Inv) = \vf(s+ s\Inv) = \vf(0)=0$. \If that
$\vf(r)=\vf(s\Inv) = \vf(s)\Inv = -\vf(s)$. If $\vf(s)=s$ \fa $s\in \wh X_+$, then $\vf(r)=-\vf(s)= -s
= s\Inv =r$. Since $r$~is arbitrary in~$\wh X_-$, the claim is proved. We now show that \er{3.4} holds
\fa $r\in\wh X_+$. In view of the discussion preceding Theorem 4.5.6, an \el\ $r\in \wh X_+$ is of the form $+q$ where $q$~is a fraction $\frac ab$
with $a,b\in\Na$ defined in (4.5.1). If $r:=+\frac ab$, $s:=+\frac cd$, $a,b,c,d\in\Na$, then $r\cdot s$, the product of $r$ and~$s$ in
$\bg(\wh X,\cdot,+\frac 11)$, the \mlv\ monoid of the field $\bg(\wh X,+,\cdot,0,+\frac11)$ is equal to $+\frac{ac}{bd}$ by~(4.5.9) where $ac$
(resp.~$bd$) is the product of~$a$ (resp.~$b$) and $c$ (resp.~$d$) in $\pz2{\Na}$. By \E\df s 4.1.7 and 2.1.7, we have $\vf(r\cdot s) = \vf(r)\cdot
\vf(s)$. \If that if $m,n\in\Na$, we have $+\frac mn = \bigl(+\frac m1\bigr)\cdot\bg(+\frac 1n)$, hence
\beq 3.5
\vf\Bgg(+\frac nm) = \vf\Bgg(+\frac m1)\cdot\vf\Bgg(+\frac 1n).
\e
We claim that if $\vf\bg(+\frac m1)=+\frac m1$ \fa $m\in\Na$, then $\vf\bg(+\frac1n) = +\frac 1n$ \fa $n\in\Na$. We assume $\vf\bg(+\frac m1)
=+\frac m1$ \fa $m\in\Na$, and show that $\vf\bg(+\frac1n)=+\frac1n$ \fa $n\in\Na$. Note that if $n\in\Na$, then we have $\bg(+\frac n1)\cdot
\bg(+\frac 1n) = +\frac nn = +\frac11$. \E\Tf $+\frac1n$ is the inverse of $+\frac n1$ in the group $(\wh X\sms0,\cdot,+\frac11)$. \If that
$\vf\bg(+\frac n1)\cdot\vf\bg(+\frac1n) = \vf\bg(+\frac11) = +\frac11$ by \E\df s 4.1.7 and 2.1.7. Hence, by our \as, $\vf\bg(+\frac n1)\cdot
\vf\bg(+\frac1n) = +\frac11$. Since $\bg(+\frac n1)\cdot\bg(+\frac1n) = + \frac11$, we obtain by \cnc ity $\vf\bg(+\frac1n) = +\frac1n$. This proves
the claim. We now prove \er{3.4} for $r:=+\frac m1$, $m\in\Na$. We recall that $\Q_{\ge0}:=\Q_{>0}\cup\{0\}$ by (4.5.3), and $\pz1{\Q_{\ge0}}$ is
a \PM, where $+$ is defined by (4.5.11) and (4.5.12), in view of \E\Pr\ 4.5.4. Since the \nog\ of the \PM\ $\pz1{\Q_{\ge0}}$ is total, Corollary
4.3.86\,(i) can be applied and the in\jc\ \hm sm $j:\pz1{\Q_{>0}} \to \bg(\wh X,\cdot,+\frac11)$ introduced in Theorem 4.3.84 \sf ies $j\bg(\frac ab)
=+\frac ab$ \fa $a,b\in\Q_{>0}$ by (4.5.16). We prove $\vf\bg(+\frac m1)=+\frac m1$ by \In\ on $m\in\Na$. Set
\[
A:= \biggl\{m\in\Na: \vf\Bgg(+\frac m1)=+\frac m1\biggr\}.
\]
We have $1\in A$, since $+\frac11$ is the \nel\ of the monoid $\bg(\wh X,\cdot,+\frac11)$, and $\vf$~is an endo\mf\ of the monoid
$\bg(\wh X,\cdot,+\frac11)$. \Mo $k\in A$ implies $k+1\in A$. Indeed, let $k\in\Na$, then $\frac{k+1}1 = \frac k1+\frac11$ in $\Q_{>0}$ by (4.5.12)
since $1\cdot1 = 1$. \Mo $+\frac{k+1}1 = j\bg(\frac{k+1}1) = j\bg(\frac k1+\frac11) = j(\frac k1)+j(\frac11) = \bg(+\frac k1)+\bg(+\frac11)$ since
$j$~is a \hm sm. Hence if $k\in A$, $\vf\bg(+\frac{k+1}1) = \vf\bg(+\frac k1) + \vf(+\frac11) \nad{k\in A}= \bg(+\frac k1)+\bg(+\frac 11)
= j\bg(\frac k1)+j\bg(\frac11) = j\bg(\frac k1+\frac11) = j\bg(\frac{k+1}1) = +\frac{k+1}1$. Thus $k+1\in A$. \If that $A\in\Na$, hence $\vf\bg(+\frac
m1) = + \frac m1$ \fa $m\in\Na$. \csq, if $m,n\in\Na$, then $\vf\bg(+\frac mn) \nde3.5 = \vf\bg(+\frac mn)\cdot \vf\bg(+\frac 1n) = \bg(+\frac m1)
\cdot\bg(+\frac1n) \nad{(4.5.19)}= +\bg(\frac m1\cdot\frac1n) \nad{(4.5.3)}= +\frac mn$. \E\Tf if $r\in\wh\Z\sms0$, then either $r=+\frac mn$ or
$r=-\frac mn$ \fs $m,n\in\Na$. In both cases $\vf(r)=r$. Since $\vf(0)=0$, $\vf$~is the identity in~$\Q$.
\endproof

The \Pf\ $\Q$ is an \ti{ordered field\/} (see \E\df\ 4.5.50 and \E\Pr\ 4.5.48).

\bns 3.3 \

\hph i,i, Let $F$ be a field, then
\bea 3.6
\th_a(x) &:= x+a, \kern-50pt&& x,a\in F, \\
\d_c(x) &:= c\cdot x, \kern-50pt&& x\in F,\ c\in F\sms0. \lb{3.7}
\e

\hph ii,, If $F$ is equipped with a total \og~$\le$, we set
\beq 3.8
F_{\ge0} \hbox{ (resp.\ $F_{>0}$)} := \{x\in F: x\ge0 \hbox{ (resp.\ $x>0$)}\}.
\e
\ens

\brm 3.4 \

\hph i,ii, The bi\jc\ selfmaps of~$F$, $\th_a$, $a\in F$ (resp.\ $\d_c$, $c\in F\sms0$) are sometimes called \ti{\tl s}\index{translations} (resp. \ti{dilations}\index{dilations}).

\hph ii,i, The inverse of the selfmap $\th_a$ is $\th_{-a}$, $a\in F$, the inverse of $\d_c$ is $\d_{\frac1c}$, $c\in F\sms0$.

\hph iii,, \E\df\ 4.5.50 states that a field is an ordered field if \te s a total \og~$\le$ on~$F$ \st $\th_a$ (resp.~$\d_c$) are strictly in\cre\
\fa $a\in F$ (resp.\ ${c\in F_{>0}}$). The inverse of~$\th_a$, denoted by~$\th_a\Inv$, is strictly in\cre. Similarly for~$\d_c$, $c>0$. Indeed, if
$a,x,y\in F$ \sf y $x<y$ then $\th_a\Inv x \le \th_a\Inv y$. Suppose $\th_a\Inv x=\th_a\Inv y$, then $x= \th_a\circ \th_a\Inv x = \th_a\circ
\th_a\Inv y = y$, a \cd ion. Similarly for $\d_c$, $c>0$.

\csq, if $F$ is an ordered field, we have \fa $a,x,y\in F$ and all $c\in F_{>0}$:
\bea 3.9
x<y &\q\ \hbox{iff }x+a<y+a, \\
x<y &\q\ \hbox{iff }cx<cy. \lb{3.10}
\e

\hph iv,, We use the notation $y\ge x$ (resp.~$y>x$) iff $x\le y$ (resp.~$x<y$), \fa $x,y\in F$.
\erm

\bpr 3.5
Let $(F,\ge)$ be an ordered field. Then the \fw\ assertions hold\/{\rm:}
\bea 3.11
{}&\hbox{\E\fa $x,y\in F_{>0}$, both $x+y$ and $x\cdot y$ belong to $F_{>0}$.}\\
&\hbox{If $a\in F\sms0$, then $a^2\in F_{>0}$. \E\Ip $1\in F_{>0}$.} \lb{3.12}
\e
\epr

\proof \

``\er{3.11}'': Let $x,y\in F_{>0}$. Then $0<x=0+x \nde3.6 = \th_x0 \nde3.9 < \th_xy \nde 3.6 = y+x = x+y$. Hence $0<x+y$ by \tr ity of~$<$.
Similarly, $0=x\cdot 0\nde3.7 =\d_x0 \nde3.10 < \d_xy \nde3.7 = x\cdot y$.

``\er{3.12}'': Let $a\in F\sms0$. Since the \og~$\le$ is total, we have either $0<a$ or $a<0$. If $0<a$, then $a\in F_{>0}$. If $a<0$, then
$0= (-a)+a \nde3.9 < (-a)+0 = -a$, hence
$-a\in F_{>0}$. Thus either $a\in F_{>0}$ or $-a\in F_{>0}$. In the former case, $a^2 = a\cdot a\in F_{>0}$ by \er{3.11}. In the latter case,
$a^2\nad\ast= (-a)\cdot(-a)\in F_{>0}$ by \er{3.11}. In $\nad\ast=$ we used $(-x)\cdot y=y\cdot(-x) = -(x\cdot y)$ and $-(-x)=x$, $x,y\in F$.
Indeed, $0\nad{(4.1.1)} = 0\cdot y = (x+(-x))\cdot y \nad{(4.1.3)} = x\cdot y + (-x)\cdot y$, $x,y\in F$. Finally, $1\ne0$ and $1=1^2$ implies $1\in
F_{>0}$.
\endproof

\bco 3.6
Let $(F,\ge)$ be an \of, let $P:=F_{>0}\cup\{0\}$ $(=F_{\ge0})$, and let $\dot P:=F_{>0}$. Then the \fw\ assertions hold\/{\rm:}
\bea 3.13
{}&\hbox{$P$ is a \sbm\ of the monoid $\pz1F$,}\\
{}&\hbox{$P$ is a \sbm\ of the monoid $\pz2F$,} \lb{3.14} \\
&\hbox{$\pz0P$ is a $($\cmt e$)$ semiring $($with unity$)$,} \lb{3.15} \\
&\hbox{The monoid $\pz1P$ is a \PM\ whose \nog\ is the \rt ion}\lb{3.16} \\
&\hbox{to~$P$ of the total \og~$\ge$,} \non \\
&\hbox{$\dot P$ is a \sbm\ of the monoid $\pz2P$,} \lb{3.17} \\
&\hbox{$\pz2{\dot P}$ is an \ag.} \lb{3.18}
\e
\eco

\proof \

``\er{3.13}'': Clearly, $0\in P=\dot P\cup\{0\}$. If $x\in P$, then $x+0 = 0+x = x\in P$. If $x,y\in\dot P$, then $x+y\in\dot P\sbs P$ by \er{3.11}.

``\er{3.14}'': $1\in\dot P\sbs P$ by \er{3.12}. If $x\in P$, then $0\cdot x = x\cdot0 = 0\in P$. If $x,y\in P$, then $x\cdot y \in\dot P\sbs P$ by
\er{3.11}.

``\er{3.15}'': In view of \er{3.13} and \er{3.14}, the \rt ions to~$P$ of the binary \op s $+$~and~$\cdot$ are binary \op s on~$P$. \E\Ip
(4.1.1)--(4.1.3) hold, hence $\pz0F$ is a \sr\ in view of \er{3.13}, \er{3.14}.

``\er{3.16}'': Since the monoid $\pz1F$ is abelian and \cnc e, so is its \rt ion to~$P$. We next show that $x+y=0$ in $P$, where $x,y\in P$,
implies $x=y=0$. If $x\in\dot P$ and $y=0$, $x+y=x\in \dot P$ contrary to the \as. Similarly, $x=0$ and $y\in\dot P$ implies $x+y\in\dot P$.
In view of \er{3.11} the case $x,y\in\dot P$ is impossible. Thus
$x=y=0$. Hence $\pz1P$ is a \hbox{\PM}. Let $\le$ denote the \rt ion to~$P$ of the total \og\ of the \of~$F$, and let $\nad+\le$ denote
the \nog\ in~$P$, defined in (3.1.4). Let $x,y\in P$ be \st $x\nad+\le y$. Then \te s $z\in P$ \st $y=x+z$. Since $0\le z$, we obtain $x=0+x
\nde3.9 \le z+x=x+z=y$. Hence $x\le y$ by \tr ity of~$\le$. Now, let $x,y\in P$ be \st $x\le y$. Then $0=x+(-x)\nde3.9 \le y+(-x) = y-x$. Thus
$y-x\in P$ and $y=y+0= y+((-x)+x) = (y+(-x))+x = (y-x)+x = x+(y-x)$. Since $y-x\in P$, we obtain $x\nad+\le y$. \E\Tf the \rt ion of the \og\ of the
field~$F$ to~$P$ is identical to the \nog\ of the \PM\ $\pz1P$.

``\er{3.17}'': $1\in\dot P$, and  $x,y\in\dot P$ implies $x\cdot y\in\dot P$ by \er{3.11}.

``\er{3.18}'': Since the binary \op~$\cdot$ on~$F$ is abelian, so is its \rt ion to~$P$. Hence $P$~is an \am. Since $\dot P$ is a \sbm\ of $\pz2P$,
$\pz2{\dot P}$ is an \am. In view of \E\df\ 4.3.2 it suffices to show that \fe $x\in \dot P$ \te s $y\in\dot P$ \st $x\cdot y=1$. Since $F$~is a
field, \te s $y\in F\sms0$ \st $x\cdot y=1$. Then either $y>0$, or $y=0$, or $-y>0$. We have $y\ne0$, otherwise $x\cdot y=0$, which is impossible
since $0\ne1$. Suppose, for \cd ion, that $-y>0$. Then $x\cdot(-y)\in\dot P$ by \er{3.11}. However, $0=x\cdot0 = x\cdot(y+(-y)) = x\cdot y +
(x\cdot(-y)) = 1+z$ where both $1$ and $z$ belong to~$\dot P$. Hence $0 = 1+z \in\dot P$ by \er{3.11}, a~\cd ion. \If that $y\in\dot P$ is the
inverse of~$x$ in $\pz2{\dot P}$. \E\Tf $\pz2{\dot P}$ is an \ag.
\endproof

\bdf3.7
A $($\cmt e$)$ \ti{\sr} $($with unity$)$ $\pz0P$ is called an \tb{ordered semifield} if the \fw\ \as s are\index{ordered semifield} \sf ied:

\hph i,ii, The additive monoid $\pz1P$ is a \PM.

\hph ii,i, The \nog\ of $\pz1P$, introduced in (3.1.4) and denoted by~$\nad+\ge$, is a total \og\ on~$P$.

\hph iii,, $\dot P:=P\sms0$ is a \sbm\ of the \mlv\ monoid $\pz2P$.

\hph iv,, The monoid $\pz2{\dot P}$ is an \ag.
\edf

\brm3.8
\E\df\ \rf{d3.7} is not standard and is motivated by the \fw\ theorem.
\erm

\bth3.8 \

{\rm(i)\hphantom{i}} Let $(F,\ge)$ be an \of. Then $F_{\ge0}$ \er{3.8} is a \sbm\ of both the additive and the \mlv\ monoids of~$F$, and the
\rt ion of~$F$ to~$F_{\ge0}$ is an \emph{ordered semifield}.

{\rm(ii)} Conversely, let $P$ be an ordered semifield. Then \te s an \of\ $(F,\ge)$ \st $F_{\ge0}=P$.
\eth

\proof
(i) is a direct con\sq\ of Corollary \rf{c3.6}.

(ii): Step 1. ``\ti{Embedding of\/ $P$ in a totally ordered group $\pz1F$}'': \ By \E\df\ \rf{d3.7}, $\pz1P$ is a \PM\ containing
$1\ne0$. We show that $P$ is \ti{infinite}. Otherwise, $P$~is a finite \Cm\ (see \E\df\ 2.1.4), hence $P$~is a group by \E\Pr\ 4.3.10\,(i). \If that
every \el~$x$ of~$P$ has an inverse $x\Inv\in P$ \sf ying $x+x\Inv = 0$. Since $P$~is a \PM, we have $x = x\Inv = 0$. \E\Ip $1+1\Inv=0$, hence
$1=0$, \cd ing $1\ne0$. Hence a nontrivial \PM\ is infinite. \Wanp apply Corollary 4.3.86 with $X:=P$, ${\qu}:=+$ and $e:=0$.

\sdim{(3.19)}
We may identify $X$ with $j(X)$, and we find a set $\wh X$ with the \fw\ \pp ies:
\bit
\refstepcounter{equation}\lb{3.19}\item[\er{3.19}]
$\wh X$ is the disjoint union of the set $X$ and the set $\wh X_-:=\{-x: x\in X\sms e\}$.
\refstepcounter{equation}\lb{3.20}\item[\er{3.20}]
\E\te s a binary \op~$\hqu$ on~$\wh X$ \st $(\wh X,\hqu,e)$ is an \ag, \sf ying $x\hqu y = x\qu y$ \fa $x,y\in X$ and $x\hqu (-x)=e$ \fa $x\in X$.
\refstepcounter{equation}\lb{3.21}\item[\er{3.21}]
\E\te s a total \og\ $\whge$ on $\wh X$ defined by $\wh x\whge\wh y$ if $\wh x\hqu(\wh y\Inv)\in X$ where $\wh y\hqu \wh y{}\Inv = e$ \fa $\wh x,
\wh y\in \wh X$. \Mo $\wh x\whge e$ iff $\wh x\h{\in}\h X$ and $\wh x\whle e$ iff $\wh x\h{\in}\h \wh X_-{\cup}\{e\}$.
\refstepcounter{equation}\lb{3.22}\item[\er{3.22}]
The \tl\ maps $\th_{\hat a}$, $\wh a\in\wh X$, are strictly in\cre.
\eit

\E\pp y \er{3.22} needs some explanation. Indeed, if $\wh x\whg \wh y$, $\wh x,\wh y\in\wh X$, then \te s $p\in X\sms e$ \st $\wh x=\wh y+p$, since
$\wh x\hqu(\wh y{}\Inv)\in X\sms e$. Then $\wh x+\wh a = \wh y+\wh a+p$, hence $\th_{\hat a}\wh x = \th_{\hat a}\wh y\hpl p$, and $\th_{\hat a}\wh x
\hpl(\th_{\hat a}\wh y)\Inv = p\in X\sms e$, which implies $\th_{\hat a}\wh x \whg \th_{\hat a}\wh y$.

We now give a more explicit description of the binary \op~$\hqu$, and of the \op\ $\INV(\wh x):=\wh x{}\Inv$, $\wh x\in \wh X$.

Let $\wh x\in\wh X$, then
\bea 3.23
{}& \INV \wh x = \bca
-x & \hbox{if }\wh x:= x\in X_+,\\
x  & \hbox{if }\wh x:=-x, \ x\in X,
\eca \\
& \wh x\hqu e = e \hqu \wh x = \wh x \qh{\fa $\wh x\in \wh X$,} \lb{3.24} \\
& x\hqu y = x\qu y \qh{\fa $x,y\in X$,} \lb{3.25} \\
& (-x)\hqu(-y) = -(x\qu y) \qh{\fa $x,y\in X$,} \lb{3.26} \\
& x\hqu (-y) = p\in X \q \hbox{if } y=x\qu p, \ x,y\in X, \lb{3.27} \\
& x\hqu (-y) = -p, \ p\in X, \q \hbox{if } y=x\qu p, \ x,y\in X, \lb{3.28} \\
& (-y)\hqu x = x\hqu (-y) \q \hbox{if } x,y\in X. \lb{3.29}
\e

\er{3.26} follows from (4.3.5) and the \cmt ity of $\hqu$. \er{3.27} and \er{3.28} follow from (4.3.110). \er{3.29} follows from the \cmt ity
of~$\hqu$.

Finally, if $(F,+,0):=(\wh X,\hqu,e)$, then $\pz1F$ is a totally ordered group
\st $F$~is the disjoint union of $P\sms0$, $\{0\}$, and $\{-x: x\in P\}$, and \st the \tl s $\th_{\hat a}$, $\wh a\in F$, are strictly in\cre.

\ssk
Step 2. ``\ti{The \mlv\ monoid of $F$}'': Our aim is to find an \ext\ of the \mlc~$\cdot$ of~$F$ which makes~$F$ a~\sr.
It turns out that \te s \ooo such \ext\ which we still denote by~$\cdot$. Indeed, we should have
$$
0\cdot x= x\cdot0 =0 \qh{\fa $x\in F$.}
$$
Hence $0= (x+(-x))\cdot y = x\cdot y+(-x)\cdot y$ \fa $x,y\in F$, thus
$$
(-x)\cdot y = -(x\cdot y) \qh{\fa $x,y\in F$.}
$$
\Mo $y\cdot(-x) = (-x)\cdot y = -(x\cdot y)$ and $(-x)\cdot(-y) = -(x\cdot(-y)) = -(-(x\cdot y)) = x\cdot y$, $x,y\in F$. \E\Tf we should have
\bea 3.30
x\cdot(-y) &= -(x\cdot y),\\
(-x)\cdot y &= -(x\cdot y), \lb{3.31} \\
(-x)\cdot(-y) &= x\cdot y, \lb{3.32} \\
0\cdot(-x) &= 0, \lb{3.33} \\
(-x)\cdot0 &= 0, \qh{for all $x,y\in F$.} \lb{3.34}
\e
\E\fa $x\in F$, we have either $x=0$, or $x\in P$, or $-x\in P$. Thus \rt ing \er{3.30}--\er{3.34} to $x,y\in P$, and replacing the \et y sign
in \er{3.30}--\er{3.34} by the defining symbol~$:=$, we obtain a \df\ of the \ext\ of the \mlc\ of~$P$. Thus we obtain (4.5.18), (4.5.19) where
$\Q_{>0}$ is replaced by $F\sms0$. Proceeding as in the proof of Theorem 4.5.6, replacing $\Q_{>0}$ by~$F\sms0$, we find that $\pz2{F\sms0}$ is an
\ag\ and that $F\sms0$ is a \sbm\ of the \am\ $\pz2F$.

\sdim{(iii)}
\ssk
Step 3. ``\ti{$F$ is an \of\ \st $F_{\ge0}=P$}'': It remains to show:
\ben
\item (4.1.1) and (4.1.2), which imply by what precedes that $\pz0F$ is a field.
\item $\th_a$, $a\in F$, and $\d_c$, $c\in F_{>0}$, are strictly in\cre, which implies that $(F,\ge)$ is an \of.
\item $P= F_{\ge0}$.
\een

Note that (4.1.1) follows from Step~2, $\th_a$ strictly in\cre, and (iii) follow from Step~1.

``(4.1.2)'': The case $x:=0$ follows from (4.1.1) and $0+0=0$. If (4.1.2) holds for $x\in F_{>0}$, $y,z\in F$, then it holds for $-x$,
$x\in F_{>0}$, $y,z\in F$. Indeed, $(-x)\cdot(y+z) = -(x\cdot(y+z)) = x\cdot(-(y+z)) = x\cdot((-y)+(-z)) = (x\cdot(-y)) + (x\cdot(-z))
= (-x)\cdot y+(-x)\cdot z$. We now suppose $x\in F_{>0}$. The case $y,z\in F_{<0}$ follows from $x\cdot(y+z) = x\cdot(-(-y)+(-(-z)))=
-(x\cdot((-y)+(-z))) \nad*= -((x\cdot(-y))+(x\cdot(-z)))= -(x\cdot(-y))+\break(-(x\cdot(-z))) = x\cdot y + x\cdot z$. In~$\nad*=$ we used $x,-y,-z
\in F_{>0}$. It remains to consider the case $x,y,-z\in F_{>0}$ in view of the \cmt ity of the \ad. We distinguish  three cases: $y+z=0$, $y+z=:p>0$,
$y+z = -p$ with $p>0$. If $y+z=0$, then $x\cdot (y+z)= x\cdot0 =0 = (x\cdot y) + (-(x\cdot y)) = x\cdot y+x \cdot(-y) = x\cdot y+x\cdot z$. If
$y+z=p>0$, then $y = p+(-z)$ with $p>0$ and $-z>0$. Then $x\cdot y=x\cdot p+ x\cdot(-z)$, hence $x\cdot(y+z) = x\cdot p = x\cdot y+ (-(x\cdot (-z)))
= x\cdot y+x\cdot z$. Finally, if $y+z =p <0$, then $y+(-p) = -z$, with $y,-p,-z \in F_{>0}$. Then $x\cdot(-z) = x\cdot(y+(-p)) = x\cdot y+ x\cdot
(-p)$. Thus $-(x\cdot(-p)) = x\cdot y+(-(x\cdot(-z)))$, and $x\cdot(y+z) = x\cdot p = -(x\cdot(-p)) = x\cdot y+ x\cdot z$. This completes the proof of
(4.1.2).

\ssk
``\ti{$\d_c$, $c>0$, is strictly in\cre}'': Let $x,y\in F$, $x<y$. Let $p\in F_{>0}$ be \st $y=x+p$. Then $c\cdot y\nad{\rm(4.1.2)}= c\cdot(x+p) =
c\cdot x+c\cdot p$. Since $c,p\in F_{>0}$ and $\pz2{F_{>0}}$ is a monoid, we have $c\cdot p\in F_{>0}$, hence $c\cdot p>0$. \E\Tf $c\cdot x <
c\cdot y$, and the proof of Theorem \rf{t3.8} is complete.
\endproof

We now consider an \ev t \df\ of \of s. We begin with \oag s. An \ag\ $(G,\qu,e)$ is called an \ti{ordered group} if a (partial) \og~$\ge$ is
defined on~$G$ \st all \tl s $\th_a$, $a\in G$, are in\cre. The \ti{\po\ cone}\index{positive cone} of $(G,\ge)$ is defined by $G_{\ge e}:=\{x\in G: x\ge e\}$ (see
\cite[pp.~287--288]{Lattice}) and the set $\{e\ge x: x\in G\}$ is denoted by $G_{\le e}$. In view of the \ti{anti\sy y} of~$\ge$, we have
$G_{\ge e}\cap G_{\le e}=\{e\}$. If $A,B$ are \ti{nonempty} subsets of~$G$, we set
\bea 3.35
A\Inv &:= \{x\Inv \in G: x\in A\}, \\
A\qu B &:= \{x\qu y\in G: x\in A,\ y\in B\}.\lb{3.36}
\e
The cone $G_{\ge e}$ is called \ti{generating}\index{generating cone} if $G = G_{\ge e}\qu G_{\le e}$ and \ti{total\/} if $G=G_{\ge e}\cup G_{\le e}$. Hence the \og~$\ge$
is a \ti{total \og} iff the \fw\ assertion holds:
\beq 3.37
G = G_{\ge e}\cup G_{\le e} \qh{and }G_{\ge e}\cap G_{\le e} = \{e\}.
\e
Note that \ev ce \rl s and total \og s are everywhere defined \rl s on~$G$ since they are \ti{\rfl}. Let $R$~be an everywhere defined \rl\ on
an \ag\ $(G,\qu,e)$. Motivated by (4.3.98), the \op~$\qu$ and \rl~R are called \ti{compatible} if the \fw\ assertion holds:\index{compatible}
\beq 3.38
x \mathrel{\rm R} y \qh{implies } \th_ax \mathrel{\rm R} \th_ay \hbox{ \fa}a,x,y\in G.
\e
Since $\th_a$'s, $a\in G$, are bi\jc\ and $\th_a\Inv = \th_{a\Inv}$, $a\in G$, we also have \fa $a,x,y\in G$:
\[
\th_ax \mathrel{\rm R} \th_ay \qh{implies }\th_{a\Inv}\th_a x \mathrel{\rm R} \th_a\Inv \th_ay,
\]
hence $\th_ax \mathrel{\rm R} \th_ay$ implies $x\mathrel{\rm R}y$. \E\Tf we have $x\mathrel{\rm R}y$ iff $\th_ax\mathrel{\rm R}\th_ay$,
$a,x,y\in G$. Taking $a:=y\Inv$, respectively $a:=y$, we find
\beq 3.39
x\mathrel{\rm R}y \qh{iff } x\qu y\Inv \mathrel{\rm R} e,\q x,y\in G.
\e

\blm 3.9
Let $(G,\qu,e)$ be an \ag\ $($see \E\df s {\rm4.3.2} and {\rm 4.3.8)}, let $\th_a$, $a\in G$ 
$($see \E\df\ {\rm 4.3.46)}, and let {\rm R}
be an everywhere defined \rl\ on~$G$. Let $K$ and $K\Inv$ be the subsets of~$G$ defined by
\bea3.40
K&:= \{x\in G: x\mathrel{\rm R}e\},\\
K\Inv &:= \{x\Inv \in G: x\in K\}, \q \hbox{if }K\ne\vn.\lb{3.41}
\e
Then the \fw\ assertions hold\/{\rm:}

\noindent If the \op~$\qu$ and the everywhere defined \rl\ $\mathrel{\rm R}$ are compatible \er{3.38}, we have
\beq 3.42
\hbox{\fa }x,y\in G, \ x\mathrel{\rm R}y \hbox{ iff }x\qu y\Inv \in K.
\e
\E\Ip $K\ne\vn$. Suppose that \er{3.42} holds, then\/{\rm:}
\bea 3.43
{}&\hbox{R is \rfl\ iff $e\in K$}, \\
&\hbox{R is \tr e iff $K\qu K\subset K$,} \lb{3.44} \\
&\hbox{R is \sy ic iff $K\sbs K\Inv$,} \lb{3.45} \\
&\hbox{R is anti\sy ic iff $K\cap K\Inv\subseteq\{e\}$.} \lb{3.46}
\e
\elm

\proof
``\er{3.42}'': follows from \er{3.39}.

``\er{3.43}'': Suppose $e\in K$. \E\fa $x\in G$ we have $x\qu x\Inv=e$, hence $x\mathrel{\rm R}x$ \fa $x\in G$. Suppose R is \rfl, then $x
\mathrel{\rm R}x$ \fa $x\in G$. Thus \te s $x\in G$ \st $x\mathrel{\rm R}x$, hence $e= x\qu x\Inv \in K$ by \er{3.42}.

``\er{3.44}'': Suppose $K\qu K\sbs K$, and let $x,y,z\in G$ be \st $x\mathrel{\rm R}y$ and $y\mathrel{\rm R}z$. Then $x\qu y\Inv$ and $y\qu z\Inv$
belong to~$K$ by \er{3.42}. \E\Tf $x\qu z\Inv = x\qu e\qu z\Inv = x\qu (y\Inv\qu y)\qu z\Inv \nad{(2.1.36)}= (x\qu y\Inv)\qu (y\qu z\Inv)\in K$.
Hence $x\mathrel{\rm R}z$, and R~is \tr e. Suppose that R~is \tr e, and let $x,y\in K$. Then $x\mathrel{\rm R}e$ and $y\mathrel{\rm R}e$ by \er{3.40}.
From \er{3.38} we infer $\th_x y\mathrel{\rm R}\th_x e$. From $\th_xe = e\qu x =x$, and $x\mathrel{\rm R}e$, we find $\th_xe\mathrel{\rm R}e$.
Using \tr ity of~R we obtain $\th_x y\mathrel{\rm R}e$. Hence $x\qu y = \th_xy\in K$ by \er{3.40}. Since $x,y$ are arbitrary \el s of~$K$, we
obtain $K\qu K \sbs K$.

``\er{3.45}'': Suppose that $K\sbs K\Inv$. Let $x,y\in G$ be \st $x\mathrel{\rm R}y$, hence $x\qu y\Inv\in K$ by \er{3.42}. Since $K\sbs K\Inv$, we
have $y\qu x\Inv = (x\qu y\Inv)\Inv\in K$. Hence $y\mathrel{\rm R}x$, and R~is \sy ic. Suppose that R~is \sy ic.
Let $x\in K$. Hence $x\mathrel{\rm R}e$ by
\er{3.40}. Then $e\mathrel{\rm R}x$ by \sy y of~R. We obtain $x\Inv = e\qu x\Inv\in K$. Thus $K\sbs K\Inv$, since $x$ is an arbitrary \el\
of~$K$.

``\er{3.46}'': Observe that $K\cap K\Inv=\vn$ iff there is no pair $(x,y)\in G\t G$ \st $x\mathrel{\rm R}y$ and $y\mathrel{\rm R}x$. In that case
\as\ (2.3.4) is (vacuously) \sf ied, hence R~is \asm. (For example, if $(E,\ge)$ is a \tos\ then the \rl~R, defined by $x\mathrel{\rm R}y$ if $x>y$,
is said to be \asm.) We now suppose $K\cap K\Inv\ne\vn$. We first assume $K\cap K\Inv=\{e\}$. If $x\mathrel{\rm R}y$ and $y\mathrel{\rm R}x$, then
$x\qu y\Inv\in K$ and $y\qu x\Inv\in K$, hence $x\qu y\Inv \in K\cap K\Inv=\{e\}$. \If that $x=y$, hence R~is \asm\ in view of (2.3.4). We now suppose
that R~is \asm. Let $z\in K\cap K\Inv$. Then $z\mathrel{\rm R}e$ and $e\mathrel{\rm R}z$. Thus $z=e$ by anti\sy y.
\endproof

\brm3.10
Let $(G,\qu,e)$ and {\rm R} be as in Lemma {\rm\rf{l3.9}}. Then

\hph i,ii, If $\qu$ and R are compatible, then R is \rfl\ and \tr e iff $K$ is a \sbm\ of $(G,\qu,e)$ and $(K,\qu,e)$ is a \Cm.

\hph ii,i, If $\qu$ and R are compatible and R is an \ev ce \rl, then the \sbm~$K$ is a \ti{subgroup} of the group~$G$. The \ev ce classes $[x]$
(\st $x\in[x]$) are of the form $x\qu K:=\{x\}\qu K$ \fa $x\in G$, and are called \ti{cosets} (see \E\Pr\ 4.5.19).

\hph iii,, If $\qu$ and R are compatible and R is an \og~$\ge$ on~$G$, then the \Cm~$K$ is a \ti{\PM}. Indeed, if $x,y\in K$ \sf y $x\qu y=e$,
then $y\Inv= e\qu y\Inv = (x\qu y)\qu y\Inv = x\qu(y\qu y\Inv) = x\qu e =x$. Hence $x\in K\cap K\Inv=\{e\}$, and $x=e$, $y=e$. Note that if $x,y\in
K$ and $x\ge y$, then $x\qu y\Inv\in K$ by \er{3.42}, hence \te s $z\in K$ \st $x\qu y\Inv=z$. Thus $x=z \qu y = y\qu z$, hence $x\stackrel\qum\ge y$
(see (3.1.4)). Conversely, if $x\stackrel\qum\ge y$, then trivially $x\ge y$. Thus \ti{the \nog\ on~$K$ is the \rt ion to~$K$
of the \og~$\ge$ on~$G$}. Observe that if the \og\ of~$G$ is compatible and total, then the \nog\ $\stackrel\qum\ge$ of~$K$ is total. Indeed, if
$x,y\in K$, then either $x\qu y\Inv\ge e$ or $e\ge x\qu y\Inv$, and either $x\stackrel\qum\ge y$ or $y\stackrel\qum\ge x$. We prove the converse in
the next lemma.
\erm

\blm3.11
Let $(G,\qu,e)$ be an infinite \ag.

\hph i,i, Suppose \te s a total \og\ $\ge$ on~$G$ \st all \tl s $\th_a$, $a\in G$, are in\cre. Then the subset~$P$ of~$G$ defined by
\beq 3.47
P:= \{x\in G: x\ge e\}
\e
\sf ies
\bit
\item $P$ is a \sbm\ of $G$,
\item $P$ is a \PM,
\item $x\ge y$ iff $x\qu y\Inv\in P$ \fa $x,y\in G$.
\eit
\E\Ip the \rt ion of $\ge$ to $P$ is the \nog\ of the \PM~$P$.

\hph ii,, Conversely, if $H$ is a subset of $G\sms e$ \sf ying
\bea 3.48
{}&\hbox{$G$ is the disjoint union of $H$, $H\Inv$ and $\{e\}$.}\\
&\hbox{\E\fa $x,y\in H$ we have $x\qu y\in H$.} \lb{3.49}
\e
Let $\ge$ denote the \rl\ on $G$ defined by
\beq 3.50
x\ge y \qh{if } x\qu y\Inv \in H\cup \{e\}, \q x,y \in G.
\e
Then $\ge$ is the only total \og\ on~$G$ \st $\qu$ and $\ge$ are compatible and \st
\beq 3.51
\{x\in G: x>e\} = H.
\e
\Mo we have
\beq 3.52
x\ge y \qh{iff }y\Inv \ge x\Inv \hbox{ \fa}x,y\in G.
\e
\elm

\proof
Part (i) follows from Lemma \rf{l3.9} and Remark \rf{r3.10}.

``(ii)'': Set
\beq 3.53
P:=H\cup \{e\}.
\e
\E\fa $x\in G$ we have $x\qu x\Inv=e \in P$, hence $x\ge x$. Thus the \rl\ defined in \er{3.50} is \ti{everywhere defined\/} and \ti{\rfl}. \E\fa
$a,x,y\in G$ we have $x\qu y\Inv \in P$ iff $(x\qu a)\qu{(y\qu a)\Inv} \nad{(4.3.5)}= (x\qu a)\qu(a\Inv\qu y\Inv) \nad{(2.1.36)}= x\qu (a\qu a\Inv)\qu
y\Inv = x\qu e \qu y\Inv = x\qu y\Inv \in P$. Hence \ti{the \op~$\qu$ and the \rl~$\ge$ are compatible} by~\er{3.38}. Clearly, $e\in P$ by \df\
of~$P$. \Mo if $x,y\in H$ then $x\qu y\in H\sbs P$ by~\er{3.49}. If $x=e$ and $y\in H$, then $x\qu y = e\qu y = y \in H\sbs P$.
Similarly, if $x\in H$ and $y=e$. Finally, $e\qu e\in P$. \If that $x\qu y\in P$ whenever $x,y\in P$. Thus $P$ is a monoid, hence
$\ge$ is \ti{\tr e} by~\er{3.44}. If $x\ge y$ and $y\ge x$, we claim $x=y$. Set $z:=x\qu y\Inv$. Then $z\in P$,
and $z\Inv = (x\qu y\Inv)\Inv \nad{(4.3.4)(4.3.5)}= y\qu x\Inv \in P$. \E\Tf $z\in P\cap P\Inv = \{e\}$ by \er{3.48}, \er{3.53}. Hence $\ge$~is
\ti{anti\sy ic} by~\er{3.46}. Hence $\ge$~is an \og\ on~$G$. By \er{3.48} if $z\in G\sms e$, then $z\in H$ or $z\in H\Inv$. Let $x,y\in G$. If
$x\ne y$, then $x\qu y\Inv\ne e$ and $y\qu x\Inv\ne e$. Set $z:=x\qu y\Inv$. Then $x\qu y\Inv\in H$ or $x\qu y\Inv\in H\Inv$. In the former case,
we obtain $x\ge y$ since $H\sbs P$. In the latter case $y\qu x\Inv=(x\qu y\Inv)\Inv\in H$ since $(P\Inv)\Inv=P$. (We need (4.3.4) and (4.3.5).)
\E\Tf in the latter case $y\ge x$. \If that \fa pairs $x,y$ in~$G$ \st $x\ne y$, we have either $x>y$ or $y>x$. Thus the \og\ is total. Next we prove
\er{3.52}. Let $x,y\in G$ be \st $x\qu y\Inv\in P$. We have to show that $y\Inv\ge x\Inv$, that is, $y\Inv \qu (x\Inv)\Inv \in P$. But $y\Inv \qu
(x\Inv)\Inv = y\Inv \qu x = x\qu y\Inv\in P$. Hence $y\Inv\ge x\Inv$. Next we prove \er{3.51}. Let $x\in G$. Then $x>e$ iff $x\ge e$ and $x\ne e$.
We have $x\ge e$ iff $x\qu e\Inv \in P$ iff $x\qu e\in P$ iff $x\in P$. Hence $x\ge e$ iff $x\in P\sms e = H$, which
proves \er{3.51}. It remains to prove the statement about ``\uq'' of the \og~$\ge$. Let $\wtge$ be a total \og\ on~$G$ \st $\qu$~and~$\wtge$ are
compatible. Then it follows from \er{3.39} and \er{3.42} that for $x,y\in G$, $x\wtge y$ iff $x \qu y\Inv \in \{z\in G: z\wtge e\}$. Since the \rl\
R \sf ies \er{3.51}, we have $\{z\in G: z\wtge e\} = \{z\in G: z\mathrel{\wt>}e\}\cup\{e\} \nde3.51 = H\cup\{e\}$. Hence $x\wtge y$ iff $x\ge y$
by \er{3.50} \fa $x,y\in G$. \E\Tf the \rl s $\ge$ and $\wtge$ are equal.
\endproof

\Wanp give an \ev t \df\ of an \of.

\bpr 3.13 \

\hph i,i, Let $F$ be an infinite field, and let $H$ be a subset of $F$ \sf ying
\beq 3.54
\hbox{$F$ is the disjoint union of $H$, $\{0\}$ and $-H$,}
\e
where $-H:= \{-x\in F: x\in H\}$. Suppose
\beq 3.55
\bal
x+y \in H \q &\hbox{\fa}x,y \in H, \hbox{ and}\\
x\cdot y \in H \q &\hbox{\fa }x,y \in H.
\eal
\e
Then $(F,\ge)$ is an \of\ with the \rl\ $\ge$ defined by
\beq 3.56
x\ge y \hbox{ if }x+(-y)\in H\cup\{0\}, \ x,y\in F.
\e
\Mo we have
\bea 3.57
{}&H = F_{>0} \hbox{ defined in \er{3.8}}, \\
&\{x\in F: x<0\} = -H. \lb {3.58}
\e

\hph ii,, Conversely, if $(F,\ge)$ is an \of, and $H:=\{x\in F: x>0\}$, then \er{3.54}--\er{3.55} hold.

\Mo we have
\beq 3.59
x>y \qh{iff }x+(-y)\in F_{>0}.
\e
\epr

\proof
``(i)'': We first show that the group $\pz1F$ equipped with the \og\ defined in \er{3.56} is a total \og\ on~$F$ \st $+$ and~$\ge$ are compatible.
To this end we apply Lemma \rf{l3.11} with $(G,\qu,e)=\pz1F$, noting that \er{3.54}, \er{3.55} imply \er{3.48}, \er{3.49}. Hence the \rl\ defined in
\er{3.56} is a total \og\ on~$F$ \st $\th_a$~are in\cre\ \fa $a\in F$. Since $\th_a$'s are bi\jn s, they are strictly in\cre. (If $y>x$ then $\th_ay
\ge \th_a x$. If $\th_ax=\th_ay$, then $x=y$, \cd ing $y>x$.) \Mo \er{3.57} follows from \er{3.51}, and \er{3.59} from \er{3.52} with $y:=0$. It
remains to show that $\d_c$ is strictly in\cre\ \fa $c>0$. Then $z:=x+(-y)\in F_{\ge0}\sms0= F_{>0}$, and $cx+(-cy)\nad{(4.4.12)}= cx+c(-y)
\nad{(4.1.2)}=c(x+(-y))= cz\in F_{>0}$ by \er{3.55}. Hence $cx>cy$. This completes the proof of part~(i) of the \Pr.

\ssk
``(ii)'': We now suppose that $(F,\ge)$ is an \of\ in the sense of \E\df\ 4.5.50. Set
\beq 3.60
x\ge y \qh{if } y\le x \q \hbox{\fa }x,y\in F,
\e
and set $H:=F_{>0}$. Let $z\in F$. Since the \og\ $\ge$ is total, we have either $z\ge 0$ or $0\ge z$. Thus if $z\ne0$, we have either $z>0$
or $0>z$. In view of (4.5.20) and \er{3.9}, we have $0>z$ iff $0+(-z) > z+(-z)$, that is, $0>z$ iff $-z>0$ \fa $z\in F$. \If that \fa
$z\in F\sms0$, we have either $z\in H$ or $-z\in H$, which proves \er{3.58}. \Mo the first part of \er{3.55} follows from \er{3.11}.

Next we prove the second part of \er{3.55}. Let $x,y\in H$. Thus $x,y\in F_{>0}$ by the \df\ of~$H$. From $x>0$ and $y>0$, we infer from (4.5.121)
with $x:=0$, $y:=y$ and $z:=x$ that $0=0\cdot x<y\cdot x=x\cdot y$. Hence \er{3.55} holds. This completes the proof of the \Pr.
\endproof

\brm 3.14 \

\hph i,i, Let $(E,\le)$ be an \os\ (see \E\df\ 1.3.2). Then the \rl~$<$ on~$E$ defined by
\beq 3.61
x<y \qh{if } x\ne y \hbox{ and }x\le y, \q x,y\in E,
\e
is \tr e. More is true, we have
\bea 3.62
x< y \hbox{ and } y\le z &\hbox{ implies }x<z, \q x,y,z\in E,\\
x\le y \hbox{ and } y< z &\hbox{ implies }x<z, \q x,y,z\in E, \lb{3.63}
\e
\Mo the \fw\ assertions are \ev t:
\bea 3.64
{}&\hbox{either }x=y, \hbox{ or }x<y, \hbox{ or }y<x,\\
&\hbox{either }x\le y, \hbox{ or }y\le x. \lb{3.65}
\e
If \er{3.64} or \er{3.65} holds, the order $\le$ is called \ti{total\/} (see \E\df\ 1.3.18).

\hph ii,, Let R be a nonempty \rl\ on a \ns~$X$. Then the \rl\ $\mathrel{\wt{\rm R}}$ on~$X$ defined by $x\mathrel{\wt{\rm R}}y$ if $y\mathrel{\rm R}
x$, $x,y\in X$ (also denoted by $\mathrel{\wt{\rm R}}$ in \cite[p.~3]{Lattice} or ${\rm R}^{-1}$ in \cite{Kelley}) is called the \ti{converse} of the
\rl~R in \cite[p.~3]{Lattice}. One verifies that if the \rl~R is \rfl\ (resp.\ \tr e, \sy ic, anti\sy ic), then so is~$\mathrel{\wt{\rm R}}$. \E\Ip
the converse of an \ev ce \rl~R is equal to~R. The converse of an \og\ (resp.\ total \og) is an \og\ (resp.\ total \og). If the set~$X$ is not a
singleton, then a total \og\ on~$X$ is not equal to its converse. One verifies that if $(G,\qu,e)$ is a nontrivial \ag\ and $\ge$~is a total \og\
on~$G$, \st $\qu$~and~$\ge$ are compatible, then so is the converse \rl~$\wtge$. Thus $(G,\qu,e)$ is also an ordered group with the \rl~$\wtge$.
If $\pz0F$ is an \of\ with \og~$\ge$, then $1>0$ by \er{3.12}. If $F$~would be an \of\ with the converse \rl~$\wtge$, then we would have
$x\cdot y\in F_{\mathrel{\wt>}0}$ whenever $x,y\mathrel{\wt>}0$ by~\er{3.11}. Thus $0>x$ and $0>y$ would imply $0>x\cdot y$. Since $0>z$ iff $-z>0$
\fa $z\in F$, we have $(-x)\cdot(-y)>0$ by~\er{3.11}. However, $(-x)\cdot(-y)= x\cdot y<0$. A~\cd ion. Thus if $(F,\ge)$ is an \of, then $(F,\wtge)$
is \ti{not}. Another \chz\ of an \of\ can be found in \cite[Corollary 25.22, p.~411]{Is}.
\erm

We next consider \ti{\of s of \qt s}. The field $\Q$ of \ra\ \nm s is the prototype of a field of \qt s. The \ra\ \nm s are ``\qt s'' $\frac ab$,
$a\in \Z$, $b\in\Z\sms0$, where $\Z$~is the domain of \ig s (see Theorem 4.5.8). We recall that a domain~$D$ is a \sr\ \st its additive monoid
$\pz1D$ is an \ag, its \mlv\ monoid $\pz2D$ is abelian, and \st $D^\t:=D\sms0$ is a \sbm\ of $\pz2D$, and $\pz2{D^\t}$ is a \Cm\ (see
(4.5.60)--(4.5.64)). The self-maps of~$D$: $\th_a$, $a\in D$, $\d_c$, $c\in D^\t$, are defined in \er{3.6}, \er{3.7}, where $F$~is allowed to be a
domain. The maps $\th_a$, $a\in D$, are bi\jc\ with inverse $\th_a\Inv=\th_{-a}$. The maps $\d_c$, $c\in D^\t$, are in\jc. We also recall that if
$(E,\ge)$ is an \os, then an in\jc\ self-map of~$E$ is strictly in\cre\ iff it is in\cre. In \E\df\ 4.5.50, a field $(F,\ge)$ is called an
\ti{\of} if \te s a total \og~$\ge$ \st both $\th_a$ and~$\d_c$ are strictly in\cre\ \fa $a,c\in F$ with $c>0$. In \E\Pr\ \rf{p3.13} we gave an \ev t
\df. Note that a field is a domain.

\bdf3.15
A domain $\pz0D$ is called an \tb{ordered domain} if \te s\index{ordered domain} a total \og~$\ge$ on~$D$ \sf ying (4.5.120) and (4.5.121) \fa $x,y,z\in D$.
\edf

\bpr3.16
Proposition {\rm \rf{p3.13}} holds if $F$~is allowed to be an infinite domain~$D$ and $(D,\ge)$ is an ordered domain.
\epr

Replacing the field $F$ by the domain $D$, and \E\df\ 4.5.50 by \E\df\ \rf{d3.15}, the proof of \E\Pr\ \rf{p3.16} is verbatim the same as the proof
of \E\Pr\ \rf{p3.13}.

In Theorem 4.5.35 it is shown that an infinite domain can be embedded in a field (note that a finite domain is already a~field). More precisely, the
set $\wh D$ consisting of \ev ce classes denoted by~$\frac ab$ where $(a,b)\in D\t D^\t$ in (4.5.81), equipped with an \ad~$\hpl$ defined in
(4.5.84) and a \mlc~$\hcd$ defined in (4.5.85) is a field $(\wh D,\hpl,\hcd,\wh 0,\wh 1)$ where $\wh0:=\frac01$ and $\wh1:=\frac11$. \Mo the map
$j:D\to\wh D$ defined by $j(x):=\frac x1$, $x\in D$, is an in\jc\ ring-\hm sm. The image of~$D$ under~$j$, $j(D)$ is a domain contained in~$\wh D$.
The fractions $\frac ab$ and $\frac cd$ are equal iff $ad=bc$. It turns out that an ordered domain $(D,\ge)$ can be embedded in an \of\
$(\wh D,\whge)$ \sf ying $\frac ab\whge \frac cd$ iff $ad\ge bc$, $a,c\in D$, $b,d\in D^\t$ and $bd>0$.

\bth3.17
Let $(D,\ge)$ be an ordered domain, and let $(\wh D,j)$ be the field of \qt s of~$D$ introduced in Theorem {\rm 4.5.35}. Then $(\wh D,\whge)$ is an
\of\ with the \og~$\whge$ defined \fa $\wh x,\wh y\in\wh D$ by
\beq 3.66
\wh x\whge \wh y \qh{if } \wh x = \wh y\hpl \wh z \hbox{ \fs} \wh z:=\frac ab \hbox{ with } a\in D_{\ge0} \hbox{ and } b\in D_{>0}.
\e
\Mo we have
\beq 3.67
\wt D_{\htg\hat0} = \bigl\{\tfrac ab\in \wh D: a,b\in D_{>0}\bigr\},
\e
where $D_{>0}:=\{a\in D: a>0\}$,
\beq 3.68
j(x)\whg j(y) \qh{\fa} x,y\in D \hbox{ \st} x>y,
\e
and the \og\ $\whge$ is the \em{only} \og\ on~$\wh D$ \st $(\wh D,\whge)$ is an \of\ \sf ying \er{3.68}.
\eth

We shall use the \fw\ lemma in the proof of Theorem \rf{t3.17}.

\blm 3.18
Let $(F,\ge)$ be an \of, and let $a>0$ and $b\in F$ be \st $a\cdot b=1$. Then $b>0$.
\elm

\proof
We know that $1>0$ by \er{3.12}. In view of \er{2.3} we have $b\ne0$ since $1\ne0$. Suppose for \cd ion that $b\not>0$. Since the \og~$\ge$ is total,
we have $b<0$. Hence by (4.5.120) $-b = 0+(-b) > b+(-b) =0$, thus $-b>0$. By \er{3.11} we have $a\cdot(-b)>0$. However, $a\cdot b+ a\cdot(-b)
\nde2.4 = a\cdot(b+(-b)) = a\cdot0 \nde2.3 = 0$, hence $a\cdot(-b) = -(a\cdot b)=-1$. \If that $-1=a\cdot(-b)>0$. By (4.5.120) again, we find
$0 = (-1)+1 > 0+1 = 1$, that is, $0>1$. A~\cd ion, since the \og~$\ge$ is total. Hence $b>0$.
\endproof

\proof[Proof of Theorem \rf{t3.17}]
We first observe that if $\wh z\in \wh D$, then we may assume that $\wh z=\frac ab$, with $a\in D$ and $b\in D_{>0}$ defined in \er{3.67}. By \df\
of~$\wh D$ \te s $(c,d)\in D\t D^\t$ \st $\wh z=\frac cd$ defined in (4.5.81). If $d\not>0$, then $d<0$ since $d\in D\sms0$ and the \og~$\ge$ is
total. As in the proof of the preceding lemma we show that $c\cdot(-d) = -(c\cdot d) = (-c)\cdot d$ since in the domain~$D$ the monoid $\pz1D$ is
a group and \er{2.3}, \er{2.4} hold. \E\Tf $\frac cd=\frac{-c}{-d}$ by (4.5.83). Thus $\frac cd=\frac ab$ where $a:=-c$ and $b:=-d>0$. In what follows
we assume that $b$, the denominator of the fraction~$\frac ab$, belongs to $D_{>0}$. We apply \E\Pr\ \rf{p3.13}\,(i) with $F:=\wh D$ and $H:=\bigl\{
\wh z\in \wh D: z=\frac ab,\ a,b\in D_{>0}\bigr\}$. Since the \og~$\ge$ of~$D$ is total, $\wh D$~is the disjoint union of $H$, $\{\wh0\}$ and
$\bigl\{\wh z\in\wh D:\wh z:=\frac ab,\ {a<0},\break b>0,\ a,b\in D^\t\bigr\}$. Note that if $a\in D$, then $a<0$ iff $-a>0$ where $-a$~is the inverse
of~$a$ in the group $\pz1D$. Indeed, we have $0=a+(-a)<0+(-a)=-a$ iff $a<0$. \If that $\bigl\{\wh z\in \wh D: z=\frac ab, \ a<0,\ b>0,\ a,b\in D^\t
\bigr\} = -H$. \csq, \as\ \er{3.54} of \E\Pr\ \rf{p3.13} is \sf ied. We next verify \er{3.55}. Let $\wh x,\wh y\in H$ be \st $\wh x=\frac ab$,
$\wh y=\frac cd$ with $a,b,c,d\in D_{>0}$. Then $\wh x+\wh y \nad{(4.5.84)}= \frac{ad+bc}{bd}$. The RHS is well-defined since $bd\in D^\t$. Indeed,
$b,d\in D_{>0} \sbs D^\t$. \Mo $bd\in D_{>0}$ by \E\Pr\ \rf{p3.13}\,(ii) with $(F,\ge):=(D,\ge)$, $H:=D_{>0}$. Then $ad,bc,bd\in H$ and $ad+bc\in H$
by \er{3.55}. \If that $\wh x+\wh y\in H$. Similarly, $\wh x\cdot \wh y\nad{(4.5.85)}= \frac{ac}{bd}$ with $ac,bd\in D_{>0}$, hence $\wh x\cdot\wh y
\in H$. Thus \er{3.55} of \E\Pr\ \rf{p3.13}\,(i) is \sf ied. \E\Tf $(\wh D,\whge)$ is an \of\ with $\whge$ defined by $\wh x\whge\wh y$ if
$\wh x=\wh y+\wh z$, $\wh z\in H\cup\{\wh 0\}$. \Mo $H = \{\wh x\in \wh D: \wh x\whg\wh 0\}$ by \er{3.57}. Thus $\wh D_{\htg\hat0} =
\{\wh x\in \wh D: \wh x\whg\wh 0\} = H =\bigl\{\wh z\in \wh D: \wh z=\frac ab,\ a,b\in D_{>0}\bigr\}$, which is \er{3.67}. Next we prove \er{3.68}.
Let $x,y\in D$ be \st $x>y$. Then \te s $p\in D_{>0}$ \st $x=y+p$. \csq, $j(x)=j(y)\hpl j(p)= j(y)\hpl \frac p1$, by (4.5.87). Since $\frac p1\in
D_{\htg\hat0}$, and $H\cup \{\wh0\} = D_{\htg\hat0} \cup\{\wt0\}$, we obtain $j(x)\whg j(y)$ from \er{3.56}. It remains to prove the ``\uq'' part
of the theorem. If $(F,\ge)$ is an \of\ and $F_{\ge0}:=\{x\in F: x\ge0\}$, we have $x>y$ iff $x+a>y+a$ \fa $a,x,y\in F$ by \er{3.9}. Clearly, $x=y$
implies $\th_ax = \th_ay$ \fa $a,x,y\in F$ (the \et y and the map $\th_a$
are compatible \fa $a\in F$). Conversely, $\th_ax=\th_ay$ implies $x=y$ \fa $a,x,y\in F$ since $\th_a$'s are in\jc. \E\Tf $x\ge y$ iff $x+a\ge y+a$
\fa $a,x,y\in F$. Choosing $a:=-y$, we find $x\ge y$ iff $x-y\ge0$ iff $x-y\in F_{\ge0}$, \fa $x,y\in F$. Thus in an \of\ (more generally, ordered
\ag) the \og~$\ge$ is completely determined by the ``\po\ cone'' $F_{\ge0}$ (\Ip if $y=0$, then we have $x\ge0$ iff $x\in F_{\ge0}$, hence the
notation $F_{\ge0}$ is consistent). In Theorem \rf{t3.17} the \og~$\whge$ is subject to the condition: $x>y$ implies $j(x)\whg j(y)$ for $x,y\in D$.
Taking $y=0$, we find $x\in D_{>0}$ implies $j(x)\in \wh D_{\htg\hat0}$, since $j(0)=\frac01=\wh0$. Thus $j(D_{>0})\sbs \wh D_{\htg\hat0}$
($=F_{>0}$). \Mo if $\wh z\in D_{\htg\hat0}$, then $\wh z=\frac ab$ with $a,b\in D_{>0}$. Hence $\wh z=\frac a1 \hcd \frac 1b = j(a)\hcd j(b)\Inv$
since $\frac1b\hcd\frac b1 = \frac bb = \frac11= \wh 1$. We have $j(a),j(b)\in \wh D_{\htg\hat0}$, and $j(b)\Inv \in \wh D_{\htg\hat0}$ by Lemma
\rf{l3.18}. Since $(\wh D,\whge)$ is an \of, we have $j(a)\hcd j(b)\Inv \in \wh D_{\htg\hat0}$ by \er{3.11}. \E\Tf $\wh D_{\htg\hat0} \sbs
j(D_{>0}) \hcd j(D_{>0})\Inv \sbs \wh D_{\htg\hat0}$. \csq, we obtain
\beq 3.69
\wh D_{\htg\hat0} = j(D_{>0}) \hcd j(D_{>0})\Inv,
\e
where $j(D_{>0})\Inv := \{\wh z\in\wh D: z=j(b)\Inv,\ b\in D_{>0}\}$, and $A\hcd B= \{\wh z\in \wh D: \wh z=\wh x\hcd \wh y,\break \wh x\in A,\ \wh b
\in B\}$, where $A,B$ are \nss s of~$\wh D$. \If that $\wh D_{\htg\hat0}$ is \ti{uniquely} determined by $D_{>0}$.
\endproof

\bco 3.19
Let $(D,\ge)$, $(\wh D,j)$ and $\whge$ be as in Theorem {\rm\rf{t3.17}}. Let $a,b\in D^\t$. Then
\beq 3.70
\frac ab \whg \wh0 \qh{iff } a\cdot b>0.
\e
Let $a,b,d \in D^\t$, $c\in D$. If $bd>0$, then
\beq 3.71
\frac ab \whg \frac cd \qh{iff }ad>bc.
\e
If $bd<0$, then
\beq 3.72
\frac ab \whg \frac cd \qh{iff }ad<bc.
\e
\eco

\bex 3.20
Prove Corollary \rf{c3.19}. (Hint: Look at the proof of \E\Pr\ 4.5.4.)
\eex

Our next goal is to prove the assertion \fw\ Exercise 4.5.51 on page 266. To this end we shall use the \fw\ lemma.

\blm 3.21
Let $\pz0\Z$ be the ring of \ig s introduced in Theorem {\rm4.5.8}. Then the \fw\ assertions hold\/{\rm:}

\hph i,i, The ring $\Z$ is a domain.

\hph ii,, \E\te s \ooo \og~$\ge$ on~$\Z$ \st $(\Z,\ge)$ is an ordered domain. \Mo we have
\beq 3.73
a\ge b \qh{if } a-b\in\N \hbox{ \fa} a,b\in\Z.
\e
\elm

\proof
(i) follows from Theorem 4.5.8\,(ii).

(ii) The ``\ex'' part follows from \E\Pr\ \rf{p3.13}\,(i) with $F:=\Z$, $H:=\Na$ (see \E\df\ 4.3.89). Note that if $z\in\Na$, then $z+(-z)=0$ by
(4.3.127). Hence $-z$~is the inverse of~$z$ in the group $\pz1\Z$, which, together with \E\df\ 4.3.89, implies that \er{3.54} is \sf ied. \Mo the
first part of \er{3.55} follows from (4.3.125). We now prove the second part of \er{3.55}. Let $z\in H$. Then $z\ne0$, $\sgn(z)={+}$, and $\{z\}=z$
where $\sgn$ is defined in (4.5.35), $|z|$ in (4.5.36). Let $a,b\in H$. Then $a\cdot b\in \Z\sms0$ with $\sgn(a\cdot b)=\sgn(a)\cdot\sgn(b) =
{+}\cdot{+} = {+}$, and $|a\cdot b| = |a|\cdot|b| = a\cdot b$. We have $a\cdot b\in\Na$ since $\pz2{\Na}$ is a monoid. Thus $a\cdot b\in H$, which
proves \er{3.55}. Note that $1\in\Na=H$. Then \er{3.73} follows from \E\Pr\ \rf{p3.13}\,(i) since $H\cup\{0\} = \Na\cup\{0\}$. Note that we just
proved the last assertion of \E\Pr\ 4.3.92, which was left as an exercise (see Exercise 4.3.94).

``\E\uq'': Let $\wtge$ be an \og\ \st $(\Z,\wtge)$ is an ordered domain. Set $\Z_{\tilde\ge0}:= \{z\in\Z: z\wtge0\}$. Since $+$ and $\wtge$ are
compatible, we have $x\wtge y$ if $x+(-y)\in\Z_{\tilde\ge0}$ \fa $x,y\in\Z$ in view of \er{3.40}, \er{3.42} where $(G,\qu,e):=\pz1\Z$. Since $x\ge y$
if $x+(-y)\in\Z_{\ge0}$ \fa $x,y\in\Z$, we need to show that $\Z_{\ge0}=\Z_{\tilde\ge0}$. By part~(ii) of \E\Pr\ \rf{p3.13} we find that
\er{3.54}--\er{3.55} hold with $H:=\Z_{\tilde\ge0}$. Set $P:=H\cup\{0\}$. Hence $P=\Z_{\tilde\ge0}$. Clearly $0\in P$, and $x+y\in P$ whenever $x,y
\in P$ by \er{3.55}, since $0+x=x+0 = x\in P$ if $x\in P$, and $x+y\in H\sbs P$ if $x,y\in H$. \If that $P$~is a \sbm\ of the group $\pz1\Z$, and
that the monoid $\pz1P$ is a \PM\ as in the proof of \er{3.16}. We claim $1\in H$. We have $(-1)\cdot(-1)=1$. Indeed, since $(-1)+1=0$, we obtain
$0\nad{(4.1.1)}= (-1)\cdot 0 = (-1)\cdot(1+(-1)) \nad{(4.1.2)} = (-1)\cdot1 + (-1)\cdot(-1) = -1 + (-1)\cdot(-1)$. Hence $(-1)\cdot(-1)=1$. Clearly,
$1\cdot1=1$. Since $1\ne0$ and $-1\ne0$, we have either $1\in\Z_{\tilde>0}$ or $-1\in\Z_{\tilde>0}$. In the former case $1\in\Z_{\tilde>0}$ and in
the latter case $1=(-1)\cdot(-1)\in\Z_{\tilde>0}$ by~\er{3.55}. We claim that $H=\Na$. We proceed by \In. Let $A:=\{n\in\Na: n\dpl1 \in H\}$. We have
$1\in A$. Suppose $n\in A$. We show that $n+1\in A$. Indeed, $(n+1)\dpl1 \nad{(2.2.3)}= (n \dpl 1) + (1\dpl1) \nad{(2.2.3)} = (n\dpl1) +1 \in A$ by
\er{3.55} since both $1$~and~$n\dpl1$ belong to~$A$. \If that $A=\Na$. Since $n\dpl1=n$, $n\in\N$, by (2.2.8) and (2.2.13). \If that $\Na\sbs H$,
hence $\N\sbs P$. From \er{3.53}, we infer that $\Z$~is the disjoint union of $\Z_{\htg0}$, $\{0\}$, and~$\Z_{\htl0}$. Suppose, for \cd ion, that
\te s $z\ne0$ \st $z\notin\Na$ and $z\in P$. Then $z\in -\Na$, hence $-z\in \Na\sbs P$. \E\Tf $z,-z\in P$, and $z+(-z)=0$. Since $P$~is a \PM, we
have $z=-z=0$, a~\cd ion. \If that $\Na
=H$, hence $\Z_{\tilde>0} = \Na = \Z_{>0}$, and $\Z_{\wtge0}=\Z_{\ge0}$. This implies ${\wtge}={\ge}$, which completes the proof of the lemma.
\endproof

\Wanp prove \E\Pr\ 4.5.48. Let $\pz0\Q$ denote the field of \ra\ \nm s introduced in \E\df\ 4.5.7. We apply \E\Pr\ \rf{p3.13}\,(i) with $F:=\Q$ and
$H:=\bigl\{q\in\Q: q=\frac ab,\ a,b\in\Na\bigr\}$ where $\frac ab$ is defined in (4.5.1). In view of the discussion \fw\ \E\Pr\ 4.5.4, \as\ \er{3.54}
is \sf ied. \Mo if $P:=H\cup\{0\}$, then $\pz1P$ is a \PM\ by \E\Pr\ 4.5.4. \E\Tf if $r,s\in H$, then $r+s\in H\cup\{0\}$. However, if $r+s=0$, then
$r=s=0$, which \cd s $r,s\in H$. \E\Tf $r+s\in H$, and the first part of \er{3.55} is \sf ied. Finally, if $r+s\in H$,
then $r=+r$, $s=-s$ and $r\cdot s = (+r)\cdot(+s) \nad{(4.5.19)}= +r\cdot s$, where $r\cdot s$ is defined in (4.5.3). \If that \er{3.55} is \sf ied,
and the conclusions of \E\Pr\ 4.5.48 hold.

\newpage
\bds3.22 \

\hph i,i, Let $\pz0X$ be a \ti{\sr} (see \E\df\ 4.1.1, see also \er{2.1}--\er{2.5}). A~subset $Y$ of~$X$ is called a \ti{sub\sr} of~$X$ (see \E\df\
4.4.33) if $Y$~is a \sbm\ of the additive monoid $\pz1X$ and of the \mlv\ monoid $\pz2X$ of~$X$. Then $\pz0Y$ is a \sr. A~sub\sr~$Y$ of the \sr~$X$
is called a \ti{subdomain} of~$X$ if the \sr\ $\pz0Y$ is a domain (see \E\df\ 4.5.30). In that case $\pz1Y$ is an \ag, $Y\sms0$ is a \sbm\ of the
monoid $\pz2Y$ and $\pz2{Y\sms0}$ is a \Cm. A~sub\sr\ of~$X$ is called a \ti{subfield\/} of~$X$ if $\pz0Y$ is a field. In other words, $Y$~is a
subfield of~$X$ if $Y$~is a subdomain of~$X$ \st $\pz2{Y\sms0}$ is an \ag. Compare this \df\ with the one given in \E\df\ 4.4.24 (see Lemma 4.4.25).

\hph ii,, Let $X,X'$ be \sr s, and let $f:X\to X'$ be a \ti{ring-\hm sm} (see \E\df\ 4.1.7). The ring-\hm sm~$f$ is called a \ti{ring-iso\mf}
(resp.\ \ti{endo\mf}, resp.\ \ti{auto\mf}) if $f$~is bi\jc\ (resp.\ if $X=X'$, resp.\ if $X=X'$ and $f$ is bi\jc).
\eds

\brm3.21
It is customary (see \cite[p.~39]{Alg}, \cite{Is}) \ti{not\/} to assume that a ring-\hm sm from a ring~$R$ to a ring~$R'$ maps the unity of~$R$ into
the unity of~$R'$, if they exist.
\erm

\blm3.23
Let $D$ be a domain, let $F$ be a field and let $j:D\to F$ be an in\jc\ ring-\hm sm. Then $j(D)$ is a subdomain of the field~$F$.
\elm

\proof
In Lemma 5.4.22 we proved that if $D$ is a field in Lemma \rf{l3.23}, then $j(D)$ is a subfield of the field~$F$. We refer the reader to parts
(a)~and~(b) to show that $j(D)$ is a subgroup of the group $\pz1F$, and that $j(D)$ is a \sbm\ of the monoid $\pz2F$. We replace part~(c) as follows:

We first prove that $j(D)\sms{j(0)}$ is a \sbm\ of the group $(F,\cdot,j(1))$. Since $0\ne1$ and $j$~is in\jc, we have $j(1)\in j(D)\sms{j(0)}$. Now
let $x,y\in j(D)\sms{j(0)}$. Then \te\ $a,b\in D\sms0$ \st $x=j(a)$, $y=j(b)$. \If that $x\cdot y= j(a)\cdot j(b)= j(a\cdot b)\in j(D)$. \Mo
$a\cdot b\ne0$ since $D$~is a domain. Hence $x\cdot y = j(a\cdot b)\in D\sms{j(0)}$. Thus $j(D)\sms{j(0)}$ is a \sbm\ of the monoid $\pz2{F\sms0}$.
Since $\pz2{F\sms0}$ is a group, it is a \Cm. \E\Tf the monoid $(j(D)\sms0,\cdot,j(1))$ is a \Cm. Indeed, let $x,y,z\in j(D)\sms0$, $a,b,c \in
D\sms0$ be \st $x=j(a)$, $y=j(b)$, $z=j(c)$. If $x\cdot z=y\cdot z$, then $j(a)\cdot j(c)= j(b)\cdot j(c)$. Hence $j(a\cdot c)=j(b\cdot c)$, and
$a\cdot c=b\cdot c$ by in\ji\ of~$j$. \If that $a=b$ since a group is a \Cm. Hence $x=j(a)=j(b)=y$. \If that $j(D)$ is a subdomain of~$F$, hence
$(j(D),+,\cdot,j(0),j(1))$ is a domain, where $+,\cdot$ are the binary \op s in~$F$.
\endproof

Let $D,F,j$ be as in Lemma \rf{l3.23}. Suppose, in \ad, that $(D,\ge)$ is an ordered domain. The map $j:D\to F$, being in\jc, is a bi\jn\ from~$D$
onto $j(D)$. Let $(E,\ge)$ be an \os\ and let $E'$ be a \ns\ \ep\ to~$E$. Let $f:E\to E'$ be a bi\jn. If \te s an \og~$\ge'$ on~$E'$ \st both
$f$~and~$f\Inv$ are in\cre, the map~$f$ is called an \tb{\ois sm} (see \E\df\ 1.3.32). If the \og~$\ge$ on~$E$ is total, then \te s an \og~$\ge'$
which makes the bi\jn~$f$ an \ois sm. Indeed, let $\ge'$ be the \rl\ on~$E'$ defined by\index{order-isomorphism}
\beq 3.74
x' \mathrel{\ge'} y' \qh{if } f\Inv(x')\ge f\Inv(y'), \q x',y'\in E',
\e
then $\ge'$ is an \og\ on $E'$ and $f$ is an \ois sm between $(E,\ge)$ and $(E',\ge')$.

``\ti{\E\tr ity}'': Let $x'\mathrel{\ge'}y'$ and $y'\mathrel{\ge'}z'$, then $f\Inv(x')\ge f\Inv(y')$ and $f\Inv(y') \ge f\Inv(z')$ implies $f\Inv(x')
\ge f\Inv(z')$. Hence $x'\mathrel{\ge'}z'$.

``\ti{Reflexivity}'': Let $x'\in E'$. Then $f\Inv(x')\ge f\Inv(x')$, hence $x'\mathrel{\ge'}x'$.

``\ti{Anti\sy y}'': Let $x',y'\in E'$ be \st $x'\ge'y'$ and $y'\ge'x'$. Then $f\Inv(x')\ge f\Inv(y')$ and $f\Inv(y')\ge f\Inv(x')$ implies
$f\Inv(x')=f\Inv(y')$, hence $x'=y'$.

``\ti{$f$ is in\cre}'': Let $x,y\in E$ with $x\ge y$. Set $x':=f(x)$, $y':=f(y)$, then $f\Inv(x')\ge f\Inv(y')$, hence $x'\mathrel{\ge'}y'$ by
\er{3.74}. Thus $f(x)\mathrel{\ge'}f(y)$.

``\ti{$f$ is strictly in\cre}'': Since $f$ is in\jc\ and in\cre.

``\ti{$f\Inv$ is in\cre}'': It is \sft\ to show that $x'>'y'$ implies $f\Inv(x')\ge f\Inv(y')$. Suppose, for \cd ion, that $f\Inv(x')\not\ge
f\Inv(y')$. Since $\ge$ is total, we have $f\Inv(x')< f\Inv(y')$. Since $f$~is strictly in\cre,
we have $x'=f(f\Inv(x'))< f(f\Inv(y')) = y'$. Hence $y'>'x'>'y'$, hence $y'\ge'x' \ge'y'$, and $x'=y'$ by anti\sy y.
A~\cd ion, since $x'>'y'$. \If that $f\Inv$ is in\cre\ (hence strictly in\cre), which completes the proof of the claim.

Summarizing, we obtain

\blm3.24
Let $(E,\ge)$ be a {\em totally ordered} set, let $E'$ be a set and let $f:E\to E'$ be a bi\jn. Then the \rl\ on~$E'$ defined in \er{3.74} is an \og\
on~$E'$, and the map $f:(E,\ge)\to(E',\ge')$ is an \ois sm. \E\Ip $f$~and~$f'$ are strictly in\cre\ and the \og~$\ge'$ is total.
\elm

\proof
In view of the discussion preceding the lemma, it is \sft\ to show that the \og~$\ge'$ is total. Let $x',y'\in E'$ with $x'\ne y'$. Then $f\Inv(x')
\ne f\Inv(y')$, since $f$~is a bi\jn. The \og~$\ge$ on~$E$ is total. Hence either $f\Inv(x')>f\Inv(y')$ or $f\Inv(x')<f\Inv(y')$. Since $f$~is
strictly in\cre, we have either $x' = f(f\Inv(x')) > f(f\Inv(y')) = y'$ or $x' = f(f\Inv(x')) < f(f\Inv(y')) = y'$.
\endproof

\brm 3.25
Let $(E,\ge)$ be a \tos, let $E'$ be a \ns\ and let $j:E\to E'$ be in\jc. Then \fe $x'\in j(E)$ \te s \ooo $x\in E$ \st $x'=j(x)$. \If that the
\rl~$\ge'$ on~$j(E)$ defined by
\bga3.75
x'\ge'y' \qh{if }j\Inv(x')\ge j\Inv(y') \hbox{ \fa}x',y'\in j(E),\\
\hbox{where $j\Inv(z')$, $z'\in E'$, is the only \el\ of~$E$ \st $j(j\Inv(z'))=z'$,}\non
\e
is well-defined. Thus replacing $E'$ by $j(E)$, we may apply Lemma \rf{l3.24} and find that $\ge'$ is a total \og\ on $j(E)$ and that $j$~is an \ois
sm between $(E,\ge)$ and $(j(E),\ge')$.
\erm

\blm3.26
Let $(F,\ge)$ be an \of, and let $D$ be a subdomain $($resp.\ subfield\/$)$. Then $(D,\ge_D)$ is an ordered domain $($resp.\ field\/$)$ where $\ge_D$
denotes the \rt ion of the \og~$\ge$ to~$D$.
\elm

\proof
We apply \E\Pr\ \rf{p3.16} in case $D$~is a subdomain and \E\Pr\ \rf{p3.13} in case $D$ is a subfield. Set $H:=\{z\in D: z>0\}$. Note that $-H =
\{z\in D: z<0\}$. Indeed, using (4.5.120) for the \of~$F$, we obtain $z<0$ iff $0 = z+(-z) < 0+(-z) = -z$ \fa $z\in F$. Since the \og~$\ge$ is total,
$F$~is the disjoint union of $F_{>0}$, $\{0\}$, and~$F_{<0}$. \If that $D$~is the disjoint union of $H$, $\{0\}$, and~$-H$, which proves \er{3.54}.
Let $x,y\in H$. Then $x,y\in D\cap F_{>0}$. Then $x+y\in D$ since $\pz1D$ is a monoid, and $x+y\in F_{>0}$, $x\cdot y\in F_{>0}$ by \er{3.11}. Hence
$x+y$ and $x\cdot y$ belong to~$H$, which proves \er{3.55}. Then $(D,\ge)$ is an ordered domain (resp.\ field) with the \og\ $\ge_D$ defined
in~\er{3.56}.
\endproof

\Wanp prove the promised theorem about ``\uq'' of the \og\ of an \of\ of \qt s.
\newpage

\bth3.27 \

\hph i,i, Let $(D,\ge)$ be an \od\ and let $(\wh D,j)$ be the embedding of~$D$ in the field of \qt s of~$D$ introduced in Theorem {\rm4.5.35}. Let
$\ge_D$ denote the \og\ on~$j(D)$ defined by
\beq 3.76
j(x) \ge_D j(y) \qh{if } x>y \hbox{ \fa}x,y\in D.
\e
Then $j(D)$ is a subdomain of the field~$\wh D$, the domain $j(D)$ equipped with the \og~$\ge_D$ is an \od, and $j:(D,\ge) \to (j(D),\ge_D)$ is an
\ois sm.

\hph ii,, There is at most one \og\ $\whge$ on~$\wh D$ \st $(\wh D,\whge)$ is an \of\ iff there is at most one \og~$\ge$ on~$D$ \st $(D,\ge)$ is an
\od.
\eth

\proof
(i) Since $j$ is an in\jc\ ring-\hm sm, $j(D)$ is a subdomain of~$\wh D$ by Lemma \rf{l3.23}. Since the \og~$\ge$ is total and $j$~is in\jc, the \rl\
on $j(D)$ defined by \er{3.76} is a total \og\ on $j(D)$ and $j:(D,\ge) \to (j(D),\ge_D)$ is an \ois sm by Remark \rf{r3.25}. It remains to prove
(4.5.120) and (4.5.121). Let $x,y\in D$ be \st $j(x)<_D j(y)$. Then $x<y$, since $j$~is an \ois sm. Let $a\in D$ and $c\in D$ be \st $j(c)>_D\wh 0$.
Then $c>0$, for the same reason. Since $(D,\ge)$ is an \od, we have $x+a>y+a$ and $c\cdot x>c\cdot y$. Then $j(x)\hpl j(a) = j(x+a)>_D j(y+a) =
j(y)\hpl j(a)$, which implies (4.5.120) for $(j(D),\ge_D)$, and $j(c)\hcd j(x) = j(c\cdot x) >_D j(c\cdot y) = j(c)\hcd j(y)$, which implies
(4.5.121) for $(j(D),\ge_D)$.

\ssk
(ii) ``\ti{If\/}'': We suppose that the \og\ $\ge$ on $D$ is the only one which makes $(D,\ge)$ an \of, and we suppose, for \cd ion, that there are
two distinct
\og s $\whge$ and~$\wtge$ on~$\wh D$ which make $(\wh D,\whge)$ and $(\wt D,\wtge)$ \of s. In view of \er{3.59}, \te s $\wh z{}'\in \wh D$ \st
$\wh z{}'\in \wh D_{\hat>\hat0}$ and $\wh z{}'\notin \wh D_{\tilde>\hat0}$, or \te s $\wh w\in \wt D$ \st $\wh w\in \wh D_{\tilde>\hat0}$ and
${\wh w\notin \wh D_{\hat>\hat0}}$. In the first case, we have $\wh z{}'\whg\wh 0$ and $\wh z{}'\wtl\wh0$, since $\wt z{}'\ne\wh0$ and the
\og~$\wtge$ is total. Similarly, in the second case, $\wh w\wtg \wt0$ and $\wh w\whl\wh0$. Note that $\wh w\wtg\wh0$
iff $-\wh w\whl0$, and $\wh w\whl \wh0$ iff $-\wh w\wtg\wh0$.
\E\Tf we may assume that \te s $\wh z\whg\wh0$ \st $\wh z\wtl\wh 0$. By (4.5.88) \te s $(a,b)\in D\t D\sms0$ \st $\wh z = j(a)\hcd j(b)\Inv$. Since
$\wh z\ne\wh 0$ and $j$~is in\jc, we have $a\ne0$, thus $a,b,j(a),j(b),j(b)\Inv \in D\sms0$. \csq, $j(b)\hcd j(b)\in \wh D\sms{\wh0}$ since
$\wh D\sms{\wh0}$ is a \sbm\ of $(\wh D,\hcd,\wh1)$. \Mo $j(b)\hcd j(b)\whg\wh0$ by \er{3.12}. \E\Tf $j(a)\hcd j(b) = j(a)\hcd \wh1\hcd j(b) =
j(a)\hcd (j(b)\Inv \hcd j(b))\hcd j(b) \nad{(2.1.36)}= (j(a)\hcd j(b)\Inv)\hcd (j(b))^2\whg\wh0$, by \er{3.11}. Thus $j(a\cdot b)=j(a)\hcd j(b)\whg
\wh0$. Since $\wh z\wtl\wh0$, we obtain in the same fashion that $j(a)\hcd j(b)\wtl\wh0$. Hence $j(a\cdot b)\whg\wh0$ and $j(a\cdot b)\wtl\wt0$. Since
$j:D\to j(D)$ is a bi\jn, \te s a map $k:j(D)\to D$ \st $j\circ k=\id_{j(D)}$ and $k\circ j=\id_D$. \E\Ip $k$~is a bi\jn. Note that $k$~is a
ring-\is sm from the domain $(j(D),\hpl,\hcd,\wh0,\wh1)$ onto the domain $\pz0D$. In view of Lemma \rf{l3.24}, \te\ \og s $\ge'$ and $\ge''$ on~$D$
\st $j:(D,\ge')\to (j(D),\whge)$ and $j:(D,\ge'')\to (j(D),\wtge)$ are \ois sms. We next show that both $(j(D),{\whge}{}')$ and $(j(D),{\whge}{}'')$
are \od s. Since $\whge$ and $\wtge$ are total \og s, the \og s $\ge'$ and $\ge''$ are also total. We next show that both
\og s \sf y (4.5.120) and (4.5.121). Let $x,y,a
\in D$ be \st $x>'y$ (resp.\ $x>''y$). Then $j(x)\whg j(y)$, hence $j(x)\hpl j(a)\whg j(y)\hpl j(a)$ since $(j(D),\wh\ge|_{j(D)})$ is an \od\ by Lemma
\rf{l3.26}. \If that $j(x+a) = j(x)\hpl j(a) \whg j(y)\hpl j(a) = j(y+a)$, hence $x+a >' y+a$ (resp.\ $x+a>''y+a$) since $j$~is an \ois sm. The proof
of $x\cdot c\ge' y\cdot c$ (resp.\ $x\cdot c>'' y\cdot c$) whenever $x>'y$ and $c'>0$ (resp.\ $x>''y$ and $c>''0$) is similar. Thus both $(D,\ge')$
and $(D,\ge'')$ are \od s. By \as\ the \og s $\ge'$~and~$\ge''$ are equal. \If that \fa $x,y\in D$, $x\ne y$, we have $j(x)\whg j(y)$ iff $x>'y$ iff
$x>''y$ iff $x\whg y$. A~\cd ion, since we have $j(a\cdot b)=j(a)\hcd j(b)\whg\wh0$ and $j(a\cdot b)=j(a)\hcd j(b) \wtl\wh0$.

``\ti{Only if\/}'': We suppose that $(D,\ge)$ and $(D,\ge')$ are distinct \od s, and let $\whge$ (resp.~$\whge'$) be the \crs\ \og s on~$\wh D$
introduced in Theorem \rf{t3.17} \er{3.66}. We shall suppose that these \og s are equal and obtain a \cd ion. As in the ``if'' part of the proof, we
find that \te s $z\in D\sms0$ \st $z>0$ and $z<'0$, \ev tly $-z>'0$. Then $j(z)\in \wh D_{\htg\hat0}$ and $j(-z)\in \wh D_{\htg'\hat0}$ by \er{3.68}.
In view of \er{3.66}, \er{3.67}, the \og s $\whge$ and $\whge'$ are equal iff $\wh D_{\htg\hat0} = \wh D_{\htg'\hat0}$. Suppose for \cd ion that
$\wh D_{\htg\hat0} = \wh D_{\htg'\hat0}$, then $\wh0 = z\hpl (-z) \in \wh D_{\htg'\wh0}$ by~\er{3.11}. A~\cd ion, since $\wh0\notin \wh
D_{\htg'\hat0}$.
\endproof

Combining Theorems \rf{t3.17}, \rf{t3.27} and Lemma \rf{l3.21}, we find that \te s \ooo \og~$\whge$ on~$\wh \Z$, the field of \qt s of the
domain~$\Z$, introduced in Theorem 4.5.35. \Mo the \og~$\whge$ is defined in \er{3.66}. \If from \E\Pr\ 4.5.48, whose proof is given above, as a
con\sq\ of \E\Pr\ \rf{p3.13}\,(i), that the field $\Q$ of \ra\ \nm s is an \of\ with the \og~$\ge$ defined in (4.5.118). In \E\Pr\ 4.5.40 we
introduced a ring-\is sm $\th:\Q\to \wh \Z$ defined in (4.5.97). By \E\Pr\ \rf{p3.2}, $\th$~is the only ring-\is sm from~$\Q$ into~$\wh \Z$. It turns
out that $\th$~is also an \ois sm.

\blm3.28
Let $(F,\ge)$ and $(F',\ge')$ be \of s and let $\vf:F\to F'$ be a ring-\is sm. Then $\vf$~is also an \ois sm iff
\beq3.77
\vf(F_{>0}) \sbs F'_{>'0'}.
\e
\elm

\proof
``\ti{Only if\/}'': The map $\vf$ is bi\jc\ and in\cre, hence strictly in\cre. \E\Tf if $x>0$, $x\in F$, then $\vf(x)>'0'$. Hence $\vf(F_{>0})\sbs
F'_{>'0'}$, and \er{3.77} holds.

``\ti{If\/}'': Let $x,y\in F$ be \st $y>x$. Then $y=x+z$ \fs $z\in F_{\ge0}$ by Lemma \rf{l3.11} and \er{3.9}. We have $z\ne0$, otherwise $y=x$, \cd
ing $y>x$. Hence $z\in F_{>0}$. \E\Tf $\vf(y) = \vf(x+z) = \vf(x)+\vf(z)$ where $\vf(z)\in F'_{>'0'}$ by \er{3.77}. Thus $\vf(y)>' \vf(x)$, hence
$\vf$~is strictly in\cre. Let $x',y'\in F'$ be \st $y'>'x'$. We claim $\vf\Inv(y')>\vf\Inv(x')$. Suppose, for \cd ion, that $\vf\Inv(y')\not>
\vf\Inv(x')$, then, since $\ge$~is a total \og, we have either $\vf\Inv(y')=\vf\Inv(x')$ or $\vf\Inv(y')<\vf\Inv(x')$. In the former case, we would
have $y'=x'$, \cd ing $y'>'x'$, and in the latter case, we would have $y'=\vf(\vf\Inv(y')) <' \vf(\vf\Inv(x'))=x'$, hence $y'<'x'$, also \cd ing
$y'>'x'$.
\endproof

\brm3.29
We could have replaced \cn\ \er{3.77} by
\beq3.78
\vf(F_{>0}) = \vf(F'_{>'0'}).
\e
Indeed, clearly \er{3.78} implies \er{3.77}. \E\oh under the \as\ of the lemma, $\vf\Inv$~is a ring-\is sm from $F'$ onto~$F$, hence $\vf\Inv
(F'_{>'0'}) \sbs F_{>0}$. Since $\vf$~is a bi\jn, we also have $F'_{>'0'}\sbs \vf(F_{>0})$, which together with \er{3.77} implies \er{3.78}.
\erm

We now apply Lemma \rf{l3.28} to the case $(F,\ge):= (\Q,\ge)$, $(F',\ge'):=(\wh Z,\whge)$ and $\vf:=\th$. We find by (4.5.119), \er{3.67} and
(4.5.97):
\beag
{}&\Q_{>0} = \biggl\{+\frac ab\in\Q: a,b\in\Na\biggr\}, \\
&\wh Z_{\htg\hat0} = \biggl\{\frac ab\in\wh Z: a,b\in\Z_+=\Na\biggr\}, \\
&\th\Bgg(+\frac ab) = \frac ab, \q a,b\in\Na.
\e
Thus $\th(\Q_{>0}) = \wh\Z_{\htg\hat0}$. \If that the \of s $(\Q,\ge)$ and $(\wh \Z,\whge)$ are \ois c, hence the \og~$\ge$ of~$\Q$ is the only one
which makes $(\Q,\ge)$ an \of.

We summarize in

\bth3.30
Let $\Q$ be the field of \ra\ \nm s introduced in \E\df\ {\rm 4.5.7} and let $\ge$ denote the \og\ on~$\Q$ introduced in {\rm(4.5.118)}. Then
$(\Q,\ge)$ is an \of\ in the sense of \E\df\ {\rm4.5.50}. \Mo the \og~$\ge$ is the only \og\ making $(\Q,\ge)$ an \of.
\eth

\bex3.31
Let $a,b,c,d,a',b',c',d'\in \Na$. Suppose $a'b=ab'$ and $c'd=cd'$. Show
\beq3.79
ad>bc \qh{iff } a'd'>b'c'.
\e
\eex

So far the field of \ra\ \nm s is the only \of\ considered in this book. By Lemma \rf{l3.26} every subfield of an \of\ is an \of. However, the field
of \ra\ \nm s does not have a proper subfield by Lemma 4.5.42. It turns out that $\Q[X]$, the set of \fp s over~$\Q$, more generally, over
a~field~$K$, is a domain (see the last sentence of Section~2 of this Appendix), and that $K[X]$ is an \od\ whenever $K$~is an \of. \csq, in view of
Theorem \rf{t3.17}, the field of \qt s $K(X)$, where $K$~is an \of, is an \of.

\bpr3.32
Let $K$ be an \of\ with \tb{\po\ cone} $K_{>0}$. Then $K[X]$ is an \od\ with \po\ cone defined in
\beq3.80
K[X]_{>0} := \{a\in K[X]\sms0: a(\deg(a))>0\},
\e
where $\deg(a)$ is defined in {\rm(5.2.11)}.

\E\Ip \fa $a,b\in K[X]$ we have
\beq3.81
a\ge b \qh{if } a+(-b)\in K[X]_{>0} \cup\{0\}.
\e
\epr

\proof
We apply the first part of \E\Pr\ \rf{p3.16}, that is, \E\Pr\ \rf{p3.13}\,(i) where $F$~is a~domain (instead of a~field). We recall that an \of~$K$
is infinite, since $\1\in K_{>0}$ and the map $n\mt n\dpl1$ from~$\N$ into $K_{\ge0}$ is in\jc. \Mo the map $j:K \to K[X]$ defined by $j(0)=0\in
K[X]$, $j(a)(0):=a\in K\sms0$, and $j(a)(k):=0$ \fa $k\in\Na$, $a\in K\sms0$, is an in\jc\ ring-\hm sm. Hence \Ip $K[X]$ is infinite. We set
$H:=K[X]_{>0}$ defined in \er{3.80}. We first observe that if $a\in K[X]\sms0$ then $\deg(a)$ is well-defined, $\deg(a)\in\N$, and $a(\deg(a))\in
K\sms0$. If $a(\deg(a))>0$, then $a\in H$. If $a(\deg(a))<0$, then $\deg(-a)=\deg(a)$ by (5.2.13), where $-a$ is the \fp\ defined by
\beq3.82
(-a)(j) = -(a(j)), \q j\in\N.
\e
\E\Tf we have $(a+(-a))(j) = a(j)+(-a)(j) = a(j)+(-(a(j))) =0$ in~$K$ \fa $j\in\N$. Hence $-a$ is the inverse of~$a$ in the group $(K[X],+,0)$. \Mo
$(-a)(\deg(-a)) = (-a)(\deg(a)) \nde3.82 = -(a(\deg(a)))$. \If that $(-a)(\deg(-a))>0$ in~$K$. Thus $-a\in H$, \ev tly $a\in -H$. \csq, if $a\in K[X]
\sms0$, then either $a\in H$, or $a\in -H$, which shows that \er{3.54} is \sf ied. It remains to verify \er{3.55}. Let $a,b\in H$. From (5.2.29),
(5.2.30) we infer that $a\cdot b\ne 0$, and $(a\cdot b)(\deg(a\cdot b)) = a(\deg(a))b(\deg(b))$ in $K\sms0$. Since $a(\deg(a))>0$ and $b(\deg(b))>0$,
we find that $a\cdot b\in H$. Finally, we show that $a+b\in H$. If $\deg(a)=\deg(b)$, then $(a+b)(\deg(a)) = a(\deg(a))+b(\deg(b))>0$, hence
$a+b\ne0$, and $\deg(a+b)$ is well-defined. By (5.2.12) we have $\deg(a+b)\le \max(\deg(a),\deg(b)) = \deg(a)$. Since $(a+b)(\deg(a+b))>0$, we have
$\deg(a+b)=\deg(a)$. Hence $(a+b)(\deg(a+b))>0$, and $a+b\in H$. If $\deg(a)\ne \deg(b)$, we may assume \wlg $\deg(a)>\deg(b)$ since $a+b=b+a$ and
the \og\ of~$\N$ is total. In that case, $b(\deg(a))=0$, hence $(a+b)(\deg(a)) = a(\deg(a))+b(\deg(a)) = a(\deg(a))+0 > 0$, since $a\in H$. \E\Ip
$a+b\ne0$, hence $\deg(a+b) \nad{(5.2.13)}= \max(\deg(a),\deg(b)) = \deg(a)$. Since $a(\deg(a+b))\ne0$, $\deg(a+b)=\deg(a)$ and $(a+b)(\deg(a+b))=
a(\deg(a))>0$. \E\Tf $a+b\in H$, and \cn\ \er{3.55} is \sf ied. Finally, \er{3.81} is just a reformulation of \er{3.56}.
\endproof

\bex3.33 \

\hph i,i, Let $a,b\in K[X]$ where $K$ is an \of. Find an explicit formulation of the \cn\ $a+(-b)\in K[X]_{>0}\cup\{0\}$ in \er{3.81} in terms of the
\cf s $\{a(j)\}_{j\in\N}$ and $\{b(j)\}_{j\in\N}$.

\hph ii,, Assume, moreover, that $a,b\in K[X]\sms0$. Define $\sgn c={+}$ if $c(\deg (c))>0$ and $-$ if $c(\deg(c))<0$. Show that
\beq 3.83
\frac ab \whg0 \hbox{ in }K(X) \qh{iff } \sgn a=\sgn b,
\e
where $\whge$ denotes the \og\ on $K(X)$ introduced in \er{3.66} where $(D,\ge):=(K[X],\ge)$ with $\ge$ defined in \er{3.81}.
\eex

\bex3.34
Let $a$ denote the \fp\ $X^2+1$ (in our notation: $\ve_2+\ve_0$) in $\Q[X]$.

\hph i,ii, Show that $a$ is irreducible in $\Q[X]$.

\hph ii,i, Let $(\Q_a,+_a,\cdot_a,0,\ve_0)$ denote the $K$-algebra \itd in \E\Pr\ 5.3.7. Use Theorem 5.3.31 and \E\Pr\ \rf{p3.5} in order to show that
there is \ti{no} \og\ on the field $(\Q_a,+_a,\cdot_a,0,\ve_0)$ which makes this field an \of.

\hph iii,, Show that the same conclusion holds if we replace $\Q$ by an arbitrary \of~$K$.
\eex

\If from Lemma 4.5.52 that the \of\ $(\Q,\ge)$ is \ti{\Ar} (see \E\df\ 4.5.53). It turns out that if $(D,\ge)$ is an \od, then \E\Pr\ \rf{p3.5}
holds, that is, \er{3.11}, \er{3.12} hold when $F$~is an \od. Indeed, \er{3.11} follows from \E\Pr\ \rf{p3.16} and the proof of \er{3.12} is word for
word identical to the proof for~$F$ being an \of. \E\Tf \E\df\ 4.5.53 could be extended to an \od. Observe that \cn s (i) and~(iii) in Lemma
4.5.54 are \ev t in an \of. Indeed, (iii) implies (i) by taking $x:=1$ if $y>1$. If $y\le1$ then $1<1+1$, hence $y<k\dpl1$ with $k:=2$.
Conversely, if $0<x<y$, $x,y\in F$, then $y\cdot x\mo < n\dpl 1$ \fs $n\in\Na$
by~(i), hence by \ml ying the in\et y by~$x$ ($>0$), we obtain $y<(n\dpl1)\cdot x$. Using \In\ on $m\in\N$, one finds in a \sr\ $(X,+,\cdot,0,\1)$:
\beq3.84
(m\dpl\1)\cdot a = m\dpl(\1\cdot a) = m\dpl a, \q a\in X,\ m\in\N.
\e
\E\Tf one finds $y< n\dpl x$. This motivates the \fw

\bdf 3.35 \
An \od\ $(D,+,\cdot,0,1,\ge)$ is called \tb{\Ar} if \fe pair $(x,y)\in D_{>0}\t D_{>0}$ \sf ying $x<y$ \te s $n\in\Na$ \st $y< n\dpl x$, where
$n\dpl x$ denotes the $n$-th \IT\ of~$x$ in the \ag\ $\pz1D$.
\edf

In view of Lemma 4.5.54, \E\df\ \rf{d3.35} is \ev t to \E\df\ 4.5.53 when $D$ is a field.

\bxs 3.36 \

\hph i,i, $(\Z,+,\cdot,0,1,\ge)$ is \Ar. Indeed, note that $\Z_{>0}=\Na$. If $a,b\in\Na$ \sf y $a<b$, then by Theorem 2.1.38 and (2.2.8), \te\
$q\in\N$ and $r\in\N$ with $0\le r<a$ \st $b=q\cdot a+r$. Then $b<q\cdot a+a = q\cdot a + 1\cdot a = (q+1)\cdot a \nad{(2.2.8)}= (q+1)\dpl a$.
\csq, $b<n\dpl a$ with $n:=q+1\in \Na$.

\hph ii,, The \od\ $(\Q[X],\ge)$ \itd in \E\Pr\ \rf{p3.32} is \ti{not\/} \Ar. Let $a:=\ve_0$ and $b:=\ve_1$, see (5.2.4). In the usual notation
$a=1$ and $b=X$. Then $\psi_0(\ve_1-\ve_0)=-1$, $\psi_1(\ve_1-\ve_0)=1$ and $\psi_i(\ve_1-\ve_0)=0$ \fa $i\ge2$, where $\psi_i$ is defined in (5.2.2).
\E\Tf $\ve_1-\ve_0 \in\Q[X]_{>0}$ by \er{3.80} and $\ve_0<\ve_1$ by \er{3.81}. However, \fa $n\in\Na$, we have
$\psi_0(\ve_1-n\dpl\ve_0) = \psi_0(\ve_1) - n\dpl \psi_0(\ve_1) \hcd 0 - n \ddt -n$,
$\psi_1(\ve_1- n\dpl \ve_0) = \psi_1(\ve_1) - n\dpl \psi_1(\ve_0)= 1-0=1$, and $\psi_i(\ve_1 - n\dpl\ve_0) = \psi_i(\ve_1) - n\dpl\psi_i(\ve_0) =
0- 0 =0$ \fa $n\ge2$. \E\Tf $\deg(\ve_1- n\dpl\ve_0) = 1$ \fa $n\in\N$, and $(\ve_1 - n\dpl\ve_0)(1) = 1 > 0$ \fa $n\ge1$. \If that $\ve_1 - n\dpl
\ve_0 \in \Q[X]_{>0}$, hence $n\dpl\ve_0<\ve_1$ \fa $n\in\Na$, and there is \ti{no} $n\in\Na$ \st $n\dpl a>b$.
\exs

\bpr3.37
Let $(D,\ge)$ be an \Ar\ domain and let $\wh D$, $j$ and $\whge$ be as in Theorem \rf{t3.17}. Then $(\wh D,\whge)$ is \Ar\ iff $(D,\ge)$ is \Ar.
\epr

\proof
``\ti{If\/}'': Let $\wh z\in \wh D_{\htg\hat0}$ and suppose $(D,\ge)$ \Ar. In view of the discussion preceding \E\df\ \rf{d3.35} it suffices to find
$n\in\Na$ \st $\wh z\whl n\dhp\wh1$. By \er{3.67} \te\ $a,b\in D_{>0}$ \st $\wh z=\frac ba$. By \as\ and \E\df\ \rf{d3.35} \te s $n\in\Na$ \st
$b< n\dpl a$, that is, \te s $p\in D_{>0}$ \st $n\dpl a=b+p$, hence by \er{3.84} $(n\dpl1)\cdot a = b+p$ where $1$~is the \nel\ of the \mlv\ monoid
$\pz2D$. \E\ml ying both sides by $a\Inv = \frac1a$, we find
$$
\frac{(n\dpl1)\cdot a}{1} \cdot \frac 1a = \frac{(n\dpl1)\cdot a\cdot 1}{1\cdot a} = \frac{b+p}a = \frac ba\hpl \frac pa = \wh z \hpl \frac pa,
$$
where $\frac pa\in \wh D_{>0}$. From \er{3.66} we infer $\wh z\whl \frac{n\dpl1}1 \hcd \frac a1\hcd \frac 1a = \frac{n\dpl1}1 = j(n\dpl1)
\nad{(2.1.45)} = n\dhp j(1) = n\dhp \wh 1$.

``\ti{Only if\/}'': Let $x,y\in D_{>0}$ be \st $x<y$. Since $\frac yx\in\wh D_{>0}$ and $(\wh D,\whge)$ is \Ar, \te s $n\in\Na$ \st $\frac yx\whl
n\dhp\wh 1$, by Lemma 4.5.54. Using (4.5.121), we find $\frac yx \hcd \frac x1 \whl (n\dhp\wh 1)\hcd \frac x1$. Hence $\frac y1\whl n\dhp \frac x1$
by \er{3.84}. But $n\dhp \frac x1= n\dhp j(x) \nad{(2.1.45)}= j(n\hpl x)$. Since $\frac y1=j(y)$ and $j:(D,\ge) \to (j(D),\wh\ge|_{j(D)})$ is an
\ois sm, we obtain $y< n\dpl x$. Hence $(D,\ge)$ is \Ar.
\endproof

\bco3.38
The \of\ $(\Q(X),\ge)$ $($field of \qt s of $(\Q[X],\ge))$ is \emph{not} \Ar.
\eco

Note that by Lemma \rf{l3.26} every subfield of an \of\ is an \of. As a direct con\sq\ of \E\df\ \rf{d3.35} we obtain

\blm3.39
Every subfield of an \Ar\ \of\ is \Ar.
\elm

In the next section we shall ``construct'' an \Ar\ field $(\R,\ge)$ \st every \Ar\ \of\ is order- and ring-\is c to a subfield of~$\R$.

\newpage
\Subsubsection{The field $\R$}\label{ass.4}

In Exercise 4.5.56 it is claimed that every \ra\ \nm~$x$ is the supremum (see \E\df s and Notations 3.1.13) of the set of \ra\ \nm s~$y$ \st
$y<x$. We also have
\beq 4.1
x = \sup\{y\in\Q_{>0}: y<x\} \qh{\fa} x\in\Q_{>0}.
\e
It turns out that a similar result holds in an \Ar\ \of\ $(K,+,\cdot,0,\1,\ge)$. We recall that, as an \of, the field~$K$ has \ch istic zero (see
\E\df\ 4.4.37). Indeed, it follows from \er{3.12}, \er{3.11} and \In\ on $n\in\Na$ that $n\dpl\1 \in K_{>0}$ \fa $n\in\Na$. In view of \E\Pr\
\rf{p3.1}, the \pf\ of~$K$, denoted by $\check K$, is infinite and by Theorem 4.5.44 $\check K$ is ring-\is c to the field~$\Q$. More precisely, it
follows from the proof of Theorem 4.5.44 with $\pz0F := (\check K,+,\cdot,0,\1)$ that if $(z,w)$ and $(z',w')\in \Z\t\Z^\t$ \sf y $zw'=wz'$, then
(4.5.109) holds. Hence the map $j:\wh\Z\to \check K$ defined by (4.5.108), where $\wh \Z$ is the field \itd in Theorem 4.5.35, makes sense. \Mo the
map~$j$ is a ring-\is sm. Note that by \E\Pr\ \rf{p3.2}, the identity is the only ring-auto\mf\ of~$\Q$ \itd in \E\df\ 4.5.7, that $\Q$ and~$\wh\Z$
are ring-\is c by \E\Pr\ 4.5.40, hence the identity is the only ring-auto\mf\ of~$\wh\Z$. \If that the map~$j$ defined by (4.5.108) is the \ti{only}
ring-\is sm from~$\wh \Z$ onto~$\check K$ and that the identity is the only ring-auto\mf\ of~$\check K$. We recall that \te s \ooo \og~$\whge$ making
$(\wh\Z,\whge)$ an \of\ and that $\wh\Z_{\htg\hat0} = \bigl\{\frac ab\in\wh\Z: a,b\in\Na\bigr\}$. See the discussions \fw\ the proof of Theorem
\rf{t3.27} and preceding Theorem \rf{t3.30}. We next show that $j$~maps $\wh\Z_{\htg\hat0}$ into $\check K_{>0}$. Let $r\in\wh\Z_{\htg\hat0}$ and
$m,n\in\Na$ be \st $r=\frac mn$. Then $j(r)=(n\dpl\1)\cdot(m\dpl\1)\mo$ belongs to $\check K_{>0}$ since both $(n\dpl\1)$ and $(m\dpl\1)$ belong
to $\wh K_{>0}$, $(m\dpl\1)\mo \in\wh K_{>0}$ by Lemma \rf{l3.18} and $(n\dpl\1)\cdot(m\dpl\1)\mo \in\wh K_{>0}$ by \er{3.11}. By Lemma \rf{l3.28}
$j$~is also an \ois sm, and by Remark \rf{r3.29} we have
\beq4.2
j(\wh\Z_{\htg\hat0}) = \check K_{>0}.
\e
In the next \Pr\ we show that between two \el s of $K_{>0}$ \te s an \el\ of~$\check K_{>0}$.

\bpr4.1
Let $(K,+,\cdot,0,\1,\ge)$ be an \Ar\ \of, and let $x,y\in K_{>0}$ with $y<x$. Then \te\ $m,n\in\Na$ \st the \fw\ assertion holds\/{\rm:}
\beq4.3
y < (n\dpl\1) \cdot (m\dpl\1)\mo < x.
\e
\epr

\proof
We apply Lemma 4.5.54\,(ii) with $y:=x-y$ and we find $m\in\Na$ \sf ying
\beq4.4
(m\dpl\1)\mo < x-y.
\e
We recall that $m\dpl\1\in K_{>0}$ and that $(m\dpl\1)\mo \in K_{>0}$ by Lemma \rf{l3.18}. Set
$$
A:= \{k\in\Na: x\le k\dpl(m\dpl\1)\mo\}.
$$
We have $(m\dpl\1)\mo < x$ by \er{4.4} since $x-y<x$. By Lemma 4.5.54\,(iii) with $x:=\break(m\dpl\1)\mo$ and $y:=x$ \te s $\ov k\in\Na$ \st $x<\ov k
\dpl(m\dpl\1)\mo$. \E\Tf $A\ne\vn$. In view of Theorem 1.3.21, the set $A$ has a least \el\ $\wh k\in\Na$. Since $(m\dpl\1)\mo < x <\break \wh k\dpl
(m\dpl\1)\mo$, we have $\wh k\ge2$, and
\beq4.5
(\wh k-1)\dpl(m\dpl\1)\mo < x < \wh k\dpl (m\dpl\1)\mo.
\e
Set $n:=\wh k-1\in\Na$. Thus we obtain $n\dpl(m\dpl\1)\mo<x$.

\E\oh we have $y-(y-x) < y +(m\dpl\1)\mo -(y-x) \nde4.4 < x-(y-x) < (n+1)\dpl(m\dpl\1)\mo - (y-x)\nde4.4 < (n+1)\dpl (m\dpl\1)\mo-(m\dpl1)\mo
\nad{\rm(2.2.3)I1,I2}= n\dpl(m\dpl\1)\mo$. \E\Tf we arrive at $y<n\dpl(m\dpl\1)\mo < x$. Since $n\dpl(m\dpl\1)\mo \nde3.84 = (n\dpl\1)\cdot\break
(m\dpl\1)\mo$ we obtain \er{4.3}.
\endproof

We recall that $(n\dpl\1)\cdot(m\dpl\1)\mo$ in \er{4.3} is the image of the \po\ \ra\ \nm~$\frac nm$ under the map~$j$ defined in (4.5.108). \Wanp
prove a \gn\ of~\er{4.1}.

\bco4.2
Let $(K,+,\cdot,0,\1,\ge)$ be an \Ar\ \of\ and let $x\in K_{>0}$. Then the \fw\ assertion holds:
\beq4.6
x=\sup\{(n\dpl\1)\cdot (m\dpl\1)\mo: n,m\in\Na, \ (n\dpl\1)\cdot(m\dpl\1)\mo < x\}.
\e
\eco

\proof
Set $A:=\{(n\dpl\1)\cdot (m\dpl\1)\mo: n,m\in\Na, \ (n\dpl\1)\cdot(m\dpl\1)\mo < x\}$.

``$A\ne\vn$'':  By Lemma 4.5.54\,(ii) with $y:=x$ \te s $m\in\Na$ such that\break $(m\dpl\1)\mo<x$. Note that $1\dpl\1=\1$ by (2.2.3)\,I1. Hence
$(1\dpl\1)\cdot(m\dpl\1)\mo \in A$.

``\ti{$x$ is an \ub\ of $A$}'': obvious.

``\ti{$x$ is the least \ub\ of $A$}'': Suppose, for \cd ion, that \te s $y\in K_{>0}$ \sf ying $0<y<x$ and $a\le y$ \fa $a\in A$. By \E\Pr\
\rf{p4.1}, \te\ $m,n\in\Na$ \st $y<(n\dpl\1)\cdot(m\dpl\1)\mo < x$. Since $y<(n\dpl\1)\cdot(m\dpl\1)\mo \in A$, we obtain a \cd ion.
\endproof

\bex4.3 \

\hph i,ii, Show that \er{4.6} does not hold in the \of\ $\Q(X)$ \itd in Corollary \rf{c3.38}.

\hph ii,i, Show that under the \as s of Corollary \rf{c4.2} the \fw\ assertion holds\dw
\beq4.7
x = \inf\{j(r)\in K_{>0}: r\in\wh\Z_{\htg\hat0}\hbox{ and }j(r)>x\}.
\e

\hph iii,, Show that if $x,y\in\Q_{\ge0}$ and $y<x$, then $y<\frac12(x+y)<x$.
\eex

We recall that one of the goals of this section, mentioned at the end of Section \ref{ass.3}, is to prove the \ex\ of an \Ar\ field $(\R,\ge)$ \st
every \Ar\ \of\ is order- and ring-\is c to a subfield of~$\R$.

In a first step we shall prove that \te s an ordered semifield $(\R_{\ge0},\ge)$ (see \E\df~\rf{d3.7}) \st the \po\ cone of every \Ar\ field can be
embedded in an appropriate sense into $\R_{\ge0}$. To this end we shall first identify the \po\ \el s of an \Ar\ field~$K$ with subsets of~$\Q_{>0}$.
So far the field $(\Q,\ge)$ is the only \Ar\ field that we have introduced. We associate with a \po\ \ra\ \nm\ the order \il
\beq4.8
(0,x):=\{y\in\Q: 0<y<x\}.
\e
Such \il s, considered as subsets of $\Q_{>0}$, have the \fw\ \pp ies. Let us denote by~$I$ such an \il. We have:

C1 \ $I\ne\vn$ and $I\ne\Q_{>0}$.

C2 \ If $a,b\in\Q_{>0}$ \sf y $a\in I$ and $b<a$, then $b\in I$.

C3 \ \E\fe $a\in I$ \te s $c\in I$ \st $a<c$.

We recall that $(\Q,\ge)$, where the \og\ is defined in (4.5.118), (4.5.119), is an \of. \E\Ip (4.5.120) and (4.5.121) hold.
Let $x\in\Q_{>0}$ and $I:=(0,x)$. Then $x\notin I$, and $\frac12x\in I$ by Exercise \rf{ex4.3}\,(iii) with $y:=0$. Thus C1 holds.
Let $a,b\in\Q_{>0}$ be as in~C2. Then $b<a<x$, hence $b<x$ by (1.3.11), and $b\in I$. Finally, let $a\in I$. Set $c:=\frac12(a+x)$.
By the same Exercise as above with $y:=a$, we obtain $a<c<x$, since $a+c=c+a$. Thus C3 holds.

\ssk
A subset of $\Q_{>0}$ \sf ying C1--C3 is usually called a (\po) cut or a (\po) Dedekind-cut (see \cite{10a}).\index{positive Dedekind cut} We shall call these subsets
\ti{ideals} of $\Q_{>0}$,\index{ideal of $\Q_{>0}$} and the set of these ideals will be denoted by $\Id(\Q_{>0})$. We define a map $i:\Q_{>0}\to \Id(\Q_{>0})$ by setting:
\beq4.9
i(x):= (0,x) \qh{\fa} x\in\Q_{>0}.
\e

\blm4.4
The map $i:\Q_{>0}\to\Id(\Q_{>0})$ \sf ies
\beq4.10
a<b \qh{iff } i(a)\sbn i(b) \qh{\fa} a,b\in\Q_{>0}.
\e
\E\Ip the map $i$ is in\jc.
\elm

\proof
``\ti{Only if\/}'': Let $0<a<b$, $a,b\in\Q$, and $x\in(0,a)$. Then $0<x<a<b$, hence $x\in i(b)$. \Mo $a<\frac12(a+b)<b$ by Exercise \rf{ex4.3}\,(iii)
(see also \er{4.23}). Hence $\frac12(a+b)\in i(b)\sm i(a)$.

``\ti{If\/}'': Let $(0,a)\sbn(0,b)$, $a,b\in\Q$. \E\te s $c\in(0,b)$ \st $c\notin(0,a)$. Thus $c\ge a$. Hence $a\le c<b$, and $a<b$ by (1.3.11).

Note that $(\Id(\Q_{>0}),\subset)$ is an \os. The last sentence of the lemma follows from Lemma 1.3.34\,(ii).
\endproof

\blm4.5
The map $i$ defined in \er{4.9} is \emph{not} sur\jc.
\elm

\proof
In Examples 4.5.59, we showed that \fa \Pn s~$p$ and all $n\in\Na\sms1$, there is no $\a\in(0,1)\sbs\Q_{>0}$ \st the sets
\beq4.11
C_{n,p} := \bigl\{y\in\Q_{>0}: y^n<\tfrac1p\bigr\}
\e
are of the form $(0,\a)$. It remains to show that the sets $C_{n,p}$ \sf y C1--C3.

``C1'': Let $a:=\frac1p$. Then $a\in\Q_{>0}$, $p<p^n$, hence $\bigl(\frac1p\bigr)^n = \frac1{p^n} <\frac1p$, and $a\in C_{n,p}$.

Let $a:=1$. Then $1^n=1>\frac1p$, hence $a\notin C_{n,p}$.

``C2'': Let $a\in C_{n,p}$, and $0<b<a$. Then $0<b^n$ and $b^n<a^n$ by Examples 4.5.59. Hence $b^n<a^n<\frac1p$. Thus $b^n<\frac1p$ and $b\in
C_{n,p}$.

``C3'': Let $a\in C_{n,p}$. Then $a^n<\frac1p$. Let $h\in\Q_{>0}$. Then $(a+h)^n-a^n \nad{(4.4.46)}= h\suml_{k=0}^{n-1} (a+h)^{n-1-k}a^k$.
Choose $h\le a$. Then $(a+h)\le 2a$, and for $0\le k\le n-1$ we have $(a+h)^{n-1-k}a^k \le 2^{n-1-k}a^{n-1-k}a^k = 2^{n-1-k}a^{n-1} \le 2^{n-1}
a^{n-1}$. Set $M:=\suml_{k=0}^{n-1} 2^{n-1}a^{n-1}$. Then $(a+h)^n-a^n \le hM$, where $M>0$. Set $\a:=\frac1p-a^n>0$, and let $h$~be also less than
$\a M\mo$. We obtain $(a+h)^n = ((a+h)^n-a^n)+a^n \le hM+a^n < \a M\mo M+a^n = \a+a^n = \frac1p-a^n+a^n = \frac1p$. Then $c:=a+h\in C_{n,p}$ and
$a<c$.
\endproof

Our goal is to prove the \ex\ of an \of\ $(\R,\ge)$ \st the ordered semifield $P:=\{x\in \R:x\ge0\}$ \sf ies
\beq 4.12
\dot P=\Id(\Q_{>0}), \qh{where $\dot P:=P\sms0$.}
\e

As in \cite{10a} we shall use small Greek letters for \el s of~$\dot P$. Since $\dot P$~is a set of subsets of $\Q_{>0}$, inclusion induces
a (partial) \ti{\og} on~$\dot P$ (see Exercise 1.3.4\,(ii)), which we also denote by~$\ge$:
\beq 4.13
\a\ge \b \qh{if }\a\supset\b \qh{(possibly \et y).}
\e
There is no ambiguity in using the same symbol $\ge$ for the \og\ on $\Q_{>0}$ and the \og\ on~$\dot P$ (subsets of $\Q_{>0}$).

\blm 4.6
Let $\a,\b\in \dot P$ with $\a\ne\b$. Then either $\a>\b$ or $\b>\a$.
\elm

\proof
If $\a\ne\b$, then either \te s $a\in\a$ \st $a\notin\b$, or \te s $b\in\b$ \st $b\notin\a$. It suffices to show that if $\g,\d\in\dot P$ \sf y
$c\in\g$, $c\notin\d$, then $\g\ge\d$. Suppose, for \cd ion, that $\g\not\ge\d$. In that case, \te s $d\in\d$ \st $d\notin\g$. Note that $c\ne d$
since $c\in\g$ and $d\notin\g$. Since the \og\ of $\Q_{>0}$ is total, we have either $d>c$ or $c>d$. In the former case, $c\in\d$ by C2 since
$d\in\d$, a~\cd ion. In the latter case, $d\in\g$ by~C2 since $c\in\g$, a~\cd ion. If $\g\ne\d$, we obtain $\g>\d$.
\endproof

We now reformulate axioms C1--C3. To this end we introduce the notion of open subset of $\Q_{>0}$.

\bdf4.7
A subset $A$ of $\Q_{>0}$, possibly empty, is called \ti{open} if \fe $x\in A$ \te\ $a,b\in\Q_{>0}$, $a<b$, \st $x\in(a,b)\subset A$.\index{open}
\edf

\bxs4.8
The empty subset of $\Q_{>0}$; $\Q_{>0}$; $(a,b)$, $a,b\in\Q_{>0}$, $a<b$; $(a,b)\cup(c,d)$, $a,b,c,d\in\Q_{>0}$, $a<b<c<d$; $\{x\in\Q_{>0}: a<x\}$,
$a\in\Q_{>0}$, are open subsets of~$\Q_{>0}$. The subset $(a,b]:=\{x\in\Q_{>0}: a<x\le b\}$ is \ti{not\/} open.
\exs

We recall (see for example \cite{18}) that if $J$~is a \ns\ and that if \fe $j\in J$ there corresponds a subset~$A_j$ of a \ns~$A$ ($j\mt A_j$ is
a map from~$J$ into $\cP(A)$, the power set of~$A$), then $\cA:=\{A_j\in\cP(A):j\in J\}$ is called an \ti{indexed family} of subsets of~$A$.
The union of all subsets of~$A$ in~$\cA$, denoted by $\bcl_{j\in J}A_j$, is defined by
\beq 4.14
\bigcup_{j\in J}A_j := \{x\in A: x\in A_j \hbox{ \fs}j\in J\}.
\e

\blm4.9
The union of a family of nonempty open subsets of $\Q_{>0}$ is a nonempty open subset of~$\Q_{>0}$.
\elm

\proof
Let $\cA:= \{O_j\sbs \Q_{>0}: j\in J\}$ be an indexed family of nonempty open subsets of $\Q_{>0}$ indexed by~$J$, and let $x\in\bcl_{j\in J}O_j$.
Then \te\ $j\in J$ and $a_j,b_j\in\Q_{>0}$, $a_j<b_j$, \st $x\in O_j$ and $x\in(a_j,b_j)\subset O_j\sbs \bcl_{l\in J}O_l$. Hence $\bcl_{j\in J}O_j$
is a nonempty open subset of~$\Q_{>0}$.
\endproof

\Wanp give another \ch ization of ideals of $\Q_{>0}$.

\blm4.10
Let $I$ be a subset of $\Q_{>0}$. Then $I$~is an \tb{ideal\/} of~$\Q_{>0}$ iff the \fw\ \cn s are \sf ied\dw

{\rm I1} \ $I\ne\vn$ and $I\ne\Q_{>0}$.

{\rm I2} \ If $a,b\in\Q_{>0}$ \sf y $a\in I$ and $b<a$, then $b\in I$.

{\rm I3} \ $I$ is open.
\elm

\proof
Clearly, C1 (resp.\ C2) is identical to I1 (resp.~I2). C3~follows from I1 and~I3. Indeed, let $I\in\dot P$ and $a\in I$. In view of~I3 \te\
$x,y\in\Q_{>0}$ \st $x<y$, and $a\in(x,y)\sbs I$. Set $c:=\frac12(a+y)$. Then $x<a=\frac12(a+a)<\frac12(a+y)<\frac12(y+y)=y$. Hence $a<c$
and~$c\in I$. Conversely, C1, C2~and~C3 imply I3. Indeed, let $I\in \dot P$ and let $a\in I$. By~C3 \te s $c\in\Q_{>0}$ \st $a<c$ and $c\in I$, and
by~C2 the set $\{y\in\Q_{>0}: 0<y<c\}\sbs I$. Since $0<\frac12a<a$ and $a<c$, we have $a\in(\frac a2,c)\sbs (0,c)\sbs I$. Since $a$~is arbitrary
in~$I$, we infer that $I$~is a nonempty open subset of $\Q_{>0}$. Hence I3 holds.
\endproof

We finally give an \ev t \ch ization of an ideal~$\a$ of~$\Q_{>0}$ (``\po'' cut) given in~\cite[p.~43]{10a}.
A~subset~$\a$ of~$\Q_{>0}$ is an ideal of~$\Q_{>0}$ if it \sf ies the \fw\ \cn s:

\ssk
L1 \ $\a\ne\vn$ and $\a\ne\Q_{>0}$.

L2 \ \E\fe $a\in\a$ and \fe $b\in\a^c$, the \cpl\ of~$\a$ in $\Q_{>0}$, we have $a<b$.

L3 \ There is no $d\in\Q_{>0}$ \st $a\le d$ \fa $a\in\a$.


\blm4.11
A subset $\a$ of $\Q_{>0}$ is an ideal of $\Q_{>0}$ iff the subset $\a$ \sf ies \cn s {\rm L1--L3}.
\elm

\proof
Let $\a$ be a subset of $\Q_{>0}$. Clearly, L1 is identical to C1 and~I1. We assume L1.

``L2 and L3 imply C2 and C3'': We first prove C2. Let $a\in\a$, and $b\in\Q_{>0}$, $b<a$. Then either $b\in\a$ or $b\in\a^c$. If $b\in\a^c$, we
obtain a~\cd ion. Indeed, from L2 we infer $a<b$ since $a\in\a$. Hence $b\in \a$. Next we prove~C3. Suppose, for \cd ion, that \te s $a\in\a$ \st
every $b\in\Q_{>0}$ greater than~$a$ does not belong to~$\a$. If $b<a$ then $b\in\a$ by~C2 since $a\in\a$. By \as\ $a\in\a$. Hence $\a =
\{x\in\Q_{>0}: 0<x\le a\}=:(0,a]$. \If that $\a$ is the greatest \el\ of~$\a$, \cd ing~L3. Hence C2 and~C3 hold.

``C2 and C3 imply L2 and L3'': We prove L2. Let $a,b\in\Q_{>0}$ be \st $a\in\a$ and $b\notin\a$. Suppose for \cd ion that $a\not< b$. Then either
$a=b$ or $a>b$ since the \og~$\le$ of~$\Q_{>0}$ is total. The case $a=b$ is not possible since $a\in\a$ and $b\notin\a$. If $a>b$ then $b\in\a$
by~C2 since $a\in\a$, a~\cd ion. Hence $a<b$. We now prove~L3. Suppose, for \cd ion, that the subset~$\a$ has a greatest \el~$a$. In view of~C3,
\te s $c\in\a$ with $a<c$. A~\cd ion.
\endproof

\brm4.12 \

\hph i,i, Let $\a\in\dot P$ and $\a^c$ be its \cpl\ in~$\dot P$. Then in view of I1 and $\a=(\a^c)^c$ we have:
\ssk

D1 \ $\a^c\ne\vn$ and $\a^c\ne\Q_{>0}$.

\ssk
\noindent\Mo the \fw\ assertion holds:

\ssk
D2 \ \E\fe $a\in\a^c$ and every $b\in\Q_{>0}$, $b>a$ implies $b\in\a^c$.

\ssk \noindent
Indeed, if $a\in\a^c$ then $d<a$ \fa $d\in\a$ by L2. Thus if $b>a$ then $b>a>d$, hence $b>d$ \fa $d\in\a$. \If that $b\notin\a$, since otherwise
$b>b$, a~\cd ion. Note that D1 and~D2 are C1 and~C2 where $\a$~is replaced by~$\a^c$ and the \rl~$<$ is replaced by its reverse (or converse).
However, the ``converse'' of~C3 (see Exercise 1.3.4\,(i)) does not hold \fa $\a\in\dot P$. Indeed, if $x\in\a^c$ then $x$~is an \ub\ of~$\a$ by~L2.
\E\oh if $x\in\a$ then $x$~is \ti{not\/}
an \ub\ of~$\a$ in view of~C3. \E\Tf $\a^c$ is the set of \ub s of~$\a$ (see \E\df s and Notations 3.1.3) denoted by $\UB(\a)$, thus
\beq 4.14n
\a^c = \UB(\a) \qh{\fe}\a\in\dot P.
\e
\If that $\a$ has a supremum iff $\a^c$ has a least \el\ $a\in\Q_{>0}$ iff $\a^c=\{x\in\Q_{>0}: x\ge a\}$ \fs $\a\in\Q_{>0}$ iff $\a=\{x\in\Q_{>0}:
x<a\}$ \fs $a\in\Q_{>0}$ iff $\a=i(a)$ \fs $a\in\Q_{>0}$. In this case we shall call~$\a$ a \ti{\pn\ ideal\/} (of~$\Q_{>0}$) with \ti{\Gn}~$a$. \If
from Lemma \rf{l4.5} that not all ideals of $\Q_{>0}$ are \pn. Note that if $\a$ is not \pn, then $\a^c$ is a nonempty open set of~$\Q_{>0}$. Indeed,
$\a^c\ne\vn$ by~C1 and since $\a^c$~has no least \el, \fe \el\ $b\in\a^c$ \te s $a\in\a^c$ \st $a<b$. \E\Tf we have $b\in(a,b+1)\sbs \a^c$, hence
$\a^c$ is open, and $\Q_{>0}$ is the \ti{disjoint union of two nonempty open subsets}, namely $\a$ and~$\a^c$.

\hph ii,, Observe that \fe $\a\in\dot P$ the \fw\ assertion holds:
\beq4.15
\hbox{\fe} a\in\a \hbox{ we have }i(a)<\a.
\e
Indeed, if $a\in\a$ (such an $a$ exists by I1), and $b\in\Q_{>0}$ \sf ies $b<a$, then $b\in\a$ by~I2. Hence we obtain $i(a)\nde4.9 = (0,a)\sbs\a$,
and $i(a)\le\a$ in view of \er{4.10}. Since $a\notin(0,a)$ and $a\in\a$, we have $i(a)<\a$.
\erm

\If from \er{4.15} that $\bcl_{a\in\a}i(a)\sbs\a$. Indeed, \fe $x\in \bcl_{a\in\a}i(a)$, \te s $a\in\a$ \st $x\in i(a)$ by \er{4.14}. Hence
$\bcl_{a\in\a} i(a)\le\a$ by \er{4.10}. \E\oh
if $x\in\a$, then \te s $b\in\a$ \st $b>x$ by~C3. Hence $x\in i(b)$ by \er{4.9}. Since $i(b)\sbs \bcl_{a\in\a}i(a)$, we obtain
$x\in\bcl_{a\in\a}i(a)$. Hence $\a\sbs\bcl_{a\in\a}i(a)$. \E\csq, we have
\beq4.16
\a = \bigcup_{a\in\a}i(a) \qh{\fa} \a\in\dot P.
\e

\Mo we have
\beq4.17
\a = \Bca_{b\in\a^c}i(b). 
\e
If $\a=(0,a)$ \fs $a\in\Q_{>0}$, then
\beq4.18
a \notin i(a),
\e
and
\beq4.19
\Bca_{b>a}i(b) = i(a)\cup \{a\}\ (=(0,a]).
\e
Indeed, \fe $c\in\a$ and every $b \in \a^c$, we have $c<b$ by~L2, hence $c\in(0,b)=i(b)$. We obtain $\a\sbs i(b)$ \fe $b\in\a^c$.
\csq, $\a\sbs\Bca_{b\in\a^c}i(b)$.
If \te s $a\in\Q_{>0}$ \st
$\a=(0,a)$, then $a\notin\a$ since $a\not<a$, and \er{4.18} holds. \Mo in this case $\a^c = \{{x\in\Q_{>0}}:x\ge a\}$, hence $b\in \a^c$ iff $b\ge a$,
and $\Bca_{b\in\a^c}i(b) = \Bca_{b\ge a}i(b) = i(a) = \a$, since the map~$i$ is in\cre\ by \er{4.10}, which proves \er{4.17} in the case where $\a$~is
\pn. \Mo in this case we have $\Bca_{b>a} i(b) \nde4.9 = \{x\in\Q_{>0}: x<b \text{ \fa} b>a\} \supset \{{x\in\Q_{>0}:}\break x\le a\}$ since $x\le
a<b$ implies $x<b$ \fa $x\in\Q_{>0}$. However, if $x>a$ then $a<\frac12(a+x)<x$ by \er{4.23}, hence $x\notin\Bca_{b>a}i(b)$, which proves \er{4.19}.
We now suppose that $\a$~is not \pn, and for \cd ion,
that \te s $a\in\a^c$ \st $a\in \Bca_{b\in\a^c}i(b)$. In view of the discussion \fw\ \er{4.14n} \te s $d\in\a^c$, $d<a$. Since $d\in\a^c$ and
$a\in\Bca_{b\in\a^c}i(b)$, we obtain $a\in i(d)$, hence $a<d$ by \er{4.9}. A~\cd ion. \E\Tf \er{4.17} holds.
\ssk

We now consider the union of a family of ideals of $\Q_{>0}$.

\bpr4.13
Let $J$ be a \ns\ and let $\{\a_j\}_{j\in J}$ be a family of ideals of $\Q_{>0}$. Then $\bcl_{j\in J}\a_j$ is either $\Q_{>0}$ or an
ideal of~$\Q_{>0}$. \E\Ip
\beq4.20
\bigcup_{n\in\Na}(0,n) = \Q_{>0}.
\e
\epr

\proof
In view of \E\Pr\ 4.5.61 \fe $x\in\Q_{>0}$, \te\ $m\in\N$ and $p\in [0,1)\break\sbs\Q_{\ge0}$ \st $x=m+p$. Set $n:=m+1$. Then $x<n\in\Na$, hence
$x\in(0,n)$, and $\Q_{>0}\sbs\bcl_{n\in\Na}(0,n)$. The other inclusion is trivial, hence \er{4.20} holds. Suppose that \te s $r\in\Q_{>0}$ \st
$r\notin\bcl_{j\in J}\a_j$. \E\Ip $\bcl_{j\in J}\a_j\ne\Q_{>0}$. By I3 and Lemma \rf{l4.9}, $\bcl_{j\in J}\a_j$ is a nonempty open subset
of~$\Q_{>0}$. Then
\cn s I1~and~I3 are \sf ied. Let $a,b\in\Q_{>0}$ \sf y $a\in\bcl_{j\in J}\a_j$ and $b<a$. Then \te s $k\in J$ \st $a\in\a_k$. Then $b\in\a_k$ by~I2
since $\a_k$~is an ideal. \csq, $b\in\bcl_{j\in J}\a_j$. Thus \cn~I2 is \sf ied, and $\bcl_{j\in J}\a_j\in\dot P$.
\endproof

\begin{prp}[\cite{10a}] \lb{p4.14}
Let $\a\in\dot P$. Then \fe $e\in\Q_{>0}$ \te\ $a\in\a$ and $b\in\a^c$ \st $b-a=e$.
\epr

\proof
By I1 \te s $x\in\a$. Set
\beq4.21
N := \{n\in\Na: x+ne\in\a^c\}.
\e

``$N\ne\vn$: Let $y\in\a^c$. Then $y>x$, hence $y-x>0$. By Lemma 4.5.54\,(iii), \te s $k\in\Na$ \st $ke>y-x$. \If that $x+ke>(y-x)+x = (y+(-x))+x =
y+(x+(-x)) = y$. Since $x+ke>y$ and $y\in\a^c$, we have $x+ke\in\a^c$ by~D2. Hence $N\ne\vn$.
\ssk

\E\te s a smallest \el\ $m$ of~$N$ by Theorem 1.3.21. If $m=1$, then $a:=x$ and $b:=x+e$ \sf y $a\in\a$, $b\in\a^c$ and $b-a=e$. If $m>1$, we set
$a:=x+(m-1)e$ and $b:=x+me$. Then $a\in\a$ by the ``minimality'' of~$m$ and $b\in\a^c$. \Mo $b-a=e$.
\endproof

In the next lemma we give another \ch ization of an \el\ of~$i(x)$, $x\in\Q_{>0}$. We defined $i(x)$ by using the \og~$\ge$ of the \ra\ \nm s induced
by the \ad. It is useful to give a \ch ization related to the \mlc\ on~$\Q_{>0}$. Recall that in the proof of~C1 for $i(x)$ we used the fact that
$\frac12x\in(0,x)$ as well as $\frac12(a+x)\in(a,x)$ \fe $a\in(0,x)$.

\blm4.15
Let $x\in\Q_{>0}$ and $i(x)=(0,x)$. Then the \fw\ assertions hold\/\dw
\bea4.22
{}&\hbox{Let $y\in\Q_{>0}$. Then $y\in i(x)$ iff $y=tx$ \fs $t\in(0,1)$.}\\
&\hbox{Let $a,b\in\Q_{>0}$ \st $a<b\le x$. Then $(1-t)a+tb\in (a,x)$ \fa $t\in(0,1)$.}\lb{4.23} \\
&\hbox{If $s,t\in(0,1)$ then $st\in(0,1)$.}\lb{4.24}
\e
\elm

\proof
We recall that $\pz2{\Q_{>0}}$ is a group, $\pz1{\Q_{\ge0}}$ is a \PM, and $0<y<x$ iff \te s $z\in\Q_{>0}$ \st $x=y+z$, \fa $x,y\in\Q_{>0}$.

``\er{4.22} \ti{Only if\/}'': Let $x,y\in\Q_{>0}$ \sf y $y<x$. Let $z\in\Q_{>0}$ be \st $x=y+z$. Let $x\Inv\in\Q_{>0}$ \sf y $x\Inv\cdot x=1$,
then $1=x\Inv\cdot x = x\Inv(y+z) \nad{\rm(4.1.2)} = x\Inv y+x\Inv z > x\Inv y$. Hence $x\Inv y<1$. Set $t:=x\Inv y$, then $0<t<1$ and $tx=
xt =x(x\Inv y)=(xx\Inv)y= 1\cdot y=y$.

``\ti{If\/}'': Let $t\in(0,1)$ and $0<x$. Then $0 = 0t\nad{(4.5.121)}< xt = tx < x$.

``\er{4.23}'': Let $a,b\in\Q_{>0}$ \sf y $a<b\le x$. Then \te\ $c>0$ and $d\ge0$ \st $b=a+c$ and $x=b+d$. Let $t\in(0,1)$. Then
$(1-t)a+tb = (1-t)a
+ t(a+c) = (1-t)a+ ta+tc = ((1-t)+t)a + tc = 1a+tc = a+tc > a>0$. Thus $(1-t)a+tb > a$. \Mo $(1-t)a < (1-t)x$ since $a<x$ and $1-t>0$. Similarly,
$tb\le tx$. Hence $(1-t)a+tb \nad{(4.5.120)}< (1-t)x+tb = tb+(1-t)x \le tx+(1-t)x = x$ since $tb\le tx$. Thus $(1-t)a+tb < x$ and \er{4.23} holds.

``\er{4.24}'': Follows from $0<st<1t= t1 <1\cdot1 =1$.
\endproof

\Wanp introduce \tb{\art\ \op s} on the set of \pn\ ideals. Since the map $i:\Q_{>0} \to\dot P$ defined by \er{4.9} is in\jc, the map $f:\Q_{>0}\to
i(\Q_{>0})$, the set of \tb{\pn\ ideals}, defined by $f(x)=i(x)$, $x\in\Q_{>0}$, is bi\jc. Proceeding as in Example 2.1.5\,(v) where $(M,\qu,e):=
\pz2{\Q_{>0}}$, we define a \mlc~$\cdot'$ on $i(\Q_{>0})$ by setting
\beq4.25
\a\cdot'\b := f(f\Inv(\a)\cdot f\Inv(\b)) \qh{\fa $\a,\b\in i(\Q_{>0})$.}
\e
Then $(i(\Q_{>0}),\cdot',i(1))$ becomes a \Cm, since $\pz2{\Q_{>0}}$ is an \ag, hence a \Cm. \Mo $f$~is a monoid-\is sm from $\pz2{\Q_{>0}}$ onto
$(i(\Q_{>0}),\cdot',i(1))$. In view of Lemma 4.3.25, $(i(\Q_{>0}),\cdot',i(1))$ is an \ag\ and (4.3.17) holds.

Since $(\Q,+,\cdot,0,1,\ge)$ is an \of, its \rt ion to $\Q_{\ge0}\nad{(4.5.9)}= \Q_{>0}\cup\{0\}$ is an ordered semifield by Theorem \rf{t3.8}\,(i).
The union on the \RHS\ of (4.5.9) being disjoint we may extend the bi\jn\ $i:\Q_{>0}\to i(\Q_{>0})$ to $\Q_{\ge0}$ by setting $i(0):=\vn$, the empty
subset of $\Q_{\ge0}$. Note that $i(0)\ne i(x)$ \fa $x\in\Q_{>0}$. Thus we may define a map $f:\Q_{\ge0} \to i(\Q_{>0})\cup\{\vn\}$ by setting:
\beq4.26
f(x):= \bca
i(x) & \hbox{for }x\in\Q_{>0},\\
\vn & \hbox{for }0\in\Q_{\ge0}.
\eca
\e
Observe that the map $f$ is a \ti{bi\jn}.

\Wanp apply the \fw\ \Pr.

\bpr4.16
Let $\pz0F$ be an ordered semifield. Suppose \te\ a set~$F'$ and a bi\jn\ $g:F\to F'$. Define the \el s of~$F'$, binary \op s and \rl s on~$F'$
by setting\dw
\bea4.27
0'&:=g(0) \hbox{ and } 1':=g(1),
\intertext{and \fa $x',y'\in F'$}
x'+'y'&:= g(g\Inv(x')+g\Inv(y')), \lb{4.28}\\
x'\cdot'y'&:= g(g\Inv(x')\cdot g\Inv(y')), \lb{4.29}\\
x'\ge'y'&\hbox{ if \ }g\Inv(x')\ge g\Inv(y'). \lb{4.30}
\e
Then the \fw\ assertions hold\/\dw
\bea4.31
{}&(F',+',\cdot',0',1',\ge') \hbox{ is an ordered semifield,}\\
&\hbox{the map $g:F\to F'$ is an order- and $($semi-$)$ring \is sm.}\lb{4.32}
\e
\epr

\proof
``$(F',+',\cdot',0',1')$ \ti{is a \sr\ with unity}'': $(F',+',0')$ and $(F',\cdot',1')$ are \am s by Example 2.1.5\,(v), \E\df s 4.1.1 and \rf{d3.7}.
In view of the \cmt ity of~$\cdot'$, it is sufficient to prove that $0'\cdot'x'=0'$, $x'\in F'$, and (4.1.2). We have
$0'\cdot'x' \nde4.29 = g(g\Inv(0')\cdot g\Inv(x')) \nde4.27 = g(0\cdot g\Inv(x')) \nad{(4.1.1)}= g(0)\nde4.27 = 0'$. Set $x:=g\Inv(x')$, $y:=g\Inv
(y')$ and $z:=g\Inv(z')$. Then $x'\cdot'(y'+'z') \nad{(2.1.14)}= g(x)\cdot' g(y+z) \nde4.29 = g(g\Inv g(x)\cdot'g\Inv (g(y+z))) =\break
g(x\cdot(y+z)) \nad{(4.1.2)}= g(xy+xz) \nad{(2.1.14)}= g(xy)+'g(xz) \nad{(2.1.14)}= (g(x)\cdot' g(y)) +' (g(x)\cdot' g(z)) = x'\cdot'y' +' x'\cdot'
z'$.

Next we prove (i)--(iv) of \E\df\ \rf{d3.7}.

(i) follows from Example 2.1.5\,(v) and \E\df\ \rf{d3.7}\,(i).

(ii) Let $x',y'\in F'$ \sf y $x'\ne y'$, and let $x:=g\Inv(x')$, $y:=g\Inv(y')$. Then by \E\df\ \rf{d3.7}\,(ii) for~$F$, we have either $g\Inv(x') >
g\Inv(y')$ or $g\Inv(y') > g\Inv(x')$, since $g\Inv(x')\ne g\Inv(y')$. Hence either $x'>'y'$ or $y'>'x'$ by~\er{4.30}.

It remains to show that the \og~$\ge'$ \sf ies (3.1.4), that is, \fa $x',y'\in F'$ \st $x'\ge'y'$ we have $x'=y'+z'$ \fs $z'\in F'$. Let $x',y'\in
F'$ \sf y $x'\ge'y'$ and set $x:=g\Inv(x')$, $y:=g\Inv(y')$. Then $x\ge y$ by \er{4.30}. \E\te s $z\in\Q_{\ge0}$ \st $x=y+z$ by (4.5.118), (4.5.119).
We obtain $x'=g(x)=g(y+z)\nde4.28 = g(y)+g(z) = y'+z'$ where $z':=g(z)\in F'$. Hence the \og s $\ge'$ and $\stackrel{+'}{\ge'}$ are identical.

(iii) By \E\df\ 4.1.1, $1\ne0$, hence $1'\ne0'$ since $g$ is in\jc. Therefore\break $1'\in F'\sms{0'}$. \Mo let $x',y'\in F'\sms{0'}$. Let
$x:=g\Inv(x')$, $y:=g\Inv(y')$. Then $g\Inv(x'\cdot' y') \nad{(2.1.11)} = g\Inv g(xy)=xy$. Since $x'\ne0'$, $y'\ne0'$, $g$~is in\jc\ and $g(0)=0'$,
we have $x\ne0$, $y\ne0$. Hence $xy\ne0$ from \E\df\ \rf{d3.7}\,(iii) for~$F$. Hence $x'\cdot' y' = g(x)\cdot' g(y) = g(xy)\ne0'$ as above. Hence
$F'\sms{0'}$ is a \sbm\ of the monoid $(F',\cdot',1')$.

(iv) Since $(F',\cdot',1')$ is an \am, so is $(F'\sms{0'},\cdot',1')$. It suffices to prove that \fe $x'\in F'$ \te s $y'\in F'$ \st $x'\cdot'y'=1'$,
in view of \E\df\ 4.3.2. Let $x:=g\Inv(x')$. Then by~(iv) for~$F$ \te s $y\in F$ \st $x\cdot y=1$. Hence $g(x)\cdot' g(y)\nad{(2.1.14)}= g(x\cdot y)
=g(1) \nde4.27 = 1'$. Set $y':=g(y)$.

``\er{4.32}'': We have $g(0)=0'$, $g(1)=1'$ by \er{4.27}, $g(x+y) = g(x)+'g(y)$ by \er{4.28} and $g(x\cdot y) = g(x)\cdot'g(y)$ by \er{4.29}. \Mo
$x\ge y$, $x,y\in F$, implies $g(x)\ge g(y)$ by \er{4.30}. Since $\ge$~is total and $g$~is bi\jc, $g\Inv$~is in\cre\ by Lemma 1.3.24\,(i).
\endproof

We now apply \E\Pr\ \rf{p4.16} with $\pz0F := \pz0{\Q_{\ge0}}$, $F':= i(\Q_{>0})\break \cup \{\vn\}$, and $g:=f$ defined in \er{4.26}.
Inspection of the proof shows that for~$g$ to be a \hbox{(semi-)}\hskip0.4pt ring \is sm (see \E\df\ 4.1.7), it is sufficient for~$g$ to be bi\jc.
Since $f|_{\Q_{>0}} : \Q_{>0} \to
i(\Q_{>0})$ is bi\jc\ we could choose for $f(0)$ any subset of $\Q_{>0}$ not belonging to $i(\Q_{>0})$. However, if we want $f$ to be an \ois sm from
$(\Q_{\ge0},\ge)$ onto $(i(A)\cup f(0),\supset)$ we need $f(0)=\vn$. Indeed, since $0<\frac1n \in \Q_{>0}$ \fa $n\in\Na$, we must have $f(0)\sbs
i\bigl(\frac1n\bigr)$ \fa $n\in\Na$. Suppose, for \cd ion, that \te s $0<a\in f(0)$. Then $(0,a)\sbs\bigl(0,\frac1n\bigr)$, hence $a<\frac1n$ by Lemma
\rf{l4.10}. Let $a=\frac pq$, $p,q\in\Na$, then $np<q$, hence $n=n1\le np<q$ \fa $n\in\Na$. But $n:=q+1\not< q$. A~\cd ion. Since $f(0)=\vn\subsetneq
i(x) = f(x)$ \fa $x\in \Q_{>0}$ and $f(y)\subsetneq f(x)$ \fa $0<y<x$, and $f:\Q_{\ge0}\to \{\vn\}\cup i(\Q_{>0})$ is a bi\jn, we find that $f$~is
an \ois sm by Lemma 1.3.34\,(i) and \E\df\ 1.3.32.

Summarizing, we obtain:

\blm4.17
Let $L:= i(\Q_{>0})\cup\{\vn\}$. Let $\ge$ denote the \rt ion to~$L$ of the \og\ on $\cP(\Q_{\ge0})$
$($the power set of $\Q_{\ge0})$ defined by \er{4.13} \fa $\a,\b\in\cP(\Q_{\ge0})$. Let $f:\Q_{\ge0}\to L$ be defined in \er{4.26}, let
$+'$~$($resp.~$\cdot')$ be the binary \op s on~$L$ defined in \er{4.28} $($resp.~\er{4.29}$)$ where $g:=f$, and let $i(1)$ be defined in \er{4.9}.
Then the \fw\ assertions hold\/\dw

\hph i,i, $(L,+',\cdot',\vn,i(1),\ge)$ is an ordered semifield,

\hph ii,, $f:\Q_{\ge0}\to L$ is a $($semi$)$ring- and \ois sm.
\elm

Our next goal is to give a \chz\ of $+'$~and~$\cdot'$ which ``naturally'' extends to~$\dot P$.

We first recall
\bea4.33
\vn+'\a &= \a+'\vn = \a \qh{\fa $\a\in L$,}\\
\vn \cdot'\a &= \a\cdot'\vn= \vn \qh{\fa $\a\in L$.}\lb{4.34}
\e
Indeed, given $\a\in L$ let $a\in\Q_{\ge0}$ be \st $a = f\Inv(\a)$. Note that $\vn = f\Inv(0)$. We have $\vn+'\a \nde4.28 = f(0+a) = f(a) = \a =
f(a+0) = \a+'\vn$. The proof of \er{4.34} is similar.

\blm4.33
Let $\a,\b \in i(Q_{>0})$. Then
\bea4.35
\a+'\b &= \bigcup_{x\in\a,y\in \b} \{x+y\},\\
\a\cdot'\b &= \bigcup_{x\in\a,y\in \b} \{x\cdot y\}.\lb{4.36}
\e
\elm

\proof
Let $a,b\in \Q_{>0}$ be \st $\a=i(a)=f(a)$, $\b=i(b)=f(b)$. Then by \er{4.28}, \er{4.29} we have
\beq4.37a
i(a)+'i(b) = i(a+b); \q i(a)\cdot' i(b)=i(a\cdot b).
\e

``\er{4.35}, \er{4.36} $\supset$'': Let $z\in\{x+y: x\in\a,y\in\b\}$. Then $z=x+y$ \fs $x\in(0,a)$, $y\in(0,b)$. Thus $x>0$, $y>0$ and $x<a$, $y<b$
by \er{4.9}. Hence $x+y<a+y=y+a < b+a=a+b$ by (4.5.120). Thus $x+y<a+b$, and $x+y\in i(a+b)$, which proves ``\er{4.35}~$\supset$''. In the case of
``\er{4.36}'' we use (4.5.121) instead of (4.5.120).

``\er{4.35} $\sbs$'': Let $z\in i(a+b)$. By \er{4.22} \te s $t\in(0,1)$ \st $z=t(a+b)=ta+tb$. By~\er{4.22} $ta\in i(a)=\a$, $tb\in i(b)=\b$, hence
$z\in\{x+y: x\in\a,y\in\b\}$.

``\er{4.36} $\sbs$'': Let $z\in i(a\cdot b)$. By \er{4.22} \te s $s\in(0,1)$ \st $z=sa\cdot b$. Set $t':=\frac12(s+1) = \frac12s+\frac12$. Then
$0<s<t'<1$ by \er{4.23} where $a:=s$, $b:=1$, $x:=1$ and $t:=\frac12$. Hence $0<sab < t'ab <ab$ by (4.5.121). Let $t''$~denote the inverse of~$t'$
in the \ag\ $\pz2{\Q_{>0}}$, then we obtain $0<s\cdot t''<t'\cdot t''=1$ by (4.5.121). \E\Tf $z=sab = s\cdot1\cdot ab = s\cdot(t''\cdot t')\cdot ab
\nad{(2.1.36)}= (s\cdot t'')\cdot(t'ab) = (s\cdot t'')\cdot(a\cdot t'\cdot b)= ((st'')a)\cdot(t'b)$. Since $s\cdot t''<1$ and $t'<1$ we have
$(st'')a \in i(a)$ and $t'b\in i(b)$ by \er{4.22}. Thus setting $x:=(st'')a$ and $b:=t'b$ we obtain $z=xy$ where $x\in\a$ and $y\in\b$, hence
$z=\bcl_{a\in\a,b\in\b}\{ab\}$.
\endproof

We now extend to $\dot P=\id(\Q_{>0})$ the binary \op s $+'$ and $\cdot'$ by using the fact that the \RHS s of \er{4.35} and \er{4.36} make sense
when $\a$ and~$\b$ are not necessarily \pn\ ideals. However, it is not obvious that the \RHS s of \er{4.35} and \er{4.36} belong to~$\dot P$. In the
next lemma, we give \ev t formulae for these RHS's, which allow us to show that indeed these RHS's belong to~$\dot P$.

\blm4.34
Let $\a,\b\in\dot P$. Then we have
\bea4.37
\bigcup_{a\in\a,b\in\b}\{a+b\} &= \bigcup_{c\in\a,d\in\b}i(c+d),\\
\bigcup_{a\in\a,b\in\b}\{a\cdot b\} &= \bigcup_{c\in\a,d\in\b}i(c\cdot d),\lb{4.38}
\e
where $+$ $($resp.~$\cdot)$ denotes the \ad\ $($resp.\ \mlc$)$ on $\Q_{>0}$.
\elm

\proof \

``$\sbs$ \ti{for \er{4.37} and \er{4.38}}'': Let $z\in \bcl_{a\in\a,b\in\b}\{a\qu b\}$ where ${\qu}:=+$ in the case of \er{4.37} and ${\qu}:=\cdot$
in the case of \er{4.38}. \E\te\ $a\in\a$, $b\in\b$ \st $z=a\qu b$. In view of~C3 \te\ $a',b'\in\Q_{>0}$ \st $a<a'\in\a$ and $b<b'\in\b$. From
(4.5.120), (4.5.121) and the \cmt ity of~$\qu$, we infer $a\qu b\le a'\qu b'$. Hence $a\qu b\in i(a'\qu b')$. Since $a'\in\a$, $b'\in \b$, we obtain
$z\in \bcl_{c\in\a,d\in\b}i(c\qu d)$.

``$\supset$'': Let $z\in\bcl_{c\in\a,d\in\b}i(c\qu d)$. Then \te\ $c'\in\a$, $d'\in\b$ \st $z< c'\qu d'$. In the case where ${\qu}:=+$, we have
$c'+d'>0$ and $z(c'+d')\mo < (c'+d')(c'+d')\mo = 1$ by (3.5.121). Set $t:= z(c'+d')\mo$. Then $z = t(c'+d') = tc'+td'$. Since $t\in(0,1)$, we have
$tc' < c'\in\a$ and $td' < d'\in\b$ by~C2, hence $tc'+td'\in \bcl_{c\in\a,d\in\b}\{c+d\}$.

In the case where ${\qu}:=\cdot$, we have $z<c'\cdot d'$, and $z=c'\cdot(z\cdot(c')\mo)$ with $c'\in\a$. Note that $z\cdot(c')\mo < (c'\cdot d')
(c')\mo$ by (3.5.121). Hence $z\cdot(c')\mo < d'\in\b$ by~C2. Setting $c'':=c'$ and $d'':=z\cdot(c')\mo$, we have $z = c''\cdot d''$ with $c''\in\a$
and $d''\in\b$, hence $z\in\bcl_{c\in\a,d\in\b}\{c\cdot d\}$.
\endproof

\blm4.35
Let $\a,\b\in\dot P$. Then both $\bcl_{a\in\a,b\in\b}\{a+b\}$ and $\bcl_{a\in\a,b\in\b}\{a\cdot b\}$ belong to~$\dot P$.
\elm

\proof
In view of \er{4.37}, \er{4.38} and \E\Pr\ \rf{p4.13}, it is \sft\ to show that\break $\bcl_{a\in\a,b\in\b}\{a\qu b\}$, ${\qu}:=+$ and ${\qu}:=\cdot$,
are not equal to $\Q_{>0}$. By~C1, \te\ $c_\a$ (resp.~$c_\b$) in~$\a^c$ and by~L2 $a\le c_\a$, $b\in c_\b$ \fa $a\in\a$, $b\in\b$. \E\Tf $a\qu b \le
c_\a\qu c_\b$ \fa $a\in\a$, $b\in\b$. Hence $i(a\qu b)\nde4.10 \sbs i(c_\a\qu c_\b)$ \fa $a\in\a$, $b\in\b$. \E\Tf $\bcl_{a\in\a,b\in\b}i({a\qu b})
\break \sbs i(c_\a\qu c_\b)$. Since $(c_\a\qu c_\b) < (c_\a\qu c_\b)+1$, we have $i(c_\a\qu c_\b) \subsetneq i((c_\a\qu c_\b)+1)$ by \er{4.10}, and
$\bcl_{a\in\a,b\in\b} i(a\qu b)\ne\Q_{>0}$.
\endproof

As a con\sq\ of Lemma \rf{l4.35} we may define two binary \op s $\hpl$ and $\hcd$ on~$\dot P$ by setting
\bea4.39
\a\hpl \b &:= \bcl_{a\in\a,b\in\b}\{a+b\}, \q \a,\b\in\dot P,\\
\a\hcd \b &:= \bcl_{a\in\a,b\in\b}\{a\cdot b\}, \q \a,\b\in\dot P,\lb{4.40}
\e
where $+$ and $\cdot$ denote the \ad\ and \mlc\ on $\Q_{>0}$ defined in (4.5.12) and (4.5.3). Since both \op s $+$ and~$\cdot$ are \cmt e, the \fw\
holds:
\beq4.41
\a\hpl\b = \b\hpl\a, \q \a\hcd\b = \b\hcd\a \qh{\fa $\a,\b\in\dot P$}.
\e

We now show that both \op s $\hpl$ and $\hcd$ are \asc e. We have for $\a,\b,\g \in\dot P$, ${\qu}:=+$ or~$\cdot$, ${\hqu}:={\hpl}$ or~${\hcd}$\,:
$(\a\hqu\b)\hqu \g = \bcl_{z\in(\a\hqu \b),\,c\in\g}\{z\h{\qu}\h c\} = \bcl_{c\in\g}\bcl_{z\in(\a\hqu\b)}\{z\h{\qu}\h c\} = \bcl_{c\in\g}
\bcl_{a\in\a,b\in\b}\{(a\h{\qu}\h b)\h{\qu}\h c\} = \bcl_{c\in\g} \bcl_{a\in\a,b\in\b}\{a\h{\qu}\h (b\h{\qu}\h c)\} =
\bcl_{a\in\a}\bcl_{b\in\b}\bcl_{c\in\g}\{a\h{\qu}\h(b\h{\qu}\h c)\} =
\bcl_{a\in\a}\bcl_{z\in(\b\hqu\g)}\{a\h{\qu} \h z\} = \a\hqu(\b\hqu\g)$. We have proved that
\bea4.42
(\a\hpl\b)\hpl\g &:= \a\hpl(\b\hpl\g) \qh{\fa} \a,\b,\g \in \dot P,\\
(\a\hcd\b)\hcd\g &:= \a\hcd(\b\hcd\g) \qh{\fa} \a,\b,\g \in \dot P.\lb{4.43}
\e

We next define $\wh1 \in i(\Q_{>0}) \sbs \id(\Q_{>0}) = \dot P$, by setting
\beq4.43a
\wh1 := i(1).
\e
Then \fe $\a\in\dot P$ we have $\a\hcd\wh1 = \a\hcd i(1) \nde4.40 = \bcl_{a\in\a,\,b\in i(1)}\{a\cdot b\} \nde4.9 = \bcl_{a\in\a}\bcl_{t\in(0,1)}
\{t\cdot a\}\break \nad{(2.2.12),(4.5.3),\er{4.22}} = \bcl_{a\in\a} i(a) \nde4.17 = \a$. Thus we have
\beq4.44
\a\hcd \wh1 = \a \qh{\fa $\a\in\dot P$.}
\e
\If from \er{4.40}, \er{4.41}, \er{4.44} that $(\dot P,\hcd,\wh1)$ is an \am. We next show that $(\dot P,\hcd,\wh1)$ is an \ag. In view of \E\df\
4.3.2 it suffices to prove \fa $\a\in\dot P$ the \ex\ of an \el\ $\b\in\dot P$ \sf ying
\beq4.45
\a\hcd \b=\wh1.
\e
If $\a:=i(a)$ \fs $a\in\Q_{>0}$, and $b\in\Q_{>0}$ \sf ies $a\cdot b=1$, then $i(a\cdot b)\nde4.26 = f(a\cdot b) = f(a)\cdot' f(b)$ by Lemma
\rf{l4.17}\,(ii). \Mo $f(a)\cdot'f(b) \nde4.36 = \bcl_{x\in i(a),y\in i(b)}\{xy\} \nde4.41 = i(a)\hcd i(b)$. Thus if $\b:=i(b)$, then $\a\hcd \b =
i(a\cdot b)=i(1)\nde4.43a =\wh1$. We now suppose that $\a\in\id(\Q_{>0})\sm i(\Q_{>0})$, and we set
\beq4.46
\b:=\{x\mo \in\Q_{>0}:x\in\a^c\},
\e
where $x\mo$ is the inverse of  $x$ in the group $\pz2{\Q_{>0}}$. Thus $x\mo\cdot x = x\cdot x\mo=1$ \fa $x\in\Q_{>0}$. We claim that $\b\in\dot P$
and that $\a\hcd\b=\wh1$. We first verify that $\b$ \sf ies C1--C3.

``C1'': Since $\a\in\id(\Q_{>0})$, \te\ $a\in\a$ and $b\in\a^c$ by~C1 for~$\a$ and $a<b$ by~L2. Then $\frac1b\in\b$. Note that $\b^c:=\{x\mo \in
\Q_{>0}: x\in(\a^c)^c=\a\}$. Thus $\frac1a\in\b^c$, and C1 holds for~$\b$.

``C2'': Let $c\in\b$ and $d\in\Q_{>0}$, $d<c$. Then $d(d\mo c\mo)<c(d\mo c\mo)$ by (4.5.121), hence $c\mo<d\mo$. Since $c\in\b$, we have $c\mo \in
\a^c$. Since $c\mo< d\mo$, we have $d\mo\in\a^c$ by~D2 in Remark \rf{r4.12}. Hence $d\in\b$ by \er{4.46}.

``C3'': Since $\a\in\id(\Q_{>0})\sm i(\Q_{>0})$, $\a^c$ has no least \el. Otherwise, \te s ${c\in\a^c}$ \st $c\le y$ \fa $y\in\a^c$, hence $\a=(0,c)
= i(c)$, contrary to the \as s. \If that $\b$~has no greatest \el. Indeed, let $\INV$ denote the self-map of~$\Q_{>0}$ defined by $\INV(x):=x\mo$.
By \E\Pr\ 4.3.26, $\INV$ is an involution, that is, $\INV\circ\INV = {\rm I}_{\Q_{>0}}$, the identity map in~$\Q_{>0}$. \E\Ip $\INV$ is bi\jc\ and
$(\INV)\Inv=\INV$. We already showed in part ``C2'' of the proof that $\INV(y)>\INV(x)$ iff $x>y$ \fa ${x,y\in\Q_{>0}}$ \st $x>y$. Thus $\INV(y) >
\INV(x)$ iff $x>y$ \fa $x,y\in\Q_{>0}$. Note that ${\b\in\INV(\a^c)}$.
We already observed that $\a^c$ has no least \el\ whenever $\a$~is not \pn. \csq, $\b$~has no greatest \el, which proves C3 for~$\b$. Thus $\b\in
\dot P$.

We show that $\a\hcd\b=\wh1$.

``$\sbs$'': $\a\hcd\b \nde4.40 = \bcl_{a\in\a,b\in\b}\{a\cdot b\} \nde4.38 = \bcl_{c\in\a,d\in\b} i(c\cdot b) = \bcl_{c\in\a,d\in\INV(\a^c)}
i(c\cdot d) = \bcl_{c\in\a,b\in\a^c} i(c\cdot b\mo)$. If $c\in\a$ and $b\in\a^c$, then $c<b$ in view of L2. Hence $c\cdot b\mo \nad{(4.5.121)}<
b\cdot b\mo=1$ \fa $c\in\a$
and $b\in\a^c$. \E\Tf \fa $c\in\a$ and $b\in\a^c$ we have $c\cdot b\mo<1$, hence $i(c\cdot b\mo)\nde4.10 \sbs i(1)$, and $\a\hcd \b = \bcl
_{c\in\a,b\in\a^c} i(c\cdot b\mo) \sbs i(1) \nde4.43a = \wh1$, thus we have proved
\beq4.47
\a\hcd \b \sbs i(1).
\e

``$\supset$'': We follow \cite[pp.~59--60]{10a}. Let $c\in i(1)$, and let $a\in\a$. By \E\Pr\ \rf{p4.14}, \te\ $a'\in\a$ and $b\in\a^c$ \st
\beq4.48
b-a' = (1-c)a.
\e
Note that $a'<b$ by L2 and $1-c>0$ since $c\in(0,1)$. We assume and prove later that $(a'\cdot c\mo)\mo \in\b$. Then $c = a'\cdot(a'\cdot c\mo)\mo \in
\bcl_{x\in\a,y\in\b} \{x\cdot y\} \nde4.40 = \a\hcd\b$. Since $c$~is arbitrary in~$i(1)$, we obtain $i(1)\sbs \bcl_{x\in\a,y\in\a^c} \{x\cdot
y\} = \a\hcd \b$.

Thus it remains to prove:

``$(a'\cdot c\mo)\mo\in\b$'': Since $a\in\a$ and $b\in\a^c$, we have $a<b$ by~L2, hence from \er{4.48} we obtain
$$
b-a' < (1-c)b.
$$
Applying (4.5.120) we obtain
$$
(b-a') + cb < (1-c)b+ cb.
$$
Since $(1-c)b+cb = ((1-c)+c)b = b = (b-a')+a'$, we obtain $(b-a')+cb < (b-a')+a'$. From \er{3.9} we infer $cb<a'$, hence $b=1b = c\mo cb < c\mo a'$.
Since $b<c\mo a'$ and $b\in\a^c$, we find $c\mo a'\in\a^c$ by~L2. Hence $(a'c\mo)\mo \in\b$ by \er{4.46}. This completes the proof that $\a\hcd\b
= i(1)=\wh1$.

We have proved the \fw

\blm4.36
The set $\dot P$ defined in \er{4.12} and endowed with the \mlc~$\hcd$ defined in \er{4.40} is an abelian \emph{group} with \nel~$\wh1$ defined
in \er{4.43a}. We will denote it by $(\dot P,\hcd,\wh1)$.
\elm

Note that the set $\dot P$ endowed with the \ad\ defined in \er{4.42} is an abelian semigroup (see \E\df\ 1.2.2). By the adjunction of the \el
\beq4.49
\wh0 := \vn, \qh{the empty subset of $\Q_{>0}$,}
\e
$\dot P$ becomes an \am\ provided the \fw\ holds:
\beq4.50
\wh0\hpl \wh0 = \wh0; \q \wh0\hpl\a = \a\hpl\wh0=\a \qh{\fa $\a\in\dot P$.}
\e
If we also extend the \mlc\ $\hcd$ to $P$ by setting
\beq4.51
\wh0 \hcd \wh0 = \wh0; \q \wh0\hcd\a = \a\hcd \wh0 = \wh0 \qh{\fa $\a\in\dot P$,}
\e
we obtain

\blm4.37
The set $P:=\dot P\cup \wh0$ endowed with the \ad~$\hpl$ defined in \er{4.39}, \er{4.50}, the \mlc~$\hcd$ defined in \er{4.40}, \er{4.51}, denoted
by $(P,\hpl,\hcd,\wh0,\wh1)$ where $\wh1$~is defined in \er{4.43a}, is a $($\cmt e$)$ \emph{semiring} $($with unity$)$. \Mo the \pp ies {\rm(iii)}
and {\rm(iv)} of \E\df\ \rf{d3.7} hold.
\elm

\proof
We use the \fw\ simple lemma whose proof is left to the reader.

\blm4.38
Let $(S,\qu)$ be an abelian semigroup $($see \E\df\ {\rm 1.2.2)} and let $\wh e\notin S$. Setting $\wh S:=S\cup \{\wh e\}$,
$$
x\hqu y:=\bca
x\qu y &\hbox{if }x,y\in S,\\
\wh e\hqu x = x\hqu\wh e = e & \hbox{if }y=\wh e,\ x\in S,\\
\wh e\hqu\wh e = \wh e &\hbox{if }x=\wh e,\ y=\wh e,
\eca
$$
we have $(\wh S,\hqu,\wh e)$ is an \am. If, moreover, $(S,\qu)$ is \cnc e, then $(\wh S,\hqu,\wh e)$ is a \Cm.
\elm

``$(P,\hpl,\wh0)$ \ti{is an \am\/}'': follows from \er{4.42}, \er{4.41}, \er{4.50} and Lemma \rf{l4.38}.

``$(P,\hcd,\wh1)$ \ti{is an \am\/}'':

``\ti{\asc ity of\/ $\hcd$}'': follows from \er{4.43} if $\a,\b,\g\in\dot P$. If $\a:=\wh0$, then $\wh0\hcd(\b\hcd\g)\nde4.51 = \wh0$ and $(\wh0\hcd
\b)\hcd\g = \wh0\hcd\g = \wh0$. If $\b$ or~$\g$ is equal to~$\wh0$, then $\a\hcd(\b\hcd\g) = \a\hcd\wh0 = \wh0$. If $\b:=\wh0$, $\a\hcd(\wh0\hcd\g)
= \a\hcd\wh0 = \wh0$, and if $\g:=\wh0$ then $(\a\hcd\b)\hcd\wh0 \nde4.51 = \wh0$.

``\ti{\cmt ity of $\hcd$}'': follows from \er{4.41} if $\a,\b\in\dot P$, and from \er{4.51} otherwise.

``\ti{\nel\/}'': follows from \er{4.44}.

``(4.1.1)'': follows from \er{4.51}.

``(4.1.2)'': If $\a:=\wh0$ then $\wh0\hcd(\b\hpl\g) \nde4.51 = \wh0 \nde4.50 = \wh0 \hpl\wh0 \nde4.51 = \wh0\hcd\b \hpl\wh0\hcd\g$. If $\b:=\wh0$ then
$\a\hcd(\wh0\hpl\g) \nde4.50 = \a\hcd\g \nde4.50 = \wh0\hpl \a\hcd\g \nde4.51 = \a\hcd\b \hpl \a\hcd\g$. The case $\g:=\wh0$ is similar.

We now suppose $\a,\b,\g\in\dot P$.

``$\a:=i(a),\ \b:=i(b),\ \g:=i(c)$'': In view of \er{4.26}, Lemma \rf{l4.17}\,(ii) we find\break $i(a)\cdot'(i(b)+' i(c)) = (i(a)\cdot' i(b))+'
(i(a)\cdot'i(c))$. In view of \er{4.35}, \er{4.36} and \er{4.39}, \er{4.40} we may replace $+'$ (resp.~$\cdot'$) by $\hpl$ (resp.~$\hcd$).

``$\sbs$'': Let $z\in \a\hcd(\b\hpl\g)$. By \er{4.40} \te\ $a\in\a$, $d\in\b\hpl\g$ \st $z=ad$. By \er{4.39} \te\ $b\in\b$ and $c\in\g$ \st $d=b+c$.
Hence $z=a(b+c) = ab+ac$. Since $a\in\a$, $b\in\b$, we have $ab\in\a\hcd\b$ by \er{4.40}, similarly $ac\in\a\hcd\g$. From \er{4.39} we infer
$ab+ac \in \a\hcd\b \hpl \a\hcd\g$.

``$\supset$'': Let $z\in\a\hcd\b \hpl \a\hcd\g$. Then \te\ $d\in\a\hcd\b$, $e\in\a\hcd\g$ \st $z= d+e$. \E\te\ $a\in\a$, $b\in\b$ \st $d=ab$ and \te\
$a'\in\a$, $c\in\g$ \st $e=a'c$. Hence $z =ab+a'c$. Set $a'':=\max(a,a')$ (see \E\df\ 2.3.17). Thus $a''\in\a$ since both $a$~and~$a'$ belong to~$\a$,
moreover $a,a'\le a''$, hence $z=ab+a'c \le a''b+a''c = a''(b+c)$. Then $b+c \in \b\hpl\g$ and $a''(b+c)\in \a\hcd(\b\hpl\g)$.

Note that we have proved that $(P,\hpl,\hcd,\wh0,\wh1)$ is a \sr. Observe that part (iv) of \E\df\ \rf{d3.7} follows from Lemma \rf{l4.36}. Since
$\wh1\in\dot P$ and $\a\hcd\b\in\dot P$ by the same lemma, we conclude that $\dot P$ is a \sbm\ of the monoid $(P,\hcd,\wh1)$. Hence part (iii) of
\E\df\ \rf{d3.7} holds. This completes the proof of Lemma \rf{l4.37}.
\endproof

We next prove part (i) of \E\df\ \rf{d3.7}.

\blm4.38a
The monoid $(P,\hpl,\wh0)$ is a \PM.
\elm

\proof \

``$\a\hpl\b = \wh0$ \ti{implies} $\a=\b=\wh0$'': Let $\a,\b\in P$ \sf y $\a\hpl\b = \wh0$. If $\a=\wh0$ then $\b=\wh0\hpl\b = \wh0$. Similarly, if
$\b=\wh0$ then $\a=\wh0$. If $\a\ne\wh0$ and $\b\ne\wh0$ then $\a,\b\in \dot P$, hence $\a\hpl\b\in\dot P$ by \er{4.39} and Lemma \rf{l4.35}. Since
$\dot P=P\sms{\wh0}$, this case is impossible.

It remains to show that if $\a,\b,\g\in P$ \sf y $\a\hpl\b = \a\hpl\g$, then $\b=\g$. If $\a=\wh0$ then $\b=\wh0\hpl\b = \wh0+\g = \g$. We now suppose
$\a\ne\wh0$ and $\b=\wh0$. We have $\a=\a\hpl\g$. If $\g=\wh0$ then $\b=\g$. We suppose $\g\ne\wh0$, thus $\a,\g\in\dot P$, hence $\a\hpl\g\in\dot P$.
We show that this case cannot occur. Indeed, if $\g\ne\vn$, then we obtain
\beq4.51a
\bigcup_{a\in\a,c\in\g}\{a+c\} = \a\hpl\g = \a.
\e
By \E\Pr\ \rf{p4.14} \te\ $c\in\g$, $\check a\in\a$, $\hat a\in\a^c$ \st $\hat a-\check a=c$. Then $\check a+c \in \a\hpl\g = \a$, hence $\check a
+ c\in\a$. \E\oh $\check a+c = \check a +(\hat a-\check a) = \hat a
\in\a^c$, a \cd ion. The case $\a\ne\wh0$, $\g=\wh0$, $\b\ne\wh0$ is similar since $\a\hpl\b = \a\hpl\g = \a$. Finally, we assume $\a\hpl\b = \a\hpl
\g$, $\a,\b,\g\in\dot P$. We suppose, for \cd ion, that $\b\ne\g$. By Lemma \rf{l4.6} we have either $\b>\g$ or $\g>\b$, where $\ge$~is defined in
\er{4.13}. We show that if $\b>\g$, then $\a\hpl\b \ne \a\hpl\g$. The case $\g>\b$ follows by interchanging $\b$ and~$\g$.

We proceed as in the proof of Theorem 134 of \cite{10a}. Since $\b>\g$, \te s $d\in\b\cap \g^c$. By~C3 \te s $d'\in\b$ \st $d'>d$. Then $d'\in\g^c$
by~D2 of Remark \rf{r4.12}\,(i). By \E\Pr\ \rf{p4.14} \te\ $\hat a\in\a^c$, $\check a\in\a$ \st $\hat a- \check a = d'-d$. \E\Tf\ $d' = d+(d'-d) =
d+(\hat a-\check a)$, hence $d'+\check a = d+((\hat a-\check a)+\check a) = d+\hat a$. We have $d'+\check a = \check a+ d'\in\a\hpl\b$ by \er{4.39}
and Lemma \rf{l4.35}. Note that $\hat a\in\a^c$ and $d\in\g^c$. We claim that $\hat a+d \notin \a\hpl\g$. Suppose, for \cd ion, that $\hat a+d
\in\a\hpl\g$. Then by \er{4.39} and Lemma \rf{l4.35} \te\ $a\in\a$, $c\in\g$ \st $\hat a+d = a+c$. Since $a\in\a$ and $\hat a\in\a^c$, $c\in\g$ and
$d\in\g^c$, we infer that $a<\hat a$ and $c<d$ by~L2. \E\Tf $a+c < \hat a+d$, a \cd ion. \If that $\a\hpl\g \ne \a\hpl \b$. \csq, $\a\hpl\b =
\a\hpl\g$ implies $\b=\g$.
\endproof

\brm4.39 \

\hph i,i, Inspection of the proof of Lemma \rf{l4.38a} shows that if $\b,\g\in\dot P$ \sf y $\b>\g$, then $\a\hpl\b > \a\hpl\g$. This is the content
of Theorem~134 of~\cite{10a}.

\hph ii,, So far, we have proved  that the \sr\ $(P,\hpl,\hcd,\wh0,\wh1)$ \sf ies \pp ies (i), (iii) and (iv) of \E\df\ \rf{d3.7}. It remains to show
that the \nog\ of the \PM\ $(P,\hpl,\wh0)$ is \ti{total\/}. An example of a \PM\ whose \nog\ is \ti{not\/} total is the \lt\ $(\Na,|)$ considered in
Section~3 of Chapter~3. The set $\Na$ inherits the \nog~$\le$ of the \PM\ $\pz1{\N}$, which is weaker than the \og~$|$. Indeed, we have $m|n$ implies
$m\le n$ by (3.1.14). This \pp y is used in the proof that $(\Na,|)$ is a \lt.
\erm

We now define a total \og\ on $P=\Id(\Q_{>0})\cup\{\vn\}$. We leave it to the reader to verify that the \rl~$\ge$ on~$P$ defined by
\beq4.52
\bca
\a\ge\b & \hbox{if $\a\supset\b$ \fa}\a,\b\in\dot P \hbox{ (see \er{4.13})},\\
\a\ge\wh0 & \hbox{\fa}\a\in P,
\eca
\e
is indeed a \ti{total\/} \og\ on~$P$ (see Lemma \rf{l4.6}).

\ssk
Our next goal is to show that the \nog\ of the \PM\ $(P,\hpl,\wh0)$ is identical to the \og~$\ge$. We first give another \chz\ of the \og~$\ge$.

\blm4.40
Let $\a,\b\in\dot P$. Then the \fw\ assertions are \ev t\dw

\hph i,ii, $\a>\b$,

\hph ii,i, \te s $z\in \a\cap \b^c$,

\hph iii,, \te s $c\in\Q_{>0}$ \st $\a > i(c)> \b$.
\elm

\proof \

``(i) \ti{implies} (ii)'': $\a>\b$ means $\a\supset\b$ and $\a\ne\b$. Thus \te s $z\in\a\cap\b^c$.

``(ii) \ti{implies} (iii)'': Let $z\in\a\cap\b^c$. By~C3 \te s $c\in\Q_{>0}$ \st $c>z$ and $c\in\a$. By~D2 (see Remark \rf{r4.12}\,(i)), $c\in\b^c$.
By~L2, $0<b<c$ \fa $b\in\b$. Hence $i(c)\ge \b$. Since $z\in i(c)\cap\b^c$ we have $i(c)\ne\b$. Finally, if $x\in i(c)$ then $x<c\in\a$. Hence
$x\in\a$ by~C2. Thus $\a\ge i(c)$. Since $c\in\a \cap(i(c))^c$, we have $\a\ne i(c)$.

``(iii) \ti{implies} (i)'': Let $c\in\Q_{>0}$ be \st $\a>i(c)>\b$. Then $\a\ge i(c)\ge\b$, hence by the \tr ity of~$\ge$, we have $\a\ge\b$. \Mo
$\a>\b$ follows from (1.3.11).
\endproof

We next show that the \og\ $\stackrel\hpl\ge$ is stronger than the \og~$\ge$.

\blm4.41
Let $\a,\b\in\dot P$ and suppose that \te s $\g \in\dot P$ \st $\a=\b\hpl\g$. Then $\a>\b$.
\elm

\proof[Proof \rm(\cite{10a})] Let $c\in\g$. By \E\Pr\ \rf{p4.14} \te\ $\check b\in\b$ and $\hat b\in\b^c$ \st $\hat b-\check b=c$. Hence
$\hat b=\check b+c \nde4.39 \in \b\hpl\g = \a$. \E\Tf $c\in\a\cap \b^c$, and $\a>\b$ by the previous lemma.
\endproof

For the sake of completeness we state the easily verified

\bco4.42
Let $\a,\b,\g\in P$ be \st $\a=\b\hpl\g$. Then $\a\ge\b$.
\eco

We now prove that the \og~$\ge$ is stronger than the \og~$\stackrel\hpl\ge$. Let $\a,\b\in P$. If $\a=\b$ then $\a=\b\hpl\wh0$ by \er{4.50}.
If $\a>\wh0$ then $\a\ne\wh0$, hence $\a\in\dot P$, and $\a=\a\hpl\wh0$ by \er{4.50}. It remains to prove

\bpr4.43
Let $\a,\b\in\dot P$ \sf y $\a>\b$. Then \te s $\g\in\dot P$ \st $\a=\b\hpl\g$.
\epr

Note that such a $\g$ is unique by Lemma \rf{l4.38a}.

\proof
In view of Lemma \rf{l4.40} \te s $c\in\Q_{>0}$ \st $\a> i(c)>\b$. We first show that \te s $\g_1\in\dot P$ \st $\a=i(c)\hpl\g_1$, and secondly, that
\te s $\g_2\in\dot P$ \st $i(c)=\b\hpl\g_2$. Then $\a= i(c)\hpl\g_1 = (\b\hpl\g_2)\hpl\g_1 = \b\hpl(\g_2\hpl\g_1)$. Set $\g:=\g_2\hpl\g_1$.

``\ti{Case} $\a >i(c)$'': Let $\a\in\dot P$ and $c\in\Q_{>0}$ \sf y $\a>i(c)$. We have $\a\supset i(c)$ by \er{4.52}. \Mo $\a = \bcl_{a\in\a}i(a)$
by~\er{4.17}. We claim that $\a = \bcl_{a\in\a,a>c}i(a)$. The inclusion $\supset$ is trivial. We prove the converse one. Since $\a> i(c)$, \te s
$z\in\a\cap i(c)^c$, that is, $z\in\a$ \st $z\ge c$. In view of~C3, \te s $z'>z$ \st $z'\in\a$. \E\Tf $z'\in\a$ and $z'>c$ by (1.3.11). Since
$\a=\bcl_{a\in\a}i(a)$, \te s $a'\in\a$ \st $z'\in i(a')$. Hence $z'<a'$, and $c<z'$ by (1.3.11). Finally, $\a = \bcl_{a\in\a}i(a) = \bg(\bcl_
{a\in\a,a\le c}{i(a)})\cup\bg(\bcl_{a\in\a,a>c}{i(a)}) \nad{C2,\er{4.10}} = i(c)\cup \bcl_{a\in\a,a>c}{i(a)}$. Since $c<z'\in i(a')\sbs \a$, we
have $i(c)\nde4.10 \subsetneq i(z')\sbs \bcl_{a\in\a,a>c}{i(a)}$. Thus $\a\sbs \bg(\bcl_{a\in\a,a>c}{i(a)}) \cup \bg(\bcl_{a\in\a,a>c}{i(a)}) =
\bcl_{a\in\a,a>c}{i(a)}$.

We have thus proved:
\beq4.53
\a = \bcl_{a\in\a,a>c}{i(a)} \qh{whenever $\a\in\dot P$ and $\a>i(c)$, $c\in\Q_{>0}$.}
\e
We recall that if $a,b\in\Q_{>0}$, then $i(a+b) \nde4.28 = i(a)+' i(b) \nad{\er{4.35},\er{4.39}}= i(a)\hpl i(b)$. We have $\a \nde4.53 =
\bcl_{a\in\a,a>c}i(a) = \bcl_{a\in\a,a>c}(i(c)+i(d))$ where $d\in\Q_{>0}$ \sf ies $a=c+d$, by \E\Pr\ 4.5.4, since $c<a$. We set $a-c:=d$. \E\Tf
$\a = \bcl_{a\in\a,a>c}i(a) = \bcl_{a\in\a,a>c} (i(c) \hpl i(a-c))$. We claim that if $D$~is a \nss\ of~$\Q_{>0}$ \sf ying
\beq4.54
\hbox{\te s $M\in \Q_{>0}$ \st $d\le M$ \fa $d\in D$,}
\e
then the \fw\ assertions hold:
\bga4.55
\bigcup_{d\in D}i(d) \in \dot P,\\
\bg(\bigcup_{d\in D} i{(d)}) \hpl i(e) = \bigcup_{d\in D}(i(d)\hpl i(e)) \qh{\fa}e\in\Q_{>0}.\lb{4.56}
\e

We first show that \er{4.55}, \er{4.56} allow us to complete the proof of the case ``$\a>i(c)$''.

We set $D:=\{a-c\in\Q_{>0}: a\in\a \hbox{ and }a>c\}$. $D$~is not empty since \te s $z'\in\a$ \st $z'>c$. \Mo by~L1 the set $\a^c$ is not empty.
Let $\wh a\in\a^c$. By~L2 we have $a<\wh a$ \fa $a\in\a$, hence $a-c<\wh a-c$ \fa $a\in D$. Indeed, $\wh a=a+\wh d$ \fs $\wh d\in\Q_{>0}$,
$a=c+(a-c)$, hence $\wh a=(c+(a-c))+\wh d = c+((a-c)+\wh d)$. Thus $a<\wh a$ since $(a-c)+\wh d>0$. \If that $d = a-c<\wh a-c$ \fa $d\in D$.
By \er{4.54}, \er{4.55} $\bcl_{a\in\a,a>c}i(a-c) \in\dot P$ and by \er{4.56} $\bg(\bcl_{a\in\a,a>c}{i(a-c)}) \hpl i(c) = \bcl_{a\in\a,a>c}\bigl(i(a-c)
\hpl i(c)\bigr) = \bcl_{a\in\a,a>c}i((a-c)+c) = \bcl_{a\in\a,a>c}i(a)\nde4.53 =\a$. Setting $\g_1:=\bcl_{a\in\a,a>c}i(a-c)$, we find $\a = i(c)+\g_1$.

It remains to prove \er{4.55}, \er{4.56}.

``\er{4.55}'': We have $d\le M$ \fa $d\in D$. Hence $i(d)\le i(M)$ \fa $d\in D$ by \er{4.10}. By \E\Pr\ \rf{p4.13} it suffices to show that
$\bcl_{d\in D}i(d) \ne \Q_{>0}$. We have $i(d)\le i(M)< i(M+1)$
\fa $d\in D$. Hence $\bcl_{d\in D}i(d) \le i(M) < i(M+1)$. Hence $i(M+1) \not\le \bcl_{d\in D}i(d)$, and $\bcl_{d\in D}i(d)\ne \Q_{>0}$. \If that
$\bcl_{d\in D}i(d)\in\dot P$.

``\er{4.56}'': We first prove the inclusion $\sbs$. Let $e\in\Q_{>0}$ and $z\in\bg(\bcl_{d\in D}{i(d)}) \hpl i(e)$. Then by \er{4.39} $z=x+e'$ where
$x\in i(d')$ \fs $d'\in D$ and $e'\in i(d)$. Thus by \er{4.39} $z\in i(d')\hpl i(e) = i(d'+e) \sbs \bcl_{d\in D} i(d+e) = \bcl_{d\in D}(i(d)\hpl
i(e))$, which proves the inclusion $\sbs$ for ``\er{4.56}''. Next we prove the converse inclusion for \er{4.56}. Let $z\in\bcl_{d\in D}(i(d)\hpl
i(e)) = \bcl_{d\in D} i(d+e)$. Then \te s $d'\in D$ \st $z\in i(d'+e)$. Then by \er{4.22} \te s $t\in(0,1)$ \st $z=t(d'+e) = td'+te$. By \er{4.22}
$td'\in i(d')$. Clearly, $i(d')\sbs \bcl_{d\in D} i(d)$, hence $z = td'+te' \in \bg(\bcl_{d\in D}{i(d)}) \hpl i(e)$, which proves the inclusion
$\supset$. \E\Tf \er{4.56} holds, and the proof of ``Case $\b>i(c)$'' is complete.

``\ti{Case} $i(c)> \b$'': Observe that if $\b=i(b)$ \fs $b\in\Q_{>0}$, then $i(b)<i(c)$, hence $b<c$ by \er{4.10}, and \te s $d\in\Q_{>0}$ \st $c=b+d$
by \E\Pr\ 4.5.4. From \er{4.26}, \er{4.28} we find that $i(c)=i(b)+' i(d)$. We set $\g_2:= i(d) = i(c-b)$.

Let $\b\in\dot P\sm i(\Q_{>0})$ be \st $i(c)>\b$. We have $\b=\bcl_{b\in\b}i(b)$ by \er{4.17} and $\b=\Bca_{\hat b\in\b^c}i(\wh b)$ by \er{4.19}.
Then we have
\beq4.57
i(b) \sbs \b \sbs i(\wh b) \sbs i(c)
\e
\fa $b\in\b$, $\wh b\in\b^c$ \sf ying $\wh b<c$. We set
\beq4.58
B:=\{\wh b\in\b^c: \wh b<c\}.
\e
Then $B$ is not empty by Lemma \rf{l4.40} since we assumed $i(c)>\b$.

Next we denote by $\g$ the subset of $\Q_{>0}$ defined by
\beq4.59
\g:=\bigcup_{\hat b\in B}i(c- \wh b).
\e
We show that $\g\in\dot P$. Indeed, we have $(c-\wh b)+\wh b = c$, hence $c-\wh b< c$ \fa $\wh b\in B$. \E\Tf $\g\in\dot P$ in view of \er{4.55}
since \er{4.54} holds. We claim
\beq4.60
\b\hpl\g \le i(c).
\e
Indeed, let $z\in\b\hpl\g$. \E\te\ $x\in\b$ and $y\in\g$ \st $z=x+y$ by \er{4.39}. By \er{4.59} \te s $\wh b\in B$ \st $y\in i(c-\wh b)$, that is,
$y<c-\wh b$. Hence $z = x+y = y+x \nad{(4.5.120)} < (c-\wh b)+x$. Since $x\in\b$ and $\wh b\in\b^c$, we have $x<\wh b$ by~L2. Hence $(c-\wh b)+x
= x+(c-\wh b) < \wh b + (c-\wh b) = c$. \If by \tr ity that $z<c$, thus $z\in i(c)$, which proves the claim.

Next we prove
\beq4.65
i(c)\le \b\hpl\g,
\e
that is, $a\in \b\hpl\g$ \fa $a\in(0,c)$. We first prove
\beq4.66
i(c') \le \b\hpl\g \qh{\fa $c'\in B$},
\e
where $B$ is defined in \er{4.58}. We recall that $B\ne\vn$ by Lemma \rf{l4.40}. Let $c'\in B$ we have $c'<c$. Let $a\in(c',c)$, that is, $a-c'>0$
and $a<c$. By \E\Pr\ \rf{p4.14} \te\ $b\in\b$ and $\wh b\in\b^c$ \st $\wh b-b = a-c'$. (Note that $\wh b>b$ by~L2 hence $\wh b-b>0$.)

We claim that
\beq4.67
c'-b = a-\wh b \ \hbox{ and }\ a-\wh b>0.
\e
Indeed, using also \ng\ \ra\ \nm s, we have $c'-b:= c'+(-b) = c'+(-\wh b)+\wh b+(-b) = c'+(-\wh b)+(\wh b-b) = (\wh b-b) + (c'+(-\wh b)) =
a+(-c') + c' +(-\wh b) = a+0+(-\wh b) = a+(-\wh b) =: a-\wh b$. \Mo we have $a-\wh b = c'-b >0$ by~L2 since $b\in\b$ and $c'\in B\sbs \b^c$.
We leave it to the reader to prove \er{4.67} using only non\ng\ \ra\ \nm s, where $x-y$ is defined for $x\ge y$ by $x=(x-y)+y$, $x,y\in\Q_{>0}$.

We obtain
\beq4.68
c' = b+ (c'-b) = b+(a-\wh b).
\e
We have $b\in\b$ and claim that $a-\wh b\in\g$. Since $a<c$, we have $0\nde4.67 < a-\wh b < c-\wh b$, hence $a-\wh b\in i(c-\wh b)$. Since $\wh b\in
\b^c$, $i(c-\wh b)\sbs\g$ in view of \er{4.59}. \E\Tf $a-\wh b\in\g$. From $b\in\b$, $a-\wh b\in\g$ and \er{4.39} we infer $c'\in\b\hpl\g$. We recall
that $\b\hpl\g \in \dot P$, hence by~C2 for $\b\hpl\g$ we have $d\in\b\hpl\g$ \fa $d<c'$. Since $c'\in B$ is arbitrary, we have proved \er{4.66}.

Now let $a\in(0,c)$. We show that $a\in\b\hpl\g$. In view of \er{4.66} we may assume $a\notin\b^c$, that is, $a\in\b$. Let $d\in B$. Then $d\in\b^c$
and $d<c$. Hence $a<d$ by~L2 and $d<c$.
Thus $a\in i(d)$ where $d\in B$. Since $i(d)\le \b\hpl\g$ by \er{4.66}, we obtain $a\in\b\hpl\g$ \fa $a\in(0,c)$. Hence \er{4.65} holds.
In view of \er{4.60} we obtain $i(c)=\b\hpl\g$. \E\Tf the proof of \E\Pr\ \rf{p4.43} is complete.
\endproof

Combining Lemmata \rf{l4.37} and \rf{l4.38a}, we find that the set $(P,\hpl,\hcd,\wh0,\wh1)$, introduced in Lemma \rf{l4.37}, is a \sr\ \sf ying \cn s
(i), (iii) and (iv) of \E\df\ \rf{d3.7}. \If from Corollary \rf{c4.42} and \E\Pr\ \rf{p4.43} that the \nog\ $\stackrel\hpl\ge$ of the \hbox{\PM}\
$(P,\hpl,\wh0)$ defined in (3.1.4) is identical to the total \og~$\ge$ on~$P$ defined in \er{4.52}. \E\Tf \cn\ (ii) of \E\df\ \rf{d3.7} is \sf ied,
hence $(P,\hpl,\hcd,\wh0,\wh1)$ is an ordered semifield. \csq, by Theorem \rf{t3.8}\,(ii) \te s an \of\ whose \po\ cone is equal to~$P$.

We have thus proved

\bth4.30
\E\te s an \of\ $(\R,\ge)$ whose \po\ cone $\R_{\ge0}$ $(:=\{x\in\R:\break x\ge0\})$ is the ordered semifield $(P,\hpl,\hcd,\wh0,\wh1)$ introduced in
Lemma \rf{l4.37}. \Mo the \og\ of the \PM\ $(P,\hpl,\wh0)$ defined in {\rm(3.1.4)} is identical to the \og~$\ge$ defined in \er{4.52}.
\eth

Our next goal is to investigate \pp ies and \ev t \df s of the field~$\R$. It turns out that the field~$\R$ is \Ar. Note that the proof of Lemma
4.5.54 was left as an exercise, and we recall that a proof is given in the discussion preceding \E\df~\rf{d3.35}.

\bpr4.31
The \of\ $(\R,\ge)$ is \emph{\Ar}.
\epr

\proof
We want to prove that \cn\ (i) of Lemma 4.5.54 is \sf ied, that is, given $\a\in\dot P=P\sms{\wh0}$ \te s $n\in\Na$ \st $n\stackrel\hpl\cdot \wh1$,
the $n$-th \IT\ of~$\wh1$ in the monoid $(P,\hpl,\wh0)$, is greater than~$\a$. We recall that $\wh1\in\dot P$ by \er{3.12}. We also recall that the
map $f:\Q_{\ge0} \to i(\Q_{>0})\cup\wh0$ (where $\wh0$ is the empty subset of~$\Q_{\ge0}$), defined in \er{4.26} is a monoid-\is sm from
$(\Q_{\ge0},+,0)$ onto $(i(\Q_{>0})\cup\wh0,+',\wh0)$ by Lemma \rf{l4.17}\,(ii). Let $\a\in\dot P:=P\sms{\wh0}$. Then $\a<\a\hpl\a$ by (3.1.4). By
Lemma \rf{l4.40} \te s $c\in\Q_{>0}$ \st $\a<i (c)<\a\hpl\a$. If $c\le\frac11$, then $c<\frac11+\frac11 = 2\dpl\frac11$, the second \IT\ of~$\frac11$
in the monoid $(\Q_{\ge0},+,\frac11)$ (see \E\Pr\ 4.5.4). If $c>\frac11$, then by Lemma 4.5.52 \te s $n\in\Na$ \st $c<\frac n1\cdot\frac11
\nad{(4.5.3)}= \frac n1$. We claim that $\frac n1 = n\dpl\frac11$ \fa $n\in\Na$, where $n\dpl\frac11$ is the $n$-th \IT\ of~$\frac11$ in the monoid
$(\Q_{\ge0},+,\frac11)$. We use \In\ on $n\in\Na$. We have $1\dpl\frac11 = \frac11$ by (2.2.3)\,I1. \Mo $(k+1)\dpl\frac11 = (k\dpl\frac11)+
1\dpl\frac11$ \fa $k\in\Na$ by (2.2.3)\,I2. We have $\frac11 = 1\dpl\frac11$. We show that $l\dpl\frac11= \frac l1$ implies $(l+1)\dpl\frac11 =\frac
{l+1}1$ \fa $l\in\Na$. Indeed, $(l+1)\dpl\frac11 = l\dpl\frac11 + 1\dpl\frac11 = \frac l1+\frac11 \nad{(4.5.2)} = \frac{l\cdot1+1\cdot1}{1\cdot1}
=\frac{l+1}1$. Set $A:=\{k\in\Na: k\dpl\frac11 = \frac k1\}$. Then $1\in A$ and $k\in A$ implies $k+1\in A$. \If that $A=\Na$, which proves the claim.
\E\Tf $c<n\dpl\frac11$, and $i(c)\subsetneq i(n\dpl\frac11)$ by \er{4.10}. We recall that $i(n\dpl\frac11) = f(n\dpl\frac11)$ by \er{4.26}. By Lemma
2.1.22 we find that $f(n\dpl\frac11) = n\stackrel{+'}\cdot f(\frac11) = n\stackrel{+'}\cdot i(\frac11) = n\stackrel{+'}\cdot \wh 1$. \Mo we have
$i(r)+'i(s) = i(r)\hpl i(s)$ \fa $r,s\in\Q_{>0}$ by \er{4.35}, \er{4.37} and \er{4.39}. \E\Tf the monoid $(i(\Q_{>0})\cup\{\wh0\},+',\wh0)$
is equal to the monoid $(i(\Q_{>0})\cup\{\wh0\},\hpl,\wh0)$. Hence $i(n\dpl1)= n\stackrel{\hpl}\cdot\wh1$. \csq, $\a<i(c) < i(n\dpl\frac11) =
n\stackrel{\hpl}\cdot i(\frac11) = n\stackrel{\hpl}\cdot\wh1$. \If that $\a<n\stackrel{\hpl}\cdot\wh1$, which shows that \cn~(i) of Lemma 4.5.54 is
\sf ied. Thus the \of\ $(\R,\ge)$ is \Ar.
\endproof

We now consider \ti{\ad al \pp ies} of \Ar\ \of s. Let $(K,+,\cdot,{\bf0},{\bf1})$ be an \Ar\ \of\ and let $K_{\ge0}$ denote the \crs\ ordered
semifield. Let $\check K$ denote the \pf\ of~$K$ (see \E\df\ 4.4.29 and \E\Pr\ 4.4.30). We know that, as an \of, the field $K$ has \ch istic~$0$. Let
$j:\Q\to \check K$ denote the ring-\is sm defined in (4.5.108). Then, if $r\in\Q_{>0}$ and $r=\frac nm$ with $n,m\in\Na$, we have
\beq4.69
j(r) = (n\dpl{\bf1})\cdot (m\dpl{\bf1})\mo, \q m,n\in\Na.
\e
We recall that this \df\ makes sense since the \fw\ holds by (4.5.109):
\beq4.70
(n\dpl{\bf1})\cdot(m\dpl{\bf1})\mo = (l\dpl{\bf1})\cdot(k\dpl{\bf1})\mo \qh{iff }nk = ml.
\e
We shall use the notation
\beq4.71
r{\bf1} := j(r) \qh{\fa $r\in\Q$},
\e
in order to emphasize the role played by~${\bf1}$, the \nel\ of the \am\ $(K,\cdot,{\bf1})$. \E\Ip if $r=\frac11$, we have $r{\bf1} =
(1\dpl{\bf1})\cdot(1\dpl{\bf1})\mo = {\bf1}\cdot({\bf1})\mo = {\bf1}$. Hence $1\cdot{\bf1} = {\bf1}$. \Mo if $r,s\in\Q_{>0}$, then $(r+s){\bf1} =
j(r+s) = j(r)+j(s) = r{\bf1}+s{\bf1}$ and $(r\cdot s){\bf1} = j(rs) = j(r)\cdot j(s) = r{\bf1}\cdot s{\bf1}$. Similarly, $0{\bf1} = j(0)= {\bf0}$. All
these formulae are a rewriting in terms of the notation $r{\bf1}$ of the fact that the map $j:\Q \to\check K$, the \Pf\ of~$K$, is a ring-\hm sm.
Since $j:\Q \to \check K$ is a bi\jn, we have
\beq4.72
\check K:= \{r{\bf1}: r\in\Q\}.
\e

Observe that we did not use the fact that the field~$K$ is an \Ar\ \of. We only used the fact that $K$~has \ch istic~$0$. \If from Example 5.1.2\,(ii)
that the pair $(\check K,K)$ endowed with the map $(\la,x)\mt\la x$ from $\check K\times K$ into~$K$ defined by $\la\cdot x$, where $\cdot$~is the
\mlc\ in~$K$, is a \ti{\vs} over the field~$\check K$ (see \E\df\ 5.1.1). Since \fa $\la\in\check K$, \te s \ooo $r\in\Q$ \st $\la=j(r)$, the \mlc\
of the vector $x\in K$ by the scalar $\la\in\check K$ \sf ies $\la x=j(r)\cdot x = (r{\bf1})\cdot x$. \E\Ip $\la\cdot{\bf1} = j(r)\cdot{\bf1} = j(r)
=r{\bf1}$. Setting
\beq4.73
rx := j(r)\cdot x, \q r\in\Q,\ x\in K,
\e
we find that the map $(r,x)\mt rx$ from $\Q\t K$ into~$K$ \sf ies (5.1.1).

We recall that $j$ is the only ring-\is sm from $\Q$ onto $\check K$. Indeed, if $j'$ is such an \is sm, then the \cm\ $(j')\Inv\circ j$ is a
ring-auto\mf\ of the field~$\Q$. Since the identity is the only ring-auto\mf\ of~$\Q$ by \E\Pr\ \rf{p3.2}, we have $j=j'$. We obtain

\bpr 4.31a
Let $(K,+,\cdot,{\bf0},{\bf1})$ be a field of \ch istic~$0$. Then $K$ is a \vs\ over the field~$\Q$ where the \mlc\ by scalars is defined by
\er{4.73}.

The \pf\ of $K$ denoted by $\check K$ is a \lss\ of~$K$ of \dm\ one, and $\bf1$ is a basis of $(\Q,\check K)$.
\epr

\proof
We check (i)--(iii) of \E\df\ 5.1.1.

\hph i,ii, $(K,+,\Q)$ is an \ag\ since $K$ is a ring.

\hph ii,i, $\Q$ is a field by Theorem 4.5.6.

\hph iii,, Let $r,s\in\Q$ and $x,y\in K$. Then

``I0'': $0x = j(0)\cdot x = {\bf0}\cdot x \nad*= {\bf0}$.

``I1'': $1x = j(1)\cdot x = {\bf1}\cdot x \nad*= x$.

``I2'': $(r+s)x = j(r+s)\cdot x = (j(r)+j(s))\cdot x \nad*= j(r)\cdot x+j(s)\cdot x = rx+sx$, \ $r,s\in\Q$.

``I3'': $(r\cdot s)x = j(rs)\cdot x = (j(r)\cdot j(s))\cdot x = j(r)\cdot(j(s)\cdot x) = j(r)\cdot(sx) = r(sx)$.

``I4'': $r{\bf0} = j(r)\cdot{\bf0} \nad*= {\bf0}$.

``I5'': $r(x+y) = j(r)\cdot(x+y) \nad*= j(r)\cdot x+j(r)\cdot y = rx+ry$.

In $\nad*=$ we used the fact that $(K,+,\cdot,{\bf0},{\bf1})$ is a \sr.

Let $\check K$ be the \pf\ of~$K$, let $r\in\Q$ and $x,y\in \check K$. Since $j:\Q\to \check K$ is a ring-\is sm, \te\ \ooo $r\ (\hbox{resp.}~s)\in\Q$
\st $j(r)=x$ (resp. $j(s)=y$). Then if $t\in\Q$, we have $tx = j(t)\cdot x = j(t)\cdot j(r) = j(t\cdot r)\in \check K$. \Mo $x+y = j(r)+j(s) =
j(r+s)\in \check K$. Hence $\check K$~is a \lss\ of~$K$ in view of \E\df\ 5.2.6. By Lemma 4.4.25 $\check K$~is a field. Clearly, it has \ch istic~$0$.
As a subfield of~$K$, we have ${\bf0},{\bf1} \in \check K$. Hence, if $x\in\check K$, then $x=j(r) = j(r)\cdot{\bf1} = r{\bf1}$ \fs $r\in\Q$.
Conversely, \fa $s\in\Q$, $s{\bf1} = j(s)\cdot{\bf1} = j(s)\in \check K$. Thus $\check K:=\{r{\bf1}: r\in\Q\}$. Note that $\Q$~is a \vs\ over~$\Q$ by
what precedes, since $\Q$ has \ch istic~$0$. The identity being the only ring-auto\mf\ of~$\Q$, we have $j=\id$. Hence the \mlc\ by scalars is given
by $rs = \id(r)\cdot s = rs$, $r,s\in\Q$. \If that the \Op\ identity~$\id$, in the \vs\ $(\Q,\Q)$ is a bi\jc\ \ti{linear} \Op\ since we also trivially
have $\id(r+s) = r+s = \id(r)+\id(s)$ \fa $r,s\in\Q$. \E\Tf the \dm\ of $(\Q,\Q)$ is one by \E\df\ 4.1.13. Then {\bf1} is a basis of this \vs, since
$r=r{\bf1}$ \fa $r\in\Q$. Similarly, the map $j:\Q\to\check K$ is a bi\jc\ linear map, hence $\check K$ is a one-\dm al \lss\ of the \vs~$K$
over~$\Q$, with {\bf1} as a basis. Clearly, $r{\bf1}$ is also a basis of $\check K$ provided $r\ne0$.
\endproof

\bco4.32
Every \of\ $(K,\ge)$ can be viewed as a vector field over~$\Q$.

Let $\check K$ denote the \pf\ of~$K$ and $j:\Q\to\check K$ be the map defined in \er{4.69}. Then $j:(\Q_{>0},\ge)\to(\check K,\ge)$ is an \ois sm.
\E\Ip we have
\beq4.74
r<s \qh{iff }r{\bf1}= j(r) < j(s) = s{\bf1}, \q s,r\in\Q.
\e
\eco

\proof
An \of\ $(F,+.\cdot,{\bf0},{\bf1},\ge)$ has \ch istic~$0$ by (4.5.126). (Use \er{3.11}, \er{3.12}, (2.2.1)\,I1, (2.2.1)\,I2 and \In\ on $n\in\Na$.)
Thus the first assertion of the corollary follows from \E\Pr\ \rf{p4.31a}. The second assertion follows from Lemma \rf{l3.28} with $(F,\ge):=
(\Q,\ge)$, $(F',\ge'):=(\check K,\ge)$, and $\vf:=j$, provided that $j(r)=0$ whenever $r>0$ \fa $r\in\Q_{>0}$. Let $r>0$ and $m,n\in\Na$ be \st
$r=\frac nm$. Then $j(\frac nm) =\break (n\dpl{\bf1})\cdot(m\dpl{\bf1})\mo$. \E\fa $n,m\in\Na$, $n\dpl{\bf1}>0$, $m\dpl{\bf1}>0$ by what precedes. \Mo
$(m\dpl{\bf1})\mo>0$ \fa $m\in\Na$ by Lemma \rf{l3.18}. By \er{3.11} $(n\dpl{\bf1})\cdot(m\dpl{\bf1})\mo>0$ \fa $n,m\in\Na$.
\endproof

We now assume that the field $(K,+,\cdot,{\bf0},{\bf1},\ge)$ is an \Ar\ \of\ (or simply \Ar\ field). \If from \E\Pr\ \rf{p4.1} that given $x,y>0$
\st $y<x$, \te s $r\in\Q_{>0}$ \sf ying $y< r{\bf1}<x$. Note that if $y<\bf0$ and ${\bf0}<x$ then $y< {\bf0}<x$ with ${\bf0} = 0{\bf1}$. \Mo if
$y<x<{\bf0}$, then we have ${\bf0}<-x<-y$ by \er{3.9}, hence \te s $r\in\Q_{>0}$ \st $-x<r <-y$. \csq, we have $y=-(-y) < -r <-(-x)= x$. Note that
${\bf0}=0{\bf1} = (r+(-r)){\bf1} = r{\bf1} + (-r){\bf1}$ \fa $r\in\Q$. Hence
\beq4.75
-(r{\bf1}) = (-r){\bf1} \qh{\fa $r\in\Q$.}
\e
Thus the \fw\ holds:
\beq4.76
\hbox{\E\fe pair $x,y\in K$ with $x<y$, \te s $z\in\check K$ \st $x<z<y$.}
\e

\begin{dfn}[{\cite[p. 23]{Nrs}}]\lb{d4.33}
Let $(F,+,\cdot,{\bf0},{\bf1},\ge)$ be an \of, and let $H$ be a subfield of~$F$. The subfield~$H$ is said to be \ti{dense} in~$F$ if \fe pair $x,y
\in F$ with $x<y$, \te s an \el\ $z\in H$ \st $x<z<y$.\index{dense in an \of}
\edf

\bpr4.34
Let $(K,\ge)$ be an \of. Then the \fw\ assertions are \ev t\dw

\hph i,i, $K$ is \Ar,

\hph ii,, The \pf\ of $K$ is dense in $(K,\ge)$.
\epr

\proof \

``(i)$\Rightarrow$(ii)'': follows from \er{4.76}.

``(ii)$\Rightarrow$(i)'': Let $x\in K_{>0}$. Then $x={\bf0}+x \nde3.9 < {\bf1}+x = x+{\bf1}$. Then \te s $r\in\Q$ \st $x<r{\bf1} <x+{\bf1}$ by~(ii).
Since $(\Q,\ge)$ is \Ar, \te s $n\in\Na$ \st $r<n\dpl\frac11$. Hence $j(r)< j(n\dpl\frac11)$ by Corollary \rf{c4.32}. \Mo $j(r\dpl\frac11) =
n\dpl j(\frac11) = n\dpl{\bf1}$ by (2.1.45) since $j$ is a ring-\hm sm. Hence $x<n\dpl{\bf1}$, and $(K,\ge)$ is \Ar.
\endproof

We now prove the result promised at the end of Section~3.

\begin{thm}[{\cite[p. 413]{Is}, \cite[p. 241]{Alg}}]\lb{t4.35}
Every \Ar\ \of\ is order- and ring-\is c to a subfield of the field~$\R$.
\eth

In a first step we construct an in\jc\ \hm sm from the ordered semifield $K_{\ge0}$ into the ordered semifield $(\R_{\ge0},\ge)$. In a second (easy)
step we extend this \hm sm to~$K$.

Capital (resp.\ small) letters will be used to denote \el s of the ordered semifield $(K_{\ge0},+,\cdot,{\bf0},{\bf1})$ (resp.\ $\pz0{\Q_{\ge0}}$).

We define a map $i_K: K_{\ge0} \to\R_{\ge0}$ by setting:
\beq4.77
\bca
{}& i_K(X):= \{r\in\Q_{>0}: r{\bf1}<X\} \hbox{ for }X\in K_{>0},\\
& i_K({\bf0}) := \vn, \hbox{ the empty subset of }\Q_{\ge0}.
\eca
\e
We recall that $r{\bf1}$ is defined in \er{4.71}. If $K:=\Q$, the map $i_K$ is just the map $f:\Q_{\ge0}\to\R_{\ge0}$ defined in \er{4.26}.

\blm4.36a
Let $(K,\ge)$ be an \Ar\ \of\ and let $i_K: K_{>0}\to\cP(\Q_{>0})$ be defined by \er{4.77}. Then $i_K(X) \in \R_{>0}$ \fa $X\in K_{>0}$.
\elm

\proof \

``C1'': By Lemma 4.5.54\,(ii) \te s $m\in\Na$ \st $(m\dpl{\bf1})\mo <X$, and by \er{4.3} \te s $r\in\Q_{>0}$ \st $(m\dpl{\bf1})\mo< r{\bf1} <X$. \Mo
$X\notin i_K(X)$, hence $i_K(X)\ne \Q_{>0}$.

``C2'': If $r\in i_K(X)$, and $s\in\Q_{>0}$ \sf ies $s<r$, then $r{\bf1} <X$ by \er{4.77}, $s{\bf1}<r{\bf1}$ by \er{4.74}, hence $s{\bf1}<X$, and
$s\in i_K(X)$.

``C3'': Let $r\in i_K(X)$. Then $r{\bf1} < X$. By \er{4.76} \te s $s\in\Q_{>0}$ \st $r{\bf1} < s{\bf1} <X$. Hence $r<s$ by \er{4.74}, and $s\in
i_K(X)$. \csq, $i_K(X)$ is an ideal of $\Q_{>0}$, hence an \el\ of $\R_{>0}$ in view of Theorem \rf{t4.30}.
\endproof

\blm4.37a
Let $(K,+,\cdot,{\bf0},{\bf1},\ge)$ be an \Ar\ field, and let $r\in\Q_{>0}$, $X,Y \in K_{>0}$. Then the \fw\ assertions hold\dw

\hph i,ii, $i_K(r{\bf1}) = i(r)$\sd

\hph ii,i, $s\in i_K(rX)$ iff $r\mo s\in i_K(X)$ \fa $s\in\Q_{>0}$\sd

\hph iii,, $X<Y$ iff $i_K(X) < i_K(Y)$\sd

\hph iv,, $i_K(X) \hpl i_K(Y) \le i_K(X+Y)$\sd

\hph v,i, $i_K(X) \hcd i_K(Y) \le i_K(X\cdot Y)$\sd

\hph vi,, \fa $e\in\Q_{>0}$ \te s $x\in\Q_{>0}$ \sf ying
\beq4.78
x{\bf1} <X < (x+e){\bf1}.
\e
\elm

The \og\ used on the \RHS\ of (iii) is the \og\ on~$\R$ defined in \er{4.52}, and the binary \op s $\hpl$ (resp.~$\hcd$) are defined in \er{4.39}
(resp.~\er{4.40}).

\proof \

(i) Let $s\in\Q_{>0}$. Then by \er{4.74}, $s{\bf1} < r{\bf1}$ iff $s<r$. Thus $s\in i_K(r{\bf1})$ iff $s{\bf1}< r{\bf1}$ iff $s<r$ iff $s\in i(r)$.

(ii) Let $s\in\Q_{>0}$. Then $s\in i_K(rX)$ iff $s{\bf1} < rX$. Since $s{\bf1} = j(s)\cdot{\bf1}$, $rX=j(r)\cdot X$ by \er{4.73}, and $j(r)\cdot X
= j(r{\bf1})\cdot X = (rj(1))\cdot X = (r{\bf1})\cdot X$ we have $s\in i_K(rX)$ iff $j(s)\cdot {\bf1}< j(r)\cdot X$ iff $(j(r)\mo \cdot j(s)) \cdot
{\bf1} < {\bf1}\cdot X = X$ iff $(r\mo s){\bf1} = (j(r)\mo j(s))\cdot {\bf1}=j(r)\mo j(s) <X$ iff $r\mo s \in i_K(X)$.

(iii) The proof of (iii) is similar to the proof of \er{4.10}.

(iv), (v) Set ${\qu} :=+$ (resp.~$\cdot$) in case (iv) (resp.~(v)), and ${\hqu}:={\hpl}$ (resp.~$\hcd$) in case (iv) (resp.~(v)). Let $u\in i_K(X)
\hqu i_K(Y) \nad{\er{4.39},\er{4.40}} = \bcl_{a\in\a,b\in\b}\{a\qu b\}$ where $\a:=i_K(X)$ and $\b:=i_K(Y)$. Then \te\ $r\in\a$ and $s\in\b$ \st
$u=r\qu s$.
We have $j(u) = j(r)\hqu j(s)$. Since $r\in i_K(X)$, $j(r) = j(r)\cdot {\bf1} = r\cdot{\bf1} <X$. Similarly, $j(s)<Y$. In view of \er{3.9}, \er{3.10}
and the \cmt ity of~$\qu$, we obtain $j(u) = j(r)\qu j(s) < X\qu j(s) < X\qu Y$. Hence $u{\bf1} = j(u){\bf1} < X\qu Y$. \If that $u\in i_K(X\qu Y)$.

(vi) \E\te s $s \in \Q_{>0}$ \st $0<s{\bf1} <X$ by \er{4.3} where $y:=\frac12 X$ and $x:=X$. We distinguish two cases: $X\le e{\bf1}$ and
$e{\bf1} <X$.

``$X\le e{\bf1}$'': Then $X\le e{\bf1} \nde 3.9 < s{\bf1} + e{\bf1} = (s+e){\bf1}$. Set $x:=s$.

``$e{\bf1} < X$'': By \E\df\ 4.5.53 \te s $k\in\Na$ \st $X < k\dpl e{\bf1}$. Set $A:=\{n\in\Na: X<n\dpl e{\bf1}\}$. Since $A\ne\vn$, the set $A$
as a subset of~$\Na$ has a least \el~$k$. We have $k>1$ since $1\dpl e{\bf1} = e{\bf1} <X$. Set $n:=k-1>0$. Then $n\dpl e{\bf1}
<X <(n+1)\dpl e{\bf1} = n\dpl e{\bf1} + 1\dpl e{\bf1} = n\dpl e{\bf1}+ e{\bf1}$. Moreover, $n\dpl e{\bf1} \nde4.73 =\break n\dpl j(e)\nad{(2.1.45)}
= j(n\dpl e) \nde 4.71 = (n\dpl e){\bf1}$. Set $x:= n\dpl e$, the $n$-th \IT\ of~$e$ in the \PM\ $\pz1{\Q_{\ge0}}$. Since $e>0$, we have $x>0$.
Thus \er{4.78} holds.
\endproof

\blm4.38b
Let $K$, $X$ and $Y$ be as in Lemma \rf{l4.37a}, and let $u\in\Q_{>0}$. Then\dw

\hph i,i, If\/ $u{\bf1} < X+Y$, then \te\ $r,s\in\Q_{>0}$ \st $r{\bf1} < X$, $s{\bf1} < Y$ and $u = r+s$.

\hph ii,, If\/ $u{\bf1} < X\cdot Y$, then \te\ $r,s\in\Q_{>0}$ \st $r{\bf1} < X$, $s{\bf1} < Y$ and $u = rs$.
\elm

\proof \

(i) Let $u{\bf1} < X+Y$. By \er{4.76} \te s $z\in\Q_{>0}$ \st $u{\bf1} < z{\bf1} < X+Y$. Similarly, \te s $z'\in\Q_{>0}$ \st $z{\bf1} < z'{\bf1} <
X+Y$. Set $h:=z'-z\in\Q_{>0}$. By Lemma \rf{l4.37a}\,(vi) \te\ $x,y\in\Q_{>0}$ \st $x{\bf1} < X < (x+\frac12 h){\bf1}$ and $y{\bf1} <Y < (y+\frac12
h){\bf1}$. \E\Tf $(z+h){\bf1} < X+Y < ((x+y)+h){\bf1}$, hence $(z+h){\bf1} < ((x+y)+h){\bf1}$, and $z+h < (x+y)+h$. \csq, $z<x+y$ by \er{3.9}. Since
$u{\bf1}<z{\bf1}$, hence $u<z$ by \er{4.74}, we obtain $u<x+y$. We use the notation $\frac c{a+b} := c(a+b)\mo$, $a,b,c\in\Q_{>0}$. Thus $\frac u{x+y}
<1$ and setting $r:=\frac u{x+y}x$, $s:=\frac u{x+y}y$, we find $r<x$, $s<y$, $r{\bf1} < x{\bf1} < X$, $s{\bf1} < y{\bf1} <Y$. Finally, $u=u1 =
u\frac{x+y}{x+y} = u\frac x{x+y} + u\frac y{x+y} = x\frac u{x+y} + y\frac u{x+y} = r+s$, where $r\h {<}\h x$, $s\h{<}\h y$, hence $r{\bf1}\h{<}\h
x{\bf1} \h{<}\h X$ and $s{\bf1} \h{<}\h y{\bf1} \h{<}\h Y$. Thus $u\h{=}\h r+s$ and $r{\bf1} \h{<}\h X$, $s{\bf1} \h{<}\h Y$.

(ii) (\cite[p.~414]{Is}) Let $u\in\Q_{>0}$ be \st $u{\bf1}<X\cdot Y$. Since $Y>0$, we have $Y\mo>0$ by Lemma \rf{l3.18}. Hence $Y\mo\cdot u{\bf1}
\nde3.10 < Y\mo\cdot(X\cdot Y) = Y\mo\cdot(Y\cdot X) = (Y\mo\cdot Y) \cdot X = {\bf1}\cdot X = X$. By \er{4.3} \te s $r\in\Q_{>0}$ \st $Y\mo \cdot
u{\bf1} < r{\bf1} < X$. Set $s:=r\mo u$. We have $r\mo>0$ by the same lemma and $s>0$ by \er{3.11}. Hence $rs = r(r\mo u) = (rr\mo)u = 1u = u$. Since
$u=rs$ and $r{\bf1} < X$ it remains to show that $s{\bf1} < Y$, \ev tly $(r\mo u){\bf1} <Y$. From $Y\mo\cdot u{\bf1} < r{\bf1}$ we infer $u{\bf1} <
(r{\bf1})\cdot Y$ by \er{3.10}. Then $(r\mo u){\bf1} = j(r\mo u) = j(r\mo)\cdot j(u) = (r\mo {\bf1})\cdot(u{\bf1}) < (r\mo{\bf1})\cdot(r{\bf1}) \cdot
Y = {\bf1}\cdot Y = Y$. \E\Tf $s{\bf1} = (r\mo u){\bf1} <Y$.
\endproof

As a direct con\sq\ of Lemmata \rf{l4.37a} and \rf{l4.38b}, we obtain

\bco4.39
Under the \as s of Lemma \rf{l4.37a}, the \fw\ holds\dw

\hph i,ii, $i_K(X \hpl Y) \ge i_K(X+Y)$.

\hph ii,i, $i_K(X \hcd Y) \ge i_K(X\cdot Y)$.

\hph iii,, The map $i_K : K_{\ge0} \to \R_{\ge0}$ is an in\jc\ order- and $($semi$)$ring-\hm sm.
\eco

Then Theorem \rf{t4.35} follows from Corollary \rf{c4.39}\,(iii) and

\bpr4.40
Let $(F,\ge)$ and $(F',\ge')$ be \of s, and let $\vf: F_{\ge0}\to F'_{\ge'0'}$ be an in\jc\ \sr\ \hm sm.

Define a map $\ov \vf : F\to F'$ by setting
\beq4.79
\ov\vf(x) = \bca
\vf(x) &\hbox{if } x\ge0,\\
-\vf(-x) & \hbox{if } x<0.
\eca
\e
Then $\ov\vf$ is an in\jc\ order- and ring-\hm sm from $F$ into~$F'$.
\epr

\bex4.41
Prove \E\Pr\ \rf{p4.40}.
\eex

\E\pp ies of proper subfields of the field~$\R$ will be investigated in the next section. We just mention that the \isc\ of subfields of~$\R$
containing the ideal $C_{n,p}$, $n\ge2$, $p$~prime, introduced in Lemma \rf{l4.5} is a well-defined subfield of~$\R$ by Lemma 4.4.28. We denote this
field by $\Q[C_{n,p}]$. Since $C_{n,p}$ is not a \pn\ ideal, $\Q[C_{n,p}]$ is not equal to~$\Q$. It will be shown in the next sections that
$\Q[C_{2,p}]$ is \ct y infinite and that $\R$ is un\ct e, hence $\Q[C_{2,p}]\ne\R$. We conclude this section by giving an exercise which reveals
an important \pp y of the field~$\R$. No proper subfields share this \pp y.

\bex4.42
Let $A$ be a \nss\ of the \tos\ $(\Q_{>0},\ge)$. If \te s $m\in\Q_{>0}$ \st $x\le m$ \fa $x\in A$, then the set $A$ is called \ti{bounded above}
and $m$~is called an \ub\ for~$A$. We denote the set of \ub s of~$A$ by $\UB(A)$. Clearly, $\UB(A)\ne\vn$ iff $A$ is bounded above. If the set~$A$
possesses a greatest \el~$\ov m$ then $\ov m\in\UB(A)$. Conversely, if $A\cap \UB(A)\ne\vn$, then $A$~possesses a greatest \el~$\ov m$ and
$A\cap \UB(A)=\{\ov m\}$. We now suppose that the set~$A$ is bounded above but does not have a greatest \el. Then \fe $x\in A$ and every $y\in\UB(A)$,
we have $x<y$.

\hph i,i, Show that if $\UB(A)$ has a least \el~$\ov m$, then $x<\ov m$ \fa $x\in A$, $i(x)<i(\ov m)$ \fa $x\in A$, $\Q_{>0}\sm\UB(A)$ is the \pn\
ideal $i(\ov m)$, and $i(\ov m)=\bcl_{x\in A}i(x)$.

\hph ii,, Suppose that $\UB(A)$ has no least \el. Show that $\Q_{>0}\sm\UB(A)$ is an ideal of~$\Q_{>0}$, hence an \el~$\xi$ of~$\R_{>0}$.
\eex

Suppose that $\O$ is a \nss\ of $(\R_{>0},\ge)$ which is bounded above. We will show in the next section that the set~$\O$ has a least \ub, also
called supremum, see \E\df s and Notations 3.1.13. \E\Ip if $\O$ \sf ies C1--C3 where $(\Q_{>0},\ge)$ is replaced by $(\R_{>0},\ge)$, that is,
$\O$~is an ideal of~$\R_{>0}$, then \te s $\xi\in\R_{>0}$ \st $\O=\{\eta\in \R_{>0}:\eta<\xi\}$. Thus  in $\R_{>0}$, every ideal is \pn.

\newpage
\Subsubsection{Order completeness and connectedness}\label{ass.5}

In the previous section we showed that the field~$\R$ as well as all of its subfields are \Ar\ \of s and that every \Ar\ \of\ is \is c to a subfield
of~$\R$. It turns out that \ti{no proper} subfield of~$\R$ is \is c to~$\R$. \Mo the \og\ of~$\R$ can be defined by means of its algebraic \sc\ and
the identity is the only auto\mf\ of~$\R$. We shall show that every \po\ \el\ of~$\R$ has a \po\ $n$-th root, $n\ge2$. \E\Ip every \po\ \ra\ \nm\
possesses a \po\ $n$-th root in~$\R$.

We first show that the \of~$\R$ is \tb{order-complete}.

\bth5.1
Every \nss\ of $\R$ that is bounded above has a least \ub.
\eth

An \of\ \sf ying this \pp y is usually called a \tb{complete \of}.\index{complete \of} \E\ev t \df s of order-\cp ness will be investigated. It turns out that a \cp\ \of\
is \Ar. The \fw\ lemma will be used in the proof of Theorem \rf{t5.1}. We recall that the least \ub\ for a \nss~$A$ of an \os\ $(E,\ge)$ is also
called the supremum of~$A$ and is denoted by $\sup A$ (see \E\df s and Notations 3.1.13).

\blm5.2
Let $(E,\ge)$, $(E',\ge')$ be \os s and let $\vf:E\to E'$ be an \ois sm $($see \E\df\ {\rm1.3.32)}. If $A$~is a \nss\ of~$E$ having a supremum, then
$\vf(A)$~is a \nss\ of~$E'$ having a supremum and
\beq5.1
\sup \vf(A) = \vf(\sup A).
\e
\elm

\proof
Let $x\in A$. Since $\vf$~is sur\jc, \fa $y\in A$ \te s $x\in A$ \st $y=\vf(x)$. Since $\sup A$ is an \ub\ for~$A$, we have $x\le\sup A$. Since
$\vf$~is in\cre, we have $\vf(x)\le' \vf(\sup A)$. \If that $\vf(\sup A)$ \ti{is an \ub\ for} $\vf(A)$. We now show that $\vf(\sup A)$ is the least
one. Let $y\in E'$ be an \ub\ for $\vf(A)$, and let $z\in E$ be \st $\vf(z)=y$. We then have $\vf(x)\le' y$, and $\vf(x)\le'\vf(z)$ \fa $x\in A$.
Since $\vf$~is bi\jc\ and $\vf\Inv$~is in\cre, we infer that $x\le z$ \fa $x\in A$. \If that $z$~is an \ub\ for~$A$, hence $\sup A\le z$ since
$\sup A$ is the least \ub\ for~$A$. \E\Tf $\vf(\sup A)\le' \vf(z)=y$. \csq, $\vf(\sup A)$ is \ti{the least \ub} for $\vf(A)$, which proves \er{5.1}.
\endproof

\bex5.3
Give an example of a \lt\ $(X,\lor,\land)$ (see \E\df\ 3.3.1) and of a strictly in\cre\ self-map $\vf$ of~$X$ \st \te\ $\ov x,\ov y\in X$,
$\ov x\ne\ov y$, \sf ying $\vf(\ov x\lor\ov y) \ne \vf(\ov x)\lor \vf(\ov y)$.
\eex

\proof[Proof of Theorem \rf{t5.1}]
Let $\cA$ denote a \nss\ of~$\R$ which is bounded above, and let $M\in\R$ denote an \ub\ for~$\cA$. We first consider the case where \te s $\ov\xi
\in\cA$ \st $\ov\xi>0$. (We use the symbols $0,1$ instead of $\wh0,\wh1$.) Then $0<\ov\xi \le M$, as well as $0<\ov\xi\le M'$ \fa \ub s~$M'$
for~$\cA$. Using the notation $\UB(\cA)$ for the set of all \ub s for~$\cA$, we have $\UB(\cA)\sbs \R_{>0}$. We claim $\UB(\cA) =\UB(\cA \cap
\R_{>0})$. Indeed, the inclusion $\subset$ follows from the fact that $\cA\cap\R_{>0}\sbs \cA$, and the converse inclusion $\supset$ follows from the
fact that if $\eta\in\cA \cap \R_{\le0}$, then $\eta\le 0<\ov \xi$. We set $\cA_{>0}:= \cA\cap\R_{>0}$ and observe that $\sup\cA$ exists iff $\sup
\cA_{>0}$ exists, in which case $\sup\cA = \sup\cA_{>0}$. In view of Theorem \rf{t4.30}, $\cA_{>0}$ is a \nss\ of ideals of~$\Q_{>0}$ which is
bounded above.

\bex5.4
Let $X$ be a \ns\ and let $\cP(X)$ be its power set. Let $\le$ denote the \og\ on $\cP(X)$ defined by $A\le B$, $A,B\in\cP(X)$ if $A\sbs B$ (possibly
\et y). Let $\cB$~be a nonempty family of subsets of~$X$. Show that
\beq5.2
\bigcup_{B\in\cB}B = \sup\cB \qh{in} (\cP(X),\ge).
\e
\eex

In view of \E\Pr\ \rf{p4.13}, $\bcl_{\a\in\cA_{>0}}\a$, the union of all ideals belonging to the family $\cA_{>0}$, is either $\Q_{>0}$ or
an ideal of~$\Q_{>0}$. Since $\R$ is \Ar\ by \E\Pr\ \rf{p4.31}, \te s $m\in\Q_{>0}$ \st $M<i(m)<M\hpl1$. \If that $\bcl_{\a\in\cA_{>0}}\a \le i(m)
\ne\Q_{>0}$. Thus $\bcl_{\a\in\cA_{>0}}\a \in \R_{>0}$. \Mo \fe $\b\in\cA_{>0}$, we have $\b\le\bcl_{\a\in\cA_{>0}}\a$, hence $\bcl_{\a\in\cA_{>0}}
\a$ is an \ub\ for $\cA_{>0}$. Finally, if $\b\in\R_{>0}$ is an \ub\ for~$\cA_{>0}$, then $\a\le\b$ \fa $\a\in\cA_{>0}$. \E\Tf $\bcl_{\a\in\cA_{>0}}\a
\sbs \b$, hence $\bcl_{\a\in\cA_{>0}}\a = \sup\cA_{>0} = \sup \cA$.

We now consider the case where $\cA\sbs \R_{\le0}$. Clearly, $0$~is an \ub\ for~$\cA$. By \as\ \te s $\ov\xi\in\cA$. Thus $\ov\xi\le 0$. Let $\th_\a$,
$\a\in\R$, denote the map defined in \er{3.6}. We have $\th_\a\circ\th_{-\a} = \id_\R = \th_{-\a}\circ\th_\a$ \fa $\a\in\R$, hence $\th_\a$ is a
bi\jn\ of~$\R$. By \er{3.9} $\th_a$~is an \ois sm \fa $\a\in\R$. Set $\a:=(-\ov\xi)+1$. Then $\th_\a \ov\xi = \ov\xi+(-\ov\xi)+1 = 1 > 0$. Hence
$1\in(\th_\a\cA)\cap\R_{>0}$. In view of the preceding part of the proof $\sup\th_\a\cA$ exists. Since $\cA= \th_{-\a}\circ(\th_\a\cA)$, $\sup\cA$
exists by Lemma~\rf{l5.2}. This completes the proof of Theorem \rf{t5.1}.
\endproof

We proved in Examples 4.5.59 that given a \Pn~$p$ and $n\in\Na\sms 1$, there is no \ra\ \nm\ $x\in(0,1)$ \sf ying $x^n=\frac1p$. Note that the \cn\
$x<1$ is necessary for ``\ex'' since if $x>1$, then $x^n > 1 > \frac1p$ \fa $n\ge2$. The situation is essentially different if we allow $x$ to be an
\el\ of~$\R_{>0}$. In that case, we reformulate the \eq\ $x^n=\frac1p$ as follows: Find $\xi\in\R_{>0}$ \st
\beq5.3
\xi^n = i(\tfrac1p),
\e
where $i(\frac1p):= \{y\in\Q_{>0}: y<\frac1p\}$, the \pn\ ideal of $\Q_{>0}$ generated by~$\frac1p$.

\E\eq\ \er{5.3} is a special case of the \eq
\beq5.4
\xi^n = \a, \q \a\in\R_{\ge0}, \q \xi\in\R_{\ge0}.
\e

\bth5.5
\E\fe $n\in\Na\sms1$ and every $\a\in\R_{\ge0}$, \te s \ooo $\xi\in\R_{\ge0}$ \st \er{5.4} holds.
\eth

\brm5.6 \

\hph i,ii, We recall that $\xi^n$ in \er{5.4} is the $n$-th \IT\ of~$\xi$ in the \am\ $\pz2{\R_{\ge0}} := (P,\hcd,i(1))$. We have $\xi^0:=1$,
$\xi^1:=\xi$ and $\xi^{n+1}:=\xi^n\cdot \xi$ \fa $n\in\N$.

\hph ii,i, If we denote by $f_n$ the self-map of $\R_{\ge0}$ defined by
\beq5.5
f_n(\xi):=\xi^n, \q n\in\Na,
\e
we find that $f_n$ is \ti{strictly in\cre} for $n\ge1$ in view of identity (4.4.46). \E\Tf \eq\ \er{5.4} possesses \ti{at most one} \so\ \fa $n\ge1$.

\hph iii,, Clearly $0^n=0$ and $1^n=1$ \fa $n\ge1$.
\erm

We next show that it suffices to consider $\xi$ in the \il\ of $\R_{>0}$ defined by
\beq5.6
(0,1):=\{\eta\in\R_{>0}: \eta<1\}.
\e
We still denote by $f_n$ the \rt ion of~$f_n$ to $(0,1)$. The \fw\ lemma allows us to \rt\ ourselves to the case $\a\in(0,1)$.

\blm5.7
Let $n\ge2$ and $\a\in(0,1)$. Suppose that \te s $\eta\in(0,1)$ \st $\eta^n=\a$. Then $\xi:=\eta\mo$, the inverse of~$\eta$ in the group
$\pz2{\R_{>0}}$, \sf ies
\beq5.7
\xi^n = \a\mo.
\e
\elm

Note that for $\a\in\R_{>0}$, we have $\a<1$ iff $1<\a\mo$. Indeed, $\a<1$ iff $1=\a\mo\a < \a\mo 1=\a\mo$, and $1<\a\mo$ iff $\a=\a1 < \a\a\mo=1$,
by \er{3.10} and Lemma \rf{l3.18}.

\proof[Proof of Lemma \rf{l5.7}]
It suffices to prove given $\eta\in\R_{>0}$:
\beq5.8
(\eta^n)\mo = (\eta\mo)^n \qh{\fa $n\ge1$.}
\e
We proceed by \In\ on $n\in\Na$, and we set $A:=\{n\in\Na: \hbox{\er{5.8} holds}\}$. Clearly $1\in A$. Suppose $n\in\Na$. Then we have
$(\eta^{n+1})\mo = (\eta^n\cdot\eta)\mo \nad{(4.3.5)}= \eta\mo \cdot(\eta^n)\mo \nad{n\in A}= \eta\mo\cdot (\eta\mo)^n = (\eta\mo)^{n+1}$. Hence
$n+1\in A$. Since $A$~is \iv, we have $A=\Na$.
\endproof

We now return to the proof of Theorem \rf{t5.5} with $\a\in(0,1)$. We give a proof based on the ``\ti{connectedness}'' of the \il\ $(0,1)$ defined in
\er{5.6}. We suppose, for \cd ion, that there is no $\xi\in(0,1)$ \st $\xi^n=\a$. We define two subsets $O_1$ and~$O_2$ of $(0,1)$ by setting
\bea5.9
O_1 &:= \{\xi\in(0,1): \xi^n < \a\},\\
O_2 &:= \{\xi\in(0,1): \xi^n > \a\}. \lb{5.10}
\e
Since the \og\ of $\R$ is total, $O_1$ and~$O_2$ are \ti{disjoint\/} ($O_1\cap O_2=\vn$) and $(0,1)=O_1\cup O_2$. We shall prove that both
$O_1$ and~$O_2$ are nonempty and that both $O_1$ and~$O_2$ are \ti{open} in the sense of the \fw\ \df.

\bdf5.8
A subset $U$ (possibly empty or equal to $(0,1)$) of the \il\ $(0,1)$ defined in \er{5.6} is called \ti{open}\index{open} if \fe $\d\in U$ \te\ $\b,\g\in(0,1)$
with $\b<\g$ \st $\d\in(\b,\g):=\{\rho\in(0,1): \b<\rho<\g\}$ and $(\b,\g)\sbs U$.
\edf

The \fw\ \pp y of the \il\ $(0,1)$ will lead to a \cd ion.

\blm5.9
The set $(0,1)$ defined in \er{5.6} is \emph{``connected''},\index{connected} that is, $(0,1)$ is not the disjoint union of two nonempty open subsets of $(0,1)$.
\elm

We first prove Lemma \rf{l5.9}, which is a con\sq\ of the order-\cp ness of the field~$\R$.

\proof[Proof of Lemma \rf{l5.9}]
Suppose, for \cd ion, that \te\ $O_1,O_2$ nonempty open subsets of~$(0,1)$ \st $O_1\cup O_2 = (0,1)$ and $O_1\cap O_2 = \vn$. Let $\xi_1\in O_1$ and
$\xi_2\in O_2$. Then $\xi_1\ne \xi_2$ since $O_1\cap O_2=\vn$. By possibly exchanging
the roles of $O_1$ and~$O_2$, we may assume that $\xi_1<\xi_2$. Note that $A:=O_1\cap[\xi_1,\xi_2]$ is nonempty since $\xi_1\in A$, and \ba\
by~$\xi_2$. \E\Tf $\sup A$ exists by Theorem \rf{t5.1}. We now show that both cases $\sup A\in O_1$ and $\sup A\in O_2$ lead to a~\cd ion. Indeed,
if $\sup A\in O_1$ then $\sup A\ge \xi_1$ since $\xi_1\in A$, and $\sup A\le \xi_2$ since $\xi_2$~is an \ub\ for~$A$ and $\sup A$~is the least \ub\
for~$A$. \If that $\sup A\in A$. However, since $\sup A\in O_1$ which is open and $\sup A<\xi_2\notin O_1$, \te s $\eta\in(0,\xi_2-\xi_1)$ \st
$[\xi_1,\xi_1+\eta) \sbs O_1\cap[\xi_1,\xi_2] = A$. A~\cd ion, since $\xi_1+\frac12\eta \in O_1\cap[\xi_1,\xi_2]$, and $\xi_1+\frac12\eta >\sup A$.

We now suppose that $\sup A\in O_2$. Since $O_2$ is open, \te s $\eta>0$ \st $\oz (\sup A)-\eta,\sup A \sbs O_2$. However, there must exist
$\a\in((\sup A)-\eta,\sup A)$ \st $\a\in O_1$. Otherwise, $\sup A$ would be less than or equal to $(\sup A)-\eta$, which is impossible. \E\Tf
$x\in O_1\cap O_2$, a~\cd ion since $O_1\cap O_2=\vn$.
\endproof

In the rest of the proof we use the notation $f$ instead of $f_n$. There is no risk of confusion with~$f$ defined in \er{4.26}.

In view of what precedes, the proof of Theorem \rf{t5.5} will be complete provided we show that both subsets $O_1$ and~$O_2$ defined in \er{5.9} and
\er{5.10} are nonempty and open. To this end, we first observe that Lemma 4.5.57 holds if we replace the \il\ $[0,1]:=\{x\in\Q_{\ge0}:\break
0\le x\le 1\}$ by the \il\ $[0,1]:=\{\eta\in\R_{\ge0}: 0\le \eta\le 1\}$. Indeed, the proof of this lemma uses identity (4.4.46) which holds in an
arbitrary (\cmt e) ring. The rest of the proof is verbatim the same as for the case of $\Q_{\ge0}$. Note that we already used (4.5.130) for the proof
of ``\uq''. We now use (4.5.131). Note that $M\in\Na$ in (4.5.131).

``\ti{$O_1$ and $O_2$ are not empty}'': Let $h\in(0,1)\sbs \R_{\ge0}$. Then $f(h)=f(h)-0 = f(h)-f(0)\le M(h-0)=Mh$. We want $Mh=\frac\a2<\a$. \E\Tf
$h:=M\mo \cdot\frac\a2 \le \frac\a2$ since $M\in\Na$ and $M\mo\le1$. Note that $h\in(0,1)$, $f(M\mo\cdot\frac\a2)=f(h)=Mh = \frac\a2< \a$. \E\Tf
$M\mo\cdot\frac\a2 \in O_1$. Similarly, $f(1)-f(1-h) \le M(1-(1-h))= Mh$. We want $f(1-h) = \frac12(\a+1)> \a$. We have $f(1-h)
\ge f(1)-Mh = 1-Mh$. If $1- Mh= \frac12(\a+1)$, then
$f(1-h)>\a$. \E\Tf we choose $h>0$ \st $Mh=1-\frac12(\a+1) = \frac12(1-\a)$, hence $h=M\mo \frac12(1-\a)$. Note that $1-\a$, $\frac12(1-\a)$ and
$M\mo\frac12(1-\a)$ belong to $(0,1)$, hence $1-M\mo\frac12(1-\a)\in(0,1)$. \E\Tf $f(1-M\mo\frac12(1-\a)) = f(1-h)>\a$, hence $1-M\mo\frac12(1-\a)
\in O_2$.

\ssk
``\ti{$O_1$ is open}'': Let $\eta\in(0,1)$ be \st $f(\eta)\in(0,\a)$. In view of (4.5.130) it suffices to find $\ov\si\in(\eta,1)$ \st $f(\ov\si)<\a$.
Indeed, in this case we have $0\h{<}\h f(\rho) \h{<}\h f(\eta) \h{<}\h f(\ov\si) \h{<}\h \a$ \fa $0<\rho<\eta<\ov\si$, hence $\eta\in(\rho,\ov\si)$
and \fa $\xi\in(\rho,\ov\si)$ we have $0< f(\xi)<\a$. We know that $f(\eta)<\a$ and $f(\si)-f(\eta) \le M(\si-\eta)$, thus $f(\si)\le
f(\eta)+M(\si-\eta)$. Since $\frac12(\a-f(\eta)) >0$, we have \fa $\o>0$ \sf ying
\beq5.11
\o \le M\mo \tfrac12(\a-f(\eta)),
\e
$M(\eta+\o-\eta) = M\o \le \frac12(\a-f(\eta))$. Setting $\si:= \eta+\o$, we find $f(\eta+\o) \le f(\eta)+\break \frac12(\a-f(\eta)) =
\frac12(f(\eta)+\a)
< \frac12(\a+\a) = \a$. Hence $f(\eta+\o)<\a$ provided $\o>0$ \sf ies \er{5.11}. Choosing moreover, if necessary, $\o\in(0,1-\eta)$ which is possible
since $\eta<1$, we find that $\ov\si:= \eta+\o$ \sf ies
\beq5.12
\ov\si \in (\eta,1) \qh{and } f(\ov\si)<\a.
\e
This completes the proof that $O_1$ is open.

\ssk
``\ti{$O_2$ is open}'': Inspection of the proof that $O_1$ is open shows that we could choose every \el\ of $(0,1)$ instead of~$\a$. \E\Ip we could
replace $\a$ by $1-\a$. \Mo we only used the \fw\ \pp ies of the \f\ $f:[0,1]\to[0,1]$:  $f(0)=0$, $f(1)=1$, (4.5.130) and (4.5.131). Note that if we
define $g: [0,1]\to\R_{\ge0}$ by setting $g(y):=1-f(1-y)$, $y\in[0,1]$, we find that $g$ is a self-map of $[0,1]$ \sf ying $g(0)=0$, $g(1)=1$,
(4.5.130) and (4.5.131). \E\Tf $A:=\{y\in[0,1]: g(y)<1-\a\}$ is open. We now show that $A=O_2$. We have $A=\{y\in[0,1]: 1-f(1-y)<1-\a\} = \{y\in[0,1]:
\break\a<f(1-y)\} \nad{x:=1-y}= \{x\in[0,1]: f(x)>\a\} = O_2$. This completes the proof that $O_2$ is open. \E\Tf the proof of Theorem \rf{t5.5} is
complete.
\endproof

\bex5.10
Show that the \f\ $g$ introduced in the proof that $O_2$ is open \sf ies (4.5.130) and (4.5.131).
\eex

We now show that a \cp\ \of\ is \Ar. To this end we will use the \fw\ lemma.

\blm5.11
Let $(F,+,\cdot,0,1,\ge)$ be an \of. Let $\d_c:F\to F$, $c>0$, denote the bi\jn\ of~$F$ introduced in \er{3.7} and let $\d_{-1}:F\to F$ denote the
bi\jn\ of~$F$ defined by
\beq5.13
\d_{-1}(x):= -x,\q x\in F.
\e
Let $A$ be a \nss\ of~$F$ that is bounded above \st $\sup A$ exists. Then so is $\d_c(A)$ and the \fw\ holds\dw
\beq5.14
\sup \d_c(A) = c \sup A.
\e
\Mo $\d_{-1}(A)$ is nonempty, bounded below, $\inf \d_{-1}(A)$ exists and
\beq5.15
\inf\d_{-1}(A) = -\sup A.
\e

Conversely, if $B$ is nonempty, bounded below and $\inf B$ exists, then $\d_{-1}(B)$ is non\-empty, bounded above, $\sup\d_{-1}(B)$ exists and
\beq5.16
\sup\d_{-1}(B) = -\inf B.
\e
\elm

We recall that the infimum of a set is defined in Notations and \E\df s 3.1.13.

\proof \

``\er{5.14}'': Let $A\sbs F$ be nonempty and let $M\in F$ be \st $x\le M$ \fa $x\in A$. Then $cx\le cM$ \fa $x\in A$ by \er{3.10}. \E\Ip if $\sup A$
exists, we have $x\le\sup A$ \fa $x\in A$, hence $cx\le c\sup A$ \fa $x\in A$, hence $c\sup A$ is an \ub\ for $\d_c(A)$. If $y\in F$ is an \ub\ for
$\d_c(A)$, then $cx\le y$ \fa $x\in A$, hence $x= c\mo cx\le c\mo y$ \fa $x\in A$ by \er{3.10} since $c\mo>0$ by Lemma \rf{l3.18}. \E\Tf $c\mo y$
is an \ub\ for~$A$, hence $\sup A\le c\mo y$, since $\sup A$ is the least \ub\ for~$A$. By \er{3.10} again we obtain $c\sup A\le y$. \If that
$c\sup A$ is the least \ub\ for $\d_c(A)$, hence \er{5.14} holds.

``\er{5.15}'': We use
\beq5.17
x\le y \hbox{ \ iff \ } {-y}\le -x \hbox{ \quad \fa}x,y\in F.
\e
We recall that \er{5.17} ``only if'' follows from \er{3.9} with $a:=(-y)+(-x)$ and \er{5.17} ``if'' from \er{3.9} with $a:=x+y$. Thus if $x\le M\in F$
\fa $x\in A$, then $-M\le \d_{-1}x$ \fa $x\in A$ and $\sup A\le M$. \If that $-M$ is a \lo\ for~$-A$ and $-\sup A$ is the greatest \lo\ for~$-A$ since
$-M\le -\sup A$.

The proof of ``\er{5.16}'' is similar and \Tf omitted.
\endproof

\bco5.12
Let $F$ be an \of. Then the \fw\ assertions are \ev t\dw

\hph i,ii, Every \nss\ of $F$ that is bounded above has a supremum.

\hph ii,i, Every \nss\ of $F$ that is bounded below has an infimum.

\hph iii,, Every \nss\ of $F$ that is bounded has a supremum and an infimum.

\hph iv,, $F$ is a \cp\ \of.
\eco

\Wanp show that for an \of\ the \pp y to be \cp\ is stronger than the \pp y to be \Ar.

\bpr5.13
Let $(F,+,\cdot,{\bf0},{\bf1},\ge)$ be a \cp\ \of. Then $F$ is \Ar.
\epr

\proof
In view of \E\df\ 4.5.53 it suffices to prove \pp y (ii) of Lemma 4.5.54 (see also the discussion preceding \er{3.36}). Suppose, for \cd ion, that
\te s $\ov y\in F_{>0}$ \st $\ov y\le(m\dpl {\bf1})\mo$ \fa $m\in\Na$. Set
\beag
A&:= \{y\in F: y=(m\dpl{\bf1})\mo \hbox{ \fa}m\in\Na\},\\
B&:= \{z\in F: z=((n+n)\dpl{\bf1})\mo \hbox{ \fa}n\in\Na\}.
\e
We have $B\sbs A$ since $n+n=m\in \Na$ \fa $n\in\Na$. Note that ${\bf0} < (m\dpl{\bf1})\mo$ \fa $m\in\Na$, since ${\bf1}\in F_{>{\bf0}}$,
$m\dpl1 \in F_{>{\bf0}}$ \fa $m\in\Na$ by \In\ on $m\in\Na$ and $(m\dpl1)\mo \in F_{>\bf0}$ by Lemma \rf{l3.18}. \E\Tf $\bf0$ is a \lo\ for~$A$, hence
also for~$B$. By Corollary \rf{c5.12} $\inf A$ and $\inf B$ exist. Since $\bf0$~is a \lo\ both for~$A$ and for~$B$, we have ${\bf0}\le \inf A$,
${\bf0}\le \inf B$. Since $B\sbs A$, if $x\le a$ \fa $a\in A$, then $x \le b$ \fa $b\in B$. Hence $\LB(A)\sbs \LB(B)$, and $\inf A\le \inf B$.
\E\Tf ${\bf0} \le \inf A\le \inf B$. Note that $B = \d_{({\bf1}+{\bf1})\mo}(A)$. This follows from
\beq5.18
({\bf1}+{\bf1})\mo \cdot (m\dpl{\bf1})\mo = ((m+m)\dpl{\bf1})\mo \qh{\fa}m\in\Na.
\e
We now prove \er{5.18}. Let $m\in\Na$. Then $({\bf1}\h{+}\h{\bf1})\mo\cdot(m\h{\dpl}\h{\bf1})\mo \nad{(4.3.5)} =
((m\h{\dpl}\h{\bf1})\cdot({\bf1}\h{+}\h{\bf1}))\mo\break
\nad{(4.1.2)} = (((m\dpl{\bf1})\cdot{\bf1}) + (m\dpl{\bf1})\cdot{\bf1})\mo = ((m\dpl{\bf1})+(m\dpl{\bf1}))\mo \nad{(2.2.3)\,\rm I2} =
((m+m)\dpl{\bf1})\mo$.

We now show that
\beq5.19
\inf \d_c(A) = c\inf A \qh{\fa $c\in F_{>0}$.}
\e
We first claim that
\beq5.20
\d_{-1}(\d_c(A)) = \d_c(\d_{-1}(A)) \qh{\fa $c\in F_{>0}$.}
\e
Indeed, $\d_{-1}(\d_c(A)) = \{x\in F: x=-cy,\ y\in A\} = \{x\in F: x=c\cdot (-y), \ y\in A\}$ since ${\bf0} \nad{(4.1.1)}= c\cdot{\bf0} = c\cdot
(y+(-y)) \nad{(4.1.2)}= c\cdot y + c\cdot(-y)$. Hence $c(-y) = -c\cdot y$, $y\in A$, which proves \er{5.20}.

We have $\inf \d_c(A) \nde5.16 = -\sup\d_{-1}(\d_c(A)) \nde5.20 = -\sup \d_c(\d_{-1}(A)) \nde5.14 = -c\sup(\d_{-1}(A)) \break \nde5.16 = -c(-\inf A)
\nde5.20 = c(-(-\inf A)) = c\inf A$, which proves \er{5.19}.

\csq, since $B= \d_{({\bf1}+{\bf1})\mo}(A)$, we obtain $\inf B = \inf\d_{({\bf1}+{\bf1})\mo}(A) =\break ({\bf1}+{\bf1})\mo \inf A$, \ev tly $\inf A =
({\bf1}+{\bf1})\cdot \inf B = \inf B+ \inf B$. \csq, we have
\beq5.21
{\bf0} \le \inf A \le \inf B \hbox{ \ and \ } \inf A = \inf B+\inf B.
\e
Thus $0\le \inf B+\inf B \le 0+\inf B$, hence $\inf B\le 0$ by \er{3.9}, which implies $\inf B=0$ and $\inf A=0$. A~\cd ion, since $(m\dpl{\bf1})\mo
\ge \ov y>0$ \fa $m\in\Na$ entails $\inf\limits_{m\in\Na}(m\dpl{\bf1})\mo\ge \ov y>0$. This completes the proof of the \Pr.
\endproof

\Wanp prove assertions mentioned at the beginning of this section.

\bpr5.14 \

\hph i,ii, No proper subfield of the field $\R$ is a \cp\ \of.

\hph ii,i, Every \cp\ \of\ is order- and ring-\is c to the field $(\R,\ge)$.

\hph iii,, The identity is the only ring-auto\mf\ of the field~$\R$.
\epr

\proof \

(i) Let $H$ denote an order-\cp\ subfield of the field $\R$ and let $x\in\R$. If $x=0$, then $x\in H$. Let $x\in\R_{>0}$. We claim that $x\in H$. Let
$A:=\{y\in\check\R: 0<y<x\}$, where $\check \R$ is the \pf\ of~$\R$ (that is the set of the form $r1$, where $r\in\Q$ and $1$~is the \nel\ of the
\mlv\ monoid $\pz2\R$). By \E\Pr\ \rf{p4.31}, $(\R,\ge)$ is \Ar. \csq, we have $x=\sup A$ in $(\R,\ge)$ by Corollary \rf{c4.2}, and \te s $r\in\Q$ \st
$x<r1<x+1$ by \E\Pr\ \rf{p4.1}. Hence the subset~$A$ of $(\check\R,\ge)$ is bounded above by~$r1$. Since the \pf\ of~$H$ is equal to $\check R$ by
\E\Pr\ 4.4.30\,(iii), we have $A\sbs H$ and $r1\in H$. Since $(H,\ge)$ is a \ti{\cp} \of\ by \as, \te s $\ov h\in H$ and \st $y\le\ov h$ \fa $y\in A$,
and $\ov h\le h$ \fa $h\in H$, $h$~\ub\ for~$A$. Since $\ov h\in H\sbs \R$, $\ov h$~is an \ub\ for~$A$ belonging to~$\R$, hence $x\le\ov h$, $x$~being
the least \ub\ for~$A$ in~$\R$. Suppose,
for \cd ion, that $\ov h<x$. Since $\R$~is \Ar\ and both $\ov h$ and~$x$ belong to~$\R$, \te s $s1\in \check K$ ($s\in\Q_{>0}$) \st $\ov h< s1 <x$.
\E\Tf $s1\in A$ by \df\ of~$A$ and $s1>\ov h$. However, $y\le \ov h$ \fa $y\in A$. A~\cd ion. \csq, $x=\ov h\in H$. We conclude the proof by noting
that if $x<0$ then $-x\in\R_{>0}$, hence $-x\in H$ by what precedes, and $x = -(-x)\in H$.

(ii) Let $(F,\ge)$ be a \cp\ \of. Then $(F,\ge)$ is \Ar\ by \E\Pr\ \rf{p4.31}, and by Theorem \rf{t4.35} $(F,\ge)$ is order- and ring-\is c to a
subfield~$H$ of~$\R$. By Lemma \rf{l5.2} the field $(H,\ge)$ is \cp, hence $H= \R$ by part~(i).

(iii) Let $\Phi:\R\to \R$ be a ring-auto\mf. We first show that $\Phi$~maps $\R_{>0}$ into itself. Let $x\in\R_{>0}$. By Theorem \rf{t5.5} \te s
$y\in\R_{>0}$ \st $x=y^2$. Then $\Phi(x) = \Phi(y^2) = \Phi(y\cdot y) = \Phi(y)\cdot \Phi(y)\in\R_{>0}$ by~\er{3.12}. By Lemma \rf{l3.28} $\Phi$~is
an \ois sm of $(\R,\ge)$, and by Lemma 4.4.31\,(ii) $\Phi$~maps~$\check \R$ into itself and is a ring-auto\mf\ of~$\check \R$. Since $\check \R$ is
ring-\is c to~$\Q$ and the identity is the only ring-auto\mf\ of~$\Q$ (see the discussion at the beginning of Section~\ref{ass.4}), the identity is
the only ring-auto\mf\ of~$\check \R$. Let $x\in\R_{>0}$. Then $x = \sup\{y\in\check \R: 0 < y < x\}$ by~\er{4.6}, hence $\Phi(x) = \Phi(\sup
\{y\in\check \R: 0<y<x\}) \nde5.1 = \sup\Phi(\{y\in\check \R: 0<y<x\}) = \sup\{y\in\check \R: 0<y<x\} \nde4.6 = x$. If $x<0$ then $-x>0$, and
$\Phi(x) = \Phi(-(-x)) = -\Phi(-x) = -(-x)= x$, which completes the proof of (iii) and of \E\Pr\ \rf{p5.14}.
\endproof

Our next goal is to give an example of an \ti{\Ar} field \st the identity is \ti{not\/} the only ring-auto\mf. We begin with an exercise.

\bex5.15
Let $\pz0F$ be a field, and let $H$ be a subset of~$F$ \sf ying:
\bea5.22
{}&0,1\in H.\\
&\hbox{\E\fa} x,y\in H,\ x+y \hbox{ and } x\cdot y\in H. \lb{5.23}\\
&\hbox{\E\fa} x\in H,\ \hbox{\te s } y\in H \hbox{ \st}x+y=0. \lb{5.24}\\
&\hbox{\E\fa} x\in H\sms0,\ \hbox{\te s } y\in H \hbox{ \st}x\cdot y=1. \lb{5.25}
\e
Show that $H$ is a \ti{subfield\/} of $F$.
\eex

\bxa5.16
We use the notation $(\R,+,\cdot,{\bf0},{\bf1})$ for the field introduced in Theorem \rf{t4.30}. As an \of, the field~$\R$ has \ch istic~$0$
by (4.5.126). By \E\Pr\ \rf{p4.31a}, the field $\R$ can be viewed as a \vsf0$\Q$ where the \mlc\ by scalars is defined by \er{4.73}. Let $p$~be
a \Pn, and let $\sqrt p$ denote the \po\ \sqr\ of~$p$, that is, the unique $\xi\in\R_{>0}$ \sf ying $\xi\cdot\xi=p$. Such an \el\ is introduced in
Theorem \rf{t5.5}. We define a subset $H$ of the \vs\ $(\Q,\R)$ by setting
\beq5.26
H:= \{a{\bf1} + b\sqrt p: a,b\in\Q\}
\e
where $a{\bf1}$ and $b\sqrt p$ are defined by \er{4.73}. We use the notation $\pz0\Q$ for the field of \ra\ \nm s, \Ip the same symbol~$+$ for the
\ad\ in $\Q$, $\R$ and $(\Q,\R)$. We omit the dot for the \mlc\ of two \ra\ \nm s, and for \el s of $\Q\sms{0,1}$ we use small letters $a,b,c,\dots
$. We observe that ${\bf0}\in H$ since ${\bf0} = {\bf0}+{\bf0} = 0{\bf1} + 0\sqrt p$. Similarly, ${\bf1}\in H$ since ${\bf1} = {\bf1}+{\bf0} = 1
{\bf1}+0\sqrt p$. \Mo if $z:=a{\bf1} + b\sqrt p$ and $z':=a'{\bf1}+b'\sqrt p$, then $z+z' = (a{\bf1}+ b\sqrt p) + (a'{\bf1}+b'\sqrt p)
\nad{(2.1.36)} = a{\bf1} + (b\sqrt p+ a'{\bf1})+ b'\sqrt p = a{\bf1} + (a'{\bf1}+b\sqrt p) + b'\sqrt p \nad{(2.1.36)} = (a{\bf1}+a'{\bf1}) +
(b\sqrt p+b'\sqrt p) \nde4.73 = (j(a)\cdot{\bf1}+j(a')\cdot{\bf1}) + (j(b)\cdot\sqrt p + j(b')\cdot\sqrt p) \nad{(4.1.3)}= (j(a)+j(a'))\cdot{\bf1}
+ (j(b)+j(b'))\cdot\sqrt p = j(a+a')\cdot{\bf1} + j(b+b')\sqrt p$, since $j:\Q\to \check R$ is a ring-\hm sm.

\E\Tf $z+z'\in H$. Similarly, we have $cz = j(c)\cdot z = j(c)\cdot(a{\bf1}+b\sqrt p)\nad{(4.1.2),\er{4.73}} = j(c)\cdot (j(a)\cdot{\bf1}) + j(c)\cdot
(j(b)\cdot\sqrt p) = (j(c)\cdot j(a))\cdot{\bf1} + (j(c)\cdot j(b))\cdot\sqrt p = j(ca)\cdot{\bf1} + j(cb)\cdot\sqrt p$. Hence $cz\in H$. In view of
\E\df\ 5.2.6, $H$~is a \lss\ of $(\Q,\R)$, hence a \vs\ over~$\Q$ by \E\Pr\ 5.2.7. Note that \er{5.24} holds since $(a{\bf1}+b\sqrt p) + ((-a){\bf1}
+(-b)\sqrt p) = (a+(-a)){\bf1} + (b+(-b))\sqrt p = 0{\bf1} + 0\sqrt p = {\bf0}+{\bf0} = {\bf0}$. Since $H=\spn\{{\bf1},\sqrt p\}$ in the \vs\
$(\Q,H)$ (see \E\df\ 5.1.5\,(iii)), the pair is a basis of $(\Q,H)$ iff ${\bf1}$ and $\sqrt p$ are \li\ iff $\sqrt p\notin\spn\{{\bf1}\}$ (see \E\df\
5.1.24). We claim that there is no $a\in\Q$ \st $\sqrt p = a{\bf1}$. From Example 4.5.54 we know that there is no $a\in\Q$ \st $a^2=p$. Proceeding as
in Lemma \rf{l5.7} we find:
\beq5.27
\hbox{there is no }b\in\Q_{>0} \hbox{ \st} b^2=p.
\e
Suppose, for \cd ion, that $\sqrt p = r{\bf1}$ \fs $r\in\Q$. Then $i(p) = \sqrt p\cdot\sqrt p = r{\bf1}\cdot r{\bf1} = (j(r)\cdot{\bf1})\cdot(j(r)
\cdot{\bf1}) = j(r)\cdot j(r)\cdot{\bf1} = j(r^2)\cdot{\bf1} = i(r^2)$. Hence $p=r^2=|r|^2$, a \cd ion.

\If that $\{{\bf1},\sqrt p\}$ is a basis of $(\Q,H)$, hence $H$ is a \ti{two-\dm al\/} \vs\ over~$\Q$. We next prove that $z\cdot z'\in H$ \fa
$z,z'\in H$. We have $z\cdot z' = (a{\bf1} + b\sqrt p)\cdot\break (a'{\bf1}+b'\sqrt p) \nad{(4.1.2),(4.1.3)}= aa'{\bf1}\cdot{\bf1} + ab'{\bf1}\cdot
\sqrt p + ba'\sqrt p\cdot{\bf1} + bb'\sqrt p\cdot\sqrt p = (aa'+pbb'){\bf1} + (ab'+ba')\sqrt p \in H$.

Observe that we proved \er{5.22}--\er{5.24}. It remains to prove \er{5.25}. Given $z\in H\sms0$ we want to find $z'\in H$ \st $z\cdot z'={\bf1}$. We
have $(a,b)\ne(0,0)\in\Q\t\Q$, and we need $aa'+p\cdot bb'=1$ and $ab'+ba'=0$. Setting $x:=a'$, $y:=b'$ we find a system of \eq s \er{4.6} with
$a:=a$, $b:=pb$, $c:=b$, $d:=a$, $f:=1$, $g:=0$. We verify \cn\ (4.4.7), namely
\beq5.28
a^2 \ne pb^2.
\e
We consider four cases: (i)~$a\ne0$, $b\ne0$; (ii)~$a\ne0$, $b=0$; (iii)~$a=0$, $b\ne0$; (iv)~$a=b=0$. Case~(iv) is excluded by the \cn\ $z\ne{\bf0}$.
Case (iii): $a^2 = 0^2 = 0 \ne pb^2$ since $b\ne0$ implies $b^2\ne0$ and $pb^2\ne0$. Hence \er{5.28} holds. Case~(ii): $b=0$ implies $pb^2=0$ and
$a\ne0$ implies $a^2\ne0$. Hence \er{5.28} holds. Case~(i): we have $b\ne0$, hence $b^2\ne0$, $a=1$, and $a^2(b^2)\mo = (ab\mo)(ab\mo)\in\Q_{>0}$.
\E\Tf $(ab\mo)^2\ne p$, otherwise $p$~would have a \po\ \ra\ \sqr, which is impossible by \er{5.27}. \If that all \cn s of Exercise \rf{ex5.15} are
\sf ied, \Tf $H$~is a subfield of the \Ar\ field $(\R,\ge)$, hence $(H,\ge)$ is \Ar\ in view of \cn\ 4.5.54\,(i). We finally show that the field~$H$
possesses a ring-auto\mf\ different from the identity. Let $\Phi$ be a ring-auto\mf\ of~$H$. By Lemma 4.4.31\,(ii) $\Phi$~maps the \pf~$\check H$
of~$H$ into itself. Since $\check H$ is ring-\is c to the field~$\Q$, which has only the identity as ring-\is sm, we have $\Phi(a{\bf1}) = a1$ \fa
$a\in\Q$. \Mo $\Phi(\sqrt p)$ \sf ies the \eq\ $\Phi(\sqrt p)^2=p{\bf1}$. Indeed, $\Phi(\sqrt p)\cdot\Phi(\sqrt p) = \Phi(\sqrt p\cdot\sqrt p) =
\Phi(p{\bf1}) = p{\bf1}$. There is \ooo $\xi\in\R_{>0}$ \sf ying $\xi^2=p{\bf1}$, namely $\xi:=\sqrt p$. \E\oh $-\sqrt p$ \sf ies $(-\sqrt p)^2 =
(-\sqrt p)\cdot(-\sqrt p) = \sqrt p\cdot\sqrt p = p{\bf1}$.

We claim that $\sqrt p$ and $-\sqrt p$ are the only roots of the \pl\ map $\xi\mt \xi^2-p{\bf1}$ from~$\R$ into~$\R$. Let $\ov\xi\in\R\sms0$ \sf y
$\ov\xi{}^2-p{\bf1} = {\bf0}$, and let $\ov \rho\in\R$ be defined by $\ov\rho:= \ov\xi\cdot(\sqrt p)\mo$. We have $\ov\xi = \ov\rho\cdot\sqrt p$,
hence $\ov\xi{}^2=\ov\rho{}^2\cdot p{\bf1}$ and ${\bf1}\cdot p{\bf1} = \ov\rho{}^2\cdot p{\bf1}$. By \cnc ity ($p{\bf1}\ne0$), we obtain $\ov\rho{}^2
={\bf1}$. \E\Tf $\ov\rho{}^2$ is a \ru{} in~$\R$, hence ${\bf1}$ and $-\bf1$ are the only roots of $\ov\rho{}^2-\bf1$ by \E\Pr\ 4.4.42. \If that
$\ov\xi\in\R$ \sf ies $\ov\xi{}^2 = p\bf1$ iff $\ov\xi=\sqrt p$ or $\ov\xi=-\sqrt p$. \If that $\Phi(\sqrt p)$ is either equal to~$\sqrt p$ in which
case $\Phi$~is the identity or $\Phi(\sqrt p)=-\sqrt p$. We finally show that the map $\Phi:H \to H$ defined by
\beq5.29
\Phi(a{\bf1} + b\sqrt p) = a{\bf1} - b\sqrt p, \q a,b\in\Q,
\e
is an \ti{auto\mf} of $H$ (different from id).

First observe that the pair $({\bf1},-\sqrt p)$ is a basis of~$H$ by \E\Pr\ 5.1.30 since $(\Q,H)$ is two-\dm al and $\spn({\bf1},-\sqrt p)=H$.
One easily verifies that $\Phi({\bf0})=\bf0$, $\Phi({\bf1})={\bf1}$, $\Phi$~is a \lop\ of $(\Q,H)$, the \rt ion of~$\Phi$ to~$\check \R$ is
$\id_{\check \R}$. It remains to show that $\Phi(z\cdot z') = \Phi(z)\cdot \Phi(z')$ and $\Phi$~is bi\jc.

``$\Phi(z)\cdot \Phi(z')=\Phi(z\cdot z')$'': $\Phi(a{\bf1}+b\sqrt p)\cdot\Phi(a'{\bf1} + b'\sqrt p) = (a{\bf1}-b\sqrt p)\cdot(a'{\bf1}-b'\sqrt p)
= (aa'+pbb'){\bf1} - (ab'+a'b)\sqrt p$. \E\oh $z\cdot z' = (a{\bf1}+b\sqrt p)(a'{\bf1}+b'\sqrt p) = (aa'+bb'){\bf1} + (ab'+a'b)\sqrt p$. Thus
$\Phi(z)\cdot\Phi(z')= \Phi(z\cdot z')$.

``\ti{$\Phi$ is bi\jc}'': Since $(\Q,H)$ is finite-\dm al, it suffices to show that $\Phi$ is in\jc, hence in view of the linearity of~$\Phi$ that the
null-space of~$\Phi$ is~$\{{\bf0}\}$. We have $\Phi(z)=\bf0$ iff $a{\bf1}+b\sqrt p = \bf0$ iff $a=0$, $b=0$ iff $z:=a{\bf1}+b\sqrt p=\bf0$ since
$({\bf1},\sqrt p)$ are \li.

We conclude this example by some remarks.

The \nm\ of auto\mf s of $H$ is equal to the \dm\ of~$H$ as a \vs\ over~$\Q$. The auto\mf\ $\Phi$ maps the root $\sqrt p$ into the root $-\sqrt p$.
If we denote by~$H_p$ (resp.\ $H_{p'}$, $p'$~prime, $p'\ne p$) the \crs\ \vs s spanned by $({\bf1},\sqrt p)$ (resp.\ $({\bf1},\sqrt{p'})$) in the \vs\
$(\Q,\R)$, then $H_p$ and $H_{p'}$ are $\Q$-linear \is sms. Observe that both $H_p$ and~$H_{p'}$ are \ti{splitting} fields for the \fp s $X^2-p$
(resp.\ $X^2-p'$) belonging to $\Q[X]$.
\exa

\bex5.17
Let $\Phi$ denote the auto\mf\ of the subfield~$H$ defined by \er{5.29}. Let $H_+:=\{\eta\in H: \eta>0\}$ and $H_-:=\{\eta\in H: \eta<0\}$. Set
$\wt H_+:=\Phi(H_+)$ and $\wt H_-:=\Phi(H_-)$. Show that

\hph i,ii, $H$ is the disjoint union of $\wt H_+$, $\{0\}$ and $\wt H_-$, and $\wt H_-=-\wt H_+$.

\hph ii,i, $\wt H_+ + \wt H_+ \sbs \wt H_+$ and $\wt H_+ \cdot \wt H_+\sbs \wt H_+$.

\hph iii,, If $\ge'$ denotes the \og\ of~$H$ defined by $\xi\ge'\eta$ if $\xi-\eta \in \wt H_+$ \sf ies $\xi\ge'\eta$ iff $\Phi\Inv(\xi) \ge
\Phi\Inv(\eta)$, $\xi,\eta\in H$.

\hph iv,, $(H,\ge')$ is an \of\ \ti{distinct\/} from $(H,\ge)$.
\eex
\goodbreak

\brm5.18
We use notations \er{5.5}, \er{5.6} introduced in Remark \rf{r5.6}. Let $n\in\Na$.

\hph i,i, \If from the strict in\cre ness of~$f_n$, $f_n(0)=0$, $f_n(1)=1$ and Theorem \rf{t5.5} (or from the proof of Theorem \rf{t5.5}) that
$f_n$ is a bi\jn\ of~$(0,1)$. We denote its inverse by $g_n$, thus we have
\beq5.30
g_n \circ f_n = \id_{(0,1)}, \q f_n\circ g_n = \id_{(0,1)}.
\e
Let $h_1:(0,1)\to(0,1)$ be the map defined by
\beq5.31
h_1(y):=1-y, \q y\in(0,1).
\e
We have
\beq5.32
h_1\circ h_1 = \id_{(0,1)},
\e
hence $h_1$ is a bi\jn\ (involution) of $(0,1)$.

We next consider the \cm\ of the maps $f_n$, $h_1$ and $g_n$ by setting:
\beq5.33
h_n:= g_n\circ h_1\circ f_n, \q n\in\Na\sms1.
\e
Note that $g_1 \circ h_1 \circ f_1 = h_1$, since $f_1=\id_{(0,1)}$. As a \cm\ of bi\jn s, $h_n$~is a bi\jn\ of $(0,1)$. \Mo since $f_n,g_n$ are
strictly in\cre\ and $h_1$~is strictly de\cre, the bi\jn\ $h_n$ is strictly de\cre. It turns out that \fa $n\ge3$ and all $x\in(0,1)\cap\Q_{>0}$,
$h_n(x)\notin\Q_{>0}$. Indeed, if $y_n:=h_n(x)$ for $x\in(0,1)$ would belong to $(0,1)\cap\Q_{>0}$, we would have $y_n^n = f_n\circ(g_n\circ h_1\circ
f_n) = (f_n\circ g_n)\circ h_1\circ f_n = 1-x^n$, that is, $x^n+y_n^n=1$. Since \te\ $a,b,c,d\in\Na$ \st $x=\frac ab$, $y_n=\frac cd$, we would have
$(\frac ab)^n+ (\frac cd)^n = \frac11$. Using \In\ on $n\in\Na$, one finds that $(\frac ab)^n = \frac{a^n}{b^n}$, $(\frac cd)^n= \frac{c^n}{d^n}$.
\If that $a^nd^n + b^nc^n = b^nd^n$. Setting $k:=ad$, $l:=bc$ and $m:=bd$, we find $k,l,m\in\Na$ \sf ying $k^n+l^n = m^n$. As a con\sq\ of the famous
highly nontrivial Fermat--Wiles theorem (see~\cite[p.~474]{Is}), we find $n\le2$. A~\cd ion.

Note that $(0,1)_\Q := \{x\in\Q: 0<x<1\}$ is \ct y infinite. Indeed, we have
\beq5.34
A \sbs (0,1)_\Q \sbs \Q
\e
where $A:=\{\frac1n\in\Q_{>0}:n\in\Na\}$. Since $A$~is \ep\ to~$\N$ as well as $\Q$ by Theorem 4.5.44, by \E\Pr\ 1.4.27\,(A)(i) and~(iv) we find that
$(0,1)_\Q$ is \ct y infinite. In view of the same \Pr\ we find that $h_n((0,1)_\Q)$ is \ct y infinite.
It will be shown in the next section that $(0,1)$ is \ti{un\ct e}. Note that $(0,1)_\Q$ and $h_n((0,1)_\Q)$, $n\ge3$,
are two disjoint \ct y infinite subsets of~$(0,1)$. We claim that the union $(0,1)_\Q \cup h_n((0,1)_\Q)$ is also \ct y infinite. Indeed, let $\vf_1:
\N \to (0,1)_\Q$ and $\vf_2:\N \to h_n((0,1)_\Q)$ be bi\jn s. Define $\vf:\N \to (0,1)_\Q \cup h_n((0,1)_\Q)$ by setting $\vf(2m):=\vf_1(m)$,
$\vf(2m+1):=\vf_2(m)$, $m\in\N$. One verifies that $\vf$~is a bi\jn. Since an un\ct e set is not \ep\ to a \ct y infinite set by \df\ (see \E\df\
1.4.21), we have $(0,1)\ne(0,1)_\Q \cup h_n((0,1)_\Q)$. We know that $(0,1)_\Q$ is order-dense in $(0,1)$ since $(\R,\ge)$ is \Ar. A~natural question
to ask is whether $h_n((0,1)_\Q)$ is order-dense in $(0,1)$, that is, given $\xi,\eta\in(0,1)$, $\xi<\eta$, does \te s $r\in(0,1)_\Q$ \st $\xi<
h_n(r{\bf1})<\eta$\,? Let $\a_n,\b_n\in(0,1)$ be \st $\a_n:=h_n\Inv(\xi)$, $\b_n:=h_n\Inv(\eta)$. We have $\b_n<\a_n$, $\b_n,\a_n\in(0,1)$. Since
$(\R,\ge)$ is \Ar, \te s $r\in\Q_{>0}$ \st $\b_n < r{\bf1} < \a_n$. Hence $\xi = h_n(\a_n) < h_n(r{\bf1}) < h_n(\b_n) = \eta$. \E\Tf $h_n((0,1)_\Q)$
is \ti{order-dense} in $(0,1)$.

\hph ii,, The map $g_n:(0,1)\to(0,1)$, $n\ge2$, is strictly in\cre\ (analogue of (4.5.130)) by Lemma 1.3.34, but does not \sf y the analogue of
(4.5.131). Suppose, for \cd ion, that \te s $M\in\Na$ \st $g_n(y)-g_n(x)\le M(y-x)$ \fa $0\le x\le 1$, $x\in\R$, then setting $x:=0$ we find $g_n(y)
\le My$, hence $y\le M^ny^n$ \fa $y\in(0,1)$. \E\ml ying both sides by $y\mo(M^n)\mo$ we find $(M^n)\mo \le y^{n-1}$. Since $n\ge2$ and $y\le1$, we
have $y^{n-1}\le y$, hence $0<(M^n)\mo\le y$, $y\in(0,1)$. Choosing $y:=((M+1)^n)\mo$ we obtain a \cd ion.
\erm

We consider \pp ies of the map $g_n:\R_{\ge0} \to \R_{\ge0}$:
\beq 5.35a
g_n:= f_n\Inv \q \hbox{where $f_n$ is defined in \er{5.5}.}
\e

\blm5.19
The map $g_n$ defined above \sf ies \fa $x,y\in\R_{\ge0}$\dw
\bea5.35
g_n(x) &< g_n(y) \hbox{ if } x<y,\\
g_n(x\cdot y) &= g_n(x)\cdot g_n(y), \lb{5.36}\\
g_n(x+y) &\le g_n(x)+g_n(y). \lb{5.37}
\e
\elm

\proof \

``\er{5.35}'': Follows from Lemma 1.3.34.

``\er{5.36}'': From \E\Pr\ 2.2.3\,I5 with $M:=\R_{\ge0}$, ${\qu}:={\cdot}$ and $e:=1$, we find\break $(g_n(x)\cdot g_n(y))^n = (g_n(x))^n\cdot
(g_n(y))^n$. Hence $g_n(x)\cdot g_n(y) = g_n\bigl((g_n(x)\cdot g_n(y))^n\bigr) = g_n\bigl((g_n(x))^n\cdot (g_n(y))^n\bigr) = g_n(x\cdot y)$.

``\er{5.37}'': The case $x:=0$ or $y:=0$ is trivial. We assume $x>0$ and $y>0$. Then $g_n(x)>0$ and $g_n(y)>0$ by \er{5.35}. Set
$$
t:= \frac{g_n(x)}{g_n(x)+g_n(y)},
$$
hence
$$
1-t:= \frac{g_n(y)}{g_n(x)+g_n(y)}.
$$
We have $0<t,1-t<1$, hence $0<t^n<t$ and $0<(1-t)^n<1-t$, by using \In. \E\Tf $t^n+(1-t)^n < t+(1-t)=1$, and $g_n(t^n+(1-t)^n) < g(1)=1$ by \er{5.35}.
We obtain
$$
g_n\biggl(\biggl(\frac{g_n(x)}{g_n(x)+g_n(y)}\biggr)^n + \biggl(\frac{g_n(y)}{g_n(x)+g_n(y)}\biggr)^n\biggr)\le1.
$$
\E\ml ying both sides by $g_n(x)+g_n(y)>0$, and using $g_n(x)+g_n(y) = g_n\bigl((g_n(x)+g_n(y))^n\bigr)$, we obtain
\bmlg
g_n\biggl((g_n(x)+g_n(y))^n \cdot \biggl(\frac{g_n(x)}{g_n(x)+g_n(y)}\biggr)^n + (g_n(x)+g_n(y))^n \cdot \biggl(\frac{g_n(y)}{g_n(x)+g_n(y)}\biggr)^n
\biggr) \\ \le g_n(x)+g_n(y).
\e
Note that $(g_n(x)+g_n(y))^n \cdot \bigl(\frac{g_n(x)}{g_n(x)+g_n(y)}\bigr)^n = (g_n(x))^n = x$, similarly $(g_n(x)+g_n(y))^n \cdot
\bigl(\frac{g_n(y)}{g_n(x)+g_n(y)}\bigr)^n = (g_n(y))^n = y$. Thus $g_n(x+y)\le g_n(x)+g_n(y)$.
\endproof

Using the more ``appealing'' usual notations
\bea5.38
\root n \of x &:= g_n(x), \q n\in\Na\sms1, \qquad \root 1 \of x := x, & \hskip-50pt x\in\R_{\ge0}, \\
\sqrt x &:= g_2(x), & \hskip-50pt x\in\R_{\ge0}, \lb{5.39}
\e
we have the \fw\ analogue of (4.5.131).

\bco5.20
Let $0<x<y$, $x,y\in\R$. Then
\beq5.40
\root n \of y - \root n \of x \le \root n \of {y-x}  \qh{for $n\ge2$.}
\e
\eco

\proof
Let $z\in\R_{>0}$ be \st $y=x+z$. Then $\root n \of y = \root n \of {x+z} \nde5.37 \le \root n \of x + \root n \of {y-x}$. Hence $\root n \of y
+ (-\root n \of x) \le \root n \of x + \root n \of {y-x} + (-\root n \of x) = \root n \of {y-x}$.
\endproof

\brm5.21
Formula (4.5.131) (resp.\ \er{5.40}) expresses an important \pp y of the in\cre\ maps $f$ (resp.\ $\root m \of{\phantom{x}}$). In order to remove
the \cn\ of in\cre ness, it is useful to introduce the notion of \av\ of an \el\ of $(\R,\ge)$.
\erm

\bdf5.22
Let $(F,\ge)$ be an \of. Given $x\in F$ we set
\beq5.41
|x| := \bca
x &\hbox{if } x\in F_{\ge0},\\
-x & \hbox{if }x\in F_{<0}.
\eca
\e
Then $|x|$ is called the \it{\av} of $x$.
\edf

It turns out that both the map $x\mt|x|$ in $(F,\ge)$ and the map $x\mt \root m \of {|x|}$ in $(\R,\ge)$ are valuations.

\begin{dfn}[{\cite[p.~191]{A2}}] \lb{d5.23}
Let $K$ be a field and let $(F,\ge)$ be an \of. A~map $\vf :K \to F$ is called an \tb{$F$-valued valuation}\index{F-valued valuatuion@$F$-valued valuation} of the field~$K$, if \fa $x,y\in K$
we have
\bea5.42
{}&\vf(x)\ge0 \hbox{ \ and \ }\vf(y)=0 \hbox{ iff }y=0,\\
&\vf(x\cdot y)= \vf(x)\cdot \vf(y), \lb{5.43} \\
&\vf(x+y) \le \vf(x)+\vf(y), \lb{5.44} \\
&\hbox{\te s $z\in K\sms0$ \st}\vf(z)\ne1. \lb{5.45}
\e
\edf

\bex5.24 \

\hph i,ii, Let $\vf:K \to F$ be a valuation. Prove:
\bea5.46
{}& \vf(1)=1, \q \vf(-1)=1, \q \vf(x)=\vf(-x) \hbox{ \fa} x\in K.\\
& |\vf(x)-\vf(y)| \le \vf(x-y) \hbox{ \fa} x,y\in K. \lb{5.47}
\e

\hph ii,i, Let $K$ be an \of. Set $\vf(x):= |x|$, $x\in K$. Show that $\vf$ is a $K$-valued valuation of~$K$.

\hph iii,, Let $K:=\Q$ or $\R$ and let $\vf(x):=\root m \of{|x|}$, $x\in K$, $m\ge2$. Show that $\vf$ is an $\R$-valued valuation of~$K$.

\hph iv,, Let $K:=\R$ and let $m\in\Na$. Then $\root 2m \of{x^{2m}} = |x|$ \fa $x\in\R$.
\eex

\brm5.25
In a group $(G,\qu,e)$ with \mlv\ notation ${\qu}:={\cdot}$, the $n$-th \IT\ of an \el\ $x\in G$ is usually called the $n$-th \ti{power} of~$x$
(see for example \cite[p.~17]{Alg}) and is denoted by~$x^n$.

\ti{\E\ng\ powers} are defined by $x^{-n}:= (a\mo)^n$, $n\in\Na\sms1$, where $a\mo$ denotes the inverse of~$x$ in $\pz2G$ (see \cite[p.~18]{Alg}).
If $(G,\qu,e):=\pz2{\R_{>0}}$, $\root n \of x$ may be called a \ti{fractional power} of~$x$, and is denoted by $x^\frac1n$. In what follows we define
\po\ \ra\ powers of~$x$. For example $x^\frac2n$ could be defined by $(x^2)^\frac1n$. However, if we want to use the notation $x^\frac ab$, $a,b\in
\Na$, we need to show that $x^\frac ab= x^\frac{a'}{b'}$ whenever $ab'=ba'$.
\erm

This motivates the \fw\ lemma.

\blm5.26
Let $x\in\R_{>0}$ and let $a,b,a',b'\in\Na$. Let $\xi,\xi'\in\R_{>0}$ \sf y
\beq5.48
\xi^b = x^a \hbox{ \ and \ }{\xi'}^{b'} = x^{a'}.
\e
Then $\xi=\xi'$ if $ab' = ba'$.
\elm

\proof
Let $z\in\R_{>0}$, $c,d\in\Na$. Then $(z^c)^d = z^{(cd)}$ by \E\Pr\ 2.2.3\,I3 where $M:=\R_{>0}$, ${\qu}:={\cdot}$ and $e:=1$. From \er{5.48} we infer
$\xi^{ba'} = (\xi^b)^{a'} = (x^a)^{a'} = x^{aa'} = x^{a'a} = ({\xi'}^{b'})^a = {\xi'}^{b'a} = {\xi'}^{ab'} = {\xi'}^{ba'}$. Hence $\xi^{ba'} =
{\xi'}^{ba'}$. Thus $\xi = \root ba' \of{\xi^{ba'}} = \root ba' \of {{\xi'}^{ba'}} = \xi'$.
\endproof

\If from Lemma \rf{l5.26}, \er{5.38} and (4.5.2) that
\beq5.49
\root b \of {x^a} = \root b' \of {x^{a'}} \hbox{ \ if \ } \frac ab = \frac{a'}{b'}, \q a,a',b,b'\in\N, \ x\in\R_{>0}.
\e
Note that if $\frac ab= \frac1n$, then $\root b \of {x^a} = \root n \of x$.

\E\Tf we may define $x^\frac ab\in\R_{>0}$ by setting
\beq5.50
x^\frac ab := \root b \of {x^a}, \q a,b\in\Na, \ x\in\R_{>0}.
\e
Note that $x\in\R_{>0}$ implies $x^a \in\R_{>0}$, hence $\root b \of {x^a}\in\R_{>0}$.

We next \es\ some basic algebraic \pp ies of \ra\ powers of an \el\ of~$\R_{>0}$.

\blm5.27
Let $a,b,n\in\Na$, $r,s\in\Q_{>0}$ and $x,y\in\R_{>0}$. Then\dw
\bea5.51
\root n \of {xy} &= \root n \of x \cdot \root n \of y,\\
x^\frac ab &= \bigl(\root b \of x\bigr)^a, \lb{5.52} \\
x^{rs} &= (x^r)^s, \lb{5.53} \\
\biggl(\frac 1x\biggr)^r &= \frac1{x^r}, \lb{5.54} \\
x^{r+s} &= x^r \cdot x^s. \lb{5.55}
\e
\elm

\proof \

``\er{5.51}'': Set $x':=\root n\of x$, $y':=\root n\of y$. Then $x=(x')^n$, $y=(y')^n$. Then $\root n \of {xy} = \root n \of {(x')^n(y')^n} \nad*=
\root n \of{(x'y')^n} = x'y' = \root n \of x \cdot \root n \of y$. In $\nad*=$ we use \E\Pr\ 2.2.3\,I5 with $M:=\R_{>0}$, ${\qu}:={\cdot}$ and
$e:=1$.

``\er{5.52}'': We use \In\ on $a\in\Na$. Let $x\in\R_{>0}$, $b\in\Na\sms1$, and let $A:=\{{c\in\Na}: x^\frac cb=\bigl(\root b\of x\bigr)^c\}$. We have
$1\in A$. Suppose $a\in A$. Then $x^\frac{a+1}b \nde5.50 {:=} \root b\of {x^{a+1}} = \root b \of {x^a x} \nde5.51 = \root b \of {x^a}\root b \of x
\nad{a\in A,\er{5.50}} = \bigl(\root b \of x\bigr)^a \root b \of x = \bigl(\root b \of x\bigr)^{a+1}$ by (2.2.3)\,I2. Hence $a+1\in A$, and $A=\Na$.

``\er{5.53}'': For convenience we recall
\bea5.56
(x^n)^\frac1n &= x = \bigl(x^\frac1n\bigr)^n,\\
(x^a)^\frac1b &= x^\frac ab = \bigl(x^\frac1b\bigr)^a. \lb{5.57}
\e
Let $r:=\frac ab$, $s:=\frac cd$. Set $x':=x^\frac1{bd}$. We have $x^{\frac ab\frac cd} \nad{(4.5.3)}= x^\frac{ac}{bd} \nde5.56 = ({x'}^{bd})^
\frac{ac}{bd} \nde5.57 = \bigl(({x'}^{bd})^\frac1{bd}\bigr)^{ac}\break \nde5.56 = (x')^{ac}$.

\Mo $\bigl(x^\frac ab\bigr)^\frac cd \nde5.56  = \bigl(({x'}^{bd})^\frac ab\bigr)^\frac cd$, and $({x'}^{bd})^\frac ab \nad{(2.2.12)}= ({x'}^{db})
^\frac ab \nad{(2.2.3)\,{\rm I3},(2.2.8),(2.2.12)} = (({x'}^d)^b)^\frac ab \break \nde5.57 = \bigl((({x'}^d)^b)^\frac1b\bigr)^a \nde5.56 = ({x'}^d)^a
= {x'}^{da} = {x'}^{ad} = ({x'}^a)^d$. \E\Tf $\bigl(x^\frac ab\bigr)^\frac cd = (({x'}^a)^d)^\frac cd$. Finally, $(({x'}^a)^d)^\frac cd \nde5.57 =
\bigl((({x'}^a)^d)^\frac1d\bigr)^c \nde5.56 = ({x'}^a)^c = ({x'})^{ac}$. \csq, $x^{\frac ab\frac cd} = \bigl(x^\frac ab\bigr)^\frac cd$, and \er{5.53}
holds.

``\er{5.54}'': Since $x\in\R_{>0}$, $\frac 1x$ belongs to $\R_{>0}$ by Lemma \rf{l3.18} since $x\cdot\frac1x=1$ in the group $\pz2{\R_{>0}}$. We have
to show $x^r\cdot\bigl(\frac1x\bigr)^r=1$, that is, $x^\frac ab\cdot\bigl(\frac1x\bigr)^\frac ab=1$. Set $y=\frac 1x$, thus $x\cdot y=1$. By \er{5.53}
we have $1^\frac ab = (x\cdot y)^\frac ab= x^\frac ab\cdot y^\frac ab$. It suffices to show $1^\frac ab=1$. We have $1^\frac ab \nde5.57 =
(1^a)^\frac 1b$. By (2.2.3)\,I4, $1^a=1$. \Mo $\root b\of 1 = 1$ since $1^b=1$.

``\er{5.55}'': Since $\frac ab+\frac cd \nad{(4.5.12)} = \frac{ad+bc}{bd}$, we have to show:
\beq5.58
x^\frac{ad+bc}{bd} = x^\frac ab \cdot x^\frac cd.
\e
We have $x^\frac{ad+bc}{bd} = (x^{ad+bc})^\frac1{bd} \nad{(2.2.3)\,\rm I2}= (x^{ad}\cdot x^{bc})^\frac1{bd} \nde5.51 = (x^{ad})^\frac1{bd}\cdot
(x^{bc})^\frac1{bd}$. As in the proof of \er{5.53} we have $x^{ad} = (x^a)^d$ and $x^{bc} = x^{cb} = (x^c)^b$. \Mo $(y^\frac1b)^ \frac1d \nde5.53 =
y^{\frac1b\frac1d} \nad{(4.5.3)} = y^\frac1{bd} \nad{(2.2.12)} = y^\frac1{db} = \bigl(y^\frac 1d\bigr)^\frac1b$. \E\Tf we have $(x^{ad})^\frac1{bd}
= \bigl(((x^a)^d)^\frac1d\bigr)^\frac1b = (x^a)^\frac1b = x^\frac ab$. Similarly, $(x^{bc})^\frac1{bd} = (x^{cb})^\frac1{bd} =
\bigl(((x^c)^b)^\frac1b\bigr)^\frac1d = x^\frac cd$. Hence \er{5.58} and \er{5.55} hold.
\endproof

\bex5.28
Let $x,y\in\R_{>0}$ and $r\in\Q_{>0}$. Show
\beq5.59
(x\cdot y)^r = x^r \cdot y^r.
\e
\eex

\bpr5.29
Let $x\in\R_{>0}$. Let $\psi_x$ denote the map from~$\Q$ into $\R_{>0}$ defined by
\beq5.60
\psi_x(r) := \bca
x^r &\hbox{for } r>0,\\
1 &\hbox{for } r=0, \\
(x^{-r})\mo & \hbox{for } r<0,
\eca
\e
where $x^r$ is defined in \er{5.50} and $y\mo$ denotes the inverse of $y\in\R_{>0}$ in the \ag\ $\pz2{\R_{>0}}$.

Then

\hph i,i, $\psi_x: \pz1{\Q} \to \pz2{\R_{>0}}$ is a \emph{monoid-\hm sm}.

\hph ii,, If $x>1$ $($resp.\ $x<1)$, then $\psi_x$ is strictly in\cre\ $($resp.\ de\cre$)$ \wrt the \og s of the fields $\Q$ and~$\R$. \E\Ip if
$x\ne1$ then $\psi_x$~is in\jc.
\epr

\proof \

(i) Since $\psi_x(0)=1$ by \er{5.60}, it suffices to prove
\beq5.61
\psi_x(r+s) = \psi_x(r)\cdot \psi_x(s), \q r,s\in\Q.
\e
In this part of the proof we skip the index $x$. We claim
\beq5.62
\psi(-t) = (\psi(t))\mo, \q t\in\Q.
\e
Indeed, if $t\ge0$ then \er{5.62} follows from  \er{5.60}. If $t<0$ then $t=-|t|$, hence $\psi(-|t|)\nde5.60 = \psi(|t|)\mo$. \E\Tf $\psi(-t)
= \psi(|t|) = \bigl(\psi(|t|)\mo\bigr)\mo = (\psi(-|t|))\mo = (\psi(t))\mo$. Note that
\er{5.61} holds if $r=0$ or $s=0$ since $\psi(0)=1$. In view of the \cmt ity of $+$~and~$\cdot$, it suffices to consider three
cases: $r,s>0$, $r,s<0$ and $r>0,s<0$.

``$r>0,\ s>0$'': Follows from \er{5.55}.

Note
\beq5.63
(\psi(r)\cdot\psi(s))\mo \nad{(4.3.5)} = \psi(s)\mo \cdot \psi(r)\mo = \psi(r)\mo \cdot \psi(s)\mo, \q r,s>0.
\e

``$r<0,\ s<0$'': $-(r+s)=(-r)+(-s) >0$, hence $\psi(r+s) = \psi(-(-(r+s))) = \psi(-((-r)+(-s))) \nde5.62 = (\psi((-r)+(-s)))\mo \nde5.55 =
(\psi(-r)\cdot \psi(-s))\mo \nde5.63 = (\psi(-r))\mo \cdot (\psi(-s))\mo \nde5.62 = \psi(r)\cdot\psi(s)$.

``$r>0,\ s<0$ and $r\ge|s|$'': Set $t:=r-|s|$, i.e.\ $r=|s|+t$, $t\ge0$. Then $r=(-s)+t$ or $r+s=t$. We obtain $\psi(r+s)=\psi(t)$. \E\oh $\psi(r)
\nde5.55 = \psi(-s)\cdot\psi(t)$ and $\psi(s) = \psi(-(-s)) \nde5.62 = (\psi(-s))\mo$, hence $\psi(r)\cdot \psi(s) = \psi(s)\cdot \psi(r) =
(\psi(-s))\mo \cdot (\psi(-s)\cdot\psi(t)) = \bigl((\psi(-s))\mo \cdot \psi(-s)\bigr) \cdot \psi(t) = 1\cdot\psi(t) = \psi(t)$.

``$r>0,\ s<0$ and $|s|>r$'': Set $t:=|s|-r$, hence $t=(-s)+(-r) = -(r{+}s)$, and $r+s=-t$. Thus $\psi(r+s)= \psi(-t)$. \E\oh $\psi(t)\cdot\psi(r)
\nde5.62 = \psi(|s|)= \psi(-s) = \psi(s)\mo$. \E\Tf $\psi(r)\cdot\psi(s) = \psi(-t) = \psi(r+s)$.
\ssk

(ii) We claim
\beq5.64
\psi_x(t)>1 \q\hbox{whenever $x>1$ and }t\in\Q_{>0}.
\e
We recall that if $x,y\in\R_{>0}$ \sf y $y<x$, then $y^a<x^a$ \fa $a\in\Na$ (see Remark \rf{r5.6}). \E\Ip $1=1^a < x^a$ if $a\in\N$ and $x>1$. We now
use the notation $f_n$ for the self-map of~$\R_{>0}$ defined by
\beq5.65
f_n(x):= x^n, \q x\in\R_{>0},\ n\in\Na.
\e
We just mentioned that the map $f_n$ is strictly in\cre, hence in\jc. By Theorem~\rf{t5.5} $f_n$ is bi\jc. We denote its inverse by~$g_n$. Hence
\beq5.66
g_n \circ f_n = \id_{\R_{>0}} \qh{and } f_n \circ g_n = \id_{\R_{>0}}.
\e
By Lemma 1.3.34\,(i) $g_n$ is strictly in\cre. \E\Tf if we use the notation $\root n\of y$ for $g_n(y)$, $y\in\R_{>0}$, we obtain $\root b\of{x^a}
>\root b\of1$, $x\in\R_{>0}$, $a,b\in\Na$, provided $x>1$. Using Lemma \rf{l5.26} we showed that the notation $x^r:=\root b\of{x^a}$, $r:=\frac ab$,
is justified for $r\in\Q_{>0}$. Since $1^b=1$, we infer $\root b\of 1 = 1$, hence $x^r>1$ whenever $r\in\Q_{>0}$, $x\in\R_{>0}$, which proves the
claim. Now let $r,s\in\Q$ with $r>s$. Set $t:=r-s$, that is, $r=s+t = t+s$. We obtain $\psi_x(r) = \psi_x(t+s) \nde5.55 = \psi_x(t)\cdot
\psi_x(s)$. Since $\psi_x(t)>1$ by \er{5.64} and $\psi_x(r)=1\cdot\psi_x(r)$, we infer from \er{3.10}:
\beq5.67
\psi_x(r) > \psi_x(s) \hbox{ if } r>s,\ r,s\in\Q.
\e
Finally, let $x\in\R_{>0}$, $x<1$. By \er{3.10} $x\mo>1$ and $\psi_\frac1x(r) > \psi_\frac1x(s)$. We claim
\beq5.68
\psi_\frac1x(t) = \psi_x(-t),\q x\in\R_{>0},\ t\in\Q_{>0}.
\e
Hence $0<\psi_x(-t)\mo<1$. Indeed, if $y\in\R_{>0}$, $y>1$, and $y\cdot y\mo = y\mo\cdot y=1$, we obtain from \er{3.10} and Lemma \rf{l3.18}:
$1 = y\cdot y\mo > 1\cdot y\mo = y\mo >0$. \csq, we have $\psi_x(r) = \psi_x(t)\cdot \psi_x(s)$, hence
\beq5.68a
\psi_x(r) < \psi_x(s) \hbox{ if } r>s,\ r,s\in\Q, \ x\in\R_{>0}, \ x<1.
\e

It remains to prove the claim. Let $x\in\R_{>0}$ and $t\in\Q_{>0}$. Let $x\mo\in\R_{>0}$ \sf y $x\cdot x\mo = x\mo\cdot x=1$. To prove
\beq5.69
(x\mo)^t = (x^t)\mo
\e
let $t:=\frac ab$, $a,b\in\Na$. Then $x\mo = \bigl(\bigl(x^\frac1b\bigr)^b\bigr)\mo \nde5.8 = \bigl(\bigl(x^\frac1b\bigr)\mo\bigr)^b$. Hence
$(x\mo)^\frac1b = \bigl(x^\frac1b\bigr)\mo$. \E\Tf $(x\mo)^t = (x\mo)^\frac ab = \bigl((x\mo)^\frac1b\bigr)
^a = \bigl(\bigl(x^\frac1b\bigr)\mo\bigr)^a \nde5.8 = \bigl(\bigl(x^\frac1b\bigr)^a\bigr)\mo = \bigl(x^\frac ab\bigr)\mo = (x^t)\mo$.

This completes the proof of \E\Pr\ \rf{p5.29}.
\endproof

In \E\Pr\ 2.2.3 we gave a list of \pp ies of \IT s of an \el\ of an \ti{\am\/} $(M,\qu,e)$. We proved that the $n$-fold \IT, $n\in\N$, of an \el\
$x\in M$ \sf ies (2.2.3)\,I1--I5. In \E\Pr\ 4.3.97 we considered $z$-fold \IT s, $z\in\Z$, of an \el\ $x\in G$, where $(G,\qu,e)$ is an \ti{\ag}.
\E\pp ies are listed in (4.3.143)\,SI0--SI2, (4.3.143)\,SI4--SI7, and in (4.5.46)\,SI3, where SI stands for signed \IT s. These \IT s are usually
called \ti{powers} if the \mlv\ notation $(G,\cdot,e)$ is used (see \cite[pp.~36--37]{Alg}). The \rl\ between \po\ and \ng\ powers is given in SI17,
namely $-z\cdot x = z\cdot x\Inv$. \E\Ip if $z\in\Na$, then the \LHS\ is defined by the \RHS, i.e.\ $-n\dqu x:=n\dqu x\Inv$, $n\in\Na$.

Note that if $n:=1$, then $(-1)\dqu x = 1\dqu x\Inv = x\Inv$. The usual notation for powers of~$x$~\el\ of a group $\pz2G$ is $x^z:= z\ddt x$,
$z\in\Z$. \E\Ip $x\mo = (-1)\dqu x= x\Inv$. The notation $x\mo$ for the inverse of an \el\ of a ``\mlv'' group is standard. In the case of an ``\ad''
group the standard notation for $(-1)\dpl x$ is $-x$ (see \E\df\ 4.3.2).

\ssk
We, finally, consider the special case $\pz2G := \pz2{\R_{>0}}$ where \ra\ powers $x^r$, $x\in\R_{>0}$, $r\in\Q$, can be defined. \E\el s of~$\Q$ will
be denoted, as in \E\df\ 4.5.7, by~$0$, $\frac ab$, $-\frac ab$, where $a,b\in\Na$, $\frac ab$ defined in (4.5.1) and $-\frac ab:= \bigl(\frac
ab\bigr)\Inv$, where $\frac ab+\bigl(\frac ab\bigr)\Inv = 0$. \Mo if $i:\Z\to\Q$ denotes the map defined by
\beq5.70
i(0):=0, \q i(n):=\frac n1 \qh{and } i(-n):=-\frac n1, \q n\in\Na,
\e
one finds that $i$ is an in\jc\ (semi)ring-\hm sm from the ring $\pz0\Z$ into the field $\pz0\Q$ (which is also a ring). Thus $\Z$~is imbedded
in~$\Q$ (see Section~\ref{ass.2}), and we may \ti{identify} $z\in\Z$ with $i(z)\in\Q$. Note that, strictly speaking, we identify $\frac ab$ with
$+\frac ab$ in \E\df\ 4.5.7.

\bex5.30
Show that the map $i$ defined in \er{5.70} is an in\jc\ ring-\hm sm.
\eex

It turns out that the \ra\ powers of an \el\ of the group $\pz2{\R_{>0}}$ \sf y \pp ies (4.3.143)\,SI1--7 where $(G,\qu,e)$ (resp.\ $z,w,\Z$) are
replaced by $\pz2{\R_{>0}}$ (resp.\ $r,s,\Q$). Note that $\pz2{\R_{>0}}$ is a subgroup of the \mlv\ group of the field~$\R$.
\newpage

\bpr5.31
Let $\pz2{\R_{>0}}$ be the \mlv\ group of \po\ \el s of~$\R$. Given $r\in\Q$ and $x\in\R_{>0}$, let $x^r:=\psi_x(r)$, defined in \er{5.60}. Then \fa
$r,s\in\Q$ and $x,y\in\R_{>0}$ we have
\beq5.71
\aligned
{}&\bca
{\rm RI0} \q\ x^0 =1,\\
{\rm RI1} \q\ x^1 =x,
\eca\\
&\bca
{\rm RI2} \q\ x^{r+s} = x^r\cdot x^s,\\
{\rm RI3} \q\ (x^r)^s = x^{rs},
\eca \\
&\bca
{\rm RI4} \q\ 1^r = 1,\\
{\rm RI5} \q\ (x\cdot y)^r = x^r\cdot y^r,\\
{\rm RI6} \q\ (x\mo)^r = (x^r)\mo, \\
{\rm RI7} \q\ x^{-r} = (x\mo)^r,
\eca
\endaligned
\e
where $x\mo$ denotes the inverse of~$x$ in $\pz2{\R_{>0}}$, that is, $x\cdot x\mo = x\mo \cdot x =1$.
\epr

\proof \

RI0: $x^0 = \psi_x(0) \nde5.60 = 1$.

RI1: $x^1 = x$ (see Remark \rf{r5.6}).

RI7: follows from \er{5.69}.

RI6: follows from \er{5.68} for $t\in\Q_{>0}$, and from \er{5.60} for $t=0$. Let $r\in\Q$, $r<0$. Then $r=-|r|$, $|r|\in\Q_{>0}$. We have
$(x\mo)^{-|r|} \nad*= ((x\mo)\mo)^{|r|} = x^{|r|} = x^{-r}$, where $\nad*=$ is $\nde5.68 = $ with $x:=\frac1x=x\mo$.

RI4: If $r:=0$ then $1^0 = 1$ by RI0. If $r:=\frac ab$, $a,b\in\Na$, then $1^{\frac ab} \nde5.52 = (1^{\frac 1b})^a \nad*= 1^a \nad*= 1$.
In $\nad*=$ we used $1^b \nad{(2.2.27)} = 1$, $1^a =1$, and \er{5.38}. If $r<0$ then $r=-|r|$, hence $1^{|r|}=1$, and $1^r = 1^{-|r|} \nde5.68 =
(1^{|r|})\mo = 1\mo =1$, since $1\cdot1=1$.

RI5: If $r=0$ then $(x\cdot y)^0 = 1 = 1\cdot1 = x^0\cdot y^0$. If $r:=\frac ab$, $a,b\in\Na$, then $(x\cdot y)^{\frac ab} \nde5.57 =((|x\cdot y)^a)
^{\frac1b} \nad*= (x^a\cdot y^a)^{\frac1b}$ where $x^a$ (resp.~$x^b$) is the $a$-th (resp.~$b$-th) \IT\ of~$x$ in the group $\pz2{\R_{>0}}$. In
$\nad*=$ we applied (2.2.3)\,I5 with $(M,\qu,e):=\pz2{\R_{>0}}$, $E:=\N$. From \er{5.51} we infer $(x^a\cdot y^a)^{\frac1b} = (x^a)^{\frac1b}\cdot
(y^a)^{\frac1b} = x^r\cdot y^r$. If $r:=-\frac ab$, $a,b\in\Na$, we have $(x\cdot y)^{-\frac ab} \nde5.62 = \bigl((x\cdot y)^{\frac ab}\bigr)\mo
= (x^{|r|})\mo \cdot (y^{|r|})\mo \nde5.62 = x^r\cdot y^r$.

RI2: follows from \er{5.61}.
\endproof

\bex 5.32
Prove \er{5.71}\,RI3.
\eex

We now turn to the study of order \pp ies of the maps $t\mt x^t$ and $x\mt x^t$, $t\in\Q$, $x\in\R_{>0}$. To this end we shall use the \fw\ lemma.

\blm5.33
Let $(E,\ge)$ and $(E',\ge')$ be \os s, and let $\vf:E\to E'$ be bi\jc\ and de\cre. Then\dw

\hph i,ii, $\vf$ is strictly de\cre\sd

\hph ii,i, if $(E,\ge)$ is totally ordered, then $\vf\Inv$ is strictly de\cre\sd

\hph iii,, if $A$ is a \ns\ of~$E$ \st $\sup A$ $($resp.\ $\inf A)$ exists, then $\inf \vf(A)$ $($resp.\ $\sup A)$ exists and we have $\inf\vf(A) =
\vf(\sup A)$ $($resp.\ $\sup \vf(A)= \vf(\inf A))$.
\elm

\bex5.34
Prove Lemma \rf{l5.33}. (Hint: imitate the proofs of Lemmata 5.3.34, \rf{l5.2} and \rf{l5.11}.)
\eex

In view of \E\Pr\ \rf{p5.29} the map $t\mt x^t$, $x\in\R_{>0}$, from $\Q$ into $\R_{>0}$ is strictly in\cre\ (resp.\ de\cre) whenever $x>1$ (resp.\
$x<1$). More is true.

\blm5.35
Let $r,s,t\in\Q$, $x,y\in\R_{>0}$. Then
\bea5.72
{}& x^r < y^r \hbox{ if } x<y,\\
&\sup_{s<r} x^s = x^r = \inf_{s>r}x^s\hbox{ if }x>1, \lb{5.73}\\
&\inf_{s<r} x^s = x^r = \sup_{s>r}x^s\hbox{ if }x<1, \lb{5.74}\\
&0< x^t-x^s = x^{\min(t,s)}(x^{(t-s)}-1) \hbox{ if }x>1,\ t>s, \lb{5.75}\\
&0< x^s-x^t = x^{\max(t,s)}(x^{(s-t)}-1) \hbox{ if }x<1,\ t>s. \lb{5.76}
\e
\elm

\proof
Let $r:=\frac ab$, $s:=\frac cd$, $a,b,c,d\in\Na$.

``\er{5.72}'': Set $x':=x^\frac 1b$, $y':=y^\frac1b$. Then $x'=g_b(x)$, $y':=g_b(y)$, where $g_b$ is the inverse of~$f_b$ defined in~\er{5.5}.
We have $x'<y'$ by \er{5.35}, hence $x^r = (x')^a < (y')^a = y^r$ by (4.4.46).

``\er{5.73}'': We first prove: let $r>1$, then
\beq5.77
\hbox{there is no $u\in\Q_{>0}$ \st $r^n\le u$ \fa $n\in\Na$.}
\e
Since $r>1$, \te s $h\in\Q_{>0}$ \st $r=1+h$. Hence if $n\ge2$ we have $(1+h)^n \nad{(4.4.113)} = \suml_{k=0}^n \frac{n!}{k!(n-k)!}\,1^{n-k}h^k =
\frac{n!}{0!n!}\,1^{n}h^0 + \frac{n!}{1!(n-1)!}\,1^{n-1}h^1 + \suml_{k=2}^n \frac{n!}{k!(n-k)!}h^k > 1+nh$, since $0!= \frac{n!}{n!} =1$,
$\frac{n!}{(n-1)!}=n$, $1^{n-k}=1$ for $k\in[0,n]$, $h^0=1$ and $\suml_{k=2}^n \frac{n!}{k!(n-k)!}h^k> 0$. Suppose, for \cd ion, that \er{5.77} does
not hold, then \te s $u$ \st $u\ge r^n > 1+nh > nh$, $n\ge2$. By Lemma 4.5.52 \te s $\ov n\in\Na$ \st $\ov nh>u$, hence $(\ov n+1)h = \ov nh+h > \ov
nh>u$. Since $\ov n+1\ge2$, we obtain $u>r^{\bar n+1} > (\ov n+1)h$. Hence $(\ov n+1)h > u > (\ov n+1)h$, a \cd ion.

In\et y $(1+h)^n \ge 1+nh$, $n\in\Na$, $h\in\Q_{>0}$, is usually called \ti{Bernoulli's in\et y}.

We next prove: Let $r\in\Q_{>0}$ satisfy $r<1$. Then
\beq5.78
\hbox{there is no $e\in\Q_{>0}$ \st $e<r^n$ \fa $n\in\Na$.}
\e
Otherwise, $(r^n)\mo\h {<}\h e\mo$ \fa $n\in\Na$. Since $(r^n)\mo \nad{\rm RI7}= (r\mo)^n$, we would have ${(r\mo)^n\h {<}\h e\mo}$ \fa $n\in\Na$
with $r\mo>1$ \cd ing \er{5.77}.

We claim
\beq5.79
\inf_{r>0} x^r = 1 \qh{if }x>1.
\e
By \E\Pr\ \rf{p5.29}\,(ii) we have $1=x^0 < x^r$ \fa $r>0$. Hence $1$~is a \lo\ for the set $A:=\{x^r\in\R: r>0\}$. Since $A$~is \bb, \te s a greatest
\lo\ of~$A$, denoted by $\inf A$,
in view of Theorem \rf{t5.1} and Corollary \rf{c5.12}. Suppose, for \cd ion, that $\ov x:=\inf A>1$. Then $x^r\ge \ov x$ \fa $r>0$,
in particular \fa $n\in\Na$, $x^\frac1n\ge \ov x$, $n\in\Na$. By \E\Pr s \rf{p4.31} and \rf{p4.1} \te s $\ov s\in\Q_{>0}$ \st $\ov x>\ov s>1$. Hence
$x \ge \ov x{}^n>\ov s{}^n$ \fa $n\in\Na$, \cd ing \er{5.79}.

\Wanp prove the second \et y of \er{5.73}. Let $s,r\in\Q$ with $s>r$. Set $h:=s-r\in\Q_{>0}$. We have $\inf\limits_{h>0}x^h =1$ by \er{5.79}. \Mo
$x^s = x^{r+h} \nde5.71 = x^r\cdot x^h$. Setting $A:=\{x^h: h>0\}$, we have $\inf A=1$, hence $x^r = x^r \cdot \inf A \nde5.14 = \inf\{x^r\cdot x^h:
{h>0}\} = \inf\{x^{r+h}: h>0\} = \inf \{x^s: s>r\} = \inf\limits_{s>r} x^s$.

We now complete the proof of \er{5.73}. Let $s<r$ and $h:=r-s>0$. Then $s=r+(-h)$, $x^s=x^{r+(-h)} \nde5.71 = x^r\cdot x^{-h} \nde5.71 = x^r \cdot
(x^h)\mo$. Let $A$ be as above, then $1=\inf\limits_{h>0} x^h = \inf A$, and $1\nde5.71 = 1\mo = (\inf A)\mo$.

Set $A\mo := \{x^{-h}: h>0\}$. We claim
\beq5.80
\sup A\mo = (\inf A)\mo.
\e
Indeed, we have $x^{-h}\nde5.71 = (x\mo)^h$, hence the map $h\mt x^{-h}$ is de\cre\ since $0<x\mo<1$. \E\Tf $\sup A\mo = (\inf A)\mo$ by Lemma
\rf{l5.33}\,(iii). \csq, $\sup A\mo = (\inf A)\mo = 1\mo \nde5.71 = 1$. Finally, $x^r = x^r\cdot 1 = x^r\cdot \sup\limits_{h>0}x^{-h} \nde5.14 =
\sup\limits_{h>0} x^r\cdot x^{-h} = \sup\limits_{h>0} x^{r-h} = \sup\limits_{s<r}x^s$ since $s=r-h < r$ whenever $h>0$.

\ssk
``\er{5.74}'': Set $y:=x\mo\in \R_{>0}$. Then $y>1$ and $\sup\limits_{s<r} y^s \nde5.73 = y^r \nde5.73 = \inf\limits_{s>r}y^s$. Since the map $z\mt
z\mo$ from $\R_{>0}$ into $\R_{>0}$ is bi\jc\ and de\cre, we infer from Lemma \rf{l5.33} that $y^r=\sup\limits_{s<r} y^s = \sup\limits_{s<r}(x\mo)^s
= (\inf\limits_{s<r}x^s)\mo$, $y^r = \inf\limits_{s>r}y^s = \inf\limits_{s>r}(x\mo)^s =
(\sup\limits_{s>r} x^s)\mo$. Thus $\inf\limits_{s<r}x^s = (y^r)\mo = ((x\mo)^r)\mo = (x^{-r})\mo = x^r$. Similarly, $\sup\limits_{s>r} x^s = (y^r)\mo
= x^r$, and \er{5.74} holds.

\ssk
``\er{5.75}, \er{5.76}'': In view of \E\Pr\ \rf{p5.29}\,(ii) it suffices to prove the \et y in \er{5.75}, \er{5.76}. Let $t>s$, then
$\min(t,s)=s$, $\max(t,s)=t$, and $x^t-x^s \nde5.71 = x^s\cdot x^{t-s} - x^s = x^s(x^{t-s}-1)$ if $x>1$, $x^s-x^t = x^t\cdot x^{s-t} - x^t =
x^t(x^{s-t}-1)$ if $x<1$.
\endproof

\brm5.36
In \E\df\ 2.2.17 we introduced the \op\ ``\epc'' $x^t$ \fa $x,t\in\N$. The case $x=t=0$ is not interesting. We proved \pp ies \er{5.71}\,R0--R5 in
\E\Pr\ 2.2.18. In \E\Pr\ \rf{p5.29} we were able to define $x^t$ \fa $x\in\R_{>0}$ and $t\in\Q$, and to prove \er{5.71}\,R0--R7 (the proof of~R5 is
Exercise \rf{ex5.28}) in Lemma \rf{l5.27}. We found that the map $\psi_x:\Q\to \R_{>0}$ defined in \er{5.60} is an in\jc\ \hm sm from the group
$\pz1{\Q}$ into the group $\pz2{\R_{>0}}$ if $x\ne1$. We showed in particular cases that $\psi_x$ is not sur\jc\ (for example if $x=\frac 1p$,
$p$~prime). In fact, \fa $x\in\R_{>0}\sms1$ $\psi_x$~is not bi\jc, since $\Q$~is \ct y infinite and $\R$~is un\ct e (see the next section).
However, using Lemma \rf{l5.35}, we will be able to find a map denoted by $\log_x$ from $\R_{>1}$ (${} :=\{y\in\R: y>1\}$) into $\R_{>0}$, \sf ying
${\log_x} \circ \psi_x =\id_{\Q_{>0}}$.
\erm

\bpr5.37
Let $x,\a\in\R_{>1}$ and let
$$
A:=\{s\in\Q_{>0}: x^s<\a\}.
$$
Then $A$ is an ideal of $\Q_{>0}$ denoted by $\log_x\a$ and called the \tb{logarithm of the \nm~$\a$ to the base~$x$}. By \df\ $\log_x1:=0$.
\epr \index{logarithm of the \nm~$\a$ to the base~$x$}

\proof \

``$\Q_{>0}\sm A\ne\vn$'': Let $x\in\R$, $x>1$. Since $(\R,\ge)$ is \Ar\ by \E\Pr\ \rf{p4.31}, it follows from \E\Pr\ \rf{p4.1} that \te s $u\in\Q$
\sf ying $1 < j(u) < x$ where $j:\Q\to\R$ is defined in \er{4.69}. Suppose, for \cd ion, that $\Q\sm A=\vn$, then, since the \og\ of~$\R$ is total,
we would have $x^s<\a$ \fa $s\in\Q_{>0}$, in particular $x^n<\a$ \fa $n\in\Na$. \Mo $j(u)\nde5.72 < x^n < \a$ \fa $n\in\Na$, \cd ing \er{5.77}. Hence
$\Q\sm A \ne\vn$.

\ssk
``$A\ne\vn$'': Suppose, for \cd ion, that $A=\vn$. Then we would have $\a\le x^s$ \fa $s\in\Q_{>0}$. However, $\inf\limits_{s>0} x^s = 1$ by
\er{5.79}. \E\Tf $\a$ is a \lo\ for the set $\{x^s: s>0\}$ and $1$~is the greatest \lo\ of this set, we obtain $\a\le1$, \cd ing $\a>1$. Hence
$A\ne\vn$.

\ssk
``C2'': Let $s\in A$ and $s'\in\Q_{>0}$ \sf y $s'<s$. Then $x^{s'}<x^s$ by \E\Pr\ \rf{p5.29}\,(ii). Hence $x^{s'}<\a$ since $x^s<\a$, and $s'\in A$.

\ssk
``C3'': Let $s\in A$ and let $h\in\Q_{>0}$. Then $x^s < x^{s+h}$ as above. \Mo $\inf\limits_{h>0} x^h=1$. We conclude by using the \fw\ useful lemma.

\blm5.38
Let $B$ be a \nss\ of~$\R$. Then $\xi\in\R$ is the supremum of~$B$ $($resp.\ infimum of~$B)$ iff the \fw\ \cn s hold\/\dw
\bea5.81
{}& \xi\ge\eta \ (\hbox{resp.\ }\xi\le\eta) \hbox{ \fa}\eta\in B.\\
& \hbox{\E\fe $\ve\in\R_{>0}$ \te s $\b\in B$ \st $\b+\ve>\xi$ $($resp.\ $\b<\xi+\ve)$.} \lb{5.82}
\e
\elm

See Lemma \rf{l6.27}.

\bex5.39
Prove Lemma \rf{l5.38}.
\eex

Returning to the proof of ``C3'' we find that \fe $\ve>0$ \te s $\ov h\in\Q_{>0}$ \st $x^{\bar h}< 1+\ve$. We have $x^s<\a$ and we would like to have
$x^{s'}= x^{s+\bar h} = x^s \cdot x^{\bar h} < \a$. Thus if $x^s\cdot x^{\bar h} < x^s +x^s\cdot \ve < \a$, we would have $x^{s'}\in A$ provided
$x^s\cdot\ve < \a-x^s$, that is, $\ve< x^{-s}\cdot(\a-x^s) = \a\cdot x^{-s}-1$. Note that $\a\cdot x^{-s}-1>0$ ($\a x^{-s}>1$ \ev t to $\a>x^s$).
Thus choosing $\ve:=\frac12(\a x^{-s}-1)$ $(<\a x^{-s}-1)$ we obtain $x^{s'}<\a$ hence $s'\in A$.
\endproof

\blm5.40
Let $x\in\R_{>1}$ and let $\log_x$ denote the map from $\R_{\ge1}$ into~$\R_{\ge0}$ introduced in \E\Pr\ \rf{p5.37}. Then
\beq5.84
\log_x(x^t) = t \qh{\fa} t\in\Q_{\ge0}.
\e
\elm

\proof
If $t:=0$ then $x^0=1$ and $\log_x1 = 0$. Let $t\in\Q_{>0}$. Then $\log_x(x^t) = \{s\in\Q_{>0}: x^s<x^t\}$. By \E\Pr\ \rf{p5.29}\,(ii) we have
$x^s<x^t$ whenever $s<t$. The converse is also true. Indeed, if $x^s<x^t$, then $s\ne t$ and either $s<t$ or $t<s$. The case $t<s$ is impossible since
$x^t<x^s$. \E\Tf $\log_x x^t = \{s\in\Q_{>0}: s<t\} = i(t)$, the \pn\ ideal defined in \er{4.9}, which is identical to~$t$ in~$\Q_{>0}$.
\endproof

In the next section it will be proved that the map $\psi_x$ defined in \er{5.60} can be extended to~$\R$, and that this extended map is bi\jc\ with
the map $\log_x$ as inverse.

\ssk
We conclude this section with an ``\ex'' result concerning \pl\ \eq s in~$\R$, where the connectedness of \il s of~$\R$ plays an essential role.

\blm5.41
Let $a,b\in\R$, $a<b$, and let $f:[a,b]\to\R$ \sf y
\bga5.85
f(a)\ne f(b), \\
\hbox{\te s $M\in\R_{>0}$ \st } |f(x)-f(y)| \le M|x-y| \hbox{ \fa} x,y\in[a,b]. \lb{5.86}
\e
Then
\beq5.87
\bigl(\min(f(a),f(b)), \max(f(a),f(b))\bigr) \sbs f\bigl((a,b)\bigr).
\e
\elm

\brm5.42 \

\hph i,i, Lemma \rf{l5.41} is a weaker version of what is called Bolzano--Weierstrass ``Intermediate Value Theorem'' (see \cite[p.~34]{Nrs}). \E\cn\
\er{5.86} is useful but not necessary.

\hph ii,, If, in \ad\ to \er{5.85}, \er{5.86}, the map $f$ is required to be \ti{in\cre}, then $f(a)=\min(f(a),f(b))$, $f(b)=\max(f(a),f(b))$ and
$f([a,b])\sbs [f(a),f(b)]$, hence $f([a,b])=[f(a),f(b)]$. \Mo if $f$ is \ti{strictly in\cre}, then $f:[a,b]\to[f(a),f(b)]$ is a bi\jn.
\erm

\proof[Proof of Lemma \rf{l5.41}]
We first consider the case $a:=0$ and $b:=1$. For convenience, we still use letters $a$~and~$b$ instead of 0~and~1.
Set $\a:=\min(f(a),f(b))$ and $\b:=\max(f(a),f(b))$. Then $\a<\b$ since $f(a)\ne f(b)$. Suppose, for \cd ion, that \te s $\g\in(\a,\b)$ \st $\g\notin
f((a,b))$. Set $O_1:=\{x\in(a,b): f(x)<\g\}$ and $O_2:=\{x\in(a,b): f(x)>\g\}$. Since the \og\ of~$\R$ is total, we have $(a,b)=O_1\cup O_2$ and
$O_1\cap O_2=\vn$. The proof that $O_1$~and~$O_2$ are open is similar to the one given in the case of Theorem \rf{t5.5}.
We next show that both $O_1$ and~$O_2$ are not empty. We show that
$O_1\ne\vn$. We consider two cases: $f(a)<f(b)$ and $f(b)<f(a)$. In the first case
$f(a)=\a<\g<\b$. Let $h>0$ be \st $a+h<b$. We have $f(a+h) = f(a+h) - f(a) + f(a) \le |{f(a+h)-f(a)}|+f(a) \nde5.86 \le Mh+f(a)$. Choosing $\ov h
<M\mo(\g-f(a))$, we obtain $f(a+\ov h)<\g$. Thus $a+\ov h\in(a,b)$ and $f(a+\ov h)<\g$, hence $a+\ov h\in O_1$. We now choose $h>0$ \st $b-h>a$, hence
$f(b-h) = f(b-h) - f(b)+ f(b) \ge f(b) - |f(b-h)-f(b)| \ge f(b) - Mh$. Choosing $\ov h < M\mo(f(b)-\g)$, we find $b-\ov h \in (a,b)$ and $f(b-\ov h)
>\g$, hence $b-\ov h\in O_2$. The proof of the case $f(b)<f(a)$ is similar, and \Tf omitted. \If that $(a,b)$ is not connected, \cd ing the
connectedness of the \il\ $(a,b)$ of~$\R$. We proved the case $a=0$, $b=1$ in Lemma \rf{l5.9}. The proof of the general case is similar, hence
omitted.

In Section \ref{ass.8} the notion of open sets and open subsets of~$\R$ will be considered. We claim that the general case can be deduced from the
case $a:=0$, $b:=1$. Indeed, define a map $\vf:[0,1]\to[a,b]$ by setting $\vf(t):=(1-t)a+tb$, $t\in[0,1]$. Then $\vf(0)=a$, $\vf(1)=b$, $\vf(t)
\in[a,b]$, $t\in[0,1]$. Set $\psi:[a,b] \to[0,1]$ by setting $\psi(s):=\frac{s-a}{b-a}$, $s\in[a,b]$. Then $\psi(a)=0$, $\psi(b)=1$, $\psi(s)\in
[0,1]$ for $s\in[a,b]$. \Mo one verifies that $\psi\circ\vf(t)=t$, $t\in[0,1]$, and $\vf\circ\psi(s) = s$, $s\in[a,b]$. Thus $\vf$~and~$\psi$ are
bi\jc\ and $\psi=\vf\Inv$, $\vf=\psi\Inv$. In \ad, we have $0<\vf(t)-\vf(t') \le (b-a)(t-t')$ for $0<t'<t \le1$, $0<\psi(s)-\psi(s') \le (b-a)\mo
(s-s')$ for $a<s'<s\le b$, hence $|\vf(t)-\vf(t')| \le (b-a)|t-t'|$ for $t,t'\in[0,1]$, and $|\psi(s)-\psi(s')| \le (b-a)\mo|s-s'|$ for $s,s'\in
[a,b]$. \csq, if $f:[a,b]\to\R$ \sf ies \er{5.85}, \er{5.86} then $\wt f:[0,1]\to\R$ defined by $\wt f:=f\circ \vf$ \sf ies \er{5.85}, \er{5.86}
with $a:=0$, $b:=1$ and $M:=(b-a)M$. Indeed, $\wt f(0) = f\circ\vf(0) = f(a)$, $\wt f(1) = f\circ\vf(1) = f(b)$, $\wt f(0) = f(a)\ne f(b) = \wt f(1)$
and $|\wt f(t)-\wt f(t')| = |f\circ\vf(t) - f\circ\vf(t')| \le M|\vf(t)-\vf(t')| \le M(b-a)|t-t'|$, $t,t'\in[0,1]$. By Lemma \rf{l5.41} with $a:=0$,
$b:=1$, $f:=\wt f$ and $M:=(b-a)M$, we obtain $\bigl(\min(\wt f(0),\wt f(1)), \max(\wt f(0),\wt f(1))\bigr) \sbs \wt f((0,1)) = f\circ\vf((0,1))$.
Since $\vf$~is an \ois sm from $[0,1]$ onto $[a,b]$ and $\vf(0)=a$, $\vf(1)=b$, we have $\vf((0,1)) = (a,b)$, which proves the claim.

This completes the proof of Lemma \rf{l5.41}.
\endproof

\newpage
\bpr5.43
Let $f:\R\to \R$ be the \pl\ map defined by
\beq5.88
f(x) := \sum_{k=0}^n a_kx^k, \q x\in\R,
\e
where $a_k\in\R$ \fa $k\in\N$ \st $k\le n$. Suppose $n=2m+1$ \fs $m\in\N$, $a_n\ne0$ and $a_0\ne0$. Then \te s $\ov x\in\R$ \st
\beq5.89
f(\ov x)=0.
\e
\epr

\brm5.44 \

\hph i,i, If $a_0=0$, then $\ov x:=0$ \sf ies \er{5.89}. The case $m=0$ is trivial, since $\ov x:=-a_1\mo\cdot a_0$ \sf ies \er{5.89}. If $n=2m$,
$m\in\Na$, then $x^{2m}+1 = (x^m)^2+1 \ge1$ \fa $x\in\R$.

\hph ii,, Setting $g(x):= a_n\mo \cdot f(x)$, $x\in\R$, we have $g(x) = \suml_{k=0}^n b_kx^k$ where $b_k = a_n\mo a_k$, $1\le k\le n$. Note that
$b_n=1$, $b_0\ne0$ and $g(x)=0$ iff $f(x)=0$, $x\in\R$.  \E\Tf we may assume \wlg (wlog) that $a_n:=1$.
\erm

\proof
We first show that \te s $R\in\R_{>1}$ \st
\beq5.90
f(-R)<0 \qh{and } f(R)>0.
\e
Indeed, in this case, setting $M_1:=\suml_{k=0}^{n-1}|a_k|$ we find (using \er{5.43}, \er{5.44} and \In)
$|\suml_{k=0}^{n-1}a_kx^k| \le \suml_{k=0}^{n-1}|a_k|\,|x^k| \le \suml_{k=0}^{n-1} |a_k|\,|x|^{n-1} \le M_1 |x|^{n-1}$ for $|x|>1$. Thus
$f(x)\ge\break x^n - M_1x^{n-1} = x^{n-1}(x-M_1) >0$ if $x>M_1+1$, and $f(x) \le -|x|^n + M_1|x|^{n-1} =
-|x|^{n-1}(|x|-M_1) <0$ if $x<-(M_1+1)$. Setting $R:=M_1+1$ we obtain \er{5.90}, hence \er{5.85} with $a:=-R$ and $b:=R$. Next we prove \er{5.86}. Let
$x,y\in\R$ be \st $|x|,|y|\le R$. Then $|f(x)-f(y)| = \bigl|\suml_{k=1}^n a_k(x^k-y^k)\bigr| \le \suml_{k=1}^n |a_k|\,|x^k-y^k|$. From (4.4.46) we
obtain $|x^k-y^k| \le |x-y|\suml_{l=0}^{n-1} |y|^{n-1-l}|x|^l \le |x-y| n R^{n-1}$, $1\le k\le n$. \If that $|f(x)-f(y) \le nR^{n-1}M_1|x-y|$ for
$|x|,|y|\le R$. Hence \er{5.86} holds with $ M=nR^{n-1}M_1$. \csq, in view of Lemma \rf{l5.41} \te s $\ov x\in[-R,R]$ \st $f(\ov x)=0$ by \er{5.87}.
\endproof

\brm5.45 \

\hph i,ii, \E\Pr\ \rf{p5.43} guarantees the \ex\ of a root of~$f$, but does not give a formula or an algorithm for computing it.

\hph ii,i, The \il\ $[-R,R]$ gives a rough ``location''of the root. Better estimates can be found for specific \f s~$f$.

\hph iii,, Lemma \rf{l5.41} can be used to prove the \ex\ of an $n$-th root without \rt ion on~$n$. Indeed, given $\xi\in\R_{>0}$, setting $f(x):=
x^n-\xi$ we find $f(0)=-\xi<0$ and $f(x)>0$ provided $x^n>\xi$. \E\Ip if $\xi<1$, then $f(1)>0$.
\erm

\newpage
\Subsubsection{Decimal expansions and logarithms}\label{ass.6}

The first part of this section is concerned with \de s of \po\ \el s of the \of\ $(\R,\ge)$ introduced in Theorem 4.30. We shall prove and generalize
assertions of Exercises 4.5.63, 4.5.65, 4.5.66, Remark 4.5.67. Some of these results will be used in the next part of this section.

\bdf6.1
An \el\ $r$ of the field of \qt s $(\wh\Z,\wh+,\wh\cdot,\wh0,\wh1)$ of the domain\break $\pz0\Z$ (see \E\Pr\ 4.5.40) is called an \tb{$a$-decimal \nm}
if $r=\frac z{a^n}$ \fs $z\in\Z$, $n\in\Na$ and $a\in\Na\sms1$. The set of $a$-\dn s will be denoted by~$\Q_a$ (not standard notation). The
$a$-\dn~$r$ is called \ti{dyadic} if $a:=2$ and simply \ti{decimal\/} if $a:=10$.\index{a-decimal \nm@$a$-decimal \nm}\index{dyadic}
\edf

We shall identify \el s $z$ of the domain $\pz0\Z$ with \el s of its field of \qt s $(\wh\Z,\wh+,\wh\cdot,\wh0,\wh1)$ of the form $\frac z1$ (see
\E\Pr\ 4.5.40). This is possible since the map
\beq6.1
z \mt \frac z1, \q z\in\Z,
\e
is a well-defined in\jc\ ring-\hm sm. \E\Tf we shall also use notations $0$ and~$1$ instead of $\wh0$ and~$\wh1$ when no ambiguity arises.

\bex6.2
Show that (2.4.38), (2.4.39), (2.4.41), (2.4.42) hold whenever $m,n,p$ belong to an \of\ $(F,\ge)$ and $a,b$ are maps from $[0,n]$ into~$F$,
$n\in\N$.
\eex

\blm6.3
The set $\Q_a$ is a subdomain $($see \E\df\ \rf{d3.22}$)$ of the field~$\wh\Z$, and is \ct y infinite.
\elm

\proof \

(i) ``$\Q_a$ is a \ti{sub\sr} of $\wh\Z$ (see \E\df\ 4.4.33)'': Note that $\wh0\in\wh\Z$ belongs to~$\Q_a$ since $\wh0\nad{(4.5.86)}{:=} \frac01
\nad{(4.5.83)}=\frac0a$. Let $r,s\in\Q_a$. We have $r=\frac z{a^n}$, $s=\frac w{a^m}$ \fs $z,w\in\Z$, $n,m\in\Na$. Then $r\hpl s\nad{(4.5.84)}=
\frac{za^m+a^nw}{a^n\cdot a^m}$. Set $z':=za^m+a^nw \in\Z$, $n':=n+m\in\Na$. Then $r\hpl s= \frac{z'}{a^{n'}}\in\Q_a$. Hence $\Q_a$ is a \sbm\ of
$(\wh\Z,\wh+,\wh0)$.

Note that $\wh1\in\wh Z$ belongs to $\Q_a$ since $\wh1=\frac11=\frac aa$. \Mo if $r,s\in\Q_a$, then
$$
r\hcd s = \frac z{a^n} \hcd \frac w{a^m} \nad{(4.5.85)} = \frac{zw}{a^{n+m}}\in\Q_a.
$$
Hence $\Q_a$ is a \sbm\ of $(\Z,\wh\cdot,\wh1)$.

(ii) ``$(\Q_a,\wh+,\wh\cdot,\wh0,\wh1)$ \ti{is a domain}'': $(\Q_a,\wh+,\wh0)$ is an abelian semiring by~(i) and \E\df\ 1.2.2. Let $r:=\frac z{a^n}\in
\Q_a$ and $s:=\frac{-z}{a^n}\in\Q_a$. Then $r\hpl s= \frac{za^n+(a^n(-z))}{a^n\cdot a^n} = \frac0{(a^n)^2}=\wh0$. Hence $s$ (resp.~$r$) is the inverse
of~$r$ (resp.~$s$) in the monoid $(\Q_a,\wh+,\wh0)$. Thus $(\Q_a,\wh+,\wh0)$ is an \ti{\ag}. \Mo $\Q_a\sms{\wh0}$ is a \sbm\ of the monoid
$(\Q_a,\wh\cdot,\wh1)$. Indeed, if $r,s\in\Q_a\sms{\wh0}$, then as above $r\hcd s= \frac{zw}{a^{n+m}}$ where $z,w\in\Z\sms0$. We have $z\cdot w
\in\Z\sms0$ since $\pz0\Z$ is a domain, hence $r\hcd s\in\Q_a\sms{\wh0}$. \If that $\Q_a$~is a domain, hence a subdomain of the field $\wh\Z$ by
\E\df\ \rf{d3.22}.

(iii) ``$\Q_a$ \ti{is \ct y infinite}'': The map $j_1$ from~$\N$ into $\Q_a$ defined by $j_1(k):=\frac{ka}a$ is in\jc, since $\frac{ka}a=\frac{la}a$,
$k,l\in\N$, implies $ka^2=la^2$, hence $k=l$ since $a^2\ne0$. \E\Tf $\Q_a$ is infinite by Theorem 1.4.18\,(iv) and \E\Pr\ 1.4.27\,(A)(i). \Mo
$\Q_a$~is included in $\wh\Z$, hence $\Q_a$ is \ct y infinite by \E\Pr\ 1.4.27\,(A)(iv).
\endproof

As a subdomain of the \of\ $(\wh\Z,\ge,0)$ the domain $\Q_a$ inherits the \og\ of~$\wt\Z$ by \rt ing it to~$\Q_a$. In what follows we shall be only
interested in $\Q_{a,\ge0}$ the set of non\ng\ $a$-\dn s. One verifies that $(\Q_{a,\ge0},\wh+,\wh0)$ is a \PM, itself being a \sbm\ of the ordered
semifield $\pz0{\Q_{\ge0}}$. Note that if $r,s\in\Q_{a,\ge0}$ and $r\ge s$, then \te s $p\in\Q_{\ge0}$ \st $r=s+p$. Since $p=r\hpl(-s)\in\Q_a$, we
have $p\in\Q_{a,\ge0}$ hence the \og\ of $\Q_{a,\ge0}$ is the \nog\ of the \PM\ $(\Q_{a,\ge0},\wh+,\wh0)$. Since $\Q_{a,\ge0}$ is a subset
of~$\Q_{\ge0}$ we may invoke \E\Pr\ 4.5.61 to show that \fe $r\in\Q_{a,\ge0}$ \te s \ooo pair $(c,p)\in\N\t[0,1)_\Q$ \sf ying $r=c+p$. Then $c$
(resp.~$p$) is called the integral (resp.\ fractional)\index{integral part}\index{fractional part} part of~$r$, denoted by $\Int(r)$ (resp.\ $\Fr(r)$). Observe that \fe $c\in\Na$, we have
$c=\frac c1 \nad{(4.5.2)}= \frac{ac}c$, hence $c$ is an $a$-\dn.

\bdf6.4
A \po\ \ra\ \nm\ $x(\in\Q_{>0})$ is said to be a \ti{\po\ $a$-decimal fraction} if $x\in\Q_a\cap [0,1)_\Q$,\index{positive $a$-decimal fraction} where $[0,1)_\Q := \{r\in\Q: 0\le r<1\}$.
\edf

\bns6.5
Let $a\in\Na\sms1$. Then
\beq6.2
\bca
D_a := \{l\in\N: 0\le l\le a-1\}\\
\dot D_a:= D_a\sms0.
\eca
\e
\ens

Compare \er{6.2} with (2.5.1).

\blm6.6
Let $x$ a \po\ $a$-decimal fraction. Then \te\ \ooo $N\in\Na$ and \ooo map $A:[1,N]\to D_a$, where $[1,N]:= \{l\in\Na: 1\le l\le N\}$, \st the \fw\
holds\dw
\bea6.3
{}&x = \sum_{l=1}^N A_la^{-l}, \\
&A_N\ne0. \lb{6.4}
\e
\elm

\proof
By \df\ \te\ $k'\in\Na$ and $n'\in\Na$ \st $x=\frac{k'}{a^{n'}}$ and $a^{n'}>k'$. Let $M:=\{l\in\N: a^l \hbox{ divides }k'\}$. Then $0\in M$ since
$a^0=1$ divides~$k'$. Since $a\ge2$ we have $a^p<a^q$ \fa $p,q\in\N$, $p<q$, by (2.2.33). We claim that $M$~is \ba. Indeed, set $c_l:=a^l-1$,
$l\in\N$. Then $c_0=0$ and the map $l\mt c_l+1$ from $\N$ into~$\N$ is strictly in\cre\ by (2.2.33), hence the same holds for the map $l\mt c_l$.
\If from Lemma 2.1.32 where $(\wt E,\wt e,\wt S):= (\N,0,S)$ that \te s \ooo $\ov l\in\N$ \st $k'\in [c_{\bar l},c_{\bar l+1}]$. Thus if $l\in M$,
then $a^l|k'$, and $a^l\le k'$ by (3.1.14). \E\Tf $a^l < c_{\bar l+1} = a^{\bar l+1} - 1 < a^{\bar l+1}$. From (2.2.35) we infer $l<\ov l+1$, which
proves the claim. By Theorem 1.3.38, \te s \ooo $l\in\N$, denoted by~$\wt l$, \st the \fw\ holds:
\beq6.5
a^{\tilde l} | k' \q\hbox{and}\q a^{\tilde l+1} \nmid k'.
\e
\Mo we have $a^{n'}>k'$. Setting
\beq6.6
k:= \frac{k'}{a^{\tilde l}}\in\Na,
\e
we find $a^{n'}> ka^{\tilde l} \ge a^{\tilde l}$. Hence $n'>\wt l$ by (2.2.35). Setting
\beq6.7
n:=n'-\wt l \in \Na,
\e
we find
\beq6.8
x = \frac k{a^n} \qh{\fs} k,n\in\Na
\e
\sf ying
\beq6.9
\phi_a(k) \ne0 \qh{and } a^n>k,
\e
where $\phi_a(k)$ is defined in (4.1.20).

Indeed, $x = \frac{k'}{a^{n'}} = \frac{ka^{\tilde l}}{a^{n'}} = \frac k{a^{n'-\tilde l}} = \frac k{a^n}$. \Mo $\frac k{a^n} = \frac{k'}{a^{n'}} < 1$
implies $k<a^n$. Finally, $a^{\tilde l+1}\nmid k'$ implies $a^{\tilde l+1}\nmid ka^{\tilde l}$ implies $a\nmid k$, that is, $\phi_a(k)\ne0$ by
(4.1.19).

We now show that the \rp ation of $x$ by means of \er{6.8}, \er{6.9} is unique. Suppose, for \cd ion, that \te\ $k,k'',n,n''\in\Na$ \st
$\frac k{a^n} = \frac{k''}{a^{n''}}$ with $\phi_a(k)\ne0$ and $\phi_a(k'')\ne0$. Then $k=k''$ and $n=n''$. Suppose, for \cd ion, that \te s $n>n''$
(resp.\ $n''>n$). Then \te s~$p$ (resp.~$p'$) in~$\Na$ \st $n=n''+p$ (resp.\ $n''=n+p'$). Then since $ka^{n''} = a^nk''$, we obtain $ka^{n''} =
a^{n''}a^pk''$ (resp.\ $ka^n a^{p'}=a^nk''$). Dividing by $a^{n''}$ (resp.~$a^n$) we obtain $k=a^p k''$ (resp.\ $k''= a^{p'}k$). Hence $\phi_a(k)=0$
(resp.\ $\phi_a(k'')=0$). A~\cd ion.

Summarizing, we proved that given a \po\ $a$-decimal fraction~$x$ \te s \ooo pair $(k,n)\in\Na\t\Na$ \sf ying \er{6.8} and~\er{6.9}. We conclude the
proof by using Lemmata 2.5.3 and 2.5.4 where $m:=k$. We obtain \ex\ and \uq\ of $N\in\N$ and of a map $B:[0,N] \to D_a$ \sf ying $B(N)\ne0$ and
$k = \suml_{l=0}^N B_la^l$. Dividing by~$a^n$ we obtain $x= a^{-n}\suml_{l=0}^N B_la^l = \suml_{l=0}^Na^{-n}B_la^l = \suml_{l=0}^N B_la^{l-n}$,
together with $\phi_a(k)\ne0$. Note that $\phi_a(k) = \phi_a\bigl(\suml_{l=0}^N B_la^l\bigr) = \phi_a(B_0) \nad{(4.1.22)}= B_0$ since $\phi_a(a^l)
\nad{(4.1.21)}= 0$ whenever $l>1$, and $\phi_a(a^0)= \phi_a(1) =1$ by (4.1.23). \E\Tf $B_0\ne0$.

\Mo we have $k<a^n$, that is, $\suml_{l=0}^NB_la^l < a^n$. Since $B_la^l\ge0$, $l\in[0,N]$, we obtain $B_Na^N \le \suml_{l=0}^NB_la^l$. Since $B_N
\ne0$, we have $B_N\ge1$, hence $a^N \le B_Na^N \le \suml_{l=0}^NB_la^l <a^n$. Thus $N<n$ by (2.2.35). We have $\suml_{l=0}^NB_la^{l-n} =
\suml_{l\in[0,N]}B_la^{l-n}$. Setting $p:=n-N \in\Na$ we find that $\suml_{l\in[0,N]}B_la^{l-n} = \suml_{i\in[p,n]}B_{n-i}a^{-i}$. Setting
\beq6.10
A_l:=\bca
0, & 0\le l\le p-1,\\
B_{n-l}, & p \le l \le n,
\eca
\e
we finally find $x=\suml_{l=0}^N A_la^{-l}$, which is \er{6.3}. \Mo $A_n=B_0\ne0$, which is \er{6.4}.
\endproof

\bxs6.7 \

\hph i,i, $a:=10$ (see Notation 2.5.12).
$$
\begin{array}{lll}
k':=2350\ & n':=4\ & x:=\dfrac{2350}{10^4} = \dfrac{2350}{10'000} \\[8pt]
k:= 235 & n:=3 & x =\dfrac{235}{1000} \\[8pt]
\multispan3{$x = 2\cdot10\mo + 3\cdot10^{-2} + 5\cdot10^{-3} = 0{,}235.$}
\end{array}
$$
Note that $x=\dfrac{47}{200}$.

\hph ii,, $a:=2$, $x:=\dfrac{34}{64}$,
$$
\begin{array}{ll}
k':=100010\ (34) \q & n':= 1000000\ (64)\\
k:= 10001 \ (17) \q & n:=100000 \ (32)\\
\multispan2{$x=1\cdot2\mo + 0\cdot2^{-2} + 0\cdot2^{-3} + 0\cdot2^{-4} + 1\cdot2^{-5}.$}
\end{array}
$$
Note that if $x:=0{,}235$, then $10\cdot x=2{,}35$, hence $2$~is the integral part of~$2{,}35$. \Mo $100\cdot x = 23{,}5$, hence $3$~is the last digit
of~$23$, the integral part of~$23{,}5$. We recall that $\phi_{10}(n)$, $n\in\N$, defined in (4.1.20), is the last digit of~$n$.
\exs

\blm6.8
Let $N\in\Na$, $a\in\Na\sms1$, and let $A$~denote a map from $[1,N]$ into $D_a$. Then the \fw\ holds\dw
\bea6.11
{}&\sum_{l=1}^N A_la^{-l} \in \Q_{a,\ge0} \cap [0,1)_\Q,\\
&A_l = \phi_a\Bigl(\Int\bigl(a^l\cdot \sum_{i=1}^N A_ia^{-i}\bigr)\Bigr),\ l\in[1,N], \lb{6.12}
\e
where $\phi_a$ is defined in $(4.1.20)$.
\elm

\proof \

``\er{6.11}'': Note that $A_la^{-l}\in\Q_{\ge0}\cap \Q_a$ \fa $l\in[1,N]$. Since $\pz1{\Q_{\ge0}}$ is a \hbox{\PM} and $\pz1{\Q_a}$ is a group by
Lemma \rf{l6.3}, $\suml_{l=1}^N A_la^{-l} \in \Q_{\ge0}\cap\Q_a =: \Q_{a,\ge0}$. It remains to show that $\suml_{l=1}^N A_la^{-l}<1$. Since $0\le
A_l\le a-1 $ \fa $l\in[1,N]$, we have $\suml_{l=1}^N A_la^{-l} \nad{(2.4.41)}\le \suml_{l=1}^N(a-1)a^{-l} \nad{(2.4.35)}= (a-1)\suml_{l=1}^N a^{-l}$.
We recall that in every field $\pz0X$ the \fw\ formula holds
\beq6.13
\sum_{l=1}^N u^l + \frac{u^{N+1}}{1-u} = \frac u{1-u} \q \hbox{\fa}u\in X\sms1
\e
(see for example the proof of \E\Pr\ 4.5.60). \E\Ip if $(X,\ge)$ is an \of, we have
\beq6.14
0 < \frac u{1-u} - \sum_{l=1}^N u^l = \frac{u^{N+1}}{1-u} \q \hbox{\fa} 0<u<1.
\e
\E\Tf $\suml_{l=1}^N a^{-l} < \frac{a\mo}{1-a\mo} = \frac1{a-1}$. Hence $\suml_{l=1}^NA_la^{-l} \le(a-1)\suml_{l=1}^Na^{-l} < (a-1)\cdot \frac1{a-1}
=1$.

``\er{6.12}'': Set $x:=\suml_{i=1}^NA_i a^{-i}$. We consider three cases: (i)~$N=1=l$, (ii)~$N>1$, $l=N$, (iii)~$N>1$, $l<N$.

(i) $x=A_1a\mo$, hence $a\cdot x = A_1 = \Int(A_1) = \phi_a(\Int(A_1)) = A_1$ by (4.1.29) since $A_1<a$.

(ii) $a^N\cdot x = \suml_{i=1}^N A_ia^{N-i}$. Note that $A_ia^{N-1}\in\N$ \fa $i\in[1,N]$, hence $a^N\cdot x=\break \Int(a^N\cdot x)$.
$\suml_{i=1}^N A_ia^{N-i} \nad{(2.4.32)}= \suml_{i=1}^{N-1}A_ia^{N-i} + A_N$. Note that $a^{N-i} = a\cdot a^{(N-1)-i}$ for $i\in[1,N-1]$. Hence
$\phi_a\bigl(\suml_{i=1}^N A_ia^{N-i}\bigr) = \phi_a\bigl(a\suml_{i=1}^{N-1} A_ia^{(N-1)-i}+ A_N\bigr) \nad{(4.1.21)}= \phi_a(A_N) \nad{(4.1.22)} =
A_N$ since $A_N<a$. Thus \er{6.12} holds.

(iii) $a^lx = \suml_{i=1}^N A_i a^{l-i} = \suml_{i=1}^l A_i a^{l-i}+ \suml_{i=l+1}^N A_i a^{l-i}$. Note that $A_i a^{l-i}\in\N$ \fa $i\in[1,l]$,
hence $\suml_{i=1}^l A_i a^{l-i}\in\N$. \Mo $\suml_{i=l+1}^N A_i a^{l-i}= \suml_{i=1}^{N-l} A_i a^{-i} < 1$ by \er{6.11}. Hence $\Int(a^lx) =
\suml_{i=1}^l A_i a^{l-i}$ by \E\Pr\ 4.5.61. Hence as above $\phi_a\bigl(\suml_{i=1}^l A_i a^{l-i}\bigr) = A_l$. Thus \er{6.12} holds.
\endproof

We now consider an \ti{infinite} \sq~$A$ of \el s of~$D_a$, i.e.\ a map from~$\Na$ into~$D_a$. We define a \sq\ $\zb sn{\Na}$ of \el s of
$\Q_{a,\ge0}$ by setting
\beq6.15
s_n:= \sum_{l=1}^n A_l a^{-l}, \q A_l\in D_a, \q n\in\Na.
\e
\E\fe $p\in\Na$ we have $s_{n+p} = \suml_{l=1}^{n+p} A_l a^{-l} = \suml_{l=1}^n A_l a^{-l} + \suml_{l=n+1}^{n+p} A_l a^{-l}$. Since
$\suml_{l=n+1}^{n+p}A_l a^{-l}\ge0$, we obtain
\beq6.16
s_n \le s_{n+p} \q \hbox{\fa} n,p\in\Na.
\e
Thus the \sq\ $\zb sn{\Na}$ is in\cre.

\E\oh we have $s_n<1$ \fa $n\ge1$ by \er{6.11}. If $i$ denotes the map defined in \er{4.9}, we find $i(s_n)\le i(1)$ by \er{4.10}. We have $i(s_n)
\supset \bcl_{k=1}^n i(s_k)$ by \er{6.16}. Thus $\bcl_{n\in\Na}i(s_n)\sbs i(1)$ (possibly \et y). Since $i(1)\ne\Q_{>0}$, we find by \E\Pr\
\rf{p4.13} that \te s $\a\in\R_{>0}$ \st
\bea6.17
\a&= \bigcup_{n\in\Na} i(s_n),\\
\a&\sbs i(1).\lb{6.18}
\e
Observe that $\{i(s_n)\in\R_{\ge0}: n\in\Na\}$ is a set of \el s of $\R_{>0}\cup\{0\}$, and that $\bcl_{n\in\Na}i(s_n)$ is equal to $\sup
\{i(s_n)\in\R_{\ge0}: n\in\Na\}$ in $(\R,\ge)$. Hence
\beq6.19
\a=\sup\{i(s_n): n\in\Na\}.
\e

Recall that $\sup\{\g\}=\g$, $\g\in\R$, in view of \E\df s and Notations 3.1.13. \E\Ip if $A_i:=0$ \fa $i\in\Na$, then $\a=0$ in \er{6.19}.

\newpage
\bdf6.9
Let $a\in \Na\sms1$ and let $\a\in\R_{\ge0}$. A~\sq\ $\zb Al{\Na}$ of \el s of~$D_a$ is said to be a \ti{\de\ to base~$a$} of~$\a$, if \er{6.19}
holds. A~\de\ to base~$a$ is called \ti{terminating} (see \cite[p.~6]{Niven}) if \te s $N\in\Na$ \st $A_l=0$ \fa $l>N$, and \ti{non-terminating}
otherwise.\index{decimal expansion}
\edf

In view of Lemma \rf{l6.6} every \po\ $a$-decimal fraction~$x$ has a terminating \de\ to base~$a$. However, the \nm~$x$ has also a non-terminating
\de\ to base~$a$.

\blm6.10
Let $x$ \sf y \er{6.3}, \er{6.4}. Let $B:\N\to D_a$ be defined by
\beq6.20
B_l :=
\bca
A_l, &\hbox{if }1\le l<N, \ N>1,\\
A_N-1, &\hbox{if }l=N,\\
a-1, &\hbox{if }N<l.
\eca
\e
Set
\beq6.21
t_n:=\sum_{l=1}^n B_l\cdot a^{-l}, \q n\ge1.
\e
Then the \sq\ $\zb tn\N$ is in\cre\ and
\beq6.22
i(x) = \sup\{i(t_n): n\in\Na\} \hbox{ in }(\R,\ge).
\e
\elm

\proof
Note that $B_l\in[0,a-1]$, $l\in\Na$, hence the \sq\ $\zb tn\N \sbs \Q_{\ge0}$ is in\cre\ and bounded by~$1$ in view of \er{6.11}. We have
\beq6.23
x = t_n+a^{-n}, \q n\ge N+1.
\e
We use
\beq6.24
(a-1)\sum_{l=N+1}^n a^{-l} = a^{-N} - a^{-n}, \q n\ge N+1.
\e
We first prove \er{6.24}: $\suml_{l=1}^n a^{-l} \nad{(2.4.32)} = \suml_{l=1}^N a^{-l} + \suml_{l=N+1}^n a^{-l}$. Hence $\suml_{l=N+1}^n a^{-l}
= \suml_{l=1}^n a^{-l} - \suml_{l=1}^N a^{-l} \nde6.13 = \frac{a\mo}{1-a\mo} - \frac{a^{-(n+1)}}{1-a\mo} - \bigl(\frac{a\mo}{1-a\mo} -
\frac{a^{-(N+1)}}{1-a\mo}\bigr) = \frac{a^{-(N+1)}-a^{-(n+1)}}{1-a\mo}= \frac{a^{-N}-a^{-n}}{a-1}$.
We next prove \er{6.23}. If $N=1$ and $n\ge2$, we have $t_n+a^{-n} = (A_1-1)a\mo + \suml_{l=2}^n (a-1)a^{-l} + a^{-n}\nad{\er{6.20},\er{6.24}}=
(A_1-1)a\mo + a\mo - a^{-n} + a^{-n} \nde6.21 = A_1a\mo = x$. If $N>1$ and $n\ge N+1$, we have $c_n+a^{-n} \nad{\er{6.20},\er{6.24}} =
\suml_{l=1}^{N-1} A_la^{-l} + (A_N-1)a^{-N} + a^{-N} + a^{-n} + a^{-n} \nde6.21 = x$. This completes the proof of \er{6.23}. Set
\beq6.25
M:= \{a^{-n}\in\Q: n\in\Na\}.
\e

We claim
\beq6.26
\inf M=0 \hbox{ in }(\Q,\ge).
\e
Indeed, since $a\ge2$, we have $0<a\mo\le\frac12<1$, hence by \er{5.78} there is no $e\in\Q_{>0}$ \st $e<a^{-n}$ \fa $n\in\Na$. \Mo $0<a^{-n}$ \fa
$n\in\Na$, hence $\Q_{\le0}$ is the set of all \lo s for~$M$. Clearly, $0$~is the greatest \el\ of~$M$. \E\Tf \er{6.26} holds. Applying Lemma
\rf{l5.11} with $F:=\Q$ and $A:=M$, we find by \er{5.15}:
$$
\sup\{-a^{-n}\in\Q: n\in\Na\} = -\inf M=0.
$$
Since the \sq\ $\{-a^{-n}\}_{n\in\N}$ is in\cre, we also have
\beq6.27
\sup\{-a^{-n}\in\Q: n\ge m\}=0 \qh{\fa} m\in\N.
\e
Note that $t_n=x+(-a^{-n})$, $n\ge N+1$, by \er{6.23}. Thus $t_n=\th_x(-a^{-n})$, $n\ge1$, where $\th_x$ is defined in \er{3.6}. \If from
Lemma \rf{l5.2} with $E=E'=\Q$, $f=\th_x$ that $\sup\{t_n:n\ge N+1\} \nde5.1 = \th_x(\sup \{\a^{-n}:n\ge N+1\})\nde6.28 = \th_x(0)=x$. Since the \sq\
$\{t_n\}_{n\ge1}$ is in\cre, we obtain \er{6.22}.
\endproof

We now turn to the case where $x\in\R_{\ge0}$, $x<1$. We need a \gn\ of \E\Pr\ 4.5.61.

\bpr6.11
\E\fe\ $\a\in\R_{\ge0}$ \te s \ooo pair $(n,\rho)\in\N\t\zo 0,1 $, with $\zo0,1 :=\{\eta\in\R_{\ge0}: \eta<1\}$, \sf ying
\beq6.28
\a = f(n)+\rho,
\e
where $f(n)$ is defined in \er{4.26}.

The \nm\ $n\in\N$ is usually called the \tb{integral part\/} of~$\a$, and $\rho\in\zo0,1 $ the \tb{fractional part\/} of~$\a$.
\epr

The notations $[x]$, $\lfloor x\rfloor$, $\Int(x)$ are used for the integral part of~$x$, and $\{x\}$, $\Fr(x)$ for the fractional part of~$x$.

We shall adopt the notations
\beq6.29
\Int(\a):=n \q \hbox{and} \q \Fr(\a):=\rho,
\e
whenever $\a=f(n)+\rho$ as in \er{6.28}.

\proof
We identify $r$ and $f(r)$, $r\in\Q_{\ge0}$. By \E\Pr s \rf{p4.31} and \rf{p4.1} \te s $r\in\Q_{>0}$ \st $\a<r<\a+1$. Set $B:=\{m\in\N: m\le\a\}$.
Then $B\ne\vn$ and $B$~is \ba\ since $m\le\a<r$ \fa $m\in B$. By Theorem 1.3.38 \te s $n\in\N$ \st $n\le \a<n+1$. Note that $\a-n < (n+1)-n =1$,
hence $\a-n\in\zo0,1 $. Let $n'\in\N$ and $\rho'\in\zo0,1 $ \sf y $n'+\rho'=\a=n+\rho$. If $n=n'$ then $\rho=\rho'$ by \cnc ity. If $n>n'$ then
$\rho'-\rho= n-n'\ge1$, and $\rho'>\rho+1>1$ \cd ing $\rho'\in\zo0,1 $. The case $n'>n$ is similar. Thus \te s \ooo pair $(n,\rho)\in\N\t\zo0,1 $
\st $\a=n+\rho$. Set $(\a):=n$ and $[\a]:=\rho$.
\endproof

\begin{lem}[\cite{Rham}] \lb{l6.12}
Every \el\ $\a$ of\/ $\R_{>0}$ \sf ying $\a<1$ admits a non-terminating \de\ to base~$a$.
\elm

\proof
Let $\a\in\R_{>0}$, $\a<1$. In view of Lemma \rf{l6.10}, it suffices to consider $\a\in\R_{>0}\sm \Q_{a,\ge0}$ \st $\a<1$. We define \rc vely a \sq\
$\zb xn\N$ of \el s of $\zo0,1 _\R$ by setting:
\bea6.30
x_{n+1}&:=\Fr(i(a)\cdot x_n), \q n\in\N,\\
x_0&:= x. \lb{6.31}
\e
We also set
\beq6.32
A_{n+1}:=\Int(i(a)\cdot x_n), \q n\in\N.
\e
The \sq\ $\zb xn\N$ is well-defined in view of Theorem 1.1.7 where $(E,e,S):=(\N,0,S)$, $F:=\zo0,1 _\R$ and $f(y):=\Fr(\a\cdot y)$, $y\in\zo0,1 _\R$.

We claim that
\beq6.33
x = \sum_{l=1}^n A_la^{-l} + x_n\cdot a^{-n} \qh{\fa} n\in\Na.
\e
We use \In\ on $n\in\Na$. The case $n:=1$ holds since $\suml_{l=1}^1 A_la^{-l}+ x_1a\mo = A_1a\mo+x_1a\mo$. Note that $x_1=\Fr(ax_0)=\Fr(ax)
\nde6.29 =ax - \Int(ax)\nde6.32 = ax-A_1$. Hence $x_1a\mo = x-A_1a\mo$. Hence \er{6.33} holds for $n:=1$. We next assume that \er{6.33} holds for~$n$,
and prove that \er{6.33} holds for $n+1$. We have $\suml_{l=1}^{n+1} A_la^{-l} + x_{n+1}a^{-(n+1)} \nad{(2.4.32),\er{6.30}} = \suml_{l=1}^n A_la^{-l}
+ A_{n+1}a^{-(n+1)} + \Fr(ax_n)a^{-(n+1)} \nde6.33 = x-x_na^{-n} + A_{n+1}a^{-(n+1)} + \Fr(ax_n)a^{-(n+1)} \nad{\er{6.32},\er{6.28}} = x-x_na^{-n}
+ ax_n\cdot a^{-(n+1)} = x-x_n\cdot a^{-n}+ x_n\cdot a^{-n} = x$. This completes the proof of the claim.

From \er{6.30} we have $0\le x_n<1$ and from \er{6.33} we infer
\beq6.34
0 < x-\sum_{l=1}^n A_la^{-l} < a^{-n} \qh{\fa } n\in\Na.
\e
We set
\beq6.35
\b_n:= x - \sum_{l=1}^n A_la^{-l} \in \R, \q n\ge1.
\e
Clearly, $\b_n\ge0$, $n\ge1$, and $\zb \b n\Na$ is de\cre, since $\bigl\{\suml_{l=1}^nA_la^{-l}\}_{n\ge1}$ is in\cre.

We claim that
\beq6.36
\inf\{\b_n: n\ge1\} = 0 \q \hbox{in }(\R,\ge).
\e
Indeed, if not, \te s $\ve\in\R_{>0}$ \st $\b_n\ge\ve$ \fa $n\in\N$. However, $\b_n<a^{-n}$ for $n\ge1$ by \er{6.34}. Hence $\ve<a^{-n}$, $n\ge1$.
Since $(\R,\ge)$ is \Ar, \te s $e\in\Q_{>0}$ \st $e<\ve$. \E\Tf $e<a^{-n}$, $n\ge1$. In view of \er{5.78}, $e=0$. A~\cd ion. This proves the claim.

Note that
\beq6.37
x+(-\b_n) = \sum_{l=1}^n A_la^{-l}, \q n\ge1.
\e
As in the proof of Lemma \rf{l6.10}, we find $\sup\{-\b_n: n\ge1\}=0$ (by Lemma \rf{l5.11} with $F:=\R$), and $\sup\bigl\{\suml_{l=1}^n A_la^{-l}:
n\in\Na\} = \sup\{\th_x(-\b_n): n\in\Na\} = \th_x(0) = x$ (by Lemma \rf{l5.2} with $E=E'=\R$). This completes the proof of Lemma \rf{l6.12}.
\endproof

From now on we shall use the more usual notations for $\sup$ (resp.~$\inf$) of \ti{in\cre} (resp.\ \ti{de\cre}) \sq s of \el s of~$\R$, more
generally of an \os\ $(E,\ge)$:
\bea6.38
\sup_{n\ge m}\a_n &:= \sup\{\a_n\in\R: n\ge m\}, \q m\in\Na,\\
\inf_{n\ge m}\b_n &:= \inf\{\b_n\in\R: n\ge m\}, \q m\in\Na. \lb{6.39}
\e

\bex6.13
Let $(E,\ge)$ be an \os, and let $\zb\a n\Na$ (resp.\ $\zb\b n\Na$) be an in\cre\ (resp.\ de\cre) \sq\ of \el s of $(E,\ge)$. Let $m\in\Na$.
Suppose that $\sup\limits_{n\ge m}\a_n$ (resp.\ $\inf\limits_{n\ge m}\b_n$) exists. Show that the \fw\ holds:
\bea6.40
\sup_{n\ge m}\a_n &= \sup_{n\ge k}\a_n \qh{\fa}k\in\N,\\
\inf_{n\ge m}\b_n &= \inf_{n\ge k}\b_n \qh{\fa}k\in\N.\lb{6.41}
\e
One says that the supremum (resp.\ infimum) of an in\cre\ (resp.\ de\cre) \sq\ only ``depends'' on the tail of the \sq.
\eex

\bth6.14
Let $\a\in\R_{>0}$, $\a<1$. Then

\hph i,i, If $\a\notin\Q_a$, then $\a$ admits a unique non-terminating \de\ to the base~$a$ \itd in Lemma \rf{l6.12}.

\hph ii,, If $\a\in\Q_a$, then $\a$ admits a unique terminating \de\ to the base~$a$ \sf ying \er{6.3} and \er{6.4}, and a unique non-terminating \de\
to the base~$a$ \itd in \E\df\ \rf{d6.9}.
\eth

\proof
The ``\ex'' part of (i) follows from \E\Pr\ \rf{p6.11}. The ``\ex'' part of~(ii) follows from Lemma \rf{l6.6} and \E\Pr\ \rf{p6.11}.

``\ti{\E\uq}'': The \uq\ of a terminating expansion is proved in Lemma \rf{l6.6}. It remains to prove that if $\a$~admits two distinct expansions,
one of them is terminating and the other one is not. \E\Tf we assume that $A$ (resp.~$A'$) are \sq s of \el s of~$D_a$ indexed by $n\in\Na$, and that
the \fw\ holds:
\bea6.42
{}&\sup_{n\ge1}\sum_{l=1}^n A_la^{-l} = \a = \sup_{n\ge1} \sum_{l=1}^n A_l'a^{-l},\\
&\hbox{\E\te s $\ov n\in\Na$ \st}A_{\bar n} \ne A'_{\bar n}. \lb{6.43}
\e
Set $M:=\{n\in\Na: t_n\ne t'_n\}$ where
\beq6.44
t_n:=\sum_{l=1}^n A_la^{-l}, \q t'_n:=\sum_{l=1}^n A'_la^{-l}, \q n\ge1.
\e
By \er{6.43} $M\ne\vn$. \E\te s $\ov n\in\Na$ \st $\ov n\le m$ \fa $m\in M$, in view of Theorem 1.3.30\,(iv). We obtain
\beq6.45
A_{\bar n}a^{-\bar n} + \sup_{n\ge{\bar n}+1} t_n = A'_{\bar n}a^{-\bar n} + \sup_{n\ge{\bar n}+1}t'_n.
\e
Indeed, if $\ov n=1$ then $\sup\limits_{n\ge1} t_n \nde6.42 = \sup\limits_{n\ge2}t_n = \sup\limits_{n\ge2}\bigl(A_1a\mo + \suml_{l=2}^nA_la^{-l}\bigr)
\nad*= A_1a\mo + \sup\limits_{n\ge2}\suml_{l=2}^n A_la^{-l}$. In $\nad*=$ we used, as in the proof of \E\Pr\ \rf{p6.11}, Lemma \rf{l5.2} with $f(y):=
\th_{A_1a\mo}(y)$, $y\in\R$. Clearly, the same holds for $t'_n$ instead of~$t_n$, and \er{6.45} follows.

If $\ov n>1$, we have $\sup\limits_{n\ge1}t_n = \sup\limits_{n\ge\bar n+1}t_n = \sup\limits_{n\ge\bar n+1}\bigl(\suml_{l=1}^{\bar n}A_la^{-l}
+ \suml_{l=\bar n+1}^n A_la^{-l}\bigr) = \suml_{l=1}^{\bar n}A_la^{-l} + \sup\limits_{n\ge\bar n+1}\suml_{l=\bar n+1}^n A_la^{-l}$. The same holds
for~$t_n'$. Note that $\suml_{l=1}^{\bar n-1}A_la^{-l} = \suml_{l=1}^{\bar n-1}A_l' a^{-l}$. Hence, by \cnc ity, we obtain \er{6.45}.

Subtracting $A_n'a^{-n}$ from both sides of \er{6.45} we find
\beq6.46
(A_{\bar n}-A'_{\bar n})a^{-\bar n} + \sup_{n\ge \bar n+1} \sum_{l=\bar n+1}^n A_la^{-l} = \sup_{n\ge \bar n+1} \sum_{l=\bar n+1}^n A'_la^{-l}.
\e

Interchanging $t_n$ and $t'_n$ if necessary, we may assume that $A_{\bar n} > A'_{\bar n}$, that is, $A_{\bar n}\ge A'_{\bar n}\h{+}\h1$. Then the
\LHS\ of \er{6.46} is greater than or equal to $a^{-\bar n}$. \E\oh we have $\suml_{l=\bar n+1}^n A_l'a^{-l} \le \suml_{l=\bar n+1}^n(a-1)a^{-l}
=(a-1)\suml_{l=\bar n+1}^na^{-l} \nde6.24 = a^{-\bar n} - a^{-n} \le a^{-\bar n}$ \fa $n\ge\ov n+1$. Hence $\sup\limits_{n\ge\bar n+1}
\suml_{l=\bar n+1}^n A'_la^{-l} \le a^{-\bar n}$. Thus we obtain from \er{6.46} $a^{-\bar n} \le\break (A_{\bar n}\h {-}\h A'_{\bar n})a^{-\bar n}$,
which implies $A_{\bar n}-A'_{\bar n}= 1$. From \er{6.46} again we infer $a^{-\bar n} + \suml_{l=\bar n+1}^nA_la^{-l}\le a^{-\bar n}$, $n\h{\ge}\h
\ov n+1$. \E\Tf $\suml_{l=\bar n+1}^n A_la^{-l}\h{=}\h0$, $n\ge\ov n+1$. \E\Ip $A_ka^{-k} \le \suml_{l=\bar n+1}^nA_la^{-l}=0$, $n\ge\ov n+1$,
$n\ge\ov n+1$, hence $A_k=0$, $k\ge\ov n+1$. Thus $\{t_n\}_{n\ge 1}$ is a \tm\ \de\ of~$\a$, and we have
\beq6.47
\a = \sum_{l=1}^{\bar n}A_la^{-l}, \q A_{\bar n}=1.
\e

Finally, we infer from \er{6.46}:
$$
a^{-\bar n} = \sup_{n\ge{\bar n}+1}\sum_{l={\bar n}+1}^n A'_la^{-l}.
$$
Since $\suml_{l=\bar n+1}^n(a-1)a^{-l}\nde6.24 = a^{-\bar n}-a^{-n}$, $n\ge\ov n+1$, we obtain
\beq6.48
A'_l = a-1, \q l\ge \ov n+1.
\e
Indeed, we have $A'_l\le a-1$, $l\ge\ov n-1$. Suppose, for \cd ion, that \te s $\ov l \ge \ov n+1$ \st $A'_{\bar l} < a-1$, hence $a-1 \ge A'_{\bar l}
a^{-l}$. We would have $\suml_{l=\bar n+1}^n(a-1)a^{-l}- \suml_{l=\bar n+1}^n A'_la^{-l} \ge a^{-\bar l}$, $n\ge \ov n+1$, and from \er{6.24}
\beq6.49
a^{-\bar l} + \sum_{l=\bar n+1}^n A'_la^{-l} \le a^{-\bar n}, \q n\ge\ov n+1.
\e
However, $\sup\limits_{n\ge\bar n+1}\bigl(a^{-\bar l} + \suml_{l=\bar n+1}^nA'_la^{-l}\bigr) \nad*= a^{-\bar l} + \sup\limits_{n\ge\bar n+1}
\suml_{l=\bar n+1}^n A'_la^{-l} \nde6.47 = a^{-\bar l} + a^{-\bar n}$. \E\oh $\sup\limits_{n\ge\bar n+1}\bigl(a^{-\bar l} +
\suml_{l=\bar n+1}^nA'_la^{-l}\bigr)\le a^{-\bar n}$ by \er{6.49}. A~\cd ion, since $a^{-\bar l}+a^{-\bar n}\not\le a^{-\bar n}$. In $\nad*=$ we used
Lemma \rf{l5.2} as in the proof of \E\Pr\ \rf{p6.11}. \E\Tf \er{6.48} holds.
We observe that the expansion $\{t'_n\}$ is the non-\tm\ expansion of~$\a$ \itd in \E\Pr\ \rf{p6.11}. Hence the proof of Theorem \rf{t6.14} is
complete.
\endproof

\bco6.15
The set $[0,1]:=\{x\in\R: 0\le x\le 1\}$ is \emph{un\ct e}. \E\Ip the field $\R$ is un\ct e.
\eco

\proof
The power set of $\Na$ is un\ct e by Cantor's Theorem 1.4.26. In view of \E\Pr\ 1.4.27\,(ii), it suffices to find an \ti{in\jc} map from $\cP(\Na)$
into $[0,1]$. Let $\Phi:\cP(\Na)\to[0,1]$ be defined by setting
\beq6.50
\Phi(M) := \sup_{n\in\Na} \sum_{l=1}^n A_l10^{-l}, \q M\in\cP(\Na),
\e
where
$$
A_l:=\bca 2 &\hbox{for }l\in M,\\
3 &\hbox{for }l\notin M. \eca
$$
Note that the \RHS\ of \er{6.50} is an \el\ of $[0,1]$ by \er{6.11} and Theorem \rf{t5.1}. Suppose that $\Phi(M)=\Phi(M')$, $M,M'\in\cP(\Na)$. Set
$A'_l:=2$ if $l\in M'$ and $A'_l:=3$ if $l\notin M'$, $t_n:=\suml_{l=1}^N A_l 10^{-l}$, $t'_n:=\suml_{l=1}^N A'_l10^{-l}$, $n\ge1$. Then both
$\zb tn\N$ and $\zb{t'}n\N$ are non-\tm\ \de s to the base~$10$ of $\a:=\Phi(M) = \Phi(M')$. By Theorem \rf{t6.14}, $t_n=t_n'$ \fa $n\in\N$, hence
using \In\ on $n\in\Na$ we find that $A_n=A_n'$ \fa $n\in\Na$. Hence $M=\{m\in\Na: A_m=2\}=\{m\in\Na: A_m'=2\}=M'$, and $\Phi$ is in\jc.
\endproof

Our next goal is to prove the assertion in Remark 4.5.67. To this end and for the remaining part of this section, we first collect some \pp ies of
bounded in\cre\ (resp.\ de\cre) \sq s of \el s of $(\R,\ge)$. These \sq s will be indexed by \el s of~$\Na$ but the results are also valid for
indices in an arbitrary \os\ of \nn s.

\blm6.16
Let $\zb \si n\Na$ denote a bounded in\cre\ \sq\ of \el s of $(\R,\ge)$, and let $\a:=\supl_{n\ge1}\si_n$. The \fw\ holds\dw

\hph i,ii, If $j:\Na\to\Na$ is a strictly in\cre\ map, then the \sq\ $\{\si_{j(k)}\}_{k\in\Na}$ is also bounded and in\cre\ and
\beq6.51
\a = \sup_{k\ge1} \si_{j(k)}.
\e
The \sq\ $\{\si_{j(k)}\}_{k\in\Na}$ is called a \emph{sub\sq} of the \sq\ $\{\si_n\}_{n\in\N}$.

\hph ii,i, If $\b\in\R$ \sf ies $\si_n\le\b$ \fa $n\in\Na$, then
\beq6.52
\a\le\b.
\e

\hph iii,, Let $\g\in\R$ and let $\rho_n:=\si_n+\g$ \fa $n\in\Na$. Then $\zb \rho n\Na$ is bounded, in\cre\ and \sf ies
\beq6.53
\sup_{n\ge1}\rho_n = \a+\g.
\e

\hph iv,, Let $\g\in\R_{>0}$ and let $\rho_n:=\g\cdot \si_n$ \fa $n\in\Na$. Then $\zb\rho n\Na$ is bounded, in\cre\ and \sf ies
\beq6.54
\sup_{n\ge1}\rho_n = \g\cdot\a.
\e
\elm

\bex6.17 \

\hph i,i, Prove Lemma \rf{l6.16}.

\hph ii,, Find and prove the analogues of (i)--(iv) when $\{\si_n\}_{n\ge1}$ is a bounded de\cre\ \sq\ of \el s of~$(\R,\ge)$.
\eex

We now \ch ize the \el s $x\in\R_{>0}\sm \Q_a$ \sf ying $x<1$ which possess  a \ti{periodic} \de\ to base~$a$ (see \E\df\ 4.1.23 and Lemma
4.1.24). Let $\zb Al\Na$ denote the \sq\ defined in \er{6.32}. Let $T\in\Na$ be a period\index{period} of the \sq\ $\zb Al\Na$, then
\beq6.55
A_{l+kT} = A_l \qh{\fa}l,k\in\Na.
\e
\If that $x=\sup\limits_{n\ge1} \suml_{l=1}^n A_la^{-l}$ by \E\df\ \rf{d6.9}. \Mo $\supl_{n\ge1}\suml_{l=1}^n A_la^{-l}\nde6.51 = \break\supl_{k\ge1}
\suml_{l=1}^{kT}A_la^{-l}$. We claim that \fa $k\in\Na$
\beq6.56
\sum_{l=1}^{kT} A_la^{-l} = \biggl(\sum_{l=1}^T A_la^{-l}\biggr) \cdot \biggl(\sum_{m=0}^{k-1}a^{-mT}\biggr).
\e
We use \In\ on $k\in\Na$. Set $M:=\{k\in\Na: \hbox{\er{6.56} holds}\}$. We have $0\in M$ since $\suml_{m=0}^0 a^{-mT} \nad{(2.4.29)} = a^{-0T} =
a^{-0} =1$. Suppose $k\in M$, we show that $k+1\in M$: $\suml_{l=1}^{(k+1)T} A_la^{-l} \nad{(2.4.32)} = \suml_{l=1}^{kT} A_la^{-l} +
\suml_{l=kT+1}^{kT+T}A_la^{-l} \nad{k\in M,(2.4.31)} = \bigl(\suml_{l=1}^T A_la^{-l}\bigr)\cdot\bigl(\suml_{m=0}^{k-1} a^{-mT}\bigr)
+\break \suml_{l=1}^T A_{l+kT}a^{-(l+kT)}$.

Note that $A_{l+kT}=A_l$ by \er{6.55} and $a^{-(l+kT)} = a^{-l}\cdot a^{-kT}$. Indeed, $a^{l+kT}\nad{(2.2.25)}= a^l\cdot a^{kT}$, hence $a^{-(l+kT)}
= (a^l\cdot a^{kT})\mo$ in $\Q_{>0} = (a^l)\mo (a^{kT})\mo = a^{-l}a^{-kT}$. \E\Tf $\suml_{l=1}^T A_{l+T}a^{-(l+T)} \nad{(2.4.39)}= a^{-kT}
\suml_{l=1}^T A_la^{-l}$, $\suml_{l=1}^{(k+1)T}A_la^{-l} = \bigl(\suml_{l=1}^T A_la^{-l}\bigr)\cdot\bigl(\suml_{m=0}^{k-1}a^{-mT}+a^{-kT}\bigr)$, and
\er{6.56} holds for $k:=k+1$. The proof of the claim is complete. \E\Tf we obtain $\supl_{n\ge1}\suml_{l=1}^n A_la^{-l} = \supl_{k\ge1}\bigl(
\suml_{l=1}^T A_la^{-l}\bigr)\cdot \bigl(\suml_{m=0}^{k-1} a^{-mT}\bigr) \nde6.54 = \bigl(\suml_{l=1}^T A_la^{-l}\bigr)\cdot \supl_{k\ge1}
\suml_{m=0}^{k-1} a^{-mT} = \bigl(\suml_{l=1}^T A_la^{-l}\bigr)\cdot \supl_{l\ge0}\suml_{m=0}^l (a^{-T})^m \nad*=
\bigl(\suml_{l=1}^T A_la^{-l}\bigr) \cdot
\frac1{1-a^{-T}} = \bigl(\suml_{l=1}^T A_la^{-l}\bigr)\cdot\frac{a^T}{a^T-1} = \bigl(\suml_{l=1}^T A_la^{T-l}\bigr)\cdot \frac1{a^T-1}$. In
$\nad*=$ we used (4.5.134), \er{3.10} and the fact that $a^{-T}<1$. Indeed, $a^T \nad{(2.2.36)}> 2^T \nad{(2.2.33)}\ge 2>1$.

\If that if $x\in\R_{>0}$, $x<1$, and \te s $\ov l\in\Na$ \st \er{6.55} holds where $\zb Al\Na$ is defined in \er{6.32}, then
\beq6.57
x = \frac B{a^T-1},
\e
where
\beq6.58
B := \sum_{l=1}^T A_la^{T-l} \in\Na.
\e
Note that $B\in\N$ since $A_l,a^{T-l}\in\N$ for $l\in[1,T]$, and $\sum$ is the \cme \ad\ in $\pz1\N$. \Mo $B\ne0$ since otherwise $x=0\cdot
\frac1{a^T-1} = 0$, \cd ing $x\in\R_{>0}$. Since $x<1$, we also have
\beq6.59
B < a^T-1.
\e
\E\Tf $x\in\Q_{>0}$ and by Lemma 4.5.46 \te s \ooo pair $(N,D) \in\Na\t\Na$ \st
\beq6.60
\frac B{a^T-1} = \frac ND, \q \gcd(N,D) = 1, \q N<D.
\e
Note that $N<D$ since $x<1$.

From $BD = N(a^T-1)$, $\gcd(D,N)=1$, and Euclid's Lemma 3.1.30, we infer that $D|a^T-1$. Thus \te s $k\in\Na$ \st $a^T = 1+kD$. We claim that
\beq6.57a
\gcd(a,D)=1.
\e
Suppose for \cd ion that \te s $d\in\Na\sms1$ \st $d|a^T$ and $d|D$. Then $0\nad{(4.1.21),(4.1.23)}= \phi_d(1+dkD) \nad{(4.1.21),(4.1.23)}=1$,
a \cd ion. We have thus proved that a \ti{necessary} \cn\ for a fraction $\frac ND$, $N,D\in\Na$, $N<D$, $\gcd(D,N)=1$, to possess a $T$-periodic
\de\ to base $a\in\Na\sms1$, is \er{6.57a}.

We now show that \cn\ \er{6.57a} is also a \ti{\sft\/} \cn. To this end we shall use the \fw\ result of Number Theory.

\ssk
\ti{Let $a,D\in\Na\sms1$ \sf y \er{6.57a}, and let
\beq6.58a
M:= \{T\in\Na: D|a^T-1\}.
\e
Then \te s $\ov T\in M$ \st}
\beq6.59a
M = \{m\ov T: m\in\Na\}.
\e

\ssk
We postpone the proof of \er{6.59a}.

Now, let $x:=\frac ND$, $N,D\in\Na$, $N<D$ and $\gcd(D,N)=1$. Suppose that \er{6.57a} holds. Let
\beq6.60a
x = \sup_{n\ge1} \sum_{l=1}^n A_la^{-l}
\e
be the \de\ of $x$ to base $a$. Let $\ov T$ be as in \er{6.59a} and let $k\in\Na$ be \st $a^{\bar T}-1 = kD$. Then $(a^{\bar T}-1)\frac ND = kN$.
Hence $a^{\bar T}\cdot\frac ND = kN+\frac ND$. Since $kN\in\N$ and $\frac ND<1$, we have
\beq6.61
\frac ND = \Fr\biggl(a^{\bar T}\cdot \frac ND\biggr).
\e
\Mo $a^{\bar T}\cdot \frac ND = \supl_{n\ge1} \suml_{l=1}^n A_la^{\bar T-l} = \supl_{n\ge\bar T+1}\suml_{l=1}^n A_la^{\bar T-l} = \supl_{n\ge\bar T+1}
\bigl(\suml_{l=1}^{\bar T}A_la^{\bar T-l}+ \suml_{l=\bar T+1}^n A_la^{\bar T-l}\bigr) = \suml_{l=1}^{\bar T}A_la^{\bar T-l}+ \supl_{m\ge1}
\suml_{l=1}^m B_la^{-l}$, where $B_l:=A_{l+\bar T}$, $l\in\Na$. Note that $\suml_{l=1}^{\bar T}A_la^{\bar T-l}\in\N$. \Mo we have $0<x -
\suml_{l=1}^{\bar T}A_la^{-l} < a^{-\bar T}$ by \er{6.34}. Hence $a^{\bar T}x - \suml_{l=1}^{\bar T}A_la^{\bar T-l}<1$. \If that
\beq6.62
\sup_{m\ge1} \sum_{l=1}^m B_la^{-l} = \Fr(a^{\bar T}x).
\e
From \er{6.60a}, \er{6.61}, \er{6.62} and Theorem \rf{t6.14}, we obtain $A_{l+T}=A_l$ \fa $l\in\Na$, and by using \In\ on $k\in\Na$, we arrive at
\er{6.55}.

It remains to prove \er{6.59a}. Observe that $D|a^T-1$ \fs $T\in\Na$ iff $\phi_D(a^T)=1$ where $\phi_D$~is defined in (4.1.20). Indeed, if $D|a^T-1$
\te s $k\in\Na$ \st $a^T-1=kD$. Hence $a^T = 1+kD$, and $\phi_D(1+kD)\nad{(4.1.21),(4.1.23)}= 1$. Conversely, if $\phi_D(a^T)=1$, then \te s $q\in\Na$
\st $a^T = q\cdot B+1$ by (4.1.19). Hence $a^T-1 = qD$, and $D|a^T-1$. Since $D\ge2$, it follows from \E\Pr\ 4.1.18 that $(\N_D,+_D,\cdot_D,0,1)$
(where $\N_D:=\{m\in\N: 0\le m<D\}$) is a \sr\ and that $\phi_D:\pz0\N \to (\N_D,+_D,\cdot_D,0,1)$ is a \sr-\hm sm. \E\Ip $\phi_D(a^T) = T
\mathrel{\lower3pt\hbox{$\stackrel{\cdot_D}\cdot$}}$ the $T$-fold \IT\ of~$\phi_D(a)$ in the monoid $(\N_D,\cdot_D,1)$ by (2.1.45). We recall that
$\N_D^\t$~is the set of invertible \el s of the (finite) \am \ $(\N_D,\cdot_D,1)$ (see \E\df\ 4.3.2) and that $(\N^\t_D,\cdot_D,1)$ is a group (in
particular a \Cm), by Lemma 4.3.6. We also recall that the monoid $(\N_D,+_D,0)$ is a finite abelian (cyclic) group by \E\Pr\ 4.1.15. Thus in view of
\E\df\ 4.3.11, the \sr\ $(\N_D,+_D,\cdot_D,0,1)$ is a (\cmt e) ring (with unity). We recall that the invertible \el s of $(\N_D,\cdot_D,1)$ are
called units of the ring~$\N_D$ (see \E\df\ 4.3.34) and that $(\N_D,\cdot_D,1)^\t = \{x\in\N_D: \gcd(x,D)=1\}$. \Mo since $(\N_D,\cdot_D,1)^\t$ is
a \Cm\ \st $\#(\N_D)=D\ge2$, the \fw\ holds by Lemma 4.1.22: \ti{\E\fe $x\in\Na$, \te s a unique $\ov T\in\Na$ \st $\ov T
\mathrel{\lower3pt\hbox{$\stackrel{\cdot_D}\cdot$}}x=1$, $k\mathrel{\lower3pt\hbox{$\stackrel{\cdot_D}\cdot$}}x\ne1$ for $1\le k<\ov T$, and
$m\ov T\mathrel{\lower3pt\hbox{$\stackrel{\cdot_D}\cdot$}}=1$ \fe $m\in\Na$}. Finally, observe that the \fw\ holds:
\beq6.63
\gcd(a,D)=1 \qh{iff }\gcd(\phi_D(a),D)=1.
\e
Indeed, from $a=q\cdot D+\phi_D(a)$, $\phi_D(a)<D$, \fs $q\in\N$, by (4.1.19), we infer that $b\in\Na\sms1$ divides~$a$ and~$D$ iff it divides
$\phi_D(a)$ and~$D$ (see (3.1.31), (3.1.32)). Note that if $\gcd(a,D)=1$ then $\phi_D(a)\in\N_D$ and $\phi_D(a)\in(\N_D,\cdot_D,1)^\t$. Hence
\er{6.59a} holds, since $\phi_D(a^T) = T\mathrel{\lower3pt\hbox{$\stackrel{\cdot_D}\cdot$}}\phi_D(a)$ \fa $T\in\Na$.

We conclude by observing that $\ov T$ is the minimal period\index{period!minimal} of the \de\ of the \ti{reduced\/} fraction $\frac ND$ (that is, $\gcd(N,D)=1$). Indeed,
if the expansion of~$x$ is \hbox{$T$-periodic}, then $\gcd(a,D)=1$ by \er{6.57a}. Hence $\gcd(\phi_D(a),D)=1$, $\phi_D(a)\in (\N_D,
\mathrel{\lower3pt\hbox{$\stackrel{\cdot_D}\cdot$}},1)^\t$, and there is no $T\in\Na$ less than $\ov T$ \st $T
\mathrel{\lower3pt\hbox{$\stackrel{\cdot_D}\cdot$}}\phi_D(a) = \phi_D(a^T) = 1$,
\ev tly there is no $T<\ov T$ \st $\phi_D(a^T)=1$. We have thus proved

\bth6.18
Let $a\in\Na\sms1$ and let $x\in\R_{>0}\sm\Q_a$ \sf y $x<1$. Then the \fw\ assertions hold.

\hph i,i, Let $\supl_{n\ge1}\suml_{l=1}^n A_la^{-l}$ denote the \de\ of~$x$ to base~$a$ \itd in Lemma \rf{l6.12}. If the expansion is periodic, then
$x\in\Q_{>0}$, and if $x=\frac ND$, $N,D\in\Na$, $N<D$, $\gcd(D,N)=1$, then $\gcd(a,D)=1$.

\hph ii,, If $x=\frac ND$, $N,D\in\Na$, $N<D$, $\gcd(D,N)=1$, and $\gcd(a,D)=1$, then the \de\ to base~$a$ of~$x$ is $\ov T$-periodic, where
$\ov T$~is the smallest \po\ \ig~$T$ \sf ying\dw\ $D$~divides $a^T-1$.
\eth

\brm6.19 \

A \po\ $a$-decimal fraction~$y$ (see \E\df\ \rf{d6.4}) does \ti{not\/} \sf y $\gcd(a,D)=1$, $a\in\Na\sms1$, whenever $y=\frac ND$, $N<D$,
$N,D\in\Na$, $\gcd(D,N)=1$. Indeed, suppose for \cd ion that $\gcd(a,D)=1$, then we have $\phi_D(a^{m\bar T})=1$ \fs $\ov T\in\Na$ and \fa
$m\in\Na$. \E\oh from $\frac ND=\frac k{a^n}$, $k,n\in\Na$, we infer $kD=Na^n$, hence $D|Na^n$, and by Euclid's Lemma, $D|a^n$. \E\Tf $D|a^{n\bar T}$,
hence $\phi_D(a^{n\bar T})=0$, a~\cd ion. Note that if $A_l=a-1$ \fa $l\in\Na$, then $\supl_{n\ge1}\suml_{l=1}^n A_la^{-l}=1\notin\zo0,1 $.
\erm

In the proof of \er{6.62} we used (and proved) a special case of the \fw\ result, which we rewrite as a lemma.

\blm6.19
Let $x\in\R_{>0}$ be \st $x=\Fr(x)$, and let $\supl_{n\ge1}\suml_{l=1}^n A_la^{-l}$ be its \de\ to base $a\in\Na\sms1$ \itd in Lemma \rf{l6.12}.
Then \fa $m\in\Na$ the \fw\ holds\dw
\beq6.68
i(a^m)\cdot x \nad*= \sum_{l=1}^m A_la^{m-l} + \sup_{n\ge1}\sum_{l=1}^n B_la^{-l},
\e
where
\beq6.69
B_l:= A_{l+m} \qh{\fa}l\ge1.
\e
\Mo
\bea6.70
\Int(i(a^m)\cdot x) &\nad*= i\biggl(\sum_{l=1}^m A_la^{m-l}\biggr),\\
\Fr(i(a^m)\cdot x) &= \sup_{n\ge1} \sum_{l=1}^n B_la^{-l}. \lb{6.71}
\e
\elm

In $\nad*=$ we identify an \el\ $n$ of~$\Na$ with $i(n)$, and $0\in\N$ with $0\in\R$. In the proof of \er{6.62} we had $x\in\Q_{>0}$.

\bex6.20
Prove Lemma \rf{l6.19}.
\eex

In view of Lemma \rf{l6.19}, Theorem \rf{t6.18}, and recalling the \df\ of \ti{eventual periodicity} given in Remark 4.5.67, we obtain

\blm6.21
Let $a\in\Na\sms1$, and let $x\in\R_{>0}\sm\Q_a$ be \st $x=\Fr(x)$. If $x$ possesses an \tb{eventually periodic} \de\ to base~$a$, then $x\in\Q_{>0}$.\index{eventually periodic \de}
\elm

\proof
Since the expansion is eventually periodic, \te s $m\in\Na$ \st the \sq\ $\zb Bl\Na$ defined in \er{6.69} is periodic. By Theorem \rf{t6.18}
$\supl_{n\ge1}\suml_{l=1}^n B_la^{-l}\in\Q_{\ge0}$, and by \er{6.70} $\suml_{l=1}^m A_la^{m-l}\in\Q_{\ge0}$, hence by \er{6.68} $i(a^m)\cdot x
\in\Q_{\ge0}$. \E\Tf $x=i(a^{-m})\cdot i(a^m)\cdot x\in\Q_{\ge0} \cap \R_{>0}=\Q_{>0}$.
\endproof

As a direct con\sq\ of Lemma \rf{l6.19} we find that the \de\ of $x\in\Q_{>0}$, $x<1$, to base~$a$ is eventually periodic iff \te s $m\in\Na$ \st
$\Fr(a^mx)$ has a periodic expansion. Indeed, the expansion of~$x$ is eventually periodic iff \te s $m\in\Na$ \st the \sq\ $\zb Bl\Na$ defined in
\er{6.69} is periodic, and iff $\Fr(a^mx)$ has a periodic expansion \fs $m\in\Na$ in view of \er{6.71}. The next lemma shows that this \cn\
is \sf ied \fe $x\in \Q_{>0}\sm\Q_a$, $x<1$.

\newpage
\blm6.22
Let $a\in\Na\sms1$ and let $x\in\Q_{>0}\sm\Q_a$, $x<1$, be \st
\bga6.72
x=\frac ND, \q N,D\in\Na, \q \gcd(D,N)=1,\\
\gcd(a,D)>1. \lb{6.73}
\e
Then \te s $m\in\Na$ \st
\bga6.74
\Fr(a^mx) = \frac PQ, \q P\in\Na,\ Q\in\Na\sms{1,2},\\
\gcd(P,Q)=1 \qh{and } \gcd(a,Q)=1. \lb{6.75}
\e
\elm

We only prove the case $a:=10$. The proof of the general case is left as an exercise.

\proof \

``\ti{Case} $a:=10$'': The \Pn s dividing 10 are 2~and~5. Let $\a,\b\in\N$ be \st $2^\a | D$, $2^{\a+1}\nmid D$, and $5^\b|D$, $5^{\b+1}\nmid D$.
Set $R:=2^\a 5^\b$. Since $\gcd(a,D)>1$, we have $\a+\b\ge1$. Note that $R|D$ by Lemma 3.1.29, since $2^\a|D$, $5^\b|D$ and $\gcd(2^\a,5^\b)=1$.
Set $Q:=\frac DR \in\Na$. Observe that $\gcd(a,Q)=1$. Indeed, if $\gcd(a,Q)>1$, then either $2|Q$, or $5|Q$, or both. Since $D=QR$, we would have
$2^{\a+1}|D$, or $5^{\b+1}|D$, or both, \cd ing the \df\ of $\a$ and~$\b$. We obtain
\beq6.76
x = \frac N{2^\a5^\b Q} \qh{where }\gcd(a,Q)=1.
\e
\Mo if $\a>0$ then $2\nmid N$, otherwise $2|N$ and $2|D$ \cd ing $\gcd(D,N)=1$. Similarly, $5\nmid N$. Set
\beq6.77
m:= \max(\a,\b).
\e
Then
\beq6.78
a^mx = \frac{2^\g 5^\d N}Q \qh{where }\g,\d\in\N,
\e
since $\g:=m-\a$ and $\d:=m-\b$.

Note that $Q>1$, otherwise $a^mx = 2^\g5^\d N$, hence $x = \frac{2^\g 5^\d N}{a^m}$, and $x$~would belong to $\Q_a$, \cd ing the \as s of the lemma.
In view of the division algorithm (Theorem 2.1.38) and (2.2.8), \te\ $k,P\in\N$ \sf ying $2^\g 5^\d N = kQ+P$ and $P<Q$. From \er{6.78} we infer
\beq6.79
a^m x = k+\frac PQ.
\e
Since $k\in\N$ and $\frac PQ\in\zo0,1 _\Q$, we obtain from \E\Pr\ 4.5.61 that
\bea6.80
\frac PQ&= \Fr(a^mx),\\
P &= 2^\g5^\d N - kQ. \lb{6.81}
\e
We claim that $\gcd(P,Q)=1$. Indeed, suppose for \cd ion that \te s a~\Pn\ $p$ dividing $P$~and~$Q$, then $p|2^\g5^\d N$ by \er{6.81}. Since $p|Q$,
$p\ne2$ and $p\ne5$. Hence $\gcd(p,2^\g5^\d)=1$. By Euclid's Lemma (Lemma 3.1.30) we obtain $p|N$. Since $p|Q$, and $D=QR$, we obtain $p|D$.
A~\cd ion, since $\gcd(D,N)=1$. \E\csq, $\Fr(a^mx)=\frac PQ$ \sf ies $\gcd(P,Q)=1$ and $\gcd(a,Q)=1$.
\endproof

\If from \er{6.71} and Theorem \rf{t6.18} that $\{B_l\}_{l\ge1}$ is periodic, and from the discussion preceding Lemma \rf{l6.22} that $\{A_l\}
_{l\ge1}$ is eventually periodic.

\bex6.23
Prove Lemma \rf{l6.22} for $a=2^k5^l$, $k,l\in\N$, $k+l\ge1$.

Hint: Show first that $x=\frac N{2^\a5^\b Q}$ for $Q\in\Na$ \fs $\a,\b\in\N$ \sf ying $\a+\b\ge1$, $2\nmid N$, $5\nmid N$ and $(a,Q)=1$. Set
\beag
\g:=\biggl[\frac \a k\biggr]\ &\hbox{if }k|\a \hbox{ and }\biggl[\frac \a k\biggr]+1 \hbox{ if }k\nmid \a,\\
\d:=\biggl[\frac \b l\biggr]\ &\hbox{if }l|\b \hbox{ and }\biggl[\frac \b l\biggr]+1 \hbox{ if }l\nmid \b.
\e
Set $m:=\max(\g,\d)$. Show that \er{6.74}, \er{6.75} hold. Then prove the general case.
\eex

\bxs6.24
Let $a:=10$ and $x\in\Q_{>0}\sm \Q_{10}$, $x<1$ be \st \er{6.72} holds.

\noi \hph I,I, $\gcd(10,D)=1$. Let $\ov T$ be as in Theorem \rf{t6.18}.

\vskip2pt
\hph i,ii, $x=\frac ND=\frac13$: \ $3|10^1-1$; $\ov T=1$, $x = 0{,}3\,3\,3\,\dots$

\vskip2pt
\hph ii,i, $x=\frac23$: \ $3|10^1-1$; $\ov T=1$, $x= 0{,}6\,6\,6\,\dots$

\vskip2pt
\hph iii,, $x=\frac17$: \ $7|10^6-1$, $7\nmid10^n-1$ for $1\le n<6$, $\ov T=6$, $x=0{,}142857\,142857\,\dots$

\vskip2pt
\hph iv,, $x=\frac27$: \ $x=0{,}285714\,285714\,\dots$

\vskip2pt
\hph v,i, $x=\frac1{11}$: \ $11|10^2-1$; $\ov T=2$, $x=0{,}09\,09\,09\,\dots$

\vskip2pt
\hph vi,, $x=\frac9{11}$: \ $11|10^2-1$; $\ov T=2$, $x=0{,}81\,81\,81\,\dots$

\vskip2pt
\noi \hph II,, $\gcd(10,D)>1$.

\vskip2pt
\hph i,ii, $x=\frac16 = \frac1{2^1\cdot 5^0\cdot 3}$; $m=1$, $x=0{,}16666\dots$

\vskip2pt
\hph ii,i, $x=\frac56 = \frac5{2^1\cdot 5^0\cdot 3}$; $m=1$, $x=0{,}83333\dots$

\vskip2pt
\hph iii,, $x=\frac1{60} = \frac1{2^2\cdot 5^1\cdot 3}$; $m=2$, $x=0{,}01666\dots$

\vskip2pt
\hph iv,, $x=\frac1{3000} = \frac1{2^3\cdot 5^3\cdot 3}$; $m=3$, $x=0{,}000333\dots$
\exs

We now return to the study of the map $\psi_x:\Q \to\R_{>0}$ \itd\ in \E\Pr\ \rf{p5.29} \er{5.60}. From now on we assume
\beq6.82
x>1,
\e
and we omit the subscript $x$.

For convenience we recall some \pp ies of the map (\f) $\psi$.
\bea6.83
{}&\hbox{Positivity: } \psi(r)>0 \hbox{ \fa} r\in\Q.\\
&\hbox{Strict in\cre ness: } \psi(r)<\psi(s) \hbox{ \fa} r,s\in\Q,\ r<s. \lb{6.84}\\
&\hbox{\E\hm sm: } \psi(0)=1 \hbox{ and }\psi(r+s)= \psi(r)\cdot \psi(s) \hbox{ \fa} r,s\in\Q. \lb{6.85}\\
&\hbox{Order-continuity: } \sup_{s<r}\psi(s) = \psi(r) = \inf_{s>r}\psi(s). \lb{6.86}
\e

Our next goal is to \es\ the \ti{convexity} of the \f~$\psi$.

\bdf6.25
A \f\ $F:\Q\to\R$ (resp.\ $F:\R\to\R$) is called \ti{convex} (resp.\ \ti{concave}) if\index{convex, concave \f}
\beq6.87
F((1-t)r+ts) \le \hbox{(resp.\ $\ge$)}\, (1-t)F(r)+tF(s)
\e
\fa\ $r,s\in\Q$ (resp.~$\R$) and all $t\in\Q$ (resp.~$\R$) \sf ying $0\le t\le 1$.
\edf

In a first step we show that the \f\ $\psi$ \sf ies \E\df\ \rf{d6.25} for \ti{dyadic}~$t$'s (see \E\df\ \rf{d6.1}). In a second step we use the
order-continuity of~$\psi$ and Lemma \rf{l6.12} to prove the convexity of~$\psi$.

\ssk
Clearly, every \f\ $F:\Q \hbox{\ (resp.\ $\R$)}\to\R$ \sf ies \er{6.87} for $t=0$ and $t=1$. We first show that $\psi$~\sf ies \er{6.87} for $t:=\frac
12$. We use the in\et y
\beq6.88
a\cdot b \le \tfrac12 a^2+\tfrac12 b^2 \qh{\fa} a,b\in\R_{\ge0}.
\e
In\et y \er{6.88} follows from
$$
a^2-2ab+b^2 = (a-b)^2\ge0 , \q a,b\in\R.
$$
In\et y \er{6.88} is a special case of the useful in\et y
\beq6.89
a\cdot b \le \frac \ve 2\,a^2 + \frac1{2\ve}\,b^2, \q a,b,\ve\in\R_{>0}.
\e
(Replace $a,b$ by $\sqrt\ve\, a$, $\frac1{\sqrt\ve}b$.)

\blm6.26
Let $G:\Q \hbox{ $($resp.\ $\R)$}\to\R_{>0}$ \sf y
\beq6.90
G(r+ s) = G(r)\cdot G(s) \qh{\fa $r,s\in\Q$ $($resp.\ $\R)$.}
\e
Then
\beq6.91
G((1-t)r+ts) \le (1-t)G(r) + tG(s)
\e
\fa $r,s\in\Q$ $($resp.\ $\R)$ and all $t\in\Q_2\cap (0,1)$.
\elm

\proof[Proof \rm(inspired by {\cite[p.~5]{Gamma}})]
Set $A_n:=\{k\in\Na: 1\le k\le 2^n\}$, $n\in\Na$. We have to show that \er{6.91} holds \fa $r,s\in\Q$ (resp.~$\R$), all $t:=\frac k{2^n}$, where
$k\in A_n$, $n\in\Na$.

\ssk
``\ti{Case} $n:=1$'': Let $r,s\in\Q$ (resp.\ $\R$). Then $G(\frac12r+\frac12s)\nde6.90 = G(\frac12r)\cdot G(\frac12s) \nde6.88 \le
\frac12(G(\frac12r))^2 + \frac12(G(\frac12s))^2 = \frac12G(\frac12r)\cdot G(\frac12r) + \frac12G(\frac12s)\cdot G(\frac12s) \nde6.90 =
\frac12G(\frac12r+\frac12r) + \frac12G(\frac12s+\frac12s) = \frac12G(r)+\frac12G(s)$.

\ssk
``\ti{Case} $n>1$'': A \f\ $G$ \sf ying $G(\frac12(x_1+x_2))\le \frac12(G(x_1)+G(x_2))$, $x_1,x_2\in\R$ is called \ti{weakly convex} in~\cite{Gamma},
and this \pp y is shown to be generalized in the \fw\ way:

Let $x:A_n\to\Q$ (resp.\ $\R$) denote an arbitrary finite \sq\ of \el s of~$\Q$ (resp.~$\R$). Then we have
\beq6.92
G\biggl(\frac1{2^n}\sum_{l=1}^{2^n} x_l\biggr) \le \frac1{2^n} \sum_{l=1}^{2^n} G(x_l), \q n\in\Na.
\e
We proceed by \In\ on $n\in\Na$. The case $n:=1$ follows from ``Case $n:=1$''. We assume that \er{6.92} holds for $n\in\Na$ and show that \er{6.92}
holds for $n:=n+1$. We have
\bmlg
G\Bgg(\frac1{2^{n+1}}\sum_{l=1}^{2^{n+1}}x_l) = G\Bgg(\frac12\,\frac1{2^n}\sum_{l=1}^{2^n}x_l + \frac12\,\frac1{2^n} \sum_{l=2^n+1}^{2^{n+1}}x_l)\\
\nad{t=\frac12} \le \frac12\biggl( G\Bgg(\frac1{2^n}\sum_{l=1}^{2^n}x_l) + \frac12\,G\Bgg(\frac1{2^n}\sum_{l=1}^{2^n} x_{2^n+l})\biggr) \\
\nde6.92 \le \frac12\,\frac1{2^n} \sum_{l=1}^{2^n} G(x_l) + \frac12\,\frac1{2^n}\sum_{l=1}^{2^n} G(x_{2^n+l})
\nad{\rm(2.4.31),(2.4.32)} = \frac1{2^{n+1}}\sum_{l=1}^{2^{n+1}} G(x_l).
\e
Hence \er{6.92} holds for $n:=n+1$. \csq, \er{6.92} holds \fa $n\in\Na$.

Finally, given $r,s\in\Q$ (resp.\ $\R$), $n\in\Na$, $k\in A_n$, setting $x_l:=\frac1{2^n}r$, $1\le l\le k<2^n$, and $x_l:=\frac1{2^n}s$, $k+1\le l
\le2^n$, we find from \er{6.92}:
\beq6.93
G\Bgg(\frac k{2^n}\,r + \frac{2^n-k}{2^n}\,s) \le \frac k{2^n}\,G(n) + \frac{2^n-k}{2^n}\,G(s).
\e
Note that if $\wt t\in\Q_2\cap(0,1)$, then \te\ $n\in\Na$ and $k\in A_n\sms{2^n}$, \st $\wt t=\frac k{2^n}$. \Mo $1-\wt t = \frac{2^n-k}{2^n}\in
\Q_2\cap(0,1)$. Setting $1-t:=\wt t$, we have $t,1-t\in\Q_2\cap(0,1)$ and obtain \er{6.91}.
\endproof

Our next goal is to prove that the \f~$\psi$ \sf ies \er{6.87} \fa $r,s\in\Q$ and all $t\in(0,1)_\Q$. Since $\psi$ \sf ies \er{6.85}, it \sf ies \cn\
\er{6.90} \fa $r,s\in\Q$, hence, by Lemma \rf{l6.26}, it also \sf ies \er{6.91} \fa $r,s\in\Q$ and $t\in\Q_2\cap (0,1)$. By Lemma \rf{l6.12}, given
$t\in(0,1)_\Q$, \te s an in\cre\ \sq\ $\zb tn\Na$ of \el s of $\Q_2\cap(0,1)$ \st $\supl_{n\ge1}t_n=t$. From \er{6.91} we infer \fa $r,s\in\Q$:
\beq6.94
\psi(r+t_n(s-r)) \le \psi(r) + t_n(\psi(s)-\psi(r)), \q n\in\Na.
\e
We assume
\beq6.95
r<s.
\e

The case $r=s$ is trivial and the case $s<r$ will be treated later on. Since $\psi$ is strictly in\cre, $\psi(s)-\psi(r)\in\R_{>0}$. \E\Tf
$\supl_{n\ge1} t_n(\psi(s)-\psi(r)) \nde6.54 = t(\psi(s)-\psi(r))$, $\supl_{n\ge1} \psi(r)+t_n(\psi(s)-\psi(r)) \nde6.53 = \psi(r) + \supl_{n\ge1}
t_n(\psi(s)-\psi(r)) = \psi(r) + t(\psi(s)-\psi(r)) = (1-t)\psi(r) + t\psi(s)$. \If that $(1-t)\psi(r)+t\psi(s)$ is an \ub\ for
$\{\psi(r)+t_n(\psi(s)-\psi(r)): n\in\Na\}$, hence from \er{6.94} we infer
\beq6.96
\psi(r+t_n(s-r)) \le \psi(r) + t(\psi(s)-\psi(r)), \q n\in\Na.
\e
Note that $r+t(s-r) = \supl_{n\ge1} r+t_n(s-r)$ by \er{6.54}, \er{6.53}.

We  recall that if $\zb \si n\Na$ is a bounded in\cre\ \sq\ of \el s of~$\R$ and $\vf:\R\to\R_{>0}$ is in\cre, then $\si_m\le \supl_{n\ge1}\si_n$
\fa $m\in\Na$ and $\vf(\si_m) \le \vf\bigl(\supl_{n\ge1}\si_n\bigr)$, $m\ge1$. Hence by \er{6.52} $\supl_{m\ge1}\vf(\si_m)\le \vf\bigl(
\supl_{n\ge1}\si_n\bigr)$. \E\Ip if $\vf:=\psi :\Q\to\R_{>0}$, and $\si_n:=r+t_n(s-r)$, we obtain $\supl_{m\ge1}\psi(r+t_m(s-r)) \le \psi(r+t(s-r))$.
\Mo $\supl_{m\ge1}\psi(r+t_m(s-r)) \le \psi(r+t(s-r))$. \Mo in order to \es\ \er{6.87} from \er{6.56} we need $\supl_{m\ge1}\psi(r+t_m(s-r))\nad!=
\psi(r+t(s-r))$. Thanks to \pp y \er{6.86}, \et y holds as a con\sq\ of the \fw\ lemma.

\newpage
\blm6.27
Let $\vf:\Q\to\R$ be in\cre\ and \sf y \fa $r\in\Q$\dw
\beq6.97
\sup_{s<r,s\in\Q}\vf(s) = \vf(r).
\e
\elm

Let $\{s_n\}_{n\ge1}$ be a bounded in\cre\ \sq\ of \el s of~$\Q$, \sf ying
\beq6.98
\sup_{n\ge1} s_n = \a\in\Q.
\e
Then we have
\beq6.99
\vf(\a) = \sup_{n\ge1}\vf(s_n).
\e

\ssk
Before proving \er{6.99} we give a \chz\ of the supremum of a subset of $(\Q,\ge)$. See Lemma \rf{l5.38} for $(\R,\ge)$ instead of $(\Q,\ge)$.
Let $B$~be a \nss\ of~$\Q$. Then $\xi\in\Q$ is the supremum of~$B$ iff the \fw\ \cn s hold:
\bea6.100
{}&\xi\ge\eta \hbox{ \fa}\eta\in B.\\
&\hbox{\E\fe $\ve\in\Q_{>0}$ \te s $\b\in B$ \st }\b+\ve>\xi. \lb{6.101}
\e
Since the proof of Lemma \rf{l5.38} was left as an exercise, we give it now. \E\cn\ \er{6.100} means that $\xi$~is an \ub\ for the set~$B$. We first
show that if $\xi'\in\Q$ \sf ies $\xi'\ge\eta$ \fa $\eta\in B$, that is, $\xi'$~is an \ub\ of~$B$, then $\xi\le\xi'$, that is, $\xi$~is an \ub\
for~$B$. Indeed, \fe $\ve\in\Q_{>0}$ \te s $\b\in B$ \st $\b+\ve>\xi$, hence $\xi'+\ve\nde3.10 \ge \b+\ve>\xi$, and $\xi'+\ve>\xi$ by (1.3.12).
Thus $\xi'+\ve>\xi$ \fa $\ve\in\Q_{>0}$. Suppose, for \cd ion, that $\xi>\xi'$. Set $\a:=\xi-\xi'>0$. Then we obtain $\xi>\a>0$ \fa $\xi>0$, \Ip
for $\ve:=\frac \a2$. A~\cd ion, since $\a=2\ve$ and $\ve\not> 2\ve$. \E\Tf $\xi\le\xi'$. Observe that we used the fact that $(\Q,\ge)$ is \Ar.
We next prove that if $\xi$~is the least \ub\ for~$B$, then \er{6.101} holds. Suppose, for \cd ion, that \er{6.101} does not hold. Then \te s
$\ov\ve>0$ \st \fa $\b\in B$ we have $\b+\ov \ve\le \xi$. Then $\b\le\b+\frac{\ov \ve}2 \le \xi-\frac{\ov\ve}2$ \fa $\b\in B$. \csq, $\xi-\frac
{\ov\ve}2$ is an \ub\ for~$B$ which is less than~$\xi$. A~\cd ion. Summarizing, we proved

\blm6.28
Let $B$ be a \nss\ of~$\R$. Then $\xi\in\Q$ is the supremum of~$B$ if \er{6.100} and \er{6.101} hold.
\elm

Note that Lemma \rf{l5.38} can be proved in the same fashion.

\proof[Proof of Lemma \rm\rf{l6.27}.]
We have to show that supremum $\{\vf(s_n): n\in\Na\}=\vf(\a)$ in $(\Q,\ge)$. By \er{6.98} $\vf(s_n)\le\vf(\a)$ \fa $n\in\Na$, since $\vf$~is in\cre.
Thus \er{6.100} holds with $\xi:=\vf(a)$. In view of Lemma \rf{l6.28} it remains to show that \er{6.101} holds with $B:=\{\vf(s_n):n\in\Na\}$ and
$\xi:=\vf(\a)$, that is, \fe $\ve\in\Q_{>0}$ \te s $n\in\Na$ \st $\vf(s_n)+\ve>\vf(\a)$. By \er{6.98} and Lemma \rf{l6.28}, \fe $\ve\in\Q_{>0}$ \te s
$t\in\Q$, $t<\a$, \st
\beq6.102
\vf(t)+\ve > \vf(\a).
\e
We claim that \fe $t<\a$ \te s $n\in\Na$ \st $s_n>t$. Suppose, for \cd ion, that \te s $t<\a$ \st \fa $n\in\Na$ there is no $n\in\N$ \st $s_n>t$,
\ev tly that $s_n\le t$ \fa $n\in\Na$. In this case $\a=\supl_{n\ge1}s_n \le t$ since $t$ is an \ub\ for $\{s_n:n\in\Na\}$ in~$\Q$. A~\cd ion, since
$t<\a$. \csq, \fe $t<\a$ \te s $n\in\Na$ \st $s_n>t$, hence \fe $\ve\in\Q_{>0}$ \te s $n\in\Na$ \st
\beq6.103
\vf(s_n)+\ve \ge \vf(t)+\ve > \vf(\a).
\e
\E\Tf \er{6.101} holds and by Lemma \rf{l6.28} \er{6.93} holds.
\endproof

We now return to the proof of Lemma \rf{l6.26}. In view of the discussion preceding Lemma \rf{l6.27} we know that \er{6.96} holds, we have
$\supl_{n\ge1} \psi(r+t_n(s-r)) \le \psi(r)+t(\psi(r)-\psi(s)$, for $r,s\in\Q$, $r<s$. Since $r+t(s-r) = \supl_{n\ge1} r+t_n(s-r)$ and \er{6.86}
holds, we infer from Lemma \rf{l6.27} that $\supl_{n\ge1}\psi(r+t_n(s-r)) = \psi(r+t(s-r))$, which implies \er{6.87} with $F:=\psi$, and $r<s$.
We now remove \cn\ \er{6.82}. Let $r>s$, $r,s\in\Q$. Set $\wt s:=r$, $\wt r:=s$ and $\wt t:=1-t$, $t\in(0,1)$. Note that $\wt s>\wt r$ and $\wt t\in
(0,1)$. In view of \er{6.87} we obtain $\psi((1-\wt t)\wt r + \wt t\wt s) \le (1-\wt t)\psi(\wt r) + \wt t\psi(\wt s)$. Since $1-\wt t=t$, $\wt t =
1-t$, we obtain $\psi((1-t)r+ts) \le (1-t)\psi(r)+t\psi(s)$, which is \er{6.87}.

Summarizing, we have proved

\bpr6.29
The map $\psi_x:\Q\to\R$, $x\in\R_{>1}$, \itd\ in \E\Pr\ \rf{p5.29} is \emph{convex} in the sense of \E\df\ \rf{d6.25}, where $t\in[0,1]_\Q$.
\epr

\brm6.30 \

\hph i,i,
In \cite{Gamma} a proof of \E\Pr\ \rf{p6.29} is given, which does not use \pp y \er{6.86} and does not invoke Lemma \rf{l6.12}. This proof is based
on the observation that the \fw\ \pp y \er{6.104} implies convexity: Let $G:\Q\to\R$ \sf y
\beq6.104
G\Bgg(\frac1n \sum_{l=1}^n x_l) \le \frac1n \sum_{l=1}^n G(x_l)
\e
\fa finite \sq s $\{x_l\}_{l\in[1,n]}$ of \el s of~$\Q$ \fa $n\in\Na$. Then \er{6.91} holds \fa $r,s\in\Q$ and all $t\in[0,1]_\Q$. Indeed, given
$r,s\in\Q$, $t:=\frac ab$, $a,b\in\Na$, $a<b$, setting $n:=b$, $x_l:=r$ \fa $l\in[a+1,b]$, $x_l:=s$ \fa $l\in[1,a]$, we find $G((1-\frac ab)r+\frac ab
s) \le \frac1b \suml_{l=a+1}^b G(r) + \frac1b\suml_{l=1}^a G(s) = (1-\frac ab)G(r) + \frac abG(s)$. Hence \er{6.91} holds \fa $r,s\in\Q$ and $t\in
(0,1)_\Q$. Suppose now that \er{6.104} holds for $n:=2$, that is,
\beq6.105
G(\tfrac12 r+\tfrac12 s) \le \tfrac12G(r) + \tfrac12 G(s) \qh{\fa $r,s\in\Q$.}
\e
Then, one can prove that if \er{6.104} holds for~$n$, then it holds for~$2n$. Hence by \In\ \er{6.104} holds \fa $n$ of the form $2^k$, $k\in\Na$.
The proof goes as follows. Given $\{x_l\}_{l\in[1,2n]}$, set $r:=\frac1n \suml_{l=1}^n x_l$ and $s:=\frac1n \suml_{l=n+1}^{2n}x_l$, apply \er{6.105}.
Finally, we observe that \er{6.104} holds for $n+1$, then it also holds for~$n$. Indeed, given $\{x_l\}_{l\in[1,n]}\sbs\Q$, we may adjoin the \el\
$x_{n+1}:=\frac1n \suml_{l=1}^n x_l$ to the \sq\ $\{x_l\}_{l\in[1,n]}$ and apply \er{6.104} for $n:=n+1$ to this new \sq. We obtain
$G(x_{n+1}) = G\bg(\frac{nx_{n+1}+x_{n+1}}{n+1}) = G\bg(\frac1{n+1} \suml_{l=1}^{n+1}x_l) \le \frac1{n+1}\suml_{l=1}^{n+1} G(x_l) = \frac1{n+1}
\suml_{l=1}^n G(x_l) + \frac1{n+1}G(x_{n+1})$. Hence $\frac n{n+1}G\bg(\frac1n \suml_{l=1}^n x_l) \le \frac1{n+1}\suml_{l=1}^n G(x_l)$, which implies
\er{6.104} for $n:=n$. We conclude by using the \fw\ lemma.
\erm

\blm6.31
Let $A$ be a nonempty bounded subset of $\N$ with $\#(A)\ge2$. Let $a$ denote its least \el\ and $b$ its greatest \el. Suppose that \fa $n\in A$,
$n>a$, we have $n-1\in A$. Then $A=\{n\in\N: a\le n\le b\}$.
\elm

\bex6.32
Prove Lemma \rf{l6.31}.
\eex

We have thus proved that if $G:\Q\to\R$ \sf ies \er{6.105} then \er{6.104} holds \fa $m\ge2$. Indeed, given $m>2$, \te s $k\in\Na$ \st $2^k\ge m$,
hence \er{6.104} holds for $n=2^k$ as well as \fa $n<2^k$, in particular for $n:=m$. \E\Tf if $G$ \sf ies \er{6.90}, then $G$~is convex in the sense
of \E\df\ \rf{d6.25}.

\hph ii,, The reason why we gave another proof of the convexity of~$\psi$ is that in the sequel we shall present a proof of the convexity of a map
$G:\R\to\R$ where $t\in[0,1]_\R$, which is almost verbatim the same as for the case $t\in[0,1]_\Q$.
\ssk

The next \Pr\ allows us to strengthen \pp y \er{6.86} of the \f~$\psi$.

\bpr6.32
Let $F:\Q (\hbox{resp.\ }\R) \to \R$ \sf y $F((1-t)r+ts) \le (1-t)F(r)+tF(s)$ \fa $t\in[0,1]_\Q$ $($resp.\ $[0,1]_\R)$ and all $r,s\in\Q$
$($resp.~$\R)$. Then \fa $a,b,c\in\Q$ $($resp.~$\R)$ \sf ying $a<b<c$, we have\dw
\beq6.106
\frac{F(b)-F(a)}{b-a} \le \frac{F(c)-F(a)}{c-a} \le \frac{F(c)-F(b)}{c-b}.
\e
\epr

\proof
Set $t:=\frac{b-a}{c-a}\in(0,1)$. Then $1-t=\frac{c-b}{c-a}$, and $b=(1-t)a+bc$. Indeed,\break $(1-t)a+tc = \frac{c-b}{c-a}a + \frac{b-a}{c-a}c =
\frac{bc-ba}{c-a} = b$. Hence it follows from the \as\ that\break $F(b) = F((1-t)a+tc) \le (1-t)F(a)+tF(c)$. \E\Tf $F(b)-F(a) \le t(F(c)-F(a)) =
\frac{b-a}{c-a}(F(c)-F(a))$, which implies the first in\et y in \er{6.106}. \E\oh $F(c)-F(b) = (1-t)F(c) + t(F(c)-F(b)) \ge (1-t)(F(c)-F(a)) =
\frac{c-b}{c-a} (F(c)-F(a))$, which implies the second in\et y in \er{6.106}.
\endproof

\bco6.33
Let the map $\psi_x$ be as in \E\Pr\ \rf{p6.29}. Let $s,r,u\in\Q$ \sf y ${s<r<u}$. Then \fa $h\in\Q_{>0}$\dw
\beq6.107
0< \psi_x(r)-\psi_x(s) \le \frac{\psi_x(u+h)-\psi_x(u)}h\cdot (r-s).
\e
\E\Ip
\beq6.108
\psi_x(r)-\psi_x(s) \le (\psi_x(u+1)-\psi_x(u))\cdot(r-s).
\e
\eco

\proof
We have $s<r<u<u+h$. The first in\et y in \er{6.107} follows from the strict in\cre ness of~$\psi_x$. From \er{6.106} with $a:=s$, $b:=r$ and $c:=u$
we infer $\frac{\psi_x(r)-\psi_x(s)}{r-s} \le \frac{\psi_x(u)-\psi_x(r)}{u-r}$, and from \er{6.106} with $a:=r$, $b:=u$ and $c:=u+h$ we obtain
$\frac{\psi_x(u)-\psi_x(r)}{u-r} \le \frac{\psi_x(u+h)-\psi_x(u)}{u+h-u} = \frac{\psi_x(u+h)-\psi_x(u)}h$.
\endproof

\bex6.34 \

\hph i,i, Let $a,b\in\Q$, $a<b$, and let $\vf:(a,b)\to\R$ be strictly in\cre\ and \sf y
\beq6.109
\vf(r)-\vf(s) \le C(r-s), \q u< s<r< b,
\e
\fs $C\in\R_{>0}$. Show that the \fw\ holds:
\beq6.110
\sup_{a<s<r<b}\vf(s) = \vf(r) = \inf_{a<r<u<b}\vf(u) \qh{\fa $r\in(a,b)$.}
\e

\hph ii,, Let $a:=-1$, $b:=1$ and $\vf(r):=0$ for $-1<r\le0$, $\vf(r):=\sqrt r$ for $0<r<1$. Show that \er{6.110} holds, and that there is no
$C\in\R_{>0}$ \st \er{6.109} holds. \csq, Corollary \rf{c6.33} strictly improves \pp y \er{6.86} of the \f~$\psi_x$.
\eex

\Wanp prove one of the main results of this section, namely the \ex\ of monoid- and order-\is sms between the additive group $(\R,+,0)$ and the
\mlv\ group $(\R_{>0},\cdot,1)$. It turns out that these \is sms are \ext s of the maps $\psi_x$, $x\in\R_{>1}$.

For the reader's convenience we recall that if $A,B$ are \nss s of~$\R$ which are \ba\ and \sf y $A\sbs B$, then
\beq6.111
\sup A\le \sup B.
\e
However, if $\sup A=\sup B$ then $A$ is not necessarily equal to $B$.

Indeed, $\sup A$ and $\sup B$ exist by Theorem \rf{t5.1}, and every \ub\ for~$B$ is an \ub\ for~$A$ since $A\sbs B$. \E\Ip $\sup B$ the least \ub\
for~$B$ is an \ub\ for~$A$, which is greater than or equal to $\sup A$ the least \ub\ for~$A$. Let $A:=\{y\in\R: 0\le y<1\}$ and $B:=\{y\in\R:
0\le y\le 1\}$. Then $1$ is the greatest \el\ of~$B$, hence $1=\sup B$. We have $\sup A\le1$ by \er{6.111}, and $0\le\sup A$ since $0\in A$. Suppose,
for \cd ion, that \te s $\ov a\in\R$ \st $\ov a:=\sup A<1$. Then $0<\ov a<\frac12(\ov a+1)<1$. A~\cd ion, since $\frac12(a+1)>\ov a\ge a$ \fa $a\in
A$, and $\frac12(\ov a+1)>\ov a$, hence $\ov a\ne \sup A$.

\blm6.35
Let $z\in\R$ and set
\beq6.112
A_z:= \{\psi_x(r)\in\R: r<z,\ r\in\Q\}.
\e
Then the \fw\ holds\dw
\bga6.113
A_z \hbox{ is a nonempty \ba\ subset of }\R,\\
A_z\sbs A_{z'} \hbox{ \fa } z'<z,\ z'\in\R. \lb{6.114}
\e
\elm

\proof \

``\er{6.113}'': Let $z\in\R$. By \E\Pr s \rf{p4.31} and \rf{p4.34}, \te s $r\in\Q$ \st $z-1< r<z$. Hence $r\in A_z$.

``\er{6.114}'': Let $z\in\R$ and let $u\in A_z$. Then by \er{6.112} \te s $r\in\Q$ \st $r<z$, $u=\psi_x(r)$. Let $z'\in\R$ be \st $z<z'$. Then
$r<z<z'$ implies $r<z'$. Since $u=\psi_x(r)$, we have $u\in A_{z'}$ by \er{6.113}.
\endproof

\bdf6.36
\E\fa $x\in\R_{>1}$ and $z\in\R$ set
\beq6.115
\wh\psi_x(z):= \sup\{\psi_x(r)\in\R: r<z,\ r\in\Q\}.
\e
\edf

Note that $\sup\{\psi_x(r)\in\R: r<z\}= \sup A_z$, which is well-defined in view of Theorem \rf{t5.1} and Lemma \rf{l6.35}. \Mo $\wh\psi_x(z)>0$
\fa $z\in\R$ since $\psi_x(r)>0$ \fa $r\in\Q$.

\blm6.37
The \f\ $\wh\psi_x:\R \to\R_{>0}$ \sf ies
\beq6.120
\wh\psi_x(z)=\psi_x(z) \qh{\fa $z\in\Q$,}
\e
thus $\wh\psi_x$ extends $\psi(x)$. \E\Ip $\wh\psi_x(0)=1$ and $\wh\psi_x(1)=x$.
\beq6.121
\wh\psi_x(z) \le \wh\psi_x(z') \qh{\fa $z,z'\in\R$ \st $z<z'$.}
\e
\elm

\proof
\er{6.120} follows from \er{6.86} and \er{6.121} from \er{6.114} and \er{6.111}.
\endproof

\bex6.38
Define $\ov F:\R\to\R$ by setting
\beq6.122
\ov F(r):= \bca
r &\hbox{for } r<\sqrt2,\\
1+\sqrt2 &\hbox{for } r=\sqrt2,\\
2+r &\hbox{for }r>\sqrt2,
\eca
\e
and $F:\Q\to\R$ by setting $F(r)=\ov F(r)$ \fa $r\in\Q$.

\hph i,i, Show that $F:\Q\to\R$ is order-continuous, i.e.\ \sf ies \er{6.86} with $\psi$ replaced by~$F$, and is strictly in\cre.

\hph ii,, Set $\wh F(z):= \supl_{r<z}F(r)$ \fa $z\in\R$. Show that $\wh F(r)=\ov F(r)= F(r)$ \fa $r\in\Q$, $\wh F(\sqrt2)=\sqrt2 < \ov F(\sqrt2)$.
Hence $\wh F$~is \ti{not\/} order-continuous.
\eex

\bpr6.39
The \f\ $\wh\psi_x:\R\to\R_{>0}$ is order-continuous.
\epr

\proof
We have to show that \fa $z\in\R$,
\beq6.123
\sup_{y<z,y\in\R} \wh\psi_x(y) = \wh\psi_x(z) = \inf_{y>z,y\in\R}\wh\psi_x(y).
\e
In what follows we use the letters $r,s,t,u$ for \el s of~$\Q$, and $w,x,y,z$ for \el s of~$\R$. Since $\wh\psi_x$~is in\cre\ by \er{6.121}, we have
$\wh\psi_x(y)\le \wh\psi_x(z)\le \wh\psi_x(w)$ for $y<z<w$, hence $\supl_{y<z}\wh\psi_x(y) \le \wh\psi_x(z)\le \inf\limits_{w>z}\wh\psi_x(w)$. In
view of Lemma \rf{l5.38} we have to show that \fe $\ve\in\R_{>0}$ \te\ $y<z$ and $w>z$ \st
\beq6.124
\wh\psi_x(z) - \wh\psi_x(y) < \ve \qh{and } \wh\psi_x(z)-\wh\psi_x(w)<\ve.
\e
Exercise \rf{ex6.38} shows us that order-continuity of~$\psi$ is \ti{not\/} \sft\ for proving \er{6.124}. We invoke \er{6.107} of Corollary
\rf{c6.33}. By \E\Pr s \rf{p4.31} and \rf{p4.34} \te s $u\in\Q$ \st $z+1<u<z+2$. Set $M':= \psi_x(u+1)-\psi_x(u)\in\R_{>0}$, and let $M\in\Q_{>0}$ \st
$M'<M<M'+1$. Then we obtain
\beq6.125
\psi_x(r) - \psi_x(s) \le M(r-s), \q s<r<u.
\e
Now given $\ve\in\R_{>0}$, let $\ve'\in\Q$ be \st $0<\frac12\ve <\ve' <\ve$, and set
\beq6.126
e:= \ve' M\mo \in\Q.
\e
We make the \ad al \as\ on~$z$ by assuming
\beq6.127
z\in\R_{>0}.
\e
This \as\ will be removed later on. \Wanp apply \E\Pr\ \rf{p4.14} with $\a:=z$ and $e:=e$, we find $r,s\in\Q_{>0}$ \st $s<z<r$ \sf ying
$r-s=\ve'M\mo$.

\E\Tf by \er{6.125} we obtain
\beq6.128
\psi_x(r) - \psi_x(s) \le \ve'< \ve.
\e
Noting that $\wh\psi_x(r) = \psi_x(r)$, $\wh\psi(s)=\psi_x(s)$ by \er{6.120} since $r,s\in\Q_{>0}$, we obtain $\wh\psi_x(r) < \wh\psi_x(z) <
\wh\psi_x(s)$ by \er{6.121} since $r<z<s$. \Mo we have $(\wh\psi_x(r)-\wh\psi_x(z)) +\break (\wh\psi_x(z)-\wh\psi_x(s)) = \wh\psi_x(r) -
\wh\psi_x(s) = \psi_x(r) - \psi_x(s) < \ve$. Since both terms of the \LHS\ are non\ng\ we obtain
\beq6.129
\wh\psi_x(r) - \wh\psi_x(z) < \ve \qh{and } \wh\psi_x(z) - \wh\psi_x(s) < \ve,
\e
which completes the proof of the lemma whenever $z\in\R_{>0}$.

We now remove this \ad al \as\ $z>0$. Recall that we have to prove \er{6.124}.

``\ti{Case} $z=0$'': If $z=0$, or more generally, if $z\in\Q$, \er{6.124} follows from \er{6.86}. Indeed, by Lemma \rf{l5.38}, \fe $z\in\R_{>0}$
\te\ $y\in\Q$ and $w\in\Q$ \st $\psi_x(z)-\psi_x(y) < \ve$ and $\psi_x(z) - \psi_x(w) < \ve$. \Mo by \er{6.120} we have $\psi_x(z) = \wh\psi_x(z)$,
$\psi_x(y)=\wh\psi_x(y)$ and $\psi_x(w) = \wh\psi_x(w)$, which completes the proof.

``\ti{Case} $z<0$'': Note that we used \as\ \er{6.127} in order to be able to apply \E\Pr\ \rf{p4.14}. Setting $z':=|z|=-z>0$ we find $r',s'\in
\Q_{>0}$ \st $r'<z'<s'$ and $\psi_x(s')-\psi_x(r')<\ve$ by \er{6.128}. Setting $s:=-s'$, $r:=-r'$ we obtain $\wh\psi_x(-s) - \wh\psi_x(z)<\ve$ and
$\wh\psi_x(z)-\wh\psi_x(-r)<\ve$ where $-r<z<-s$. This completes the proof of the \Pr.
\endproof

It is convenient to introduce the \fw

\bnt6.40
We denote by $a_n\ua a$ (resp.\ $a_n\da a$) a \sq\ $\zb an\Na$ of \el s of~$\R$, which is in\cre\ (resp.\ de\cre) and bounded, with $a\in\R$ as
supremum (resp.\ infimum).
\ent

\blm6.41 \

\hph i,i, Let $a_n\ua a$, $b_n\ua b$. Then we have
\bea6.130
a_n+b_n &\ua a+b,\\
a_n\cdot b_n &\ua a\cdot b \hbox{ if $a_n>0$, $b_n>0$ \fa }n\in\Na. \lb{6.131}
\e

\hph ii,, \Mo
\bea6.132
a_n \ua a &\hbox{ iff }-a_n\da -a,\\
a_n \da a &\hbox{ iff }-a_n\ua -a. \lb{6.133}
\e
\elm

\proof[Proof of \er{6.130}]
We have $a_n\le a_m$, $b_n\le b_m$, hence $a_n+b_n \nde3.9 \le a_m+b_m$, $n\le m$, $n,m\in\Na$. \Mo $a_n+b_n \le a+b$, $n\in\Na$. By Lemma \rf{l5.38}
\fe $\ve \in \R_{>0}$ \te s $\ov n\in\Na$ (resp.\ $\ov m\in\Na$) \st  $a-a_{\ov n}<\frac\ve2$, $b-b_{\ov m}<\frac\ve2$. Hence $a-a_n \le a-a_{\ov n}
\le \frac\ve2$, $n\ge\ov n$, $b-b_m \le b-b_{\ov m} \le \frac\ve2$, $m\ge\ov m$. Setting $\ov k:=\max(\ov n,\ov m)$, we obtain $(a+b)-(a_n+b_n) =
(a-a_n)+ (b-b_n)\le\frac\ve2 + \frac\ve2=\ve$, whenever $n\ge\ov k$. Hence $a_n+b_n\ua a+b$ by Lemma \rf{l5.38}.
\endproof

\bex6.42
Prove \er{6.131}--\er{6.133}.
\eex

\blm6.43
\E\fe $z\in\R$ \te s a \sq\ $\zb rn\Na$ of \el s of~$\Q$ \sf ying $r_n\ua z$.
\elm

\proof
If $z\in\Q$ then take $r_n:=z$ \fa $n\in\Na$. If $z\in\R_{>0}\sm\Q$, then $z=m+\rho$, $m\in\N$, and $\rho\in(0,1)_\R$ by \er{6.28}. \E\fe $a\in\Na
\sms1$ \te s $\wt s_n\ua\rho$, $\wt s_n\in\Q_a\cap(0,1)_\Q$, $n\in\Na$, by Lemma \rf{l6.12}. Hence $m+\wt s_n\ua m+\rho=z$ by \er{6.53} and
$r_n:=m+\wt s_n\in\Q$, $n\in\Na$. If $z\in\R_{<0}\sm\Q$, then $-z\in\R_{>0}\sm\Q$, and $-z=m+\rho$, $m\in\N$, $\rho\in(0,1)_\R$, and $z=(-m-1)+(1-\rho
)$. Since $1-\rho\in(0,1)_\R$, \te\ $s_n'\in(0,1)_\Q$ \st $s_n'\ua 1-\rho$. Hence $-(m+1)+s_n' \ua -(m+1)+1-\rho = -(m+\rho)= z$. Set $r_n:=(m+1)
+s_n'$, $n\in\Na$.
\endproof

\bpr6.44
The \f\ $\wh\psi_x$ \sf ies
\beq6.134
\wh\psi_x(y+z) = \wh\psi_x(y)\cdot \wh\psi_x(z), \q y,z\in\R.
\e
\epr

\proof
By Lemma \rf{l6.43} \te\ $r_n,s_n\in\Q$ \st $r_n\ua y$ and $s_n\ua z$. By \er{6.130} $r_n+s_n\ua y+z$, hence $\wh\psi_x(y_n+s_n) \ua \wh\psi_x(y+z)$,
$\wh\psi_x(r_n)\ua \wh\psi_x(y)$, $\wh\psi_x(s_n)\ua\wh\psi_x(z)$ by \E\Pr\ \rf{p6.39}. By \E\Pr\ \rf{p5.29}, we have \er{5.61} with $r:=r_n$,
$s:=s_n$, $n\in\Na$. Note that $\psi_x(r_n),\psi_x(s_n)>0$ for $n\in\Na$, hence $\psi_x(r_n)\cdot \psi_x (s_n)\ua \wh\psi_x(y)\cdot\wh\psi(z)$
by \er{6.131}. \E\Tf \er{6.134} holds.
\endproof

\Wanp prove one of the main results of this section.

\bth6.50
Let $x\in\R$, $x>1$ and let $\wh\psi_x:\R \to\R_{>0}$ be the \f\ defined in \er{6.115}. Then the \fw\ assertions hold.

\hph i,vi, The \f\ $\wh\psi_x$ is order-continuous, i.e.\ \sf ies \er{6.123}, and is the only order-continuous \ext\ of the \f\ $\psi_x:\Q\to
\R_{>0}$ \itd\ in \E\Pr\ \rf{p5.29}.

\hph ii,v, $\wh\psi_x(y+z) = \wh\psi_x(y)\cdot\wh\psi_x(z)$ \fa $y,z\in\R$.

\hph iii,i, The \f\ $\wh\psi_x:(\R,\ge)\to(\R_{>0},\ge)$ is strictly in\cre\ and convex $($in the sense of \E\df\ \rf{d6.25} with $t\in[0,1]_\R)$.

\hph iv,i, $\inf\limits_{z\in\R}\wh\psi_x(z) = 0$ and $\{\wh\psi_x(z): z\in\R\}$ is not \ba.

\hph v,ii, The \f\ $\wh\psi_x$ is sur\jc.

\hph vi,i, The \f\ $\wh\psi_x$ is an \ois sm from $(\R,\ge)$ onto $(\R_{>0},\ge)$.

\hph vii,, The \f\ $\wh\psi_x$ is a monoid-\is sm from the group $(\R,+,0)$ onto the group $(\R_{>0},\cdot,1)$.
\eth

\proof \

(i) Follows from \E\Pr\ \rf{p6.39} and \er{6.120}. The ``\uq'' \pp y follows from \er{6.120} and \er{6.115}.

(ii) Follows from \E\Pr\ \rf{p6.44}.

(iii) ``\ti{Strict in\cre ness}'': Let $y<z$, $y,z\in\R$, and set $h:=z-y>0$. Then $\wh\psi_x(z) = \wh\psi_x(y+h) = \wh\psi_x(y)\cdot\wh\psi_x(h)$.
Let $r\in\Q$ be \st $0<r<h$. Then $\wh\psi_x(h)\ge \wh\psi_x(r)=\psi_x(r)$. By \E\Pr\ \rf{p5.29}\,(ii) $\psi_x(r)>\psi_x(0)\nde5.60 = 1$. Hence
$\wh\psi_x(h)>1$, and $\wh\psi_x(z)=\wh\psi_x(y)\cdot\wh\psi_x(h) > \wh\psi_x(y)\cdot1$ by \er{3.10}. \E\Tf $\wh\psi_x(y)< \wh\psi_x(z)$.

``\ti{Convexity}'': From Lemma \rf{l6.26} we find \er{6.91} with $G:=\wh\psi_x$ as a con\sq\ of~(ii). Given $t\in(0,1)$, we find by Lemma \rf{l6.12}
a~\sq\ $\zb tn\N$, $t_n\in\Q_2\cap(0,1)$ \st $t_n\ua t$. Proceeding as in the proof of Lemma \rf{l6.26} we find that \er{6.87} holds for
$F:=\wh\psi_x$.

(iv) Suppose, for \cd ion, that $\{\wh\psi_x(z): z\in\R\}$ is \ba. Then \te s $C'\in\R_{>0}$ \st $\wh\psi_x(n)\le C'$ \fa $n\in\Na$. Let $C\in\Q_{>0}$
\sf y $C'<C<C'+1$. Let $r\in\Q_{>0}$ \sf y $1<r<x$. In view of~(ii) and using \In\ on $n\in\Na$ we find $\wh\psi_x(n) = (\wh\psi_x(1))^n$ \fa $n\in
\Na$. Note that $\wh\psi_x(0)\nad{\rm(i)}= \psi_x(0)=1$ and $\wh\psi_x(1)\nad{\rm(i)}= x$, hence $\wh\psi_x(1)\nad{\rm(iii)}> r>1$. \Mo $\wh\psi_x(n)
= (\wh\psi_x(1))^n > r^n$ \fa $n\in\Na$ by using \In\ on~$n$. Hence $r^n<C$, $C\in\Q_{>0}$, \fa $n\ge1$, \cd ing \er{5.77}. Hence $\{\wh\psi_x(z):
z\in\R\}$ is not \ba.

Suppose, for \cd ion, that $\infl_{z\in\R}\wh\psi_x(z) = \a\in\R_{>0}$. \E\Ip $\wh\psi_x(-n)\ge\a$ \fa $n\in\Na$. We have $\wh\psi_x(-n)\nad{\rm(i)}=
\psi_x(-n) \nde5.60 = (\psi_x(n))\mo\ge\a$ \fa $n\in\N$. Hence $\psi_x(n)\nde3.10 \le \a\mo$, and $(\wh\psi_x(1))^n = \wh\psi_x(n) = \psi_x(n)
\le \a\mo$. As above we find $r\in\Q$, $r>1$, and $C\in\Q$ \st $\a\mo < C < \a\mo+1$, $r^n\le C$ \fa $n\ge1$. A~\cd ion.

(v) Since $\wh\psi_x$ is strictly in\cre\ by (iii) and convex in the sense of \E\df\ \rf{d6.25} where $t\in(0,1)_\R$ by~(iii), it follows from the
proof of Corollary \rf{c6.33} that \er{6.108} holds if $s,r,u\in\R$, $s<r<u$. \E\Ip \fe $n\in\Na$, setting $s:=-n$, $r:=n$ and $u:=n+1$, we find
$0<\wh\psi_x(r)-\wh\psi_x(s)\le M$ \fs $M\in\R_{>0}$ depending on~$u$ but not on $r,s$. We may apply Lemma \rf{l5.41} with $a:=-n$, $b:=n$, $n\in\Na$,
$f$~the \rt ion of~$\wh\psi_x$ to $[-n,n]$, since $\wh\psi_x(-n) < \wh\psi_x(n)$ and $0<\wh\psi_x(z)-\wh\psi_x(y) \le M(z-y)$, $a<y<z<b$. We find that
\beq6.135
\{w\in\R: \wh\psi_x(-n) < w <\wh\psi_x(n)\} \sbs \wh\psi_x((-n,n))
\e
where $(-n,+n):=\{w\in\R: -n<w<n\}$. Therefore, given $t\in\R$, \te s $n\in\N$ \st $-n<t<n$, since $\R$~is \Ar, hence \te s $w\in\R$ \st $\wh\psi_x
(w)=t$ in view of \er{6.135}.

(vi) Since $\wh\psi_x$ is strictly in\cre, it is in\jc, hence bi\jc\ by~(v). By Lemma 1.3.34, $\wh\psi_x$ is an \ois sm from $(\R,\ge)$ into
$(\R_{>0},\ge)$.

(vii) Since $\wh\psi_x(0)=1$ and $\wh\psi_x$ \sf ies (ii), $\wh\psi_x$ is a monoid-\hm sm from $\pz1\R$ into $\pz2{\R_{>0}}$. Since $\wh\psi_x$
is bi\jc, it follows from Lemma 2.1.8\,(i) that $\wh\psi_x$ is a monoid-\is sm.
\endproof

\bnt6.51
The usual notation for $\wh\psi_x(z)$ is $x^z$ for $x\in\R_{>1}$ and $z\in\R$.
\ent

\begin{exe}[see Lemma \rf{l9.15}]\lb{ex6.52} \

\hph i,i, Let $x_1,x_2\in\R_{>1}$ and $z\in\R$. Show
\beq6.136
(x_1\cdot x_2)^z = (x_1^z)\cdot(x_2^z).
\e

\hph ii,, Let $x\in\R_{>1}$ and $z_1,z_2\in\R$. Show
\beq6.137
x^{z_1\cdot z_2} = (x^{z_1})^{z_2}.
\e
\eex

\bdf6.53
Let $x\in\R_{>1}$ and let $\wh\psi_x :\R \to \R_{>0}$ be the map defined in \er{6.115}. The map $\log_x : \R_{>0}\to\R$ is by \df\ the \ti{inverse}
of the map~$\wh\psi_x$, which exists by Theorem \rf{t6.50}\,(vi),\,(vii). It is called \ti{logarithm to base }$x$.
\edf

We obtain \fa $x\in\R_{>1}$, $z\in\R$,
\beq6.138
x^{\log_xz} = z, \q \log_x x^z=z.
\e

\bex6.54 \

\hph i,ii, Show that $\log_x$ defined in \E\df\ \rf{d6.53} extends $\log_x$ defined in \E\Pr\ \rf{p5.37}.

\hph ii,i, $\log_x(y\cdot z) = \log_x y + \log_xz$, $y,z\in\R_{>0}$.

\hph iii,, $\log_x:\R_{>0}\to\R$ is \ti{concave}, i.e.\ \sf ies
\beq6.139
\log_x((1-t)y + tz) \ge (1-t)\log_xy + t\log_xz
\e
\fa $y,z\in\R_{>0}$ and all $t\in[0,1]$.
\eex

\bnt6.55
Let $(\R_{>0},\cdot,\cdot_x',1,x)$ be the field \itd\ in \E\Pr\ 4.4.2\,(ii), where $f:=\psi_x$. We shall denote this field by~$x^\R$.
\ent

\bex6.56 \

\hph i,ii, Show that if $y',z'\in x^\R$, then the binary \op\ \sf ies
\beq6.140
y' \cdot_x' z' = x^{(\log_xy')\cdot(\log_xz')}.
\e

\hph ii,i, Let $x,x'\in\R_{>1}$ and $y\in\R_{>0}$. Show that
\beq6.141
\log_xy = (\log_xx')\cdot(\log_{x'}y).
\e

\hph iii,, Show that $\wh\psi_x$ \sf ies for $h,h',h''\in\R_{>0}$
\beq6.142
\frac{\wh\psi_x(h)-1}h \le \frac{\wh\psi_x(h')-1}{h'} \le \frac{\wh\psi_x(h'')-1}{h''},
\e
whenever $h<h'<h''$, and that $\infl_{h>0} \frac{\wh\psi_x(h)-1}h$ exists.
\eex


\begin{exe}[see Section \ref{ass.9}]\lb{ex6.57} \

\hph i,ii, Show that the \pf\ (see \E\df\ 4.4.29) of the field $x^\R$, $x\in\R_{>1}$, is~$x^\Q$.

\hph ii,i, Show that $x^\R$ is an \of\ when endowed with the \og~$\ge_x$ defined by $\a\ge_x\b$, $\a,\b\in x^\R$ if $\log_x\a \le \log_x\b$.

\hph iii,, Show that $(x^\R,\le_x)$ is complete and that the \pf\ of~$x^\R$ is dense in~$x^\R$.
\eex

\bqu6.58
Let $(F,\qu_0,\qu_1,e_0,e_1)$ be an \of\ with \og~$\ge$. Suppose that \te\ $x\in F_{> e_1}$ and a monoid-\is sm $f_x$ from $(F_{\ge e_0},\qu_0,e_0)$
onto $(F_{\ge e_1},\qu_1,e_1)$ \sf ying $f_x(e_1)=x$. Suppose, moreover, that the \fw\ holds\dw
\beq6.143
\hbox{the \pf\ of $x^F$ is dense in $x^F$ $($see \E\df\ \rf{d4.33}$)$.}
\e
Is $(F,\qu_0,\qu_1,e_0,e_1)$ a complete \of?
\equ

\newpage
\Subsubsection{Metrics and valuations}\label{ass.7}

In the first part of this section we introduce some \fd\ notions belonging to Analysis, which will play an important role in the study of \pp
ies of the field~$\R$. We begin with the notion of \tb{limit} of a \sq\ of \el s of~$\R$, which extends the notion of supremum (resp.\ infimum)
of a bounded in\cre\ (resp.\ de\cre) \sq\ in~$\R$. Using Notation \rf{n6.40} and recalling Lemma \rf{l5.38} we have $a_n \ua a$ iff $a_n\le a_m\le a$
\fa $n,m\in\N$, $n\le m$, and \fe $\ve\in\R_{>0}$ \te s an $N\in\N$ \st $a_N+\ve>a$. Note that these two \cn s imply the \fw\ one:
\beq7.1
\hbox{\E\fe $\ve\in\R_{>0}$ \te s an $N\in\N$ \st $|a-a_n|<\ve$ \fa}n\ge N.
\e
Indeed, $|a-a_n|\nde5.41 = a-a_n$ since $a_n<a$, and if $a-a_N<\ve$ then $a-a_n<\ve$ \fa $n\ge N$ since the \sq\ $a_n$ is in\cre. Similarly, one
verifies that $a_n\da a$ implies \er{7.1}.

\bds7.1 \

\hph i,i, A \sq\ $\zb an\N$ of \el s of~$\R$ is said to \ti{be convergent\/}, or to \ti{converge to a limit\/}~$a$ if \er{7.1} holds. Otherwise, the
\sq\ $\zb an\N$ is said to be \ti{divergent}.

\hph ii,, Let $\zb an{\N\ (\text{resp.\ }\Na)}$ be a \sq\ of \el s of a \ns~$A$. Then the \sq\ $\{a_{j(n)}\}_{n\in\N\ \text{(resp.\ }\Na)}$ is called
a \ti{sub\sq} of the \sq\ $\zb an{\N\ (\text{resp.\ }\Na)}$ if $j$~is a \ti{strictly in\cre} self-map of~$(\N,\ge)$ (resp.\ $(\Na,\ge)$).
\eds

\blm7.2 \

\hph i,ii, A \cvg\ \sq\ has at most one limit.

\hph ii,i, A \cvg\ \sq\ is bounded in $(\R,\ge)$.

\hph iii,, Every sub\sq\ of a \cvg\ \sq\ to a limit~$a$ converges to the same limit~$a$.

\hph iv,, An in\cre\ \sq\ in~$\R$ which is \ba\ $($resp.\ below$)$ \cg s to its supremum $($resp.\ infimum$)$.
\elm

\proof\

(i) Suppose that $a,a'\in\R$ are limits of a \sq\ $\zb an\N$. Then $|a-a'| = |(a-a_n)+(a_n-a')| \nde5.44 \le |a-a_n|+|a_n-a'|\nde5.46 = |a-a_n| +
|a'-a_n|$ \fa $n\in\N$. \E\fe $\ve\in\R_{>0}$ \te\ $N\in\N$ and $N'\in\N$ \st $|a-a_n|<\ve$ for $n\ge N$ and $|a'-a_n|<\ve$ for $n\ge N'$. \E\Tf \fe
$\ve>0$, $|a-a'|<2\ve$ for $n\ge\max(N,N')$. Thus $|a-a'|<2\ve$ \fe $\ve\in\R_{>0}$. Recall that $|a-a'|\ge0$. Suppose that $\a:=|a-a'|>0$. Then,
choosing $\ve:=\frac\a2$ we obtain $\a<\a$, a~\cd ion. Hence $|a-a'|=0$, $a-a'=0$ by \er{5.41}, and $a=a'$.

(ii) Let $\zb an\N$ be a \sq\ of \el s of~$\R$ converging to a limit $a\in\R$. From \er{7.1} with $\ve:=1$ we infer that \te s $N\in\N$ \st $|a-a_n|
<1$ \fa $n\ge N$. If $N=0$ then  $|a-a_n|<1$ \fa $n\in\N$. Hence $|a_n|=\mathopen|-a_n\mathclose|=|a+(-a_n)-a| = |(a-a_n)+(-a)| \le |a-a_n|+
\mathopen|-a\mathclose| =  |a-a_n|+|a| \le 1+|a|
=: M \in\R_{\ge0}$ \fa $n\in\N$. Hence $-M\le a_n\le M$ \fa $n\in\N$, and $\{a_n\in\R: n\in\N\}$ is a bounded subset of $(\R,\ge)$. If $N\in\Na$ then
$\{a_n\in\R: n\in\N\} = \{a_n\in\R: 0\le n\le N-1\} \cup \{a_n\in\R: n\ge N\}$. The first set is finite. We leave it to the reader to show that
a finite subset of $(\R,\ge)$ is bounded. One shows as in the case where $N=0$ that the second set is bounded. The union of two bounded sets
of $(\R,\ge)$ is easily shown to be bounded in $(\R,\ge)$.
\endproof

\bex7.3
Prove Lemma \rf{l7.2}\,(iii) and (iv).
\eex

\bnt7.4
We will use the notation $a_n\to a$, $\lim a_n=a$ or $\lim\limits_{n\ge0}a_n=a$, instead of the usual notation $\lim\limits_{n\to\iy}a_n=a$
for a \sq\ converging to a limit~$a$.
\ent

The ``constant'' \sq\ $a_n:=a\in\R$ \fa $n\in\N$ converges to~$a$. The bounded \sq\ $a_n:=(-1)^n$, $n\in\N$, has a sub\sq\ converging to~$1$ and
another sub\sq\ converging to~$-1$, hence this \sq\ is divergent in view of Lemma \rf{l7.2}. We now give an example of a \sq\ of \po\ \ra\ \nm s
converging to~$\sqrt2$ which is not the same as the \de\ of~$\sqrt2$. Note that $\sqrt2-1\in(0,1)$, hence $\sqrt2-1$ has a \de\ by Theorem
\rf{t6.14}\,(i) which we denote by~$s_n$, hence $\supl_{n\ge1}1+s_n = 1+(\sqrt2-1) = \sqrt2$ by \er{6.53}. We call this \sq\ the \de\ of~$\sqrt2$.
Since $1+s_n\ua\sqrt2$, we have $1+s_n\to\sqrt2$.

\bxa7.5
Let $\a\in\R_{>0}$  and let $f:\R_{>0}\to\R$ be defined by
\beq7.2
f(x):=\frac12\biggl(x+\frac\a x\biggr), \q x\in\R_{>0}.
\e
Note that if $x\in\R_{>0}$ then $\frac1x\in\R_{>0}$ by Lemma \rf{l3.18}, and $\frac12\bigl(x+\frac\a x\bigr)\in\R_{>0}$ by \er{3.9}, \er{3.10}.
Hence $f$ \ti{maps $\R_{>0}$ into $\R_{>0}$}. Applying Theorem 1.1.7 with $(E,e,S):=(\N,0,S)$ where $Sx:=x+1$, $x\in E$, $a\in\R_{>0}$, $F:=
\R_{>0}$, and setting $x_n:=\phi(n)$ \fa $n\in\N$, we find a \sq\ $\zb xn\N$ of \el s of $\R_{>0}$ \sf ying
\bea7.3
x_{n+1}&=f(x_n) \qh{\fa $n\in\N$,}\\
x_0&=a. \lb{7.4}
\e
We shall prove that $x_n\to\sqrt\a$. One shows as above that if $\a\in\Q_{>0}$ and $x\in\Q_{>0}$, then $f(x)\in\Q_{>0}$. \Mo if $\a:=2$, $a:=1$,
then $x_1:=\frac32$ and $x_2=1+\frac5{12}$ which is not a \dn.

We first show that if $x_n\to L$ \fs $L\in\R_{>0}$, then $L=\sqrt\a$. To this end we use the \fw
\exa

\blm7.6
Suppose $a_n\to L_1$, $b_n\to L_2$, where $a_n,b_n$ are \sq s of \el s of~$\R$, $L_1,L_2\in\R$. Then we have\dw
\bea7.5
{}&\hbox{If $a_n\in\R_{>0}$ \fa $n\in\N$, then $L_1\ge0$.}\\
&a_n\pm b_n \to L_1\pm L_2. \lb{7.6}\\
&a_n\cdot b_n \to L_1\cdot L_2. \hbox{ \E\Ip if $a_n=a$ \fa $n\in\N$, then $a_n\cdot b_n \to a\cdot L_2$.} \lb{7.7}\\
&\hbox{If $a_n\in\R_{>0}$ \fa $n\in\N$ and if $L_1>0$, then $\frac1{a_n} \to \frac1{L_1}$.} \lb{7.8}
\e
\elm

\proof[Proof of \er{7.8}]
We have to show that given $\ve\in\R_{>0}$ \te s an $N\in\N$ such that\break $\bigl|\frac1{a_n} - \frac1{L_1}\bigr|<\ve$ whenever $n\ge N$, $n\in\N$.
Note that $\bigl|\frac1{a_n}-\frac1{L_1}\bigr| = \bigl|\frac{L_1-a_n}{a_nL_1}\bigr| = \frac1{a_n}\cdot \frac1{L_1}\cdot |L_1-a_n|$, since $a_n,L_1>0$.
From $a_n\to L_1$, we infer that \te s $N_1\in\N$ \st $|a_n-L_1| < \frac12L_1$ for $n\ge N_1$. Let $n\ge N_1$, then $a_n>0$ and $a_n = L_1+a_n-L_1
= |L_1+a_n-L_1| \nde5.47 \ge \bigl||L_1|-|a_n-L_1|\bigr| = L_1 - |a_n-L_1| > L_1-\frac12L_1 = \frac12L_1$. Hence $a_n>\frac12L_1$, and $\frac1{a_n}
<\frac2{L_1}$. \If that $\bigl|\frac1{a_n}-\frac1{L_1}\bigr| \le \frac1{L_1}\cdot \frac2{L_1}|a_n-L_1|$. Given $\ve\in\R_{>0}$ we may find $N\in\N$
\st $N>N_1$ and $|a_n-L_1|<\frac12L_1^2\ve$ \fa $n\ge N$, since $a_n\to L_1$. \If that $\frac1{a_n}\to\frac1{L_1}$.
\endproof

\bex7.7
Prove \er{7.5}--\er{7.7} (see the proof of Lemma \rf{l7.35} and of \E\Pr\ \rf{p8.12}).
\eex

\Wanp show that $x_n\to\sqrt\a$. Set $a_n:=x_{n+1}$ and $b_n:=\frac12\bg(x_n+\frac\a{x_n})$, $n\in\N$. Suppose that \te s $L\in\R_{>0}$ \st
$x_n\to L$. Note that the \sq\ $\zb an\N$ is a sub\sq\ of the \sq\ $\zb xn\N$ (see \E\df\ \rf{d7.1}\,(ii)). Indeed, $x_{n+1}=x_{j(n)}$, $n\in\N$,
where $j(n):=n+1$, $n\in\N$. Observe that $j$~is the map~$S$ in \E\df\ 1.1.2, hence strictly in\cre. We have $a_n\to L$ by Lemma \rf{l7.2}\,(iii). \Mo
$\frac1{x_n}\to \frac1L$ by \er{7.8}, $\frac\a{x_n} = \a\cdot \frac1{x_n}\to\frac\a L$ by~\er{7.7}, $x_n+\frac\a {x_n}\to L+\frac\a L$ by~\er{7.6},
hence $b_n\to\frac12\bigl(L+\frac\a L\bigr)$ by~\er{7.7}. Since $a_n=b_n$ \fa $n\in\N$, and $a_n\to L$, $b_n\to \frac12\bg(L+\frac\a L)$, we obtain
$L=\frac12\bg(L+\frac\a L)$ by Lemma \rf{l7.2}, hence $L^2=\a$. Since $L>0$, we have $L=\sqrt\a$ by Theorem \rf{t5.5}.

\ssk
We now investigate the \ti{\cg nce} of the \sq\ $\zb xn\N$.

``\ti{Case} $x_0=\sqrt\a$'': We have $x_n=\sqrt\a$ \fa $n\ge1$. Indeed, $f(\sqrt\a)\nde7.2 = \frac12\bg(\sqrt\a+\frac\a{\sqrt\a}) \nad*=
\frac12(\sqrt\a +\sqrt\a)=\sqrt\a$. In $\nad*=$ we used Theorem \rf{t5.5}. Set $A:=\{n\in\N: x_n=\sqrt\a\}$. Clearly, $0\in A$. \Mo if $n\in A$ then
$n+1\in A$. Hence $A=\N$, and $x_n=\sqrt\a=f(x_n)$ \fa $n\in\N$.

``$x_0\in\R_{>0}$'': We recall that $f$ is a self-map of~$\R_{>0}$. Since $x_0\in\R_{>0}$, we have $x_n>0$ \fa $n\in\Na$ by Theorem 1.1.7. Set
$y_n:=\frac\a{x_n}$, $n\in\N$. Then $y_n\in\R_{>0}$ and $x_{n+1}\nde7.2 = \frac12(x_n+y_n)$ \fa $n\in\N$. We have $\bigl(\frac12(x_n+y_n)\bigr)^2
= \frac12\bg(\frac12x_n^2+\frac12y_n^2) + \frac12(x_ny_n) \nde6.88 \ge \frac12(x_ny_n) + \frac12(x_ny_n) = x_ny_n$. Thus $\frac12(x_n+y_n)\ge
\sqrt{x_ny_n}$, the \ti{\art\ mean} of $x_n,y_n$ is greater than or equal to the \ti{geometric mean} of $x_n,y_n$. Observe that $\frac12(x_n+y_n)
=\sqrt{x_ny_n}$ iff $x_n=y_n$. Indeed, if $x_n=y_n$, then $\frac12(x_n+y_n) = x_n = \sqrt{x_n^2} = \sqrt{x_ny_n}$. Conversely, if $\frac12(x_n+y_n)
= \sqrt{x_ny_n}$, then $x_n^2+2x_ny_n+y_n^2 = 4x_ny_n$, hence $(x_n-y_n)^2=0$, and $x_n=y_n$. Note that $x_n=y_n$ iff $x_n=\frac\a{x_n}$ iff
$x_n=\sqrt\a$. We claim
\beq7.9
x_{n+1}\ge \sqrt\a \qh{\fa}n\in\N.
\e
Indeed, we have $x_{n+1}=\frac12(x_n+y_n)\ge \sqrt{x_ny_n} = \sqrt{x_n\frac\a{x_n}} = \sqrt\a$. \Mo in view of what precedes the \fw\ holds: \E\fa
$n\in\N$
\beq7.10
x_n=\sqrt\a \hbox{ \ iff \ } x_{n+1}=x_n=y_n=\sqrt\a.
\e
\E\Ip if $x_0<\sqrt\a$ then $x_1>\sqrt\a$.

In order to study the \cg nce of the \sq\ $\zb xn\N$, it is \sft\ to consider the \sq\ $\zb zn\N$ defined by
\bea7.11
z_n&:= x_{n+1}, \q n\in\N, \\
z_0&>\sqrt\a. \lb{7.12}
\e
To this end we observe
\beq7.13
z_{n+1} = z_n - \frac{z_n^2-\a}{2z_n}, \q n\in\N,
\e
since the \RHS\ of the \et y is equal to $\frac{2z_n^2-z_n^2+\a}{2z_n}$. From \er{7.9} and \er{7.13} we infer
\beq7.14
z_{n+1} \le z_n, \q n\in\N.
\e
\E\oh from \er{7.9} we find that the \sq\ $\zb zn\N$ is bounded, hence $z_n\da \ov z$ \fs $\ov z\ge\sqrt\a$ in view of Corollary \rf{c5.12} and the
fact that $\ov z$~is the greatest \lo\ for the set $\{z_n\in\R_{>0}:n\in\N\}$. Since $z_n\da \ov z$ implies $z_n\to \ov z$, we obtain from the
discussion \fw\ Exercise \rf{ex7.7} that $\ov z=\sqrt\a$. \E\Ip $x_n\to\sqrt\a$ \fa $x_0\in\R_{>0}$. \If that $\sqrt\a$ is the only \el\ of $\R_{>0}$
\sf ying $f(\xi)=\xi$, $\xi\in\R_{>0}$. Such \el~$\sqrt\a$ is usually called a \ti{fixed point\/} of the self-map~$f$ of~$\R_{>0}$.

\ssk
Limits of \sq s of \el s of $\R$ can be defined by means of the notion of {\bf distance \f} or {\bf metric}. This notion plays an important role
in Analysis. This notion is also used in Algebra. An important example is the \ti{Hamming} metric defined on finite-\dm al \vs s over finite fields
(see \cite[p.~474]{Lidl}).

\bdf7.8
A \ti{metric} on a \ns~$X$\index{metric} is a map $d:X\t X\to\R_{\ge0}$ \sf ying the \fw\ axioms. \E\fa $x,y,z\in X$ we have\dw
\begin{align}
d(x,y)&=0 \qh{iff $x=y$}, \tag{M1}\\
d(x,y)&=d(y,x), \tag{M2}\\
d(x,y)&\le d(x,z)+d(z,y).\tag{M3}
\end{align}
A set $X$ endowed with a metric $d$ is called a \ti{metric space}, and is denoted by $(X,d)$.
\edf

\bxs7.9 \

\hph i,ii, If $X:=\R$, and
\beq7.15
d(x,y):=|x-y| \qh{\fa $x,y\in\R$,}
\e
then $d$ is a metric on $\R$.

\hph ii,i, If $X$ is a \ns, and
\beq7.16
d(x,y):=0 \hbox{ if }x=y \qh{and }d(x,y):=1 \hbox{ if }x\ne y,
\e
then $d$ is a metric on $X$, called the \ti{discrete} metric.\index{metric!discrete}

\hph iii,, If $(X,d)$ is a metric space and $Y$ is a \nss\ of~$X$, then $d_Y$, the \rt ion of~$d$ to $Y\t Y$, is a metric on~$Y$. The metric space
$(Y,d_Y)$ is called a (metric) \ti{subspace} of the space $(X,d)$.

\hph iv,, If $X:=\R_{>0}$ and $a\in\R_{>1}$, then the map $d:X\t X\to \R_{\ge0}$ defined by
\beq7.18
d(x,y):=\Bigl|\log_a \frac xy\Bigr|, \q x,y\in X,
\e
is a metric on $X$ (see \E\df\ \rf{d6.53}).
\exs

\bex7.10
Show that the maps $d$ defined in (i)--(iv) are metrics.
\eex

\begin{prt}The Euclidean metric on $\R^2$. \lb{p7.11}

Let $x:=\binom{x_1}{x_2}$, $y:=\binom{y_1}{y_2}$ and $z:=\binom{z_1}{z_2}$ be \el s of~$\R^2$. Let $d:\R^2\to\R_{\ge0}$ be the map defined by
\beq7.19
d(x,y):= \sqrt{(x_1-y_1)^2 + (x_2-y_2)^2}.
\e
We claim that $d$ is a metric on $\R^2$ called the \tb{Euclidean metric}.\index{metric!Euclidean}

``(M1)'': If $x=y$, that is, $x_1=y_1$ and $x_2=y_2$, then $d(x,x)=\sqrt{0^2+0^2} = \sqrt{0} =0$. Conversely, if $d(x,y)=0$, then $(x_1-y_1)^2
+(x_2-y_2)^2=0^2=0$. Hence $(x_1-y_1)^2=0$ and $(x_2-y_2)^2=0$ (see (2.1.9) for $(\R_{\ge0},+,0)$). \E\Tf $x_1=y_1$ and $x_2=y_2$.

``(M2)'': Follows from $(x_1-y_1)^2=(y_1-x_1)^2$ and $(x_2-y_2)^2=(y_2-x_2)^2$.

``(M3)'': Let $a,b\in\R^2$ and set $\|c\|_2:=\sqrt{c_1^2+c_2^2}$ for $c\in\R^2$. We first claim that
\beq7.20
\Biggl|\sum_{i=1}^2a_ib_i\Biggr| \le \|a\|_2\|b\|_2.
\e
Note that $\|c\|_2=0$ iff $c=0$. If $a=0$ or $b=0$ then \er{7.20} is trivial. We suppose $a\ne0$ and $b\ne0$. We have
\bmlg
\biggl|\sum_{i=1}^2\frac{a_i}{\|a\|_2}\cdot \frac{b_i}{\|b\|_2}\biggr| \le
\sum_{i=1}^2\biggl|\frac{a_i}{\|a\|_2}\cdot \frac{b_i}{\|b\|_2}\biggr|
= \sum_{i=1}^2\frac{|a_i|}{\|a\|_2}\cdot \frac{|b_i|}{\|b\|_2} \nde6.88 \le
\sum_{i=1}^2 \frac12\,\frac{|a_i|^2}{\|a\|_2^2} + \sum_{i=1}^2\frac12\,\frac{|b_i|^2}{\|b\|_2^2}\\
= \frac12\,\frac1{\|a\|_2^2}\sum_{i=1}^2|a_i|^2 + \frac12\,\frac1{\|b\|_2^2}\sum_{i=1}^2|b_i|^2
=\frac12\,\frac{\|a\|_2^2}{\|a\|_2^2} + \frac12\,\frac{\|b\|_2^2}{\|b\|_2^2} =1.
\e
Hence
\bmlg
\biggl|\sum_{i=1}^2a_ib_i\biggr| = \|a\|_2\|b\|_2\cdot\frac1{\|a\|_2\|b\|_2}\biggl|\sum_{i=1}^2a_ib_i\biggr|
=\|a\|_2\|b\|_2\biggl|\frac1{\|a\|_2\|b\|_2} \sum_{i=1}^2a_ib_i\biggr|\\
= \|a\|_2\|b\|_2\biggl|\sum_{i=1}^2\frac{a_ib_i}{\|a\|_2\|b\|_2}\biggr| \le \|a\|_2\|b\|_2.
\e
In what precedes we used (2.4.38), (2.4.39), (2.4.41) with $a,b:[0,n]\to\R$ and $m\in\R$.
We now prove:
\beq7.21
\|a+b\|_2 \le \|a\|_2+\|b\|_2, \q a,b\in\R^2.
\e
If $a+b=0$ then \er{7.21} is trivial. We assume $a+b\ne0$. We have
\bmlg
\|a+b\|_2^2 = \sum_{i=1}^2(a_i+b_i)^2 = \sum_{i=1}^2|a_i+b_i|\,|a_i+b_i|\le \sum_{i=1}^2(|a_i|+|b_i|)|a_i+b_i| \\ =
\sum_{i=1}^2|a_i|\,|a_i+b_i| + \sum_{i=1}^2 |b_i|\,|a_i+b_i|
\nde7.20 \le \Bgg(\sum_{i=1}^2|a_i|^2)^{1/2}\Bgg(\sum_{i=1}^2|a_i+b_i|^2)^{1/2}\\ + \Bgg(\sum_{i=1}^2|b_i|^2)^{1/2}\Bgg(\sum_{i=1}^2|a_i+b_i|^2)^{1/2}
=\|a\|_2\|a+b\|_2 + \|b\|_2\|a+b\|_2.
\e
Dividing both sides by $\|a+b\|_2$ we find \er{7.21}.

Finally, setting $a:=x-z$, $b:=z-y$ in \er{7.21}, we find $\|x-y\|_2 = \|(x-z)+(z-y)\|_2 = \|a+b\|_2 \le \|a\|_2+\|b\|_2 = \|x-z\|_2
+\|z-y\|_2$, which proves (M3).
\end{prt}

\bex7.12 \

\hph i,i, Show that
\beq7.22a
d(x,y) := |x-y|^{1/2}, \q x,y\in\R,
\e
is a metric on $\R$. Find $x,y,z\in\R$ \st the in\et y (M3) is strict.

\hph ii,, Find $x,y,z\in\R^2$ \st the in\et y in \er{7.21} is strict.
\eex

\brm7.13 \

\hph i,i, In\et y (M3) is usually called the \ti{triangle in\et y}. In $\R^2$ if $x,y,z$ are the vertices of a triangle, (M3)~says
that the sum of the lengths of two sides is greater than the length of the other side. This geometric \cn\ on $x,y,z\in\R^2$ is \ev t to the \fw\
algebraic one: $y-x$ and $z-x$ are \ti{\li\/} in the \vs\ $(\R,\R^2)$. The length of the side~$xy$ is equal to $\|y-x\|_2$. Note that the \et y in
(M3) holds iff $y-x$ and $z-x$ are linearly dependent.

\hph ii,, Proceeding as in \rf{p7.11} one shows that the map $d:\R^n\t\R^n \to\R$, $n\in\N$, $n\ge3$, defined by
\beq7.22
d(x,y):= \sqrt{\sum_{i=1}^n (x_i-y_i)^2}, \q x,y\in\R^n,
\e
is a metric, called the \tb{Euclidean metric} on $\R^n$.
\erm

\bex7.14
Show that the \fw\ maps from $\R^n$ into $\R_{\ge0}$ are metrics on $\R^n$:
\bea7.23
d(x,y)&:=\sum_{i=1}^n |x_i-y_i|,\\
d(x,y)&:=\max_{1\le i\le n} |x_i-y_i|. \lb{7.24}
\e
\eex

We now turn to the notion of \ti{\cvg\ \sq} in a \ms. Observe that if we replace in \E\df\ \rf{d7.1} $|a-a_n|$ by $d(a,a_n)$ where $d$~is the metric
on~$\R$ defined in \er{7.15}, and if we set $c_n:=d(a,a_n)$ \fa $n\in\N$, then $a_n\to a$ iff $c_n\to 0$. Thus the \sq\ $\zb an\N$
converges to $a\in\R$ iff the \sq\ $\{d(a,a_n)\}_{n\in\N}$ of \el s of~$\R_{\ge0}$ converges to~$0$ in the sense of \E\df\ \rf{d7.1}. These
considerations motivate the \fw\ \df.

\bdf7.15
A \sq\ $\zb xn\N$ in a \ms\ $(X,d)$ is said to be \ti{\cvg\/} if \te s an $x\in X$ \st \fa $\ve>0$ there is an $n\in\N$ \st $d(x,x_n)<\ve$ \fa
$n\ge N$. We write $x_n \nad d\to x$ and $x$~is called a \ti{limit\/} (\wrt the metric~$d$) of the \sq\ $\zb xn\N$.\index{limit of a \sq}
\edf

\bdf7.16
A \nss\ $Y$ of a \ms\ $(X,d)$ is called \ti{bounded\/} (\wrt the metric~$d$) if \te s an $M\in\R_{\ge0}$ \st
\beq7.25
d(x,y)\le M \qh{\fa $x,y\in Y$.}
\e
\edf

The notions of order-boundedness and metric-boundedness in $(\R,\ge)$ are not necessarily \ev t. In Lemma \rf{l7.24} we will introduce a non-trivial
metric~$\wt d$ on~$\R$ which \sf ies $\wt d(x,y)\le 1$ \fa $x,y\in\R$. Obviously, all subsets of~$\R$ are bounded \wrt the metric~$\wt d$.

\bex7.17
Show that Lemma \rf{l7.2} holds for \sq s in a \ms\ $(X,d)$.
\eex

\bxs7.18
We introduced on $\R$ three metrics, namely in \er{7.15}, \er{7.16} and \er{7.22a}, which we denote \rsp\ by $d,d',d''$. We now compare the notions
of \cg nce of \sq s induced by these metrics. Let $\zb xn\N$ be a \sq\ in~$\R$ and let $x\in\R$. Suppose that $x_n \nad{d'} \to x$. By \df, given
$\ve:=\frac12$ \te s an $N\in\N$ \st $|x-x_n|<\frac12$ \fa $n\ge N$. By \er{7.16} we have $|x-x_n|=0$ \fa $n\ge N$, that is, $x_n=x$ for $n\ge N$.
Such a \sq\ $\{x_n\}$ is called \ti{eventually} constant. \If that \fe $\ve\in\R_{>0}$, $d(x,x_n)=d(x,x)=0<\ve$ whenever $n\ge N$. \E\Tf $x_n\nad
d\to x$. Similarly $x_n\nad{d''}\to x$. More generally, if $(X,\wt d)$ is a \ms, then $x_n\nad{d'}\to x$ implies $x_n\nad{\tilde d}\to x$.

Now let $x_n:=\frac1n$, $n\in\N$, and $x=0$. We have $x_n\da 0$ hence $x_n\nad d\to0$. However, $d'(x,x_n)=1$ \fa $n\in\Na$, hence
$d'(x,x_n)\to1\ne0$. \E\Tf $\zb xn\N$ does \ti{not\/} converge to~$0$ \wrt the discrete metric~$d'$. One says that the \cg nce \wrt the metric~$d'$
is \ti{strictly stronger} than the \cg nce \wrt the metric~$d$.

We now compare metrics $d$ and $d''$. It turns out that \fe $x\in\R$ and every \sq\ $\zb xn\N$ in~$\R$, we have
\beq7.26
x_n\nad d\to x \qh{iff \ } x_n\nad{d''}\to x.
\e
Setting $c_n:=d(x,x_n)\ge0$, $n\in\N$, we have to show that $c_n\nad d\to0$ iff $\sqrt{c_n}\nad d \to0$.

``\ti{If\/}'': Given $\ve:=1$ \te s an $N_1\in\N$ \st $\sqrt{c_n}\le1$ \fa $n\ge N_1$. Note that $c_n= \sqrt{c_n}\cdot \sqrt{c_n}$ \fa $n\in\N$,
hence $0\le c_n\le 1\cdot\sqrt{c_n}$ \fa $n\ge N_1$. By \df, given $\ve\in\R_{>0}$ \te s an $N_2\in\N$ \st $\sqrt{c_n}\le \ve$ \fa $n\ge N_2$. Hence
for $n\ge N:=\max(N_1,N_2)$ we have $c_n\le \ve$. Thus $d''(x,x_n)\to0$ implies $d(x,x_n)\to0$, that is, $x_n\nad d\to x$ if $x_n\nad{d''}\to x$.

``\ti{Only if\/}'': Given $\ve\in\R_{>0}$ \te s an $N\in\N$ \st $c_n<\ve^2$ \fa $n\ge N$, since $c_n\nad d\to 0$. \If that given $\ve\in\R_{>0}$,
we have $\sqrt{c_n}<\sqrt{\ve^2}=\ve$ \fa $n\ge \N$ since the self-map $y\mt\sqrt y$ of $\R_{\ge0}$ is strictly in\cre. Hence $\sqrt{c_n}\to0$,
and $x_n\nad{d''}\to x$.
\exs

We now consider {\bf pointwise \cg nce} of \sq s of \el s of the \vs\ $(\R,\R^I)$ \itd in Examples 5.1.2 where $I$~is a \ns. The \co\ maps of~$\R^I$
(see (5.1.5)) are defined by
\beq7.27
\psi_i(\bX):=\bX(i), \q \bX\in\R^I,\ i\in I.
\e

\bdf7.19
A \sq\ $\zb \bX n\N$ of \el s of~$\R^I$ is said to converge \ti{pointwise} to $\bY\in\R^I$ if the \fw\ holds:\index{pointwise convergence}
\beq7.28
\psi_i(\bX_n) \to \psi_i(\bY) \qh{\fa $i\in I$}.
\e
\edf

Note that \fe $i\in I$ $\{\psi_i(\bX_n)\}_{n\in\N}$ is a \sq\ in~$\R$. We shall use the notation
\beq7.29
\bX_n \lm \bY
\e
for \er{7.28}.

Our goal is to investigate whether \te s a metric $D$ on~$\R^I$ \st
\beq7.30
\bX_n \lm \bY \qh{iff }D(\bX_n,\bY)\to0.
\e
We first consider the case where $I:=[1,l]$, $l\in\Na$ (see (1.3.27) where $(E,\le):=(\N,\le)$). Then $\R^{[1,l]}$ is usually denoted by~$\R^l$.

\newpage
\bpr7.20
Let $\zb\bX n\N$ be a \sq\ of \el s of~$\R^l$ and let $\bY\in\R^l$. Then the \fw\ assertions are \ev t\/\dw

\hph i,i, $\bX_n \lm \bY$,

\hph ii,, $\bX_n \nad d\to \bY$ \fa metrics $d$ \itd in \er{7.22}, \er{7.23} and \er{7.24}.
\epr

We will use the \fw\ lemma in the proof of \E\Pr\ \rf{p7.20}.

\blm7.21
Let $\zb an\N$ be a \sq\ in $\R_{\ge0}$ \sf ying $a_n\to0$. Let $\zb bn\N$ be a \sq\ in~$\R$ \sf ying
\beq7.30a
|b_n|\le a_n \qh{\fa $n\ge N$ \fs $N\in\N$}.
\e
Then $b_n\to0$.
\elm

\proof[Proof of Lemma \rf{l7.21}]
Let $\ve\in\R_{>0}$, and let $N_1\in\N$ be \st $0\le a_n<\ve$ \fa $n\ge N_1$ and set $N_2:=\max(N_1,N)$. Then $|b_n-0| = |b_n| \le a_n = |a_n-0| <
\ve$ \fa $n\ge N_2$. Then use (1.3.11).
\endproof

\bex7.22
Let $\bZ\in\R^l$. Prove the \fw\ in\et ies:
\bea7.31
{}& \max_{1\le i\le l} |\bZ(i)| \le \|\bZ\|_2 \le \sum_{i=1}^l |\bZ(i)|,\\
& \sum_{i=1}^l|\bZ(i)| \le \sqrt l\cdot\|\bZ\|_2 \le l\cdot \max_{1\le i\le l}|\bZ(i)|.\lb{7.32}
\e
Hint: Use \er{7.20} with $l$ instead of $2$, and use $\a^2\le\a$ if $0<\a\le1$.
\eex

\proof[Proof of \E\Pr\ \rf{p7.20}]\

\ti{Step} 1: ``\ti{Reduction to the case $\bY=\bz$}.'' Note that \er{7.28} is \ev t to: $\psi_i(\bc_n){\to}\psi_i(\bz)$ \ti{\fa} $1\le i\le l$
where $\bc_n:=\bX_n-\bY$. Indeed, let $1\le i\le l$. Then $\psi_i(\bX_n) \to \psi_i(\bY)$  iff $\psi_i(\bX_n)-\psi_i(\bY)\to \psi_i(\bY) -
\psi_i(\bY) = 0$ by Lemma \rf{l7.6}. \Mo in view of the linearity of~$\psi_i$, we have $\psi_i(\bX_n)-\psi_i(\bY) = \psi_i(\bX_n -\bY)$.
Hence $\psi_i(\bX_n)\to\psi_i(\bY)$ iff $\psi_i(\bc_n)\to0 = \psi_i(\bz)$.
Similarly, we note that $d(\bX_n,\bY) = d(\bX_n-\bY,\boldsymbol0)$ for the three metrics above. For example,
$$
\sum_{i=1}^l |\psi_i(\bX_n)-\psi_i(\bY)| = \sum_{i=1}^l |\psi_i(\bX_n-\bY)| = \sum_{i=1}^l |\psi_i(\bX_n-\bY) - \psi_i(\boldsymbol0)|.
$$
Set $\bZ_n:=\bX_n-\bY$.

\ti{Step} 2: ``\ti{\E\ev ce between pointwise \cg nce and \cg nce \wrt the metric~$d$ defined in \er{7.23}}.'' Suppose that
$\suml_{i=1}^l |\psi_i(\bZ_n)|\to0$. We have $|\psi_k(\bZ_n)| \le \suml_{i=1}^l |\psi_i(\bZ_n)|\to 0$ for $1\le k\le l$. Then $|\psi_k(\bZ_n)| \to 0$,
$1\le k\le l$, by Lemma \rf{l7.21}. Note that $|\psi_k(\bZ_n)-0| = |\psi_k(\bZ_n)|$, hence $\psi_k(\bZ_n)\to0$ for $1\le k\le l$. Hence $\bZ_n$ \cg s
pointwise to~$\boldsymbol0$. Conversely, suppose that $|\psi_k(\bZ_n)|\to0$, $1\le k\le l$. Then \fe $\ve\in\R_{\ge0}$ and every $k\in[1,l]$ \te s
$N_k\in\N$ \st $|\psi_k(\bZ_n)| \le \frac1l \ve$ \fa $n\ge N_k$. Set $N:=\sup\{N_k\in\N: 1\le k\le l\} = \max\{N_k\in\N: 1\le k\le l\}$. Then
$\suml_{k=1}^l |\psi_k(\bZ_n)|\nad{(2.4.41)}\le \suml_{k=1}^l \frac1l \ve \nad{(2.4.40)}= l\cdot \frac1l \ve = \ve$ \fa $n\ge N$. Hence by \er{7.1}
$d(\bZ_n,\boldsymbol0)\to0$, where $d$~is defined in \er{7.23}.

\ti{Step} 3: ``\ti{\E\ev ce between the \cg nce \wrt the metrics defined in \er{7.22}, \er{7.23} and \er{7.24}}.'' This \ev ce is a direct con\sq\
of Lemma \rf{l7.21}, Exercise \rf{ex7.22}, and the fact that if $a_n\to0$, $c\in\R_{>0}$, then $ca_n\to0$.
\endproof

We now turn to the case $I:=\N$. To this end we first introduce a metric on~$\R$ that is bounded in $(\R,\ge)$ and that induces the same notion of
\cg nce as the one induced by the metric defined in \er{7.15}.

\blm7.23
Let $\vf:\R_{\ge0}\to\R_{\ge0}$ be a \f\ \sf ying\dw

\def\labelenumi{(\rm G\theenumi)}
\sdim{(G1)}
\ben
\item \ $\vf(0)=0$ and $\vf(x)>0$ \fa $x\in\R_{>0}$,
\item \ $\vf$ is in\cre,
\item \ $\vf(x+y) \le \vf(x)+\vf(y)$ \fa $x,y\in\R_{\ge0}$,
\item \ $\vf(x)\le M$ \fs $M\in\R_{>0}$ and all $x\in\R_{>0}$.
\een
Then $\wt d:\R \to \R_{\ge0}$ defined by
\beq7.34
\wt d(a,b):=\vf(|a-b|), \q a,b\in\R,
\e
is a \emph{metric} on $\R$ \sf ying \fa $a,b,c\in\R$
\beq7.35
\wt d(a+c,b+c) = \wt d(a,b),
\e
and
\beq7.36
\wt d(a,b)\le M.
\e
\elm

\proof
``\ti{$\wt d$ is a metric}'': (M1) Let $a,b\in\R$. Then if $a=b$, we have $|a-b|=0$ hence $\wt d(a,b)\nde7.32 = \vf(0) \nad{\rm (G1)} = 0$. Conversely,
if $\wt d(a,b)=0$ then $\vf(|a-b|)=0$, hence $|a-b|=0$ by~(G1), and $a=b$.

(M2) $\wt d(a,b) = \vf(|a-b|) \nde5.41 = \vf(|b-a|) = \wt d(b,a)$.

(M3) Let $c\in\R$. We have to show $\wt d(a,b)\le \wt d(a,c)+\wt d(c,b)$, $a,b,c\in\R$. We have $|a-b| = |a+((-c)+c)-b| \nad{(2.1.36)} =
|(a+(-c)) + (c-b)| \nad*\le |a-c|+|c-b|$. \E\Tf $\wt d(a,b) \nad{\rm (G2)}\le \vf(|a-c|+|a-b|) \nad{\rm (G3)}\le \vf(|a-c|)+\vf(|c-b|) = \wt d(a,c)
+\wt d(c,b)$. In~$\nad*\le$ \pp y (M3) for the metric defined in \er{7.15} is used.

\E\pp ies \er{7.35}, \er{7.36} easily follow from \er{7.34} and (G4).
\endproof

\blm7.24
Let $\vf:\R_{\ge0} \to \R_{\ge0}$ be defined by
\beq7.37
\vf(x):=\min(x,1), \q x\in\R_{\ge0},
\e
where
\beq7.38
\min(x,y) := \bca
x &\hbox{if }x\le y,\\
y &\hbox{if }y<x,
\eca
\q x,y\in\R.
\e
Then $\vf$ \sf ies {\rm (G1)--(G4)}, and the metric $\wt d$ defined in \er{7.34} \sf ies \er{7.35}, \er{7.36} with $M:=1$, and
\beq7.39
a_n\nad{\wt d}\to a \qh{iff \ }a_n\nad{d}\to a
\e
\fa \sq s $\zb an\N$ in~$\R$ and all $a\in\R$.
\elm

\proof
``(G1), (G2), (G4)'' are easy to verify.

``(G3)'': \E\wlg we may suppose $x\le y$. If $x+y\le1$, then $x,y\le1$, hence $\vf(x+y) = x+y = \vf(x) + \vf(y)$. Suppose now that $x+y>1$. Then
$\vf(x+y)=1$. If $y>1$, then $\vf(y)=1$ and $\vf(x+y)=1\le \vf(x)+1 = \vf(x)+\vf(y)$. If $y\le1$, then $x\le y\le 1$, hence $\vf(x)=x$, $\vf(y)=y$,
and $\vf(x+y)=1<x+y= \vf(x)+\vf(y)$.

``\er{7.39}'': Setting $c_n:=a_n-a$, $n\in\N$, we have to show $\vf(|c_n|)\to0$ iff $|c_n|\to0$. Note that $\vf(x)\nde7.37 \le x$, $x\in\R_{\ge0}$,
hence if $|c_n|\to0$,
then $\vf(|c_n|)\le |c_n|$, hence $\vf(|c_n|)\to0$ by Lemma \rf{l7.21}. \E\oh if $\vf(|c_n|)\to0$ then there is an $N\in\N$ \st $\vf(|c_n|)\le
\frac12<1$ \fa $n\ge N$. If $n\ge N$ then $\vf(|c_n|)=|c_n|$ by \er{7.37}, hence $|c_n|=\vf(|c_n|)\to0$.
\endproof

\bpr7.25
Let $\bX,\bY\in\R^\N$, let $\psi_i$ be as in \er{7.27}, and let $\wt d$ be as in Lemma \rf{l7.24}. Set
\beq7.40
c_l:= \sum_{i=0}^l \frac1{2^{i+1}} \,\wt d(\psi_i(\bX),\psi_i(\bY)), \q l\in\N.
\e
Then the \fw\ holds\dw

\hph i,i,
\beq7.41
0 \le c_l \le  c_{l+k} \le 1, \q l,k\in\N.
\e
In view of \er{7.41}, $\supl_{l\ge0}c_l$ exists by Theorem \rf{t5.1}. Set
\beq7.42
D(\bX,\bY) := \sup_{l\ge0}c_l.
\e

\hph ii,, $D$ is a metric on $\R^\N$ \sf ying
\bea7.43
{}&D(\bX+\bZ,\bY+\bZ) = D(\bX,\bY), &&\hskip-60pt \bX,\bY,\bZ\in\R^\N,\\
{}& D(\bX,\bY)\le1, &&\hskip-60pt \bX,\bY\in\R^\N. \lb{7.44}
\e
\epr

\proof \

\ti{Step} 1. \ \E\fe $l\in\N$ we set
\beq7.45
D_l(u,v) := \sum_{i=0}^l \frac1{2^{i+1}}\,\wt d(\psi_i(u),\psi_i(v)), \q u,v\in\R^{[0,l]}.
\e

``\ti{$D_l$ is a metric on $\R^{[0,l]}$}'':

(M1) Let $u,v\in\R^{[0,l]}$. If $u=v$ then $\wt d(\psi_i(u),\psi_i(v)) = 0$ for $0\le i\le l$, since $\wt d$ is a metric on~$\R$. Then
$D_l(u,v) = \suml_{i=0}^l 0 \nad{(2.4.39),(2.4.40)} = (l+1)\cdot 0=0$. If $u\ne v$ then \te s $i_0\in[0,l]$ \st $\psi_{i_0}(u)-\psi_{i_0}(v)\ne0$,
hence $|\psi_{i_0}(u)-\psi_{i_0}(v)|>0$, and $\wt d(\psi_{i_0}(u),\psi_{i_0}(v)) = \vf(|\psi_{i_0}(u)-\psi_{i_0}(v)|)>0$ by~(G1). Hence $\frac1{2^{i_0}}
\wt d(\psi_{i_0}(u),\psi_{i_0}(v))>0$. We have $\wt d(\psi_{i}(u),\psi_{i}(v))\ge0$ \fa $i\in[0,l]$ since $\wt d$~is a~metric, hence $\frac1{2^i}
\wt d(\psi_{i}(u),\psi_{i}(v))\ge0$ \fa $i\in[0,l]$, and $\suml_{i=0}^l \frac1{2^i}\wt d(\psi_{i}(u),\psi_{i}(v))>0$ by (2.4.42).

(M2) follows from $\wt d(\psi_i(u),\psi_i(v)) = \wt d(\psi_i(v),\psi_i(u))$ \fa $i\in[0,l]$ and $u,v\in\R^{[0,l]}$ since $\wt d$~is a~metric
on~$\R$.

(M3) Let $u,v,w\in\R^{[0,l]}$. Then $\wt d(\psi_i(u),\psi_i(v)) \le \wt d(\psi_i(u),\psi_i(w))+\wt d(\psi_i(w),\psi_i(v))$ \fa $i\in[0,l]$ since
$\wt d$~is a~metric on~$\R$. \E\Tf $D_l(u,v) = \suml_{i=0}^l \frac1{2^{i+1}}\,\wt d(\psi_i(u),\psi_i(v)) \break\nad{\er{3.10},(2.4.41)} \le
\suml_{i=0}^l \Bg({\frac1{2^{i+1}} \,\wt d(\psi_i(u),\psi_i(w)) + \frac1{2^{i+1}}\,\wt d(\psi_i(w),\psi_i(v))}) \nad{(2.4.38)}=
\suml_{i=0}^l \frac1{2^{i+1}} \wt d(\psi_i(u),\psi_i(w)) + \suml_{i=0}^l \frac1{2^{i+1}} \wt d(\psi_i(w),\psi_i(v)) = D_l(u,w)+D_l(w,v)$.

We claim that \fa $l\in\N$, $u,v,w\in\R^{[0,l]}$ we have:
$$\dsl{
(\rf{7.43}_l) \hfill D_l(u+w,v+w) = D_l(u,v),\hfill  \cr
(\rf{7.44}_l) \hfill D_l(u,v) \le1. \hfill
}
$$

Indeed, $\suml_{i=0}^l \frac1{2^{i+1}}\wt d(\psi_i(u+w),\psi_i(v+w)) \nde7.34 = \suml_{i=0}^l \frac1{2^{i+1}}\vf(|\psi_i(u+w)-\psi_i(v+w))| \nad*=
\break\suml_{i=0}^l \frac1{2^{i+1}}\vf(|(\psi_i(u)+\psi_i(w))-(\psi_i(v)+\psi_i(w))| = \suml_{i=0}^l \frac1{2^{i+1}}\vf(|\psi_i(u)-\psi_i(v))|)$.
In $\nad*=$ we used the \ti{linearity} of~$\psi_i$ (see Examples 5.1.7\,(ii)), $i\in[0,l]$. Then $(\rf{7.43}_l)$ is proved.

We now prove $(\rf{7.44}_l)$. $\suml_{i=0}^l \frac1{2^{i+1}}\wt d(\psi_i(u),\psi_i(v)) \nad*\le \suml_{i=0}^l \frac1{2^{i+1}} <1$ \fa $l\in\N$. In
$\nad*\le$ we used $0\le \wt d(\psi_i(u),\psi_i(v)) \le1$ by Lemma \rf{l7.24}, \er{5.10} and (2.4.41).

\ssk
\ti{Step} 2. \ \E\fa $l\in\N$ and all $\bZ\in\R^\N$ we define an \el\ $P_l\bZ \in \R^{[0,l]}$ by setting
\beq7.46
P_l\bZ(i) := \bca
\psi_i(\bZ) & \hbox{for } 0\le i\le l,\\
0 &\hbox{for } i>l.
\eca
\e
We also define \fa $l\in\N$ and all $u\in\R^{[0,l]}$ an \el\ $J_lu\in\R^\N$ by setting
\beq7.47
J_lu(i) := \bca
u(i) &\hbox{for }0\le i\le l, \\
0 &\hbox{for } i>l.
\eca
\e
We clearly have
\beq7.48
P_l(J_lu) = \psi_l(u), \q l\in\N \hbox{ and } u\in\R^{[0,l]}.
\e
Observe that \fa $\bX,\bY\in\R^\N$ we have
\beq7.49
c_l = D_l(P_l\bX,P_l\bY), \q l\in\N.
\e
Thus $c_l\ge0$, $l\in\N$, by (M1) for the metric $D_l$. \Mo
\bmlg
c_{l+1} = D_{l+1}(P_{l+1}\bX, P_{l+1}\bY) = \suml_{i=0}^{l+1} \frac1{2^{i+1}}
\wt d(\psi_i(P_{l+1}\bX),\psi_i(P_{l+1}\bY))\\ \nad{(2.4.32)} = \suml_{i=0}^{l}\frac1{2^{i+1}}\wt d(\psi_i(P_{l+1}\bX),\psi_i(P_{l+1}\bY))+
\frac1{2^{l+2}}\wt d(\psi_{l+1}(P_{l+1}\bX),\psi_{l+1}(P_{l+1}\bY)) \\ \nad*\ge \suml_{i=0}^{l}\frac1{2^{i+1}}\wt d(\psi_i(P_{l+1}\bX),
\psi_i(P_{l+1}\bY)).
\e
In $\nad*\ge$ we used $\frac1{2^{l+1}}\wt d(\cdot,\cdot)\ge0$, since $\wt d$~is a metric. \Mo $\psi_i(P_{l+1}\bX) =
\psi_i(P_l\bX)$ if $i<l+1$. Similarly, if $\bX$~is replaced by~$\bY$. Hence $c_{l+1} \ge \suml_{i=0}^l \frac1{2^{i+1}}\wt d(\psi_i(P_l\bX),
\psi_i(P_l\bY)) = c_l$. Using \In\ one finds that $c_l\le c_{l+k}$, $l,k\in\N$. Finally, $c_l = D_l(P_l(\bX),P_l(\bY)) \le 1$ by
$(\rf{7.44}_l)$. This completes the proof of \er{7.41}.

\ssk
\ti{Step} 3. \ In view of \er{7.42} and \er{7.49} we have
\beq7.50
D(\bX,\bY) = \sup_{l\ge0} D_l(P_l\bX,P_l\bY), \q \bX,\bY \in\R^\N.
\e
From \er{7.50} and $(\rf{7.44}_l)$ we infer \er{7.44}, since $D(\bX,\bY)$ is the least \ub\ for the set $\{D_l(P_l\bX,P_l\bY)\in\R: l\in\N\}$,
which is bounded by~$1$.

We also prove \er{7.43}. From $(\rf{7.43}_l)$ we obtain $D_l(P_l\bX+P_l\bZ,P_l\bY+P_l\bZ) = D_l(P_l\bX,P_l\bY)$ for $\bX,\bY,\bZ\in\R^\N$ and
$l\in\N$. Note that $P_l\bX+P_l\bZ = P_l(\bX+\bZ)$ and $P_l\bY+P_l\bZ = P_l(\bY+\bZ)$ (see Example 5.1.7\,(ii)). Then $D(\bX+\bZ,\bY+\bZ)
\nde7.50 = \supl_{l\ge0} D_l(P_l(\bX+\bZ),P_l(\bY+\bZ)) \nad{(\rf{7.43}_l)} = \supl_{l\ge0} D_l(P_l\bX,P_l\bY) \nde7.50 = D(\bX,\bY)$.

It remains to show that $D$ is a metric on $\R^\N$.

(M2): $D(\bX,\bY) = \supl_{l\ge0} D_l(P_l\bX,P_l\bY) = \supl_{l\ge0}D_l(P_l\bY,P_l\bX) = D(\bY,\bX)$ since $D_l$~is a metric.

(M1): $D(\bX,\bY) = \supl_{l\ge0} D_l(P_l\bX,P_l\bY) \ge D_1(P_1\bX,P_1\bY) = \frac12\wt d(P_1\bX,P_1\bY) \ge0$.

If $\bX=\bY$ then $P_l\bX = P_l\bY$ \fa $l\in\N$. Then $D(\bX,\bY) = \supl_{l\ge0} D_l(P_l\bX,P_l\bY) = \supl_{l\ge0}0 = 0$. If $\bX\ne\bY$,
there is $i_0\le l$ \st $\psi_{i_0}(\bX)\ne \psi_{i_0}(\bY)$, then $\wt d(\psi_{i_0}(\bX),\psi_{i_0}(\bY)) > 0$ and $\suml_{i=0}^{i_0}
\frac1{2^{i+1}}\wt d(\psi_{i}(\bX),\psi_{i}(\bY)) \ge \frac1{2^{i_0+1}}\wt d(\psi_{i_0}(\bX),\psi_{i_0}(\bY)) > 0$. Hence
$$
D(\bX,\bY) \ge \suml_{i=0}^{i_0} \frac1{2^{i+1}}\,\wt d(\psi_{i}(\bX),\psi_{i}(\bY)) > 0.
$$

(M3) Let $\bX,\bY,\bZ\in\R^\N$. Then $D_l(P_l\bX,P_l\bY) \nad*\le D_l(P_l\bX,P_l\bZ) + D_l(P_l\bZ,P_l\bY) \nde7.50 \le D(\bX,\bZ) +
D(\bZ,\bY)$ \fa $l\in\N$. In~$\nad*\le$ we used the fact that $D_l$~is a metric. Finally, $D(\bX,\bY) = \supl_{l\ge0} D_l(P_l\bX,P_l\bY)
\le D(\bX,\bZ)+D(\bZ,\bY)$. This completes the proof of \E\Pr\ \rf{p7.25}.
\endproof

\bth7.26
Let $D$ denote the metric on $\R^\N$ \itd in \E\Pr\ \rf{p7.25}. Then we have
\beq7.51
\bX_n \lm \bX \qh{iff \ }D(\bX_n,\bX)\to0
\e
\fe \sq\ $\zb\bX n\N$ of \el s of $\R^\N$ and every $\bX$ in~$\R^\N$.
\eth

\proof \

``\ti{Reduction to the case $\bX:=\boldsymbol0$}'': $\bX_n\lm\bX$ iff $\psi_i(\bX_n)\to \psi_i(\bX)$ iff $|\psi_i(\bX_n)-\psi_i(\bX)|\to0$ \fa
$i\in\N$ iff $|\psi_i(\bX_n-\bX)-0| \to0$ iff $|\psi_i(\boldsymbol e_n)-\psi_i(\boldsymbol0)|\to0$ \fa $i\in\N$ where $\boldsymbol e_n:=\bX_n-\bX$.
$D(\bX_n,\bX)\to0$ iff $D(\bX_n+(-\bX),\bX+(-\bX))\to0$ by \er{7.43} iff $D(\boldsymbol e_n,\boldsymbol 0)\to0$.

``\ti{If\/}'': Suppose $D(\be_n,\bz)\to0$. Let $k\in\N$. Then $D_k(P_k\be_n,P_k\bz) \nde7.50 \le D(\be_n,\bz)$ \fa $n\in\N$. \E\Tf
$D_k(P_k\be_n,P_k\bz)
\nde7.45 = \suml_{i=0}^k \frac1{2^{i+1}}\wt d(\psi_i(P_k\be_n),0) \nad{(2.4.41),\er{7.34}}\ge \frac1{2^{k+1}}\vf(|\psi_k(P_k\be_n)|)$ \fa $n\in\N$,
since $\wt d(\cdot,\cdot)\ge0$. Note that $\psi_k(P_k\be_n) = \psi_k(\be_n)$ \fa $n\in\N$. Given $\ve\in\R_{>0}$, $\ve<1$, there is an $N\in\N$ \st
$D(\be_n,\bz)<\frac1{2^{k+1}}\ve$ \fa $n\ge N$. \E\Tf $\vf(|\psi_k(\be_n)|)<\ve<1$ \fa $n\ge N$. From \er{7.37} we infer $|\psi_k(\be_n)|<\ve$. \E\Tf
$\psi_k(\be_n)\to0$. Since $k$~is arbitrary, we obtain $\be_n\lm \bz$.

``\ti{Only if\/}'': Suppose $\be_n\lm \bz$, hence $\psi_i(\be_n)\to0$ \fa $i\in\N$. Let $0<\ve<1$. We have to show that \te s an $N'\in\N$ \st
$D(\be_n,\bz)<\ve$ \fa $n>N'$. Note that $\suml_{i=k}^l \frac1{2^{i+1}}\wt d(\psi_i(\be_n),0) \le \suml_{i=k}^l\frac1{2^{i+1}}$ since $\wt d(\cdot,
\cdot)\le 1$ \fa $l,k\in\N$ \st $k\le l$. Recall that $\suml_{i=k}^l \frac1{2^{i+1}} < \frac1{2^k}$ \fa $l\in\N$. Hence $\suml_{i=k}^l
\frac1{2^{i+1}}\wt d(\psi_i(\be_n),0)\nad{\er{3.10},(2.4.41),(2.4.31)} \le \suml_{i=k}^l \frac1{2^{i+1}} < \frac1{2^k} = 2^{-k}$ \fa $l,k\in\N$ \sf
ying $k\le l$.

We claim that \te s $\ov k\in\Na$ \st $2^{-\bar k} < \frac\ve3$. Suppose, for \cd ion, that $2^{-m}\ge\frac\ve3$ \fa $m\in\Na$. Then $2^m \le
\frac3\ve$ \fa $m\in\Na$. Let $r\in\Q_{>0}$ be \st $\frac3\ve < r <\frac3\ve+1$. Such an~$r$ exists by \E\Pr\ \rf{p4.31} and \E\Pr\ \rf{p4.1}. Then
$2^m\le r$ \fa $m\in\Na$, \cd ing \er{5.77} since $2>1$. \E\Tf \te s $\ov k\in\Na$ \st
\beq7.52
\sum_{i=\bar k}^l \frac1{2^{i+1}}\,\wt d(\psi_i(\be_n),0)<\frac\ve3
\e
\fa $n>N'$ and all $l\in\N$, $l\ge \ov k$.

We now estimate $\suml_{i=0}^{\bar k-1}\frac1{2^{i+1}}\wt d(\psi_i(\be_n),0)$. Note that $\frac1{2^{i+1}}\le1$ for $0\le i\le \ov k-1$, hence
$
\suml_{i=0}^{\bar k-1}\frac1{2^{i+1}}\,\wt d(\psi_i(\be_n),0) \le \suml_{i=0}^{\bar k-1}\wt d(\psi_i(\be_n),0).
$
We want $\suml_{i=0}^{\bar k-1} \wt d(\psi_i(\be_n),0)<\frac\ve3 <\frac13$, hence since $\wt d(\cdot,\cdot)\ge0$ we need $\wt d(\psi_i(\be_n),0)<1$.
If it is the case, then $\wt d(\psi_i(\be_n),0)= |\psi_i(\be_n)|$ by \er{4.33}, \er{4.36}, and $\suml_{i=0}^{\bar k-1} \wt d(\psi_i(\be_n),0) =
\suml_{i=0}^{\bar k-1}|\psi_i(\be_n)| = \ov d(\bX_n,\bz)$ where $\bX_n\in\R^{[0,\bar k-1]}$, $\bX_n(i)= \psi_i(\be_n)$, $0\le i\le \ov k-1$, and
$\ov d$~is the metric on $\R^{[0,\bar k-1]}$ defined in \er{7.23}. From \E\Pr\ \rf{p7.20}, we obtain the \ex\ of $\wt N\in\N$ \st $\suml_{i=0}^{\bar k
-1} |\psi_i(\be_n)|<\frac\ve3$ ($<\frac13$) \fa $n\ge \wt N$.

Finally, setting $N=\max(N',\wt N)$, we have
$$
\suml_{i=0}^l \frac1{2^{i+1}}\wt d(\psi_i(\be_n),0) \nad{(2.4.32)} = \suml_{i=0}^{\bar k-1}
\frac1{2^{i+1}}\wt d(\psi_i(\be_n),0) + \suml_{i=\bar k}^l \frac1{2^{i+1}}\wt d(\psi_i(\be_n),0) \le\frac\ve3 +\frac\ve3 = \frac{2\ve}3
$$
\fa $l\ge N$. Hence we obtain $D_l(P_l\be_n,P_l\bz)<\frac{2\ve}3$ \fa $l\ge N$.
Since $D_l(P_l\be_n,P_l\bz) = c_l$, $l\in\N$, we obtain
\beq7.53
D(\be_n,\bz) = \supl_{l\ge0} D_l(P_l\be_n,P_l\bz) \le \frac{2\ve}3 < \ve.
\e
This completes the proof of the ``only if'' part of the proof, hence the proof of Theorem \rf{t7.26} is complete.
\endproof

\bex7.27
Let $\wh\vf(x) = \frac x{1+x}$, $x\in\R_{\ge0}$.

\hph i,ii, Show that $\wh\vf$ \sf ies (G1)--(G4), and let $\wh d$ denote the metric on~$\R$ defined by \er{7.34} with $\vf:=\wh\vf$.

\hph ii,i, Show that $a_n\nad{\hat d}\to a$ iff $a_n\to a$ \fe \sq\ $\zb an\N$ and $a_n,a\in\R$.

\hph iii,, Setting $\wh D(\bX,\bY) := \supl_{l\ge0} \suml_{i=0}^l \frac1{2^{i+1}} \wh d(\psi_i(\bX),\psi_i(\bY))$, show that $\wh D$~is a metric on
$\R^\N$ \sf ying \er{7.43}, \er{7.44} where $D:=\wh D$.

\hph iv,, Show that $\bX_n \nad D\to \bY$ iff $\bX_n \nad{\hat D}\to \bY$ \fa \sq s $\zb \bX n\N$ and $\bY$ in~$\R^\N$.
\eex

\brs7.28 \

\hph i,ii, If in \E\df\ \rf{d7.19} the set $I$ is \ct e, then \te s a bi\jn\ $b:\N\to I$, and we may replace $\psi_i$ in \er{7.40} by~$\psi_{b(i)}$.
We obtain a metric~$D_b$ on~$\R^I$ \sf ying $\bX_n\lm\bY$ iff $D_b(\bX_n,\bY)\to0$ \fa \sq s $\bX_n$ and $\bY$ in~$\R^I$. \E\Ip we may choose for~$I$
any infinite subset of~$\Q$.

\hph ii,i, It turns out that if $I:=[0,1]_\R$ then there is \ti{no} metric~$D$ on~$\R^I$ \st $\bX_n\lm\bY$ iff $\bX_n\nad D\to \bY$ \fa \sq s
$\zb \bX n\N$ and $\bY\in\R^I$. In this case the pointwise \cg nce is said to be \ti{not metrizable} (see \cite{Kelley}, \cite{Choquet}).

\hph iii,, Observe that all metrics $d$ \itd so far in \er{7.15}, \er{7.16}, \er{7.19}, \er{7.22a}, \er{7.22}, \er{7.23}, \er{7.24}, \er{7.34},
\er{7.37}, \er{7.42} \sf y
\beq7.54
d(\bX+\bZ,\bY+\bZ) = d(\bX,\bY) \qh{\fa $\bX,\bY,\bZ\in\R^I$}
\e
(note that we may identify $\R$ with $\R^{[0]}$) where $+$ denotes the pointwise \ad\ on~$\R^I$ (see Examples 5.1.2). Observe that $(\R^I,+,\bz)$ is
an \ag. \Mo we have for the metrics mentioned above
\beq7.55
d(-\bX,-\bY) = d(\bX,\bY) \qh{\fa $\bX,\bY\in\R^I$.}
\e
\ers

\bex7.28 \

\hph i,i, Prove \er{7.54} and \er{7.55}.

\hph ii,, Let $d$ be the metric on $\R_{>0}$ defined in \er{7.18}. Recall that $\pz2{\R_{>0}}$ is an \ag. Prove
\beq7.56
d(x\cdot z,y\cdot z) = d(x,y) \qh{\fa $x,y,z\in\R_{>0}$,}
\e
and find $x,y,z\in\R_{>0}$ \st $d(x+z,y+z)\ne d(x,y)$.
\eex

So far we have \itd only five metrics on $\R$, namely the metric $|x-y|$ (which we shall call the \ti{usual\/} metric on~$\R$) \er{7.15},
the discrete metric \er{7.16}, the metric $|x-y|^{1/2}$ \er{7.22a}, and the metrics $\min(|x-y|,1)$ \er{7.37} and $\frac{|x-y|}{1+|x-y|}$ (see
Exercise \rf{ex7.27}).

\blm7.29
Let $d$ be a metric on $\R$ \sf ying \er{7.54} and let $\vf:\R\to\R$ be defined by
\beq7.57
\vf(x):= d(x,0), \q x\in\R.
\e
Then \pp ies \er{5.42}, \er{5.44} of \E\df\ \rf{d5.23} hold and the \f~$\vf$ is \ti{even}, that is,
\beq7.58
\vf(-x) = \vf(x) \qh{\fa} x\in\R.
\e
\elm

\bdf7.30
A \f\ $F:\R_{\ge0} \to\R$ is called \ti{convex} (resp.\ \ti{concave}) if \er{6.87} holds \fa $r,s\in\R_{\ge0}$ and all $t\in[0,1]_\R$. (Note that
$(1-t)r+ts\in\R_{\ge0}$ whenever $r,s\in\R_{\ge0}$ and $t\in[0,1]_\R$).
\edf


\bex7.31
Let $\vf:\R_{\ge0}\to\R_{\ge0}$ \sf y $\vf(0)=0$, $\vf$~be in\cre\ and concave, then $\vf(x+y)\le \vf(x)+\vf(y)$,
$x,y\in\R_{\ge0}$.
\eex

\bco7.32
Under the \cn s of Exercise \rf{ex7.31}, the map $\wt d:\R\t\R\to\R$ defined in \er{7.34} is a \emph{metric} on~$\R$ \sf ying \er{7.54}.
\eco

\bex7.33
Prove the \fw\ assertions:

\hph i,ii, Let $f:\R_{\ge0}\to\R$ be concave (resp.\ convex). Then $-f:\R_{\ge0}\to\R$ defined by $(-f)(x):=-f(x)$, $x\in\R_{\ge0}$, is convex (resp.\
concave).

\hph ii,i, Let $f_1,f_2 : \R_{\ge0}\to \R$ be concave. Then $f_1{\land} f_2: \R_{\ge0}\to\R$ defined by $(f_1{\land} f_2)(x):=
\min(f_1(x),f_2(x))$, $x\in \R_{\ge0}$, is concave.

\hph iii,, Let $f:\R_{\ge0}\to\R$ be the \pw\ limit of a \sq\ $\zb fn\N$ of concave (resp.\ convex) \f s on $\R_{\ge0}$, then $f$~is concave (resp.\
convex).
\eex

\bxs7.34 \

\hph i,ii, Let $f(x):=ax+b$, $a,b\in\R$. Then $f$~is both \ti{concave} and \ti{convex}. \E\Ip $f(x):=x$ \sf ies the \as s of Corollary \rf{c7.32}, and
$d(x,y):=|x-y|$, $x,y\in\R$, is the usual metric on~$\R$.

\hph ii,i, Observe that every \f\ $F:\R \hbox{\,(resp.\ $\R_{\ge0}$)} \to\R$ \sf ies \er{6.87} for $t=0$, $t=1$ and $r=s$. \Mo if \er{6.87} holds for
$r<s$, then \er{6.87} also holds for $s<r$, since $1-(1-t) = t$.

Let $d$ be the discrete metric on $\R$. Then the \f~$\vf$ defined in \er{7.57} \sf ies $\vf(0)=0$ and $\vf(x)=1$ for $x>0$. We show that $\vf$~is
\ti{concave}. Indeed, if $r:=0$, $s>0$ and $t\in(0,1)$, we have $\vf((1-t)0+ts)= \vf(ts)=1 = (1-t)\vf(0)+t\vf(s)$. If $r>0$, $s>0$ and $t\in(0,1)$,
then $\vf((1-t)r+ts)=1=(1-t)+t = (1-t)\vf(r)+t\vf(s)$.

\hph iii,, Let $d$ be the metric on $\R$ defined in \er{7.37}, and let $\vf(x):=\min(x,1)$, $x\in\R_{\ge0}$. Both $\vf_1(x):=x$ and
$\vf_2(x):=1$, $x\in\R_{\ge0}$, are concave on~$\R_{\ge0}$, hence $\vf:=\vf_1{\land}\vf_2$ (see Exercise \rf{ex7.33}\,(ii)) is \ti{concave} on
$\R_{\ge0}$.

\hph iv,, Later on, it will be easily shown, by making use of the notion of \ti{\dv} of a \f, that the \f s $f(x):=x^m$,
$x\in\R_{\ge0}$, $m\in\Na$, are convex on $\R_{\ge0}$ and that the \f\ $f(x):=\frac x{1+\la x}$, $\la>0$, is concave on $\R_{\ge0}$. We now show that
the \f\ $f(x)=x^2$ is convex without using the notion of \dv. We first prove \er{6.87} for $t:=\frac12$. Indeed, let $r,s\in\R_{\ge0}$ and
$t:=\frac12$, then $(\frac12r+\frac12s)^2 = \frac14(r+s)^2 = \frac14(r^2+s^2+2rs)$. Since $2rs \le r^2+s^2$, we obtain $(\frac12r+\frac12s)^2 \le
\frac14 (2r^2+2s^2) = \frac12r^2+\frac12s^2$ (see \er{6.88}). By Lemma \rf{l6.26} we find that \er{6.87} holds for all dyadic $t$'s in~$(0,1)$. By
Lemma \rf{l6.12} \te s a \sq\ $t_n\ua t$, $t\in(0,1)$, $t_n\in\Q_2\cap(0,1)$. Hence $t_n\to t$ and $1-t_n\to 1-t$.

Given $r,s\in\R_{\ge0}$ we have $((1-t_n)r+t_ns)^2 \le (1-t_n)r^2+t_ns^2$, $n\in\N$. Then $(1-t_n)r \to (1-t)r$, $t_ns\to ts$, $(1-t_n)r+t_ns \to
(1-t)r+ts$. Set $u_n:=(1-t_n)r+t_ns$, $n\in\N$, and $v:= (1-t)r+ts$. We have $(u_n-v)^2 = (u_n-v)(u_n-v) = |u_n-v|\,|u_n-v|$,  $n\in\N$. Since
$u_n-v \to u-v$, \te s $M\in\R\ge0$ \st $|u_n-v|\le M$. Hence $(u_n-v)^2 \le M|u_n-v|$. Since $|u_n-v|\to0$, we also have $(u_n-v)^2\to0$ by Lemma
\rf{l7.21}. \E\Tf $a_n:=((1-t_n)r+t_ns)^2 \to ((1-t)r+ts)^2=:a$. Set $b_n:=(1-t_n)r^2+t_ns^2$. We have $b_n\to b$ where $b:=(1-t)r^2+ts^2$. From
$a_n\to a$, $b_n\to b$ and $a_n\le b_n$, $n\in\N$, we infer from the \fw\ lemma that $a\le b$, hence \er{6.87} holds \fa $r,s\in\R_{\ge0}$ and
$t\in[0,1]_\R$.
\exs

\blm7.35
Let $\zb xn\N$ and $\zb yn\N$ be \sq s in~$\R$ \sf ying $x_n\le y_n$ \fa $n\in\N$. Suppose that \te\ $x,y\in\R$ \st $x_n\to x$ and $y_n\to y$. Then
$x\le y$.
\elm

\proof
Set $z_n=y_n-x_n$, $n\in\N$. We have $z_n\ge0$, $n\in\N$.
We claim that $y_n-x_n\to y-x$. We could use \er{7.6} but since we did not prove it, we do
it now. We have $(y_n-x_n)-(y-x) = y_n+(-x_n)+(x-y) = (y_n-y)+(x-x_n)$, hence $|(y_n-x_n)-(y-x)| \le |y_n-y| + |x-x_n|$ \fa $n\in\N$. Given $\ve>0$,
\te $N_1,N_2\in\N$ \st $|y_n-y|<\frac\ve2$, $n\ge N_1$, and $|x-x_n|<\frac\ve2$, $n\ge N_2$. Set $N:=\max (N_1,N_2)$, then $|(y_n-x_n)-(y-x)|
<\frac\ve2 + \frac\ve2=\ve$. Hence $z_n:=y_n-x_n \to y-x$. We claim that $y-x\ge0$. Suppose, for \cd ion, that $z:=y-x<0$. Choosing $\ve:=
\frac12(y-x)>0$, we find for $n\ge N$, $z=(z-z_n)+z_n \ge z-z_n$, since $z_n\ge0$. However, since $z_n\le|z_n|$, hence $-z_n\ge -|z_n|$, we find
$z-z_n \ge z-|z_n| > z-\ve = \frac12(y-x) > 0$. Hence $z>\frac12(y-x)>0$, a~\cd ion since $z<0$.
\endproof

\bex7.36
Justify all steps in the proof given above of the convexity of $f(x)=x^2$, $x\ge0$.
\eex

Recall that the map $f:\R_{\ge0} \to\R_{\ge0}$ defined by $f(x):=x^2$ is a bi\jn, and that $g:=f\Inv$ is denoted by $g(x):=\sqrt x=x^{1/2}$ in
Theorem \rf{t5.5}.
In\et y~(G3) is easy to prove for $x\in\R_{\ge0}$. Indeed, from $(\sqrt x+\sqrt y)^2 =(\sqrt x)^2+(\sqrt y)^2+2\sqrt x\,\sqrt y\ge x+y$, we obtain
$\sqrt x+\sqrt y = \sqrt{(\sqrt x+\sqrt y)^2}\ge \sqrt{x+y}$ \fa $x,y\in\R_{\ge0}$. We used the in\cre ness of the map $x\mt \sqrt x$, $x\ge0$.
The \ti{concavity} of the map $g(x)=\sqrt x$, $x\in\R_{\ge0}$, is a con\sq\ of the \fw\ lemma.

\blm7.37
Let $F:\R\,\hbox{\rm(resp.\ $\R_{\ge0}$)} \to \R\,\hbox{\rm(resp.\ $\R_{\ge0}$)}$ be convex $($resp.\ concave$)$ and bi\jc. Let $G:=F\Inv$.
Then $G$~is
concave $($resp.\ convex$)$ if $F$ is in\cre, and $G$~is convex $($resp.\ concave$)$ if $F$~is de\cre.
\elm

\proof
Let $x,y\in\R$ (resp.\ $\R_{\ge0}$) and $t\in(0,1)_\R$ \sf y $F((1-t)x+ty) \le \hbox{(resp.\ $\ge$)}\break (1-t)F(x)+tF(y)$. Set $u:=F(x)$, $v:=F(y)$.
Then $x=G(u)$, $y=G(v)$ and $F((1-t)G(u)+tG(v)) \le \hbox{(resp.\ $\ge$)}\,(1-t)u+tv$. If $F$~is in\cre, then $F$~is strictly in\cre\ in view of
the bi\ji\ of~$F$, and $G$~is also in\cre\ since the \og\ of~$\R$ is total. We obtain
$$
(1-t)G(u) + tG(v) = G\circ F((1-t)G(u) + tG(v)) \le G((1-t)u + tv).
$$
Given $u,v\in\R$ (resp.\ $\R_{\ge0}$), set $x:=G(u)$, $y:=G(v)$, then \er{6.87} holds, hence\break $G((1-t)u+tv)\ge (1-t)G(u)+tG(v)$. The case $F$~is
de\cre\ is similar.
\endproof

\brm7.39
An important class of metrics on~$\Q$ \sf ying \er{7.54} is the class of metrics ``induced'' by a \ti{\vl} of the field~$\Q$ (see \E\df\ \rf{d5.23}).
These metrics are of the form $d(x,y)=\vf(x-y)$, $x,y\in\Q$, where the map $\vf:\Q\to\R$ \sf ies \er{5.42}--\er{5.45}. Note that $\vf(|x|)=\vf(x)$
\fa $x\in\Q$. Observe that the usual metric on~$\Q$ and the metric defined in \er{7.22a} are of this form.  In~\cite{Ostr} (see also~\cite{A2})
Ostrowski gave a \chz\ of \ti{all\/} valuations
of~$\Q$. Examples of valuations of~$\Q$ having no \ext\ to~$\R$ are the \hbox{$p$-adic} valuations $\vf_p$ where $p$~is an arbitrary \Pn.
Assume $\psi$~is a map from~$\Q$ taking values
in~$\R$, and \sf ying \er{5.42}, \er{5.43}. We have $\psi(y)=0$ iff $y=0$ by \er{5.42}. \Mo $\psi(1) = \psi(1\cdot1) \nde5.43 = (\psi(1))^2$, hence
$\psi(1)=1$ since $\psi(1)\ne0$. Since $(\psi(-1))^2 \nde5.43 = \psi((-1)\cdot(-1)) = \psi(1)=1$, and $\psi(-1)\ge0$ by~\er{5.42}, we have $\psi(-1)
=1$, and $\psi(-1)=\psi(1)$. Finally, $\psi(-x)=\psi((-1)\cdot x) \nde5.43 = \psi(-1)\cdot\psi(x)= 1\cdot\psi(x)= \psi(x)$ \fa $x\in\Q$. \E\Tf the
map~$\psi$ is even, and since $\psi(0)=0$, the map~$\psi$ is \cp ly determined by its values on~$\Q_{>0}$. \E\fe $r\in\Q_{>0}$ \te s \ooo pair
$(c,d)\in\Na\t\Na$ \sf ying $r=\frac cd$ and $\gcd(c,d)=1$ in view of Lemma 4.5.46. Hence $\psi(r)\cdot\psi(\frac d1)\nde5.43 = \psi(\frac cd \cdot
\frac d1) = \psi(\frac c1)= \psi(c)>0$. \If that $\psi(r) = \psi(c)\cdot(\psi(d))\mo$ since $\psi(d) = \psi(\frac d1)\in \R_{>0}$ and $\pz2{\R_{>0}}$
is an \ag. \csq, $\psi(r)$ is \cp ly determined by $\psi(c)$ and~$\psi(d)$, $c,d\in\Na$. If $r=1$, then $\psi(r)=1$. We consider three cases:
$c>1$ and $d=1$, $c=1$ and $d>1$, $c,d>1$.

``$c>1$, $d=1$'': \ In view of \E\Pr\ 2.4.18, we have $c = \prodl_{p\in\cP, p|c}p^{m_p}$, where $m_p\in\Na$ \sf ies $p^{m_p}|c$ and $p^{m_p+1}\nmid
c$ ($m_p$~is uniquely defined). We now apply \E\Pr\ 2.4.18\,(iv) with $(X,\qu,e) := \pz2{\Na}$, $(\wt X,\tqu,\wt e) := \pz2{\R_{>0}}$, and
$I:= \{p\in\cP: p|c\}$. Note that \er{5.43} implies (2.4.57). We obtain $\psi(c) = \prodl_{p\in\cP, p|c} \psi(p^{m_p})$ in~$\R_{>0}$. Applying
Lemma 2.1.22, we obtain
\beq7.59
\psi(c)=\prod_{p\in\cP, p|c} (\psi(p))^{m_p},
\e
where $m_p$ is the \ml icity of~$p$ in the \fc\ (3.2.4) of~$c$ in powers of primes.

``$c=1$, $d>1$'': Setting $s:=\frac d1\in\Q_{>0}$ we infer from \er{7.59} that $\psi(s)=\prodl_{p\in\cP,p|d}\psi(p)^{n_p}$ where $n_p\in\Na$ \sf ies
$p^{n_p}|d$ and $p^{n_p+1}\nmid d$. Since $1:=\psi(1) = \psi(\frac d1\cdot \frac1d) \nde5.43 = \break \psi(\frac d1)\cdot\psi(\frac1d)$, we obtain
$\psi(\frac1d)=\psi(\frac d1)\mo$, the inverse of $\psi(\frac d1)$ in the \ag\ $\pz2{\R_{>0}}$. Thus $\psi(\frac1d) = (\psi(s))\mo =
\Bigl(\prodl_{p\in\cP,p|d}\psi(p)^{n_p}\Bigr)\mo$. Note that $\INV$, the self-map of $\R_{>0}$ defined by $\INV(x):=x\mo$, $x\in\R_{>0}$, is a
monoid \is sm of the group $\pz2{\R_{>0}}$, hence $\Bigl(\prodl_{p\in\cP,p|d}\psi(p)^{n_p}\Bigr)\mo =
\prodl_{p\in\cP,p|d}\bigl(\psi(p)^{n_p}\bigr)\mo$
in view  of \E\Pr\ 2.4.18\,(iv). \E\Tf we obtain
\beq7.60
\psi\Bigl(\frac1d\Bigr) = \prod_{p\in\cP,p|d}\bigl(\psi(p)^{n_p}\bigr)\mo,
\e
where $n_p$ is the \ml icity of~$p$ in the \fc\ (3.2.4) of~$d$ in powers of primes.

``$c>1$ and $d>1$'': Since $r=\frac cd = \frac c1\cdot\frac 1d$, we have $\psi(r) = \psi(\frac c1\cdot\frac1d) \nde5.43 = \psi(\frac c1)\cdot
\psi(\frac1d)$, hence we obtain from \er{7.59}, \er{7.60}:
\beq7.61
\psi(r) = \Bigl(\prod_{p\in\cP,p|c}\psi(c)^{m_p}\Bigr) \cdot \Bigl(\prod_{p\in\cP,p|d}\psi(d)^{n_p}\Bigr)\mo.
\e
Observe that in \er{7.61}, the sets $\{p\in\cP:p|c\}$ and $\{p\in\cP: p|d\}$ are \ti{disjoint\/} since $\gcd(c,d)=1$. Indeed, the set $\{p\in\cP:
p|d\}$ (resp.\ $\{p\in\cP: p|c\}$) is empty whenever the set $\{p\in\cP: p|c\}$ (resp.\ $\{p\in\cP: p|d\}$) is nonempty.

We now consider two cases where the map~$\psi$ is defined by means of formulae \er{7.59}--\er{7.61}.

``$\psi(p):=p$, $p\in\cP$'': If $r:=c$, $c>1$, then $c\nad{(3.2.4)} = \prodl_{p\in\cP,p|c}p^{m_p}$, hence $\psi(c) = \prodl_{p\in\cP,p|c}
(\psi(p))^{m_p} = \prodl_{p\in\cP,p|c}p^{m_p} = c$.

If $r:=\frac1d$, $d>1$, then $\psi(r) \nde7.60 = \prodl_{p\in\cP,p|d}(\psi(p)^{r_p})\mo = \Bigl(\prodl_{p\in\cP, p|d}p^{n_p}\Bigr)\mo
\nad{(3.2.4)}= d\mo = \frac1d = r$.

Finally, if $r:=\frac cd$, $c,d>1$ and $\gcd(c,d)=1$, we have $\psi(r) \nde7.61 = c\cdot d\mo = \frac c1\cdot \frac1d = \frac cd = r$. \E\Tf
$\psi(r)=r$ \fa $r\in\Q_{>0}$. Since $\psi(0)=0$ and $\psi(1)=1$, we obtain $\psi(r) = \psi(|r|) = |r|$ \fa $r\in\Q$. One verifies that $\psi$
\sf ies \er{5.42}--\er{5.45}. Hence $d(x,y):=|x-y|$, $x,y\in\Q$, the induced metric~$d$ is the \rt ion to~$\Q$ of the usual metric of~$\R$.

We now consider the case of a family of maps $\vf_p$, $p\in\cP$, defined by \er{7.59}--\er{7.60} and by
\beq7.62
\vf_p(p') := \bca
\frac1p & \hbox{if } p'=p,\\
1 & \hbox{if } p'\ne p,\ p'\in\cP.
\eca
\e
If $r=\frac c1$, $c>1$, we infer from \er{7.59} that $\vf_p(c) = \prodl_{p'\in\cP,p'|c} (\vf_p(p'))^{m_{p'}}$, where $m_{p'}$ is the \ml icity of~$p'$
in the \fc\ (3.2.4) of~$c$.

If $p\nmid c$, then $\vf_p(c)=1$ in view of (2.4.36). If $p|c$, then $\vf_p(c) = \prodl_{p'\in\cP\sms p,p'|c}\vf_p(p')^{m_{p'}} \cdot
\prodl_{p'\in\{p\}}\vf_p(p')^{m_{p'}}$ by (2.4.59). The first factor is equal to $\prodl_{p'\in\cP\sms p,p'|c}1^{m_{p'}}=1$ by (2.4.34) since
$1^{m_{p'}}=1$. \E\oh the second factor is equal to $\vf_p(p)^{m_p}$ by (2.4.9), hence $\vf_p(c)=(\frac1p)^{m_p}$ by \er{7.62}. In view of (3.2.4)
with $n:=c$, we find that $c\cdot\vf_p(c)= \prodl_{p'\in\cP\sms p,p'|c}{p'}^{m_{p'}}$. Note that $p$~does not divide the \RHS\ of this \et y, which we
denote by~$c'$. We obtain
\beq7.63
c\cdot \vf_p(c)= c',\q c'\in\Na \hbox { \st} p\nmid c'.
\e
Observe that \er{7.63} also holds when $p\nmid c$ since in this case $\vf_p(c)=1$. We also showed
\beq7.64
(\vf_p(c))\mo = \bca
1 &\hbox{if }p\nmid c,\\
p^{m_p} &\hbox{if }p|c.
\eca
\e

Now let $c=c_1c_2$, $c_1,c_2\in\Na$. Then $c_1c_2\vf_p(c_1c_2)=c'$, $p\nmid c'$, $c_i\vf_p(c_i)=c_i'$, $p\nmid c_i$, $i=1,2$. \E\Tf $c_1c_2\vf(c_1)
\vf(c_2) = c_1'c_2'$. Since $\vf_p(c_1c_2)\in\Na$, we obtain $c_1c_2 = c'(\vf_p(c_1c_2))\mo$. Hence $\vf_p(c_1c_2)\mo \vf_p(c_1)\vf_p(c_2)c' =
c_1'c_2'$, and
\beq7.65
\vf_p(c_1c_2)\mo c' = \vf_p(c_1)\mo \vf_p(c_2)\mo c_1'c_2'.
\e
By \er{7.64} $\vf_p(c_1c_2)\mo = p^k$, $\vf_p(c_1)\mo \vf_p(c_2)\mo= p^l$ \fs $k,l\in\N$. From Euclid's Lemma 3.1.30 we infer that $p\nmid c_1'c_2'$.
\Mo from $p^kc' = p^l(c_1'c_2')$ we infer $k=l$. Indeed, if $k>l$ then, dividing both sides by~$p^l$, we find $p^{k-l}c' = c_1'c_2'$, hence
$p^{k-l}|c_1'c_2'$, a~\cd ion. The case $k<l$ is similar. \If that $\vf_p(c_1c_2)\mo = \vf_p(c_1)\mo \cdot\vf_p(c_2)\mo$, hence we obtain
\beq7.66
\vf_p(c_1c_2) = \vf_p(c_1)\vf_p(c_2).
\e
If $r:=\frac1d$, $d>1$, $d=\prodl_{p\in\cP,p|d}p^{n_p}$, then
$$
\vf_p(d) \nde7.64 = \bca
1 &\hbox{if }p\nmid d,\\
p^{n_p} &\hbox{if }p|d.
\eca
$$
Then
$$
\vf_p\Bg(\frac1d) \!\nad{\er{7.60},\er{7.62}} = \!\!\!\!\prod_{p'\in\cP,p'|d}\!(\vf_p(p')^{n_{p'}})\mo =
\biggl(\prod_{p'\in\cP,p'|d}\!\vf_p(p')^{n_{p'}}\biggr)\mo \nde7.59 = (\vf_p(d))\mo  = \bca
1 &\hbox{\hskip-8pt if }p\nmid d,\\
p^{n_p} &\hbox{\hskip-8pt if }p|d.
\eca
$$
Hence we have
\beq7.67
\vf_p\biggl(\frac1d\biggr)\mo = \bca
1 &\hbox{if }p\nmid d,\\
p^{-n_p} &\hbox{if }p|d.
\eca
\e
Thus if $p|d$ we obtain
$$
\frac1d \,\vf_p\biggl(\frac1d\biggr)\mo = \frac1{\prodl_{p'\in\cP}p'{}^{n_{p'}}} \cdot \frac1{p^{-n_d}} = \frac1{\prodl_{p'\in\cP\sms p,p'|d}
p'{}^n_{p'}} = \frac1{d'},
$$
where $d'\in\Na$, $p\nmid d'$.

\E\Tf the analogue of \er{7.63} is
\beq7.68
\frac1d \,\vf_p\biggl(\frac1d\biggr) = \frac1{d'}\,.
\e
If $p\nmid d$, then $\vf_p(d)\mo=1 =\vf_p(d)$, and $d'=d$. Hence \er{7.68} also holds in this case.

\Wanp prove that \fa $r\in\Q_{>0}$ and all $p$~prime, the \fw\ holds:
\beq7.69
r \vf_p(r) = r',
\e
where
\beq7.70
r' = \frac ab \q \hbox{\fs $a,b\in\N$ \sf ying $p\nmid a$, $p\nmid b$,}
\e
and
\beq7.71
\vf_p(r)= p^z \qh{\fs $z\in\Z$.}
\e
Let $p\in\cP$ and let $r\in\Q_{>0}$. Let $c,d\in\Na$ be \st $r=\frac cd$ and $\gcd(c,d)=1$. The \nm s $c,d$ are well and uniquely determined by
Lemma 4.5.46. We distinguish three cases:

\hph i,ii, $p\nmid c$ and $p\nmid d$,

\hph ii,i, $p|c$ and $p\nmid d$,

\hph iii,, $p\nmid c$ and $p|d$.

(i) If $c=d=1$, then $r=1$, $\F(r)=1$ and \er{7.69}--\er{7.71} hold with $a=b=1$, $z=0$. If $c>1$, $d=1$, then $r=\frac c1=c$, and $\vf_p(c)
\nde7.59 = \prodl_{p'\in\cP,p'|c}\vf_p(p') \nde7.62 = \prodl_{p'\in\cP,p'|c} 1 \nad{(2.4.36)} = 1$. Then \er{7.69}--\er{7.71} hold with $a=c$,
$d=1$, $z=0$. If $c=1$, $d>1$, then $r=\frac1d$. We have $\prodl_{p'\in\cP,p'|d}\vf_p(p') = \prodl_{p'\in\cP,p'|d}1=1$. \Mo $\vf_p(\frac1d)
\nad{\er{7.60},\er{7.62}} = \prodl_{p'\in\cP,p'|d}\bigl((\vf_p(p'))^{n_{p'}}\bigr)\mo =1$. Then \er{7.69}--\er{7.71} hold with $a=1$, $b=d$,
$z=0$. If $c>1$, $d>1$, then $\vf_p(\frac1d) \nde7.61 = 1\cdot 1\mo =1$, and \er{7.69}--\er{7.71} hold with $a=c$, $b=1$, $z=0$.

(ii) Since $p|c$, we have $c>1$ and $p\nmid d$ (since $\gcd(c,d)=1$). Thus $\vf_p(\frac cd) \nde7.61 =\break \vf_p(c)\cdot 1\mo = \vf_p(c) \nde7.64 =
\frac1{p^{m_p}}$. Then \er{7.69}--\er{7.71} hold with $a=c\cdot\frac1{p^{m_p}}$, $b=d$, $z=-m_p$. (\E\Ip if $c=p$, $d=1$, then $r=p$ and $\vf_p(r)
=\frac1p$.)

(iii) If $c>1$, $p\nmid c$, $p|d$ (hence $d>1$), then $\vf_p(\frac cd) = 1\cdot \vf_p(\frac1d) = (\vf_p(d))\mo = p^{n_p}$. Then \er{7.69}--\er{7.71}
hold with $a=c$, $b=d\cdot\frac1{p^{n_p}}$, $z=n_p$.

Summarizing, we find that
\beq7.72
z:=\bca
0 &\hbox{if } p\nmid c \hbox{ and } p\nmid d,\\
-m_p &\hbox{if } p|c \hbox{ and } p\nmid d,\\
n_p &\hbox{if } p\nmid c \hbox{ and } p|d,
\eca
\e
where
\beq7.73
m_p \hbox{ (resp.\ $n_p$) is the multiplicity of $p$ in (3.2.4) with $n:=c$ (resp.\ $d$).}
\e
\erm

Note that \fe $p\in\cP$ the set $H_p := \{p^w\in\Q: w\in\Z\}$ is a subgroup of the \ag\ $\pz2{\Q_{>0}}$. As a group $\pz2{H_p}$ is cyclic (see \E\df\
4.3.77). \Mo \fe $p\in\cP$, the set $K_p:= \{r'\in\Q: r'=\frac ab,\ a,b\in\Na,\ p\nmid a,\ p\nmid b\}$ is a subgroup of $\pz2{\Q_{>0}}$. Indeed,
$1=\frac11 \in H_p$, for all $p\in\cP$. \Mo if $r'\in H_p$, then $(r')\mo \in H_p$, since $r'=\frac ab$, $p\nmid a$, $p\nmid b$, implies $(r')\mo =
\frac ba$, $p\nmid b$, $p\nmid a$. Finally, if $r',s'\in H_p$, where $r'=\frac ab$, $s'=\frac cd$, $p\nmid a$, $p\nmid b$, $p\nmid c$, $p\nmid d$,
then $r'\cdot s' = \frac ab\cdot \frac cd = \frac{ac}{bd}$. We have $p\nmid ac$ by Euclid's Lemma 3.1.30. Similarly $p\nmid bd$. Thus $K_c$ is a
subgroup of $\pz2{\Q_{>0}}$ by Lemma 4.3.3. Clearly, $H_p\cap K_p = \{1\}$. Indeed, suppose $p^w\in K_p$, $z\in\Z_+$, then $p^w=\frac ab$,
$a,b\in\Na$, $p\nmid a$, $p\nmid b$. Hence $a=\frac ab \cdot b = p^w\cdot b$. Hence $p|p^w$, $p^w|p^w \cdot b$, hence $p|a$, a \cd ion. Similarly, if
$z\in\Z_-$, we find $p|b$, a~\cd ion. \csq, if $r_1,r_2\in H_p$, $s_1,s_2\in K_p$ and $r_1s_1 = r_2s_2$, then $r_1=r_2$ and $s_1=s_2$. Indeed, from
$r_1s_1 = r_2s_2$, we find $r_1r_2\mo = s_2s_1\mo \in H_p\cap K_p=\{1\}$. Hence $r_1r_2\mo=1$ and $s_2s_1\mo =1$, and $r_1=r_2$, $s_1=s_2$.
Using the fact that the map $z\mt p^w$ is in\jc, and rewriting \er{7.69} as $r = p^{-z}r'$ \fe $r\in\Q_{>0}$ where $p^{-z}\in H_p$ and $r'\in K_p$,
we find that \fe $p\in\cP$ and every $r\in\Q_{>0}$ the \fw\ holds:
\bga7.74
\hbox{\te s \ooo $z\in\Z$ \st $rp^z\in K_p$.}\\ \hbox{\Mo $\vf_p(r) = p^z$ where $z$ \sf ies \er{7.72}.}\non
\e
\Wanp prove that \fe $p\in\cP$ and all $r,s\in\Q$ we have
\beq7.75
\vf_p(rs) = \vf_p(r)\vf_p(s).
\e
Indeed, \fa $r,s\in\Q_{>0}$
we have $r\vf_p(r)=r'\in K_p$, $s\vf_p(s)=s'\in K_p$ and $rs\vf_p(rs)=t'\in K_p$. Hence $rs\vf_p(r)\vf_p(s) = r's'\in K_p$. Thus $(\vf_p(r)
\vf_p(s))\mo r's' = (\vf_p(rs))\mo t'$. Note that $\vf_p(r), \vf_p(s), \vf_p(rs)\in H_p$. From the ``\uq'' part of \er{7.74} we infer $(\vf_p(r)
\vf_p(s))\mo = (\vf_p(rs))\mo$, hence \er{7.75} holds since $(x\mo)\mo=x$ \fa $x\in\Q_{>0}$.

We have thus proved that the maps $\vf_p$ \sf y \er{5.43} \fa $r,s\in\Q_{>0}$. Setting $\vf_p(-r):=\vf_p(r)$, $r\in\Q_{>0}$, we find that \er{7.75}
holds \fa $r,s\in\Q\sms0$, since $\vf_p(rs) = \vf_p(|rs|) = \vf_p(|r|\,|s|) = \vf_p(|r|)\cdot \vf_p(|s|) = \vf_p(r)\vf_p(s)$. Setting $\vf_p(0):=0$,
we find that \er{7.75} holds \fa $r,s\in\Q$, since $0r=0s=0$ \fa $r,s\in\Q$.

It turns out that the maps $\vf_p$ also \sf y \er{5.44}. Let $r_1=(\vf_p(r_1))\mo r_1'$,
$r_2=(\vf_p(r_2))\mo r_2'$, $r_1',r_2'\in K_p$. We have $r_i'=\frac{a_i}{b_i}$, $a_i,b_i\in\Na$, $p\nmid a_i$, $p\nmid b_i$, $i=1,2$, hence
\beq7.76
r_1+r_2 = \frac{(\vf_p(r_1))\mo a_1}{b_1} + \frac{(\vf_p(r_2))\mo a_2}{b_2} = \frac{(\vf_p(r_1))\mo a_1b_2 + (\vf_p(r_2))\mo a_2b_1}{b_1b_2}\,.
\e
By \er{7.74} there is \ooo $z_i\in\Z$ \st $\vf_p(r_i)=p^{z_i}$ for $i=1,2$.

We consider three cases: (i) $\vf_p(r_1)\mo > \vf_p(r_2)\mo$, (ii)~$\vf_p(r_2)\mo > \vf_p(r_1)\mo$ and (iii)~$\vf_p(r_1)\mo = \vf_p(r_2)\mo$.

(i) $\vf_p(r_2)>\vf_p(r_1)$ hence $z_2>z_1$, $z_1,z_2\in\Z$, $w\mt p^w$ is strictly in\cre. We obtain
$$
r_1+r_2 \nde7.76 = \frac{p^{-z_1}a_1b_2 + p^{-z_2}a_2b_1}{b_1b_2} = \frac{p^{z_2-z_1}a_1b_2 + a_2b_1}{b_1b_2}\cdot p^{-z_2}.
$$
Since $z_2-z_1\ge1$, we have $p|p^{z_2-z_1}a_1b_2$, and since $p\nmid a_2b_1$, $p$~does not divide $p^{z_2-z_1}a_1b_2+a_2b_1$. \Mo $b_1b_2\in K_p$,
hence $\vf_p(r_1+r_2) = p^{z_2} = \vf_p(r_2)$. Hence if $\vf_p(r_2)>\vf_p(r_1)$, then $\vf_p(r_1+r_2)= \max(\vf_p(r_1),\vf_p(r_2))$.

(ii) Exchanging the role of $r_1$ and $r_2$, we obtain from case~(i) $\vf_p(r_1+r_2)=\vf_p(r_1)=\max(\vf_p(r_1),\vf_p(r_2))$.

(iii) We have $r_1+r_2 = \frac{a_1b_2+a_2b_1}{b_1b_2}\cdot p^{-z_1}$ and $z_1=z_2$. From \er{7.75} we infer $\vf_p(r_1+r_2) = \vf_p(a_1b_2+a_2b_1)
\cdot(\vf_p(b_1b_2))\mo \cdot\vf_p(p^{-z_1})$. (Recall that $\vf_p(t)\vf_p(t\mo)= \vf_p(t\cdot t\mo)=\vf_p(1)=1$ \fa $t\in\Q_{>0}$ by \er{7.75}.)
Since $b_1b_2\in K_p$, we have $\vf_p(b_1b_2)=1 = (\vf_p(b_1b_2))\mo$. \Mo $\vf_p(p^{-z_1}) = (\vf_p(p^{z_1}))\mo = ((\vf_p(p))^{z_1})\mo =
\bigl((\frac1p)^{z_1}\bigr)\mo = \bigl(\frac1{p^{z_1}}\bigr)\mo = p^{z_1} = \vf_p(r_1) = \vf_p(r_2) = \max(\vf_p(r_1),\vf_p(r_2))$. From \er{7.71},
\er{7.72}, \er{7.73}, we find that $\vf_p(a_1b_2+a_2b_1)$ equals~$1$ if $p\nmid a_1b_2+a_2b_1$, and is less than~$1$ $({}=p^{-m_p})$ if
$p|a_1b_2+a_2b_1$. Hence $\vf_p(r_1+r_2) \le \max(\vf_p(r_1),\vf_p(r_2))$. \If that \fa $r,s\in\Q_{>0}$ and all $p\in\cP$ we have
\beq7.76a
\vf_p(r+s) \le \max(\vf_p(r),\vf_p(s)) \le \vf_p(r)+\vf_p(s).
\e

We now consider the general case $r,s\in\Q$. Note that the cases $r=0$, $s=0$, $r+s=0$ are trivial since $\vf_p(0)=0$. If $r,s\in\Q_{<0}$, then
$\vf_p(r+s) = \vf_p((-|r|)+(-|s|)) =\break \vf_p(-(|r|+|s|)) = \vf_p(|r|+|s|)$. Hence \er{7.76a} holds if $r,s\in\Q_{<0}$. If $r\in\Q_{>0}$ and $s\in
\Q_{<0}$, it is \sft\ to consider the case $r>|s|$. Indeed, if $|s|>r$ then $\vf_p(r+s) = \vf_p(r-|s|) = \vf_p(|s|-r)$. The case $r\in\Q_{<0}$ and
$s\in\Q_{>0}$ is similar. \E\Tf it is \sft\ to consider the case $\vf_p(r_1-r_2)$ where $r_1,r_2\in\Q_{>0}$ and $r_1>r_2$. In that case we have to
replace the $+$~sign by the $-$~sign in \er{7.76}. Then in case~(i) we observe that $p^{z_2-z_1}a_1b_2 - a_2b_1 \in\Na$ and $p\nmid p^{z_2-z_1}
a_1b_2 - a_2b_1$. Hence $\vf_p(r_1-r_2) = p^{z_2} = \vf_p(r_2) = \max(\vf_p(r_1),\vf_p(r_2))$. In case~(ii) we have
$$
r_1-r_2 = \frac{a_1b_2 - p^{z_1-z_2}a_2b_1}{b_1b_2} \cdot p^{-z_1}.
$$
Note that $a_1b_2 - p^{z_1-z_2}a_2b_1\in\Na$ and $p\nmid a_1b_2 - p^{z_1-z_2}a_2b_1$. Hence $\vf_p(r_1-r_2) = p^{z_1} = \vf_p(r_1) \le
\max(\vf_p(r_1),
\vf_p(r_2))$. Finally, in case (iii), we have $r_1-r_2 = \frac{a_1b_2-a_2b_1}{b_1b_2}\cdot p^{-z_1}$ with $z_1=z_2$. Then $a_1b_2-a_2b_1\in\Na$ and
$p\nmid a_1b_2-a_2b_1$, hence $\vf_p(r_1-r_2) = \vf_p(r_1) = \vf_p(r_2)$, hence $\vf_p(r_1-r_2) \le \max(\vf_p(r_1),\vf_p(r_2))$. \E\csq, we have
proved that \er{7.76a} holds \fa $r,s\in\Q$.

Since the maps $\vf_p$ \sf y all \cn s \er{5.42}--\er{5.45} (note that $\vf_p(p)=\frac1p$, $p\in\cP$), these maps are $\Q$- (resp.~$\R$-) \ti{valued
valuations} of the field~$\Q$. The valuations $\vf_p$ \sf y $\vf_p(n)\le1$ \fa $n\in\Na$. Indeed, $\vf_p(1)=1$, $\vf_p(n)\le1$ implies $\vf_p(n+1)
\le1$ \fa $n\in\Na$. A~valuation \sf ying $\vf_p(n)\le1$ \fa $n\in\N$ is usually called \ti{non-\Ar} (see \cite[p.~194]{A2}). If $p=2$ then
$\vf_2(n)=1$ if $n$~is odd and $\vf_2(2^k)=2^{-k}$, $k\in\N$. It is shown in \cite[p.~194]{A2} that an $\R$-valued valuation on an \of~$K$ is
non-\Ar\ iff \er{7.76a} holds.

\If from \er{7.76a} that the induced metric defined by $d(x,y):=\vf_p(|x-y|)$, $x,y\in\Q$, \sf ies instead of~(M3) the stronger \cn:
$$
d(x,y) \le \max(d(x,z),d(z,y)) \qh{\fa $x,y,z\in\Q$.} \leqno{({\rm M}3')}
$$
A metric $d$ on a \ns\ $X$ \sf ying \cn s (M1), (M2) and (M$3'$) is sometimes called an \ti{ultrametric}.\index{ultrametric} Note that if a valuation~$\vf$ on~$K$ \sf
ies \er{7.76a}, then \fa $x,y\in K$ \st $\vf(x)>\vf(y)$ we have $\vf(x+y)=\vf(x)$. Indeed, suppose for \cd ion that $\vf(x+y)<\vf(x)$. Then $\vf(-y)
=\vf(y) < \vf(x)$. From \er{7.76a} we obtain $\vf(x)= \vf(x+y+(-y)) \le \max(\vf(x+y),\vf(-y)) < \vf(x)$, hence $\vf(x) < \vf(x)$, a~\cd ion
(see \cite[p.~195]{A2}). Thus, if a metric~$d$ is induced by a valuation \sf ying \er{7.76a}, then \fa $x,y,z\in K$ we have either $d(y,x) =
d(y,z)$, or $d(y,x)\ne d(y,z)$ and $d(y,x)=d(x,z)$.

\bex7.40
\E\fe $p$ prime, give an example of $r,s\in\Q_{>0}$ \st $r\ne s$, $\vf_p(r) = \vf_p(s) \ne1$, and $\vf_p(r+s)< \max(\vf_p(r),\vf_p(s))$.
\eex

In order to state Ostrowski's theorem we need to define $x^z$ for $x\in\Q_{>0}$ and $z\in\R$. Our next goal is to define~$x^z$ \fa $x\in\R_{>0}$ and
all $z\in\R$. Recall that for $x\in\R$, $x<1$, and $z\in\Q$ we have $x^z = \bigl((x\mo)\mo)^z \nad{\er{5.71}\,\rm RI6}= \bigl((x\mo)^z\bigr)\mo$ since
$x\mo>1$.

This formula motivates the \fw

\bdf7.41
Let $x\in\R_{>0}$, $z\in\R$. Set
\beq7.77
x^z := \bca
\wh\psi_x(z) &\hbox{for }x>1,\\
1 &\hbox{for }x=1,\\
\bigl(\wh\psi_{x\mo}(z)\bigr)\mo &\hbox{for }x<1.
\eca
\e
\edf

We first show that \E\df\ \rf{d7.41} of~$x^z$ is consistent with previous ones. We claim that \E\df\ \rf{d7.41} is consistent with Notation
\rf{n6.51} and with the \df\ $x^r:=\psi_x(r)$, $x\in\R_{>0}$, $r\in\Q$, in \E\Pr\ \rf{p5.31}. Indeed, we have $\wh\psi_x(r)=\psi_x(r)$ whenever
$x\in\R_{>1}$ \fa $r\in\Q$, since $\psi_x(r)=\supl_{s<r,s\in\Q}\psi_x(s) = \wh\psi_x(r)$ in view of Theorem \rf{t6.50}\,(i). Thus if $x\in\R_{>1}$
and $r\in\Q$ we have $x^r\nde7.77 = \wh\psi_x(r) = \psi_x(r) = x^r$ defined in \E\Pr\ \rf{p5.31}. If $x:=1$, $r\in\Q$, then $1^r\nde7.77 = 1
\nad{\er{5.71}\,\rm RI4} = \psi_x(r) = 1^r$ defined in \E\Pr\ \rf{p5.31}. Finally, if $x<0$, $r\in\Q$, $x^r \nde7.77 = (\wh\psi_{x\mo}(r))\mo =
(\psi_{x\mo}(r))\mo \nad{\er{5.71}\,\rm RI6} = \bigl((\psi_x(r))\mo\bigr)\mo = \psi_x(r) = x^r$ defined in \E\Pr\ \rf{p5.31}.

The next lemma is an \ext\ of \er{5.71}\,RI6 and \er{5.71}\,RI7 to the case $x\in\R_{>0}$ and $y\in\R$. In order to emphasize that the exponent~$z$
belongs to~$\R$ we shall use Greek letters $\rho$~and~$\si$ instead of $r$~and~$s$.

\blm7.42
\E\fa $x\in\R_{>0}$ and $\rho\in\R$ the \fw\ holds\dw
\beq7.78
(x\mo)^\rho = x^{-\rho} = (x^\rho)\mo.
\e
\elm

\proof
We first prove
$$
(x\mo)^\rho = (x^\rho)\mo, \q x\in\R_{>0},\ \rho\in\R. \leqno\rm\er{7.78}(i)
$$
The case $x=1$ is trivial since $1^\rho \nde7.77 = 1$ \fa $\rho\in\R$, and $1\mo=1$.

``$x<1$'': \ $x^\rho \nde7.77 = ((x\mo)^\rho)\mo$, hence $(x^\rho)\mo = \bigl(((x\mo)^\rho)\mo\bigr)\mo = (x\mo)^\rho$.

``$x>1$'': \ Since $x\mo<1$, we have $(x\mo)^\rho \nde7.77 = \bigl(((x\mo)\mo)^\rho\bigr)\mo = (x^\rho)\mo$.
$$
(x^\rho)\mo = (x^{-\rho}), \q x\in\R_{>0},\ \rho\in\R. \leqno\rm\er{7.78}(ii)
$$

``$x<1$'': \ $x^{-\rho} \nde7.77 = ((x\mo)^{-\rho})\mo$. Since $x\mo>1$, we have $(x\mo)^{\rho+\si} \nde6.134 = (x\mo)^\rho \cdot (x\mo)^\si$,
$\rho,\si\in\R$. Setting $\si:=-\rho$, we obtain $1 \nde7.77 = (x\mo)^0 = (x\mo)^\rho \cdot (x\mo)^{-\rho}$. \E\Tf $(x\mo)^\rho = ((x\mo)^{-\rho})
\mo = x^{-\rho}$.

The case $x=1$ is trivial, since $1^\rho \nde7.77 = 1$ for $\rho\in\R$.

``$x>1$'': $(x^\rho)\mo \nde6.120 = (\supl_{\si<\rho}x^\si)\mo$. Set $A_\rho := \{x^\si\in\R_{>0}: \si<\rho,\ \si\in\R\}$. Thus $x^\rho \nde6.124 =
\sup A_\rho$. The self-map of~$\R_{>0}$ defined by $\INV(x):= x\mo$, the inverse of~$x$ in the \ag\ $\pz2{\R_{>0}}$, \sf ies $\INV\circ\INV =
{\rm I}_{\R_{>0}}$, is bi\jc\ and strictly de\cre\ since $0<x<y$ implies $0<y\mo<x\mo$.

In view of Lemma \rf{l5.33}\,(iii), we have $\inf \INV(A_\rho)$ exists, and $\inf \INV(A_\rho)=\INV(\sup A_\rho)$. Thus $(x^\rho)\mo = \inf
\{(x^\si)\mo: \si<\rho,\ \si\in\R\}$. We have $(x^\si)\mo = x^{-\si}$, $\si\in\R$, by \er{7.78}\,(ii). Hence $(x^\rho)\mo = \inf\{x^{-\si}: \si<\rho,
\ \si\in\R\} = \inf\{x^{\si'}: -\si'<\rho,\ \si'\in\R\} = \inf\{x^{\si'}: {-\rho<\si'},\ \si'\in\R\} \nde6.123 = x^{-\rho}$. \E\Tf $(x^\rho)\mo=
x^{-\rho}$.

This completes the proof of Lemma \rf{l7.42}.
\endproof

\blm7.43
\E\fa $x\in\R_{>0}$ and all $r,s\in\R$ \pp ies \er{5.71}\,{\rm RI0--RI4} hold.
\elm

\proof \

``$x=1$'': RI4 follows from \er{7.77}. RI0, RI1 follow from RI4. RI2: $1^{r+s}=1= 1\cdot1 = 1^r\cdot 1^s$. RI3: $(1^r)^s=1^s=1$.

\ssk
``$x\ne1$'': \
RI0 and RI1: $x^0 \nde7.77 = \wh\psi_x(0) \nde6.120 = \psi_x(0) \nde5.60 = 1$; $x^1 = \wh\psi_x(1) = \psi_x(1) \nde5.60 = x^1 \nde5.71 =x$.

RI2 follows from \er{6.134} if $x>1$. Suppose $x<1$. Then $x\mo>1$, hence $(x\mo)^{r+s} = (x\mo)^r \cdot (x\mo)^s$, and $(x\mo)^{r+s}
\nad{\er{7.78}\rm(i)} = (x^{r+s})\mo$, $(x\mo)^r = (x^r)\mo$ and $(x\mo)^s= (x^s)\mo$. Thus $(x^{r+s})\mo = (x^r)\mo\cdot(x^s)\mo$,
$(x^s\cdot x^r)\mo = (x^r\cdot x^s)\mo$. Finally, $x^{r+s} = ((x^{r+s})\mo)\mo =((x^r\cdot x^s)\mo)\mo = x^r\cdot x^s$.

RI3: We first observe that the case $x=1$ follows from RI4. Indeed, $(1^r)^s = 1^s =1$. We next observe that the case $x<1$ follows from the case
$x>1$. Indeed, if $\a<1$ then $x:=\a\mo>1$, hence $(\a^r)^s = ((x\mo)^r)^s \nde7.78 = ((x^r)\mo)^s \nde7.78 = ((x^r)^s)\mo \nad{\rm RI3,x>1} =
(x^{rs})\mo \nde7.78 = (x\mo)^{rs} = \a^{rs}$.

We also observe that the case $x>1$ and $r=0$ or $s=0$ is trivial. Indeed, $(x^0)^s = 1^s = 1 = x^0 = x^{0\cdot s}$. \Mo $(x^r)^0 = 1 = x^0
= x^{r\cdot0}$.

Finally, we show that if $x>1$, cases $r{>}0,\ s{<}0$; $r{<}0,\ s{>}0$; $r{<}0,\ s{<}0$ follow from the case $r{>}0,\ s{>}0$.

``$x>1$, $r>0$, $s<0$'': \ $(x^r)^s = (x^r)^{-|s|} = ((x^r)\mo)^{|s|} = ((x\mo)^r)^{|s|} = (x\mo)^{r|s|} = x^{-(r|s|)} = x^{r(-|s|)} = x^{rs}$.

``$x>1$, $r<0$, $s>0$'': \ $(x^r)^s = (x^{-|r|})^s = ((x^{|r|})\mo)^s = ((x^{|r|})^s)\mo = (x^{|r|s})\mo = x^{-(|r|s)} = x^{-(|r|)s} = x^{rs}$.

``$x>1$, $r<0$, $s<0$'': \ $(x^r)^s = (x^{-|r|})^{-|s|} = ((x^{|r|})\mo)^{-|s|} = (((x^{|r|})\mo)\mo)^{|s|} = x^{|r|\,|s|} = x^{|rs|} = x^{rs}$.

It remains to prove the case $x>1$, $r>0$ and $s>0$. To this end we use the \fw\ formula, which will be proved in the next lemma:
\beq7.79
\log_x(y^\rho) = \rho(\log_x y) \qh{\fa $x\in\R_{>1}$, $y\in\R_{>0}$ and $\rho\in\R$.}
\e
Applying \er{7.79} with $x:=x$, $y:=x$ and $\rho:=\rho\si\in\R_{>0}$, we find $\log_x(x^{\rho\si})=\rho\si\log_xx \nde6.138 =
\rho\si$. Applying \er{7.79} with $x:=x$, $y:=x^\rho\in\R_{>1}$ (by Theorem \rf{t6.50}\,(iii) and $x^0=1$), and $\rho:=\si$, we find
$\log_x((x^\rho)^\si)  = \si\log_x (x^\rho) = \si\rho \log_xx = \rho\si$. Since $\log_x:\R_{>0}\to\R$ is bi\jc\ as inverse of the bi\jc\ map
$\wh\psi_x:\R \to \R_{>0}$ (see Theorem \rf{t6.50} and \E\df\ \rf{d6.53}), the map $\log_x$ is in\jc. Hence from $\log_x((x^\rho)^\si) = \log_x
(x^{\rho\si})$, we infer $(x^\rho)^\si = x^{\rho\si}$.
\endproof

\newpage
\blm7.44
Let $x\in\R_{>1}$ and let $\log_x:\R_{>0}\to\R$ be the map defined in \E\df\ \rf{d6.53}. Then

\hph i,i, $\log_x$ is an \ois sm from $(\R_{>0},\ge)$ onto $(\R,\ge)$,

\hph ii,, \er{7.79} holds,

\noi as well as
\bea7.80
\log_x(y\cdot z) &= \log_xy +\log_xz &&\hskip-50pt y,z\in\R_{>0},\\
\log_x(y\mo) &= -\log_x y &&\hskip-50pt y\in\R_{>0}. \lb{7.81}
\e
\elm

\proof \

``\er{7.80}'': \ By Theorem \rf{t6.50} \te\ $\a,\b\in\R$ \st $y=\wh\psi_x(\a)$, $z=\wh\psi_x(\b)$. Then $\wh\psi_x(\a+\b)\nde6.134 = \wh\psi_x(\a)
\cdot \wh\psi_x(\b)$. Since $\log_x:\R_{>0}\to\R$ is the inverse of the\break map $\wh\psi_x:\R\to\R_{>0}$, we obtain $\log_xy + \log_xz = \a+\b =
(\wh\psi_x)\Inv \wh\psi_x(\a+\b) =\break (\wh\psi_x)\Inv(\wh\psi_x(\a)\cdot\wh\psi_x(\b)) = \log_x(y\cdot z)$.

``\er{7.81}'': \ From $y\cdot y\mo=1$ and $\log_x1=0$ (since $x^0=1$) and from \er{7.80} we infer $0=\log_x1 = \log_x(y\cdot y\mo) = \log_xy +
\log_x(y\mo)$.

``\er{7.79}'': \ $x\in\R_{>1}$, $y\in\R_{>0}$. Case $\rho:=n\in\Na$. To prove: $\log_x y^n = n\log_xy$, $n\in\Na$. We use \In\ on $n\in\Na$. Clearly,
\er{7.79} holds for $n:=1$. Note that $y^n>1$ \fa $n\in\Na$ by \er{3.10} and \In. Suppose that \er{7.79} holds for $\rho:=n$. Then $\log_xy^{n+1} =
\log_x(y^n\cdot y) \nde7.80 = \log_x(y^n) + \log_xy = n\log_xy +\log_xy= (n+1)\log_xy$.

Case $\rho:=\frac1m$, $m\in\Na$: \ Set $z=y^{\frac1m}\in\R_{>0}$ (see \er{5.50}). We have $y=(y^{\frac1m})^m$ by \er{5.56}. Hence $\log_xy =
\log_x(y^{\frac1m})^m = m\log_x(y^{\frac1m})$.

Case $\rho:=\frac nm$: \ We have $y^{\frac nm} \nde5.57 = (y^{\frac1m})^n$, hence $\log_x y^{\frac nm} = \log_x(y^{\frac1m})^n = n\log_x(y^{\frac1m})
= \frac nm\log_xy$.

``(i)'': \ We recall that the map $\wh\psi_x:(\R,\ge) \to (\R_{>0},\ge)$ is an \ois sm by Theorem \rf{t6.50}\,(vi). Since $(\R,\ge)$ is totally
ordered, $\log_x:= (\wh\psi_x)\Inv$ is also an \ois sm in view of Lemma 1.3.34.

Case $y\in\R_{>1}$ and $\rho\in\R_{>0}$: \ By \E\Pr\ \rf{p6.11}, Lemma \rf{l6.12} and Lemma \rf{l6.16}\,(iii), \te s a \sq\ $\zb rn\N$ in~$\Q_{>0}$
\st $r_n\ua\rho$ (see the proof of Lemma \rf{l6.43}). By Theorem \rf{t6.50}\,(vi) we have $\wh\psi_y(r_n)\ua \wh\psi(\rho)$. By~(i) $\log_xy^{r_n}
\ua \log_x(y^\rho)$. \E\oh $\log_xy^{r_n} = r_n\log_xy$, $n\in\N$. Since $y\in\R_n$, we have $\log_xy>0$, hence $r_n\log_xy \ua \rho\log_xy$ by Lemma
\rf{l6.16}\,(iv). \csq, $\log_x(y^\rho) = \supl_{n\ge1}\log_xy^{r_n} = \supl_{n\ge1}r_n\log_xy = \rho\log_xy$.

We next remove \cn\ $y>1$.

Case $0<y\le1$ and $\rho\in\R_{>0}$: \ The case $y=1$ is trivial. If $y\in\R_{>0}$, $y<1$, then $\a:=y\mo\in\R_{>1}$. Hence $\log_x\a^\rho = \rho
\log_x\a$. We have $\log_x\a = \log_x(y\mo) \nde7.81 = -\log_xy$. \Mo $\log_x(\a^\rho) = \log_x((y\mo)^\rho) \nad{y\mo>1}= \rho \log_x(y\mo)\nde7.81 =
-\rho\log_xy$. Hence $\rho\log_x y = -\log_x(\a^\rho) = -\log_x((y\mo)^\rho) \nde7.78 = -\log_x((y^\rho)\mo) \nde7.81 = -(-\log_xy^\rho) =
\log_xy^\rho$.

The case $y\in\R_{>0}$, $\rho:=0$ is trivial.

Case $y\in\R_{>0}$, $\rho:=-|\rho|=0$: \ $\log_xy^\rho = \log_xy^{-|\rho|} \nde7.78 = \log_x(y\mo)^{|\rho|} = |\rho|\log_xy\mo \nde7.81 = -|\rho|
\log_xy = \rho\log_xy$.
\endproof

\Wanp prove an \ext\ of Theorem \rf{t6.50} to the case $x\in\R_{>0}$, $x<1$.

\newpage
\bth7.45
Let $\a\in\R_{>0}$, $\a<1$, and let $\wh\psi_\a :\R\to\R_{>0}$ be defined by $\wh\psi_\a (\rho):=\a^\rho$, $\rho\in\R$ $($see \E\df\ \rf{d7.41}$)$.
Then the \fw\ assertions hold\dw

\hph i,ii, The \f\ $\wh\psi_\a $ is a monoid-\is sm from the \ag\ $\pz1\R$ onto the \ag\ $\pz2{\R_{>0}}$.

\hph ii,i, The \f\ $\wh\psi_\a $ is strictly de\cre\ as well as its inverse $\log_\a:= \smash{\wh\psi_\a} \Inv$.

\hph iii,, The \f\ $\wh\psi_\a $ \sf ies \er{5.71} {\rm RI0--RI4, RI6, RI7} with $x:=\a$, $r:=\rho$, $s:=\si$, $\rho,\si\in\R$.

\hph iv,, $(y\cdot z)^\rho = y^\rho\cdot z^\rho$ \fa $y,z,\rho\in\R_{>0}$.
\eth

\proof \

(i) \ Setting $x:=\a\mo$, we have $x\in\R_{>1}$ and $\wh\psi_\a(\rho) = \a^\rho = (x\mo)^\rho \nde7.78 = (x^\rho)\mo = (\wh\psi_x(\rho))\mo$.
Setting $\d_{-1}(y)=-y$, $y\in\R$, we find that $\wh\psi_\a$ is the \cm\ of two monoid-\is sms. Indeed, $\wh\psi_\a = \wh\psi_x\circ\d_{-1}$ where
$\d_{-1}:\pz1\R \to \pz1\R$ \sf ies $\d_{-1}\circ\d_{-1} = \id_\R$, $\d_{-1}(y+z)=\d_{-1}(y)+\d_{-1}(z)$, $\d_{-1}(0)=0$, and where $\wh\psi_x:
\pz1\R$ to $\pz2{\R_{>1}}$ is a monoid-\is sm by Theorem \rf{t6.50}\,(vii). \E\Tf assertion~(i) follows from Lemma 2.1.8.

(ii) \ From the proof of part~(i) we know that the map $\wh\psi_\a$ is bi\jc. \Mo the self-map $\d_{-1}$ of $(\R,\ge)$ is strictly de\cre, and the
map $\wh\psi_x$ is strictly in\cre\ by Theorem \rf{t6.50}\,(vi). Hence $\wh\psi_\a$ is strictly de\cre\ as easily seen. Since $(\R,\ge)$ is totally
ordered, $\smash{\wh\psi_\a}\Inv$ is also strictly de\cre\ in view of Lemma \rf{l5.33}\,(ii).

(iii) \ \E\pp ies \er{5.71}\,RI0--RI4 follow from Lemma \rf{l7.43}, and RI6, RI7 from Lemma \rf{l7.42}.

(iv) \ Let $x\in\R_{>1}$. Then $\log_x(y\cdot z)^\rho \nde7.79 = \rho\log_x(y\cdot z)$ since $y,z,y\cdot z\in\R_{>0}$. \Mo $\log_x(y\cdot z) = \log_xy
+\log_xz$ by \er{7.80}. Hence $\log_x(y\cdot z)^\rho = \rho\log_xy + \rho\log_xz \nde7.79 = \log_xy^\rho + \log_xz^\rho \nde7.80 = \log_x(y^\rho \cdot
z^\rho)$. Since $\log_x$ is in\jc, we find $(y\cdot z)^\rho= y^\rho\cdot z^\rho$.
\endproof

We now \es\ some basic \pp ies of the self-maps of $(\R_{\ge0},\ge)$, $x\mt x^\rho$, $\rho\in\R_{>0}$.

\bpr7.46
Let $x\in\R_{>0}$ and $\rho\in\R_{>0}$. Set $f_\rho(x):=x^\rho$ $($see \E\df\ \rf{d7.41}$)$. Then the \fw\ assertions hold\dw

\hph i,i, \E\fe $\rho\in\R_{>0}$, $f_\rho$ is an \ois sm of $(\R_{>0},\ge)$ onto $(\R_{>0},\ge)$.

\hph ii,, The inverse of the self-map $f_{\rho}$ is $f_{\rho\mo}$.
\epr

\proof \

(i) ``\ti{Sur\ji}'': \ Let $z\in\R_{>0}$, $y\in\R_{>1}$. Suppose \te s $x\in\R_{>0}$ \st $x^\rho=z$. Then $\log_yx^\rho = \log_yz$, hence $\rho\log_y
x = \log_yz$ by \er{7.79}, hence $\log_yx=\rho\mo \log_yz$ and $x = y^{\log_yx} = y^{\rho\mo\log_yz}$. Set
\beq7.82
x:=y^{\rho\mo\log_yz}.
\e
Then $\log_yx^\rho = \rho\log_yx = \rho\cdot\rho\mo\log_yz = \log_yz$, hence $\log_yx^\rho = \log_yz$. Since $\log_y$ is in\jc, we obtain $x^\rho=z$.

``\ti{Strict in\cre ness}'': \ Let $x,y,\rho\in\R_{>0}$ be \st $x<y$. We have to show that $x^\rho<y^\rho$. Let $z\in\R_{>1}$. Then $\log_zx< \log_zy$
by Theorem \rf{t7.45}\,(ii). Hence $\rho\log_zx< \rho\log_zy$ by~\er{3.10}. \E\Tf by \er{7.79}, $\log_zx^\rho < \log_zy^\rho$. Since $\log_z$ is an
\ois sm, we find $x^\rho < y^\rho$.

(ii) \ From Lemma 1.3.34, we infer that $f_\rho$ is in\jc, hence bi\jc, and that $f_\rho\Inv$ is also strictly in\cre, since $(\R_{>0},\ge)$ is
totally ordered.
\endproof

\bco7.47
Let $x\in\R_{\ge0}$ and $\rho\in\R_{>0}$. Set
\beq7.83
f_\rho(x) := \bca
x^\rho &\hbox{if }x>0,\\
0 &\hbox{if }x=0.
\eca
\e
Then assertions {\rm(i)} and {\rm(ii)} of \E\Pr\ \rf{p7.46} hold with $\R_{>0}$ replaced by $\R_{\ge0}$.
\eco

\proof
Observe that $0<x$, $0<x^\rho$ \fa $\R_{>0}$ and $\rho\in\R_{>0}$.
\endproof

\bpr7.48
\E\fe $\rho\in\R_{>0}$, the self-map of $\R_{\ge0}$ defined in \er{7.83} is order-continuous, that is, the \fw\ assertions hold\dw
\bea7.84
 f_\rho(0) &= \inf_{y>0,y\in\R} f_\rho(y),\\
 \sup_{y<x,y\in\R} f_\rho(y) &= f_\rho(x) = \inf_{x<z,z\in\R}f_\rho(z) \qh{\fa $x\in\R_{>0}$.}\lb{7.85}
\e
\epr

\proof

``\er{7.84} \ti{and the second \et y in \er{7.85}}'': \ In \er{7.84}, $f_\rho(y)>0$ and in \er{7.85} $f_\rho(z)>f_\rho(x)$ since $f_\rho$ is strictly
in\cre. Hence $\infl_{y>0,y\in\R}f_\rho(y)$ exists. Set $f_\rho(0^+):=\infl_{y>0,y\in\R}f_\rho(y)$. Suppose, for \cd ion, that $f_\rho(0^+)>0$. Set
$\a:=f_\rho(0^+)\in\R_{>0}$. Then there is no $x\in\R_{>0}$ \st $f_\rho(x)=\frac\a2$, since $f_\rho(x)>\frac\a2$ \fa $x\in\R_{>0}$. The proof of
the second \et y in \er{7.85} is similar.
\endproof

\bex7.49 \

\hph i,ii, Prove the first \et y in \er{7.85}.

\hph ii,i, Show that in \er{7.84} and \er{7.85}, $y>0$, $y\in\R$ (resp.\ $y<x$, $y\in\R$; $x<z$, $z\in\R$) can be replaced by $y>0$, $y\in\Q$
(resp.\ $y<x$, $y\in\Q$; $x<z$, $z\in\Q$).

\hph iii,, Show that every in\cre\ bi\jc\ self-map of $(\R_{\ge0},\ge)$ \sf ies \er{7.84}, \er{7.85}.
\eex

Our next goal is to state and prove Ostrowski's theorem.

\begin{thm}[Ostrowski \cite{Ostr}, {\cite[pp.~201--202]{A2}}]\lb{t7.51}
Let $\vf:\Q\to\R_{\ge0}$ be an $\R$-valued valuation $($see \E\df\ \rf{d5.23}$)$. Then either $\vf(t)=|t|^\rho$, $t\in\R$, \fs $\rho\in\R_{>0}$,
$\rho\le1$, or $\vf(t)=(\vf_p(t))^\si$, $t\in\R$, \fs prime~$p$ and some $\si\in\R_{>0}$.\index{Ostrowski's theorem}
\eth

In what follows valuations will always be $\R$-valued. Theorem \rf{t7.51} is an obvious con\sq\ of the \fw\ lemmata.

\blm7.52 \

\hph i,i, \E\fe $p$ prime and every $\si\in\R_{>0}$ the map $\vf:\Q \to \R_{\ge0}$ defined by
\beq7.86
\vf(t):=(\vf_p(t))^\si, \q t\in\R,
\e
is a \vl\ of $\Q$ \sf ying
\beq7.87
\vf(n)\le1 \q\hbox{\fa}n\in\N.
\e

\hph ii,, \E\fe $\rho\in\R_{>0}$, $\rho\le1$, the map $\vf:\Q\to\R_{\ge0}$ defined by
\beq7.88
\vf(t) := |t|^\rho, \q t\in\R,
\e
is a \vl\ of $\Q$ \sf ying
\beq7.89
\hbox{\E\te s $n\in\N$ \st $\vf(n)>1$.}
\e
\elm

\bdf7.53
A \vl\ of $\Q$ is called \ti{\Ar} if it \sf ies \er{7.89} and \ti{non-\Ar} if it \sf ies \er{7.87}.\index{valuation!\Ar}\index{valuation!non-\Ar}
\edf

\blm7.54
If $\vf:\Q\to\R_{\ge0}$ is a non-\Ar\ \vl\ of~$\Q$, then \te\ $p$~prime and $\si\in\R_{>0}$ \st \er{7.86} holds.
\elm

\blm7.55
If $\vf:\Q\to\R_{\ge0}$ is an \Ar\ \vl\ of~$\Q$, then \te s $\rho\in\R_{>0}$, $\rho\le1$, \st \er{7.88} holds.
\elm

\proof[Proof of Lemma \rf{l7.52}\,\rm(i)]
Let $\vf:\Q\to\R_{\ge0}$ be defined as in \er{7.84}. Note that $\vf_p(t)\in\R_{\ge0}$, hence $(\vf_p(t))^\si$ is well defined in view of \er{7.83}.
We already know from Remark \rf{r7.39} that $\vf_p$ is a \vl\ on~$\Q$ \sf ying \er{7.62}, \er{7.76a} and $\vf_p(n)\le1$ \fa $n\in\N$ (note that
$\vf_p(0)=0$). We now show that if $\si\in\R$, $(\vf_p)^\si$ is a non-\Ar\ \vl\ of~$\Q$.

``\er{5.42}'' follows from \er{7.83}.

``\er{5.43}'': \ $(\vf_p(r\cdot s))^\si = (\vf_p(r)\cdot\vf_p(s))^\si = \vf_p(r)^\si\cdot\vf_p(s)^\si$ by Lemma \rf{l7.43} whenever $\vf_p(r)\ne0$,
$\vf_p(s)\ne0$, and by \er{7.83} otherwise.

``\er{5.44}'': \ From \er{7.76a} and $\vf_p(r)\ge0$ we infer $(\vf_p(r+s))^\si \le (\max(\vf_p(r),\vf_p(s)))^\si$ by \E\Pr\ \rf{p7.46}. If $\vf_p(r)
\le\vf_p(s)$ then $(\vf_p(r))^\si \le(\vf_p(s))^\si$, hence $(\max(\vf_p(r),\vf_p(s)))^\si = \max((\vf_p(r))^\si,(\vf_p(s))^\si) \le (\vf_p(r))^\si
+(\vf_p(s))^\si$.

``\er{5.46}'': \ $(\vf_p(p))^\si = \bigl(\frac1p\bigr)^\si\ne1$. \If that $(\vf_p)^\si$ is a \vl\ of~$\Q$. \Mo $(\vf_p(0))^\si = 0^\si \nde7.83 = 0$,
$(\vf_p(1))^\si = 1^\si= 1$, $(\vf_p(n+1))^\si \le \max(\vf_p(n)^\si,1)$, hence by using \In\ we find that \er{7.87} holds.

This completes the proof of Lemma \rf{l7.52}\,(i).
\endproof

\proof[Proof of Lemma \rf{l7.52}\,\rm(ii)]
We already know that $t\mt|t|$, $t\in\Q$, is a \vl\ of~$\Q$. Let $\vf:\Q\to\R$ be as in \er{7.88}.

``\er{5.42}'' follows from \er{7.83}.

``\er{5.43}'' follows from Lemma \rf{l7.43} and from \er{7.83}.

``\er{5.45}'': \ $\vf(z)=|z|^\rho > |1|^\rho = 1^\rho = 1$ by \E\Pr\ \rf{p7.46}\,(i) and Lemma \rf{l7.43}.

``\er{7.89}'': \ $\vf(z)>1$.

``\er{5.44}'' follows from
\beq7.90
(a+b)^\a \le a^\a + b^\a \qh{\fa $a,b\in\R_{\ge0}$ and $\a\in\R_{>0}$, $\a<1$.}
\e

\ti{Proof of \er{7.90}.} The in\et y is trivial if $a=0$ or $b=0$ since $0^\a=0$. We assume $a,b\in\R_{>0}$.
Set $t:=\frac b{a+b}$, thus $1-t=\frac a{a+b}$. Then $t,1-t \in
(0,1)_\R$, hence $t^\a\ge t^1$, $(1-t)^\a\ge(1-t)^1$. Indeed, we have $s^\b < s^\g$ for $s\in(0,1)_\R$, $\b,\g\in\R$, $\b>\g$ by Theorem
\rf{t7.45}\,(ii). \E\Tf $1=(1-t)+t\le (1-t)^\a+t^\a$, hence $1\le (\frac a{a+b})^\a+(\frac b{a+b})^\a$, and
\bmlg
(a+b)^\a \le (a+b)^\a\cdot
\bigl(a\cdot(a+b)\mo\bigr)^\a + (a+b)^\a \bigl(b\cdot(a+b)\mo\bigr)^\a \nad*= (a+b)^\a
\bigl(a^\a \cdot((a+b)\mo)^\a\bigr)\\ + (a+b)^\a \bigl(b^\a \cdot((a+b)\mo)^\a\bigr) \nde7.78 = (a+b)^\a\bigl(a^\a ((a+b)^\a)\mo\bigr) +
(a+b)^\a \bigl(b^\a ((a+b)^\a)\mo\bigr) = a^\a + b^\a.
\e
In $\nad*=$ we used Theorem \rf{t7.45}\,(iv). \hfill$\Box$

\E\csq, $|a+b|^\a \nad*\le (|a|+|b|)^\a \le |a|^\a+|b|^\a$ which proves \er{5.44}. In $\nad*\le$ we used \E\Pr\ \rf{p7.46}\,(i). This completes the
proof of Lemma \rf{l7.52}\,(ii), hence of Lemma~\rf{l7.52}.
\endproof

\proof[Proof of Lemma \rf{l7.54} \rm(see {\cite[pp.\ 194, 202--203]{A2}})]\

\ti{Step} 1. \ ``$\vf(r+s)\le\max(\vf(r),\vf(s))$, $r,s\in\Q$'': \ From \er{5.44} using \In, we find $\vf\bigl(\suml_{l=1}^na_l\bigr) \le \suml_{l=1}^n
\vf(a_l)$ \fa $a_i\in\Q_{>0}$, $1\le i\le n$, and all $n\in\Na\sms1$. Let $a,b\in\Q_{>0}$ and let $n\in\Na\sms1$. Then we have
\bmlg
(\vf(a+b))^n \nad*= \vf((a+b)^n) \nad{\rm(4.4.113)}= \vf\Bgg(\sum_{k=0}^n \frac{n!}{k!(n-k)!}\,a^{n-k}b^k) \\\le
\sum_{k=0}^n \vf\Bgg(\frac{n!}{k!(n-k)!}\,a^{n-k}b^k) \nde5.44 = \sum_{k=0}^n \vf\Bgg(\frac{n!}{k!(n-k)!})\vf(a^{n-k}b^k)
\nad{**}\le \sum_{k=0}^n 1\cdot\vf(a^{n-k}b^k) \\\nde5.43 = \sum_{k=0}^n \vf(a^{n-k})\vf(b^k) \nad*= \sum_{k=0}^n \vf(a)^{n-k} \vf(b)^k.
\e
Set $M:=\max(\vf(a),\vf(b))$. Note that $M>0$ and $\vf(a+b)>0$ by \er{5.42}. Thus\break $(\vf(a+b))^n \le \suml_{k=0}^n M^{n-k}M^k = \suml_{k=0}^n M^n
\nad{\rm(2.4.40)} = (n+1)M^n$ \fa $n\in\N$. In~$\nad*=$ we used $\vf(c^n) = (\vf(c))^n$, $n\in\Na$, $c\in\Q$, which follows from \er{5.43} and \In.
In~$\nad{**}=$ we used $\frac{n!}{k!(n-k)!}\in\Na$ for $0\le k\le n$, $n\in\Na\sms1$, \er{7.86}, \er{5.42}, \er{3.10} and (2.4.41). \E\Tf we have
\beq7.91
\Bgg(\frac{\vf(a+b)}M)^n \le n+1 \qh{\fa $n\in\Na\sms1$.}
\e
Setting $\g:=\frac{\vf(a+b)}M \in \R_{>0}$, we find $\g^n\le n+1$ \fa $n\ge2$. We claim that $\g\le1$. Suppose, for \cd ion, that \te s $\d\in
\R_{>0}$ \st $\g=1+\d$. Then $\g^n\nad{\rm(4.4.54)}\ge 1^n+n1^{n-1}\d + \frac{n(n-1)}2 1^{n-2}\d^2 = 1+n\d + \frac{n(n-1)}2\d^2$ \fa $n\ge2$.
However, \te s $\ov n\ge2$ \st
\beq7.92
\ov n\d + \frac{\ov n(\ov n-1)}2 \d^2 > \ov n.
\e
Indeed, it suffices to show that \te s $\ov n\in\Na$ \st $\frac{\ov n(\ov n-1)}2\d^2 > \ov n$, \ev tly, $\frac{\ov n-1}2 > \frac2{\d^2}$ or
$\ov n>1+\frac4{\d^2}$. Since $(\R,\ge)$ is \Ar, \te s $\ov n\in\Na$ \st $\ov n\cdot1 > 1+\frac4{\d^2}$. \If that $\g\le1$, hence $\vf(a+b)
\le \max(\vf(a),\vf(b))$.

\ssk
\ti{Step} 2. \ ``\ti{\E\te s $p\in\cP$ \st $\{n\in\Na: \vf(n)<1\} = \{kp: k\in\Na\}$}'': \ Set $A:=\{a\in\Z: \vf(z)<1\}$.
Then $0\in A$, $z,w\in A$ implies $z+w\in A$,
${-z\in A}$ since ${\vf(0)\nde4.46 = 0}$; $\vf(z)<1$, $\vf(w)<1$ implies $\vf(z+w)\le \max(\vf(z),\vf(w))<1$; $\vf(-z)={\vf(z)<1}$. Set $A_+:=\{z\in
A: z>0\}$, $A_-:=\{z\in A: z<0\}$. Note that $A_+\ne\vn$ (hence $A_-\ne\vn$) by~\er{5.45}. \Mo $A$~is the disjoint union of $A_+$,
$\{0\}$ and~$A_-$. Note that $A_+\cup\{0\}$ is a \sbm\ of the monoid $(\N,+,0)$. If $x,y\in A_+\cup\{0\}$ then $x-y = x+(-y)\in A$, since $x,-y,
x+(-y)\in A$. If $x\ne y$ then either $x-y>0$ or $-(x-y)>0$. Hence $\vf(x-y) = \vf(-(x-y))<1$. If $x-y>0$ then $x-y\in A_+$, hence $A_+\cup\{0\}$
is a \pn\ monoid in view of Lemma 3.1.15. If $y-x>0$ then $A_+\cup\{0\}$ is a \pn\ monoid. Thus \te s $a\in A_+$ \st $A_+=\{ka: k\in\Na\}$.

We claim that $a$~is a \Pn. Indeed, suppose for \cd ion that \te\ $b,c\in\Na\sms1$ \st $a=bc$. If $b=c$, we have $\vf(b)^2 \nde5.43 = \vf(b^2) =
\vf(a)<1$, hence $\vf(b)<1$. \E\Tf \te s $k\in\Na$ \st $b=ka$, hence $a=(ka)^2=k^2a^2$, and $\vf(a) = \vf(k^2a^2) = \vf(k^2)\vf(a^2) \le \vf(a^2)
= \vf(a)^2 < \vf(a)$, a~\cd ion. If $b\ne c$ then $\vf(a)=\vf(b)\vf(c)$. Since $\vf(a)<1$ and $\vf(b),\vf(c)\le1$, we have $\vf(b)$ and/or
$\vf(c)<1$. Suppose $\vf(b)<1$. \E\Tf $b=ka$, $k\in\Na$, hence $a=kac$. Dividing by~$a$, we find $kc=1$, thus $k=c=1$. A~\cd ion, since $c\in\Na\sms1$.
\If that $A_+=\{kp:k\in\Na\}$ \fs prime~$p$.

\ssk
\ti{Step} 3. \ ``$\vf(t)=\vf_p(t)^\si$, $\si\in\R_{>0}$, $t\in\Q$'': \ Let $t\in\Q_{>0}$. In view of \er{7.74} \te s \ooo $z\in\Z$ \st $tp^z=t'$
where $t'=\frac ab$, $a,b\in\Na$, $p\nmid a$, $p\nmid b$. Thus $\vf(t')=\vf(\frac a1\cdot \frac1b) = \vf(a)\cdot(\vf(b))\mo$, by \er{5.42},
\er{5.43}. Since $p\nmid a$ and $p\nmid b$, we have $\vf(a)\ge1$, $\vf(b)\ge1$ by Step~2. Hence $\vf(a)=\vf(b)=1$ by \er{7.87}. \E\Tf $\vf(t')=1$,
and $1=\vf(tp^z) \nde5.43 = \vf(t)\vf(p^z)$. Since $\vf:\R_{\ge0} \to \R_{\ge0}$ \sf ies \er{5.43} and $\vf(1)=1$, the map~$\vf$ is a monoid-endo\mf\
of the \am\ $(\R_{\ge0},\cdot,1)$. Note that the map $z\mt p^z$ from $\pz1\Z$ into $\pz2{\R_{\ge0}}$ is also a monoid-\hm sm. \If from Lemma
2.1.8\,(ii) that $\vf(p^z)=\vf(p)^z$, $z\in\Z$. \E\Tf we obtain
\beq7.93
\vf(t) = (\vf(p))^{-z} \qh{\fs $z\in\Z$.}
\e
\E\oh we have from \er{7.83} $\vf_p(t)\vf_p(p^z) = \vf_p(t') = 1$. Since $\vf_p(p)\nde7.62 = \frac1p$, we obtain $\vf_p(t) = (\vf_p(p^z))\mo =
((\vf_p(p))^z)\mo$, hence
\beq7.94
\vf_p(t) = p^z.
\e
If $z:=0$ then $\vf(t)= p^{-0} =1$, $\vf_p(t)=1$, hence $\vf(t)=(\vf_p(t))^\si$ \fe $\si\in\R_{>0}$ since $1^z\nde7.77 = 1$. If $z\ne0$, and if
$\vf(t) = (\vf_p(t))^\si$ \fs $\si\in\R_{>0}$, then we would have $\log_x(\vf(t)) = \log_x(\vf_p(t))^\si$ \fs $x\in\R>0$, hence from \er{7.93},
\er{7.94}: $-z\log_x\vf(p) = \si z\log_xp$. Thus we should have $\si = -\frac{\log_x\vf(p)}{\log_xp}$. Choosing $x:=p$, we find
\beq7.95
\si = \log_p\vf(p)\mo = \log_p \vf\Bgg(\frac1p).
\e
Note that $\si\in\R_{>0}$ since $\frac1p < 1 <p$.
\endproof

\bex7.55 \

\hph i,i, Verify that \er{7.86} holds with $\si$ as in \er{7.95}.

\hph ii,, Does $\frac{\log_x\vf(p)}{\log_x p}$ depend on $x\in\R_{>1}$?
\eex

\proof[Proof of Lemma \rf{l7.55}.] \

\ti{Step} 1. \ It suffices to prove that the \fw\ holds:
\bml7.96
\vf(b)\le \max\bigl(1,\vf(a)^{\log_xb/\log_xa}\bigr) \hbox{ \fa $a,b\in\Na\sms1$}\\
\hbox{and \fa $x\in\R_{>0}$ \sf ying $x<a,x<b$.}
\e
Indeed, by \er{7.89} \te s $\ov n\in\Na\sms1$ \st $\vf(\ov n)>1$. \If that $\vf(n)>1$ \fa $n\in\Na\sms1$, otherwise there would exist $\wt n>1$ \st
$\vf(\wt n)\le1$. Then by \er{7.96} we would have $\vf(\ov n)\le1$, a~\cd ion. From \er{7.96} we infer that
\beq7.97
\vf(a) \le \vf(b)^{\log_xb/\log_xa}
\e
holds \fa $x\in\R_{>0}$, $x<a$, $x<b$, since $\log_xb/\log_xa>0$.

Since $(\log_xb)\mo>0$, it follows from \er{7.97} and \E\Pr\ \rf{p7.46}\,(i) that
$$
\vf(b)^{1/\log_xb}\le \vf(a)^{1/\log_xa}.
$$
Interchanging the roles of $a$ and $b$, we find
$$
\vf(a)^{1/\log_xa}\le \vf(b)^{1/\log_xb}.
$$
Hence we obtain
\beq7.98
\vf(a)^{1/\log_xa} = \vf(b)^{1/\log_xb}.
\e
Setting $\rho:=\frac{\log_x\vf(b)}{\log_xb}$ with $x<\vf(b)$, which is possible since $\vf(b)>1$, we have $\rho\in\R_{>0}$. Thus $\log_xb^\rho =
\rho\log_xb = \log_x\vf(b)$, and by in\ji\ of~$\log_x$, we find $\vf(b)=b^\rho$. Then from \er{7.98} we infer $\log_x\vf(a) = \log_xa \cdot
\log_xb^\rho \cdot \frac1{\log_xb} = \log_xa \cdot\rho \log_xb\cdot \frac1{\log_xb} = \log_xa\cdot\rho = \log_xa^\rho$. Hence $\vf(a)=a^\rho$.
We required $x\in\R_{>0}$, $x<\vf(a)$, which is possible since $\vf(a)>1$, that is, $x\in\R_{>0}$, $x<a$, $x<\vf(a)$, $x<b$, $x<\vf(b)$.
If $t\in\Q_{>0}$ and $t=\frac ab$, we find that $\vf(t)=\vf(a)\cdot\vf(b)\mo = a^\rho\cdot b^{-\rho}=t^\rho$.

Finally, since $2^\rho=\vf(2) = \vf(1+1)\nde5.44 \le \vf(1)+\vf(1) = 1+1 =2$, we find that $\rho\le1$ in view of \E\Pr\ \rf{p7.46}\,(i). Clearly, if
$t\in\Q\sms0$ then $\vf(t)=\vf(|t|)=|t|^\rho$, and $|0|^\rho \nde7.83 = 0$.

\ssk
\ti{Step} 2. We show that \er{7.96} is a con\sq\ of the in\et y
\bml7.99
\Bgg(\frac{\vf(b)}{M^{\log_xb/\log_xa}})^n \le a\,\frac{\log_xb}{\log_xa}\,n+a,\q n\in\Na,\\
\hbox{where $a,b\in\Na\sms1$, $x\in\R_{>0}$, $x<a$, $x<b$, and $M:=\max(1,\vf(a))$.}
\e
Proceeding as in the proof of Lemma \rf{l7.54}, we show that if $\g\in\R_{>0}$ \sf ies
\beq7.100
\g^n \le \b n+\a \qh{\fa $n\in\Na$}
\e
\fs $\b\in\R_{>0}$ and some $\a\in\R_{\ge0}$, then $\g\le1$. Suppose, for \cd ion, that \te s $\d\in\R_{>0}$ \st $\g=1+\d$. Then by (4.4.54) we have
\beq7.101
\g^n \ge 1+n\d+ \frac{n(n-1)}2\,\d^2, \q n\ge2.
\e
Since $(\R,\ge)$ is \Ar, we can find $\ov n>2$ \st $\ov n\d>\a$ and $\ov n\d^2>2\b+\d^2$. Then $1+\ov n\d + \frac{\ov n(\ov n-1)}2\d^2 > \b\ov n+\g$,
a~\cd ion. Setting $\g:=\frac{\vf(b)}{M^{\log_xb/\log_xa}}$, $\b:=a\frac{\log_xb}{\log_xa}$ and $\g:=a$, we find that $\g\le1$, since $\vf(b),
\log_xb, \log_xa \in\R_{>0}$, $M\in\R_{\ge1}$, $a\in\R_{\ge0}$, hence $\g,\b\in\R_{>0}$ and $\a\in\R_{\ge0}$. Finally, observe that $\g\ge1$ implies
\er{7.96}.

\ssk
\ti{Step} 3. \ We now prove \er{7.99}. Note that $(\vf(b))^n = \vf(b^n)$, $n\in\Na$, by \er{5.43} and using \In\ on $n\in\Na$. The idea is to use
the ``decimal'' \rp ation of~$b^n$ to base~$a$. In view of Lemma 2.5.3, \fe $n\in\Na$, \te\ $N_n\in\N$ and $c\en n_k\in\N$ \sf ying $c\en n_k<a$
\fa $k\in\N$, $k\le N_n$ \st the \fw\ holds:
\beq7.102
b^n = \sum_{k=0}^{N_n} c\en n_ka^k, \q c\en n_{N_n}\ne0.
\e
We have $a^{N_n} \le c\en n_{N_n}a^{N_n} \le b^n$, hence $N_n\log_xa \le n\log_xb$ \fe $x\in\R_{>0}$, $x<a$, $x<b$. \E\Tf
\beq7.103
{N_n} \le n\,\frac{\log_xb}{\log_xa}, \q n\in\Na.
\e
We next estimate $\vf(b^n)$:
$$
\vf(b^n) = \vf\Bgg(\sum_{k=0}^{N_n}c\en n_ka^k) \nde5.44 \le \sum_{k=0}^{N_n}\vf(c\en n_ka^k) \nde5.43 = \sum_{k=0}^{N_n}\vf(c\en n_k)a^k.
$$
Note that $\vf(z)=\vf(|z|) = \vf\bigl(\suml_{k=1}^{|z|}1\bigr) \nde5.44 \le \suml_{k=1}^{|z|}\vf(1) = \suml_{k=1}^{|z|}1 = |z|$ \fa $z\in\Z\sms0$.
\Mo $\vf(0)=0$. Since $|c\en n_k|\le a$, we have $\vf(c\en n_k) \le a$, $0\le k\le {N_n}$, $n\in\Na$. \E\Tf $\vf(b^n) \le \suml_{k=0}^{N_n}a \vf(a)^k
\le a\suml_{k=0}^{N_n}M^k = a(N_n+1)M^{N_n}$. From \er{7.103} we obtain
$$
\vf(b)^n = \vf(b^n) \le a\Bgg(\frac{\log_xb}{\log_xa}\,n+1) M^{\frac{\log_xb}{\log_xa}n},
$$
hence \er{7.99} holds. This completes the proof of Step~3 and of Lemma \rf{l7.55}.
\endproof

\bex7.56
Let $\vf\en\a:\R \to \R_{\ge0}$, $\a\in\R_{>0}$, $\a\le1$, be defined by $\vf\en\a(t):=|t|^\a$, $t\in\R$. Show that $\vf\en\a$ is a \vl\ of the
field~$\R$ \fe $\a\in\oz0,1 _\R$.
\eex

We next give a \chz\ of the usual metric on~$\R$.

\bpr7.57
The usual metric $d$ defined in \er{7.15} is the only metric \sf ying the \fw\ \cn s\dw
\begin{align}
{}&d(x+z,y+z) = d(x,y) \q\hbox{\fa}x,y,z\in\R, \tag{M4}\\
&d(0,1) = 1, \tag{M5}\\
&d(x,z)=d(x,y)+d(y,z) \q\hbox{\fa $x,y,z\in\R$ \st $x\le y\le z$.} \tag{M6}
\end{align}
\epr

\proof \

(i) The usual metric on $\R$ \sf ies (M1)--(M6).

``(M6)'': \ $d(x,z) \nde7.15 = |z-x| = z-x = (z-y)+(y-x) = |y-x|+|z-y| \nde7.15 = d(x,y)+d(y,z)$.

(ii) ``\E\uq'': \ Set $\vf(t):=d(0,t)$, $t\in\R$. Then $\vf(t)\ge0$ \fa $t\in\R$, and $\vf(t)=0$ iff~$t=0$. \Mo $\vf(-t)=\vf(t)$ \fa $t\in\R$.
Indeed, $\vf(-t)=d(0,-t) \nad{\rm(M4)}= d(0+t,(-t)+t) = d(t,0) \nad{\rm(M2)} = d(0,t)$. We have to show that $\vf(t)=t$ \fa $t\in\R_{\ge0}$. Let
$t,s\in\R_{\ge0}$. Then $0\le t\le t+s$, hence $\vf(t+s) = d(0,t+s) \nad{\rm(M6)} = d(0,t)+d(t,t+s) \nad{\rm(M4)} = d(0,t)+d(0,s) = \vf(t)+\vf(s)$.
\E\Tf the map $\vf|_{R_{\ge0}}$ \sf ies the ``\f al \eq''
\beq7.104
\vf(t+s) = \vf(t)+\vf(s), \q t,s\in\R_{\ge0},
\e
together with the \cn s
\bea7.105
\vf(t)\ge0,& \q t\in\R_{\ge0},\\
\vf(0) = 0,& \q \vf(1)=1.\lb{7.106}
\e
\If from \er{7.104}, \er{7.105} that $\vf$ is in\cre. \E\eq\ \er{7.104} (for $t,s\in\R$) is usually called \ti{Cauchy's \f al \eq}. This \eq\ is
called \ti{\f al\/} since the ``unknown''~$\vf$ is a~\f, in contrast to an algebraic \eq\ such as $x^3+ax^2+bx+c=0$, $a,b,c\in\R$, where the
``unknown''~$x$ is a \nm, for example an \el\ of~$\R$. A~\so\ to \eq\ \er{7.104} is \ti{by \df} a~\f\ (map) $\vf:\R_{\ge0} \to\R$ \sf ying
\er{7.104}. Every \so~$\vf$ to~\er{7.104} \sf ies $\vf(0)=0$, since $\vf(0)=\vf(0+0)=\vf(0)+\vf(0) =2\vf(0)$. Note that $f(t):=t$, $t\in\R_{\ge0}$,
\sf ies \er{7.104}--\er{7.106}. The question is to know whether the map~$f$ is the \ti{only} \f\ \sf ying \er{7.104}--\er{7.106}. The answer is yes.
For a proof we follow \texttt{http://www.mathsolympiad.org.nz}.

We first show that if $\vf:\R_{\ge0}\to\R$ \sf ies \er{7.104}--\er{7.106}, then $\vf(t)=t$ \fa $t\in\Q_{\ge0}$. To this end, it suffices to show that
if $\vf$~\sf ies \er{7.104}, then \fa $a\in\R$ and all $k\in\N$
\beq7.107
\vf(ka) = k\vf(a).
\e
Set $A:=\{k\in\N: \hbox{\er{7.107} holds}\}$. Then $0\in A$ since $\vf(0a)=\vf(0)=0=0\vf(a)$. If $n\in A$ then $n+1\in A$. Indeed, $\vf((n+1)a) =
\vf(na+a) \nde7.104 = \vf(na)+\vf(a) \nad{n\in A}= n\vf(a)+\vf(a) = (n+1)\vf(a)$. \E\Tf $A=\N$. \E\Ip choosing $a:=m\in\Na$ and using $\vf(1)=1$, we
find $1=\vf(1)=\vf(\frac mm) = \vf(m\frac1m) \nde7.107 = m\vf(\frac1m)$. Hence $\vf(\frac1m)=\frac1m$. \Mo \fa $n\in\N$ we have $\vf(\frac nm) =
\vf(n\frac1m) = n\vf(\frac1m) = n\frac1m = \frac nm$. \E\csq, $\vf(t)=t$ \fa $t\in\Q_{\ge0}$. \Wanp show that $\vf(t)=t$ \fa $t\in\R_{\ge0}$ if
\er{7.105} holds. Clearly, $\vf$~is in\cre\ by \er{7.104}. Suppose, for \cd ion, that \te s $\ov t\in\R_{>0}$ \st $\vf(\ov t)\ne\ov t$. If
$\vf(\ov t)>\ov t$, \te s $r\in\Q_{>0}$ \st $\vf(\ov t)>r>\ov t$ since $(\R,\ge)$ is an \Ar\ \of. Note that $r=\vf(r)$, hence $\vf(\ov t)>\vf(r)=r>
\ov t$, \cd ing the in\cre ness of~$\vf$. If $\vf(\ov t)<\ov t$, \te s $s\in\Q$ \st $\vf(\ov t)<s<\ov t$, hence $\vf(\ov t)<\vf(s)=s=\ov t$, \cd ing
the in\cre ness of~$\vf$.
\endproof

Note that if $\vf:\R_{\ge0} \to\R$ \sf ies \er{7.104}, \er{7.105} and $\vf(1)=a\in\R_{>0}$, then the \f\ $\wt\vf(t):=a\mo \vf(t)$, $t\in\R$, \sf ies
\er{7.104}--\er{7.106}, hence $\vf(t)=at$, $t\in\R$.

\bex7.58 \

\hph i,i, Let $\phi$ be an order- and monoid-\is sm from $\pz1{\R_{\ge0}}$ onto $\pz2{\R_{\ge1}}$ \sf ying $\F(1)=x\in\R_{>1}$. Show that $\F(t)=x^t$
\fa $t\in\R_{\ge0}$. (Hint: Make use of the \f\ $\log_x\circ\,\F$.)

\hph ii,, Show that the map $t\mt x^t$, $x\in\R_{>1}$, $t\in\R$, is the only order- and monoid-\is sm from $\pz1\R$ onto $\pz2{\R_{>0}}$.
\eex

\bex7.60 \

\hph i,ii, Show that the maps $x\mt x^\rho$, $\rho\in\R_{>0}$, are the only order- and monoid-\is sms from $\pz2{\R_{\ge1}}$ onto $\pz2{\R_{\ge1}}$.

\hph ii,i, Show that the maps $x\mt x^\rho$, $\rho\in\R_{>0}$, are the only order- and monoid-\is sms from $\pz2{\R_{>0}}$ onto $\pz2{\R_{>0}}$.

\hph iii,, Show that the maps $t\mt |t|^\a$, $\a\in\oz0,1 $, are the only \vl s~$\vf$ of~$\R$ \st the \rt ion of~$\vf$ to $\R_{\ge0}$ is in\cre.
\eex

We conclude this section by mentioning that a metric on a~set~$X$ can be ``transported'' to a~set~$X'$ \ep\ to~$X$ (see Example 2.1.5\,(v), \E\Pr\
4.4.2). Indeed, if $(X,d)$ is a \ms, and $X'$~is a set \ep\ to~$X$, then \te s a bi\jn\ $f:X\to X'$. Setting $d'(a',b'):= d(f\Inv(a'),f\Inv(b'))$,
$a',b'\in X'$, we find that $d'$~is a metric on~$X'$, and that $d_{X'}(f(u),f(v))=d(u,v)$, $u,v\in X$. See Example \rf{xa7.9}\,(iv).

\bdf7.61
Let $(X,d)$, $(X',d')$ be \ms s. A map $f:X\to X'$ \sf ying
\beq7.108
d'(f(x),f(y)) = d(x,y), \q x,y\in X,
\e
is called an \tb{isometry}.\index{isometry}
\edf

Note that an isometry is in\jc. If $f$~is bi\jc, then $f$~is called a \ti{sur\jc\ isometry}.

\bex7.61a
Let $f$ be a bi\jn\ from~$X$ into~$X'$. Suppose that \te s a metric~$d$ on~$X$. Set $d'(a',b'):= d(f\Inv(a'),f\Inv(b'))$ \fa $a',b'\in X'$. Show that
$d'$~is a metric on~$X'$ and $f$~is a sur\jc\ isometry between $(X,d)$ and $(X',d')$.
\eex

\bex7.62
Let $f$ be an isometry from a \ms\ $(X,d)$ into a \ms\ $(X',d')$. Show that if $\{x_n\}_{n\ge1}$ is a \sq\ of \el s of~$X$, and $a$~is an \el\ of~$X$,
then
\beq7.109
x_n\nad d\to a \hbox{ implies }f(x_n) \nad{d'}\to f(a).
\e
If, moreover, $f$ is a sur\jc\ isometry, show that $f\Inv$ is an isometry and that the \fw\ holds:
\beq7.110
x_n \nad d\to a \hbox{ iff }f(x_n)\nad{d'}\to f(a).
\e
\eex

\bex7.63
What are all isometries from $(\R,d)$ into $(\R,d)$ where $d$ is the usual metric of~$\R$? Are these isometries sur\jc?
\eex

%% file: app2.tex
\newpage
\Subsubsection{Metric compactness, connectedness and completeness}\label{ass.8}

The main goal of this section is to present some \fd\ and useful \pp ies of the field~$\R$, that are related to the notion of ``\cg nce of
\sq s'' \itd in the previous section. We first consider an \ev t \df\ of ``limit'' (see \E\df s \rf{d7.1}, \E\df\ \rf{d7.8} and Example
\rf{xa7.9}\,(i)) not involving the notion of metric. To this end we introduce the notions of upper and lower limits.

\bdf8.1
Let $\zb an\N$ be a \sq\ of \el s of~$\Q$ (resp.~$\R$). An \el\ $b\in\Q$ (resp.~$\R$) \sf ying \fe $\ve\in\Q_{>0}$ (resp.~$\R_{>0}$) the \fw\ two
\cn s:
\bea8.1
\hbox{the set }&\{n\in\N: a_n\ge b+\ve\} \hbox{ is finite},\\
\hbox{the set }&\{n\in\N: a_n+\ve > b\} \hbox{ is infinite}, \lb{8.2}
\e
is called a \ti{superior limit\/} (or a \ti{limit superior}) or an \ti{upper limit\/} of the \sq\ $\zb an\N$.\index{limit!superior}\index{limit!upper}
\edf

\blm8.2
A \sq\ $\zb an\N$ in $\Q$ $($resp.~$\R)$ possesses at most one upper limit, and in this case the \sq\ $\zb an\N$ is \ba.
\elm

\proof
Suppose $b,b'\in\Q$ (resp.~$\R$) be \ut s of a \sq\ $\zb an\N$ in~$\Q$ (resp.~$\R$), and suppose, for \cd ion, that $b\ne b'$. Interchanging the roles
of $b$~and~$b'$, we may assume \wlg that $b>b'$. By \cn\ \er{8.1} for~$b'$ with $\ve:={\frac12 (b-b')>0}$, we find that the set $\{n\in\N:
a_n \ge \frac{b+b'}2 \}$ is finite, since $\frac{b+b'}2 = b'+\ve$. By \cn\ \er{8.2} for~$b$ with the same~$\ve$, we find that the set $\{n\in\N:
a_n >\frac{b+b'}2 \}$ is infinite, since $\frac{b+b'}2=b-\ve$. A~\cd ion. \E\Tf $b=b'$. We next show that the \sq\ $\zb an\N$ is \ba. Setting $\ve:=1$
in \cn~\er{8.1} for~$b$, we find that either the set $A:=\{n\in\N: a_n\ge b+1\}$ is empty, hence $a_n<b+1$ \fa $n\in\N$, or the set~$A$ is nonempty
and finite. Set $\wt a:=\max\limits_{n\in A}a_n$. Then $a_n\le \wt a$ \fa $n\in\N$. Hence the \sq\ $\zb an\N$ is \ba.
\endproof

\bnt8.3
In view of the previous lemma we may speak of \ti{the} \ut\ of a \sq\ $\zb an\N$ in~$\Q$ (resp.~$\R$) if it exists. Usual notations for \ut s are
$\limu_{n\to\iy}a_n$ or $\limsup\limits_{n\to\iy}a_n$. We will adopt the notation $\limu_{n\ge0}a_n$ $($or $\limu_{n\ge1}a_n$ if the index~$n$
belongs to~$\Na)$. We will also write $\limu a_n$ when no confusion arises.
\ent

\blm8.4
Let $\zb an\N$ be a \cg nt \sq\ of \el s of~$\Q$ $($resp.~$\R)$ with limit~$a$ $($see \E\df s \rf{d7.1}, \rf{d7.8} and \rf{d7.15}$)$. Then $a=\limu
_{n\ge0}a_n$.
\elm

\proof
By \E\df\ \rf{d7.1}, \fe $\ve\in\Q_{>0}$ (resp.\ $\R_{>0}$) the set $A:=\{n\in\N: a_n\ge a+\ve\}\cup\{n\in\N: a_n\le a-\ve\}$ is finite. In view of
Theorem 1.4.18\,(i), the subset $\{n\in\N: a_n\ge a+\ve\}$ of~$A$ is finite. Hence \cn\ \er{8.1} holds with $b:=a$. \E\fe $\ve\in\Q_{>0}$ (resp.\
$\R_{>0}$), the set~$\N$ is the disjoint union of the sets $B:=\{n\in\N: a_n+\ve>a\}$ and $C:=\{n\in\N: a_n+\ve\le a\}$. Since $C=A$, the set~$C$ is
finite. Suppose, for \cd ion, that the set~$B$ is finite. From Theorem 2.3.23 we infer that the set $B\cup C$ is finite. Since $B\cup C=\N$, we obtain
a~\cd ion, since $\N$ is infinite by Theorem 1.4.18\,(iv). \If that \cn\ \er{8.2} holds, hence $a=\limu_{n\ge0}a_n$.
\endproof

A bounded \sq\ of \el s of~$\Q$ may not have an \ut\ in~$\Q$, but has always an \ut\ in~$\R$, as we will see below.

\bex8.5
Give an example of a bounded \sq\ of \el s of~$\Q$ having no \ut\ in~$\Q$ but having an \ut\ in~$\R$.
\eex

\bpr8.6
Let $\zb an\N$ be a bounded \sq\ of \el s of~$\R$. Then $\limu_{n\ge0}a_n$ exists and \sf ies
\beq8.3
\limu_{n\ge0}a_n = \inf_{k\ge0}\sup_{n\ge k}a_n.
\e
\epr

\proof
Let $\check M\in\R$ (resp.\ $\hat M\in\R$) denote a \lo\ (resp.\ \ub) for the \sq\ $\zb an\N$. \E\fe $k\in\N$ set
\beq8.4
\wh a_k:=\sup_{n\ge k}a_n.
\e
The \ex\ of $\wh a_k$, $k\in\N$, follows from the order-\cp ness of the field~$\R$ (see Theorem \rf{t5.1}). Note that
\beq8.5
\wh a_k \ge \wh a_l \qh{\fa $k,l\in\N$ \st $k\le l$.}
\e
Indeed, if $k<l$ then $\{a_n\in\R: n\ge k\} \supset \{a_n\in\R: n\ge l\}$. Hence an \ub\ for $\{a_n\in\R: n\ge k\}$ is an \ub\ for
$\{a_n\in\R: n\ge l\}$. \E\Tf $\supl_{n\ge k}a_n \ge \supl_{n\ge l}a_n$. \Mo $\check M\le \infl_{n\ge0}a_n \le \infl_{n\ge k}a_n \le
\supl_{n\ge k}a_n = \wh a_k$ \fa $k\in\N$. Thus we have
\beq8.6
\wh a_k \ge \check M \qh{\fa $k\in\N$.}
\e
\E\Tf $\infl_{k\ge0}\supl_{n\ge k}a_n$ exists by Corollary \rf{c5.12}. Set $\ov a:= \infl_{k\ge0}\supl_{n\ge k}a_n$. We now show that $\ov a$~\sf ies
\er{8.1} and \er{8.2}.

``\er{8.1}'': \ Note that
\beq8.7
a_k \le \wh a_k \qh{\fa $k\in\N$,}
\e
since $\wh a_k \nde8.4 = \supl_{n\ge k}a_n \ge a_k$, $k\in\N$. We have $\ov a= \infl_{k\in\N}\wh a_k$. By \er{5.82} \fe $\ve\in\R_{>0}$ \te s
$\ov k\in\N$ \st $\wh a_{\bar k} < \ov a+\ve$. Since the \sq\ $\zb{\wh a}k\N$ is de\cre\ by \er{8.5}, we infer that \fe $\ve>0$ \te s $\ov k\in\N$ \st
$a_l\nde8.7 \le \wh a_l\le \wh a_k <\ov a+\ve$ \fa $l\ge\ov k$. Hence $\{l\in\N: a_l\ge \ov a+\ve\} \sbs \{l\in\N: l<\ov k\}$. Since $\{l\in\N:
l<\ov k\}$ is finite (possibly empty), we find that $\{l\in\N: a_l\ge \ov a+\ve\}$ is finite by Theorem 1.4.18\,(i). Thus \er{8.1} holds with
$b:=\ov a$.

\ssk
``\er{8.2}'': \ Suppose, for \cd ion, that \te s $\ve\in\R_{>0}$ \st $\{n\in\N: a_n+\ve>\ov a\}$ is finite. Then \te s $N\in\Na$ \st $a_n+\ve\le
\ov a$ \fa $n>N$. \csq, $a_n<\ov a-\ve$, $n>N$, hence $\supl_{n>N}a_n \le \ov a-\ve$ (since $\ov a-\ve$ is an \ub\ for the set~$A:=\{a_n\in\R:
n\in\N,\ n>N\}$ and $\sup A$ is the least \ub\ for~$A$). However, $\supl_{n>N}a_n = \wh a_{N+1} \ge \infl_{k\ge0}\wh a_k = \ov a$, hence we find
$\ov a\le \ov a-\ve$, a~\cd ion. Thus, \er{8.2} holds with $b:=\ov a$. In view of \E\df\ \rf{d8.1} and Lemma \rf{l8.2}, we find that \er{8.3} holds.
\endproof

We recall that if $\zb an\N$ is a \cg nt \sq\ in~$\R$ with limit~$a$, then the \sq\ $\zb{-a}n\N$ \cg s to~$-a$ (see \er{7.7} with $a_n:=1$ \fa
$n\in\N$, or note that ${|(-a)-(-a_n)| = |{-a}+a_n| = |a_n-a|}$, $n\in\N$). If the \sq\ $\zb an\N$ is not \cg nt, but is bounded, $\limu a_n$ exists
as well as $\limu(-a_n)$, since $|{-a_n}|=|a_n|\le M$, $n\ge0$, \fs $M\in\R_{>0}$. From \er{8.5} we infer $\limu(-a_n) = \infl_{k\ge0} \supl_{n\ge k}
(-a_n)\nde5.16 = \infl_{k\ge0}(-\infl_{n\ge k} a_n) \nde5.15 = -\supl_{k\ge0}\infl_{n\ge k}a_n$. Set
\beq8.8
c:=\sup_{k\ge0}\inf_{n\ge k}a_n.
\e

\bex8.7
Show that $c$ defined in \er{8.8} \sf ies:
\bea8.9
\hbox{the set }&\{n\in\N: a_n+\ve \le c\} \hbox{ is finite},\\
\hbox{the set }&\{n\in\N: a_n < c+\ve\} \hbox{ is infinite}.\lb{8.10}
\e
\eex

\bdf8.8
Let $\zb an\N$ be a \sq\ of \el s of~$\Q$ (resp.~$\R$). An \el\ $c\in\Q$ (resp.~$\R$) \sf ying \fe $\ve\in\Q_{>0}$ (resp.\ $\R_{>0}$) both \cn s
\er{8.9} and \er{8.10}, is called an \ti{inferior limit\/} (or \ti{limit inferior}) or \ti{lower limit} of the \sq\ $\zb an\N$.\index{limit!inferior}\index{limit!lower}
\edf

\bns8.9
Usual notations for an inferior limit or lower limit of $\zb an\N$ are $\liml_{n\to\iy}a_n$ and $\liminf\limits_{n\to\iy}a_n$. We will use the
notation $\liml_{n\ge0}a_n$ or simply $\liml a_n$.
\ens

\bpr8.10
Let $\zb an\N$ be a \sq\ of \el s of~$\Q$ $($resp.~$\R)$. Then the \fw\ assertions hold\dw

\hph i,ii, The \sq\ $\{a_n\}$ possesses at most one \lg, and in this case the \sq\ $\{a_n\}$ is \bb.

\hph ii,i, If $a_n\to a\in \Q$ $($resp.\ $\R)$ then $a=\liml a_n=\limu a_n$.

If $\{a_n\}$ is a \emph{bounded} \sq\ in~$\R$, then the \fw\ assertions hold\dw

\hph iii,, $\liml(-a_n) = -\limu a_n$\sd \ $\limu (-a_n)=-\liml a_n$,

\hph iv,, $\liml a_n = \supl_{k\ge0}\infl_{n\ge k} a_n$,

\hph v,i, $\liml a_n \le \limu a_n$.
\epr

\proof \

``(i), (ii)'': \ The proofs of (i) and (ii) are similar to the proofs of Lemma \rf{l8.2} and of Lemma \rf{l8.4}. \E\Tf they are omitted.

``(iii)'': \ Replace $a_n$ by $-a_n$ in \er{8.9}, \er{8.10}.

``(iv)'': \ Follows from (iii) and \er{8.3}.

``(v)'': \
Setting $\hat a_k:=\supl_{n\ge k}a_n$, $\check a:=\infl_{n\ge k}a_n$, we find $\hat a_k\ge a_k$, $k\in\N$, by \er{8.7}. Similarly, $\check
a_k\le a_k$, $k\in\N$, hence $\check a_k\le\hat a_k$, $k\in\N$. We recall that the \sq\ $\zb{\hat a}k\N$ is de\cre\ by~\er{8.5}. Similarly, one
finds that the \sq\ $\zb{\check a}k\N$ is in\cre. \E\Tf for $l\ge k$, $l,k\in\N$, we have $\hat a_l\ge\check a_l\ge\check a_k$, hence $\infl_{l\ge
k}\hat a_l \ge \check a_k$. Since $\zb{\hat a}k\N$ is de\cre, $\infl_{l\ge k}\hat a_l = \infl_{l\ge0}\hat a_l$, we obtain $\limu a_n\ge\check a_k$
\fa $k\in\N$. Since $\limu a_n$ is an \ub\ for the set $\{\check a_k\in\R: k\in\N\}$, we obtain $\limu a_n\ge \supl_{k\in\N}\check a_k = \liml a_n$.
\endproof

The next theorem provides us with an alternate \df\ of a limit of a \sq\ in~$\R$ without using the notion of metric.

\bth8.11
Let $\zb an\N$ be a \sq\ in $\R$, and let $a\in\R$. Then the \sq\ $\zb an\N$ \cg s to~$a$ $($see \E\df\ \rf{d7.1}$)$ iff the \fw\ holds\dw

\hph i,i, \te\ $\check M,\hat M \in \R$ \st $\check M\le a_n\le \hat M$, $n\in\N$,

\hph ii,, $\liml_{n\ge0} a_n = a = \limu_{n\ge0}a_n$.
\eth

\proof \

``Only if\/'': \ From Lemma \rf{l8.4} we infer $\limu a_n=a$, hence from Lemma \rf{l8.2} the \sq\ $\zb an\N$ is \ba. Note that the \sq\ $\zb{-a}n\N$
\cg s to~$-a$. \E\Tf similarly we obtain $\limu(-a_n)=-a$, hence $\liml a_n=a$ by \E\Pr\ \rf{p8.10}\,(iv). From $\limu(-a_n)=-a$, we infer that the
\sq\ $\zb{-a}n\N$ is \ba, hence $\zb an\N$ is \bb.

``If\/'': \ From $\limu a_n=a$, we obtain \er{8.1} with $b:=a$, and from Exercise \rf{ex8.7} we obtain \er{8.9} with $c:=a$. In view of~\er{8.1},
\fe $\ve\in\R_{>0}$ \te s $N_1\in\N$ \st $a_n-a<\ve$ \fa $n\ge N_1$, and in view of~\er{8.9} \te s $N_2\in\N$ \st $a-a_n<\ve$ \fa $n\ge N_2$. Setting
$N:=\max(N_1,N_2)$, we arrive at $|a-a_n|<\ve$ \fa $n\ge N$, that is, $a_n\to a$.
\endproof

\bpr8.12
Let $\zb an\N$, $\zb bn\N$ be bounded \sq s in~$\R$. Then the \fw\ assertions hold\dw

\hph i,v, If $a_n\le b_n$ \fa $n\in\N$, then
\beq8.11
\limu a_n\le \limu b_n\sd \q \liml a_n\le \liml b_n.
\e

\hph ii,i, $\limu(a_n+b_n) \le \limu a_n+\limu b_n$\sd \ $\liml(a_n+b_n) \ge \liml a_n+\liml b_n$.

\hph iii,, If \fa $n\in\N$, $c_n:= a_{n+m}$ \fs $m\in\N$, then
\beq8.12
\limu c_n=\limu a_n\sd \q \liml c_n=\liml a_n.
\e

\hph iv,, If $b_n\to b$ then
\beq8.13
\bca
\limu(a_n+b_n) = \limu a_n + b,\\
\liml(a_n+b_n) = \liml a_n + b.
\eca
\e

\hph v,i, If \fs $k\in\N$, $a_{n+k}\ge0$ \fa $n\in\N$, and if $\limu a_n=0$, then $a_n\to 0$.
\epr

\proof \

``(i)'': \ Since $a_l\le b_l$ for $l\in\N$, we have $\supl_{n\ge k}a_n \le \supl_{n\ge k}b_n$ \fa $k\in\N$. Indeed, we have $a_i\le b_i\le
\supl_{l\ge k}b_l$ \fa $i\in\N$ \sf ying $i\ge k$ \fa $k\in\N$. Hence $\supl_{i\ge k}a_i\le \supl_{l\ge k}b_l$ \fa $k\in\N$, which proves the claim.
Using the notation \itd in \er{8.4} we obtain $\wh a_k\le \wh b_k$ \fa $k\in\N$. In view of \er{8.5} we have $\wh a_i\le \wh b_k$ \fa $i\le k$.
\E\Tf $\infl_{i\ge0}\wh a_i \le \infl_{i\le k}\wh a_i \le \wh b_k$ \fa $k\in\N$, hence $\infl_{l\ge0}\wh a_l\le \wh b_k$ \fa $k\in\N$. We claim
$\infl_{l\ge0}\wh a_l \le \infl_{k\ge0}\wh b_k$. Indeed, \fa $k\in\N$, we have  $-\wh b_k\le -\infl_{l\ge0}\wh a_l \nde5.15 =
\supl_{l\ge0}(-\wh a_l)$. From \er{6.52} we infer $\supl_{k\ge0}(-\wh b_k) \le \supl_{l\ge0}(-\wh a_l)$. Hence $\infl_{l\ge0}\wh a_l = -\supl_{l\ge0}
(-\wh a_l) \le -\supl_{k\ge0}(-\wh b_k) = \infl_{k\ge0} \wh b_k$. \csq, $\infl_{l\ge0}\wh a_l \le \infl_{k\ge0}\wh b_k$, which proves the claim.
From \er{8.3} we obtain the part of \er{8.11}. For the second part, use $\limu a_n = -\liml(-a_n)$.

``(ii)'': \ We have $\supl_{n\ge k}(a_n+b_n) \le \supl_{n\ge k}a_n+ \supl_{n\ge k}b_n$ \fa $k\in\N$. Indeed, \fa $k\in\N$ and \fa $n\ge k$ we have
$a_n+b_n \nde3.9 \le (\supl_{n\ge k}a_k) + b_n = b_n+\supl_{n\ge k}a_n \nde3.9 \le \supl_{n\ge k}b_n + \supl_{n\ge k}a_n = \supl_{n\ge k}a_n +
\supl_{n\ge k}b_n$. Using \er{6.52} and \er{8.4} we obtain $\supl_{n\ge k}(a_n+b_n)\le \wh a_k+\wh b_k$, $k\in\N$. Set $\wh c_k:=\supl_{n\ge k}
(a_n+b_n)$, and from \er{8.5}, \er{8.3} we infer that $\wh c_k\da \limu(a_n+b_n)$, $\wh a_k\da \limu a_n$ and $\wh b_k\da \limu b_n$. \If that
$\limu(a_n+b_n) \le \wh a_k+\wh b_k$ \fa $k\in\N$. As in the proof of~(i) we find $\limu(a_n+b_n) \le \infl_{k\ge0}(\wh a_k+\wh b_k)$. \E\fa $k,l
\in\N$ \st $k\le l$ we have $\wh a_k\ge \wh a_l$, $\wh b_k\ge \wh b_l$ by~\er{8.5}, hence $\wh a_k+\wh b_k \ge \wh a_l+\wh b_l$ in view of~\er{3.10}.
\Mo the \sq s $\zb{\wh a}k\N$, $\zb{\wh b}k\N$ and $\zb{\wh a_k+\wh b}k\N$ are \bb. Hence $\wh a_k\da\infl_{k\ge0}\wh a_k$, $\wh b_k\da\infl_{k\ge0}
\wh b_k$ and $\wh a_k+\wh b_k \da \infl_{k\ge0}(\wh a_k+\wh b_k)$. We conclude the proof of~(ii) by using the \fw\ result.

{\it Let $\zb sn\N$, $\zb tn\N$ be bounded \sq s of~$\R$, and let $\a,\b\in\R$. If\/ $s_n\ua \a$ $($resp.\ $s_n\da\a)$ and $t_n\ua\b$ $($resp.\
$t_n\da\b)$, then}
\beq8.14
s_n+t_n \ua \a+\b \q \hbox{(resp.\ } s_n+t_n\da \a+\b).
\e

Observe that since $\infl_{n\ge0}t_n = -\supl_{n\ge0}(-t_n)$, it suffices to prove the first assertion of \er{8.14}. \E\fe $\ve\in\R_{\ge0}$ \te\
$N_1,N_2\in\N$ \st $\a-s_n<\frac\ve 2$ \fa $n\ge N_1$ and $\b-t_n<\frac\ve2$ \fa $n\ge N_2$. Setting $N:=\max(N_1,N_2)$ we find that $\a+\b
-(s_n+t_n) < \ve$ \fa $n\ge N$. Hence $s_n+t_n \ua \a+\b$.

\E\Tf $\infl_{k\ge0}(\wh a_k+\wh b_k) = \infl_{k\ge0}\wh a_k + \infl_{k\ge0}\wh b_k$, that is, $\limu(a_n+b_n) = \limu a_n + \limu b_n$. For the
second part, use $\liml a_n = -\limu(-a_n)$.

``(iii)'': \ $\limu c_n = \limu a_{n+m} = \infl_{k\ge0}\supl_{n\ge k}a_{n+m} = \infl_{k\ge0}\supl_{l\ge k+m}a_l = \infl_{p\ge m}\supl_{l\ge p}a_l
\nad*= \infl_{p\ge0}\supl_{l\ge p}a_l = \limu a_n$. In~$\nad*=$ we used \er{6.51} for $\inf \si_n$ instead of $\sup\si_n$ (recall $\infl_{n\ge0}
\si_n = -\supl_{n\ge0}(-\si_n)$).

``(iv)'': \ We have $\limu b_n=b$ by Lemma \rf{l8.4}. Hence
\beq8.15
\limu(a_n+b_n) \nad{\rm(ii)} \le \limu a_n + \limu b_n = \limu a_n + b.
\e
\E\oh from \er{7.1} we find that \fe $\ve\in\R_{>0}$ \te s $N\in\N$ \st $b_n > b-\ve$ \fa $n\ge N$. Hence \fa $n\ge0$ we have $a_{n+N} + b_{n+N}
\ge a_{n+N}+b-\ve$. \E\Tf $\limu(a_n+b_n) \nad{\rm(iii)}= \limu(a_{n+N}+b_{n+N}) \nad{\rm(i)}\ge \limu(a_{n+N}+b-\ve) = \infl_{k\ge0}\supl_{n\ge k}
(a_{n+N}+b-\ve) \nde6.53 = \infl_{k\ge0}(\supl_{n\ge k}a_{n+N}+b-\ve) \nad*= \infl_{k\ge0} \bigl((\supl_{n\ge k} a_{n+N}) +
(b-\ve)\bigr)=\bigl(\infl_{k\ge0}
\supl_{n\ge k}a_{n+N}\bigr)+b-\ve \nad{\rm(iii)}= \limu a_n+b-\ve$. In~$\nad*=$ we used Exercise \rf{ex6.17}\,(ii). Thus \fe $\ve\in\R_{>0}$, we have
$\limu(a_n+b_n) \ge (\limu a_n)+b-\ve$. Choosing $\ve:=\frac1{k+1}$ \fa $k\in\N$, we find $\limu(a_n+b_n)\ge(\limu a_n)+b$. Otherwise, we would have
$\limu(a_n+b_n) < \limu a_n+b$, hence $\limu a_n+b = \limu(a_n+b_n)+\a$ \fs $\a\in\R_{>0}$. In this case, we would have $-\a\ge-\frac1{k+1}$, \ev tly
$\frac1{k+1}\ge\a>0$ \fa $k\in\N$, which is impossible since $\R$ is \Ar. Finally, from $\limu(a_n+b)\ge \limu a_n+b$ and \er{8.15}, we obtain
\er{8.13}.


For the second parts of (iii) and (iv), use $\liml a_n = -\limu(-a_n)$.

``(v)'': \ Setting $a_n:=0$ and $b_n:=a_n$, $n\in\N$, in \er{8.11}, we find $a_n\to0$, hence $0=\liml a_n\nad*= \liml a_{n_k} \nad{**} \limu a_{n_k}
\nad{*{*}*}= \limu a_n=0$. Hence $\liml a_n=0=\limu a_n$, and $a_n\to0$. In~$\nad*=$ we used (iii),~(v), in~$\nad{**}=$ we used (v)~in \E\Pr\
\rf{p8.10} and $\nad{*{*}*}=$ we used (iii).
\endproof

\bex 8.13
Show that $\limu(a_n\cdot b_n) = a\cdot \limu b_n$ if $a_n\to a$ with $a_n,a>0$ \fa $n\in\N$, and $\zb bn\N$ is a bounded \sq\ in~$\R$.
\eex

\brm8.14
In an \os\ $(X,\ge)$ an \et y $A=B$ with $A,B\in X$ is proved by showing $A\le B$ and $B\le A$. In~$(\R,\ge)$ it often happens that one of these
in\et ies, say $B\le A$, is a con\sq\ of the family of in\et ies $B\le A+\ve$ \fa $\ve\in\R_{>0}$ or $B\le A+\frac1{k+1}$ \fa $k\in\N$ since
$(\R,\ge)$ is \Ar. The proof of \E\Pr\ \rf{p8.12}\,(iv) is an example thereof.
\erm

In the next lemmata we investigate some \pp ies of sub\sq s of bounded \sq s in~$\R$.

\blm8.15 \

\hph i,ii, {\rm (see \cite[pp.\ 48, 49]{Choquet})} Every \sq\ $\zb an\N$ in $\R$ has a strictly in\cre\ or a de\cre\ sub\sq.

\hph ii,i, Every bounded \sq\ $\zb an\N$ in $\R$ has a \cvg\ sub\sq.

\hph iii,, Let $\zb an\N$ be a \sq\ in $\R$ that is \bb\ and \sf ies $\liml a_n = \inf a_n$. Then there is a sub\sq\ $\{a_{n_k}\}_{k\in\N}$
converging to $\infl_{n\ge0}a_n$. Such a sub\sq\ is sometimes called a \emph{minimizing} \sq. Similarly, a \sq\ $\zb an\N$ in~$\R$ that is \ba\ and
\sf ies $\limu a_n=\sup a_n$, possesses a sub\sq\ converging to $\supl_{n\ge0}a_n$.
\elm

\proof \

(i) A term $a_k$ of the \sq\ $\zb an\N$ is called a \ti{peak\/} for the index $k\in\N$ (see~\cite{Choquet}) if $a_k\ge a_n$ \fa $n>k$. We distinguish
three cases. The set of peaks is either empty, or finite, or infinite.

\ti{Step} 1. ``There are no peaks''. \E\fe $n\in\N$ set
\beq8.16
A_n:=\{l\in\N: l>n \hbox{ and }a_l>a_n\}.
\e
Since there are no peaks, the sets $A_n$, $n\in\N$, are nonempty. By Theorem 1.3.21 the sets $A_n$ possess a least \el, which we denote by~$\si(n)$.
We next apply Theorem 1.1.7 with $(E,e,S):=(\N,0,S)$ where $S(n):=n+1$, $n\in\N$, $F:=\N$, $a:=0$, and $f(n):=\si(n)$, $n\in\N$. Denoting by~$\F$
the self-map of~$\N$ defined in (1.1.7), (1.1.8), we find that
\beq8.17
\F(n+1) = \si(\F(n)), \ n\in\N, \q \hbox{and }\F(0)=0.
\e
Since $\si(n)\in A_n$, we have $\F(n+1)>\F(n)$ \fa $n\in\N$, hence $\{a_{\F(n)}\}_{n\in\N}$ is a sub\sq\ of $\zb an\N$. \Mo since $\si(n)\in A_n$,
we have $a_{\F(n+1)} > a_{\F(n)}$, $n\in\N$, hence the sub\sq\ $\{a_{\F(n)}\}_{n\in\N}$ is strictly in\cre.

\ti{Step} 2. ``The set of peaks is nonempty and finite''. Let $m\in\Na$ denote the \nm\ of peaks (i.e.\ the \ca\ of the set of peaks), and let
$k$ denote the index of the greatest peak, which exists by Theorem 1.3.38 since a finite set is \ba. Indeed, by \E\df\ 1.4.17, a \nf\ is \ep\ to an
ordered \il\ $[0,n]$ of $(\N,\ge)$ \fs $n\in\N$. Since $[0,n]$ is nonempty and bounded, its image under a bi\jn\ is bounded by Theorem 1.3.47.
\If that the \sq\ $\zb cn\N$ in~$\R$ defined by $c_n:=a_{n+k+1}$, $n\in\N$, has no peaks. Therefore, in view of Step~1 of the proof, the \sq\
$\zb cn\N$ possesses a sub\sq, which is strictly in\cre\ or de\cre. Note that $(n+1)+k+1 > n+k+1$ \fa $n\in\N$, hence $\zb cn\N$ is a sub\sq\ of the
\sq\ $\zb an\N$.

\ti{Step} 3. ``The set of peaks is infinite''. We first show that an infinite subset~$A$ of~$\N$ is unbounded. By \E\df\ 1.4.17 an infinite set is
nonempty and not finite. If there would exist $M\in\N$ \st $a\le M$ \fa $a\in A$, then $A$~would be a subset of the \il\ $[0,M]$, which is finite by
\E\df\ 1.4.17. However, a subset of a finite set is finite by Theorem 1.4.18. A~\cd ion, since $A$~is infinite. \E\Tf it follows from Theorem 1.3.48
that \te s an \ois sm $\rho:\N\to A$. Thus the \sq\ $\{a_{\rho(n)}\}_{n\in\N}$ \sf ies $\rho(n+1)>\rho(n)$, $n\in\N$. \Mo since $a_{\rho(n)}$ is a
peak \fa $n\in\N$, we have $a_{\rho(n+1)} \ge a_{\rho(n)}$, $n\in\N$. Hence this sub\sq\ is de\cre.

\ssk (ii) Follows from (i) and Lemma \rf{l7.2}\,(iv).

\ssk (iii) \E\fe $k\in\N$ set $A_k:=\{l\in\N: l>k \hbox{ and } a_l-\inf a_n < \frac1{k+1}\}$. We claim that the sets $A_k$ are nonempty. \E\fe
$k\in\N$, let $\zb{c^{(k)}}n\N$ be the \sq\ defined by $c\en k_n:=a_{n+k+1}$, $n\in\N$. In view of \E\Pr\ \rf{p8.12}\,(iii) and~(v), we have
$\liml_{n\ge0}c\en k_n = \liml_{n\ge0}a_n$ \fa $k\in\N$. Since $\liml_{n\ge0} a_n=\infl_{n\ge0}a_n$, we have $\liml_{n\ge0}c\en k_n = \infl_{n\ge0}
a_n$. Then \fe $k\in\N$ \te s $\ov l\in\N$ \st $c\en k_{\bar l} - \frac1{k+1} < \liml_{n\ge0}c\en k_n$ in view of \er{8.10}, hence \st
$c\en k_{\bar l} - \frac1{k+1} < \infl_{n\ge0}a_n$, \ev tly \st $c\en k_{\bar l} - \inf a_n < \frac1{k+1}$. Since $c\en k_{\bar l} = a_{\bar l+k+1}$,
we find that $a_{\bar l+k+1} \in A_k$ \fa $k\in\N$. Setting $\si(k):=\ov l+k+1$, where $\ov l$ depends on~$k$, we have $0\le a_{\si(k)}- \inf a_n
<\frac1{k+1}$ \fa $k\in\N$, hence $a_{\si(k)}\to \inf a_n$.
\endproof

\bex8.16
Show that (i) of Lemma \rf{l8.15} holds if we replace in\cre\ (resp.\ de\cre) by de\cre\ (resp.\ in\cre). (Hint: use $a=-(-a)$, $a\in\R$.)
\eex

\bpr 8.17
Let $\zb an\N$ be a bounded \sq\ in~$\R$. Then the \fw\ holds\dw

\hph i,ii, \E\te s a sub\sq\ of $\zb an\N$ converging to $\limu a_n$.

\hph ii,i, \E\te s a sub\sq\ of $\zb an\N$ converging to $\liml a_n$.

\hph iii,, If $\{a_{n_k}\}_{k\in\N}$ is a sub\sq\ of $\zb an\N$, then $\liml a_n\le \liml a_{n_k} \le \limu a_{n_k} \le \limu a_n$.

\E\Ip if $a_{n_k}\to L$, then $\liml a_n\le L \le \limu a_n$.
\epr

\bxa8.18
Let $a_n\in\R$, $n\in\N$, be defined by $a_{3k}=-1$, $a_{3k+1}=0$, $a_{3k+2}=+1$ \fa $k\in\N$. Then $|a_n|\le1$, $n\in\N$, $\liml a_n=-1$,
$\limu a_n=+1$ and $\lim\limits_{k\to\iy}a_{3k+1}=0$.
\exa

\proof[Proof of \E\Pr\ \rf{p8.17}] \

(iii) ``$\limu a_{n_k}\le \limu a_n$'': \ Let $\si:\N\to \N$ be a strictly in\cre\ \sq. Then \fe $k\in\N$ we have
\beq8.18
\sup_{n\ge k}a_{\si(n)} \le \sup_{n\ge k} a_n.
\e
Indeed, the set $A_k:=\{a_{\si(n)}\in\R: n\ge k\}$ is a subset of the set $B_k:=\{a_n\in\R: n\ge k\}$ \fe $k\in\N$. Since $B_k$ is \ba, $\sup B_k$
exists as well as $\sup A_k$. \Mo from $A_k\sbs B_k$ we infer that $\sup B_k$ is an \ub\ for~$A_k$ since $\sup B_k$ is an \ub\ for~$B_k$, and
$\sup A_k\le \sup B_k$ since $\sup A_k$ is the least \ub\ for~$A_k$. \csq, we have
\beq8.19
\sup_{n\ge k}a_{\si(n)} \le \sup_{n\ge k} a_n \ \hbox{\fa $k\in\N$,\q and }\sup_{n\ge k}a_{\si(n)}\da\limu a_{\si(n)}, \ \sup_{n\ge k}a_n \da
\limu a_n,
\e
by Lemma \rf{l7.2}\,(iv).

Note that $\supl_{n\ge k}a_{\si(n)} \to \limu a_{\si(n)}$, and $\supl_{n\ge k}a_n \to \limu a_n$. By Lemma \rf{l8.4} we have $\limu_{k\ge0}
\supl_{n\ge k}a_{\si(n)} = \limu a_{\si(n)}$, and $\limu_{k\ge0}\supl_{n\ge k}a_n = \limu a_n$. From \er{8.11} and \er{8.5} we obtain
$\limu_{k\ge0}\supl_{n\ge k}a_{\si(n)} \le \limu_{k\ge0}\supl_{n\ge k}a_n$, hence $\limu a_{\si(n)}\le \limu a_n$.

``$\liml a_{\si(n)}\ge \liml a_n$'': \ Follows from what precedes, and $\liml a_n = -\limu(-a_n)$.

\csq, we obtain
$$
\liml a_n \le \liml a_{\si(n)} \le \limu a_{\si(n)} \le \limu a_n.
$$
If, moreover, $a_{n_k}\to L$, then $\liml a_{n_k}=L = \limu a_{n_k}$ by Theorem \rf{t8.11}, hence
$$
\liml a_n \le \liml a_{n_k} = L = \limu a_{n_k} \le \limu a_n,
$$
which completes the proof of part (iii) of \E\Pr\ \rf{p8.17}.

\ssk
We next prove (ii). Note that $\infl_{n\ge0}a_n = \supl_{k=0}\infl_{n\ge k}a_n \le \supl_{k\ge0}\infl_{n\ge k}a_n = \liml a_n$. We distinguish two
cases: (1)~$\liml a_n = \inf a_n$ and (2)~$\liml a_n> \inf a_n$.

(1) is a con\sq\ of Lemma \rf{l8.15}\,(ii).

(2) The set $A:=\{l\in\N:a_l<\liml a_n\}$ is not empty. Indeed, if we set $\a:=\frac12(\liml a_n+\inf a_n)$, and $\ve:=\frac12(\liml a_n-\inf a_n)>0$,
\te s $k\in\N$ \st the \fw\ holds:
\beq8.20
\inf a_n < a_k < \a = (\inf a_n)+\ve < \liml a_n,
\e
by Lemma \rf{l6.28} and $\inf a_n = -\sup(-a_n)$ (see \er{5.16}). If $A$~is finite, we define a \sq\ $\zb cn\N$ in~$\R$ by setting
\beq8.21
c_n := a_{n+\#(A)+1}, \q n\in\N.
\e
Then $\zb cn\N$ is a sub\sq\ of $\zb an\N$, since $(n+1)+\#(A)+1 > n+\#(A)+1$ \fa $n\in\N$, and $\liml c_n=\liml a_n$ by \er{8.12}. \Mo $c_k \ge
\limu a_n$ \fa $k\in\N$. We may thus apply part~(1) with $a_n:=c_n$, $n\in\N$. Hence \te s a sub\sq\ of $\zb cn\N$ converging to $\liml c_n =
\lim a_n$. However, a sub\sq\ of a sub\sq\ of $\zb an\N$ is a sub\sq\ of $\zb an\N$, since the \cm\ of two in\cre\ self-maps of~$\N$ is also strictly
in\cre. It remains to prove the case where $A$~is infinite.  As in the proof of part~(iii) of Lemma \rf{l8.15}, we may assume that \te s a strictly
in\cre\ self-map $\rho$ of~$\N$ \st $\{n\in\N: a_{\rho(n)}\in\R\} = A$. Setting $d_k:=a_{\rho(k)}$, $k\in\N$, we find that the \sq\ $\zb dk\N$ is
a sub\sq\ of $\zb an\N$ \st $d_k<\liml a_n$, $k\in\N$. In view of \er{8.11} with $a_n:=d_n$ and $b_n:=\a$ \fa $n\in\N$, where $\a:=\liml a_n$, we find
$\limu d_n\le \limu \a = \lim\a = \a$ since $b_n\to\a$, hence $\a=\limu\a$ by Lemma \rf{l8.4}. Hence $\limu d_k\le \liml a_n$. \E\oh since $\zb dn\N$
is a sub\sq\ of $\zb an\N$, we have $\liml d_k \ge \liml a_n$ by part~(iii) of \E\Pr\ \rf{p8.17}. \E\Tf $\liml a_n \le \liml d_k \nad*\le \liml d_k
\le \liml a_n$. In~$\nad*\le$ we used \E\Pr\ \rf{p8.10}\,(v). \csq\ $\liml d_k = \limu d_k= \liml a_n$. From Theorem \rf{t8.11} we infer $d_k \to
\liml a_n$. This completes the proof of part~(ii) of \E\Pr\ \rf{8.17}.

\ssk
(i) Part (i) follows from $\limu a_n = -\liml(-a_n)$ (see \E\Pr\ \rf{p8.10}\,(iii)) and from part~(ii).
\endproof

\bex8.19 \

\hph i,i, Let $\zb un\N$, $\zb vn\N$ and $\zb wn\N$ be \sq s in~$\R$ and let $\a\in\R$. Suppose $u_n\to\a$, $v_n\to\a$ and $u_n\le w_n\le v_n$,
$n\in\N$. Show that $w_n\to\a$.

\hph ii,, Let $\zb an\N$ be a bounded \sq\ in~$\R$. Suppose $\liml a_n > \inf a_n$, and let $\zb bn\N$ be a \sq\ \sf ying $\inf a_n\le b_k \le
\liml a_n$ \fa $k\in\N$. Show that $\supl_{m\ge0}\infl_{k\le m} a_m \le b_l \le \liml a_n$ \fa $l\in\N$, and $b_k \to \liml a_n$.
\eex

\brm8.20
In the proofs of preceding lemmata, \Pr s or theorems, we used the argument of the form ``since $\liml a_n = -\limu(-a_n)$, $\infl_{n\ge0}a_n =
-\supl_{n\ge0}(-a_n)$, we infer \dots''. For example, in the proof of Theorem \rf{t8.11}, \E\Pr\ \rf{p8.12}\,(v), \E\Pr\ \rf{p8.17}. Since
$a_n\le b_n$ implies $-b_n \le -a_n$, replacing $a_n$ by~$-a_n$, and $b_n$ by~$-b_n$ amounts to replace $\le$~by~$\ge$. In what follows, we shall
say that the proof of \dots\ follows by ``duality''.
\erm

As a first application of Lemma \rf{l8.15}\,(ii) or \E\Pr\ \rf{p8.17} we consider a first improvement of Lemma \rf{l5.41}.

\bdf8.21
Let $X$ be a \ns\ and let $f$~be a map from~$X$ into~$\R$. The map~$f$ is said to \ti{have an absolute} or \ti{a global minimum} at $c\in X$ (resp.\
\ti{an absolute} or \ti{a global maximum} at $d\in X$) if\index{global minimum}\index{global maximum}
\beq8.22
f(c) \le f(y) \hbox{ (resp.\ } f(d)\ge f(y)) \q \hbox{\fa} y\in X.
\e
\edf

\blm8.22
Let $a,b\in\R$, $a<b$, and let $f:[a,b]\to\R$ \sf y \er{5.86}, then \te s $c\in[a,b]$ \st $f$~has a global minimum at~$c$.
\elm

\proof
A necessary \cn\ for the map $f$ to have a global minimum is that the set $f([a,b])$ be \bb. We first show that $f([a,b])$, thanks to \cn\ \er{8.6},
is bounded. Indeed, \fa $y\in[a,b]$ we have $|f(y)-f(a)| \le M|y-a|$, hence $|f(y)| = |f(y)-f(a)+f(a)| \le |f(y)-f(a)|+|f(a)| \le M|y-a|+|f(a)|
\nad*\le M(b-a)+|f(a)|$. We recall that $\vf(x):=|x|$, $x\in\R$, is a \vl\ (see Exercise \rf{ex5.24}\,(ii)), hence $|x+y|\le|x|+|y|$, $x,y\in\R$.
\E\oh since $a\le y\le b$, we have $0=a-a\le y-a\le b-a$ by~\er{3.9} and $|y-a|=y-a\le b-a$ by \df\ of~$|\,\cdot\,|$. Thus $\nad*\le$ follows from
\er{3.9} and \er{8.10}. \csq, we obtain \fa $y\in[a,b]$:
\beq8.23
-M(b-a) - |f(a)| \le f(y)\le M(b-a)+|f(a)|.
\e

Thanks to the order-\cp ness of $(\R,\ge)$, the set $f([a,b])$ has an infimum (see Corollary \rf{c5.12}), which we denote by~$\a$. Then $\a$~is not
only a \lo\ for $f([a,b])$, but $\a$~is also the greatest \lo\ for $f([a,b])$. \If from Lemma \rf{l5.38} that \fe $\ve\in\R_{>0}$ \te s $y_\ve\in
f([a,b])$ \st $y_\ve-\a<\ve$. \E\Ip \fe $n\in\N$, \te s $z_n\in f([a,b])$ \st the \fw\ holds:
\beq8.24
0 \le z_n-\a < \frac1{n+1} \qh{\fa $n\in\N$.}
\e

\E\fe $n\in\N$ \te s $x_n\in[a,b]$ \st
\beq8.25
f(x_n) = z_n \qh{\fa $n\in\N$.}
\e
We recall that $\frac1{n+1}\da0$, hence $\frac1{n+1}\to0$, and $\liml\frac1{n+1} = 0 = \limu\frac1{n+1}$. From $0\le z_n-\a$ we have $0\nde8.11 \le
\liml(z_n-\a)$ and from $z_n-\a\le\frac1{n+1}$, we obtain $\limu(z_n-a)\nde8.11 \le 0$. Hence $0\le \liml(z_n-\a)\le \limu(z_n-\a)\le0$, which implies
$z_n-\a \to 0$, or $z_n\to\a$.

Since $x_n\in[a,b]$, we have $a\le x_n\le b$, hence the \sq\ $\{x_n\}_{n\ge0}$ is bounded.

From Lemma \rf{l8.15}\,(ii) or \E\Pr\ \rf{p8.17} we find that \te s $c\in\R$ and a sub\sq\ $\{x_{n_k}\}_{k\in\N}$ \st $x_{n_k}\to c$. Since $a\le
x_{n_k}\le b$ \fa $k\in\N$, we have $a=\limu a \nde8.11 \le \limu x_{n_k} \nde8.11 \le \limu c=c$. Since $\limu x_{n_k}=\lim x_{n_k}$ by Lemma
\rf{l8.4}, we obtain $\limu x_{n_k}=\lim x_{n_k} =c$. \csq, $a\le c\le b$, that is, $c\in[a,b]$.

Finally, we have $x_{n_k}\to c\in[a,b]$, $f(x_{n_k})\to\a$ by Lemma \rf{l7.2}\,(iii). It remains to show that
\beq8.26
f(\lim x_{n_k}) = \lim f(x_{n_k}).
\e
It turns out that \er{8.26} follows from \cn\ \er{5.85}. Indeed, $0\le |f(x_{n_k}) - f(c)| \le M|x_{n_k}-c|\to0$, \Tf $0=\limu 0 \nde8.11 \le
\limu |f(x_{n_k})-f(c)| \nde8.11 \le \limu M|x_{n_k}-c| = \lim M|x_{n_k}-c|=0$. Hence $\limu|f(x_{n_k})-\a|=0$ and $f(x_{n_k})\to\a$ by \E\Pr\
\rf{p8.12}\,(iv). \If that the map~$f$ has a global minimum at $c\in[a,b]$.
\endproof

\bco8.23
Under the \as s of Lemma \rf{l8.22}, the map $f$ has a global maximum at some $d\in[a,b]$.
\eco

\proof
Observe that the map $g(x):=-f(x)$, $x\in[a,b]$, \sf ies \er{5.85} and use ``duality'' (see Remark~\rf{r8.20}).
\endproof

Combining Lemma \rf{l5.41}, Lemma \rf{l8.22} and Corollary \rf{c8.23}, we obtain the \fw

\bpr8.24
Under the \as s of Lemma \rf{l8.22}, \te\ $c,d\in[a,b]$ \st
\beq8.27
f([a,b]) = [f(c),f(d)].
\e
\epr

\brm8.25
\E\Pr\ \rf{p8.24} is a strict \ext\ of Lemma \rf{l5.41}. Indeed, let
$$
a:=-1; \q b:=\tfrac12; \q f(t):= \bca
1+t &\hbox{for }t\in[-1,0],\\
1-t &\hbox{for }t\in(0,\frac12],
\eca
$$
then $f$ \sf ies \er{5.85} with $M=1$ (the proof is left to the reader) and $f([-1,\frac12]) = [0,1]$. The map~$f$ has an absolute minimum (resp.\
maximum) at $c:=-1$ (resp.\ $d:=0$). Thus \er{8.27} holds. However, from Lemma \rf{l5.41} we obtain $\bigl(\min(f(-1),f(\frac12)),\max(f(-1),
f(\frac12))\bigr)\break\sbs f((-1,\frac12))$, that is, $(0,\frac12) \sbs(f(0),f(\frac12)) = (0,\frac12) \ne (0,1) = f((0,\frac12))$, which proves the
claim.

\Mo the case $f(a)=f(b)$ is not included in Lemma \rf{l5.41}. If $a:=-1$, $b:=1$ and
$$
g(t):=\bca
1+t &\hbox{for } t\in[-1,0],\\
1-t &\hbox{for } t\in \oz 0,1 ,
\eca
$$
then $f([-1,1]) = [0,1]$ by \E\Pr\ \rf{p8.24}.
\erm

\proof[Proof of \E\Pr\ \rf{p8.24}]
By Lemma \rf{l8.22} and Corollary \rf{c8.23} \te\ $c,d\in[a,b]$ \st
\beq8.28
f(c) \le y \le f(d) \qh{\fa $y\in[a,b]$.}
\e
If $f(c)=f(d)$ then $f([a,b]) = [c,d] = \{c\} = \{d\}$. If $f(c)<f(d)$ then, by Lemma \rf{l8.22} with $a:=\min(c,d)$, $b:=\max(c,d)$, we obtain
$(f(c),f(d)) \sbs f((c,d))$, hence $(f(c),f(d)) \sbs f([a,b])$. Since $c,d\in[a,b]$, we also have $[f(c),f(d)] \sbs f([a,b])$. \E\oh we have
$f([a,b]) \sbs [f(c),f(d)]$ by \er{8.28}. \csq, $f([a,b]) = [f(c),f(d)]$.
\endproof

It turns out that Lemma \rf{l8.22}, Corollary \rf{c8.23} and \E\Pr\ \rf{p8.24} admit a further improvement. We need some \df s.

\bdf8.26
A \ms\ $(X,d)$ is called \tb{\sql\ compact} if every \sq\ $\zb xn\N$ in~$X$ has a \cvg\ sub\sq\ (see \E\df\ \rf{d7.15}).\index{sequentially!compact}
\edf

\bdf8.27
Let $(X,d)$ be a \ms. A~map $f:X\to\R$ is called \ti{\sql\ lower-semicontinuous} (s-l.s.c.) (resp.\ \ti{\sql\ upper semi\ctn} --- s-u.s.c.) if \fe \sq
\ $\zb xn\N$ in~$X$ converging to some \el\ $x\in X$, we have\index{sequentially!lower-semi\ctn}\index{sequentially!upper-semi\ctn}
\beq8.29
f(x)\le \liml f(x_n) \q \hbox{(resp.\ } f(x)\ge \limu f(x_n)).
\e
\edf

Observe that if $(X,d)$ is a \ms\ and $f:X\to\R$ is both s-u.s.c.\ and s-l.s.c., then \fe \sq\ $\zb xn\N$ in~$X$ converging to a limit~$x$, we have
$f(x_n)\to f(x)$. Indeed,
\beq8.30
f(x) \le \liml f(x_n) \nad*\le \limu f(x_n) \le f(x),
\e
hence $f(x) = \liml f(x_n) = \limu f(x_n)$, and $f(x_n) \to f(x)$ by Theorem \rf{t8.11}. In $\nad*\le$ we used \E\Pr\ \rf{p8.10}\,(v).

\bdf8.28
Let $(X,d_X)$ and $(Y,d_Y)$ be \ms s. A~map $f:X\to Y$ is called \tb{\sql\ \ctn\ at $x\in X$}, if \fe \sq\ $\zb xn\N$ in~$X$ \sf ying $x_n\nad{d_X}
\to x$, we have $f(x_n)\nad{d_Y}\to f(x)$. A~map $f:X\to Y$ that is \sql\ \ctn\ at every $x\in X$ is called \ti{s-\ctn}.\index{sequentially!continuous}
\edf

\bex8.29 \

\hph i,ii, Let $f,g:(\R,d)\to\R$. Show that if $f(x):= y_0\in\R$ \fa $x\in X$, and $g(x):=x$ \fa $x\in\R$,
then $f,g$ are s-\ctn. Show that if $f,g:(\R,d) \to\R$ are s-\ctn, then $f+g$, $f\cdot g$, $-f$ with
$$
\bal
(f+g)(x) &:= f(x)+g(x),\\
(f\cdot g)(x) &:= f(x)\cdot g(x),\\
(-f)(x) &:= -f(x)
\eal
$$
\fa $x\in X$, are s-\ctn. \E\Ip \pl\ maps from $\R$ into~$\R$ are s-\ctn.

\hph ii,i, Let $f:\R\to\R$ be an in\cre\ map. Show that $f$~is s-\ctn\ at $x\in\R$ iff $\supl_{y<x} f(y) = f(x) = \infl_{z>x}f(z)$. See \er{6.123}.

\hph iii,, Let $A$ be a subset of $\R$ and let
$$
\1_A(x) := \bca
0 &\hbox{if } x\notin A,\\
1 &\hbox{if } x\in A.
\eca
$$
Show that $\1_A$ is s-\ctn\ iff $A:=\vn$ or $A:=\R$. Show that if $A:=\Q$ or $A:=\R\sm \Q$ then $\1_A$~is s-\ctn\ at \ti{no} $x\in\R$. Show that
$\1_A$~is s-l.s.c.\ if $A:=(0,1)$, and s-u.s.c.\ if $A:=[0,1]$.
\eex

\blm8.30
Let $(X,d)$ be a \sql\ compact \ms, and let $f:X\to\R$ be \sql\ lower-semi\ctn\ $($resp.\ \sql\ upper-semi\ctn$)$. Then \te s $c\in X$ $($resp.\
$d\in X)$ \st
\beq8.31
f(x)\ge c \q (\hbox{resp.\ $f(x)\le d)$ \fa} x\in X.
\e
\elm

\proof
Since $f=-(-f)$ and $-f$ is s-u.s.c.\ iff $f$ is s-l.s.c., it suffices to prove the case where $f$~is s-l.s.c. Suppose, for \cd ion, that $f(X)$ is
not \bb. In that case, \fe $n\in\N$ \te s $x_n\in X$ \st $f(x_n)\le -n$. However, \te s a sub\sq\ $\{x_{n_k}\}_{k\in\N}$ converging to some $\ov x
\in X$, in view of the sequential compactness of $(X,d)$. \E\Tf $f(x_{n_k}) \to f(\ov x)$. By Lemma \rf{l7.2} \te s $M\in\R_{>0}$ \st
$f(x_{n_k})\ge -M$ \fa $k\in\N$. We have $-M\le f(x_{n_k})\le -n_k$, \ev tly, $n_k \le M$ \fa $k\in\N$, a~\cd ion in view of Lemma 2.1.32.
\endproof

Let $\zb an\N$ be defined by $a_n:=\frac1{n+1}$, $n\in\N$. We already know that $\infl_{n\ge0}a_n = 0 = \lim a_n$ since $(\R,\ge)$ is \Ar. The set
$\{a_n: n\in\N\}\cup\{0\}$ is \sql\ compact, since every sub\sq\ \cg s to~$0$ by Lemma \rf{l7.2}\,(iii). The set $A:=\{a_n:n\in\N\}$ is thus
\ti{not\/} \sql\ compact, since $\zb an\N$ is not \cvg\ in~$A$.

\bdf8.31
A \nss\ $A$ of a \ms\ $(X,d)$ is called \tb{\sql\ closed\/} if \fe \sq\ $\zb an\N$ in~$A$ converging to an \el\ $\ov a\in X$, we have $\ov a\in A$.
\E\Ip $X$~is \sql\ closed in $(X,d)$. The empty set is \sql\ closed.\index{sequentially!closed}
\edf

\blm8.32
Let $(X,d)$ be a \sql\ compact \ms. Then \fe non\-empty subset $A$~of~$X$, $(A,d)$ is \sql\ compact iff $A$~is \sql\ closed.
\elm

\bex8.33 \

\hph i,i, Prove Lemma \rf{l8.32}.

\hph ii,, Let $X:= [x_1,x_2] \cup [x_3,x_4]$ with $x_1,x_2,x_3,x_4\in\R$, $x_1<x_2<x_3<x_4$. Show that $(X,d)$ as a subset of $(\R,d)$, where $d$~is
the usual metric on~$\R$, is \sql\ compact. Let $f(x):=x$, $x\in X$. Then, clearly, $|f(x)-f(y)|\le M|x-y|$ with $M:=1$. Show that $f(X)\ne
[f(x_3),f(x_4)]$ as in \E\Pr\ \rf{p8.24}.
\eex

As a con\sq\ the set $[a,b]$ has a \pp y not shared by the set~$X$ above. We will see that the missing \pp y is the ``connectedness'' of~$X$. Before
doing
this, we show that \cn\ \er{5.86} in \E\Pr\ \rf{p8.24} can be weakened. Indeed, let $a,b\in\R$, $a<b$, and $f:[a,b]\to\R$ be s-\ctn. Then $f$~is both
s-l.s.c.\ and s-u.s.c. Indeed, if $\zb xn\N$ is a \sq\ in~$[a,b]$ \st $x_n\to x$ \fs $x\in[a,b]$, then $f(x_n)\to f(x)$. \csq, $f(x)=\lim f(x_n)
\nad*= \liml f(x_n) \le \liml f(x_n)$, hence $f(x)\le \liml f(x_n)$. In~$\nad*=$ we used \E\Pr\ \rf{p8.10}\,(ii). \If that $f$ is s-l.s.c. The proof
that $f$~is s-u.s.c.\ can be done by duality (see Remark \rf{r8.20}).

We recall that for a map $f:[a,b]\to\R$ \cn\ \er{5.86} implies s-continuity. This \cn\ can be extended to maps between two \ms s.

\bdf8.34
Let $(X,d_X)$ and $(Y,d_Y)$ be \ms s. A~map $f:X\to Y$ is called \tb{\Hc} of order $\a\in\R_{>0}$, $\a\le1$, if \te s $M\in\R_{>0}$ \st\index{continuous!H\"older-}
\beq8.32
d_Y(f(x_1),f(x_2)) \le M(d_X(x_1,x_2))^\a \qh{\fa $x_1,x_2\in X$.}
\e
If $\a=1$ then the \f~$f$ is called \tb{\Lc}.\index{continuous!Lipschitz-}
\edf

\bex8.35 \

\hph i,ii, Let $f:X\to Y$ and $\a\in\oz0,1 $ be as in \E\df\ \rf{d8.34}. Show that $f$~is s-\ctn.

\hph ii,i, Let $X:=[0,1]_\R$ and let $f_\a(x):=x^\a$, $x\in[0,1]$. Show that $f_\a:[0,1]\to\R$ is \Hc\ of order $\b\in(0,\a]$.

\hph iii,, Let $f:[0,1]\to\R$ be defined by $f(0):=0$ and $f(x):=x\log_{10}x$, $x\in\oz 0,1 $. Is there $\a\in\oz0,1 _\R$ \st $f$~is \Hc\ of
order~$\a$, \wrt the usual metric on~$\R?$ Is $f$ \sql\ \ctn$?$ Is the \rt ion of~$f$ to $[a,1]$, $a\in(0,1)_\R$, \Hc\ \fs $\a\in\oz0,1 ?$ $($See
Section {\rm\ref{ass.9}.)}
\eex

Our next goal is to prove the \fw\ \ext\ of \E\Pr\ \rf{p8.24}.

\bth8.36
Let $a,b\in\R$, $a<b$, and let $f:[a,b]\to\R$ be s-\ctn. Then \te\ $c,d\in\R$ \st the \fw\ holds\dw
$$
f([a,b]) = [f(c),f(d)].
$$
\eth

We will give a proof based on the \fw\ lemma.

\blm8.37
Let  $(X_i,d_{(i)})$, $i=1,2$, be \ms s, and $X:=X_1\t X_2$ be the product of the sets $X_1$ and~$X_2$. For $x:=(x_{(1)},x_{(2)})$, $x_{(i)}\in X_i$,
$i=1,2$, set $\pi_i(x):=x_{(i)}$, $i:=1,2$, the \pj s of~$x$ on~$X_i$, $i=1,2$ $($see {\rm(4.2.17))}. Then the \fw\ \as s hold\dw

\hph i,ii, \E\fa $x,y\in X$ set
\beq8.33
d(x,y):= \max_{i=1,2} d_{(i)}(x_{(i)},y_{(i)}).
\e
Then $d:X\to\R_{\ge0}$ is a metric on~$X$, and for $i=1,2$, $\pi_i:X\to X_i$ is a sur\jc\ \Lc\ map with Lipschitz constant~$M$ equal to~$1$.

\hph ii,i, If $(X,d)$ is \sql\ compact, then so are $(X_1,d_{(1)})$ and $(X_2,d_{(2)})$.

\hph iii,, If $(X_i,d_{(i)})$, $i=1,2$, are \sql\ compact, then so is $(X,d)$.
\elm

\bex8.38
Prove assertions (i), (ii) of Lemma \rf{l8.37}.
\eex

\proof[Proof of \rm(iii)]
Let $\zb xn\N$ be a \sq\ in $X$. Then $\{\pi_1(x_n)\}_{n\in\N}$ is a \sq\ in~$X_1$. Since $(X_1,d_{(1)})$ is s-compact, \te s a strictly in\cre\
self-map of~$\N$ denoted by~$\si_1$ \st $\{\pi_1(x_{\si_1(n)})\}_{n\in\N}$ \cg s to some \el\ $\ov x_{(1)}\in X_1$. \If that $\{\pi_2(x_{\si_1(n)})
\}_{n\in\N}$ is a \sq\ in~$X_2$. Since $(X_2,d_{(2)})$ is s-compact, \te s a strictly in\cre\ self-map of~$\N$ denoted by~$\si_2$ \st $\{\pi_2
(x_{\si_2\circ \si_1(n)})\}_{n\in\N}$ \cg s to some \el\ $\ov x_{(2)}\in X_2$. We claim that the \sq\ $\{x_{\si_2\circ\si_1(n)}\}_{n\in\N}$ in~$X$
\cg s to $\ov x:=(\ov x_{(1)},\ov x_{(2)})\in X$ \wrt the metric~$d$. Note that $\si:=\si_2\circ \si_1$ is a strictly in\cre\ self-map of~$\N$, thus
$\{x_{\si(n)}\}_{n\in\N}$ is a sub\sq\ of the \sq\ $\zb xn\N$. \Mo the \sq\ $\{\pi_1(x_{\si(n)})\}_{n\in\N}$ is a sub\sq\ of the \sq\
$\{\pi_1(x_{\si_1(n)})\}_{n\in\N}$, hence by Lemma \rf{7.2}\,(iii), $\pi_1(x_{\si(n)})\h{\nad{d_{(1)}}\to}\h \ov x_{(1)}$. We also have
$\pi_2(x_{\si(n)})\h{\nad{d_{(2)}}\to}\h \ov x_{(2)}$. \E\Tf $d_{(1)}(\pi_1(x_{\si(n)}),\ov x_{(1)})\h{\to}\h 0$, and
$d_{(2)}(\pi_2(x_{\si(n)}),\ov x_{(2)})\h{\to}\h 0$.
Setting $\ov x:=(\ov x_{(1)},\ov x_{(2)})\in X$, we obtain
$$
0\le d(x_{\si(n)},\ov x) \nde8.33 = \max_{i=1,2} d_{(i)}(\pi_i(x_{\si(n)}),\ov x_{(i)}) \le \sum_{i=1}^2 d_{(i)}(\pi_i(x_{\si(n)}),\ov x_{(i)}).
$$
\csq,
\bmlg
0 = \limu_{n\ge0} 0 \nde8.11 \le \limu_{n\ge0}\sum_{i=1}^2 d_{(i)}(\pi_i(x_{\si(n)}),\ov x_{(i)}) \\ \nde8.13 \le
\limu_{n\ge0} d_{(1)}(\pi_1(x_{\si(n)}),\ov x_{(1)}) + \limu_{n\ge0}d_{(2)}(\pi_2(x_{\si(n)}),\ov x_{(2)}) = 0+0 =0.
\e
By \E\Pr\ \rf{p8.12}\,(v) we find $d(x_{\si(n)},\ov x) \to 0$. Hence $(X,d)$ is \sql\ compact.
\endproof

\Wanp prove Theorem \rf{t8.36}.

\proof[Proof of Theorem \rf{t8.36}]
Let $X:=[a,b]$, $a,b\in\R$, $a<b$, and let $d(x,y):=|x-y|$, $x,y\in[a,b]$. We recall that $(X,d)$ is \sql\ compact. Indeed, if $\zb xn\N$ is a \sq\
in~$[a,b]$, i.e.\ $a\le x_n\le b$, $n\in\N$, it is bounded, \Tf by Lemma \rf{l8.15}\,(ii) \te s $\ov x\in\R$ and a sub\sq\ $\{x_{n_k}\}_{k\in\N}$ \st
$x_{n_k}\to \ov x$. We have $a=\limu_{k\ge0}a \nde8.11 \le \limu_{k\ge0}x_{n_k} = \ov x = \limu_{k\ge0} x_{n_k} \nde8.11 \le \limu_{k\ge0} b = b$.
Hence $a\le\ov x\le b$, which proves the claim.
Since $f:X\to\R$ is s-\ctn, it is both s-lower- and s-upper-semi\ctn, hence by Lemma \rf{l8.30} \te\ $c,d\in[a,b]$ \st $f(c) \le f(x) \le f(d)$ \fa
$x\in[a,b]$. We now show that $[f(c),f(d)]\sbs f([a,b])$. If $f(c)=f(d)$ then clearly $[f(c),f(d)]= \{f(c)\} \sbs f([a,b])$. We now suppose that
$f(c)<f(d)$, hence $c\ne d$. Set $\check c:=\min(c,d)$ and $\hat d:=\max(c,d)$. Set $K:=[\check c,\hat d]$ and $d_K(x,y):=|x-y|$, $x,y\in K$. We
observe that $(K,d_K)$ is s-compact. Indeed, $K$~is \sql\ closed in $([a,b],d)$, since for every \sq\ $\zb yn\N$ in~$K$ \sf ying $y_n\to \ov
y\in[a,b]$, we have $\ov y\in K$ (for example, as a~con\sq\ of~\er{8.11}). By Lemma \rf{l8.32}, $(K,d_K)$~is s-compact. Clearly, $f(c),f(d)\in
f([a,b])$. Let $t\in(f(c),f(d))$, which is not empty since $f(c)<f(d)$. Set
\bea8.34
K_\le &:= \{x\in K: f(x)\le t\},\\
K_\ge &:= \{x\in K: t\le f(x)\}. \lb{8.35}
\e
Then $c\in K_\le$ and $d\in K_\ge$. We claim that both $K_\le$ and $K_\ge$ are nonempty \sql\ closed subsets of~$K$. Indeed, if $\zb yn\N \sbs
K_\le$ and $|y_n-\ov y|\to0$, we have $f(y_n)\to f(\ov y)$ since $f$~is s-\ctn, hence $f(\ov y) = \lim\limits_{n\ge0} f(y_n) = \limu_{n\ge0} f(y_n)
\le \limu_{n\ge0} t = t$. Hence $\ov y\in K_\le$, and $K_\le$~is closed. Similarly, one shows that $K_\ge$~is s-closed. Hence both $K_\le$
and~$K_\ge$ are s-closed subsets of $(K,d_K)$. By Lemma \rf{l8.32} $(K_\le,d_K)$, $(K_\ge,d_K)$ are s-compact. \Mo $K_\le \cup K_\ge = K$ by
\er{8.34}, \er{8.35}. Note that $K_\le \cap K_\ge = \{x\in K: f(x)=t\}$. We want to show that $K_\le \cap K_\ge \ne \vn$. The idea of the proof is
to consider the distance \f\ $g(x,y):=|x-y|$ between $x\in K_\le$ and $y\in K_\ge$. Thus $g$~is a map from $K_\le \t K_\ge$ into $\R_{\ge0}$. By
Lemma \rf{l8.37}, the \ms\ $K_\le \t K_\ge$ endowed with the metric defined in \er{8.33} is s-compact. We recall that the metric on~$K_\le$
(resp.~$K_\ge$) is the \rt ion to $K$ of the usual metric on~$\R$. \E\Tf if $x,x'\in K_\le$ and $y,y'\in K_\ge$, then $d_{K_\le}(x,x') = |x-x'|$
and $d_{K_\ge}(y,y') = |y-y'|$. \Mo $|g(x,y) - g(x',y')| = \bigl||x-y|-|x'-y'|\bigr|$. Using $\bigl||a|-|b|\bigr| \le |a-b|$, $a,b\in\R$, we find
$|g(x,y)-g(x',y')| \le |x-y-x'+y'| = |{(x-x')}+(y'-y)| \le |x-x'| + |y'-y| \le 2\max(|x-x'|,|y-y'|) = 2d((x,y),(x',y'))$, $x,x'\in K_\le$ and $y,y'\in
K_\ge$. \If that the map $g:K_\le \t K_\ge \to \R$ is \Lc\ with constant $M:=2$, from the s-compact \ms\ $(K_\le \t K_\ge,d)$ into $\R_{\ge0}$ with
the usual metric. By Lemma \rf{l8.30} \te s $(\ov a,\ov b)\in K_\le\t K_\ge$ \st $g(\ov a,\ov b) \le g(x,y) = |x-y|$ \fa $x\in K_\le$ and $y\in
K_\ge$. Suppose, for \cd ion, that $g(\ov a,\ov b)>0$. Set $\ov e:=\frac12(\ov a+\ov b)$. Since $\ov a$~and~$\ov b$ belong to $K$, we have
$\ov e\in K = K_\le \cup K_\ge$. If $\ov e\in K_\le$, then $g(\ov e,\ov b) = |\frac12(\ov a+\ov b )-\ov b| = \frac12|\ov a-\ov b|
< g(\ov a,\ov b)$, a~\cd ion. If $\ov e\in K_\ge$, $g(\ov a,\ov e) = |\ov a-\frac12(\ov a+\ov b)| = \frac12|\ov a-\ov b| < g(\ov a,\ov b)$, a~\cd ion.
\E\Tf $g(\ov a,\ov b)=0$, hence $\ov a=\ov b\in K_\le \cap K_\ge$ and $f(\ov a)=f(\ov b)=t$. Since $t= f(\ov a) = f(\ov b) \in(f(c),f(d))$ is
arbitrary, $(f(c),f(d)) \sbs f([a,b])$, hence $[f(c),f(d)]\sbs f([a,b])$. Since $f([a,b])\sbs [f(c),f(d)]$, we obtain $f([a,b])= [f(c),f(d)]$.
\endproof

So far we used the notion of \ti{open} subset in~$\Q_{>0}$ (see \E\df\ \rf{d4.7}) in Lemma \rf{l4.10}, and of \ti{open} subset of the \il~$(0,1)$
in~$\R$ (see \E\df\ \rf{d5.8}) in Lemma~\rf{l5.9}. We recall that the notion of \ti{\il s}, as subsets of an \os\ $(E,\ge)$, is \itd in Notation
1.3.41, (1.3.25)--(1.3.30). Note that if $a=b$ in~(1.3.30), $(a,b)$ is the empty set. Since, by \df, a bounded set is nonempty (see \E\df\ 1.3.36),
the bounded \il s of~$\Q$ and~$\R$ are of the form $[a,a]$ ($=\{a\}$), $(a,b)$, $\zo a,b $, $\oz a,b $ and $[a,b]$ with $a<b$. If $E:=\Q$ (resp.~$\R$)
and if $a,b\in E$, $a<b$, we set
\bga8.36
\rho_{a,b}(r):= (1-r)a+rb = a+r(b-a),\\
\si_{a,b}(s):= \frac{s-a}{b-a}\,, \lb{8.37}
\e
with $r\in[0,1]_E$, $s\in[a,b]_E$. The map $\rho_{a,b}$ is a bi\jn\ from the \il\ $[0,1]_E$ onto the \il\ $[a,b]_E$, and $\si_{a,b}$ is the inverse
of~$\rho_{a,b}$. We also have
\bea8.38
0 &< \rho_{a,b}(r_2) - \rho_{a,b}(r_1) \le (b-a)(r_2-r_1), \q 0\le r_1<r_2\le 1,\\
0 &< \si_{a,b}(s_2) - \si_{a,b}(s_1) \le (b-a)\mo (s_2-s_1), \q a\le s_1<s_2\le b. \lb{8.39}
\e

Clearly, the \rt ion $\rho_{a,b}$ to $(0,1)$ is a bi\jn\ between $(0,1)$ and $(a,b)$, similarly for $\si_{a,b}$. The maps defined in \er{8.36},
(8.37) allow us in some cases to ``transfer'' \pp ies from \il s $[0,1]$ (resp.~$(0,1)$) to \il s $[a,b]$ (resp.~$(a,b)$). See, for example, the proof
of Lemma \rf{l5.41}. Inspection of \E\df s \rf{d4.7} and~\rf{d5.8} motivates the \fw\ one.

\bds8.39 \

\hph i,ii, Let $A$ be a \nss\ of~$\Q$ (resp.~$\R$). An \el\ $x\in A$ is called an \ti{interior point\/} (or \ti{\el\/}) of~$A$ if \te\ $\a,\b\in\Q$
(resp.~$\R$), $\a<\b$, \st the \fw\ assertion holds:\index{interior point}
\beq8.40
x\in (\a,\b) \sbs A.
\e
The set of interior points of $A$ is called the \ti{interior} of the set~$A$ and is denoted by $\kl A$ (or~$A^\circ$).

\hph ii,i, A \nss\ of $\Q$ (resp.\ $\R$) is called \ti{open} if it is equal to its interior, i.e.\ $A=\kl A$.

\hph iii,, The empty set $\vn$ is open by \df.
\eds

\bnt8.40
We will denote by $\cO(\Q)$ (resp.~$\cO(\R)$) the collection of open subsets of~$\Q$ (resp.~$\R$).
\ent

It turns out that $\cO(\Q)$ (resp.\ $\cO(\R)$) is a topology on the set~$\Q$ (resp.~$\R$).

\begin{dfn}[\cite{Du}, \cite{JvN}] \lb{d8.41}
A \tb{topology} on a \ns~$X$ is a collection of subsets of~$X$, denoted by~$\cO$, with the \fw\ \pp ies:\index{topology}

\bit
\item[$(O_1)$] \ $\vn$ and $X$ belong to $\cO$\sd
\item[$(O_2)$] \ $\cO\sms\vn$ is closed under taking unions\sd
\item[$(O_3)$] \ $\cO$ is closed under taking \ti{finite} intersections.
\eit
The pair $(X,\cO)$ is called a \ti{\tps}. A~subset $U$ of~$X$ is said to be open if $U\in\cO$ and \ti{closed\/} if its complement in~$X$ is open.

A \tps\ $(X,\cO)$ is said to be \ti{Hausdorff\/} if \fe two distinct points $x_1,x_2\in X$ the \fw\ assertion holds:\index{Hausdorff \tps}
\bit
\item[$(O_4)$] \ \te\ disjoint open sets $U_1,U_2\in\cO$ \st $x_1\in U_1$ and $x_2\in U_2$.
\eit
\edf

We are mainly interested in \Hts s.

\brm8.42
Let $X:=\Q$ (resp.\ $\R$). Observe that if $U\in\cO(X)\sms\vn$ then, in view of \er{8.40}, \fe $x\in U$ \te s an open \il\ $I:=(\a,\b)$, $\a<\b$,
depending on~$x$ \sf ying $x\in I\sbs X$. Clearly, if $y\in I$ then $y\in I\sbs I$, hence $I\in\cO(X)$. Let $\cB(X)$ denote the subcollection of
$\cO(X)\sms\vn$ consisting of \il s $I=(\a,\b)$, $\a<\b$, $\a,\b\in X$. The \fw\ assertions hold:
\bea8.41
{}&\hbox{\E\fe $x\in X$ \te s $I\in\cB(X)$ \st $x\in I$.}\\
&\hbox{\E\fe non-disjoint pair $I',I''\in\cB(X)$, \te s $I\in\cB(X)$ \sf ying}\lb{8.42}\\
&\q \hbox{$I\sbs I'\cap I''$.}\non\\
&\hbox{\E\fe two distinct points $x_1,x_2\in X$ \te\ $I_1,I_2\in\cB(X)$ \st }\lb{8.43}\\
&\q \hbox{$x_i\in I_i$, $i=1,2$, and $I_1\cap I_2=\vn$.}\non
\e
\erm

\bex8.43
Prove \er{8.41}--\er{8.43}.
\eex

\blm8.44
Let $X$ be a \ns\ and let $\cB$ be a collection of subsets of~$X$ \sf ying \er{8.41}, \er{8.42} and \er{8.43} with $\cB(X):=\cB$. \E\fe \nss~$A$
of~$X$ set\dw
\beq8.44
A^\circ := \{x\in A: \hbox{\te s $I\in \cB(X)$ \st} x\in I\sbs A\}.
\e
Let $\cO[\cB]$ denote the collection of subsets $A$ of~$X$ consisting of the empty subset of~$X$, and the subsets $A$~of~$X$ \sf ying $A=A^\circ$
whenever $A\ne\vn$. Then $(X, \cO[\cB])$ is a \emph{Hausdorff \tps} \sf ying $\cB\sbs \cO[\cB]$.
\elm

\proof \

``$(O_1)$'': \ $\vn\in\cO[\cB]$ by \df. We now show $X\in\cO[\cB]$. \E\fe $x\in X$ \te s $I\in\cB$ \st $x\in I\sbs X$ by \er{8.41}. Hence $x\in \kl X$
\fe $x\in X$, hence $X=\kl X$ and $X\in\cO[\cB]$.

``$(O_2)$'': \ Let $\zb U\a J$ denote a family of \nss s of~$X$ belonging to $\cO[\cB]$. Let $x\in\bcl_{\a\in J}U_\a$. By \df\ \te s $\ov\a\in J$ \st
$x\in U_{\bar \a}$. Since $U_{\bar \a}\in\cO[\cB]$ and $U_{\bar \a}\ne\vn$, we have $U_{\bar \a}=\kl U_{\bar \a}$, hence \te s $I\in\cB$ \st
$x\in I\sbs U_{\bar a}$. Since $U_{\bar a}\sbs \bcl_{\a\in J}U_\a$, we conclude that $x\in \Bg(\bcl_{\a\in J}U_\a)^\circ$. Since $x$~is an arbitrary
\el\ of $\bcl_{\a\in J}U_\a$, we find that $\bcl_{\a\in J}U_\a\sbs \Bg(\bcl_{\a\in J}U_\a)^\circ$. From \er{8.44} we infer that $\Bg(
\bcl_{\a\in J}U_\a)^\circ \sbs \bcl_{\a\in J}U_\a$, \Tf the sets $\Bg(\bcl_{\a\in J}U_\a)^\circ$ and $\bcl_{\a\in J}U_\a$ are equal, which implies
$\bcl_{\a\in J}U_\a \in \cO[\cB]$.

``$(O_3)$'': \ Let $\zb UiJ$ be a finite (nonempty) family of subsets of~$X$ belonging to~$\cO[\cB]$. If $U_i=\vn$ \fs $i:=i_0\in J$, then $\Bca_{i\in
J}U_i = \vn\in\cO[\cB]$. We suppose that $U_i\ne\vn$ \fa $i\in J$. \Mo if $\Bca_{i\in J'}U_i = \vn$ \fs strict subset~$J'$ of~$J$, then $\Bca_{i\in J}
U_i= \vn \in \cO[\cB]$. \E\Tf we may assume that $\Bca_{i\in J'}U_i \ne\vn$ \fe \nss\ of~$J$. We first prove the case where $\#(J)=2$. Let $U_{i_0},
U_{i_1}\in \cO[\cB]\sms\vn$ be \st $U_{i_0}\cap U_{i_1}\ne\vn$ and $U_{i_0}=U_{i_0}^\circ$, $U_{i_1}=U_{i_1}^\circ$. We have to prove
$$
U_{i_0}^\circ \cap U_{i_1}^\circ = (U_{i_0}^\circ \cap U_{i_1}^\circ)^\circ,
$$
whenever $U_{i_0}\ne\vn$, $U_{i_1}\ne\vn$, $U_{i_0}\cap U_{i_1}\ne\vn$, $U_{i_0}=\kl U_{i_0}$, $U_{i_1}=\kl U_{i_1}$.

We have $U_{i_0}^\circ\cap U_{i_1}^\circ = U_{i_0}\cap U_{i_1}\ne\vn$. Let $x\in U_{i_0}^\circ\cap U_{i_1}^\circ$. By \er{8.44} \te\ $I_0,I_1\in\cB$
\st $x\in I_0\sbs U_{i_0}$ and $x\in I_1\sbs U_{i_1}$. Hence $x\in I_0\cap I_1 \sbs U_{i_0}\cap U_{i_1}$. By \er{8.44} we find that $x\in (U_{i_0}
\cap U_{i_1})^\circ$. Since $x$~is an arbitrary \el\ of $U_{i_0}^\circ \cap U_{i_1}^\circ$, we infer that $U_{i_0}^\circ\cap U_{i_1}^\circ \sbs
(U_{i_0}\cap U_{i_1})^\circ$. \Mo since $U_{i_0} = \kl U_{i_0}$ and $U_{i_1} = \kl U_{i_1}$, we have $U_{i_0}^\circ\cap U_{i_1}^\circ \sbs
(U_{i_0}^\circ\cap U_{i_1}^\circ)^\circ$. Finally, it follows from \er{8.44} that $(U_{i_0}^\circ\cap U_{i_1}^\circ)^\circ \sbs U_{i_0}^\circ\cap
U_{i_1}^\circ$, hence $U_{i_0}^\circ \cap U_{i_1}^\circ = (U_{i_0}^\circ \cap U_{i_1}^\circ)^\circ$,
that is, $U_{i_0}\cap U_{i_1} = (U_{i_0}\cap U_{i_1})^\circ$. \csq, $U_{i_0}\cap U_{i_1} \in \cO[\cB]$, which
completes the proof of the case $\#(J)=2$. The case $\#(J)=n$, $n\in\Na\sms1$, follows by using \In\ on~$n$.

``$\cB\sbs \cO[\cB]$'': \ Let $I\in\cB$, then $I\ne\vn$. Let $x\in I$, then $x\in I=I\in \cB$, hence by \er{8.44} $x\in \kl I$. \E\Tf $I=\{x\in I\}
\sbs \{x\in \kl I\}= \kl I$. Since $I^\circ \sbs I$, we obtain $I=I^\circ$, and $I\in\cO[\cB]$. Thus $\cB\sbs \cO[\cB]$.

``$(O_4)$'': \ Let $x_1,x_2\in X$, $x_1\ne x_2$. By \er{8.43} \te\ $I_1,I_2\in\cB$ \st $x_i\in I_i$, $i=1,2$, and $I_1\cap I_2=\vn$. Since $I_i\in
\cO[\cB]$, $i=1,2$, $(O_4)$~holds.
\endproof

\bxa8.45
Let $X:=\Q$ (resp.\ $\R$), and let $\cB(\Q):=\{(a,b)\sbs \Q \hbox{ \st}a<b\}$ (resp.\ $\cB(\R):=\{(a,b)\sbs \R \hbox{ \st}a<b\}$). Then
\er{8.41}--\er{8.43} hold. Indeed,

``\er{8.41}'': \ \E\fe $x\in\Q$ (resp.\ $\R$) $x\in (x-1,x+1)$ in $\Q$ (resp.\ $\R$).

``\er{8.42}'': \ Let $I_1:=(a,b)$, $a<b$, $a,b\in X$, and let $I_2:=(c,d)$, $c<d$, $c,d\in X$. Suppose that $I_1\cap I_2\ne\vn$. Then \te s $x\in X$
\st $x\in(a,b)\cap(c,d)$, that is, $a<x<b$ and $c<x<d$. \E\Tf we must have $\max(a,c)<x<\min(b,d)$. Suppose that $y\in X$ \sf ies
$\max(a,c)<y<\min(b,d)$,
for example $y:=\frac12(\max(a,c)+\min(b,d))$. Then we have $a\le \max(a,c) <y< \min(b,d)\le b$, and $c\le \max(a,c) <y< \min(b,d)\le d$, hence
$a<y<b$ and $c<y<d$. Set $I:=(\max(a,c),\min(b,d))$. If $y\in I$, we have $y\in I_1$ and $y\in I_2$, hence $I\sbs I_1\cap I_2$.

``\er{8.43}'': \ Let $x_1,x_2\in X$, $x_1\ne x_2$. By exchanging the roles of $x_1$~and~$x_2$ if necessary, we may assume $x_1<x_2$. Set $r:=x_2
-x_1>0$, and $I_i:=(x_1-\frac12r,x_i+\frac12r)$, $i=1,2$. Then $x_i\in I_i$, $i=1,2$, and $I_1\cap I_2=\vn$. Indeed, if $x\in I_1\cap I_2$ then
$x<x_1+\frac12r$ and $x>x_2-\frac12r$, which implies $x_2-\frac12r < x_1+\frac12r$, hence $r=x_2-x_1<r$. A~\cd ion.
\exa

\brm8.46 \

\hph i,i, Recall that in \E\df\ \rf{d4.7} we \itd the notion of open subsets of~$\Q_{>0}$, and in \E\df\ \rf{d5.8} the notion of open subsets of the
\il\ $(0,1)$ in~$\R$.

In these \df s we made use of the notion of interior point without explicitly mentioning it. In \E\df\ \rf{d8.39} the notion of interior point
defined in \er{8.40} is explicit. These notions are formulated in terms of the total \og\ of~$\Q$ and~$\R$. \E\Ip the collection of (nonempty) open
\il s $(\a,\b)$ denoted by $\cB(\Q)$ (resp.\ $\cB(\R)$) is also a notion formulated in terms of the \og\ of~$\Q$ and~$\R$. Lemma \rf{l8.44} shows
that a topology $\cO(\Q)$ (resp.~$\cO(\R)$) can be \itd by means of~$\cB(\Q)$ (resp.~$\cB(\R)$).

The notion of open subset of $\Q$ (resp.\ $\R$) in the sense of \E\df\ \rf{d8.41} obviously depends on the topology. One verifies that the collection
of all subsets of a \ns~$X$, \Ip $\Q$~or~$\R$, is a topology on~$X$. In this case every subset of~$X$ is open. Given a \tps\ $(X,\cO)$ and a \nss\
$A$ of~$X$, an \el\ $x\in A$ is called an interior point of \te s $O\in\cO$ \st $x\in O\sbs A$. If $X:=\Q$ (resp.~$\R$) and $\cO:=\cO(\Q)$ (resp.\
$\cO(\R)$), then the notion of interior point \itd\ in \E\df\ \rf{d8.39} (see~\er{8.40}) coincides with the ``topological'' notion of interior point.
Indeed, if \er{8.40} holds then $x\in(a,b)\sbs A$, $a<b$, hence $x\in O\sbs A$, $O\in\cO(X)$ since $(a,b)\in \cB(X)\sbs \cO(X)$. Conversely, if
$x\in O\sbs A$, $O\in\cO(X)$, then $O=\kl O$, hence \te s $(a,b)$ \st \er{8.40} holds.

The topologies $\cO(\Q)$ on $\Q$ and $\cO(\R)$ on $\R$ are sometimes called \ti{order topologies} or topologies induced by the \og\ of the \tos s
$\Q$~and~$\R$.

\ssk
\hph ii,, Observe that for \nss s of $\Q_{>0}$, we have two \df s of ``open set''. Indeed, if $U$~is a \nss\ of~$\Q_{>0}$, $U$~may be open in the
sense of \E\df\ \rf{d4.7}, or open, as a subset of~$\Q$, in the sense of \E\df\ \rf{d8.41}, with $\cO:=\cO(\Q)$. However, if $U$~is open in the sense
of \E\df\ \rf{d4.7}, \fe $x\in U$ \te\ $\a,\b\in\Q$, $a<\b$, \st $x\in(\a,\b)\sbs U$. Thus in view of Lemma \rf{l8.44} and Example \rf{xa8.45},
$U=U^\circ$, hence $U\in\cO(\Q)$. Conversely, if $U\sbs \Q_{>0}$, $U\ne\vn$ and $U\in\cO(\Q)$, then \fe $x\in U$ \te\ $\a',\b'\in\Q$, $\a'<\b'$,
\st $x\in(\a',\b')\sbs U$. Since $U\sbs \Q_{>0}$, we have $U=U\cap \Q_{>0}$. Thus $(\a',\b')\sbs \Q_{>0}$. \E\Tf $\a'\ge0$ and $\b'>\a'>0$, hence
$\b'\in\Q_{>0}$. Setting $\a:=\frac12(\a'+x)$, we find $\a\in\Q_{>0}$ and $\a<x<\b'$. Thus $x\in(\a,\b')\sbs U$ with $\a<\b'$, $\a,\b'\in\Q_{>0}$.
\csq, $U$~is open in the sense of \E\df\ \rf{d4.7}.
\erm

\bex8.47
Let $U$ be a \nss\ of $(0,1)_\R$. Show that $U$ is open in the sense of \E\df\ \rf{d5.8} iff $U\in\cO(\R)$.
\eex

\bdf8.48
Let $(X,\cO)$ be a Hausdorff \tps\ and let $A$~be a \nss\ of~$X$. Set
\beq8.45
\cO_A := \{U\sbs A: U= A\cap O \hbox{ \fs}O\in \cO\}.
\e
\edf

\bpr8.49
Let $X$, $\cO$ and $\cO_A$ be as in \E\df\ \rf{d8.48}. Then the \fw\ assertions hold\dw

\hph i,i, $(A,\cO_A)$ is a Hausdorff \tps.

\hph ii,, $\cO_A\sbs \cO$ iff $A\in\cO$.
\epr

\bex8.50
Prove \E\Pr\ \rf{p8.49}.
\eex

\bdf8.51
The topology $\cO_A$ defined in \er{8.45} is called the \tb{relative topology}.\index{relative topology}
\edf

\Wanp introduce the notion of \tb{connectedness}\index{connectedness} for subsets of a Hausdorff \tps. Recall that such a notion has been \itd in Lemma \rf{l5.9} for the
subset $(0,1)$ of~$\R$.

\bdf8.52
Let $(H,\cO)$ be a \Hts. The set $X$ is called \ti{connected\/} (\wrt the topology~$\cO$) if it is not the disjoint union of two \nss s belonging
to~$\cO$. A~\ti{subset\/} $A$ of~$X$ is called \ti{connected\/} (\wrt the topology~$\cO$) if it is connected \wrt the relative topology~$\cO_A$.
\edf

We recall that an \of\ $(K,\ge)$ is \Ar\ iff its \Pf\ is order-dense in $(K,\ge)$ (see \E\Pr~\rf{p4.34}). \Mo by Theorem \rf{t4.35}, an \Ar\ field
$(K,\ge)$ is ring- and \ois c to a subfield of the field $(\R,\ge)$. Let $\cB(K)$ denote the collection of subsets of~$K$ defined as in Example
\rf{xa8.45} with~$\Q$ (resp.~$\R$) replaced by~$K$. The proof of \er{8.41}--\er{8.43} for $\cB(K)$ is verbatim the same as for $\Q$ or~$\R$. \E\Tf
we may introduce the order topology $\cO(K)$ on~$K$. Observe that both $\Q$~and~$\R$ are \Ar\ fields. Other examples are given in Example \rf{xa5.16}.
In view of \E\Pr\ \rf{p5.14} an \Ar\ field $(K,\ge)$ is \ois c to $(\R,\ge)$ iff $(K,\ge)$ is order-complete. In what follows we will show that this
\pp y can be replaced by the connectedness of the \tps\ $(K,\cO(K))$.

\newpage
\bth8.53
Let $(K,\ge)$ be an \Ar\ \of. Let $\cB(K)$ denote the collection of \nss s of~$K$ defined by
\beq8.46
\cB(K):= \{(a,b)\sbs K: a,b\in K \hbox{ \st} a<b\}.
\e
Then the collection of subsets of $K$ denoted by $\cO_K$ consisting of arbitrary unions of \el s of~$\cB(K)$ and of the empty subset of~$K$ is a
Hausdorff topology on the field~$K$ called the \emph{order topology} of~$K$.\index{order topology} \Mo the \fw\ assertions are \ev t\dw

\hph i,ii, \fs $a,b\in K$, $a<b$, the \il\ $(a,b)$ is connected \wrt the topology~$\cO_K$\sd

\hph ii,i, \fa $a',b'\in K$, $a'<b'$, the \il\ $(a',b')$ is connected \wrt $\cO_K$\sd

\hph iii,, the \tps\ $(K,\cO_K)$ is connected\sd

\hph iv,, the \of\ $(K,\ge)$ is ring- and \ois c to the field $(\R,\ge)$.
\eth

\proof \

``\ti{$\cO_K$ is a topology on $K$}'': \ The proof of \er{8.41}--\er{8.43} for $\cB(K)$ is verbatim the same as for $\cB(\Q)$ and $\cB(\R)$ in
Example \rf{xa8.45}, \Tf it will be omitted. Note that if $A$~is a \nss\ of~$K$ and if $A$~is the union of \el s of~$\cB(K)$, then $A=A^\circ$.
Indeed, let $A=\bcl_{j\in J}I_j$ with $I_j:=(a_j,b_j)$, $a_j<b_j$, $a_j,b_j\in K$ \fa $j\in J$. Then if $x\in A$, \te s $j_0\in J$ \st $x\in I_{j_0}$,
and $I_{j_0}\sbs \bcl_{j\in J}I_j = A$, hence $I_{j_0}\sbs A$. \E\Tf every $x\in A$ is an interior point, and \cn\ $A=A^\circ$ (defined in \er{8.44})
in Lemma \rf{l8.44} is \sf ied. Thus $A\in \cO[\cB(K)]$. Setting $\cO_K:= \cO[\cB(K)]$, we find that $\cO_K$~is a topology on the field~$K$, and that
$I\in \cO_K$ \fa $I\in \cB_K$.
\endproof

In the proof of ``$\rm (i)\imp(ii)$'' we will use the \fw\ lemma.

\blm8.54
Let $(X,\cO)$ and $(X',\cO')$ be \Hts s, and let ${f:X\h{\to}\h X'}$ be sur\jc. Suppose that the inverse image under~$f$ of every open set of~$X'$
is an open set of~$X$, and that $(X,\cO)$ is connected. Then so is $(X',\cO')$.
\elm

\proof
We suppose that $(X',\cO')$ is not connected and we show that $(X,\cO)$ is not connected. Let $U_1',U_2'$ be two nonempty open subsets of~$X'$ \sf
ying $X'=U_1'\cup U_2'$ and $U_1'\cap U_2'=\vn$. Then $X=f\mo(X') = f\mo(U_1'\cup U_2') = f\mo(U_1')\cup f\mo(U_2')$ and $f\mo(U_1')\cap
f\mo(U_2') = f\mo(U_1'\cap U_2') = f\mo(\vn)=\vn$. Since $f\mo(U_i')$, $i=1,2$, are nonempty open subsets of~$X$ whose intersection is empty and \sf y
$f\mo(U_1')\cup f\mo(U_2')=X$, we find that $(X,\cO)$ is not connected.
\endproof

\proof[Proof of $\rm(i)\imp(ii)$.]
We recall that the map $\rho_{a,b}:[0,1]\to[a,b]$, $a,b\in K$, $a<b$, defined in \er{8.36} is a strictly in\cre\ bi\jn\ with inverse $\si_{a,b}:
[a,b]\to[0,1]$ defined in \er{8.37}, also strictly in\cre. \E\Ip $\rho_{a,b}$ maps open sub\il s $(\a,\b)$ of~$[0,1]$ with $0\le\a<\b\le1$ onto
open sub\il s of~$[a,b]$. Similarly for $\si_{a,b}$. We define a map $f:(a,b)\to(a',b')$ by setting
\beq8.47
f(s):= (\rho_{a',b'}\circ \si_{a,b})(s), \q s\in(a,b).
\e
Then $f$ is bi\jc, in particular sur\jc. We want to apply Lemma \rf{l8.54} with
$$
(X,\cO):=\bigl((a,b),(\cO_K)_{(a,b)}\bigr), \q (X',\cO'):=\bigl((a',b'),(\cO_K)_{(a',b')}\bigr)
$$
where $(\cO_K)_{(a,b)}$ is the relative topology on $(a,b)$ defined in \E\df\ \rf{d8.51}. Similarly for $(\cO_K)_{(a',b')}$. Since $(a,b)$ (resp.\
$(a',b')$) belongs to~$\cO_K$, we have $(\cO_K)_{(a,b)}=\{{(\a,\b)\sbs K}: a\le \a<\b\le b\}$, similarly for $(\cO_K)_{(a',b')}$, by \E\Pr\
\rf{p8.49}. By \as\break $\bigl((a,b),(\cO_K)_{(a,b)}\bigr)$ is connected. In view of Lemma \rf{l8.54}, in order to prove that\break
$\bigl((a',b'),(\cO_K)_{(a',b')}\bigr)$
is connected, it suffices to show that $f\mo(U')\in (\cO_K)_{(a,b)}$ \fe $U'\in (\cO_K)_{(a',b')}$. Note that since $(a',b')\in\cO_K$, an \el\
$U'\in(\cO_K)_{(a',b')}$ is the union of \il s $(\a,\b)$ contained in $(a',b')$. \E\Tf $f\mo(U') = f\mo\bigl(\bcl_{j\in J}(\a_j,\b_j)\bigr) =
\bcl_{j\in J}f\mo((\a_j,\b_j)) = \bcl_{j\in J} (\rho_{a',b'}\circ \si_{a,b})\mo ((\a_j,\b_j)) = \bcl_{j\in J} (\si_{a,b})\mo (\rho_{a',b'})\mo
((\a_j,\b_j)) =\break \bcl_{j\in J} (\si_{a,b})\mo ((\a'_j,\b'_j)) = \bcl_{j\in J} (\a_j'',\b_j'')$, where $(\a_j',\b_j')\sbs(0,1)$,
$(\a_j'',\b_j'')\sbs(a,b)$ \fa $j\in J$. Since $(\a_j'',\b_j'')\in\cO_K$ and $(\a_j'',\b_j'')\sbs(a,b)$ \fa $j\in J$, we have $f\mo(U') =
\bcl_{j\in J}(a_j'',\b_j'')\in(\cO_K)_{(a,b)}$. \If that $\bigl((a',b'),(\cO_K)_{(a',b')}\bigr)$ is a connected \Hts, hence in view of \E\df\
\rf{d8.51}, $(a',b')$ is a connected subset of $(K,\cO_K)$.
\endproof

In the proof of ``$\rm(ii)\imp(iii)$'' we will use the \fw\ lemma.

\blm8.55
Let $(X,\cO)$ be a \Hts\ and let $\zb AiI$ be a family of connected subsets of~$X$. If $\Bca_{i\in I}A_i\ne\vn$, then $\bcl_{i\in I}A_i$ is connected
in $(X,\cO)$.
\elm

\proof
Let $x\in\Bca_{i\in I}A_i$. \E\fe $i\in I$, let $O_1,O_2\in\cO$ be \st $\bcl_{i\in I}A_i = O_1\cup O_2$ and $O_1\cap O_2=\vn$. \E\fe $i\in I$,
$O_1\cap A_i$ and $O_2\cap A_i$ belong to $\cO_{A_i}$. We have $(O_1\cap A_i)\cup(O_2\cap A_i) = (O_1\cup O_2)\cap A_i = X\cap A_i = A_i$ and
$(O_1\cap A_i)\cap(O_2\cap A_i) = (O_1\cap O_2)\cap A_i = \vn \cap A_i = \vn$. Since $A_i$ is connected, we have \fe $i\in I$, $A_i\cap O_1=\vn$ or
$A_i\cap O_2 = \vn$ but not both $A_i\cap O_1 = \vn$ and $A_i\cap O_2 = \vn$ since $A_i\ne\vn$. Hence \fe $i\in I$ either $A_i\cap O_1 \ne \vn$
and $A_i\cap O_2 = \vn$, or $A_i\cap O_1 = \vn$ and $A_i\cap O_2 \ne \vn$. \E\Tf \fe $i\in I$, either $A_i\sbs O_1$ or $A_i\sbs O_2$. Let $I_1:=
\{i\in I: A_i\sbs O_1\}$ and $I_2:=\{i\in I: A_i\sbs O_2\}$. We have $I=I_1\cup I_2$ and $I_1\cap I_2=\vn$. However, we have $x\in A_i$ \fa $i\in I$.
If $x\in A_i$, $i\in I_1$, then $x\in A_i\sbs O_1$, and if $x\in A_i$, $i\in I_2$, then $x\in A_i\sbs O_2$. Since $O_1\cap O_2=\vn$, we have either
$A_i\sbs O_1$ \fa $i\in I$, or $A_i\sbs O_2$ \fa $i\in I$. \If that either $\bcl_{i\in I}A_i=O_1$ and $O_2=\vn$ or $\bcl_{i\in I}A_i=O_2$ and
$O_1=\vn$. \csq, $\bcl_{i\in I}A_i$ is connected in $(X,\cO)$.
\endproof

\proof[Proof of $\rm (ii)\imp (iii)$]
We apply Lemma \rf{l8.55} with $I:=\Na$, $A_n:=(-n,n)$, $n\in\Na$. From (ii) we have $A_n\in\cO_K$, $n\in \Na$. \Mo $0\in(-1,1) = \Bca_{n\in\Na}
A_n$ and $K=\bcl_{n\in\Na}A_n$ since $(K,\ge)$ is \Ar. Indeed, \fe $y\in K$, we have $- y\in K$, hence $|y|\in K$. If $|y|=0$ then $y=0$ and
$y\in A_1$ \fa $n\in\Na$. If $|y|>0$ then in view of \E\df\ 4.5.53 \te s $k\in\Na$ \st $|y|<k$. Hence $y\in(-k,k)\sbs \bcl_{n\in\Na}A_n$. \E\Tf
$K\sbs \bcl_{n\in\Na}A_n$. Clearly, $\bcl_{n\in\Na}A_n\sbs K$ since $A_n\sbs K$ \fa $n\in\Na$. Hence $K=\bcl_{n\in\Na}A_n$, and $(K,\cO_K)$ is
connected by Lemma \rf{l8.55}.
\endproof

\advance\abovedisplayskip by-2pt
\advance\belowdisplayskip by-2pt
\proof[Proof of $\rm(iii)\imp(iv)$]
Suppose that $(K,\cO_K)$ is connected, and we want to prove~(iv). In view of \E\Pr\ \rf{p5.14}\,(ii) it suffices to show that $(K,\ge)$ is a
complete \of. By Corollary \rf{c5.12}, it suffices to show that every \nss\ of~$K$ that is \ba\ has a supremum. Let $A$~be an arbitrary \nss\ of~$K$
that is \ba. We suppose, to obtain a \cd ion, that $A$~has no supremum. \E\Tf we suppose that $A$~has no greatest \el, and that the set of \ti{\ub s}
for~$A$, denoted by~$B$, has no least \el. Note that since $A$ is \ba, the set~$B$ is not empty. \Mo $B\in\cO_K$. Indeed, we will prove that every
\el\ $y\in B$ is an \ipt\ of~$B$. Note that if $y\in B$ and $y'\in K$ \sf ies $y<y'$, then \fe $x\in A$ we have $x<y<y'$, hence $y'$~is an \ub\
for~$A$. Thus, in particular $[y,y+1)\sbs B$. \E\oh since $B$~has no least \el, \te s $z\in B$ \st $z<y$. \E\Tf since $z$~is an \ub\ for~$A$, we
obtain $(z,y]\sbs B$. \csq, we find that $y\in (z,y+1)\sbs B$, hence $y\in \kl B$. Since $y$~is arbitrary in~$B$, we have $B=\kl B$, hence $B\in
\cO_K$ by Lemma \rf{l8.44}. So far we have shown that $B$~is a \ti{nonempty open subset\/} of $(K,\cO_K)$.

We next prove that
\beq8.48
A\sbs B^c.
\e
If $x\in A$ then $x\notin B$. Otherwise $x\in B$, and $x'\le x$ \fa $x'\in A$, hence $x$~would be the greatest \el\ of~$A$, \cd ing the \as\
``$\sup A$ does not exist''. \csq, $B^c$~\ti{is not empty}, since $A\ne\vn$.

In view of the connectedness of $K$, we have:
\beq8.49
B^c \notin \cO_K.
\e
Indeed, otherwise we would have $K=B\cup B^c$, $B\cap B^c=\vn$, $B,B^c$ \nos s of~$K$, \cd ing the connectedness of $(K,\cO_K)$.

We next claim
\beq8.50
\hbox{\E\fa $y\in K$ and $x\in B^c$:\ $y<x$ implies }y\in B^c.
\e
We already observed that if $u\in B$ and $v\in K$, $v\ge u$, then $v\in B$, since $B$~is the set of \ub s for~$A$. \csq, if $u\in B$ and $v\in B^c$,
then $v<u$. Indeed, $v\ne u$ since $B\cap B^c=\vn$, and if $v>u$ then $v\in B$, \cd ing $v\in B^c$. \E\Tf if $u\in B$ and $v\in B^c$ then $v<u$.
\csq, if $x\in B^c$ and $y\in K$ \sf ies $y<x$, then $y\notin B$, and by what precedes we would have $x<y$, \cd ing $y<x$. The claim is proved.

Since $B^c\notin \cO_K$ by \er{8.49}, \te s $\ov x\in B^c$ \st $\ov x$~is not an \ipt\ of~$B^c$. We have $(\ov x-1,\ov x]\sbs B^c$ by~\er{8.50}.
Since $\ov x\notin {\kl B}{}^c$, $[\ov x,z)\not\sbs B^c$ \fa $z\in K$ \sf ying $z>\ov x$, otherwise we would have $\ov x\in(\ov x-1,z)$. \If that
$\ov x$~is the greatest \el\ of~$B^c$, hence $B^c=(\lea,\ov x]$ and $B=(\ov x,\rightarrow)$. By~\er{8.48} we have $A\sbs B^c=(\lea,\ov x]$. If
$\ov x\in A$ then $\ov x$ would be the greatest \el\ of~$A$, hence $\ov x=\sup A$, \cd ing the \as\ $\sup A$ does not exist. We claim that \fe
$y\in K$ \st $y<\ov x$, \te s $x\in A$ \sf ying $y\le x<\ov x$. Otherwise, \te s $\ov{\ov x}\in K$ \st $\ov{\ov x}<\ov x$ and
$(\ov{\ov x},\ov x)\not\sbs A$. Then every \el\ $y$ of $(\ov{\ov x},\ov x)$ is an \ub\ for~$A$, hence $y\in B=(\ov x,{\to})$, a~\cd ion. \If that
if $y<\ov x$, \te s $x\in A$ \sf ying $y<x<\ov x$. Since every $z\in K$ \sf ying $z>\ov x$ is an \ub\ for~$A$, we infer that $\ov x=\sup A$, a
\cd ion.
\endproof

\proof[Proof of $\rm(iv)\imp(i)$]
In view of (iv) \te s an order- and ring-\is sm between $(\R,\ge)$ and $(K,\ge)$. We denote by~$\vf$ this \is sm. \E\Tf $\vf$~maps every open \il\
$(a,b)$ of~$\R$, $a<b$, into an \il\ $(\vf(a),\vf(b))$ of~$(K,\ge)$. Similarly for~$\vf\Inv$. Since every open set~$U$ of $(\R,\cO(\R))$ is the union
of \il s $(a,b)$ of~$\R$ contained in~$U$, the image of~$U$ under~$\vf$ is an open subset of $(K,\cO_K)$ and the image of an open set~$U'$ in
$(K,\cO_K)$ is an open set of~$(\R,\ge)$ under~$\vf\Inv$. In view of Lemma \rf{l5.9} the set $(0,1)$ in~$\R$ is connected \wrt $\cO(\R)$, then we
claim that the set $(\vf(0),\vf(1))$ in $(K,\cO_K)$ is connected. Indeed, suppose, for \cd ion, that \te\ two \nos s of $(\vf(0),\vf(1))$ denoted by
$U_1',U_2'$, \sf ying $(\vf(0),\vf(1))=U_1'\cup U_2'$ and $U_1'\cap U_2'=\vn$. Then one verifies that the subsets of~$(0,1)$, $\vf\Inv(U_1'),
\vf\Inv(U_2')$ are \nos s of~$(0,1)$, \sf ying $\vf\Inv(U_1')\cup \vf\Inv(U_2')=(0,1)$ and $\vf\Inv(U_1')\cap \vf\Inv(U_2')=\vn$. A~\cd ion, since
$(0,1)$ is connected. \E\Tf the subset of~$K$, $(\vf(0),\vf(1))$ is connected in $(K,\cO_K)$.

This completes the proof of the theorem.
\endproof

\advance\abovedisplayskip by2pt
\advance\belowdisplayskip by2pt

\bex8.56
Show that all open \il s (see (1.3.26), (1.3.30)) of a complete \of\ are connected \wrt the order topology.
\eex

Our next goal is to \ch ize \ti{all\/} connected subsets of a complete \of\ endowed with the order topology.

\ssk
We first prove that it suffices to consider the case where the complete \of\ is the field $(\R,\ge)$. Indeed, noticing that a \cp\ \of\ $(K,\ge)$ is
\Ar\ by \E\Pr\ \rf{p5.13}, we find from the part ``$\rm (iv)\imp(i)$'' of the proof of Theorem \rf{t8.53}, that \te s a bi\jn\ $\vf:(\R,\ge)\to
(K,\ge)$ that maps open subsets of~$\R$ \wrt the order topology $\cO(\R)$ into open subsets of~$K$ \wrt the order topology~$\cO_K$, and that
$\vf\Inv$ maps \el s of~$\cO_K$ into \el s of~$\cO(\R)$.

\bdf8.57
Let $(X,\cO)$, $(X',\cO')$ be \tps s. A~bi\jc\ map $\vf:X\to X'$ that maps open sets of $(X,\cO)$ into open sets of $(X',\cO')$, and \st $\vf\Inv$
maps open sets of $(X',\cO')$ into open sets of $(X,\cO)$ is usually called a \tb{homeo\mf}\index{homeo\mf} (or a topological \is sm), and the \tps s $(X,\cO),
(X',\cO')$ are called \ti{homeomorphic}.
\edf

\blm8.58 \

\hph i,i, The \cm\ of \he sms and the inverse of a \he sm is a \he sm.

\hph ii,, A \tps\ $(X,\cO)$ \he c to a Hausdorff $($resp.\ connected\/$)$ \tps\ $(X',\cO')$ is a Hausdorff $($resp.\ connected\/$)$ \tps.
\elm

\bex8.59 \

\hph i,i, Prove Lemma \rf{l8.58}.

\hph ii,, Show that the subsets $\zo0,1 $ and $\oz0,1 $ of $(\R,\cO(\R))$ endowed with their relative topology are \he c.
\eex

\brm8.60
Since for a \ns\ $X$, the power set $\cP(X)$ (i.e., the set of all subsets of~$X$) is a Hausdorff topology on~$X$ (called the \ti{discrete topology}),
the identity map from a \tps\ $(X,\cO)$, $\cO\ne\cP(X)$, maps \el s of~$\cO$ into \el s of~$\cP(X)$, and does not map \el s of~$\cP(X)$ into \el s
of~$\cO$. Hence the identity is not a \he sm in this case.
\erm

The next lemma shows that the connected subsets of $(\R,\cO(\R))$ are nonempty \il s of $(\R,\ge)$ (see (1.3.25)--(1.3.30)).

\blm8.61 \

\hph i,i, Let $a\in\R$. Then $\{a\}$ is connected.

\hph ii,, Let $A$ be a connected subset of $\R$ \wrt $\cO(\R)$. If $A$~contains two distinct points $a$~and~$b$, $a<b$, then $[a,b]\sbs A$.
\elm

\proof
(ii) Suppose, for \cd ion, that \te s $c\in(a,b)$ \st $c\notin A$. Let $U_1:=(\lea,c)\cap A$ and $U_2:=(c,\to)\cap A$. Then $U_1,U_2$ are two open
subsets of~$A$ \wrt the relative topology~$\cO_A$, since $(\lea,c)$ and $(c,\to)$ are open subsets of $(\R,\cO(\R))$ by Theorem \rf{t8.53}. Note that
$a\in U_1$ and $b\in U_2$. \Mo $U_1\cup U_2=A$ and $U_1\cap U_2=\vn$, \cd ing the connectedness of~$A$. Hence $c\in[a,b]$. Since $c$~is arbitrary
in $(a,b)$, we have $(a,b)\sbs A$, hence $[a,b]\sbs A$ since $a,b\in A$.
\endproof

\bco8.62
The connected subsets of $(\R,\cO(\R))$ are \il s.
\eco

\bex8.63
Prove Corollary \rf{c8.62}.
\eex

In view of Theorem \rf{t8.53} it remains to show that \fa $a,b\in\R$, $a<b$, the \il s $[a,b]$, $\zo a,b $, $\oz a,b $, $\oz \lea,b $ and
$[a,\to)$ are connected in $(\R,\cO(\R))$.

\blm8.64
Suppose that \fa $a,b\in\R$, $a<b$, the \il\ $[a,b]$ is connected in $(\R,\cO(\R))$. Then \fa $a,b\in\R$, $a<b$, the \il s $\zo a,b $, $\oz a,b $,
$[a,\to)$, $\oz \lea,b $ are connected in $(\R,\cO(\R))$.
\elm

\proof
Using the fact that the field $(\R,\ge)$ is \Ar, we find that $[a,b)=\break\bcl_{n\ge n_0}[a,b-\frac1n]$, $\oz a,b = \bcl_{n\ge n_1}[a+\frac1n,b]$ \fs
$n_0,n_1\in\Na$, $[a,\to)= \bcl_{n\in\Na}[a,a+n]$, $\oz\lea,b = \bcl_{n\in\Na}[b-n,b]$. Then the conclusion follows from Lemma \rf{l8.55}.
\endproof

It remains to show that the \il s $[a,b]$, $a,b\in\R$, $a<b$, are connected in $(\R,\cO(\R))$. Using, as in the proof of Theorem \rf{t8.53}, the
\ois sm $\si_{a,b}:[0,1]\to[a,b]$, we find that the \tps s $([0,1],\cO(\R)_{[0,1]})$ and $([a,b],\cO(\R)_{[a,b]})$ are \he c. As a con\sq, it suffices
to show that the subset $[0,1]$ of $(\R,\cO(\R))$ is connected. 

We now present three different proofs of the connectedness of~$[0,1]$ in $(\R,\cO(\R))$. The first one, based on the order-completeness of $(\R,\ge)$,
is verbatim the same as the proof of Lemma \rf{l5.9} with $(0,1)$ replaced by~$[0,1]$. The verification of this assertion is left to the reader. The
second proof is based on the sequential compactness of $[a,b]$ \wrt the usual metric in~$\R$. Indeed, from the proof of Theorem \rf{t8.36} we may
extract the \fw\ lemma.

\blm8.64a
Let $K$ be a \sql\ compact subset of~$\R$ endowed with the usual metric, and let $K_1,K_2$ be two \sql\ closed subsets of~$K$ \sf ying $K=K_1\cup
K_2$. Then $K_1\cap K_2\ne\vn$.
\elm

\proof[Sketch of the proof]
We define $g:K_1\t K_2\to\R$ by setting $g(x,y):=|x-y|$, $x\in K_1$, $y\in K_2$. One verifies that $g:K_1\t K_2\to\R$ is \Lc\ whenever $K_1\t K_2$
is endowed with the metric~$d$ defined in~\er{8.33} and $d_i(x_{(i)},y_{(i)}) = |x_{(i)}-y_{(i)}|$, $x_{(i)},y_{(i)}\in\R$, $i=1,2$. Since
$(K_1\t K_2,d)$ is \sql\ compact by Lemma \rf{l8.32}, \te\ $\ov a\in K_1$, $\ov b\in K_2$ \st $g(\ov a,\ov b)\le g(x,y)$ \fa $(x,y)\in K_1\t K_2$ in
view of Lemma \rf{l8.30}. From $K=K_1\cup K_2$ we infer as in the proof of Theorem \rf{t8.36} that $g(\ov a,\ov b)=0$. \If that $\ov a=\ov b$, hence
$K_1\cap K_2\ne\vn$.
\endproof

\ti{Second proof of the connectedness of {$[a,b]$}, $a<b$}.
Let $U_1$ and~$U_2$ be two nonempty relatively open subsets of~$[a,b]$ \sf ying $[a,b] = U_1\cup U_2$. Suppose, for \cd ion, that $U_1\cap U_2=\vn$.
We first state and prove two lemmata which will be used in the proof.

For the convenience of the reader we recall some facts from set theory concerning \cpl s of sets (see for example \cite[pp.~3,~4]{Kelley}).
Let $A,B$ and~$X$ be sets. Then the set $X\sm A:=\{x\in X:x\notin A\}$ \sf ies
\bea8.51
X\sm (X\sm A) &= A\cap X,\\
X\sm(A\cup B) &= (X\sm A)\cap (X\sm B), \lb{8.52}\\
X\sm(A\cap B) &= (X\sm A)\cup (X\sm B). \lb{8.53}\\
\intertext{\E\Ip if $A,B\sbs X$, then}
X\sm(X\sm A) &= A, \lb{8.54}\\
X\sm A = X\sm B &\hbox{ implies }A=B. \lb{8.55}
\e

\Wanp prove

\blm8.65
Let $(E,\cO)$ be a \tps, let $F$~be a nonempty \emph{closed} subset of~$E$, and let $A$~be a \nss\ of~$F$. Then $A$~is relatively closed in~$F$
$(F\sm A\in\cO_F)$ iff $A$~is closed in $(E,\cO)$ $(A^c = E\sm A\in\cO)$.
\elm

\proof \

``\ti{Only if\/}'': \ By \E\df\ \rf{d8.41} if $A$ is \rv ly closed in $(F,\cO_F)$ then $F\sm A\in\cO_A$, that is, \te s $O\in\cO$ \st $F\sm A=
F\cap O$. Note that $F = F\cap E = F\cap (O\cup O^c) = (F\cap O)\cup (F\cap O^c)$ and $(F\cap O)\cap(F\cap O^c)=\vn$. Hence $F\sm(F\cap O) =
F\cap O^c$. \E\oh $F\sm(F\sm A)=A$ by \er{8.54}. Since $F\cap O=F\sm A$, we obtain $A=F\cap O^c$ from \er{8.55}. \E\ev tly, $A^c = (F\cap O^c)^c =
F^c \cup O^{cc} = F^c \cup O$. Since $F$~is closed in~$E$, we infer $F^c$ open in~$E$ as well as~$O$, hence $A^c = F^c\cup O\in\cO$ by~$O_2$.

``\ti{If\/}'': \ If $A$ is closed in $(E,\cO)$, then $A^c\in\cO$. Since $F$~is closed in $(E,\cO)$, then $F^c\in\cO$. \E\Tf $(A\cap F)^c = A^c\cup
F^c \in \cO$. Hence $A\cap F = ((A\cap F)^c)^c = \cO$ by~$O_2$. Hence $A\cap F = ((A\cap F)^c)^c$ is closed in~$(E,\cO)$.
\endproof

\blm8.66
Let $A$ be a closed subset of $(\R,\cO(\R))$ and let $d$~denote the usual metric of~$\R$. Then $A$~is \sql\ closed in $(\R,d)$ $($see \E\df\
\rf{d8.31}$)$.
\elm

\proof
Suppose, for \cd ion, that this is not the case. Then \te s a \sq\ $\zb xn\N \sbs A$ and $x\notin A$ \st $x_n\to x$. Since $x\in A^c\in\cO(\R)$,
$x$~is an \ipt\ of~$A^c$, hence \te s $\ve\in\R_{>0}$ \st $(x-\ve,x+\ve)\sbs A^c$. \E\oh since $x_n\to x$, \te s $N\in\N$ \st if $n\ge N$ then
$|x_n-x|<\ve$. For such indices $n\in\N$ we have $x-\ve<x_n<x+\ve$, hence $x_n\in(x-\ve,x+\ve)\sbs A^c$, \cd ing $x_n\in A$.
\endproof

We now return to the second proof of the connectedness of $[a,b]$. Let $U_1,U_2$ be nonempty \rv ly open subsets of~$[a,b]$ \sf ying
\beq8.56
[a,b] = U_1\cup U_2 \hbox{ and }U_1\cap U_2 =\vn.
\e
\If that $U_1=[a,b]\sm U_2$, $U_2=[a,b]\sm U_1$, hence both $U_1$ and~$U_2$ are closed \wrt $\cO_{[a,b]}:=\cO(\R)_{[a,b]}$ (see \E\df\ \rf{d8.48}).
\Mo $[a,b]^c=(\lea,a)\cup(b,\to)$. Clearly, every $x\in(\lea,a)$ (resp.\ $(b,\to)$) is an \ipt, hence both \il s $(\lea,a)$ and $(b,\to)$ as well as
their union belong to $\cO(\R)$. Since $[a,b]=([a,b]^c)^c$, $[a,b]$ is closed in $(\R,\cO(\R))$.
In view of Lemma \rf{l8.65}, the sets $U_1$, $U_2$ and $[a,b]$ are closed in $(\R,\cO(\R))$, hence \sql\ closed in $(\R,d)$ by Lemma \rf{l8.66}.
Since $([a,b],d)$ is \sql\ closed, $[a,b]$ is \sql\ compact by \E\Pr\ \rf{p8.17}. We may apply Lemma \rf{l8.64a} with $K:=[a,b]$, $K_i:=U_i$,
$i=1,2$, and we find $U_1\cap U_2\ne\vn$. A~\cd ion. Hence $[a,b]$ is connected in $(\R,\cO(\R))$. \hfill $\Box$

\ssk
We now turn to the third proof of the connectedness of~$[a,b]$. To this end, we need the notion of (topological) \ti{closure}\index{topological closure} of a subset of a \tps.

Let $(X,\cO)$ be a \tps. Using De Morgan formulae (see \cite[p.~5]{Kelley}) we obtain the \fw\ \pp ies of closed subsets of $(X,\cO)$. Denoting
by~$\cC$ the collection of all closed subsets of $(X,\cO)$ we have:
\bit
\item[$\cC_1$] $\vn$ and $X$ belong to $\cC$;
\item[$\cC_2$] $\cC$ is closed under taking \isc s;
\item[$\cC_3$] $\cC$ is closed under taking \ti{finite} unions.
\eit

Given a subset $A$ of $X$, and noticing that the collection of all closed subsets of~$X$ containing~$A$ is nonempty since $X\supset A$,
and that the collection of all open subsets of~$X$ contained in~$A$ is nonempty since $\vn\sbs A$,
we may define~$\ov A$, the \ti{closure} ($\cO$-closure) of the subset~$A$ of $(X,\cO)$
and the \ti{interior} ($\cO$-interior) of~$A$ by setting
\beq8.57
\ov A := \bigcap_{B\in \cC,B\supset A} \qh{and} \kl A= \bigcup_{B\in\cO,B\sbs A}B.
\e

\bpr8.67
Let $(X,\cO)$ be a \tps, and let $A,B$ be subsets of $X$. Then\dw
$$
{\rm (i)}\ \ov\vn=\vn; \q\q {\rm (ii)}\ A\sbs\ov A; \q\q {\rm (iii)}\ \ov{\ov A}=\ov A; \q\q{\rm (iv)}\ \ov{A\cup B}=\ov A\cup \ov B.
$$
\Mo $A\in\cC(\cO)$ iff $A=\ov A$ and

\hph v,ii, $A\sbs B \hbox{ implies } \ov A\sbs \ov B$.

\hph vi,i, ${\ov A}^c = (A^c)^\circ$.

\hph vii,, $\kl A$ defined in \er{8.57} coincides with $\kl A$ defined in \er{8.44}.
\epr

\bex8.68 \

\hph i,ii, Prove \E\Pr\ \rf{p8.67}.

\hph ii,i, Suppose that $X$ is a \ns\ and that $F:\cP(X)\to\cP(X)$ \sf ies \pp ies (i)--(iv). \Mo define $\cC_F$ by setting:
\beq8.58
\hbox{\fe}A\sbs X, \q A\in\cC_F \hbox{ if }A=F(A).
\e
Show that the collection $\cC_F$ \sf ies $\cC_1$--$\cC_3$. Setting $\cO_F:= \{A^c\sbs X: A\in\cC_F\}$, show that $(X,\cO_F)$ is a \tps. Show that
$\cO_F=\cO$ if $(X,\cO)$ is a \tps, and $F(A) = \ov A$ defined in \er{8.57}.

\hph iii,, Show that every singleton $\{x\}$ in a \Hts\ is closed, that is, $X\sms x$ is open.
\eex

\Wanp give the promised third proof of the connectedness of~$[a,b]$. This proof is based on the connectedness of the \il\ $(a,b)$ and of~$\R$.
Observe that $[a,b]$, $a,b\in \R$, $a<b$, is the closure of the \il\ $(a,b)$ \wrt the topology $\cO(\R)$, i.e.\ $[a,b]=\ov{(a,b)}$. Indeed, we already
know that $[a,b]$ is the \cpl\ in~$\R$ of an open set of $(\R,\cO)$. Hence $[a,b]\in\cC$. \Mo $(a,b)\notin\cC$. Otherwise, we would have $(a,b)^c
\in\cO$, and $\R=(a,b)\cup(a,b)^c$, $(a,b)\cap(a,b)^c=\vn$, \cd ing the connectedness of~$\R$. Since $\zo a,b = \{a\}\cup(a,b)$ and $\oz a,b =
(a,b)\cup\{b\}$, the only subsets of~$\R$ strictly containing $(a,b)$ and contained in~$[a,b]$ are $\zo a,b $ and $\oz a,b $. None of them are
closed. Indeed, $\zo a,b ^c = (\lea,a)\cup [b,\to)$. Observe that $b$~is not an \ipt\ of $\zo a,b ^c$ since there is no \il\ $(\a,\b)$, $\a,\b\in\R$,
$\a<\b$, contained in $\zo a,b ^c$ which contains $\{b\}$. Thus $\zo a,b \notin\cC$. The proof that $\oz a,b \notin\cC$ is
similar. \E\Tf the \il\ $[a,b]$ is the smallest (\wrt inclusion) closed \il\ containing the \il\ $(a,b)$. Hence $[a,b]=\ov{(a,b)}$. We now apply the
\fw

\begin{lem}[see e.g.\ {\cite[p.~54]{Kelley}}]\lb{l8.69}
Let $(X,\cO)$ be a \tps. Then the closure of a connected subset~$Y$ of~$X$ is connected.
\elm

\proof
Suppose, for \cd  ion, that \te\ nonempty \rv ly open subsets $A,B$ of~$\ov Y$ \sf ying $\ov Y= A\cup B$ and $A\cap B=\vn$. (In notation $A,B \in
\cO_{\bar Y}$.) Then $A= \ov Y\sm B$, $B=\ov Y\sm A$ are also closed in $(\ov Y,\cO_{\bar Y})$. Note that $Y=\ov Y\cap Y$ since $Y\sbs \ov Y$. \E\Tf
$A\cap Y$ and $B\cap Y$ belong to $\cO_Y:=(\cO_{\bar Y})_Y$. \Mo from $\ov Y = A\cup B$, we infer $Y=Y\cap\ov Y = Y\cap(A\cup B) = (A\cap Y)\cup
(B\cap Y)$, and from $A\cap B=\vn$, we infer $(A\cap Y)\cap(B\cap Y) = (A\cap B)\cap Y = \vn\cap Y=\vn$. Hence $A\cap Y = Y\sm(B\cap Y)$ and
$B\cap Y = Y\sm(A\cap Y)$. Hence both $A\cap Y$ and $B\cap Y$ are open in $(Y,\cO_Y)$, and \sf y $(A\cap Y)\cap(B\cap Y)=\vn$. Since $(Y,\cO_Y)$ is
connected, one of them must be empty. Interchanging $A$~and~$B$ if necessary, we may assume $B\cap Y=\vn$. \E\Tf $Y=A\cap Y$, hence $Y\sbs A$. By
\E\Pr\ \rf{p8.67}\,(ii), we have $\ov Y\sbs\ov A$ in~$(X,\cO)$. Since $\ov Y$~is closed in $(X,\cO)$ and $A$~is closed in $(\ov Y,\cO_{\bar Y})$,
we infer from Lemma \rf{l8.65} that $A$~is closed in~$(X,\cO)$. Thus $\ov A=A$ by \E\Pr\ \rf{p8.67}. \E\Tf $\ov Y\sbs A$, hence $\ov Y=A$ and
$B=\vn$, since $\ov Y=A\cup B$ and $A\cap B=\vn$. A~\cd ion since $B\ne\vn$ by \as.
\endproof

\bex8.70
Let $(X,\cO)$ be a \tps\ and let $Y,Z$ be \nss s of~$X$ \sf ying $Y\sbs Z\sbs \ov Y$. Show that if $X$~is connected then so is~$Z$. (Hint: Apply
Lemma \rf{l8.69} to~$Z$ with the \rv\ topology.)
\eex

We summarize some of the previous results concerning \Ar\ fields and connectedness in the next theorem.

\bth8.71
Let $(K,+,\cdot,0,1,\ge)$ be an \Ar\ \of\ endowed with the order topology $\cO(K)$. Then the \fw\ assertions are \ev t\dw

\hph i,ii, $(K,\ge)$ is order-complete\sd

\hph ii,i, $K$ is connected\/\sd

\hph iii,, $(0,1)$ and $[0,1]$ are connected\/\sd

\hph iv,, all \il s of $(K,\ge)$ are connected.

\noi \Mo if $A$ is a subset of $K$, $A\ne K$, then $A$~is connected iff $A$~is an \il.
\eth

\brm8.72
Observe that if we delete an \el\ $x$ of~$\R$ from~$\R$, then $\R\sms x$ is not connected. Similarly, if we delete $c\in(a,b)$ from $[a,b]$, then
$[a,b]\sms c$ is not connected. However, if we delete $a$~or~$b$ from $A:=[a,b]$ (resp.\ $\zo a,b $, $\oz a,b $) then $A\sms a$ (resp.\ $A\sms b$)
is connected. Note that $\{a\},\{b\}$ are subsets of $[a,b]\sm(a,b)$.
\erm

\bdf8.73
Let $(X,\cO)$ be a \Hts\ and let $A$~be a \nss\ of~$X$. If $\kl A$ denotes the interior of~$A$, i.e.\ $\{x\in A: \hbox{\te s }O\in\cO \hbox{ \sf
ying}\break{x\in O\sbs A}\}$, then $\ov A\sm\kl A$ is called the (\ti{topological\/}) \ti{boundary} of~$A$ in~$(X,\cO)$ and is some\-times denoted
by~$b(A)$ or $\partial A$.
\edf

\bex8.74 \

\hph i,i, Let $(X,\cO)$ be a \Hts\ and let $A$~be a \nss\ of~$X$. Prove the \fw\ assertions:
\beag
{}&b(A)= \ov A\cap \ov{(X\sm A)}, \q X\sm b(A) = {\kl A}\, \cup \nad{\lower3pt\hbox{$\hskip8pt\scriptstyle\circ$}}{(X\sm A)};\\
&\ov A = A\cup b(A), \q \kl A= A\sm b(A);\\
&A \hbox{ is closed iff }b(A)\sbs A;\\
&A \hbox{ is open iff }A\cap b(A)=\vn.
\e

\hph ii,, If $A$ is connected in $(\R,\cO(\R))$, is $A\sm b(A)$ connected?
\eex

\brm8.75
Let $a,b,c\in\R$ be \st $a<b<c$. Then $(a,b)$ and $(b,c)$ are two nonempty disjoint connected subsets of $(\R,\cO(\R))$. Clearly, the set $(a,b)
\cup(b,c)$ is not connected in $(\R,\cO(\R))$. However, the subset of~$\R$ $(a,b)\cup\zo b,c $, which is the disjoint union of two connected sets,
\ti{is} connected. Indeed, $(a,b)\cup\zo b,c = \oz a,b \cup(b,c) = (a,c)$ is connected. Note that $\oz a,b \cap \zo b,c =\{b\}\ne\vn$.
\erm

\newpage
More generally, we have

\blm8.76
Let $(X,\cO)$ be a \tps\ and let $A,B$ be nonempty disjoint connected subsets of~$X$ \sf ying $\ov A\cap B\ne\vn$ or $A\cap \ov B\ne\vn$. Then
$A\cup B$ is connected.
\elm

\bex8.77
Prove Lemma \rf{l8.76} (see e.g.\ \cite[pp.\ 52--54]{Kelley}).
\eex

\bex8.78
Let $(X,\cO)$ be a \Hts\ and let $x\in X$. Prove the \fw\ assertions:

\hph i,ii, The union of all connected subsets of $(X,\cO)$ containing $\{x\}$ is a closed connected subset of $(X,\cO)$, called the (connected)
component of~$x$, and sometimes denoted by~$C(x)$.

\hph ii,i, If $(X,\cO):=(\R,\cO(\R))$ then $C(x)=\R$ \fa $x\in\R$.

\hph iii,, If $(X,\cO):=(\R\sms0,\cO_{\R\sms0})$ then $C(x)=\R_{<0}$ (resp.\ $\R_{>0}$) whenever $x\in\R_{<0}$ (resp.\ $x\in\R_{>0}$).

\hph iv,, If $(X,\cO):=(\Q,\cO(\Q))$ then $C(x)=\{x\}$ \fa $x\in\Q$.
\eex

So far we considered on the field $\R$ the order topology. In order to define \cvg\ \sq s in~$\R$, we \itd the notion of metric, in particular the
usual metric. It turns out that a metric~$d$ on a \ns~$X$ induces a natural topology, which we denote by~$\cO(d)$. \E\Ip if $(X,d)$ and $(X',d')$
are \ms s and $f$~is a map from $X$ into~$X'$, the notion of sequential continuity of~$f$ defined in \E\df\ \rf{d8.28} can be formulated in terms of
topological notions such as open sets, closed sets, and so on. We will also give an example where the pointwise \cg nce of a \sq\ of \f s (maps)
cannot be described by a metric, but by a~topology.

The key object which makes the link between a metric and a topology on a set is the notion of ball defined as follows.

\bdf8.79
Let $(X,d)$ be a \ms. The set
\beq8.59
B(x;r) := \{y\in X: d(x,y)<r\}
\e
is called a \tb{ball of center $x$ and radius} $r\in\R_{>0}$.\index{ball of center $x$ and radius~$r$}
\edf

Note that if the \ms\ is $(\R,d)$ with its usual metric, then we have: $B(x;r)= (x-r,x+r)$. Conversely, if $a,b\in\R$, $a<b$, then $(a,b) =B(x;r)$
with $x:=\frac12(a+b)$ and $r:=\frac12(b-a)$.

\blm8.80
Let $(X,d)$ be a \ms\ and let
$$
\cB(d) := \{B(x;r)\sbs X: x\in X \hbox{ and } r\in\R_{>0}\}.
$$
Then $\cB(d)$ \sf ies \er{8.41}--\er{8.43} with $\cB :=\cB(d)$.
\elm

\proof \

``\er{8.41}'': $x\in B(x,r)$ \fa $x\in X$ and all $r\in\R_{>0}$.

``\er{8.42}'': Let $x',x''\in X$, $r',r''\in\R_{>0}$. We claim that if $B(x',r')\cap B(x'',r'')\ne\vn$, then \te\ $x\in X$ and $r\in\R_{>0}$ \st the
\fw\ holds:
\beq8.60
B(x;r) \sbs B(x';r') \cap B(x'';r'').
\e
Indeed, since $B(x';r') \cap B(x'';r'')\ne\vn$, \te s $x\in B(x';r') \cap B(x'';r'')$. \E\Ip we have $d(x,x')<r'$, hence $r'-d(x,x')>0$. Similarly,
$r''-d(x,x'')>0$. Set $r:=\min(r'-d(x,x'),r''-d(x,x''))$. Then $r>0$. If $y\in B(x,r)$ then $d(x',y) \le d(x',x) + d(x,y) < d(x',x)+r < d(x',x)
+r'-d(x,x') = r'$. Hence $y\in B(x';r')$. Similarly, $d(x'',y) < d(x'',y)+r''-d(x'',y) = r''$, hence $y\in B(x'';r'')$.

``\er{8.43}'': Suppose that \te\ $z\in X$, $r_1,r_2\in\R_{>0}$ \st $z\in B(x_1;r_1)\cap B(x_2;r_2)$. Then $d(x_1,x_2)\nad{\rm(M2)}\le d(x_1,z) +
d(z,x_2)\nad{\rm(M2),\er{8.59}}\le r_1+r_2$. Thus if $r_1>0$, $r_2>0$ \sf y $r_1+r_2 = \frac12d(x_1,x_2)$, there is no $z\in B(x_1;r_1) \cap
B(x_2;r_2)$, hence $B(x_1;r_1) \cap B(x_2;r_2)=\vn$.
\endproof

\bpr8.81
Let $(X,d)$ be a \ms\ and let $\cO(d)$ denote the collection of subsets of~$X$ consisting of the empty subset of~$X$ and of the \nss s $A$ of~$X$
\sf ying
\beq8.61
\hbox{\fe $x\in A$ \te s $r\in\R_{>0}$ \st $B(x,r)\sbs A$.}
\e
Then $(X,\cO(d))$ is a \Hts\ and $B(x,r)\in \cO(d)$ \fa $x\in X$ and all $r\in\R_{>0}$.
\epr

\proof
\E\Pr\ \rf{p8.81} is a con\sq\ of Lemma \rf{l8.44} and Lemma \rf{l8.80}.
\endproof

\bex8.82 \

\hph i,i, Show that if $X:=\R$ then the order topology is the same as $\cO(d)$ where $d$~is the usual metric on~$\R$.

\hph ii,, Show that on $X:=\R$ $\cO(d) = \cO(d')$, where $d'(x,y)=\a d(x,y)$, $\a,x,y\in\R$, $\a>0$.
\eex

\bdf8.83
Let $(X,\cO)$ be a \Hts, and let $\zb xn\N$ be a \sq\ in~$X$. The \sq\ $\zb xn\N$ is said to be \ti{\cvg\/} to a \ti{limit\/} $x\in X$, written
$x_n\nad{\cO}\to x$, if \fe open set $O\in\cO$ with $x\in O$ \te s $N\in\N$ \st $x_n\in O$ \fa $n\ge N$.
\edf

\bex8.84 \

\hph i,i, Show that a \sq\ $\zb xn\N$ in a \Hts\ possesses at most one limit.

\hph ii,, Let $(X,d)$ be a \ms\ and let $\cO(d)$ be the induced topology. Show that $x_n \nad{\cO}\to x$ \fs $x\in A$ iff $x_n \nad{d}\to x$.
\eex

\bdf8.85
A \Hts\ $(X,\cO)$ is called \ti{metrizable} if \te s a metric~$d$ on~$X$ \st $\cO=\cO(d)$.\index{metrizable}
\edf

\bpr8.86
Let $(X,\cO)$ be a Hausdorff \tps\ and let $A$ be a \nss\ of~$X$. Then

\hph i,i, if $A$ is closed then $A$ is \sql\ closed\sd

\hph ii,, if $(X,\cO)$ is metrizable and $A$ is \sql\ closed then $A$ is closed.
\epr

\proof \

(i) The proof is similar to the proof of Lemma \rf{l8.66}. It suffices to replace $\cO(\R)$ by~$\cO$, $B(x;\ve)$ by $O\in\cO$ \st $x\in O$ and
$|x-x_n|<\ve$ by $x_n\in O$.

(ii) To prove $A^c\in\cO$, that is, \fe $x\in A^c$ \te s $\ve\in\R_{>0}$ \st $B(x;\ve)\sbs A^c$, where $B(x;\ve)$ is defined in \er{8.59}, and $\cO
(d)=\cO$. Let $x\in A^c$ and suppose that there is no $\ve\in\R_{>0}$ \st $B(x;\ve)\sbs A^c$. Then \fe $n\in\Na$ \te s $x_n\in B(x;\frac1n)$ \st
$x_n\in(A^c)^c=A$. \If that $x_n\nad d\to x$, hence $x\in A$ since $A$~is \sql\ closed, \cd ing $x\in A^c$.
\endproof

\bdf8.87
Let $(X,\cO)$ be a \Hts\ and let $A\sbs X$, $A\ne\vn$, and $A\ne X$. Then the set $\ov A{}\en s \sbs X$ defined by
\beq8.62
\ov A{}\en s := \{x\in X: \hbox{\te s a \sq\ $\zb xn\N$ in $A$ \st}x_n\nad{\cO}\to x\}
\e
is called the \tb{sequential closure} of~$A$.
\edf

\bpr8.88
Let $(X,\cO)$ be a Hausdorff \tps. Then \fe non\-empty subset $A$ of~$X$ we have $\ov A{}\en s \sbs \ov A$, and if $(X,\cO)$ is metrizable,
then $\ov A{}\en s = \ov A$.
\epr

\proof \

``$\ov A{}\en s \sbs \ov A$'': \
Let $x\in \ov A{}\en s $ and let $\zb xn\N$ in~$A$ be \st $x_n\nad{\cO}\to x$. Suppose, for \cd ion, that $x\notin\ov A$. Since $\ov{\ov A}=\ov A$ by
\E\Pr\ \rf{p8.67}\,(ii), we infer that $\ov A$ is closed by the same \Pr. Hence $(\ov A)^c\in\cO$. By \df\ \te s $O\in\cO$ \st
$x\in O\sbs A^c$.
Since $x_n\nad{\cO}\to x$, \te s $N\in\N$ \st $x_n\in O$ for $n\ge N$ in view of \E\df\ \rf{d8.83}. Hence $x_N\in O\sbs A^c$, a~\cd ion, since
$x_N\in A$.

Let $d$ be a metric \st $\cO(d)=\cO$.

``$\ov A\sbs \ov A{}\en s$'': \ We first show that \fe $x\in \ov A$ we have $A\cap B(x;\ve)\ne\vn$ \fe $\ve\in\R_{>0}$. Indeed,
let $x\in\ov A$ be arbitrary, and suppose, for \cd ion,
that \te s $\ve\in\R_{>0}$ \st $A\cap B(x;\ve)=\vn$, then $A\sbs B(x;\ve)^c$. Since $B(x;\ve)^c$ is closed, we have $B(x;\ve)^c \nad*= \ov{B(x;\ve)
^c}$, hence $x\in \ov A \nad*\sbs B(x;\ve)^c$. In~$\nad*=$ we used \E\Pr\ \rf{p8.67}. \If that $x\in B(x;\ve)^c$, a~\cd ion, since $x\in B(x;\ve)$.
We next show that $x\in\ov A{}\en s$. \E\fe $n\in\Na$ \te s $x_n\in B(x;\frac1n)$. \E\Tf $x_n\nad d\to x$, and $x\in\ov A{}\en s$.
\endproof

\Wanp reformulate the \df\ of the sequential continuity of a map $f:X_1\to X_2$ where $(X_1,d_1)$, $(X_2,d_2)$ are \ms s. We recall that $f$~is
s-\ctn\ at $x\in X_1$ if \fe \sq\ $\zb xn\N$ \st $x_n\nad{d_1}\to x$, we have $f(x_n)\nad{d_2}\to f(x)$. If $(X_1,\cO_1)$ and $(X_2,\cO_2)$ are
\Hts s and $f:X_1\to X_2$, the map is also called \ti{\sql\ \ctn} at $x\in X_1$ if $x_n\nad{\cO_1}\to x$ implies $f(x_n)\nad{\cO_2}\to f(x)$ \fe
\sq\ $\zb xn\N$ and $x\in X$. The map~$f$ is called \ti{s-\ctn} if it is s-\ctn\ at~$x$ \fa $x\in X_1$.

\bth8.89
Let $(X_1,\cO_1)$, $(X_2,\cO_2)$ be \tps s, and let $f:X_1\to X_2$. If $(X_1,\cO_1)$ is metrizable then the \fw\ assertions are \ev t\dw

\hph i,ii, $f$ is \sql\ \ctn\sd

\hph ii,i, $f(\ov A) \sbs \ov{f(A)}$ \fe nonempty subset $A$ of $X_1$\sd

\hph iii,, $f\mo(B) \in \cC_1$ \fe set $B\in\cC_2\sms\vn$\sd

\hph iv,, $f\mo (B) \in\cO_1$ \fe set $B\in\cO_2\sms\vn$.

\noi Here $\cC_i$ denotes the collection of closed subsets of $(X_i,\cO_i)$, $i=1,2$.
\eth

{\allowdisplaybreaks For the sake of clarity we will use the notation $f\mo[A]$ instead of~$f\mo(A)$ with\break ${f\mo[\vn]:=\vn}$.

We leave it to the reader to verify the \fw\ \ass s (see \cite{Kelley}). Let $f:X_1\to X_2$ and let $A,A'$ (resp.\ $B,B'$) be \nss s of~$X_1$
(resp.~$X_2$). Then:
\bea8.63
{}&f\mo[B\sm B'] = f\mo[B] \sm f\mo[B'];\\
&f\mo[B\cup B'] = f\mo[B]\cup f\mo[B']; \q f\mo[B]\sbs f\mo[B'] \hbox{ if }B\sbs B';\lb{8.64}\\
&f\mo[B\cap B'] = f\mo[B]\cap f\mo[B']; \lb{8.65}\\
&A \sbs f\mo[f(A)]; \lb{8.66}\\
&f(f\mo[B]) \sbs B; \lb{8.67}\\
&A=f\mo[f(A)] \hbox{ if $f$ is in\jc, and $B= f(f\mo[B])$ if $f$ is sur\jc.} \lb{8.68}
\e

Note that $B\sbs B'$ implies $f\mo[B] \sbs f\mo[B']$. Indeed, $B\sbs B'$ iff $B\cup B'=B'$, hence $f\mo[B']=f\mo[B\cup B'] \nde8.64 = f\mo[B]\cup
f\mo[B']\supset f\mo[B]$. \E\Tf $f\mo[B]\sbs f\mo[B']$.}

\proof[Proof of Theorem \rf{t8.89}] \

``$\rm(i)\imp(ii)$'': \ Let $x\in\ov A$. By \E\df\ \rf{d8.87} and \E\Pr\ \rf{p8.88}, \te s a \sq\ $\zb xn\N$ in~$A$ \st $x_n\nad d\to x$ with
$\cO(d)=\cO_1$. By~(i) $f(x_n)\nad{\cO_2}\to f(x)$, hence $f(x)\in \ov{f(A)}$ in view of \E\Pr\ \rf{p8.86}\,(i).

``$\rm (ii)\imp(iii)$'': \ We claim that it suffices to prove
\beq8.69
\ov{f\mo[B]} \sbs f\mo[\ov B] \hbox{ \fa nonempty closed subsets of }X_2.
\e
Indeed, if $B\in \cC_2$ then $B=\ov B$ by \E\Pr\ \rf{p8.67}. Hence $\ov{f\mo[B]} \sbs f\mo[\ov B]$ by~\er{8.69}, and $f\mo[\ov B] = f\mo[B]$. \csq,
$\ov{f\mo[B]} \sbs f\mo[B] \sbs \ov{f\mo[B]}$ by \E\Pr\ \rf{p8.67}\,(ii). \E\Tf $f\mo[B] = \ov{f\mo[B]}$, and $f\mo[B]\in\cC_1$ by the same \Pr.
It remains to prove \er{8.69}. Let $B\sbs X_2$, $B\ne\vn$, and set $A:=f\mo[B]\sbs X_1$. We have $f(A) \sbs B$ by~\er{8.67}, hence
$\ov {f(A)}\sbs \ov B$ by \E\Pr\ \rf{p8.67}\,(v). From~(ii) we have $f(\ov A)\sbs \ov{f(A)}$, hence $f(\ov A)\sbs\ov B$. By~\er{8.66} with
$A:=\ov A$ we obtain $\ov A\sbs f\mo[f(\ov A)]$, hence $\ov{f\mo[B]} = \ov A
\sbs f\mo[f(\ov A)] \sbs f\mo[\ov B]$ by \er{8.64}, and $\ov{f\mo[B]} \sbs f\mo[\ov B]$, which is~\er{8.69}.

``$\rm(iii)\imp(iv)$'': \ Let $B\in\cO_2$. Then $X_2\sm B\in\cC_2$ by \E\df\ \rf{d8.41}. Then ${f\mo[X_2\sm B]\in\cC_1}$ by~(iii), and
$X_1\sm f\mo[B] \nde 8.63 = f\mo[X_2\sm B]\in\cC_1$ since $f\mo[X_2]=X_1$. \E\Tf $f\mo[B]\nde8.54 = X_1\sm(X_1\sm f\mo[B])\in\cO_1$.

``$\rm(iv)\imp(i)$'': \ Let $x \in X_1$ be arbitrary, and let $\zb xn\N$ be an arbitrary \sq\ in~$X_1$ \sf ying $x_n\nad{\cO_1}\to x$. We have
to show that $f(x_n)\nad{\cO_2}\to f(x)$ assuming~(iv).
In view of \E\df\ \rf{d8.83} we have to show that given $U\in\cO_2$ \st $f(x)\in U$, \te s $N\in\N$ \st $f(x_n)\in U$ \fa $n\ge N$. Set
$O:=f\mo[U]$, $U\in\cO_2$, $f(x)\in U$. Then $f(O)\sbs U$ by~\er{8.67}. Since $f(x)\in U$, we have $x\in O$. \Mo in view of $x_n\nad{\cO_1}\to x$,
\te s $N\in\N$ \st $x_n\in O$ \fa $n\ge N$. For these indices~$n$ we have $x_n\in O$, hence $f(x_n)\in f(O)\sbs U$. \E\Tf $f(x_n)\nad{\cO_2}\to f(x)$.

The proof of Theorem \rf{t8.89} is complete.
\endproof

\brm8.91
Inspection of the proof of Theorem \rf{t8.89} shows that the metrizability of $(X_1,\cO_1)$ is only required for the proof of ``$\rm(i)\imp(ii)$''. It
turns out that \te\ stronger \cn s which imply~(ii) without requiring the metrizability of $(X_1,\cO_1)$. For the convenience of the reader we state
one of them, and for the proof of the claim we refer the reader to~\cite{Du}. We need the \fw\ \df. A~\ti{directed set\/}\index{directed set} is a (partially) \os\
$(I,\ge)$ \st \fa $i,j\in I$ \te s $k\in I$ with $k\ge i$ and $k\ge j$. Note that a join-\slt\ (see \E\df\ 3.3.8) is a directed set, but not
conversely. Examples of \dis s are $(\N,\ge)$, and the join-\slt\ consisting of all finite subsets of~$\N$, (partially) ordered by
inclusion~$\supset$. Indeed, if $A,B$ are finite subsets of~$\N$, then $A\cup B$ is finite, and $A\cup B\supset A$, $A\cup B\supset B$. A~\ti{net\/}
in a \ns~$X$ is a map from a \dis~$I$ into~$X$. A~\sq\ in a set~$X$ is sometimes denoted by $\zb xi\N$. A~net in~$X$ is denoted by $\zb xiI$
in~\cite{Du}. An example of net in $(\R,\cO(\R))$, which is not a \sq, is defined as follows. Let~$I$ denote the second example of \dis\ given above,
and let $\zb yn\N$ be a \sq\ in~$\R$. Then for $A\in I$, set $x_A:=\suml_{j\in A}y_j$ if $A\ne\vn$, and $x_A:=0$ if~$A=\vn$, according to~(2.4.6).
\erm

A net $\zb xiI$ in a \Hts\ $(X,\cO)$ \ti{converges} to $a\in X$, written $x_i\nad{\cO}\to a$, if \fe open set $O\in\cO$ containing~$a$, \te s $l\in I$
\st $x_k\in I$ \fa $k\ge l$. \Mo if $x_i\nad{\cO}\to a$ and $x_i\nad{\cO}\to b$ with $b\in X$, then $a=b$, and $a$~is called the \ti{limit of the net
$\zb xiI$}.

\If from Theorem 2.1.2 in \cite{Du} and \E\df\ on page~28 of~\cite{Du} that Theorem \rf{t8.89} holds without the requirement that $(X_1,\cO_1)$ be
metrizable, provided (i)~is replaced by
\bml8.70
\hbox{\E\fe $a\in X_1$ and every net $\zb xiI$ in $(X,\cO)$ converging to $a$,}\\
\hbox{the net $\{f(x_i)\}_{i\in I}$ in $X_2$ converges to $f(a)$.}
\e
The raison for which (i) does not imply (ii) is that the closure $\ov A$ of~$A$ in $(X_1,\cO_1)$ is not always equal to $\ov A{}\en s$
(see~\cite{Du}). However, the \fw\ holds:
\beq8.71
\hbox{\E\fe $a\in \ov A$ \te s a net $\zb xiI$ converging to $a$.}
\e
(See Theorem 2.1.3 in \cite{Du}.) It should be mentioned that the proof of \er{8.66} makes use of the axiom of choice \cite{Du}.

\goodbreak
\bex8.92 \

\hph i,i, Show that the net \itd in Remark \rf{r8.91} with $y_n\in\R_{\ge0}$, $n\in\N$, is \cvg\ iff \te s $M\in\R_{\ge0}$ \st
\beq8.72
\sum_{n=0}^m y_n \le M \hbox{ \fa}m\in\N.
\e
Show that if \cn\ \er{8.72} is \sf ied then $\zb xiI \nad{\cO(\R)}\to \supl_{n\ge0}y_n$.

\hph ii,, Observe that the proofs of implications $\rm(ii)\imp(iii)$ and $\rm(iii)\imp(iv)$ do not require the metrizability of $(X_1,\cO_1)$. Show
that (iv) implies (iii) without the metrizability of $(X_1,\cO_1)$. Does the implication $\rm(iii)\imp(ii)$ need the metrizability of $(X_1,\cO_1)$?
\eex

Theorem \rf{t8.89} and Remark \rf{r8.91} motivate the \fw\ \df.

\bdf8.93
Let $(X_i,\cO_i)$, $i=1,2$, be \Hts s and $f:X_1\to X_2$. Then the map $f$ is called \ti{\ctn},\index{continuous map} if \fa \ns s $U\in\cO(X_2)$, $f\mo[U]\in\cO_1$.
\edf

In the next \Pr\ we state simple but useful observations concerning \ctn\ maps.

\bpr8.94 \

\hph i,ii, Let $(X_i,\cO_i)$, $i=1,2,3$, be \tps s, and let $f:X_1\to X_2$, $g:X_2\to X_3$. If $f$~and~$g$ are \ctn, so is $g\circ f:X_1\to X_3$.

\hph ii,i, Let $(X,\cO)$ be a \tps\ and let $A$ be a nonempty subspace of~$X$. Then the inclusion map $j_A:A\to X$ defined by $j_A(x):=x$, $x\in A$,
is \ctn\ from $(A,\cO_A)$ into $(X,\cO)$.

\hph iii,, If $f:(X_1,\cO_1)\to(X_2,\cO_2)$ is \ctn, and $A$~is a \nss\ of~$X_1$, then $f|_A$, the \rt ion of~$f$ to~$A$, is \ctn\ from $(A,\cO_A)$
into $(X_2,\cO_2)$.

\hph iv,, If for $i=1,2$ $(X_i,d_i)$ are \ms s with balls $B_i(x;r)$, and ${f:X_1\to X_2}$, then $f$~is \ctn\ iff \fe
$x\in X_1$ and every $\ve\in\R_{>0}$ \te s a $\d\in\R_{>0}$ \st $d_2(f(x),f(x'))<\ve$ \fa $x'\in X_1$ with $d_1(x,x')<\d$.
\E\Ip if $f$~is \Lc\ $($resp.\ \Hc\ of order $\a\in(0,1))$, then $f$~is \ctn.

\hph v,i, Let $f:(X,\cO)\to(\R,\ge)$. Then $f$~is \ctn\ from $(X,\cO)\to(\R,\cO(\R))$ iff \fe $a\in\R$ $\{x\in X:f(x)<a\}\in\cO$ and $\{x\in X:
f(x)>a\}\in\cO$.
\epr

\proof \

(i) \ We have $\xymatrix{X_1 \ar[r]^f \ar@/_/[rr]_{g\circ f}&X_2 \ar[r]^g &X_3}$. Let $O\in\cO_3$. Since $g$~is \ctn, $g\mo[O]\in\cO_2$. Since $f$~is
\ctn, $f\mo[g\mo[O]]\in\cO_1$. Note that \fe $A\sbs X_3$ $f\mo[g\mo[A]] = (g\circ f)\mo[A]$. Indeed, $(g\circ f)\mo[A] = \{x\in X_1: (g\circ f)(x)
\in A\} = \{x\in X_1: {g(f(x))\in A}\} = \{x\in X_1: f(x)\in g\mo[A]\} = \{x\in X_1: x\in f\mo[g\mo[A]]\}$. Thus $(g\circ f)\mo[O]= f\mo[g\mo[O]] \in
\cO_1$.

(ii) \ To prove: $j_A\mo(O)\in\cO_A$ \fe $O\in\cO$. By \df, $j_A\mo(O):= \{x\in A: j_A(x)\in O\}$. Hence $j_A\mo(O) = \{x\in A: x\in O\} = O\cap A\in
\cO_A$.

(iii) \ We have $f|_A(x)=f(x)$ \fa $x\in A$. Hence $f|_A(x) = f(j_A(x))$ \fa $x\in A$, hence $f|_A = f\circ j_A$. Since both $j_A$ and~$f$ are \ctn,
so is $f|_A$ in view of~(i).

(iv) ``\ti{Only if\/}'': \ Suppose that $f$ is \ctn. Then \fe $x\in X_2$ and every $\ve\in\R_{>0}$, $B_2(f(x);\ve)\in \cO(d_2)$ by \E\Pr\  \rf{p8.81}.
Since $f:(X_1,\cO(d_1))\to (X_2,\cO(d_2))$, we have $x\in f\mo[B_2(f(x);\ve)]$, and since $f\mo[B_2(f(x);\ve)]\in\cO(d_1)$, \te s $B_1(x;\d)$ with
$\d\in \R_{>0}$ \st $B_1(x;\d)\sbs f\mo[B_2(f(x);\ve)]$. Hence $f(B_1(x;\d)) \sbs f(f\mo[B_2(f(x);\ve)])\nde8.67 \sbs B_2(f(x);\ve)$. \csq, $d_2(f(x),
f(x'))<\ve$, $x'\in X_1$, whenever $d_1(x,x')<\d$.

``\ti{If\/}'': \ We first show that $f$ is \sql\ \ctn. Let $x\in X_1$ be arbitrary, and let $\zb xn\N$ in~$X_1$ be \st $x_n\nad{d_1}\to x$. Let
$\ve\in\R_{>0}$. To prove: \te s $N\in\N$ \st $d_2(f(x),f(x_n))<\ve$ provided $n\ge N$. By \as\ \te s a $\d\in\R_{>0}$ \st $d_1(x,x_n)<\d$ implies
$d_2(f(x),f(x_n))<\ve$. Since $x_n\nad{d_1}\to x$, \te s $N\in\N$ \st $d(x,x_n)<\d$ for $n\ge N$, hence also $d(f(x),f(x_n))<\ve$. Thus $f(x_n)
\nad{d_2}\to f(x)$. Since $f$~is \sql\ \ctn, $f$~is \ctn\ by Theorem \rf{t8.89} and \E\df\ \rf{d8.93}.

Let $f$ \sf y $d_2(f(x),f(y)) \le M(d_1(x,y))^\a$ \fs $M\in\R_{>0}$ and $\a\in\oz 0,1 $ \fa $x,y\in X_1$. Then $d_2(f(x),f(y))<\ve$ provided
$M(d_1(x,y))^\a <\ve$, that is, $d(x,y) < (M\mo \ve)^{\frac1\a} =: \d$.

(v) ``\ti{Only if\/}'': \ Let $a\in\R$ be arbitrary. Note that $\{x\in X: f(x)<\a\} = f\mo[(\lea,a)]$, and that $(\lea,a)\in\cO(\R)$, since \fe
$y\in\R$, $y<a$, we have $(y-h,y+h) \sbs (\lea,a)$ whenever $0<h<a-y$. Similarly, $(a,\to)\in \cO(\R)$. Then $f\mo[(\lea,a)]\in\cO(\R)$ and
$f\mo[(a,\to)]\in\cO(\R)$ by the continuity of~$f$.

``\ti{If\/}'': \ Let $a,b\in\R$, $a<b$. Then $(a,b)=(\lea,b)\cap(a,\to)\in\cO(\R)$ by~$(O_3)$ (see \E\df~\rf{d8.41}). Let $O\in \cO(\R)$. Then \te\
$a_i,b_i\in\R$, $a_i<b_i$ and an index set~$I$ \st $O=\bcl_{i\in I}(a_i,b_i)$. Then $f\mo[0] = f\mo\bigl[\bcl_{i\in I}(a_i,b_i)\bigr] = \bcl_{i\in I}
f\mo[(a_i,b_i)]\in\cO(\R)$ in view of~$(O_2)$.
\endproof

\brm8.95
We already mentioned (see \E\df\ \rf{d8.57}) that \he sms between \tps s could be called topological \is sms, since the \cm\ and the inverse of
\he sms are \he sms. See
\E\df\ 2.1.7 and Lemma 2.1.8 for the case of monoid-\is sms (resp.\ monoid-\hm sms). Since the \cm\ of two \ctn\ maps is a \ctn\ map, and since a
bi\jc\ map~$f$ \st both $f$~and~$f\Inv$ are \ctn\ is a \he sm, it is natural to call a \ctn\ map \ti{topological \hm sm} between \tps s. However, by
contrast to the case of a monoid-\hm sm, the inverse of a bi\jc\ \ctn\ map is not necessarily \ctn. For example, if $X_1=X_2:=\R$, $\cO_1:=\cP(\R)$,
and $\cO_2:=\cO(\R)$, $f(x)=x$ \fa $x\in\R$, then $f:X_1\to X_2$ is \ctn, but $f\mo$~is not since $\cO_2\notin \cP(\R)$.
\erm

In Remark \rf{r7.28}\,(iii) we mentioned that the pointwise \cg nce of \sq\ of \f s from the \il\ $[0,1]$ in~$\R$ into~$\R$ is not ``metrizable''.
More precisely and more generally, let $A$~be a \ns\ and let $\R^A$ denote the set of all maps from~$A$ into~$\R$. A~\sq\ $\zb fn\N$ in $\R_A$ is
said to be pointwise \cvg\ to an \el\ $f\in\R^A$, written $f_n\imp f$, if the \fw\ holds:
\bml8.72n
\hbox{\E\fe $a\in A$, the \sq\ $\{f_n(a)\}_{n\in\N}$ in $\R$ \cg s to $f(x)$}\\
\hbox{in the sense of \E\df\ \rf{d7.1}.}
\e
See also \E\df\ \rf{d7.19}, \E\Pr\ \rf{p7.20}.

The pointwise \cg nce on $\R^A$ is called \ti{metrizable} if
\beq8.73
f_n\nad{D}\to f \ \hbox{iff}\ f_n\imp f \q \hbox{\fa $f, f_n \in \R^A$ and }n\in\N
\e
holds for some metric $D$ on $\R^A$, and \ti{not metrizable} otherwise.

\bpr8.96
Let $A$ be a nonempty \emph{\ct e} set. Then the pointwise \cg nce on~$\R^A$ is metrizable.
\epr

Observe that the \df\ of $f_n\imp f$ uses the notion of \cg nce of \sq s in~$\R$ \itd in \E\df\ \rf{d7.1}.

\proof
See \E\Pr\ \rf{p7.20} if $A$ is nonempty and finite. See Theorem \rf{t7.26} if $A:=\N$.

\looseness=-1
Let $A$ be a \ct y infinite set and let $\vf:\N\to A$ be a bi\jn. In view of Theorem \rf{t7.26}, \te s a metric~$D$ on~$\R^\N$ \st $f_n\imp f$
in~$\R^\N$ iff $D(f_n,f)\to 0$ \fa $\zb fn\N$ and $f$ in~$\R^\N$. The idea is to ``transport'' the metric~$D$ from~$\R^\N$ to~$\R^A$. To this end we
define a map $\vf^*:\R^\N \to \R^A$ by setting $(\vf^*(\bX))(a) := \bX(\vf(a)) = (\bX\circ \vf)(a)$ \fa $a\in A$, and \fa $\bX\in\R^\N$. We also
define $((\vf\Inv)^*(f))(i) := f(\vf\Inv(i)) = {(f\circ\vf\Inv)(i)}$ \fa $i\in\N$, and \fa $f\in\R^A$. Then $(\vf\Inv)^*$ is a map from~$\R^A$
into~$\R^\N$. We have $((\vf\Inv)^*\circ \vf^*)(\bX) = (\vf\Inv)^*(\vf^*(\bX)) = (\vf\Inv)^*(\bX\circ\vf) = (\bX\circ\vf)\circ\vf\Inv =
\bX\circ(\vf\circ \vf\Inv) = \bX\circ\id_{\R^\N}=\bX$ \fa $\bX\in\R^\N$. Hence $(\vf\Inv)^*\circ\vf^*= \id_{\R^\N}$.
Similarly, we have $\vf^*\circ(\vf\Inv)^* = \id_{\R^A}$. Hence $\vf^*$ and $(\vf\Inv)^*$ are bi\jn s and
$(\vf^*)\Inv = (\vf\Inv)^*$. In view of Exercise \rf{ex7.61a} with $(X,d):=(\R^\N,D)$, $Y:=\R^A$, and $f:=\vf^*$, we find that $D_A:\R^A\t\R^A\to
\R_{\ge0}$ defined by $D_A(f,g):= D((\vf^*)\Inv f,(\vf^*)\Inv g)$, $f,g\in\R^A$, is a metric on~$\R^A$. \Mo $\vf^*$~is a sur\jc\ isometry from
$(\R^\N,D)$ onto $(\R^A,D_A)$. \E\Ip \fe $\zb fn\N$ and $f$~in~$\R^A$, we have $D_A(f_n,f)\to0$ iff ${D((\vf^*)\Inv f_n,(\vf^*)\Inv f)\to0}$.
By Theorem
\rf{t7.26} ${D((\vf^*)\Inv f_n,(\vf^*)\Inv f)\to0}$ iff $((\vf^*)\Inv f_n)(i) \to ((\vf^*)\Inv f)(i)$ \fa ${i\in\N}$. Since $(\vf^*)\Inv f_n(i) =
f_n(\vf\Inv(i))$, $n,i\in\N$, and $(\vf^*)\Inv f(i) = f(\vf\Inv(i))$, $i\in\N$, we obtain $f_n(\vf\Inv(i)) \to f(\vf\Inv(i))$ \fa $i\in\N$
iff $D_A(f_n,f)\to0$. \Mo since $\vf$~is bi\jc, there is \ooo $i\in\N$ \st $a=\vf(i)$. \E\Tf $f_n\imp f$ iff ${D_A(f_n,f)\to 0}$.
\endproof

It turns out that the pointwise \cg nce on~$\R^A$ is \ti{not\/} metrizable if $A:=[a,b]$, $a,b\in\R$, $a<b$. We are not yet in a position to prove
this assertion. However, we only show that \te s a nonmetrizable topology on~$\R^A$ denoted by $\cO(\R^A)$ \sf ying
\beq8.75
f_n \nad{\cO(\R^A)}\to f \qh{iff } f_n\imp f \hbox{ \fa} \zb fn\N \hbox{ and all }f\in\R^A.
\e
To this end we first prove that if a topology~$\cO$ on~$\R^A$ makes all \co\ maps $\psi_a$, $a\in A$, from $(\R^A,\cO)$ into $(\R,\cO(\R))$
\ti{\ctn} then \er{8.75} holds with $\cO(\R^A):=\cO$.

\blm8.97
Let $(X,\cO)$, $(X',\cO')$ be \Hts s and let $\vf:X\to X'$ be \ctn\ at~$x$. Then \fa \sq s $x_n\nad\cO\to x$ we have $\vf(x_n)\nad{\cO'}\to \vf(x)$.
\elm

\proof
Given $O'\in\cO'$ containing $\vf(x)$, we have $\vf\mo[O'] = \{y\in X: \vf(y)\in O'\}$, hence $x\in\vf\mo[O']$. \Mo $\vf\mo[O']\in\cO$ by the
continuity of~$f$ at~$x$. Since $x_n\nad\cO\to x$, \te s an $N\in\N$ \st $x_n\in\vf\mo[O']$ \fa $n\ge N$. Hence $\vf(x_n)\in\vf(\vf\mo[O'])\nde8.67
\sbs O'$ \fa $n\ge N$. \E\Tf $\vf(x_n)\nad{\cO'}\to \vf(x)$.
\endproof

\brm8.98 \

\hph i,ii, We recall that if $(X,\cO)$ in Lemma \rf{l8.97} is metrizable, the \ti{sequential\/} continuity of~$\vf$ at \ti{every} $x\in X$ implies the
continuity of~$\vf$ in view of Theorem \rf{t8.89}.

\hph ii,i, Under the \as s of Lemma \rf{l8.97}, the same conclusion holds if we replace \sq s by \ti{nets} (see \cite[Theorem 2.1.2]{Du}).

\hph iii,, Assuming the axiom of choice, if every \cvg\ net $x_i\to x$ in $(X,\cO)$ implies \cg nce of~$\vf(x_n)$ to $\vf(x)$ in $(X',\cO')$, then
$\vf$~is \ctn\ (same reference as above).
%
\erm

We now show that \te s a topology on~$\R^A$ making all \co\ maps \ctn. This topology is a particular case of the so-called product topology.

Let $A$ be an infinite set. Suppose that \te s a topology $\cO$ on~$\R^A$ making all the \co\ maps \ctn, then in view of Theorem \rf{t8.89} and
\E\df\ \rf{d8.93} the sets $V_{\{a\}}(\a;\b):= \psi\mo_a[(\a,\b)]$, $a\in A$, $\a,\b\in\R$, $\a<\b$, belong to~$\cO$. By~$(O_3)$ finite \isc s of such
sets also belong to~$\cO$. Let $\cA$ denote the set of all \nfs s of~$\N$. \E\fe $B\in\cA$ and \fa $\bA,\bB\in\R^B$ \sf ying $\bA(i)<\bB(i)$,
$i\in B$, we define $V_B(\bA;\bB)$ by setting
\beq8.76
V_B(\bA;\bB) := \{g\in \R^A: \bA(i)<\psi_i(g)<\bB(i), \ i\in B\}.
\e
Observe that $V_B(\bA;\bB) = \Bca_{i\in B} V_{\{i\}}(\bA(i);\bB(i))$. Hence $V_B(\bA;\bB)\in\cO$ by~$(O_3)$.

It turns out that the family of subsets of~$\R^A$ of the form $V_B(\bA;\bB)$ \sf ies \cn s \er{8.41}--\er{8.43}. Indeed, we have:

``\er{8.41}'': \ \E\fe $f\in\R^A$ and every $B\in\cA$ $f\in V_B(\bA;\bB)$ with $\bA(i):=f(i)-1$, $\bB(i):=f(i)+1$ \fa $i\in B$.

``\er{8.42}'': \ Let $V_B(\bA;\bB)$ and $V_{B'}(\bA';\bB')$ with $B,B'\in\cA$, $\bA,\bB\in\R^B$, $\bA',\bB'\in \R^{B'}$, be \st $\bA(i)<\bB(i)$,
$\bA'(j)<\bB'(j)$ \fa $i\in B$ and all $j\in B'$. If $B\cap B'=\vn$ and $B'':=B\cup B'$, then $V_B(\bA;\bB)\cap V_{B'}(\bA';\bB') = V_{B''}(\bG;\bD)$
with
$$
\bG(i):=\bca
\bA(i) &\hbox{for } i\in B,\\
\bA'(i) &\hbox{for } i\in B',
\eca
\q\q
\bD(j):=\bca
\bB(j) &\hbox{for } j\in B,\\
\bB'(j) &\hbox{for } j\in B'.
\eca
$$
If $B\cap B'\ne\vn$ and $V_B(\bA;\bB)\cap V_{B'}(\bA';\bB')\ne\vn$, then $V_B(\bA;\bB)\cap V_{B'}(\bA';\bB') = V_{B\cup B'}(\bG;\bD)$ with
\beag
\bG(i)&:=\bca
\bA(i) &\hbox{for } i\in B\sm(B \cap B'),\\
\bA'(i) &\hbox{for } i\in B'\sm(B \cap B'),\\
\max(\bA(i),\bA'(i)) &\hbox{for }i\in B\cap B',
\eca
\\
\bD(j)&:=\bca
\bB(j) &\hbox{for } j\in B\sm(B \cap B')=B\sm B',\\
\bB'(j) &\hbox{for } j\in B'\sm(B \cap B')=B'\sm B,\\
\min(\bB(j),\bB'(j)) &\hbox{for }j\in B\cap B'.
\eca
\e
Note that if $i\in B\cap B'$ and $V_B(\bA;\bB) \cap V_{B'}(\bA';\bB') \h{\ne}\h\vn$, then $(\bA(i),\bB(i))\cap (\bA'(i),\bB'(i))\h{\ne}\h\vn$, hence
$\max(\bA(i),\bA'(i)) < \min(\bB(j),\bB'(j))$.

``\er{8.43}'': \ Let $f,g\in\R^A$ be \st $f\ne g$. Then \te s $i\in A$ \st $f(i)\ne g(i)$. \E\wlg we may assume $f(i)<g(i)$ (otherwise interchange
$f$~and~$g$). Then if $W_1:= V_{\{i\}}\bigl(f(i)-1;\frac12(f(i)+g(i))\bigr)$, $W_2:= V_{\{i\}}\bigl(\frac12(f(i)+g(i));g(i)+1\bigr)$, then
$f\in W_1$, $g\in W_2$ and $W_1\cap W_2=\vn$. In view of Lemma \rf{l8.44} \te s a Hausdorff topology, denoted by $\cO(\R^A)$, \sf ying $\cB(\R^A)\sbs
\cO(\R^A)$ with
\beq8.77
\cB(\R^A) := \{V_B(\bA;\bB): B\in \cA,\ \bA,\bB\in \R^B,\ \bA(i)<\bB(i),\ i\in B\}.
\e
Then $\cO(\R^A)$ consists of the empty set and all \ti{unions} of \el s of $\cB(\R^A)$.

\ssk
We now prove the continuity of the maps $\psi_a$, $a\in A$, from $(\R^A,\cO(\R^A))$ into $(\R,\cO(\R))$. Let $a\in A$, and let $U\in\cO(\R)\sms\vn$.
We have to show that $\psi_a\mo[U] \in \cO(\R^A)$. \E\te s a family of open \il s $(\a_i,\b_i)$, $i\in I$, of~$\R$ \st $U=\bcl_{i\in I}(\a_i,\b_i)$.
Then $\psi_a\mo[U] = \psi_a\mo\bigl[ \bcl_{i\in I}(\a_i,\b_i)\bigr] \nad*= \bcl_{i\in I}\psi_a\mo[(\a_i,\b_i)] = \bcl_{i\in I} V_{\{i\}}(\a_i;\b_i)
\in\cO(\R^A)$ by~$(O_2)$. In~$\nad*=$ we used Morgan's law~\cite{Kelley}. \E\Tf the topology $\cO(\R^A)$ on~$\R^A$ makes all \co\ maps \ctn.

\begin{dfn}[{\cite[pp.\ 46--48]{Kelley}}]\label{d8.99}
Let $(X,\cO)$ be a \tps. A subfamily $\cB$ of the topology~$\cO$ is said to be a \ti{base} for~$\cO$ if every nonempty member of~$\cO$ is the union
of members of~$\cB$. A~subfamily $\cS$ of~$\cO$ is said to be a \ti{subbase} for the topology~$\cO$ if the collection of all finite \isc s of
members of~$\cS$ is a base for~$\cO$.
\edf

\brm8.100
In Lemma \rf{l8.44} $\cB(X)$ is a base for the topology $\cO[\cB]$. \If from the proof that only \cn s \er{8.41} and \er{8.42} are used. If \er{8.43}
also holds, then $\cO[\cB]$ is a Hausdorff topology.
\erm

\bxa8.101 \

\hph i,i, The \il s of the form $(\lea,b)$ and $(a,\to)$, $a,b\in\R$, form a subbase of the order topology $\cO(\R)$ of~$\R$.

\hph ii,, The sets $V_B(\bA;\bB)$ in \er{8.76} form a base of the topology $\cO(\R^A)$. The sets $V_{\{i\}}(\bA;\bB)$ form a subbase of the topology
$\cO(\R^A)$.
\exa

\bex8.102
Show that if $\cO$ is a topology on~$\R^A$ making all \co\ maps \ctn, then $\cO(\R^A)\sbs \cO$. Thus the topology $\cO(\R^A)$ is the smallest
(weakest) topology on~$\R^A$ making all \co\ maps \ctn.
\eex

\If from Lemma \rf{l8.97} that $f_n\nad{\cO(\R^A)}\to f$ implies $f_n\imp f$ \fa $f,f_n\in\R^A$, $n\in\N$. We now prove the converse. Assuming
$f_n\imp f$, and given $O\in\cO(\R^A)$ with $f\in O$, we have to find an $N\in\N$ \st $f_n\in O$ \fa $n\ge N$ in view of \E\df\ \rf{d8.83}.
Let $O\in\cO(\R^A)\sms\vn$. Then $O$~is the union of sets of the form $V_B(\bA;\bB)$. If $f\in O$ then \te s at least one \el\ of the base of the
topology $\cO(\R^A)$ containing~$f$. Thus \te\ $B\in\cA$ and $\bA,\bB\in\R^B$ \st $\bA(i)<\psi_i(f)<\bB(i)$, $i\in B$, by~\er{8.41} since
$f(i) = \psi_i(f)$, $i\in B$. We have to show that \te s an $N\in\N$ \st $f_n\in V_B(\bA;\bB)$ \fa $n\ge N$. Setting $\rho_i:= \min(\bB(i)-f(i),
f(i)-\bA(i))$, we have
\beq8.78
\rho_i>0 \hbox{\q and }\q \bA(i) < f(i)-\rho_i < f(i) < f(i)+\rho_i < \bB(i)
\e
\fa $i\in B$. Hence if $\bG(i):= f(i) - \rho_i$, and $\bD(i):= f(i)+ \rho_i$, then we obtain
\beq8.79
f \in V_B(\bG;\bD) \sbs V_B(\bA;\bB).
\e
Given $\ve\in\R_{>0}$, $\ve<\min\limits_{i\in I}\rho_i$, \te s $N_i$, $i\in \N$, \st $|f(i)-f_n(i)|<\ve$ \fa $n>N_i$. Setting $N:=
\max\limits_{i\in I}N_i$, we find that \fa $i\in I$ we have $|f(i)-f_n(i)|<\ve$, hence $f_n(i) = f_n(i)-f(i)+f(i) \le |f_n(i)-f(i)| + f(i) <
f(i)+\ve \le f(i)+\rho_i$. Similarly, $f_n(i) \ge f(i)- |f_n(i)-f(i)| \ge f(i)-\ve \ge f(i)-\rho_i$, $i\in I$. \E\Tf for $n\ge N$ we find
$f_n \in V_B(\bG;\bD) \sbs V_B(\bA;\bB)$.

\E\Tf $f_n \nad{\cO(\R^A)}\longrightarrow f$. \csq, we proved that
\beq8.80
f_n \nad{\cO(\R^A)}\longrightarrow f \qh{iff  }f_n\imp f \qh{\fa} f,f_n\in\R^A,\ n\in\N.
\e

We now show that the topology $\cO(\R^A)$ is \ti{not\/} metrizable whenever $A$~is un\ct e. The idea is to show that the map $\1_A\in \R^A$ defined
by $\1(i):=1$ \fa $i\in A$ belongs to the closure of the subset~$H$ of~$A$ defined by
\beq8.81
H := \{\1_B \in \R^A: B\in\cA\}
\e
where $\1_B(i):=1$ for $i\in B$ and 0 otherwise. We will show that $\1_A\in \ov H$. Then assuming, for \cd ion, that $\cO(\R^A)$ is metrizable, we
show that $\1\notin \ov H{}\en s$ defined in \E\df\ \rf{d8.87}. Hence $\cO(\R^A)$ is not metrizable in view of \E\Pr\ \rf{p8.88}. We first prove that
$\1_A\in\ov H$ whenever $A$~is infinite (the case $A$~finite will be treated later). To this end we use the \fw\ lemma.

\blm8.103
Let $(X,\cO)$ be a \Hts\ and let $Y$~be a \nss\ of~$X$ \st $Y\ne X$. Then~$\ov Y$, the closure of~$Y$ defined in \er{8.57}, \sf ies
\beq8.81a
\ov Y = \{x\in X: \hbox{every $O\in\cO$ containing $x$ intersects }Y\}.
\e
\elm

\proof
$\ov Y{}^c = \Bigl(\Bca_{Z\supset Y,Z^c\in\cO}Z\Bigr)^c = \bcl_{Z^c\sbs Y^c,Z^c\in\cO}Z^c \nad*= \bcl_{O\sbs Y^c,O\in\cO}O = \kl{(Y^c)}$.
In $\nad*=$ we used \er{8.44}. Clearly, $\kl Y\cap \kl{(Y^c)}=\vn$ since $\kl Y\sbs Y$ and $\kl{(Y^c)}\sbs Y^c$. Set $K:=X \sm \bigl(\kl Y\cup
\kl{(Y^c)}\bigr)$. In view of \er{8.44} we find that $K=\{x\in X: \hbox{every }O\in\cO$ containing~$x$ intersects $Y$ and~$Y^c\}$. \If that $X$~is
the disjoint union of~$\kl Y$, $\kl{(Y^c)}$ and~$K$. Since $\kl{(Y^c)}=\ov Y{}^c$ we find that $\ov Y=\kl Y\cup K$ and $\kl Y\cap K=\vn$. (\E\Ip
we find that $K=\ov Y\sm \kl Y =: b(Y)$ by \E\df\ \rf{d8.73}.) \csq, \fe $x\in X$ we have $x\in \ov Y$ iff either $x\in\kl Y$ or $x\in K$. \E\Tf
if $x\in\ov Y$ and $O\in\cO$ with $x\in O$, then $O\cap\{x\} = \{x\}\sbs Y$ if $x\in\kl Y$, and $O\cap Y\ne\vn$ if $x\in K$.

Conversely, if $x\in X$ and \fe $O\in\cO$ \st $x\in O$, we have $O\cap Y\ne\vn$, then $x\notin \kl{(Y^c)}$, thus $x\in\kl{(Y^c)}{}^c = (\ov Y{}^c)^c
= \ov Y$.
\endproof

\blm8.104
$\1_A\in \ov H$ where $\ov H$ is the closure of the set~$H$ defined in \er{8.78} in the \tps\ $(\R^A,\cO(\R^A))$.
\elm

\proof
In view of the preceding lemma, it is sufficient to show that \fe set $O\in\cO(\R^A)$ containing~$\1$ we have $O\cap H\ne\vn$. By \df\ $O$~is the
union of sets of the form $V_B(\bA;\bB)$ defined in~\er{8.76}, \te\ $B'\in \cA$ and $\bA',\bB'\in\R^{B'}$ \st $\1\in V_{B'}(\bA';\bB')$. Note that
\fa $B\in\cA$, $\1\in V_B(\bA;\bB)$ iff $\bA(i)<1<\bB(i)$ \fa $i\in B$. \If that $\1_{B'}\in V_{B'}(\bA';\bB')$ and $\1_{B'}\in H$. Hence $O\cap H
\supset V_{B'}(\bA';\bB')\cap H\supset \{\1_{B'}\}\ne\vn$.

We now suppose, for \cd ion,  that \te s a metric $D_A$ \st $\cO(\R^A)=\cO(D_A)$. We infer that $\ov H=\ov H{}\en s$ (see \E\df\ \rf{d8.87}) by
\E\Pr\ \rf{p8.88}. However, we will prove that there is \ti{no} \sq\ $\zb gn\N$ in~$H$ \sf ying $D(g_n,\1)\to0$. Indeed, if $D(g_n,\1)\to0$, then
$g_n\imp\1$, hence $g_n(i)\to1$ \fa $i\in A$. Set $B_n:= \{i\in A: g_n(i)\ne0\}$, $n\in\N$. By~\er{8.81} $B_n$~is finite \fa $n\in\N$. Set $C:=
\bcl_{n\in\N}B_n$. As a \ct e union of \ti{finite} sets, the set~$C$ is \ct e. \csq, \fa $i\in C^c$ and all $n\in\N$, $g_n(i)=0$. \E\Tf $g_n(i)\to0$
\fa $i\in C^c$. Observe that $C^c$ is un\ct e since the union of two \ct e sets is \ct e, and $A=C\cup C^c$ is un\ct e. \E\Ip $C^c\ne\vn$. Thus if
$g\in \ov H{}\en s$, \te s $i\in A$ \st $g(i)=0$. \If that $\1\notin \ov H{}\en s = \ov H$, \cd ing \E\Pr\ \rf{p8.88}.
\endproof

We have thus proved

\bth8.105
Let $A$ be an infinite set and let $\R^A$ denote the set of all maps from~$A$ into~$\R$. Then \te s a topology on~$\R^A$, denoted by~$\cO(\R^A)$,
\sf ying the \fw\ \ass s. \E\fa $f,f_n\in\R^A$ and all $n\in\N$ \er{8.75} holds. \Mo the topology $\cO(\R^A)$ is \emph{not} metrizable whenever
$A$~is un\ct e.
\eth

We now introduce on a subset of $\R^A$ a metric which \ch izes the notion of \ti{uniform} \cg nce of \sq s of \f s. This notion of \cg nce plays an
important role in Analysis since limits of uniformly \cvg\ \sq s of \ctn\ \f s are \ctn.

\bnt8.106
Let $A$ be a \ns. We denote by $B(A;\R)$ the set of all bounded real-valued maps from~$A$ into~$\R$, that is,\index{bounded real-valued map}
\bml8.82
B(A;\R) := \{g:A\to\R: \hbox{\te s $M\in\R_{>0}$ depending on $g$}\\
\hbox{\st $|g(x)|\le M$ \fa}x\in A\}.
\e
Let $D_{\sup}$ denote the metric defined by
\beq8.83
D_{\sup}(f,g) := \sup_{x\in A}|f(x)-g(x)|.
\e
\ent

The \RHS\ of \er{8.83} is well defined since the subset of~$\R_{\ge0}$ $\{|f(x)-g(x)|: {x\in A}\}$ is \ba. Indeed, for $f,g\in B(A,\R)$ \te\ $m_f,m_g,
M_f,M_g$ in~$\R$ \st $m_f\le f(x)\le M_f$, $m_g\le g(x)\le M_g$ \fa $x\in A$. Thus $|f(x)|\le \max\{|m_f|,|M_f|\}$ \fa $x\in A$ and similarly for~$g$.
Setting $M_1:=\max\{|m_f|,|M_f|\}$, $M_2:=\max\{|m_g|,|M_g|\}$, we find $|f(x)-g(x)|= |f(x)+(-g(x))| \le |f(x)|+|g(x)|\le M_1+M_2$ \fa $x\in A$. Note
that in the case where $A$~is infinite we use the \ti{order-\cp ness} of~$\R_{\ge0}$. We now verify (M1)--(M3) of \E\df\ \rf{d7.8}.

``(M1)'': \ $D_{\sup}(f,g)=0$ iff $\supl_{x\in A}|f(x)-g(x)|=0$ iff $|f(x)-g(x)|=0$ \fa $x\in A$ since $|f(x)-g(x)|\ge0$ \fa $x\in A$. \E\Tf $f(x)=
g(x)$ \fa $x\in A$, and $f=g$.

``(M2)'': \ follows from $|f(x)-g(x)| = \left|-1\right|\,|f(x)-g(x)| = \left|-f(x)-(-g(x))\right| = \left|-f(x)+g(x)\right|=|g(x)-f(x)|$ \fa $x\in A$.

``(M3)'': \ Let $f,g,h\in B(A;\R)$. Then $|f(x)-g(x)| = \bigl|(f(x)-h(x))+(h(x)-g(x))\bigr| \le |f(x)-h(x)| + |h(x)-g(x)|$, $x\in A$. \E\Tf $|f(x)-
g(x)| \nad{\er{8.83},\er{3.8}}\le \supl_{y\in A}|f(y)-h(y)| + |h(x)-g(x)| \le \supl_{y\in A}|f(y)-g(y)| + \supl_{y\in A}|h(y)-g(y)|$. Hence
$\supl_{x\in A}|f(x)-g(x)| \le D_{\sup}(f,h)+D_{\sup}(h,g)$ since $\supl_{x\in A}|f(x)-g(x)|$ is the least \ub\ for the set $\{|f(x)-g(x)|: x\in A\}$.

Thus $D_{\sup}$ is a metric on $B(A;\R)$.

\blm8.107
$D_{\sup}$ defined in \er{8.83} is a metric on the set $B(A;\R)$ defined in \er{8.82}. \Mo \fa $f,f_n\in B(A;\R)$, $n\in\N$, the \fw\ holds\dw
\beq8.84
f_n\nad{D_{\sup}}\to f \qh{implies }f_n\imp f.
\e
\elm

\proof
Suppose that \fe $\ve\in\R_{>0}$ \te s an $N\in\N$ \st $D_{\sup}(f_n,f)<\ve$ \fa $n\ge N$, then \fa $x\in A$ we have
$$
|f_n(x)-f(x)| \le \sup_{y\in A} |f_n(y)-f(y)| = D_{\sup}(f_n,f)<\ve.
$$
\E\Tf $f_n\imp f$. Note that $N$ only depends on~$\ve$ and not on $x\in A$. By contrast, $f_n\imp f$ means that \fe $x$ and \fe $\ve\in\R_{>0}$
\te s $N=N(x,\ve)\in\N$ \st $|f_n(x)-f(x)|<\ve$ \fa $n\ge N$.
\endproof

\bdf8.108
A \sq\ $\zb fn\N$ in $B(A;\R)$ is said to \cg\ to $f\in B(A;\R)$ \ti{uniformly} if $D_{\sup}(f_n,f)\to 0$.\index{uniform \cg nce}
\edf

\bex8.109
Let $A:=[0,1]$ in $\R$ and let $f_n(x):=x^n$, $x\in A$ and $n\in\N$. Let $f(x):=0$ for $x\in\zo0,1 $ and $f(1):=1$. Show that $f_n\imp f$ but
$D_{\sup}(f_n,f)\not\to0$. Show that if $A:=\zo0,1 $ then $f_n(x)\to0$ \fa $x\in A$, but $f_n$ does \ti{not\/} \cg\ to~$0$ uniformly.
\eex

We now consider an important example.

\bxa8.110
Let $A:=K$ with $(K,d)$ being a \sql\ compact \ms, for example $K:=[a,b]$, $a,b\in\R$, $a<b$. \If from Lemma \rf{l8.30} and \er{8.30} that \sql\
\ctn\ real-valued \f s on a \sql\ compact \ms\ are bounded. Since $K$~is a \ms, it follows from Theorem \rf{t8.89} that such \f s are also \ctn.
\exa

\bnt8.111
We denote by $C(K)$ the set of all \ctn\ real-valued \f s on a \sql\ compact \ms\ $(K,d)$.
\ent

We just showed that $C(K)$ is a subset of $B(A;\R)$. More is true.

\blm8.111
$C(K)$ is a \emph{closed} subset of the \ms\ $(B(K;\R),D_{\sup})$. In other words, ``uniform limits'' of \ctn\ \f s on a \sql\ compact \ms\ are \ctn.
\elm

\proof
Let $f_n\in C(K)$, $n\in\N$, and let $f\in B(K;\R)$ be \st $f_n\nad{D_{\sup}}\to f$. To prove $f\in C(K)$. By \E\Pr\ \rf{p8.94}\,(iv) it suffices to
show that \fe $x\in K$ and every $\ve\in\R_{>0}$ \te s $\d\in\R_{>0}$ \st $|f(x)-f(y)|<\ve$ whenever $d(x,y)<\d$. Let $x\in K$ be arbitrary, and let
$\ve\in\R_{>0}$. Since $f_n\nad{D_{\sup}}\to f$, \te s an $N\in\N$ \st $D_{\sup}(f_n,f)<\frac\ve3$ \fa $n\in\N$. Let $y\in K$. Then $|f(x)-f(y)|
= |(f(x)-f_n(x))+(f_n(x)-f_n(y))+(f_n(y)-f(y))| \le D_{\sup}(f,f_n) + |f_n(x)-f_n(y)| + D_{\sup}(f,f_n) \le \frac23\ve + |f_n(x)-f_n(y)|$.
Since $f_n\in C(K)$, \te s $f\in\R_{>0}$ \st $|f_n(x)-f_n(y)|< \frac\ve3$ whenever $d(x,y)<\d$. Hence $|f(x)-f(y)|<\ve$ whenever $d(x,y)<\d$.
\endproof

\E\ctn\ real-valued \f s on a \sql\ compact \ms\ are uniformly \ctn.

\bdf8.113
Let $(X_i,d_i)$, $i=1,2$, be \ms s. A~map $f:X_1 \to X_2$ is called \tb{uniformly \ctn} if \fe $\ve\in\R_{>0}$ there is a $\d\in\R_{>0}$ \st \fa
$x,y\in X_1$ we have $d_2(f(x),f(y))\le\ve$ whenever $d_1(x,y)<\d$.\index{uniformly continuous}
\edf

\bxa8.114
If $f$ \sf ies a \cn\ of the form $d_2(f(x),f(y)) \le M(d_1(x,y))^\a$ with $\a\in(0,1]$ \fa $x,y\in X_1$, then given $\ve>0$, $\d$~can be chosen
equal to $(\frac1M \ve)^{\frac1a}$.
\exa

\bth8.115
Let $(K,d_1)$ be a \sql\ compact \ms\ and let $(Y,d_2)$ be a \ms\ $($for example $\R$ with the usual metric$)$. Let $f:K\to Y$ be \ctn. Then $f$~is
\emph{uniformly} \ctn.
\eth

\proof
Suppose, for \cd ion, that $f$ is not uniformly \ctn. Then \te s an $\ve\in\R_{>0}$ \st \fa $n\in\N$ \te\ $x_n,y_n\in K$ \sf ying
$$
d_1(x_n,y_n)\le \tfrac1n \qh{and \ } d_2(f(x_n),f(y_n))\ge \ve.
$$
Then \te s an $a\in K$ and a sub\sq\ $\{x_{n_k}\}_{k\in\N}$ \st $x_{n_k}\nad{d_1}\to a$. Since $d_1(x_{n_k},y_{n_k})\le \frac1{n_k}$ for $k\in\N$,
we obtain $d_1(a,y_{n_k}) \le d_1(a,x_{n_k})+d_1(x_{n_k},y_{n_k})\le d_1(a,y_{n_k})+\frac1{n_k}$ \fa $k\in\N$. \E\Tf $\limu_{k\ge0} d_1(a,y_{n_k})
\le \limu_{k\ge0} d_1(a,x_{n_k}) + \limu_{k\ge0}\frac1{n_k} = \lim\limits_{k\ge0} d_1(a,x_{n_k})+0 = 0$. \E\Tf $d_1(a,y_{n_k})\to0$. \Mo
$d_2(f(x_{n_k}),f(y_{n_k})) \le d_2(f(x_{n_k}),f(a)) + d_2(f(a),f(y_{n_k}))$. In view of the continuity of~$f$ at~$a$, \te s an $N\in\N$ \st
$d_2(f(x_{n_k}),f(a))\le \frac\ve3$ and $d_2(f(y_{n_k}),f(a))\le\frac\ve3$ whenever $k\ge N$. \csq, for $k\ge N$ we obtain $\ve\le d_2(f(x_{n_k}),
f(y_{n_k}))\le\frac23\ve$. A~\cd ion.
\endproof

Our next goal is to present an \ev t \df\ of sequential compactness in~$\R$, involving only the topology induced by the metric.

\begin{dfn}[\cite{Kelley}]\lb{d8.116}
Let $B$ be a \ns. A~family $\cA$ of sets is called a \tb{cover} $B$ if $B$~is a subset of $\bcl_{A\in\cA}A$. A~\tb{subcover} of~$\cA$ is
a subfamily of~$\cA$ which is also a cover of~$B$.\index{cover of a set}\index{subcover}
\edf

\bdf8.117
Let $(X,\cO)$ be a \Hts. A nonempty set~$X$ is called \tb{compact\/} if every (open) cover of~$X$ \sf ying $\cA\sbs\cO$ has a finite subcover. A~\nss\
$B$ of~$X$ is called \ti{compact\/} if $B$~is compact when endowed with the \rv\ topology $\cO_B$ (see \E\df~\rf{d8.51}).\index{compact set}
\edf

\bex8.118
Let $\R$ be endowed with the order topology $\cO(\R)$ (see Remark \rf{r8.46}\,(i)). Prove the \fw\ \ass s:

\hph1,, Every finite subset of $\R$ is compact.

\hph2,, Let $\zb xn\N$ be an infinite \sq\ in $(\R,\cO(\R))$ converging to some $x\in\R$. Then the set $\{x_n:n\in\N\}\cup\{x\}$ is compact and the
set $\{x_n:n\in\N\}\sms x$ is \ti{not\/} compact.

\hph3,, The subsets $\N$ and $\Z$ are \ti{not\/} compact.

\hph4,, Every compact subset of $(\R,\cO(\R))$ is bounded (\wrt the usual \og\ of~$\R$).
\eex

\begin{thm}\lb{t8.119}
The $($closed\/$)$ \il\ $[a,b]$ in $\R$ with $a,b\in\R$, $a<b$, is \emph{compact}.
\eth

\proof[Proof \rm(\cite{Choquet,Du})]
Let $\cA:= \{\o_i\in\cO_{[a,b]}: i\in I\}$ be an open cover of~$[a,b]$. Let $A:=\break \{y\in[a,b]: [a,y]$ can be covered by a union of finitely
many sets in~$\cA\}$. We claim that $A=[a,b]$. We first show that

``$[a,d) \sbs A$ \fs $d\in\R$, $a<d\le b$'': \ Indeed, \te s $j\in I$ \st $a\in\o_j$. By \E\df\ \rf{d8.51} \te s $O\in\cO(\R)$ \st $\o_j=O\cap[a,b]$.
Since $a\in\o_j$, we have $a\in O$. By the \df\ of $\cO(\R)$, \te\ $c,d\in\R$ with $c<d$ \st $a\in(c,d)\sbs O$. Hence we have
$\zo a,d = (c,d)\cap[a,b]\sbs O\cap[a,b] = \o_j$. \E\Tf $[a,d']\sbs \o_j$ \fa $d'\in\R$ \sf ying $0<d'<d$. \E\Ip $\zo a,d \sbs A$, which proves
the claim.

Since $A\ne\vn$ and $b$ is an \ub\ for~$A$, \te s $\ov d\in\oz a,b $ \st $\ov d=\sup A$, in view of the order-\cp ness of the \il\ $[a,b]$ in
$(\R,\ge)$.

We claim

``$\ov d=b$'': \ Suppose, for \cd ion, that $\ov d<b$. \E\te s $l\in I$ \st $\ov d\in\o_l$. Let $e\in(a,\ov d)\cap \o_l$ and let $\cB$~be a finite
subcollection of~$\cA$ which covers $[a,\ov d]$. Then the subcollection of~$\cA$ consisting of~$\cB$ and~$\o_l$ covers $[a,\ov e]$ where $\ov e\in
(\ov d,b)\cap \o_l$. \csq, $\ov d\ne\sup A$. A~\cd ion. \E\Tf $\ov d=b$, hence $A=[a,b]$.
\endproof

We recall that in view of Theorem \rf{t6.14} every \el\ of~$\R$ between $0$~and~$1$ is the supremum of an in\cre\ \sq\ of decimal \nm s, hence of
\ra\ \nm s between $0$~and~$1$. \If that $\ov{[0,1]}{}\en s_\Q$, the sequential closure of $[0,1]_\Q$, is equal to $[0,1]_\R$. Since the
order topology $\cO(\R)$ is the same as the topology induced by~$d$, the usual metric of~$\R$, that is, $\cO(\R)=\cO(d)$, it follows from \E\Pr\
\rf{p8.88} that $\ov{[0,1]}{}\en s_\Q = \ov{[0,1]}_\Q$, the closure of $[0,1]_\Q$ in $(\R,\cO(d))$.


\bdf8.120
Let $A$ be a \nss\ of a \Hts\ $(X,\cO)$. The set $A$ is said to be \ti{dense} in~$X$ if $\ov A=X$. The space $(X,\cO)$ is called \tb{separable} if
\te s a \ti{\ct e} dense subset of~$X$.\index{separable topological space}
\edf

\bth8.121
Let $(X,\cO)$ be a compact metrizable \tps. Then $X$ is separable.
\eth

\proof
Let $D$ be a metric \st $\cO=\cO(D)$. \E\fe $x\in X$ \te s a ball $B(x;r)$ with $r\in\R_{>0}$ \st $x\in B(x;r)$. In view of the compactness of~$X$,
\fe $n\in\Na$ \te s a finite cover of~$X$ consisting of balls $B_j\en n(x_j\en n;\frac1n)$, $j\in[1,N\en n]$, \fs $N\en n\in\Na$. Set
$A:=\bcl_{n\in\Na}\bcl_{j\in [1,N\en n]}\{x_j\en n\}$. As a \ct e union of finite sets, the set~$A$ is \ct e. Let $O\in\cO(d)$. \E\te s $B(x;r)\sbs O$
\fs $x\in X$ and $r\in\R_{>0}$. Let $n\in\Na$ be \st $\frac1n<\frac r3$. Such an $n\in\Na$ exists since $\R$~is \Ar. Since $\bcl_{j\in[1,N\en n]}
B_j\en n(x_j\en n;\frac1n) = X$, \te\ $j\in[1,N\en n]$ and $y\in X$ \st $y\in B(x;\frac13r) \cap B_j\en n(x_j\en n;\frac1n)$. Then $d(x,x_j\en n)
\le d(x,y)+d(y,x_j\en n) < \frac13r + \frac1n < \frac23r <r $. \csq, $B(x;r)\cap A\ne\vn$,
hence $O\cap A\ne\vn$. From Lemma \rf{l8.103} we infer that $\ov A=X$, hence $X$~is separable.
\endproof

\brm8.122
The notion of separability for a \tps\ is a topological \pp y or topological invariant (see~\cite{JvN}). More precisely, if $(X,\cO)$ and $(X',\cO')$
are \tps s (see \E\df\ \rf{d8.41}), and $f:X\to X'$ is a bi\jn\ \st both $f$~and~$f\Inv$ map open sets into open sets (as in Remark~\rf{r8.95}), then
a \pp y of a subset of~$X$ or of a \sq\ in~$X$ is called \ti{topological\/} if it is preserved by the map~$f$. By \df\ the notion of \ti{openness}
for a~set (the \pp y to be open) is topological. One verifies that the \pp y to be Hausdorff (see \E\df\ \rf{d8.41}) is topological.
\erm

\bex8.123
Show that the \fw\ \pp ies and notions are topological: Closedness, Closure, Interior, Compactness, Base of a topology, Separability, \E\cvg\ \sq.
\eex

\bpr8.124
The \tps\ $(\R,\cO(\R))$ is separable.
\epr

\proof
We show that $\ov\Q$, the closure of~$\Q$ in $(\R,\cO(\R))$, is equal to~$\R$. Note that $\cO(\R)=\cO(d)$, the topology induced by the usual
metric $d$ on~$\R$, hence by \E\Pr\ \rf{p8.88} we have $\ov \Q=\ov \Q\en s$. It remains to show that every \el\ of $\R\sm\Q$ is a limit of \sq s of
\el s of~$\Q$ \wrt the usual metric on~$\R$, since every \el\ of~$\Q$ is the limit of a constant \sq\ in~$\Q$. Let $x\in\R_{>0}\sm\Q$. By Lemma
\rf{l6.43} \te s a \sq\ $\zb xn\Na$ \st $x_n\ua x$. Hence by Lemma \rf{l7.2}\,(iv) $x_n\to x$. If $y\in \R_{<0}\sm\Q$, then $-y\in\R_{>0}\sm\Q$,
and $-y=\lim x_n$ \fs \sq\ $\zb xn\N$ in $\R_{>0}\sm\Q$. Hence $-x_n=(-1)\cdot x_n \to (-1)(-y)= y$ by Lemma \rf{l7.6} \er{7.7}. \If that $\ov \Q =
\ov\Q\en s= \R$. One completes the proof by noting that $\Q$~is \ct e (see Theorem 4.5.44).
\endproof

We now introduce an important notion, which is \ti{not\/} topological, namely the notion of Cauchy \sq.

\bdf8.125
A \sq\ $\zb xn\N$ in a \ms\ $(M,d)$ is said to be a \tb{\Cs} if \fe $\ve\in\R_{>0}$ \te s an $N\in\N$ \st $d(x_n,x_m)<\ve$ \fa $n,m\in\N$ \sf ying
$n\ge N$ and $m\ge N$.\index{Cauchy \sq}
\edf

\bxa8.126
Let the field $\Q$ be endowed with the metric $d(x,y):=|x-y|$, $x,y\in\Q$, and let $\zb xn\N$ be an infinite in\cre\ \sq\ in~$\Q$ \ba\ by
some $r\in\Q$. By Theorem \rf{t5.1} and by Lemma \rf{l5.38} \te s $x\in\R$ \st $x_n\ua x$, hence $x_n\to x$ in~$\R$ by \er{7.1} and \E\df\ \rf{d7.1}.
\E\fa $n,m\in\N$ we have
$x_n-x_m = (x_n-x)+(x-x_m)$, hence $|x_n-x_m| \le |x_n-x|+|x-x_m|$. Given $\ve\in\R_{>0}$ \te s~$N$ in~$\N$ \st $|x_n-x|<\frac\ve2$ for
$n\ge N$ and $|x-x_m|<\frac\ve2$ for $m\ge N$. Hence $|x_n-x_m|<\ve$ \fa $m,n\in\N$ \sf ying $m\ge N$ and $n\ge N$. Thus $\zb xn\N$ is a \Cs\ in
$(\Q,d)$ and in $(\R,d)$. Observe that the \sq\ $\zb xn\N$ \cg s to~$x$ in~$\R$, but is not \cvg\ in~$\Q$ if $x\notin\Q$.
\exa

\bth8.127
Let the field $\R$ be endowed with the usual metric. Then \fe \sq\ $\zb xn\N$ in~$\R$, the \sq\ $\zb xn\N$ is \cvg\ iff it is a \Cs.
\eth

\proof
Let $\zb xn\N$ be a \Cs\ in~$\Q$ (resp.~$\R$). If we set $\ve:=1$, \te s an $N\in\N$ \st $|x_n-x_m|\le 1$ \fa $n,m\ge N$. Setting $m:=N$ we
obtain $|x_n| = |x_n-x_m+x_m| \le |x_n-x_m|+|x_m| \le 1+|x_m|$ \fa $n\ge N$. \E\Tf $|x_n|\le \max\bigl(\max\limits_{1\le i\le N}|x_i|,1+|x_N|)$ \fa
$n\in\N$. \If that \Cs s in $(\Q,d)$ (resp.\ $(\R,d)$) are \ti{bounded\/} in $(\Q,\ge)$ (resp.\ $(\R,\ge)$).

If the \sq\ $\zb xn\N$ is a \sq\ in~$\R$, we find that $\limu x_n$ and $\liml x_n$ exist in~$\R$ and $\liml x_n\le \limu x_n$ by
\E\Pr s \rf{p8.6} and~\rf{p8.10}. \E\Tf we obtain $0\le \limu x_n-\liml x_n = (\limu x_n-x_k)+(x_k-x_l)+(x_l-\liml x_n)$ \fa $k,l\in\N$. \If that
$0\le \limu x_n - \liml x_n \le |\limu x_n - x_k| + |x_k-x_l| + |x_l-\liml x_n|$ \fa $k,l\in\N$.

In view of \E\Pr\ \rf{p8.17} \te\ sub\sq s $x_{\si(k)}$ and $x_{\rho(l)}$ converging to $\limu x_n$ and $\liml x_n$, \rsp. Thus given $\ve\in\R_{>0}$,
\te\ $N_1,N_2,N_3 \in \N$ \st $|\limu x_n-x_{\si(k)}|<\frac\ve3$ for $k\ge N_1$, $|x_{\rho(l)}-\liml x_n|<\frac\ve3$ for $l\ge N_2$. \csq,
$0\le \limu x_n-\liml x_n < \frac{2\ve}3 + |x_{\si(k)}-x_{\rho(l)}|$ \fa $k,l\ge \max(N_1,N_2)$. Since the \sq\ $\zb xn\N$ is a \Cs, \te s an $N_3\in
\N$ \st $|x_m-x_{m'}|<\frac\ve3$ provided $m,m'\ge N_3$. Setting $N_4:=\max(N_1,N_2)$, $N:=\max(\rho(N_4),\si(N_4),N_3)$, and using the strict in\cre
ness of the maps $\rho,\si$, we find $0\le \limu x_n-\liml x_n <\ve$. Since $\ve\in\R_{>0}$ is arbitrary and since the field $(\R,\ge)$ is \Ar, we
find that $0\le \limu x_n-\liml x_n=0$, hence $\limu x_n= \liml x_n$. \csq, in view of Theorem \rf{t8.11} we have: $\lim x_n$ exists and is equal to
$\limu x_n$. We have thus proved that every \Cs\ in the \ms\ $(\R,d)$ is \ti{\cvg\/}.

Proceeding as in Example \rf{xa8.126} we find that in $(\Q,d)$ (resp.\ $(\R,d)$) every \cvg\ \sq\ is a \Cs\ since $|x_k-x_l| \le |x_k-x| + |x-x_l|$
with $x=\lim x_n$.
\endproof

\bex8.128
Show that if a \Cs\ in a \ms\ $(M,d)$ possesses a \cvg\ sub\sq, then the whole \sq\ is \cvg.
\eex

\bdf8.129
A \ms\ $(M,d)$ is called (metrically) \ti{\cp} if every \Cs\ is \cvg. A \cp\ separable \ms\ is called a \tb{\Ps}. A~\Hts\ $(X,\cO)$ is called
\ti{Polish} if the topology~$\cO$ is induced by a metric~$d$ \st $(X,d)$ is Polish (see~\cite{Du}).\index{metric \cp ness}\index{Polish space}
\edf

\Ps s play an important role in the theory of measures on \tps s (see for example \cite[Theorem E.21]{JvN}).

We recall that the field $\R$ is a \cp\ \of\ by Theorem \rf{t5.1}, and that every \cp\ \of\ is order- and ring-\is c to the field $(\R,\ge)$ by
\E\Pr\ \rf{p5.14}\,(ii). From \E\Pr\ \rf{p4.31} we infer that a \cp\ \of\ is \Ar, hence an \Ar\ \of\ is order- and ring-\is c to the field $(\R,\ge)$
iff it is \ti{order-\cp} by Theorem \rf{t4.35} and \E\Pr\ \rf{p5.14}\,(i). In Theorem \rf{t8.53} it is shown that an \Ar\ \of\ is order-\cp\ iff it
is \ti{connected\/} \wrt the order topology $\cO(K)$. In the next theorem it will be shown that the connectedness can be replaced by the \ti{metric
\cp ness} of $(K,\ge)$ \wrt the metric $d_j(x,y):=j(|x-y|_K)$, $x,y\in K$, with $j:K\to\R$ in\jc\ order- and ring-\hm sm \itd\ in Theorem \rf{t4.35},
and $|\cdot|_K :K\to K$ the \av\ in~$K$ (see \E\df\ \rf{d5.22}). Hence $j\circ |\cdot|:K\to\R_{\ge0}$ is an $\R$-valued \vl\ on~$K$. We leave it to
the reader to verify these \ass s.

Inspired by \cite[pp.~47--50]{Choquet} we use the notion of ``nested \il s''.\index{nested \il s}
A~nest of \il s in an \Ar\ \of\ $(K,\ge)$ is a \sq\ of \il s $I_n:=[a_n,b_n]$, $a_n,b_n\in K$, $a_n<b_n$, $n\in\N$, \sf ying (i)~$I_m\subseteq
I_n$ whenever $m>n$, $m,n\in\N$, and\break (ii)~$b_n-a_n\da 0$.

\bth8.130
Let $(K,\ge)$ be an \Ar\ field, and let $d_j$ be the metric on~$K$ defined above. Then the \fw\ assertions are \ev t\dw

\hph1,, The field $(K,d_j)$ $($metrically$)$ is \cp.

\hph2,, For every nest of \il s $\zb In\N$ \te s \ooo $x\in K$ \st $\{x\}= \Bca_{n\in\N}I_n$.

\hph3,, The field $(K,\ge)$ is order- and ring-\is c to the field $(\R,\ge)$, \Ip $(K,\ge)$ is order-\cp.
\eth

\proof \

``(1) \ti{holds iff $(3)$ holds}'': \ We first recall that $j\circ |\cdot|_K$ is a real-valued \vl, hence $d_j$ is a metric on~$K$. We next observe
that $j(|x|_K)=|j(x)|$ \fa $x\in K$ with $|\cdot|$ the \av\ in~$\R$, since $j(x)\ge0$ (resp.~$\le0$) iff $x\ge0$ ($x\le0$) and $j(-x)=-j(x)$
($-x\in K$!). \E\Tf $d_j(x,y)=|j(x)-j(y)|$, $x,y\in K$. Note that $j(K)$ is a subfield of the field~$K$ (use Lemma \rf{l3.23}). \If that $j:K\to j(K)$
is an order- and ring-\is sm. We denote by~$j\Inv$ the inverse of~$j$. Note that $j\Inv$ also is an order- and ring-\is sm. Hence we obtain for
$\a,\b\in j(K)\sbs\R$ and $x:=j\Inv(\a)$, $y:=j\Inv(\b)$, $d_j(x,y)=|\a-\b|=d(\a,\b)$ with $d$~the usual metric on~$\R$. \If that $j\Inv$ is a sur\jc\
isometry from $(j(K),d)$ onto $(K,d_j)$, hence $j$~is a sur\jc\ isometry from $(K,d_j)$ onto~$(j(K),d)$. We next observe that if $\zb xn\N$ is a \Cs\
in~$(K,d_j)$, then $\{j(x_n)\}_{n\in\N}$ is a \Cs\ in $(j(K),d)$. Indeed, if \fe $\ve\in\R_{>0}$ \te s an $N\in\N$ \st $d_j(x_n,x_m)<\ve$ \fa $m,n
\ge N$, then $|j(x_n)-j(x_m)|=d_j(x_n,x_m)<\ve$ for $n,m\ge N$. Likewise, one shows that if $\zb\a n\N$ is a \Cs\ in $(j(K),d)$, then $\{j\Inv
(\a_n)\}_{n\in\N}$ is a \Cs\ in~$(K,d_j)$. We claim that if $(K,d_j)$ is \cp, then so is $(j(K),d)$.

Indeed, let $\zb\a n\N$ be a \Cs\ in $(j(K),d)$. Then $\{j\Inv(\a_n)\}_{n\in\N}$ is a \Cs\ in $(K,d_j)$, which is metrically \cp. Hence, \te s $\ov x
\in K$ \st $x_n:=j\Inv(\a_n)$, $n\in\N$, \cg s to~$\ov x$, that is, $d_j(x_n,\ov x)\nad d\to 0$. Set $\ov\a:=j(\ov x)\in j(K)\sbs \R$. Then
$d(\a_n,\ov\a) = d_j(x_n,\ov x)\nad d\to 0$, hence $\a_n\nad d\to \ov\a\in j(K)$. \If that $(j(K),d)$ is \cp. Likewise, one shows that if $(j(K),d)$
is \cp, then $(K,d_j)$ is \cp, hence $(j(K),d)$ is not \cp\ iff $(K,d_j)$ is not \cp. In view of Theorem \rf{t8.127}, $(\R,d)$ is \cp. \E\oh if
$F$~is a proper subfield of~$\R$, that is, $F\ne\R$, then $(F,d)$ is \ti{not\/} \cp. Indeed, suppose for \cd ion that $F$~is \cp\ and let $\a\in\R\sm
F$. From the first part of Section~\ref{ass.6} we know that \te s a \sq\ $\zb rn\N$ \st $r_n\in\Q$ and $d(r_n,\a)\nad d\to 0$. \E\Tf $\zb rn\N$ is
a \Cs\ in $\Q\sbs F$, which does not \cg\ in~$(F,d)$. A~\cd ion. \If that $(j(K),d)$ is \cp\ iff $K=\R$ iff $j$~is sur\jc. \E\Tf if $(K,d_j)$ is \cp,
then $K$~is order- and ring-\is c to $(K,\ge)$, and $j$~is the only order- and ring-\is sm in view of \E\Pr\ \rf{p5.14}\,(iii). \E\oh if $(K,d_j)$
is not \cp\ and $j:K\to j(K)$ is an arbitrary in\jc\ order- and ring \hm sm, then $j(K)\ne\R$, hence $(K,d_j)$ is \ti{not\/} order- and ring-\is c
to $(\R,\ge)$. This \cp s the proof of $(1)\iff(3)$.

``(1)$\imp$(2)'': \ Let $\{[a_n,b_n]\}_{n\in\N}$ be a nest of \il s on $(K,\ge)$. We assume~(1). Hence from the first part of the proof, we infer that
\te s an \ti{order-\is sm} $\vf:(K,\ge) \to (\R,\ge)$. One verifies that $\{[\vf(a_n),\vf(b_n)]\}_{n\in\N}$ is a nest of \il s on~${(\R,\ge)}$. \E\Tf
if we can prove that $\Bca_{n\in\N}[\vf(a_n),\vf(b_n)]=\{y\}$ \fs $y\in\R$, then we obtain $\Bca_{n\in\N}[a_n,b_n] = \Bca_{n\in\N}[\vf\Inv(\vf(a_n)),
\vf\Inv(\vf(b_n))] = \Bca_{n\in\N}\vf\Inv([\vf(a_n),\vf(b_n)]) =\break \vf\Inv\bigl(\Bca_{n\in\N}[\vf(a_n),\vf(b_n)] = \{\vf\Inv(y)\}$, which is~(2).
\E\Tf it suffices to prove that (1) implies~(2) for $(K,\ge):=(\R,\ge)$.

We have for $n<m$, $n,m\in\N$:
\beq gw
\ga
{}[a_m,b_m] \sbs [a_n,b_n]\\
\begin{picture}(0,-10)(0,0)
\put(-100,0){\line(1,0){200}}\put(-93,-3){$\bullet$}\put(-95,-10){$a_0$}
\put(-73,-3){$\bullet$}\put(-38,-3){$\bullet$}\put(12,-3){$\bullet$}\put(56,-3){$\bullet$}\put(89,-3){$\bullet$}
\put(-75,-10){$a_n$}\put(-40,-10){$a_m$}\put(10,-10){$b_m$}\put(54,-10){$b_n$}\put(87,-10){$b_0$}
\end{picture}
\ega
\e
\ssk

``$a_n\le a_m$'': \ otherwise $a_m<a_n<b_n$, hence $a_m\notin[a_n,b_n]$ \cd ing \er{gw}. Likewise we have ``$b_n\ge b_m$''. Hence $\zb an\N$ (resp.\
$\zb bn\N$) is an in\cre\ (resp.\ de\cre) \sq\ in $(\R,\ge)$.

``$a_n\le b_0$'': \ follows from $a_n<b_n\le b_0$. Likewise we have ``$a_0\le b_n$''.
In view of the \cp ness of $(\R,\ge)$ we obtain
$$
a_n\ua \ov a:=\sup_{n\ge0}a_n, \q b_n\da \ov b:=\inf_{n\ge0}b_n.
$$

``$a_n\le b_m$ \ti{if\/} $n\ge m$, $n,m\in\N$'': \ $a_n\in[a_n,b_n]\sbs [a_m,b_m]$ if $n\ge m$. Hence $a_n\in[a_m,b_m]$, and $a_n\le b_m$.

``$\ov a\le b_m$ \ti{if\/} $m\in\N$'': \ $\ov a=\supl_{n\ge0}a_n = \supl_{n\ge m}a_n\le b_m$ since $\supl_{n\ge m}a_n$ is the least \ub\ for the set
$\{a_n\in\R: n\ge m\}$ \fa $m\in\N$.

``$\ov a\le \ov b$'': \ since $\ov b=\infl_{n\ge 0}b_n = \infl_{n\ge m}b_n$.

``$[\ov a,\ov b]\sbs [a_n,b_n]$ \ti{\fa}$n\in\N$'': \ Let $x\in[\ov a,\ov b]$. Then $x\ge\ov a$ and $x\le \ov b$, hence $x\in[a_n,b_n]$ \fa $n\in\N$.

``$[\ov a,\ov b]\sbs \Bca_{n\in\N}[a_n,b_n]$: \ by \df\ of $\Bca_{n\in\N}[a_n,b_n]$.

``$\ov a=\ov b$'': \ Since $(b_n-a_n)\da 0$, \fe $\ve\in\R_{>0}$ \te s $N\in\N$ \st $0\le b_n-a_n\le\ve$ \fa $n\in\N$. Since $[\ov a,\ov b]\sbs
\Bca_{n\in\N}[a_n,b_n]$, we have $|\ov b-\ov a|\le\ve$ \fa $\ve\in\R_{>0}$. Hence $|\ov b-\ov a|=0$, and $\ov a=\ov b$.

``(2)$\imp$(3)'': \ Let $A$ be an infinite subset of~$K$ which is \ba. Suppose, for \cd ion, that the set~$A$ does not have a least \ub. We construct
\rc vely a \sq\ of nested \il s $[a_n,b_n]$, $n\in\N$, \st $a_n$ is not an \ub\ for~$A$, and $b_n$~is an \ub\ for~$A$. Since $A$ is not empty and has
no least \ub, \te s $a_0\in A$ which is not an \ub\ for~$A$. Let $b_0$~be an \ub\ for~$A$. Such $b_0$ exists since the set~$A$ is \ba. Then $a_0<b_0$.
Suppose that \te\ $n\in\N$ and $a_n,b_n\in K$ \sf ying $a_n<b_n$, $a_n$~is not an \ub\ for~$A$, and $b_n$~is an \ub\ for~$A$. Set $c_n:=\frac12
(a_n+b_n)$. \Mo set $a_{n+1}:=a_n$ and $b_{n+1}:=c_n$ if $c_n$~is an \ub\ for~$A$, $a_{n+1}:=c_n$ and $b_{n+1}:=b_n$ otherwise. We claim that $a_{n+1}
<b_{n+1}$, $a_n\le a_{n+1}$, $b_n\ge b_{n+1}$, $a_{n+1}$ is not an \ub\ for~$A$, $b_{n+1}$~is an \ub\ for~$A$ and $b_{n+1}-a_{n+1}=\frac12(b_n-a_n)$.
Indeed:

\noi --- case ``$c_n$ is an \ub\ for~$A$'': $a_{n+1}=a_n<c_n=b_{n+1}$; $a_{n+1}=a_n\ge a_n$; $b_{n+1}=c_n<b_n$; and $b_{n+1}-a_{n+1}=c_n-a_n=\frac12
(b_n-a_n)$. \Mo $b_{n+1}=c_n$ is an \ub\ for~$A$, and $a_{n+1}=a_n$ is not an \ub\ for~$A$.

\noi --- case ``$c_n$ is not an \ub\ for~$A$'': $a_{n+1}=c_n<b_n=b_{n+1}$; $a_{n+1}=c_n>a_n$; $b_{n+1}=b_n\le b_n$; and $b_{n+1}-a_{n+1}=b_n-c_n=
\frac12(b_n-a_n)$. \Mo $b_{n+1}=b_n$ is an \ub\ for~$A$, and $a_{n+1}=c_n$ is not an \ub\ for~$A$.

Using \In\ on~$n$ we find that the \sq\ $\zb an\N$ is in\cre, the \sq\ $\zb bn\N$ is de\cre, $b_{n+1}-a_{n+1}< b_n-a_n$ \fa $n\in\N$ and
$b_n-a_n=2^{-n} (b_0-a_0) \to 0$ since $K$ is \Ar. By~(2) \te s $\ov c\in\R$ \st $\{\ov c\} = \Bca_{n\in\N}[a_n,b_n]$. We claim that $\ov c$ is an
\ub\ for~$A$. If not, \te s $x\in A$ \st $\ov c<x$. However, $x\le b_n$ \fa $n\in\N$ since $b_n$~is an \ub\ for~$A$ \fa $n\in\N$. \Mo $\ov
c\in[a_n,b_n]$ \fa $n\in\N$. Hence $a_n\le \ov c<x\le b_n$, and $0<x-\ov c\le b_n-a_n\da0$, which is impossible. \E\Tf $\ov c$ is an \ub\ for~$A$.

Finally,we show that $\ov c$ is the least \ub\ for~$A$. If not, \te s $y\in A$ \st $y<\ov c$ and $y$~is an \ub\ for~$A$. \Mo since $a_n$ is not an
\ub\ for~$A$, \fe $n\in\N$ \te s $y_n\in A$ \st $a_n<y_n$. Hence $a_n<y_n\le y$ \fa $n\in\N$. \E\Tf $\supl_{n\ge0}a_n \le y$ and $\ov c=\supl_{n\ge0}
a_n\le y$. A~\cd ion, since $y<\ov c$. Thus $\ov c$~is the least \ub\ for~$A$ \cd ing the \as\ there is no least \ub\ for~$A$. \csq, $\sup A$
exists.
\endproof


We now give an example which shows that the notion of \Cs\ is \ti{not\/} topological. Let $X:=\oz 0,1 $ in~$\R$, let $Y:=\R_{\ge1}$, and let $f:X
\to Y$ be defined by $f(x):=\frac1x$, and $g:Y\to X$ be defined by $g(y):=\frac1y$. Then $g\circ f=\id_X$ and $f\circ g=\id_Y$, hence $f$~is bi\jc\
and $f\Inv=g$. Both $f$~and~$g$ are \sql\ \ctn\ in view of Lemma \rf{l7.6} \er{7.8}, and \ctn\ by Theorem \rf{t8.89} and \E\df\ \rf{d8.93}. Set
$a_n:=\frac1{n+1}$ for $n\in\N$. Then $\lim\limits_{n\ge0}a_n=0$ in~$\R$. Hence $\zb an\N$ is a \Cs\ in~$(\R,d)$, and also in~$(X,d)$. However, the
\sq\ $f(a_n)=n+1$ is not a \Cs\ in~$\R_{\ge1}$ since it is not \ba. \E\Tf $f$~is a \he sm which does not preserve \Cs s. It turns out that uniformly
\ctn\ \f s preserve \Cs s.

\bpr8.129
Let $(X,d)$, $(X',d')$ be \ms s and let $f:X\to X'$ be uniformly \ctn. If $\zb xn\N$ is a \Cs\ in~$X$ then so is $\{f(x_n)\}_{n\in \N}$ in~$X'$.
\epr

\proof
For convenience of the reader we first consider the case where the map~$f$ is \Lc. Let $\zb xn\N$ be a \Cs\ in~$(X,d)$. To prove: \fe $\ve\in\R_{>0}$
\te s an $N\in\N$ \st $d(f(x_n),f(x_m))<\ve$ whenever $n,m\ge N$. Let $M\in\R_{>0}$ be \st $d'(f(x),f(y))\le Md(x,y)$ \fa $x,y\in X$. We know that
\fe $\d\in\R_{>0}$ \te s an $N\in\N$ \st $d(x_n,x_m)<\d$ whenever $n,m\ge N$. Hence given $\ve\in \R_{>0}$, setting $\d:=\frac\ve M$, we find that
\te s an $N\in\N$ \st $d(x_n,x_m)< \frac\ve M$ whenever $n,m\ge N$. Thus we obtain $d'(f(x_n),f(x_m))\le Md(x_n,x_m)<\ve$ whenever $n,m\ge N$. \If
that $\{f(x_n)\}_{n\in\N}$ is a \Cs\ in $(X',d')$.

We now consider the general case. In view of the uniform continuity of the map~$f$, \fe $\ve\in\R_{>0}$ \te s $\g\in\R_{>0}$ \st $d'(f(x),f(y))<\ve$
whenever $d(x,y)<\g$. Let $\ve\in\R_{>0}$ and let $N\in\N$ be \st $d(x_n,x_m)<\g$ whenever $n,m>N$. Thus we obtain $d'(f(x_n),f(x_m))<\ve$  whenever
$n,m>N$. Hence $\{f(x_n)\}_{n\in\N}$ is a \Cs.
\endproof


\bpr8.130a
Let $(X,d)$ be a \cp\ \ms, let $(X',d')$ be  a \ms\ and let $f:X\to X'$ be a \hm sm. If $f\Inv$ is \emph{uniformly} \ctn, then $X'$~is \cp.
\epr

\proof
Let $\zb{x'}n\N$ be any \Cs\ in $(X',d')$. Then $\{f\Inv(x'_n)\}_{n\in\N}$ is a \Cs\ in $(X,d)$ in view of \E\Pr\ \rf{p8.129}. Since $(X,d)$
is \cp, \te s $y\in X$ \st $f\Inv(x'_n)\nad d\to y\in X$. \E\Tf $x_n'= f(f\Inv(x'_n)) \nad{d'}\to f(y)$. Hence $(X',d')$ is \cp.
\endproof

We now present some basic examples of Polish spaces.

\blm8.132
Let $N\in\Na$ and let $\R^N$ denote the $N$-\dm al real \vs\ \itd in Example {\rm5.1.2}. Let~$d_2$ denote the Euclidean metric on~$\R^N$ \itd in
\er{7.22} with $n:=N$. Then the $N$-\dm al Euclidean space $(\R^N,d_2)$ is \cp.
\elm

\proof
Let $\zb\bX n\N$ be a \Cs\ in $(\R^N,d_2)$, and let $\psi_k:\R^N \to\R$, $k\in[1,N]$, denote the \co\ map defined in (5.1.5). \E\fa $\bY,\bZ\in\R^N$
we have $\psi_k(\bY)-\psi_k(\bZ) = \psi_k(\bY-\bZ)$, hence $(\psi_k(\bY)-\psi_k(\bZ))^2 = (\psi_k(\bY-\bZ))^2 \le \suml_{i=1}^N (\psi_i(\bY-\bZ))^2
= (d_2(\bY,\bZ))^2$. \If that $|\psi_k(\bY)-\psi_k(\bZ)| \le d_2(\bY,\bZ)$ for $\bY,\bZ\in\R^N$, hence the maps $\psi_k:(\R^N,d_2)\to(\R,d)$ are \Lc,
hence uniformly \ctn. From \E\Pr\ \rf{p8.129} we infer that the \sq s $\{\psi_k(\bX_n)\}_{n\in\N}$, $k\in[1,N]$, are \Cs s. By Theorem \rf{t8.127}
\fa $k\in[1,n]$ \te\ $z_k\in\R$ \st $\psi_k(\bX_n)\to z_k$. Set $\bZ:=\suml_{i=1}^N z_k\ve_i\in \R^N$ with $\ve_i$ defined in~(5.1.15). Since
\fa $k\in[1,n]$ $\psi_k(\bX_n)\to\psi_k(\bZ)$, we have $\bX_n\imp \bZ$ by~\er{7.29}. From \E\Pr\ \rf{p7.20} we obtain $d_2(\bX_n,\bZ)\to0$, hence
$(\R^N,d_2)$ is a \ti{\cp} \ms.
\endproof

\bex8.133
Let $\rho\in\R_{>0}$, $\rho<1$ and let $D:\R^N\t \R^N \to\R$, $N\in\Na$, be defined by
\beq8.86
D(\bX,\bY) := \sum_{i=1}^N |\psi_i(\bX)-\psi_i(\bY)|^\rho, \q \bX,\bY\in\R^N.
\e
Show that $D$ is a metric on $\R^N$ and that $(\R^N,D)$ is \cp. (Hint: use Theorem \rf{t7.51}.)
\eex

Proceeding as in the proof of Lemma \rf{l8.132} we might show that the space $\R^N$ endowed with the metrics defined in \er{7.23} and \er{7.24}
is \cp. Indeed, this is the case for the important class of metrics on~$\R^N$ ``induced'' by a norm.

\bdf8.136
Let $X$ be a real \vs. A~\tb{norm} on~$X$ is a map $\|\cdot\|$ from~$X$ into $\R_{\ge0}$ \sf ying the \fw\ axioms. \E\fa $u,v\in X$ and all $\la\in\R$
we have:\index{norm}
\ben
\def\labelenumi{\rm(N\theenumi)}
\sdim{(N1)}
\item $\|u\|\ge0$ and $\|u\|=0$ iff $u=0$;
\item $\|\la u\| = |\la|\,\|u\|$;
\item $\|u+v\| \le \|u\|+\|v\|$.
\een
A real \vs\ is said to be a \ti{real normed space} if it is endowed with a norm.
\edf

\brm8.137 \

\hph i,i, If $X:=\R$ and $\|\cdot\|$ is a norm on the \vs\ $(\R,\R)$, then $\|\cdot\|$ is a \vl\ (see \E\df\ \rf{d5.23}) on~$\R$.

\hph ii,, If $\vf(u):=d_2(0,u) = \bigl(\suml_{i=1}^n (u_i)^2\bigr)^{1/2}$, then (N1), (N2) are clearly \sf ied, and (N3) follows from \er{7.21}.
Indeed, $\vf(u+v)=d_2(0,u+v) \nde7.54 = d_2(-v,(u+v)+(-v)) = d_2(-v,u) \nad{\rm(M3)} \le d_2(-v,0) + d_2(0,u) \nde7.54 = d_2(0,v) + d_2(0,u)
= \vf(v)+\vf(u) = \vf(u)+\vf(v)$.
\erm

\bex8.138 \

\hph i,i, Prove (N1)--(N3) for the metrics defined in \er{7.23} and \er{7.24}.

\hph ii,, Show that if $\|\cdot\|$ is a norm on a real \vs~$X$, then $d(u,v)$, $u,v\in X$, defined by
\beq8.87
d(u,v):=\|u-v\|
\e
is a metric on $X$ (the metric ``induced'' by the norm~$\|\cdot\|$), and that the \fw\ holds:
\beq8.88
d(0,u) = \|u\|, \q u\in X.
\e
\eex

\bnt8.139
For $\bX\in\R^N$, $N\in\Na$, we set
\beq8.89
\|\bX\|_2 := d_2(0,\bX).
\e
\ent

The norm $\|\cdot\|_2$ on $\R^N$ is called the \ti{Euclidean norm} on~$\R^N$, and $(\R^N,\|\cdot\|_2)$ is called the $N$-\ti{\dm al Euclidean normed
space}.\index{Euclidean norm}

\blm8.140
Let $(\R^N,\|\cdot\|_2)$ be the $N$-\dm al Euclidean normed space. Then \fa $R\in\R_{>0}$, the sets
\beq8.90a
\{\bX\in\R^N: \|\bX\|_2\le R\} \q\hbox{and}\q \{\bX\in\R^N: \|\bX\|_2= R\}
\e
endowed with the metric~$d_2$ are \sql\ compact.
\elm

\proof \

\ti{Step} 1: \ Let $[-R,R]^N$ denote the set $\{\bX\in\R^N: -R\le \bX(i)\le R \hbox{ \fa}i\in[1,N]\}$. We claim that $[-R,R]^N$ is \sql\ compact
in~$(\R^N,d_2)$. We proceed by \In\ on $N\in\Na$. 
Let $N:=1$ and let $\zb xn\N$ in $[-R,R]$. By Lemma \rf{l8.15}\,(ii) \te\ $\ov x\in\R$ and a sub\sq\ $\{x_{n_k}\}_{k\in\N}$ \st $x_{n_k}\to \ov x$.
Define \sq s $y_k:=-R$ and $z_k:=R$ \fa $k\in\N$. By \E\Pr\ \rf{p8.10}\,(ii), we have $\liml_{k\ge0}x_{n_k} = \limu_{k\ge0}x_{n_k} = \ov x$ and
$\liml_{k\ge0}y_k = -R$, $\limu_{k\ge0}z_k= R$. We find $-R = \liml_{k\ge0}y_k \nde8.11 \le \liml_{k\ge0} x_{n_k} = \ov x = \limu_{k\ge0} x_{n_k}
\nde8.11 \le \limu_{k\ge0}y_k = R$. Hence $|\ov x|\le R$. Since the \sq\ $\zb xn\N$ in $[-R,R]$ is arbitrary, the \ms\ $[-R,R]$
endowed with the metric $d_2(x,y) = ((x-y)^2)^{1/2} = |x-y|$, $x,y\in [-R,R]$ is \sql\ compact.
Suppose that the claim holds for~$N$,
and let $\zb \bX n\N$ be a \sq\ in~$\R^{N+1}$ \sf ying $\bX_n\in[-R,R]^{N+1}$. \E\te\ a sub\sq\ $\{\bX_{\si(k)}\}_{k\in\N}$ and $\bY\in[-R,R]^N$ \st
$\bX_{\si(k)}(i) \to \bY(i)$ \fa $i\in[1,N]$. Note that 
$-R \le
\bX_{\si(k)}(N+1)\le R$ \fa $k\in\N$. \E\Tf \te\ a sub\sq\ $\{\bX_{\rho(\si(k))}(N+1)\}_{k\in\N}$ and $\a\in[-R,R]$ \st $\bX_{\rho(\si(k))}(i)\to
\bY(i)$ \fa $i\in[1,N]$, and $\bX_{\rho(\si(k))}(N+1)\to\a$. Setting $\bZ(N+1):=\a \in[-R,R]$ and $\bZ(i):=\bY(i)$ for $i\in[1,N]$, we obtain
$\bZ\in[-R,R]^{N+1}$ and $\bX_{\rho(\si(k))}(i)\to \bZ(i)$ \fa $i\in[1,N+1]$. By \E\Pr\ \rf{p7.20} we have $d_2(\bX_{\rho(\si(k))},\bZ)\to 0$. Hence
$([-R,R]^{N+1},d_2)$ is \sql\ compact, which proves the claim.

\ti{Step} 2: \ We claim that the sets defined in \er{8.90a} are closed subsets of the \ms\ $([-1,1]^N,d_2)$. Indeed, we have $\suml_{i=1}^N (\bZ(i))^2
= \|\bZ\|_2^2 \le R^2$ \fa $\bZ\in\R^N$ \st $\|\bZ\|_2 \le R$ (since $\|\bZ\|_2^2 \le\|\bZ\|_2R \le R^2$). \If that $(\bZ(i))^2\le R^2$, hence
$|\bZ(i)|\le R$ \fa $i\in[1,N]$. \E\Tf $\bZ\in[-R,R]^N$. \Mo if a \sq\ $\zb \bZ n\N$ in~$\R^N$ \sf ies $(d_2(\bz,\bZ_n))^2 = R^2$ \fa $n\in\N$, and
\te s $\bw\in\R^N$ \st $d_2(\bZ_n,\bw)\to0$, we claim that $d_2(\bz,\bw)=R$, that is, $\|\bw\|_2\nde8.89 = R$. For the convenience of the reader, we
recall that if $(M,d)$ is a \ms, then as a direct con\sq\ of (M2),~(M3) we have
\beq8.90
|d(x,y)-d(x,z)| \le d(y,z) \qh{\fa $x,y,z\in M$.}
\e
\E\Tf $|d_2(\bz,\bw)-R| = |d_2(\bz,\bw) - d_2(\bz,\bZ_n)| \le d_2(\bw,\bZ_n) \nad{\rm(M2)} = d_2(\bZ_n,\bw)\to 0$. \If that $|d_2(\bz,\bw)-R| < \ve$
\fa $\ve\in\R_{>0}$. Thus $|d_2(\bz,\bw)-R|=0$, and $\|\bw\|_2 = d_2(\bz,\bw)=R$.
\csq, the set $\{\bX\in\R^N: \|\bX\|_2=1\}$ is a \sql\ closed subset of $([-R,R],d_2)$.

We already proved that the set $\{\bZ\in\R^N: \|\bZ\|_2\le R\}$ is a subset of the
set $[-R,R]^N$. We now show that this set is s-closed in $([-R,R]^N,d_2)$. Let $\zb \bZ n\N$ be a \sq\ in $(\R^N,d_2)$ \sf ying $\|\bZ_n\|\le R$ \fa
$n\in\N$ and $\|\bZ_n-\bw\|_2\to0$ \fs $\bw\in\R^N$. We have to show that $\|\bw\|_2\le R$. From \er{8.90} with $x:=0$, we infer \fa $\bY,\bZ\in\R^N$
$\bigl|\|\bY\|_2-\|\bZ\|_2\bigr| = |d_2(\bz,\bY)-d_2(\bz,\bZ)| \le d_2(\bY,\bZ) = \|\bY-\bZ\|_2$. \If that the map $\bZ\mt \|\bZ\|_2$ from
$(\R^N,d_2)$ into $(\R,d)$ is \Lc, hence \sql\ \ctn. We obtain $\|\bw\|_2 = \lim\limits_{n\ge0}\|\bZ_n\|_2 \nad*= \limu_{n\ge0}\|\bZ_n\|_2 \nad{**}\le
\limu_{n\ge0}R=R$. In $\nad*=$ we used \E\Pr\ \rf{p8.10} and in $\nad{**}\le$ \E\Pr\ \rf{p8.12}. Note that the \cvg\ \sq\ $\{\|\bZ_n\|_2\}_{n\in\N}$
is bounded in~$\R$. \If that the set $\{\bX\in\R^N: \|\bX\|_2\le R\}$ is a s-closed subset of the set $([-R,R]^N,d_2)$. We complete the proof of
Lemma \rf{l8.140} by applying Lemma \rf{l8.32}. Since the proof of this lemma was left as an exercise, we do it now. Let $(X,d)$ be a \sql\ compact
\ms\ and let $A$~be a nonempty \sql\ closed subset of~$A$. Then \te\ a sub\sq\ $\{x_{n_k}\}_{k\in\N}$ in~$A$ and $x\in X$ \st $x_{n_k}\nad d\to x$.
Since $A$~is \sql\ closed in $(X,d)$ we have $x\in A$.
\endproof

\bth8.141
Under the \as\ of Lemma \rf{l8.30}, \te s $\check c\in\R^N$ $($resp.\ $\hat c\in\R^N)$ \st
\beq8.91
f(\check c)=\inf f(X) \q(\hbox{resp.} \ f(\hat c)=\sup f(X)).
\e
\eth

\proof
In view of Lemma \rf{l8.30} $f(X)$ is \bb\ (resp.\ above). We consider the first case. Set $\check\a:=\inf(X)\in\R$ which exists by Theorem \rf{t5.1}
and Corollary \rf{c5.12}. By Lemma \rf{l6.28} \te s a \sq\ $\zb xn\N$ in~$X$ \st $f(x_n)\to\check\a$ in~$\R$. Since $(X,d)$ is \sql\ compact, \te\
$\check c \in X$ and a sub\sq\ $\{x_{n_k}\}_{k\in\N}$ \st $x_{n_k}\nad d\to\check c$ in~$X$. In view of the sequential lower-semicontinuity of~$f$
we obtain $f(\check c)\le\liml f(x_{n_k})\nad*= \lim\limits_{k\ge0}f(x_{n_k}) = \inf f(X)$. In~$\nad*=$ we used \E\Pr\ \rf{p8.10} noting that the \sq\
$\{f(x_n)\}_{n\in\N}$ is bounded in~$\R$. \E\Tf $f(\check c)\le \inf f(X)$. Since $\inf X\le f(\check c)$, we obtain \et y. The proof of the second
case is similar, \Tf it is omitted.
\endproof

\blm8.142
Let $N\in\Na$ and let $\R^N$ be as in Lemma \rf{l8.132}. Let $\|\cdot\|_2$ be the Euclidean norm and let $\|\cdot\|$ be a norm on~$\R^N$. Then \te\
$c,C\in\R_{>0}$ \st
\beq8.92
c\|\bX\|_2 \le \|\bX\| \le C\|\bX\|_2 \qh{\fa $\bX\in\R^N$.}
\e
\elm

\proof
Let $\bE_i\in\R^N$ be defined by $\bE_i(j):=\d_{ij}$, $i,j\in[1,N]$. Let $\bX=\suml_{i=1}^n \bX(i)\bE_i\in\R^N$ with $\bX(i)\in\R$, $i\in[1,N]$. Then
we have $\|\bX\| = \bigl\|\suml_{i=1}^N \bX(i)\bE_i\bigr\| \nad*\le \suml_{i=1}^N \|\bX(i)\bE_i\| = \suml_{i=1}^N |\bX(i)|\,\|\bE_i\| \nad{**}\le
\bigl(\suml_{i=1}^N |\bX(i)|^2\bigr)^{1/2} \bigl(\suml_{i=1}^N \|\bE_i\|^2\bigr)^{1/2} = C\|\bX\|_2$ with $C:=\bigl(\suml_{i=1}^N \|\bE_i\|^2\bigr)
^{1/2}$. In~$\nad*\le$ we used (N3) and \In, in~$\nad{**}\le$ we used \er{7.20} and Remark \rf{r7.13}\,(ii). This \cp s the proof of the second
in\et y in \er{8.92}. \E\fa $\bX,\bY\in\R^N$ we obtain $\|\bX-\bY\| \le C\|\bX-\bY\|_2$, hence $\bigl|\|\bX\|-\|\bY\|\bigr| \nad*\le \|\bX-\bY\| \le
C\|\bX-\bY\|_2$, where in~$\nad*\le$ we used (N3). Indeed, $\|\bX\| = \|(\bX-\bY)+\bY\| \le \|\bX-\bY\|+\|\bY\|$, hence $\|\bX\| - \|\bY\| \le
\|\bX-\bY\|$ \fa $\bX,\bY\in\R^N$. Interchanging $\bX$ and $\bY$ leads to $\|\bY\|-\|\bX\| \le \|\bY-\bX\| = \|-(\bX-\bY)\| = \|\bX-\bY\|$. From
$\bigl|\|\bX\|-\|\bY\|\bigr| \le C\|\bX-\bY\|_2$, $\bX,\bY\in\R^N$, we infer that the map $f:(\R^N,d_2) \to (\R_{\ge0},d)$ defined by $f(\bX):=
\|\bX\|$, is \Lc, hence \sql\ \ctn. Clearly, the \rt ion of~$f$ to the set $\Si:=\{\bX\in\R^N: \|\bX\|_2=1\}$ is also s-\ctn. In view of Lemma
\rf{l8.140} with $R:=1$ and Theorem \rf{t8.141}, \te s $\bY\in\Si$ \st $f(\bY)=\inf f(\Si)$. Setting $c:=\inf f(\Si)$, we find $c=f(\bY)=\|\bY\|
\ge0$. We have $\|\bY\|>0$, otherwise $\|\bY\|=0$, hence $\bY=\bz$ by~(N1) \cd ing $\bY\in\Si$ since $\|\bY\|_2=1 \ne0=\|\bz\|_2$. \csq, $c>0$.
Clearly, \er{8.92} holds for $\bX:=\bz$. Suppose $\bX\ne\bz$, then $\|\bX\|\ne0$ and $\bigl\|\|\bX\|\mo\bX\bigr\| = \|\bX\|\mo\|\bX\|=1$. Hence
$\|\bX\|\mo\bX\in\Si$. Thus $c = f(\bY) \le f(\bZ)\le \|\bZ\|$ \fa $\bZ\in\Si$. \E\Ip $c\le \bigl\|\|\bX\|_2\mo \bX\bigr\|$, hence $c\le\|\bX\|_2\mo
\|\bX\|$, and $c\|\bX\|_2\le \|\bX\|$ by \er{3.10}. This \cp s the proof of the first in\et y in~\er{8.92}.
\endproof

\brm8.143
If $\|\cdot\|'$ is another norm on $\R^N$, we have \er{8.92} as well as $c'\|\bX\|_2\le \|\bX\|\le C'\|\bX\|_2$ \fa $\bX\in\R^N$. \csq, applying Lemma
\rf{l8.142} with $\|\cdot\|'$ instead of~$\|\cdot\|$, we obtain
\beq8.93
C\mo c'\|\bX\| \le \|\bX\|' \le C'c\mo \|\bX\|, \q \bX\in\R^N,
\e
where $c'$ and $C'$ are the \crs\ ``constants'' in \er{8.92} for $\|\cdot\|'$.
\erm

\bdf8.144
Let $X$ be a real \vs, and let $\|\cdot\|$ and $\|\cdot\|'$ be norms on~$X$. Then the norms $\|\cdot\|$ and $\|\cdot\|'$ are said to be \ti{\ev t\/}
if \te\ $c,C\in\R_{>0}$ \st $c\|\bX\| \le \|\bX\|' \le C\|\bX\|$ \fa $\bX\in X$.\index{equivalent norms}
\edf

\brm8.145 \

\hph i,i, On the \vs\ $\R^N$ all norms are \ev t.

\hph ii,, All preceding results on the normed real \vs\ $(\R^N,\|\cdot\|_2)$ can be ``transported'' to an $N$-\dm al normed real \vs\ $(X,\|\cdot\|)$
by means of the linear \is sm~$L$ from $\R^N$ into~$X$ \itd in \E\df\ 5.1.13. Indeed, we have $L\bX = L\bigl(\suml_{i=1}^N \bX(i)\bE_i\bigr) =
\suml_{i=1}^N L(\bX(i)\bE_i) = \suml_{i=1}^N \bX(i)L\bE_i$ \fa $\bX\in\R^N$. Hence $\|L\bX\| = \bigl\|\suml_{i=1}^N \bX(i)L\bE_i\bigr\| \le
\suml_{i=1}^N \|\bX(i)L\bE_i\| = \suml_{i=1}^N |\bX(i)|\,\|L\bE_i\| \le \bigl(\suml_{i=1}^N |\bX(i)|^2\bigr)^{1/2}\bigl(\suml_{i=1}^N \|L\bE_i\|^2
\bigr)^{1/2}$. Setting $C:=\bigl(\suml_{i=1}^N \|L\bE_i\|^2\bigr)^{1/2}\in\R_{\ge0}$ we obtain $\|L\bX\| \le C\|\bX\|_2$, $\bX\in\R^N$, and using the
linearity of~$L$,
\beq8.94
\|L\bX-L\bY\| \le C\|\bX-\bY\|_2
\e
\If that the map $L$ is \Lc\ from $(\R^N,d_2)$ onto $(X,d)$ with
\beq8.95
d(u,v):=\|u-v\|, \q u,v\in X.
\e
Observe that the map $\vf:\R^N\to\R_{\ge0}$ defined by $\vf(\bX):=\|L\bX\|$ is a norm on~$\R^N$. Indeed, $\vf(\bX)\ge0$ by~(N1) for~$\|\cdot\|$ and
$\vf(\bX)=0$ iff $\|L\bX\|=0$ iff $L\bX=0$ iff $\bX=\bz$. Thus $\vf$~\sf ies~(N1). \Mo $\vf(\la \bX) = \|L(\la\bX)\| = \|\la L\bX\| = |\la|\,\|L\bX\|
= |\la|\vf(x)$ for $\la\in\R$, $x\in\R^N$. Hence $\vf$ \sf ies~(N2). Finally, $\vf(\bX+\bY) = \|L(\bX+\bY)\| = \|L\bX+L\bY\| \le \|L\bX\| + \|L\bY\|
= \vf(\bX)+\vf(\bY)$, hence $\vf$~\sf ies~(N3).

In view of Lemma \rf{l8.142} \te s $c\in\R_{>0}$ \st $c\|\bX\|_2 \le \|L\bX\|$, $\bX\in\R^N$. \csq, if $\bY=L\bX$ then $\bX = L\Inv \bY$ \fa $\bX\in
\R^N$, hence \fa $\bY\in\R^N$ we have $\bX=L\Inv\bY$. \E\Tf $c\|L\Inv\bY\|_2 \le \|\bY\|$, and $\|L\Inv\bY\|\le c\mo\|\bY\|$ \fa $\bY\in\R^N$. \csq,
as above, $L\Inv$ from $(X,d)$ onto $(\R^N,d_2)$ is \Lc. \E\Tf if $\zb un\N$ is a \Cs\ in $(X,\|\cdot\|)$, then $\zb{L\Inv u}n\N$ is a \Cs\ in $(\R^N,
d_2)$, by \E\Pr\ \rf{p8.129}, hence a \cvg\ \sq\ in~$(\R^N,d_2)$ by Theorem \rf{t8.127}. Let $\bY$ be the limit in $(\R^N,d_2)$ of the \sq\
$\zb{L\Inv u}n\N$ in $(\R^N,d_2)$. Set $u:=L\bY$. Then $u_n=L(L\Inv u_n)\to L\bY$. We obtain $u_n\nad d\to u$ which shows that $(X,d)$ is \cp.
\erm

\blm8.146
Let $(X,\|\cdot\|)$ be a real $N$-\dm al \vs\ endowed with a norm $\|\cdot\|$. Let $L:\R^N\to X$ be a linear \is sm $($see \E\df\ $5.1.13)$. Then the
linear map $L:(\R^N,d_2) \to (X,d)$ is \Lc\ as well as~$L\Inv$.
\elm


\bth8.147
Let $(X,\|\cdot\|)$ be a normed $N$-\dm al real \vs\ with $N\in\Na$. Then the \fw\ assertions hold.

\hph i,ii, The \ms\ $(X,d)$ is \emph{\cp} where $d$ is the metric induced by the norm~$\|\cdot\|$.

\hph ii,i, All norms on $X$ are \emph{\ev t}.

\hph iii,, A \nss\ $A$ of $(X,d)$ is \sql\ compact iff \te s $R\in\R_{>0}$ \st $\|u\|\le R$ \fa $u\in A$ and $A$~is closed in~$(X,d)$.

\hph iv,, The \ms\ $(X,d)$ is separable and connected.
\eth


\proof \

\hph i,ii, follows from Lemma \rf{l8.146} and \E\Pr\ \rf{p8.130a}.

\hph ii,i, follows from Remark \rf{r8.143}.

\hph iii,, \ti{Step} 1: \ Case $X:=\R^N$, $\|\cdot\|:=\|\cdot\|_2$.

``\ti{If\/}'' follows from Lemma \rf{l8.140} and Lemma \rf{l8.32}.

``\ti{Only if\/}'': \ Let $A$ be a nonempty closed subset of~$\R^N$ contained in $B:=\{\bX\in\R^N: \|\bX\|_2\le R\}$ \fs $R\in \R_{>0}$. To prove:
$A$~is \sql\ compact. In view of Lemma \rf{l8.140} the set~$B$ is \sql\ compact. Let $\zb\bX n\N$ be an arbitrary \sq\ in~$A$. Since $A\sbs B$,
there are an $\ov{\bX}\in B$ and a sub\sq\ $\{\bX_{n_k}\}_{k\in\N}$ \st $\bX_{n_k}\nad{d_2}\to \ov{\bX}$. Since $A$~is closed, we have
$\ov{\bX}\in A$.

\ti{Step} 2: \ \ti{General case.} Let $L:\R^N\to X$ be as in Lemma \rf{l8.146}. Then \te\ $C_1,C_2\in\R_{>0}$ \st $\|L\bu - L\bv\|\le
C_1\|\bu-\bv\|_2$, $\bu,\bv\in\R^N$, and $\|L\Inv x-L\Inv y\|_2 \le C_2\|x-y\|$, $x,y\in X$.

``\ti{If\/}:'' \ Let $\zb xn\N$ be a \sq\ in~$A$, and let $\bu_n=L\Inv x_n\in\R^N$. Then $\|\bu_n\|_2 = \|L\Inv x_n\|_2 = \|L\Inv x_n-L\Inv 0\|_2
\le C_2\|x_n-0\| = C_2\|x_n\|\le C_2R$ since $x_n\in A$ \fa $n\in\N$. Thus $\|\bu_n\|_2 \le C_2R$, $n\in\N$. Applying Lemma \rf{l8.140} with
$R:=C_2R$ we find that \te\ $\bu\in\R^N$ and a sub\sq\ $\{\bu_{n_k}\}_{k\in\N}$ \st $\bu_{n_k}\nad{d_2}\to\bu$. \E\Tf $x_{n_k} = L\bu_{n_k}\nad d \to
L\bu \in X$, since $L$~is \Lc. We now used the closedness of~$A$ to infer that $\ov x:=L\bu\in A$. Hence $A$~is \sql\ compact.

``\ti{Only if\/}'' is a con\sq\ of the \fw\ observation. Let $f$ be a \ctn\ bi\jn\ from a \ms\ $(X_1,d_1)$ onto a \ms\ $(X_2,d_2)$. If $X_1$~is \sql\
compact then so is~$X_2$.

\hph iv,, We first show that $\Q^N$ is dense in $(\R^N,d_2)$.

``$N:=1$'': \ $0\in\Q$; let $x\in\R_{>0}$. Since $(\R,\ge)$ is \Ar, \te s $\ov m\in\Na$ \st $x<\ov m\cdot 1$. Let $A$ denote the set $\{m\in\Na:
x<m\cdot 1\}$. Since $A\ne\vn$ and is \bb\ by~$1$, it has a least \el\ denoted by~$m'$. Clearly, $m'\ge1$. Set $n:=m'-1$. Then $n\in\N$ and we have
$n\le x< n+1$. The \ig~$n$ is called the integral part of~$x$ and $x-n$ its fractional part. Note that $\a:=x-n\in(0,1)_\R$. By Theorem \rf{t6.14},
\te s an in\cre\ \sq\ $\zb sm\N$ in~$(0,1)$ \st $\supl_{m\ge0}s_m = \a$. Then $x=\supl_{m\ge0}(n+s_m)$, hence $n+s_m\to x$ with $n+s_m\in\Q$. If $x\in
\R_{<0}$ then $-x=\wt n+\wt s_m$, $\wt n\in\N$, $\wt s_m\in\Q$, \st $\wt n+\wt s_m \to -x$. Hence $-\wt n-\wt s_m\to -x$ since the map $x\mt -x$
in~$\R$ is \ctn\ ($|-x-(-y)|=|x-y|$). \E\Tf \fe $x\in\R$ \te s a \sq\ in~$\Q$ converging to~$x$. \E\Tf $\ov Q=\ov Q\en s = \R$, and $\Q$~is dense
in~$\R$.

``$N>1$'': \ Let $\bu\in\R^N$. Then \fe $i\in[1,N]$ \te s a \sq\ $\{\bu_n(i)\}_{n\in\N}$ in~$\Q$ \st $\bu_n(i)\to \bu(i)$. By \E\Pr\ \rf{p7.20}
we have $\bu_n\nad{d_2}\to \bu$. \If that $\Q^N$ is dense in $(\R^N,d_2)$.

Next, we claim that $\Q^N$ is \ct e. Indeed, $\Q$~is \ct e, \te s a bi\jn~$\vf$ from~$\N$ onto~$\Q$. Define $\Phi:\N^N \to \Q^N$ by setting $\Phi(\bn)
(i):= \vf(\bn(i))$, $\bn\in\N^N$, $i\in[1,N]$. One verifies that $\Phi$ is a bi\jn\ with inverse $\Phi\Inv(\bn)(i)=\vf\Inv(\bn(i))$, $i\in[1,N]$.
\E\Tf $\Q^N$ and~$\N^N$ are \ep. \Mo $\N^N$~is \ct e, see Theorem 3.2.15 and Exercise 3.3.60. \If that $\Q^N$ is \ct e, hence $(\R^N,d_2)$ is
separable.

Finally, we show that $(X,d)$ is separable. Since $L:(\R^N,d_2) \to (X,d)$ is \ctn, we have $L(\ov{\Q^N}) = \ov{L(\Q^N)}$ where $\ov{\Q^N}$ (resp.\
$\ov{L(\Q^N)}$) is the closure of~$\Q^N$ (resp.\ $L(\Q^N)$) in $(\R^N,d_2)$ (resp.\ $(X,d)$). Indeed, $L(\ov{\Q^N}) \sbs \ov{L(\Q^N)}$ by Theorem
\rf{t8.89}, hence $X=L(\R^N) \sbs \ov{L(\Q^N)} \sbs X$. \If that $L(\Q^N)$ is dense in $(X,d)$ and $L(\Q^N)$ is \ct e since $L$~is a bi\jn.

``\ti{Connectedness}'': \ Suppose, for \cd ion, that \te\ $O_1$ and~$O_2$ \nos s of $(X,d)$ \sf ying $O_1\cap O_2 = \vn$. Let $x_1\in O_1$ and
$x_2\in O_2$. Let $\vf:[0,1]\to X$ be defined by $\vf(t):=(1-t)x_1 + tx_2$, $t\in[0,1]$. Then for $0\le s\le t\le 1$ we have
$\|\vf(t)-\vf(s)\| = \|(s-t)(x_1-x_2)\| = |t-s|\,\|x_1-x_2\|$. \E\Tf $\vf$~is in\jc, since $t\ne s$ implies $\|\vf(t)-\vf(s)\|>0$, hence $\vf(t)\ne
\vf(s)$. \Mo $\vf$~is \Lc\ from $[0,1]$ endowed with the usual metric of~$\R$ into $(X,d)$. Note that $\vf(0)=x_1$ and $\vf(1)=x_2$. Let $U_1$
and~$U_2$ be subsets of $[0,1]$, closed \il\ of~$\R$, defined by $U_i:=\{t\in[0,1]: \vf(t)\in O_i\}$, $i=1,2$. Since $\vf$~is \ctn, $U_1$~and~$U_2$
are open subsets of $[0,1]$ endowed with the \rv\ topology. Since $\vf$~is in\jc\ and $O_1\cap O_2=\vn$, we have $U_1\cap U_2=\vn$. Since
$X=O_1\cup O_2$, $\vf\Inv[X] = \vf\Inv[O_1\cup O_2] = \vf\Inv[O_1] \cup \vf\Inv[O_2] = U_1\cup U_2$.
We claim that both $U_1$ and~$U_2$ are not empty. Indeed, since $d(x_1,x_2) = \|x_1-x_2\| > 0$,
the balls $B(x_1;\frac13d(x_1,x_2))$  and $B(x_2;\frac13d(x_1,x_2))$ are disjoint. Since $\wt U_1:= \vf\Inv[B(x_1;\frac13d(x_1,x_2))]$ is open,
contains~$0$ and does not contain~$1$, we have $\wt U_1\sbs U_1$. Since $\wt U_1\sms0\ne\vn$, \te s $t_1\in(0,1)$ \st $\vf(t_1)\in U_1$. Likewise one
shows that \te s $t_2\in(0,1)$ \st $\vf(t_2)\in U_2$. Since $[0,1]$ is connected in the \rv\ topology by, we obtain a \cd ion. \If that $(X,d)$ is
connected. The proof of Theorem \rf{t8.147} is complete.
\endproof

\brs8.148 \

\hph i,ii, We claim that the sets $\{u\in V: \|u\|\le R\}$, $R\in\R_{>0}$, are \sql\ compact in the \ms\ $(V,d)$. Indeed, from the in\et y
$\bigl|\|x\|-\|y\|\bigr| \le \|x-y\|$ we infer that the map $\vf$: $x\mt \|x\|$ from~$V$ into~$\R$ is \Lc, hence \ctn. \Mo $\{u\in V: \|u\|\le R\}
= \vf\Inv\bigl[[0,R]\bigr]$, $[0,R]$~is closed in~$\R$, hence $\vf\Inv\bigl[[0,R]\bigr]$ is closed in~$V$. The claim follows from Theorem \rf{t8.147}.
If $V:=\R$ then the set $[-R,R]$ is not only \sql\ compact, but also compact in view of Theorem \rf{t8.119}. It turns out that in a \ms\ a~set~$K$
is compact iff it is \sql\ compact (see \cite[Theorem 2.3.1]{Du} or \cite[Theorem D.10]{JvN}).

\hph ii,i, A \sql\ compact \ms\ $(K,d)$ is \ti{\cp}. Indeed, every \Cs\ in $(K,d)$ possesses a \cvg\ sub\sq, hence the whole \sq\ is \cvg. \E\oh if
a set $(K,d)$ is compact then \fe $n\in\Na$, the collection of all open balls of radius $\frac1n$ is an open cover of~$K$, \Tf a finite subset of this
collection is also a cover of~$K$.

A \ms\ $(M,d)$ possessing this \pp y is called \ti{totally bounded\/}.\index{totally bounded \ms} According to \cite[pp.\ 76, 77]{Du}, Hausdorff (1914) \itd\ this notion and
proved that a \ms\ is compact iff it is \cp\ and totally bounded. This \ev ce is also a part of the theorems mentioned above in~\cite{Du}.

\hph iii,, From what precedes and Theorem \rf{t8.121} we infer that a (\sql) compact \ms\ is a \Ps. \E\Tf a \Ps\ is compact iff it is totally bounded.
\ers

In the next theorem we show that the \ms\ $(\R^N,D)$ \itd\ in \E\Pr\ \rf{p7.25} is a \Ps.

\bth8.149
Let $(\R^\N,D)$ be the \ms\ with the metric defined in \er{7.42}, then\dw

\hph i,ii, $(\R^\N,D)$ is a \cp\ \ms.

Let $\R^\N_0$ $($resp. $\Q^\N_0)$ be defined by
$\R^\N_0$ $($resp.\ $\Q^\N_0):=\{\bX\in\R^\N$ $($resp.\ $\Q^\N)$ \sf ying $\bX(i):=\bz$
\fa $i\ge N$ \fs $N\in\N$ depending on~$\bX\}$.

\hph ii,i, $\R^\N_0$ is dense in $(\R^\N,D)$.

\hph iii,, $\Q^\N_0$ is dense in $(\R^\N,D)$.

\hph iv,, $\Q^\N_0$ is \ct y infinite.
\eth

\proof \

\hph i,ii, Let $\zb \bX n\N$ be a \Cs\ in $(\R^\N,D)$. Let $\psi_i:\R^\N\to \R$, $i\in\N$, be the $i$-th \co\ map, that is, $\psi_j(\bY):=\bY(j)$ \fa
$\bY\in\R^\N$ and $j\in\N$. Let $\ve\in\R_{>0}$. \E\te s $M\in\N$ \st $D(\bX_n,\bX_m)<\ve$ \fs $n,m\ge M$. \E\Tf $\frac1{2^{i+1}}|\psi_i(\bX_n) -
\psi_i(\bX_m)|<\ve$, $n,m\ge M$, hence $|\psi_i(\bX_n) - \psi_i(\bX_m)|\le 2^{i+1}\ve$, $n,m\ge M$. Set $\ve':=2^{i+1}\ve$. \E\te s $M'\in\N$ \st
$|\psi_i(\bX_n) - \psi_i(\bX_m)|<\ve'$, $n,m\ge M'$. \If that $\{\psi_i(\bX_n)\}_{n\in\N}$ is a \Cs\ in~$\R$. Since $(\R,D)$ is \cp, \te s $\a_i\in\R$
\st $\psi_i(\bX_n)\to \a_i$. Define $\bX\in\R^\N$ by setting $\bX(j):=\a_j$ \fa $j\in\N$. Since $\bX_n\lm\bX$, we obtain $D(\bX,\bX_n)\to0$ by
Theorem \rf{t7.26}. \csq, since the \Cs\ $\zb \bX n\N$ was arbitrary, the \ms\ $(\R^\N,D)$ is \cp.

\hph ii,i, Let $\bX\in\R$ and let $P_l \bX(i)\nde7.46 {:=}\bca \psi_i(\bX) & 0\le i\le l, \\ 0 & i>l,\eca$ \fa $l\in\N$. Then we have $P_l\bX \in
\R^\N_0$ \fa $l\in\N$ and $P_l\bX(i)\to \bX(i)$ \fa $i\in\N$ since \fa $i\in\N$, $P_l\bX(i)=\bX(i)$ \fa $l\ge i$. \E\Tf $\ov{\R^\N_0} \nad*=
\ov{R^\N_0}{}\en s = \R^\N$. In~$\nad*=$ we used \E\Pr\ \rf{p8.88}.

\hph iii,, Let $\bX\in\R^\N_0$. Then \te s $N\in\N$ \st $\bX(i)=0$ for $i>N$. Let $\bu\in\R^\N$ be defined by $\bu(i):=\bX(i)$, $0\le i\le N$. Then
$\bu\in\R^\N$. \If from the proof of the separability of $(X,d)$ in Theorem \rf{t8.147} that \te s a \sq\ $\zb \bu n\N$ in~$\Q^\N$ \st $\bu_n(i)\to
\bu(i)$, $0\le i\le N$. \E\Tf setting $\bX(i):=\bu(i)$, $0\le i\le N$, and $\bX(i):=0$, $i>N$, we find that $\bX_N\lm\bX$ in $(\R^\N,D)$, hence
$D(\bX_n,\bX)\to0$ by Theorem \rf{t7.26}. Since $\bX_n\in\Q^\N_0$, we have $\R^\N_0\sbs (\ov{\Q^\N_0})\en s = \ov{\Q^\N_0}\sbs \ov{\R^\N_0}= \R^\N$.
\E\Tf $\R^\N =\ov{\R^\N_0} =\ov{\ov{\Q^\N_0}} = \ov{\Q^\N_0} \sbs \R^\N$. \If that $\ov{\Q^\N_0} = \R^\N$, hence $\Q^\N_0$ is dense in~$\R^\N$.

\hph iv,, Consider the \sq\ $\zb \ve i\N$ in $\Q^\N$ defined by
$\ve_i(j) :=\bca  1 & \hbox{if }i=j,\\0 & \hbox{if }i\ne j,\eca$
$i,j\in\N$. Then \fe $N\in\N$, $\spn_{i\in[0,N]}\{\ve_i\}$ in $(\R,\Q^\N)$ is a $N$-\dm al \vs\ over the field~$\Q$ with basis $\{\ve_i: i\in[0,N]\}$.
Set $V_N:=
\spn_{0\le i\le N} \{\ve_i\}$. We have $V_N\sbs V_M$ if $N<M$, $M\in\N$, and $\Q^\N_0 = \bcl_{N\in\N}V_N$. \If from the proof of Theorem \rf{t8.147}
that $V_N$ is \ep\ to~$\Q^{[0,N]}$. \If from Remark 3.3.31 that $\R^\N_0$ is \ep\ to~$\N$, hence \ct y infinite.
\endproof

\brs8.150 \

\hph i,i, Observe that the map $\vf:\R^\N \to \R_{\ge0}$ defined by $\vf(\bX):=D(\bX,\bz)$, $\bX\in\R^\N$, \sf ies (N1) and (N3), but not (N2), hence
$D$~is \ti{not\/} a norm on~$\R^\N$,

\hph ii,, \E\fa $R\in\R_{\ge0}$ the sets $[-R,R]^\N$ and $[0,R]^\N$ are (\sql) compact subsets of $(\R^\N,D)$ (see \cite{Du}, \cite{JvN} and
\cite{Kelley}, pp.~238--239: D. The diagonal process and sequential compactness).
\ers

\bex8.151
Show that the topology on $\R^\N$ \itd in \er{8.77} is the same as the topology induced by the metric~$D$ of Theorem \rf{t8.149}.
\eex

We now consider infinite \dm al \ti{normed\/} real \vs s. It turns out that such spaces with \ct y infinite bases are \ti{not\/} \cp\ \wrt the induced
metric. This is a con\sq\ of Baire's lemma. As a motivation we recall that a singleton in~$\Q$ or~$\R$ endowed with the usual metric is a closed
subset with empty interior. Observe that the set~$\Q$ is a \ct e union of such sets since $\Q$~is \ct e. \E\oh every infinite \ct e subset of~$\R$ is
not equal to~$\R$ since $\R$~is un\ct e, but it turns out that every \ct e subset of~$\R$ has empty interior. We present a form of Baire's Lemma.

\def\namespec{Baire's Lemma}
\begin{tspc} \lb{l8.149}
Let $(X,d)$ be a \emph{\cp} \ms\ and let $\zb An\N$ be a \sq\ of closed subsets of~$X$ with empty interior, that is, $A_n=\ov A_n$ and $\kl A_n=\vn$
\fa $n\in\N$. Then the set $\bcl_{n\in\N}A_n$ has empty interior.\index{Baire's Lemma}
\end{tspc}

We give an application of this lemma before giving its proof. Let $(\R,V)$ be a real \vs\ having a \ct y infinite basis $\zb ei\N$. Set $V_n:=\spn
\limits_{0\le i<n}\{e_i\}$, $n\in\N$. Then \fa $n\in\N$ $V_n$~is an $(n+1)$-\dm al \lss\ of~$V$, $V_n\subsetneq V_{n+1}$ and $V=\bcl_{k\in\N}V_k$.
Let $\|\cdot\|$ be a norm on~$V$, and let $d$~be the induced metric on~$V$. We claim that $\kl V_n=\vn$ \fa $n\in\N$. Suppose, for \cd ion, that \te s
an $n\in\N$ \st $\kl V_n\ne \vn$. Then \te\ an $n\in\N$ and an $R\in\R_{>0}$ \st $B(x;R)\sbs V_n$. Note that $B(0;R) = \{z-x\in V:z\in B(x;R)\}
\sbs V_n$ since $z,x
\in V_n$ and $V_n$~is a \lss\ of~$V$. \E\Tf $B(0;R)\sbs V_n$. \If that $V\sbs V_n$. Indeed, $0\in V_n$. \Mo if $x\in V\sms0$, then $\la x\in V\sms0$
\fa $\la\in\R\sms0$. Choosing $\la:= \frac R2\|x\|\mo$ which is well-defined since $\|x\|\ne0$, we obtain $\frac R2\|x\|\mo x\in V\sms0$. \csq,
$\|\la x\|=\frac R2\|x\|\mo \|x\|= \frac R2<R$. Hence $\la x\in B(0;R)$, and $\la x\in V_n$. \E\Tf $V\sms0\sbs V_n$, a~\cd ion since $V_{n+1}\sms0
\sbs V\sms0 \sbs V_n\sms0$. \If that $\kl V_n=\vn$ \fa $n\in\N$.

We next show that $V_n$ is closed in $(V,d)$ \fa $n\in\N$. Let $n\in\N$ and let $\zb xk\N$ be a \sq\ in~$V_n$ \st $x_k\nad d\to x$ \fs $x\in V$. We
prove that $x\in V_n$. The \sq\ $\zb xk\N$ is a \Cs\ in $(V_n,d)$ since \fe $\ve\in\R_{>0}$, there is an $N\in\N$ \st $\|x-x_k\|<\frac\ve2$ \fa
$k\ge N$. \If that $d(x_k,x_l) \le d(x_k,x)+d(x,x_l) = \|x_k-x\| + \|x_l-x\|< \ve$ \fa $k,l\ge N$. In view of Theorem \rf{t8.147}\,(i) \te s $x'\in
V_n$ \st $x_n\nad d\to x'$. \E\Tf $d(x,x') \le d(x,x_k)+d(x_k,x')<\ve$ \fa $k\ge N$. Hence $0\le d(x,x')<\ve$ \fa $\ve\in\R_{>0}$, hence $d(x,x')=0$,
and $x=x'\in V_n$. \If that $V_n$ is closed in~$(V,d)$.

If $(V,d)$ were \cp, the interior of $\bcl_{k\in\N}V_k$ would be empty in view of Baire's Lemma, a~\cd ion, since $\bcl_{k\in\N}V_k = V$ and $\kl V
= V\ne\vn$.

\blm8.150
Let $(X,\cO)$ be a \Hts, and let $A$ be a subset of~$X$. Then the \fw\ holds\dw

\hph i,i, $\kl A=\vn$ iff $A^c$ is dense in $X$,

\hph ii,, $A=\ov A$ and $\kl A=\vn$ iff $A^c$ is open and dense in $X$.
\elm

\proof \

(i) \ From \E\Pr\ \rf{p8.67}\,(vi) with $A^c:=A$ we obtain $\kl A=(\ov{A^c})^c$ since $(A^c)^c=A$. If $\kl A=\vn$, we have $\ov{A^c}=X$, that is,
$A^c$~is dense in~$X$ by \E\df\ \rf{d8.120}. Conversely, if $A^c$ is dense in~$X$, then $\ov{A^c}=X$, and $\vn = X^c= (\ov{A^c})^c = \kl A$.

(ii) \ In view of (i) it suffices to observe that $A^c$~is open iff $A=(A^c)^c$ is closed.
\endproof

The Baire Lemma can be reformulated as follows.

\def\namespec{Baire's Lemma {\rm (second version)}}
\begin{tspc} \lb{l8.151}
Let $(X,d)$ be a \cp\ \ms\ and let $\zb Bn\N$ be a \sq\ of open and dense subsets in~$X$. Then $\Bca_{n\in\N}B_n$ is dense in~$X$.\index{Baire's Lemma}
\end{tspc}

\brm8.152
In the proof of Baire's Lemma we will use the notation $\ov{B(x;r)}$ for the closure of the ball $B(x;r)$ in a \ms\ $(X,d)$. Some authors use the
notation $\ov{B(x;r)}$ for the ``closed ball'' defined by $\{y\in X: d(x,y)\le r\}$. Since the map $y\mt d(x,y)$ from~$X$ into $\R_{\ge0}$ is
(Lipschitz)-\ctn\ it follows that the (sequential) closure of $B(x;r)$ is contained in the ``closed ball'' with center~$x$ and radius~$r$. However,
if $d$~is the discrete metric and $X$~contains at least two \el s we have
$$
B(x;\tfrac12) = \{x\} \ne X = \{y\in X: d(x,y)\le 1\}.
$$
Thus the closure of $B(x;1)$ is not the \crs\ ``closed ball''.
\erm

\proof[Proof of Lemma \rf{l8.151}]
We have to show that $\Bca_{n\in\N}B_n$ is dense in $X$, that is, $\ov{\Bca_{n\in\N}B_n} = X$. In view of Lemma \rf{l8.103} it suffices to show that
every nonempty open subset of~$X$ intersects $\Bca_{n\in\N}B_n$, that is, \fe $\o\in\cO(d)\sms\vn$, we have $\bigl(\Bca_{n\in\N}B_n\bigr) \cap
\o\ne\vn$. We claim that it suffices to find two \sq s $\zb xn\N$ in~$X$ and $\zb rn\N$ in~$\R_{>0}$ \sf ying \fa $n\in\N$:
\bga 8.96
B(x_{n+1};r_{n+1}) \sbs \ov{B(x_{n+1};r_{n+1})} \sbs B(x_{n};r_{n})\cap B_n,\\
0< r_{n+1}\le \tfrac12 r_n \lb{8.97}\\
\intertext{and}
x_0\in\o, \lb{8.98}\\
\ov{B(x_0;r_0)} \sbs \o. \lb{8.99}
\e
Indeed, suppose that \er{8.96}--\er{8.99} hold \fa $n\in\N$, then we obtain
\bea 8.100
{}&\hbox{the \sq\ $\zb xn\N$ is a \Cs\ in $(X,d)$,}\\
&x\in X, \hbox{ the limit $x$ of the \sq\ $\zb xn\N$ \sf ies}\lb{8.101}\\
&\hskip200pt  x \in B(x_{n};r_{n})\hbox{ \fa $n\in\N$, }\non \\
&x \in \Bigl(\bigcap_{n\in\N}B_n\Bigr) \cap \o. \lb{8.102}
\e

``\er{8.100}'': \ From \er{8.96} we infer \fa $n\in\N$ $x_{n+1}\in B(x_{n+1};r_{n+1}) \sbs B(x_n;r_n)$, hence $d(x_n,x_{n+1})< r_n$. Let $m:=n+p$,
$p\in\Na$. Then $d(x_n,x_{n+p}) \nad{\rm(M3)}\le \suml_{k=0}^{p-1} d(x_{n+k},x_{n+k+1}) \le \suml_{k=0}^{p-1} r_{n+k} \nde8.97 \le r_n\cdot
\suml_{k=0}^{p-1}\frac1{2^k}\le 2r_n$. From \er{8.97} we have $0<r_n<\frac1{2^n}r_0$ \fa $n\in\Na$. \E\Tf \fa $m>n\ge0$ we obtain $d(x_n,x_m) <
\frac1{2^{n-1}}r_0$. \If that \fe $\ve\in\R_{>0}$ \te s an $N\in\Na$ \st $m>n>N$ implies $d(x,x_m)<\ve$.

``\er{8.101}'': \ Setting $O_n:=B(x_n;r_n)$, $n\in\N$, we find from \er{8.96} that $O_{n+1}\sbs O_n$ \fa $n\in\N$. By \E\Pr\ \rf{p8.67}\,(ii) we have
$O_n\sbs \ov{O_n}$ \fa $n\in\N$. \E\Tf \fa $p\in\N$, $x_{n+p}\in
B(x_{n+p};r_{n+p})\sbs O_{n+p}\sbs O_n \sbs \ov{O_n} = \ov{B(x_n;r_n)}$. Given $n\in\N$ set $y_p:=x_{n+p}$, $p\in\N$. Then $y_p\in \ov{B(x_n;r_n)}$
\fa $p\in\N$. Since $\zb xn\N$ is a \Cs\ in $(X,d)$, and since $(X,d)$ is \cp, \te s $x\in X$ \st $x_n\nad d\to x$. Note that $\zb yp\N$ is a sub\sq\
of the \sq\ $\zb xn\N$ hence $y_p\nad d\to x$. Since $y_p\in \ov{B(x_n;r_n)}$ \fa $p\in\N$ and since $\ov{B(x_n;r_n)}$ is closed by \E\Pr\ \rf{p8.67}
(hence \sql\ closed) in $(X,d)$, we obtain $x\in\ov{B(x_n;r_n)}$ \fa $n\in\N$, which proves \er{8.101}.

``\er{8.102}'': \ Since $\ov{B(x_n;r_n)} = \ov{O_n} \sbs \ov{O_0}$ (by \E\Pr\ \rf{p8.67}\,(v)) and $\ov O\sbs\o$ by \er{8.95}, we have $x\in\o$. \Mo
since $\ov{B(x_{n+1};r_{n+1})}\nde8.96 \sbs B_n$ \fa $n\in\N$, we have $x\in B_n$ \fa $n\in\N$, hence $x\in \bigl(\Bca_{n\in\N}B_n\bigr) \cap\o$.

We now prove the \ex\ of \sq s \sf ying \er{8.96}--\er{8.99}.

``\er{8.98}, \er{8.99}'': \ Since $\o\in\cO(d)\sms\vn$ \te s $x_0\in\o$. \er{8.95} is a con\sq\ of the \fw\ observation. Let $(X',d')$ be a \ms, let
$O'$ be a \nos\ of~$X'$, and let $x'\in O'$. Then \te s $s'\in\R_{>0}$ \st $B(x';s')\sbs O'$. Clearly, $B(x';\frac12s') \sbs B(x',s')$. Then
$\ov{B(x';\frac12s')} \sbs \{y'\in X': d'(x',y')\le \frac12s'\} \sbs B(x';s')$ since $\frac12s'< s'$. Setting $r':=\frac12s'$ we have $\ov{
B(x';r')}\sbs O'$. \E\Tf \te\ $x_0\in X$ and $r_0\in\R_{>0}$ \st \er{8.98} and \er{8.99} hold.

``\er{8.96}, \er{8.97}'': \ We proceed by \In\ on $n\in\N$.

We first show that \te\ $x_1\in X$ and $r_1\in\R_{>0}$ \st \er{8.96}--\er{8.97} hold for $n:=0$. Note that $B(x_0;r_0)\cap B_0$ is open since
both $B(x_0;r_0)$ and~$B_0$ are open, and $B(x_0;r_0)\cap B_0$ is not empty by Lemma \rf{l8.103} since $B(x_0;r_0)$ is a nonempty open set and $B_0$
is dense in~$X$. \If that \te\ $x_1\in X$ and $s_1\in\R_{>0}$ \st $B(x_1;s_1)\sbs B(x_0;r_0)\cap B_0$. Setting $r_1:=\min(\frac12r_0,\frac12s_1)$, we
find that $B(x_1;r_1)\sbs \ov{B(x_1;r_1)}\sbs B(x_0;r_0)\cap B_0$ and $0<r_1\le \frac12r_0$.

We now assume that \er{8.96}, \er{8.97} hold for $n:=k\in\N$, and we prove that \er{8.96}, \er{8.97} hold for $n:=k+1$. Then $B(x_k;r_k)\cap B_k$
is open and nonempty since both $B(x_k;r_k)$ and~$B_k$ are open and $B_k$~is dense in~$X$. \E\Tf \te\ $x_{k+1}\in X$ and $s\in\R_{>0}$ \st
$B(x_{k+1};s) \sbs B(x_k;r_k)\cap B_k$. Setting $r_{k+1}:= \min(\frac12s,\frac12r_k)$, we find that \er{8.96}, \er{8.97} hold for $n:=k+1$.
This \cp s the proof of Lemma \rf{l8.151}.
\endproof

We now show that the space $B(A;\R)$ \itd in \er{8.82}, endowed with the metric $D_{\sup}$ (see \er{8.83}), is \cp.

\bth8.154
The \ms\ $(B(A;\R),D_{\sup})$ is \cp.
\eth

\proof
Let $\zb fn\N$ be a \Cs\ in $(B(A;\R),D_{\sup})$, that is, \fe $\ve\in\R_{>0}$ \te s an $N\in\N$ \st \fa $n,m\ge N$
\beq8.103
D_{\sup} (f_n,f_m)< \ve.
\e
Since $|f_n(x)-f_m(x)| \le \supl_{y\in A} |f_n(y)-f_m(y)| = D_{\sup}(f_n,f_m)$, the \sq\ $\zb fn\N$ in~$\R$ is a \Cs\ \fa $x\in A$. Since $(\R,d)$ is
\cp, the \sq\ $\zb fn\N$ \cg s to some limit in~$\R$ denoted by $f(x)$ \fa $x\in A$. Thus $x\mt f(x)$ is a map from~$A$ into~$\R$. We claim that
$f$~is bounded, that is, $f\in B(A;\R)$, and that $D_{\sup}(f_n,f)\to 0$. Let $\ve:=1$ in \er{8.103}. Then $|f_n(x)-f_m(x)|\le 1$ \fa $x\in A$ and
\fa $n,m\ge N$ \fs $N\in\N$. \If that $|f_n(x)| = |f_n(x)-f_N(x)+f_N(x)| \le |f_n(x)-f_N(x)|+|f_N(x)| \le 1+|f_N(x)|$ \fa $x\in A$ and all $n\ge N$.
Since the map $f_N$ is bounded, \te s $M\in\R_{\ge0}$ \st $|f_N(x)|\le M$ \fa $x\in A$. \csq, $|f_n(x)|\le 1+M$ \fa $x\in A$ and all $n\ge N$. \Mo
$|f_n(x)|\le \supl_{y\in A}|f_n(y)|$ \fa $n\in[0,N]$. \E\Tf $|f_n(x)| \le M'$ \fa $x\in A$ and $n\in[0,N]$ with $M':=\max\limits_{0\le n\le N}
\supl_{x\in A}|f_n(x)|$. \If that $|f_n(x)|\le 1+M+M'$, $n\in\N$, $x\in A$, hence $f\in B(A;\R)$.

Set $\vf_n(t):= |f_n(x)-t|$, $n\in\N$, $t\in\R$, $x\in A$. We have $|\vf_n(t)-\vf_n(s)| = \bigl||f_n(x)-t|-|f_n(x)-s|\bigr| \le \bigl|
(f_n(x)-t) - (f_n(x)-s)\bigr| = |t-s|$, $s,t\in\R$. Thus \fa $n\in\N$ the maps $\vf_n:\R\to\R_{\ge0}$ are \Lc, hence \sql\ \ctn. \E\Tf \fa $n\ge N$
we have $|f_n(x)-f(x)| = |f_n(x)-\lim\limits_{m\ge0} f_m(x)| = \lim\limits_{m\ge0} |f_n(x)-f_m(x)| \nad*= \limu_{m\ge0} |f_n(x)-f_m(x)|
\nad{**}\le \limu_{m\ge0}\ve = \ve$ \fa $x\in A$. In $\nad*=$ (resp.~$\nad{**}\le$) we use \E\Pr\ \rf{p8.10} (resp.~\rf{p8.12}). \csq, \fa $\ve\in\R
_{\ge0}$ \te s $N\in\N$ \st $D_{\sup}(f_n,f)<\ve$, that is, $f_n\nad{D_{\sup}}\to f$. \If that $(B(A;\R),D_{\sup})$ is \cp.
\endproof

\bco8.155
Let $C(K)$ denote the set \itd in Notation \rf{n8.111}, endowed with the metric $D_{\sup}$, then $(C(K),D_{\sup})$ is a \cp\ \ms.
\eco

\proof
Let $\zb fn\N$ be a \Cs\ in $(C(K),D_{\sup})$. Then $\zb fn\N$ is a \Cs\ in $B(K;\R)$ since $C(K)\sbs B(K;\R)$ (see Example \rf{xa8.110}). By Theorem
\rf{t8.154} \te s $f\in B(K;\R)$ \st $f_n\nad{D_{\sup}}\to f$. By Lemma \rf{l8.111} $f\in C(K)$, hence $(C(K),D_{\sup})$ is \cp.
\endproof

\begin{dfn}[{\cite[p. 158]{Du}}]\lb{8.157}
A normed real \vs, that is, \cp\ for the metric induced by the norm, is called a \ti{real \Bs}.
\edf

In view of Theorem \rf{t8.147} every finite \dm al normed real space is a \Bs. Let $A$ be a \ns\ and let $(\R,\R^A)$ denote the real \vs\ \itd\ in
Example 5.1.2\,(i). Note that the set $B(A;\R)$ defined in \er{8.82} is a \lss\ of $(\R,\R^A)$ (see \E\df\ 5.2.6) since $\la f+g\in B(A;\R)$ whenever
$f,g\in B(A;\R)$ and $\la\in\R$. Hence $(\R,B(A;\R))$ is a real \vs\ by Lemma 5.2.7. Let $\ve_i:A\to\R$, $i\in A$, be defined by
$$
\ve_i(t):=\bca 1 &\hbox{for }t=i,\\ 0 &\hbox{for }t\ne i.\eca \q t\in\R.
$$
Then $\ve_i \in B(A;\R)$ \fa $i\in A$. \Mo
if $B$ is a \nfs\ $B$ of~$A$ and $\la:B\to R$, then
$\suml_{j\in B} \la(j)\ve_j=\bz$, where $\bz(t):=0$ for $t\in A$, implies $\la(j)=0$ \fa $j\in B$.

If $A$ is finite and if $\|\cdot\|$ is any norm on $B(A;\R)$, hence $(B(A;\R),\|\cdot\|)$ is a $\#(A)$-\dm al real \Bs, then
$\R^A=B(A;\R)$. If $A$~is infinite, then $D_{\sup}$ is a metric on $B(A;\R)$. Setting $\|f\|_\iy:= D_{\sup}(f,\bz)$, $f\in B(A;\R)$, we find that
$\|\cdot\|_\iy$ is a norm on $B(A;\R)$. Indeed, $\|f+g\|_\iy = \supl_{t\in A}|(f+g)(t)+0| = \supl_{t\in A}|f(t)+g(t)|$. Since $|f(t)+g(t)| \le
|f(t)|+|g(t)|$, $t\in A$, we find $|f(t)+g(t)| \le \|f\|_\iy+ \|g\|_\iy$. $\|\la f\|_\iy = \supl_{t\in A}|\la f(t)| = \supl_{t\in A} |\la|\,|f(t)|
= |\la| \supl_{t\in A}|f(t)|$ by~\er{5.14} if $\la\ne0$, and clearly equal to~$0$ if $\la=0$.

Finally, $\|f\|_\iy=0$ if $f=0$. Conversely, if $\|f\|_\iy = 0$, then $\supl_{t\in A}|f(t)|=0$, and \fe $s\in A$ we have $0\le |f(s)| \le
\supl_{t\in A}|f(t)|=0$. Hence $f=0$. Thus $\|\cdot\|_\iy$ is a norm on $B(A;\R)$, usually called the \ti{supremum norm}. In view of Theorem
\rf{t8.154}, $(B(A;\R),\|\cdot\|_\iy)$ is a \ti{real \Bs}.

Let $C(K)$ denote the set \itd in Notation \rf{n8.111}. We have $C(K)\sbs B(K;\R)$ and \fa $f,g\in C(K)$, all $\la\in\R$: $\la f+g\in C(K)$. Hence
$C(K)$ is a real \vs\ and the \rt ion of $\|\cdot\|_\iy$ to $C(K)$, still denoted by $\|\cdot\|_\iy$, is a norm on~$C(K)$. In view of Corollary
\rf{c8.155}, $(C(K),\|\cdot\|_\iy)$ \ti{is a real \Bs}.

We next consider the \Bs\ $(B(\N;\R),\|\cdot\|_\iy)$, that is, the space of all bounded \sq s in~$\R$. The usual notation for this space is~$l^\iy$.

\bpr8.158
The space $(l^\iy,\|\cdot\|_\iy)$ is not separable.
\epr

\proof
Let $\cP(\N)$ be the power set of~$\N$ and let $\vf:\cP(\N)\to l^\iy$ be defined by
$$
(\vf(M))(i) := \bca 1 &\hbox{if }i\in M,\\ 0 &\hbox{if }i\notin M, \eca
$$
where $M\in\cP(\N)$. By Cantor's Theorem 1.4.26, the set $\cP(\N)$ is un\ct e. \E\oh the map~$\vf$ is in\jc, that is, if $\vf(M_1)=\vf(M_2)$, $M_1,
M_2\in\cP(\N)$, we find $M_1\sbs M_2$ and $M_2\sbs M_1$. Indeed, let $j\in M_1$, then $\vf(M_1)(j)=1$. Hence $\vf(M_2)(j)=1$, and $j\in M_2$.
Interchanging $M_1$ and~$M_2$, we find that $M_1=M_2$. Then the range of~$\vf$ in~$l^\iy$, denoted by $R(\vf)$, is un\ct e. Note that if $a,b\in
R(\vf)$ with $a\ne b$ then $\|a-b\|_\iy=1$. \Mo if $B_a$ (resp.~$B_b$) denotes the open ball in~$l^\iy$ of center~$a$ (resp.~$b$) and of radius
$\frac13$, then we have $B_a\cap B_b = \vn$ iff $a\ne b$. Now suppose that $H$~is a dense subset of~$l^\iy$, then in view of Lemma \rf{l8.103}, we
obtain $H\cap B_a\ne\vn$ \fa $a\in R(\vf)$. Thus $H$ cannot be \ct e.
\endproof

\blm8.159
Let $Y$ be a \nss\ of a separable \ms\ $(X,d)$. Then $(Y,d)$ is separable.
\elm

\proof
Let $\zb xn\N$ be a \sq\ in~$X$ \st $A:=\bcl_{n\in\N}\{x_n\}$ is dense in~$X$. Given $m\in\Na$ and $n\in\N$ \te s an \el\ of~$Y$, denoted by $y_{mn}$,
\sf ying $d(x_n,y_{mn})< \frac1m$. Set $H:= \bcl_{(m,n)\in\Na\t\N}\{y_{mn}\}$. Since $\Na\t\N \sbs \N\t\N$ is \ct e by (3.3.60), the set~$H$ is \ct e.
Clearly $H\sbs Y$. Let $y\in Y$ be arbitrary and let $\ve\in\R_{\ge0}$. Since $(\R,\ge)$ is \Ar, \te s $m\in\Na$ \sf ying $\frac1m < \frac\ve2$.
Since $A$~is dense in~$X$, \te s $n\in\N$ \st $d(y,x_n)<\frac1m$. Then $d(y,y_{mn}) \le d(y,x_n)+ d(x_n,y_{mn}) < \frac1m+\frac1m = \frac2m<\ve$.
\E\Tf $H$~is dense in~$(Y,d)$, and $(Y,d)$ is separable.
\endproof

\bpr8.160
The real \Bs\ $(B(A;\R),\|\cdot\|_\iy)$ is \emph{not} separable if $A$ contains a \ct y infinite set.
\epr

\proof
Let $A$ be a set containing a \ct y infinite set~$H$. Suppose, for \cd ion, that $(B(A;\R),\|\cdot\|_\iy)$ is separable. We can embed the space
$B(H;\R)$ into the space $B(A;\R)$ by means of the in\jc\ map $j:B(H;\R) \to B(A;\R)$ defined by setting
\beq8.105
j(f)(t) := \bca f(t) & t\in H,\\ 0 & t\notin H. \eca
\e
Note that $\|j(f)\|_\iy = \|f\|_\iy$ \fa $f\in B(H;\R)$. \If from Lemma \rf{l8.159} that the space $(j(B(H;\R)),\|\cdot\|_\iy)$ is separable.
Since the map $j:(B(H;\R),\|\cdot\|_\iy) \to (j(B(H;\R),\|\cdot\|_\iy))$ is a sur\jc\ isometry, hence a \he sm \wrt the metric $D_{\sup}$, the space
$(B(H;\R),\|\cdot\|_\iy)$ is separable (see Exercise \rf{ex8.123}). We now use the bi\jc\ map $\vf^*:\R^\N \to \R^H$ \itd in the proof of \E\Pr\
\rf{p8.96} with $A:=H$. Note that the \rt ion of~$\vf^*$ to $B(\N;\R)$ is a sur\jc\ isometry from $(B(\N;\R),\|\cdot\|_\iy)$ onto $(B(H;\R),
\|\cdot\|_\iy)$. Hence $(B(\N;\R),\|\cdot\|_\iy)$ is separable, \cd ing \E\Pr\ \rf{p8.158}.
\endproof

\bth8.161
Let $(K,d)$ be a $($\sql$)$ compact \ms. Then the space $(C(K),\|\cdot\|_\iy)$ is separable, hence a \Ps.
\eth

For a proof we refer the reader to \cite[\E\Pr\ 2.8]{JvN}, see also \cite{Zabczyk}.

In the last part of this section we present some \Ps s consisting of real-valued \sq s.

We use the notations: $\bX,\bY$ for \el s of~$\R^\N$, $x_k$~for~$\bX(k)$, $k\in\N$ (thus $\bX$~is another notation for $\zb xn\N$), $\bz(k):=0$ and
$\1(k):=1$, $k\in\N$. Endowed with the pointwise \ad~$+$ and \mlc~$\cdot$:
$$
\bal
(\bX+\bY)(k)&:=\bX(k)+\bY(k),\\
(\bX\cdot \bY)(k)&:= \bX(k)\cdot\bY(k), \q k\in\N,
\eal
$$
$(\R^\N,+,\bz)$ (resp.\ $(\R^\N,\cdot,\1)$ is an \ag\ (resp.\ \am). Note that (4.4.1), (4.4.2) and (4.4.3) hold, \Tf $(\R^\N,+,\cdot,\bz,\1)$ is a
\ti{\cmt e semiring with unity} (see \E\df\ 4.1.1).

Note that the map $(\la,\bX)\mt \la\bX$ from $\R\t\R^\N$ into~$\R^\N$ defined by $(\la\bX)(k) := \la x_k$, $\la\in\R$, $k\in\N$, \sf ies (5.1.1),
hence $(\R,\R^\N)$ becomes a real \vs, and $(\R,\R^\N,\cdot)$ an \asc e, \cmt e $\R$-algebra with unity (see \E\df\ 5.2.18). We set $\be_k(l):=
\d_{k,l}$, $k,l\in\N$. Since $\suml_{k=0}^n \la_k\be_k=\bz$ implies $\la_k=0$, $k\in[0,n]$, \fa $n\in\N$, this algebra is \ti{infinite} \dm al.

\If from Theorem \rf{t8.149} that the set $\R^\N$ endowed with the metric~$D$ \itd in \E\Pr\ \rf{p7.25} is a \cp\ \ms\ and that the \ct y infinite
subset $\Q_0^\N$ is dense in $(\R^\N,D)$. Hence $(\R^\N,D)$ is a \ti{\Ps}, but \ti{not\/} a normed space. Indeed,
\bmlg
D(\bz,\1)\nad{\er{7.27},\er{7.37},\er{7.40},\er{7.42}}= \sup_{l\in\N}\sum_{i=0}^l \frac1{2^{i+1}}\min(|0-1|,1)=1\\
=\sup_{l\in\N}\sum_{i=0}^l \frac1{2^{n+1}}\min(|0-2|,1) = D(\bz,2\1) = D(2\bz,2\1) \ne2D(\bz,\1).
\e
\E\Tf the metric $D$ is not induced by a norm (see Exercise \rf{ex8.138}\,(ii)).

We next introduce the \fw

\bns8.163 \

$c:=\{\bX\in\R^\N: \lim\limits_{n\ge0}|x_n-\a|=0 \hbox{ \fs }\a\in\R\}$;

$c_0:= \{\bX\in c: \lim\limits_{n\ge0}x_n=0\}$;

$c_{00} := \R^\N_0$ \itd in Theorem \rf{t8.149}\,(ii).
\ens

We recall that $l^\iy:=B(\N;\R)$ and \fa $\bX\in l^\iy$, $\|x\|_\iy:=D_{\sup}(\bX,\bz)$.

We have $c_{00} \subsetneq c_0 \subsetneq c \subsetneq l^\iy \subsetneq \R^\N$. \Mo $l^\iy$ is a \lss\ of the real \vs\ $(\R,\R^\N)$. The inclusions
$\subsetneq $ can be replaced by ``is a \lss\ of''.

We leave it as a (good) exercise to verify that $c$~is a \ti{closed\/} subset of $(l^\iy,\|\cdot\|_\iy)$. Note that since $(l^\iy,\|\cdot\|_\iy)$ is
a \ms, it suffices to show that $c$~is \ti{\sql} closed in $(l^\iy,\|\cdot\|_\iy)$.

Observe that \fe $\bX\in c$ \te\ unique $\bY\in c_0$ and $\la\in\R$ \st $\bX=\bY+\la\1$ with $\1(k):=1$, $k\in\N$. Indeed, $\la=\lim x_n$. \Mo
$c_0$~is closed in $(c,\|\cdot \|_\iy)$. \If that both $(c,\|\cdot\|_\iy)$ and $(c_0,\|\cdot\|_\iy)$ are real \Bs s.

\begin{dfn}[{\cite[p.~239]{Nrs}}]\lb{d8.164}
An \ti{\asc e} $\R$-algebra $(V,\cdot,\|\cdot\|)$ is called a \tb{real Banach algebra} if $(V,\|\cdot\|)$ is a real \Bs\ and $\|a\cdot b\|\le\|a\|
\,\|b\|$ \fa $a,b\in V$. The Banach algebra $(V,\cdot,\|\cdot\|)$ is called \ti{\cmt e} if $a\cdot b=b\cdot a$ \fa $a,b\in V$, and \ti{unital\/} if
\te s a unity~$e$, that is, $e\cdot a=a\cdot e=a$ \fa $a\in V$.\index{real Banach algebra}
\edf

One verifies that $(c,\cdot,\|\cdot\|_\iy)$ as well as $(C(K),\cdot,\|\cdot\|_\iy)$ with $C(K)$ defined in Notation \rf{n8.111} are \cmt e unital
Banach algebras with $(u\cdot v)(x):=u(x)v(x)$ for $x\in K$, $u,v\in C(K)$. However, $(c_0,\cdot,\|\cdot\|_\iy)$ is a \cmt e $\R$-algebra
\ti{without\/} unity. We recall that $\Q^\N_0$ is \ct e and dense in $(\R_0^\N,\|\cdot\|_\iy)=c_{00}$. Given $\bX\in c_0$ and $m\in\N$ we set
$\bX_m(k)=\bX(k)$ for $0\le k\le m$ and $0$~for $k>0$. One verifies that $\|\bX-\bX_m\|_\iy\to 0$ in~$\R$, hence $\ov{c_{00}} = \ov{c_{00}}{}\en s=
c_0$ where $\ov{\phantom n}$ (resp.\ $\ov{\phantom n}{}\en s$) denotes the closure (resp.\ sequential closure) in the \Bs\ $(l^\iy,\|\cdot\|_\iy)$.
\E\Tf
$(c_0,\|\cdot\|_\iy)$ is \cp\ (see the proof of Corollary \rf{c8.155}). Hence $(c_0,\|\cdot\|_\iy)$ \ti{is a \Ps}. Finally, we show that
$(c,\|\cdot\|_\iy)$ is also a \Ps. Using the notation $X\approx Y$ for two sets $X$~and~$Y$ for which \te s a bi\jn\ $f:X\to Y$, we have
$E:=\{\la\1+\bY:\la\in\Q,\bY\in\Q_0^\N\} \approx \Q\t\Q_0^\N \approx \N\t\N \approx N$. If $\bX\in c$, then $\bX=\la\1+\bY$ with $\la\in\R$ and
$\bY\in c_0$, hence \te s a \sq\ $\{\la_n\1+\bY_n\}_{n\in\N}$ in~$E$ \st $\la_n\to\la$ in~$\R$ and $\bY_n\to \bY$ in $(l^\iy,\|\cdot\|_\iy)$. \E\Tf
$\|(\la\1+\bY)-(\la_n\1+\bY_n)\|_\iy \le \|(\la-\la_n)\1\|_\iy + \|\bY-\bY_n\|_\iy = |\la-\la_n|+\|\bY-\bY_n\|_\iy \to0$ since $\|\1\|_\iy=1$. Since
the closure of~$E$ in $(l^\iy,\|\cdot\|_\iy)$ is \cp\ and since $E$~is \ct e, $(c,\|\cdot\|_\iy)$ \ti{is a \Ps}.

\brm8.165
\ti{Cantor's construction of the field of real \nm s.}

Observe that $c_0$ is not only a \lss\ of the ring $(c,+,\cdot,\bz,\1)$, but it is an \ti{ideal\/} (see \E\df\ 5.4.37). Indeed, if $\bX\in c$ and
$\bY\in c_0$, then $\bX\cdot \bY\in c_0$. Let $\sim$ denote the \ev ce \rl\ on~$c$ defined by $\bX\sim\bY$ if $\bX-\bY\in c_0$, let~$[\bX]$ denote the
\crs\ \ev ce class containing~$\bX$, and let $c/c_0$ denote the set of \ev ce classes. Then $(c/c_0,\hpl,[\bz])$ with $[\bX]\hpl[\bY] = [\bX+\bY]$
is known to be an \ag. If $\bX:=
\la\1+x_0$, $\bY:=\mu\1+y_0$ with $\la,\mu\in \R$, $x_0,y_0\in c_0$, then $\bX\cdot\bY = (\la\mu)\1 + (\la x_0+\mu y_0+x_0\cdot y_0)$ since $\1
\cdot\bZ = \bZ$ \fa $\bZ\in c$. Note that $\la x_0+\mu y_0+ x_0\cdot y_0\in c_0$, hence $[\bX\cdot \bY] = [(\la\mu)\1] = [\la\1\cdot \mu\1]$.

Note that $\la$ and $\mu$ are uniquely defined by $\bX$ and~$\bY$. By setting $[\bX]\hcd[\bY]:=[\bX\cdot\bY]$, $\bX,\bY\in c$, it turns out that
$(c/c_0,\hpl,\hcd,[\bz],[\1])$ is a \cmt e ring with unity. Defining $\vf:c/c_0\to\R$ by setting $\vf([\bX]) = \vf(\la\1+\bY):=\la\in\R$, we find that
$\vf$~is a ring-\is sm from $c/c_0$ into~$\R$. \E\Tf $(c/c_0,\hpl,\hcd,[\bz],[\1])$ is a field, and $\vf(\bX)=\lim\limits_{n\ge0}x_n$.
\erm

In view of Theorem \rf{t8.127} one may replace \cvg\ \sq s by \ti{\Cs s} in~$\R$. It turns out that it suffices to consider \Cs s in~$\Q$ \wrt the
\hbox{$\Q$-\vl} in~$\Q$ defined by $r \mt |r|$, $r\in\Q$ (see \E\df\ \rf{d5.23}). This is \ti{Cantor's approach}.

In \cite[Section 11.2]{Alg}, the field $(\Q,\ge)$ is replaced by an order field $(K,\ge)$, and \Cs s by \ti{\fd\ \sq s} defined as follows:

A \sq\ $\bX\in K^\N$ is said to be a \ti{\fd\ \sq} if \fe $\ve\in K_{>0}$ \te s an $N\in\N$ \st $|x_m-x_n|<\ve$ \we $m,n>N$, $m,n\in\N$.\index{fundamental sequence}

We denote the set of all \fd\ \sq s by $c(K)$. Clearly, $\bz,\1$ belong to~$c(K)$. A \ti{null \sq} $\bX\in K^\N$ is a \fd\ \sq\ \sf ying: \fe
$\ve\in K_{>0}$ \te s $N\in\N$ \st $|x_n|<\ve$ \we $n>N$. We denote the set of null \sq s by $c_0(K)$. Clearly, $\bz\in c_0(K)$,
but $\1\notin c_0(K)$.

It turns out that $(c(K),+,\cdot,\bz,\1)$ is a \ti{\cmt e ring with unity}, that $c_0(K)$ is a \cmt e ring without unity, and that $c_0(K)$ is an
\ti{ideal in the ring} $c(K)$. \E\Tf by \E\Pr\ 4.5.19 $c(K)/c_0(K)$ is a \cmt e ring with unity~[\1]. \Mo \fe $[\bX]\in c(K)/c_0(K)\ne[\bz]$, \te s
$\bY\in c(K)$
\st $[\bX]\hcd[\bY]=[\1]$. \E\Tf $(c(K)/c_0(K),\hpl,\hcd,[\bz],[\1])$ is a \ti{field\/} by \E\df\ 4.4.1, which we denote by~$\wh K$.
The map $j:K\to \wh K$ defined by $j[a] = [\bX]$ with $x_n(k):=a$ \fa $k,n\in\N$ and $a\in K$ is an in\jc\ \hm sm, $j(K)$ is a subfield
of~$\wh K$, hence $(\wh K,j)$ is a field \ext.

A \ti{\fd\ \sq} $\bX$ is called \ti{\po} if \te s an $\ve\in K_{>0}$ and an $N\in\N$ \st $x_n>\ve$ \fa $n>N$.

It turns out that $\bX$ is \po\ iff $\bX+\bY$ is \po\ \fs $\bY\in c_0(K)$ iff $\bX+\bY$ is \po\ \fa $\bY\in c_0$, in which case $[\bX]$ is called \po.
\Mo the \rl\ $\wh\ge$ on~$\wh K$ defined by $[\bX]\wh\ge[\bY]$, $x,y\in K$, if either $[\bX]=[\bY]$ or $[\bX\hat- \bY]$ is \po, makes $(\wh K,\wh\ge)$
an \ti{\of}.

If we denote by $|\cdot|_K$ the \av\ in $(\wh K,\wh \ge)$, then the \fw\ holds:

\ti{Every \fd\ \sq\ $\zb \xi n\N$ is \cvg\ in~$\wh K$}, that is, \te s \ooo $\xi\in \wh K$ \st \fe $\ve\in\wh K_{\hat>\hat0}$ \te s an $N\in\N$ \st
$|\xi\hat- \xi_n|_{\hat K}\wh\le\ve$ \we $n\ge N$. Finally, if $(K,\ge)$ is \Ar, then so is $(\wh K,\wh\ge)$. \E\Ip if $(K,\ge):=(\Q,\ge)$ then
$(\wh K,\wh \ge)$ is order- and ring-\is c to $(\R,\ge)$. This is Cantor's construction of real \nm s.

We recall that if $(K,\ge)$ is an \of, then $K[X]$ is a domain (see Theorem 5.2.24). \Mo $K[X]$ is an \od\ with the \og\ defined in \E\Pr\ \rf{p3.32}.
By Example \rf{xa3.36}\,(ii) this \og\ is \ti{not\/} \Ar\ (see \E\df\ \rf{d3.35}), hence the \crs\ \of\ of \qt s $K[X]$ is \ti{not\/} \Ar\ (see
\E\Pr\ \rf{p3.37} and Corollary \rf{c3.38}).

\bex8.166 \

\hph i,i, Find an \el\ $e$ of $\Q[X]$ \sf ying $0<e<\frac1n \1$ \fa $n\in\N$, with $\ge$ defined in \er{3.81} and $K:=\Q$ endowed with its usual \og.

\hph ii,, Is $\{\frac1n e\}_{n\in\N}$ a \fd\ \sq\ in $(\Q[X],\ge)$? (See \cite[pp.\ 412--413]{Is}.)
\eex

\brm8.166
In the construction of the field $(\wh K,\wh\ge)$ (\cite[Sections 11, 12]{Alg}) the field~$K$ is an \of\ not necessarily \Ar. For example $(K,\ge):=
(\Q[X],\ge)$. However, in that case $(\wh K,\wh \ge)$ cannot be order-\cp, since if it were the case, $(\wh K,\wh \ge)$ would be \Ar\ by \E\Pr\
\rf{p5.13}. Since the map $j:K\to \wh K$ is order-\is c from $K$ onto~$j(K)$, $(j(K),\wh\ge)$ would be \Ar\ by Lemma \rf{l3.39}, as well as $(K,\ge)$,
a~\cd ion.
\erm

We now turn to the case where $(K,\vf)$ is a field $K$ (not necessarily ordered) with an $F$-valued \vl\ where $F$~is an \Ar\ \of. Thanks to Theorem
\rf{t4.35} \te s an in\jc\ order- and ring-\is sm $j:F\to\R$. \If that the map $d_j:K\t K\to\R$ defined by $d_j:=j\circ\vf$ is a \ti{metric} on~$K$
\sf ying \er{7.43} where $\R^\N$~is replaced by~$K^\N$. If we replace $\ve\in K$ in the \df\ of \fd\ \sq s by $\ve\in\R_{>0}$, we find that \fd\ \sq s
are \Cs s. In view of \cite[p.~197]{Alg}, the proof of the \fw\ assertions carry over verbatim. The set of \Cs s endowed
with the pointwise \ad~$+$ and \mlc~$\cdot$ form a \cmt e ring with unity~($=\1$), the set of null-\sq s forms an ideal of this ring, hence the
``\qt\ ring'' forms a \cmt e ring with unity by \E\Pr\ 4.5.19. \Mo the nonzero \el s of this ``\qt\ ring'' are units (see \E\df\ 4.3.34), hence the
``\qt\ ring'' denoted by~$\wh K$ is a \ti{field\/} in view of \E\df\ 5.4.1.

Finally, let $i:K\to\wh K$ be the map defined by $i(a):=[\bX]$ with $\bX\in K^\N$ \st $x_k=\bX(k)=a$ \fa $k\in\N$ and $a\in K$, and $[\bX]$ is the
\ev ce class containing~$a$. It turns out that the map~$i$ is an in\jc\ ring-\hm sm, hence $(\wh K,i)$ is a \ti{field \ext}. So far for the algebraic
\sc\ of~$\wh K$.

Since the field $K$ is a \ms\ with metric $d_j$, we may invoke Hausdorff's theorem on the completion of a \ms\ (see \cite[Theorem 2.5.1 and p.~78
\S2.5]{Du}). The proof makes use of the \fw\ in\et ies.

Let $(M,d)$ be a \ms\ and let $a,b,c,d\in M$. Then we have
\bga8.107
|d(a,b)-d(c,d)| \le d(a,c)+d(b,d), \\
|d(a,c)-d(b,c)| \le d(a,b). \lb{8.108}
\e
Let $es$ denote the set of all \Cs s in $(M,d)$. The set $es$ is not empty, since \fe $a\in M$ the constant \sq\ $i_0(a)$ defined by $i_0(a):=a$ \fa
$k\in\N$ belongs to~$cs$. The map $i_0:K\to cs$ defined by $a\mt i_0(a)$ is \ti{in\jc}.

For the convenience of the reader we formulate these results as a theorem.

\bdf8.167
Let $\pz0F$ be a field of \ch istic~$0$, endowed with a real-valued \vl~$\F$ \sf ying \er{5.42}--\er{5.45}. A~\sq\ $\zb an\N$ is said to be a \tb{\fd\
\sq} if, \fe $\ve\in\R_{>0}$, \te s an \el~$N$ of~$\N$ \st
\beq8.106a
\F(a_m-a_n)<\ve \qh{\fa} m,n\in\N, \ m,n\ge N.
\e
The \sq\ $\zb an\N$ is said to \ti{\cg} to an \el\ $\ov a\in F$ if \fe $\ve\in\R_{>0}$ \te s an $N\in\N$ \st
\beq8.107a
\F(a_n-\ov a)<\ve \qh{\fa} n\in\N,\ n\ge N.
\e
Two \sq s $\zb an\N$, $\zb bn\N$ in $F$ are said to be \ti{\ev t\/} if \fe $\ve\in\R_{>0}$ \te s $N\in\Na$ \st
\beq8.108a
\F(a_n-b_n)<\ve \qh{\fa} n\ge N.
\e
A field with a real-valued \vl\ is said to be \ti{\cp}, if every \fd\ \sq\ \cg s to some \el\ of~$F$.
\edf

\bth8.168
Let $(F,\th)$ be a field $F$ with an $\R$-valued \emph{\vl} and let $d$ be the \emph{metric} defined by $d(x,y):=\th(x-y)$, $x,y\in F$. Then \te\
a field~$\wh F$, an $\R$-valued \vl~$\wh \th$ on~$\wh F$ with \crs\ metric~$\wh d$, and a map $i:F\to \wh F$ \st $(\wh K,i)$ is a \emph{field \ext},
that is, $i$~is an in\jc\ ring-\hm sm. \Mo $i:(F,d)\to(\wh F,\wh d)$ is an \emph{isometry} whose range is \emph{\sql\ dense} in $(\wh F,\wh d)$. If
$(F,d)$ is \cp, then $i$~is sur\jc.
\eth

\bxs8.169 \

\hph i,ii, $F:=\Q$, $\F:=|\cdot|^\rho$, $\rho\in\oz 0,1 $ (see Theorem \rf{t7.51}). Let $a\in\Q$, then $j(a)$ consists of all \sq s in~$\Q$ converging
to~$a$, and $\wh a\in\wh F$ consists of all \sq s in~$\Q$ converging to some \el\ $\a\in\R$ with $\a$~depending on~$\wh a$ but not on
$\rho\in\oz0,1 $.

\hph ii,i, $F:=\Q$, $\F(a):=(\F_p(a))^\si$, $\si\in\R_{>0}$ (see Theorem \rf{t7.51}). Let $\a\in\Q\sms0$ then, according to \cite[Section 6.3]{Nrs},
$\wh F$~can be identified with $p$-adic fields. These fields are cornerstones of the $p$-adic analysis.\index{p-adic analysis@$p$-adic analysis} They are used in Number Theory (see also
\cite[Chapter~18]{Alg} and references therein).

\hph iii,, Let $F$ be the set of all $2\t2$ matrices of the form $\bmk \begin{smallmatrix}a&-b\\b&a\end{smallmatrix}\emk$ with $a,b\in\Q$.
Setting
$\bz:=\bmk\begin{smallmatrix} 0&0 \\ 0&0 \end{smallmatrix}\emk$,
$\1:=\bmk\begin{smallmatrix} 1&0 \\ 0&1 \end{smallmatrix}\emk$ and
$J:=\bmk\begin{smallmatrix} 0&-1 \\ 1&0 \end{smallmatrix}\emk$, we find that $F$ endowed with the matrix \ad\ and \mlc\ forms not only a \cmt e ring
with unity~$\1$, but a field (note that
$\bmk\begin{smallmatrix} a&-b \\ b&a \end{smallmatrix}\emk = \bmk\begin{smallmatrix} c&-d \\ d&c \end{smallmatrix}\emk$ if $a=c$, $b=d$).
The field~$F$ can also be viewed as a $2$-\dm al \vs\ over~$\Q$ with
$\la \bmk\begin{smallmatrix} a&-b \\ b&a \end{smallmatrix}\emk:= \bmk\begin{smallmatrix} \la a&-\la b \\ \la b&\la a \end{smallmatrix}\emk$,
$\la\in\Q$, as well as an \asc e \cmt e \hbox{$\Q$-algebra} (see \E\df\ 5.2.18). Note that $I,J$ form a basis of the \vs\ and\break $I^2=I$,
$IJ=JI=J$, $J^2=-I$. The field~$F$ is ring-\is c to the subfield $\{a+ib: a,b\in\Q\}$ of the complex field~$\C$. \E\nm s of the form $a+ib$,
$a,b\in\Q$, are called Gauss \nm s (see \cite[p.~56]{Alg}). Clearly $F$~is not ordered since $J^2=-I$ ($i^2=-1$). If $z:=a+ib$ and $\ov z:=a-ib$
(the complex conjugate of~$z$), then $z\ov z=a^2+b^2=\ov zz$. The modulus of~$z$ defined by $\sqrt{a^2+b^2} = \sqrt{z\ov z}$, written $|z|$, is an
\Ar\ \vl\ on~$F$ since $|z|\ge0$, $|z\ov z|= |z|\,|\ov z|$, $|z+w|\le |z|+|w|$ and $|2+0i|=2\ne1$. It turns out that $\wh F=
\{\bmk\begin{smallmatrix} x&-y \\ y&x \end{smallmatrix}\emk:{x,y\in\R}\}$. \Mo if $\F(\bmk\begin{smallmatrix} a&-b \\ b&a \end{smallmatrix}\emk
:=\sqrt{a^2+b^2})$, $a,b\in\Q$, then $\wh\F (\bmk\begin{smallmatrix} x&-y \\ y&x \end{smallmatrix}\emk) = \sqrt{x^2+y^2}$, $x,y\in\R$.

\hph iv,, $F$ is an algebraic \nm\ field. We refer the reader to \cite[Section 18.6]{A2}.
\exs

\proof[Proof of Theorem \rf{t8.168}]
From Remark \rf{r8.165} we know that if $c(F)$ denotes the \cmt e ring with unity of \fd\ \sq s in $(F,\th)$ and $c_0(F)$ denotes the ideal of
null-\sq s in $(F,\th)$, then $\wh F:=c(F)/c_0(F)$ is a field. \Mo if the map $j:F\to \wh F$ is defined by $j(a):=[x]$ with $\bX_n(k):=a$ \fa
$n,k\in\N$ and $a\in F$, then $(\wh F,j)$ is a field \ext\ of the field~$F$.

The next observation is crucial. Let $\bX,\bX'\in c(F)$. Then $\bX-\bX'\in c_0(F)$ iff $\bX$ and~$\bX'$ are \ev t (see \E\df\ \rf{d8.167}). Indeed,
since the \vl~$\F$ is real-valued, $\bX-\bX'\in c_0(F)$ means that $\lim\limits_{n\ge0}\F(x_n-x_n')=0$. \Mo if $d(a,b):=\F(a-b)$, $a,b\in F$, $d$~is
a metric on~$F$, \sf ying $d(a+c,b+c)=d(a,b)$, hence $d(-a,-b)=d(a,b)$, $a,b,c\in F$ (see Remark \rf{r7.28}). \E\Ip $\F(x_n-x_n') = d(x_n,x_n')$,
$n\in\N$, and $\lim\limits_{n\ge0} d(x_n,x'_n)=0$.

It turns out that the \rl\ $\sim$ on $c(F)$, defined by $\bX\sim\bY$, $\bX,\bY\in c(F)$, if\break ${\lim d(x_n,y_n)=0}$, is an \ev ce \rl\ on an
\ti{arbitrary} \ms. Indeed, \ti{\sy y} is an immediate con\sq\ of (M2) in \E\df\ \rf{d7.8}. \ti{Reflexivity} follows from (M1) and the fact that the
constant \sq\ $\zb an\N$ in~$\R_{\ge0}$ with $a_n:=0$ \fa $n\in\N$ \cg s to~$0$. From (M1), (M3) we find for $\bX,\bY,\bZ\in c(F)$:
$$
0\le d(x_n,z_n) \le d(x_n,y_n)+d(y_n,z_n), \q n\in\N.
$$
If $\bX\sim\bY$ and $\bY\sim\bZ$, given $\ve\in\R_{>0}$, \te\ $N_1,N_2\in\N$ \st $d(x_n,y_n)<\ve/2$, $d(y_n,z_n)<\ve/2$ \fa $n\ge N:=\max(N_1,N_2)$.
\E\Tf $d(x_n,z_n)<\ve$ \we $n>N$, hence $\bX\sim\bZ$. Since $\bX,\bY,\bZ$ are arbitrary in~$c(F)$, the \rl~$\sim$ is \ti{\tr e}. We denote the \ev ce
class containing $\ov x$ by~$[\ov x]$, and by $i_0(a)$ the constant \sq\ in~$F$ \sf ying $i_0(a)(k)=a\in F$ \fa $k\in\N$. Note that $i_0(a)\in c(F)$
and the map $i_0:K\to c(K)$ is in\jc.

\Wanp to use Hausdorff's proof of his ``completion theorem'' (see  Theorem 2.5.1 in \cite{Du} and historical notes on page~78). This proof can be
found in many books. We refer the reader to \cite[pp.~640--641]{JvN}. For the convenience of the reader we give a statement of the theorem as well as
a sketch of its proof. Let $(X,d)$ be a \ms, let $c(X)$ be the set of \Cs s, let $\sim$ be the \ev ce \rl\ on $c(X)$ \itd above. Set $\wh X:=c(X)
/{\sim}$. For every pair $\bX,\bY\in c(X)$ we observe that $\lim d(x_n,y_n)$ exists in $\R_{\ge0}$. Indeed, from \er{8.107} we obtain
$|d(x_n,y_n)-d(x_m,y_m)| \le d(x_n,x_m)+ d(y_n,y_m)$, $m,n\in \N$. Since $\zb xn\N$ and $\zb yn\N$ are \Cs s, so is the \sq\ $\{d(x_n,y_n)\}_{n\in\N}$
in~$\R$, which is metrically \cp. Set $\a:=\lim d(x_n,y_n)$, $\b:=\lim d(x_n',y_n')$, $\bX',\bY'\in c(F)$. If $\bX\sim \bX'$, $\bY\sim\bY'$, then we
have
$$
0\le |d(x_n,y_n)-d(x_n',y_n')| \nde8.107 \le d(x_n,x_n')+d(y_n,y_n') \to 0.
$$
Hence $0\le |\a-\b|\to0$, and $\a=\b$. We set $\wh d([\bX],[\bY]) = \lim d(x_n,y_n)$. It turns out that $\wh d$ is a metric on~$\wh X$.

Next one shows that $i:X\to\wh X$ defined by $i(a):=[i_0(a)]$, $a\in X$, is an \ti{isometry} with \ti{dense range} in~$X$. Finally, one shows that
$(\wh X,\wh d)$ is a \ti{\cp} \ms.

Returning to the case $(X,d):=(F,d)$ with $d(a,b):=\th(a-b)$, $a,b\in F$, we set $\wh F:=\wh X=X/{\sim} = c(F)/c_0(F)$, $\wh\th([\bX]):=\wh d([\bX],
[\bz])$, $\bX,\bz\in c(F)$ and $i(a):=[i_0(a)]$, $a\in F$. It remains to show that (i)~$\wh\th$ is a $\R$-\vl\ on~$\wh F$, \sf ying $\wh d(\xi,\eta)=
\wh\th(\xi-\eta)$, $\xi,\eta\in\wh F$; (ii)~$\wh d(i(a),i(b)) = d(a,b)$, $a,b\in F$; and (iii)~$i(F)=\wh F$ iff $F$~is \cp.

\hph i,ii, We know that $\wh d$ is a metric on $F$ and that $\wh d([\bX],[\bY]) = \lim d(x_n,y_n)$ \fa $\bX,\bY\in c(F)$. \Mo $d(x,y)=\th(x-y)$,
$x,y\in F$. \If that $d(x+z,y+z)=d(x,y)$, $x,y,z\in F$, since $(x+z) - (y+z)=x-y$. \Mo $d(-x,-y)=\th(-x+y)=\th(-(x-y))=\th(x-y)$. \E\Ip $\th(x)=
d(x,0)$, $\th(x-y)=d(x-y,0)$ and $d(xy,0)=\th(xy) = \th(x)\cdot\th(y) = d(x,0)\cdot d(y,0)$, $x,y\in F$.

\Wanp prove (i). Note that \fe $\xi\in c(F)/{\sim}$ \te s $\bX\in c(F)$ \st $\xi=[\bX]$. We have $\wh\th([\bX]):=\wh d([\bX],\bz)$, $\bX\in F$.

\er{5.42} follows from (M1) in \E\df\ \rf{d7.8}.

\er{5.43}: $\wh\th([\bX]\hcd[\bY]) = \wh\th([\bX\cdot\bY]) = \wh d([\bX\cdot\bY],\bz) = \lim d(x_ny_n,0) = \lim (d(x_n,0)d(y_n,0)) = \lim d(x_n,0)
\cdot\lim d(y_n,0) = \wh d([\bX],\bz)\cdot\wh d([\bY],\bz) = \wh \th([\bX])\cdot\wh \th([\bY])$, $\bX,\bY\in c(F)$.

\er{5.44}: $\wh\th([\bX]{\hpl}[\bY]) = \wh\th([\bX+\bY]) = \wh d(\bX+\bY,\bz) = \lim d((\bX+\bY)(n),0) = \lim {d(x_n{+}y_n,0)} = \limu d(x_n,-y_n)
\le \limu d(x_n,0)+\limu d(0,-y_n) = \limu d(x_n,0)+\limu d(y_n,0) = \lim d(x_n,0)+\lim d(y_n,0) = \wh\th([\bX])+\wh \th([\bY])$.

We used \E\Pr s \rf{p8.10} and \rf{p8.12}.

\er{5.45}: We have $1+1\notin\{0,1\}$ in $\pz0F$ since $F$~is a field of \ch istic~$0$. \E\Tf $\wh\th([1+1]) = \wh d([1+1],\bz) = \lim d(1+1,0) =
d(1+1,0) \ne 0$.

\hph ii,i, We have $\wh d(i(a),i(b)) = \wh d([i_0(a)],[i_0(b)]) = \lim d(a,b) = d(a,b)$.

\hph iii,, If $i(F)=\wh F$ and $i$ is an isometry, then $i$~is a bi\jn\ since an isometry is in\jc, moreover $i\Inv$~is also a sur\jc\ isometry. Since
$(\wh F,\wh d)$ is \cp, isometries are \Lc, hence uniformly \ctn,
hence (\sql) \ctn. From \E\Pr\ \rf{p8.130a} we infer that $(F,d)$ is \cp. \E\oh if $(F,d)$ is \cp, we claim
that $i(F)=\wh F$. Indeed, clearly $i(F)\sbs \wh F$. Let $\xi\in\wh F$. Since $i(F)$ is \sql\ dense in~$\wh F$, \te s a \sq\ $\zb an\N$ in~$F$ \st
$\wh d(i(a_n),\xi)\to0$. \E\Tf $\{i(a_n)\}_{n\in\N}$ is a \Cs\ in $(i(F),\wh d)$. Since $i$~is an isometry, $\zb an\N$ is a \Cs\ in $(F,d)$, hence
since $(F,d)$ is \cp\ \te s $a\in F$ \st $d(a_n,a)\to0$. Hence $\wh d(i(a_n),i(a))\to0$. In view of the \uq\ of the limit we obtain $\xi=i(a)$, which
implies $\xi\in\wh F$, hence $\xi\in i(F)$. Since $\xi$~is arbitrary in~$\wh F$ we obtain $\wh F\sbs i(F)$, hence $\wh F=i(F)$.
\endproof

Note that if $F,\th,d,\wh F,i$ are as in Theorem \rf{t8.168}, and if $f$~is a bi\jn\ from $\wh F$ onto a set~$G$, and the metric $d_f$ on~$G$ is
defined by $d_f(g_1,g_2):= d(f\Inv(g_1),f\Inv(g_2))$, then $f:(\wh F,\wh d) \to (G,d_f)$ is a sur\jc\ isometry. \If that the \ms\ $(G,d_f)$ is \cp\
and the map $f\circ i:F\to G$ is an isometry with dense range in $(G,d_f)$. \E\Tf $(G,f\circ i)$ is another (if $f\ne\id$) \ti{completion} of the \ms\
$(F,d)$. One could call the completion used in the proof of Theorem \rf{t8.168} the Hausdorff completion of $(F,d)$.

\E\oh if $({\wh F}',i')$ with $i(F)$ dense in~${\wh F}'$, is a completion of the \ms\ $(F,d)$, then \te s a unique sur\jc\ isometry $j:(\wh F,\wh d)
\to ({\wh F}',{\wh d}')$.
$$
\xymatrix{{}&(\wh F,\wh d)\ar[dd]^j\\
(F,d)\ar[ru]^i\ar[rd]_{i'}&\\
&({\wh F}',{\wh d}')}
$$
This \ass\ is a con\sq\ of the next lemma.

\blm8.170
Let $(M_i,d_i)$, $i=1,2$, be \ms s and let $A$ be a dense subset of $M_1$.

\hph i,i, Suppose that $f,g$ are \emph{\ctn} maps from $(M_1,d_1)$ into $(M_2,d_2)$. If $f|_A=g|_A$, then $f=g$.

\hph ii,, Suppose that $f$ is a map from~$A$ into~$M_2$. If $(M_2,d_2)$ is \emph{\cp} and $f$~is \emph{uniformly \ctn}, then \te s a uniformly
\ctn\ map $\ov f:(M_1,d_1)\to(M_2,d_2)$ \st $\ov f|_A =f$.
\elm

\proof \

\hph i,i, Let $x\in M_1\sm A$. Since $A$ is \ti{dense} in~$M_1$, \te s a \sq\ $\zb xn\N$ in~$A$ \st $d(x_n,x)\to0$. Since $f,g$ are \ti{\ctn}, $f(x)=
f(\lim x_n)= \lim f(x_n) = \lim g(x_n) = g(x)$.

\hph ii,, Let $x\in M_1\sm A$, and $\zb xn\N$ as above. Then $\zb xn\N$ is a \Cs\ in $(A,d_1)$. Since $f$ is \ti{uniformly \ctn},
$\{f(x_n)\}_{n\in\N}$ is a \Cs\ in $(M_2,d_2)$ by \E\Pr\ \rf{p8.129}. Since $(M_2,d_2)$ is \ti{\cp}, \te s $y\in M_2$ \st\break ${d_2(f(x_n),y)\to0}$.
Let $\zb {x'}n\N$ in~$A$ be \st $d(x'_n,x)\to0$. Then the \sq\ $\zb zn\N$ defined by $z_{2n}:=x_n$ and $z'_{2n+1}:=x'_n$, $n\in\N$, also \sf ies
$d(z_n,x)\to0$. Indeed, given $\ve\in\R_{>0}$ \te\ $N_1,N_2\in\N$ \st $d(z_{2n},x)<\ve$ and $d(z_{2n+1},x)<\ve$ \fa $n\in\N$, $n\ge\max(N_1,N_2)$.
\E\Tf $\zb zn\N$ is a \Cs\ in $(A,d_1)$, hence $\{f(z_n)\}_{n\in\N}$ is a \Cs\ in $(M_2,d_2)$ by \E\Pr\ \rf{p8.129}. The sub\sq\
$\{f(z_{2n})_{n\in\N}$ \cg s to~$y$, hence the whole \sq\ \cg s to~$y$ (see Exercise \rf{ex8.128}), hence the sub\sq\ $\{f(z_{2n+1})\}_{n\in\N}=
\{f(x'_n)\}_{n\in\N}$ \cg s to~$y$. \E\Tf we may define $\ov f:M_1\to M_2$ by setting $\ov f(x):=y$.

We now prove that $\ov f$ is uniformly \ctn. In view of the uniform continuity of~$f$, the \fw\ holds: given $\ve\in\R_{>0}$, \te s $\d\in\R_{>0}$
\st $d_2(f(a),f(b))< \frac\ve2$ \we $a,b\in A$ \sf y $d_1(a,b)<\d$. Since $A$~is dense in $(M_1,d_1)$, given $x,y\in M_1$, \te\ \sq s $\zb xn\N$
and $\zb yn\N$ in~$A$ \st $d_1(x_n,x)\to0$ and $d_1(y_n,y)\to0$. We have by \er{8.107} $|d_1(x_n,y_n)-d_1(x,y)| \le d_1(x_n,x) + d_1(y_n,y)$ \fa
$n\in\N$, hence
$$
0\le \limu|d_1(x_n,y_n)-d_1(x,y)| \le \limu(d_1(x_n,x)+d_1(y_n,y)) \le \limu d_1(x_n,x)+\limu d_1(y_n,y) =0.
$$
Hence $d_1(x_n,y_n)\to0$. \E\te s $N\in\N$ \st $d_1(x_n,y_n)<\frac\ve2$ \fa $n\ge N$. \E\Tf $d_2(f(x_n),f(y_n))<\d$ \fa $n\ge N$.
Since $f(x_n)\nad{d_2}\to \ov f(x)$ and $f(y_n)\nad{d_2}\to \ov f(y)$, we obtain as above $d_2(f(x_n),f(y_n)) \to d_2(\ov f(x),\ov f(y))$. Hence
$d_2(\ov f(x),\ov f(y)) \le \frac\ve2<\ve$. Hence $\ov f$ is uniformly \ctn.
\endproof

\bex8.171
Show that if $f$ in Lemma \rf{l8.170}\,(ii) \sf ies $d_2(f(a),f(b)) \le Cd_1(a,b)^\a$ \fs $C\in \R_{>0}$ and $\a\in \oz0,1 _\R$ \fa $a,b\in A$, then
the \f~$\ov f$ \sf ies the same estimate. \Mo show that if $f$ is an isometry, then so is~$\ov f$.
\eex

\brm8.172 \

\hph i,ii, If in Lemma \rf{l8.170} $M_1:=[0,1]_\R$, $d_1$ the usual metric in~$\R$, $M_2:=\zo 0,1 _\R$ with $d_2:=d_1$, $A:=M_2$ and $f(x):=x$,
$x\in A$, then $f$~is uniformly \ctn, however, $f$~has no \ctn\ \ext. Hence the \cp ness of $(M_2,d_2)$ cannot be omitted.

\hph ii,i, If $M_1:=[0,1]_\R$, $M_2:=\R$, $d_1,d_2$ usual metric in~$\R$, $A:=\oz 0,1 _\R$, and $f(x):=\frac1x$, $x\in A$, then both $(M_1,d_1)$ and
$(M_2,d_2)$ are \cp, $f$~is \ctn, however, $f$~has no \ctn\ \ext. Note that $f$~is not uniformly \ctn.

\hph iii,, Note that in Section~\ref{ass.6} we extended the \f~$\psi_x$ by using the order \cp ness of $(\R,\ge)$ and the Lipschitz continuity
of~$\psi_x$ on bounded \il s (induced from the convexity of~$\psi_x$), for proving the continuity of~$\psi_x$.
\erm

\bex8.173
Let $\R_{\ge0}\t\R_{\ge0}$ be endowed with the Euclidean metric and $\R$ with the usual metric.

\hph i,ii, Let $f:\R_{\ge0}\t\R_{\ge0}\to\R$ be defined by
$$
f(x,y):= \bca
\dfrac{xy}{x^2+y^2} & \hbox{for } (x,y)\in(\R_{\ge0}\t\R_{\ge0})\sms{(0,0)}\\
0  &\hbox{for } (x,y)=(0,0).
\eca
$$
Show that $f|_{(\R_{\ge0}\t\R_{\ge0})\sms{(0,0)}}$ is \ctn, $f|_{\R_{\ge0}\t\{0\}}$ and $f|_{\{0\}\t\R_{\ge0}}$ are \ctn, $f$~is not \ctn\ at $(0,0)$.

\hph ii,i, Let $f:\R_{\ge0}\t\R_{\ge0}\to\R_{\ge0}$ be defined by $f(x,y):=xy$. Show that $f$~is \ctn\ but not uniformly \ctn. Show that
$f|_{[0,1]\t[0,1]}$ is uniformly \ctn.

\hph iii,, Let $\R_{\ge0}\t\R_{\ge0}\to\R_{\ge0}$ be defined by $f(x,y):=x+y$. Show that $f$~is uniformly \ctn.

\hph iv,, Let $f:\R\t\R\sms{(0,0)}$ be defined by $f(x,y):=\frac xy$. Show that $f$ is \ctn.
\eex

\def\namespec{Final Remark}
\begin{rspc}\lb{r8.174}
According to \cite[Section 8.8]{Nrs}, Ostrowski proved in 1918 that every \cp, \cmt e field with \Ar\ \vl\ is \is c to the field $\R$ or~$\C$. The
field~$\C$ can be viewed as a two-\dm al \asc e \cmt e $\R$-algebra with unity~$e$ (see \E\df\ 5.2.18) and with basis $\{e,i\}$ \st $i\cdot i=-e$.
Note that if $z:=\a e+\b i$, $\a,\b\in\R$, $\ov z:=\a e-\b i$, then $z\cdot\ov z=\a^2+\b^2$, hence $\frac1{\a^2+\b^2}\ov z$ is the inverse of~$z$
provided $z\in\C^\t:=\C\sms0$. \E\Ip $\C$~is an $\R$-\ti{division algebra}, that is, an $\R$-algebra \st the two \eq s: $ax=0$ and $bx=0$ have only
the trivial \so\ $x:=0$ provided $a\ne0$, $b\ne0$. In~1940 Heinz Hopf, using topological arguments, was able to prove that every finite \dm al \cmt e
$\R$-division algebra is at most two-\dm al. Note that \asc ity is not required. If, moreover, the algebra has a unity~$e$, then the algebra is \is c
to~$\R$ or to~$\C$. We invite the interested reader to look at Section~11 of~\cite{Nrs} for more recent developments.
\end{rspc}

\newpage
\Subsubsection{The \Ar\ field $\ee^\R$}\label{ass.9}

\advance\abovedisplayskip by-2pt
\advance\belowdisplayskip by-2pt
In what follows, we consider an \of, ring- and \ois c to the field $(\R,\ge)$, \st the \el s of its \pf, apart from the \nel\ of the additive group,
are ir\ra\ \nm s.

We first introduce a family of fields ring- and \ois c to $(\R,\ge)$, denoted by~$\a^\R$ with $\a\in\R_{>1}$. Proceeding in the same way we obtain
a subfield of~$\a^\R$ denoted by~$\a^\Q$, and we show that $\a^\Q$~is the \pf\ of the field~$\a^\R$. We observe that a \sft\ \cn\
on~$\a$ for the \pf~$\a^\Q$ \sf y the \cn\ mentioned above is that $\a$~be \tc\ (see \E\df\ \rf{d9.3}). Finally, we introduce the famous Euler
\nm~$e$, and give a sketch of the proof of its transcendence.

\subsubsection*{The fields $\a^\R$}  We recall that if $(X,\qu,e)$ is an \am\ (resp.\ \Cm, \PM, \ag) and $f$~is a bi\jn\ from~$X$ onto
a set~$X'$, then so is
$(X',\qu',e')$ with $e':=f(e)$, $x'\qu'y' := f(f\Inv(x')\qu f\Inv(y'))$ \fa $x',y'\in X'$. See Examples 2.1.5\,(v). Note that if $x\qu y=e$ for
$x,y\in X$, then $f(x)\qu' f(y)= e'$. Hence $f$~is a monoid-\is sm from $X$ onto~$X'$. See \E\df\ 2.1.7 and Lemma 2.1.8.

If $(X,\qu,e):=(\R,+,0)$, $X':=\R_{>0}$ and $f(x):=\wh\psi_\a(x)$, $\a\in\R_{>1}$, $x\in\R$, defined in \er{6.115}, then by Theorem \rf{t6.50}\,(iii),
(iv), $f$~is a bi\jn\ from $X$ onto~$X'$. Note that $(\R,+,0)$ is an \ag, and $f(x+y)=f(x)\cdot f(y)$ by Theorem \rf{t6.50}\,(ii).
\Mo $f(0)=\wh\psi_\a(0)=1$ by Theorem
\rf{t6.50}\,(vii), since $f$~is a monoid (group) \is sm from $(\R,+,0)$ onto $(\R_{>0},\cdot,1)$.

Recall that $\pz2\R$ is an \am\ and $\pz2{\R\sms0}$ is an \ag. Note that $f(1):= \wh\psi_\a(1)=\psi_\a(1) =\a$ by Lemma \rf{l6.37}. We define a binary
\op~$\odot$ on~$\R_{>0}$ by setting
\beq9.1
x'\odot y' := f(f\Inv(x') \cdot f\Inv(y')), \q x',y'\in\R_{>0}.
\e
\If that $(\R_{>0},\odot,\a)$ is an \am\ and that $(\R_{>0}\sms1,\odot,\a)$ is an \ag. We recall that $f$~is a bi\jn\ from~$\R$ onto $\R_{>0}$ and
that $f\Inv(x')=\log_\a(x')$ in view of \E\df\ \rf{d6.53}. We find that
\beq9.2
x'\odot y' = \wh\psi_\a(\log_\a(x') \cdot \log_\a(y')), \q x',y'\in\R_{>0}.
\e
We next show that (4.1.1), (4.1.2) (hence (4.1.3)) hold for $0:=1$, $x',y',z'\in\R_{>0}$ and ${{+}:={\cdot}}$, ${\cdot}:={\odot}$. Indeed, \fa $x'\in
\R_{>0}$ we have $1\odot x' \nde9.2 = \wh\psi_\a(\log_\a 1\cdot \log_\a x') = \wh\psi_\a(0\cdot \log_\a x')= \wh\psi_\a(0) =1$. Hence $1\odot x' =
1 = x'\odot1$. \Mo $x'\odot(y'\cdot z') = \wh\psi_\a(\log_\a x'\cdot \log_\a(y'\cdot z')) \nde6.138 = \wh\psi_\a(\log_\a x'\cdot(\log_\a y' +
\log_\a z')) = \wh\psi_\a\bigl((\log_\a x'\cdot \log_\a y') + (\log_\a x'\cdot \log_\a z')\bigr) \nad*=
\bigl(\wh\psi_\a(\log_\a x' \cdot \log_\a y')\bigr)\cdot (\wh\psi_\a(\log_\a x'\cdot\log_\a z') = (x'\odot y')\cdot(x'\odot z')$,
where in~$\nad*=$ we used Theorem \rf{t6.50}\,(ii). \If that \cn s (4.1.1)--(4.1.3) are \sf ied.
\E\Tf $(\R_{>0},\cdot,\odot,1,\a)$ is a \cmt e semiring with unity~$\a$ \st every \el\ of $\R_{>0}\sms1$ is invertible (that is, \fe
$x'\in\R_{>0}\sms1$ \te s $y'\in\R_{>0}$ \st $x'\odot y'=\a$). According to \E\df\ 4.4.1, $(\R_{>0},\cdot,\odot,1,\a)$ is a field. \Mo since $f$~is
a monoid-\is sm from $\pz1\R$ (resp.\ $\pz2{\R_{>0}\sms1}$) onto $\pz2{\R_{>0}}$ (resp.\ $(\R_{>0}\sms1,\odot,\a)$) by Examples 2.1.5\,(v), the
bi\jc\ map~$f$ is a ring-\is sm from the field $\pz0\R$ onto  the field $(\a^\R,\cdot,\odot,1,\a)$ with
\beq9.3
\a^\R := \R_{>0},
\e
in view of \E\df\ 4.1.7.
\goodbreak

\advance\abovedisplayskip by2pt
\advance\belowdisplayskip by2pt

Examination of the proof given above shows that we furnished a proof of \E\Pr\ 4.4.2\,(ii), whose proof was left as an exercise. Applying \E\Pr\
4.4.2\,(ii) with $\pz0F := \pz0\Q$, with $f:=\psi_\a$, the \rt ion of~$\wh \psi_\a$ to $\Q$ \er{6.120}, and setting
\beq9.4
\a^\Q := \psi_\a(\Q),
\e
we find that $(\a^\Q,\cdot,\odot,1,\a)$ is a field ring-\is c to $\pz0\Q$.

\ssk
\setbox0=\hbox{(b)} \sz=\wd0
We recall that the field $\Q$ is a \Pf\ by Lemma 4.5.42. \Mo $\Q$~is a subfield of~$\R$ in view of \E\df\ 4.4.24. Indeed, \\
\hbox to\wd0{(a)\hfil} $\pz1\Q$ is a subgroup of the group $\pz1\R$,\\
(b) $\pz2\Q$ is a \sbm\ of the monoid $\pz2\R$,\\
\hbox to\wd0{(c)\hfil} every \el\ of $\Q\sms0$ has an inverse in $\pz2\Q$ by Lemma 4.3.3.

\noi \Mo by \er{9.3}  and by \E\Pr\ \rf{p5.29}, the map
$\psi_\a :\Q \to \a^\Q$ is a ring-\is sm. We claim that $\a^\Q$ is the \pf\ of the field~$\a^\R$. The claim follows from Lemma 4.4.31. Since its
proof was left as an exercise, we give a hint for the convenience of the reader. Concerning part~(i), we consider the \fw\ situation. Let $Y,Y'$
be \am s (resp.\ groups), let $g:Y\to Y'$
be a monoid-\is sm, and let $H$~be a \sbm\ (resp.\ subgroup) of~$Y$. We claim that $g(H)$ is a \sbm\ of~$Y'$. We have $g:(Y,\qu,e) \to (Y',\qu',e')$,
$g(e)=e'$, $g(u\qu v) = {g(u)\qu'g(v)}$, $u,v\in Y$. Since $H$~is a \sbm\ of~$Y$, we have $e\in H$, hence $e'=g(e)\in g(H)$. Let $u',v'\in g(H)$, then
\te\ $u,v\in H$ \st $u'=g(u)$, $v'=g(v)$. Hence $u'\qu'v' = g(u)\qu g(v) = g(u\qu v)\in g(H)$. Hence $g(H)$ is a \sbm\ of~$Y'$. If $H$~is an \ag, that
is, an abelian subgroup of~$Y$, then \fe $u\in H$ \te s $v\in H$ \st $u\qu v=e$. Now let $u'\in g(H)$, then \te\ $u,v\in H$ \st $u'=g(u)$ and $u\qu v
= e$. Hence $g(u)\qu' g(v)=g(u\qu v) = g(e)=e'$. Thus $u'\qu' g(v)=e'$, and $u'$~is invertible in~$g(H)$, which \cp s the proof of the claim. In case
$Y,Y'$ are fields and $g$~is a ring-\is sm, use \E\df\ 4.4.24.

Concerning part (ii) of Lemma 4.4.31, observe that $f(\check F)$ is a subfield of~$F'$ in view of part~(i). Denote by~$\wt f$ the \rt ion of~$f$
to~$\check F$. Then $\wt f$ is a ring-\is sm from~$\check F$ onto $f(\check F)$.

Let $g$ denote the inverse of $\wt f$. Then $g:f(\check F)\to \check F$ is a ring-\is sm. Suppose, for \cd ion, that $f(\check F)$ has a proper
subfield~$H$, then $H\sbs f(\check F)$ and $H\ne f(\check F)$. Since $f$~is bi\jc, $\check F$~would have a proper subfield, which would \cd\ the \as\
that $\check F$~is a \Pf. Summarizing, we have proved

\blm9.1
The field $\a^\Q$ defined in \er{9.4} is the \pf\ of the field $\a^\R$ defined in~\er{9.3}.
\elm

We next show that the field $\a^\R$ is an \of\ (see \E\df\ 4.5.50) \ois c to the \of\ $(\R,\ge)$. To this end, we first introduce a total \og~$\le'$
on~$\a^\R$ by setting:
\beq9.5
x' \le' y' \q\hbox{if} \q \wh\psi_\a\Inv(x') \le \wh\psi_\a\Inv(y'), \q x',y'\in\a^\R.
\e
with $\le$ being the \og\ of $\R$.

Note that $\wh\psi_\a\Inv = \log_\a$. Since $\a^\R = \R_{>0}$, we have $\log_\a x' \le \log_\a y'$ iff $\log_\a x'-\log_\a y'\le0$ iff
$\log_\a \frac{x'}{y'}\le 0$ iff $\frac{x'}{y'}\le 1$ iff $x'\le y'$, $x',y'\in\R_{>0}$. \If that the \og~$\le'$ on~$\R_{>0}$ is the \ti{\rt ion to
$\R_{>0}$ of the \og\ on~$\R$}. In view of Theorem \rf{t6.50}\,(iii) we find that the map $f:\R \to \a^\R$ defined by $f(x):= \wh\psi_\a(x)$,
$x\in\R$, \sf ies: $x<y$ implies $f(x)<'f(y)$ \fa $x,y\in\R$ \fa $x,y\in\R$. It also follows that $x'<y'$ implies $f\Inv(x') <f\Inv(y')$, $x',y'\in
\R_{>0}$, since $f$~is bi\jc\ and the \og\ of~$\R$ is total (see Lemma 1.3.34). \E\Tf $f$~is an \ois sm from $(\R,\ge)$ onto $(\a^\R,\ge')$. It
remains to show that $(\a^\R,\ge')$ is an \of. Clearly, the \og~$\ge'$ on~$\R_{>0}$ is total since the \og~$\ge$ on~$\R_{>0}$ is total. It remains to
verify (4.5.120) and (4.5.121).

``(4.5.120)'': \ Let $x',y',z'\in\R_{>0}$ be \st $x'<'y'$, that is, $x'<y'$. We have to show that $x'\cdot z'<' y'\cdot z'$, that is, $x'\cdot z'
<y'\cdot z'$. But this is \cn\ (4.5.121) for the \of\ $\pz0\R$.
``(4.5.121)'': \ We have to verify that if $z'>'1$, and $x'<'y'$, then $z'\odot x' <' z'\odot y'$. Assuming $z'>1$, $x'<y'$, we have to show that
$z'\odot x' < z'\odot y'$, that is, $\wh\psi_\a(\log_\a z'\cdot \log_\a x') < \psi_\a(\log_\a z'\cdot \log_\a y')$. Since $\wh\psi_\a$ is an \ois sm,
it suffices to show that
$$
(\log_\a z')\cdot(\log_\a x') < (\log_\a z')\cdot (\log_\a y'),
$$
\ev tly $(\log_\a z')\cdot(\log_\a x'-\log_\a y')<0$. Since $\a>1$ and $z'>1$, we have $\log_\a z'\nad*> 0$. Hence it suffices to show that
$\log_\a(\frac{x'}{y'}) = \log_\a x' - \log_\a y' <0$. Since $0<x'<y'$, we have $0<\frac{x'}{y'}<1$, hence $\log_\a\frac{x'}{y'} \nad*<0$. In~$\nad*>$
and $\nad*<$ we used the fact that $\log_\a$ is strictly in\cre\ and that $\log_\a1=0$.

We thus proved

\blm9.2
The field $(\a^\R,\cdot,\odot,1,\a;\ge')$ is an \of\ and the map $\wh\psi_\a : (\R,\ge) \to (\a^\R,\ge')$ is a ring- and \ois sm.
\elm

\bdf9.3
An \el\ $\b\in\R\sms0$ is called \tb{algebraic} (over the field~$\Q$), if \te\ $N\in\Na$, $a:[0,N]\to\Z$ with $a_0,a_N\in\Na$ \st
\beq9.6
\sum_{k=0}^N a_k\b^k =0,
\e
and \tb{\tc\/} otherwise.\index{algebraic}\index{transcendental}
\edf

\blm LE4
Let $\a\in\R_{>1}$. If $\a$ is \tc, then $\a^r$ is ir\ra\ \fa $r\in\Q\sms0$.
\elm

\proof
Suppose $r:=\frac mn$, $m,n\in\Na$, and $\a^r$ is \ra. Then $\a^m=(\a^r)^n\in\Q\sms0$. \E\Tf $\a^r=(\a^m)^{\frac1n}$ is algebraic, hence not \tc. If
$r=-\frac mn$, $m,n\in\Na$, and $\a^r$~is \ra, then $\a^{\frac mn}= (\a^r)\mo\in\Q\sms0$, which is impossible by the first part of the proof.
\endproof

\bpb9.5
Can we replace in Lemma \rf{lLE4} \ti{ir\ra\/} by \ti{\tc\/}?
\exb

\If from Lemma \rf{l9.2}, Theorem \rf{t5.1} and Lemma \rf{l5.2} that the field $(\a^\R,\ge')$ is order-\cp. \Mo by \E\Pr\ \rf{p4.31},
the field $(\R,\ge)$ is \Ar, and by \E\Pr\ \rf{p4.34}, the field
$(\Q,\ge)$ is order-dense in $(\R,\ge)$. Since $\wh\psi_\a$ is an \ois sm from $(\R,\ge)$ onto $(\a^\R,\ge')$, it follows from Lemma \rf{l9.2} that
$\a^\Q$~is order-dense in $(\a^\R,\ge')$.

One verifies that $d'(x',y'):=\bigl|\log_\a\frac{x'}{y'}\bigr|$, $x',y'\in\R_{>0}=\a^\R$, is a metric on~$\a^\R$ (see \er{7.18}). It turns out that
$\wh\psi_\a:(\R,d)\to(\a^\R,d')$ ($d$~usual metric on~$\R$) is a sur\jc\ isometry. \If that the \ms\ $(\a^\R,d')$ is \cp, since $\wh\psi_\a$ as well
as~$\wh\psi_\a\Inv$ are uniformly \ctn. \Mo since $\Q$~is dense in $(\R,d)$, we have $\a^\Q$ (metrically) dense in $(\a^\R,d')$. Since $\a^\Q$~is
\ct y infinite, $(\a^\R,d')$~is a Polish space.

\def\psf #1#2{\frac{\wh\psi_\a(#2)-\wh\psi_\a(#1)}{#2-#1}}
\subsubsection*{\E\dfb ility of $\wh\psi_\a$} In view of Theorem \rf{t6.50} the \f~$\wh\psi_a$ is \po, strictly in\cre\ and convex (see \E\df\
\rf{d6.25}). By \E\Pr\ \rf{p6.32}, we obtain\index{differentiability}
\beq9.7
\psf ab \le \psf ac \le \psf bc
\e
\fa $a,b,c\in\R$ \st $a<b<c$. Setting $h_1:=b-a$, $h_2:=c-a$ and $g(h):=\frac{\wh\psi_\a(a+h)-\wh\psi_\a(a)}h$ for $h\in(0,c-a)$, we find by \er{9.7}
that $g(h_1)\le g(h_2)$ for $0<h_1<h_2$, hence $g$~is \ti{in\cre}. We claim that $g$~is \ti{\bb}. To this end, let $u<a$ and substituting $a:=u$,
$b:=a$, $c:=b$ in~\er{9.7} we obtain
$$
\psf ua \le \psf ab.
$$
\E\Tf we have
\beq9.8
\psf ua \le g(h) \qh{\fa $h\in(0,c-a)$.}
\e

By Theorem \rf{t5.1} and Corollary \rf{c5.12}, $\inf\{g(h): 0<h<b-a\}$ exists. We extend the \f~$g$ by setting
\beq9.9
\ov g(h):=\bca
g(h) & \hbox{if }h\in(0,b-a),\\
\infl_{0<h<b-a}g(h) &\hbox{if }h=0.
\eca
\e

\blm9.6
The \f\ $\ov g:\zo 0,b-a \to\R$, where $\zo 0,b-a $ and $\R$ are endowed with the usual metric on~$\R$, is \emph{\sql\ \ctn\ at} $0$ $($see \E\df\
\rf{d8.28}$)$.
\elm

\proof
We have $\ov g(0)\le \ov g(h)$ \fa $h\in(0,b-a)$. \E\fe $\ve\in\R_{>0}$ \te s $\d$ \st $\ov g(0)\le \ov g(\d)\le \ov g(0)+\ve$, by the \df\ of the
infimum. Since $\ov g$~is in\cre, we have $\ov g(0)\le \ov g(h)\le \ov g(\d)\le \ov g(0)+\ve$ \fa $h\in(0,\d)$. Hence \fe $\ve\in\R_{>0}$ \te s
$\d'\in (0,b-a)$ \st $|\ov g(h)-\ov g(0)|\le\ve$ for $h\in[0,\d')$. Let $\zb xn\N$ be an arbitrary \sq\ in $(0,b-a)$ \st $x_n\to0$. Then \te s
$N\in\N$ \st $0\le x_n<\d$ \fa $n>N$. \E\Tf $\ov g(0) \le \ov g(x_n)\le \ov g(\d)\le \ov g(0)+\ve$. \If that $\ov g(x_n)\to \ov g(0)$. Since the
\sq\ $\zb xn\N$ is arbitrary, the \f~$\ov g$ is \sql\ \ctn\ at~$0$ in view of \E\df\ \rf{d8.28}.
\endproof

\blm9.7
Let $(X_i,d_i)$, $i=1,2$, be \ms s, let $f:X_1 \to X_2$ and let $x_0\in X_1$. Then the \fw\ assertions are \ev t\dw
\bea9.10
{}&\hbox{$f$ is \sql\ \ctn\ at $x_0$,}\\
&\hbox{\fe $\ve\in\R_{>0}$ \te s $\d\in\R_{>0}$ \st $f(B(x_0;\d)) \sbs B(f(x_0);\ve)$.}\lb{9.11}
\e
\elm

\proof \

``\er{9.11} \ti{implies} \er{9.10}'': \ The proof is similar to the proof of Lemma \rf{l9.6}, \Tf it will be omitted.

``\er{9.10} \ti{implies} \er{9.11}'': \ By contraposition. Suppose \te s $\ve\in\R_{>0}$ \st \fa $n\in\Na$ \te s $x_n\in B(x_0;\frac1n) \not\sbs
B(f(x_0);\ve)$. Then $x_n\to x_0$ and $d(f(x_0),f(x_n))\ge \ve$, hence $\{f(x_n)\}_{n\in\N}$ does not \cg\ to $f(x_0)$, and $f$~is not \sql\ \ctn\
at~$x_0$.
\endproof

\def\podd#1{_{\raise5pt\hbox{$\scriptstyle #1$}}}
\blm9.8 \

\hph i,i, Let $a,b,c\in\R$ with $a<c<b$ and let $\vf:(a,c)\,\hbox{$($resp.\ $(c,b))$} \to \R$. Then \te s \emph{at most one} $L\in\R$ \st the \f\
$\ov \vf: \oz a,c \,\hbox{$($resp.\ $\zo c,b )$} \to\R$ defined by
$$
\ov \vf(x) := \bca
\vf(x), & x\in(a,c) \hbox{ $($resp.\ }(c,b)),\\
L, & x=c,
\eca
$$
is \sql\ \ctn\ at $c$ \wrt the usual metric on~$\R$.

If such an $L$ exists, then $L$ is called the \emph{left-hand} $($resp.\ \emph{right-hand}$)$ \emph{limit of $\vf$ at~$c$} and is denoted by
$\lml\podd{x\to c}\vf(x)$ $($resp.\ $\lmr\podd{x\to c}\vf(x))$.

\hph ii,, Let $\vf:(a,c)\cup(c,b)\to\R$ be \st $\lml\podd{x\to c}\vf(x)$ and $\lmr\podd{x\to c}\vf(x)$ \emph{exist} and \emph{are equal}, then
$\ov\vf:(a,b)\to\R$ defined by $\ov\vf(x):=\vf(x)$, $x\in(a,c)\cup(c,b)$ and $\ov\vf(c):=\lmr\podd{x\to c}\vf(x)$ $(=\lml\podd{x\to c}\vf(x))$ is
\sql\ \ctn\ at~$c$ \wrt the usual metric on~$\R$.
\elm

\proof \

\hph i,i, Case $\vf:(a,c)\to\R$. Suppose that \te\ \sq s $\zb xn\N$ and $\zb{x'}n\N$ \st $\vf(x_n)\to L_1$ and $\vf(x'_n)\to L_2$ with $L_1,L_2\in
\R$, $L_1\ne L_2$. Since $\vf(x_n)= \ov\vf(x_n)$ and $\vf(x'_n)=\ov\vf(x'_n)$ for $n\in\N$, $\ov\vf$~is not \sql\ \ctn\ at~$c$ in view of \E\df\
\rf{d8.28}. The proof of the case $\vf:(c,b)\to\R$ is similar, hence omitted.

\hph ii,, In view of Lemma \rf{l9.7}, it suffices to prove that \fe $\ve\in\R_{>0}$ \te s $\d\in\R_{>0}$ \st $|\ov\vf(x)-\ov\vf(c)|<\ve$ \fa $x\in\R$
\sf ying $|x-c|<\d$. (Note that the \cn\ $|\ov\vf(x)-\ov\vf(c)|<\ve$ can be replaced by $|\ov\vf(x)-\ov\vf(c)|\le \ve$. Indeed, given $\ve\in\R_{>0}$
and setting $\ve':=\frac\ve2\in\R_{>0}$, \te s $\d'\in\R_{>0}$ \st $|\ov\vf(x)-\ov\vf(c)|\le\ve'=\frac\ve2<\ve$. Set $\d:=\d'$.)

In view of part (i) and Lemma \rf{l9.7} \fe $\ve\in\R_{>0}$ \te\ $\d_+,\d_-\in\R_{>0}$ \st the \fw\ holds. Setting $L:=\lml\podd{x\to c}\vf(x) =
\lmr\podd{x \to c}\vf(x)$, we have $|L-\ov\vf(x)|<\ve$, ${|L-\ov\vf(y)|<\ve}$, \we $x\in \oz c-\d_-,c $ and $y\in\zo c,c+\d_+ $. Let
$\d:=\min(\d_+,\d_-)>0$. Then ${|L-\ov\vf(x)|<\ve}$ \we $|x-c|<\d$. We conclude the proof by applying again Lemma \rf{l9.7}.
\endproof

\bpr9.9
Let $a,b\in\R$, $a<b$, and let $\vf:(a,b)\to\R$ be convex. Then the \fw\ holds. \E\fa $x\in(a,b)$,

\hph i,ii, $\lmr\podd{h\to0}$ $($resp.\ $\lml\podd{h\to0})$ $\frac{f(x+h)-f(x)}h$ exists and is denoted by $f'_+(x)$ $($resp.\ $f'_-(x))$.

\hph ii,i, The maps $x\mt f'_-(x)$ and $x\mt f'_+(x)$, $x\in(a,b)$, are in\cre\ and \sf y
\beq9.12
f'_-(x) \le f'_+(x).
\e
If $y\in(a,b)$ with $x<y$, then
\beq9.13
f'_+(x)\le f'_-(y).
\e

\hph iii,, If $f'_-(x)<f'_+(x)$, then \fa $y\in(a,b)$\dw\ $f(y)-f(x)\ge z(y-x)$ \fa $z\in[f'_-(x),f'_+(x)]$.
\epr

\advance\abovedisplayskip by-2pt
\advance\belowdisplayskip by-2pt
\unitlength0.9truemm
\bxa9.10
Let $f:\R\to\R$ be defined by $f(x):=|x|$ \fa $x\in\R$. Then $f$ is convex,
$$
f'_-(x):=\bca -1, &x\le 0,\\+1, &x>0,\eca \qquad
f'_+(x):=\bca -1, &x< 0,\\+1, &x\ge0.\eca
$$
\exa

$$
\begin{picture}(0,40)(0,0)
\put(-75,20){\vector(1,0){40}}\put(-55,0){\vector(0,1){40}}\put(-37,17){$x$}\put(-53.5,16.5){0}\put(-56,-5){$f$}
\put(-20,20){\vector(1,0){40}}\put(0,0){\vector(0,1){40}}\put(-0.9,2){$\bullet$}\put(2,1.5){$-1$}\put(-1,-5){$f'_-$}
\put(35,20){\vector(1,0){40}}\put(55,0){\vector(0,1){40}}\put(54.1,36){$\bullet$}\put(52,37){$1$}\put(54,-5){$f'_+$}
\thicklines
\put(-73,38){\line(1,-1){18}}\put(-55,20){\line(1,1){18}}
\put(-18,3){\line(1,0){18}}\put(18,37){\vector(-1,0){18}}
\put(37,3){\vector(1,0){18}}\put(55,37){\line(1,0){18}}
\end{picture}
$$

\vskip10pt

For a nice proof of \E\Pr\ \rf{p9.9} we refer the reader to the first chapter of~\cite{Gamma}, where it is also shown that a real-valued \f\
on~$(a,b)$ is convex if and only if $f$~is \ctn\ on~$(a,b)$, $f'_+(x)$, $f'_-(x)$ exist \fa $x\in(a,b)$ and $f'_+,f'_-$ are in\cre.

We return to the study of the \f~$\wh\psi_\a$. We need some \df s.

\bdf9.11
Let $a,b\in\R$, $a<b$, and let $f:(a,b)\to\R$. The \f~$f$ is called \tb{\dfb le} at~$c$ if $f'_+(c),f'_-(c)$ exist and are equal.\index{differentiable at $c$}
\edf

\blm9.12
Let $f:(a,b)\to\R$ be \dfb le at $c\in(a,b)$. Then we have\dw

\hph i,i, the \f\ $\ov g:(a,b)\to\R$ defined by
\beq9.14
\ov g(h):=\bca
\dfrac{f(c+h)-f(c)}h & \hbox{ for } h\in(a-c,0)\cup (0,b-c),\\
f'_+(c)\,(=f'_-(c)) & \hbox{ for }h=0,
\eca
\e
is \ctn\ at $0$\sd

\hph ii,, the \f\ $f$ is \ctn\ at $c$.
\elm

\proof \

\hph i,i, The \f\ $g:(a,b)\to\R$ defined by $g(h):=\frac{f(c+h)-f(c)}h$ for $h\in(a-c,b-c)\sms0$, $g(0):=f_-(c)$ $(=f'_+(c))$ is \ctn\ at~$0$ in view
of Lemma \rf{l9.8}.

\hph ii,, \E\fe $\ve\in\R$ \te s $\d_+\in\R_{>0}$ \st $\ov g$ defined in \er{9.14} \sf ies $|\ov g(h)-\ov g(0)|<\ve$ provided $0<h<\d_+$. \E\Ip if
$\ve:=1$, there is $\d_+>0$ \st
$$
\left|\frac{f(c+h)-f(c)-hf'_+(c)}h\right|<1 \qh{for $h\in(0,\d_+)$},
$$
hence $|f(c+h)-f(c)|\le h(|f'_+(c)|+1)$ for $h\in(0,\d_+)$. Thus if $h_n\to0$, $h_n\in(0,\d_+)$, we have $|f(c+h_n)-f(c)|\to0$. Since the \sq\
$\{h_n\}$ is arbitrary, $f$~is \sql\ \ctn\ at~$c$. The other case $h\in(-\d_-,0)$ is similar. The conclusion follows from Lemma \rf{l9.12}\,(ii).
\endproof

\bds9.13
If \er{9.14} holds, then we write
\bga 9.15
\lim_{h\to0} \frac{f(x+h)-f(x)}h := \ov g(0),\\
f'(x):= \lim_{h\to0}\frac{f(x+h)-f(x)}h,
\e
and $f'(x)$ is called the \tb{\dv} of~$f$ at~$c$.\index{derivative}

\goodbreak
If the \f\ $f:(a,b)\to\R$ is \dfb le at all $c\in(a,b)$, then $f$~is called \ti{\dfb le on} $(a,b)$. Finally, $f'_+(x)$ (resp.\ $f'_-(x)$) is called
\ti{right-hand\/} (resp.\ \ti{left-hand\/}) \ti{\dv} of~$f$ at~$x$. $f'_+(x)$ and $f'_-(x)$ are called \ti{one-sided \dv s} of~$f$ at~$x$.
\eds

\advance\abovedisplayskip by2pt
\advance\belowdisplayskip by2pt
In the next theorem we collect some basic \pp ies of the differentiation of \f s. Proofs can be found in every good introductory book to Analysis.

\bth9.14 \

\hph I,I, Let $I$ be a nonempty open \il\ of~$\R$, and let $f,g:I\to \R$ be \dfb le at $c\in I$. Then the \fw\ holds\dw

\hph i,ii, If $\la\in\R$ then $\la f+g$ is \dfb le at~$c$ and
$$
(\la f+g)'(c) = \la f'(c)+g'(c).
$$

\hph ii,i, $f\cdot g$ is \dfb le and
$$
(f\cdot g)'(c) = f'(c)\cdot g(c)+f(c)\cdot g'(c).
$$

\hph iii,, If $g(c)\ne0$ then $g(x)\ne0$ on an \il\ $(c-\ve,c+\ve)\sbs I$ \fs $\ve\in\R_{>0}$ since $g$~is \ctn\ at~$c$. Then $\frac fg$ is \dfb le
at~$c$ and
$$
\biggl(\frac fg\biggr)'(c)= \frac{g(c)f'(c) - f(c)g'(c)}{(g(c))^2}.
$$

\hph II,, Let $I$, $c$ and $f$ be as in {\rm(I)}. Let $J$~be an open \il\ of~$\R$ containing $f(c)$ and let\break $g:J\to\R$ be \dfb le at $f(c)$.
Then the \f\ $g\circ f:I\to\R$ is \dfb le at~$c$ and
$$
(g\circ f)'(c) = g'(f(c))\cdot f'(c).
$$
\eth

The next lemma will be used in the proof of the Main Theorem.

\blm9.15
Let $\a,\b\in\R_{>1}$ and $x,y\in\R_{>0}$. We use Notation \rf{n6.51} for the \f\ $\wh\psi_\a(x)$, that is, $\a^x:=\wh\psi_\a(x)$. Then the \fw\
holds\dw
\bga9.17a
(\a^x)^y = \a^{xy},\\
\log_\a \b \cdot \log_\b \a = 1. \lb{9.18a}
\e
\elm

\proof
\er{9.17a} is \er{5.71}\,RI3 for $\a\in\R_{>0}$, $x,y\in\Q$, and \er{7.79} for $\a\in\R_{>1}$, $x,y\in\R$ (see Lemmata \rf{l7.43}, \rf{l7.44}).

\er{9.18a} follows from $\a^{(\log_\a\b)\cdot(\log_\b\a)}\nde9.17a = (\a^{\log_\a\b})^{\log_\b\a}=\b^{\log_\b\a} = \a = \a^1$. Since $\wh\psi_\a$
is in\jc, \er{9.18a} follows.
\endproof

The next lemma is usually called the \ti{Mean Value Theorem}.

\def\namespec{Fundamental Lemma}
\begin{tspc}\lb{l9.14}
Let $a,b\in\R$, $a<b$, let $f:[a,b]\to\R$ be \ctn\ on $[a,b]$ and \dfb le on $(a,b)$. Then \te s $c\in(a,b)$ \st
\beq9.17
\frac{f(b)-f(a)}{b-a} = f'(c).
\e
\end{tspc}

\proof
Let $g:[a,b]\to\R$ denote the \ti{linear interpolation} of~$f$ on $[a,b]$, that is, the only \pl\ map of degree less than or equal to one \sf ying
$g(a)=f(a)$ and $g(b)=f(b)$. Then $g(x) = f(a)+ (x-a)\frac{f(b)-f(a)}{b-a}$ for $x\in[a,b]$. Define $F:[a,b]\to\R$ by setting
$$
F(x):=f(x)-g(x), \q x\in[a,b].
$$
\If that $F\in C[a,b]$ and $F(a)=F(b)=0$. We distinguish two cases: either \te s $\ov x\in(a,b)$ \st $F(\ov x)\ne0$ or $F(x)=0$ \fa $x\in (a,b)$.
In the first case either $F(\ov x)<0$ or $F(\ov x)>0$. If $F(\ov x)<0$ then by Theorem \rf{t8.36}, \te s $c\in(a,b)$ \st $F(x)\ge F(c)$ \fa
$x\in[a,b]$. Thus \fa $h\in(a-c,0)$,
$$
\frac{F(c+h)-F(c)}h = \frac{F(c-|h|)-F(c)}{-|h|} = \frac{F(c)-F(c-|h|)}{|h|} \le0.
$$
Since $F'_-(c)$ exists, we obtain $F'_-(c)\le0$. Similarly, we have $F'_+(c)\ge0$. Since $F'_-(c) = F'_+(c) = F'(c)$, we have $F'(c)=0$. From \er{7.6}
and Lemma \rf{l9.12} we infer $0=F'(c) = f'(c)-g'(c)$. Observe that $g'(c)= \frac{f(b)-f(a)}{b-a}$, hence \er{9.17} holds. If $F(\ov x)<0$, then set
$\wt F(x):=-F(x)$, $x\in(a,b)$. We obtain $0=\wt F{}'(\wt c) = -f'(\wt c)+g'(\wt c)$ \fs $\wt c\in(a,b)$. Hence $f'(\wt c)=\frac{f(b)-f(a)}{b-a}$.
Finally, if $F(x)=0$ \fa $x\in[a,b]$, then $f(x)=g(x)$ \fa $x\in[a,b]$. \E\Tf $f'(x)=g'(x)= \frac{f(b)-f(a)}{b-a}$ \fa $x\in(a,b)$. Hence \er{9.17}
holds.
\endproof

\subsubsection*{\E\dv s of \pl s}
Real \pl\ maps and their \dv s play an important role in the study of the \f~$\wh\psi_\a$. Real \pl\ maps are \f s~$f$ from~$\R$ into~$\R$ defined by
$f(t):=\suml_{k=0}^N a_kt^k$ with \cf s $\zb ak{[0,N]}$, $N\in\N$. If $a_N\ne0$ then $N$~is called the \ti{degree} of the \pl. We will use the \fw\
notations:
\beq9.18
e_k(t) := t^k \qh{and } \wt e_k(t):= \frac{t^k}{k!}, \q k\in\N, \ t\in\R.
\e
Given $h\in \R\sms0$, we have
$$
\frac{e_0(t+h)-e_0(t)}h = 0, \q \frac{e_1(t+h)-e_1(t)}h = 1, \q t\in\R,
$$
hence
\beq9.19
e_0'(t)=0, \q e_1'(t)=1, \q t\in\R.
\e

We claim that
\beq9.20
te'_n(t) = ne_n(t),\q n\in\Na,\ t\in\R.
\e
We use \In\ on $n\in\Na$. The case $n=1$ follows from \er{9.19} and \er{9.18}. Suppose that \er{9.20} holds for $n:=k \in \Na$. To prove
$te'_{k+1}(t) = (k+1)e_{k+1}(t)$, $t\in\R$. We have $e_{k+1}(t) = te_k(t)$, hence for $h\in\R\sms0$
$$
\frac{e_{k+1}(t+h)-e_{k+1}(t)}h = \frac{(t+h)e_k(t+h)-te_k(t)}h = t\,\frac{e_k(t+h)-e_k(t)}h + e_k(t+h).
$$
Taking $\lim\limits_{t\nad\ne\to0}$ we obtain $te'_{k+1}(t) = t^2e'_k(t) + te_k(t) = tke_k(t)+te_k(t) = (k+1)te_k(t) = (k+1)e_{k+1}(t)$. Hence
\er{9.20} holds \fa $n\in\Na$ and all $t\in\R$, which proves the claim. We also have
\beq9.21
\wt e{}'_{n+1}= \wt e_n, \q n\in\N; \q\q \wt e{}'_0=0.
\e
Indeed, from \er{9.18}, \er{9.20} we find
$$
t\wt e{}'_{n+1}(t) \nde9.18 = t\,\frac1{(n+1)!}\,e'_{n+1}(t) \nde9.20 = \frac1{(n+1)!}(n+1)e_{n+1}(t) = \frac1{n!}\,te_n(t) = t\wt e_n(t), \q t\in\R.
$$
For $t\in\R\sms0$, dividing by $t$, we find
\beq9.22
\wt e{}'_{n+1}(t) = \wt e_n(t), \q n\in\N.
\e

\ti{Case} $n\in\Na$, $t:=0$: \ $\wt e_n(t) = \frac1{n!}0^n = \frac1{n!}0 = 0$. If $h\ne0$ then $\frac{\wt e_{n+1}(h) - \wt e_{n+1}(0)}h
= \frac1{(n+1)!}\,\frac{h^{n+1}-0^{n+1}}h = \frac1{(n+1)!}h^n \to 0$ as $h\nad{\ne0}\to0$, hence $\wt e{}'_{n+1}(0)=0$, and \er{9.22} holds for
$t:=0$, $n\in\Na$.

\ti{Case} $n:=0$, $t:=0$: \ $\wt e_0(t) = \frac1{0!}e_0(t) = 1\cdot0^0 = 1$. If $h\ne0$ then $\frac{\wt e_1(h)-\wt e_1(0)}h = \frac1{1!}
\,\frac{h-0}h = 1\cdot1$. Hence $\wt e{}'_1(0)=1$, and \er{9.22} holds for $t:=0$ and $n:=0$. It is customary to set $f\en0:=f$.

This completes the proof of \er{9.21}.

\bdf9.15
If a \f\ $f:(a,b)\to\R$ is \dfb le on $(a,b)$, and if the \f\ $g:(a,b)\to\R$, defined by $g(x):=f'(x)$, $x\in(a,b)$, is also \dfb le on $(a,b)$, then
$f$~is said to be \ti{twice \dfb le} on~$(a,b)$, $g'$~is called the \ti{second \dv} of~$f$ on $(a,b)$, and $g'$~is denoted by~$f''$ or $f\en2$. More
generally, $f\en{n+1}:= (f\en n)'$, $n\in\Na$, with $f\en0:=f$, \fa $n\in\N$, if $f\en n$ exists on $(a,b)$ and $f\en n$~is \dfb le on $(a,b)$.\index{derivative!second}
\edf

With these notations, the \fw\ holds:
\beq9.23
(\wt e_m)\en n = \bca
\wt e_{m-n} & \hbox{if }m\ge n,\\
0 &\hbox{if }m<n,
\eca
\e
\fa $m,n\in\N$.

\bpr9.16
Let $f$ be a real \pl\ of degree $N\in\N$ with \cf s $\zb ak{[0,N]}$. Then we have
\bea9.24
{}&a_k=\frac1{k!}\,f\en k(0) \qh{\fa $k\in[0,N]$,}\\
&f\en{N+1}(t)=0 \qh{\fa $t\in\R$.}\lb{9.25}
\e
\epr

\proof
To prove \er{9.24}, write $f(x)=\suml_{k=0}^N \wt a_k\wt e_k(x)$, $x\in\R$, with $\wt a_k:=k!a_k$, $k\in[0,N]$. Use \In\ on $N\in\N$, \er{9.23}, and
the fact that $g_1,g_2$ \dfb le on~$\R$ implies $\la_1g_1+\la_2g_2$ \dfb le on~$\R$ with $(\la_1g_1+\la_2g_2)'(x) = \la_1g_1'(x)+\la_2g_2'(x)$,
$\la_1,\la_2,x\in\R$.

For proving \er{9.25} use \er{9.23}.
\endproof

\newpage
\def\namespec{Main Theorem}
\begin{tspc}\lb{t9.17}
Let $\a\in\R_{>1}$. Then the \f\ $\wh\psi_\a:\R\to\R$ defined in \er{6.115} \sf ies\dw

\hph i,v, $\wh\psi_\a(x)$ is \dfb le at all $x\in\R$.

\hph ii,i, $\wh\psi_\a{}'(x) = \wh\psi_\a{}'(0)\wh\psi_\a(x)$ \fa $x\in\R$.

\hph iii,, $\wh\psi_\a{}'(0)>0$.

\hph iv,, \E\te s \ooo $\a\in\R_{>1}$, denoted by~$\bf e$, sometimes called \emph{Euler's \nm}, \sf ying $\wh\psi_\ee{}'(0)=1$. The \f\ $t\mt \wh\psi
_\ee(t)$ is usually denoted by $t\mt \ee^t$ or $t\mt \exp(t)$. \Mo we have $\wh\psi_\a(0)=\log_\ee a$.

\hph v,i, The \f\ $\log:\R_{>0}\to\R$ defined by $\log:=\wh\psi_\ee{}\Inv$ is \dfb le at all $y\in\R_{>0}$ and $\log'(y)=\frac1y$ \fa $y\in\R_{>0}$.
\Mo $\log'_\a(y) = (\log_e(\a))\mo \,\frac1y$.
\end{tspc}

\brm9.18
The \nm\ $\ee$ plays an important role in Analysis, Harmonic Analysis, Probability, Number Theory, Differential \E\eq s and so on. In order to give
some examples, we need the \fw\ \df s.
\erm

\bdf9.19
Let $f:\R_{>0}\,(\hbox{resp.\ }\N) \to\R$. If \te s $L\in\R$ \sf ying: \fe $\ve\in\R_{>0}$ \te s a real \nm\ (resp.\ \po\ \ig) $N$ \st $|f(x)-L|<\ve$
\we $x>N$ (\ev tly $x\ge N$), then $f(x)$ is said to \cg\ to~$L$ as $x$~tends to infinity, written
\beq9.26
\lim_{x\to\iy} f(x) = L.
\e
\edf

We mention the \fw\ famous results.

\begin{pnt}[see for example {\cite[Chapter 7]{Lax}}]
Let $\pi(x)$ be the \nm\ of primes less than or equal to~$x\in\R_{\ge2}$. Then
\beq9.27
\lim_{x\to\iy}x^{\frac{\pi (x)}x} = \ee.
\e
\end{pnt}

\bex pnt
Show the \ev ce of Theorem \rf{t9.17} and the usual formula given in~\cite{Lax}.
\eex

A much simpler but not trivial result is
\beq9.28
\lim_{n\to\iy}\frac1n(n!)^\frac1n = \ee\mo
\e
(see \cite[p.\ 387]{Widder}), which is a simplified version of the \ti{Stirling formula} (see \cite[p.~383]{Widder}):\index{Stirling formula}
\beq9.29
\lim_{n\to\iy} \frac{n!}{(n/\ee)^n \sqrt{2\pi n}} = 1,
\e
where $\pi$ is the ratio of a circle's circumference to its diameter.

\medskip
\ti{Proof of Theorem \rf{t9.17}}

\hph i,v, \If from the convexity of $\wh\psi_\a$ and Lemma \rf{l9.6} that $(\wh\psi_\a)'_+(x)$ exists \fa $x\in\R$. We next show that
$(\wh\psi_\a)'_-(x)$ exists \fa $x\in\R$, by using Theorem \rf{t6.50}\,(ii). Indeed, let $x\in\R$, $h\in\R_{<0}$, then $|h|\in\R_{>0}$ and
\bmlg
\frac{\wh\psi_\a(x+h)-\wh\psi_\a(x)}h = \frac{\wh\psi_\a(x-|h|)-\wh\psi_\a(x)}{-|h|}\\ = \wh\psi_\a(-|h|)\cdot \wh\psi_\a(|h|)
\frac{\wh\psi_\a(x)-\wh\psi_\a(x-|h|)}{|h|} = \wh\psi_\a(-|h|) \cdot \frac{\wh\psi_\a(x+|h|)-\wh\psi_\a(x)}{|h|},
\e
since $\wh\psi_\a(-|h|)\cdot\wh\psi_\a(|h|) = \wh\psi_\a(0) = 1$, and $h=-|h|$. Let $\zb hn\N$ be an arbitrary \sq\ in $\R_{<0}$ converging to~$0$.
Then
$$
\frac{\wh\psi_\a(x+|h_n|)-\wh\psi_\a(x)}{|h_n|} \to (\wh\psi_\a)'_+(x).
$$
We claim that $\wh\psi_\a(h_n) \to \wh\psi_\a(0)=1$, that is, \fe $\ve\in\R_{>0}$ \te s $N\in\Na$ \st \fa $n\ge N$, we have $|1-\wh\psi_\a(h_n)|<\ve$.
By \er{6.123} and Lemma \rf{l5.38}, \te s $\ov h<0$ \st $1-\ve < \wh\psi_\a(\ov h)\le 1$. Since $\wh\psi_\a$ is strictly in\cre\ by Theorem
\rf{t6.50}\,(iii), we obtain $1-\ve<\wh\psi_\a(\ov h)<\wh\psi_\a(h)<1$ \fa $h\in(\ov h,0)$. Since $h_n\to0$, \te s $N\in\N$ \st $|h_n-0|<-\ov h
\in\R_{>0}$ \we $n\ge N=N(\d)$. \E\Tf if $n\ge N$, we have $\ov h<-|h_n|=h_n<0$. Hence $1-\ve < \wh\psi_\a(\ov h)<\wh\psi_\a(h_n)<1$, and
$|1-\wh\psi_\a(h_n)|<\ve$. \csq, $\wh\psi_\a(h_n)\to1$, which proves the claim.

Applying Lemma \rf{l7.6} with $a_n:=\wh\psi_\a(h_n)$ and $b_n:=\frac{\wh\psi_\a(x+|h_n|)-\wh\psi_\a(x)}{|h_n|}$, we find $a_n\cdot b_n \to a\cdot b
= (\wh\psi_\a)'_+(x)$. \If that $(\wh\psi_\a)'_-(x)$ exists and is equal to $(\wh\psi_\a)'_+(x)$ \fa $x\in\R$. By \E\df\ \rf{d9.13} the \f\
$\wh\psi_\a$ is \dfb le at every $x\in\R$ with \dv\ $\wh\psi{}'_\a(x) = (\wh\psi_\a)'_-(x) = (\wh\psi_\a)'_+(x)$.

\ssk
\hph ii,i, Let $x\in\R$ and $h\in\R\sms0$. Then
$$
\frac{\wh\psi_\a(x+h)-\wh\psi_\a(x)}h = \frac{\wh\psi_\a(x)\cdot\wh\psi_\a(h) - \wh\psi_\a(x)}h = \wh\psi_\a(x)\cdot\frac{\wh\psi_\a(h)-\wh\psi_\a
(0)}h,
$$
by Theorem \rf{t6.50}\,(ii). Then
\bmlg
\wh\psi{}'_\a(x)= \wh\psi_+{}'(x) = \lmr_{h\to0}\frac{\wh\psi_\a(x+h)-\wh\psi_\a(x)}h = \lmr_{h\to0}\wh\psi_\a(x)\cdot \frac{\wh\psi_\a(h)
-\wh\psi_\a(0)}h\\ \nad*= \wh\psi_\a(x)\cdot\lmr_{h\to0}\frac{\wh\psi_\a(h)-\wh\psi_\a(0)}h = \wh\psi_\a(x)\cdot(\wh\psi_\a(0))'_+ =
\wh\psi_\a(x)\cdot\wh\psi{}'_\a(0).
\e
In $\nad*=$ we used \er{7.7} \fe \sq\ $h_n\nad>\to 0$.

\ssk
\hph iii,, Suppose, for \cd ion, that $\wh\psi{}'_\a(0)=0$. Then $\wh\psi{}'_\a(x)=0$ \fa $x\in\R$ and by Lemma \rf{l9.14} we find $\wh\psi_\a(1)
- \wh\psi_\a(0) = (1-0)\wh\psi{}'_\a(\xi)$ \fs $\xi\in(1,\a)$. However, $\wh\psi_\a(1)=\a$ and $\wh\psi_\a(0)=1$, hence $0=\wh\psi{}'_\a(\xi)=
\a-1>0$. A~\cd ion.

\ssk
\hph iv,, We first suppose that \te s $\b\in\R_{>1}$, $\b\ne\a$, \st $\wh\psi{}'_\b(0)=1$. Then $\wh\psi_\a(x) =\a^x = (\b^{\log_\b\a})^x$, $x\in\R$.
We use the ``Chain Rule'' (see Theorem \rf{t9.14}\,(II)) with $f(x):=\b^x$ and $g(x):=(\log_\b \a)x$, $x\in\R$. We find $(f\circ g)'(x) =
f'(g(x))\cdot g'(x) = \wh\psi_\b(x) \cdot \log_\b\a$, $x\in\R$. Setting $x:=0$, we obtain $(f\circ g)'(0) = \wh\psi_\b(0)\cdot\log_\b\a= \log_\b\a$.
Since $(f\circ g)(x) = \wh\psi_\a(x)$, $(f\circ g)'(0) = \wh\psi_\a(0)$. We have $\log_\b\a = \wh\psi{}'_\a(0)>0$ by (iii). By \er{9.18} we find
$\log_\a\b = (\wh\psi{}'_\a(0))\mo$, hence $\b = \a^{(\wh\psi{}'_\a(0))\mo}$. We now define
\beq9.30
\ee:= \a^{((\wh\psi_\a)'(0))\mo}.
\e
Since $\a^0=1$ and $0<(\wh\psi{}'_\a(0))\mo$, we have $1=\a^0 < \a^{(\wh\psi{}'_\a(0))\mo}=\ee$. Hence $\ee\in\R_{>1}$. We claim that
\beq9.31
\wh\psi_\ee{}'(x) = \wh\psi_\ee(x) \qh{\fa $x\in\R$.}
\e
Indeed, $\wh\psi_\ee{}'(x) = (f\circ g)'(x)$ with $f(x):=\a^x$, $g(x):=(\wh\psi{}'_\a(0))\mo x$, $x\in\R$. Applying again the ``Chain Rule'' we find
\bmlg
\wh\psi_\ee{}'(x) \nad{\rm(ii)}= \wh\psi{}'_\a(0) \cdot \wh\psi_\a ((\wh\psi{}'_\a(0))\mo x)(\wh\psi{}'_\a(0))\mo
= \wh\psi_\a \bigl((\wh\psi{}'_\a(0))\mo x\bigr)\\ \nde9.30 = \wh\psi_\a ((\log_\a \ee)x) = \a^{(\log_\a\ee)x} \nde9.18 = (\a^{\log_\a\ee})^x = \ee^x
= \wh\psi_\ee(x)
\e
\fa $x\in\R$, which proves the claim.

By (ii) we have $\wh\psi_\ee(0)=1$. Suppose that \te s $\b\in\R_{>1}$, $\b\ne\ee $ \st $\wh\psi{}'_\b(0)=1$, then $\b^x= \ee^{(\log_\ee \b)x}$ and
$\wh\psi{}'_\b(x)= (\log_\ee\b)\wh\psi_\b(x)$, $x\in\R$, hence $1=\wh\psi{}'_\b(0) = \log_\ee\b\cdot\wh\psi_\b(0) = \log_\ee\b$. \E\Tf $\b=\ee$ since
the \f\ $\log_\ee$ is in\jc. A~\cd ion.

Finally, by \er{9.30}, $(\wh\psi{}'_\a(0))\mo = \log_\a \ee$, hence $\wh\psi{}'_\a(0) = \frac1{\log_\a\ee} = \log_\ee\a$ by~\er{9.18}.

\ssk
\hph v,i, We first observe that $\wh\psi_\a:\R \to\R_{>0}$ is a \ti{\he sm} where $\R$ and $\R_{>0}$ are endowed with the usual metric on~$\R$.
Indeed, by Lemma \rf{l7.44} $\wh\psi_\a$ is an \ois sm from ${(\R,\ge)}$ onto $(\R_{>0},\ge)$. By Lemma \rf{l9.12}\,(ii), $\wh\psi_\a$~is \ctn\
on~$\R$. As an \ois sm, $\wh\psi_\a$ maps order-\il s $(a,b)$, $a\ne b$, into order-\il s $(\wh\psi_\a(a),\wh\psi_\a(b))$, hence $\wh\psi_\a([a,b])$
is contained in $[\wh\psi_\a(a),\wh\psi_\a(b)]$. From Theorem \rf{t8.36} we infer that $f\bigl([\wh\psi_\a(a),\wh\psi_\a(b)]\bigr) =
[\wh\psi_\a(a),\wh\psi_\a(b)]$, hence $\wh\psi_\a((a,b)) = (\wh\psi_\a(a),\wh\psi_\a(b))$.

Since \nos s of $(\R,d)$ are unions of open \il s, $\wh\psi_\a(U)$, $U\sbs \cO(\R)$, belongs to $\cO(\R_{>0})$. \E\Tf $(\wh\psi_\a)\Inv$ is \ctn,
hence $\wh\psi_\a$~is a \he sm, \Ip $(\wh\psi_\a)\mo = \log_\a$ is \sql\ \ctn.

We next show that $\log_\a$ is \ti{\dfb le} on~$\R_{>0}$, that is, \fa $y\in\R_{>0}$ and \sq s $\zb hn\N$ in~$\R$ \st $h_n\to0$ and $h_n\ne0$ \fa
$n\in\N$, \te s $L=L(y)$ with $\frac{\log_\a(y+h_n)-\log_\a y}{h_n} \to L$. Since $\wh\psi_\a$ is \dfb le at $x:= (\wh\psi_\a)\Inv y$,
$$
\lim_{x_n\nad\ne \to x} \frac{\wh\psi_\a(x_n)-\wh\psi_\a(x)}{x_n-x} = \wh\psi{}'_\a(x),
$$
for every \sq\ $\zb xn\N$ \st $x_n\to x$ and $x_n\ne x$ \fa $n\in\N$. Indeed, set $\wt h_n:=x_n-x$, $n\in\N$. Setting $x_n:=(\wh\psi_\a)\Inv y_n$,
$n\in \N$, and $h_n:=x_n-x$, $n\in\N$, we obtain $x_n = \wh\psi_\a{}\Inv(y_n) \to \wh\psi_\a{}\Inv(y)=x$ and $h_n\ne0$ for $n\in\N$ in view of the
bi\ji\ of~$\wh\psi_\a$. Note that $\wh\psi_\a(x_n)=y_n$. \E\Tf
$$
\lim_{\log_\a y_n\nad\ne\to \log_\a y} \frac{y_n-y}{\log_\a y_n-\log_\a y} = \wh\psi{}'_\a(x).
$$
Since $\log_\a y_n\nad\ne\to \log_\a y$ iff $y_n\nad\ne\to y$, setting $h_n:=y_n-y$, we obtain $h_n\nad\ne\to 0$ and
$$
\lim_{h_n\nad\ne\to 0} \biggl(\frac{\log_\a(y+h_n)-\log_\a y}{h_n}\biggr)\mo = \wh\psi{}'_\a(x).
$$
Finally, since $\wh\psi_\a{}'(x)>0$ (hence $\ne0$), we may apply Lemma \rf{l7.6} \er{7.8} and obtain
\bmlg
((\psi_\a)'(x))\mo = \biggl(\lim_{h_n\nad\ne\to0}\biggl(\frac{\log_\a(y+h)-\log_a(y)}{h_n}\biggr)\mo \biggr)\mo \\=
\biggl(\biggl(\lim_{h_n\nad\ne\to0}\frac{\log_\a(y+h)-\log_a(y)}{h_n}\biggr)\mo \biggr)\mo =
\lim_{h_n\nad\ne\to 0} \frac{\log_\a(y+h)-\log_a(y)}{h_n}.
\e
Hence $(\log_\a)'(y)$ exists and is equal to $(\psi_\a'(x))\mo$.
By (ii) $(\wh\psi_\a)'(x) = (\wh\psi{}'_\a)(0) \wh\psi_\a(x) \nad{\rm(iv)}= (\log_\ee\a)\wh\psi_\a(x) = (\log_\ee\a)y$. Hence $(\log_\a)'(y)=
\frac1{(\log_\ee\a)y}$. \E\Ip
$$
\log'(y) = (\log_\ee)'y = \frac1{\log_\ee \ee}\,\frac1y = \frac1y.\eqno{\Box}
$$

We claim that \fa $n\in\N$ we have
\beq9.34n
\sum_{k=0}^n \frac1{k!}\,x^k \le \ee^x \qh{\fa $x\in\R_{>0}$.}
\e
We use \In\ on $n\in\N$. If $n:=0$ then $\frac1{0!}x^0 = \frac11 = 1\le \ee^0$. Suppose \er{9.34n} holds for $n:=l\in\N$. We show that \er{9.34n}
holds for $n:=l+1$. Set
$$
S_n(x) := \sum_{k=0}^n \frac1{k!}\,x^k, \q n\in\N,\ x\in\R.
$$
Then $S_n$ is \dfb le on $\R$ \fa $n\in\N$. \Mo taking the \dv\ of $\exp-S_{l+1}$ at~$x$, we find $\ee^x - S_l(x)\ge0$, $x\in\R_{\ge0}$, by the \iv\
hypothesis, since $S_{l+1}'(x) = S_l(x)$ \fa $x\in\R$. \E\Tf by Lemma \rf{l9.14} we obtain $\ee^x - S'_{l+1}(x) = \ee^\xi - S'_{l+1}(\xi) = \ee^\xi
-S_l(\xi) \ge0$, \fs $\xi\in(0,x)$. Thus \er{9.34n} holds for $n:=l+1$, hence \fa $n\in\N$ and all $x\in\R_{>0}$. Thus the claim is proved. Note that
\fe $x\in\R_{\ge0}$, we have $0\le S_n(x)\le S_{n+1}(x) \le \ee^x$. \E\Tf $\supl_{n\ge0}S_n(x)$ exists and $\supl_{n\ge0}S_n(x)\le \ee^x$ \fa
$x\in\R_{>0}$. In order to prove \et y we use Taylor's theorem.

\bth9.23n
Let $I$ be a nonempty open \il\ of~$\R$, let $a,b\in\R$ be \st $a<b$ and $[a,b]\sbs I$. Let $f:I\to\R$ be $(n+1)$-\dfb le on~$I$ with $n\in\N$. Then
\te s $c\in(a,b)$ \st the \fw\ holds\dw
\beq9.35n
f(b) = \sum_{k=0}^n \frac{f\en k (a)}{k!}(b-a)^k + \frac{f\en{n+1}(c)}{(n+1)!}(b-a)^{n+1}
\e
with $f\en 0(a):=f(a)$.
\eth

Note that the \RHS\ of \er{9.35n} is well-defined since $f$~is $(n+1)$-\dfb le on~$I$ and $a\in I$. Observe that if $n:=0$, \er{9.35n} becomes
$$
f(b) = f(a)(b-a)^0 + f'(c)(b-a) = f(a)+f'(c)(b-a),
$$
which is \er{9.17}. However, in the Mean Value Theorem we only assumed the \ex\ of~$f'$ on~$(a,b)$. A stronger version of Theorem
\rf{t9.23n} requires only the \ex\ and continuity of~$f\en j$, $1\le j\le n$, on~$[a,b]$ and the \ex\ of $f\en{n+1}$ on~$(a,b)$.

\proof
Let $n\in\Na$. Set
$$
\vf(x) := f(b) - \sum_{k=0}^n \frac{f\en k(x)}{k!} (b-x)^k - \frac{(b-x)^{n+1}}{(n+1)!}\cdot M
$$
\fa $x\in[a,b]$, with $M\in\R$ chosen \st $\vf(a)=0$. This is possible since $\frac{(b-a)^{n+1}}{(n+1)!}>0$. Note that $\vf(b)=0$ since $f(b)-
\frac{f\en 0(b)}{0!}(b-b)^0 = f(b)-f(b) =0$, and $(b-b)^k=0$ for $k\in[1,n+1]$. Since $f\en k$ exists on~$I$ \fa $k\in[1,n]$, $f\en l$ is \ctn\
on~$[a,b]$ for $l\in[0,n]$ by Lemma \rf{l9.12}\,(ii). \E\pl\ maps are \dfb le \fa $n\in\N$ on~$I$, \Tf \ctn\ \fa $k\in[0,n+1]$. \If that $\vf$~is
\ctn\ on~$I$, hence $\vf|_{[a,b]}$ is \ctn. Since sums and product of \dfb le \f s on~$I$ are \dfb le on~$I$ by Theorem \rf{t9.14}\,I (i), (ii),
$\vf\en k$ is \dfb le on~$I$, hence on~$(a,b)$ \fa $k\in[1,n]$. In view of the Mean Value Theorem \te s $c\in(a,b)$ \st $\vf'(c)=0$. We now compute
the \dv\ of~$\vf$ at $x\in(a,b)$. Note that the \dv\ of $x\mt \frac{f\en k(x)}{k!}(b-x)^k$ for $k\ge1$ at~$x$ is
$$
\frac{f\en{k+1}(x)}{k!} (b-x)^{k-1} - \frac{f\en k(x)}{k!}\,k(b-x)^{k-1} = \frac{f\en{k+1}(x)}{k!} (b-x)^k - \frac{f\en k(x)}{(k-1)!}(b-x)^{k-1},
$$
by Theorem \rf{t9.14}\,I(ii) and II, and for $k=0$ at $x$ is equal to $-f'(x)$. \E\Tf
$$
\vf'(x) = -f'(x) - \sum_{k=1}^n \biggl[\frac{f\en{k+1}(x)}{k!} (b-x)^{k-1} - \frac{f\en k(x)}{(k-1)!}(b-x)^{k-1}\biggr] + \frac{(b-x)^n}{n!}\,M.
$$
If $n:=1$ then
\bmlg
\vf'(x) = -f'(x) + \frac{f\en{k+1}(x)}{n!} (b-x)^n - \frac{f\en 1(x)}{0!}\,1 + \frac{(b-x)^n}{n!}\,M \\=
-f'(x) + f'(x) + \frac{(b-x)^n}{n!}(-f\en{n+1}(x)+M).
\e
Hence, since $\vf'(c)=0$ and $\frac{(b-c)^n}{n!}\ne0$, we obtain $M=f\en{n+1}(c)$, which \cp s the proof for the case $n:=1$.

\ti{Case} $n\ge2$: \ Set $A_k:=\frac{f\en k(x)}{(k-1)|}(b-x)^k$, $k\in[1,n+1]$. Note that $-f'(x)=-A_1$. Then $\vf'(x) = -A_1 - \suml_{k=1}^n
(A_{k+1}-A_k) + \frac{(b-x)^n}{n!}M$. Observe that $\suml_{k=1}^n(A_{k+1}-A_k) = A_{n+1}-A_1$ (use \In\ on $n\in\Na$). \E\Tf
$$
\vf'(x) = -\frac{f\en{k+1}(x)}{n!}(b-x)^n + \frac M{n!}(b-x)^n.
$$
Since $\vf'(c)=0$, we obtain $M=f\en{n+1}(c)$ as in the case $n:=1$. This \cp s the proof of Theorem \rf{t9.23n}.
\endproof

\bco9.24
We have
\beq9.36n
\ee = \sup_{n\ge0} \sum_{k=0}^n \frac1{k!} \q\hbox{and}\q 2<\ee<3.
\e
\eco

\proof
We apply Theorem \rf{t9.23n} with $I:=(-1,3)$, $a:=0$, $b:=1$, $f(x):=\ee^x$, $x\in[0,1]$. \E\fa $n\in\Na$ $f\en n(x)=\ee^x$. \csq, \fa $n\in\Na$
\te s $c_n\in(0,1)$ \st $\ee-\suml_{k=0}^n\frac1{k!} = \frac1{(n+1)!}\ee^{c_n}$. Since $\frac1{(n+1)!}\ee^{c_n}>0$ and $\ee^{c_n}<\ee^1$, we obtain
$$
\biggl|\ee - \sum_{k=0}^n \frac1{k!}\biggr| \nde9.34n = \ee - \sum_{k=0}^n \frac1{k!} \le \frac1{(n+1)!}\,\ee^1.
$$
\Mo $2=\suml_{k=0}^1\frac1{k!} < \suml_{k=0}^n \frac1{k!}\le \ee$ for $n>2$. Note that $\suml_{k=2}^n \frac1{k!} \le \suml_{k=2}^n \frac1{2^{k-1}}
<1$ \fa $n\ge2$, hence $\suml_{k=0}^1\frac1{k!} + \suml_{k=2}^n\frac1{k!}<3$ \fa $n\in\N$. Hence $2<\ee<3$. Actually, $\ee\in(2.718,2.719)$.
\endproof

\bco9.25
The \fw\ holds\dw
\beq9.37n
\sum_{k=0}^n \frac{x^k}{k!} \ua \ee^x \qh{\fa $x\in\R_{\ge0}$.}
\e
\eco

\proof
Set $S_n(x):=\suml_{k=0}^n \frac{x^k}{k!}$ \fa $n\in\N$ and all $x\in\R_{\ge0}$. Clearly, $S_n(x)\le S_{n+1}(x)$ \fa $n\in\N$ and all $x\in\R_{\ge0}$.
Let $M\in\Na$. Applying Theorem \rf{t9.23n} with $I:=(1,M+1)$, $a:=0$, $b:=M$ and $f(x):=\ee^x$, $x\in I$, one finds
$$
0\nde9.34 \le \ee^x - \sum_{k=0}^n \frac1{k!}\,x^k = \ee^x\,\frac1{(n+1)!}\,M^{n+1} \le \ee^M\cdot\frac{M^{n+1}}{(n+1)!}
$$
\fa $x\in[0,M]$ and $n\in\Na$. From the next lemma, we infer
$$
0\le \limu_{n\ge0}\biggl(\ee^x - \sum_{k=1}^n \frac{x^k}{k!}\biggr) \le \limu_{n\ge0}e^M\,\frac{M^{n+1}}{(n+1)!} = \ee^M \cdot \limu_{n\ge0}
\frac{M^{n+1}}{(n+1)!} = \ee^M\cdot0 =0.
$$
\If that \fe $M\in\Na$ we have $0\le \limu_{n\ge0}\bigl(\ee^x-\suml_{k=0}^n \frac{x^k}{k!}\bigr) = 0$ \fa $x\in[0,M]$. \csq, $\suml_{k=0}^n \frac{x^k}
{k!} \ua \ee^x$ \fa $x\in[0,M]$, and \fa $M\in\R_{>0}$, hence \fa $x\in\R_{>0}$.
\endproof

\blm9.26
Let $R\in\R_{>0}$. Then
\beq9.38n
\lim_{n\ge0}\frac{R^n}{n!} = 0.
\e
\elm

\proof
Since $\R$ is \Ar, \te s $\ov n\in\Na$ \st $2R<\ov n$. Note that $(\ov n+l)! = \ov n! \prodl_{k=1}^l (\ov n+k)! \ge \ov n!\,\ov n{}^l$ for $l\in\Na$.
Hence for $n:=\ov n+l$ we have $0<\frac{R^n}{n!} \le \frac{R^{\bar n}R^l}{\ov n!\, (2R)^l} \le \frac{R^{\bar n}}{n!}\bigl(\frac12\bigr)^l \to 0$ as
$n$~(hence~$l$) tends to infinity.
\endproof

\bco9.27
The \nm\ $\ee$ is not \ra.
\eco

\proof
It suffices to show that $\ee\mo$ is not \ra. Suppose, for \cd ion, that \te\ $p,q\in\Na$ \st $\ee\mo=\frac p{q+1}$. Applying Theorem
\rf{t9.23n} with $I:=(-1,2)$, $f(x):=\ee^{-x}$, $x\in I$, $a:=0$, $b:=1$, $n:=q$ and noting that $f\en k(0)=(-1)^k$ by the ``Chain rule'', we find
$(q+1)!\,\ee\mo = p\cdot q!\in\Na$. \Mo $(q+1)!\,\ee\mo = (q+1)!\suml_{k=0}^q(-1)^k\frac1{k!}+\break (q+1)!\,(-1)^{q+1}\frac1{(q+1)!}\ee^{-c}$,
$c\in(0,1)$. Since $(q+1)!\,\frac1{k!}\in\Na$ for $k\in[0,q]$, the first term belongs to~$\Z$, and $\ee^{-c}\in(0,1)$ since $\ee\in(2,3)$ and
$c\in(0,1)$. \If that the \RHS\ does not belong to~$\Na$. A~\cd ion.
\endproof

\bth9.28
The \nm\ $\ee$ is \tc.
\eth

(See Theorem 1.2 and its proof in \cite{Baker}.)

\proof
We suppose that $\ee$ is algebraic and find a \cd ion. \E\Tf we suppose
\beq9.39n
\sum_{k=0}^n q_k\cdot \ee^k =0
\e
\fs $n\in\Na$ and some \ig s $q_i\in\Z$, $i\in[0,n]$, with $q_0>0$. In \cite[p.~4]{Baker} a \cd ion is obtained by introducing the \f
\beq9.40n
I_p(t) := \int_0^t \ee^{t-s}f_p(s)\,ds, \q p\in\Na,
\e
with
\beq9.41n
f_p(s):= s^{p-1} \prod_{k=1}^n (s-k)^p, \q s\in\R_{\ge0},
\e
and the \f\ $J:\Na\to\R$ defined by
\beq9.42
J(p) := \sum_{k=0}^n q_k I_p(k).
\e
It is shown that \te s $c\in\R_{>0}$ independent of~$p$ \st
\beq9.43
|J(p)| \le c^p, \q p\in\Na,
\e
and that
\beq9.44
J(p)\in\Z\sms0, \q (p-1)! \bigm| |J(p)| \ne0 \q \hbox{for $p$ \sf ying } p>\max(q_0,n), \ p \hbox{ prime}.
\e
Since $(p-1)! \,|J(p)|\ne0$, \te s $d\in\Na$ \st $|J(p)|=d(p-1)!$, hence $(p-1)!\le |J(p)| \le c^p = cc^{p-1}$
\cd ing \er{9.38n} for $p$~large.

Since we did not introduce the notion of integral, we will use another \rp ation of the \f~$I_p$.

\ssk
\ti{Proof of estimate \er{9.43}}. \ Let $f$ denote a \pl\ of degree $m\in\N$ with real \cf s, that is,
\beq9.46
f(t) = \sum_{k=0}^m c_ke_k(t), \q t\in\R,
\e
with $c_k\in\R$, $0\le k\le m$, $c_m\ne0$ and $e_k(t):=t^k$, $k\in\N$, $t\in\R$. From \er{9.40n} we infer that if $f:=f_p$ as in \er{9.41n} then we have
\beq9.47
I_p(t) = \ee^t \sum_{j=0}^m f\en j(0) - \sum_{j=0}^m f\en j(t), \q t\in\R,
\e
by using repeated integration by parts. We will use the \RHS\ of \er{9.47} as \df\ of $I(f_p)(t)$, that is, we set
\beq9.48
I(f)(t) := \ee^t \sum_{j=0}^m f\en j(0) - \sum_{j=0}^m f\en j(t), \q t\in\R.
\e
Note that
\beq9.49
I(f)(0) = 0.
\e
For $f$ as in \er{9.46} we set
\beq9.50
\ov f(t):= \sum_{k=0}^m |c_k|e_k(t), \q t\in\R.
\e
We have
\beq9.51
f(t)\le \ov f(t) \qh{and } {-f}(t)\le \ov f(t) \qh{\fa $t\ge0$,}
\e
hence $|f(t)| \le \ov f(t)$, $t\ge0$. However, $t\mt |f(t)|$, $t\ge0$, is \ti{not\/} always a \pl. Our first goal is to prove:
\beq9.52
I(f)(t) \le t\ee^t \ov f(t), \q t\ge0,
\e
It will be convenient to use $\wt e_k(t):=\frac1{k!}\ee_k(t)$, $t\in\R$, since $(\wt e_k)\en j = \wt e_{k-j}$ for $0\le j\le k$. \Mo $\wt e_0(0)=1$
and $\wt e_k(0)=0$ for $k\in\Na$. We have
$$
I(\wt e_k)(t) \nde9.48 = \ee^t \sum_{j=0}^k \wt e_{k-j}(0) - \sum_{j=0}^k \wt e_{k-j}(t) = \ee^t - \sum_{j=0}^k \wt e_j(t)\ge0, \q t\in\R_{\ge0}.
$$
In view of Taylor's theorem \rf{t9.23n} with $a:=0$, $b:=t\in\R_{\ge0}$, $n:=k$ and $f(t):=\ee^t$ we obtain in view of \er{9.49}:
\beq9.53
I(\wt e_k)(t) = \frac{(t-0)^{k+1}}{(k+1)!}\,\ee^c \qh{\fs $c\in(0,t)$.}
\e
Since $I(\la f)=\la I(f)$ for $\la\in\R$, we have $I(e_k)(t) = k!\,I(\wt e_k) = \frac1{k+1} tt^k\ee^c$ whenever $k\in\Na$ and $t\in\R_{>0}$. Hence for
$k\in\Na$ $I(e_k)\le t \ee^te_k$, $t\ge0$, in view of \er{9.53}. \If that if $k=0$ and $t\in\R_{>0}$, we have $I(e_k)(t)=\ee^t-1$ by \er{9.47},
hence $I(e_0) \le t\ee^c \le t\ee^te_0$ by the Mean Value Theorem. \E\Tf
$$
I(f)(t) = I\biggl(\sum_{k=0}^m c_ke_k\biggr) = \sum_{k=0}^mc_kI(e_k) \le \sum_{k=0}^m|c_k|\cdot I(e_k) \le \sum_{k=0}^m |c_k| t\ee^te_k
 = t\ee^t\ov f(t), \q t>0.
$$
The case $t:=0$ follows from \er{9.49}. This \cp s the proof of \er{9.52}.

Observe that if $(-f)(t):=-f(t)$, $f$ as in \er{9.46} and $t\ge0$, we have $(-f)(t) = \suml_{j=0}^m (-c_j)e_j(t)$, hence $\ov{-f}(t) = \suml_{j=0}^m
|c_j|e_j(t) = \ov f(t)$, $t\ge0$. \E\Tf $I(-f)(t) \le t\ee^t\ov f(t)$, $t\ge0$, thus $-I(f)(t) \le t\ee^t\ov f(t)$. \If that
\beq9.53a
|I(f)(t)| \le t\ee^t\ov f(t), \q t\ge0.
\e
Observe that we obtain (2) on page 4 of \cite{Baker} for $t\ge0$ without using the notion of integral. We now apply formula \er{9.53a} with $f:=f_p$
\itd in \er{9.41n}. Using formula (5.2.28) for the product of two \pl s $f$~and~$g$, we find that $\ov{f\cdot g}\le \ov f\cdot\ov g$. Hence, using
\In, we obtain
$$
\ov f_p(x) \le x^{p-1}\prod_{k=1}^n (x+k)^p, \q x\in\R_{\ge0}.
$$
Setting $A(n):=n(2n)^n$,  $n\in\Na$, we obtain
 $\ov f_p(t) \le n^{p-1}(2n)^{np} \le A(n)^p$ for $0\le t\le n$ and $n,p\in\Na$. \csq,
$$
|J(p)| \nde9.42 \le \sum_{k=0}^n |q_k|\,|I(\ov f)(k)| \le \biggl(\sum_{k=0}^n |q_k|\biggr)A(n)^p\le\bigl(\max_{0\le k\le n}|q_k|n+A(n)\bigr),
$$
hence \er{9.43} holds. \hfill$\Box$

\ssk
\ti{Proof of estimate \er{9.44}}. \ We first use \as\ \er{9.39n} ($\ee$~is algebraic) and obtain
\bmlg
J(p)\nde9.42 = \sum_{k=0}^n q_k I_p(k) \nde9.48 = \sum_{k=0}^nq_k\biggl(\ee^k\sum_{j=0}^m f\en j(0)-\sum_{j=0}^m f\en j(k)\biggr)\\ =
\sum_{j=0}^mf\en j(0)\biggl(\sum_{k=0}^nq_k\ee^k\biggr) - \sum_{j=0}^m\sum_{k=0}^nq_kf\en j(k) = -\sum_{j=0}^m\sum_{k=0}^nq_kf\en j(k),
\e
with $m:=(n+1)p-1$. \csq, we have
\beq9.54
J(p) = -\sum_{j=0}^m\sum_{k=0}^nq_kf\en j(k).
\e
In order to compute $f\en j(k)$ we use Leibniz's rule for the $n$-th \dv\ of a product of two \pl s, that is, if $g_1,g_2$ are two \pl s then
\bml9.59
(g_1\cdot g_2)\en n(t) = \sum_{k=0}^n\binom nkg_1\en{n-k}(t)g_2\en k(t) \\= \sum_{k=0}^n\binom nk g_1\en k(t)g_2\en{n-k}(t), \q t\in\R,\ n\in\Na.
\e
Formula \er{9.59} follows from the case $n:=1$ (see Theorem \rf{t9.14}\,(ii)) and \In\ on~$\Na$. Set
$$
A:=\sum_{k=1}^n q_k\sum_{j=0}^mf\en j_p(k) \qh{and } B:=q_0\sum_{j=0}^mf\en j_p(0),
$$
we have $J=-(A+B)$, and $f\en j_p(k)\in\Z$ for $1\le k\le n$ and $0\le j\le m$.

One verifies that $f\en j_p(k)=0$ for $k\in[1,n]$ and $0\le j< p$, and that $p! \mid f_p\en p(k)$ for $k\in[1,n]$, $p\le j\le m$. \If
that $p! \mid \bigl|\suml_{j=0}^mf\en j(k)\bigr|$ for $k\in[1,n]$, hence we obtain
\beq9.60
p! \mid |A|.
\e
Since $f\en j_p(0)=0$, $0\le j<p-1$, and $p-1 < j \le m$, we obtain
\beq9.61
p!\mid \biggl|\sum_{j=0,j\ne p-1}^m f\en j_p(0)\biggr|.
\e
We also have
\beq9.62
f_p\en{p-1}(0) = (p-1)!(-1)^{np}(n!)^p.
\e
\E\Tf since $p>|q_0|$ by \er{9.44}, we obtain
\beq9.63
p!\nbmid |q_0f_p\en{p-1}(0)| \qh{and } p!\nmid |B|.
\e
In view of \er{9.23} $(p-1)!\bigm| |B|$, hence $(p-1)!\bigm|J(p)$ and $p\nmid |J(p)|$. \E\Tf $J(p)\ne0$, and \er{9.44} holds.
\endproof

We conclude this section, hence the Appendix, by some final remarks not necessarily related to Section~\ref{ass.9}.

\brs9.22 \

\hph i,ii, \ti{About connectedness}

In Theorem \rf{t8.53} we \es ed the connectedness of the field~$\R$ endowed with the order topology. Since every point~$\{x\}$ in~$\R$ is closed,
$\R\sms x$, the \cpl\ of~$\{x\}$ in~$\R$ is open and is the disjoint union of two nonempty open sets $\{y\in\R:y<x\}$ and $\{y\in\R:y>x\}$. Thus
deleting a point
from~$\R$ ``disconnects''~$\R$. Similarly, deleting a point from the open connected set $(0,1)$ ``disconnects'' it. This fact was used in the proof
of Theorem \rf{t5.5} where the \ex\ of a $n$-th root of $\xi\in(0,1)$ was \es ed.\index{disconnects@``disconnects''}

In Theorem \rf{t8.147} we showed that $\R^N$, $N\in\Na$, endowed with the topology induced by the Euclidean metric (\ev tly, by every metric induced
by a norm on~$\R^N$) is \cnt ed. However, in~$\R^N$, $N\ge2$, deleting a point \ti{does not\/} ``dis\cnt'' $\R^N$. We give a sketch of a proof for the
case $N=2$ and $\bX:=\bz$. Suppose, for \cd ion, that \te\ two \nos s $U_1,U_2$ of $\R^2\sms{\bz}$ \st $U_1\cup U_2=\R^2\sms\bz$ and $U_1\cap
U_2=\vn$. Let $\bX_i\in U_i$, $i=1,2$, and let the segment $I\sbs \R^2$ be defined by $I:=\{\bY\in\R^2: \bY=(1-t)\bX_1+t\bX_2,\ t\in[0,1]\}$. Thanks
to the openness of~$U_1$ and~$U_2$ and the continuity of the map $t\mt \bY(t)=(1-t)\bX_1+t\bX_2$ from $[0,1]$ onto~$I$, \te s $\ve\in\R_{>0}$ \st
$\bY((0,\ve))\sbs U_1$ and $\bY((1-\ve,1))\sbs U_2$. Let $O_1:=\{t\in(0,1):\bY(t)\in U_1\}$ and $O_2:=\{t\in(0,1):\bY(t)\in U_2\}$. Then both
$O_1$~and~$O_2$ are open in $(0,1)$ and $O_1\cap O_2=\vn$. If $\bz\notin I$, then $O_1\cup O_2=I$, \cd ing the \cnt edness of $(0,1)$. If $\bz\in I$,
there is no \cd ion, but if we replace $\bX_1$ by~$\bX_1'$ close enough to~$\bX_1$ ($\|\bX_1-\bX_1'\|$ small enough), $\bX_1\in U_1$, then we may
choose $\bX_1'$ \st $\bz\notin I'$. For example, if $\bX_1:=(-1,0)$ and $\bX_2:=(1,0)$, then $\frac12\bX_1+\frac12\bX_2=\bz\in I$, but if
$\bX_1':=(-1,\ve)$, $\ve\in\R_{>0}$ small enough, then $\bz\notin I':=\{(1-t)\bX_1+t\bX_2 = (2t-1,(1-t)\ve), \ t\in[0,1]\}$. We obtain a \cd ion by
taking~$I'$ instead of~$I$. As a con\sq\ we find that $\R$ endowed with order topology and $\R^2$ with the ``Euclidean'' topology are \ti{not\/} \he
c. Indeed, if they were \he c, so would be $\R\sms0$ and $\R^2\sms\bz$, which is impossible since the image of a \cnt ed set ($\R^2\sms{\bz}$)
under a \ctn\ map ($\R\sms0$) is also \cnt ed. Note that $\R^2\sm\{\bX\in\R^2: \|\bX\|_2=1\}$ is \ti{not\/} \cnt ed. Observe that instead of~$I'$
above we could have chosen the union of two
segments $I_1,I_2$ with $I_1$~joining~$\bX$ to some $\bZ\in\R^2$ and $I_2$ joining $\bZ$~to~$\bX_2$ \st $\bz\notin I_1\cup I_2$.

\ssk
\hph ii,i, \ti{About pointwise \cg nce}

In Remark \rf{r7.28}\,(ii) it is claimed that if $I:=[0,1]_\R$, then there is \ti{no} metric~$D$ on~$\R^I$ \st $\bX_n\lm \bY$ iff $\bX_n\nad D\to \bY$
\fa \sq s $\zb \bX n\N$ and $\bY$~in~$\R^I$. Our aim is to give a sketch of a proof of the claim. It suffices to consider the case $I:=(0,1)$. Indeed,
if such a metric~$D$ would exist in $\R^{(0,1)}$, then the metric $D'(\bX',\bY'):= D(\bX,\bY) + |\bX'(0)-\bY'(0)| + |\bX'(1)-\bY'(1)|$ with
$\bX',\bY'\in\R^{[0,1]}$, $\bX:=\bX'|_{(0,1)}$, $\bY:=\bY'|_{(0,1)}$, would metrize the pointwise \cg nce for real-valued \f s on~$[0,1]$.
Suppose, for \cd ion, that the pointwise \cg nce on~$(0,1)$ is metrizable, with metric~$D$ on $\R^{(0,1)}$. We denote
by~$\1_A$ the indicator \f\ of a subset~$A$ of $(0,1)_\R$ defined by
$$
\1_A(x) := \bca
0 &\hbox{if } x\notin A,\\
1 &\hbox{if } x\in A.
\eca
$$
Set $f:=\1_A$ with $A:=(0,1)_\Q$ ($=(0,1)_\R\cap \Q$). Since $A$ is \ct y infinite, \te s a bi\jn\ $\vf:\N \to(0,1)_\Q$. \E\fe $n\in\N$ set $f_n:=
\1_{\vf([0,n])}$. Then $f_n\lm f$. Let $H:=\{\1_B: B\sbs(0,1)_\Q,\ B\hbox{ finite}\}$. Then $f$~belongs to $\ov H{}\en s$, the sequential closure
\wrt $D$ of the set~$H$, see \E\df\ \rf{d8.87}.

We claim that \fe\ $n\in\N$ \te s a \sq\ of \ctn\ \f s on~$(0,1)$ with values in~$\R$, denoted by $f_{n,k}$, $k\in\N$, \st $f_{n,k}\lm f_n$. Indeed,
\fe $i\in\N$, set $x_i:=\vf(i)$, and let $y_i\in(0,x_i)$, $z_i\in(x_i,1)$. Then we have $0<y_i<x_i<z_i<1$. Set for $k\in\Na$:
$$
y_{i,k}:= x_i-\frac1k (x_i-y_i), \qquad z_{i,k}:= x_i+\frac1k (z_i-x_i).
$$
We obtain
$$
0<y_i<y_{i,k}<x_i<z_{i,k}<y_i<1, \qh{and} y_{i,k}\ua x_i,\q z_{i,k}\da x_i.
$$

Define $g_{i,k}:(0,1)\to \R$ by setting
$$
g_{i,k}(t):=\bca
0, & t\in\oz 0,y_{i,k}, \\
k\,\dfrac{t-y_{i,k}}{x_i-y_i}, & t\in(y_{i,k},x_i), \\
1, & t=x_i, \\
-k\,\dfrac{t-z_{i,k}}{z_i-x_i}, & t\in (x_i,z_{i,k}), \\
0, & t\in[z_{i,k},1).
\eca
$$

$$
\hskip-320pt
\begin{picture}(-20,38)(0,0)
\put(0,0){\line(0,1){38}\line(1,0){33}}
\put(58,0){\line(1,0){42}}
\put(33,0){\line(1,3){12.5}}
\put(58,0){\line(-1,3){12.5}}
\put(44.5,36.5){$\bullet$}
\put(19,16){\vector(1,-1){10}}\put(72,16){\vector(-1,-1){10}}
\put(16,18){$y_{i,2}$}\put(69,18){$z_{i,2}$}
\put(-1,37){$\bullet$}\put(-3,36){$1$}
\put(-1,-3.5){$0$}\put(99,-3.5){$1$}
\put(20,-3){$y_i$}\put(43.5,-3){$x_i$}\put(70,-3){$z_i$}
\put(21,-1){$\bullet$}\put(45,-1){$\bullet$}\put(71,-1){$\bullet$}\put(99,-1){$\bullet$}
\end{picture}
$$

\vskip4pt
One verifies that $g_{i,k}\in C(0,1)$, $i\in\N$, $k\in\Na$, $g_{i,k}\lm\1_{\{x_i\}}$, $i\in\N$, and $\suml_{i=0}^n g_{i,k}\lm \1_{\vf([0,n])}=f_n$.
\If that $f_n\in\ov{C(0,1)}{}\en s$, where $C(0,1)$ denotes the set of all real-valued \ctn\ \f s on~$(0,1)$.

In view of \E\Pr\
\rf{p8.88} we have $f\in\ov H$ and $H\sbs\ov{C(0,1)}$. From \E\Pr\ \rf{p8.67}\,(v) and (iii), we obtain $\ov H\sbs \ov{\ov{C(0,1)}} = \ov{C(0,1)}$.
Hence $f\in\ov{C(0,1)}$. Since $\ov{C(0,1)} = \ov{C(0,1)}{}\en s$, \te s a \sq\ of \ctn\ \f s on~$(0,1)$ converging pointwise to the \f\
$f=\1_{(0,1)_\Q}$. However, it follows from Baire's lemma \rf{l8.149} that the pointwise limit of a \sq\ of \ctn\ \f s on~$(0,1)$ is \ctn\ at least
at one point in~$(0,1)$ (actually, at a dense subset of $(0,1)$). We invite the interested reader to find a reference in the literature. Note that
the \f\ $\1_{(0,1)_\Q}$ is nowhere \ctn.
Indeed, if $x\in(0,1)_\R \sm (0,1)_\Q$ then $\1_{(0,1)_\Q}(x)=0$. Given $x\in(0,1)_\R \sm (0,1)_\Q$ \te s a \sq\ $\zb xn\N$ in~$(0,1)_\Q$ \st $x_n\ua
x$, hence $x_n\to x$. Since $\1_{(0,1)_\Q}(x_n)=1$ \fa $n\in\N$, we have $\1_{(0,1)}(x_n)=1\not\to0 = \1_{(0,1)}(x)$. Similarly, \fe $y\in(0,1)_\Q$
\te s a \sq\ $\zb yn\N$ in $(0,1)_\R \sm (0,1)_\Q$ \st $y_n\to y$. Hence $\1_{(0,1)}(y_n)=0 \not\to 1 = \1_{(0,1)}(y)$. \E\Tf we obtain a \cd ion,
which proves the claim of Remark \rf{r7.28}\,(ii).

\ssk \goodbreak
\hph iii,, \ti{About \ap ion of \ctn\ \f s}

The aim of this remark is to show that real-valued \ctn\ \f s on~$[0,1]$ (more generally on~$[a,b]$) can be uniformly \ap ed by \pl s and by piecewise
linear \f s. Weierstrass proved in 1885 that \fe \f\ $f\in C[a,b]$ \te s a \sq\ of \pl s $\zb pn\N$ \st $\|f-p_n\|_{\sup}\to0$. It is natural to think
that the \sq\ of \pl s $\zb pn\N$ of degree less than or equal to~$n$ obtained by the Lagrange interpolation to~$f$ (that is, \sf ying
$p_n(x_k\en n)= f(x_k\en n)$,
$k\in[0,n]$, with $x_k\en n := a+\frac kn(b-a)$, $n\in\Na$) \cg s to~$f$ on $[a,b]$. C.~Runge showed that if $[a,b]:=[-5,5]$ and
$f(x):=\frac1{1+x^2}$, $x\in[-5,5]$, then $\zb pn\N$ diverges for $|x|>3.63\dots$ (see \cite[pp.~78--79]{Davis}).

Motivated by probabilistic arguments, S.~Bernstein (1912) constructed an explicit \sq\ of \pl s converging uniformly to a \ctn\ \f\ $f$ on~$[0,1]$.
For the convenience of the reader we present the elegant proof of Bernstein.

Let $f\in C[0,1]$ and let $n\in\Na$. Set
\beq9.34
(B_nf)(x):= \sum_{k=0}^n f\biggl(\frac kn\biggr) r_{k,n}(x)
\e
with $r_{k,n}(x):=\dfrac{n!}{k!(n-k)!}\,x^k(1-x)^{n-k}$, $x\in[0,1]$.
Note that $r_{k,n}(x)\ge0$ and the degree of $B_nf$ is less than or equal to~$n$. Setting $e_k(x):=x^k$, $x\in[0,1]$, $k\in\N$, we have \fa $n\in\Na$:
\beq9.35
B_ne_0=e_0; \q B_ne_1=e_1; \q B_ne_2-e_2 = \frac1n(e_2-e_1).
\e
Indeed, from $(x+y)^n=\suml_{k=0}^n \binom nk x^ky^{n-k}$, $x,y\in[0,1]$, setting $y:=1-x$, we find $e_0=B_ne_0$.

\Mo differentiating and \ml ying by~$x$ the \f\
$$
x\mt (x+y)^n = \suml_{k=1}^n \binom nk x^ky^{n-k},\q x\in\R, y\in[0,1],
$$
we find
\beq9.36
nx(x+y)^{n-1} = \sum_{k=1}^n k \binom nk x^ky^{n-k}.
\e
Setting $y:=1-x$, we obtain
$$
x = \sum_{k=0}^n \biggl(\frac kn\biggr)\binom nk x^ky^{n-k},
$$
that is, $e_1=B_ne_1$. Assuming $n\ge2$ and proceeding as above, we obtain by Theorem \rf{t9.14}\,(ii):
\beq9.37
nx(x+y)^{n-1} + n(n-1)x^2 (x+y)^{n-2} = \sum_{k=0}^n k^2\binom nk x^ky^{n-k}.
\e
Setting $y:=1-x$, we arrive at
$$
\frac1n x +\biggl(1-\frac1n\biggr)x^2 = (B_ne_2)(x), \q x\in[0,1],
$$
that is,
\beq9.38
B_ne_2 - e_2 = \frac1n(e_1-e_2), \q n\ge2.
\e
Note that
$$
B_1e_2(x) = \sum_{k=0}^1 \frac{k^2}{1^2} \frac1{k!(1-k)!}\,x^k(1-x)^{1-k} = 0+ x(1-x)^0 = x = \frac1n x+ \biggl(1-\frac1n\biggr)x^2 \qh{for}n=1.
$$
\E\Tf \er{9.38} holds \fa $n\in\Na$.

\Wanp estimate $\|B_nf-f\|_{\sup}$ for $n\ge1$. Since $B_ne_0=e_0$, we have
\bmlg
|B_nf(x)-f(x)| = \biggl|\sum_{k=0}^n f\biggl(\frac kn\biggr)r_{k,n}(x) - \sum_{k=0}^n f(x)r_{k,n}(x)\biggr|\\
= \biggl|\sum_{k=0}^n \biggl(f\biggl(\frac kn\biggr)-f(x)\biggr) r_{k,n}(x)\biggr| \le \sum_{k=0}^n \biggl|f\biggl(\frac kn\biggr)-f(x)\biggr|
r_{k,n}(x), \q x\in[0,1].
\e
In view of the \ti{uniform} continuity of~$f$ on~$[0,1]$ by Theorem \rf{t8.115}, \fe $\ve\in\R_{\ge0}$ \te s $\d\in\R_{>0}$ \st $|f(y)-f(x)|\le\ve$
\we $|x-y|\le\d$, $x,y\in[0,1]$. \Mo if $|x-y|>\d$ we obtain $|f(y)-f(x)|\le |f(y)|+|f(x)| \le 2\|f\|_{\sup} \le 2\|f\|_{\sup}\d^{-2}|x-y|^2$
for $x,y\in[0,1]$. \E\Tf
\beq9.39
|f(y)-f(x)| \le \ve +2\|f\|_{\sup}\d^{-2}(x-y)^2 \qh{\fa $x,y\in[0,1]$.}
\e
\csq, $|B_nf(x)-f(x)| \le \suml_{k=0}^n |f(\frac kn)-f(x)| \binom nk x^k(1-x)^{n-k} \le  \suml_{k=0}^n\frac12\ve r_{k,n}(x) +  \suml_{k=0}^n
2\|f\|_{\sup}\d^{-2}\bigl(x^2-2\frac kn x+(\frac kn)^2\bigr) r_{k,n}(x) = \frac12\ve (B_ne_0)(x) + 2\|f\|_{\sup}\d^{-2}\bigl[x^2(B_ne_0)(x)
\break -2x(B_ne_1)x + (B_ne_2)(x)\bigr]$ for $x\in[0,1]$.

Observe that $E(x):= x^2(B_ne_0)(x) - 2x(B_ne_1)(x)+(B_ne_2)(x) = -x^2+(B_ne_2)(x)= (B_ne_2-e_2)(x) = \frac1n(e_1-e_2)(x)$. Hence $|E(x)|\le
\frac1n\|e_1-e_2\|_{\sup}$, $x\in[0,1]$. Hence \fe $\ve\in\R_{>0}$, $\|B_nf-f\|_{\sup}\le \ve + \frac1n2\|f\|_{\sup}\d^{-2}\|e_1-e_2\|_{\sup}$.
\If that $\limu_{n\ge0}\|B_nf-f\|_{\sup}\le\ve$ \fe $\ve\in\R_{>0}$, thus $\limu_{n\ge0}\|B_nf-f\|_{\sup}=0$, which implies $\|B_nf-f\|_{\sup}\to0$
since $\|B_nf-f\|_{\sup}\ge0$.

\E\pl s $\{B_nf\}_{n\in\N}$ are usually called \tb{Bernstein \pl s}. We will denote by~$B_n$, $n\in\N$, the self-maps of $C[0,1]$ defined by
$f\mt B_nf$. The \Op s $B_n$ in the \Bs\ $(C[0,1],\|\cdot\|_{\sup})$ are linear and \ti{\po}, that is, $(B_nf)(x)\ge0$ \fa $x\in[0,1]$ \we
$f(x)\ge0$ \fa $x\in[0,1]$. We will denote by~$\le$ the \ti{partial \og} on $C[0,1]$ defined by $f\le g$ if $f(x)\le g(x)$ \fa $x\in[0,1]$ and
$f,g\in C[0,1]$. As a con\sq\ of their linearity, these \Op s are \ti{in\cre}, that is, $f\le g$ implies $B_nf \le B_ng$, $f,g\in C[0,1]$, $n\in\Na$.
Notice that if $f\in C[0,1]$ then $|f|$, defined by $|f|(x):=|f(x)|$, $x\in[0,1]$, is also \ctn\ on~$[0,1]$. Indeed, if $x_n\to x$ in~$[0,1]$ then
$0\le \bigl||f(x_n)|-|f(x)|\bigr| \le |f(x_n)-f(x)|\to 0$. \E\Tf \fa $f,g\in C[0,1]$ and $n\in\Na$ $B_nf-B_ng \le |B_nf-B_ng|$ and $-(B_nf-B_ng)
\le |B_nf-B_ng|$, hence $|B_nf-B_ng|\le B_n|f-g|$.
\ers

We will also use the \fw\ important \pp y of~$C[0,1]$, namely, if $f,g\in C[0,1]$, then $f\cdot g$, defined by
$(f\cdot g)(x):=f(x)\cdot g(x)$, $x\in[0,1]$, belongs to $C[0,1]$.

\Wanp prove the \fw

\blm9.22
Let $\zb Ln\Na$ be a \sq\ of \po\ \lop s in $C[0,1]$. If $L_ne_i=e_i$, $i=0,1$, \fa $n\in\Na$ and $\|L_ne_2-e_2\|_{\sup}\to0$, then $\|L_nf-f\|_{\sup}
\to0$ \fa $f\in C[0,1]$.
\elm

\proof
Let $f\in C[0,1]$ and $n\in\Na$. We have $(L_nf)(x) - f(x) = (L_nf)(x)-f(x)e_0(x) = (L_nf)(x) - f(x)(L_ne_0)(x) = (L_nf-L_nf(x)e_0)(x)$ \fa
$x\in[0,1]$. Hence $|L_nf-f| = |L_nf - L_nf(x)e_0| \le L_n|f-f(x)e_0|$. Note that $|f-f(x)e_0|(y) = |f(y)-f(x)e_0(y)| = |f(y)-f(x)|$, hence by
\er{9.39} \fe $\ve\in\R_{>0}$ there is $\d\in\R_{>0}$ \st
$$
L_n|f-f(x)e_0| \le \ve+\frac1n \,2\|f\|_\iy\d^{-2}\|L_ne_2-e_2\|_{\sup}.
$$
Hence $\limu_{n\ge0}\|L_nf-f\|_{\sup} \le \ve$ \fe $\ve\in\R_{>0}$, and $\|L_nf-f\|_{\sup}\to0$.
\endproof

\goodbreak

As an application we consider the problem of uniform \ap ion of a \ctn\ \f\ on~$[0,1]$ by piecewise linear interpolation. Let $f\in C[0,1]$ and for
$n\in\Na$ let $h:=\frac1n$ and $x_k:=kh$, $k\in[0,n]$, $k\in\N$. Define $L_nf$, $n\in\Na$, by setting
\bga9.40
(L_nf)(x_k):=f(x_k),\q k\in[0,n],\\
(L_nf)(x):= \frac{x_{k+1}-x}{x_{k+1}-x_k}f(x_k)+\frac{x-x_k}{x_{k+1}-x_k}\,f(x_{k+1})\lb{9.41}\\
\hskip150pt \hbox{for }x\in(x_k,x_{k+1}),\ k\in[0,n-1]. \non
\e
One verifies that $L_nf\in C[0,1]$, $L_n: C[0,1]\to C[0,1]$ is linear and \po\ \fa $n\in\Na$. \Mo we have $L_ne_0=e_0$ and $L_ne_1=e_1$, $n\in\Na$.

In view of Lemma \rf{l9.22}, $\|L_nf-f\|_{\sup}\to 0$ provided $\|L_ne_2-e_2\|_\iy\to 0$. We compute $L_ne_2(x)-x^2$ on $[x_k,x_{k+1}]$.
Setting $a:=x_k$, $b:=x_{k+1}$, we find
\bmlg
x^2 - L_ne_2(x) = \frac{b-x}{b-a}\,x^2 + \frac{x-a}{b-a}\,x^2 - \biggl(\frac{b-x}{b-a}\,a^2 + \frac{b-x}{b-a}\,b^2\biggr) \\
= \frac{b-x}{b-a}(x+a)(x-a) + \frac{x-a}{b-a}(x+b)(x-b) \\
= \frac1{b-a}\bigl[(b-x)(x+a)(x-a) - (b-x)(x+b)(x-a)\bigr]\\
= \frac1{b-a}(b-x)(x-a)(x+a-x-b) = \frac1{b-a}(b-x)(x-a)(a-b) = - (b-x)(x-a).
\e
Hence $0\le (b-x)(x-a) = L_ne_2(x)-x^2 = |L_ne_2(x)-x^2|\le h^2$.
Hence $|L_ne_2(x) - x^2| \le h^2 = (\frac1n)^2$ \fa $x\in[0,1]$. \E\Tf $\|L_ne_2-e_2\|_{\sup}\le\frac1{n^2}$, and $\|L_nf-f\|_{\sup}\to0$.

For the convenience of the reader (see Corollary 3.1.2 of \cite{Davis}) let $f\in C^1[a,b]$, $a<b$, and suppose that $f''(x)$ exists \fa $x\in(a,b)$,
then \te s $c\in(a,b)$ \st
$$
f(x) - \biggl(\frac{b-x}{b-a}\,f(a) + \frac{x-a}{b-a}\,f(b)\biggr) = -(b-x)(x-a) \cdot \frac12\,f''(c).
$$
Note that if, in \ad, $f''(x)\ge0$ \fa $x\in(a,b)$, then $f$~is convex on $[a,b]$.

More generally, the \fw\ theorem due to P.~Korovkin (see \cite{Sikk}) holds.

\bth9.23
Let $\zb Ln\Na$ be a \sq\ of \po\ \lop s on $C[0,1]$. Then $\|L_nf-f\|_{\sup}\to0$ \fa $f\in C[0,1]$ iff $\|L_ne_i-e_i\|_{\sup}\to0$ for $i=0,1,2$.
\eth

The proof is left as an exercise. Note that if $L:C[0,1]\to C[0,1]$ is \po, linear and \sf ies $Le_i=e_i$, $i=0,1,2$, then $L$ is the identity \Op.

Important \pp ies of the Bernstein \Op s can be found for example in~\cite{Davis}. We only mention three of them. The first one is easy. By using the
\f\ $x\mt \frac{x-a}{b-a}$ one may replace $C[0,1]$ by $C[a,b]$, $a<b$. Secondly, if $f\in C^1[a,b]$, that is, $f\in C[a,b]$, $f'_+(a)$, $f'_-(b)$,
$f'(x)$ for $x\in(a,b)$ exist, and the \f\ $g:[a,b]\to\R$ defined by
$$
g(x):=\bca
f'_+(a), & x=a,\\
f'(x), & x\in(a,b),\\
f'_-(b), & x=b,
\eca
$$
belongs to $C[a,b]$, then $\|B_nf-f\|_{\sup}\to0$ and $\|(B_nf)'-f'\|_{\sup}\to 0$, where $B_nf$ are the Bernstein \pl s on~$[a,b]$. See Theorem 6.3.2
in~\cite{Davis}. Finally, Generalized Bernstein \pl s can be defined on the hypercube $[0,1]^N$, see \cite[p.~122]{Davis}. Non-explicit \ap ions can
be obtained by Stone--Weierstrass theorems, see for example \cite[pp.~38--39]{JvN}.

Since the advent of the computer, new methods of \ap ion of \so s to partial differential \eq s became available. \E\Ip the ``finite \el'' method
plays an important role. In~\cite{Clement} an \Op\ is introduced with \crs\ error estimates. Such an approach has been used and improved by many
authors. A~nontrivial \fd\ \gn\ to vector-valued \f s can be found in~\cite{Gawlik}.

%% file: references.tex
\section*{Epilogue}\setcounter{equation}0
\def\sct{Epilogue}\def\sbsct{Epilogue}
\def\sss{}
\newdimen\wys
\addcontentsline{toc}{section}{Epilogue}

The notion of \ti{field\/} in mathematics is an algebraic \sc\ which plays
an important role in algebra (Galois theory \cite{Alg}, \cite{Is}, \cite{Artin}),
in number theory (fields of \ra~\cite{Nrs}, algebraic~\cite{Is}, $p$-adic~\cite
{Nrs} \nm s). The fields of real and complex \nm s play a \fd\ role in analysis
(series, Riemann integral, Riemann--Stieltjes integral), in measure theory
(Borel, Lebesgue measure), in distribution theory, in \f al analysis (Banach
and Hilbert spaces). There are not only fields of \nm s, but also of matrices,
\f s, and so on. For example the field of complex \nm s can be ``viewed'' as
a field of $2\times2$ real \ti{matrices} (\cite{Nrs}, p.~69). In \cite{Rem},
p.~318, a field of meromorphic \ti{\f s} is introduced.
According to K.~Mainzer (\cite{Nrs}, p.~22), the concept of field, not yet
under this name, appears in the writings of Abel and Galois, later in Kronecker
as \ti{domain of \ra ity} (\cite{Kron}, \S1), in Dedekind \cite{Ded}. Steinitz
gave an abstract \df\ of a field in~\cite{Stein}.

\newpage

\def\sct{References}\def\sbsct{References}
\def\sss{}